\RequirePackage{pdf14}
\documentclass[10pt,leqno,oneside,openright]{memoir}
\def\OPTstockformat{\stockaiv}
\def\OPTpageformat{\pageaiv}
\usepackage{microtype}
\usepackage[utf8]{inputenc}
\usepackage[T1]{fontenc}
\usepackage[sfscaled=false,ttscaled=false,charter]{mathdesign} 
\DeclareFontFamily{U}{pseus}{\skewchar\font'60}
\DeclareFontShape{U}{pseus}{m}{n}{%
     <-6> s * [0.95] eusm5%
    <6-8> s * [0.95] eusm7%
    <8-> s * [0.95] eusm10%
}{}
\DeclareFontShape{U}{pseus}{b}{n}{%
     <-6> s * [0.95] eusb5%
    <6-8> s * [0.95] eusb7%
    <8-> s * [0.95] eusb10%
}{}
\DeclareMathAlphabet\EuScript{U}{pseus}{m}{n}
\SetMathAlphabet\EuScript{bold}{U}{pseus}{b}{n}
\renewcommand{\mathscr}{\EuScript}
\DeclareSymbolFont{eulersymbols}  {U}{zeuex}{m}{n}
\DeclareMathSymbol{\euinfty}{\mathord}{eulersymbols}{153}

\makeatletter
\def\fvm@ratio{0.9}%
\def\fvs@ratio{0.9}%
\def\fve@ratio{0.9}%
\input{mdsffont.def}
\input{mdttfont.def}
\makeatother
\usepackage{mathtools}

\usepackage{amsthm}
\usepackage{thmtools}
\usepackage{stmaryrd}
\mathtoolsset{centercolon,mathic}
\PassOptionsToPackage{dvipsnames,svgnames}{xcolor}
\usepackage{tikz-cd}
\usepackage[style=authoryear-comp,maxnames=10,backend=biber,%
            hyperref=true]{biblatex}
\usepackage{enumitem}
\usepackage{stackengine}
\setstackEOL{\\}

\usetikzlibrary{decorations.pathreplacing}
\pgfdeclarearrow{
  name = mdbchto,
  setup code = {
    \pgfarrowssettipend{1.5\pgflinewidth}
    \pgfarrowssetbackend{-2.55\pgflinewidth}
    \pgfarrowssetlineend{-.25\pgflinewidth}
    \pgfarrowssetvisualbackend{-.27437\pgflinewidth}
    \pgfarrowsupperhullpoint{1.5\pgflinewidth}{0.25\pgflinewidth}
    \pgfarrowsupperhullpoint{-1.88\pgflinewidth}{3.05\pgflinewidth}
    \pgfarrowsupperhullpoint{-2.55\pgflinewidth}{2.31\pgflinewidth}
  },
  drawing code = {
    \pgfsetdash{}{0pt}%
    \pgfpathmoveto{\pgfpoint{1.5\pgflinewidth}{0.25\pgflinewidth}}%
    \pgfpathlineto{\pgfpoint{-1.88\pgflinewidth}{3.05\pgflinewidth}}%
    \pgfpathlineto{\pgfpoint{-2.55\pgflinewidth}{2.31\pgflinewidth}}%
    \pgfpathlineto{\pgfpoint{-0.5\pgflinewidth}{0.5\pgflinewidth}}%
    \pgfpatharctoprecomputed%
      {\pgfpoint{-.94118\pgflinewidth}{0\pgflinewidth}}
      {48.576}
      {-48.576}
      {\pgfpoint{-0.5\pgflinewidth}{-0.5\pgflinewidth}}
      {.66681\pgflinewidth}{.66681\pgflinewidth}{1}{1}
    \pgfpathlineto{\pgfpoint{-2.55\pgflinewidth}{-2.31\pgflinewidth}}%
    \pgfpathlineto{\pgfpoint{-1.88\pgflinewidth}{-3.05\pgflinewidth}}%
    \pgfpathlineto{\pgfpoint{1.5\pgflinewidth}{-0.25\pgflinewidth}}%
    \pgfpathclose%
    \pgfusepathqfill
  }
}
\pgfdeclarearrow{
  name = mdbchbar,
  setup code = {
    \pgfarrowssettipend{0.5\pgflinewidth}
    \pgfarrowssetbackend{-0.5\pgflinewidth}
    \pgfarrowsupperhullpoint{0.5\pgflinewidth}{2.0345\pgflinewidth}
    \pgfarrowsupperhullpoint{-0.5\pgflinewidth}{2.0345\pgflinewidth}
  },
  drawing code = {
    \pgfsetdash{}{0pt}%
    \pgfpathmoveto{\pgfpoint{0.5\pgflinewidth}{2.0345\pgflinewidth}}%
    \pgfpathlineto{\pgfpoint{-0.5\pgflinewidth}{2.0345\pgflinewidth}}%
    \pgfpathlineto{\pgfpoint{-0.5\pgflinewidth}{-2.0345\pgflinewidth}}%
    \pgfpathlineto{\pgfpoint{0.5\pgflinewidth}{-2.0345\pgflinewidth}}%
    \pgfpathclose%
    \pgfusepathqfill
  }
}
\makeatletter
\pgfdeclarearrow{
  name = implies,
  setup code = {
    \pgfmathsetlength{\pgfutil@tempdimb}{.5\pgflinewidth-.5*\pgfinnerlinewidth}%
    \pgfarrowssettipend{1.5\pgfutil@tempdimb}
    \pgfarrowssetbackend{-4.948\pgfutil@tempdimb}
    \pgfarrowssetlineend{0.0\pgflinewidth}
    \pgfarrowssetvisualbackend{0.0\pgflinewidth}
    \pgfarrowsupperhullpoint{1.5\pgfutil@tempdimb}{0.25862\pgfutil@tempdimb}
    \pgfarrowsupperhullpoint{-4.2931\pgfutil@tempdimb}{4.9138\pgfutil@tempdimb}
    \pgfarrowsupperhullpoint{-4.948\pgfutil@tempdimb}{4.155\pgfutil@tempdimb}
  },
  drawing code = {
    \pgfsetdash{}{0pt}%
    \pgfmathsetlength{\pgfutil@tempdimb}{.5\pgflinewidth-.5*\pgfinnerlinewidth}%
    \pgfsetdash{}{0pt}%
    \pgfpathmoveto{\pgfpoint{1.5\pgfutil@tempdimb}{.25862\pgfutil@tempdimb}}%
    \pgfpathlineto{\pgfpoint{-4.2931\pgfutil@tempdimb}{4.9138\pgfutil@tempdimb}}%
    \pgfpathlineto{\pgfpoint{-4.948\pgfutil@tempdimb}{4.155\pgfutil@tempdimb}}%
    \pgfpathlineto{\pgfpoint{0\pgfutil@tempdimb}{0\pgfutil@tempdimb}}%
    \pgfpathlineto{\pgfpoint{-4.948\pgfutil@tempdimb}{-4.155\pgfutil@tempdimb}}%
    \pgfpathlineto{\pgfpoint{-4.2931\pgfutil@tempdimb}{-4.9138\pgfutil@tempdimb}}%
    \pgfpathlineto{\pgfpoint{1.5\pgfutil@tempdimb}{-.25862\pgfutil@tempdimb}}%
    \pgfpathclose%
    \pgfusepathqfill%
    \pgfsetfillcolor{white}%
    \pgfpathmoveto{\pgfpoint{1.5\pgfutil@tempdimb}{.25862\pgfutil@tempdimb}}%
    \pgfpathlineto{\pgfpoint{1.5\pgfutil@tempdimb}{4.9138\pgfutil@tempdimb}}%
    \pgfpathlineto{\pgfpoint{-4.2931\pgfutil@tempdimb}{4.9138\pgfutil@tempdimb}}%
    \pgfpathclose%
    \pgfusepathqfill%
    \pgfsetfillcolor{white}%
    \pgfpathmoveto{\pgfpoint{-4.2931\pgfutil@tempdimb}{-4.9138\pgfutil@tempdimb}}%
    \pgfpathlineto{\pgfpoint{1.5\pgfutil@tempdimb}{-4.9138\pgfutil@tempdimb}}%
    \pgfpathlineto{\pgfpoint{1.5\pgfutil@tempdimb}{-.25862\pgfutil@tempdimb}}%
    \pgfpathclose%
    \pgfusepathqfill
  }
}
\makeatother

\tikzset{>=mdbchto}
\tikzset{every picture/.style={line width=0.58pt}}
\tikzcdset{arrow style=tikz,%
  tikzcd implies/.tip={implies},%
  invisible/.style={/tikz/draw=none}}

\usepackage{hyperxmp}
\PassOptionsToPackage{hyphens}{url}
\PassOptionsToPackage{hidelinks}{hyperref}
\usepackage{bookmark}

\hypersetup{%
  pdftitle={Pursuing Stacks},
  pdfauthor={Alexander Grothendieck},
  pdflang={en},
  bookmarksopen=true,
  bookmarksopenlevel=0,
  bookmarksnumbered=true,
  hypertexnames=false,
  linktocpage=true,
  plainpages=false,
  colorlinks,
  citecolor=DarkGreen,
  linkcolor=DarkGreen,
  urlcolor=DarkGreen,
  breaklinks
}
\DeclareFieldFormat{eprint:mr}{%
  \mkbibacro{MR}\addcolon\space
  \ifhyperref
    {\href{http://www.ams.org/mathscinet-getitem?mr=#1}{\nolinkurl{#1}}}
    {\nolinkurl{#1}}}
\DeclareFieldFormat{eprint:eudml}{%
  \mkbibacro{EUDML}\addcolon\space
  \ifhyperref
    {\href{https://eudml.org/doc/#1}{\nolinkurl{#1}}}
    {\nolinkurl{#1}}}
\DeclareFieldFormat{eprint:euclid}{%
  Project Euclid\addcolon\space
  \ifhyperref
    {\href{http://projecteuclid.org/euclid.lnl/#1}{\nolinkurl{#1}}}
    {\nolinkurl{#1}}}

\SetLabelAlign{parright}{\parbox[t]{\labelwidth}{\raggedleft#1}}
\SetLabelAlign{pushright}{\makebox[{\ifdim\labelwidth<\width\width\else\labelwidth\fi}]%
  {\noindent\raggedleft\hfill #1}}
\setlist[description]{style=sameline,font=\normalfont,leftmargin=2.5em,align=pushright}
\firmlists

\declaretheorem[numbered=no]{theorem}
\declaretheorem[numberwithin=section,name=Theorem]{theoremnum}

\declaretheorem[numbered=no]{lemma}
\declaretheorem[numberwithin=section,name=Lemma]{lemmanum}

\declaretheorem[numbered=no]{proposition}
\declaretheorem[numberwithin=section,name=Proposition]{propositionnum}

\declaretheorem[numbered=no]{observation}
\declaretheorem[numbered=no]{corollary}
\declaretheorem[numberwithin=section,name=Corollary]{corollarynum}

\declaretheorem[numbered=no,style=definition]{definition}
\declaretheorem[numberwithin=section,style=definition,name=Definition]{definitionnum}

\declaretheorem[numbered=no,style=definition]{example}
\declaretheorem[numbered=no,style=definition]{remark}
\declaretheorem[numberwithin=section,style=definition,name=Remark]{remarknum}

\declaretheorem[numbered=no,style=definition]{remarks}
\declaretheorem[numberwithin=section,style=definition,name=Remarks]{remarksnum}

\declaretheorem[numbered=no,style=definition]{question}
\declaretheorem[numberwithin=section,style=definition,name=Question]{questionnum}
\renewcommand*{\thequestionnum}{\arabic{questionnum}}
\declaretheorem[numbered=no,style=definition]{comments}
\ignoretheorems{definition,definitionnum,example,examples,%
  remark,remarks,remarknum,remarksnum,%
  question,questionnum,comments}
\declaretheorem[numbered=no,name=Homotopy lemma]{homotopylemma}
\declaretheorem[numbered=no,name=Homotopy lemma reformulated]{homotopylemmareformulated}
\declaretheorem[numbered=no,name=Comparison lemma for homotopy intervals]{comparisonlemmaforHI}
\declaretheorem[numbered=no,name=``Scholie'']{scholie}

\usepackage{marginnote}

\DeclareMathSymbol{.}{\mathpunct}{letters}{"3A}
\DeclareMathSymbol{\decimalperiod}{\mathord}{letters}{"3A}

\newcommand{\thetitlepage}{{%
  \thispagestyle{empty}
  \calccentering{\unitlength}
  \begin{adjustwidth*}{\unitlength}{-\unitlength}
    \centering
    {\huge \textls[250]{PURSUING STACKS}}\par
    \vspace{\baselineskip}
    {\Large (\kern1pt\emph{À la poursuite des Champs})}\par
    \vspace{5\baselineskip}
    First episode:\par
    \vspace{\baselineskip}
    {\large \textls[200]{THE MODELIZING STORY}}\par
    \vspace{\baselineskip}
    (\kern0.4pt\emph{histoire de modèles})\par
    \vspace{5\baselineskip}
    by
    \vspace{5\baselineskip}
    {\Large Alexander Grothendieck}\par
    \vspace{\baselineskip}
    1983\par
  \end{adjustwidth*}
  \vfill\clearpage}}

\newcommand{\thecopyrightpage}{{%
  \raggedright
  Pursuing Stacks (À la poursuite des Champs)\par
  First episode: The modelizing story (histoire de modèles)\par
  by Alexander Grothendieck (* March 28, 1918 in Berlin; \textdagger{}
  November 13, 2014 Lasserre, Ariège)\par
  \vspace*{\baselineskip}\par
  \noindent The 2015 edition (including all source files) by
  \emph{the scrivener} is in the public domain.\par
  \url{https://github.com/thescrivener/PursuingStacks}\par
  \vspace*{\baselineskip}\par
  \noindent  This edition is an extended version by Mateo Carmona with the collaboration of Ulrik Buchholtz.
  Remarks, comments, and corrections are welcome.\par
  \url{https://agrothendieck.github.io/}\par
  \vspace*{5\baselineskip}\par
  Typeset with \LaTeX{} and the fabulous \emph{memoir} class.
  Composed with Bitstream Charter for the text type and Paul
  Pichaureau's \texttt{mathdesign} package for the math type. All
  diagrams and illustrations were made with Till Tantau's Ti{\itshape
    k}Z (Ti{\itshape k}Z ist \emph{kein} Zeichenprogramm) and the
  \texttt{tikz-cd} package.\par
  \clearpage}}

\setrmarg{0pt plus 1fil}
\setpnumwidth{2.55em}

\OPTstockformat
\OPTpageformat
\settypeblocksize{*}{24.5pc}{*}
\setlrmargins{1in}{*}{*}
\setulmarginsandblock{1in}{2in}{*}
\setheadfoot{\baselineskip}{2\baselineskip}
\setheaderspaces{*}{2\baselineskip}{*}
\setmarginnotes{2pc}{12pc}{\baselineskip}
\checkandfixthelayout
\footnotesinmargin

\sideparmargin{right}
\newlength{\extrawidth}
\newenvironment{widematter}{\begin{adjustwidth*}{0mm}{-\extrawidth}}%
  {\end{adjustwidth*}}

\addtopsmarks{companion}{}{
  \createmark{chapter}{both}{shownumber}{}{\enspace}
  \createmark{section}{right}{shownumber}{\S}{}
}
\strictpagecheck
\newcounter{mastersection} 
\hangsecnum
\counterwithout{section}{chapter}
\makeatletter
\@addtoreset{section}{mastersection}
\newcommand*{\kernifitalic}[1]{%
  \ifx\f@shape\my@test@it
    \kern#1
  \else\relax\fi
}
\newcommand*{\my@test@it}{it}
\def\hangsection{\@dblarg\@hangsection}
\def\@hangsection[#1]#2{%
  \leavevmode\par\needspace{3\baselineskip}\everypar{}\noindent%
  \setcounter{footnote}{0}%
  \refstepcounter{section}%
  \llap{\normalfont\bfseries\large\thesection\quad}%
  \phantomsection%
  \markright{\memUChead{\S\thesection\enspace #1}}%
  \addcontentsline{toc}{section}{\numberline{\thesection}#2}%
  \checkoddpage%
  \marginpar{\if@twoside\ifoddpage\relax\else{$\phantom{1}$\\}\fi\fi\itshape #2}%
  \ignorespaces%
}
\def\presectionfill{\noindent\needspace{8\baselineskip}\hfill}
\makeatother

\newcommand{\DontChangeNextSectionNumber}{%
  \edef\sectioncountervalue{\numexpr\the\value{section}-1\relax}%
  \refstepcounter{mastersection}%
  \setcounter{section}{\sectioncountervalue}%
}
\setaftersubsecskip{-1em}
\maxsecnumdepth{subsection}

\makeatletter
\def\namedlabel#1#2{\begingroup
  #2%
  \def\@currentlabel{#2}%
  \phantomsection\label{#1}\endgroup
}
\def\phantomlabel#1#2{\begingroup
  \def\@currentlabel{#2}%
  \phantomsection\label{#1}\endgroup
}
\makeatother

\makeatletter
\newcommand{\namedpspage}[2]{%
  \phantomlabel{p:#1}{#2}%
  \checkoddpage%
  \marginpar{%
    \if@twoside\ifoddpage\relax\else\noindent\hfill\fi\fi%
    [p.\ #2]}}
\makeatother
\newcommand{\pspage}[1]{\namedpspage{#1}{#1}}
\newcommand{\scrcomment}[1]{\marginpar{\color{DarkRed}{[#1]}}}
\newcommand{\scrcommentinline}[1]{%
  {\color{DarkRed}{[#1]}}}
\newcommand{\ondate}[1]{\namedlabel{date:#1}{#1}}
\newcommand{\alsoondate}[1]{#1}

\mathchardef\Hat="705E
\mathchardef\Vee="7014
\newcommand*{\uphat}{^{\vbox to 0.8ex{\hbox{\large$\Hat$}\vss}}}
\newcommand*{\adjuphat}{^{\vbox to 0.8ex{\hbox{\!\large$\Hat$}\vss}}}
\newcommand*{\upvee}{^{\vbox to 0.8ex{\hbox{\large$\Vee$}\vss}}}

\makeatletter
\newcommand*\wt[2][0.1ex]{%
        \begingroup
        \mathchoice{\wt@helper{#1}{#2}{\displaystyle}{\textfont}}
                   {\wt@helper{#1}{#2}{\textstyle}{\textfont}}
                   {\wt@helper{#1}{#2}{\scriptstyle}{\scriptfont}}
                   {\wt@helper{#1}{#2}{\scriptscriptstyle}{\scriptscriptfont}}%
        \endgroup
        #2%
}
\newcommand*\wt@helper[4]{%
        \def\currentfont{\the#41}%
        \def\currentskewchar{\char\the\skewchar\currentfont}%
        \setbox\tw@\hbox{\currentfont#2\currentskewchar}%
        \dimen@ii\wd\tw@
        \setbox\tw@\hbox{\currentfont#2{}\currentskewchar}%
        \advance\dimen@ii-\wd\tw@
        \rlap{\raisebox{-#1}{$\m@th#3\kern\dimen@ii\widetilde{\phantom{#2}}$}}%
}
\makeatother

\newcommand*{\starsbreak}{\fancybreak{\raisebox{3pt}[18pt][6pt]*\kern1em
    *\kern1em\raisebox{3pt}*}}
\newcommand*{\preslice}[2]{\vphantom{#1}_{#2\backslash}{#1}} 
\newcommand*{\cst}[1]{\mathrm{#1}} 
\newcommand*{\bN}{{\ensuremath{\mathbf{N}}}}

\newcommand*{\bZ}{{\ensuremath{\mathbf{Z}}}}

\newcommand*{\bP}{{\ensuremath{\mathbb{P}}}}

\newcommand*{\bQ}{{\ensuremath{\mathbf{Q}}}}

\newcommand*{\bR}{{\ensuremath{\mathbb{R}}}}
\newcommand*{\bC}{{\ensuremath{\mathbf{C}}}}
\newcommand*{\bB}{{\ensuremath{\mathbf{B}}}}
\newcommand*{\bF}{{\ensuremath{\mathbb{F}}}}
\newcommand*{\bG}{{\ensuremath{\mathbb{G}}}}
\newcommand*{\bI}{{\ensuremath{\mathbb{I}}}}
\newcommand*{\bJ}{{\ensuremath{\mathbb{J}}}}
\newcommand*{\bK}{{\ensuremath{\mathbb{K}}}}
\newcommand*{\bL}{{\ensuremath{\mathbb{L}}}}
\newcommand*{\bM}{{\ensuremath{\mathbb{M}}}}
\newcommand*{\scrA}{{\ensuremath{\mathscr{A}}}}
\newcommand*{\scrAab}{{\ensuremath{\mathscr{A}\subab}}}
\newcommand*{\scrB}{{\ensuremath{\mathscr{B}}}}
\newcommand*{\scrC}{{\ensuremath{\mathscr{C}}}}
\newcommand*{\scrD}{{\ensuremath{\mathscr{D}}}}
\newcommand*{\scrE}{{\ensuremath{\mathscr{E}}}}
\newcommand*{\scrF}{{\ensuremath{\mathscr{F}}}}
\newcommand*{\scrH}{{\ensuremath{\mathscr{H}}}}
\newcommand*{\scrL}{{\ensuremath{\mathscr{L}}}}
\newcommand*{\scrM}{{\ensuremath{\mathscr{M}}}}
\newcommand*{\scrN}{{\ensuremath{\mathscr{N}}}}
\newcommand*{\scrO}{{\ensuremath{\mathscr{O}}}}
\newcommand*{\scrOX}{{\ensuremath{\mathscr{O}_{\mathscr X}}}}
\newcommand*{\scrOk}{{\ensuremath{\mathscr{O}_k}}}
\newcommand*{\scrP}{{\ensuremath{\mathscr{P}}}}
\newcommand*{\scrQ}{{\ensuremath{\mathscr{Q}}}}
\newcommand*{\scrU}{{\ensuremath{\mathscr{U}}}}
\newcommand*{\scrW}{{\ensuremath{\mathscr{W}}}}
\newcommand*{\scrWA}{{\ensuremath{\mathscr{W}_A}}}
\newcommand*{\scrWB}{{\ensuremath{\mathscr{W}_B}}}
\newcommand*{\scrWoo}{{\ensuremath{\mathscr{W}_\euinfty}}}
\newcommand*{\scrWz}{{\ensuremath{\mathscr{W}_0}}}
\newcommand*{\scrX}{{\ensuremath{\mathscr{X}}}}
\newcommand*{\Ahat}{{\ensuremath{A\uphat}}}
\newcommand*{\Ahatc}{{\ensuremath{A\uphat\subc}}}
\newcommand*{\Ahata}{{\ensuremath{A\uphat\suba}}}
\newcommand*{\Ahatas}{{\ensuremath{A\uphat\subas}}}
\newcommand*{\Ahatab}{{\ensuremath{A\uphat\subab}}}
\newcommand*{\AhatM}{{\ensuremath{A\uphat_\scrM}}}
\newcommand*{\AhatN}{{\ensuremath{A\uphat_\scrN}}}
\newcommand*{\Ahatk}{{\ensuremath{A\uphat_k}}}
\newcommand*{\Ahatfp}{{\ensuremath{A\uphat\subfp}}}
\newcommand*{\Ahattotas}{{\ensuremath{A\uphat_{\mathrm{tot.as}}}}}
\newcommand*{\Ahatlocas}{{\ensuremath{A\uphat_{\mathrm{loc.as}}}}}
\newcommand*{\Bhat}{{\ensuremath{B\uphat}}}
\newcommand*{\Bhatas}{{\ensuremath{B\uphat\subas}}}
\newcommand*{\Bhatab}{{\ensuremath{B\uphat\subab}}}
\newcommand*{\BhatM}{{\ensuremath{B\uphat_\scrM}}}
\newcommand*{\Bhatk}{{\ensuremath{B\uphat_k}}}
\newcommand*{\Bhatc}{{\ensuremath{B\uphat\subc}}}
\newcommand*{\Chat}{{\ensuremath{C\uphat}}}
\newcommand*{\PampM}{{\ensuremath{P\supamp_\scrM}}}
\newcommand*{\PampN}{{\ensuremath{P\supamp_\scrN}}}
\newcommand*{\Simplex}{{\text{\raisebox{-0.1ex}{\includegraphics[height=1.7ex]{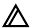}}}}}
\newcommand*{\Square}{{\text{\raisebox{-0.1ex}{\includegraphics[height=1.5ex]{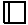}}}}}
\newcommand*{\Globe}{{\text{\raisebox{-0.2ex}{\includegraphics[height=1.8ex]{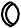}}}}}
\newcommand*{\Simplexop}{{\ensuremath{\Simplex\!\op}}}
\newcommand*{\Simplexhat}{{\ensuremath{\Simplex\adjuphat}}}
\newcommand*{\Simplexhatab}{{\ensuremath{\Simplex\subab\adjuphat}}}
\newcommand*{\Simplextildehat}{{\ensuremath{\wt[0.2ex]{\Simplex}\adjuphat}}}
\newcommand*{\Simplexfhat}{{\ensuremath{{{\Simplex\!}^{\mathrm f}}\uphat}}}
\newcommand*{\sor}{\mathbin{\cup}}
\newcommand*{\sand}{\mathbin{\cap}}
\newcommand*{\suba}{_{\mathrm{a}}}
\newcommand*{\subm}{_{\mathrm{m}}}
\newcommand*{\subc}{_{\mathrm{c}}}
\newcommand*{\subd}{_{\mathrm{d}}} 
\newcommand*{\subs}{_{\mathrm{s}}} 
\newcommand*{\subab}{_{\mathrm{ab}}}
\newcommand*{\subas}{_{\mathrm{as}}}
\newcommand*{\subfp}{_{\mathrm{fp}}}
\newcommand*{\subpt}{_{\mathrm{pt}}}
\newcommand*{\supamp}{^{\&}}
\newcommand*{\supab}{^{\mathrm{ab}}}
\newcommand*{\supDP}{^{\mathrm{DP}}}
\newcommand*{\Ab}{{\ensuremath{(\mathrm{Ab})}}}
\newcommand*{\Hot}{{\ensuremath{(\mathrm{Hot})}}}
\newcommand*{\Hotab}{{\ensuremath{(\mathrm{Hotab})}}}
\newcommand*{\HotA}{{\ensuremath{(\mathrm{Hot}_A)}}}
\newcommand*{\HotB}{{\ensuremath{(\mathrm{Hot}_B)}}}
\newcommand*{\HotM}{{\ensuremath{(\mathrm{Hot}_M)}}}
\newcommand*{\HotW}{{\ensuremath{(\mathrm{Hot}_\scrW)}}}
\newcommand*{\HotAW}{{\ensuremath{(\mathrm{Hot}_A^\scrW)}}}
\newcommand*{\Hotz}{{\ensuremath{\mathrm{Hot}_0}}}
\newcommand*{\Hotabz}{{\ensuremath{\mathrm{Hotab}_0}}}
\newcommand*{\Sets}{{\ensuremath{(\mathrm{Sets})}}}
\newcommand*{\pSets}{{\ensuremath{(\mathrm{Sets}^\bullet)}}}
\newcommand*{\Gsets}{{\ensuremath{(G\textup{-sets})}}}
\newcommand*{\Sssets}{{\ensuremath{(\mathrm{Ss~sets})}}}
\newcommand*{\Spaces}{{\ensuremath{(\mathrm{Spaces})}}}
\newcommand*{\Ord}{{\ensuremath{(\mathrm{Ord})}}}
\newcommand*{\Preord}{{\ensuremath{(\mathrm{Preord})}}}
\newcommand*{\kMod}{{\ensuremath{(k\textup{-Mod})}}}
\newcommand*{\Rings}{{\ensuremath{(\mathrm{Rings})}}}
\newcommand*{\Cat}{{\ensuremath{(\mathrm{Cat})}}}
\newcommand*{\oCat}{{\ensuremath{\overline{\vphantom{W}\smash{(\mathrm{Cat})}}}}}
\newcommand*{\cTop}{{\ensuremath{(\mathrm{Top})}}}
\newcommand*{\Topoi}{{\ensuremath{(\mathrm{Topoi})}}}
\newcommand*{\Simplextilde}{{\ensuremath{\wt[0.2ex]{\Simplex}}}}
\newcommand*{\Simplexf}{{\ensuremath{\Simplex\!^{\mathrm{f}}}}}
\newcommand*{\smashSimplexf}{{\ensuremath{\smash{\Simplex\!^{\mathrm{f}}}}}}
\newcommand*{\Simplexprimef}{{\ensuremath{{\Simplex\!'}^{\mathrm{f}}}}}
\newcommand*{\Simplextildef}{{\ensuremath{\wt[0.2ex]{\Simplex}\!^{\mathrm{f}}}}}
\newcommand*{\iSimplexIntoCat}{{\ensuremath{i:\Simplex\hookrightarrow\Cat}}}
\newcommand*{\piz}{{\ensuremath{\pi_0}}}
\newcommand*{\bpi}{{\ensuremath{\boldsymbol\pi}}}
\newcommand*{\bpiz}{{\ensuremath{\boldsymbol\pi_0}}}
\newcommand*{\Isim}{\mathrel{\underset{\bI}{\sim}}}
\newcommand*{\Izsim}{\mathrel{\underset{\bI_0}{\sim}}}
\newcommand*{\Jsim}{\mathrel{\underset{\bJ}{\sim}}}
\newcommand*{\Ksim}{\mathrel{\underset{\bK}{\sim}}}
\newcommand*{\Wsim}{\mathrel{\underset{W}{\sim}}}
\newcommand*{\hsim}{\mathrel{\underset{h}{\sim}}}
\newcommand*{\Rsim}{\mathrel{\underset{R}{\sim}}}

\newcommand*{\oo}{\ensuremath{\euinfty}}

\newcommand*{\op}{^{\mathrm{op}}}

\newcommand*{\restrto}{\ensuremath{\mathbin{\!\upharpoonright}}}
\newcommand*{\equ}{\mathrel{\approx}}
\newcommand*{\equeq}{\mathrel{\approxeq}}
\newcommand*{\tosim}{\xrightarrow{\raisebox{-2.5pt}[0pt][0pt]{$\scriptstyle\sim$}}}
\newcommand*{\fromsim}{\xleftarrow{\raisebox{-2.5pt}[0pt][0pt]{$\scriptstyle\sim$}}}
\newcommand*{\tosimeq}{\xrightarrow{\raisebox{-1.5pt}[0pt][0pt]{$\scriptstyle\simeq$}}}
\newcommand*{\toequ}{\xrightarrow{\raisebox{-1.5pt}[0pt][0pt]{$\scriptstyle\approx$}}}
\newcommand*{\fromequ}{\xleftarrow{\raisebox{-1.5pt}[0pt][0pt]{$\scriptstyle\approx$}}}
\newcommand*{\UW}{{\ensuremath{\mathrm UW}}}
\newcommand*{\GalQQ}{{\ensuremath{\mathrm{Gal}_{\overline{\bQ}/\bQ}}}}
\newcommand*{\Last}{\mathbin{\overset{\mathrm L}{*}}}
\newcommand*{\Loast}{\mathbin{\overset{\mathrm L}{\oast}}}
\newcommand*{\hatotimes}{\mathbin{\hat{\otimes}}}
\DeclareMathOperator{\Add}{Add}
\DeclareMathOperator{\Addinf}{Addinf}

\DeclareMathOperator{\Ker}{Ker}
\DeclareMathOperator{\Coker}{Coker}
\DeclareMathOperator{\End}{End}
\DeclareMathOperator{\Aut}{Aut}
\DeclareMathOperator{\Ext}{Ext}
\DeclareMathOperator{\bExt}{\mathbf{Ext}}
\DeclareMathOperator{\Tor}{Tor}
\DeclareMathOperator{\Spec}{Spec}
\DeclareMathOperator{\Qucoh}{Qucoh}
\DeclareMathOperator{\Bil}{Bil}
\DeclareMathOperator{\Kar}{Kar}
\DeclareMathOperator{\KarAdd}{KarAdd}
\DeclareMathOperator{\KarAddinf}{KarAddinf}
\DeclareMathOperator{\Hom}{Hom}
\DeclareMathOperator{\Homadd}{Homadd}

\DeclareMathOperator{\Homcont}{Homcont}
\DeclareMathOperator{\Maps}{Maps}
\DeclareMathOperator{\Mod}{Mod}
\DeclareMathOperator{\RHom}{RHom}
\DeclareMathOperator{\RGamma}{R\Gamma}
\DeclareMathOperator{\RH}{RH}
\DeclareMathOperator{\LH}{LH}
\DeclareMathOperator{\Lpi}{L\pi}
\DeclareMathOperator{\tH}{\widetilde H}
\DeclareMathOperator{\LtH}{L\widetilde H}
\DeclareMathOperator{\bBiadd}{\mathbf{Biadd}}
\DeclareMathOperator{\bHom}{\mathbf{Hom}}
\DeclareMathOperator{\RbHom}{R\mathbf{Hom}}
\DeclareMathOperator{\bHomex}{\mathbf{Homex}}
\DeclareMathOperator{\oHom}{\overline{Hom}}
\DeclareMathOperator{\bHomadd}{\mathbf{Homadd}}
\DeclareMathOperator{\bHomaddinf}{\mathbf{Homaddinf}}
\DeclareMathOperator{\bHomaddkar}{\mathbf{Homaddkar}}
\DeclareMathOperator{\bHomaddinfkar}{\mathbf{Homaddinfkar}}
\DeclareMathOperator{\bHommultinf}{\mathbf{Hommultinf}}

\DeclareMathOperator{\bAsph}{\mathbf{Asph}}
\DeclareMathOperator{\Int}{Int}
\DeclareMathOperator{\Lie}{Lie}
\DeclareMathOperator{\bInt}{\mathbf{Int}}
\DeclareMathOperator{\bEnd}{\mathbf{End}}
\DeclareMathOperator{\Crib}{Crib}
\DeclareMathOperator{\Comp}{Comp}
\DeclareMathOperator{\Homtp}{Homtp}
\DeclareMathOperator{\Cont}{Cont}
\DeclareMathOperator{\Contr}{Contr}
\DeclareMathOperator{\HotOf}{Hot}
\DeclareMathOperator{\HotabOf}{Hotab}
\DeclareMathOperator{\sKM}{sK\mathscr M}
\DeclareMathOperator{\AbOf}{Ab}
\DeclareMathOperator{\Alg}{Alg}
\DeclareMathOperator{\Aff}{Aff}
\DeclareMathOperator{\WAsph}{\mathscr{W}\mathrm{-Asph}}
\DeclareMathOperator{\WzAsph}{\mathscr{W}_0\mathrm{-Asph}}
\DeclareMathOperator{\WooAsph}{\mathscr{W}_\euinfty\mathrm{-Asph}}
\DeclareMathOperator{\WprimeAsph}{\mathscr{W}'\mathrm{-Asph}}
\newcommand*{\HotOfzc}{\ensuremath{\HotOf(0\textup{-}\mathrm{conn})}}
\newcommand*{\HotabOfzc}{\ensuremath{\HotabOf(0\textup{-}\mathrm{conn})}}
\newcommand*{\zc}{\ensuremath{(0\textup{-}\mathrm{conn})}}
\newcommand*{\D}{{\ensuremath{\mathrm{D}}}} 
\newcommand*{\Lotimes}{\mathbin{\overset{\mathrm{L}}{\otimes}}} 
\newcommand*{\DP}{{\ensuremath{\mathrm{DP}}}} 
\newcommand*{\Ch}{{\ensuremath{\mathrm{Ch}}}} 
\DeclareMathOperator{\Sd}{Sd}

\DeclareMathOperator{\Sh}{Sh}
\DeclareMathOperator{\Top}{Top}
\DeclareMathOperator{\Ob}{Ob}
\DeclareMathOperator{\Fl}{Fl}
\DeclareMathOperator{\bFl}{\mathbf{Fl}}
\DeclareMathOperator{\Pro}{Pro}
\DeclareMathOperator{\Proj}{Proj}
\DeclareMathOperator{\UlProj}{UlProj}

\DeclareMathOperator{\id}{id}

\DeclareMathOperator{\Imm}{Im}
\DeclareMathOperator{\inv}{inv}
\DeclareMathOperator{\pr}{pr}
\DeclareMathOperator{\Sym}{Sym}
\DeclareMathOperator{\Symbin}{Symbin}
\DeclareMathOperator{\Wh}{Wh}
\DeclareMathOperator{\bWh}{\mathbf{Wh}}
\DeclareMathOperator{\forg}{forg}
\DeclarePairedDelimiter\abs{\lvert}{\rvert} 
\DeclarePairedDelimiter\angled{\langle}{\rangle} 
\DeclarePairedDelimiterX\set[2]{\{}{\}}
  {#1 \mathrel{}\mathclose{}\delimsize|\mathopen{}\mathrel{} #2}

\newcommand*{\eqdef}{\overset{\textup{def}}{=}}

\DeclareMathOperator{\Gr}{Gr}
\DeclareMathOperator{\B}{B}
\DeclareMathOperator{\Br}{Br}
\DeclareMathOperator{\Lf}{Lf}
\DeclareMathOperator{\Rf}{Rf}
\DeclareMathOperator{\Ens}{Ens}
\DeclareMathOperator{\Alb}{Alb}

\begin{document}

\frontmatter
\thetitlepage
\thecopyrightpage

\renewcommand*{\contentsname}{Contents}
\setcounter{tocdepth}{2} 
\phantomsection%
\pdfbookmark[0]{Contents}{contents}%
\tableofcontents*



\chapter{Preface}

Since the month of March last year, so nearly a year ago, the greater part of my energy has been devoted to a work of reflection on the \emph{foundations of non-commutative (co)homological algebra}, or what is the same, after all, of \emph{homotopical algebra}. These reflections have taken the concrete form of a voluminous stack of typed notes, destined to form the first volume (now being finished) of a work in two volumes to be published by Hermann, under the overall title \emph{``Pursuing Stacks''}. I now foresee (after successive extensions of the initial project) that the manuscript of the whole of the two volumes, which I hope to finish definitively in the course of this year, will be about 1500 typed pages in length. These two volumes are moreover for me the first in a vaster series, under the overall title \emph{``Mathematical Reflections''}, in which I intend to develop some of the themes sketched in the present report.

Since I am speaking here of work which is actually now being written up and is even almost finished, the first volume of which will doubtless appear this year and will contain a detailed introduction, it is undoubtedly less interesting for me to develop this theme of reflection here, and I will content myself with speaking of it only very briefly. This work seems to me to be somewhat marginal with respect to the themes I sketched before, and does not (it seems to me) represent a real renewal of viewpoint or approach with respect to my interests and my mathematical vision of before 1970. If I suddenly resolved to do it, it is almost out of desperation, for nearly twenty years have gone by since certain visibly fundamental questions, which were ripe to be thoroughly investigated, without anyone seeing them or taking the trouble to fathom them. Still today, the basic structures which occur in the homotopical point of view in topology are not understood, and to my knowledge, after the work of Verdier, Giraud and Illusie on this theme (which are so many beginnings still waiting for continuations\dots) there has been no effort in this direction. I should probably make an exception for the axiomatisation work done by Quillen on the notion of a category of models, at the end of the sixties, and taken up in various forms by various authors. At that time, and still now, this work seduced me and taught me a great deal, even while going in quite a different direction from the one which was and still is close to my heart. Certainly, it introduces derived categories in various non-commutative contexts, but without entering into the question of the essential internal structures of such a category, also left open in the commutative case by Verdier, and after him by Illusie. Similarly, the question of putting one’s finger on the natural “coefficients” for a non-commutative cohomological formalism, beyond the stacks (which should be called 1-stacks) studied in the book by Giraud\scrcomment{cf. \textcite{Giraud1971}}, remained open – or rather, the rich and precise intuitions concerning it, taken from the numerous examples coming in particular from algebraic geometry, are still waiting for a precise and supple language to give them form.

I returned to certain aspects of these foundational questions in 1975, on the occasion (I seem to remember) of a correspondence with Larry Breen ([three] letters from this correspondence will be reproduced as an appendix to Chap. I of volume 1, ``History of Models'', of Pursuing Stacks). At that moment the intuition appeared that $\infty$-groupoids should constitute particularly adequate models for homotopy types, the $n$-groupoids corresponding to \emph{truncated} homotopy types (with $\pi_i = 0$ pour $i > n$). This same intuition, via very different routes, was discovered by Ronnie Brown and some of his students in Bangor, but using a rather restrictive notion of $\infty$-groupoid (which, among the 1-connected homotopy types, model only products of Eilenberg-Mac Lane spaces). Stimulated by a rather haphazard correspondence with Ronnie Brown\scrcomment{cf. \textcite{Kunzer2015}}, I finally began this reflection, starting with an attempt to define a wider notion of $\infty$-groupoid (later rebaptised stack in $\infty$-groupoids or simply “stack”, the implication being: over the 1-point topos), and which, from one thing to another, led me to Pursuing Stacks. The volume ``History of Models'' is actually a completely unintended digression with respect to the initial project (the famous stacks being temporarily forgotten, and supposed to reappear only around page 1000\dots).

This work is not completely isolated with respect to my more recent interests. For example, my reflection on the modular multiplicities $\widehat{M}_{g, \nu}$ and their stratified structure renewed the reflection on a theorem of van Kampen in dimension > 1 (also one of the preferred themes of the group in Bangor), and perhaps also contributed to preparing the ground for the more important work of the following year. This also links up from time to time with a reflection dating from the same year 1975 (or the following year) on a ``De Rham complex with divided powers''\scrcomment{cf. \textcite{AG76IHES}}, which was the subject of my last public lecture, at the IHES in 1976; I lent the manuscript of it to I don’t remember whom after the talk, and it is now lost. It was at the moment of this reflection that the intuition of a ``schematisation'' of homotopy types germinated, and seven years later I am trying to make it precise in a (particularly hypothetical) chapter of the History of Models.

The work of reflection undertaken in Pursuing Stacks is a little like a debt which I am paying towards a scientific past where, for about fifteen years (from 1955 to 1970), the development of cohomological tools was the constant Leitmotiv in my foundational work on algebraic geometry. If in this renewal of my interest in this theme, it has taken on unexpected dimensions, it is however not out of pity for a past, but because of the numerous unexpected phenomena which ceaselessly appear and unceremoniously shatter the previously laid plans and projects – rather like in the thousand and one nights, where one awaits with bated breath through twenty other tales the final end of the first.

\bigskip
\noindent\hfill\emph{Translation of a section in ``Esquisse d'un Programme'' of 1984 \scrcomment{\textcite{LochakSchneps2000}}}\par


\mainmatter


\chapter{Take-off}
\label{ch:I}

\par\hfill Les Aumettes \ondate{19.2.}1983\namedpspage{L.1}{1}\par

Dear Daniel\scrcomment{Daniel Quillen},

\hangsection{The importance of innocence.}\label{sec:1}%
Last year Ronnie Brown from Bangor\scrcomment{cf.\ \textcite{Kunzer2015}} sent me a heap of reprints and
preprints by him and a group of friends, on various foundational
matters of homotopical algebra. I did not really dig through any of
this, as I kind of lost contact with the technicalities of this kind
(I was never too familiar with the homotopy techniques anyhow, I
confess) -- but this reminded me of a few letters I had exchanged with
Larry Breen in 1975\footnote{These letters are reproduced as an 
``appendix'' at the end of this chapter.}, where I had developed an outline of a program for
a kind of ``topological algebra'', viewed as a synthesis of
homotopical and homological algebra, with special emphasis on topoi --
most of the basic intuitions in this program arising from various
backgrounds in algebraic geometry. Some of those intuitions we
discussed, I believe, at IHES eight or nine years before\scrcomment{\textcite{tapisQ}}, then, at a time
when you had just written up your nice ideas on axiomatic homotopical
algebra,\scrcomment{\textcite{Quillen1967}} published since in
Springer's Lecture Notes. I write you under the assumption that you
have not entirely lost interest for those foundational questions you
were looking at more than fifteen years ago. One thing which strikes
me, is that (as far as I know) there has not been any substantial
progress since -- it looks to me that an understanding of the basic
structures underlying homotopy theory, or even homological algebra
only, is still lacking -- probably because the few people who have a
wide enough background and perspective enabling them to feel the main
questions, are devoting their energies to things which seem more
directly rewarding. Maybe even a wind of disrepute for any
foundational matters whatever is blowing nowadays\footnote{When making this suggestion about there being a ``wind of disrepute for any foundational matters whatever'', I little suspected that the former friend to whom I was communicating my ponderings as they came, would take care of providing a most unexpected confirmation. As a matter of fact, this letter never got an answer, nor was it even read! Upon my inquiry nearly one year later, this colleague appeared sincerely surprised that I could have expected even for a minute that he might possibly read a letter of mine on mathematical matters, well knowing the kind of ``general nonsense'' mathematics that was to be expected from me\dots}!  In this respect,
what seems to me even more striking than the lack of proper
foundations for homological and homotopical algebra, is the absence I
daresay of proper foundations for topology itself! I am thinking here
mainly of the development of a context of ``tame'' topology, which (I
am convinced) would have on the everyday technique of geometric
topology (I use this expression in contrast to the topology of use for
analysts) a comparable impact or even a greater one, than the
introduction of the point of view of schemes had on algebraic
geometry\footnote{For some particulars about a program of ``tame topology'', I refer to ``Esquisse d'un Programme'', sections 5 and 6, which is included in Réflexions Mathématiques 1.}. The psychological drawback here I believe is not anything
like messiness, as for homological and homotopical algebra (as for
schemes), but merely the inrooted inertia which prevents us so
stubbornly from looking innocently, with fresh eyes, upon things,
without being dulled and imprisoned by standing habits of thought,
going with a familiar context -- \emph{too} familiar a context! The
task of working out the foundations of tame topology, and a
corresponding structure theory for ``stratified (tame) spaces'', seems
to me a lot more urgent and exciting still
than\namedpspage{L.1prime}{1'} any program of homological,
homotopical or topological algebra.

\hangsection{A short look into purgatory.}\label{sec:2}%
The motivation for this letter was the latter topic
however. Ronnie Brown and his friends are competent algebraists and
apparently strongly motivated for investing energy in foundational
work, on the other hand they visibly are lacking the necessary scope
of vision which geometry alone provides\footnote{I have to apologise for this rash statement, as later correspondence made me realise that ``Ronnie Brown and his friends'' do have stronger contact with ``geometry'' than I suspected, even though they are not too familiar with algebraic geometry!}. They seem to me kind of
isolated, partly due I guess to the disrepute I mentioned before -- I
suggested to try and have contact with people such as yourself, Larry
Breen, Illusie and others, who have the geometric insight and who
moreover, may not think themselves too good for indulging in
occasional reflection on foundational matters and in the process help
others do the work which should be done.

At first sight it has seemed to me that the Bangor group\footnote{The ``Bangor group'' is made up by Ronnie Brown and Tim Porter as the two fixed points, and a number of devoted research students. Moreover Ronnie Brown is working in close contact with J.L. Loday and J. Pradines in France.} had indeed
come to work out (quite independently) one basic intuition of the
program I had envisioned in those letters to Larry Breen -- namely
that the study of $n$-truncated homotopy types (of semisimplicial
sets, or of topological spaces) was essentially equivalent to the
study of so-called $n$-groupoids (where $n$ is any natural
integer). This is expected to be achieved by associating to any space
(say) $X$ its ``fundamental $n$-groupoid'' $\Pi_n(X)$, generalizing
the familiar Poincar\'e fundamental groupoid for $n=1$. The obvious
idea is that $0$-objects of $\Pi_n(X)$ should be the points of $X$,
$1$-objects should be ``homotopies'' or paths between points,
$2$-objects should be homotopies between $1$-objects, etc. This
$\Pi_n(X)$ should embody the $n$-truncated homotopy type of $X$, in
much the same way as for $n=1$ the usual fundamental groupoid embodies
the $1$-truncated homotopy type. For two spaces $X,Y$, the set of
homotopy-classes of maps $X\to Y$ (more correctly, for general $X,Y$,
the maps of $X$ into $Y$ in the homotopy category) should correspond
to $n$-equivalence classes of $n$-functors from $\Pi_n(X)$ to
$\Pi_n(Y)$ -- etc. There are very strong suggestions for a nice
formalism including a notion of geometric realization of an
$n$-groupoid, which should imply that any $n$-groupoid (or more
generally of an $n$-category) is relativized over an arbitrary topos
to the notion of an $n$-gerbe (or more generally, an $n$-stack), these
become the natural ``coefficients'' for a formalism of non-commutative
cohomological algebra, in the spirit of Giraud's thesis.

But all this kind of thing for the time being is pure heuristics -- I
never so far sat down to try to make explicit at least a definition of
$n$-categories and $n$-groupoids, of $n$-functors between these
etc. When I got the Bangor reprints I at once had the feeling that
this kind of work\namedpspage{L.2}{2} had been done and the homotopy category
expressed in terms of \oo-groupoids. But finally it appears this is
not so, they have been working throughout with a notion of
\oo-groupoid too restrictive for the purposes I had in mind (probably
because they insist I guess on strict associativity of compositions,
rather than associativity up to a (given) isomorphism, or rather,
homotopy) -- to the effect that the simply connected homotopy types
they obtain are merely products of Eilenberg-MacLane spaces, too bad!
They do not seem to have realized yet that this makes their set-up
wholly inadequate to a sweeping foundational set-up for homotopy. This
brings to the fore again to work out the suitable definitions for
$n$-groupoids -- if this is not done yet anywhere. I spent the
afternoon today trying to figure out a reasonable definition, to get a
feeling at least of where the difficulties are, if any. I am guided
mainly of course by the topological interpretation. It will be short
enough to say how far I got. The main part of the structure it seems
is expressed by the sets $F_i$ ($i\in\bN$) of $i$-objects, the source,
target and identity maps
\begin{align*}
  \cst s_1^i, \cst t_1^i &: F_i \to F_{i-1} \quad (i\ge 1) \\
  \cst k_1^i &: F_i \to F_{i+1} \quad (i\in\bN) \\
\intertext{and the symmetry map (passage to the inverse)}
  \inv_i &: F_i \to F_i \quad (i\ge 1),
\end{align*}
satisfying\marginpar{\emph{notation:}\\ $d_a=\cst k_1^i(a)$,\\ $\check u =
  \inv_i(u)$} some obvious relations: $\cst k_1^i$ is right inverse to
the source and target maps $\cst s_1^{i+1}, \cst t_1^{i+1}$, $\inv_i$
is an involution and ``exchanges'' source and target, and moreover for
$i\ge2$
\begin{align*}
  \cst s_1^{i-1} \cst s_1^i
  &= \cst s_1^{i-1} \cst t_1^i \;
    \Bigl(\eqdef \cst s_2^i : F_i \to F_{i-2} \Bigr) \\
  \cst t_1^{i-1} \cst s_1^i
  &= \cst t_1^{i-1} \cst t_1^i \;
    \Bigl(\eqdef \cst t_2^i : F_i \to F_{i-2} \Bigr) ;
\end{align*}
thus the composition of the source and target maps yields, for $0\le
j\le i$, just \emph{two} maps
\begin{equation*}
  \cst s_\ell^i, \cst t_\ell^i : F_i \to F_{i-\ell} = F_j \quad (\ell = i-j).
\end{equation*}
The next basic structure is the composition structure, where the
usual compoisition of arrows, more specifically of $i$-objects
($i\ge1$) $v \circ u$ (defined when $\cst t_1(u)=\cst s_1(v)$) must be
supplemented by the Godement-type operations $\mu * \lambda$ when
$\mu$ and $\lambda$ are ``arrows between arrows'', etc. Following this
line of thought, one gets the composition maps
\begin{equation*}
  (u,v) \mapsto (v *_\ell u) : (F_i, \cst s_\ell^i) \times_{F_{i-\ell}}
  (F_i,\cst s_\ell^i) \to F_i,
\end{equation*}
the composition of $i$-objects for $1\le \ell\le i$, being defined
when the $\ell$-target of $u$ is equal to the $\ell$-source of $v$,
and then we have
\begin{equation*}
  \left.\begin{aligned}
    \cst s_1^i(v *_\ell u) &=  \cst s_1^i(v) *_{\ell-1} \cst s_1^i(u) \\
    \cst t_1^i(v *_\ell u) &=  \cst t_1^i(v) *_{\ell-1} \cst s_1^i(u)
    \end{aligned}\right\}
  \quad \text{$\ell\ge2$ i.e. $\ell-1\ge1$}
\end{equation*}
and\namedpspage{L.2prime}{2'} for $\ell=1$
\begin{align*}
  \cst s_1(v *_1 u) &= s_1(u) \\
  \cst t_1(v *_1 u) &= t_1(v)
\end{align*}
(NB\enspace the operation $v *_1 u$ is just the usual composition $v \circ
u$).

One may be tempted to think that the preceding data exhaust the
structure of \oo-groupoids, and that they will have to be supplemented
only by a handful of suitable \emph{axioms}, one being
\emph{associativity} for the operation
$\overset{i}{\underset{\ell}{*}}$, which can be expressed essentially
by saying that that composition operation turns $F_i$ into the set of
arrows of a category having $F_{i-\ell}$ as a set of objects (with the
source and target maps $\cst s_\ell^i$ and $\cst t_\ell^i$, and with
identity map $\cst k_\ell^{i-\ell} : F_{i-\ell} \to F_i$ the
composition of the identity maps $F_{i-\ell} \to F_{i-\ell+1} \to
\dots \to F_{i-1} \to F_i$), and another being the Godement relation
\begin{equation*}
  (v' *_\alpha v) *_\nu (u' *_\alpha u) = (v' *_\nu u') *_\alpha (v
  *_\nu u)
\end{equation*}
(with the assumptions $1\le\alpha\le\nu$, and $u,u',v,v''in F_i$ and
\begin{equation*}
  \left\{\begin{aligned}
      \cst t_\alpha(u)=\cst s_\alpha(u') \\
      \cst t_\alpha(v)=\cst s_\alpha(v') \\
    \end{aligned}\right.\quad
  \cst t_\nu(u) = \cst s_\nu(v) = \cst s_\nu(v') = \cst t_\nu(u')
\end{equation*}
implying that both members are defined), plus the two relations
concerning the inversion of $i$-objects ($i\ge1$) $u \mapsto \check
u$,
\begin{multline*}
  u *_1 \check u = \id_{\cst t_1(u)} , \qquad \check u *_1 u =
  \id_{\cst s_1(u)} , \qquad
  (\check v *_\ell \check u) = \text{?} \quad(\ell\ge2)
\end{multline*}
It just occurs to me, by the way, that the previous description of
basic (or ``primary'') data for an \oo-groupoid is already incomplete
in some rather obvious respect, namely that the symmetry-operation
$\inv_i : u \mapsto \check u$ on $F_i$ must be complemented by $i-1$
similar involutions on $F_i$, which corresponds algebraically to the
intuition that when we have an $(i+1)$-arrow $\lambda$ say between two
$i$-arrows $u$ and $v$, then we must be able to deduce from it another
arrow from $\check u$ to $\check v$ (namely $u\mapsto\check u$ has a
``functorial character'' for variable $u$)? This seems a rather
anodine modification of the previous set-up, and is irrelevant for the
main point I want to make here, namely: that for the notion of
\oo-groupoids we are after, all the equalities just envisioned in this
paragraph (and those I guess which will ensure naturality by the
necessary extension of the basic involution on $F_i$) should be
replaced by ``homotopies'', namely by $(i+1)$-arrows between the two
members. These arrows should be viewed, I believe, as being part of
the data, they appear here as a kind of ``secondary'' structure. The
difficulty which appears now is to work out the natural coherence
properties\namedpspage{L.3}{3} concerning this secondary structure. The first thing I
could think of is the ``pentagon axiom'' for the associativity data,
which occurs when looking at associativities for the compositum (for
$\overset{i}{\underset{\ell}{*}}$ say) of four factors. Here again the
first reflex would be to write down, as usual, an \emph{equality} for
two compositions of associativity isomorphisms, exhibited in the
pentagon diagram. One suspects however that such equality should,
again, be replaced by a ``homotopy''-arrow, which now appears as a
kind of ``ternary'' structure -- before even having exhausted the list
of coherence ``relations'' one could think of with the respect to the
secondary structure! Here one seems caught at first sight in an
infinite chain of ever ``higher'', and presumably, messier structures,
where one is going to get hopelessly lost, unless one discovers some
simple guiding principle for shedding some clarity into the mess.

\hangsection{``Fundamental \texorpdfstring{\oo}{oo}-groupoids'' as
  objects of a ``model category''?}\label{sec:3}%
I thought of writing you mainly because I believe that,
if anybody, you should know if the kind of structure I am looking for
has been worked out -- maybe even \emph{you} did? In this respect, I
vaguely remember that you had a description of $n$-categories in terms
of $n$-semisimplicial sets, satisfying certain exactness conditions,
in much the same way as an ordinary category can be interpreted, via
its ``nerve'', as a particular type of semisimplicial set. But I have
no idea if your definition applied only for describing $n$-categories
with strict associativities, or not\footnote{Definitely only for \emph{strict} associativity.}.

Still some contents in the spirit of your axiomatics of homotopical
algebra -- in order to make the question I am proposing more seducing
maybe to you! One comment is that presumably, the category of
\oo-groupoids (which is still to be defined) is a ``model category''
for the usual homotopy category; this would be at any rate one
plausible way to make explicit the intuition referred to before, that
a homotopy type is ``essentially the same'' as an \oo-groupoid up to
\oo-equivalence. The other comment: the construction of the
fundamental \oo-groupoid of a space, disregarding for the time being
the question of working out in full the pertinent structure on this
messy object, can be paraphrased in any model category in your sense,
and yields a functor from this category to the category of
\oo-groupoids, and hence (by geometric realization, or by
localization) also to the usual homotopy category\footnote{This idea is taken up again in section 12. The statement made here is a little rash, as for existence and uniqueness (in a suitable sense) of this functor. Compare note $(^{17})$ below.}. Was this functor
obvious beforehand? It is of a non-trivial nature only when the model
category is \emph{not} pointed -- as a matter of fact the whole
construction can be carried out canonically, in terms of a ``cylinder
object'' $I$ for the final object $e$ of the model category, playing
the role of the unit segment.\namedpspage{L.3prime}{3'}
It's high time to stop this letter -- please excuse me if it should
come ten or fifteen years too late, or maybe one year too early. If
you are not interested for the time being in such general nonsense,
maybe you know someone who is \ldots

\bigskip

Very cordially yours

\bigbreak

\presectionfill\ondate{20.2.}\namedpspage{L.4}{4}\par

I finally went on pondering about a definition of
\oo-groupoids, and it seems to me that, after all, the topological
motivation does furnish the ``simple guiding principle'' which
yesterday seemed to me to be still to be discovered, in order not to
get lost in the messiness of ever higher order structures. Let me try
to put it down roughly.

\hangsection{A bit of ordering in the mess of ``higher order
  structures''.}\label{sec:4}%
First I would like to correct somewhat the rather indiscriminate
description I gave yesterday of what I thought of viewing as
``primary'', secondary, ternary etc.\ structures for an
\oo-groupoid. More careful reflection conduces to view as the most
primitive, starting structure on the set of sets $F_i$ ($i\in\bN$), as
a skeleton on which progressively organs and flesh will be added, the
mere data of the source and target maps
\begin{equation*}
  \cst s_i, \cst t_i : F_i \rightrightarrows F_{i-1}
  \quad(i \ge 1),
\end{equation*}
which it will be convenient to supplement formally by corresponding
maps $\cst s_0, \cst t_0$ for $i=1$, from $F_0$ to
$F_{-1} \eqdef \text{one-element set}$. In a moment we will pass to a
universal situation, when the $F_i$ are replaced by the corresponding
``universal'' objects $\boldsymbol F_i$ in a suitable category stable
under finite products, where $F_{-1}$ will be the final element. For
several reasons, it is not proper to view the inversion maps
$\inv_i : F_i \to F_i$, and still less the other $i-1$ involutions on
$F_i$ which I at first overlooked, as being part of the primitive or
``skeletal'' structure. One main reason is that already for the most
usual $2$-groupoids, such as the $2$-groupoid whose $0$-objects are
ordinary ($1$-)groupoids, the $1$-objects being \emph{equivalences} between
these (namely functors which are fully faithful and essentially
surjective), and the $2$-objects morphisms (or ``natural
transformations'') between such, there is \emph{not}, for an
$i$-object $f: C \to C'$, a natural choice of an ``inverse'' namely of
a quasi-inverse in the usual sense. And even assuming that such
quasi-inverse is chosen for every $f$, it is by no means clear that
such choice can be made involutive, namely such that
$(f\upvee)\upvee = f$ for every $f$ (and not merely $(f\upvee)\upvee$
isomorphic to $f$). The maps $\inv_i$ will appear rather, quite
naturally, as ``primary structure'', and they will not be involutions,
but ``pseudo-involutions'' (namely involutions ``up to homotopy''). It
turns out that among the various functors that we will construct, from
the category of topological spaces to the category of \oo-groupoids
(the construction depending on arbitrary choices and yielding a large
bunch of mutually non-isomorphic functors, which however are
``equivalent'' in a sense which will have to be made precise) -- there
are choices neater than others, and some of these will yield in the
primary structure maps $\inv_i$ which are actual involutions and
similarly for the other pseudo-involutions, appearing in succession as
higher order structure. The possibility of such neat and fairly
natural choices had somewhat misled me\namedpspage{L.4prime}{4'} yesterday.

What may look less convincing though at first sight, is my choice to
view as non-primitive even the ``degeneration maps'' $\cst k_1^i : F_i
\to F_{i+1}$, associating to every $i$-objects the $i+1$-object acting
as an identity on the former. In all cases I have met so far, these
maps are either given beforehand with the structure (of a $1$-category
or $2$-category say), or they can be uniquely deduced from the
axioms. In the present set-up however, they seem to me to appear more
naturally as ``primary'' (not as primitive) structure, much in the
same way as the $\inv_i$. Different choices for associating an
\oo-groupoid to a topological space, while yielding the same base-sets
$F_i$, will however (according to this point of view) give rise to
different maps $\cst k_1^i$. The main motivation for this point of
view comes from the fact that the mechanism for a uniform construction
of the chain of ever higher order structures makes a basic use of the
source and target maps only and of the ``transposes'' (see below), and
(it seems to me) not at all of the degeneration maps, which in this
respect rather are confusing the real picture, if viewed as
``primitive''.
The degeneration maps rather appear as typical cases of primary
structure, probably of special significance in the practical handling
of \oo-groupoids, but not at all in the conceptual machinery leading
up to the construction of the structure species of ``\oo-groupoids''.

Much in the same way, the composition operations
$\overset{i}{\underset{1}{*}}$ are viewed as primary, not as
primitive or skeletal structure. Their description for the fundamental
\oo-groupoid of a space -- for instance the description of composition
of paths -- depends on arbitrary choices, such as the choice of an
isomorphism (say) between $(I,1) \amalg_e (I,0)$ and $I$, where $I$ is
the unit interval, much in the same way as the notion of an inverse of
a path depends on the choice of an isomorphism of $I$ with itself,
exchanging the two end-points $0$ and $1$. The operations
$\overset{i}{\underset{2}{*}}$ of Godement take sense only once the
composition operations $*_1$ are defined -- they are ``secondary
structure'', and successively the operations $*_3, \ldots, *_i$ appear
as ternary etc.\ structure. This is correctly suggested by the
notations which I chose yesterday, where however I hastily threw all
the operations into a same pot baptized ``primary structure''!

\hangsection{Jumping over the abyss!}\label{sec:5}%
It is about time though to come to a tentative precise definition of
description of the process of stepwise introduction of an increasing
chain of higher order structure. This will be done by introducing a
canonical sequence of categories and functors
\[ C_0 \to C_1 \to C_2 \to \dots \to C_n \to C_{n+1} \to \cdots,\]
where $C_n$ denotes the category harbouring the ``universal''\namedpspage{L.5}{5}
partial structure of a would-be \oo-groupoid, endowed only with its
``structure of order $\le n$''. The idea is to give a direct inductive
construction of this sequence, by describing $C_0$, and the passage
from $C_n$ to $C_{n+1}$ $(n\ge0$), namely from an $n$-ary to
$(n+1)$-ary structure. As for the meaning of ``universal structure'',
once a given structure species is at hand, it depends on the type of
categories (described by the exactness properties one is assuming for
these) one wants to take as carriers for the considered structure, and
the type of exactness properties one assumes for the functor one
allows between these. The choice depends partly on the particular
species; if it is an algebraic structure which can be described say by
a handful of composition laws between a bunch of base sets (or
base-\emph{objects}, when looking at ``realizations'' of the structure not
only in the category \Sets), one natural choice is to take categories
with finite products, and functors which commute to these. For more
sophisticated algebraic structures (including the structure of
category, groupoid or the like), which requires for the description of
data or axioms not only finite products, but also some fiber products,
one other familiar choice is to take categories with finite inverse
limits, and left exact functors. Still more sophisticated structures,
when the description of the structure in terms of base objects
requires not only some kind of inverse limits, but also more or less
arbitrary direct limits (such as the structure of a comodule over an
algebra, which requires consideration of tensor products over a ring
object\footnote{Another important example is the structure of a ``torsor'' under a group $G$ (torsor = principal homogeneous space). When this group $G$ is fixed, the corresponding classifying topos $\B_G$ is the natural purely algebraic substitute for the familiar ``classifying space'' for the discrete group $G$}.), still more
stringent conditions will have to be imposed upon categories and
functors between these, for the structure to make a sense in these
categories, and the functors to transform a structure of this type in
one category into one of same type in another. In most examples I have
looked up, everything is OK taking categories which are topoi, and
functors between these which are inverse image functors for morphisms
of topoi, namely which are left exact and commute with arbitrary
direct limits. There is a general theorem for the existence of
universal structures\footnote{Such a theory was developed in a seminar I gave at Buffalo in 1973}, covering all these cases -- for instance there
is a ``classifying topos'' for most algebro-geometric structures,
whose cohomology say should be viewed as the ``classifying
cohomology'' of the structure species considered. In the case we are
interested in here, it is convenient however to work with the smallest
categories $C_n$ feasible -- which amounts to being as generous as
possible for the categories one is allowing as carriers for the
structure of an \oo-groupoid, and for the functors between these which
are expected to carry an \oo-groupoid into an \oo-groupoid. What we
will do is define ultimately a structure of an \oo-groupoid in a
category $C$, as a sequence of objects $F_i$ ($i\in\bN$), endowed with
some structure to be defined,\namedpspage{L.5prime}{5'} assuming merely that in $C$
finite products of the $F_i$ exist, plus certain finite inverse limits
built up with the $F_i$'s and the maps $\cst s_\ell^i$, $\cst
t_\ell^i$ between them (the iterated source and target maps). It
should be noted that the type of $\varprojlim$ we allow, which will
have to be made precise below, is fixed beforehand in terms of the
``skeletal'' or ``primitive'' structure alone, embodied by the family
of couples $(\cst s_i, \cst t_i)_{i\in\bN}$. This implies that the
categories $\bC_i$ can be viewed as having \emph{the same set of
  objects}, namely the objects $\boldsymbol F_i$ (written
in bold\scrcomment{was: \emph{underlined}} now to indicate their
universal nature, and including as was said before $\boldsymbol F_{-1} =
\text{final object}$), plus the finite products and iterated fiber
products of so-called ``standard'' type. While I am writing, it
appears to me even that the finite products here are of no use (so we
just drop them both in the condition on categories which are accepted
for harbouring \oo-groupoids, and in the set of objects of the
categories $C_i$). Finally, the common set of objects of the
categories $\bC_i$ is the set of ``standard'' iterated fiber
products of the $\boldsymbol F_i$, built up using only the primitive
structure embodied by the maps $\cst s_i$ and $\cst
t_i$ (which I renounce henceforth to underline!).
This at the same time gives, in principle, a precise definition of
$\bC_0$, at least up to equivalence -- it should not be hard anyhow
to give a wholly explicit description of $\bC_0$ as a small
category, having a countable set of objects, once the basic notion of
the standard iterated fiber-products has been explained.

Once $\bC_0$ is constructed, we will get the higher order
categories $\bC_1$ (primary), $\bC_2$ etc.\ by an inductive
process of \emph{successively adding arrows}. The category $\bC_\oo$
will then be defined as the direct limit of the categories
$\bC_n$, having the same objects therefore as $\bC_0$, with
\[\Hom_\oo(X, Y) = \varinjlim_n \Hom_n(X, Y)\]
for any two objects. This being done, giving a structure of
\oo-groupoid in any category $C$, will amount to giving a functor
\[\bC_\oo \to C\]
commuting with the standard iterated fiber-products. This can be
reexpressed, as amounting to the same as to give a sequence of objects
($F_i$) in $C$, and maps $s_\ell^i, t_\ell^i$ between these,
satisfying the two relations I wrote down yesterday (page \hyperref[p:L.2]2) (and
which of course have to be taken into account when defining $C_0$ to
start with, I forgot to say before), and such that
``standard'' fiber-products (I'll drop the qualification ``iterated''
  henceforth!) defined in terms of these data should
exist in $C$, plus a bunch of maps between these fiber-products (in
fact, it will suffice to give such maps with target among the
$F_i$'s), satisfying certain relations embodied in the structure of
the category $\bC_\oo$. I am convinced that this bunch of maps
(namely the maps stemming from arrows in $\bC_\oo$) not only is
infinite, but cannot either be generated in the\namedpspage{L.6}{6} obvious
sense by a finite number, nor even by a finite number of infinite
series of maps such as $\cst k_\ell^i$, $\inv_i$,
$\overset{i}{\underset{\ell}{*}}$, the compatibility arrows in the
pentagon, and the like. More precisely still, I am convinced that none
of the functors $\bC_n \to \bC_{n+1}$ is an equivalence, which
amounts to saying that the structures of increasing order form a
strictly increasing sequence -- at every step, there is actual extra
structure added. This is perhaps evident beforehand to topologists in
the know, but I confess that for the time being it isn't to me, in
terms uniquely of the somewhat formal description I will make of the
passage of $\bC_n$ to $\bC_{n+1}$\footnote{That $\overset{i}{\underset{2}{*}}$ is of order $\ell$ is heuristically clear, but will require a proof none the less}. This step theoretically is all
that remains to be done, in order to achieve an explicit construction
of the structure species of an \oo-category (besides the definition of
standard fiber-products) -- without having to get involved, still less
lost, in the technical intricacies of ever messier diagrams to write
down, with increasing order of the structures to be added\ldots

\hangsection[The topological model: hemispheres building up the \dots]%
{The topological model: hemispheres building up the
  \texorpdfstring{\textup(}{(}tentative\texorpdfstring{\textup)}{)}
  ``universal \texorpdfstring{\oo-\textup(co\textup)}{oo-(co)}groupoid''.}%
\label{sec:6}%
In the outline of a method of construction for the structure species,
there has not been any explicit mention so far of the topological
motivation behind the whole approach, which could wrongly give the
impression of being a purely algebraic one. However, topological
considerations alone are giving me the clue both for the description
of the so-called standard fiber products, and of the inductive step
allowing to wind up from $\bC_n$ to $\bC_{n+1}$? The heuristics
indeed of the present approach is simple enough, and suggested by the
starting task, to define pertinent structure on the system of sets
$F_i(X)$ of ``homotopies'' of arbitrary order, associated to an
arbitrary topological space. In effect, the functors
\[ X \mapsto F_i(X) : \Spaces \to \Sets
\]
are representable by spaces $D_i$, which are easily seen to be
$i$-cells. The source and target maps $s_\ell^i, t_\ell^i :
F_i(X) \rightrightarrows F_{i-1}(X)$ are transposed to maps, which I may
denote by the same letters,
\[s_\ell^i, t_\ell^i : D_{i-1} \rightrightarrows D_i.\]
Handling around a little, one easily convinces oneself that all the
main structural items on $F_*(X)$ which one is figuring out in
succession, such as the degeneracy maps $k_\ell^i$, the inversion maps
$\inv_i$, the composition $v\cdot u = v *_i u$ for $i$-objects, etc.,
are all transposed of similar maps which are defined between the cells
$D_i$, or which go from such cells to certain spaces, deduced from
these by gluing them together -- the most evident example in this
respect being the composition of paths, which is transposed of a map
from the unit segment $I$ into $(I,1)\amalg_e (I,0)$, having
preassigned values on the endpoints of $I$ (which correspond in fact
to the images of the two maps $s_1^1,t_1^1 : D_0
= \text{one point} \rightrightarrows D_1 =
I$).\namedpspage{L.6prime}{6'} In a more suggestive way, we could say from this
experiment that the family of discs $(D_i)_{i\in\bN}$, together with
the maps $s,t$ and a lot of extra structure which enters into the
picture step by step, is what we would like to call a
\emph{co-\oo-groupoid in the category} \cTop{} of topological spaces
(namely an \oo-groupoid in the dual category $\cTop\op$), and
that the structure of \oo-category on $F_*(X)$ we want to describe is
the transform of this co-structure into an \oo-groupoid, by the
contravariant functor from \cTop{} to \Sets{} defined by $X$. The
(iterated) amalgamated sums in \cTop{} which allow to glue together the
various $D_i$'s using the $s$ and $t$ maps between them, namely the
corresponding fibered products in $\cTop\op$, are indeed
transformed by the functor $h_X$ into fibered products of \Sets.
The suggestion is, moreover, that if we view our co-structure in \cTop{}
as a co-structure in the subcategory of Top, say $B_\oo$, whose
objects are the cells $D_i$ and the amalgamated sums built up with
these which step-wise enter into play, and whose arrows are all those
arrows which are introduced step-wise to define the co-structure, and
all compositions of these -- that this should be \emph{the universal
  structure} in the sense dual to the one we have been contemplating
before; or what amount to the same, that the corresponding
\oo-groupoid structure in the dual category $B_\oo\op$ is
``universal'' -- which means essentially that it is none other than
the universal structure in the category $\bC_\oo$ we are
after. Whether or not this expectation will turn out to be correct (I
believe it is\footnote{(Added \alsoondate{23.2.83} I don't believe it now any more - and I do not really care - compare comments in section 11}), we should be aware that, while the successive
introduction of maps between the cells $D_i$ and their ``standard''
amalgamated sums (which we will define precisely below) depends at
every stage on arbitrary choices, the categories $\bC_n$ and their
limit $\bC_\oo$ do not depend on any of these choices; assuming the
expectation is correct, this means that up to (unique) isomorphism,
the category $B_\oo$ (and each of the categories $B_n$ of which it
appears as the direct limit) is independent of those choices -- the
isomorphism between two such categories transforming any one choice
made for the first, into the corresponding choice made for the
second. Also, while this expectation was of course the crucial
motivation leading to the explicit description of $\bC_0$ and of
the inductive step from $\bC_n$ to $\bC_{n+1}$, this description
seems to me a reasonable one and in any case it makes a formal sense,
quite independently of whether the expectation proves valid or not.

\hangsection{Gluing hemispheres: the ``standard amalgamations''.}%
\label{sec:7}%
Before pursuing, it is time to give a more complete description of the
primitive structure on $(D_i)$, as embodied by the maps $s,t$, which I
will now denote by\footnote{\emph{NB}\enspace it is more natural to consider
  $\varphi^+$ as ``target'' and $\varphi^-$ as ``source''.}
\[\varphi_i^+,\varphi_i^- :  D_i \rightrightarrows D_{i+1}. \]
It appears that these maps are injective, that their images make
up\namedpspage{L.7}{7} the boundary $S_i = \dot D_{i+1}$ of $D_{i+1}$, more
specifically these images are just two ``complementary'' hemispheres
in $S_i$, which I will denote by $S_i^+$ and $S_i^-$. The kernel of
the pair $(\varphi_i^+, \varphi_i^-)$ is just $S_{i-1} = \dot D_i$,
and the common restriction of the maps $\varphi_i^+, \varphi_i^-$ to
$S_{i-1}$ is an isomorphism
\[ S_{i-1} \simeq  S_i^+ \cap S_i^-.\]
This $S_{i-1}$ in turn decomposes into the two hemispheres $S_{i-1}^+,
S_{i-1}^-$, images of $D_{i-1}$. Replacing $D_{i+1}$ by $D_i$, we see
that the $i$-cell $D_i$ is decomposed into a union of $2i+1$ closed
cells, one being $D_i$ itself, the others being canonically isomorphic
to the cells $S_j^+,S_j^-$ ($0\le j\le i$), images of $D_j\to D_n$ by
the iterated morphisms
\[\varphi_{n,j}^+,\varphi_{n,j}^- :  D_j \rightrightarrows D_n. \]
This is a cellular decomposition, corresponding to a partition of
$D_n$ into $2n+1$ open cells $D_n$, $S_j^+ = \varphi_{n,j}^+(D_j)$,
$S_j^- = \varphi_{n,j}^-(D_j)$. For any cell in this decomposition,
the incident cells are exactly those of strictly smaller dimension.

When introducing the operation $\overset{n}{\underset{\ell}{*}}$ with
$\ell=n-j$, it is seen that this corresponds to choosing a map
\[ D_n \to (D_n, S_j^+) \amalg_{D_j} (D_n, S_j^-),\]
satisfying a certain condition *, expressing the formulas I wrote down
yesterday for $s_1$ and $t_1$ of $u *_\ell v$ -- the formulas
translate into demanding that the restriction of the looked-for map of
$D_n$ to its boundary $S_{n-1}$ should be a given map (given at any
rate, for $\ell\ge2$, in terms of the operation $*_{\ell-1}$, which
explains the point I made that the $*_\ell$-structure is of order just
above the $*_{\ell-1}$-structure, namely (inductively) is of order
$\ell$ \ldots). That the extension of this map of $S_{n-1}$ to $D_n$
does indeed exist, comes from the fact that the amalgamated sum on the
right hand side is contractible for obvious reasons.

This gives a clue of what we should call ``standard'' amalgamated sums
of the cells $D_i$. The first idea that comes to mind is that we
should insist that the space considered should be contractible,
excluding amalgamated sums therefore such as
\[
\begin{tikzcd}[row sep=tiny]
  & \bullet \ar[dr] \\ \bullet \ar[ur] \ar[rr] & & \bullet
\end{tikzcd} \quad\text{or}\quad
\begin{tikzcd}
  \bullet \ar[r, bend left] \ar[r, bend right] & \bullet
\end{tikzcd}\]
which are circles. This formulation however has the inconvenience of
not being directly expressed in combinatorial terms. The following
formulation, which has the advantage of being of combinatorial nature,
is presumably equivalent to the former, and gives at any rate (I expect) a large
enough notion of ``standardness'' to yield for the corresponding
notion of \oo-category enough structure for whatever one will ever
need. In any case, it is\namedpspage{L.7prime}{7'} understood that the ``amalgamated
sum'' (rather, finite $\varinjlim$) we are considering are of the most
common type, when $X$ is the finite union of closed subsets $X_i$,
with given isomorphisms
\[ X \simeq  D_{n(i)},\]
the intersection of any two of these $X_i\cap X_j$ being a union of
closed cells both in $D_{n(i)}$ and in $D_{n(j)}$. (This implies in
fact that it is either a closed cell in both, or the union of two
closed cells of same dimension $m$ and hence isomorphic to $S_m$, a
case which will be ruled out anyhow by the triviality condition which
follows.\footnote{This is nonsense, as one sees already in the following picture where $f \cap f'' = S_0$: [picture of disk divided in three]}) The
triviality or ``standardness'' condition is now expressed by demanding
that the set of indices $I$ can be totally ordered, i.e., numbered in
such a way that we get $X$ by successively ``attaching'' cells
$D_{n(i)}$ to the space already constructed, $X(i-1)$, by a map from a
sub-cell of $D_{n(i)}$, $S_j^\xi \to X(i-1)$ ($\xi\in\{\pm1\}$), this
map of course inducing an isomorphism, more precisely \emph{the}
standard isomorphism, $\varphi_{n(i)}^\xi : S_j^\xi \simeq  D_j$ with
$S_j^\xi$ one of the two corresponding cells $S_j^+, S_j^-$ in some
$X_{i'}\simeq  D_{n(i')}$. The dual translation of this, in terms of
fiber products in a category $C$ endowed with objects $F_i$
($i\in\bN$) and maps $s_1^i,t_1^i$ between these, is clear: for a
given set of indices $I$ and map $i\mapsto n(i) : I \to \bN$, we
consider a subobject of $\prod_{i\in I} D_{n(i)}$, which can be
described by equality relations between iterated sources and targets
of various components of $u=(u_i)_{i\in I}$ in $P$, the structure of
the set of relations being such that $I$ can be numbers, from $1$ to
$N$ say, in such a way that we get in succession $N-1$ relations on
the $N$ components $u_i$ respectively ($2\le i\le N$), every relation
being of the type $f(u_i) = g(u_i')$, with $f$ and $g$ being iterated
source of target maps, and $i'<i$. (Whether source or target depending
in obvious way on the two signs $\xi,\xi'$.)

\bigbreak

\presectionfill\ondate{21.2.}\namedpspage{L.8}{8}\par

\hangsection{Description of the universal \emph{primitive} structure.}%
\label{sec:8}%
Returning to the amalgamated sum $X = \bigcup_i X_i$, the cellular
decompositions of the components $X_i \simeq  D_{n(i)}$ define a
cellular decomposition of $X$, whose set of cells with incidence
relation forms a finite ordered set $K$, finite union of a family of
subsets $(K_i)_{i\in I}$, with given isomorphisms
\[ f_i : K_i \simeq  \cst J_{n(i)} \quad (i\in I),\]
where for every index $n\in\bN$, $\cst J_n$ denotes the ordered set of
the $2n+1$ cells $S_j^\xi$ ($0\le j\le n-1$, $\xi\in\{\pm1\}$), $D_n$
of the relevant cellular decomposition of $D_n$. We may without loss
of generality assume there is no inclusion relation between the $K_i$,
moreover the standardness condition described above readily translates
into a condition on this structure $K$, $(K_i)_{i\in I}$, $(f_i)_{i\in
  I}$, and implies that for $i,i'\in I$, $K_i\to K_{i'}$, is a
``closed'' subset in the two ordered sets $K_i,K_{i'}$ (namely
contains with any element $x$ the elements smaller than $x$), and
moreover isomorphic (for the induced order) to some $\cst J_n$. Thus
the category $B_0$ can be viewed as the category of such ``standard
ordered sets'' (with the extra structure on these just said), and the
category $C_0$ can be defined most simply as the dual category
$B_0\op$. (NB\enspace the definition of morphisms in $B_0$ is clear I guess
\ldots) I believe the category $B_0$ is stable under amalgamated sums
$X \amalg_Z Y$, provided however we insist that the empty structure
$K$ is \emph{not} allowed -- otherwise we have to restrict to
amalgamated sums with $Z \ne \emptyset$. It seems finally more
convenient to exclude the empty structure in $B_0$, i.e.\ to exclude
the final element from $\bC_0$, for the benefit of being able to
state that $\bC_0$ (and all categories $\bC_n$\footnote{This is nonsense again for $n \leq 1$, see P.S. at the end of section 13. Even this P.S. is still inaccurate. Compare comments section 18.} ) are stable under binary amalgamated sums,
and that the functors $\bC_n\to C$ (and ultimately $\bC_\oo\to
C$) we are interested in are just those commuting to arbitrary binary
amalgamated sums (without awkward reference to the objects $\boldsymbol F_i$
and the iterated source and target maps between them\ldots).

\hangsection[The main inductive step: just add coherence arrows! The \dots]%
{The main inductive step: just add coherence arrows! The abridged
  story of an \texorpdfstring{\textup(}{(}inescapable and
  irrelevant\texorpdfstring{\textup)}{)} ambiguity}\label{sec:9}%
The category $\bC_0$ being fairly well understood, it remains to
complete the construction by the inductive step, passing from $\bC_n$
to $\bC_{n+1}$. The main properties I have in mind therefore, for the
sequence of categories $\bC_n$ and their limit $\bC_\oo$, are the
following two.
\begin{enumerate}[label=(\Alph*)]
\item\label{it:8.A} For any $K\in\Ob(\bC_\oo)$ ($=\Ob(\bC_0)$), and any two
  arrows in $\bC_\oo$
  \[ f,g : K \rightrightarrows F_i,\]
  with $i\in\bN$, and such that either $i=0$, or the equalities
  \begin{equation}\label{eq:8.1}
    s_1^i\, f = s_1^i\, g , \quad
    t_1^i\, f = t_1^i\, g\tag{1}
  \end{equation}
  hold (case $i\ge1$), there exists $h : K \to F_{i+1}$ such that
  \begin{equation}\label{eq:8.2}
    s_1^{i+1}\, h = f, \quad
    t_1^{i+1}\, h = g.\tag{2}
  \end{equation}
\item\label{it:8.B}
  For any $n\in\bN$, the category $\bC_{n+1}$ is deduced from $\bC_n$
  by keeping the same objects, and just adding new arrows $h$ as
  in \ref{it:8.A}, with $f,g$ arrows in $\bC_n$.
\end{enumerate}
The\namedpspage{L.8prime}{8'} expression ``deduces from'' in \ref{it:8.B} means that we are
adding arrows $h: K\to F_i$ (each with preassigned source and target
in $\bC_n$), with as ``new axioms'' on the bunch of these uniquely
the two relations \eqref{eq:8.2} of \ref{it:8.A}, the category $\bC_{n+1}$
being deduced from $\bC_n$ in an obvious way, as the solution of a
universal problem within the category of all categories where binary
amalgamated products exist, and ``maps'' between these being functors
which commute to those fibered products\footnote{Same mistake as the one noticed in the previous note. The fibered products exists in $\bC$ only, and \emph{these} should be preserved by the functors under consideration. Thus the ``universal problem'' has to b rephrased somewhat\dots}. In practical terms, the
arrows of $\bC_{n+1}$ are those deduced from the arrows in $\bC_n$ and
the ``new'' arrows $h$, by combining formal operations of
composing arrows by $v \circ u$, and taking (binary) amalgamated
products of arrows.\footnote{This has to be corrected -- amalgamated
  sums exist in $\bC$, only -- and \emph{those} should be
  respected.}

NB\enspace Of course the condition \eqref{eq:8.1} in \ref{it:8.A} is necessary for the
existence of an $h$ satisfying \eqref{eq:8.2}. That it is sufficient
too can be viewed as an extremely strong, ``universal'' version of
coherence conditions, concerning the various structures introduced on
an \oo-groupoid. Intuitively, it means that whenever we have two ways
of associating to a finite family $(u_i)_{i\in I}$ of objects of an
\oo-groupoid, $u_i\in F_{n(i)}$, subjected to a standard set of
relations on the $u_i$'s, an element of some $F_n$, in terms of the
\oo-groupoid structure only, then we have automatically a ``homotopy''
between these built in in the very structure of the \oo-groupoid,
provided it makes at all sense to ask for one (namely provided
condition \eqref{eq:8.1} holds if $n \ge 1$). I have the feeling
moreover that conditions \ref{it:8.A} and \ref{it:8.B} (plus the
relation $\bC_\oo = \varinjlim \bC_n$) is all what will be ever
needed, when using the definition of the structure species, -- plus of
course the description of $\bC_0$, and the implicit fact that the
categories $\bC_n$ are stable under binary fiber products and the
inclusion functors commute to these.\footnote{Inaccurate; see above}
Of course, the category which really interests us is $\bC_\oo$, the
description of the intermediate $\bC_n$'s is merely technical --
the main point is that there should exist an increasing sequence
$(\bC_n)$ of subcategories of $\bC_\oo$, having the same objects (and
the ``same'' fiber-products), such that $\bC_\oo$ should be the
limit (i.e., every arrow in $\bC_\oo$ should belong to some $\bC_n$),
and such that the passage from $\bC_n$ to $\bC_{n+1}$
should satisfy \ref{it:8.B}. It is fairly obvious that these
conditions alone do by no means characterize $\bC_\oo$ up to
equivalence, and still less the sequence of its subcategories $\bC_n$.
The point I wish to make though, before pursuing with a proposal
of an explicit description, is (firstly) that \emph{this ambiguity is in the nature of things}. Roughly saying, two different mathematicians,
working independently on the conceptual problem I had in mind, assuming
they both wind up with some explicit definition, will almost certainly
get non-equivalent definitions -- namely with non-equivalent
categories of (set-valued, say) \oo-groupoids! And, secondly and as
importantly,\namedpspage{L.9}{9} that \emph{this ambiguity however is an
  irrelevant one}. To make this point a little clearer, I could say
that a third mathematician, informed of the work of both, will readily
think out a functor or rather a pair of functors, associating to any
structure of Mr.\ X one of Mr.\ Y and conversely, in such a way that
by composition of the two, we will associate to an $X$-structure ($T$
say) another $T'$, which will not be isomorphic to $T$ of course, but
endowed with a canonical \oo-equivalence (in the sense of Mr.\ X) $T
\underset{\oo}{\simeq } T'$, and the same on the Mr.\ Y side. Most
probably, a fourth mathematician, faced with the same situation as the
third, will get his own pair of functors to reconcile Mr.\ X and Mr.\
Y, which very probably won't be equivalent (I mean isomorphic) to the
previous one. Here however, a fifth mathematician, informed about this
new perplexity, will probably show that the two $Y$-structures $U$ and
$U'$, associated by his two colleagues to an $X$-structure $T$, while
not isomorphic alas, admit however a canonical \oo-equivalence between
$U$, and $U'$ (in the sense of the $Y$-theory). I could go on with a
sixth mathematician, confronted with the same perplexity as the
previous one, who winds up with another \oo-equivalence between $U$
and $U'$ (without being informed of the work of the fifth), and a
seventh reconciling them by discovering an \oo-equivalence between
these equivalences. The story of course is infinite, I better stop
with seven mathematicians, a fair number nowadays to allow themselves
getting involved with foundational matters \ldots
There should be a mathematical statement though resuming in finite
terms this infinite story, but in order to write it down I guess a
minimum amount of conceptual work, in the context of a given notion of
\oo-groupoids satisfying the desiderata \ref{it:8.A} and \ref{it:8.B}
should be done, and I am by no means sure I will go through this, not
in this letter anyhow.

\hangsection{Cutting down redundancies -- or: ``l'embarras du
  choix''.}\label{sec:10}%
Now in the long last the explicit description I promised of $\bC_{n+1}$
in terms of $\bC_n$. As a matter of fact, I have a handful
to propose! One choice, about the widest I would think of, is: for
every pair $(f,g)$ in $\bC_n$ satisfying condition \eqref{eq:8.1}
of \ref{it:8.A}, add one new arrow $h$. To avoid set-theoretic
difficulties though, we better first modify the definition of $\bC_0$
so that the set of its objects should be in the universe we are
working in, preferably even it be countable. Or else, and more
reasonably, we will pick one $h$ for every isomorphism class of
situations $(f,g)$ in $\bC_n$. Another restriction to avoid too
much redundancy -- this was the first definition actually that flipped
to my mind the day before yesterday -- is to add a \emph{new} $h$ only
when there is no ``old'' one, namely in $\bC_n$, serving the same
purpose. Then it came to my mind that there is a lot of redundancy
still, thus there would be already\namedpspage{L.9prime}{9'} infinitely many
operations standing for the single operation $v \overset{i}{\circ} u$
say, which could be viewed in effect in terms of an arbitrary
$n$-sequence ($n\ge2$) of ``composable'' $i$-objects $u_1=u, u_2=v,
u_3, \dots, u_n$, just ``forgetting'' $u_3, \dots, u_n$! The natural way to meet this ``objection'' would be
to restrict to pairs $(f,g)$ which cannot be factored non-trivially
through another objects $K'$ as
\[\begin{tikzcd}
  K \ar[r] & K' \ar[r, shift left, "f'"] \ar[r, shift right, swap,
  "g'"] & F_i \end{tikzcd}.\]
But even with such restrictions, there remain a lot of redundancies --
and this again seems to me in the nature of things, namely that there
is no really natural, ``most economic'' way for achieving condition
\ref{it:8.A}, by a stepwise construction meeting condition
\ref{it:8.B}. For instance, in $\bC_1$ already we will have not
merely the compositions $v \overset{i}{\circ} u$, but at the same time
simultaneous compositions
\begin{equation}
  \label{eq:10.star}
  u_n \circ u_{n-1} \circ \dots \circ u_1 \tag{*}
\end{equation}
for ``composable'' sequences of $i$-objects ($i\ge1$), without
reducing this (as is customary) to iteration of the binary composition
$v \circ u$.
Of course using the binary composition, and more generally iteration
of $n'-ary$ compositions with $n'<n$ (when $n\ge3)$, we get an
impressive bunch of operations in the $n$ variables $u_1, \dots, u_n$,
serving the same purpose as \eqref{eq:10.star}. All these will be tied
up by homotopies in the next step $\bC_2$. We would like to think
of this set of homotopies in $\bC_2$ as a kind of ``transitive
system of isomorphisms'' (of associativity), now the transitivity
relations one is looking for will be replaced by homotopies again
between compositions of homotopies, which will enter in the picture
with $\bC_3$, etc. Here the infinite story is exemplified by the
more familiar situation of the two ways in which one could define a
``$\otimes$-composition with associativity'' in a category, starting
either in terms of a binary operation, or with a bunch of $n$-ary
operations -- with, in both cases the associativity isomorphisms being
an essential part of the structure. Here again, while it is generally
(and quite validly) felt that the two points of view are equivalent;
and both have their advantages and their drawbacks, still it is not
true, I believe, that the two categories of algebraic structures
``category with associative $\otimes$-operation'', using one or the
other definition, are equivalent\footnote{I guess they are not equivalent, even when restricting to objects which are groupoids (i.e. so-called $\Gr$-categories).}. Here the story though of the relation between the two
notions is a finite one, due to the fact that it is related to the
notion of $2$-categories or $2$-groupoids, instead of \oo-groupoids as
before\ldots

Thus I don't feel really like spending much energy in cutting down
redundancies, but prefer working with a notion of \oo-groupoid which
remains partly indeterminate, the main features being embodied in the
conditions \ref{it:8.A} and \ref{it:8.B} and in the description of
$\bC_0$, without other specification.\namedpspage{L.10}{10}

\hangsection[Returning to the topological model (the canonical
functor \dots]%
{Returning to the topological model \texorpdfstring{\textup(}{(}the
  canonical functor from spaces to
  ``\texorpdfstring{\oo}{oo}-groupoids''\texorpdfstring{\textup)}{)}.}%
\label{sec:11}%
One convenient way for constructing a category $\bC_\oo$ would be
to define for every $K,L\in\Ob(\bC_0) = \Ob(B_0)$ the set
$\Hom_\oo(K,L)$ as a subset of the set $\Hom(\abs L,\abs K)$ of
continuous maps between the geometric realizations of $L$ and $K$ in
terms of gluing together cells $D_i$, the composition of arrows in
$\bC_\oo$ being just composition of maps. This amounts to defining
$\bC_\oo$ as the dual of a category $\bB_\oo$ of topological
descriptions. It will be sufficient to define for every cell $D_n$ and
every subset $\Hom_\oo(D_n,\abs K)$ of $\Hom(D_n,\abs K)$, satisfying
the two conditions:
\begin{enumerate}[label=(\alph*)]
\item\label{it:11.a}
  stability by compositions $D_n \to \abs K \to \abs{K'}$, where $K\to
  K'$ is an ``allowable'' continuous map, namely subjected only to the
  condition that its restriction to any standard subcell $D_{n'}
  \subset \abs K$ is again ``allowable'', i.e., in $\Hom_\oo$.
\item\label{it:11.b}
  Any ``allowable'' map $S_n \to \abs K$ (i.e., whose restrictions to
  $S_n^+$ and $S_n^-$ are allowable) extends to an allowable map
  $D_{n+1}\to\abs K$.
\end{enumerate}

Condition \ref{it:11.a} merely ensures stability of allowable maps
under composition, and the fact that $\bB_\oo$ (endowed with the
allowable maps as morphisms) has the correct binary amalgamated sums,
whereas \ref{it:11.b} expresses condition \ref{it:11.a} on $\bC_\oo$.
These conditions are satisfied when we take as $\Hom_\oo$
subsets defined by tameness conditions (such as piecewise linear for
suitable piecewise linear structure on the $D_n$'s, or differentiable,
etc.). The condition \ref{it:11.b} however is of a subtler nature in
the topological interpretation and surely not met by such sweeping
tameness requirements only! Finally, the question as to whether we can
actually in this way describe an ``acceptable'' category $\bC_\oo$,
by defining sets $\Hom_\oo$, namely describing $\bC_\oo$ in terms
of $\bB_\oo$, seems rather subsidiary after all. We may think of
course of constructing stepwise $\bB_\oo$ via subcategories $\bB_n$,
by adding stepwise new arrows in order to meet condition
\ref{it:11.b}, thus paraphrasing condition \ref{it:8.B} for passage
from $\bC_n$ to $\bC_{n+1}$ -- but it is by no means clear that
when passing to the category $\bB_{n+1}$ by composing maps of $\bB_n$
and ``new'' ones, and using amalgamated sums too, there might not
be some undesirable extra relations in $\bB_{n+1}$, coming from the
topological interpretation of the arrows in $\bC_{n+1}$ as maps. To
say it differently, universal algebra furnishes us readily with an
acceptable sequence of categories $\bC_n$ and hence $\bC_\oo$,
and by the universal properties of the $\bC_n$ in terms of
$\bC_{n+1}$, we readily get (using arbitrary choices) a contravariant
functor $K \mapsto \abs K$ from $\bC_\oo$ to the category of
topological spaces (i.e., a co-\oo-groupoid in \cTop), but it is by no
means clear that this functor is \emph{faithful} -- and it doesn't really
matter after all!\namedpspage{L.10prime}{10'}

\hangsection{About replacing spaces by objects of a ``model
  category''.}\label{sec:12}%
I think I really better stop now, except for one last comment. The
construction of a co-\oo-groupoid $D_*$ in \cTop, giving rise to the
fundamental functor
\[ \cTop \longrightarrow (\text{\oo-groupoids}), \quad X \mapsto (\Hom(D_i, X))_i,
\]
generalizes, as I already alluded to earlier, to the case when \cTop{}
is replaced by an arbitrary ``model category'' $M$ in Quillen's sense. Here
however the choices occur not only stepwise for the primary,
secondary, ternary etc.\ structures, but already for the primitive
structures, namely by choice of objects $D_i$ ($i\in\bN$) in $M$, and
source and target maps $D_i \rightrightarrows D_{i+1}$. These choices
can be made inductively, by choosing first for $D_0$ the final object,
or more generally any object which is fibrant and trivial (over the
final objects), $D_{-1}$ being the initial object, and defining
further $S_0 = D_0 \amalg_{D_{-1}} D_0 = D_0 \amalg D_0$ with obvious
maps $\psi_0^+,\psi_0^- : D_0 \to S_0$, and then, if everything is
constructed up to $D_n$ and $S_n = (D_n,\varphi_{n-1}^+)
\amalg_{D_{n-1}} (D_n,\varphi_{n-1}^-)$, defining $D_{n+1}$ as any
fibrant and trivial\footnote{An object of $M$ is called fibrant (resp. trivial) if the map from it to the final object is fibrant (resp. a weak equivalence).} object together with a cofibrant map
\[ S_n \to D_{n+1}, \]
and $\varphi_n^+, \varphi_n^-$ as the compositions of the latter with
$\psi_n^+,\psi_n^-: D_n \rightrightarrows S_n$.
Using this and amalgamated sums in $M$, we get our functor
\[ \bB_0 = \bC_0\op \to M , \quad K\mapsto \abs K_M,\]
commuting with amalgamated sums, which we can extend stepwise through
the $\bC_n\op$'s to a functor $\bB_\oo = \bC_\oo\op \to M$,
provided we know that the objects $\abs K_M$ of $M$ ($K\in\Ob\bC_0$)
obtained by ``standard'' gluing of the $D_n$'s in $M$, are again
fibrant and trivial -- and I hope indeed that your axioms imply that,
via, say, that if $Z \to X$ and $Z\to Y$ are cofibrant and $X,Y,Z$ are
fibrant and trivial, then $X \amalg_Z Y$ is fibrant and trivial\footnote{This seems doubtful, unless all objects of $M$ are fibrant (which is true for the most familiar cases I was having in mind). But even assuming this, we still need a unicity statement (in a suitable sense) for the functors $\bB_0 \to $ thus obtained, in order to be sure that the corresponding functors $M \to \Hot$ are all canonical  isomorphic. These kind of questions may be viewed as closely related to the question of existence (and unicicity) of ``test functors'' for a given test category $A$ (here, the category of ``standard hemispheres'' $\Globe$) into a given asphericity or contractibility structure, as discussed in section 90 (without yet getting there any clear-cut handy existence and unicicity theorems as are to be hoped for). In my notes, the set-up of asphericity and contractibility structures, which has been worked out in the months following, has been gradually replacing Quillen's approach to homotopy models. (I'll have to come back in due course to the question of the relationship between the two approaches, which deserves to be understood.)}\ldots

Among the things to be checked is of course that when we localize the
category of \oo-groupoids with respect to morphisms which are ``weak
equivalences'' in a rather obvious sense (NB\enspace the definition of the
$\Pi_i$'s of an \oo-groupoid is practically trivial!), we get a
category equivalent to the usual homotopy category \Hot. Thus we get a
composed functor
\[ M \to (\text{\oo-groupoids}) \to \Hot,\]
as announced. I have some intuitive feeling of what this functor
stands for, at least when $M$ is say the category of semisimplicial
sheaves, or (more or less equivalently) of $n$-gerbes or \oo-gerbes\footnote{The word ``gerb'' here stands as a ``translation'' of the French word ``gerbe'', as used in Giraud's book on non-commutative cohomological algebra. With the terminology of next section, we should call it rather a ``stack of groupoids'' or ``$\Gr$-stack'' (with a specification by $n$ or $\infty$ if needed!)} on
a given topos: namely it should correspond to the operation of
``integration'' or ``sections'' for $n$-gerbes (more generally for
$n$-stacks) over a topos -- which is indeed \emph{the} basic operation
(embodying non-commutative cohomology objects of the topos) in
``non-commutative homological algebra''.

I guess that's about it for today. It's getting late and time to go to
bed! Good night.

\bigbreak

\presectionfill\ondate{22.2.}1983\namedpspage{L.11}{11}\par

\hangsection[An urgent reflection on proper names: ``Stacks''
  and \dots]%
  {An urgent reflection on proper names: ``Stacks'' and
  ``coherators''.}\label{sec:13}%
It seems I can't help pursuing further the reflection I started with
this letter! First I would like to come back upon terminology. Maybe
to give the name of $n$-groupoids and \oo-groupoids to the objects I
was after is not proper, for two reasons: a) it conflicts with a
standing terminology, applying to structure species which are
frequently met and deserve names of their own, even if they turn out
to be too restrictive kind of objects for the use I am having in mind
-- so why not keep the terminology already in use, especially for
two-groupoids, which is pretty well suited after all; b) the structure
species I have in mind is not really a very well determined one, it
depends as we saw on choices, without any one choice looking
convincingly better than the others -- so it would be a mess to give
an unqualified name to such structure, depending on the choice of a
certain category $\bC_\oo = \bC$. I have been thinking of the
terminology \emph{$n$-stack} and \emph{\oo-stack} (stack = ``champ'' in French), a
name introduced in Giraud's book (he was restricting to champs =
1-champs), which over a topos reduced to a one-point space reduces in
his case to the usual notion of a category, i.e., $1$-category. Here
of course we are thinking of ``stacks of groupoids'' rather than
arbitrary stacks, which I would like to call (for arbitrary order
$n\in\bN$ or $n=\oo$) $n$-Gr-stack -- suggesting evident ties with the
notion of Gr-categories, we should say Gr-$1$-categories, of Mme Hoang
Xuan Sinh. One advantage of the name ``stack'' is that the use it had
so far spontaneously suggests the extension of these notions to the
corresponding notions over an arbitrary topos, which of course is what
I am after ultimately. Of course, when an ambiguity is possible, we
should speak of $n$-$\bC$-stacks -- the reference to $\bC$
should make superfluous the ``Gr'' specification. Thus $n$-$\bC$-stacks
are essentially the same as $n$-$\bC$-stacks over the
final topos, i.e., over a one point space. When both $\bC$ and
``Gr'' are understood in a given context, we will use the terminology
$n$-stack simply, or even ``stack'' when $n$ is fixed throughout. Thus
it will occur that in certain contexts ``stack'' will just mean a
usual groupoid, in others it will mean just a category, but when $n=2$
it will not mean a usual $2$-groupoid, but something more general,
defined in terms of $\bC$.

The categories $\bC=\bC_\oo$ described before merit a name too
-- I would like to call them ``\emph{coherators}'' (``\emph{coh\'ereurs}'' in
French). This name is meant to suggest that $\bC$ embodies a
hierarchy of coherence relations, more accurately of coherence
``homotopies''. When dealing with stacks, the term
\emph{$i$-homotopies} (rather than $i$-objects or $i$-arrows) for the
elements of the $i$\textsuperscript{th} component $F_i$ of a stack
seems to me the most suggestive -- they will of course be denoted
graphically by arrows\namedpspage{L.11prime}{11'} (such as $h: f\to g$ in the
formulation of \ref{it:8.A} yesterday). More specifically, I will call
coherator any category \emph{equivalent} to a category $\bC_\oo$ as
constructed before. Thus a coherator is stable under binary fiber
products,\footnote{false, see PS} moreover the $F_i$ are recovered up
to isomorphism as the indecomposable elements of $\bC$ with respect
to amalgamation. However, in a category $\bC_\oo$, the objects
$\boldsymbol F_i$ have non-trivial automorphisms -- namely the ``duality
involutions'' and their compositions (the group of automorphisms of
$\boldsymbol F_i$ should turn out to be canonically isomorphic to
$(\pm1)^i$),\footnote{This is false, as seen below} in other words by the mere
category structure of a coherator we will not be able to recover the
objects $\boldsymbol F_i$ in $\bC$ up to unique isomorphism. Therefore, in
the structure of a coherator should be included, too, the choice of
the basic indecomposable objects $\boldsymbol F_i$ (one in each isomorphism
class), and moreover the arrows $\cst s_1^i,\cst t_1^i : \boldsymbol F_i
\rightrightarrows \boldsymbol F_{i-1}$ for $i\ge1$ (a priori, only the
\emph{pair} $(\cst s_1^i,\cst t_1^i)$ can be described intrinsically
in terms of the category structure of $\bC$, once $\boldsymbol F_i$ and
$\boldsymbol F_{i-1}$ are chosen\ldots). But it now occurs to me that we
don't have to put in this extra structure after all -- while the
$\boldsymbol F_i$'s separately do have automorphisms, the system of objects
$(\boldsymbol F_i)_{i\in\bN}$ and of the maps $(\cst s_1^i)_{i\ge1}$ and $(\cst
t_1^i)_{i\ge1}$ has only the trivial automorphism (all this of course
is heuristics, I didn't really prove anything -- but the structure of
the full subcategory of a $\bC_\oo$ formed by the objects $\boldsymbol F_i$ seems
pretty obvious\ldots). To finish getting convinced that the
mere category structure of a coherator includes already all other
relevant structure, we should still describe a suitable intrinsic
filtration by subcategories $\bC_n$. We define the $\bC_n$
inductively, $\bC_0$ being the ``primitive structure'' (the arrows
are those deducible from the source and target arrows by composition
and fiber products), and $\bC_{n+1}$ being defined in terms of
$\bC_n$ as follows: add to $\Fl \bC_n$ all arrows of $\bC$ of
the type $h: K \to \boldsymbol F_i$ ($i\ge1$) such that $\cst s_1^i\, h$ and
$\cst t_1^i\,h$ are in $\bC_n$, and the arrows deduced from the
bunch obtained by composition and fibered products.\footnote{This description is dubious, as it may give categories $\bC_n$ which are larger that the ones we want.} In view of
these constructions, it would be an easy exercise to give an intrinsic
characterization of a coherator, as a category satisfying certain
internal properties.

I was a little rash right now when making assertions about the
structure of the group of automorphisms of $\boldsymbol F_i$ -- I forgot that
two days ago I pointed out to myself that even the basic operation
$\inv_i$ upon $\boldsymbol F_i$ need not even be involutions (Nor
  even automorphisms)! However, I just checked that if in the inductive
construction of coherators $\bC_\oo$ given yesterday, we insist on
the most trivial irredundancy condition (namely that we don't add a
``new'' homotopy $h: f\to g$ when there is already an old one), then
any morphism $h: \boldsymbol F_i \to \boldsymbol F_i$ such that $\cst s\,f=\cst s$
and $\cst t\,f=\cst t$, is the identity -- and that\namedpspage{L.12}{12}
implies inductively that an automorphism of the system of $\boldsymbol F_i$'s
related by the source and target maps $\cst s_1^i,\cst t_1^i$ is the
identity. Thus it is correct after all, it seems, that the category
structure of a coherator implies all other structure relevant to
us\footnote{It seems dubious however that the mere category structure of $\bC$ will allow us to recover the ``primitive'' subcategory $\bC_0$, and it looks safer to add the latter as an extra structure to $\bC$}.

I do believe that the description given so far of what I mean by a
coherator, namely something acting like a kind of pattern in order to
define a corresponding notion of ``stacks'' (which in turn should be
the basic coefficient objects in non-commutative homological algebra,
as well as a convenient description of homotopy types) embodies some
of the essential features of the theory still in embryo that wants to
be developed. It is quite possible of course that some features are
lacking still, for instance that some extra conditions have to be
imposed upon $\bC$, possibly of a very different nature from mere
irredundancy conditions (which, I feel, are kind of irrelevant in this
set-up). Only by pushing ahead and working out at least in outline the
main aspects of the formalism of stacks, will it become clear whether
or not extra conditions on $\bC$ are needed. I would like at least
to make a commented list of these main aspects, and possibly do some
heuristic pondering on some of these in the stride, or afterwards. For
today it seems a little late though -- I have been pretty busy with
non-mathematical work most part of the day, and the next two days I'll
be busy at the university. Thus I guess I'll send off this letter
tomorrow, and send you later an elaboration (presumably much in the
style of this unending letter) if you are interested. In any case I
would appreciate any comments you make -- that's why I have been
writing you after all! I will probably send copies to Ronnie Brown,
Luc Illusie and Jean Giraud, in case they should be interested (I
guess at least Ronnie Brown is). Maybe the theory is going to take off
after all, in the long last!

\bigskip

Very cordially yours

\bigskip

PS (25.2.)\enspace I noticed a rather silly mistake in the notes of two days
ago, when stating that the categories $\bC_n$ admit fiber products:
what is true is that the category $\bC_0$ has fiber products (by
construction, practically), and that these are fiber products also in
the categories $\bC_n$ (by construction equally), i.e., that the
inclusion functors $\bC_0 \to \bC_n$ commute to fiber
products. Stacks in a category $C$ correspond to functors $\bC_\oo\to
C$ whose \emph{restriction to} $\bC_0$ commutes with fiber
products. I carried the mistake along in the yesterday notes -- it
doesn't really change anything substantially. I will have to come back
anyhow upon the basic notion of a coherator\ldots

\chapter*{Appendix: Three letters to Larry Breen}
\addcontentsline{toc}{chapter}{Appendix}
\label{ch:AppchI}

In this appendix, I am including three letters to Larry Breen\scrcomment{Lawrence Breen}, dated 5.2, 
17.2 and 17-19.7.1975. These letters are written in French, and the first 
two are reproduced textually, whereas the third appears here in English translation.
I am grateful to Ronnie Brown who took the trouble last year to make such a translation 
(from a hardly legible copy of the handwritten letter to Larry Breen) with the assistance of Larry Breen himself and J. L. Loday. I am now using his translation rather than the original letter, which no printer could possibly decipher!
Also, for the second letter I am using a typed copy made in 1945 by Larry Breen (who presumably had difficulties too deciphering the handwriting). Thanks are due to him for his interest and patience (with someone like me, very unknowledgeable in standard homotopy techniques), which appeared in the letters I got from him in response and his verbal explanations on related matters, as well as in the trouble he took in retyping the letter and sending me copies of my own letters to him, and allowing me to reproduce them in this volume of Pursuing stacks. My only present contribution to this set of letters is adding a few comments (in the notes), adding subtitles (with numbers 1 to 18), and correcting some inaccuracies in the English translation (due mainly to my handwriting\dots). 
Also, I skipped the beginning of the first letter, which doesn't seem of general relevance.

The first two letters are an attempt to explain to Larry Breen (who has a wide background in algebraic geometry and homological and homotopical algebra) some of the main points of the programme I had in mind around the notions of $n$-categories and $n$-stacks (which is what I am supposed to be pursuing now in my present work ``Pursuing stacks''). They were written under the impetus of the new intuition (new to me at any rate) which then had just appeared to me, namely that (non strict) $n$-groupoids should model (in a suitable sense) $n$-truncated homotopy types. The third letter, written in answer to a number of questions in Larry Breen's response to the first two, is of a wider scope. A large part of the letter outlines (very sketchily) some main points of a duality program (including a cohomological formulation of ``geometric'' local and global class field theory), which emerged by the end of the fifties and appears here for the first time in print. The later part of the letter gives also some hints about the need of a framework of ``tame topology'' suitable for writing up a ``dévissage theory'' of stratified spaces, and for working with étale tubular neighbourhoods, for the common purpose of coming to grips with a suitable notion of ``finite homotopy type'' of a ``tame'' topological space or of a scheme say, in terms of the inverted system of ``indexed homotopy types'' corresponding to all equi-singular stratifications.

\bigskip

\par\hfill Villecun \ondate{5.2.}1975\par

Cher Breen

\dots

\setcounter{section}{0}

\hangsection{Examples de $2$-catégories de Picard.}\label{sec:app1}%
\dots Pour tout le reste de ta lettre, elle mériterait un lecteur plus averti, aussi, pour qu'elle ne soit pas entièrement perdue au monde, je vais l'envoyer à Illusie ! J'ai néanmoins constaté, avec intérêt, ton intérêt â demi refoulé pour des $2$-catégories de Picard, $n$-catégories et autre faune de ce genre, et ton espoir que je te prouverai peut-être que ces animaux sont tout à fait indispensables pour faire des maths sérieuses dans telle circonstance. J'ai bien peur que cet espoir ne soit dé\c{c}u, je crois que jusqu'à maintenant on a toujours pu d'en tirer en éludant de tels objets et l'engrenage dans lequel ils pourraient nous entraîner. Est-ce nécessairement une raison pour continuer à les éluder ? Les situations où on a l'impression ``d'éluder'' en effet me semblent en tous cas devenir toujours plus nombreuses - et si on s'abstenait de tirer une situation complexe et chargée de mystère au clair, chaque fois qu'on ne serait pas \emph{forcé} de le faire pour des raisons techniques provenant de la math déjà faite, - il y aurait sans doute beaucoup de parties des maths aujourd'hui reputées ``sérieusses'' qui n'auraient jamais été developpées (Il n'est pas dit non plus que le mode s'en trouverait plus mal...). Ton commentaire (que j'ai également entendu chez Deligne) que la classification d'objets géométriques relativement merdiques se réduit finalement à des invariants cohomologiques essentiellement ``bien connus'' et relativement simples n'est pas non plus convainquant; n'est ce pas négliger la différence entre la \emph{compréhension d'un objet géométrique}, et la détermination de sa ``classe à isomorphisme (ou équivalence) près'' ?

Tu me demandes des exemples ``convainquants'' de $2$-catégories de Picard. Voici quelques exemples, en vrac (je ne sais s'ils seraient convainquants !):

\begin{enumerate}
 \item[1)]\label{it:app1.1} Si $L$ est un lien\footnote{For the notion of a ``lien'' (or ``tie''), which is one of the main ingredients of the non-commutative cohomology panoply of Giraud's theory, I refer to his books (Springer, Grundlehren 179, 1971). A \emph{Picard category} is a groupoid endowed with an operation $\otimes$ together with associativity, unity and commutativity data for this operation, which make it resemble to a commutative group. A \emph{``Champ de Picard''} (or ``Picard stack'') is defined accordingly, by relativizing over an arbitrary space or topos (replacing the groupoid by a stack of groupoids over this topos). The necessary ``general nonsense'' on these notions is developed rather carefully in an exposé of Deligne in SGA 4 (SGA 4 XVIII 1.4). In this letter to Larry Breen, I am asuming ``known'' the notion of an $n$-stack (for $n = 3$ at any rate), and the corresponding notion of (strict) \emph{Picard $n$-stack}, which should be describable (as was explained in Deligne's notes in the case $n = 1$) by an $n$-truncated chain complex in the category of abelian sheaves on $X$ (viewed mainly as an object of the relevant derived category). The ``strictness'' condition on usual Picard stacks refers to the restriction that the commutativity isomorphism within an object $L \otimes L'$, when $L = L'$, should reduce to the identity. It is assumed (without further explanation) that the condition carries over in a natural way to Picard $n$-stacks, in such a a way as to allow an interpretation of these by truncated objects in a suitable derived category, as hinted above.} de centre $Z$ sur le topos $X$, les \emph{gerbes liées par $Z$} forment une $2$-catégorie de Picard stricte, représentée par le complexe $\RGamma_X(Z)$ tronqué en degré $2$, dont les objets de cohomologie non triviaux sont les $\mathrm H^i(X, Z)$, $0 \leq i \leq 2$. Les gerbes liées par $L$ forment un \emph{pseudo-$2$-torseur} sous le gerbe précédente, qui est un $2$-torseur (i.e. non vide) si et seule si une certaine obstruction dans $\mathrm H^3(X, Z)$ est nulle. Pour comprendre cette classe de notre point de vue, il y a lieu de passer aux $2$-champs correspondants: le $2$-champs de Picard strict des $Z$-gerbes sur des objets variables de $X$, et le $2$-champ des $L$-gerbes sur des objets variables. Ce dernier est bel et bien un $2$-torseur sous le champ précédent, or la classification de ces $2$-torseurs (à $2$-équivalence près) se fait par le $\mathrm H^3(X, Z)$, (tout comme les $Z$-$L$-gerbes peuvent être interprétées comme des torseurs sous la $Z$-$L$-champ de Picard strict des $Z$-torseurs, et sont classifiées  par le $\mathrm H^2(X, Z)$). On voit déjà, bien sûr, poindre ici l'oreille de la $3$-catégorie de Picard stricte des $2$-gerbes liées par $Z$, ou (de fa\c{c}on équivalente) des $2$-torseurs sous le $2$-champ de Picard strict des $Z$-$L$-gerbes; cette $3$-catégorie de Picard stricte étant décrite par $\RGamma_X(Z)$ tronqué en dimension $3$, ayant comme invariants de cohomologie non triviaux les $\mathrm H^i(X, Z)$ ($0 \leq i \leq 3$). Quant au $3$-champ de Picard correspondant, il est décrit par une résolution injective de $Z$ tronqué en degré $3$, alors que le $2$-champ de Picard précédent se décrivait en tronquant en degré $2$.
    
\item[2)]\label{it:app1.2} Si $M$ et $N$ sont deux faisceaux abéliens sur $X$, les \emph{champs de Picard} (N.B. $1$-champs !) \emph{d'invariants $M$ et $N$} forment eux-même une $2$-catégorie de Picard stricte, représentée sans doute par le complexe $\RHom(X(M), N)$\footnote{When $M$ is any abelian sheaf on a topos, the ``MacLane resolution'' $X(M)$ is a certain canonical left resolution of $M$ by sheaves of $\mathbf{Z}$-modules which are ``free'', and more specifically, which are finite direct sums of sheaves of the type $\mathbf{Z}^{(T)}$, where $T$ is any sheaf of the type $M^n$ (finite product of copies of $M$). This canonical construction was introduced by MacLane (for abelian groups), and gained new popularity in the French school of algebraic geometry and homological algebra in the late sixties, because it gives a very handy way to relate the $\Ext^i(M, N)$ invariants (when $N$ is another abelian sheaf on $X$) to the ``spacial'' cohomology of $M$ (i.e. of the induced topos $X_{/M}$) with coefficients in $N$.} tronqué en degré $2$, dont les invariants de cohomologie non triviaux sont donc ``le drôle de $\Ext^2$'' de ma lettre à Deligne, et les honnêtes $\Ext^i(M, N)$ ($0 \leq i \leq 2$). Les champs de Picard stricts forment une sous-$2$-catégorie de Picard pleine, représentée par $\RHom(M, N)$ tronqué en degré $2$, d'invariants les $\Ext^i(M; N)$ ($0 \leq i \leq 2$). Bien sûr, $\Ext^2$ donne les $0$-objets à équivalence près, $\Ext^1$ les automorphismes à isomorphisme près de l'objet nul, $\Ext^0$ les automorphismes de l'automorphisme identique audit\dots Je n'ai pas réfléchi à une bonne interprétation géométrique de la $n$-catégorie de Picard associée à $\RHom(M, N)$ tronqué en degré $n$, et encore moins bien sûr pour $\RHom(X(M), N)$, mais sans doute il faut regarder dans la direction des $n$-champs de Picard. 
    
\item[3)]\label{it:app1.3} Soit $G$ un Groupe sur $X$, opérant sur un faisceau abélien $N$. Les \emph{champs en $\Gr$-catégories sur $X$ liés par $(G, N)$} forment une $2$-catégorie de Picard, dont les invariants sont $\mathrm H^3(\B_G mod X, N)$, $\mathrm H^2(\B_G mod X, N)$ et $Z^1(G, N)$ (groupe des $1$-cocycles de $G$ à coefficients dans $N$) - je te laisse le soin de deviner quel est le complexe qui le décrit ! J'ai écrit il y a quelques mois à Deligne\scrcomment{cf.\ \textcite{LGD7874}} à ce sujet, et l'ai prié de t'envoyer une copie de la lettre.
    
\item[4)]\label{it:app1.4} Soit $X$ un topos localement annelé, on peut considérer les \emph{Algèbres d'Azumaya sur $X$} (i.e. les Algèbres localement isomorphes à une algèbre de matrices d'ordre $n$, $n \geq 1$) comme les objets d'une $2$-catégorie de Picard, où la catégorie $\bHom(A, B)$, pour $A$ et $B$ des Algèbres d'Azumaya, est la catégorie des ``trivialisations'' de $A^{\circ} \otimes B$, i.e. des couples $(E, \emptyset)$, $E$ un Module localement libre et $\emptyset$ un isomorphisme $\bEnd(E) \simeq A^{\circ}\otimes B$. Il faut travailler un peu pour définir les accouplements $\bHom(A, B) \times \bHom(B, C) \to \bHom(A, C)$; l'opération $\otimes$ dans la $2$-catégorie de Picard à construire est bien sûr le produit tensoriel d'Algèbres, et l'opération ``puissance $-1$'' est le passage à l'algèbre opposée. On vérifie qu'en associant à toute Algèbre d'Azumaya la $1$-gerbe de ses trivialisations, on trouve un $2$-$\otimes$-foncteur de la $2$-catégorie de Picard (dite ``de Brauer'') dans celle des $1$-gerbes liées par $\mathbf{G}_m$, qui est $2$-fidèle. Les invariants de la première sont donc les groupes $\mathrm H^2(X, \mathbf{G}_m)_{\Br}$, $\mathrm H^1(X, \mathbf{G}_m)$ et $\mathrm H^0(X, \mathbf{G}_m)$, où dans le premier terme l'indice $\Br$ désigne le sous-groupe du $\mathrm H^2$ formé des classes de cohomologie provenant d'Algèbres d'Azumaya. On aurait envie de parler du $2$-champ de Picard des Algèbres d'Azumaya sur des objets variables de $X$, mais c'est bien une $2$-catégorie de Picard fibrée sur $X$, mais pas tout à fait un $2$-champ (whatever that means), san doute - la condition de $2$-recollement (whatever that means) ne doit pas être satisfaite - sinon il n'y aurait pas d'indice $\Br$ au $\mathrm H^2$...
\end{enumerate}

\hangsection{Théorème de Lefschetz (faible) en termes de Champs.}\label{sec:app2}%
La considération des $n$-catégories de Picard strictes (qui s'imposent à nous pas à pas dans un contexte essentiellement ``commutatif'') me semblent la clef du passage de l'algèbre homologique ordinaire (``commutative''), en termes de complexes, à une algèbre homologique non commutative, du fait qu'elles donnent une interprétation géométrique correcte des ``complexes tronqués à l'ordre $n$'' (en tant qui objets de catégories dérivées), donc, essentiellement (par passage à la limite sur $n$) des complexes tout courts. L'idée naïve qui se présente est alors que les ``complexes non commutatifs'' (qui seraient les objets-fantômes d'une algèbre homologique non commutative) sont peut-être ce qui reste des $n$-catégories de Picard (strictes) quand on oublie leur caractère additif, i.e. leur structure de Picard - c'est à dire qu'on ne retient que la $n$-catégorie ! (Quand on se place sur un topos $X$, on s'intéresse donc aux $n$-champs sur $X$...) A vrai dire, cette idée est venue d'abord d'une autre direction, quand il s'est agi en géométrie algébrique de démontrer des \emph{théorèmes de Lefschetz} à coefficients discrets en cohomologie étale, dans le cas d'une variété projective disons et de toute section hyperplane, ou d'une variété quasi-projective et de presque toute section hyperplane (pour ne mentionner que le cas global le plus simple), sous les hypothèses de profondeur cohomologique ``le plus naturelles'' (en fait, essentiellement des conditions nécessaires et suffisantes de validité du dit théorème). Dans le cas commutatif, les techniques de dualité nous suggèrent très clairement quels sont les meilleurs énoncés possibles, cf. l'exposé de Mme Raynaud dans SGA 2\scrcomment{\textcite{SGA2}}. Mais ces techniques ne valent qu'en se restreignant à des coefficients premiers aux caractéristiques, alors que des démonstrations directes plus géométriques (développées dans SGA 2 avant le développement du formalisme de la cohomologie étale) donnaient des résultats très voisins pour le $\mathrm H^0$ et le $\mathrm H^1$ (ou le $\pi_0$ et le $\pi_1$, si on préfère) sans telles restrictions, du moins dans le cas propre (i.e. projectif, au lien de quasi-projectif). En fait, ce sont les ``résultats les meilleurs possibles'' eux-mêmes, énoncés comme conjectures dans SGA 2 dans l'exposé cité de Mme Raynaud, qui sont démontres ultérieurement par elle dans sa thèse\scrcomment{\textcite{Raynaud1975}}. Ce qui est remarquable de notre point de vue, c'est que les énoncés les plus forts se présentent le plus naturellement sous forme d'énoncés sur des \emph{1-champs} sur le site étale de la variété algébrique considérée - la notion de ``profondeur $\geq i$ (pour $i = 1, 2, 3$) s'énon\c{c}ant aussi le plus naturellement en termes de champs. Non seulement cela, mais alors même qu'on voudrait ignorer la notion technique de champ et travailler exclusivement en termes de $\mathrm H^0$ et $\mathrm H^1$ en utilisant à bloc le formalisme cohomologique non commutatif de Giraud, pour démontrer disons un théorème de bijectivité $\pi_1(Y) \to \pi_1(X)$ (ce qui est le résultat le plus profond établi dans la thèse de Mme Raynaud), il semble bien qu'on n'y arrive pas, faute à ce formalisme d'avoir la souplesse nécessaire. En fait, il faut utiliser comme ingrédients techniques, de fa\c{c}on essentielle, les trois théorèmes suivantes directement pour les $1$-champs ``de torsion'' (i.e. où les faisceaux en groupes d'automorphismes sont de ind-torsion): a) théorème de changement de base pour une morphisme propre, b) théorème de changement de base par un morphisme lisse c) théorème de ``propreté cohomologique générique'' pour un morphisme de type fini $f: X \to S$, $S$ intègre (disant que l'on peut trouver dans $S$ un ouvert $U \neq \emptyset$ tel que pour \emph{tout} changement de base $S' \to S$ se factorisant par $u$, la formule de changement de base est vraie). (Pour b) et c), il faut faire des hypothèses que les faisceaux d'automorphismes sont premiers aux charactéristiques, et dans c) ne servent que dans la version ``générique'' du théorème de Lefschetz). C'est avec en vue de telles applications que Giraud a pris la peine dans son bouquin (si je ne me trompe) de démontrer a), b) (et c) ?) dans le contexte des 1-champs et de leurs images directes et inverses. Mais du même coup il dévient clair que le contexte ``naturel'' des théorèmes de changement de base en cohomologie étale, des théorèmes du type de Lefschetz (dits ``faibles'') sur les ``sections hyperplanes'', tout comme de la notion de profondeur qui y joue un rôle crucial, doit être celui des $n$-champs. Et que le développement hypothétique de ce contexte ne risque pas de se réduire à une jonglerie purement formelle et absolument bordélique avec du ``general nonsense'', mais qu'on se trouvera aussitôt confronté à des tests ``d'utilisabilité'' aussi sérieux que la démonstration des théorèmes de changement de base et ceux du type de Lefschetz (qui même dans le contexte commutatif ne sont pas piqués de vers...). [] pour variantes analytiques complexes etc.

Je ne sais si ces commentaires te ``passent par dessus la tête'' à ton tour, ni si elles te donnent l'impression qu'il aurait peut-être des choses intéressantes à tirer au clair. Si cela t'intéresse, je pourrais expliciter sous forme un peu plus systématique quelques ingrédients d'une hypothétique algèbre homologique non commutative et les liens de celle-ci à l'algèbre homologique commutative. Plus mystérieux pour moi (et pour cause, vu mon ignorance en homotopie) seraient les relations entre celle-là et l'algèbre homotopique, i.e. les structures semi-simpliciales, et je n'ai que des commentaires assez vagues à faire en ce sens (*). Par ailleurs, je te rappelle que même l'algèbre homologique commutative n'est pas, il s'en faut, dans un état satisfaisant, pour autant que je sache, vu qu'on ne sait\footnote{Reflecting on the ``right'' version of the provisional Verdier notion of a triangulated category (which was supposed to describe adequately the relevant internal structure of the derived categories of abelian categories) is part of my present program for the notes on Pursuing stacks, and will be the main task in one of the chapters of volume two. For some indications along these lines, see also section 69 (sketching the basic notion of a ``\emph{derivator}'').} toujours pas quelle est la ``bonne'' notion de catégorie triangulée. Or il me semble bien clair que ce n'est pas une question purement académique - même si on a pu se passer de le savoir jusqu'ici (en se bornant comme Monsieur Jourdain à ``faire de la prose sans le savoir'' - en travaillant sur des catégories de complexes, éventuellement filtrés, sans trop se demander quelles structures il y a sur ces catégories...).

Bien cordialement à toi

\hangsection{Les $n$-groupoïdes comme types d'homotopie tronqués.}\label{sec:app3}%
(*) P.S. Réflexion faite, j'ai quand même envie de te mettre un peu en appétit, en faisant ces ``quelques commentaires assez vagues''. Il s'agit du yoga qu'une (petite) $n$-catégorie ou groupoïdes (à $n$-équivalence près) ``est essentiellement la même chose'' qu'une ensemble semi-simplicial pris à homotopie près et où on néglige les $\pi_i$ pour $n + 1 \geq i$ (où, si tu préfères, ``où on a tué les groupes d'homotopie en dimension $\geq n + 1$). Voici des éléments heuristiques pour ce yoga. Si $K_\bullet$ est un ensemble simplicial (il peut être prudent de le prendre de Kan) on lui associe une $n$-catégorie $C_n(K_\bullet)$, dont les $0$-objets sont les $0$-simplexes, les $1$-objets sont les chemins (ou homotopies) entre $0$-simplexes, les $2$-objets sont les homotopies entre chemins (à extrémités fixées) etc. Pour les $n$-objets, cependant, on ne prend pas les homotopies entre homotopies de fourbis, mais classes d'équivalence de homotopies (modulo la relation d'homotopie) entre homotopies. La composition des $i$-objets ($i \geq 1$) se définit de fa\c{c}on évidente, on notera qu'elle n'est pas strictement associative, mais associative modulo homotopie. Donc la $n$-catégorie qu'on obtient n'est pas ``stricte'' - et on prévoit pas mal d'emmerdement pour définir de fa\c{c}on raisonnable une $n$-catégorie pas stricte (dans la description des compatibilités pour les ``données d'associativité). La mise sur pied du yoga qui suit pourrait constituer un fil d'Ariadne pour la définition en forme des $n$-catégories (pas strictes), les $n$-foncteurs entre elles (pas non plus stricts, et pour cause), les $n$-équivalences etc, au même titre que le yoga initial ``une $n$-catégorie est une catégorie ou les $\Hom$ et leurs accouplements de composition sont des $(n-1)$-catégories et des accouplements entre telles''. Cette $n$-catégorie $C_n(K_{\bullet})$ dépend fonctoriellement de $K_{\bullet}$, tout morphisme simplicial $K_{\bullet} \to K'_{\bullet}$ définit un $n$-foncteur $C_n(K_{\bullet} \to C_n(K'_{\bullet}))$ ; en fait, cela doit en dépendre même $n$-fonctoriellement, vu qu'on voit (en s'inspirant de ce qui précède et l'application à des ensembles semi-simpliciaux de la forme $\Hom_{\bullet}(K_{\bullet}, K'_{\bullet})$) que les ensembles semi-simpliciaux forment eux-mêmes les $0$-objets d'une $n$-catégorie, quel que soit $n$\dots

En fait, $C_n(K_{\bullet})$ est un $n$-groupoïde, i.e. une $n$-catégorie où toute $i$-flèche ($1 \leq i \leq n$) (= $i$-objet) est une ``équivalence'' i.e. admet un quasi-inverse (donc un inverse si la $n$-catégorie est ``réduite''). Si $C$ est une telle $n$-catégorie i.e. un $n$-groupoïde, et $X$ un $0$-objet de $C$, il s'impose de désigner par $\pi_i(C, x)$ ($0 \leq i \leq n$) successivement : l'ensemble des classes de $0$-objets à équivalence près de $1$-objets (ou $1$-flèches) $x \to x$ (c'est un groupe, pas nécessairement commutatif), l'ensemble des classes modulo équivalence des $2$-flèches $1_x \to 1_x$, où $1_x$ est la $1$-flèche identique de $x$ (c'est un groupe commutatif $\pi_2(C, x)$, ainsi que les groupes qui vont suivre), l'ensemble des classes modulo équivalence de $3$-flèches $1_{1_x} \to 1_{1_x}$, etc. Ces groupes forment, comme de juste, des ``systèmes locaux'' sur l'ensemble des $0$-objets de $C$, et modulo le grain de sel habituel, les $\pi_i(C, x)$ ne dépendent que de la ``composante connexe'' du $0$-objet $x$ i.e. de sa classe modulo équivalence de $0$-objets. Ceci dit, si $C$ est de la forme $C_n(K_{\bullet})$, il résulte pratiquement des définitions que l'on a des isomorphismes canoniques $\pi_i(K_{\bullet}, x) \simeq \pi_i(C_n(K_{\bullet}))$ pour $0 \leq i \leq n$, qui pour $x$ variable peuvent s'interpréter comme des isomorphismes de systèmes locaux. Il s'ensuit que pour une application semi-simplicial $f: K_{\bullet} \to K'_{\bullet}$, le $n$-foncteur correspondant $C_n(K_{\bullet}) \to C_n(K'_{\bullet})$ est une $n$-équivalence si et seule si $f$ induit un isomorphisme sur les $\pi_0$ et sur les $\pi_i$ en tout point ($1 \leq i \leq n$). On serait plus heureux de pouvoir dire à la place ``et de plus un homomorphisme surjectif pour $i = n + 1$, car c'est, il me semble, cela qu'il faudrait pour espérer pouvoir conclure que la catégorie localisée de la catégorie des ensembles semi-simpliciaux, obtenue en inversant les flèches ``qui induisent des isomorphismes sur les $\pi_i$ pour $0 \leq i \leq n$ (ou encore, ``en négligeant'' les ensembles semi-simpliciaux $n$-connexes), est équivalente à la catégorie localisée de la catégorie des $n$-catégories, où on rend inversibles les $n$-équivalences ? Quoi qu'il en soit, ces petites bavures devraient disparaître lorsqu'on ``stabilise'' en faisant augmenter $n$. A ce propos, on voit que le foncteur ``troncature en dimension $n$'' de la théorie homotopique (consistant à tuer les groupes d'homotopie à partir de la dimension $n + 1$) s'interprète dans la langage des $n$-catégories par l'opération faisant passer d'une $N$-catégorie ($N > n$) à une $n$-catégorie, en conservant tels quels les $i$-objets ($0 \leq i \leq n-1$) et leur composition ($1 \leq i \leq n-1$), et en rempla\c{c}ant les $n$-objets par les classes de $n$-objets ``à équivalence près'', avec la composition obtenue par passage au quotient. De même, le foncteur d'inclusion évident en théorie homotopique, consistant à regarder un ensemble semi-simplicial ``où on a négligé les $\pi_i$ pour $i \geq n + 1$'' comme un ensemble semi-simplicial (dans la catégorie homotopique) qui se trouve avoir des $\pi_i$ nuls pour $i \geq n + 1$, se traduit par le foncteur allant des $n$-catégories vers les $N$-catégories, obtenue en ajoutant à une $n$-catégorie des $i$-flèches ($n + 1 \leq i \leq N$) identiques exclusivement. (Ainsi, un ensemble est regardé comme une catégorie ``discrète'', une catégorie comme une $2$-catégorie où les $\bHom(A, B)$, $A$ et $B$ des $0$-objets, sont des catégories discrètes, etc\dots).

Bien entendu, rien n'empêche de considérer aussi la notion de $\infty$-catégorie, à laquelle celle de $n$-catégorie est comme la notion d'ensemble semi-simplicial tronqué à celle d'ensemble semi-simplicial. Sauf erreur, la localisée de la catégorie des $\infty$-catégories, pour les flèches de $\infty$-équivalence, est équivalente à ``la catégorie homotopique'', localisée de la catégorie des ensembles semi-simpliciaux, ou du moins une sorte de complétée de celle-là. Dans cette optique, le tapis consistant à interpréter une $\infty$-catégorie de Picard stricte (i.e. quelque chose qui ressemble à un groupe abélien de la catégorie des $\infty$-catégories) comme donnée (à $\infty$.équivalence près) par un complexe de chaînes regardé comme un objet d'une catégorie dérivée, est à relier au tapis de Dold-Puppe, interprétant ces derniers comme des groupes abéliens semi-simpliciaux.

Pour se donner confiance dans ce yoga général, on peut essayer d'interpréter en termes de $n$-catégories ou $\infty$-catégories des constructions familières en homotopie. Ainsi, l'espace des lacets $\Omega (K_\bullet, x)$ correspond manifestement à la $(n-1)$-catégorie $\bHom(x, x)$ formée des $i$-flèches de $C$ $(1 \leq i \leq n)$ dont la 0-origine et la 0-extrémité sont $x$, réindexées en les appelant $(i-1)$-flèches. Je n'aper\c{c}ois pas à vue de nez un joli candidat pour la suspension en termes de $n$-catégories. Par contre le $\Hom_\bullet (K_\bullet, K_\bullet')$ doit correspondre au $\bHom(C, C')$, qui est une $n$-catégorie quand $C$, $C'$ en sont. La ``fibre homotopique'' d'une application semi-simpliciale $f: K_\bullet \to K_\bullet'$ (transformée d'abord, pour les besoins de la cause, en une fibration de Serre par le procédé bien connu de Serre-Cartan) correspond sans doute à l'opération bien familière de produit $(n+1)$-fibré (du moins les cas $n = 0, 1$ sont bien familiers !) $C \times_{C'} C''$ pour des $n$-foncteurs $c \to C'$ et $C'' \to C'$, dans le cas où $C''$ est la $n$-catégorie ponctuelle, donc la donnée de $C'' \to C'$ correspond à la donnée d'un 0-objet de $C'$. Les espaces $K(\pi, n)$ ont une interprétation évidente comme $n$-gerbes liées par $\pi$. Enfin, on voit aussi poindre l'analogue du dévissage de Postnikoff d'un ensemble semi-simplicial - mais la fa\c{c}on dont je l'entrevois (vue ma prédilection pour les topos) passe par la notion de topos classifiant d'un $n$-groupoïde (généralisant de fa\c{c}on évidente le topos classifiant d'un groupe). En termes de cette notion, on peut, il me semble, interpréter un $n$-groupoïde en termes d'un $(n-1)$-groupoïde (savoir son tronqué), \emph{muni} d'une $n$-gerbe sur le topos classifiant, liée par $\pi_n$ (``fordu'' bien sûr par l'action du $\pi_1$...).

\hangsection{Relativisation sur un topos.}\label{sec:app4}%
Bien sûr, il faut relativiser encore tout le yoga qu'on vient de décrire, au dessus d'un topos quelconque $X$. Il s'agirait donc de mettre en relation et d'identifier, dans un certaine mesure, d'une part l'algèbre homotopique sur $X$, d'autre part l'algèbre catégorique sur $X$ construite en termes de la notion de $n$-champ en groupoïdes ($n \geq 0$ fini ou infini). On espère que la notion d'image inverse de faisceau semi-simplicial par un morphisme de topos $f: X \to X'$ (qui est évidente) correspond à la notion évidente d'image directe de $n$-champs; et inversement, la notion évidente d'image directe de $n$-champs par $f$ devrait correspondre à une notion plus subtile d'image directe $\Lf_* (K_\bullet)$ d'un faisceau semi-simplicial, construit sansa doute dans l'esprit des foncteurs dérivés à partir de la notion naïve (mais on hésite s'il faut mettre $\Lf_*$ ou $\Rf_*$)... Les dévissages à la Postnikoff doivent avoir encore une interprétation remarquablement simple en termes de $n$-champs. Comparer à la remarque de Giraud qu'un 1-champ en groupoïdes sur $X$ peut s'identifier au couple d'un faisceau $\pi_0$ sur $X$, et d'une 1-gerbe sur le topos induit $X_{/\pi_0}$ (dont le lien, comme de juste, devrait être note $\pi_1$ !). D'ailleurs, dans le cas des 1-champs en groupoïdes, la traduction de ces animaux en termes de topos classifiants au dessus de $X$ est, je crois, développé en long et en large dans Giraud (il parle, si je me rappelle bien, d'``extensions'' du topos $X$). L'extension (si j'ose dire) de ce tapis aux $n$-champs ne devrait pas poser de problème.

Remords : tâchant de préciser heuristiquement la notion de topos classifiant d'un $n$-champ en groupoïdes (ou plus particulièrement, d'un $n$-groupoïde) pour $n \geq 2$, je vois que je n'y arrive pas à vue de nez. (Bien sûr, il suffirait (procédant de proche en proche) de savoir définir un topos classifiant raisonnable pour une \emph{$n$-gerbe}, liée par un faisceau abélien $\pi_n$). Donc je ne sais comment décrire le dévissage de Postnikoff en termes de $n$-champs, sauf pour $n \leq 2$. Ceci est lié à la question d'une description directe des groupes de cohomologie d'un $n$-groupoïde $C$ (ou d'un $n$-champ), à coefficients disons dans un système local commutatif, de fa\c{c}on que pour $C = C_n(K_\bullet)$, $K_\bullet$ un ensemble semi-simplicial dont les $\pi_i$ pour $i \geq n + 1$ sont nuls, on trouve les groupes de cohomologie correspondants de $K_\bullet$. Peut-on le faire en associant à $C$, de fa\c{c}on convenable, un ensemble semi-simplicial ``nerf'' de $C$ ?

Bien entendu, si on réussit à définir un topos classifiant pour $C$, celui-ci devrait être homotope à $K_\bullet$ ci-dessus, donc avoir les mêmes invariants homotopiques $\pi_i$ et cohomologiques $H^i$ ; itou pour les champs. La définition habituelle du topos classifiant, dans le cas $n = 1$, a bien cette vertu. Cas particulier typique de problème de la définition du topos classifiant : pour $\pi$ un groupe commutatif, trouver un topos canonique (fonctoriel en $\pi$ bien sûr...) ayant le type d'homotopie de $K(\pi, n)$, et qui généralise la définition du topos classifiant pour $n = 1$ (topos des ensembles où $\pi$ opère). On frémit à l'idée que les topos pourraient ne pas faire l'affaire, et qu'il y faille des ``$n$-topos'' !! (J'espère bien que ces animaux n'existent pas...)

\hangsection{Ingrédients principaux vers une ``Algèbre Topologique''.}\label{sec:app5}%
La théorie ``d'algèbre homologique non commutative'' que j'essaie de suggérer pourrait se définir, vaguement, comme l'étude parallèle des notions suivantes et de leurs relations des notions suivantes et de leurs relations multiples: a) espaces topologiques, topos, b) ensembles semi-simpliciaux, faisceaux semi-simpliciaux etc. c) $n$-catégories (notamment $n$-groupoïdes), $n$-champs (notament $n$-champs en groupoïdes) etc. d) complexes de groupes abéliens, de faisceaux abéliens. (Les ``etc'' réfèrent surtout aux structures supplémentaires qu'on peut envisager sur les objets du type envisagé...). C'est donc de l'algèbre avec la présence constante de motivations provenant de l'intuition topologique. Si une telle théorie devait voir le jour, il lui faudrait bien un nom, je me demande si ``algèbre topologique'' ne serait pas le plus adéquat (``algèbre homologique non commutative'' ne peut guère aller à la longue, pour des raisons évidentes). Ce qui est aujourd'hui parfois désigné sous ce [] n'est guère qu'un bric à brac de notions (telles que anneau topologique, corps topologique, groupe topologique etc) qui ne forment guère un corps de doctrine cohérent - il ne s'impose donc pas que cela accapare un nom qui servirait mieux d'autres usages. (Comparer le nouvel usage du terme ``géométrie analytique'' introduit par Serre, et qui ne semble guère avoir rencontré de résistance.)

Re-salut, et au plaisir de te lire

\bigskip

\par\hfill Villecun le \ondate{17.2.}1975\par

Cher Breen\scrcomment{Translation in \textcite{Gletters}},

\hangsection{Champs essentiellement localement constants et types d'homotopie.}\label{sec:app6}%
Encore un ``afterthought'' à une lettre-fleuve sur le yoga homotopique. Comme tu sais sans doute, à un topos $X$ on associe canoniquement un pro-ensemble simplicial, donc un ``pro-type d’homotopie'' en un sens convenable. Dans le cas où $X$ est "localement homotopiquement trivial", le pro-objet associé est essentiellement constant en tant que pro-objet dans la catégorie homotopique, donc $X$ définit un objet de la catégorie homotopique usuelle, qui est son "type d’homotopie". De même, si $X$ est ``localement homotopiquement trivial en dim $\leq n$'', il définit un type d’homotopie ordinaire ``tronqué en dim $\leq n$'' - construction familière pour $i = 0$ ou $1$, même à des gens comme moi qui ne connaissent guère l’homotopie !

Ces constructions sont fonctorielles en $X$. D’ailleurs, si $f: X \to Y$ est un morphisme de topos, Artin-Mazur ont donné une condition nécéssaire et suffisante \emph{cohomologique} pour que ce soit une "équivalence d’homotopie en dim $\leq n$'' : c’est que $\mathrm H^i(Y, F) \tosim \mathrm H^i(X, f^*(F))$ pour $i \leq n$, et tout faisceau de groupes \emph{localement constant} $F$ sur $Y$, en se restreignant de plus à $i \leq 1$ dans le cas non commutatif. Ce critère, en termes de $n$-gerbes ``localement constantes'' $F$ sur $Y$, s’interprète par la condition que $F(Y) \to F(X)$ est une $n$-équivalence pour tout tel $F$ et $i \leq n$. Il est certainement vrai que ceci équivaut encore au critère suivant 
\begin{enumerate}
\item[(A)]\label{it:App6.A} Pour tout $n$-champ ``localement contant'' $F$ sur $Y$, le $n$-foncteur $F(Y) \to f^*(F)(X)$ est une $n$-équivalence;
\end{enumerate}
ou encore à
\begin{enumerate}
\item[(B)]\label{it:App6.B} Le $n$-foncteur $F \to f^*(F)$ allant de la $n$-catégorie des $(n-1)$-champs localement constants sur $Y$ dans celle des $(n-1)$-champs localement constants sur $X$, est une $n$-équivalence.
\end{enumerate}

En d’autres termes, les constructions sur un topos $X$ qu'on peut faire en termes de $(n-1)$-champs \emph{localement} constants ne dépendent que de son "(pro)-type d’homotopie $n$-tronquée", et le définissent. Dans le cas où $X$ est localement homotopiquement trivial en dim $\leq n$, donc définit un type d’homotopie $n$-tronqué ordinaire, on peut interpréter ce dernier comme un $n$-groupoïde $C_n$, (défini à $n$-équivalence près). En termes de $C_n$, les $(n-1)$-champs localement constants sur $X$ doivent s’identifier aux $n$-foncteurs de la $n$-catégorie $C_n$ dans la $n$-catégorie $(n-1)-\Cat$ de toutes les $(n-1)$-catégories. Dans le cas $n = 1$ ceci n’est autre que la théorie de Poincaré de la classification des revêtements de $X$ en termes du ``groupoïde fondamental'' $C_1$ de $X$. Par extension, $C_n$ mérite le nom de \emph{$n$-groupoïde fondamental de $X$}, que je propose de noter $\Pi_n(X)$. Sa connaissance induite donc celle des $\pi_i(X)$ $(0 \geq i \geq n)$ et des invariants de Postnikoff de tous les ordres jusqu’à $\mathrm H^{n+1}(\Pi_{n-1}(X), \pi_n)$.

Dans le cas d’un topos $X$  quelconque, pas nécessairement localement homotopiquement trivial en dim $\leq n$, on espère pouvoir interpréter les $(n-1)$-champs localement constants sur $X$ en termes d’un $\Pi_n(X)$ qui sera un pro-n-groupoïde. Ça a été fait en tous cas, plus ou moins, pour $n = 1$ (du moins pour $X$ connexe); le cas où $X$ est le topos étale d’un schéma est traité in extenso dans SGA 3\scrcomment{\textcite{SGA3}}, à propos de la classification des tores sur une base quelconque. 

Dans le cas $n = 1$, on sait qu'on récupère (à équivalence près) le $1$-groupoïde $C_1$ à partir de la $1$-catégorie $\bHom(C_1, \Ens)$ de ces foncteurs dans $\Ens = 0-\Cat$ (i.e. des ``systèmes locaux'' sur $C_1$ qui est un topos, dit "multigaloisien") comme la catégorie des ``foncteurs fibres'' sur le dit topos, i.e. la catégorie opposée à la catégorie des points de ce topos (lequel n’est autre que le \emph{topos classifiant} de $C_1$). Pour préciser pour $n$ quelconque la façon dont le $n$-type d’homotopie d’un topos $X$ (supposé localement homotopiquement trivial en dim $\leq n$, pour simplifier), i.e. son $n$-groupoïde fondamental $C_n$, s’exprime en termes de la $n$-catégorie des ``$(n-1)$-systèmes locaux sur $X$'' i.e. des $(n-1)$-champs localement constants sur $X$, et par là élucider complètement l’énoncé hypothétique $(B)$ ci-dessus, il faudrait donc expliciter comment un $n$-groupoïde $C_n$ se récupère, à $n$-équivalence près, par la connaissance de la $n$-catégorie $C_n = n-\bHom (C_n, (n-1)-\Cat)$ des $(n-1)$-systèmes locaux sur $C_n$. On aurait envie de dire que $C_n$ est la catégorie des ``$n$-foncteurs fibres'' sur $C_n$, i.e. des $n$-foncteurs $C_n \to (n-1)-\Cat$ ayant certaines propriétés d’exactitude (pour $n = 1$, c’était la condition d’être les foncteurs image inverse pour un morphisme de topos, i.e. de commuter aux $\varprojlim$ quelconques et aux $\varinjlim$ finies ...)  C’est ici que se matérialise la peur, exprimée dans ma précédente lettre, qu’on finisse par tomber sur la notion de $n$-topos et morphismes de tels ! $C_n$ serait un topos (appelé le "$n$-topos classifiant du $n$-groupoïde $C_n$), $(n-1)-\Cat$ serait le $n$-topos ``ponctuel'' type, et $C_n$ d’interprète modulo $n$-équivalence comme la $n$-catégorie des ``$n$-points'' du $n$-topos classifiant $C_n$. Brr !

Si on espère encore pouvoir définir un bon vieux $1$-topos classifiant pour un $n$-groupoïde $C_n$, comme solution d’un problème universel, je ne vois guère que le problème universel suivant : pour tout topos $T$, considérons $\bHom(\Pi_n(T), C_n)$. C’est une $n$-catégorie, mais prenons en la $1$-catégorie tronquée $\tau_1 \bHom(\Pi_n(T), C_n)$. Pour $T$ variable, on voudrait $2$-représenter le $2$-foncteur contravariant $\Top^{\circ} \to 1-\Cat$ par un topos classifiant $\B = \B_{C_n}$, donc trouver un $\Pi_n(\B) \xrightarrow{\phi} C_n$ $2$-universel en le sens que pour tout $T$, le foncteur
$$
\bHom_{\Top}(T, \B) \xrightarrow {u \mapsto \phi \circ \Pi_n(u)} \tau_1 \bHom(\Pi_n(T), C_n)
$$
soit une équivalence. Pour $n = 1$ on sait que le topos classifiant de $C_1$ au sens usuel fait l'affaire, mais pour $n = 2$ déjà, je doute que ce problème universel ait une solution. C'est peut-être lié au fait que le ``théorème de Van Kampen'', qu'on peut exprimer en disant que le $2$-foncteur $T \to \Pi_1(T)$ des topos localement $1$-connexes vers les groupoïdes transforme (à $1$-équivalence près) sommes amalgamées (et plus généralement commute aux $2$-limites inductives), n'est sans doute plus vrai pour le $\Pi_2(T)$. Ainsi, si $T$ est un espace topologique réunion de deux fermés $T_1$ et $T_2$, il n'est sans doute plus vrai  que la donnée d'un $1$-champ localement constant sur $T$ ``équivaut à'' la donnée d'un $1$-champ localement constant $F_i$ sur $T_i$ $(i = 1, 2)$ et d'une équivalence entre les restrictions de $F_1$ et $F_2$ à $T_1 \cup T_2$ (alors que l'énoncé analogue en termes de $0$-champs, i.e. de revêtements, est évidemment correct).

\hangsection{Intégration de champs et cohomologie.}\label{sec:app7}%
L'énoncé (B) plus haut rend clair comment expliciter la cohomologie d'un $n$-groupoïde $C_n$. Si $C_n = \Pi_n(X)$, et si $F$ est un $(n-1)$-champ localement constant sur $X$, $e^X_{n-1}$ est le $(n-1)$-champ ``final'', on a une $(n-1)$-équivalence de $(n-1)$-catégories
$$
\Gamma_X(F) = F(X) \simeq \bHom(e^X_{n-1}, F)
$$
qui montre que le foncteur $\Gamma_X$ ``intégration sur $X$'' sur les $(n-1)$-champs localement constants, qui inclut la cohomologie (non commutative) localement constante de $X$ en dim $\leq n-1$, s'interprète en termes de ``$(n-1)$-systèmes locaux'' sur le groupoïde fondamental comme un $\bHom(e^{C_n}_{n-1}, F)$ où maintenant $F$ est interprété comme un $n$-foncteur 
$$
C_n \xrightarrow{F} (n-1)-\Cat
$$
et $e^{C_n}_{n-1}$ est le $n$-foncteur constant sur $C_n$, de valeur la $(n-1)$-catégorie finale.

Pour interpréter ceci en notation cohomologique, il faut que j'ajoute, comme ``remords'' à la lettre précédente, l'interprétation explicite de la cohomologie non commutative sur un topos $X$, en termes d'intégration de $n$-champs sur $X$. Soit $F$ un $n$-champ de Picard strict sur $X$, il est donc défini par un complexe de cochaines $L'$ sur $X$
$$
0 \to L^0 \to L^1 \to L^2 \to ... \to L^n \to 0
$$
concentré en degrés $0 \leq i \leq n$ (défini à isomorphisme unique près dans la catégorie dérivée de $\Ab(X)$). Ceci dit, les $\mathrm H^i(X, L')$ (hypercohomologie) \emph{pour $0 \leq i \leq n$} s'interprètent comme $\mathrm H^i(X, L') = \pi_{n-i}\Gamma_X(F)$.

Si on s'intéresse à tous les $\mathrm H^i$ (pas seulement pour $i \leq n$) on doit, pour tout $N \geq n$, regarder $L'$ comme un complexe concentré en degrés $0 \leq i \leq N$ (en prolongeant $L'$ par des $0$ à droite).Le $N$-champ de Picard strict correspondant n'est plus $F$ mais $C^{N-n}F$, où $C$ est le foncteur ``espace classifiant'', s'interprétant sur les $n$-catégories de Picard strictes comme l'opération consistant à ``translater'' les $i$-objets en des $(i+1)$-objets, et à rajouter un unique $0$-objet; il se prolonge aux $n$-champs de Picard ``de fa\c{c}on évidente'', on espère, de fa\c{c}on à commuter aux opérations d'image inverse de $n$-champs. On aura donc pour $i \leq N$
$$
\mathrm H^i(X, L') = \pi_{N-i}\Gamma_X(C^{N-n}F) \quad i \leq N.
$$
Ceci posé, il s'impose, pour tout $n$-champ de Picard strict $F$ sur $X$, de poser 
$$
\boxed{
\mathrm H^i(X, F) = \pi_{N-i}\Gamma_X(C^{N-n}F) \quad~\text{si}~\quad N \geq i,n
}
$$
ce qui ne dépend pas du choix de l'entier $N \geq Sup(i, n)$ [{\textbf{NB}} On a un morphisme canonique de $(n-1)$-groupoïdes,
$$
C(\Gamma_X F) \to \Gamma_X(C F),
$$
comme le montrent les constructions évidentes en termes de complexes de cochaines, et on voit de même que celui-ci induit des isomorphismes pour les $\pi_i$ pour $1 \leq i \leq n+1$.]

{\textbf{NB}} On voit en passant que pour un $n$-champ en groupoïdes $F$ sur $X$, si on se borne à vouloir définir les $\mathrm H^i(X, F)$ pour $0 \leq i \leq n$, on n'a pas besoin sur $F$ d'une structure de Picard, car il suffit de poser
$$
\mathrm H^i(X, F) = \pi_{n-i}(\Gamma_X(F)) \quad 0 \leq i \leq n.
$$
Si d'autre part $F$ est un $n$-$\Gr$-champ (i.e. muni d'une loi de composition $F \times F \to F$ ayant les propriétés formelles d'une loi de groupe) le $(n+1)$-``champ classifiant'' est défini, et on peut définir $\mathrm H^i(X, F)$ pour $i \leq n+1$ par
$$
\mathrm H^i(X, F) = \pi_{n+1-i}(\Gamma_X(CF))
$$
en particulier
$$
\mathrm H^{n+1}(X, F) = \pi_0(\Gamma_X(CF)) =~\text{sections de}~CF~\text{à équivalence près}.
$$

Mais on ne peut former $CCF = C^2F$ et définir $\mathrm H^{n+2}(X, F)$, semble-t-il \emph{que} si $CF$ est lui-même un $\Gr$-$(n+1)$-champ, ce qui ne sera sans doute le cas que si $F$ est un $n$-champ de Picard strict...

Venons en maintenant au cas où $F$ est un $n$-champ \emph{localement constant} sur $X$, donc défini par un $(n+1)$-foncteur
$$
C_{n+1} \xrightarrow{F} n-\Cat.~\text{de Picard strictes}.
$$
Alors, posant pour $0 \leq i \leq n$
$$
\mathrm H^i(C_ {n+1}, F) = \pi_{n-1}(\bHom(e_n^{C_{n+1}}, F)),
$$
``on a fait ce qu'il fallait'' pour que l'on ait un isomorphisme canonique
$$
\mathrm H^i(C_{n+1}, F) \simeq \mathrm H^i(X, F),
$$
(valable en fait sans structure de Picard sur $F$...). Il s'impose, pour tout $\infty$-groupoïde $C$ et tout $(n+1)$-foncteur
$$
C \xrightarrow{F} n-\Cat.~\text{de Picard strictes}.
$$
de définir les $\mathrm H^i(C, F)$, pour tout $i$, par
$$
\mathrm H^i(C, F) = \pi_{N-i}\bHom(e_N^C, C^{N-n}F)
$$
où on choisit $N \geq Sup (i, n)$. Si $F$ n'a qu'une $\Gr$-structure (pas nécessairement de Picard) on peut définir encore les $\mathrm H^i(C, F)$ pour $i \leq n+1$ par
$$
\mathrm H^i(C, F) = \pi_{n+1-i}\bHom(e_{n+1}^C, CF).
$$
Dans le cas $C = C_{n+1} = \Pi_{n+1}(X)$, il doit être vrai encore (en vertu de (A) plus haut), que cet ensemble est canoniquement isomorphe à $\mathrm H^{n+1}(X, F) = \pi_0 \Gamma_X(CF)$ (c'est vrai et bien facile pour $n = 0$). Décrire la flèche canonique entre les deux membres de 
$$
\mathrm H^{n+1}(X, F) \simeq \mathrm H^{n+1}(\Pi_{n+1}X, F) \quad ?
$$
Si on veut réexpliciter (A) et (B), en termes du yoga (C), on arrive à la situation suivante:

On a un $(n+1)$-foncteur entre $(n+1)$-groupoïdes
$$
f_{n+1}: C_{n+1} \to D_{n+1}
$$
induisant par troncature un $n$-foncteur
$$
f_n: C_n \to D_n
$$
On doit avoir alors:
\begin{enumerate}
    \item[(A')] $f_n$ est une $n$-équivalence si et seule si le $n$-foncteur $\phi \to \phi \circ f_n$
    $$
    f_n^*: \bHom(D_n, (n-1)-\Cat) \to \bHom(C_n, (n-1)-\Cat)
    $$
\end{enumerate}
allant des $(n-1)$.systèmes locaux sur $D_n$ (ou $D_{n+1}$, c'est pareil) vers les $(n-1)$-systèmes locaux sur $C_n$, est une $n$-équivalence.
\begin{enumerate}
    \item[(B')] $f_n$ est une $n$-équivalence si et seule si pour tout $n$-système local $F$ sur $D_{n+1}$,
    $$
    F: D_{n+1} \to n-\Cat,
    $$
    le $n$-foncteur induit par $f_{n+1}$
    $$
    \underbrace{\bHom(e_n^{D_{n+1}}, F)}_{\hspace*{-5mm}\Gamma_{D_{n+1}(F)}\hspace*{-5mm}}\; \to \underbrace{\bHom(e_n^{D_{n+1}}, f^*_{n+1} F)}_{\hspace*{-5mm}\Gamma_{C_{n+1}(F)}\hspace*{-5mm}}\;
    $$
    est une $n$-équivalence.
\end{enumerate}

\hangsection{Les trois approches vers la cohomologie d'un topos.}\label{sec:app8}%
La construction de la cohomologie d'un topos en termes d'intégration des champs ne fait aucun appel à la notion de complexe de faisceaux abéliens, encore moins à la technique des résolutions injectives. On a l'impression que dans son esprit, via la définition (qui reste à expliciter !) des $n$-champs, elle s'apparenterait plutôt aux calculs ``Cechistes'' en termes d'hyperrecouvrements. Or ces derniers se décrivent à l'aide d'une petite dose d'algèbre semi-simpliciale. Si oui, cela ferait essentiellement trois approches distinctes pour construire la cohomologie d'un topos : 
\begin{enumerate}
\item[a)] point de vue des complexes de faisceaux, des résolutions injectives, des catégories dérivées (\emph{algèbre homologique commutative});
\item[b)] point de vue Cechiste ou semi-simplicial (\emph{algèbre homotopique});
\item[c)] point de vue des $n$-champs (algèbre catégorique, ou \emph{algèbre homologique non-commutative}).
\end{enumerate}
Dans a) on ``résoud'' les coefficients, dans b) on résoud l'espace (ou topos) de base, et dans c) en apparence on ne résoud ni l'un ni l'autre. 

Bien cordialement,

\bigskip

\par\hfill Villecun \ondate{17/19} July 1975\par

Dear Larry,

\hangsection{Question about a non-abelian Dold-Kan theorem.}\label{sec:app9}%
I am happy to finish by receiving an echo to my long letter and even a beginning to a constructive approach to a theory of the type I envisaged. The construction which you propose for the notion of a non-strict $n$-category, and of the nerve of the functor, has certainly the merit of existing, and of being a first precise approach, but otherwise can be subject to some evident criticism: it is very technical, unintuitive (yet at the level of $1-\Cat$, etc, and even of $2-\Cat$, everything is so clear ``you just follow your nose...''). And finally the absence of a definition of a functor sending (semi-)simplicial sets to $n$-groupoids. This functor correspond to a geometric intuition so clear that a theory which does not include it seems to me kind of a joke! Perhaps in trying to write down (like a sort of list of Christmas presents!) in a complete and explicit enough way the notions which one would like to have at ones disposal, and the relations (functor, equivalence, etc.) which should link them, one would arrive finally at a kind of axiomatic description sufficiently complete which should either give the key to a explicit \emph{ad hoc} construction, or should permit at least to enunciate and prove a theorem of existence and uniqueness\footnote{As was seen in section 9, ``uniqueness'' here has to b understood in a considerably wider sense than I expected, when writing this letter to Larry Breen. It now appears that the whole theory of stacks of groupoids will depend on the choice of a ``coherator'' $\bC$, as seen in section 13.} for a theory of the required type.

Otherwise, not having understood the idea of Segal in your last letter (which I have generously sent to Illusie\dots), I do not see how you define the Picard $n$-categories - but this matters little. As far as ``strict'' Picard $n$-categories are concerned, all I ask of them is that they finally form an $(n + 1)$-category $(n + 1)$-equivalent to that of chain complexes of length $n$. Agreed? I thank you for having rectified in my mind a big blunder, due to my great ignorance of algebraic topology and homotopy - I was in fact of the impression that $H$-spaces satisfying conditions of associativity and commutativity strict enough (say equivalent to an $\Omega^i X$ with $i$ arbitrarily large) correspond to commutative topological groups (inspired by several analogies\dots). Thus I am entirely in agreement with your observations on p.5.

On the other hand, I am still intrigued by the following question: is there an analogue of the ``tapis'' of Dold-Puppe\footnote{Tim Porter pointed out to me that ``Dold-Puppe'' is an inaccuracy name for this basic theorem, which should be called \emph{Dold-Kan theorem.}} for semi-simplicial groups (\emph{not necessarily commutative}) and what form should it take? To tell the truth I consider the yoga 
$$
\text{simplicial sets} \leftrightarrow \infty-\text{groupoids}
\leqno{(*)}
$$
as being essentially the ultimate ``set theoretic'' version of Dold-Puppe, which I would deduce from (*) by making explicit solely the fact that the abelian groups in $\infty-\Cat$ are ``nothing else'' than the chain complexes in $\Ab$. One should therefore first determine what should be the groups in $\infty-\Cat$. I can tell you what these are in $1-\Cat$ . this will be discussed at length in the book of Mme. Sinh\scrcomment{\textcite{GCS}}, I think in the chapter ``\emph{strict} $\Gr$-categories'' (i.e. the isomorphisms of associativity, for unity and inverse $X X^{-1} \simeq 1$ are \emph{identities}). One can make explicit for example how (\emph{via} the fact that a $\Gr$-category is $\Gr$-equivalent to a strict $\Gr$-category) the calculation with the $\Gr$-categories reduces to a very algebraic calculation with the \emph{strict} $\Gr$-category, by a kind of ``calculus of fractions'' (by choice, left or right) of the type which is used in giving the construction of derived categories. In any case, here is the explicit formulation of the structures (groups in $1-\Cat$) in terms of the theory of groups ($1$-categories in $\Gr$\footnote{AS was pointed out to me by Ronnie Brown, this structure was already well-known to J.H.C. Whitehead, under the name of ``crossed module'', and extensive use and extensive generalizations of this notion (in quite different directions from those I was having in mind, in terms of $\Gr$-stacks over an arbitrary topos) have been made by him and others. With respect to the question on next page, of generalizing this notion of ``non-commutative chain complex'' from length one to length two, Ronnie says there is a work in preparation by D. Conduché ``Modules croisés généralisés de longueur 2''.}). The structure is described by a quadruplet $(L_1, L_0, d, \theta)$ with 
$$
L_1 \xrightarrow{d} L_0
$$
a homomorphism of ordinary groups, 
$$
\theta: L_0 \to \Aut_{\Gr}(L_1)
$$
an operation of $L_0$ on $L_1$, with the following two axioms:
\begin{enumerate}
    \item[(a)] $d$ commutes with the operation of $L_0$, when $L_0$ acts on $L_1$ via $\theta$ and on itself by inner automorphisms:
    $$
    d(\theta(x_0) x_1) = \text{int}(x_0) d(x_1)
    $$
    \item[(b)] $\theta (d()x_1) = \text{int}(x_1)$.
\end{enumerate}
These properties imply that $\text{Im}d$ is normal in $L_0$ (hence $\Coker d = \pi_0$ is defined) and $\pi_1 = \Ker d$ is central in $L_1$, and finally that $L_0$ operates on $L_1$ leaving $\pi_1$ invariant, and it operates \emph{via} $\pi_0$. The principal cohomological invariant of this situation is evidently the Postnikoff-Sinh invariant
$$
\alpha \in \mathrm H^3 (\pi_0, \pi_1).
$$
I have met these animals - without even looking for them - in many situations, which I will not list now (I came across them recently \emph{a propos} the classification of ``ordinary'' formal groups over a perfect field, in terms of \emph{affine} algebraic groups, and \emph{commutative} formal groups, related by the strict $\Gr$-structures of this type (except that one has to use this formalism in an arbitrary topos (not merely in $\Sets$)) - to make explicit the yoga that ``the transcendent character of a formal group is concentrated essentially in the commutative formal groups'', discovered it seems by Dieudonné\dots). The question which I wish to raise is the generalisation to groups in $n-\Cat$, where I expect to find a non-commutative chain complex
$$
L_n \to L_{n-1} \to \dots \to L_1 \to L_0 \to 1
$$
with supplementary structures doubtless of the type of $\theta$, but what are they? It is understood that the topological significance of such structures is that they express exactly the ``truncated homotopy type in dim $\leq n$'' of topological groups, or equivalently the homotopy type in dim $\leq n + 1$ of pointed connected topological spaces\dots). Have you candidates to propose ?

\starsbreak

\hangsection{The ``six operations'' and homology.}\label{sec:app10}%
Your reflections on biduality and homology, however formal, tie in with a crowd of developments, of which only some exists at present, and others would demand considerable work still. Here are the reminiscences which your naive questions bring to mind:
(A) The formalism of the $\Rf_!$, $\Rf^!$ (combined with $\Rf_*$, $\Lf^*$, $\Lotimes$ and $\RHom$, ``the six operations'') carries implicitly in itself the definition of homology and the essential identity between homology and cohomology. One now has this formalism for quasi-coherent sheaves on schemes - seminar Hartshorne (Springer L.N. 20) - for the topological spaces and arbitrary sheaves of coefficients - Verdier, \emph{exposé} Bourbaki (SNLM 300) - and for the étale cohomology of schemes for ``discrete'' coefficients (``$\ell$-adic'' or torsion) prime to the residual characteristic (SGA 5), finally, for coherent sheaves on analytic spaces (Verdier-Ruget). (The formalism remains to be developed in the crystalline context, and in the characteristic 0 in the context of stratified modules with singularities, à la Deligne, with perhaps - over the field $\mathbf{C}$ - the introduction of additional Hodge structures, finally in the context of motives; I am convinced that it exists about anywhere - maybe, wherever there is a formalism of a cohomological nature.)
    
    Working in étale cohomology on a separated scheme of finite type over a field $k$, say, with a ring of coefficients $\Lambda$ of torsion prime to the characteristic, the complex of sheaves $f^!(\Lambda_e)$ (where $\hat{f}: X \to \Spec k = e$) plays the role of \emph{complex of singular chains on $X$ with coefficients in $\Lambda$}, and $\Rf_!$ $(f^! \Lambda_e)$ plays the role of a \emph{homology $\mathrm{H}_*(X/e)$}, \emph{vis a vis} of course, of coefficients on $e$ which are complexes of $\Lambda$-modules. You can easily justify this assertion with the help of the ``global duality theorems'', by one or two tricks which I spare you here.

{\textbf REMARKS}.
\begin{enumerate}
    \item[(1)] There is no need to truncation, it works in all dimensions.
    \item[(2)] This is related (at least as far as the philosophy is concerned) to he fact that for the various types of coefficients (under conditions of ``constructibility'') one has a theorem of ``biduality'', at least if one allows resolution of singularities (but Deligne has told me I believe that he knows a proof without that), with values in a ``dualizing complex'' $K_e$ (on $e$), $K_X$ (on $X$). If for example $\Lambda$ is ``self-dualising'' (or Gorenstein) for example $\Lambda = \mathbf{Z}/n \mathbf{Z}$, one can take $K_e = \Lambda$, therefore the dualising complex $K_X = f^! (K_e)$ is nothing else than the ``complex of singular chains with coefficients in $\Lambda$.
    \item[(3)] One can do the same thing for coefficients such as $\mathbf{Z}_{\ell}$ (Jouanolou, thesis non published\scrcomment{\textcite{JOU69}}, I fear!)
    \item[(4)] This works also for $f: X \to S$ finitely presented separated if $f$ has the properties of ``cohomological local triviality'' (properties ``local upstairs'') for example $f$ \emph{smooth}; one finds that $\mathrm{H}_*(X/S) = \Rf_! f^!(\Lambda_S)$.
\end{enumerate}

\starsbreak

(B) Artin-Mazur have studied in a spirit close to yours the \emph{autoduality} of the Jacobian of a relative curve $X/S$. It is necessary to ask them for precise results, perhaps it works say if $X/S$ is proper and flat or relative dimension 1 - in any case it is OK on a discrete valuation ring with smooth \emph{generic} fibre. The special fibre could be very wild. (I have used their results in SGA 7 to prove, in the case of Jacobians, a duality conjecture on the group of connected components associated to the Neron models of abelian varieties dual one to another\dots). 
Towards the end of the 50's (beginning of 60's?), when the grand cohomological stuff ($f^!$, $f^!$, étale cohomology, etc.) just came out from darkness, the course given by Serre on the theory of Rosenlich and Lang on generalised jacobians and the geometric class field theory (see Serre's book) and later the ``geometric'' theory of \emph{local} class field theory making use of pro-algebraic groups (see his article on this subject), made me reflect on the cohomological formulations of these and other results, which should be of a ``geometric'' nature, such that the ``arithmetic'' results over an arbitrary base field (or residue field) $k$ (finite, for example) follow immediately by descent from the ``geometric'' case of base field $\overline{k}$. I exchanged letters with Serre - I don't know if I can find copies - but I recall that I sketched projects for some ambitious enough theories on generalised residues, generalised local jacobians, etc.,  in at least three different directions. But I have never, in spite of numerous attempts, succeeded in mobilising someone for developing one of these programmes. Here a few words on them: 

\starsbreak

\hangsection{Complex of generalised jacobians.}\label{sec:app11}%
(C) In the situation where $X$ is of finite type over a \emph{field} $k$, construction of a complex of generalised jacobians $J_{* X/k}$ (of length equal to dim $X$).
    
    This is a complex of affine commutative pro-algebraic groups on $k$, with the exception of $J_0$ if I remember well, ($J_0$ had as abelian part the abelian part of $\Alb_{X/k}$, the usual generalised jacobian). It's construction, inspired by the residual complex, passes by generalised jacobians (in an appropriate cohomological sense) of the localisation $\Spec_{\mathbf{O}_{X, x}}$ of $X$ at its different point.
N.B. $\mathbf{H}_O (J_*)$ was the ``generalised Jacobian'' of $X$, i.e. there existed a homeomorphism $X \to \mathbf{H}_0$, which was universal for homomorphisms of $X$ into commutative locally proalgebraic groups. For $X$ connected, $\mathbf{H}_O$ is an extension of $\mathbf{Z}$ by an appropriate proalgebraic group. It is possible that, at first, I restricted to the case of $X$ smooth.

The cohomology role of this complex was that of a complex of \emph{homology}
$$
\mathrm{H}^i(X, G_X) \simeq \Ext^i(J_{* X/k}, G)
\leqno{(*)}
$$
but for which coefficients? I believe I took arbitrary commutative algebraic groups $G$ but worked with the Zariski topology (malédiction !). Even in the case of discrete $G$, I considered the Zariskian $\mathrm{H}^i$, this gives slightly stupid cohomology groups, evidently. I realised that one should work ultimately in étale cohomology, and that the construction of the $(J_i)_{X/k}$ will evidently be modified accordingly. As for the significance of the $\Ext^i$ (hypercohomology), at a moment where Serre had developed the formalism for proalgebraic groups, one was not too fearful of taking it in the category of such objects - and in the sense of a ``derived category'' which at that moment had never yet been explicitly defined and studied. (We have, after all, somewhat progressed since those days!). I have the impression, in view of these antique cogitations, heuristic as they were, that it should now be possible to develop at present such a theory of $J_{* X/k}$, in cohomology fppf, giving a formula (*) without limitation on the degree $i$ of the cohomology. (N. B. But $J_*$ evidently no longer stops in dim $X = n$ but in dim $2n$. It is nevertheless possible that the components $J_i$ might be of dim 0 for $i > n$).

I believe that the construction of the $J_*$ does not commute with base change, but merely does so in the derived category sense.

\starsbreak

\hangsection{Global ``geometric'' class field theory as a cohomological duality formula. Serre duality and the ``Lang trick''.}\label{sec:app12}%
(D) Let $X/k$ be a smooth scheme (for simplicity) over a field $k$, separated and of finite type, or relative dimension $d$, and $n$ an integer $> 0$. If $n$ is prime to the characteristic and if $F$ is a sheaf of coefficients on $X$ which is annihilated by $n$, the global duality tells us that $\Rf_!(F)$ and $\Rf_*(\RHom(F, \mu_n^{\otimes d}))$ ($\mu_n$ = sheaf of n-th roots of unity = $\Ker(G_m \xrightarrow{n} G_m)$) are dual to each other with values in $(\mathbf{Z}/n\mathbf{Z})_k$, for example $\Rf_!(\mathbf{Z}/n\mathbf{Z})$ and $\Rf_*(\mu_n^{\otimes d})$, or $\Rf_!(\mu^{\otimes}_m)$ and $\Rf_*(\mathbf{Z}/n\mathbf{Z})$, are dual to each other - at least with a shift of amplitude $2d$ in dimension. (As $\mathbf{Z}/n\mathbf{Z}$ is injective over itself, this gives in fact perfect duality
    $$
    \mathrm{R}^if_!(F) \times \mathrm{R}^{2d-i}f_*(\RHom(F, \mu_n^{\otimes d})) \to \mathbf{Z}/n\mathbf{Z}.)
    $$
    If now one no longer assumes $n$ prime to he characteristic, for example $n$ is a power of $p$ = characteristic of $k > 0$), it seems that everything collapses: to start with, one no longer knows (for $d > 1$) by what to replace $\mu_n^{\otimes d}$\dots
The extraordinary miracle is that for $d = 1$, i.e. $X$ a smooth curve, everything continues to work perfectly, provided one states things with care! The first verifications are made for example with $F = \mathbf{Z}/p\mathbf{Z}$, $\mu_p$, or $\alpha_p$, with $X$ complete - one finds it's O.K. by virtue essentially of the autoduality of the jacobian. One can make these examples more sophisticated on taking \emph{twisted} coefficients, and $X$ not complete - one convinces oneself this works always! Simply, it is necessary to note that here the $\mathrm{R}^if_*(F)$, $\mathrm{R}^if_!(F)$ have a ``continuous'' structure (they are essentially poalgebric groups). This corresponds to the well known phenomenon in class field theory that the structure of $\pi_{1\text{ab}}$ of $X$, when $X$ is not complete, is \emph{continuous} - hence same holds for $\mathrm{H}^1(X, \mathbf{Z}/p^n\mathbf{Z})$ say.

By the way, I point out for you that Serre once proposed (without ever writing it down, I think) a theory of duality for \emph{commutative unipotent} algebraic groups, \emph{modulo radical isogeny}, duality with values in $\mathbf{Q}/\mathbf{Z}$ (or $\mathbf{Q}_p/\mathbf{Z}_p$). He found that if (when $k$ is algebraically closed, say) $G$ is such a group, then $G' = \Ext^1(G, \mathbf{Q}/\mathbf{Z})$ can canonically be given a structure of quasi-algebraic group (i.e. defined modulo radical isogeny), doubtless in a unique manner provided it verifies some functorial properties, and on requiring that for $G = \mathbf{G}_a$ one finds that $\Ext^1(\mathbf{G_a}, \mathbf{Q}/\mathbf{Z}) \simeq \mathbf{G}_a$ with the usual structure. Let $\Delta G = G' = \Ext^1(G, \mathbf{Q}/\mathbf{Z})$. One finds $G \simeq \Delta\Delta G$ i.e. $\Delta$ is an authentic autoduality! I call $\Delta$ \emph{Serre duality}. It surely goes over to ind-progroups on an arbitrary base field (not necessarily algebraically closed) in the case $p > 0$. Moreover, for finite étale groups, it is $\Ext^0(G, \mathbf{Q}/\mathbf{Z})$ (pontrjagin duality) which gives a perfect duality. One could screw together, in an appropriate derived category, Serre duality and Pontrjagin duality, by taking $G \mapsto \Delta G = \RHom(G, \mathbf{Q}/\mathbf{Z})$: one calls this (``cohomological'') Serre duality. This will be a magnificent autoduality, if one puts oneself in a derived category where the $\mathbf{H}^i$ of the envisaged complexes are (up to passing to the limit) extensions of étale groups by connected unipotent groups. Now one gets only such complexes, by ``integrating'' finite coefficients $F$ on $X$ by $\Rf_!$ or $\Rf_*$. This being said, by passing to the limit in the initial formulation (or equivalently by replacing the $(\mathbf{Z}/n\mathbf{Z})_k$, previously considered, by $(\mathbf{Q}/\mathbf{Z})_k$ on $k$, and forming $f^!(\mathbf{Q}/\mathbf{Z})_k = (\mu_\infty)_X$) the duality formula takes the form
$$
\Delta(\Rf_!(F)) \simeq \Rf_*(\text{DF}[2]) \quad \text{``shift'' of dimension}
$$
where $D$ is the ``Cartier duality'' $\RHom(F, \mu_\infty)$ (or $\RHom(F, \mathbf{G}_m)$ if one prefers?), and $\Delta$ is the Serre duality: cohomology with proper supports and with arbitrary supports are exchanged by duality, when one takes upstairs Cartier duality, and downstairs Serre duality.

The validity of the duality formula is not open to doubt - the principal work for establishing it consist certainly in a careful description of the category of coefficients with which one is working, as well on $X$ as on $k$, and of the functors $D$ and $\Delta$. As the definition of an arrow is immediate, once the building of the machine has been accomplished, the validity of the formula should result without difficulty from the usual ``dévissages'' which allow one to verify the duality in the particular standard cases $F = \mathbf{Z}/p\mathbf{Z}$, $\mu_p$, $\alpha_p$ on a smooth, complete $X$. (N.B. the case of coefficients prime to the characteristic is already known.) Let us make explicit what the formula of duality says for $\mathrm{R}^1f_*(\mathbf{G}_X)$, where $G$ is a finite group étale on $k$ (the most important case being $G = (\mathbf{Z}/p^m\mathbf{Z})_k$); one recovers Serre's description of ``geometric class field theory'' in terms of extensions by $G$ of a generalised jacobian of $X$. Thus, the duality formula can be understood as a cohomological version, considerably enriched, of geometric class field theory. When the base field $k$ is finite, to retrieve the class field theory in the classical form, one can use ``the trick of Lang'' (on the relation between the ``arithmetic'' $\pi_1$ of a smooth, connected commutative algebraic group $J$ on $k$ and its $\mathrm{H}^0(k, J) = J(k)$: the $\pi_1^{\text{ar}}(J)$ classifies the isogenies above $J$ with kernel a constant group $\pi_1^{\text{ar}}(J) \simeq \mathrm{H}^0(k, J)$) - i its cohomological form, which may be stated: 
$$
\Delta_0\RGamma_K(J^*) \simeq \RGamma_K(\Delta J^*[1]),
$$
where $\Delta$ is Serre duality, $\Delta_0$ Pontrjagin duality for the totally disconnected topological abelian groups (duality with values in $\mathbf{Q}/\mathbf{Z}$), $J^*$ a complex of algebraic ind-progroups on $k$. Taking account of this ``Lang duality formula'' and applying $\RGamma_K$ to the formula of duality for geometric class fields, one gets the ``duality formula of arithmetic class field theory'':
$$
\Delta_0 (\mathrm{H}_!(X, F)) \simeq \mathrm{H}^*(X, D(F) [3])
$$
(isomorphism of totally disconnected topological groups).

Another remark: when $F$ is not an ``étale sheaf'', but has a continuous structure such as $\alpha_p$, one must be careful in the definition of $\Rf_!(F)$, for $X$ non complete, starting from the compactification $\widetilde{X}$; thus, if $F$ comes from an ``admissible'' sheaf $\widetilde{F}$ on $\widetilde{X}$, one must have an exact triangle
$$
\begin{tikzcd}[row sep=tiny]
  & \Rf_!(\hat{\widetilde{F}}) \ar[dl] \\ \Rf_!(F) \ar[rr] & & \mathrm{R}\widetilde{f}_*(\widetilde{F}) \ar[ul],
\end{tikzcd}
$$
where $\hat{\widetilde{F}}$ is the \emph{formal completion} of $\widetilde{F}$ along $\widetilde{X} - X$ (a finite number of points\dots). It is here, unless I am mistaken, that appears the link with local class field theory, in its cohomological version, on which I am going now to say a few words.

\starsbreak

\hangsection{Case of local ``geometric'' class field theory.}\label{sec:app13}%
(E) \textbf{Local class field theory as a duality formula}
    
    Let $V$ be a complete discrete valuation ring with residue field $k$ - assume either that $k$ has been lifted to $k \subset V$ (and therefore $V \simeq k[[T]]$) or that $k$ is perfect of characteristic $p > 0$. In order to fix ideas, and to be sure that I'm on solid ground, I consider at first on $K$ (= the field of fractions of $V$) \emph{finite} coefficients $F$ (as on $X$ previously) and I consider the objects $\mathrm{H}^1(K, F)$, or $\RGamma_K(F)$. The main work to be done consists in defining an adequate category of coefficients over $k$ (perhaps the same one as in (D)) and a functor
    $$
    F \mapsto \mathrm R \underline{\Gamma}_K(F)
    $$
    with values in the category of such coefficients, in such a manner that the following isomorphisms holds.
    $$
    \RGamma_K(F) \simeq \RGamma(\mathrm R \underline{\Gamma}_K(F)).
    $$
    This correspond to the intuition (acquired directly from elementary examples) according to which for $k$ algebraically closed, say, the $\mathrm{H}^0(K, F)$, $\mathrm{H}^1(K, F)$\dots are endowed with a structure of $k$-algebraic group (ind-pro\dots). In this construction, the ring scheme of Witt vectors over $k$ (introduced by Serre) and the ``Greenberg functor'' (associating to a $V$-scheme a $k$-prescheme) will play an essential role.
    
This being done, the duality formula will be formally stated as in (D) above: 
$$
\Delta \mathrm R \underline{\Gamma}_K(F) \simeq \mathrm R \underline{\Gamma}_K(\text{DF}[1])
$$
where $D$ stands for Cartier duality, $\Delta$ for Serre duality. When the residue field is finite, it becomes (via ``Lang's trick'' mentioned previously)
$$
\Delta_0 \RGamma_k(F) \simeq \RGamma_K(\text{DF}[2])
$$
$\Delta_0$ standing for Pontrjagin duality. The formula contains local geometric class field theory à la Serre, and arithmetical local class field theory in its classical form.

\textbf{Remarks}.

(a) If $F$ is prime to the residue characteristic the formula is very easy to prove and well known. It may be considered a very special case of the ``induction formula'' for a morphism $i: s \mapsto S$, in the duality formalism:
$$
i^!(D_S(F)) = D_S(i^*(F))
$$
(we take here the inclusion of $p = \Spec(k)$ in $S = \Spec(V)$). Thus the work to be done concerns the $p$-primary coefficients, for $p =$ characteristic $k > 0$. The most subtle case is that of unequal characteristic.

(b) The functor $\RGamma$ may be obtained by composing $\mathrm{R}j_*$ (where $j: U = \Spec(K) \to \Spec(V) = S$ is the inclusion) with a cohomological version of the ``Greenberg functor''.

(c) In (D) and (E), I restricted myself to finite coefficients $F$ - it's for those that I am sure of what I assert. But it is certainly true that the duality formula is even richer, that something may still be asserted for example for $F$ a not necessarily finite group scheme, for example an abelian scheme (with a few degenerate fibres in the case of (D)?), but I have never entirely clarified this question, even on a heuristic basis. I vaguely recall a formula which should be contained in the formalism (say if $k$ is algebraically closed): for $F$ an abelian scheme on $K$, $F$ the dual abelian scheme and $G'$ the pro-algebraic group over $k$ attached ``à la Greenberg'' to its Néron model, then one has
$$
\mathrm{H}^1(K, F) \stackrel{?}{\simeq} \Ext^1_{k-\text{grp}}(G', \mathbf{Q}/\mathbf{Z})
$$
(N.B. without any guarantee.) In principle, the previously mentioned duality conjecture concerning Néron models of SGA 6 should come out of the local duality machine.

(d) You may ask Deligne if he didn't dive into questions (D) and (E) lately.

\hangsection{Comments on fppf cohomology versus crystalline cohomology.}\label{sec:app14}%
(F) \textbf{Significance and limitations of the fppf topology}
    
    Since the attempts of Serre to find a ``Weil chomology'' by using the cohomology of a scheme with coefficients not only discrete $\mathbf{Z}/p^n\mathbf{Z}$ ($n \to \infty$) or $\mu_{p^n}$ ($n \to \infty$), but also continuous (for example $W_n$, $n \to \infty$), which give good results for recovering a correct $\mathrm{H}^1$, during numerous years I have come upon the impression, which I have tried in vain to make precise, that a correct ``$p$-adic'' ``Weil cohomology'', in the case $p > 0$ and $k$ of characteristic $p$, should come, in one way or another, from the fppf cohomology, for finite coefficients for example, or more general coefficients, e.g. algebraic groups over $k$. The construction in (B) of the local jacobian complex was, of course, related to this hope: the homology might reveal what is hidden to us in cohomology! For some time now, one has at ones disposal the formalism of crystalline cohomology, and one knows (Berthelot) that (at least for $X$ projective and smooth) it has the correct properties. If one uses that as a kind of standard by which to ``measure'' the other cohomologies, one finds that the part of the crystalline cohomology $\mathrm{H}^i_{\text{cris}}(X)$ which could be described in terms of fppf cohomology of $X$ with coefficients in algebraic $k$-groups is a small part of $\mathrm{H}^i$ only; more precisely, using the very rigid supplementary structure of the $\mathrm{H}^i$ (modules of finite type on the ring $W(k)$ of Witt vectors) which comes from the existence of the Frobenius homomorphism (an isogeny), $\mathrm{H}^i \xrightarrow{F} \mathrm{H}^i$ (semi-linear), one finds that one keeps always in the part ``of slope $\leq 1$ (although the possible slopes vary between 0 and $i$\dots). This explains why for $i = 1$ one can obtain via fppf a correct $\mathrm{H}^i$, although for $\mathrm{H}^2$ already all the attempts have been unfruitful.
In truth, one conjectures that \emph{all} the part of slope $\leq 1$ in $\mathrm{H}^i_{\text{cris}}$ comes from fppf. But I have completely lost contact with these questions - people such as Mazur, Kats, Messing - and of course Deligne - should be knowledgeable as to the present states of these questions.

\starsbreak

\hangsection{The homotopical trinity: ``spaces'', ``$\infty$-groupoids'', (essentially locally constant) ``coefficients categories''\dots}\label{sec:app15}%
Your question 7 seems to indicate that there is a misunderstanding on your part on the significance of the ``homotopy type'' of $X$, for $X$ a topos (for example the étale topos of a scheme). Doubtless you must be confusing the homotopical algebra which one can perform on $X$, using semi-simplicial sheaves, stacks of all kinds, the relations between these - and the other point of view according to which $X$ (with its very rich structure of topos) virtually disappears so as to become no more than a pale element of a ``homotopical category'' (or pro-homotopical), deduced from the topos by a very thorough process of ``localisation''. At first sight, all that still remains with poor stripped $X$, are the $\pi_i$ - and its cohomology groups with constant coefficients - or at the worst twisted constant coefficients. When one digs more into this definition of ``what is left to this poor $X$'' one falls precisely on \emph{the locally constant $n$-stacks} (as an $f: X \to X'$ which is a homotopy equivalence induces a $(n + 1)$-equivalence between the categories of locally constant $n$-stacks on $X$ and on $X'$) - which of course contain the abelian chain complexes of length $n$ of sheaves with locally constant cohomology sheaves, and the hyper-cohomology of these. It is thus that one arrives at this triangle of objects which mutually determine each other
\[\begin{tikzcd}[column sep=-3em,every arrow/.append style={<->}]
 & \begin{tabular}{@{}c@{}}
    topos (or topological space\\
    or semi-simplical complex)\\
    \textit{modulo} $n$-homotopy
   \end{tabular}
   \arrow[dl]\arrow[dr]&\\
 \begin{tabular}{@{}l@{}}
    $n$-groupoids\\
    (up to $n$-equivalence)
   \end{tabular}
   \arrow[rr]&
 & \begin{tabular}{@{}l@{}}
    ``special'' $n$-categories\\
    (up to $n$-equivalence)
   \end{tabular}
\end{tikzcd}
\]
One says that an $n$-category $E_n$ is ``special'' (or \emph{$n$-galois}) if it is $n$-equivalent to the category of locally constant $(n-1)$-stacks on an appropriate topological space (or a topos), or, what should be equivalent, if it is $n$-equivalent to the category of $n$-functors $G_n \to n-\Cat$, where $G_n$ is an $n$-groupoid. If $X$, $G_n$, $E_n$ correspond in this way, one calls $G_n$ \emph{the fundamental $n$-groupoid} of $X$, or of $E_n$, or says that $E_n$ is \emph{the category of local $(n-1)$-systems} on $X$, or on $G_n$, or that $X$ is \emph{the geometric realisation} of $G_n$ or of $E_n$. In analogy with the familiar case $n = 1$, it should be possible to interpret $G_n$ as the full sub-$n$-category of $\Hom_n(E_n, n-\Cat)$ formed by the $n$-functors $E_n \to n-\Cat$ satisfying certain exactness properties (one feels like saying: which commute with finite $\varprojlim$ and arbitrary $\varinjlim$); but this raises the disquieting vision of $n$-limits in $n$-categories. (N. B. The case $n = 2$ begins to become familiar to us\dots). It is prudent in all of this to suppose that $X$ is ``locally homotopically trivial'', which ensures the pro-simplicial set which Artin-Mazur associate to it (with the help of nerves of hyper-coverings) is essentially constant in the ordinary homotopy category - thus $X$ defines a homotopy type in the usual sense. This is surely \emph{not} the case for the étale topos of a scheme. In such case, the fundamental $n$-groupoid should be conceived as a \emph{pro-$n$-groupoid} (nothing surprising in that, in view of the familiar theory of $\pi_1$), and $E_n$ as an (ind)-$n$-category (the ind-structure will correspond to the exigencies of local triviality for a variable $n$-stack, relative to coverings more and more fine on $X$).

\starsbreak

\hangsection{The ``six operations'' for $P$-constructible sheaves (for a given equi-singular stratification $P$).}\label{sec:app16}%
I nevertheless understand your instinctive resistance to conceive this extreme stripping of a beautiful topos $X$, to the point of retaining only the meagre homotopy type. Even more, I am persuaded that going to the root of this instinctive resistance, one arrives at a generalisation and deepening of the notion of ``homotopy type'', and to bring new grist to the mill of the development of a good homotopical yoga. Here is what I have in mind.

Let us speak first of sheaves (of sets, or of modules, etc.) instead of stacks, for simplicity, and place ourselves in the étale topos of a scheme. The locally constant sheaves - modulo a supplementary condition of finiteness which is sufficiently anodyne - form the easiest of the \emph{constructible} sheaves, for the definition of which they serve as models. Supposing $X$ coherent (= quasi-coherent and quasi-separated), then the general constructible sheaves are those for which there exists a finite partition $X = \cup_{i \in I} X_i$ of $X$ into ``cells'' or ``strata'' $X_i$, each locally closed and constructible, such that the restriction of $F$ to every $X_i$ is locally constant (also a finiteness condition\dots). Thus the category of constructible sheaves on $X$ (which gives back the category of all sheaves on passing to a category of ind-objects\dots) may itself be thought of as an inductive limit of categories associated to finer and finer partitions on $X$. One can then, for such a fixed partition $P$, set out to study the category of sheaves (or complexes of sheaves, or stacks) which are ``$P$-constructible'' (or, more generally, which are ``locally constant'' on every $X_i$). These categories will not have truly satisfying structures unless they are stable for the usual operations - such as $\RHom$, or $\mathrm{R}j_* j^*$ where $j: X_i \to X$ is a ``cell'' of the partition, etc. In fact, if $X$ is excellent and one has resolution of singularities at ones disposal, one knows that the torsion constructible sheaves (under the proviso of being prime to the characteristic) are stable for all these operations - but not for a finite partition of $X$ fixed once and for all. To have such a finer stability, it is necessary to make some very strict hypotheses of \emph{``equi-singularity''} on the given stratifications of $X$, along the strata. I think nonetheless that a refinement of known techniques will show that $X$ admits arbitrarily fine stratifications having these properties of equi-singularity (and with the $X_i$ regular and connected, but this does not matter for our present purpose).

By way of example, suppose that there are just two strata, the closed one $X_0$, and $X_1 = X \textbackslash X_0$. According to Artin's devissage, giving oneself a sheaf $F$ on $X$ is equivalent to giving a sheaf $F_0 = i^*_0(F)$ on $X_0$, a sheaf $F_1(= i_1^* F)$ on $X_1$, and a homomorphism $F_0 \to i^*_0 i^*_{1^*}(F_1) = \varphi(F_1)$, where $i_0$, $i_1$ are the inclusions $X_0 \xrightarrow{i_0} X \xleftarrow{i_1} X_1$. In order that $F$ should be $P$-constructible, it is necessary and sufficient that $F_0$ and $F_1$ should be locally constant (plus some accessory finiteness conditions\dots), on $X_0$ and $X_1$ respectively. Then (by virtue of the hypothesis of equi-singularity) the same will be true of $\varphi(F_1)$, and the category of sheaves in which we are interested can be expressed entirely in terms of the category of locally constant sheaves on $X_0$ and $X_1$, i.e. of the mere homotopy type of $X_0$ and $X_1$, except that we must make explicit the nature of the left exact functor $\phi$. I think tht this should be possible, in the context of \emph{schemes} in which I am placed (technically rather sophisticated), on introducing an ``étale tubular neighbourhood'' of $X_0$ in $X_1$ (which is a very interesting topos, but not associated to a scheme). But this technical construction is only a paraphrase of an extraordinary simple topological intuition, which I will make explicit, supposing, to fix the ideas, that the base field is $\mathbf{C}$ and so one may work with locally compact spaces in the usual sense. The topological idea behind the hypothesis of equi-singularity  that there exists a \emph{tubular neighbourhood} $T$ of $X_0$ in $X$ retracting onto $X_0$ and such that the pair $(X_0, T)$ over $X_0$ should be a locally trivial bundle, i.e. that $T \textbackslash X_0$ is locally trivial over $X_0$. In fact if $\partial T$ is the ``boundary'' of $T$, which also should be a locally trivial bundle on $X_0$, then $T$ over $X$ is the conic bundle (= bundle where fibres are cones) ($\simeq (\partial T \times I) \amalg_{\partial T} X_0$ where $I = [0, 1]$, $\partial T \to \partial T \times I$ is defined by $x \mapsto (x, 1)$, and $\partial T \to X_0$ is the projection) then $\overset{\circ}{T} = T \textbackslash X_0 \simeq \partial T \times [0, 1 [$ is $X_0$-homotopic to $\partial T$. If $X_0$ and $X_1$ are non singular, then so also will be $\overset{\circ}{T}$ and $\partial T$, which are then topologically smooth fibrations on $X_0$. Moreover, putting $\widetilde{T}_1 = X_1 \textbackslash \overset{\circ}{T}$, the inclusion $\widetilde{X}_1 \to X_1$ is a homotopy equivalence, and $X$ can be recovered, up to homeomorphism, from the diagram of spaces
\[\begin{tikzcd}
  (\overset{\circ}{T} = T \textbackslash X_0 \simeq) \partial T \ar[r,swap,"\text{inclusion}","j"']\ar[d,swap,"\text{fibration}","p"'] & \widetilde{X}_1(\simeq X_1) \\
  X_0 &
\end{tikzcd}\]
as an amalgamated sum. In terms of this diagram of spaces, the above functor $\phi$ interprets immediately as
$$
\phi(F_1) \simeq p_* j^*(\widetilde{F}_1)
$$
where $F_1 \to \widetilde{F}_1$ is the restriction from $X_1$ to $\widetilde{X}_1$ (which is an equivalence of categories for the envisaged (locally constant) sheaves). Giving $F = (F_0, F_1, u:F_0 \to \phi F_1)$ can then also be made explicit as giving
$$
F_0, \widetilde{F}_1, \widetilde{u}: p^*(F_0) \to j^*(\widetilde{F}_1)
$$
where $F_0(\widetilde{F}_1)$ are locally constant sheaves on $X_0$ (respectively $X_1$). It is necessary to recall that here $p$ is a real fibration, and $j$ is an inclusion (in practice, for the case $X_0$, $X_1$ smooth, the inclusion of the boundary in a manifold with boundary).

If you prefer, one can also take the diagram which is less pretty (but a little more canonical)
\[\begin{tikzcd}
  \overset{\circ}{T} \ar[r,"j'"]\ar[d,swap,"p'"] & X_1 \\
  X_0 &
\end{tikzcd}\]
coming essentially to the same thing, as it is formed from spaces homotopic to the preceding one. One can even replace $X_0$ by $T$ ($X_0$ being a deformation retract of it) and write 
\[\begin{tikzcd}
  \overset{\circ}{T} \ar[r]\ar[d,swap,"p''"] & X_1 \\
  T &
\end{tikzcd}\]
where ``literally'' $p''$ is now an inclusion, but ``morally'', it is a \emph{fibration} with very pretty fibres (notably compact of finite dimension, and moreover non-singular varieties - this is much better than that which is given by the yoga of Cartan-Serre ``every continuous mapping is equivalent to a fibration''\dots). This last diagram however has the advantage of being amenable to a purely algebraic, direct construction, in the context of schemes, once one has developed the construction of étale tubular neighbourhoods\footnote{Tim Porter pointed out to me that work on étale tubular neighbourhoods was done by D.A. Cox: ``Algebraic tubular neighbourhoods I, II'', Math. Scand. 42 (1978) 211-228, 229-242. I've not seen yet this work, and can't say therefore whether it meets the rather precise expectations I have for a theory of tubular neighbourhoods, for the needs of a dévissage theory of stratified schemes (or, more generally, stratified topoi)}.

\hangsection{Relation with ``indexed homotopy types'' and with ``dévissage'' of stratified spaces. The ``fine homotopy type'' of a (tame) space or of a scheme.}\label{sec:app17}%
The point I wish to come to, is that the consideration f sheaves (or complexes thereof, or $n$-stacks\dots) which are $P$-constructible on an $X$, where $P$ is a given ``equi-singular'' stratification, reduces in our particular cases to the knowledge of a diagram of ordinary \emph{homotopy types} (or pro-types, if one comes back to the étale topology)
\[\begin{tikzcd}
  \Delta \ar[r,"j"]\ar[d,swap,"p"] & X_1 \\
  X_0 &
\end{tikzcd}\]
by taking local coefficients systems (or locally constant $n$-stacks) on the vertices $X_0$, $X_1$, which are related to each other by a homomorphism of compatibility of the type $p^*(F_0) \to j^*(F_1)$. It should be an amusing exercise (which I have not yet done) to verify and to make explicit how the ``six operations'' on sheaves (either on $X$, or on a subspace which is a union of strata of $X$) can be expressed in this dictionary, in the case, let us say, of non-singular strata (otherwise, there will be a difficulty with the dualising complexes, which one would prefer to have as objects in our category), and to reestablish the known formulae involving these operations. But it appears probable that, to carry out this transcription well, it would be necessary, rather than considering a diagram of type
\[\begin{tikzcd}
  \bullet \ar[r]\ar[d] & \bullet \\
  \bullet &
\end{tikzcd}\]
in the homotopical category formed from the category of semi-simplicial sets, to consider the category of diagrams of semi-simplicial sets, and to pass from these to the homotopical category of fractions\footnote{This is the typical game embodied in the ``derivator'' associated to the theory $\Hot$ of usual homotopy types (compare section 69).}

I have recently more or less made explicit, while thinking on the foundations of ``tame topology'', (i.e. where one eliminates from start all wild phenomena) how an equi-singular stratification, say with non singular strata, of a compact ``tame space'', gives rise canonically to a diagram of space which are manifolds with boundary, the arrows of the diagram being essentially locally trivial fibrations of manifolds with boundary on the others (with fibres which are compact manifolds with boundary), \emph{and} the inclusion of the boundary in these manifolds with boundary (in fact, one finds slightly more general inclusions, certain boundaries which appear being endowed with an ``elementary'' cellular decomposition, i.e. the closed strata are again manifolds with boundary which are glued together along common parts of the boundary; and it is also necessary to consider the inclusions of these pieces one in another\dots), and can be reconstituted from this diagram by gluing\footnote{Some more details on this program are outlines in \emph{``Esquisse d'un Programme''} (section 5), in \emph{Réflexions Mathématiques 1}}. In other words, one has a canonical devissage description of tame compact spaces $X$, eventually endowed with equi-singular stratifications with non-singular strata, in terms of finite diagrams of a precise nature made out o manifolds with boundary. When we are interested in sheaves (or complexes of sheaves, or $n$-stacks) which are $P$-constructible on $X$, where $P$ is such a fixed stratification, these may be described in terms of the envisaged diagram, of which only the ``homotopy type'' is to be retained. One foresees that the six operations on these sheaves can be translated in an ad hoc manner to this homotopical context. Finally, if instead of having only one compact tame space $X$, one has, let us say, a tame morphism $f: X \to Y$ of such objects, then by choosing equi-singular stratifications on $X$ and $Y$ adapted to $f$ (the strata of $X$ being in particular locally trivial fibrations on those of $Y$\dots), one should find a ``morphism'' from the diagram of manifolds with boundary expressing $X$ int that expressing $Y$ (with natural morphisms which essentially reduce to fibrations of compact manifolds with boundaries on other such objects) in such a way that the four operations $\Rf_*$, $\Rf_!$, $\mathrm{L}g^*$, $g^!$ between $P_{X^-}$ and $P_{Y^-}$ constructible sheaves on $X$ and $Y$ (or on locally closed sub-spaces $X'$, $Y'$ which are union of strata, such that $f$ induces $g: X' \to Y'$) can be expressed in terms of standard operations between the mere homotopy types. Finally, all these constructions, still partially hypothetical (there is work on the foundations to be done!) should be able to be paraphrased in the framework of excellent schemes, by making use of the machinery of étale tubular neighbourhoods. In one or other case, the ``fine homotopy type'' of a tame space, respectively of an excellent scheme, is defined by passage to the limit from ``$P$-homotopy type'' associated to finer and finer equi-singular stratifications $P$ (with non-singular strata).

This ``fine'' homotopy type would embody the knowledge, not only of sheaves or locally constant $n$-stacks, but (via a passage to the inductive limit) the knowledge of \emph{all of them}. And it would depend, in a suitable sense, functorially on $X$. In the case of a scheme of finite type on an algebraically closed field $k$ say, the strongest cohomological and homotopical \emph{finiteness theorem} would be expressed precisely in terms of a fine homotopy type, and would say that \emph{the ordinary homotopy types which are their constituents are essentially ``finite polyhedra''} - and even compact manifolds with boundary - or in more precise fashion, their profinite completions (in the sense of Artin-Mazur) prime to the characteristic $p$ of $k$ are those of such polyhedra. One sees clearly how to begin on such programme in characteristic 0, but one foresees supplementary amusement, or even mystery, in the case $p > 0$, for the varieties which, even birationally, resist being lifted to characteristic 0!

\starsbreak

\hangsection{How far do ``essentially constant'' \emph{abelian} coefficients determine a homotopy type?}\label{sec:app18}%
From these essentially geometric thoughts, I could not at this moment draw up a precise programme for developing adequate algebraic structures to express them. I restrict myself to several marginal remarks.

For a long time I have been intrigued by the idea of a ``linearisation'' of an (ordinary) homotopy type, i.e. questions of the type: if $X$ is a homotopy type, how much cohomological information of the type: cohomology of $X$ with variable coefficients $M$ (constant or twisted constant), multiplicative structure $\mathrm{H}^i(X, M) \times \mathrm{H}^j(X, N) \to \mathrm{H}^{i + j}(X, M \otimes N)$, then eventually other cohomology operations - is it necessary to have to reconstruct entirely the homotopy type? (say, in this preliminary pre-derived category approach, assuming given the fundamental group $\pi_1$, and therefore the category of constant twisted coefficients (= $\pi_1$-modules), the functors $\mathrm{H}^i(-, M)$ over these, together with the structure of cohomological functors relative to exact sequences, the structure of cup-product, etc. - related by certain formal properties?) Once one has at one's disposal the language of derived categories: the sub-category of the derived category of abelian complexes on $X$, formed from complexes the sheaves of cohomology of which are locally constant on $X$, with its triangulated structure and its multiplicative structure $\Lotimes$ (and eventually $\RHom$\dots) gives a more satisfying candidate for hoping to recover the homotopy type. I don't really know if this suffices the recovery indeed\footnote{I was informed by knowledgeable people soon later that the answer is well known to be negative, by working with ``rational homotopy types'' (the cohomology of which is made up with vector spaces over $\mathbf{Q}$). It is well known indeed that a 1-connected rational homotopy type is \emph{not} known from its rational cohomology ring alone, which contains already all the information I was contemplating. At last this is so if we assume that $\mathrm{H}^i(X)$ is of finite dimension over $\mathbf{Q}$ for all $i$. But is there a counterexample still when $X$ is a homotopy type ``of finite type''?}, but on the other hand I have no doubt that on pursuing ``linearisation'' to the end, that is to say by going to the \emph{non-abelian} framework, and working with the $(n + 1)$-category (without any supplementary structure on it!) of locally constant $n$-stacks of constructible sheaves on $X$, for all $n$, one manages to reconstruct the homotopy type via its fundamental $\infty$-groupoid, as explained in my previous letter and recalled in this one. (This signifies in particular that all the possible and imaginable cohomology operations are already included in the data furnished by such a system of $n$-categories\dots).

Similarly, the more elaborate homotopy type, which are related to certain finite diagrams, which one can associate to certain types of stratification $P$ of tame topological spaces $X$, let's say, should correspond in as perfect a fashion to the $(n + 1)$-category of $n$-stacks on $X$ which are locally constant on each of the strata of $P$ (say: which are subordinated to $P$). If the above description of homotopy types by the ``locally constant derived category'' was valid indeed, one would expect to recover here the mixed homotopy type from the corresponding sub-category of the derived category of abelian sheaves on $X$, provided by the complexes which have locally constant cohomology on each of the strata - with also the operations $\Lotimes$, $\RHom$, plus in case of need, the four operations $\mathrm{R}g_!$, $\mathrm{R}g_*$, $\mathrm{L}g^*$, $g^!$ for the induced $g: Z' \to Z''$ of the various locally closed unions of strata\dots The problem here is that we don't at present even know what is a triangulated category, not any more than what is its non-commutative version, described probably more simple and more fundamentally: a ``homotopical category'' with operations of taking ``fibres'' and ``cofibres''\footnote{This ``problem'' is met with by the notion of a ``derivator'', which ``was in the air'' already by the late sixties, but was never developed (instead even derived categories became tabu in the seventies\dots).}.

It is surely time that I finish this ``lettre-fleuve'', which is becoming more and more vague. Just one question: what is this marvellous formula of Bloch-Quillen to which you allude, of which I have never heard, and which makes my mouth water?

\bigskip

Very cordially yours, 

\begin{flushright}
Alexandre Grothendieck
\end{flushright}



\setcounter{chapter}{1}

\chapter{Test categories and test functors}
\label{ch:II}

\centerline{\itshape Reflections on homotopical algebra}\pspage{1}

\bigbreak

\presectionfill\ondate{27.2.}83.\par

The following notes are the continuation of the reflection started in
my letter to Daniel Quillen written previous week (19.2 -- 23.2),
which I will cite by \hyperref[ch:I]{(L)} (``letter''). I begin with
some corrections and comments to this letter.

\setcounter{section}{13}
\hangsection[The unnoticed failure \dots]{The unnoticed failure
  \texorpdfstring{\textup(}{(}of the present foundations of topology
  for expressing topological intuition\texorpdfstring{\textup)}{)}.}%
\label{sec:14}%
Homotopical algebra can be viewed as being concerned mainly with the
study of spaces of continuous maps between spaces and the algebraic
analogons of these, with a special emphasis on homotopies between such
homotopies between homotopies etc. The kind of restrictive properties
imposed on the maps under consideration is exemplified by the typical
example when demanding that the maps should be extensions, or
liftings, of a given map. Homotopical algebra is not directly suited
for the study of spaces of homeomorphisms, spaces of immersions, of
embeddings, of fibrations, etc.\ -- and it seems that the study of
such spaces has not really yet taken off the ground. Maybe the main
obstacle here lies in the wildness phenomena, which however, one
feels, makes a wholly artificial difficulty, stemming from the
particular way by which topological intuition has been mathematically
formalized, in terms of the basic notion of topological spaces and
continuous maps between them. This transcription, while adequate for
the homotopical point of view, and partly adequate too for the use of
analysts, is rather coarsely inadequate, it seems to me, in most other
geometrico-topological contexts, and particularly so when it comes to
studying spaces of homeomorphisms, immersions, etc.\ (in all those
questions when ``isotopy'' is replacing the rather coarse homotopy
relation), as well as for a study of stratified structures, when it
becomes indispensable to give intrinsic and precise meaning to such
notions as tubular neighborhoods, etc.  For a structure theory of
stratifications, it turns out (somewhat surprisingly maybe) that even
the somewhat cumbersome context of topoi and pretopoi is better suited
than topological spaces, and moreover directly applicable to
unconventional contexts such as \'etale topology of schemes, where the
conventional transcription of topological intuition in terms of
topological spaces is quite evidently breaking down. To emphasize one
point I was making in (\hyperref[ch:I]{L}, p.~\ref{p:L.1}), it
seems to me that this breakdown is almost as evident in isotopy
questions or for the needs of a structure theory of stratified
``spaces'' (whatever we mean by ``spaces'' \ldots). It is a matter of
amazement to me that this breakdown has not been clearly noticed, and
still less overcome by working out a more suitable transcription of
topological objects and\pspage{2} topological intuition, by the people
primarily concerned, namely the topologists. The need of eliminating
wildness phenomena has of course been felt repeatedly, and (by lack of
anything better maybe, or rather by lack of any attempt of a
systematic reflection on what \emph{was} needed) it was supposed to be
met by the notion of piecewise linear structures. This however was
falling from one extreme into another -- from a structure species with
vastly too many maps between ``spaces'', like a coat vastly too wide
and floating around in a million wild wrinkles, to one with so few
(not even a quadratic map from $\bR$ to $\bR$ is allowed!) that it
feels like too narrow a coat, bursting apart on all edges and
ends. The main defect here, technically speaking, seems to me the fact
that numerical piece-wise linear functions are not stable under
multiplication, and as a geometric consequence of this, that when
contracting a compact p.l.\ subspace of a (compact, say) piecewise
linear space into one point, we do not get on the quotient space a
natural p.l.\ structure. This alone should have sufficed, one might
think, to eliminate the piecewise linear structure species as a
reasonable candidate for ``doing topology'' without wildness
impediments -- but strangely enough, it seems to be hanging around
till this very day!

\hangsection{Overall review of standard descriptions of homotopy
  types.}\label{sec:15}%
But my aim here is not to give an outline of foundations on ``tame
topology'', but rather to fill some foundational gaps in homotopy
theory, more specifically in homotopical algebra. The relative success
in the homotopical approach to topology is probably closely tied to
the well known fact (Brouwer's starting point as a matter of fact,
when he introduced the systematics of triangulations) that every
continuous map between triangulated spaces can be approximated by
simplicial maps. This gave rise, rather naturally, to the hope
expressed in the ``Hauptvermutung'' -- that two homeomorphic
triangulated spaces admit isomorphic subdivisions, a hope that finally
proved a delusion. With a distance of two or three generations, I
would comment on this by saying that this negative result was the one
to be expected, once it has become clear that neither of the two
structure species one was comparing, namely topological spaces and
triangulated spaces, was adequate for expressing what one is really
after -- namely an accurate mathematical transcription, in terms of
``spaces'' of some kind or other, of some vast and deep and misty and
ever transforming mass of intuitions in our psyche, which we are
referring to as ``topological'' intuition\footnote{For another way of viewing the role of the Hauptvermutung, namely ``rehabilitating'' it as being essentially valid, see note $(^6)$ at the end of Esquisse d'un Programme (in RM 1).}. There \emph{is} something
positive though, definitely, which can be viewed as an extremely
weakened version of the Hauptvermutung, namely the fact that
topological spaces on the one hand, and semi-simplicial sets on the
other, give rise, by a suitable ``localization''\pspage{3} process
(formally analogous to the passage from categories of chain complexes
to the corresponding derived categories), to
essentially ``the same'' (up to equivalence)
``homotopy category''. One way of describing it is via topological
spaces which are not ``too wild'' as objects (the CW-spaces),
morphisms being homotopy classes of continuous maps. The other is via
semi-simplicial sets, taking for instance Kan complexes as objects,
and again homotopy classes of ``maps'' as morphisms. The first
description is the one most adapted to direct topological intuition,
as long as least as no more adequate notion than ``topological
spaces'' is at hand. The second has the advantage of being a purely
algebraic description, with rather amazing conceptual simplicity
moreover. In terms of the two basic sets of algebraic invariants of a
space which has turned up so far, namely cohomology (or homology) on
the one hand, and homotopy groups on the other, it can be said that
the description via topological spaces is adequate for direct
description of neither cohomology nor homotopy groups, whereas the
description via semi-simplicial sets is fairly adequate for
description of cohomology groups (taking simply the abelianization of
the semi-simplicial set, which turns out to be a chain complex, and
taking its homology and cohomology groups). The same can be said for
the alternative algebraic description of homotopy types, using cubical
complexes instead of semi-simplicial ones, which were introduced by
Serre as they were better suited, it seems\footnote{As a matter of fact, it is pretty obvious today how to get the Leray spectral sequence directly in the semisimplicial set-up, without using cubes.}, for the study of
fibrations and of the homology and cohomology spectral sequences
relative to these. One somewhat surprising common feature of those two
standard algebraic descriptions of homotopy types, is that neither is
any better adapted for a direct description of homotopy groups than
the objects we started with, namely topological spaces. This is all
the more remarkable as it is the homotopy groups really, rather than
the cohomology groups, which are commonly viewed as \emph{the} basic
invariants in the homotopy point of view, sufficient, e.g., for test
whether a given map is a ``weak equivalence'', namely gives rise to an
isomorphism in the homotopy category. 

It is here of course that the
point of view of ``stacks'' (``champs'' in French) of
\hyperref[ch:I]{(L)} (previously called ``\oo-groupoids'' in the
beginning of the reflections of \hyperref[ch:I]{(L)}) sets in. These
presumably give rise to a ``category of models'' and from there, to the usual homotopy category, in much the same way as topological spaces or simplicial (or cubical)
complexes, thus yielding a third main description of
homotopy types, and corresponding wealth of algebraic-geometric
intuitions. Moreover, stacks are ideally suited for expressing the
homotopy groups, in an even more direct way than simplicial
complexes\pspage{4} allow description of homology and cohomology
groups. As a matter of fact, the description is formally analogous,
and nearly identical, to the description of the homology groups of a
chain complex -- and it would seem therefore that that stacks (more
specifically, Gr-stacks) are in a sense the \emph{closest possible
non-commutative generalization of chain complexes}, the homology groups
of the chain complex becoming the homotopy groups of the
``non-commutative chain complex'' or stack.

It is well understood, since Dold-Puppe\footnote{Ronnie Brown pointed out to me that ``Dold-Puppe'' is an accurate name here, as this theorem has been found independently by Dold and Kan. I have corrected this inaccuracy in later occurrences in my typescript, writing ``Dold-Kan''instead.}, that chain complexes form a
category equivalent to the category of abelian group objects in the
category of semi-simplicial sets, or equivalently, to the category of
semi-simplicial abelian groups. By this equivalence, the homology
groups of the chain complex are identified with the \emph{homotopy}
groups of the corresponding homotopy type. As for the homology and
cohomology groups of this homotopy type, their description in terms of
the chain complex we started with is kind of delicate (I forgot all
about it I am afraid!). A fortiori, when a homotopy type is described
in terms of a stack, i.e., a ``non-commutative chain complex'', there
is no immediate way for describing its homology or cohomology groups
in terms of this
\clearpage
structure\footnote{This statement is mistaken - it appeared soon after that the cohomology of a stack can be expressed in terms of the ``primitive structure'' (just face operations) much in the same way as for simplicial or cubical complexes - it is even simpler, because in each dimension there are only two face operators (whose formal difference will be the boundary operator). It is at this point that the idea of \emph{``hemispherical complexes''} (or $h$ ``sets'') comes in very natural, as the kind of objects that remain, when a ``stack'' is stripped from all it's structure except the primitive structure, namely a family of objects $(F_i)$ together with the ``face operations'' $(s_i, t_i)$ $(i \geq 1)$, namely the ``source'' and ``target'' maps - plus, preferably, also the degeneracies $k_i$ $(i \geq 1)$, one in each dimension. These may be viewed too as the objects in $\hat{\Globe}$, where $\Globe$ is the category of \emph{``standard hemispheres''}, very similar to the categories $\Simplex$ and $\square$ of standard simplices resp. cubes, and even simpler than these (being based on bigons, rather than on triangles or squares). In the sequel, we often refer to this category $\Globe$, and to the ``trinity'' of the three most beautiful ``test categories'' (in fact, they are even ``contractors'') $\Globe$, $\Simplex$, $\square$, without having introduced formally $\Globe$ anywhere in the notes. This category $\Globe$ has been first introduced in the thesis of \emph{David W. Jones} (a student of Ronnie Brown) \emph{``Poly $T$-complexes''}, and he calls the corresponding notion of complexes ``\emph{globular} complexes'' (instead of ``hemispherical''). That these do modelize homotopy types, more accurately that $\Globe$ is a (strict) test category, will follow from the general criterion of section 31, much in the same way as for $\Simplex$ and $\square$ (see section 34 for these). It is clear that for the study of stacks, hemispherical complexes will be better suited than semisimplicial or cubical ones. It shouldn't be hard anyhow to describe functors going from one type of complexes to the other, in such a way that the homotopy types modelled by two complexes which correspond should be canonically isomorphic. In terms of the theory of ``test-functors'' which is to follow, this amounts essentially to the question of finding an explicit test functor $\Simplex \to \hat{\Globe}$ (and similarly in the five other cases). This reminds me of the problem of giving handy existence and unicity theorems for test functors, dwelled upon somewhat in section 90, but by no means settle yet\dots}. 

What probably should be done, is to define first a
\emph{nerve} functor from stacks to semi-simplicial sets (generalizing
the familiar nerve functor defined on the category of categories), and
define homology invariants of a stack via those of the associated
semi-simplicial set (directly suited for calculating these).

\hangsection[Stacks over topoi as unifying concept for homotopical and
\dots]{Stacks over topoi as unifying concept for homotopical and
  cohomological algebra.}\label{sec:16}%
These reflections on the proper place of the notion of a stack which in
standard homotopy algebra are largely a posteriori -- the clues they
give are surely not so strong as to give an imperative feeling for the
need of developing this new approach to the homotopy category. Rather,
the imperative feeling comes from the intuitions tied up with
non-commutative cohomological algebra over topological spaces, and
more generally over topoi, in the spirit of Giraud's thesis,\scrcomment{\cite{Giraud1971}} where a
suitable formalism for non-commutative $K^i$'s for $i=0$, $1$ or $2$
is developed. He develops in extenso the notion of stacks, we should
rather say now $1$-stacks, over a topos, constantly alluding (and for
very understandable reasons!) to the notion of a $2$-stack, appearing
closely on the heels of the $1$-stacks. Keeping in mind that
$0$-stacks are just ordinary sheaves of sets, on the space or the
topos considered, the hierarchy of increasingly higher and more sophisticated notions of
$0$-stacks, $1$-stacks, $2$-stacks, etc., which will have to be
developed over an arbitrary topos, just parallels the hierarchy of
corresponding notions over the one-point\pspage{5} space, namely
sets (= $0$-stacks), categories (or $1$-stacks), $2$-categories, etc.
Among these structures, those generalizing groupoids among categories,
namely Gr-stacks of various orders $n$, play a significant role,
especially for the description of homotopy types, but equally for a
non-commutative ``geometric'' interpretation of the cohomology groups
$\mathrm H^i(X,F)$ of arbitrary dimension (or ``order''), of a topos $X$ with
coefficients in an abelian sheaf $F$. The reflections in \hyperref[ch:I]{(L)} therefore
were directly aimed at getting a grasp on a definition of such
Gr-stacks, and whereas it seems to me to have come to a concrete
starting point for such a definition, a similar reflection for
defining just stacks rather than Gr-stacks is still lacking. This is
one among the manifold things I have in mind while sitting down on the
present reflections.

Thus \emph{$n$-stacks, relativized over a topos to ``$n$-stacks over $X$'',
are viewed primarily as the natural ``coefficients'' in order to do
(co)homological algebra of dimension $\le n$ over $X$}. The
``integration'' of such coefficients, in much the same spirit as
taking objects $\RGamma_*$ (with $\RGamma$ the derived functor of the
sections functor $\Gamma$) for complexes of abelian sheaves $F_j$ on
$X$, is here merely the trivial operation of taking sections, namely
the ``value'' of the $n$-stack on the final object of $X$ (or of a
representative site of $X$, if $X$ is described in terms of a
site). The result of integration is again an $n$-stack, whose homotopy
groups (with a dimension shift of $n$) should be viewed as the
cohomology invariants $\mathrm H^i(X,F_*)$, where $F_*$ now stands for the
$n$-stack rather than for a complex of abelian sheaves\footnote{Compare with my letter to Breen (I App.7)}. In my letters
(two or three) to Larry Breen in 1975, I develop some heuristics along
this point of view, with constant reference to various geometric
situations (mainly from algebraic geometry), providing the
motivations. The one motivation maybe which was the strongest, was the
realization that the classical Lefschetz theorem about comparison of
homology and homotopy invariants of a projective variety, and a
hyperplane section -- once it was reformulated suitably so as to get
rid of non-singularity assumptions, replaced by suitable assumptions
on cohomological ``depth'' -- could be viewed as comparison statements
of ``cohomology'' with coefficients in more or less arbitrary
stacks\footnote{Compare with (I App.2)}. This is carries through completely, within the then existing
conceptual framework restricted to $1$-stacks, in the thesis\scrcomment{\textcite{Raynaud1975}} of Mme
Raynaud, a beautiful piece of work. There seems to me to be
overwhelming evidence that her results (maybe her method of proof
too?) should generalize in the context of non-commutative
cohomological algebra of arbitrary dimension, with a suitable property
of ind-finiteness as the unique\pspage{6} restriction on the
coefficient stacks under consideration.

Technically speaking, \oo-stacks are the common denominator of
$n$-stacks for arbitrary $n$, in much the same way as $n$-stacks
appear both as the next-step generalization of $(n-1)$-stacks (the
former forming a category which admits the category of $(n-1)$-stacks
as a full subcategory), and as the most natural ``higher order
structure'' appearing on the category of all $(n-1)$-stacks (and on
various analogous categories whose objects are $(n-1)$-stacks subject
to some restrictions or endowed with some extra structure). I'll have
to make this explicit in due course. For the time being, when speaking
of ``stacks'' or ``Gr-stacks'', it will be understood (unless
otherwise specified) that we are dealing with the infinite order
notions, which encompass the finite ones.

Working out a theory of stacks over topoi, as the natural foundation
of non-commutative cohomological algebra, would amount among others to
write Giraud's book within this considerably wider framework. Of
course, this mere prospect wouldn't be particularly exciting by
itself, if it did not appear as something more than grinding through
an unending exercise of rephrasing and reproving known things,
replacing everywhere $n=0$, $1$ or $2$ by arbitrary $n$. I am
convinced however that there is a lot more to it -- namely the
fascination of gradually discovering and naming and getting acquainted
with presently still unknown, unnamed, mysterious structures. As is
the case so often when making a big step backwards for gaining new
perspective, there is not merely a quantitative change (from $n\le2$
to arbitrary $n$ say), but a qualitative change in scope and depth of
vision. One such step was already taken I feel by Daniel Quillen and
others, when realizing that homotopy constructions make sense not only
in the usual homotopy category, or in one or the other categories of
models which give rise to it, but in more or less arbitrary
categories, by working with semi-simplicial objects in these. The step
I am proposing is of a somewhat different type. \emph{The notion of a
stack here appears as the unifying concept for a synthesis of
homotopical algebra and non-commutative cohomological algebra}. This (rather than merely furnishing us
with still another description of homotopy types, more convenient for
expression of the homotopy groups) seems to me the real ``raison
d'\^etre'' of the notion of a stack, and the main motivation for
pushing ahead a theory of stacks.\pspage{7}

\DontChangeNextSectionNumber
\renewcommand{\thesection}{\arabic{section}bis}
\hangsection[Categories as models for homotopy types. First glimpse
\dots]{Categories as models for homotopy types. First glimpse
  upon an ``impressive bunch'' \texorpdfstring{\textup(}{(}of
  modelizers\texorpdfstring{\textup)}{)}.}\label{sec:16bis}%
One last comment before diving into more technical
matters. Without even climbing up the ladder of increasing
sophistication, leading up to $\infty$-stacks, there is on the very first
step, namely with just usual categories, the possibility of describing
arbitrary homotopy types. Namely, there are two natural, well-known ways to associate to a
category $C$ (I mean here a ``small'' category, belonging to the given
universe we are working in) some kind of topological object, and hence
a homotopy type. One is by associating to $C$ the topos $\Chat$
or $\Top(C)$, the sheaves on which are just the presheaves on $C$, namely arbitrary contravariant functors from $C$ to
\Sets. The other is through the nerve functor\footnote{The definition of the nerve functor, realizing a full embedding of $\Cat$ into $\Simplexhat$, is given for the first time (I believe) in a Bourbaki talk of mine, which at that time wasn't concerned at all with the topological interpretation of small categories, but with some rather formal prerequisites for the operation of ``passage to quotient'' of an object by a ``pre-equivalence relation'', in any category (with a view of applying it in the category of schemes\dots)}\scrcomment{See FGA $n^\circ$212 in \textcite{FGA}}, associating to $C$
a semi-simplicial set -- and hence, if this suits us better, a
topological space, by taking the geometrical realization. By a
construction of Verdier, any topos and therefore $\mathrm{Top}(C)$
gives rise canonically to a pro-object in the category of semi-simplicial
sets, and hence by ``localization'' to a pro-object in the homotopy
category (namely a ``pro-homotopy type'' in the terminology of
Artin-Mazur). In the same way, the nerve $\cst N(C)$ gives rise to a
homotopy type -- and of course the two (pro)-homotopy types thus obtained coincide essentially, and may be called \emph{the}
homotopy type of $C$. When $C$ is a groupoid, we get merely a
$1$-truncated homotopy type, namely with homotopy groups $\pi_n$ which
vanish for $n\ge2$, or equivalently, with connected components
(corresponding of course to connected components of $C$) $K(\pi,1)$
spaces. This had led me at one moment in the late sixties to
hastily surmise that even for arbitrary $C$, we got merely such sums
of $K(\pi,1)$ spaces (namely, that the homotopy type of $C$ does not
change when replacing $C$ by the universal enveloping groupoid,
deduced from $C$ by making formally invertible all its arrows).
As Quillen pointed out to me, this is definitely not so -- indeed,
using categories, we get (up to isomorphism) arbitrary homotopy
types. This is achieved, I guess, using the left adjoint functor $\cst
N'$ from the inclusion functor
\[ \begin{tikzcd} \Cat \ar[r, hook, "\cst \bN"] & \Sssets\end{tikzcd}\]
which is fully faithful (the adjoint functor being therefore a
localization functor), and showing that for a semisimplicial set $K$,
the natural map
\[ K \to \cst N\,\cst N' (K) \]
is a weak equivalence; or what amounts to the same, that the set of
arrows in \Sssets{} by which we localize in order to get \Cat\ (namely
those transformed into invertible arrows by $\cst N'$) is made up with
weak equivalences only\footnote{This is well known to be false, and the mistake is corrected later (section 24). The next statement, about $\Cat$ modelizing the homotopy category, is true however and has been fairly well known (due to Fritsch and Latch, Ronnie Brown tells me). In these notes, it is going to come out in section 34, under the form ``$\Simplex$ is a test category'', using the general criterion for test categories obtained in section 31.}. This would imply
that we may reconstruct the usual homotopy category, up to
equivalence, in terms of \Cat, by\pspage{8} just localizing \Cat{}
with respect to weak isomorphism, namely functors $C \to C'$ inducing
an isomorphism between the corresponding homotopy types. Pushing a
little further in this direction, one may conjecture that \Cat{} is a
``closed model category'', whose weak equivalences are the functors
just specified, and whose cofibrations are just functors which are
injective on objects and injective on arrows\footnote{As seen a little later, this tentative description of a suitable class of cofibrations is inadequate (see beginning of section 30). An adequate description of a closed model structure in $\Cat$ was found by \emph{Thomason} (see comments on this in section 87).}.

The fact that categories are suitable objects for defining arbitrary homotopy types is quite remarkable, and one moreover it seems which has
not found its way still into the minds of topologists or homotopists,
with only few exceptions I guess. These objects are extremely simply
and familiar to most mathematicians; what is somewhat more
sophisticated is the process of localization towards homotopy types,
or equivalently, the explicit description of weak equivalences, within
the framework of usual category theory. This would amount more or less
to the same as describing the homotopy groups of a category, which
does not seem any simpler than the same task for its nerve. As for the
cohomology invariants, which can be interpreted as the left derived
functors of the $\varprojlim_C$ functor, or rather its values on
particular presheaves (for instance constant presheaves), they are of
course known to be significant, independently of any particular
topological interpretation, but they are not expressible in direct
terms. (The most common computation for these is again via the nerve
of $C$.)

This situation suggests that for any natural integer $n\ge1$, the
category of $n$-stacks can be used as a category of models for the
usual homotopy category, in particular any $n$-stack gives rise to a
homotopy type, and up to equivalence we should get any homotopy type
in this way (for instance, through the $n$-category canonically
associated to any $1$-category giving rise to this homotopy type). The
homotopy types coming from $n$-Gr-stacks, however, should be merely
the $n$-truncated ones, namely those whose homotopy groups in
dimension $>n$ are zero. Moreover, $n$-Gr-stacks appear as (in some respects) the most
adequate algebraic structures for expressing $n$-truncated homotopy
types, the latter being deduced from the former, presumably, by the
same process of localization by weak equivalences. Moreover, in the
context of $n$-Gr-stacks, the notions of homotopy groups and of weak
equivalences are described in a particularly obvious way. Thus,
passing to the limit case $n=\oo$, it is \oo-Gr-stacks? rather than general \oo-stacks which
appear as the neatest algebraic description of homotopy types.

These\pspage{9} reflections suggest that there should be a rather
impressive bunch of algebraic structures, each giving rise to a model
category for the usual homotopy category, or in any case yielding this
category by localization with respect to a suitable notion of ``weak
equivalences''. The ``bunch'' is all the more impressive, if we
remember that the notion of stack (dropping now the qualification $n$,
namely assuming $n=\oo$) is not really a uniquely defined one, but
depends on the choice of a ``\emph{coherator}'', namely (mainly) a category
$C$ satisfying certain requirements, which can be met in a vast
variety of ways, presumably. The construction of coherators is
achieved in terms of universal algebra, which seems here the
indispensable Ariadne's thread not to get lost in overwhelming
messiness. The natural question which arises here (and which do not
feel though like pursuing\footnote{Despite this ``feeling'', it turns out that the whole first volume of ``Pursuing Stacks'' is spent on pursuing this general question of understanding homotopy models in general, while the stacks are being forgotten for a thousands pages or so\dots}) is to give in terms of universal algebra
some kind of characterization, among all algebraic structures, of
those which give rise in some specified way (including the known
cases) to a category of models say for usual homotopy theory.

\renewcommand{\thesection}{\arabic{section}} \hangsection[The
Artin-Mazur cohomological criterion for weak \dots]{The Artin-Mazur
  cohomological criterion for weak equivalence.}\label{sec:17}%
When referring (section \ref{sec:15}) to the notion of a stack as a unifying
concept for homotopical algebra and non-commutative cohomological
algebra, I forgot to mention one significant observation of
Artin-Mazur along those lines (messy unification), namely that (for
ordinary homotopy types) weak equivalences (namely maps inducing
isomorphisms for all homotopy groups) can be characterized as being
those which induce isomorphisms on cohomology groups of the spaces
considered not only for constant coefficients, but also for arbitrary
twisted coefficients on the target space, including also the
non-commutative $\mathrm H^0$ and $\mathrm H^1$ for twisted
(non-commutative) group coefficients. This is indeed the basic
technical result enabling them, from known results on \'etale
cohomology of schemes (including non-commutative $\mathrm H^1$'s) to
deduce corresponding information on homotopy types. Maybe however that
the observation has acted rather as a dissuasion for developing higher
non-commutative cohomological algebra, as it seemingly says that the
non-commutative $\mathrm H^1$, plus the commutative $\mathrm H^i$'s,
was all that was needed to recover stringent information about
homotopy types. In other words, there wasn't too little in Giraud's
book, but rather, too much!

\bigbreak

\presectionfill\ondate{28.2.}\pspage{10}\par

\hangsection[Corrections and contents to letter. Bénabou's lonely
\dots]{Corrections and contents to letter. Bénabou's lonely
  approach.}\label{sec:18}%
I still have to correct a number of ``\'etourderies'' of
\hyperref[ch:I]{(L)}. The most persistent one, ever since page~\ref{p:L.8} of
that letter, is about fiber products in the coherator $\bC_\oo$, or,
equivalently, amalgamated sum in the dual category $\bB_\oo$. The
``correction'' I added in the last PS (section \ref{sec:13}) is
still incorrect, namely it is not true even in the subcategory $\bB_0$
of $\bB$ that arbitrary amalgamated sums exist. I was mislead by the
interpretation of elements of $\bB_0$ in terms of (contractible)
spaces, obtained inductively by gluing together discs $D_n$
($n\in\bN$) via subdiscs, corresponding to the cellular subdivisions
of the discs $D_n$ considered p.~\ref{p:L.7}. I was implicitly
thinking of amalgamated sums of the type
\[ K \amalg_L M,\]
where $L\to K$ and $M\to K$ are \emph{monomorphisms}, corresponding to
the geometric vision of embeddings -- in which case the usual
amalgamated sum in the category of topological spaces is indeed
contractible, which was enough to make me happy. But I overlooked the
existence of morphisms in $\bB_0$ which are visibly no
monomorphisms, such as $K \amalg_L K \to K$ the codiagonal map, when
$L \to K$ is a strict inclusion of discs. In any case, I will have to
come back upon the description of the categories $\bB_0$ and its
dual $\bC_0$ and upon the definition of coherators, after the
heuristic introduction \hyperref[ch:I]{(L)}. It will be time then too to correct the
mistaken description of $\bC_{n+1}$ in terms of $\bC_n$, which I
propose on p.~\ref{p:L.11prime}, yielding probably much too big a subcategory of
$\bC_\oo$ -- the correct definition should make explicit use of the
total set $A$ of ``new'' arrows, by which $\bC_\oo$ is described in
terms of universal algebra via $\bC_0$. Another \'etourdie the day
before, (section 7), is the statement that in a standard amalgamated sum,
any intersection of two maximal subcells is a subcell -- which is seen
to be false in the standard example.\footnote{Drawing of two
  hemispheres intersecting in a sphere.}

Of lesser import is the misstatement (in section 2),
stating that the associativity relation for the operation $*_\ell$
should by replaced by a homotopy arrow
$(\lambda * \mu) * \nu \to \lambda * (\mu * \nu)$. This is OK for
$\ell=1$ (primary composition), but already for $\ell=2$ does not make
sense as stated, because the sources and the targets of the two compositions to be compared are not equal, but only homotopic. Here the
statement should be replaced by one making sense, with a homotopy
``making commutative'' a certain square, and accordingly for cubes,
etc.\ for higher order compositions $\lambda *_\ell \mu$, to give
reasonable meaning to associativity. Anyhow, such painstaking
explicitations of particular coherence properties (rather, coherence
homotopies) is kind of ruled out by the sweeping axiomatic description
of the kind of structure species we want for a ``stack'', at least, I
guess, in\pspage{11} a large part of the development of the theory of
Gr-stacks. A systematic study of particular sets of homotopies is
closely connected of course to an investigation into irredundancy
conditions which can be figured out for a coherator $\bC$. This is
indeed an interesting topic, but I decided not to get involved in
this, unless I am really forced to!

The basic notion which has been peeling out in the reflection \hyperref[ch:I]{(L)} is
of course the notion of a coherator. Concerning terminology, it
occurred to me that the dual category $\bB_\oo$ to $\bC_\oo$ is
more suggestive in some cases, for instance because of the topological
interpretation attached to its objects, and (more technically) because
of formal analogy of the role of this category, for developing
homotopical algebra, with the category of the standard (ordered)
simplices\footnote{Finally, it turns out more suitable still to replace this category by the smaller one of standard hemispheres $\Globe$}. Both mutually dual objects $\bB_\oo$, $\bC_\oo$ seem
to me to merit a name, I suggest to call them respectively \emph{left}
and \emph{right coherators}, or simply \emph{coherators} of course
when for a while it is understood on which side of the mirror we are
playing the game.

One last comment still before taking off for a heuristic voyage of
discovery of stacks! I just had a glance at Bénabou's\scrcomment{\textcite{Benabou1967}} exposé in
1967 of what he calls ``bicategories'' (Springer Lecture Notes
n\textsuperscript{\b o}~47, p.~1--77). These are none else, it appears,
than \emph{non-associative} $2$-categories, namely $2$-stacks in the
terminology I am proposing (but not $2$-Gr-stacks -- namely it is a
particular case of a general notion of \oo-stack which has still to be
developed). The most interesting feature of this expos\'e, it seems to
me, is the systematic reference to topological intuition, notably of
the structure of various diagrams. His terminology, referring to
elements of $F_0$, $F_1$, $F_2$ respectively as $0$-cells, $1$-cells
and $2$-cells, is quite suggestive of an idea of topological
realization of a $2$-stack -- it is not clear from this expos\'e
whether B\'enabou has worked out this idea, nor whether he has made a
connection with Quillen's ideas on axiomatics of homotopy theory,
which appeared the same year in the same series\footnote{From a letter of Bénabou which I got later, I understand that he got discouraged from pursuing his topological approach to $n$-categories, because of the complete lack of interest of anybody else in this approach. I was glad to hear from him last year that my present work has stimulated him into starting a seminar on this topic of stacks, using my notes as a leading thread.}. In any case, in the
last section of his expos\'e, he deals with his bicategories formally
as with topological spaces, much in the same spirit as the one I was
contemplating since around 1975, and which is now motivating the
present notes. While there is no mention of B\'enabou's ideas in my
letters to Larry Breen, it is quite possible that on the unconscious
level, the little I had heard of his approach on one or two casual
occasions, had entered into reaction with my own intuitions, coming
mainly from geometry and cohomological algebra, and finally resulted
in the program outlined in those letters.

\hangsection{Beginning of a provisional itinerary
  \texorpdfstring{\textup(}{(}through
  stacks\texorpdfstring{\textup)}{)}.}\label{sec:19}%
I\pspage{12} would like now to write down a provision itinerary of
the voyage ahead -- namely to make a list of those main features of a
theory of stacks which are in my mind these days. I will write them
down in the order in which they occur to me -- which will be no
obligation upon me to follow this order, when coming back on those
features separately to elaborate somewhat on them. This I expect to
do, mainly as a way to check whether the main notions and intuitions
introduced are sound indeed, and otherwise, to see how to correct
them.

\begin{enumerate}[label=\arabic*\textsuperscript{\b{o}})]
\item Definition of the categories $\bB_0$ and its dual $\bC_0$, and
  formal definition of of left and right coherators. This definition will
  still be a provisional one, and will presumably have to be adjusted
  somewhat to allow for the various structures we are looking for in
  the corresponding category of Gr-stacks.
\item Relation between the category of Gr-stacks and the category of
  topological spaces, via two adjoint functors, the ``topological
  realization functor'' $F_* \mapsto \abs{F_*}$, and the ``singular
  stack functor'' $X \mapsto \bF_*(X)$. The situation should be
  formally analogous to the corresponding situation for
  semi-simplicial sets versus topological spaces, the role of the
  category $\cst S_*$ of standard ordered simplices being taken by the
  left coherator $\bB$ we are working with. The main technical
  difference here is that the category of Gr-stacks is not just the
  category of presheaves on $\bB$, but the full subcategory defined
  by the requirement that the presheaves considered should transform
  ``standard'' amalgamated sums into fibered products. As a matter of
  fact, the topological realization functor $F_* \mapsto \abs{F_*}$
  could be defined in a standard way on the whole of $\bB\uphat$, in
  terms of any functor
  \begin{equation}
    \label{eq:19.star}
    \bB \mapsto \Spaces \tag{*}
  \end{equation}
  (by the requirement that the extension of this functor to $B\uphat$
  commutes with arbitrary direct limits). A second difference with
  the simplicial situation lies in the fact that the only really
  compelling choice for the functor \eqref{eq:19.star} is it's
  restriction to $\bB_0$, in terms of the cells $\cst D_n$ and the
  standard ``half-hemisphere maps'' between these (\hyperref[ch:I]{L},
  section 6). The
  extension of this functor to $\bB$ is always possible, due to the
  inductive construction of $\bB$ and to the interpretation of
  elements of $\bB_0$ as \emph{contractible} spaces, via the
  functor \eqref{eq:19.star}$_0$; but it depends on a bunch of
  arbitrary choices. To give precise meaning to the intuition that
  these choices don't really make a difference, and that the choice of
  coherator we are starting with doesn't make much of a difference
  either, will need some elaboration on the notion of\pspage{13} \emph{equivalences
  between Gr-stacks}, which will have to be developed at
  a later stage.
\end{enumerate}

The idea just comes to my mind whether the exactness condition
implying standard amalgamated sums, defining the subcategory of stacks
within $\bB\uphat$, cannot be interpreted in terms of some more
or less obvious topology on $\bB$ turning $\bB$ into a site, as
the subcategory of corresponding \emph{sheaves}. This would mean that
the category (Gr-stacks) is in fact a topos, with the host of
categorical information and topological intuition that goes with such
a situation. In this connection, it is timely to recall that the
related categories \Cat{} of ``all'' categories, and (Groupoids) of
``all'' groupoids, are definitely \emph{not} topoi (if my recollection
is correct -- it isn't immediately clear to me why they are not). This
seems to suggest that, granting that (Gr-stacks) is indeed a topos,
that this would be a rather special feature of the structure species
of infinite order we are working with (as one ward so to say, among a
heap of others, for conceptual sophistication!), in contrast to the
categories (Gr-$n$-stacks) with finite $n$, which presumably are not
topoi. (For a definition of Gr-$n$-stacks in terms of Gr-stacks,
namely Gr-\oo-stacks, see below.)

A related question is whether the category (Gr-stacks) is a model
category for the usual homotopy category, the pair of adjoint functors
considered before satisfying moreover the conditions of Quillen's
comparison theorem\footnote{I am relying here on Quillen's set-up as the only axiomatic set-up for axiomatic homotopy theory which was known to me at that moment. In the course of my reflection, different set-ups are peeling out, which seem more directly suited for finding ``modelizers'' for the specific category $\Hot$ (notably, the so-called \emph{asphericity structures} and, still more importantly, the \emph{contractibility structures}). Thus, the notion of a model category in Quillen's sense is somewhat fading into the background after section 23, where one main new idea gets emphasized, namely the description of a canonical functor $$M \to \Hot$$ for any category $M$ endowed with a structure of a (locally contractible) ``site''\dots}. The obvious idea that comes to mind here, in order
to define the model structure on (Gr-stacks), is to take as ``weak
equivalences'' the maps which are transformed into weak equivalences
by the topological realization functor (which should be readily
expressible in algebraic terms), for cofibrations the monomorphisms,
and defining fibrations by the Serre-Quillen lifting property with
respect to cofibrations which are weak equivalences (with the
expectation that we even get a ``closed model category'' in the sense
of Quillen). Here it doesn't look too unreasonable to expect the same
constructions to work in each of the categories (Gr-$n$-stacks),
$n\ge1$, as well as in the categories ($n$-stacks) without Gr, which
are still to be defined though.

Coming back upon the question of a suitable topology on $\bB$, the
idea that comes to mind immediately of course is to define covering
families of an object $K$ of $\bB$, i.e., of $\bB_0$, in terms
of the components which occur in the description of $K$ as iterated
amalgamated sum of cells $D_n$. A quick glance (too quick a glance?)
seems to show this is indeed a topology, and that the sheaves for this
topology are what we expect.

\bigbreak

\presectionfill\ondate{5.3.}\pspage{14}\par

\hangsection{Are model categories sites?}\label{sec:20}%
I had barely stopped writing last Tuesday, when it became clear that
the ``quick glance'' had been too quick indeed. As a matter of fact,
the ``topology'' I was contemplating on $\bB$ in terms of covering
families does not satisfy the conditions for a ``site'' -- that
something was fishy first occurred to me through the heuristic
consequence, that the functors ``$n$-th component''
\[F_* \mapsto F_n\]
from stacks to sets are fiber functors, or equivalently, that direct
limits in the category of Gr-stacks can be computed componentwise --
now this is definitely false, for much the same reasons why it is
false already in the category \Cat{} of categories. However, it
occurred to me that the latter category \Cat{}, although not a topos
(for instance because, according to Giraud's paper on descent theory
of 1965,\scrcomment{\cite{Giraud1964}} the implications for epimorphisms
\[
\begin{tikzcd}[row sep=tiny,column sep=huge,arrows=Rightarrow]
  & \text{effective} \ar[dr] & \\
  \text{effective universal} \ar[ur] \ar[dr] & & \text{just epi} \\
  & \text{universal} \ar[ur]
\end{tikzcd}
\]
are strict), there \emph{is} a very natural topology, turning it into
a site, namely the one where a family of morphisms namely functors
\[ A_i \to A\]
is covering if and only if the corresponding family in \Sssets
\[ \text{Nerve}(A_i) \to \text{Nerve}(A)\]
is covering, i.e., if{f} every sequence of composable arrows in $A$
\[ a_0 \to a_1 \to a_2 \to a_3 \to \dots \to a_n\]
can be lifted to one among the $A_i$'s. As a matter of fact, this
condition (where it suffices to take $n=2$, visibly) is equivalent
(according to Giraud) to the condition that the family be
``universally effectively epimorphic'', i.e., covering with respect to
the ``canonical topology'' of \Cat{}\footnote{5.3. This is definitely false, see §24 below}. This suggests a third ward for associating
to a category $A$ a topology-like structure, namely the topos of all
sheaves over $\Cat_{/A}$, the site of all categories over $A$,
endowed with the topology induced by the canonical topology of \Cat{}
(which is indeed the canonical topology of $\Cat_{/A}$). It should be
an easy exercise in terms of nerves to check that the homotopy type (a
priori, a pro-homotopy type) associated to this site is just ``the''
homotopy type of $A$, defined either as the homotopy type of
Nerve$(A)$, or as the (pro)homotopy type of the topos $\Ahat$ of
all presheaves over $A$. This of course parallels the similar familiar
fact in the category \Sssets{} of all semisimplicial sets, namely that
for such a ss~set $K$, the homotopy type of $K$ can be viewed also as
the homotopy type of the induced topos $\Sssets_{/K}$ of all ss~sets
over $K$. The difference in the two cases is that in the second
case\pspage{15} the category \Sssets{} (and hence any induced
category) is already a topos, namely equivalent to the category of
sheaves on the same for the canonical topology, whereas on the
category \Cat{} this is not so.

These reflections suggest that in most if not all categories of models
encountered for describing usual homotopy types, there is a natural
structure of a site on the model category $M$
(\emph{presumably}\footnote{It turns out that the relevant topology is not always the canonical one - however the intuition which appears here is going to be basic in the whole ``Modelizing Story''.} the
one corresponding to the canonical topology), with the property that
the homotopy type of any object $X$ in $M$ is canonically isomorphic
to the (pro)homotopy type of the induced site $M_{/X}$, or what
amounts to the same of the corresponding topos
\[ (M_{/X})^\sim = M^\sim_{/X}.\]
I would expect definitely this to be the case for each of the
categories (Gr-stacks), ($n$-Gr-stacks) (although this is not really a
model category in the sense of Quillen, as it gives rise only to
$n$-truncated homotopy types \ldots), (stacks) and ($n$-stacks),
corresponding to an arbitrary choice of coherator, defining the notion
of a (Gr-stack), or of a stack (for the latter and relations between
the two, see below).

\hangsection[Further glimpse upon the ``bunch'' of possible model
\dots]{Further glimpse upon the ``bunch'' of possible model categories
  and a relation between \texorpdfstring{$n$}{n}-complexes and
  \texorpdfstring{$n$}{n}-stacks.}\label{sec:21}%
I would like to digress a little more, to emphasize still about the
vast variety of algebraic structures giving rise to model categories
for the usual homotopy category, or at any rate suitable for
expressing more or less arbitrary homotopy types. Apart from stacks,
where everything is still heuristics for the time being, we have
noticed so far three examples of such structures, namely as complexes, cubical complexes, and categories. There are a few familiar variants of the two
former, such as the ``simplicial complexes'' (in contract to
semi-simplicial ones), namely presheaves on the category of non-empty
finite sets, which can be interpreted as ss~complexes enriched with
symmetry operations on each component -- and there is the
corresponding variant in the cubical case. It would be rather
surprising that there were not just as good model categories, as the
more habitual ones\footnote{It will turn out later that these kinds of complexes, as well as the multicomplexes envisioned below (even the ``mixed'' ones, namely semisimplicial or simplicial for some indices, and cubical (with or without symmetries), or even hemispherical, for some others\dots) do indeed model homotopy types. This is a consequence of the criterion for test categories (section 31), and of the fact that a product of test categories is again a test category (cf. section 74).} -- all the more as the singular complex (simplicial
or cubical) is naturally endowed with this extra structure, which one
generally chooses to forget. More interesting variants are the
$n$-multicomplexes (simplicial or cubical, with or without
symmetries), defined by contravariant functors to sets in $n$
arguments rather than in just one, where $n\ge1$. These complexes are
familiar mainly, it seems, because of their connections with product
spaces and the K\"unneth-Eilenberg-Zilber type relations. It is
generally understood that to such a multicomplex is
associated\pspage{16} the corresponding ``diagonal'' complex, which is
just a usual complex and adequately describes the ``homotopy type'' of
the multicomplex. So why bother with relatively messy kinds of models,
when just usual complexes suffice! Here however the point is not to
get the handiest possible model categories (whatever our criteria of
``handyness''), but rather to get an idea of the variety of algebraic
structures suitable for defining homotopy types, and perhaps to come
to a clue of what is common to all these. Moreover, I feel the
relation between ordinary complexes and $n$-multicomplexes, is of much
the same nature as the relation between just ordinary categories,
perfectly sufficient for describing homotopy types, and $n$-categories
or $n$-stacks. This reminds of course of Quillen's idea of
of defining $n$-categories or $n$-semisimplicial complexes, in much the same way as categories can be
defined (via the Nerve functor) in terms of usual ss~complexes. I hit
again upon multisimplices (without symmetry), when trying to reduce to
a minimum the category $\bB_0$ of what should be called ``standard''
amalgamated sums of the cells $D_n$, where my tendency initially
(\hyperref[ch:I]L, p.~\ref{p:L.7}) had rather been to be as generous
as possible, in order to be as stringent as imaginable for the
completion condition \ref{it:8.A} of (\hyperref[ch:I]L, p.~%
\ref{p:L.8}).  Now it turns out that the coherence relations which
seem to have been written down so far (and the like of which
presumably will suffice to imply full completeness of coherence
relations, in the sense of \ref{it:8.A}) make use only of very
restricted types of such amalgamated sums, expressible precisely in
terms of multisimplices. This I check for instance on the full list of
data and axioms for B\'enabou's ``bicategories'' namely $2$-stacks, in
his 1967 Midwest Category Seminar expos\'e (already referred to). I'll
have to come back upon this point with some care, which gives also a
pretty natural way for getting Quillen's functor from $n$-stacks to
$n$-ss~complexes.

\hangsection{Oriented sets as models for homotopy types.}\label{sec:22}%
All these examples of possible models for homotopy types can be viewed as
generalizations of usual complexes, or of usual categories; I would
like to give a few others which go in the opposite direction -- they
may be viewed as particular cases of categories. One is the (pre)order
structure, which may be viewed as a category structure when the map
$\Fl{} \to \Ob{} \times \Ob{}$ defined by the source and target maps
is injective. Such category is equivalent (and hence homotopic) to the
category associated with the corresponding \emph{ordered} set (when
$x\ge y$ and $y\ge x$ imply $x=y$). Ordered sets are more familiar I
guess as model objects for describing combinatorially a topological
space, in terms of a ``cellular subdivision'' by compact\footnote{The notes are a little floppy here - the ``cells'' are going to be compact only under suitable restriction on the ordered set $K$, for instance if $K$ is finite - however, it is clear I want also to consider arbitrary ordered sets.} subsets or
``cells'' (``strata'' would be a more\pspage{17} appropriate term),
which actually need not be topological cells in the strict sense, but
rather conical (and hence contractible) spaces, each being
homeomorphic to the cone over the union of all strictly smaller strata
(this union is compact). The ordered set associated to such a
(conically) stratified space $X$ is just the set of strata, with the
inclusion relation, and it can be shown that there is a perfect
dictionary between the topological objects (at least in the case of
finite or locally finite stratifications), and the corresponding
(finite or locally finite) ordered sets, via a ``topological
realization functor''
\[ X \mapsto \abs X\]
from ordered sets to (conically stratified) topological spaces. As a
matter of fact, when $K$ is finite, $K$ is endowed with a canonical
triangulation (the so-called barycentric subdivision), the
combinatorial model of which (of ``maquette'') is as follows: the
vertices are in one-to-one correspondence with the elements of $K$
(they correspond to the vertices of the corresponding cones), and the
combinatorial simplices are the ``flags'', or subsets of $K$ which are
totally ordered for the induced order. This is still OK when $K$ is
only locally finite (restricting of course to subsets of $K$ which are
finite, when describing simplices), but in any case we can define (via
infinite direct limits in the category \Spaces of all topological
spaces) the geometric realization of $K$, together with an
interpretation of $K$ as the ordered set of strata of $\abs K$. As a
matter of fact, we get a canonical isomorphism (nearly tautological)
\[ \abs K \simeq \abs{\text{Nerve}(K)} \]
where in the second member, we have written shortly $K$ for the
category defined by $K$. I did not reflect whether it was reasonable
to expect that the categories \Preord{} and \Ord{} of preordered and
ordered sets are model categories, or even closed model categories, in
Quillen's precise sense -- but it is clear though that using ordered
sets we'll get practically any homotopy type, in any case any homotopy
type which can be described in terms of locally finite
triangulations. However it should be noted that the inclusion functors
into \Cat{} or \Sssets
\[ \Preord \hookrightarrow \Ord \hookrightarrow
\Cat \hookrightarrow \Sssets,\]
while giving the correct results on homotopy types, do not satisfy
Quillen's general conditions on pairs of adjoint functors between
model categories -- namely the adjoint functor say
\[ \Cat \to \Preord, \]
associating to a category $A$ the set $\Ob A$ with the obvious
relation, does \emph{not} commute to formation of homotopy types, as
we see in the trivial case when $\Ob A$ is reduced to one point \ldots

There\pspage{18} is another amusing interpretation of the homotopy
type associated to any preordered set $K$, via the topological space
whose underlying set is $K$ itself, and where the closed sets are the
subsets $J$ of $K$ such that $x\in J$, $y\subseteq x$ implies $y\in
J$. This is a highly non-separated topology $\tau$ (except when the
preorder relation is the discrete one), where an arbitrary union of
closed subsets is again closed. I doubt its singular homotopy type to
make much sense, however its homotopy type as a topos does, and
(possibly under mild local finiteness restrictions) it should be the
same as the homotopy type of $K$ just envisioned. Thus, a sheaf on the
topological space $K$ can be interpreted via its fibers as being just
a \emph{covariant} functor
\[ K \to \Sets \]
(NB\enspace the open sets of $K$ are just the closed sets of $K\op$, namely
for the opposite order relation, and thus every $x\in K$ has a
smallest neighborhood, namely the set $K_{\ge x}$ of all $y\in K$ such
that $y\ge x$), or what amounts to the same, a sheaf on the topos
$(K\op)\uphat$ defined by the opposite order. Hence the derived
functors of the ``sections'' functor, when working with abelian
sheaves on $(K,\tau)$, i.e., the cohomology of the topological space
$(K,\tau)$, can be interpreted in terms of the topos associated to
$K\op$. This suggests that the definition I gave of the topology of
$K$ was awkward and maybe it is indeed (although it is the more
natural one in terms of incidence relations between open strata of
$\abs K$), and that we should have called ``open'' the sets I called
``closed'' and vice-versa, or equivalently, replace $K$ by $K\op$ in
the definition I gave of a topology on the set $K$. But as far as
homotopy types are concerned, it doesn't make a difference, namely the
homotopy types associated to $K$ and $K\op$ are canonically
isomorphic. This can be seen most simply on the topological
realizations, via a homeomorphism (not only a homotopism)
\[ \abs K \simeq \abs{K\op},\]
coming from the fact that the maquettes (by which I mean the
combinatorial model for a triangulation, which H. Cartan time ago called
``sch\'ema simplicial'' \ldots) of the two spaces are canonically
isomorphic, because the ``flags'' of $K$ are $K\op$ are the same. A
similar argument due to Quillen using the nerves shows that for any
category $K$ (not necessarily ordered), $K$ and $K\op$ are homotopic,
although $K\uphat$ and $(K\op)\uphat$ are definitely not
equivalent, i.e., not ``homeomorphic''.

To come back to the decreasing cascade of algebraic structure suitable
for describing homotopy types, we could go down one more step still,
to the category of ``maquettes'' (Maq), namely sets $S$ together with
a family $K$ of finite subsets (the simplices of the set of
vertices\pspage{19} $S$), such that one-point subsets are simplices
and a subset of a simplex is a simplex. This category, via the functor
$(S,K) \mapsto K$, is equivalent to the full subcategory of \Ord,
whose objects are those ordered sets $K$, such that for every $x\in
K$, the set $K_{\le x}$ be isomorphic to the ordered set of non-empty
subsets of some finite set (or ``simplex''). Here the question whether
this category is a model category in the technical sense (of Quillen) doesn't
really arise\footnote{It does make sense however to ask whether this category is a ``modelizer'' (in the precise sense of section 28) - I believe it is, but didn't check it yet.}, because this category doesn't even admit finite products
-- rien \`a faire!\footnote{\alsoondate{5.3.} indeed, the notion of a maquette is
  \emph{not} an algebraic structure species!}

\hangsection[Getting a basic functor $M\to\Hot$ from a site
structure \dots]{Getting a basic functor
  \texorpdfstring{$M\to\Hot$}{M->(Hot)} from a site structure
  \texorpdfstring{$M$}{M} \texorpdfstring{\textup(}{(}altering
  beginning of a systematic
  reflection\texorpdfstring{\textup)}{)}.}\label{sec:23}%
It may be about time to get back to stacks, still I can't help going
on pondering about algebraic structures as models for homotopy
types. If we have any algebraic structure species, giving rise to a
category $M$ of set-theoretic realizations, the basic question here
doesn't seem so much whether $M$ is a model category for a suitable
choice of the three sets of arrows (fibrations, cofibrations, weak
equivalences), but rather how to define a natural functor
\begin{equation}
  \label{eq:23.star}
  M \to \Hot, \tag{*}
\end{equation}
where \Hot{} is the category of usual homotopy types, and see whether
via this functor \Hot{} can be interpreted as a category of fractions
(or ``localization'') of $M$ -- namely, of course, by the operation of
making invertible those arrows in $M$ which are transformed into
isomorphisms in \Hot{}. In any case, if we have such a natural functor,
the natural thing to do is to call those arrows ``weak
equivalences''. If we want $M$ to be a category of models, various
examples suggest that the natural thing again is to take as
cofibrations the monomorphisms, and then (expecting that the model
categories we are going to meet will be closed model categories) to
define fibrations by the Serre-Quillen lifting property with respect
to cofibrations (=monomorphisms) which are weak equivalences. This
being done, it becomes meaningful to ask if indeed $M$ is a category
of models.

Now the reflections of the beginning of today's notes (section \ref{sec:20})
suggest a rather natural way for describing a functor
\eqref{eq:23.star}, which makes sense in fact, in principle, for any
category $M$, namely: endow $M$ with its canonical topology (unless a
still more natural one appears at hand\footnote{As we'll see, in the all-important case $\bB = \Cat$, the canonical topology is definitely \emph{not} the right one!} -- I am not sure there is any
better one in the present context),\footnote{\alsoondate{5.3.} the canonical
  topology is \emph{not} always suitable, see \S24 below.} and assume
that for every $X\in\Ob M$, the pro-homotopy type of the induced site
$M_{/X}$ is essentially constant, i.e., can be identified with an
object in \Hot{} itself. We then get the functor \eqref{eq:23.star} in
an obvious way. It then becomes meaningful to ask whether this functor
is a localization functor.

When $M$ is defined in terms of an algebraic structure species,
it\pspage{20} admits both types of limits, without finiteness
requirement -- and we certainly would expect indeed at least existence
of finite \emph{and} infinite direct sums in $M$, if objects of $M$
were to describe arbitrary homotopy types. However, in view of the
special exactness properties of \Hot{}, which are by no means autodual,
we will expect moreover direct sums in $M$ to be ``universal and
disjoint'', in Giraud's sense. This condition, which characterizes to
a certain extent categories which at least mildly resemble or parallel
categories such as \Sets, \Spaces{} and similar categories, whose
objects more or less express ``shapes'' -- this condition at once
rules out the majority of the most common algebraic structures, such
as rings, groups, modules over a ring or anything which yields for $M$
an abelian category, etc. If we describe an algebraic structure
species in terms of its universal realization in a category stable
under finite inverse limits, then such a structure species can be
viewed as being defined by such a category $\bC$, and its
realizations in any other such $C$ as the left-exact functors $\bC\to C$
(the universal realization of the structure within $\bC$
corresponding to the identity functor $\bC\to \bC$). In terms of
the dual category $\bB$, associating to every element in $\Ob{\bB}
= \Ob{\bC}$ the covariant functor $\bC\to\Sets$ it
represents, we get a fully-faithful embedding
\[\bB\hookrightarrow M,  \]
by which $\bB$ can be interpreted as the category of the
(set-theoretic) realizations of the given structure which are of
``finite presentation'' in a suitable sense (in terms of a given
family of generators of $\bB$ namely cogenerators of $\bC$,
considered as corresponding to the choice of ``base-sets'' for the
given structure species -- such choice however being considered as a
convenient way merely to describe the species in concrete terms\dots).
If I remember correctly, $M$ can be deduced from $\bB$, up to
equivalence, as being merely the category of Ind-objects of $\bB$,
i.e., the inclusion functor above yields an equivalence of categories
\[ \text{Ind}(\bB) \toequ M, \]
which implies, I guess, that the exactness properties of $M$ mainly
reflect those of $\bB$. Thus I would expect the condition we want
on direct sums in $M$ to correspond to the same condition for finite
direct sums in $\bB$, not more not less. Thus the algebraic
structure species satisfying this condition should correspond exactly
to small categories $\bB$, stable under finite direct limits, and
such that finite sums in $\bB$ are disjoint and universal. This
condition is presumably necessary, if we want the functor
\eqref{eq:23.star} from $M$ to \Hot{} to be defined and to be a
localization functor -- a condition which it would be nice to
understand directly in\pspage{21} terms of $\bB$, and
(presumably) in terms of the canonical topology of $\bB$, which
should give rise to a localization functor
\[ \bB \to \Hot_{\text{ft}}, \]
where the subscript ft means ``finite type'' -- granting that the
notion of homotopy types of finite type (presumably the same as
homotopy types of finite triangulations, or of finite CW space) is a
well-defined notion\footnote{The consideration of homotopy types ``of finite presentation'' turns out to be irrelevant at this stage of the reflection, leading up mainly to the notion of ``test categories'' and ``test functors'', which are independent of ``finite type'' kind of conditions.}. As usual, it is in $\bB$, not in $\bC$,
that geometrical constructions take place which make sense for
topological intuition. More specifically, it seems that in the cases
met so far, there \emph{are} indeed privileged base-sets for the
structure species considered (such as the ``components'' or a
semisimplicial or cubical complex, or of a stack, etc.), indexed by
the natural integers or $n$-tuples of such integers, and which
``correspond'' to topological cells of various dimensions. Moreover,
some of the basic structural monomorphic maps between these objects of
$\bB$ define cellular decompositions of the topological spheres
building these topological cells. These objects and ``boundary maps''
between them define a (non-full, in general) subcategory of $\bB$,
say $\bB_\oo$, which looks like the core of the category $\bB$,
from which the topological significance of $\bB$ is springing. It
is in terms of $\bB_\oo$ that the ``correspondence'' (vaguely
referred to above) with topological cells and spheres takes a precise
meaning. Namely, associating to any object of $\bB_\oo$ the ordered
set of its subobjects (within $\bB_\oo$ of course, not $\bB$),
the (stratified) topological realization of this ordered set is a
cell, the family of subcells of smaller dimension (namely different
from the given one) is a cellular subdivision of the sphere bounding
this cell\footnote{For the time being, these considerations of cellular structures turn out to be irrelevant - compare the comments in section 26. The reflection on such cellular structures is taken up in section 91, $8^\circ$, but later again it is kind of superseded by the formalism of ``integratos'' and ``cointegrators'' (in part $V$ of the notes, on abelianization of homotopy types) which makes sense in $\hat{A}$ for an arbitrary small category $A$}.

\bigbreak
\presectionfill\ondate{7.3.}\par

\hangsection{A bunch of topologies on \texorpdfstring{\Cat}{(Cat)}.}%
\label{sec:24}%
Yesterday there occurred to me a big ``\'etourderie'' again of the day
before, in connection with a reflection on a suitable ``natural'' site
structure on the category \Cat{} -- namely when asserting that for a
family of morphisms, i.e., functors in \Cat{}
\begin{equation}
  \label{eq:24.star}
  A_i \to A, \tag{*}
\end{equation}
in order for the corresponding family to be ``covering'' namely
epimorphic in the category (a topos, as a matter of fact) \Sssets,
namely for the corresponding families of mappings of sets
\begin{equation}
  \label{eq:24.starstar}
  \Fl_n(A_i) \to \Fl_n(A) \tag{**}
\end{equation}
to be epimorphic, it was sufficient that his condition be satisfied
for $n=2$ (which, according to Giraud, just means that the family
\eqref{eq:24.star} is covering for the canonical topology of
\Cat{}). This is obviously false -- morally, it would mean that an
$n$-simplex is ``covered'' in a reasonable sense by its
sub-$2$-simplices, which is pretty absurd. To\pspage{22} give
specific positive statements along these lines, let's for any
$N\in\bN$, denote by $T_N$ the topology on \Cat{} for which a family
\eqref{eq:24.star} is covering if{f} the families
\eqref{eq:24.starstar} are for $n\le N$ -- which means also that the
corresponding family of $N$-truncated nerves
\[ \text{Nerve}(A_i)[N] \to \text{Nerve}(A)[N] \]
is epimorphic. For $N=1$, this means that the family is ``universally
epimorphic'', for $N=2$, that it is ``universally effectively
epimorphic'', i.e., covering for the canonical topology of \Cat{}
(Giraud, loc.\ cit.\ page~28). It turns out that this decreasing
sequence of topologies on \Cat{} is strictly decreasing -- as a matter
of fact, denoting by $\Delta[N]$ the category of standard, ordered
simplices of dimension $\le N$, and using the inclusions
\[ \Delta[N] \hookrightarrow \Cat \hookrightarrow
\Delta[N]\uphat \]
(for $N\ge 2$ say), an immediate application of the ``comparison
lemma'' for sites shows that we have an equivalence
\[ (\Cat, T_N)^\sim \simeq {\Delta[N]}\uphat,\]
and hence, for any object $A$ in \Cat{}, namely a category $A$, we get
an equivalence
\[ (\Cat, T_N)_{/A}^\sim \simeq (\Delta[N]_{/A})\uphat,\]
and hence the homotopy type of the first hand member is \emph{not}
described by and equivalent to the homotopy type of the whole
Nerve$(A)$ object, but rather by its $N$-skeleton, which has the same
homotopy and cohomology invariant in dimension $\le N$, but by no
means for higher dimensions. This shows that among the topologies
$T_N$, none is suitable for recovering the homotopy type of objects of
\Cat{} in the way contemplated two days ago (page 15); the one topology
which \emph{is} suitable is the one which may be denoted by $T_\oo$,
and which is the one indeed which I first contemplated (page 14),
before the mistaken idea occurred to me that it was the same as the
canonical topology.

A very similar mistake occurred earlier, when I surmised that the left
adjoint functor $N'$ to the inclusion or Nerve functor $N$ from \Cat{}
to $\Sssets = \Simplexhat$ had the property that for any object
$K$ in \Sssets, $K$ and $N'(K)$ had the same homotopy type. Looking
up yesterday the description in Gabriel-Zisman of this functor, this
recalled to my mind that it factors (via the natural restriction
functor) through the category ${\Delta[2]}\uphat$, hence any
morphism $K \to K'$ inducing an isomorphism on the $2$-skeletons
(which by no means implies that it is a homotopy equivalence) induces
an isomorphism $N'(K) \tosim N'(K')$, and a fortiori a
homotopy equivalence. This shows that, even if it should be true that
\Cat{} is\pspage{23} a model category in Quillen's sense, the
situation with the inclusion functor
\[ N : \Cat \to \Sssets \]
and the left adjoint functor $N'$ is by no means the one of Quillen's
comparison theorem, where the two functors play mutually dual roles
and both induce equivalences on the corresponding localized
homotopical categories. Here only $N$, not $N'$, induces such
equivalence. This is analogous to the situation, already noted before,
of the inclusion functors
\[ \Preord \hookrightarrow \Cat
\quad\text{or}\quad
\Ord \hookrightarrow \Cat,\]
which (presumably) by localization induce equivalences on the associate homotopy
categories, but the left adjoint functors do not share this property.

The topology $T_\oo$ just considered on \Cat{} as a ``suitable''
topology for describing homotopy types of objects of \Cat{}, was of
course directly inspired by the semi-simplicial approach to homotopy
types, via simplicial complexes, namely via the two associated
inclusion functors
\[ \Simplex \hookrightarrow \Cat \hookrightarrow
\Simplexhat,\]
the second functor associating to every category $A$ the ``complex''
of its ``simplicial diagrams'' $a_0 \to a_1 \to \dots \to a_n$. If we
had been working with cubical complexes rather than ss ones for
describing homotopy types, this would give rise similarly to two
functors
\[ \square \hookrightarrow \Cat \hookrightarrow
\square\uphat,\]
where $\square$ is the category of ``standard cubes'' and face and
degeneracy maps between them, and where the second inclusion
associates to every category $A$ the cubical complex of its ``cubical
diagrams'', namely commutative diagrams in $A$ modelled after the
diagram types $\square_n[1]$, the $1$-skeleton of the standard
$n$-cubes $\square_n$, with suitable orientations on its edges
(indicative of the direction of corresponding arrows in $A$). The
natural idea would be to endow \Cat{} with the topology, $T_\oo'$ say,
induced by $\square\uphat$, namely call a family \eqref{eq:24.star}
covering if{f} for every $n\in\bN$, the corresponding family of maps
of sets
\[ \text{Cub}_n(A_i) \to \text{Cub}_n(A)\]
is epimorphic. This topology appears to be coarser than $T_\oo$ (i.e.,
there are fewer covering families), and the comparison lemma gives now
the equivalence
\[ (\Cat,T_\oo')^\sim \simeq \square\uphat,\]
which shows that definitely the topology is \emph{strictly} coarser
than $T_\oo$, as it gives rise to a non-equivalent topos,
$\square\uphat$ instead of $\Simplexhat$.

These reflections convince me 1)\enspace there are indeed topologies
on \Cat{}, suitable for describing the natural homotopy types of
objects of \Cat{} namely of categories, and 2)\enspace that there is
definitely\pspage{24} no privileged choice for such a topology. We
just described two such, but using multicomplexes (cubical or
semisimplicial) rather than simple complexes should give us infinitely
many others just as suitable, and I suspect now that there must be a
big lot more of them still!

\hangsection{A tentative equivalence relation for topologists.}\label{sec:25}%
A fortiori, coming back to the intriguing question of characterizing
the algebraic structure species suitable for describing homotopy types
in the usual homotopy category \Hot{}, and for recovering \Hot{} as a
category of fractions of the category $M$ of all set-theoretic
realizations of this species, it becomes clear now that we cannot hope
reasonably for a ``natural'' topology on $M$, distinguished among all
others, giving rise to the wished-for functor
\begin{equation}
  \label{eq:25.star}
  M \to \Hot \tag{*}
\end{equation}
in the way contemplated earlier (section \ref{sec:23}). Here it occurs to me
that anyhow, if concerned mainly with defining the functor
\eqref{eq:25.star}, we should consider that two topologies $T,T'$ on
$M$ such that $T \ge T'$ and hence giving rise to a morphism of topoi
\begin{equation}
  \label{eq:25.starstar}
  M_T^\sim \to M_{T'}^\sim \tag{**}
\end{equation}
(the direct image functor associated to this morphism being the
natural inclusion functor, when considering $T$-sheaves as particular
cases of $T'$-sheaves) are ``\emph{equivalent}'' -- maybe we should
rather say ``\emph{Hot-equivalent}'' -- if for any object $A$ in $M$, the
induced morphism
\[ (M_T^\sim)_{/A} \to (M_{T'}^\sim)_{/A} \]
is a homotopy equivalence. This can be viewed as an intrinsic property
of the morphism of topoi \eqref{eq:25.starstar}, of a type rather
familiar I guess to people used to the dialectics of \'etale
cohomology, where a very similar notion was met and given the name of
a ``globally acyclic morphism''. The suitable name here would be
``globally aspheric morphism'' which is a reinforcement of the former,
in the sense of being expressible in terms of isomorphism relations in
cohomology with arbitrary coefficient sheaves on the base,
\emph{including non-commutative} coefficient sheaves\footnote{For a \emph{proper} map of paracompact spaces, this condition is equivalent to the fibers being ``aspheric'' (i.e. weakly contractible).}. The relation just introduced between two topologies $T,T'$
on a category $M$ makes sense for any $M$ (irrespective of the
particular way $M$ was introduced here), it is not yet an equivalence
relation though -- so why not introduce the equivalence relation it
generates, and call this ``Hot-equivalence'' -- unless we find a
coarser, and cleverer notion of equivalence, deserving this name. The
point which, one feels now, should be developed, is that this notion
of equivalence should be the coarsest we can find out, and which still
implies that to any Hot-equivalence class of topologies on $M$ there
should be canonically associated a functor \eqref{eq:25.star}, which
should essentially be ``the'' common value of all the similar
functors, associated to the topologies $T$ within\pspage{25} this
class. Maybe even it could be shown that this equivalence class can be
recovered in terms of the corresponding functor \eqref{eq:25.star}, in
the same way as (according to Giraud) a ``topology'' on $M$ can be
recovered from the associated subtopos of $M\uphat$, namely the
associated category of sheaves on $M$ (in such a way that the set of
``topologies'' on $M$ can be identified with the set of ``closed
subtopoi'' of $M\uphat$). The very best one could possibly hope for
along these lines would be a one-to-one correspondence between
isomorphism classes of functors \eqref{eq:25.star} (satisfying certain
properties?), and the so-called ``Hot-equivalence classes'' of
topologies on $M$.

Whether or not this tentative hope is excessive, when it comes to the
(still somewhat vague) question of understanding ``which algebraic
structure species are suitable for expressing homotopy types'', it
might not be excessive though to expect that in all such cases, $M$
should be equipped with just \emph{one} ``natural'' Hot-equivalence
class of topologies on $M$, which moreover (one hopes, or wonders)
should be expressible directly in terms of the intrinsic structure of
$M$, or, what amounts to the same, in terms of the full subcategory
$\bB$ of objects of ``finite presentation'', giving rise to $M$ via
the equivalence
\[ M \simeq \text{Ind}(\bB). \]
As we just saw in the case $M = \Cat$, the so-called
``canonical topology'' on $M$ need not be within the natural
Hot-equivalence class -- and I am at a loss for the moment to give a
plausible intrinsic characterization of the latter, in terms of the
category $M=\Cat$\footnote{This point of view of ``$\Hot$-equivalence classes'' of topologies is going to lead to the notions of test categories and test functors, which have proved so handy that I didn't feel later a need to come back to the former point of view, in terms of the understanding we got with the test notions.}.

\hangsection{The dawn of test categories and test functors\dots}%
\label{sec:26}%
What comes to mind though is that the categories such as $\Simplex$,
$\Square$ and their analogons (corresponding to multicomplexes rather
than monocomplexes, for instance) can be viewed as (generally not
full) subcategories of $M$ (in fact, even of the smaller category
\Ord{} of all ordered sets). The topologies we found on $M$ were in
fact associated in an evident way to the choice of such
subcategories. As was already felt by the end of the reflection two
days ago (p.~\ref{p:21}),
these subcategories (denoted there by $\bB_\oo$)
have rather special features -- they are associated to simultaneous
cellular decompositions of spheres of all dimensions -- and it is this
feature, presumably, that makes the associated ``trivial'' algebraic
structure species, giving rise to the category of set-theoretic
realizations $\bB_\oo\uphat$, eligible for ``describing
homotopy types''. In the typical example $M = \Cat$ though,
contrarily to what was suggested by the end of section \ref{sec:23}
(when thinking mainly of
the rather special although important case when $M$ is expressible as
a category $B\uphat$, for some category $B$ such as\pspage{26}
$\Simplex$, $\Square$ etc.), there is no really privileged choice of
such subcategory $\bB_\oo$ -- we found indeed a big bunch of such,
the ones just as good as the others. The point of course is that the
corresponding topologies on $M$, namely induced from the canonical
topology on $\bB_\oo\uphat$ by the canonical functor
\[ M \to \bB_\oo\uphat,\]
are Hot-equivalent for some reason or other, which should be
understood. The plausible fact that emerges here, is that the
``natural'' Hot-equivalence classes of topologies on $M$ is
associated, in the way just described, to a class (presumably an
equivalence class in a suitable set for suitable equivalence
relation\ldots) of subcategories $\bB_\oo$ in $M$. The question of
giving an intrinsic description of the former, is apparently reduced
to the (possibly more concrete) one, of giving an intrinsic
description of a bunch of subcategories $B=\bB_\oo$ of $M$. This
description, one feels, should both \namedlabel{step:26.a}{a)}\enspace
insist on intrinsic properties of $B$, independent of $M$, namely of
the structure species one is working with, and
\namedlabel{step:26.b}{b)}\enspace be concerned with the particular
way in which $B$ is embedded in $M$, which should by no means be an
arbitrary one.

The properties \ref{step:26.a} should be, I guess, no more no less
than those which express that the ``trivial'' algebraic structure
species defined by $B$, giving a category of set-theoretic
realizations $M_B=B\uphat$, should be ``suitable for describing
homotopy types''. The examples at hand so far suggest that in this
case, the canonical topology on the topos $M_B$ is within the natural
Hot-equivalence class, which gives a meaning to the functor
\[ M_B = B\uphat \to \Hot,\]
indeed it associates to any $a \in \Ob{B\uphat}$ the homotopy type
of the induced category $B_{/a}$ of all objects of $B$ ``over
$a$''\footnote{We use here the canonical equivalence $B\uphat_{/a} \simeq (B_{/a})\uphat$.}.
Thus a first condition on $B$ is that \emph{this functor should be a
``localization functor''}, identifying \Hot{} with a category of
fractions of $M_B = B\uphat$. This does look indeed as an extremely
stringent condition on $B$, and I wonder if the features we noticed in
the special cases dealt with so far, connected with cellular
decompositions of spheres, have any more compulsive significance than
just giving some handy \emph{sufficient} conditions (which deserve to
be made explicit sooner or later!) for ``eligibility'' of $B$ for
recovering \Hot{}. Beyond this, one would of course like to have a
better understanding of what it \emph{really} means, in terms of the internal
structure of $B$, that the functor above from $B\uphat$ to \Hot{} is
a localization functor.

Once\pspage{27} this internal condition on $B$ is understood, step
\ref{step:26.b} then would amount to describing, in terms of an arbitrary
``eligible'' algebraic structure species expressed by the category
$M$, of what we should mean by ``eligible functors''\footnote{A little later we will say ``test functor'' instead of ``eligible functor'', and call $B$ a ``test category'' (in a provisional sense, which will be still strengthened somewhat later, cf. sections 39 or 44).}
\[ B \to M,\]
giving rise in the usual way to a functor
\[ M \to B\uphat.\]
In any case, the latter functor defines upon $M$ an induced topology,
$T$ say, and the comparison lemma tells us that if either $B \to M$ or
$M\to B\uphat$ is fully faithful, then the topos associated to $M$
is canonically equivalent to $B\uphat$ (using this comparison lemma
for the functor which happens to be fully faithful). From this follows
that the functor (*) $M^{\sim} \to \Hot$ defined by
$T$ is nothing but the compositum 
\[ M^\sim \tosimeq B\uphat \to \Hot,\]
and hence a localization functor. In other words, when $B$ is a
category satisfying the condition seen above,\footnote{a ``test
  category'', as we will say} then \emph{any} functor
$B \to M$ satisfying one of the two fullness conditions above yields a
corresponding description of \Hot{} as a localization of $M^\sim$. What
is still lacking though is a grasp on when two such functors $B \to
M$, $B' \to M$ define essentially ``the same'' functor $M \to
\Hot$, or (more or less equivalently) two \Hot{}-equivalent
topologies $T,T'$ on $M$; \emph{is it enough, for instance, that they give rise to the same notion of ``weak equivalences''} (namely morphisms in
$M$ which are transformed into an isomorphism of \Hot{})? And moreover,
granting that this equivalence relation between certain full
subcategories (say) $B$ of $M$ is understood, how to define, in terms
of $M$, a ``natural'' equivalence class of such full subcategories,
giving rise to a canonical functor $M\to\Hot$?\footnote{This is now fairly well understood with the notion of an \emph{asphericity structure} (see notably section 77), and of a \emph{``canonical modelizer''} (the idea of which first occurs in section 50).}

Recalling that the algebraic structure considered can be described in
terms of an arbitrary small category $\bB$ where arbitrary finite
direct limits exist (namely $\bB$ is the full subcategory of $M$ of
objects of finite presentation), it seems reasonably to assume that
indeed
\[B \hookrightarrow \bB \quad\text{(a full embedding)},\]
and the question transforms into describing a natural equivalence
class of such full subcategories, in (more or less) any small category
$\bB$ where finite direct limits exist, and where moreover there
exist such full subcategories $B$. Also, we may have to throw in some
extra conditions on $\bB$, such as the condition that direct sums
be ``disjoint'' and ``universal'' already contemplated before.

Maybe\pspage{28} I was a little overenthusiastic, when observing
for any full embedding of a category $B$ in $M$ (let's call the
categories $B$ giving rise to a localization functor $B\uphat\to
\Hot$ \emph{homotopy-test} categories, or simply \emph{test}
categories) we get a localization functor
\[ M^\sim = (M, T_B)^\sim \to \Hot,\]
where $T_B$ is the topology on $M$ corresponding to the full
subcategory $B$. After all, there is a long way in between $M$ itself
and the category of sheaves $M^\sim$ -- and what we want is to get
\Hot{} as a localization of $M$ itself, not of $M^\sim$. It is not even
clear, without some extra assumptions, that the natural functor from
$M$ to $M^\sim$ is fully faithful, namely that $M$ can be identified
with a full subcategory of the category $M^\sim$ we've got to localize
to get \Hot{}. We definitely would like this to be true, or what
amounts to the same, that the functor $M\to B\uphat$ defined by
$B\to M$ should be fully faithful -- which means also that the full
subcategory $B$ of $M$ is ``generating by strict epimorphisms''\footnote{See digression in section 90 on strictly generating families.}, namely
that for every $K$ in $M$, there exists a strictly epimorphic family
of morphisms $b_i \to K$, with sources $b_i$ in $B$. This
interpretation of full faithfulness of $M \to B\uphat$ is OK when
$B \to M$ is fully faithful, a condition which I gradually put into
the fore without really compelling reason, except that in those
examples I have in mind and which are \emph{not} connected with the
theory of stacks of various kinds, this condition is satisfied
indeed. 

Apparently, with this endless digression on algebraic models
for homotopy types, stacks (which I am supposed to be after, after
all) are kind of fading into the background! Maybe we should after all
forget about the fully faithfulness condition on either $B\to M$ or
$M\to B\uphat$, and just insist that the compositum
\[ M \to M^\sim \to B\uphat \to \Hot\]
(which can be described directly in terms of the topology $T_B$ on $M$
associated to the functor $B\to M$\footnote{And even without using this topology, via the obvious functor $M \to B\uphat$!}) should be a localization
functor. I guess that for a given $M$ or $\bB$, the mere fact that
there should exist a test category $B$ and a functor $B\to M$, or
$B\to\bB$, having this property is already a very strong condition on
the structure species considered, namely on the category $\bB$ which
embodies this structure. It possibly means that the corresponding
functor
\[ \bB \to \Hot\]
factors through a functor
\[ \bB \to \Hot_{\mathrm{f.t.}} \quad \text{(f.t.\ means ``finite
  type'')}\]
which is itself a localization functor. It is not wholly impossible,
after all, that this condition on a functor $B \to \bB$ ($B$ a test
category) is so stringent, that all such functors (for variable $B$)
must be\pspage{29} already ``equivalent'', namely define
Hot-equivalent topologies on $\bB$ (or $M$, equivalently), and hence
define ``the same'' functor $M\to\Hot$ or
$\bB\to\Hot_{\mathrm{f.t.}}$.

\hangsection{Digression on ``geometric realization'' functors.}\label{sec:27}%
All this is pretty much ``thin air conjecturing'' for the time being
-- quite possibly the notion of a ``test category'' itself has to be
considerably adjusted, namely strengthened, as well as the notion of a
``test functor'' $B \to \bB$ or $B \to M$ -- some important features
may have entirely escaped my attention. The one idea though which may
prove perhaps a valid one, it that a suitable localization functor
\begin{equation}
  \label{eq:27.star}
  M \to \Hot \tag{*}
\end{equation}
may be defined, using either various topologies on $M$ (related by a
suitable ``Hot-equivalence'' relation), or various functors $B \to M$
or $B \to \bB$ of suitable ``test categories'' $B$, and how the two
are related. I do not wish to pursue much longer along these lines
though, but rather put now into the picture a third way still for
getting a functor \eqref{eq:27.star}, namely through some more or less
natural functor
\begin{equation}
  \label{eq:27.starstar}
  M \to \Spaces, \quad K \mapsto \abs K, \tag{**}
\end{equation}
called a ``geometric (or topological) realization functor''\footnote{This reflection on geometric realization is going to stop short, and the use of the catgory of topological spaces is fading out altogether.}. There is
a pretty compelling choice for such a functor, in the case of
(semisimplicial or cubical) complexes or multicomplexes of various
kinds, and accordingly for the subcategories \Cat, \Preord, \Ord{} or
\Sssets, using geometric realization of semisimplicial complexes. In
the case of the considerably more sophisticated structure of Gr-stacks
though (or the related structure of stacks, which will be dealt
with in much the same way below), although there is a pretty natural
choice for geometric realization on the subcategory $\bB_\oo$ of $M$
embodying the ``primitive structure'' (namely the structure of an
\oo-graph, see below also); it has been seen that the extension of
this to a functor on the whole of $M$ (via its extension to the left
coherator defining the structure species, which we denoted by $\bB$ at
the beginning of these notes, but which is not quite the $\bB$
envisioned here) is by no means unique, that it depends on a pretty
big bunch of rather arbitrary choices. This indeterminacy now appears
as quite in keeping with the general aspect of a (still somewhat
hypothetical) theory of algebraic homotopy models, gradually emerging
from darkness. It parallels the corresponding indeterminacy in the
choice of an ``eligible'' topology on $M$ (call these topologies the
test topologies), or of a test functor $B \to M$. What I would like
now to do, before coming back to stacks, is to reflect a little still
about the relations between such choice of a ``geometric realization
functor'', and test topologies or test functors relative to $M$.

\bigbreak

\presectionfill\ondate{8.3.}\pspage{30}\par

\hangsection{The ``inspiring assumption''. Modelizers.}\label{sec:28}%
While writing down the notes yesterday, and this morning still while
pondering a little more, there has been the ever increasing feeling
that I ``was burning'', namely turning around something very close,
very simple-minded too surely, without quite getting hold of it
yet. In such a situation, it is next to impossible just to leave it at
that and come to the ``ordre du jour'' (namely stacks) -- and even the
``little reflection'' I was about to write down last night (but it was
really too late then to go on) will have to wait I guess, about the
``geometric realization functors'', as I feel it is getting me off
rather, maybe just a little, from where it is ``burning''!

There was one question flaring up yesterday
(section \ref{sec:26}) which I nearly
dismissed as kind of silly, namely whether two localization functors
\begin{equation}
  \label{eq:28.star}
  M \to \Hot \tag{*}
\end{equation}
obtained in such and such a way were isomorphic (maybe even
canonically so??) provided they defined the same notion of ``weak
equivalence'', namely arrows transformed into isomorphisms by the
localization functors. Now this maybe isn't so silly after all, in
view of the following\par
\noindent \textbf{Assumption}: The category of equivalences of $\Hot$
with itself, and of natural isomorphisms (possibly even any morphisms)
between such, is equivalent to the one point category.

This means 1)\enspace any equivalence $\Hot \tosimeq \Hot$ is
isomorphic to the identity functor, and 2)\enspace any automorphism of
the identity functor (possibly even any endomorphism?) is the
identity.

Maybe these are facts well-known to the experts, maybe not -- it is
not my business here anyhow to set out to prove such kinds of
things. It looks pretty plausible, because if there was any
non-trivial autoequivalence of \Hot, or automorphism of its identity
functor, I guess I would have heard about it, or something of the sort
would flip to  my mind. It would not be so if we abelianized \Hot\
some way or other, as there would be the loop and suspension functors,
and homotheties by $-1$ as automorphisms of the identity functor.

This assumption now can be rephrased, by stating that \emph{a localization
functors \eqref{eq:28.star} from any category $M$ into \Hot\ is well
determined, up to a unique isomorphism, when the corresponding class
$W \subset \Fl(M)$ of weak equivalences is known}, in positive response
to yesterday's silly question!

Such situation \eqref{eq:28.star} seems to me to merit a name. As the
work ``model category'' has already been used in a somewhat different
and more sophisticated sense by Quillen, in the context of homotopy, I
rather use another\pspage{31} one in the situation here. Let's call
a ``modelizing category'', or simply a ``\emph{modelizer}''
(``mod\'elisatrice'' in French), any category $M$, endowed with a set
$W \subset \Fl(M)$ (the weak equivalences), satisfying the obvious
condition:

\noindent\parbox[t]{0.1\textwidth}{(Mod)\par}
\parbox[t]{0.9\textwidth}{\vspace*{-11pt}%
  \begin{enumerate}[,label=\alph*)]
  \item\label{it:Mod.a}
    $W$ is the set of arrows made invertible by the localization
    functor $M \to W^{-1}M$, and
  \item\label{it:Mod.b}
    $W^{-1}M$ is equivalent to \Hot,
  \end{enumerate}}
\footnote{Condition \ref{it:Mod.a} will eventually be weakened (cf. section 41).}or equivalently, there exists a localization functor
\eqref{eq:28.star} (necessarily unique up to unique isomorphism) such
that $W$ be the set of arrows made invertible by this functor.

Let $(M,W)$, $(M',W')$ be two modelizers, a functor $F: M\to M'$ is
called \emph{model-preserving}, or a \emph{morphism} between the
modelizers, if it satisfies either of the following equivalent
conditions:
\begin{enumerate}[label=(\roman*)]
\item $F(W) \subset W'$, hence a functor $F_{W,W'}: W^{-1}M \to
  W'^{-1}M'$, and the latter is an equivalence.
\item The diagram
\[  \begin{tikzcd}[column sep=tiny]
    M \ar[dr] \ar[rr, "F"] & & M'\ar[dl] \\ & \Hot &
  \end{tikzcd}\]
  is commutative up to isomorphism (where the vertical arrows are the
  ``type functors'' associated to $M,M'$ respectively).
\end{enumerate}
When dealing with a modelizer $(M,W)$, $W$ will be generally
understood so that we write simply $M$. When $M$ is defined in terms
of an algebraic structure species, the task arises to find out whether
(if any) there exists a \emph{unique} $W \subset \Fl(M)$ turning $M$
into a modelizer, and if not so, if we can however pinpoint one which
is a more natural one, and which we would call ``canonical''.

Here is a diagram including most of the modelizers and
model-preserving functors between these which we met so far\footnote{I didn't prove yet that $\Ord$ and $\Preord$ are modelizers! I suppose they are\dots} (not
included however those connected with the theory of ``higher'' stacks
and Gr-stacks, which we will have to elaborate upon later on):
\[
\begin{tikzcd}[row sep=tiny,column sep=small]
  & & & \Simplexhat = \text{(ss~sets)}\ar[dr, "\xi"] & \\
  \Ord\ar[r, hook] & \Preord\ar[r, hook] & \Cat\ar[ur, hook, "\alpha"]
  \ar[dr, hook, swap, "\beta"] & & \Cat \\
  & & & \square\uphat = \text{(cub.\ sets)} \ar[ur, swap, "\eta"] &
\end{tikzcd}\]
where the two last functors, with values in \Cat, are the two obvious functors,
obtained from
\begin{equation}
  \label{eq:28.starstar}
  i_A : \Ahat \to \Cat, \quad F\mapsto A_{/F} \tag{**}
\end{equation}
by particularizing to $A=\Simplex$ or $\Square$. As we noticed before,
the four first among these six functors admit left adjoints, but
except for the first, these adjoints are \emph{not} model
preserving. The two last functors, and more\pspage{32} generally
the functor \eqref{eq:28.starstar}, admit right adjoints, namely the
functor
\[ j_A: \Cat\to\Ahat,\]
where $j_A(B) = (a \mapsto \Hom(A_{/a},B))$ ($B\in\Ob\Cat$). It should
be noted that the functors $j_\Simplex, j_\Square$ are \emph{not} the
two functors $\alpha,\beta$ which appear in the diagram above, the
latter are associated to the familiar functors
\begin{equation}
  \label{eq:28.starstarstar}
  \begin{tikzcd}[row sep=tiny]
    \Simplex\ar[dr] & \\ & \Cat & \\ \Square\ar[ur] &
  \end{tikzcd}\tag{***}
\end{equation}
(factoring in fact through \Ord), associating to every ordered
simplex, or to each multiordered cube, the corresponding $1$-skeleton
with suitable orientations on the edges, turning the vertices of this
graph into an ordered set; while the two former are associated to the
functors deduced from \eqref{eq:28.starstar} by restricting to $A
\subset \Ahat$, namely
\[ n \mapsto \Simplex_{/\Simplex_n}\]
in the case of $\Simplex$, and accordingly for $\Square$. The values of
these functors, contrarily to the two preceding ones, are
\emph{infinite} categories, and they cannot be described by (i.e.,
``are'' not) (pre)ordered sets. If however we had defined the
categories $\Simplex,\Square$ in terms of iterated boundary operations
only, excluding the degeneracy operations (which, I feel, are not
really needed for turning $\Simplexhat$ and $\Square\uphat$ into
modelators), we would get indeed \emph{finite} ordered sets, namely
the full combinatorial simplices or cubes, each one embodied by the
ordered set of all its facets of all possible dimensions.

Contrarily to what happens with the functors $\alpha,\beta$, I feel
that for the two functors $\xi,\eta$ in opposite direction, not only
are they model preserving, but the right adjoint functors
$j_\Simplex,j_\Square$ must be model preserving too, and we will have to
come back upon this in a more general context.

We could amplify and unify somewhat the previous diagram of
modelizers, by introducing multicomplexes, which after all can be as
well ``mixed'' namely partly semisimplicial, partly cubical. Namely,
we may introduce the would-be ``test categories''\footnote{That these are test categories indeed follows from section 74, and it is equally true that the natural functors from these into $\Cat$ (used below in order to describe the functors $\alpha_{p, q}$) are test functors (cf. remark after the proof of prop.2 in section 82), and therefore $\alpha_{p, q}$ is indeed modelizing.}
\[ \Simplex^p \times \Simplex^q = T_{p,q}\quad (p,q\in\bN, p+q\ge1)\]
giving rise to the category ${T_{p,q}}\uphat$ of
$(p,q)$-multicomplexes ($p$ times simplicial, $q$ times cubical). We
have a natural functor (generalizing the functors
\eqref{eq:28.starstarstar})
\[ T_{p,q} \to \Ord {( \hookrightarrow \Cat)},\]
associating to a system of $p$ standard simplices and $q$ standard
cubes (of variable dimensions), the \emph{product} of the $p+q$
associated ordered sets. We get this way a functor
\[\alpha_{p,q} : \Cat \to {T_{p,q}}\uphat\]
which\pspage{33} presumably (as I readily felt yesterday, cf.\
first lines p.~\ref{p:24})
is not any less model-preserving than the functors
$\alpha,\beta$ it generalizes. Of course, taking $A=T_{p,q}$ above, we
equally get a natural functor
\[i_{p,q} : {T_{p,q}}\uphat \to \Cat \]
admitting a right adjoint $j_{p,q}$, and both functors I feel must be
model preserving.

\hangsection[The basic modelizer \Cat.  Provisional definition of test
\dots]{The basic modelizer \texorpdfstring{\Cat}{(Cat)}.  Provisional
  definition of test categories and elementary
  modelizers.}\label{sec:29}%
It is time now to elaborate a little upon the notion of a test
category, within the context of modelizers. Let $A$ be a small
category, and consider the functor \eqref{eq:29.starstar}
\begin{equation}
  \label{eq:29.starstar}
  i_A : \Ahat \to \Cat, \quad F \mapsto A_{/F}.\tag{**}
\end{equation}
Whenever we have a functor $i: M \to M'$, when $M'$ is equipped with a
$W'$ turning it into a modelizer, there is (if any) just one $W
\subset \Fl(M)$ turning $M$ into a modelizer and $i$ into a morphism
of such, namely $W = i_{\Fl{}}^{-1}(W')$. In any case, we may define
$W$ (``weak equivalences'') by this formula, and get a functor
\[ W^{-1}M \to W'^{-1}M',\]
which is an equivalence if{f} $(M,W)$ is indeed a modelizer and $i$
model preserving. We may say shortly that $i:M\to M'$ is model
preserving, even without any $W$ given beforehand. Now coming back to
the situation \eqref{eq:29.starstar}, the understanding yesterday was
to call $A$ a \emph{test category}, to express that the canonical
functor \eqref{eq:29.starstar} is model preserving. (In any case,
unless otherwise specified by the context, we will refer to arrows in
$\Ahat$ which are transformed into weak equivalences of \Cat\ as
``\emph{weak equivalences}''.) It may well turn out, by the way, that
we will have to restrict somewhat still the notion of a test category.

In any case, \emph{the basic modelizer, in this whole approach to homotopy
models, is by no means the category \Sssets\ (however handy) or the
category \Spaces\ (however appealing to topological intuition), but
the category \Cat\ of ``all'' (small) categories}. In this setup, the
category \Hot\ is most suitably defined as the category of fractions
of \Cat\ with respect to ``weak equivalences''. These in turn are most
suitably defined in cohomological terms, via the corresponding notion
for topoi -- namely a morphism of topoi
\[ f : X \to X' \]
is a ``weak equivalence'' or homotopy equivalence, if{f} for every
locally constant sheaf $F'$ on $X'$, the maps
\[ \mathrm H^i(X',F') \to \mathrm H^i(X, f^{-1}(F') ) \]
are isomorphisms whenever defined -- namely for $i=0$, for $i=1$ if
moreover $F'$ is endowed with a group structure, and for any $i$ if
$F'$ is moreover commutative (\emph{criterion of Artin-Mazur}). Accordingly,
a functor\pspage{34}
\[f : A \to A'\]
between small categories (or categories which are essentially small,
namely equivalent to small categories) is called a weak equivalence,
if{f} the corresponding morphism of topoi
\[ f\uphat : \Ahat \to {A'}\uphat \]
is a weak equivalence.

Coming back to test categories $A$, which allow us to construct the
corresponding modelizers $\Ahat$, our point of view here is
rather that \emph{the test categories are each just as good as the others},
and $\Simplex$ just as good as $\Square$ or any of the $T_{p,q}$ and not
any better! Maybe it's the one though which turns out the most
economical for computational use, the nerve functor $\Cat \to
\Simplexhat$ being still the neatest known of all model preserving
embeddings of \Cat\ into categories $\Ahat$ defined by
modelizers\footnote{However $\Globe$}
. Another point, still more important it seems to me, is
that the natural functor
\[\Topoi \to \Pro\Hot\]
defined by the \v Cech-Cartier-Verdier process, and which allows for
another description of weak equivalences of topoi, namely as those
made invertible by this functor, are directly defined via
semi-simplicial structures of simplicial structures (of the type
``nerve of a covering'')\footnote{I am not sure that the passage through semisimplicial formalism, for describing a $\check{C}$ech-Cartier-Verdier functors, is the most efficient way. I am going to come back upon this in a later chapter in volume 2 of Pursuing Stacks (presumably chapter X: Back to topoi) - compare remarks in section 35.}.

Modelizers of the type $\Ahat$, with $A$ a test category, surely
deserve a name -- let's call them \emph{elementary modelizers}, as
they correspond to the case of an ``elementary'' or ``trivial''
algebraic structure species, whose set-theoretical realizations can be
expressed as just \emph{any} functors
\[ A\op \to \Sets, \]
without any exactness condition of any kind; in other words they can
be viewed as just diagrams of sets of a specified type, with specified
commutativity relations. A somewhat more ambitious question maybe is
whether on such a category $M=\Ahat$, namely an elementary
modelizer, there cannot exist any other modelizing structure. In any
case, the one we got is intrinsically determined in terms of $M$,
which is a topos, by the prescription that an arrow $f: a \to b$
within $M$ is a weak equivalence if and only if the corresponding
morphism for the induced topoi $M_{/a}$ and $M_{/b}$ is a weak
equivalence (in terms of the Artin-Mazur criterion of section \ref{sec:17}).

A more crucial question I feel is whether the right adjoint functor
$j_A$ to $i_A$ in \eqref{eq:29.starstar}
(cf.\ p.~\ref{p:33}) is equally
model preserving, whenever $A$ is a test category. This, as we have
seen, is \emph{not} automatic, whenever we have a model preserving
functor between modelizers, whenever this functor admits an adjoint
functor. In more general terms still, let
\[
\begin{tikzcd}
  M \ar[r, bend left, "i"] & M' \ar[l, bend left, "j"] 
\end{tikzcd}\]
be a pair of adjoint functors, with $M,M'$ endowed with a
``saturated''\pspage{35} set of arrows $W,W'$. Then the following
are equivalent:
\begin{enumerate}[label=(\alph*)]
\item\label{it:29.a}
  $i(W)\subset W'$, $j(W')\subset W$, and the two corresponding
  functors
  \[  W^{-1}M \rightleftarrows W'^{-1}M' \]
  are quasi-inverse of each other, the adjunction morphisms between
  them being deduced from the corresponding adjunction morphisms for
  the pair $(i,j)$.
\item\label{it:29.b}
  $W=i^{-1}(W')$, and for every $a'\in M'$, the adjunction morphism
  \[ ij(a') \to a'\]
  is in $W'$.
\item[(b')]\label{it:29.bprime}
  dual to \ref{it:29.b}, with roles of $M$ and $M'$ reversed.
\end{enumerate}
In the situation we are interested in here, $M = \Ahat$ and
$M'=\Cat$, we know already $W=i^{-1}(W')$ by definition, and hence all
have to see is whether for any (small) category $B$, the functor
\begin{equation}
  \label{eq:29.T}
  i_Aj_A(B) \to B \tag{T}
\end{equation}
is a weak equivalence. This alone will imply that not only $i_{W,W'}$,
but equally $j_{W',W}$ is an equivalence, and that the two are
quasi-inverse of each other. (NB\enspace even without assuming beforehand that
$W,W'$ are saturated, \ref{it:29.b} (say) implies \ref{it:29.bprime}
and \ref{it:29.a}, provided we assume on $W'$ the very mild saturation
condition that for composable arrows $u',v'$, if two among
$u',v',v'u'$ are in $W'$, so is the third; if we suppose moreover that
$M'$ is actually saturated namely made up with all arrows made
invertible by $M' \to W'^{-1}M'$, then condition \ref{it:29.b} implies
that $M$ is saturated too -- which ensures that if $(M',W')$ is
modelizing, so is $(M,W)$.)

It is not clear to me whether for every test category $A$, the
stronger condition \eqref{eq:29.T} above is necessarily
satisfied. This condition essentially means that for any homotopy
type, defined in terms of an arbitrary element $B$ in \Cat\ namely a
category $B$, we get a description of this homotopy type by an object
of the elementary modelizer $\Ahat$, by merely taking
$j_A(B)$. This condition seems sufficiently appealing to me, for
reinforcing accordingly the notion of a test-category $A$, and of an
elementary modelizer $\Ahat$, in case it should turn out to be
actually stronger\footnote{We'll see later it \emph{is} stronger.}. Of course, any category equivalent to an elementary
modelizer $\Ahat$ will be equally called by the same name. It
should appear in due course whether this is indeed the better suited
notion. One point in its favor already is that it appears a lot more
concrete.

Another natural question, suggested by the use of simplicial or
cubical multicomplexes, is whether the product of two test categories
is again a test category -- which might furnish us with a way to
compare directly the description of \Hot\ by the associated elementary
modelizers, without\pspage{36} having to make a detour by the
``basic'' modelizer \Cat\ we started with. But here it becomes about
time to try and leave the thin air conjecturing, and find some simple
and concrete characterization of test categories, or possibly some
reinforcement still of that notion, which will imply stability under
the product operation.

Here it is tempting to use semi-simplicial techniques though, by lack
of independent foundations of homotopy theory in terms of the
modelizer \Cat. Thus we may want to take the map of semi-simplicial
sets corresponding to \eqref{eq:29.T} when passing to nerves, and
express that a)\enspace this is a fibration and b)\enspace the fibers
of this fibration are ``contractible'' (or ``aspheric''), which
together will imply that we have a weak equivalence in \Sssets. Or we
may follow the suggestion of Quillen's set-up, working heuristically
in \Cat\ as though we actually know it is a model category, and
expressing that the adjunction morphism in \eqref{eq:29.T}, which is a
functor between categories, is actually a ``fibration'', and that its
fibers are ``contractible'' namely weakly equivalent to a one-point
category. In any case, a minimum amount of technique seems needed
here, to give the necessary clues for pursuing.

\bigbreak
\presectionfill\ondate{14.3.}\par

\hangsection{Starting the ``asphericity game''.}\label{sec:30}%
Since last week when I stopped with my notes, I got involved a bit
with recalling to mind the ``Lego-Teichm\"uller construction game''
for describing in a concrete, ``visual'' way the Teichm\"uller groups
of all possible types and the main relationships between them, which I
had first met with last year. This and other non-mathematical
occupations left little time only for my reflections on homotopy
theory, which I took up mainly last night and today. The focus of
attention was the ``technical point'' of getting a handy
characterization of test categories. The situation I feel is beginning
to clarify somewhat. Last thing I did before reading last weeks' notes
and getting back to the typewriter, was to get rid of a delusion which
I was dragging along more or less from the beginning of these notes,
namely that our basic modelizer \Cat, which we were using as a giving
the most natural definition of \Hot\ in our setting, was a ``model
category'' in the sense of Quillen, more specifically a ``closed model
category'', where the ``weak equivalences'' are the homotopy
equivalences of course, and where cofibrations are just monomorphisms
(namely functors injective on objects and on arrows) -- fibrations
being defined in terms of these by the Serre-Quillen lifting
property. Without even being so demanding, it turns out still that
there is no reasonable structure of a model category on \Cat, having
the correct weak equivalences, and such that the standard ``Kan
inclusions'' of the following two ordered sets\pspage{37}
\[
\begin{tikzcd}[row sep=small,column sep=tiny]
  & b & \\ a\ar[ur]\ar[rr] & & c
\end{tikzcd} \quad\text{and}\quad
\begin{tikzcd}[row sep=small,column sep=tiny]
  & b\ar[dr] & \\ a\ar[rr] & & c
\end{tikzcd}\]
into
\[\Simplex_2 = a \to b \to c\]
be cofibering. Namely, for a category $bC$, to say that it is
``fibering'' (over the final category $e$) with respect to one
or the other monomorphism, means respectively that every arrow in
$\bC$ has a left respectively a right inverse -- the two together mean
that $\bC$ is a groupoid. But groupoids are definitely \emph{not}
sufficient for describing arbitrary homotopy types, they give rise
only to sums of $K(\pi,1)$ spaces -- thus contradicting Quillen's
statement that homotopy types can be described by ``models'' which are
both fibering and cofibering!

The feeling however remains that any elementary modelizer, namely one
defined (up to equivalence) by a test category, should be a closed
model category in Quillen's sense -- I find it hard to believe that
this should be a special feature just of semi-simplicial complexes!

While trying to understand test categories, the notion of
\emph{asphericity} for a morphism of topoi
\[ f:X\to Y\]
came in naturally -- this is a natural variant of the notion of
$n$-acyclicity (concerned with commutative cohomology) which has been
developed in the context of \'etale topology of schemes in SGA~4. It
can be expressed by ``localizing upon $Y$'' the Artin-Mazur condition
that $f$ be a weak equivalence, by demanding that the same remain true
for the induced morphism of topoi
\[ X \times_Y Y' \to Y' \]
for any ``localization morphism'' $Y'\to Y$. In terms of the
categories of sheaves $E,F$ on $X,Y$, $Y'$ can be defined by an object
(equally called $Y'$) of $F$, the category of sheaves on $Y'$ being
$F_{/Y'}$, and the fiber-product $X'$ can be defined likewise by an
object of $E$, namely by $X'=f^*(Y')$, hence the corresponding
category of sheaves is $E_{/f^*(Y')}$. In case the functor $f^*$
associated to $f$ admits a left adjoint $f_!$ (namely if it commutes
to arbitrary inverse limits, not only to finite ones), the category
$E_{/f^*(Y')}$ can be interpreted conveniently as $E_{f_!/Y'}$ (or
simply $E_{/Y'}$ if $f_!$ is implicit), whose objects are pairs
\[ (U,\varphi), \quad U\in\Ob E, \quad \varphi : f_!(U) \to Y',\]
with obvious ``maps'' between such objects. For the time being I am
mainly interested in the case of a morphism of topoi defined by a
functor between categories, which I will denote by the same symbol
$f$:
\[ f: C' \to C \quad\text{defines $f$ or $f\uphat : {C'}\uphat
  \to C\uphat$.}\]
Using\pspage{38} the fact that in the general definition of
asphericity it is enough to take $Y'$ in a family of generators of the
topos $Y$, and using here the generating subcategory $C$ of $C\uphat$,
we get the following criterion for asphericity of $f\uphat$: it
is necessary and sufficient that for every $a \in \Ob C$, the induced
morphism of topoi
\[ {C'}\uphat_{/a} \equeq ({C'_{/a}})\uphat \to C\uphat_{/a}
\simeq (C_/a)\uphat\]
be a weak equivalence, i.e., that the natural functor $C'_{/a}\to
C_{/a}$ be a weak equivalence. But it is immediate that $C_{/a}$ is
``contractible'', i.e., ``aspheric'', namely the ``map'' from $C_{/a}$
to the final category is a weak equivalence (this is true for any
category having a final object). Therefore we get the following

\noindent{\bfseries Criterion of asphericity for a
  functor}\phantomlabel{lem:asphericitycriterion}{asphericity criterion}
$f: C'\to C$ between categories: namely it is n.\ and s.\ that for any
$a\in\Ob C$, $C'_{/a}$ be aspheric.

Let's come back to a category $A$, for which we want to express that
it is a test-category, namely for any (small) category $B$, the natural
functor
\begin{equation}
  \label{eq:30.star}
  i_Aj_A(B) \to B \tag{*}
\end{equation}
is a weak equivalence. One immediately checks that for any $b\in\Ob
B$, the category $i_Aj_A(B)_{/b}$ over $B_{/b}$ is isomorphic to
$i_Aj_A(B_{/b})$. Hence we get the following

\noindent\textbf{Proposition.} For a (small) category $A$, the following are
equivalent:
\begin{enumerate}[label=(\roman*)]
\item\label{it:30.i}
  $A$ is a test category, namely \eqref{eq:30.star} is a weak
  equivalence for any category $B$.
\item\label{it:30.ii}
  \eqref{eq:30.star} is aspheric for any category $B$.
\item\label{it:30.iii}
  $i_Aj_A(B)$ is aspheric for any category $B$ with final element.
\end{enumerate}

This latter condition, which is the most ``concrete'' one so far,
means also that the element
\[ F = j_A(B) = (a \mapsto \Hom(A_{/a}, B)) \in \Ahat\]
is an aspheric element of the topos $\Ahat$, namely the induced
topos $\Ahat_{/F}$ is aspheric (i.e., the category $A_{/F}$ is
aspheric), \emph{whenever} $B$ has a final element.

\hangsection[The end of the thin air conjecturing: a criterion for
  test \dots]{The end of the thin air conjecturing: a criterion for
  test categories.}\label{sec:31}%
For the notion of a test category to make at all sense, we should make
sure in the long last that $\Simplex$ itself, the category of standard
simplices, is indeed a test category. So I finally set out to prove at
least that much, using the few reflexes I have in semi-simplicial
homotopy theory. A proof finally peeled out it seems, giving clues for
handy conditions in the general case, which should be
\emph{sufficient} at least to ensure that $B$ is a test category, but
maybe not quite necessary. I'll try now to get it down explicitly.

Here\pspage{39} are the conditions I got: 
\begin{enumerate}[label=(T~\arabic*)]
\item\label{it:31.T1}
  $A$ is aspheric.
\item\label{it:31.T2}
  For $a,b\in\Ob A$, $A_{/{a\times b}}$ is aspheric i.e. for every $a \in\Ob A$, $A_{/a} \to A$ is aspheric (NB\enspace $a\times
  b$ need not exist in $A$ but it is in any case well defined as an
  element of $\Ahat$).
\item\label{it:31.T3}
  There exists an \emph{aspheric} element $I$ of $\Ahat$, and
  two subobjects $e_0$ and $e_1$ of $I$ which are final elements of
  $\Ahat$, such that $e_0 \sand e_1 (\eqdef e_0 \times_I e_1) =
  \varnothing_{\Ahat}$, the initial or ``empty'' element of
  $\Ahat$.
\end{enumerate}

In case when $A=\Simplex$ (as well as in the cubical analogon
$\Square$), I took $I=\Simplex_1$ which is an element of $A$ itself, and
moreover $A$ has a final element (which is a final element of
$\Ahat$ therefore) $e$, thus $e_0$ and $e_1$ defined by well
defined arrows in $A$ itself, namely $\delta_0$ and $\delta_1$. But it
does not seem that these special features are really relevant. In any
case \emph{intuitively $I$ stands for the unit interval, with
  endpoints $e_0,e_1$}. If $F$ is any element in $\Ahat$, the
standard way for trying to prove it is aspheric would be to prove that
we can find a ``constant map'' $F \to F$, namely one which factors
into
\[ F \to e \xrightarrow c F\quad\text{($e$ the final element of $\Ahat = \hat{A}$)}\]
for suitable $c : e \to F$ or ``section'' of $F$, which be
``homotopic'' to the identity map $F \to F$. When trying to make
explicit the notion of a ``homotopy'' $h$ between two such maps, more
generally between two maps $f_0,f_1 : F \rightrightarrows G$, we hit
of course upon the arrow $h$ in the following diagram, which should
make it commutative
\begin{equation}
  \label{eq:31.D}
  \begin{tikzcd}
    & F \times I \ar[dd, dashed, swap, "h"] & \\
    F \simeq F \times e_0 \ar[ur, hook] \ar[dr, swap, "f_0"] & &
    F \simeq F \times e_1 \ar[ul, hook] \ar[dl, "f_1"] \\
    & G &
  \end{tikzcd}.
  \tag{D}
\end{equation}
This notion of a homotopy is defined in any category $\mathscr A$
where we've got an element $I$ and two subobjects $e_0,e_1$ which are
final objects. Suppose we got such $h$, and we know moreover (for
given $F,G,f_0,f_1,h$) that the two inclusions of $F \times e_0$ and
$F\times e_1$ into $F \times I$ are weak equivalences, and that $f_0$
is a weak equivalence (for a given set of arrows called ``weak
equivalences'', for instance defined in terms of a ``topology'' on
$\mathscr A$, in the present case the canonical topology of the topos
$\Ahat$), then it follows (with the usual ``mild saturation
condition'' on the notion of weak equivalence) that $h$, and hence
$f_1$ are weak equivalences. Coming back to the case
$F=G, f_0=\id_F, f_1=$ ``constant map'' defined by a $c : e \to F$, we
get that this constant map $f_1$ is a weak equivalence. Does this
imply that $F \to e$ is equally a weak equivalence? This is not quite
formal for general $(\mathscr A,W)$, but it is true though in the case
$\mathscr A = \Ahat$ and with the usual meaning of
``weak\pspage{40} equivalence'', in this case it is true indeed
that if we have a situation of inclusion with retraction $E \to F$ and
$F \to E$ ($E$ need not be a final element of $\mathscr A$), such that
the compositum $p : F \to E \to F$ (a projector in $F$) is a weak
equivalence, then so are $E \to F$ and $F \to E$. To check this, we
are reduced to checking the corresponding statement in \Cat, in fact
we can check it in the more general situation with two topoi $E$ and
$F$, using the Artin-Mazur criterion. (We get first that $E\to F$ is a
weak equivalence, and hence by saturation that $F\to E$ is too.)

Thus the assumptions made on $F\in\Ob{\Ahat}$ imply that $F \to
e$ is a weak equivalence, i.e., $A_{/F} \to A$ is a weak equivalence,
and if we assume now that $A$ satisfies \ref{it:31.T1} namely $A$ is
aspheric, so is $F$.

We apply this to the case $F = j_A(B) = (a \mapsto \Hom(A_{/a},B))$,
where $B$ is a category with final element. We have to check (using
\ref{it:31.T1} to \ref{it:31.T3}):
\begin{enumerate}[label=(\alph*)]
\item\label{it:31.a}
  The inclusions of $F\times e_0, F\times e_1$ into $F\times I$
  are weak equivalences (this will be true in fact for any
  $F\in\Ob{\Ahat}$),
\item\label{it:31.b}
  there exists a ``homotopy'' $h$ making commutative the previous
  diagram \eqref{eq:31.D}, where $G=F$, $f_0=\id_F$, and where $f_1:
  F\to F$ is the ``constant map'' defined by the section $c : e\to F$
  of $F$, associating to every $a\in A$ the \emph{constant} functor
  $e_{a,B} : A_{/a}\to B$ with value $e_B$ (a fixed final element of
  $B$), thus $e_{a,B}\in F(a) = \Hom(A_{/a},B)$ (and it is clear that
  his is ``functorial in $a$'').
\end{enumerate}
Then \ref{it:31.a} and \ref{it:31.b} will imply that $F$ is aspheric
-- hence $A$ is a test-category by the criterion \ref{it:30.iii} of
the proposition above.

To check \ref{it:31.b} we do not make use of \ref{it:31.T1} nor
\ref{it:31.T2}, nor of the asphericity of $I$. We have to define a
``map''
\[ h : F \times I \to F, \]
i.e., for every $a \in \Ob A$, a map (functorial for variable $a$)
\[ h(a) : \Hom(A_{/a},B) \times \Hom(a,I) \to \Hom(A_{/a},B)\]
(two of the $\Hom$'s are in \Cat, the other is in $\Ahat$). Thus,
let
\[ f:A_{/a}\to B, \quad u: a\to I,\]
we must define
\[ h(a)(f,u) : A_{/a} \to B \]
a functor from $A_{/a}$ to $B$, depending on the choice of $f$ and
$u$. Now let, for any $u\in\Hom(a,I)$, $u:a\to I$, $a_u$ be defined in
$\Ahat$ as
\[a_u = u^{-1}(e_0) = (a,u)\times_I e_0,\]
viewed as a subobject of $a$, and hence $C_u=A_{/a_u}$ can be viewed as a
subcategory of $C=A_{/a}$, namely the full subcategory of those objects
$x$ over $a$, i.e., arrows $x\to a\in A$, which factor through $a_u$
namely such that the compositum $x\to a\to I$ factors through
$e_0$. This subcategory is clearly\pspage{41} a ``crible'' (or ``sieve'') in
$A_{/a}$, namely for an arrow $y\to x$ in $A_{/a}$, if the target $x$
is in the subcategory, so is the source $y$. This being so, we define
the functor
\[ f' = h(f,u) : A_{/a}=C \to B\]
by the conditions that
\begin{align*}
  f' &| (A_{/a_u} = C_u) = f | C_u\\
  f' &| (C \setminus C_u) = \text{constant functor with value $e_B$.}
\end{align*}
(where $C\setminus C_u$ denotes the obvious \emph{full} subcategory of
$C$, complementary to $C_u$). This defines $f'$ uniquely, on the
objects first, and on the arrows too because the only arrows left in
$C$ where we got still to define $f'$ are arrows $x\to y$ with $x$ in
$C_u$ and $y\in C\setminus C_u$ (because $C_u$ is a crible), but then
$f'(y)=e_B$ and we have no choice for $f'(x)\to f'(y)$! It's trivial
checking that this way we get indeed a \emph{functor} $f':C \to B$,
thus the map $h(a)$ is defined -- and that this map is functorial with
respect to $a$, i.e., comes from a map $h : F\times I \to F$ as we
wanted. The commutativity of \eqref{eq:31.D} is easily checked: for
the left triangle, i.e., that the compositum $F\simeq F\times e_0 \to
F\times I \xrightarrow h F$ is the identity, it comes from the fact
that if $u$ factors through $e_0$, then $C_u=C$ hence $f'=f$; for the
right triangle, it comes from the fact that if $u$ factors through
$e_1$, then $C_u=\varnothing$ (here we use the assumption $e_0 \sand e_1
= \varnothing_{\Ahat}$), and hence $f'$ is the constant functor $C
\to B$ with value $e_B$. This settles \ref{it:31.b}.

We have still to check \ref{it:31.a}, namely that for any
$F\in\Ob{\Ahat}$, the inclusion of the objects $F\times e_i$ into
$F\times I$ are weak equivalences, or what amounts to the same, that
the projection
\[ F\times I \to F\]
is a weak equivalence -- this will be true in fact for any object $I$
of $\Ahat$ which is aspheric. Indeed, we will prove the stronger
result that $F \times I\to F$ is \emph{aspheric}, i.e., that the
functor
\[A_{/F\times I} \to A_{/F} \]
is aspheric, We use for this the criterion of
p.~\ref{p:38}, which here
translates into the condition that for any $a$ in $A$ (such that we
got an $a\to F$, i.e., such that $F(a)\ne\emptyset$, but never mind),
the category $A_{/a\times I}$ is aspheric, i.e., the lement $a\times
I$ of $\Ahat$ is aspheric. Again, \emph{as $I$ is aspheric}, we
are reduced to checking that $a\times I \to I$ is aspheric, which by
the same argument (with $F,I$ replaces by $I,a$) boils down to the
condition that $A_{x\times b}$ is aspheric for any $b$ in $A$. Now
this is just condition \ref{it:31.T2}, we are through.

\bigbreak

\presectionfill\ondate{15.3.}\pspage{42}\par

\hangsection{Provisional program of work.}\label{sec:32}%
I definitely have the feeling to be out of the thin air -- the
conditions \ref{it:31.T1} to \ref{it:31.T3} look to me so elegant and
convincing, that I have no doubts left they are ``the right ones''! A
lot of things come to mind what to do next, I'll have to look at them
one by one though. Let me make a quick provisional planning.
\begin{enumerate}[label=\arabic*)]
\item\label{it:32.1}
  Have a closer look at the conditions \ref{it:31.T1} to
  \ref{it:31.T3}, to see how far they are necessary for $A$ to be a
  test category in the (admittedly provisional) sense I gave to this
  notion last week and yesterday, and to pin down the feeling of these
  being just the right ones.
\item\label{it:32.2}
  Use these conditions for constructing lots of test categories,
  including all the simplicial and cubical types which have been used
  so far.
\item\label{it:32.3}
  Check that these conditions are stable under taking the product
  of two or more test categories, and possibly use this fact for
  comparing the homotopy theories defined by any two such categories.
\item\label{it:32.4}
  Look up (using \ref{it:31.T1} to \ref{it:31.T3}) if an
  elementary modelizer $\Ahat$ is indeed a ``closed model
  category'' in Quillen's sense, and maybe too get a better feeling of
  how far apart \Cat\ is from being a closed model category\footnote{As seen a little later, this tentative description of a suitable class of cofibrations is inadequate (see beginning of section \ref{sec:30}). An adequate description of a closed model structure in $\Cat$ was found by \emph{Thomason} (see comments on this in section \ref{sec:87})}. Visibly
  there \emph{are} some natural constructions in homotopy theory which
  do make sense in \Cat.
\item\label{it:32.5}
  Using the understanding of test categories obtained, come back
  to the question of which categories $\Ahat$ associated to
  algebraic structure species can be viewed as modelizers, and to the
  question of unicity or canonicity of the modelizing structure $W
  \subset \Fl(M)$.
\end{enumerate}

That makes a lot of questions to look at, and the theory of stacks I
set out to sketch seems to be fading ever more into the background! It
is likely though that a better general understanding of the manifold
ways of constructing the category \Hot\ of homotopy types will not be
quite useless, when getting back to the initial program, namely
stacks. Quite possibly too, on my way I will have to remind myself of
and look up, in the present setting, the main structural properties of
\Hot, including exactness properties and generators and
cogenerators. This also reminds me of an intriguing foundational
question since the introduction of derived categories and their
non-commutative analogs, which I believe has never been settled yet,
namely the following:
\begin{enumerate}[label=\arabic*),resume]
\item\label{it:32.6}
  In an attempt to grasp the main natural structures associated to
  derived categories, and Quillen's non-commutative analogons
  including \Hot, try to develop a comprehensive notion of a
  ``triangulated category'', without the known drawbacks of Verdier's
  provisional notion.\pspage{43}
\end{enumerate}

\hangsection[Necessity of conditions T1 to T3, and transcription in
\dots]{Necessity of conditions T1 to T3, and transcription in terms of
  elementary modelizers.}\label{sec:33}%
For the time being, I'll use the word ``test category'' with the
meaning of last week, and refer to categories satisfying the
conditions \ref{it:31.T1} to \ref{it:31.T3} as \emph{strict} test
categories. (NB\enspace $A$ is supposed to be essentially small in any
case.)

First of all, the conditions \ref{it:31.T1} and \ref{it:31.T3} are
\emph{necessary} for $A$ to be a test category. For \ref{it:31.T1}
this is just the ``concrete'' criterion \ref{it:30.iii} of yesterday
(page 38), when $B$ is the final element of \Cat. For \ref{it:31.T3},
we get even a \emph{canonical choice} for $I,e_0,e_1$, namely starting
with the ``universal'' choice in \Cat:
  \\[\baselineskip]%
  \hspace*{3em}$\cst I = \Simplex_1$, \quad%
\parbox[t]{0.7\textwidth}{$\cst e_0$ and $\cst e_1$ the two subobjects of $\Simplex_1$ in \Cat\
  (or in \Ord) corresponding to the two unique sections of $\cst I$
  over the final element $\cst e$ of \Cat,}
  \\[\baselineskip]%
i.e., in terms of $\cst I$ as an ordered set $\{0\} \to \{1\}$, $\cst
e_0$ and $\cst e_1$ are just the two subobjects defined by the two
vertices $\{0\}$ and $\{1\}$ (viewed as defining two one-point ordered
subsets of $\cst I=\Simplex_1$). We now apply $j_A$ to get
\[ I = j_A(\cst I) = j_A(\Simplex_1), \quad e_i = j_A(\cst e_i) \quad
\text{for $i=0,1$.} \]
As $\cst I\to \cst e$ is a weak equivalence so is $I\to e_{\Ahat}
= j_A(\cst e)$, and hence (as $e_{\Ahat}$ is aspheric by
\ref{it:31.T1}) $I$ is aspheric. As $\cst e_0 \sand \cst e_1 =
\varnothing_\Cat$ and $j_A$ commutes with inverse limits, and with sums
(and in particular transforms the initial element $\varnothing$ of \Cat\
in the initial element $\varnothing$ of $\Ahat$), it follows that
$e_0 \sand e_1 = \varnothing$.

It's worthwhile having a look at what this object just constructed is
like. For this end, let's note first that for any category $C$, we
have a canonical bijection, functorial in $C$
\[ \Hom(C , \Simplex_1) \tosim \Crib(C) \simeq \text{set of all
  subobjects of $e_{C\uphat}$,}\]
by associating to any functor $f : C \to \Simplex_1$ the full
subcategory $f^{-1}(\cst e_0)$ of $C$, which clearly is a ``crible'' (= sieve)
in $C$. Thus we get, for $a\in\Ob A$
\[ I(a) \simeq \Crib(A_{/a}) \quad \text{(subobjects of $a$ in
  $\Ahat$).}\]
More generally, we deduce from this, for any $F\in\Ob{\Ahat}$:
\[ \Hom(F,I) = \Gamma(I_F = I\times F/F) = \text{set of all subobjects
  of $F$,} \]
in other words, the object $I$ is defined intrinsically in the topos
$\Ahat$ (up to unique isomorphism) as the ``\emph{Lawvere
  element}'' representing the functor $F \mapsto$ set of subobjects of
$F$, generalizing the functor $F \mapsto \mathfrak P(F)$ in the
category of sets (namely sheaves over the one-point topos),
represented by the two-point set $\{0,1\}$. The condition $e_0\sand e_1
=\varnothing_\scrA$ is automatic in any topos \scrA, provided \scrA\ is
not the ``empty topos'', namely corresponding to a category of sheaves
equivalent to a one-point category.

There\pspage{44} is a strong temptation now to diverge from test
categories, to expand on intrinsic versions of conditions
\ref{it:31.T1} to \ref{it:31.T3} for any topos, and extract a notion
of a (strictly) modelizing topos, generalizing the ``elementary
modelizers'' $\Ahat$ defined by strict test categories $A$. But
it seems more to the point for the time being to look more closely to
the one condition, namely \ref{it:31.T2}, which does not appear so far
as necessary for $A$ being a test category -- and I suspect it is
\emph{not} necessary indeed, as no idea occurred to me how to deduce
it\footnote{It will become clear later that it isn't necessary indeed, (cf. end of section \ref{sec:39})}. (I have no idea of how to make a counterexample though, as I don't
see any other way to check a category $A$ is a test category, except
precisely using yesterday's criterion via \ref{it:31.T1} to
\ref{it:31.T3}.) Still, I want to emphasize about the fact that
\ref{it:31.T2} \emph{is indeed a very natural condition}. In this connection,
it is timely to remember that in the category \Hot\scrcomment{See letter of 23.5.83 in\ \textcite{Kunzer2015}}, finite direct and
inverse limits exist\footnote{This of course is coarsely false - the only types of limits (in the usual technical sense) which exist in $\Hot$, as is well known, are finite products and arbitrary sums. I was confusing with the more sophisticated notion of ``homotopy limits'', which are   altogether a different story!} (and even infinite ones, I guess, but I feel I'll
have to be a little careful about these\ldots). The existence of such
limits, in terms of the description of $\Hot$ as a category of fractions
of \Cat, doesn't seem at all a trivial fact, for the time being I'll
admit it as ``well known'' (from the semisimplicial set-up, say), and
probably come back upon this with some care later. Now if $(M,W)$ is
any modelizer, hence endowed with a localization functor
\[ M \to \Hot,\]
it is surely not irrelevant to ask about which limits this functor
commutes with, and study the case with care. Thus in no practical
example I know of does this functor commute with binary amalgamated
sums or with fibered products without an extra condition on at least
one of the two arrows involved in $M$ -- a condition of the type that
one of these is a monomorphism or a cofibration (for co-products), or
a fibration (for products). However, in all cases known, it seems that
the functor commutes to (finite, say) sums and products. For sums, it
really seems hard to make a sense out of a localizing functor $M\to
\Hot$, namely playing a ``model'' game, without the functor commuting
at the very least to these! In this respect, it is reassuring to
notice that for any category $A$, the associated functors $i_A,j_A$
between $\Ahat$ and \Cat\ do indeed both commute with (arbitrary)
sums -- which of course is trivial anyhow for $i_A$ (commuting to
arbitrary direct limits), and easily checked for $j_A$ (which
apparently does not commute to any other type of direct limits, but of
course commutes with arbitrary inverse limits). Now for products too,
it is current use to look at products of ``models'' for homotopy
types, as models for the product type -- so much so that this fact is
surely tacitly used everywhere, without any feeling of a need to
comment. It seems too that any category $M$ which one has looked at
so far for possible use as a category of ``models'' in one sense or
other\pspage{45} for homotopy types, for instance set-theoretic
realizations of some specified algebraic structure species, or
topological spaces, and the like, do admit arbitrary direct and
inverse limits, and surely sums and products therefore, so that the
question of commutation of the localization functor to these arises
indeed and is felt to be important. Possibly so important even, that
the notion I introduced of a ``modelizer'' should take this into
account, and be strengthened to the effect that the canonical functor
$M\to\Hot$ should commute at least to finite sums and products, and
possibly even to infinite ones (whether the latter will have to be
looked up with care). I'll admit provisionally that in \Hot, finite
sums and products can be described in terms of the corresponding
operations in \Cat, namely that the canonical functor (going with our
very definition of \Hot\ as a localization of \Cat)
\[\Cat \to \Hot\]
commutes with finite (presumably even infinite) sums and
products\footnote{This is indeed easily checked, and will be done in part VIII (volume 2 of Pursuing Stacks).}. This is indeed reasonably, in view of the fact that the
Nerve functor
\[\Cat \to \Simplexhat = \Sssets\]
does commute to sums and products.

In the case of a test category $A$ and the corresponding elementary
modelizer $\Ahat$, the corresponding localizing functor is the
compositum
\[ \Ahat \xrightarrow{i_A} \Cat\to\Hot,\]
which therefore commutes with finite (presumably even infinite) sums
automatically, because $i_A$ does. Commutation with finite
\emph{products} though does not look automatic. The property of
commutation with final elements is OK and is nothing but condition
\ref{it:31.T1}, which we saw is necessary for $A$ to be a test
category. Thus remains the question of commutation with binary
products, which boils down to the following condition, for any two
elements $F$ and $G$ in $\Ahat$:
\[ i_A(F\times G) \to i_A(F) \times i_A(G) \quad\text{should be a weak
  equivalence,}\]
i.e.,
\begin{equation}
  \label{eq:33.star}
  A_{/F\times G} \to A_{/F} \times A_{/G} \quad\text{a weak equivalence.}
  \tag{*}
\end{equation}
This now implies condition \ref{it:31.T2}, as we see taking $F=a$,
$G=b$, in which case the condition \eqref{eq:33.star} just means
asphericity of $a\times b$ in $\Ahat$, namely \ref{it:31.T2}. To
be happy, we have still to show that conversely, \ref{it:31.T2}
implies \eqref{eq:33.star} for any $F,G$. As usual, it implies even
the stronger condition that the functor in question is aspheric, which
by the standard criterion (page 38) just means that the categories
$A_{/a\times b}$ (for $a,b$ in $A$, and such moreover that $F(a)$ and
$G(b)$ non-empty, but never mind) are aspheric.

Everything turns out just perfect -- it seems worthwhile to summarize
it in one theorem:\pspage{46}

\begin{theorem}
Let $A$ be an essentially small category,
and consider the composed functor
\[ m_A : \Ahat \xrightarrow{i_A} \Cat \to \Hot,\]
where $i_A(F)=A_{/F}$ for any $F\in\Ob{\Ahat}$, and where
$\Cat\to\Hot$ is the canonical functor from \Cat\ into the localized
category with respect to weak equivalences. The following conditions
are equivalent:
\begin{enumerate}[label=(\roman*),font=\normalfont]
\item\label{it:33.i}
  The functor $m_A$ commutes with finite products, and is a
  localization functor, i.e., induces an equivalence $W_A^{-1}\Ahat
  \toequ\Hot$, where $W_A$ is the set of weak equivalences in
  $\Ahat$, namely arrows transformed into weak equivalences by
  $i_A$.\footnote{(Added 9.4.1983) Presumably, this condition (i) is too weak as stated, we'll have to add that $j_A = i_A^*$ is equally modelizing, cf. remark 3 in section 65, D).}
\item\label{it:33.ii}
  The functor $i_A$ and its right adjoint $j_A$ define functors
  between $W_A^{-1}\Ahat$ and $\Hot=W^{-1}\Cat$ \textup{(}thus $j_A$
  should transform weak equivalences of \Cat\ into weak equivalences
  of $\Ahat$\textup{)} which are \emph{quasi-inverse} to each other, the
  adjunction morphisms for this pair being deduced from the adjunction
  morphisms for the pair $i_A,j_A$. Moreover, the functor $m_A$
  commutes with finite products.
\item\label{it:33.iii}
  The category $A$ satisfies the conditions \textup{\ref{it:31.T1}} to
  \textup{\ref{it:31.T3}} \textup{(}page 39\textup{)}.
\end{enumerate}
\end{theorem}

I could go on with two or three more equivalent conditions, which
could be expressed intrinsically in terms of the topos $\Ahat$
and make sense (and are equivalent) for any topos, along the lines of
the reflections of p.~\ref{p:43} and of yesterday. But
I'll refrain for the time being!

In the proof of the theorem above, I did not make use of
semi-simplicial techniques nor of any known results about \Hot, with
the only exception of the assumption (a fact, I daresay, but not
proved for the time being in the present framework, without reference
to semi-simplicial theory say) that the canonical functor
$\Cat\to\Hot$ commutes with binary products. We could have avoided
this assumption, by slightly changing the statement of the theorem,
the condition that $m_A$ commute with finite products in \ref{it:33.i}
and \ref{it:33.ii} being replaced by the assumption that $i_A$ commute
with finite products ``up to weak equivalence'', as made explicit in
\eqref{eq:33.star} above for the case of binary products (which are
enough of course).

I suspect that the notion of a test category in the initial, wider
sense will be of no use any longer, and therefore I will reserve this
name\pspage{47} to the \emph{strict} case henceforth, namely to the
case of categories satisfying the equivalent conditions of the theorem
above\footnote{Finally, it turns out later (section 39) that it is better to keep the notion of a ``strict test category'' as we had it before (implying T2), whereas the notion of a ``test category'' will be slightly stronger than before. For an overall review of the test notions (using the usual notion of weak equivalence in the basic modelizer $\Cat$), see section 44)}. Accordingly, I'll call ``\emph{elementary modelizer}'' any
category \scrA\ equivalent to a category $\Ahat$, with $A$ a test
category. Such a category will be always considered as a modelizer, of
course, with the usual notion of weak equivalence $W \subset
\Fl(\scrA)$, namely of a ``map'' $F \to G$ in \scrA\ such that the
corresponding morphism for the induced topoi is a weak
equivalence. The category \scrA\ is an elementary modelizer,
therefore, if{f} it satisfies the following conditions:
\begin{enumerate}[label=\alph*)]
\item\label{it:33.a}
  \scrA\ is equivalent to a category $\Ahat$ with $A$ small,
  which amount to saying that \scrA\ is a topos, and has sufficiently
  many ``essential points'', namely ``points'' such that the
  corresponding fiber-functor $\scrA \to \Sets$ commutes with
  arbitrary products -- i.e., there exists a conservative family of
  functors $\scrA\to\Sets$ which commute to arbitrary direct and
  inverse limits. (Cf.\ SGA~4 IV~7.5.)
\item\label{it:33.b}
  The pair $(\scrA,W)$ is modelizing (where $W \subset \Fl(\scrA)$ is
  the set of weak equivalences), i.e., the category $W^{-1}\scrA$ is
  equivalent to \Hot.
\item\label{it:33.c}
  The canonical functor $\scrA\to\Hot$ (or, equivalently, $\scrA\to
  W^{-1}\scrA$) commutes with finite products (or, equivalently, with
  binary products -- that it commutes with final elements follows
  already from \ref{it:33.a}, \ref{it:33.b}).
\end{enumerate}
As I felt insistently since yesterday, there is a very pretty notion
of a ``modelizing topos'' generalizing the notion of an elementary
modelizer, where \scrA\ is a topos but not necessarily equivalent to
one of the type $\Ahat$ -- but where however the aspheric (=
``contractible'') objects form a generating family (which generalizes
the condition $\scrA\simeq\Ahat$, which means that the
$0$-connected projective elements of \scrA\ form a generating
family). I'll come back upon this notion later I guess -- it is not
the most urgent thing for the time being\ldots

\hangsection{Examples of test categories.}\label{sec:34}%
I will now exploit the handy criterion \ref{it:31.T1} to
\ref{it:31.T3} for test categories, for constructing lots of such. In
all cases I have in mind at present, the verification of
\ref{it:31.T1} and \ref{it:31.T3} is obvious, they are consequences
indeed of the following stronger conditions:
\begin{enumerate}[label=(T~\arabic*')]
\item\label{it:34.T1prime}
  $A$ has a final element $e_A$.
  \addtocounter{enumi}{1}
\item\label{it:34.T3prime}
  There exists an element $I = I_A$ in $A$, and two sections
  $e_0,e_1$ of $I$ over $e=e_A$
  \[ \delta_0,\delta_1 : e \to I,\]
  such that the corresponding subobjects $e_0,e_1$ of $I$ satisfy
  \[ e_0 \sand e_1 = \varnothing_{\Ahat},\]
  namely that for any $a\in\Ob A$, if $p_a : a \to e$ is the
  projection, we have
  \[ \delta_0 p_a \ne \delta_1p_a.\]
\end{enumerate}

Of\pspage{48} course, the element $I$ (playing the role of unit interval with
endpoints $e_0,e_1$) is by no means unique, for instance it can be
replaced by any cartesian power $I^n$ ($n\ge1$), provided it is in $A$
-- or in the case of $A = \Simplex$ we can take for $I$ any $\Simplex_n$
($n\ge1$) and for $\delta_0,\delta_1$ any two \emph{distinct} maps
from $e=\Simplex_0$ to $\Simplex_n$, instead of the usual choice
$\Simplex_1$ for $I$, and corresponding $\delta_0$ and $\delta_1$. This
high degree of arbitrariness in the choice of $I$ should be no
surprise, this was already one striking feature in Quillen's theory of
model categories.

There remains the asphericity condition \ref{it:31.T2} for the
categories $A_{/a\times b}$, with $a,b\in\Ob A$, which is somewhat
subtler. There is one very evident way though to ensure this, namely
assuming
\begin{enumerate}[label=(T~\arabic*')]
  \addtocounter{enumi}{1}
\item\label{it:34.T2prime}
  If $a,b$ are in $A$, so is $a\times b$, i.e., $A$ is stable under
  binary products.
\end{enumerate}

In other words, putting together \ref{it:34.T1prime} and
\ref{it:34.T2prime}, we may take categories (essentially small)
\emph{stable under finite products}. When such a category satisfies
the mild extra condition \ref{it:34.T3prime} above, it is a test
category! This is already an impressive bunch of test categories. For
instance, take any category $C$ with finite products and for which
there exist $I,e_0,e_1$ as in \ref{it:34.T3prime} -- never mind
whether $C$ is essentially small, for instance any ``non-empty'' topos
will do (taking for $I$ the Lawvere element for instance, if no
simpler choice comes to mind). Take any subcategory $A$ (full or not)
stable under finite products and containing $I$ and
$\Simplex_0,\Simplex_1$ and where the $C$-products are also $A$-products
-- for instance this is OK if $A$ is full. Then $A$ is a test
category.

The simplest choice here, the smallest in any case, it to take the
subcategory made up with all cartesian powers $I^n$ ($n\ge0$). We may
take the full subcategory made up with these elements, but instead we
may take, still more economically, only the arrows
\[ I^n \to I^m\]
where $m$ components $I^n\to I$ are each, either a projection $\pr_i$
($1\le i\le n$) or of the form $\delta_i(p_{I^n})$ ($i\in\{0,1\}$).
The category thus obtained, up to unique isomorphism, does not depend
on the choice of $(C,I,e_0,e_1)$, visibly -- it is, in a sense,
\emph{the} smallest test category satisfying the stronger conditions
\ref{it:34.T1prime} to \ref{it:34.T3prime}. Its elements can be
visualized as ``standard cubes'': One convenient way to do so is to
take $C = \Ord$ (category of all ordered sets), $I = \Simplex_1 =
(\{0\}\to\{1\})$, $e_0$ and $e_1$ as usual, thus $A$ can be
interpreted as a subcategory of \Ord, but this embedding is not full
(there are maps of $I^2$ to $I$ in \Ord\ which are not in $A$, i.e.,
do not ``respect the cubical structure''). We may also take
$C=\Spaces$, $I=$ unit interval $\subset\bR$, $e_0$ and $e_1$ defined
by the endpoints $0$ and $1$,\pspage{49} thus the elements of $A$
are interpreted as the standard cubes $I^n$ in $\bR^n$; the allowable
maps between them $I^n\to I^m$ are those whose components $I^n\to I$
are either constant with value $\in\{0,1\}$, or one of the projections
$\pr_i$ ($1\le i\le n$). We may denote this category of standard cubes
by $\square$, but recall that the cubes here have symmetry operations,
they give rise to the notion of ``(fully) cubical complexes''
$(K_n)_{n\ge0}$, in contrast with what might be called ``semicubical
complexes'' without symmetry operations on the $K_n$'s, in analogy
with the case of simplicial complexes, where there are likewise
variants ``full'' and ``semi''. To make the distinction, we better
call the corresponding test category (with symmetry operations, hence
with more arrows) $\widetilde\square$ rather than $\square$, and
accordingly for $\widetilde\Simplex$, $\Simplex$. On closer inspection, it
seems to me that even apart the symmetry operations, the maps between
standard cubes allowed here are rather plethoric -- thus we admit
diagonal maps such as $I\to I\times I$, which is I guess highly
unusual in the cubical game. It is forced upon us though if we insist
on a test category stable under finite products, for easier checking
of condition \ref{it:31.T2}.

A less plethoric looking choice really is
\[ A = \text{\emph{non-empty} finite sets, with arbitrary maps}\]
or, equivalently, the full subcategory $\widetilde\Simplex$ formed by
the standard finite sets $\bN\sand{[0,n]}$, which we may call
$\widetilde\Simplex_n$ in contrast to the $\Simplex_n$'s (viewed as being
endowed with their natural total order, and with correspondingly more
restricted maps). This gives rise to the category $\widetilde\Simplex$
of ``(fully) simplicial complexes (or sets)''. The elements of $A$, or
of $\widetilde\Simplex$, can be interpreted in the well-known way as
\emph{affine simplices}, and simplicial maps between these.

It is about time now to come to the categories $\Simplex$ and $\square$
and check they are test categories, although they definitely do not
satisfy \ref{it:34.T2prime}, thus we are left with checking the
somewhat delicate condition \ref{it:31.T2}. It is tempting to dismiss
the question by saying that it is ``well-known'' that the categories
\[\Simplex_{/\Simplex_n\times\Simplex_m}\quad(n,m\in\bN)\]
are aspheric -- but this I feel would be kind of cheating. Maybe the
very intuitive homotopy argument already used yesterday (section 31)
will do. In general terms, under the assumptions \ref{it:31.T1},
\ref{it:31.T3}, we found a sufficient criterion of asphericity for an
element $F$ of $\Ahat$, which we may want to apply to the case
$F=a\times b$, with $a$ and $b$ in $A$.\pspage{50}

Now let's say that $e_{\Ahat}$ is a \emph{deformation retract of
  $F$} ($F$ any element in $\Ahat$) if the identity map of $F$ is
homotopic (with respect to $I$) to a ``constant'' map of $F$ into
itself. It is purely formal, using the diagonal map $I\to I\times I$
in $\Ahat$, that if (more generally) $F_0$ is a deformation
retract of $F$, and $G_0$ a deformation retract of $G$, then
$F_0\times G_0$ is a deformation retract of $F\times G$. Thus, a
\emph{sufficient} condition for \ref{it:31.T2} to hold is the
following ``homotopy-test axiom'':

\noindent\parbox[t]{0.3\textwidth}{\textbf{Condition}
  \namedlabel{cond:TH}{(T~H)}:\\[\baselineskip]%
  \hspace*{3em}(1)\phantomlabel{cond:TH1}{(T~H~1)}}%
\parbox[t]{0.7\textwidth}{For any $a\in A$, $e=e_{\Ahat}$ is a
  deformation retract of $a$ with respect to $(I,e_0,e_1)$, namely
  there exists a section $c_a : e\to a$ (hence a constant map $u_a =
  p_ac_a : a\to a$) and a homotopy $h_a:a\times I\to a$ from $\id_a$
  to $u_a$, i.e., a map $h_a$ such that
  \[h_a(\id_a \times \delta_0) = \id_a , \quad
  h_a(\id_a\times\delta_1) = u_a.\]}

We have to be cautious though, it occurs to me now, not to make a
``vicious circle'', as the homotopy argument used yesterday for
proving asphericity of $f$ makes use of the fact that $F\times I\to F$
is a weak equivalence. To check that $F\times I\to F$ is a weak
equivalence, and even aspheric, we have seen though that it amounts to
the same to prove $a\times I\to a$ is a weak equivalence for any $a$
in $A$, i.e., $a\times I$ is aspheric for any $a$ in $A$. This, then,
was seen to be a consequence of the assumption that all elements
$a\times b$ are aspheric (for $a,b$ in $A$) -- but this latter fact is
now what we want to prove! Thus it can't be helped, we have to
complement the condition \ref{cond:TH} above (which I will call
\ref{cond:TH1}, by the extra condition\phantomlabel{cond:TH2}{(T~H~2)}
\begin{equation}
  \label{eq:34.2}
  \text{For any $a$ in $A$, $a\times I$ is aspheric,}
  \tag{2}
\end{equation}
which, in case $I$ itself is in $A$, appears as just a particular case
of \ref{it:31.T2}, which now, it seems, we have to check directly some
way or other.

\bigbreak

\presectionfill\ondate{19.3.}\par

\hangsection[The notion of a modelizing topos. Need for revising the
\dots]{The notion of a modelizing topos. Need for revising the
  \v Cech-Verdier-Artin-Mazur construction.}\label{sec:35}%
I really feel I should not wait any longer with the digression on
modelizing topoi that keeps creeping into my mind, and which I keep
trying to dismiss as not to the point or not urgent or what not! At
least this will fix some terminology, and put the notion of an
elementary modelizer into perspective, in terms of a wider class of
topoi.

A topos $X$ is called \emph{aspheric} if the canonical ``map'' from
$X$ to the ``final topos'' (corresponding to a one-point space, and
whose category of sheaves is the category \Sets) is aspheric. It is
equivalent to say that $X$ is $0$-connected (namely non-empty, and not
decomposable non-trivially into a direct sum of two topoi), and that
for any constant sheaf of groups $G$ on\pspage{51} $X$, the
cohomology group $\mathrm H^i(X,G)$ ($i\ge1$) is trivial whenever defined,
namely for $i\ge1$ if $G$ commutative, and $i=1$ if $G$ is not
supposed commutative. The $0$-connectivity is readily translated into
a corresponding property of the final element in the category of
sheaves $\Sh(X)=\scrA$ on $X$. An element $U$ of \scrA\ is called
\emph{aspheric} if the induced topos $X/U$ or $\scrA_{/U}$ is
aspheric. An arrow $f : U' \to U$ in \scrA\ is called a \emph{weak
  equivalence} resp.\ \emph{aspheric}, if the corresponding ``map'' or
morphism of the induced topoi
\[ X/U'\to X/U \quad\text{or}\quad \scrA_{/U'}\to\scrA_{/U} \]
has the corresponding property. Of course, in the notations it is
often convenient simply to identify objects $U$ of \scrA\ with the
induced topos, and accordingly for arrows. Clearly, if $f$ is aspheric,
it is a weak equivalence, more specifically, $f$ is aspheric if{f} it
is ``universally a weak equivalence'', namely if{f} for any base
change $V\to U$ in \scrA, the corresponding map $g: V'=V\times_U U'\to
V$ is a weak equivalence. Moreover, it is sufficient to check this
property when $V$ is a member of a given generating family of \scrA.

An interesting special case is the one when \scrA\ admits a family of
aspheric generators, or equivalently a generating full subcategory $A$
whose objects are aspheric in \scrA; we will say in this case that the
topos $X$ (or \scrA) is ``\emph{locally aspheric}''. This condition is
satisfied for the most common topological spaces (it suffices that
each point admit a fundamental system of contractible neighborhoods),
as well as for the topoi $\Ahat$ associated to (essentially
small) categories $A$ (for such a topos, $A$ itself is such a
generating full category made up with aspheric elements of \scrA). In
the case of a locally aspheric topos endowed with a generating full
subcategory $A$ as above, a map $f:U'\to U$ is aspheric if{f} for any
$a\in\Ob A$, and map $a \to U$, the object $U' \times_U a$ is aspheric
(because for a map $g: V'\to V$ with $V$ aspheric, $g$ is a weak
equivalence if{f} $V'$ is aspheric, which is a particular case of the
corresponding statement valid for any map of topoi).

Using this criterion, we get readily the following:

\begin{propositionnum}\label{prop:35.1}
  Let \scrA\ be a locally aspheric topos, $A$ a generating full
  subcategory made up with aspheric elements. The following conditions
  on \scrA\ are equivalent:
  \begin{description}[leftmargin=2em]
  \item[\namedlabel{it:35.a}{a)}]
    Any aspheric element of \scrA\ is aspheric over $e$ (the final
    object).
  \item[\namedlabel{it:35.b}{b)}]
    The product of any two aspheric elements of \scrA\ is again
    aspheric.
  \item[\namedlabel{it:35.aprime}{a')}, \namedlabel{it:35.bprime}{b')}]
    as {\normalfont\ref{it:35.a}}, {\normalfont\ref{it:35.b}} but with elements restricted to be
    in $\Ob A$.
  \end{description}
\end{propositionnum}

These conditions, in case $\scrA = \Ahat$ are nothing but
\ref{it:31.T2} for $A$. \emph{It turns out they imply already}
\ref{it:31.T1}. More generally, let's call a topos\pspage{52}
\emph{totally aspheric} if is locally aspheric, and satisfies more the
equivalent conditions above. It then turns out:
\begin{corollary}
  Any totally aspheric topos is aspheric.
\end{corollary}

This follows readily from the definition of asphericity, and the
\v{C}ech computation of cohomology of $X$ in terms of a family
$(a_i)_{i\in I}$ covering $e$, with the $a_i$ in $A$. As the mutual
products and multiple products of the $a_i$'s are all aspheric (a
fortiori acyclic for any constant coefficients), the \v Cech
calculation is valid (including of course in the case of
non-commutative $\mathrm H^1$) and yields the desired result.

The condition for a topos to be totally aspheric, in contrast to
\emph{local} asphericity, is an exceedingly strong one. If for
instance the topos is defined by a topological space, which I'll
denote by $X$, then it is seen immediately that any two non-empty open
subsets of $X$ have a non-empty intersection, in other words $X$ is an
\emph{irreducible space} (i.e., not empty and not the union of two
closed subsets distinct from $X$). It is well known on the other hand
that any irreducible space is aspheric, and clearly any non-empty open
subset of an irreducible space is again irreducible. Therefore:
\setcounter{corollarynum}{1}
\begin{corollarynum}\label{cor:35.1.2}
  A topological space is totally aspheric \textup{(}i.e., the
  corresponding topos is t.a.\textup{)} if{f} it is irreducible. In
  case $X$ is Hausdorff, means also that $X$ is a one-point space.
\end{corollarynum}

Let's now translate \ref{it:31.T3} in the context of general topoi. We
get:

\begin{propositionnum}\label{prop:35.2}
  Let $X$ be a topos. Then the following two conditions are
  equivalent:
  \begin{enumerate}[label=\alph*),font=\normalfont]
  \item\label{it:35.P2.a}
    The ``Lawvere element'' $L_X$ of $\scrA = \Sh(X)$ is aspheric
    over the final object $e_X=e$.
  \item\label{it:35.P2.b}
    There exists an object $I$ in \scrA\ which is aspheric over
    $e$, and two sections $\delta_0,\delta_1$ of $I$ \textup(over
    $e$\textup), such that $\Ker(\delta_0,\delta_1)=\varnothing$, i.e.,
    such that $e_0 \sand e_1 = \varnothing$, where $e_0,e_1$ are the
    subobjects of $I$ defined by $\delta_0,\delta_1$.
  \end{enumerate}
\end{propositionnum}

I recall that the Lawvere element is the one which represents the
functor
\[ F \mapsto \text{set of all subobjects of $F$}\]
on \scrA, it is endowed with two sections over $e$, corresponding to
the two ``trivial'' subobjects $\varnothing, e$ of $e$, and the kernel
of this pair of sections is clearly $\varnothing$, thus \ref{it:35.P2.a}
$\Rightarrow$ \ref{it:35.P2.b}. Conversely, by the homotopy argument
already used (section \ref{sec:31}),
we readily get that \ref{it:35.P2.b}
implies that\pspage{53} the projection $L\to e$ is a weak
equivalence. As the assumption \ref{it:35.P2.b} is stable by
localization, this shows that for any $U$ in \scrA, $L_U \to U$ is
equally a weak equivalence, and hence $L$ is aspheric over $e$.

In proposition~\ref{prop:35.2} we did not make any assumption of the
topos $X$ being locally aspheric, let alone totally aspheric. The
property of total asphericity, and the one of prop.~\ref{prop:35.2} about
the existence of a handy substitute for the unit interval, seem to be
independent of each other. The case of topological spaces is
instructive in this respect. Namely:
\setcounter{corollarynum}{0}
\begin{corollarynum}\label{cor:35.2.1}
  Assume $X$ is a topological space, and $X$ is $0$-connected. Then
  the Lawvere element $L_X$ in $\scrA=\Sh(X)$ is $0$-connected, except
  exactly when $X$ is irreducible, in which case $L_X$ decomposes into
  two connected components.
\end{corollarynum}

Indeed, $L$ is disconnected if{f} there is a direct summand $L'$ of
$L$ containing $e_0$ and $\ne L$, or equivalently, if for any open $U
\subset X$, and $U\to L$, namely an open subset $U_0$ of $U$, we can
associate a \emph{direct summand} $U'=L'(U,U_0)$ of $U$ containing
$U_0$, functorially for variable $U$, and such that
\[ L'(U,\emptyset) = \emptyset.\]
The functoriality in $U$ means that for $V$ open in $U$ and $V_0=U_0
\sand V$, we get
\begin{equation}
  \label{eq:35.star}
  L'(V,V_0) = V \sand L'(U,U_0).
  \tag{*} 
\end{equation}
This implies that $L'(U,U_0)$ is known when we know $L'(X,U_0)$,
namely
\[ L'(U,U_0) = U \sand L'(X,U_0).\]
As $L'(X,U_0)$ is a direct summand of $X$ containing $U_0$, and $X$ is
connected, we see that
\[ \text{$L'(X,U_0) = X$ if $U_0\ne\emptyset$,} \quad
L'(X,\emptyset)=\emptyset; \]
and hence
\begin{equation}
  \label{eq:35.starstar}
  \text{$L'(U,U_0) = U$ if $U_0\ne\emptyset$,} \quad
  \text{$L'(U,U_0) = \emptyset$ if $U_0=\emptyset$,}
  \tag{**}
\end{equation}
thus the association $(U,U_0) \mapsto L'(U,U_0)$ is uniquely
determined, and it remains to see whether the association
\eqref{eq:35.starstar} is indeed functorial, i.e., satisfies
\eqref{eq:35.star} for any open $V \subset U$, with $V_0=V\sand
U_0$. It is OK if $V_0\ne\emptyset$, or if $V_0=\emptyset$ and
$U_0=\emptyset$. If $V_0=\emptyset$ and $U_0\ne\emptyset$, it is OK
if{f} $V$ is empty, in other words, if both open subsets $U_0$ and $V$
of $U$ are non-empty, so must be their intersection $V_0$. But this
means that $X$ is irreducible.

\begin{corollarynum}\label{cor:35.2.2}
  A topological space $X$ cannot be totally aspheric (i.e.,
  irreducible) and satisfy the condition of prop.~\textup{\ref{prop:35.2}}.
\end{corollarynum}

Because for a 0-connected topos, this latter condition implies that
$L_X$ is equally 0-connected, which contradicts
corollary~\ref{cor:35.2.1}.

Thus, surely total asphericity for a topos \emph{does not} imply
the\pspage{54} condition of prop.~\ref{prop:35.2} -- the simplest
counterexample is the final topos, corresponding to one-point
spaces (or any irreducible space). But I don't expect either that condition of
prop.~\ref{prop:35.2}, even for a locally aspheric topos $X$, and even
granting $X$ is aspheric moreover, implies that $X$ is totally
aspheric. The positive result in cor.~\ref{cor:35.2.1} gives a hint
that the condition of prop.~\ref{prop:35.2} may be satisfied for locally
aspheric topological spaces which are sufficiently far away from those
awful non-separated spaces (including the irreducible ones) of the
algebraic geometry freaks -- possibly even provided only $X$ is
locally Hausdorff.

Now the conjunction of the conditions of prop.~\ref{prop:35.1} and
\ref{prop:35.2}, namely total asphericity plus existence of a
monomorphism of $e \amalg e$ into an aspheric element, seems an
extraordinarily strong assumption -- so strong indeed that no
topological space whatever can satisfy it! I feel like calling a topos
satisfying these conditions a ``\emph{modelizing topos}''. For a topos
of the type $\Ahat$, this just means it is an ``elementary
modelizer'', i.e., $A$ is a test category, which are notions which, I
feel, are about to be pretty well understood. I have not such feeling
yet for the more general situation though. For instance, there is
another kind of property of a topos, which from the very start of our
model story, one would have thought of pinpointing by the name of a
``modelizing topos''. Namely, when taking $W \subset\Fl{\scrA}$ to be
weak equivalences in the sense defined earlier, we demand that
$(\scrA,W)$ should be a modelizer in the general sense
(cf.\ p.~\ref{p:31}),
with some exactness reinforcement
(cf.\ p.~\ref{p:45}), namely that
$w^{-1}\scrA$ should be equivalent to \Hot, and $W$ saturated, namely
equal to the set of maps made invertible by the localization functor
-- and moreover I guess that this functor should commute with (at
least finite) sums and products. Of course, we would even expect a
little more, namely that the ``canonical functor''
\[\scrA \to \Hot\]
can be described, by associating to every $F\in\Ob{\scrA}$ the
(pro-)homotopy type associated to it by the Verdier-Artin-Mazur
process. One feels there is a pretty juicy bunch of intimately related
properties for a topos, all connected with the ``homotopy model
yoga'', and which one would like to know about.

This calls to my mind too the question of understanding the
Verdier-Artin-Mazur construction in the present setting, where
homotopy types are thought of in terms of categories as preferential
models, rather than semi-simplicial sets. This will be connected with
the question, which\pspage{55} has been intriguing me lately, of
understanding the ``structure'' of an arbitrary morphism of topoi
\[ X\to\Top(B),\]
where $\Top(B)$ is the topos associated to an arbitrary (essentially
small) category $B$ -- a situation, it seems, which generalizes the
situation of a ``fibered topos'' over the indexing category $B$, in
the sense of SGA~4 IV. The construction of Verdier in SGA~4 V
(appendix) corresponds apparently to the case of categories $B$ of the
type $\Simplex_{/F}$ with $F\in\Ob{\Simplexhat}$, as brought near by
the standard \v Cech procedure. When however the topos $X$ is already
pretty near a topos of the type $\Top(B)$, for instance if it is such
a topos, to describe its homotopy type in terms of a huge inverse
system of semi-simplicial objects of $B\uphat$, rather than just
take $B$ and keep it as it is, seems technically somewhat prohibitive,
at least in the present set-up (and with a distance of twenty years!).

\hangsection[Characterization of a particular type of test functors
  with \dots]{Characterization of a particular type of test functors
  with values in \texorpdfstring{\Cat}{(Cat)}.}\label{sec:36}%
Now back to test categories, and more specifically to her majesty
$\Simplex$ -- we still have to check it is a test category indeed! The
method suggested last Tuesday
(p.~\ref{p:50}) does work indeed, without any
reference to ``well-known'' facts from the semi-simplicial theory. To
check \ref{cond:TH1}, we take of course $I=\Simplex_1$ and $e_0,e_1$ as
usual, we then have for every $n\in\bN$ to construct a homotopy
\[ h_n : \Simplex_n \times \Simplex_1 \to \Simplex_n,\]
where the product and the arrow can be interpreted as one in \Cat\ or
even in \Ord\ of course (of which $\Simplex$ is a full subcategory), and
there is a unique such $h$ if we take for $u_a$ the ``constant''
endomorphism of $\Simplex_n$ whose value is the \emph{last} element of
$\Simplex_n$. There remains only \ref{cond:TH2} -- namely to prove that
the objects $\Simplex_n\times\Simplex_1$ in $\Simplex$ are aspheric
($n\ge0$). For this we cover $\Simplex_n'=\Simplex_n\times\Simplex_1$ by
maximal representable subobjects, namely maximal flags of the ordered
product set (these are I guess what are called the ``shuffles'' in
semi-simplicial algebra). It then turns out that $\Simplex_n'$ can be
obtained by successive gluing of flags, the intersection of each flag
we add with the sup or union of the preceding ones being just a
subflag, hence representable and aspheric. Thus we only have to make
use of the ``Mayer-Vietoris'' type easy lemma:
\begin{lemma}
  Let $U,V$ be two subobjects of the final object $e$ of a topos (and covering $e$), assume
  $U,V$ and $U\sand V$ aspheric, then so is $e$, i.e., so is the topos.
\end{lemma}

This finishes the proof of $\Simplex$ being a test category. Of course,
one cannot help thinking of the asphericity of \emph{all} the products
$\Simplex_m\times\Simplex_n$ as just meaning asphericity of the product
ordered set (namely of the corresponding category), which follows for
instance from the fact that it\pspage{56} has a final object. But
we should beware the vicious circle, as we implicitly make the
assumption that for an object $C$ of \Cat\ (at least an object $C$
such as $\Simplex_m\times\Simplex_n$, if $C$ is aspheric, then so is the
corresponding element $\alpha(C)$ in $\Simplexhat$:
\[ \alpha(C) : \Simplex_n \mapsto \Hom(\Simplex_n,C),\]
namely the \emph{nerve} of $C$. We expect of course something stronger
to hold, namely that for any $C$, $\alpha(C)$ has the same homotopy
type as $C$ -- in fact that we have a canonical isomorphism between
the images of both in \Hot. This relationship has still to be
established, as well as the similar statement for $\beta : \Cat \to
\square\uphat$, and accordingly with $\widetilde\Simplex$ and
$\widetilde\square$. For all four categories (where the simplicial
cases are known to be test categories already), we have, together with
the category $A$ (the would-be test category), a functor
\[ A \to \Cat,\quad\text{say $a\mapsto\Simplex_a$}\]
(notation inspired by the case $A=\Simplex$), hence a functor
\[\Cat \xrightarrow\alpha \Ahat, \quad \alpha(C) = (a \mapsto
\Hom_{\Cat}(\Simplex_a, C)),\]
and the question is to compare the homotopy types of $C$ and
$\alpha(C)$, the latter being defined of course as the homotopy type
of $A_{/\alpha(C)} \simeq A_{/C}$, the category of all pairs $(a,p)$
with $a\in\Ob A$ and $p:\Simplex_a\to C$, arrows between pairs
corresponding to \emph{strictly} commutative diagrams in \Cat{}
\[\begin{tikzcd}[row sep=tiny,column sep=small]
  \Simplex_{a'} \ar[rr,"\Simplex_f"] \ar[dr, swap, "p'"] & & \Simplex_a \ar[dl,"p"] \\
  & C &
\end{tikzcd}.\]
\begin{proposition}
  Let the data be as just said, assume moreover that the categories
  $\Simplex_a$ \textup($a\in\Ob A$\textup) have final elements, that $A$
  is aspheric and that $\alpha(\Simplex_1)$ is aspheric over the final
  element of $\Ahat$ \textup(the latter automatic if $A$ is a test
  category, $A\to\Cat$ is fully faithful and $\Simplex_1$ belongs to its
  essential image\textup). Then we can find for every $C$ in \Cat\ a
  "map" in \Cat\ \textup(i.e., a functor\textup)
  \begin{equation}
    \label{eq:36.star}
    \varphi : A_{/C} \to C,
    \tag{*}
  \end{equation}
  functorial in $C$ for variable $C$, and this map being aspheric, and
  a fortiori a weak equivalence \textup(hence induces an isomorphism between
  the homotopy types\textup).
\end{proposition}

This implies that the compositum $C \mapsto A_{/C}$:
\[\Cat \xrightarrow{\alpha} \Ahat \xrightarrow{i_A} \Cat\]
transforms weak equivalences into weak equivalences, and that the
functor deduced from it by passage to the localized categories
$W^{-1}\Cat = \Hot$ is isomorphic to the identity functor of
\Hot.\pspage{57}

To define a functor \eqref{eq:36.star}, functorial in $C$ for variable
$C$, we only have to choose a final element $e_a$ in each
$\Simplex_a$. For an element $(a,p)$ of $A_{/C}$, we define
\[\varphi(a,p) = p(e_a),\]
with evident extension to arrows of $A_{/C}$ (NB\enspace Here it is important
that the $e_a$ be \emph{final} elements of the $\Simplex_a$, initial
elements for a change wouldn't do at all!), thus we get a functor
$\varphi$, and functoriality with respect to $C$ is clear. Using the
standard \hyperref[lem:asphericitycriterion]{criterion of asphericity
  of a functor} $\varphi$ we have to check that the categories
\[ (A_{/C})_{/x}\]
for $x\in\Ob C$ are aspheric, but one checks at once that the category
above is canonically isomorphic to $A_{/C'}$, where $C'=C_{/x}$. Thus
asphericity of $\varphi$ for arbitrary $C$ is equivalent to
asphericity of $A_{/C}$ when $C$ has a final element, or equivalently,
to asphericity of the element $\alpha(C)$ in $\Ahat$ for such
$C$. Now let $e_C$ be the final element of $C$, and consider the
unique homotopy
\[ h : \Simplex_1 \times C \to C\]
between $\id_C$ and the constant functor from $C$ to $C$ with value
$e_C$. Applying the functor $\alpha$, we get a corresponding homotopy
in $\Ahat$
\[ I \times \alpha(C) \to \alpha(C), \quad\text{where
  $I=\alpha(\Simplex_1)$,}\]
between the identity map of $F=\alpha(C)$ and a ``constant'' map $F\to
F$. As we assume that $I$ is aspherical over the final element of
$\Ahat$, hence $I\times F\to F$ is a weak equivalence, it follows
that $F\to e$ is a weak equivalence, hence $F$ is aspheric as $e$ is
supposed to be aspheric, which was to be proved.

If we apply this proposition to $A=\Simplex$ and $n\mapsto\Simplex_n$, we
are still reduced to check just \ref{cond:TH2} (which we did above
looking at shuffles), whereas \ref{cond:TH1} appears superfluous now
-- we get asphericity of the elements $\Simplex_m \times \Simplex_n$ in
$\Ahat$ as a consequence (but the homotopy argument used in the
prop.\ is essentially the same as the one used for proving
\ref{cond:TH1}). This proposition applies equally to the case when
$A=\widetilde\Simplex$ (category of non-ordered simplices
$\widetilde\Simplex_n$), here the natural functor
\[ i : \widetilde\Simplex \to \Cat, \quad a\mapsto\mathfrak P^*(a)\]
is obtained by associated to any finite set $a$ the combinatorial
simplex it defines, embodied by the ordered set $\mathfrak P^*(a)$ of
its facets of all dimensions, which can be identified with the ordered
set of all non-empty subsets of the finite set $a$. As we know already
$\widetilde\Simplex$ is a test category, all that remains to be done is
to check in this case that $I=\alpha(\Simplex_1)$ is aspheric (hence
aspheric over the final object of $\Ahat$). Now this follows
again from a homotopy argument involving the unit segment
substitute\pspage{58} in $\Ahat$. We could formalize it as
follows:
\begin{corollary}\scrcomment{Deleted in original, see Cor. \ref{cor:38.3}}
  In the proposition above, the condition that $\alpha(\Simplex_1)$ be
  aspheric over $e_{\Ahat}$ is a consequence of the following
  assumptions:
  \begin{enumerate}[label=\alph*),font=\normalfont]
  \item\label{it:36.a}
    The condition \textup{\ref{it:31.T3}} is valid in $A$\kern1pt, with the
    stronger assumption that $I$ is aspheric over the final object $e$
    of $\Ahat$, and moreover:
  \item\label{it:36.b}
    for every $a\in\Ob A$ and $u:a\to I$, let $a_u=u^{-1}(e_0)$
    \textup(subobject of $a$ in $\Ahat$\textup), $\Simplex_{a_u} =
    \varinjlim_{\text{$b$ in $A_{/a_u}$}} \Simplex_b$, and $\Simplex_{a,u}
    = \Imm(\Simplex_{a_u} \to \Simplex_a)$. We assume that
    $\Ob{\Simplex_{a,u}}$ is a \emph{crible} in $\Simplex_a$, say $C_u$.
  \end{enumerate}
\end{corollary}

Indeed, interpreting $\alpha(\Simplex_1)$ as the functor $a \mapsto
\Crib{\Simplex_a}$, we define a homotopy
\[ h : I \times F \to F\quad\text{(where $F=\alpha(\Simplex_1)$)}\]
from the identity $\id_F$ to the constant map $c:F\to F$ associating
to every crible in some $\Simplex_a$ the empty crible in the same, by
taking for every $a$ in $A$ the map
\[ (u,C_0) \mapsto C_u \circ C_0.\]
It follows that $F\to e$ is a weak equivalence, and the same argument
in any induced category $A_{/a}$ shows that this is universally so,
i.e., $F\to e$ is aspheric, O.K.

In the case above, $A=\widetilde\Simplex$, we get indeed that for
$u:a\to I$, $a_u$ is empty in $A$ and $\Simplex_{a_u} \to \Simplex_a$ is a
full embedding, turning $\Simplex_{a_u}$ into a crible in $\Simplex_a$,
thus \ref{it:36.b} is satisfied (with the usual choice
$I=\widetilde\Simplex_1$ of course).

\bigbreak

\presectionfill\ondate{21.3.}\par

\hangsection{The ``asphericity story'' told anew -- the ``key
  result''.}\label{sec:37}%
By the end of the notes two days ago, there was the pretty strong
impression of repeating the same argument all over again, in very
similar situations. When I tried to pin down this feeling, the first thing
that occurred to me was that the homotopy equivalence I was after
$A_{/C}\to C$ (in the proposition of
p.~\ref{p:56}), when concerned with a
general functor $A\to\Cat$, and the criterion obtained, was applicable
to the situation I was concerned with at the very start when defined
test categories, namely when looking at the \emph{canonical} functor
$A\to\Cat$ given by $a\mapsto A_{/a}$; and that this gave a criterion
in this case for $A$ to be a test category in the wider sense, which I
had overlooked when peeling out the characterization of test
categories (cf.\ theorem of
p.~\ref{p:46}).Finally it becomes clear that it
is about time to recast from scratch the asphericity story, and (by
one more repetition) tell it anew in a way stripped from its
repetitive features!\pspage{59}

On a more technical level, I got aware too that the last corollary
stated was slightly incorrect, because it is by no means clear that
the crible $C_u$ constructed there is functorial in $u$, i.e.,
corresponds to a map $I\to\alpha(\Simplex_1)$, this in fact has still to
be assumed (and turns out to be satisfied in all cases which turned up
naturally so far and which I looked up). But let's now ``retell the story''!

One key notion visibly in the homotopy technique used, and which needs
a name in the long last, is the notion of a \emph{homotopy interval}
(``segment homotopique'' in French). To be really outspoken about the
very formal nature of this notion and the way it will be used, let's
develop it in any category \scrA\ endowed with a subset
$W\subset\Fl(\scrA)$ of the set of arrows of \scrA, and satisfying the
usual conditions:
\begin{enumerate}[label=\alph*)]
\item\label{it:37.a}
  $W$ contains all isomorphisms of \scrA,
\item\label{it:37.b}
  for two composable arrows $u,v$ in \scrA, if two among $u,v,vu$ are
  in $W$, so is the third, and
\item\label{it:37.c}
  if $i:F_0\to F$, and $r:F\to F_0$ is a left inverse (i.e., a
  retraction), and if $p=ir:F\to F$ is in $W$, so is $r$ (and hence
  $i$ too).
\end{enumerate}

The condition \ref{it:37.c} here is the ingredient slightly stronger
than what we used to call ``mild saturation property'' of $W$, meaning
\ref{it:37.a} and \ref{it:37.b}.

We'll call \emph{homotopy interval in \scrA} (with respect to the
notion of ``weak equivalence'' $W$, which will generally be implicit
and specified by context) a triple $(I,e_0,e_1)$, where $I$ is an
object of \scrA, $e_0$ and $e_1$ two subobjects, such that the
following conditions \ref{cond:HIa} to \ref{cond:HIc} hold:
\begin{enumerate}[label=\alph*)]
\item\label{cond:HIa}
  $e_0$ and $e_1$ are final objects of \scrA,
\end{enumerate}
which implies that \scrA\ has a final object, unique up to unique
isomorphism, which we'll denote by $e_\scrA$ or simply $e$, so that
the data $e_0,e_1$ in $I$ are equivalent to giving two sections
\[ \delta_0,\delta_1 : e \to I\]
of $I$ over $e$. Note that $e_0 \sand e_1 = \Ker(\delta_0,\delta_1)$.
\begin{enumerate}[label=\alph*),resume]
\item\label{cond:HIb}
  $e_0\sand e_1 = \varnothing_\scrA$ (a strict initial element of \scrA),
\item\label{cond:HIc}
  $I\to e$ is ``universally in $W$'' or, as we'll say, is
  \emph{$W$-aspheric} or simply aspheric,
\end{enumerate}
which just means here that $I\to e$ is ``squarable'' and that for any
$F$ in $\Ob\scrA$, $F\times I \to F$ is in $W$, i.e., is a ``weak
equivalence''.

It is clear that if $(I,e_0,e_1)$ is a homotopy interval in \scrA,
then for any $F\in\Ob\scrA$, the ``induced interval'' in $\scrA_{/F}$,
namely $(I\times F,e_0\times F,e_1\times F)$ is a homotopy interval in
$\scrA_{/F}$.\pspage{60}

We now get the (essentially trivial)
\begin{homotopylemma}\label{lem:homotopylemma}
  Let $h: F\times I\to F$ be a ``homotopy'' with respect to
  $(I,e_0,e_1)$ of $\id_F$ with a constant map $c = ir: F\to F$, where
  $r: F\to e$ and $i: e\to F$ is a section of $F$ over $e$. Then $F\to
  e$ is $W$-aspheric.
\end{homotopylemma}

In the proof of this lemma, we make use of
\ref{it:37.a}\ref{it:37.b}\ref{it:37.c} for $W$, but only
\ref{cond:HIa}\ref{cond:HIc} for $(I,e_0,e_1)$, namely we don't even
need $e_0\sand e_1=\varnothing_\scrA$. This is still the case for the
proof of the
\begin{comparisonlemmaforHI}\label{lem:comparisonlemmaforHI}
  Let $L$ be an object of \scrA, squarable over $e$ and endowed with a
  composition law $x \land y$, let $\delta_i^L:e \to L$ be two sections of
  $L$ \textup($i\in\{0,1\}$\textup), which are respectively a left
  unit and zero element for the multiplication, namely the
  corresponding elements $e_0^a,e_1^a$ in any $\Hom(a,L)$ satisfy
  \[ e_0^a \land x = x, \quad e_1^a\land x = e_1^a\]
  for any $x$ in $\Hom(a,L)$. Assume moreover we got a homotopy
  interval $(I,e_0,e_1)$ and a ``map of intervals''
  \[ \varphi: I \to L\]
  \textup(in the sense: compatible with endpoints\textup). Then
  $\id_L$ is homotopic \textup(with respect to $(I,e_0,e_1)$\textup)
  to the constant map $c_L : L\to L$ defined by $\delta_1^L$, and
  hence by the previous \hyperref[lem:homotopylemma]{homotopy lemma},
  $L\to e$ is $W$-aspheric and therefore $L$ endowed with
  $e_0^L,e_1^L$ defined by $\delta_0^L,\delta_1^L$ is itself a
  homotopy interval in \scrA\ \textup(provided, at least, we know that
  $e_0^L\sand e_1^L = \varnothing$\textup).
\end{comparisonlemmaforHI}

The homotopy $h$ is simply given (for $u:a\to I$, $x:a\to L$) by
\[ h(x,u) = \varphi(u) \land x;\]
when $u$ factors through $e_0$, then $\varphi(u) = e_0^a$ and hence
$h(x,u)=x$; if it factors through $e_1$, then $\varphi(u)=e_1^a$ and
we get $h(x,u)=e_1^a$, qed.

\begin{corollary}\label{cor:ofcomparisonlemmaforHI}
  Assume in \scrA\ \textup(endowed with $W$\textup) finite inverse limits
  exist \textup(i.e., final object and fibered products exist\textup),
  and moreover that the presheaf on \scrA
  \[F \mapsto \text{set of all subobjects of $F$}\]
  is representable by an element $L$ of \scrA\ \textup(the ``Lawvere
  element''\textup). Assume moreover \scrA\ has a strict initial element
  $\varnothing_\scrA$, i.e., an initial element and that any map
  $a\to\varnothing_\scrA$ is an isomorphism. Consider the two sections
  of $L$ over $e$, $\delta_0^L$ and $\delta_1^L$, corresponding to the
  full and to the empty subobject of $e$, so that we get visibly
  $e_0^L\sand e_0^L=\varnothing$. Then, for a homotopy interval to exist
  in \scrA, it is necessary and sufficient that $L$ be $W$-aspheric
  over $e$, i.e., that $(L,e_0^L,e_1^L)$ be a homotopy interval.
\end{corollary}

By the comparison lemma, it is enough to show that for any homotopy
interval $(I,e_0,e_1)$ in \scrA, there exists a morphism of
``intervals''\pspage{61} from $I$ into $L$, using the fact of
course that the intersection law on $L$ is a composition law admitting
$\delta_0^L,\delta_1^L$ respectively as unit and as zero element
(still using the fact that the initial object is strict). But the
subobject $e_0$ of $I$, by definition of $L$, can be viewed as the
inverse image of $e_0^L$ by a uniquely defined map $I\to L$. The
induced map $e_1\to L$ corresponds to the induced subobject $e_0 \sand
e_1$ of $e_1$ which by assumption is $\varnothing_\scrA$, and hence
$e_1\to L$ factors through $e_1^L$, qed.

Of course, the case for the time being which mainly interests us is
the one when \scrA\ is a topos (more specifically even, a topos
equivalent to a category $\Ahat$, with $A$ a small category), in
which case it is tacitly understood that $W$ is the set of weak
equivalences in the usual sense. We now come to the key result:
\begin{theorem}\label{thm:keyresult}
  Let $A$ be a small category, and
  \[ i : A \to \Cat\]
  a functor, hence a functor
  \[ i^* : \Cat\to\Ahat, \quad C\mapsto(a\mapsto \Hom(i(a),C)).\]
  Consider the canonical functor $i_A : \Ahat\to\Cat$, $F\mapsto
  A_{/F}$, and the compositum
  \[\Cat \xrightarrow{i^*} \Ahat\xrightarrow{i_A} \Cat, \quad
  C \mapsto A_{/i^*(C)} \eqdef A_{/C}.\]
  Assume that for any $a\in\Ob A$, $i(a)\in\Ob\Cat$ has a final
  element $e_a$ \textup(but we don't demand that for $u: a\to b$,
  $i(u): i(a)\to i(b)$ transforms $e_a$ into $e_b$ nor even into a
  final element of $i(b)$\textup). Consider the canonical functor
  \begin{equation}
    \label{eq:key.star}
    A_{/C}\to C, \quad (a, p:i(a)\to C) \mapsto p(e_a),\tag{*}
  \end{equation}
  which is functorial in $C$, and hence defines a map between functors
  from \Cat\ to \Cat:
  \begin{equation}
    \label{eq:key.starstar}
    i_Ai^* \to \id_\Cat.\tag{**}
  \end{equation}
  \begin{enumerate}[label=\alph*),font=\normalfont]
  \item\label{it:key.a}
    The following conditions are equivalent:
    \begin{enumerate}[label=(\roman*),font=\normalfont]
    \item\label{it:key.a.i}
      For any $C$ in \Cat, \textup{\eqref{eq:key.star}} is aspheric.
    \item\label{it:key.a.ii}
      For any $C$ in \Cat, \textup{\eqref{eq:key.star}} is a weak
      equivalence, i.e., $i_Ai^*$ transforms weak equivalences into
      weak equivalences, and \textup{\eqref{eq:key.starstar}} induces
      an isomorphism of the corresponding functor $\Hot\to\Hot$ with
      the identity functor.
    \item\label{it:key.a.iii}
      The functor $i_Ai^*$ transforms weak equivalences into weak
      equivalences, and the induced functor $\Hot\to\Hot$ transforms
      every object into an isomorphic one, i.e., for any $C$ in \Cat,
      $A_{/C}$ is isomorphic to $C$ in \Hot.
    \item\label{it:key.a.iv}
      For any $C$ with a final element, $A_{/C}$ is aspheric.
    \end{enumerate}
  \item\label{it:key.b}
    The following conditions are equivalent, \emph{and they imply the
      conditions in \textup{\ref{it:key.a}} provided $A$ is
      aspheric}:\pspage{62}
    \begin{enumerate}[label=(\roman*),font=\normalfont]
    \item\label{it:key.b.i}
      For any $C$ in $\Cat$ the functor
      \begin{equation}
        \label{eq:key.starstarstar}
        A_{/C} \to A\times C\quad\text{deduced from
          \textup{\eqref{eq:key.star}} and
          $A_{/C}\xrightarrow{\textup{can}} A$}
        \tag{***}
      \end{equation}
      is aspheric.
    \item\label{it:key.b.ii}
      For any $C$ in \Cat\ with final element and any $a$ in $A$\kern1pt,
      $a\times i^*(C)$ is an aspheric element in $\Ahat$, i.e.,
      for any such $C$, $i^*(C)$ is aspheric over the final element
      $e$ of $\Ahat$.
    \item\label{it:key.b.iii}
      The element $i^*(\Simplex_1)$ of $\Ahat$ is aspheric over the
      final object $e$.
    \end{enumerate}
  \end{enumerate}
\end{theorem}

\begin{remark}
  Of course the conditions in \ref{it:key.a} imply that $A$ is
  aspheric (take $C$ to be the final category), and hence by an easy
  lemma (of below, \S\ref{sec:40}) we get that for any $C$ in $\Cat$ $A\times C\to C$ is a
  weak equivalence (and even aspheric), and hence $A_{/C}\to A\times
  C$ is a weak equivalence (because its compositum with the weak
  equivalence $A\times C\to C$ is a weak equivalence by
  assumption). It is unlikely however that \ref{it:key.a} implies the
  conditions of \ref{it:key.b}, namely asphericity (not merely weak
  equivalence) of \eqref{eq:key.starstarstar} in
  b\ref{it:key.b.i}. But the opposite implication, namely
  b\ref{it:key.b.i} + asphericity of $A$ implies a\ref{it:key.a.i}, is
  trivial, because $A\times C\to C$ is aspheric.
\end{remark}

\begin{proof}[Proof of theorem] We stated \ref{it:key.a} and
  \ref{it:key.b} in a way to get a visibly decreasing cascade of
  conditions; and moreover that the weakest in \ref{it:key.a} implies
  the strongest, or that b\ref{it:key.b.ii} implies b\ref{it:key.b.i},
  is an immediate consequence of the standard \ref{lem:asphericitycriterion}
  for a functor between categories
  (p.~\ref{p:38}). The only point which is a
  little less formal is that b\ref{it:key.b.iii} implies
  b\ref{it:key.b.ii}. But using the final object in $C$, we get a
  (unique) homotopy in \Cat\ (relative to $\Simplex_1$)
  \[ \Simplex_1 \times C \to C,\]
  between the identity map of $C$ and the constant map with value
  $e_C$ -- hence by applying $i^*$, a homotopy relative to
  $i^*(\Simplex_1)$ (viewed as an ``interval'' by taking as
  ``endpoints'' the arrows deduced from
  $\delta_0,\delta_1 : e_\Cat=\Simplex_0\rightrightarrows \Simplex_1$ by
  applying $i^*$), between the identity map of $i^*(C)$ and a constant
  map of $i^*(C)$, which by the \hyperref[lem:homotopylemma]{homotopy
  lemma} implies that $i^*(C)$ itself is aspheric over
  $e_{\Ahat}$, qed.
\end{proof}
\addtocounter{remarknum}{1}
\begin{remarknum}
  The Condition b\ref{it:key.b.iii} can be stated by saying that
  $i^*(\Simplex_1)$ is a homotopy interval in $\Ahat$. All which
  needs to be checked for this, is that condition \ref{cond:HIb} for
  homotopy intervals (p.~\ref{p:59}) namely $e_0 \sand
  e_1 = \varnothing$ is satisfied in $\Ahat$, but this follows from
  the corresponding property for $\Simplex_1$ in \Cat\ (as the functor
  $i^*$ is left exact) and from the fact that $i^*$ transforms initial
  element $\varnothing_\Cat$ into initial element $\varnothing_{\Ahat}$.
  This last fact is equivalent to $i(a)\ne\varnothing$ for any
  $a$ in $A$, which is true as $i(a)$ has a final element.
\end{remarknum}

\bigbreak

\presectionfill\ondate{22.3.}\pspage{63}\par

\hangsection{Asphericity story retold
  \texorpdfstring{\textup(}{(}cont'd\texorpdfstring{\textup)}{)}:
  generalized nerve functors.}\label{sec:38}%
Let's get back to the ``asphericity story retold'' -- I had to stop
yesterday just in the middle, as it was getting prohibitively late.

I want to comment a little about the ``\hyperref[thm:keyresult]{key
  result}'' just stated and proved. The main point of this result,
forgetting the game of givings heaps of equivalent formulations of two
kinds of properties, is that the extremely simple condition
b\ref{it:key.b.iii}, namely that $i^*(\Simplex_1)$ is aspheric over the
unit element $e_\Ahat$ of \Ahat, plus asphericity of the latter,
ensure already the conditions in \ref{it:key.a}, which can be viewed
(among others) as just stating that the compositum
\[ i_Ai^* : \Cat\to\Ahat\to\Cat\]
from \Cat\ to \Cat\ induces an autoequivalence of the localized
category \Hot, or (what amounts still to the same) a functor
$\Hot\to\Hot$ isomorphic to the identity functor. (NB\enspace that two do
indeed amount to the same follows at once from the implication
\ref{it:key.a.iv} $\Rightarrow$ \ref{it:key.a.ii} in \ref{it:key.a}.)
It is interesting to note that both properties, the stronger one that
$i^*(\Simplex_1)$ is aspheric over $e_\Ahat$, and the weaker one in
terms of properties of the compositum $i_Ai^*$, \emph{make a sense without
any reference to the extra assumption that the categories $i(a)$ have
a final element each}, nor to the corresponding map \eqref{eq:key.star}
$A_{/C}\to C$ (functorial in $C$). This suggests that there should exist a more general statement than in
the theorem, without making the assumption about final objects in the
categories $i(a)$, and without the possibility of a \emph{direct}
comparison of $A_{/C}$ and $C$ through a functor between them. Indeed
I have an idea of a statement in this respect, however for the time
being it seems that the theorem as stated is sufficiently general for
handling the situations I have in mind.

Applying the theorem to the \emph{canonical} functor $i$
\[i_0 : a\mapsto A_{/a} : A \to \Cat,\]
whose canonical extension to \Ahat\ (as a functor $\Ahat\to\Cat$
commuting to direct limits) is the functor
\[i_A: F \mapsto A_{/F} : \Ahat\to\Cat,\]
giving rise to the right adjoint
\[i_0^* = j_A : \Cat\to\Ahat,\]
the condition \ref{it:key.a.i} in \ref{it:key.a} is nothing but the
familiar condition
\[ i_Aj_A(C)\to C\quad\text{a weak equivalence for any $C$ in \Cat,}\]
which we had used for our first (or rather, second already!)
definition of so-called ``test categories''. Later on we considerably
strengthened this condition -- we now call them ``test categories in
the wide sense''. On the other hand, as we already noticed before,
here $i_0^*(\Simplex_1) = j_A(\Simplex_1)$ \emph{is nothing but the
  Lawvere element $L_A$ of \Ahat}. Thus the main content of
the\pspage{64} theorem in the present special case can be
formulated thus:
\begin{corollarynum}\label{cor:38.1}
  Assume the Lawvere element $L_A$ in \Ahat\ is aspheric over the
  final object $e_\Ahat$ of \Ahat, and that moreover the latter be
  aspheric, i.e., $A$ aspheric. Then $A$ is a test category in the
  wide sense, namely for any $C$ in \Cat, the canonical functor
  $i_Aj_A(C)\to C$ is a weak equivalence \textup(and even
  aspheric\textup -- see prop.\ on p.~\ref{p:38}, and
  also p.~\ref{p:35}, for equivalent formulations).
\end{corollarynum}

Thus we \emph{did} get after all a ``handy criterion'' sufficient to ensure
this basic test-property, which looks a lot less strong a priori than
the condition \ref{it:31.T2} (plus \ref{it:31.T1}, \ref{it:31.T3}) of total asphericity of
\Ahat.\footnote{\alsoondate{25.3.} or rather, than the set of conditions
  \ref{it:31.T1} to \ref{it:31.T3}}

But let's now come back to the more general situation of the theorem,
with a functor
\[ i : A \to \Cat\]
subjected only to the mild condition that the categories $i(a)$ (for
$a$ in $A$) have final objects. Assume condition b\ref{it:key.b.iii}
to be satisfied, namely that $i^*(\Simplex_1)$ is aspheric over
$e_\Ahat$, and therefore it is a homotopy interval in \Ahat. Let again
$L_A$ be the Lawvere element in \Ahat, we define (independently of any
assumption on $i$) a morphism of ``intervals'' in \Ahat, compatible
even with the natural composition laws (by intersection) on both
members
\[\varphi : J=i^*(\Simplex_1) \to L_A.\]
For this we remember that for $a$ in $J$, we get
\[ J(a) \simeq \Crib{i(a)} \hookrightarrow \text{set of all subobjects
  of $i(a)$ in \Cat}\]
(this bijection and the inclusion being functorial in $a$), thus if
$C$ is in $J(a)$, i.e., a crible in $i(a)$, we associate to this
    \\[\baselineskip]%
  \hspace*{3em}$\varphi(C)$, \quad%
\parbox[t]{0.7\textwidth}{subobject of $a$ in \Ahat,
  corresponding to the crible in $A_{/a}$ of all $b/a$ such that
  $i(b)\to i(a)$ factors through the crible $C\subset i(a)$;}
  \\[\baselineskip]%
it is immediate that the map $\varphi_a : J(a)\to L(a)$ thus obtained
is functorial in $a$ for variable $a$, hence a map $\varphi : J\to L$,
and it is immediately checked too that this is ``compatible with
endpoints'' -- namely when $C$ is full, respectively empty, then so is
$\varphi(C)$ (for the ``empty'' case, this comes from the fact that
the categories $i(a)$ are non-empty). Applying now the comparison lemma
for homotopy intervals (p.~\ref{p:60}), we get the
following

\begin{corollarynum}\label{cor:38.2}
  Under the general conditions of the theorem, and assuming moreover
  that $i^*(\Simplex_1)$ is aspheric over $e_\Ahat$ \textup(condition
  \textup{b\ref{it:key.b.iii})}, it follows that the Lawvere element $L_A$ of
  \Ahat\ is aspheric over $e_\Ahat$. Assume moreover that $e_\Ahat$ is
  aspheric, i.e., $A$ aspheric. Then we get:\pspage{65}
  \begin{enumerate}[label=\alph*),font=\normalfont]
  \item\label{it:38.a}
    The category $A$ is a test category in the wide sense \textup(cf.\
    cor.\ \textup{\ref{cor:38.1})}.
  \item\label{it:38.b}
    Both functors $i^* : \Cat\to\Ahat$ and $i_A:\Ahat\to\Cat$ are
    ``modelizing'', namely the set of weak equivalences in the source
    category is the inverse image of the corresponding set of arrows
    in the target category, and the functor induced on the
    localizations with respect to weak equivalences is an equivalence
    of categories.
  \item\label{it:38.c}
    Let $W$ \textup(resp.\ $W_A$\textup) be the set of weak
    equivalences in \Cat\ \textup(resp.\ in \Ahat\textup). Then the
    functor
    \[\Hot \eqdef W^{-1}\Cat \to W_A^{-1}\Ahat\]
    induced by $i^*$ is canonically isomorphic to the quasi-inverse of
    the functor in opposite direction induced by $i_A$, or
    equivalently, can.\ isomorphic to the functor in the same
    direction induced by $i_0^*=j_A:\Cat\to\Ahat$ \textup(cf.\ again
    cor.\ \textup{\ref{cor:38.1}} above for the notations\textup).
  \end{enumerate}
\end{corollarynum}

Of course \ref{it:38.a} follows from cor.\ \ref{cor:38.1}, and implies
that $i_A$ has the properties stated in \ref{it:38.b}. That the
analogous properties hold for $i^*$ too, and the rest of the
statement, i.e., \ref{it:38.c}, follows formally, using the fact that
the compositum $i_Ai^*$ is canonically isomorphic to the identity
functor once we pass to the localized category \Hot\ (using the
theorem, b\ref{it:key.b.iii} $\Rightarrow$ \ref{it:key.a}).

This corollary shows that, up to canonical isomorphism, the functor
\[\Hot\to W_A^{-1}\Ahat\quad\text{induced by $i^*:\Cat\to\Ahat$}\]
\emph{does not depend on the choice of the functor} $i:A\to\Cat$,
provided only this functor satisfies the two conditions that it takes
its values in the full subcategory of \Cat\ of all categories with
final objects, and that moreover $i^*(\Simplex_1)$ be aspheric over the
final object of \Ahat\ (plus of course the condition of asphericity on
$A$ itself). There \emph{is} of course always the \emph{canonical}
choice of a functor $i:A\to\Cat$, namely $i_A: a\mapsto A_{/a}$ which
(because of its canonicity) looks as the best choice theoretically --
and it was the first one indeed we investigated into. But in practical
terms, the categories $A_{/a}$ are (in the concrete cases one might
think of) comparatively big (for instance, infinite) and the
corresponding functor $j_A=i_A^*$ gives comparatively clumsy
``models'' in \Ahat\ for describing the homotopy types of given
``models'' in \Cat, whereas we can get away with considerably neater
models in \Ahat, using a functor $i$ giving rise to categories $i(a)$
which are a lot easier to compute with (for instance, finite
categories of very specific type). The most commonly used is of course
the \emph{nerve functor} $i^*$, corresponding to the standard
embedding of $A=\Simplex$ into \Cat\ -- and in the general case of the
theorem above, complemented by cor.\ \ref{cor:38.2},
\emph{the functor $i^*$ should be viewed as a generalized nerve
  functor}.\pspage{66}

To be completely happy, we still need a down-to-earth sufficient
criterion to ensure that $i^*(\Simplex_1)$ is aspheric over $e_\Ahat$,
in the spirit of the somewhat awkward (deleted) corollary of end of section \ref{sec:36}. The following seems quite adequate for all
cases I have in mind at the present moment:

\begin{corollarynum}\label{cor:38.3}
  Under the general conditions of the theorem on $A$ and $i:A\to\Cat$,
  assume moreover we got a homotopy interval $(I,\delta_0,\delta_1)$
  in \Ahat, that $A$ has a final object $e$ \textup(which is therefore
  also a final object of \Ahat\ so we may view $\delta_0,\delta_1$ as
  maps $e\rightrightarrows I$, i.e., elements in $I(e)$\textup), and
  let $i_!:\Ahat\to\Cat$ be the canonical extension of $i$ to \Ahat\
  \textup(commuting with direct limits\textup). Assume moreover
  $i(e)=\Simplex_0$ \textup(the final object in \Cat\textup), and that
  we can find a map in \Cat
  \[i_!(I) \to \Simplex_1\]
  compatible with $\delta_0,\delta_1$, i.e., whose compositae with
  $i_!(\delta_n) : \Simplex_0\to\Simplex_1$ for $n\in\{0,1\}$ are the
  two standard maps $\boldsymbol\delta_0,\boldsymbol\delta_1$ from
  $\Simplex_0$ to $\Simplex_1$. Then the condition
  \textup{b\ref{it:key.b.iii}} of the theorem holds, namely
  $i^*(\Simplex_1)$ is aspheric over $e$.
\end{corollarynum}

Indeed, to give a map $i_!(I)\to\Simplex_1$ in \Cat\ amounts to the same
as giving a map $I\to i^*(\Simplex_1)$ in \Ahat\ (namely $i_!$ and $i^*$
are adjoint), moreover the extra condition involving
$\delta_0,\delta_1$ just means that his map respects endpoints. Using
the composition law of intersection on $i^*(\Simplex_1) = (a \mapsto
\Crib(i(a)))$, the comparison lemma for homotopy intervals (p.~%
\ref{p:60} implies that $i^*(\Simplex_1)$ is aspheric over
$e$, qed.

In the cases I have in mind, $I$ is even an element of $A$, hence
$i_!(I)=i(I)$, moreover $i(I)=\Simplex_1$ and the map
$i_!(I)\to\Simplex_1$ above is the \emph{identity}! The choice of the
functor $i$ is in every case ``the most natural one'' (discarding
however the clumsy $i_0=i_A$, and trying to get away with categories
$i(a)$ which give the simplest imaginable description of objects and
arrows of the would-be test category $A$), and the choice of $I$
itself is still more evident -- it is \emph{the} object of $A$ (or
\emph{one} among the objects, in cases such as $A=\Simplex^n$, giving
rise to simplicial \emph{multi}complexes\ldots) which suggests most
strongly the picture of an ``interval''. Thus the one key verification
we are left with (all the rest being ``formal'' in terms of what
precedes) is the asphericity of $I$ over $e$, i.e., that all the
products $I\times a$ are aspheric.

\hangsection[Returning upon terminology: strict test categories, and
\dots]{Returning upon terminology: strict test categories, and strict
  modelizers.}\label{sec:39}%
Maybe it is time now to come back to the property of total asphericity
of \Ahat, expressed by the condition \ref{it:31.T2} on $A$, namely
that the product in \Ahat\ of any two elements in $A$ is aspheric. As
we saw, (if $A \neq \emptyset$) this implies already asphericity of $A$, i.e., of the topos
\Ahat. In our present setting, total asphericity is of special
interest only when coupled with the property \ref{it:31.T3}, which
amounts to saying that the Lawvere element $L_A$ in \Ahat\
is\pspage{67} aspheric over $e_\Ahat$, or (what now amounts to the
same) that $L_A$ is aspheric. However, when faced with the question to
decide whether a given $A$ does indeed satisfy \ref{it:31.T2}, it will
be convenient to use a system $(I,e_0,e_1)$ in \Ahat\ of which we know
beforehand it is a homotopy interval (the delicate part of this notion
being the asphericity of all products $I\times a$ for $a$ in $A$). As
part of the ``story retold'', I recall now the most natural geometric
assumption which will ensure \ref{it:31.T2}, i.e., total asphericity
of \Ahat:

\begin{proposition}
  Let $A$ be a small category, $(I, e_0,e_1)$ a homotopy interval in
  \Ahat. Assume that for any $a$ in $A$\kern1pt, there exists a homotopy
  \[h_a : I\times a\to a\]
  between $\id_a$ and a constant map $c_a:a\to a$ \textup(defined by a
  section of $a$ over the final element of \Ahat, i.e., by a map $e\to
  a$\textup). Then \Ahat\ is totally aspheric, i.e., every $a\in\Ob A$
  is aspheric over $e$ \textup(or, equivalently, the product in \Ahat\ of any
  two elements $a,b$ in $A$ is aspheric\textup).
\end{proposition}

This is a particular case of the
``\hyperref[lem:homotopylemma]{homotopy lemma}'' (p.~%
\ref{p:60}). In fact we don't even use the condition $e_0
\sand e_1 = \varnothing_\Ahat$ on the ``interval'' $(I,e_0,e_1)$, but we
easily see that this condition follows from the existence of the
homotopy and the sections of any $a$ in $A$ over $e$.

Before returning to the investigation of specific test categories, I
want to come back on some terminology. The condition on \Ahat\ that
the Lawvere element $L_A$ should be aspheric over $e$ has taken lately
considerable geometric significance, and merits a name. I will say
from now on that $A$ is a \emph{test category}, and that \Ahat\ is an
\emph{elementary modelizer}, if this condition is satisfied, and if
moreover $A$ is aspheric. This condition (which amounts to our former
\ref{it:31.T1} + \ref{it:31.T3}) is weaker than what we had lately
called a test category, as we had overlooked so far the fact that
\ref{it:31.T1} and \ref{it:31.T3} alone already imply the basic
requirement about $i_Aj_A(C)\to C$ being a weak equivalence for any
$C$, and hence $W_A^{-1}\Ahat$ being canonically equivalent to
\Hot. Thus it seemed for a while that the only handy conditions we
could get for ensuring this requirement were
\ref{it:31.T1}\ref{it:31.T2}\ref{it:31.T3}, all together (in fact, it
turned out later that \ref{it:31.T2} already implies
\ref{it:31.T1}). When all three conditions \ref{it:31.T1} to
\ref{it:31.T3} are satisfied, I'll now say that\pspage{68} $A$ is a
\emph{strict test category}, and that \Ahat\ is an \emph{elementary
  strict modelizer}. Here the notion of a ``strict modelizer'' (not
necessary an elementary one) makes sense independently, it means a
category $M$ endowed with a set $W\subset\Fl(M)$ satisfying conditions
\ref{it:37.a}\ref{it:37.b}\ref{it:37.c}
of p.~\ref{p:59}, such that $W^{-1}M$ is equivalent to
\Hot, and such moreover that the localization functor $M\to W^{-1}M$
\emph{commutes with finite products} (perhaps we should also insist on
commutation with finite sums, and possibly include also infinite sums
and products in the condition -- whether this is adequate is not quite
clear yet). The mere condition that $i_Aj_A(C)\to C$ should be a weak
equivalence for any $C$ in \Cat, or equivalently that $i_Aj_A(C)$
should be aspheric when $C$ has a final element, will be referred to
by saying that $A$ is a \emph{test category in the wide sense}. It
means little more than the fact that $(\Ahat,W_A)$ is a
modelizer. Finally, if we merely assume that \Ahat\ admits a homotopy
interval, or equivalently, that $L_A$ is aspheric over $e_\Ahat$, we
will say that $A$ is a \emph{local test category} (because it just
means that the induced categories $A_{/a}$ for $a$ in $A$ are test
categories), and accordingly \Ahat\ will be called a \emph{local
  elementary modelizer}. More generally, we call a topos \scrA\ such
that the Lawvere element $L_\scrA$ be aspheric over the final element
a \emph{locally modelizing topos}, and we call it a \emph{modelizing
  topos} if it is moreover aspheric. When the topos is locally
aspheric, i.e., admits a generating family made up with aspheric
objects of \scrA, then \scrA\ is indeed a locally modelizing topos
if{f} the final object can be covered by elements $U_i$ such that the
induced topoi $\scrA_{/U_i}$ be modelizing topoi. This terminology for
topoi more general than of the type \Ahat\ is possibly somewhat hasty,
as the relation with actual homotopy models, namely with the question
whether $(\scrA,W_\scrA)$ is a modelizer, has not been investigated
yet. Still I have the feeling the relation should be a satisfactory
one, much along the same lines as we got in the case of topoi of the
type \Ahat. We will not dwell upon this now any longer.

For completing conceptual clarification, we should still make sure
that a test category $A$ need not be a strict test category, i.e.,
need not be totally aspheric. As a candidate for a counterexample, one
would think about a category $A$ endowed with a final element and an
element $I$, together with $\delta_0,\delta_1:e\to I$, satisfying
$\Ker(\delta_0,\delta_1)=\varnothing_\Ahat$ in \Ahat, $I$ being
squarable in $A$, namely the products $I\times a$ ($a \in \Ob A$) are
in $A$ -- this alone will imply that $I$ is aspheric over $e$ in
\Ahat, hence $A$ is a test category, but it is unlikely that this
alone will imply equally total asphericity of \Ahat. We would think,
as the most ``economical'' example, one where $A$ is made up with
elements of the type $I^n$ ($n\in\bN$) namely cartesian\pspage{69}
powers of $I$, plus products $a_0\times I^n$ ($n\in\bN$), where $a_0$
is an extra element, and those maps between these all those (and not
more) which can be deduced from $\delta_0,\delta_1$ and the assumption
that the elements $I^n$ and $a\times I^n$ are cartesian products
indeed. This is much in the spirit of the construction of a variant of
the category $\square$ of ``standard cubes'', which naturally came to
mind a while ago (cf.\ p.~%
\ref{p:48}--\ref{p:49}). I very much doubt $A$
satisfies \ref{it:31.T2}, and would rather bet that $a_0 \times a_0$
in \Ahat\ is \emph{not} aspheric.

Another type of example comes to my mind, starting with a perfectly
good (namely \emph{strict}) test category $A$, and taking an induced
category $A_{/a_0}$, with $a_0$ in $A$. This is of course a test
category (it would have been enough that $A$ be just a local test
category), however it is unlikely that it will satisfy \ref{it:31.T2},
namely the induced topos be totally aspheric. This would imply for
instance that for any two non-empty subobjects of $a_0$ in \Ahat\
(namely subobjects of the final element in $\Ahat_{/a_0} \simeq
(A_{/a_0})\uphat$) have a non-empty intersection. Now this is an
exceedingly strong property of $a_0$, which is practically never
satisfied, except for the final object of $A$. Now here is a kind of
``universal'' counterexample. Take $A$ any test category and $I$ a
homotopy interval in \Ahat, thus $I$ is aspheric and hence $A_{/I}$ is
again a test category (namely locally test and aspheric), but it is
\emph{never} a strict test category, because the two standard
subobjects $e_0,e_1$ given with the structure of $I$ are non-empty,
and however their intersection is empty!

These reflections bring very near how much stronger the extra strictness
requirement \ref{it:31.T2} is for test categories, than merely the
conditions \ref{it:31.T1}, \ref{it:31.T3} without \ref{it:31.T2}.

\hangsection[Digression on cartesian products of weak equivalences in
\dots]{Digression on cartesian products of weak equivalences in
  \texorpdfstring{\Cat}{(Cat)}; 4 weak equivalences \emph{relative} to
  a given base object.}\label{sec:40}%
Yesterday I incidentally made use of the fact that if $A$ is an
aspheric element in \Cat, then for any other object $C$,
$C\times A\to C$ is aspheric (and a fortiori a weak equivalence). The
usual \ref{lem:asphericitycriterion} for a functor shows that it is
enough to prove that for any $C$ with final element, $C\times A$ is
aspheric. For this again, it is enough to prove that the projection
$C\times A\to A$ is aspheric, which by localization upon $A$ means
that the product categories $C\times (A_{/a})$ (with $a$ in $A$) are
aspheric. Finally we reduced to checking that the product of two
categories with final element is aspheric, which is trivial because
such a category has itself a final element.

The result just proved can be viewed as a particular case of the
following
\begin{proposition}
  In \Cat, the cartesian product of two weak equivalences is a weak equivalence.
\end{proposition}

(We\pspage{70} get the previous result by taking the weak
equivalences $A\to e_\Cat$ and $\id_C:C\to C$.) To prove the
proposition, we are immediately reduced to the case when one of the
two functors is an identity functor, i.e., proving the
\begin{corollarynum}\label{cor:40.1}
  If $f:A'\to A$ is a weak equivalence in \Cat, then for any $C$ in
  \Cat, $f\times \id_C: A'\times C\to A\times C$ is a weak equivalence.
\end{corollarynum}

For proving this, we view the functor $f\times\id_C$ as a morphism of
categories ``over $C$'', which corresponds to a situation of a
morphism of topoi over a third one
\[\begin{tikzcd}[baseline=(Y.base),column sep=tiny,row sep=small]
  X' \ar[dr]\ar[rr,"f"] & & X\ar[dl] \\ & |[alias=Y]| Y &
\end{tikzcd},\]
we will say that $f$ is a \emph{weak equivalence relative to $Y$}, if
not only this is a weak equivalence, but remains so by any base change
by a localization functor
\[ Y_{/U} \to U \]
giving rise to
\[ f_{/U} : X'_{/U} \to X_{/U}.\]
As usual, standard arguments prove that it is enough to take $U$ in a
set of generators of the topos $Y$. In case $Y=\Chat$, we may
take $U$ in $C$. In case moreover $X,X'$ are defined by small
categories $P,P'$ and a functor of categories over $C$
\[\begin{tikzcd}[baseline=(C.base),column sep=tiny,row sep=small]
  P' \ar[dr]\ar[rr,"F"] & & P\ar[dl] \\ & |[alias=C]| C &
\end{tikzcd},\]
this condition amounts to demanding that for any $c$ in $C$, the
induced functor
\[F_{/c} : P'_{/c} \to P_{/c}\]
be a weak equivalence. In our case $P=A\times C$, $P'=A'\times C$,
$F=f\times\id_C$, the induced functor can be identified with
$f\times\id_{C_{/c}}$. This reduces us, for proving the corollary, to
the case when $C$ has a final element. But consider now the
commutative diagram
\[\begin{tikzcd}[baseline=(A.base),row sep=small]
  A'\times C\ar[d]\ar[r,"f\times\id_C"] & A\times C\ar[d]\\
  A' \ar[r,"f"] & |[alias=A]| A
\end{tikzcd},\]
where the vertical arrows are the projections and hence, by what was
proved before, weak equivalences. As $f$ is a weak equivalence, it
follows that $f\times\id_C$ is a weak equivalence too, qed.

The proposition above goes somewhat in the direction of looking at
``homotopy properties of \Cat'' and ``how far \Cat\ is from being a
closed model category in Quillen's sense''. It is very suggestive
for\pspage{71} having a closer look at maps $f: B \to A$ in
\Cat\ which are ``universally weak equivalences'', i.e.,
$W_\Cat$-aspheric, namely such that for \emph{any} map in \Cat\ $A'\to
A$ (not only a localization $A_{/a}\to A$), the induced functor
$B'=B\times_A A' \to A'$ is a weak equivalence. This property is a lot
stronger than just asphericity, and reminds of the ``trivial
fibrations'' in Quillen's theory. The usual criterion of asphericity
for $B'\to A'$ shows that $f$ is $W$-aspheric if{f} for any $A'$ with
final element, and any functor $A'\to A$, the fiber product
$B'=B\times_A A'$ is aspheric. The feeling here is (suggested partly
by Quillen's terminology) that this property is tied up some way with
the property for $f$ to be a fibering (or cofibering?) functor, in the
sense of fibered and cofibered categories, with aspheric fibers
moreover. Presumably, fibered or cofibered categories (in the purely categorical sense), with ``base
change functors'' which are weak equivalences, will play the part of
Serre-Quillen's fibrations -- and it is still to be guessed what kind
of properties of a functor will play the part of
\emph{co}fibrations. Apparently they will have to be a lot more
stringent than just monomorphisms in \Cat, cf.\ section \ref{sec:30} (and \ref{sec:9}).

However, I feel it is not time yet to dive into the homotopy theory
properties of the all-encompassing basic modelizer \Cat, but rather
come back to the study of general (and less general) test categories.

\hangsection[Role of the ``inspiring assumption'', and of saturation
\dots]{Role of the ``inspiring assumption'', and of saturation
  conditions on ``weak equivalences''.}\label{sec:41}%
One comment still, upon the role played in the theory I am developing
of the assumption (p.~\ref{p:30}) that the category of
autoequivalences of \Hot\ is equivalent to the final category. This
assumption has been a crucial guide for putting the emphasis where it
really belongs, namely upon the set $W\subset\Fl(M)$ of weak
equivalences within a category $M$ which one would like to take in
some sense as a category of models for homotopy types -- the functors
$M\to\Hot$ following along automatically. However, in no statement
whatever I proved so far, was this assumption ever used. On the other
hand, the notion of a \emph{modelizer} introduced in the wake of the
``assumption'' (cf.\ \S \ref{sec:28}) was tacitly changed during the
reflection, by dropping altogether the condition \ref{it:Mod.a} of
(strong) saturation of $W$, namely that $W$ is just the set of arrows
made invertible by the localization functor $M\to W^{-1}M$. Instead of
this, it turned out that the saturation condition we really had at
hands and which was adequate for working, was the conditions
\ref{it:37.a} to \ref{it:37.c} I finally wrote down explicitly
yesterday (\S \ref{sec:37}). As for the strong saturation condition,
for the time being (using nothing but what has actually been proved so
far, without reference to ``well-known facts'' from homotopy theory),
it is not even clear that the basic modelizer \Cat\ is one in the
initial sense, namely that the set of weak equivalences\pspage{72} in
\Cat\ is strongly saturated. However, from the known relation between
\Cat\ and an elementary modelizer \Ahat, it follows that $W_\Cat$ is
strongly saturated if{f} $W_A\subset\Fl(\Ahat)$ is. This implies that
it is enough to prove strong saturation in \emph{one} elementary
modelizer \Ahat, to deduce it in all others, as well as in
\Cat. However, in the case at least when $A=\Simplex$, it is indeed
``well-known'' that the weak equivalences in \Ahat\ satisfy the strong
saturation condition. In terms of Quillen's set-up, it follows from
the fact that $\Simplexhat$ is a ``closed model category'', and prop.\
1 in I~5.5 of Quillen's expos\'e. The only thing which is not quite
understood, not by me at any rate, is why $\Simplexhat$ \emph{is}
indeed a closed model category -- Quillen's proof it seems relies
strongly on typical simplicial techniques. I'll have to look if the
present set-up will suggest a more conceptual proof, valid possibly
for any test category (or at least, any strict test category). I'll
have to come back upon this later. For the time being, I feel a
greater urge still to understand about the relationship between
different test categories -- also, I did not really finish with my
review of what may be viewed as the ``standard'' test categories such
as $\Simplex$, $\square$ and their variants.

\bigbreak

\presectionfill\ondate{25.3.}\par

\hangsection[Terminology revised (model preserving
functors). \dots]{Terminology revised
  \texorpdfstring{\textup(}{(}model preserving
  functors\texorpdfstring{\textup)}{)}. Submodelizers of the basic
  modelizer \texorpdfstring{\Cat}{(Cat)}.}\label{sec:42}%
In the notes last time I made clear what finally has turned out for me
to correspond to the appelation of a ``modelizer'', as prompted by the
internal logics of the situations I was looking at, in terms of the
information available to me. I should by then have added what a
\emph{model-preserving} (or \emph{modelizing}) \emph{functor} between
modelizers $(M,W)$ and $(M',W')$ has turned out to mean, which has
undergone a corresponding change with respect to what I first
contemplated calling by that name (\S \ref{sec:28}). Namely, here it
turned out that I should be more stringent for this notion, replacing
the condition $f(W)\subset W'$ by the stronger one
\[ W = f^{-1}(W'),\]
and moreover, of course, still demanding that the corresponding
functor
\[W^{-1}M \to W'^{-1}M'\]
should be an equivalence. It is by now established, with the exception
of just the two first that, all the functors occurring in the diagram
on p.~\ref{p:31} are indeed model preserving, namely it
is so for $\alpha$, $\beta$, $\xi$, $\eta$ -- and also the right
adjoints to $\xi$, $\eta$, as expected by then. It should be more or
less trivial that the first of the functors in this diagram, namely
the canonical inclusion $\Ord\to\Preord$, as well as the left adjoint
from \Preord\ to \Ord, are model preserving -- except for the fact
that it has not been yet established that these two categories are
indeed modelizers (for the natural notion of weak equivalences, induced
from \Cat)),\pspage{73} namely that the localized categories with
respect to weak equivalences are indeed equivalent to \Hot. The only
natural way one might think of this to be proved, is by proving that
the inclusion functor from either into \Cat\ induces an equivalence
between the localizations, which would imply at the same time that
this inclusion functor is indeed model-preserving, and hence that all
the functors in the diagram of p.~\ref{p:31} are model
preserving functors between modelizers. I still do believe this should
be so, and want to give below a reflection which might lead to a proof
of this.

Before, there is still one noteworthy circumstance I want to
emphasize. Namely, it occurred in a number of instances that we got in
a rather natural way several modelizing functors between two
modelizers $(M,W)$ and $(M',W')$, in one direction or the other -- for
instance from \Cat\ to \Ahat\ using different functors $i:A\to\Cat$,
or from \Ahat\ to \Cat\ using $a\mapsto A_{/a}$. It turned out that
the corresponding functors between $W^{-1}M$ and $W'^{-1}M'$ were
always canonically isomorphic when in the same direction, and
quasi-inverse of each other when in opposite directions. This is
indeed very much in the spirit of the ``inspiring assumption'' of p.~%
\ref{p:30}, that the category of autoequivalences of \Hot\
is equivalent to the unit category, which implies indeed that for two
categories $H,H'$ equivalent to \Hot, for two equivalences from $H$ to
$H'$, there is a unique isomorphism between them. Quite similarly,
there have been a number of situations (more or less summed up in the
end in the ``\hyperref[thm:keyresult]{key theorem}'' of p.~%
\ref{p:61}) when by localization we got a natural functor
$f$ from some $W^{-1}M$ to another $W'^{-1}M'$, and another $F$ which
is known already to be an equivalence, and it turns out that $f$ is
isomorphic to $F$ (and in fact, then, canonically so) if{f} for any
object $x$ in the source category, $f(x)$ and $F(x)$ are
isomorphic. This suggests that \emph{presumably, every functor
  from \Hot\ into itself, transforming every object into an isomorphic
  one, is in fact isomorphic to the identity functor}.

I want now to make a comment, implying that there are many full
subcategories $M$ of \Cat, such that for the induced notion of weak
equivalence, $M$ becomes a modelizer, and the inclusion functor a
modelizing functor -- or, what amounts to the same, that the canonical
functor
\[W_M^{-1}M \to W_\Cat^{-1}\Cat = \Hot\]
is an equivalence. To see this, let more generally $(M,W)$ be any
category endowed with a subset $W\subset\Fl(M)$, and let $h:M\to M$ be
a functor such that $h(W)\subset W$, and such that the induced functor
$W^{-1}\to W^{-1}M$ is an equivalence. Let now $M'$ be any full
subcategory of $M$ such that $h$ factors through $M'$, let $h' : M\to
M'$ be the corresponding induced functor, and $W'=W\sand
\Fl(M')$. Then it is formal that the inclusion $g:M'\to
M$\pspage{74} and $h':M\to M'$ induce functors between $W'^{-1}M'$
and $W^{-1}M$, which are quasi-inverse to each other -- hence these
two categories are equivalent. In case $(M,W)$ is a modelizer, this
implies that $(M',W')$ is a modelizer too and the inclusion functor
$g$, as well as $h'$, are model preserving.

We can apply this remark to the modelizer \Cat, and to the functor
\[ i_Ai^*:\Cat\to\Cat\]
defined by any functor $i:A\to\Cat$ satisfying the conditions of the
``\hyperref[thm:keyresult]{key theorem}'' (\S \ref{sec:37}), for instance the functor $i_A\restrto A: a
\mapsto A_{/a}$ (where $A$ is any test category). We get the following

\begin{proposition}
  Let $M$ be any full subcategory of \Cat, assume there exists a
  test-category $A$ such that for any $F$ in \Ahat, the category
  $A_{/F}$ belongs to $M$ \textup(i.e., the functor $i_A:\Ahat\to\Cat$,
  $F\mapsto A_{/F}$, factors through $M$\textup). Let $W_M$ be the set
  of weak equivalences in $M$. Then $(M,W_M)$ is a modelizer and the
  inclusion functor $M\to\Cat$ is model preserving.
\end{proposition}

The same of course will be true for any full subcategory of \Cat\
containing $M$ -- which makes an impressive bunch of modelizers
indeed! When the test category $A$ is given, one natural choice for
$M$ is to take all categories $C$ which are ``locally isomorphic to
$A$'', namely such that for any $x$ in $C$, the induced category
$C_{/x}$ be isomorphic to a category of the type $A_{/a}$, with $a$ in
$A$.

It would be tempting to apply this result to the full subcategory
\Ord\ of \Cat\ -- but for this to be feasible, would mean exactly that
there exists a test-category $A$ defined by an ordered set (or at
least ``locally ordered''). To see whether there exists indeed such an
ordered set looks like a rather interesting question -- maybe it would
give rise to algebraic models for homotopy types, simpler than those
used so far, namely simplicial and cubical complexes and
multicomplexes. It is interesting to note that if such a test category
should exist, it will \emph{not} be in any case a strict test
category. Indeed, the topos \Ahat\ associated to an ordered set $A$
can be viewed also, as we saw before (\S \ref{sec:22}), as
associated to a suitable topological space (namely $A$ endowed with a
suitable topology, the open sets being just the ``cribles'' in
$A$). But we have seen that the topos associated to a topological
space \emph{cannot} be \emph{strictly} modelizing (prop. 2 cor.\
\ref{cor:35.2.2} in section \ref{sec:35}).

This remark confirms the feeling that it was worth while emphasizing
the notion of a test category (just \ref{it:31.T1} to \ref{it:31.T3})
by a simple and striking name as I finally did, rather than bury it
behind the notion I now call a strict test category, which is
considerably more stringent and, moreover, more ``rigid''. For
instance, it is not stable under localization $A_{/a}$, whereas the
notion of a test category is -- indeed, for any aspheric $I$ in \Ahat,
$A_{/I}$ is still a test category.\pspage{75}

\hangsection[The category \protect\smashSimplexf{} of simplices without
degeneracies as a \dots]{The category
  \texorpdfstring{\protect\Simplexf}{Delta-f}{} of simplices without
  degeneracies as a weak test category -- or ``face complexes'' as
  models for homotopy types.}\label{sec:43}%
Now let's come back for a little while again to the so-called
``standard test categories'', and check how nicely the ``story
retold'' applies to them.

Not speaking about multicomplexes, there are essentially two variants
for ``categories of simplices'' as test categories. The smaller, more
commonly used one, is the category $\Simplex$ of ``ordered simplices'',
most conveniently described as the full subcategory of \Cat\ defined
by the family of simplices $\Simplex_n$ ($n\in\bN$). Here the most
natural choice for $i:\Simplex\to\Cat$ is of course the inclusion
functor. As $i$ is fully faithful and $\Simplex_1$ is in the image, it
follows that $i^*(\Simplex_1)=\Simplex_1$, and we have only to check (for
$A$ to be a test-category with ``test-functor'' $i$) that $\Simplex_1$
is aspheric over $e=\Simplex_0$, namely that all products
$\Simplex_1\times\Simplex_n$ are aspheric -- which we did. The extra
condition of the ``total asphericity criterion'' (proposition on p.~%
\ref{p:67}), namely existence of a homotopy in
$\Simplexhat$ from the identity map to a constant map, for any
$\Simplex_n$, is indeed satisfied: it is enough to define such a homotopy
in \Cat, which is trivial using the final element of $\Simplex_n$. Thus
$\Simplex$ is in fact a \emph{strict} test category.

As for $\widetilde\Simplex$, the most elegant choice theoretically is to
take the category of all non-empty finite sets, but his leads to
set-theoretic difficulties, as this category is not small - thus we
take again the standard non-ordered simplices
\[\widetilde\Simplex_n = \bN \sand {[0,n]} \quad (n\in\bN),\]
so as to get a ``reduced'' category with a countable set of
objects. This time, as $\widetilde\Simplex$ is stable under finite
products, and contains the ``interval''
$(\widetilde\Simplex_1,\delta_0,\delta_1)$ (which is necessarily then a
homotopy interval, as all elements of $\widetilde\Simplex$ are aspheric
over $e=\widetilde\Simplex_0$), the fact that $\widetilde\Simplex$ is a
strict test category is trivial. As for a test functor, the neatest
choice is the one we said before, namely associating to every finite
set the ordered set of all non-empty subsets. We thus get
      \\[\baselineskip]%
  \hspace*{3em}$\widetilde i : \widetilde\Simplex\uphat \to \Cat$ \quad%
\parbox[t]{0.7\textwidth}{factoring in fact through \Ord, as
does \\ $i:\Simplex\to\Cat$ above)}
  \\[\baselineskip]%
To prove it is indeed a test functor, the corollary
\ref{cor:38.3} to the ``key theorem'' (p.~%
\ref{p:66}) applies, taking of course
$I=\widetilde\Simplex_1$, hence
\[\widetilde i_!(I) = i(\widetilde\Simplex_1) =
\begin{tikzcd}[baseline=(A.base),row sep=-7pt,column sep=small]
  \{1\} \ar[dr] & \\ & |[alias=A]| \{0,1\} \\ \{0\}\ar[ur] &
\end{tikzcd}.\]
We\pspage{76} map this into the object $\Simplex_1$ of \Cat, by
taking $\{0\}$ into $0$, $\{1\}$ and $\{0,1\}$ into $1$, we do get
indeed a morphism compatible with endpoints, which implies that
$\widetilde i$ is a test functor.

If we denote still by $i,\widetilde i$ the functors from
$\Simplex,\widetilde\Simplex$ to the category \Ord\ of ordered sets
factoring the previous two test functors, we get a commutative diagram
of functors
\[\begin{tikzcd}[baseline=(O.base)]
  \Simplex\ar[r]\ar[d,swap,"i"] & \widetilde\Simplex\ar[d,"\widetilde i"] \\
  \Ord\ar[r,"\text{bar}"] & |[alias=O]| \Ord
\end{tikzcd},\]
where the first horizontal arrow is the inclusion functor (bijective
on objects, and injective but not bijective on arrows), and the second
is the ``barycentric subdivision'' functor, or ``flag''-functor,
associating to every ordered set the set of all ``flags'', namely
non-empty subsets which are totally ordered for the induced order
(here, all subsets, as the simplices are totally ordered).

For a while I thought there was an interesting third variant, namely
ordered simplices with \emph{strictly} increasing maps between them --
which means ruling out degeneracy operators. This feeling was prompted
of course by the fact that the face operators in a complex are enough
for computing homology and cohomology groups, which are felt to be
among the most important invariants of a complex. Equally, the
fundamental groupoid of a semi-simplicial set can be described, using
only the face operators. As a consequence, for a map between
semisimplicial complexes $K_*\to K_*'$, to check whether this is a
weak equivalence, in terms of the Artin-Mazur cohomological criterion,
depends only on the underlying map between ``simplicial face
complexes'' (namely, forgetting degeneracies). These \emph{are} indeed
striking facts, which will induce us to put greater emphasis on the
face operators than on degeneracies. It seems, though that the
degeneracies play a stronger role than I suspected, even though it is
a somewhat hidden one. In any case, as soon as we try to check (``par
acquit de conscience'') that the category \Simplexf\ of simplices with
strictly increasing maps is a test category, it turns out that it is
very far from it! Thus, as there is no map from any $\Simplex_n$ with
$n>1$ into $\Simplex_0$, it follows that in
$(\Simplexf)\uphat$, we get
\[\Simplexf_{/\Simplex_0\times\Simplex_n} = \text{discrete category with
  $n+1$ elements,}\]
thus these products are by no means aspheric, poor them! Even throwing
out $\Simplex_0$ (a barbarous thing to do anyhow!) doesn't rule out the
trouble, and restricting moreover to products $\Simplex_1\times\Simplex_n$
(to have at least a test category, if not a strict one). In any case,
$\Simplex_1$ wouldn't be of much use, because it has got no ``section''
anymore (nor does any other element of \Simplexf) -- because this would
imply that \emph{any} element of \Simplexf{} maps into it -- but for
given $\Simplex_n$,\pspage{77} only the $\Simplex_m$'s with $m\le n$
map into it.

Maybe I am only being imprisoned still by the preconception of finding
a homotopy interval in \Simplexf{} itself, rather than in
$(\Simplexf)\uphat$. After all, just applying the definition of a
test-category $A$ with test-functor $i:A\to\Cat$, all we have to care
about is whether a)\enspace $A$ is aspheric and b)\enspace
$i^*(\Simplex_1)$ is an aspheric element of \Ahat. We just got to
apply this to the case of the functor
\[ i^{\mathrm f} : \Simplexf \to \Cat\]
induced by $i:\Simplex\to\Cat$ above, taking into account of course the
extra trouble that $i^{\mathrm f}$ is no longer fully faithful.

I just stopped to look, with a big expectation that \Simplexf{}
\emph{is} a test category after all -- but it turns out it definitely
isn't! Indeed, $A_{/\Simplex_0\times I}$ (where
$I=(i^{\mathrm f})^*(\Simplex_1)$) is again a discrete two-point
category, not aspheric. Taking the canonical functor $a\mapsto A_{/a}$
from $A=\Simplexf$ to \Cat, which is the ultimate choice for checking
whether or not $A$ is a test category, finally gives the answer: it is
not. Because with $I$ now the Lawvere element, we still have that
$A_{/\Simplex_0\times I}$ is a two-point discrete category. Thus the
topos $\Ahat=(\Simplexf)\uphat$ isn't locally modelizing, i.e., it
hasn't got any homotopy interval, which at any rate is a very big
drawback I would think. The only hope which still remains, to account
for the positive features of face-complexes recalled above, is that
\Simplexf{} is at least a test category \emph{in the wide sense} (or weak test category), namely
that for any category $C$ with final element, the category
$i_Aj_A(C) = A_{/j_A(C)}$ ($=$ category of all pairs $(n,u)$, with $u$
a map of the (ordered) category $A_{/\widetilde\Simplex_n}$) is
aspheric. The first case to check is for
$C=\text{final category $\Simplex_0$}$, i.e., asphericity of $A$, next
step would be $C=\Simplex_1$, i.e., asphericity of the Lawvere element
$L_A$ (but of course \emph{not} asphericity over the final element
$e_\Ahat$!).

The question certainly deserves to be settled. If the answer is
affirmative, i.e., \Simplexf{} \emph{is} a test category in the wide
sense, then the proposition stated earlier (\S \ref{sec:42}), which clearly applies equally when $A$ is a
test category in the wide sense, implies that if $M$ is any full
subcategory of \Cat{} containing all those $C$ which are ``locally
isomorphic'' to \Simplexf, (i.e., such that for every $x\in C$, $C_{/x}$
is isomorphic to the ordered category of all subsimplices of some
simplex), then for the induced notion $W_M$ of weak equivalences, $M$
is a modelizer and the inclusion functor from $M$ into \Cat{} is
modelizing. This does not yet apply to \Ord, however, it reopens the
question whether the full subcategory of \Ord{} made up by all ordered
sets $J$ which are locally isomorphic to \Simplexf{} in the sense above
(namely for any $x\in J$, the ordered subset $J_{\le x}$ is
isomorphic\pspage{78} to the ordered set of subsimplices of some
simplex) is a modelizer. To ensure this, it would be enough to find an
ordered set $J$ satisfying the previous condition (for instance on
stemming from a ``simplicial maquette''), such that the corresponding
category is a test category, or at least a test category in the wide
sense. The first candidate that comes to my mind, is to take any
\emph{infinite} set $S$ of vertices, and take $J$ to be the ordered
set of all \emph{finite non-empty} subsets (called the simplices --
thus the elements of $S$ can be interpreted in terms of $J$ as the
minimal simplices). By the way, the category associated to $J$, in
case $S=\bN$, can be interpreted in terms of $A=\Simplexf$ as the
category $A_{/\Simplex_\oo}$, where $\Simplex_\oo$ is defined as the
filtering direct limit in $\Ahat$ of the $\Simplex_n$'s, arranged into a
direct system in the obvious way:
\[\Simplex_\oo = \varinjlim_n \Simplex_n\quad \text{in
  $(\Simplexf)\uphat$.}\]
Asphericity of $J$ looks intuitively evident, and should be easy by a
direct limit argument, as a matter of fact any filtering category (the
next best to having a final element) should be aspheric, at least if
it has a countable cofinal family of objects. The Lawvere element
$L_A$ in \Ahat{} is not aspheric though over the final object, because
when inducing over a zero simplex, we get the same contradiction as
before. As a matter of fact, I am getting aware I have been very silly
and prejudiced not to see one trivial common reason, applicable to
\Simplexf{} as to $J$, showing that they are not test categories nor
even local test categories: namely the induced categories $A_{/a}$
should be test categories too, but among these there are one-point
categories (take $a=\Simplex_0$, or a zero-simplex), and such a category
is \emph{not} a test-category!

Still, $J$ may be a test category in the wide sense -- the basic test
here, as we know already asphericity of $J$ itself, would be
(absolute) asphericity of $L_J$, or equivalently, of the category
$J_{/L_J}$, an ordered set in fact (as is the case for any category
$A_{/F}$ for $A$ defined by an ordered set and $F$ in \Ahat). This is
now the ordered set of all pairs $(K\subset T)$ of finite subsets of
$S$, with $T$ in $J$, namely $T$ non-empty, with the rule
\[ (K',T') \le (K,T) \quad\text{if{f}}\quad \text{$T'\subset T$ and
  $K'=K\sand T'$.}\]

\bigbreak

\presectionfill\ondate{26.3.}\pspage{79}\par

\hangsection{Overall review of the basic notions.}\label{sec:44}%
\renewcommand*{\thesubsection}{\alph{subsection})}%
I finally convinced myself that \Simplexf, the category of standard
ordered simplices with face operations (and no degeneracies) \emph{is}
a ``test category in the wide sense'' after all -- although definitely
not a test category, as was seen yesterday (turning out to be a
practically trivial observation). This now does rehabilitate the
notion of a test category in the wide sense, which I expected to be of
little or no interest -- much the way as previously, the notion I now
call by the name ``test category'' was rehabilitated or rather,
discovered, after I expected that the only one proper notion for
getting modelizers of the form \Ahat{} was in terms of conditions
\ref{it:31.T1} to \ref{it:31.T3} on $A$ (including the very strong
condition \ref{it:31.T2} of total asphericity). This notion I finally
called by the name ``\emph{strict} test category'', and I was
fortunately cautious enough to reserve a name too for the notion which
appeared then as rather weak and unmanageable, of test categories in
the wide sense, or, as I will say now more shortly, \emph{weak test
  categories} (by which of course I do not mean to exclude the
possibility that it be even a test category), thus getting the trilogy
of notions with strict implications
\[\text{weak test categories} \Leftarrow \text{test categories}
\Leftarrow \text{strict test categories}.\]
In order not to get confused, I will recall what exactly each of these
notions means.

\subsection{Weak test categories.}
\label{subsec:44.a}
For a given small category $A$, we
look at the functor
\begin{equation}
  \label{eq:44.1}
  i_A : \Ahat\to\Cat, \quad F \mapsto A_{/F}, \tag{1}
\end{equation}
which commutes with direct limits, and at the right adjoint functor
\begin{equation}
  \label{eq:44.2}
  j_A = i_A^*:\Cat\to\Ahat, \quad C\mapsto
  j_A(C)=(a\mapsto\Hom(A_{/a},C)). \tag{2}
\end{equation}
We get an adjunction morphism
\begin{equation}
  \label{eq:44.3}
  i_Aj_A(C) \to C\quad\text{in \Cat,}\tag{3}
\end{equation}
and another
\begin{equation}
  \label{eq:44.4}
  F \to j_Ai_AF\quad\text{in \Ahat,}\tag{4}
\end{equation}
functorially in $C$ resp.\ in $F$. Presumably, the following are
equivalent (I'll see in a minute how much I can prove about these
equivalences):
\begin{enumerate}[label=(\roman*)]
\item\label{it:44.a.i}
  The functors \eqref{eq:44.1} and \eqref{eq:44.2} are compatible
  with weak equivalences, and the two induced functors between the
  localized categories
  \begin{equation}
    \label{eq:44.5}
    W_A^{-1}\Ahat \rightleftarrows W^{-1}\Cat \eqdef \Hot\tag{5}
  \end{equation}
  are equivalences.
\item\label{it:44.a.ii}
  As in \ref{it:44.a.i}, and moreover the two equivalences are
  quasi-inverse of each other, with adjunction morphism in
  $W^{-1}\Cat=\Hot$ deduced from \eqref{eq:44.3} by localization.
\item\label{it:44.a.iii}
  As\pspage{80} in \ref{it:44.a.ii}, but moreover the adjunction
  morphism in $W_A^{-1}\Ahat$ for the pair of quasi-inverse
  equivalences in \eqref{eq:44.5} being likewise deduced from
  \eqref{eq:44.4}.
\item\label{it:44.a.iv}
  For any $C$ in \Cat, \eqref{eq:44.3} is a weak equivalences.
\item\label{it:44.a.v}
  Same as \ref{it:44.a.iv}, with $C$ restricted to having a final
  element, i.e., for any such $C$, $j_Ai_A(C)=A_{/j_A(C)}$ is
  aspheric.
\item\label{it:44.a.vi}
  For any $F$ in \Ahat, \eqref{eq:44.4} is a weak equivalence,
  moreover
  \[W_\Cat = j_A^{-1}(W_A),\]
  i.e., a map $f:C'\to C$ in \Cat{} is a weak equivalence if{f}
  $j_A(f)$ is a weak equivalence in \Ahat{} -- which means, by
  definition (more or less) that $i_Aj_A(f)$ is a weak equivalence.
\item\label{it:44.a.vii}
  The functor
  \begin{equation}
    \label{eq:44.6}
    i_Aj_A : \Cat\to\Cat,\quad C\mapsto A_{/j_A(C)},\tag{6}
  \end{equation}
  transforms weak equivalences into weak equivalences, i.e., gives
  rise to a functor $\Hot\to\Hot$, and moreover the latter respects
  final object (i.e., $A$ is aspheric).
\item\label{it:44.a.viii}
  Same as \ref{it:44.a.vii}, but restricting to weak equivalences of
  the type $C\to e$, where $C$ has a final element and $e$ is the
  final object in \Cat{} (the one-point category), plus asphericity
  of $A$.
\end{enumerate}

The trivial implications between all these conditions can be
summarized in the diagram
\[\begin{tikzcd}[baseline=(O.base),math mode=false,%
  row sep=small,arrows=Rightarrow]
  \ref{it:44.a.iii} \ar[r,Leftrightarrow]\ar[d] & \ref{it:44.a.vi} \ar[r] &
  \ref{it:44.a.vii} \ar[r] & \ref{it:44.a.viii} \ar[dl] \\
  \ref{it:44.a.ii} \ar[r]\ar[d] & \ref{it:44.a.iv} \ar[r] &
  \ref{it:44.a.v} & \\
  |[alias=O]| \ref{it:44.a.i} & , & &
\end{tikzcd}\]
the only slightly less obvious implication here is \ref{it:44.a.viii}
$\Rightarrow$ \ref{it:44.a.v}, which is seen by looking at the
commutative square deduced from $C\to e$ ($C$ in \Cat{} is in
\ref{it:44.a.v} namely with final element) by applying
\eqref{eq:44.3}
\[\begin{tikzcd}[baseline=(O.base)]
  i_Aj_A(C)\ar[d]\ar[r] & C\ar[d] \\
  A=i_Aj_A(e)\ar[r] & |[alias=O]| e
\end{tikzcd},\]
by assumption the vertical arrows are weak equivalences, and so is
$A\to e$ (because $A$ is supposed to be aspheric), therefore the same
holds for the fourth arrow left. On the other hand, an easy
asphericity argument showed us (\S \ref{sec:30}) that \ref{it:44.a.v} $\Rightarrow$
\ref{it:44.a.iii}, hence all conditions \ref{it:44.a.ii} to
\ref{it:44.a.viii} are equivalent, and they are equivalent to the
stronger form of \ref{it:44.a.iv}, say
\namedlabel{it:44.a.ivprime}{(iv')}, saying that \eqref{eq:44.3} is
\emph{aspheric} for any $C$ in \Cat. The only equivalence which is not
quite clear yet is that \ref{it:44.a.i} implies the other
conditions. But it is so\pspage{81} if we grant the ``inspiring
assumption'', implying that any autoequivalence of \Hot{} is
isomorphic to the identity functor -- in this case it is clear that
even the weaker form \namedlabel{it:44.a.iprime}{(i')} of \ref{it:44.a.i},
demanding only that $i_Aj_A$ induce an autoequivalence of \Hot{} and
nothing on either factor $i_A$, $j_A$ in \eqref{eq:44.5}, implies
\ref{it:44.a.v}. Also, when we assume moreover $A$ aspheric, it is
clear that \ref{it:44.a.i} (and even \ref{it:44.a.iprime}) implies
\ref{it:44.a.vii}, i.e., all other conditions. Thus, instead of the
assumption on \Hot, it would be enough to know that the particular
autoequivalence of \Hot{} induced by $i_Aj_A$ transforms the final
element of \Hot{} (represented by the element $e$ of \Cat) into an
isomorphic one,\footnote{\alsoondate{27.3.} But this is true for \emph{any}
  equivalence of categories -- I'm really being very dull!} which
looks like a very slight strengthening of \ref{it:44.a.i} indeed.

In any case, the conditions \ref{it:44.a.ii} to \ref{it:44.a.viii} are
equivalent, and equivalent to \ref{it:44.a.i} (or
\ref{it:44.a.iprime}) \emph{plus}\footnote{The ``plus'' is unnecessary
  visibly, see note above.} asphericity of $A$. The formally strongest
form is \ref{it:44.a.iii}, the formally weakest one (with the
exception of \ref{it:44.a.i} or \ref{it:44.a.iprime}) is
\ref{it:44.a.v}, which is also the one which looks the most concrete,
namely amenable to practical verification. This was indeed what at the
very beginning was attractive in the condition, in comparison to the
first one that came to my mind when introducing the notion of a
modelizer and of model preserving functors -- namely merely that $i_A$
should induce an equivalence between the localized categories, or
equivalently, that \Ahat{} should be a ``modelizer'' and $i_A$ should
be model-preserving. That we should be more demanding and ask for
$j_A$ to be equally model preserving crept in first rather timidly --
and I still don't know (and didn't really stop to think) if the first
one implies the other.\footnote{If fact, it is the condition that
  $j_A=i_A^*$, not $i_A$, should be model preserving, which is ``the
  right'' condition on $A$, equivalent to $A$ being a weak test
  category, if $A$ is assumed aspheric.}

In any case, as far as checking a property goes, I would consider
\ref{it:44.a.v} to be \emph{the} handy definition of a weak test
category, wheres \ref{it:44.a.iii} is the best, when it goes to making
use of the fact that $A$ is indeed a weak test category. As for
\ref{it:44.a.i}, it corresponds to the main intuitive content of the
notion, which means that homotopy types can be ``modelled'' by
elements of \Ahat, using $i_A$ for describing which homotopy type is
described by an element $F$ of \Ahat, and using $j_A$ for getting a
model in \Ahat{} for a given homotopy type, described by an object $C$
in \Cat.

\subsection{Test categories and local test categories.}
\label{subsec:44.b}
Let $A$ still be a small category. Then the following conditions are
equivalent:
\begin{enumerate}[label=(\roman*)]
\item\label{it:44.b.i}
  All\pspage{82} induced categories $A_{/a}$ (with $a$ in $A$)
  are weak test categories.
\item\label{it:44.b.ii}
  For any aspheric $F$ in \Ahat, $A_{/F}$ is a weak test category.
\item\label{it:44.b.iii}
  There exists a ``homotopy interval'' in \Ahat, namely an element $I$
  in \Ahat, \emph{aspheric over the final element} $e_\Ahat$ (i.e.,
  such that all products $I\times a$ (for $a$ in $A$) are aspheric),
  endowed with two sections $\delta_0$, $\delta_1$ ever $e=e_\Ahat$,
  i.e., with two subobjects $e_0,e_1$ isomorphic to $e$, such that
  $\Ker(\delta_0,\delta_1)=\varnothing_\Ahat$, i.e., $e_0\sand
  e_1=\varnothing_\Ahat$.
\item\label{it:44.b.iv}
  The Lawvere element $L_A$ or $L_\Ahat$ in \Ahat{} (i.e., the
  presheaf $a\mapsto$ subobjects of $a$ in $\Ahat\simeq$ cribles of
  $A_{/a}$) is aspheric over the final element $e$ (in other words, it
  is a homotopy interval, when endowed with the two sections
  $\delta_0,\delta_1$ corresponding to the full and the empty crible
  in $A$).
\item\label{it:44.b.v}
  For any $C$ in \Cat, the canonical map in \Cat
  \begin{equation}
    \label{eq:44.7}
    i_Aj_A(C) = A_{/j_A(C)} \to A \times C \tag{7}
  \end{equation}
  (with second component \eqref{eq:44.3}) is (not only a weak
  equivalence, which definitely isn't enough, but even) \emph{aspheric}.
\end{enumerate}

By the usual \hyperref[lem:asphericitycriterion]{criterion of
  asphericity} for a map in \Cat, condition \ref{it:44.b.v} is
equivalent with the condition \namedlabel{it:44.b.vprime}{(v')}: For
any $a$ in $A$ and any $C$ in \Cat{} with final element, $A_{/a\times
  j_A(C)}$ is aspheric, i.e., $a\times j_A(C)$ is an aspheric element
in \Ahat; now this is clearly equivalent to \ref{it:44.b.i} (by the
checking-criterion \ref{it:44.a.v} above for weak test categories,
applied to the categories $A_{/a}$). Thus we get the purely formal
implications
\[\begin{tikzcd}[baseline=(O.base),math mode=false,arrows=Rightarrow]
  \ref{it:44.b.ii} \ar[r] & \ref{it:44.b.i}\ar[d,Leftrightarrow]
  \ar[r] & \ref{it:44.b.iv} \ar[r] & \ref{it:44.b.iii} \\
  & |[alias=O]| \ref{it:44.b.v} & &
\end{tikzcd},\]
where \ref{it:44.b.i} $\Rightarrow$ \ref{it:44.b.iv} is obtained by
applying the criterion \ref{it:44.b.vprime} just recalled to the case
of $C=\Simplex_1$. As the condition \ref{it:44.b.iii} is clearly stable
by localizing to a category $A_{/F}$ (using the equivalence
$(A_{/F})\uphat\simeq \Ahat_{/F}$), we see that \ref{it:44.b.iii}
implies \ref{it:44.b.ii}, we are reduced (replacing $A$ by $A_{/F}$)
to the case when $F$ is the final element in \Ahat, and thus to
proving that if $A$ has a final element, then \ref{it:44.b.iii}
implies that $A$ is a weak test category, i.e., that the categories
$A_{/j_A(C)}$, with $C$ having a final element, are aspheric. This was
done by a simple homotopy argument in two steps. One step (``the
comparison lemma for homotopy intervals'' on \S \ref{sec:37}) shows that \ref{it:44.b.iii} $\Rightarrow$
\ref{it:44.b.iv}, i.e., existence of a homotopy interval implies that
the Lawvere interval is a homotopy interval, the other step
(presented in a more general set-up in the
``\hyperref[thm:keyresult]{key result}'' in same \S \ref{sec:31}) proving that \ref{it:44.b.iv} implies
asphericity of the elements $j_A(C)$ ($C$ with\pspage{83} final
object) over $e_A=e_\Ahat$.

We express the conditions \ref{it:44.b.i} to \ref{it:44.b.v} by saying
that $A$ is a \emph{local test category}, or that \Ahat{} is a
\emph{locally modelizing topos}, or (if we want to recall that this
topos is of the type \Ahat) an \emph{elementary local modelizer}. If
$A$ is moreover aspheric (i.e., $e_\Ahat$ aspheric), or what amounts
to the same, if $A$ is moreover a weak test category, we say that $A$
is a \emph{test category}, or that \Ahat{} is a \emph{modelizing
  topos}, or (to recall it is an \Ahat{} and not just any topos) an
\emph{elementary modelizer}. These are intrinsic properties on the
topos \Ahat, the first one of a local nature, the second not (as
asphericity of \Ahat{} is a global notion).

Here the question arises whether the condition for $A$ to be a weak
test category can be likewise expressed intrinsically as a property of
the topos \Ahat{} -- which we then would call a \emph{weakly
  modelizing topos}, or a \emph{weak elementary modelizer}. This
doesn't look so clear, as all conditions stated in \ref{subsec:44.a}
make use of at least one among the two functors $i_A, j_A$, which do
not seem to make much sense in the more general case, except possibly
when using a specified small generating subcategory $A$ of the given
topos -- and possibly checking that the condition obtained (if
something reasonable comes out, as I do expect) does not depend on the
choice of the generating site. This should be part of a systematic
reflection on modelizing topoi, to make sure for instance they are
modelizing indeed with respect to weak equivalences -- but I'll not
enter into such reflection for the time being.

\subsection{Strict test categories.}
\label{subsec:44.c}
The following conditions on the small category $A$ are equivalent:
\begin{enumerate}[label=(\roman*)]
\item\label{it:44.c.i}
  $A$ is a weak test category (cf.\ \ref{subsec:44.a}
  \ref{it:44.a.iii} above), and thus induces a localization functor
  \[ \Ahat \to \Hot=W^{-1}\Cat,\]
  and this functor moreover \emph{commutes with binary products}; or
  equivalently, the canonical functor $\Ahat\to W_A^{-1}\Ahat$
  commutes with binary products.
\item\label{it:44.c.ii}
  $A$ is a test category (a local test category would be enough,
  even), and moreover the topos \Ahat{} is totally aspheric, namely
  (apart having a generating family of aspheric generators, which is
  clear anyhow) the product of any two aspheric elements in \Ahat{} is
  aspheric, i.e., any aspheric element of \Ahat{} is aspheric over the
  final object.
\item\label{it:44.c.iii}
  $A$\pspage{84} satisfies the two conditions (besides $A \neq \emptyset$):
  \begin{enumerate}[label=T~\arabic*),start=2]
  \item\label{it:44.T2}
    The product in \Ahat{} of any two elements in $A$ is aspheric.
  \item\label{it:44.T3}
    \Ahat{} admits a homotopy interval (which is equivalent to saying
    that $A$ is a local test category).
  \end{enumerate}
\end{enumerate}

Of course, \ref{it:44.c.ii} implies \ref{it:44.c.iii}. The condition
\ref{it:44.c.iii} can be expressed in terms of the topos $\scrA=\Ahat$
by saying that this is a totally aspheric and locally modelizing
topos, which implies already (as we saw by the \v Cech computation of
cohomology) that the topos is aspheric, and therefore modelizing,
i.e., $A$ a test category, which is just \ref{it:44.c.ii}. Thus
\ref{it:44.c.ii} $\Leftrightarrow$ \ref{it:44.c.iii}. Total
asphericity of \scrA, and the property that \scrA{} be a local
modelizer, are expressed respectively by the two neat conditions
\ref{it:44.T2} and \ref{it:44.T3}, neither of which implies the other
even for an \Ahat, with $A$ a category with final object.

For expressing in more explicit terms condition \ref{it:44.c.i}, and
check equivalence with the two equivalent conditions \ref{it:44.c.ii},
\ref{it:44.c.iii}, we need to admit that the canonical functor from
\Cat{} to its localization \Hot{} commutes with binary products\footnote{This is indeed easily checked, and will be done in part VIII (volume 2 of Pursuing Stacks)}. This
being so, the condition of commutation of $\Ahat\to\Hot$ with binary
products can be expressed by the more concrete condition that for
$F,G$ in \Ahat{}, the canonical map in \Cat
\begin{equation}
  \label{eq:44.8}
  i_A(F\times G) \to i_A(F) \times i_A(G), \quad\text{i.e.,
    $A_{/F\times G}\to A_{/F}\times A_{/G}$}\tag{8}
\end{equation}
is a weak equivalence. In fact, the apparently stronger condition that
the maps \eqref{eq:44.8} are aspheric will follow, because the usual
criterion of asphericity for a map in \Cat{} shows that it is enough
for this to check that \eqref{eq:44.8} is a weak equivalence when $F$
and $G$ are representable by elements $a$ and $b$ in $A$, which also
means that $A_{/a\times b}$ is aspheric, namely $a\times b$ in \Ahat{}
aspheric -- which is nothing but condition \ref{it:44.T2} in
\ref{it:44.c.iii}, i.e., total asphericity of \Ahat. On the other
hand, the condition that $A$ be a weak test category means that the
elements $j_A(C)$ in \Ahat, for $C$ in \Cat{} admitting a final
object, are aspheric, or (what amount to the same when \Ahat{} is
totally aspheric) that they are aspheric over the final object
$e_\Ahat$, i.e., that the products $a\times j_A(C)$ are aspheric,
which also means $A$ is a local test category. Thus \ref{it:44.c.i} is
equivalent to \ref{it:44.c.ii}.

The equivalent conditions \ref{it:44.c.i} to \ref{it:44.c.iii} are
expressed by saying that $A$ is a \emph{strict test category}, or that
$\scrA=\Ahat$ is a \emph{strictly modelizing topos}, or also that it
is an \emph{elementary strict modelizer} (when emphasizing the topos
\scrA{} should be of the type \Ahat{} indeed). This presumably will
turn out to be the more important among the three notions of a ``test
category'' and the weak and strict variant. The two less stringent
notions, however, seem interesting in their own right. The notion of a
weak test category mainly (at present) because it turns out that the
category\pspage{85} \Simplexf{} of standard ordered simplices with
only ``face-like'' maps between them, namely strictly increasing ones
(that is, ruling out degeneracies) is a weak test category, and not a
test category. On the other hand, for any test category $A$, we can
construct lots of test categories which are not strict, namely all
categories $A_{/F}$ where $F$ is any aspheric object in \Ahat{} which
admits two ``non-empty'' subobjects whose intersection is ``empty'' --
take for instance for $F$ any homotopy interval. In the case of
$A=\Simplex$, we may take for $F$ any objects $\Simplex_n$ ($n\ge1$) in
$A$, except just the final object $\Simplex_0$.

\subsection[Weak test functors and test functors (with values in
\texorpdfstring{\Cat}{(Cat)}).]{Test functors.}

Let $A$ be a weak test category, and let
\begin{equation}
  \label{eq:44.9}
  i : A \to \Cat\tag{9}
\end{equation}
be a functor, giving rise to a functor
\begin{equation}
  \label{eq:44.10}
  i^*: \Cat\to\Ahat,\quad C\mapsto i^*(C)=(a\mapsto\Hom(i(a),C)).\tag{10}
\end{equation}
We'll say that $i$ is a \emph{weak test functor} if $i^*$ is a
morphism of modelizers, i.e., model-preserving, namely
\[
\begin{cases}
  W_\Cat=(i^*)^{-1}(W_\Ahat) \quad\text{and} \\
  W_\Cat^{-1}\Cat = \Hot \to W_\Ahat^{-1}\Ahat\text{ an equivalence.}
\end{cases}\]
As by assumption we know already that \eqref{eq:44.1}
\[ i_A : \Ahat\to\Cat, F\mapsto A_{/F}\]
is model-preserving, this implies that $i$ is a weak test functor
if{f} the compositum
\begin{equation}
  \label{eq:44.11}
  i_Ai^*:\Cat\to\Cat,\quad C\mapsto A_{i^*(C)} \eqdef A_{/C}\tag{11}
\end{equation}
is model-preserving, i.e., (essentially) induces an
\emph{autoequivalence} of \Hot.

The basic example is to take $i=i_A$, hence $i^*=j_A$ \eqref{eq:44.2},
the fact that $A$ is a weak test category can be just translated (p.~%
\ref{p:80}, \ref{it:44.a.vii}) by saying that this functor
\eqref{eq:44.11} $=$ \eqref{eq:44.6} transforms weak equivalences into
weak equivalences, and the localized functor $\Hot\to\Hot$ is an
equivalent (or, which is enough, transforms final object into final
object of \Hot). But as we say before, this weak test functor gives
rise to categories $i(A)=A_{/a}$ which are prohibitively large, and
one generally prefers working with weak test functors more appropriate
for computations -- including\pspage{86} the customary ones for
categories of simplices or cubes, such as \Simplexf, $\Simplex$,
$\widetilde\Simplex$ and the cubical analogons. If we admit the
``inspiring assumption'' that any autoequivalence of \Hot{} is
uniquely isomorphic to the identity functor, it will follow that the
autoequivalence $\Hot\to\Hot$ induced by \eqref{eq:44.11} (where again
$i$ is any weak test functor) is canonically isomorphic to the
identity. This we checked directly (as part of the
``\hyperref[thm:keyresult]{key result}''
p.~\ref{p:61} and following) when we make on $i$ the
extra assumption that each category $i(a)$ has a final object. In
practical terms, \emph{the role of a weak test functor $i$ is to
  furnish us with a quasi-inverse of}
\[W_A^{-1}\Ahat \to W_\Cat^{-1}\Cat = \Hot\]
\emph{induced by $i_A$, more handy than the one deduced from $j_A$ by
  localization, namely taking $i^*$ instead of $j_A=i_A^*$}. In other
terms, for every homotopy type, described by an object $C$ in \Cat, we
get a ready model in \Ahat, just taking
$i^*(C)=(a\mapsto\Hom(i(a),C))$.

In the theorem on p.~\ref{p:61} just referred to, one
point was that we did not assume beforehand that $A$ was a weak test
category, but we were examining two sets \hyperref[it:key.a]{(a)} and
\hyperref[it:key.b]{(b)} of mutually equivalent conditions, where
\hyperref[it:key.a]{(a)} just boils down to $i_Ai^*$ \eqref{eq:44.11}
inducing an autoequivalence of \Hot, while \hyperref[it:key.b]{(b)}
would be expressed (in the present terminology) by stating that
moreover $A$ is a test category (if we assume beforehand in
\hyperref[it:key.b]{(b)} that $A$ is aspheric, plus a little more
still on $i$). It was not quite clear by then, as it is now, that this
is actually strictly stronger. The somewhat bulky statement
essentially reduces, with the present background, to the following
\begin{proposition}
  Let $A$ be a small category, $i:A\to\Cat$ a functor, such that for
  any $a$ in $A$\kern1pt, $i(a)$ has a final object. Consider the following
  conditions:
  \begin{enumerate}[label=(\roman*),font=\normalfont]
  \item\label{it:44.prop.i}
    $A$ is a test category \textup(NB\enspace not only a weak one\textup), and
    $i$ is a weak test functor.
  \item\label{it:44.prop.ii}
    The element $i^*(\Simplex_1)$ in \Ahat{} \textup(namely $a \mapsto
    \Crib(i(a))$\textup) is aspheric over $e_\Ahat$, i.e., all
    products $i^*(\Simplex_1)\times a$ \textup($a\in\Ob A$\textup) are
    aspheric; moreover $e_\Ahat$ is aspheric, i.e., $A$ itself is
    aspheric.
  \item\label{it:44.prop.iii}
    There exists a homotopy interval $(I,\delta_0,\delta_1)$ in \Ahat,
    and a map in \Cat{} compatible with $\delta_0,\delta_1$
    \begin{equation}
      \label{eq:44.12}
      i_!(I) \to \Simplex_1,\tag{12}
    \end{equation}
    where $i_!:\Ahat\to\Cat$ is the canonical extension of
    $i:A\to\Cat$; moreover, $A$ is aspheric.
  \end{enumerate}
  The conditions \textup{\ref{it:44.prop.ii}} and \textup{\ref{it:44.prop.iii}} are
  equivalent and imply \textup{\ref{it:44.prop.i}}.
\end{proposition}

That \ref{it:44.prop.ii} implies \ref{it:44.prop.iii} is trivial by
taking $I=i^*(\Simplex_1)$; the converse, as we saw, is a corollary of
the comparison lemma for homotopy intervals (p.~%
\ref{p:60}), applied to $I\to i^*(\Simplex_1)$ corresponding
to \eqref{eq:44.12}, which is\pspage{87} a morphism of intervals,
namely compatible with endpoints. (However, in practical terms
\eqref{eq:44.12} is the more ready-to-use criterion, because in most
cases $I$ will be in $A$ and therefore $i_!(I)=i(I)$, and
\eqref{eq:44.12} will be more or less trivial, for instance because
$I$ will be chosen so that $i(I)$ is either $\Simplex_1=\{0\}\to\{1\}$
or its barycentric subdivision \begin{tikzcd}[baseline=(A.base),row sep=-12pt,column sep=small]
  \{0\} \ar[dr] & \\ & |[alias=A]| \{0,1\} \\ \{1\}\ar[ur] &
\end{tikzcd}.) The same argument shows that \ref{it:44.prop.iii}
implies that $A$ is a test category, namely the Lawvere element $L_A$
in \Ahat{} is a homotopy interval in \Ahat, namely is aspheric over
$e_\Ahat$. This being so, we only got to prove (under the assumption
\ref{it:44.prop.iii}) that $i_Ai^*$ \eqref{eq:44.11} induces an
autoequivalence of \Hot, and more specifically a functor isomorphic to
the identity functor. This we achieved by comparing directly $A_{/C}$
to $C$ by the functor $(a,p)\mapsto p(e_a)$
\begin{equation}
  \label{eq:44.13}
  A_{/C} \to C,\tag{13}
\end{equation}
where $e_a$ is the final element in $i(a)$. It is enough to check this
is a weak equivalence for any $C$, a fortiori that it is aspheric, and
the usual asphericity criterion reduces us to showing it is a weak
equivalence when $C$ has a final object, namely that in this case
$A_{/C}=A_{/i^*(C)}$ is aspheric, i.e., $i^*(C)$ is aspheric in
\Ahat. This was achieved by an extremely simple homotopy argument,
using a homotopy in \Cat
\[\Simplex_1\times C\to C\]
between the identity functor in $C$ and the constant functor with
value $e_C$ (cf.\ p.~\ref{p:62}).

I feel like reserving the appellation of a \emph{test functor}
$A\to\Cat$ (in contract to a \emph{weak} test functor) to functors
satisfying the equivalent conditions \ref{it:44.prop.ii} and
\ref{it:44.prop.iii}, implying that $A$ is a test category (not merely
a weak one!), and which seem stronger a priori than assuming moreover
that $i$ is a weak test functor (condition \ref{it:44.prop.i}). I
confess I did not make up my mind if for a test category $A$, there
may be weak test functors, with all $i(a)$ having final objects, which
are not test functors in the present sense, namely whether it may be
true that for any $C$ in \Cat{} with final object, $i^*(C)$ is
aspheric, and therefore $i^*(\Simplex_1)$ aspheric, without the latter
being aspheric over $e_\Ahat$, i.e., all products $i^*(\Simplex_1)\times
a$ ($a\in\Ob A$) being aspheric (which would be automatic though if
$A$ is even a \emph{strict} test category). The difficulty here is
that it looks hard to check the condition for arbitrary $C$ with final
object, except through the homotopy argument using relative
asphericity of $i^*(\Simplex_1)$. Of course, the distinction between
weak test functors and test functors is meaningful only as long as it
is expected that the two are actually distinct, when applied to test
categories. Besides this, the restriction to categories $i(a)$ with
final objects looks theoretically a little awkward, and it shouldn't
be too hard I believe to get rid of it, if need be. But for the time
being the only would-be test functors or\pspage{88} weak test
functors which have turned up, do satisfy this condition, and
therefore it doesn't seem urgent to clean up the notion in this
respect.

\needspace{6\baselineskip}
\subsection*{Remarks.}

\paragraph{1.} I regret I was slightly floppy when translating the
condition that $\varphi=i_Ai^*$ \eqref{eq:44.11} be model-preserving,
by the (a priori less precise) condition that it induce an
autoequivalence of the localized category \Hot{} (which is of course
meant to imply that it transforms weak equivalences into weak
equivalences, and thus does induce a functor from \Hot{} into
itself). It isn't clear that if this functor is an equivalence,
$W=\varphi^{-1}(W)$, except if we admit that $W$ is strongly
saturated, i.e., an arrow in \Cat{} which becomes an isomorphism in
\Hot, is indeed a quasi-equivalence. This, as we saw, follows from the
fact that there are weak test categories (such as $\Simplex$, which is
even a strict one), such that \Ahat{} be a ``closed model category''
in Quillen's sense -- but I didn't check yet in the present set-up
that \Simplexhat{} does indeed satisfy Quillen's condition. \emph{Which}
elementary modelizers are ``closed model categories'' remains one of
the intriguing questions in this homotopy model story, which I'll have
to look into pretty soon now.

\paragraph{2.} It should be noted that if $i:A\to\Cat$ is a test
functor, with a strict test category $A$ (such as $\Simplex$, with the
standard embedding $i$ of $\Simplex$ into \Cat{} say), whereas
$i^*:\Cat\to\Ahat$ is model preserving by definition, it is by no
means always true that the left adjoint functor
\[i_! : \Ahat\to\Cat\]
(which can be equally defined as the canonical extension of $i$ to
\Ahat, commuting with direct limits) is equally model preserving. This
seems to be in fact an extremely special property of just the
``canonical'' test functor $i_A\restrto A$.

\paragraph{3.} On the other hand, I do not know if for any small
category $A$, such that $(A, W_A=W_\Ahat)$ is a modelizer and
$i_A:\Ahat\to\Cat$ is model preserving, is a weak test
category. Assume that $j_A=i_A^*:\Cat\to\Ahat$ is model preserving,
and moreover $A$ aspheric, then $A$ is a weak test category, because
$j_A$ transforms weak equivalences into weak equivalences and we apply
criterion \ref{it:44.a.vii} (p.~\ref{p:80}). More
generally, if we have a functor $i:A\to\Cat$ such that
$i^*:\Cat\to\Ahat$ is model preserving (assuming already $(\Ahat,W_A)$
to be a ``modelizer''), I wonder whether this implies that $A$ is
already a weak test category, as it does when $i$ is the canonical
functor $i_A$ and hence $i^*=j_A$. The answer isn't clear to me even
when $i_A$, or equivalently $\varphi=i_Ai^*$ model-preserving too.



\chapter{Homotopy structures and contractibility structures}
\label{ch:III}

\presectionfill\ondate{27.3.}\pspage{89}\par

\hangsection[It's burning again! Review of some ``recurring striking
\dots]{It's burning again! Review of some ``recurring striking
  features'' of modelizers and standard modelizing functors.}%
\label{sec:45}%
The review yesterday of the various ``test notions'', turning around
test categories and test functors, turned out a lot longer than
expected, so much so to have me get a little weary by the end -- it
was clear though that this ``travail d'intendance'' was necessary, not
only not to get lost in a morass of closely related and yet definitely
distinct notions, but also to gain perspective and a better feeling of
the formal structure of the whole set-up. As has been the case so
often, during the very work of ``grinding through'', there has
appeared this characteristic feeling of getting close to something
``burning'' again, something very simple-minded surely which has kept
showing up gradually and more and more on all odds and ends, and which
still is escaping, still elusive. These is an impressive bunch of
things which are demanding pressingly more or less immediate
investigation -- still I can't help, I'll have to try and pin down
some way or other this ``burning spot''.

There seem to be recurring striking features of the modelizers met
with so far -- namely essentially \Cat{} and the elementary modelizers
\Ahat{} and possibly their ``weak'' variants. In all of them, these is
a very strong interplay between the following notions, which seem to
be the basic ones and more or less determine each other mutually:
\emph{weak equivalences} (which define the modelizing structure of the
given modelizer $M$), \emph{aspheric objects} (namely such that
$x\to e_M$ is a weak equivalence), \emph{homotopy intervals}
$(I,\delta_0,\delta_1)$, and last not least, the notion of a test
functor $A \to M$, where $A$ is a test category (or more generally
weak test functors of weak test categories into $M$). The latter so
far have been defined only when $M=\Cat$, and initially they were
viewed as being mainly more handy substitutes to $j_A$, for getting a
model-preserving functor $\Cat\to\Ahat$ quasi-inverse to the
all-important model-preserving functor $i_A:\Ahat\to\Cat$,
$F\mapsto A_{/F}$. I suspect however that their role is a considerably
more basic one than just computational convenience -- and this reminds
me of the analogous feeling I had, when first contemplating using such
a thing as (by then still vaguely conceived) ``test categories'', for
investigating ways of putting modelizing structures on categories $M$
such as categories of algebraic structures of some kind or
other. (Cf.\ notes of \hyperref[date:7.3.]{March 7}, and more
specifically par.\ \ref{sec:26} -- this was the very day, by the way,
I first had this feeling of being ``burning''\ldots)

Test\pspage{90} categories seem to play a similar role here as (the
spectra of) discrete valuation rings in algebraic geometry -- they can
be mapping into anywhere, to ``test'' what is going on there -- here
it means, they can be sent into any modelizer $(M,W)$ (at least among
the ones which we feel are the most interesting), by ``test functors''
$i:A\to M$ giving rise to a model-preserving functor $i^*:M\to\Ahat$,
allowing comparison of $M$ with an elementary modelizer $\Ahat$. As
for the all-encompassing basic modelizer \Cat, it seems to play the
opposite role in a sense, at least with respect to elementary
modelizers, \Ahat, which all admit modelizing maps
$i_A:\Ahat\to\Cat$. As a matter of fact, for given test category $A$,
i.e., a given elementary modelizer \Ahat, I see for the time being
just \emph{one} way to get a modelizing functor from it to \Cat,
namely just the canonical $i_A$. There is also a striking difference
between the exactness properties of the functors
\begin{equation}
  \label{eq:45.1}
  i^*:M\to\Ahat\tag{1}
\end{equation}
one way, which commute to inverse limits, and the functors
\begin{equation}
  \label{eq:45.2}
  i_A:\Ahat\to\Cat\tag{2}
\end{equation}
in the other direction, commuting to direct limits. Another difference
is that we should not expect that the left adjoint $i_!$ to $i^*$ be
model preserving too (with the exception of the very special case when
$M=\Cat$ and $i:A\to M=\Cat$ is the canonical functor $i_A\restrto A$,
which appears as highly non-typical in this respect), whereas the
right adjoint $j_A=i_A^*$ of $i_A$ is model preserving, this $i_A$ is
part of a pair $(i_A,j_A)$ of model preserving adjoint functors.

Of course, we may want to compare directly an arbitrary modelizer $M$
to \Cat{} by sending it into \Cat{} by a modelizing functor
$M\to\Cat$; we get quite naturally such a functor (for any given
choice of test-functor $i:A\to M$)
\begin{equation}
  \label{eq:45.3}
  i_Ai^*: M \to \Cat,\tag{3}
\end{equation}
but this functor is not likely any more to commute neither to direct
nor inverse limits, even finite ones -- and it isn't too clear that
for a modelizer $M$ which isn't elementary, we have much chance to get
a modelizing functor to \Cat{} which is either left or right
exact. However, the functors \eqref{eq:45.3} we've got, whenever
modelizing and if $M$ is a \emph{strict} modelizer (namely $M\to
W_M^{-1}M\simeq\Hot$ commutes with finite products), will commute to
finite products ``up to weak equivalence''. Also the functors $i^*$,
although not right exact definitely, have a tendency to commute to
sums, and hence the same will hold (not only up to weak equivalence)
for \eqref{eq:45.3}.

As for getting a modelizing functor $\Cat\to M$, for a modelizer
$M$\pspage{91} which isn't elementary, in view of having a standard
way for describing a given homotopy type (defined by an object $C$ in
\Cat) by a ``model'' in $M$, depending functorially on $C$, there
doesn't seem to be any general process for finding one, even without
any demand on exactness properties, except of course when $M$ is
supposed to be elementary; in this case $M=\Ahat$ we get the functors
\begin{equation}
  \label{eq:45.4}
  i^*:\Cat\to\Ahat\tag{4}
\end{equation}
associated to test functors $A\to\Cat$, which can be viewed as a
particular case of \eqref{eq:45.1}, applied to the case
$M=\Cat$. Using such functors \eqref{eq:45.4}, we see that the
question of finding a modelizing functor
\begin{equation}
  \label{eq:45.star}
  \varphi: \Cat\to M,\tag{*}
\end{equation}
for a more or less general $M$, is tied up with the question of
finding such a functor from an elementary modelizer \Ahat{} into $M$
\begin{equation}
  \label{eq:45.starstar}
  \psi:\Ahat\to M.\tag{**}
\end{equation}
More specifically, if we got a $\psi$, we deduce a $\varphi$ by
composing with $i^*$ in \eqref{eq:45.4}, and conversely, if we got a
$\varphi$, we deduce a $\psi$ by composing with the canonical functor
$i_A$ in \eqref{eq:45.2}. Maybe it's unrealistic to expect modelizing
functors \eqref{eq:45.star} or \eqref{eq:45.starstar} to exist for
rather general $M$. (Which modelizers will turn out to be really ``the
interesting ones'' will appear in due course presumably\ldots) There
is one interesting case though when we got such functors, namely when
\[M=\Spaces\]
is the category of topological spaces, and taking for $\psi$ one of
the manifold avatars of ``geometric realization functor'', associated
to a suitable functor
\begin{equation}
  \label{eq:45.starstarstar}
  r:A\to\Spaces\tag{***}
\end{equation}
by taking the canonical extension $r_!$ to \Ahat, commuting with
direct limits. This is precisely the ``highly non-typical'' case, when
we get a pair of adjoint functors $r_!,r^*$
\begin{equation}
  \label{eq:45.5}
  \begin{tikzcd}
    \Ahat\ar[r, shift left, "r_!"] & \Spaces\ar[l, shift left, "r^*"]
  \end{tikzcd}\tag{5}
\end{equation}
which are \emph{both} modelizing. The situation here mimics very
closely the situation of the pair $(i_A,j_A=i_A^*)$ canonically
associated to the elementary modelizer \Ahat, with the ``basic modelizer''
\Cat{} being replaced by \Spaces, which therefore can be considered as
another ``basic modelizer'' of sorts. In this case the corresponding
functor
\begin{equation}
  \label{eq:45.6}
  r_!r^* : \Cat\to\Spaces
  \tag{6}
\end{equation}
mimics the functor $i_Ai^*$ of \eqref{eq:45.3} (where on the left hand
side $M$ is taken to be just \Cat, and on the right \Cat{} as the
basic modelizer is replaced by its next best substitute \Spaces). Here
as in \eqref{eq:45.3}, the modelizing functor we got is neither left
nor right exact, it has a tendency\pspage{92} though to commute to
sums, as usual.

I wouldn't overemphasize the capacity of \Spaces{} to serve the
purpose of a ``basic modelizer'' as does \Cat, despite the attractive
feature of more direct (or at any rate, more conventional) ties with
so-called ``topological intuition''. One drawback of \Spaces{} is the
relative sophistication of the structure species ``topological
spaces'' it corresponds to, which is by no means an ``algebraic
structure species'', and fits into algebraic formalisms only at the
price of detours. More seriously still, or rather as a reflection of
this latter feature, only for some rather special elementary
modelizers \Ahat, namely rather special test categories $A$, do we get
a geometric realization functor $r_!:\Ahat\to\Spaces$ which can be
view as part of a pair of mutually adjoint modelizing functors,
mimicking the canonical pair $(i_A,j_A)$; still less does there seem
to be anything like a really canonical choice (although some choices
are pretty natural indeed, dealing with the standard test categories
such as $\Simplex$ and its variants). At any rate, it is still to be
seen whether there exists such a pair $(r_!,r^*)$ for some rather
general class of test categories -- this is one among the very many
things that I keep pushing off, as more urgent matters are calling for
attention\ldots

To sum up the outcome of these informal reflections about various
types of modelizing functors between modelizers, the two main types
which seem to overtower the whole picture, and are likely to be the
essential ones for a general understanding of homotopy models, are the
two types \eqref{eq:45.1} and \eqref{eq:45.2} above. The first one
$i^*$ is defined in terms of an arbitrary modelizer $M$. The second
$i_A$, with opposite exactness properties to the previous one, is
canonically attached to any test category, and maps the corresponding
elementary modelizer \Ahat{} into the basic modelizer \Cat, without
any reference to more general types of modelizers $M$. The right
adjoint of the latter, which is still model preserving, is in fact of
the type \eqref{eq:45.1} again, for the canonical test functor
$A\to\Cat$ induced by $i_A$, namely $a\mapsto A_{/a}$.

This whole reflection was of course on such an informal level, that
there was no sense at that stage to bother with distinctions between
weaker or stricter variants of the test-notion. Maybe it's about time
now to start getting a little more specific.\pspage{93}

\hangsection[Test functors with values in any modelizer: an
\dots]{Test functors with values in any modelizer: an observation,
  with an inspiring ``silly question''.}\label{sec:46}%
First thing to do visibly is to define the notion of a test-functor
\[ i : A \to M,\]
where $M$ is any modelizer. Thus $M$ is endowed with a subset
$W_M\subset\Fl(M)$, i.e., a notion of weak equivalence, satisfying the
``mild saturation conditions'' of p.\ \ref{p:59}, and moreover we
assume that $W_M^{-1}M$ is equivalent to \Hot{} -- but the choice of
an equivalence, or equivalently, of the corresponding localization
functor
\[ M \to \Hot,\]
is not given with the structure. (If we admit the ``inspiring
assumption'', there is no real choice, as a matter of fact -- but we
don't want to use this in a technical sense, but only as a guide and
motivation.)

Let's start with the weak variant -- we assume $A$ to be a weak test
category, and want to define what it means that $i$ is a \emph{weak
  test functor}. In all this game, it is understood that in case
$M=\Cat$, the notions we want to define (of a weak test functor and of
a test functor) should reduce to the ones we have pinpointed in
yesterday's notes.

The very first idea that comes to mind, is to demand merely that the
corresponding $i^*$ \eqref{eq:45.1}
\[ i^*: M\to\Ahat\]
should be modelizing, which means (I recall)
\begin{enumerate}[label=\alph*)]
\item $W_M= (i^*)^{-1}(W_A)$.
\item The induced functor $W_M^{-1}M \to W_A^{-1}\Ahat$ is an
  equivalence.
\end{enumerate}

This, I just checked, does correspond to the definition we gave
yesterday (p.\ \ref{p:85}), when $M=\Cat$. There is a very interesting
extra feature though in this special case, which appears kind of ``in
between the lines'' in the ``\hyperref[thm:keyresult]{key result}'' on
p.\ \ref{p:61}, and which I want now to state in the more general
set-up.

As usual in related situations, the notion of weak equivalence in $M$
gives rise to a corresponding notion of ``aspheric'' elements in $M$
-- namely those for which the unique map
\[x \to e_M\]
is a weak equivalence. We assume now the existence of a final object
$e_M$ in $M$, and will assume too, if necessary, that it's the image
in the localization $W_M^{-1}M= H_M$ is equally a final object. Thus,
if $x$ in $M$ is aspheric, its image in $H_M$ is a final object, and
the converse holds provided as assume $W_M$ strongly saturated, namely
any map in $M$ which becomes an isomorphism in $H_M$ is a weak
equivalence.\footnote{\alsoondate{29.3.} This condition will be verified if there
  exists a weak test functor $i:A\to M$.}

I\pspage{94} can now state the ``interesting extra feature''.
\begin{observation}
  For a functor $i:A\to M$ of a weak test category $A$ into the
  modelizer $M$ \textup(with final object $e_M$, giving rise to the
  final object in $H_M=W_M^{-1}M$\textup), and in the special case
  when $M=\Cat$, the following conditions are equivalent:
  \begin{enumerate}[label=(\roman*),font=\normalfont]
  \item\label{it:46.i}
    $i^*$ transforms weak equivalence into weak equivalences,
    i.e., induces a functor $H_M\to H_\Ahat$.
  \item\label{it:46.ii}
    $i^*$ transforms aspheric objects into aspheric objects.
  \item\label{it:46.iii}
    $i$ is a weak test functor, namely $W_M=(i^*)^{-1}(W_\Ahat)$
    \textup(a stronger version of \textup{\ref{it:46.i}}\textup)
    \emph{and} the induced functor $H_M\to H_\Ahat$ is an \emph{equivalence}.
  \end{enumerate}
\end{observation}

Here the obvious implications are of course
\[ \text{\ref{it:46.iii}} \Rightarrow \text{\ref{it:46.i}} \Rightarrow
\text{\ref{it:46.ii},}\]
the second implication coming from the fact that $i^*$ is compatible
with final objects, and that $e_\Ahat$ is aspheric. Of course,
\ref{it:46.ii} means that for any aspheric $x$ in $M$, $A_{/i^*(x)}$
is aspheric in \Cat. In case $M=\Cat$, and when moreover the elements
$i(a)$ in \Cat{} have final objects (a condition I forgot to include
in the statement of the observation above, sorry), this condition was
seen to imply \ref{it:46.iii} (cf.\ ``\hyperref[thm:keyresult]{key result}'' on
p.\ \ref{p:61}, \hyperref[it:key.a.iv]{(a~iv)} $\Rightarrow$
\hyperref[it:key.a.ii]{(a~ii)} -- indeed, it is even enough to check
that for any $C$ \emph{with final element} in \Cat, $i^*(C)$ is
aspheric. The proof moreover turns out practically trivial, in terms
of the usual \ref{lem:asphericitycriterion} for a functor between
categories. So much so that the really amazing strength of the
statement, which appears clearly when looked at in a more general
setting, as I just did, was kind of blurred by the impression of
merely fastidiously grinding through routine equivalences. We got
there at any rate quite an interesting class of functors between
modelizers (an elementary and the basic one, for the time being), for
which the mere fact that the functor be compatible with weak
equivalences, or only even take aspheric objects into aspheric ones,
implies that the functor in modelizing, namely that the functor
$H_M\to H_\Ahat$ it induces (and the very \emph{existence} of this
functor was all we demanded beforehand!) is actually an
\emph{equivalence of categories}.

The question that immediately comes to mind now, is if this ``extra
feature'' is indeed an extremely special one, strongly dependent on
the assumption $M=\Cat$ and the categories $i(a)$ having final objects
-- or if it may not have a considerably wider significance. This
suggests the still more general questions, involving two modelizers
$M,M'$, neither of which needs by elementary or by \Cat{}
itself:\pspage{95}
\begin{question}\label{q:naivequestion}
  Let
  \[f : M\to M'\]
  be a functor between modelizers $(M,W)$ and $(M',W')$, assume if
  needed that $f$ commute with inverse limits, or even has a left
  adjoint, and that inverse limits (and direct ones too, as for that!)
  exist in $M,M'$. Are there some natural conditions we can devise for
  $M$ and $M'$ (which should be satisfied for elementary modelizers
  and for the basic modelizer \Cat), plus possibly some mild extra
  conditions on $f$ itself, which will ensure that whenever $f$
  transforms weak equivalences into weak equivalences, or even only
  aspheric objects into aspheric objects, $f$ is model-preserving,
  i.e., $W_M = f^{-1}(W_{M'})$ \emph{and} the induced functor
  $H_f:H_M\to H_{M'}$ on the localizations is an equivalence of categories?
\end{question}

Maybe it's a silly question, with pretty obvious negative answer -- in
any case, I'll have to find out! The very first thing to check is to
see what happens in case of a functor
\[i^*:M\to\Ahat,\]
where \Ahat{} is a weak elementary modelizer, and where $M$ is either
\Cat{} or another weak elementary modelizer $B\uphat$, $i^*$ in any
case being associated of course to a functor
\[i:A\to M,\]
with a priori no special requirement whatever on $i$. In case
$M=\Cat$, this means looking up in the end the question we have
postponed for quite a while now, namely of how to rescue the
``\hyperref[thm:keyresult]{key result}'' of p.\ \ref{p:61}, when
dropping the assumption that the categories $i(a)$ (for $a$ in $A$)
have final objects. We finally got a strong motivation for carrying
through a generalization, if this is indeed feasible.

\bigbreak
\presectionfill\ondate{30.3.}\par

\hangsection[An approach for handling \Cat-valued test functors, and
\dots]{An approach for handling \texorpdfstring{\Cat}{(Cat)}-valued
  test functors, and promise of a ``key result'' revised. The
  significance of contractibility.}\label{sec:47}%
It had become clear that the most urgent thing to do now was to come
to a better understanding of test functors with values in \Cat, when
dropping the assumption that the categories $i(a)$ have final objects,
and trying to replace this (if it should turn out that something
\emph{is} needed indeed) by a kind of assumption which should make
sense when \Cat{} is replaced by a more or less arbitrary modelizer
$M$. I spent a few hours pondering over the situation, and it seems to
me that in the case at least when $A$ is a \emph{strict}, namely when
\Ahat{} is totally aspheric, there is now a rather complete
understanding of the situation, with a generalization of the
``\hyperref[thm:keyresult]{key result}'' of p.\ \ref{p:61} which seems
to be wholly satisfactory.

The\pspage{96} basic idea of how to handle the more general situation,
namely how to compare the categories $A_{/i^*(C)}=A_{/C}$ and $C$, and
show (under suitable assumptions) that there is a canonical
isomorphism between their images in the localized category
$W_\Cat^{-1}\Cat=\Hot$, was around since about the moment I worked out
the ``key result''. It can be expressed by a diagram of ``maps'' in
\Cat
\begin{equation}
  \label{eq:47.1}
    \begin{tikzcd}[cramped]
      A_{/C} \ar[r] & A_{\sslash C} & \\
      & A\times C \ar[u]\ar[r] & C\quad,
    \end{tikzcd}
  \tag{1}
\end{equation}
where $A_{\sslash C}$ is the fibered category over $A$, associated to
the functor
\[ A\op\to\Cat, \quad a\mapsto\bHom(i(a),C).\]
Here one should be careful with the distinction between the \emph{set}
\[ \Hom(i(a),C) = \Ob \bHom(i(a),C),\]
and the \emph{category} $\bHom(i(a),C)$, both depending
bi-functorially on $a$ in $A$ and $C$ in \Cat. The former (as a
presheaf on $A$ for fixed $C$) gives rise to $A_{/C}$, a fibered
category over $A$ with \emph{discrete} fibers, whereas the latter
gives rise to $A_{\sslash C}$, which is fibered over $A$ with fibers
that need not be discrete. Identifying a set with the discrete
category it defines, we get a canonical functor
\[ \Hom(i(a),C) \to \bHom(i(a),C),\]
which is very far from being an equivalence nor even a weak
equivalence; being functorial for varying $a$, it gives rise to the
first map in \eqref{eq:47.1}. The second is deduced from the canonical
functor
\[C \to \bHom(i(a),C),\]
identifying $C$ with the full subcategory of \emph{constant} functors
from $i(a)$ to $C$. This map is functorial in $a$, and gives rise
again to a cartesian functor between the corresponding fibered
categories over $A$, the first one (which corresponds to the constant
functor $A\op\to\Cat$ with value $C$) is just $A\times C$ fibered over
$A$ by $\pr_1$, hence the second arrow in \eqref{eq:47.1}. The third
arrow is just $\pr_2$.

It seems that, with the introduction of $A_{\sslash C}$, this is the
first time since the beginning of these reflections that we are making
use of the \emph{bicategory structure} of \Cat, namely of the notion
of a morphism or map or ``homotopy'' (the tie with actual homotopies
will be made clear below), between two ``maps'' namely (here) functors
$C'\rightrightarrows C$. There is of course a feeling that such a
notion of homotopy should make sense in a more or less arbitrary
modelizer $M$, and that the approach displayed by the diagram
\eqref{eq:47.1} may well generalize to mere general
situations\pspage{97} still, with \Cat{} replaced by such an $M$.

In the situation here, the work will consist in devising handy
conditions on $i: A\to\Cat$ and $A$ that will ensure that all three
maps in \eqref{eq:47.1} are weak equivalences, for any choice of
$C$. This will imply that the corresponding maps in \Hot{} are
isomorphisms, hence a canonical isomorphism between the images in
\Hot{} of $A_{/C}$ and $C$, which will imply that a)\enspace the
functor $i_Ai^*$ from \Cat{} to \Cat{} carries weak equivalences into
weak equivalences, and hence induces a functor
\[\Hot\to\Hot,\]
and b)\enspace that this functor is isomorphic to the identity
functor. If moreover $A$ is a weak test category, and therefore the
functor
\[W_A^{-1}\Ahat\to\Hot\]
induced by $i_A$ is an equivalence, it will follow that the functor
\[\Hot\to W_A^{-1}\Ahat\]
induced by $i^*$ is equally an equivalence, namely that $i^*$ is
indeed a test-functor.

For the map
\begin{equation}
  \label{eq:47.2}
  A\times C\to C
  \tag{2}
\end{equation}
to be a weak equivalence for any $C$ in \Cat, it is necessary and
sufficient that $A$ be aspheric (cf.\ par.\ \ref{sec:40}, page
\ref{p:69}), a familiar condition on $A$ indeed! For handling the map
\begin{equation}
  \label{eq:47.3}
  A\times C\to A_{\sslash C}\quad\text{associated to $C\to\bHom(i(a),C)$,}
  \tag{3}
\end{equation}
we'll use the following easy result (which I'll admit for the time
being):
\begin{proposition}
  Let $F$ and $G$ be two categories over a category $A$\kern1pt, and $u:F\to
  G$ a functor compatible with projections. We assume
  \begin{enumerate}[label=\alph*),font=\normalfont]
  \item\label{it:47.a}
    For any $a$ in $A$\kern1pt, the induced map on the fibers
    \[u_a : F_a \to G_a\]
    is a weak equivalence.
  \item\label{it:47.b}
    Either $F$ and $G$ are both cofibering over $A$ and $u$ is
    cocartesian, or $F$ and $G$ are fibering and $u$ is cartesian.
  \end{enumerate}
  Then $u$ is a weak equivalence.
\end{proposition}
This shows that a sufficient condition for \eqref{eq:47.3} to be a
weak equivalence, is that the functors $C\to\bHom(i(a),C)$ be a weak
equivalence, for any $a$ in $A$ and $C$ in \Cat. We'll see in the next
section\footnote{\alsoondate{9.4.} Actually, this is done only a lot later, on page
  \ref{p:121} and (for the converse) on page \ref{p:143}.} that this
amounts to demanding that the objects $i(a)$ in \Cat{} should be
``contractible'', in the most concrete sense of this expression, which
is actually stronger than just asphericity. (Earlier in these notes I
was a little floppy with the terminology, by using a few times the
word ``contractible''\pspage{98} as synonymous to ``aspheric'' as in
the context of topological spaces, or CW spaces at any rate, the two
notions do indeed coincide, finally I came to use rather the word
``aspheric'' systematically, as it fits nicely with the notion of an
aspheric morphism of topos\ldots) The most evident example of
contractible categories are the categories with final object. Thus I
have the strong feeling that the condition of contractibility of the
objects $i(a)$ in \Cat{} is ``the right'' generalization of the
assumption made in the ``key result'', namely that the $i(a)$ have
final elements. Also, it seems now likely that the numerous cases of
statements, when to check some property for arbitrary $C$, it turned
out to be enough to check it for $C$ with a final object, may well
generalize to more general cases, with \Cat{} replaced by some $M$ and
the reduction is from arbitrary $C$ in $M$ to contractible ones.

The next thing to do is to develop a little the notion of
contractibility of objects and of homotopies between maps, and to get
the criterion just announced for \eqref{eq:47.3} to be a weak
equivalence for any $C$. After this, handling the question of the
first map in \eqref{eq:47.1}
\begin{equation}
  \label{eq:47.4}
  A_{/C} \to A_{\sslash C}
  \tag{4}
\end{equation}
being a weak equivalence for any $C$, in terms of the
\ref{lem:asphericitycriterion} for a functor, will turn out pretty
much formal, and we'll finally be able to state a new version of the
``key result'' about test functors $i:A\to\Cat$, with this time twice
as many equivalent formulations of the same property. On n'arr\^ete
pas le progr\`es!

\bigbreak
\presectionfill\ondate{31.3.}\par

\hangsection[A journey through abstract homotopy notions (in terms
\dots]{A journey through abstract homotopy notions
  \texorpdfstring{\textup(}{(}in terms of a set
  \texorpdfstring{$W$}{W} of ``weak
  equivalences''\texorpdfstring{\textup)}{)}.}\label{sec:48}%
It's time now to develop some generalities about homotopy classes of
maps, the relation of homotopy between objects of a modelizer, and the
corresponding notion of contractibility. For the time being, it will
be enough to start with any category $M$, endowed with a set
$W\subset\Fl(M)$ of arrows (the ``weak equivalences''), satisfying the
mild saturation assumptions \ref{it:37.a}\ref{it:37.b}\ref{it:37.c} of
par.\ \ref{sec:37} (p.\ \ref{p:59}). On $M$ we'll assume for the time
being that there exists (at least one) homotopy interval
$\bI = (I,\delta_0,\delta_1)$ in $M$ (loc.\ sit.) which implies also
that $M$ has a final object, which I denote by $e_M$ or simply $e$,
and equally an initial element $\varnothing_M$. I'm not too sure yet
whether we'll really need that the latter be strict initial element,
as required in the definition of a homotopy interval on page
\ref{p:59} (it was used in the generalities of pages \ref{p:59} and
\ref{p:60} only for the corollary on the Lawvere
element\ldots). Whether or not will appear soon enough! I'll assume it
till I am forced to. To be safe, we'll assume on the other hand that
$M$ admits binary products.

Let\pspage{99} $X,Y$ be objects of $M$, and
\[ f,g: X\rightrightarrows Y\]
two maps in $M$ from $X$ to $Y$. One key notion constantly used lately
(but so far only when $f$ is an identity map, and $g$ a ``constant''
one -- which is equally the case needed for defining contractibility)
is the notion of an \emph{\bI-homotopy from $f$ to $g$}, namely a map
$X\times I\xrightarrow h Y$ making commutative the diagram
\[\begin{tikzcd}[baseline=(O.base),sep=small]
  & X\times I \ar[dd,"h"] & \\
  X \ar[ur,"\id_X\times\delta_0"] \ar[dr,swap,"f"] & &
  X \ar[ul, swap, "\id_X\times\delta_1"] \ar[dl,"g"] \\
  & |[alias=O]| Y &
\end{tikzcd}.\]
Let's first restate the ``\hyperref[lem:homotopylemma]{homotopy lemma}'' of page \ref{p:60}
in a slightly more complete form:
\begin{homotopylemmareformulated}\label{lem:hlr}
  Assume $f$ and $g$ are \bI-homotopic. Then:
  \begin{enumerate}[label=\alph*),font=\normalfont]
  \item\label{it:48.hlr.a}
    $\gamma(f)=\gamma(g)$, where
    \[\gamma: M \to W^{-1}M=H_M\]
    is the canonical functor
  \item\label{it:48.hlr.b}
    If $f$ is a weak equivalence, so is $g$ \textup(and conversely of course,
    by symmetry of the roles of $\delta_0$ and $\delta_1$\textup).
  \item\label{it:48.hlr.c}
    Assume $f$ is an isomorphism, and $g$ constant -- then $X\to e$
    and $Y\to e$ are $W$-aspheric.
  \end{enumerate}
\end{homotopylemmareformulated}

It is important to notice that the relation ``$f$ is \bI-homotopic to
$g$'' in $\Hom(X,Y)$ is not necessarily symmetric nor transitive, and
that it depends on the choice of the homotopy interval
$\bI=(I,\delta_0,\delta_1)$. Thus the symmetric relation from
$\bI$-homotopy is $\check\bI$-homotopy, where $\check\bI$ is the
homotopy interval ``opposite'' to \bI{} (namely with
$\delta_0,\delta_1$ reversed). The example we are immediately
interested in is $M=\Cat$, with $W=W_\Cat$ the usual notion of weak
equivalence. A homotopy interval is just an aspheric small category
$I$, endowed with two \emph{distinct} objects $e_0,e_1$. (The
condition $e_0\ne e_1$ just expresses the condition $e_0\sand
e_1=\varnothing$ on homotopy intervals -- if it were not satisfied,
\bI-homotopy would just mean equality of $f$ and $g$\ldots) For the
usual choice $I=\Simplex_1=(e_0\to e_1)$, an \bI-homotopy from $f$ to
$g$ is just a morphism between functors $f\to g$ -- the $\bI$-homotopy
relation between $f$ and $g$ is the existence of such a morphism, it
is a transitive, and generally non-symmetric relation. If we take $I$
to be a category with just two objects $e_0$ and $e_1$, equivalent to
the final category, an \bI-homotopy between $f$ and $g$ is just an
isomorphism from $f$ to $g$ -- the \bI-homotopy relation now is both
transitive and symmetric, and it is a lot more restrictive than the
previous one. If we take $I$ to be the barycentric
subdivision\pspage{100} of $\Simplex_1$, which can also be interpreted
as an amalgamated sum of $\Simplex_1$ with itself, namely
\[I = \begin{tikzcd}[baseline=(A.base),row sep=-7pt,column sep=small]
  e_0 \ar[dr] & \\ & |[alias=A]| e_2 \\ e_1\ar[ur] &
\end{tikzcd},\]
an \bI-homotopy from $f$ to $g$ is essentially a triple $(k,u,v)$,
with $k:X\to Y$ and $u:f\to k$ and $v:g\to k$ maps in $\bHom(X,Y)$;
this time, the relation of \bI-homotopy is symmetric, but by no means
transitive. Returning to general $M$, it is customary to introduce the
equivalence relation in $\Hom(X,Y)$ generated by the relation of
\bI-homotopy -- we'll say that $f$ and $g$ are \emph{\bI-homotopic in
  the wide sense}, and we'll write
\[ f \Isim g,\]
if they are equivalent with respect to this relation. As seen from the
examples above (where $M=\Cat$), this relation still depends on the
choice of the homotopy interval \bI. Let's first look at what we can
do for fixed \bI, and then how what we do depends on \bI.

If $f$ and $g$ are \bI-homotopic, then so are their composition with
any $Y\to Z$ or $T\to X$. This implies that the relation $\Isim$ of
\bI-homotopy in the wide sense is compatible with compositions. If we
denote by
\[\Hom(X,Y)_\bI\]
the quotient set of $\Hom(X,Y)$ by the equivalence relation of
\bI-homotopy in the wide sense, we get composition between the sets
$\Hom(X,Y)_\bI$, and hence a structure of a category $M_\bI$ having
the same objects as $M$, and where maps from $X$ to $Y$ are
\bI-homotopy classes (in the wide sense -- this will be understood
henceforth when speaking of homotopy classes) of maps from $X$ to $Y$
in $M$. Two objects $X,Y$ of $M$, i.e., of $M_\bI$ which are
isomorphic as objects of $M_\bI$ will be called
\emph{\bI-homotopic}. This means also that we can find maps in $M$
(so-called \emph{\bI-homotopisms } -- namely $M_\bI$-isomorphisms)
\[ f:X\to Y, \quad g:Y\to X\]
such that we get \bI-homotopy relations in the wide sense
\[gf \Isim \id_X,\quad fg \Isim\id_Y.\]
We'll say $X$ is \emph{\bI-contractible} if $X$ is \bI-homotopic to
the final object $e_M=e$ of $M$ (which visibly is also a final object
of $M_\bI$), i.e., if $X$ is a final object of $M_\bI$. In terms of
$M$, this means that there exists a section $f$ of $X$ over $e$, such
that $fp_X$ is \bI-homotopic in the wide sense to $\id_X$ (where $p_X$
is the unique map $X\to e$). In fact, if there is such a section $f$,
any other section will do too.

From the homotopy lemma \ref{it:48.hlr.a} it follows that the canonical
functor\pspage{101} $M\to W^{-1}M=H_M$ factors into
\[ M \to M_\bI \to H_M = W^{-1}M,\]
and from \ref{it:48.hlr.b} it follows that if $f,g: X\rightrightarrows Y$
are in the same \bI-class, then $f$ is in $W$ if{f} $g$ is, hence by
passage to quotient a subset
\[ W_\bI \subset \Fl(M_\bI)\]
of the set of arrows in $M_\bI$, namely a notion of weak equivalence
in $M_\bI$. It is evident from the universal property of $H_M$ that
the canonical functor $M_\bI\to H_M$ induces an isomorphism of
categories
\[H_{M_\bI} = W_\bI^{-1}M_\bI \tosim H_M=W^{-1}M.\]

It's hard at this point not to expect that $W_\bI$ should satisfy the
same mild saturation conditions as $W$, so let's look into this in the
stride (even though I have not had any use of this so far). Condition
\ref{it:37.b} of saturation, namely that if $f,g$ are composable and
two among $f,g,gf$ are in $W$, so is the third, carries over
trivially. Condition \ref{it:37.a}, namely the tautological looking
condition that $W$ should contain all isomorphisms, makes already a
problem, however. It is OK though if $W$ satisfies the following
saturation condition, which is a strengthening of condition
\ref{it:37.c} of page \ref{p:59}:
\begin{enumerate}[label=\alph*'),start=3]
\item\label{it:48.cprimeone}
  Let $f:X\to Y$ and $g:Y\to X$ such that $gf\in W$ and $fg\in W$,
  then $f,g\in W$.
\end{enumerate}
This condition \ref{it:48.cprimeone} carries over to $M_\bI$
trivially. This suggests to introduce a strengthening of the ``mild
saturation conditions'', which I intend henceforth to call by the name
of ``saturation'', reserving the term of ``\emph{strong saturation}''
to what I have previously referred to occasionally as ``saturation''
-- namely the still more exacting condition that $W$ consists of
\emph{all} arrows made invertible by the localization functor $M\to
H_M=W^{-1}M$, or equivalently, by some functor $M\to H$. Thus, we'll
say $W$ is \emph{saturated} if{f} it satisfies the following:
\begin{enumerate}[label=\alph*')]
\item\label{it:48.aprime}
  For any $X$ in $M$, $\id_X\in W$.
\item\label{it:48.bprime}
  Same as \ref{it:37.b} before: if two among $f,g,gf$ are in $W$, so
  is the third.
\item\label{it:48.cprime}
  If $f:X\to Y$ and $g:Y\to X$ are such that $gf,fg\in W$, then
  $f,g\in W$.
\end{enumerate}
Each of these conditions carries over from $W$ to a $W_\bI$ trivially.

I can't help either having a look at the most evident exactness
properties of the canonical functor $M\to M_\bI$: Thus one immediately
sees that for two maps
\[f,g : X \rightrightarrows Y=Y_1\times Y_2,\]
with components $f_i,g_i$ ($i\in\{1,2\}$), $f$ and $g$ are
\bI-homotopic in the wide sense if{f} so are $f_i$ and $g_i$ (for
$i\in\{1,2\}$). The analogous statement is valid for maps into any
product object $Y=\prod Y_i$ on a finite set of\pspage{102} indices. The
dual statement so to say, when $X$ is decomposed as a sum $X=\coprod
X_i$, is valid too, provided taking products with $I$ is distributive
with respect to finite direct sums. Thus we get that $M\to M_\bI$
commutes with finite products, and with finite sums too provided they
are distributive with respect to multiplication with any object (or
with $I$ only, which would be enough).

The notion of \bI-homotopy in the wide sense between
$f,g:X\rightrightarrows Y$ can be interpreted in terms of strict
$\bI'$-homotopy with variable $\bI'$, as follows, provided we make
some mild extra assumptions on $(M,W)$, namely:
\begin{enumerate}[label=\alph*)]
\item\label{it:48.a}
  (Just for memory) $M$ is stable under finite products.
\item\label{it:48.b}
  $M$ is stable under amalgamated sums $I \amalg_e J$ under the final
  object $e$ ($I$ and $J$ endowed with sections over $e$).
\item\label{it:48.c}
  If moreover $I$ and $J$ are aspheric over $e$, then so is $I\amalg_e J$.
\end{enumerate}

Conditions \ref{it:48.b} and \ref{it:48.c} give a means of
constructing new homotopy intervals \bK, by amalgamating two homotopy
intervals \bI{} and \bJ, using as sections of \bI{} and \bJ{} for
making the amalgamation, either $\delta_0$ or $\delta_1$, which gives
four ways of amalgamating -- of course we take as endpoints of the
amalgamated interval, the sections over $e$ coming from the two
endpoints of \bI{} (giving rise to $\delta_0$ for \bK) and \bJ{}
(giving rise to $\delta_1$ for \bK) which have not been ``used up'' in
the amalgamation. Maybe the handiest convention is to define
\emph{the} amalgamated interval $\bI lor \bJ$, without any ambiguity
of choice, as being
\[ \bI \lor \bJ = (I,\delta_1^\bI) \amalg_e
(J,\delta_0^\bJ)\qquad
\parbox[t]{0.33\textwidth}{\raggedright endowed with the two
  sections coming from $\delta_0^\bI:e\to I$, $\delta_1^\bI:e\to
  J$,}\]
and defining the three other choices in terms of this operation, by
replacing one or two among the summands \bI, \bJ{} by the ``opposite''
interval $\check\bI$ or $\check\bJ$. The operation of amalgamation of
intervals, and likewise of homotopy intervals, just defined, is
clearly associative up to a canonical isomorphism, and we have a
canonical isomorphism of intervals
\[ (\bI\lor \bJ)\check{\phantom{\mathrm i}} \simeq \check\bI \lor \check\bJ.\]
I forgot to check, for amalgamation of homotopy intervals, the
condition \ref{cond:HIb} of page \ref{p:59}, namely $e_0 \sand e_1 =
\varnothing_M$ (which has not so far played any role, anyhow). To get
this condition, we'll have to be slightly more specific in condition
\ref{it:48.b} above on $M$ of existence of the relevant amalgamations,
by demanding (as suggested of course by the visual intuition of the
situation) that $I$ and $J$ should become subobjects of the
amalgamation $K$, and their intersection should be reduced to the
tautological part $e$ of it. More relevant still for the use we have
in mind is to demand that those amalgamations should commute to
products by an arbitrary element $X$ of $M$.\pspage{103} This I'll
assume in the interpretation of $\Isim$ in terms of strict
homotopies. Namely, let
\[\Comp(\bI)\]
by the set of all homotopy intervals deduced from \bI{} by taking
amalgamations of copies of \bI{} and $\check\bI$, with an arbitrary
number $n\ge1$ of summands. Thus we get just \bI{} and $\check\bI$ for
$n=1$, four intervals for $n=2$, \ldots, $2^n$ intervals for $n$
arbitrary. It is now immediately checked that for $f,g$ in
$\Hom(X,Y)$, the relation $f\Isim g$ is equivalent to the existence of
$\bK\in\Comp(\bI)$, such that $f$ and $g$ be \bK-homotopic (in the
strict sense).

\begin{remark}
  The saturation conditions
  \ref{it:48.aprime}\ref{it:48.bprime}\ref{it:48.cprime} on $W$ are
  easily checked for the usual notion of weak equivalence for
  morphisms of topoi, and hence also in the categories \Cat{} and in
  any topos, and therefore in any category \Ahat{} (where it boils
  down too to the corresponding properties on $W_\Cat$, as
  $W_A=i_A^{-1}(W_\Cat)$). Thus it seems definitely reasonable
  henceforth to take these as the standard notion of saturation
  (referring to its variants by the qualifications ``mild'' or
  ``strong''). On the other hand, the stability conditions
  \ref{it:48.a}\ref{it:48.b}\ref{it:48.c} on $(M,W)$ are satisfied
  whenever $M$ is a topos, with the usual notion of weak equivalence
  -- the condition \ref{it:48.c} being a consequence of the more
  general Mayer-Vietoris type statement about amalgamations of topoi
  under closed embeddings of such (cf.\ lemma on page). The same
  should hold in \Cat, with a similar Mayer-Vietoris argument -- there
  is a slight trouble here for applying the precedent result on
  amalgamation of topoi, because a section $e\to C$ of an object $C$
  of \Cat, namely an embedding of the one-point category $e=\Simplex_0$
  into $C$ by choice of an object of $C$, does not correspond in
  general to a closed embedding of topoi. (In geometrical language, we
  get a ``point'' of the topos \Chat{} defined by $C$, but a point
  need not correspond to a subtopos, let alone a closed one\ldots)
  This shows the asphericity criterion for amalgamation of topoi, and
  hence also for amalgamation of categories, has not been cut out with
  sufficient generality yet. As this whole $\Comp(\bI)$ story is just
  a digression for the time being, I'll leave it at that now.
\end{remark}

More important than amalgamation of intervals, is to compare the
notions of homotopy defined in terms of a homotopy interval \bI, to
the corresponding notions for another interval, \bJ. Here the natural
idea first is to see what happens if we got a morphism of intervals
(compatible with endpoints, by definition)
\[ \bJ\to\bI.\]
It\pspage{104} is clear then, for $f,g\in\Hom(X,Y)$, that any
\bI-homotopy from $f$ to $g$ gives rise to a \bJ-homotopy; hence if
$f$ and $g$ are \bI-homotopic, they are \bJ-homotopic, and hence the
same for homotopy in the wide sense. We get thus a canonical functor
\[ M_\bI \to M_\bJ\]
which is the identity on objects, entering into a cascade of canonical
functors
\[ M \to M_\bI \to M_\bJ \to H,\]
where $H=H_M=W^{-1}M$ can be viewed as the common localization of $M,
M_\bI, M_\bJ$ with respect to the notion of weak equivalences in these
categories. We may view $M_\bJ$ as a closer approximation to $H$ than
$M_\bI$. There is of course an evident transitivity relation for the
functors corresponding to two composable morphisms of intervals
\[ \bK\to\bJ\to\bI.\]
\begin{remark}
  In order to get that $\Isim$ implies $\Jsim$, it is sufficient to
  make a much weaker assumption than existence of a morphism of
  intervals $\bJ\to\bI$ -- namely it suffices to assume that the two
  sections $\delta_i^\bI: e\to I$ are \bJ-homotopic. More generally,
  let $\sim$ be an equivalence relation in $\Fl(M)$, compatible with
  compositions and with cartesian products (this is the case indeed
  for $\Jsim$), and let $\bI=(I,\delta_0,\delta_1)$ any object $I$ of
  $M$ (not necessarily aspheric over $e$) endowed with two sections
  over $e$, such that $\delta_0\sim\delta_1$. Then the interval \bI{}
  gives rise to an equivalence relation $\Isim$ in $\Fl(M)$, whose
  definition is quite independent of $W$ -- and a priori, if $f\Isim
  g$ and $f\in W$, this need not imply $g\in W$. However, the
  condition $\delta_0\sim\delta_1$ implies immediately that the
  relation $\Isim$ implies the relation $\sim$. When the latter is
  $\Jsim$, we get moreover that $W$ is the inverse image of a set of
  arrows in $M_\bI$, i.e., $f\Isim g$ and $f\in W$ implies $g\in W$.
\end{remark}

An interesting particular case is the one when we can find a homotopy
interval $\bI_0$ in $M$, which has the property that for any other
homotopy interval \bI{} in $M$, its structural sections satisfy
\[\delta_0^\bI \Izsim \delta_1^\bI.\]
This implies that the homotopy relation $\Izsim$ is implies by all
other similar relations $\Isim$, i.e., it is the coarsest among all
relations $\Isim$. We may then view $\bI_0$ as a
``\emph{fundamental}'' or ``\emph{characteristic}'' \emph{homotopy
  interval} in $M$, in the sense that the relation $f \Wsim g$ in the
sense below, namely existence of a homotopy interval \bI{} such that
$f \Isim g$, is equivalent to $f\Izsim g$, i.e., can be check using
the one and unique $\bI_0$. In the case of $M=\Cat$, we get readily
$\Simplex_1$ as a characteristic homotopy interval.\pspage{105} More
specifically, if $\delta_0,\delta_1:e \to I$ are two sections of an
object $I$ of \Cat, i.e., two objects $e_0,e_1$ of the small category
$I$, then these are $\Simplex_1$-homotopic if{f} $e_0,e_1$ belong to the
same connected component of $I$, which is automatic if $I$ is
$0$-connected, and a fortiori if $I$ is aspheric. This accounts to a
great extent, it seems, for the important role $\Simplex_1$ is playing
in the homotopy theory of \Cat, and consequently in the whole
foundational set-up I am developing here, using \Cat{} as the basic
``modelizer''.

It is not clear to me for the time being whether it is reasonable to
expect in more or less any modelizer $(M,W)$ the existence of a
characteristic homotopy interval (provided of course a homotopy
interval exists). This is certainly the case for the elementary
modelizers met so far. Presumably, I'll have to come back upon this
question sooner or later.

We'll now see that the set of equivalence relations $\Isim$ on
$\Fl(M)$, indexed by the set of homotopy intervals \bI, is ``filtrant
d\'ecroissant'', namely that for two such relations $\Isim$ and
$\Jsim$, there is a third ``wider'' one, $\Ksim$, implied by both. It
is enough to construct a \bK, endowed with morphisms
\[ \bK\to\bI,\quad \bK\to\bJ\]
of homotopy intervals. Indeed, there is a universal choice, namely the
category of homotopy intervals admits binary products -- we'll take
thus
\[ \bK = \bI \times \bJ,\]
where the underlying object of $\bK$ is just $I\times J$, endowed with
the two sections $\delta_i^\bI \times \delta_i^\bJ$ ($i\in\{0,1\}$).

As usual, we'll denote by $\Wsim$ or simply $\sim$ the equivalence
relation on $\Fl(M)$, which is the limit or union of the the
equivalence relations $\Isim$ -- in other words
\[f\Wsim g\quad\text{if{f} exists \bI, a homotopy interval in $M$,
  with $f\Isim g$.}\]
For $X,Y$ in $M$, we'll denote by
\[ \Hom(X,Y)_W\]
the quotient of $\Hom(X,Y)$ by the previous equivalence relation. This
relation is clearly compatible with compositions, and hence we get a
category $M_W$, having the same objects as $M$, which can be equally
viewed as the filtering limit of the categories $M_\bI$,
\[ M_W = \varinjlim_\bI M_\bI,\]
where we may take as indexing set for the limit the set of all \bI's,
preordered by $\bI\le\bJ$ if{f} there exists a morphism of intervals
$\bJ\to\bI$ (it's the preorder relation opposite to the usual one on
the\pspage{106} set of objects of a category\ldots).

We now get canonical functors
\[ M \to M_\bI \to M_W \to H,\]
where $H$ can again be considered as the common localization of the
three categories $M$, $M_\bI$ (any \bI), $M_W$ with respect to the
notion of weak equivalence in each. It is clear that the set of weak
equivalences in $M_W$, say $\overline W$, is saturated provided $W$
is. Also, the canonical functor $M\to M_W$ commutes with finite
products, and also with finite sums provided formation of such sums in
$M$ commutes with taking products with any fixed element of $M$.

A map $f:X\to Y$ is called a \emph{homotopy equivalence} or a
\emph{homotopism} (with respect to $W$) if it is an isomorphism in
$M_W$, namely if there exists $g:Y\to X$ such that
\[gf\Wsim\id_X, \quad fg\Wsim\id_Y.\]
This implies that $f$ and $g$ are weak equivalences. If such an $f$
exists, namely if $X$ and $Y$ are isomorphic objects of $M_W$, we'll
say they are \emph{$W$-homotopic}, or simply \emph{homotopic}. This
implies that there exist weak equivalences $X\to Y$ and $Y\to X$, but
the converse of course isn't always true. If $Y$ is the final object
$e_M$, we'll say $X$ is \emph{$W$-contractible} (or simply
\emph{contractible}) instead of $W$-homotopic to $Y=e$. In case there
exists a characteristic homotopy interval $\bI_0$, all these
$W$-notions boil down to the corresponding $\bI_0$-notions considered
before.
\begin{remark}
  If occurs to me that in all what precedes, we never made any use
  really of the existence of a homotopy interval! The only notion we
  have effectively been working with, it seems to me, is the notion of
  ``weak homotopy interval'', by which I mean a triple
  $\bI=(I,\delta_0,\delta_1)$ satisfying merely the condition that $I$
  be aspheric over $e$ (which is the condition \ref{cond:HIc} on page
  \ref{p:59}). Such an \bI{} always exists of course, we need only
  take $I=e$ itself! In case however an actual homotopy interval (with
  $\Ker(\delta_0,\delta_1)=\varnothing_M$) does exist, to make sure that
  the notion of $W$-homotopy obtained by using all weak homotopy
  intervals is the same when using only actual homotopy intervals, we
  should be sure that for any $I$ aspheric over $e$, any two sections
  $\delta_0,\delta_1$ of $I$ over $e$ are $W$-homotopic in the initial
  meaning. This is evidently so when $M=\Cat$ and $W=W_\Cat$ (indeed,
  it is enough that $I$ be $0$-connected, instead of\pspage{107} being
  aspheric), and in the case when $M$ is one among the usual
  elementary modelizers. The question of devising notions which still
  make sense when there is no homotopy interval in $M$ isn't perhaps
  so silly after all, if we remember that the category of
  semi-simplicial face complexes (namely without degeneracies) is
  indeed a modelizer, but it hasn't hot a homotopy interval. However,
  for the time being I feel it isn't too urgent to get any more into this.
\end{remark}

\hangsection{Contractible objects. Multiplicative intervals.}\label{sec:49}%
I would like now to elaborate a little on the notion of a contractible
element $X$ in $M$, which (I recall) means an object, admitting a
section $e_0:e\to X$, such that the constant map
\[c = p_Xe_0 : X\to e\to X\]
is homotopic to $\id_X$, i.e., there exists a homotopy interval $I$
such that $\id_X\Isim c$ (\bI-homotopy in the wide sense).

If $X$ is contractible, then the constant map $c:X\to X$ is a weak
equivalence (as it is homotopic to $\id_X$ which is) and hence by the
saturation condition \ref{it:48.cprime} (in fact the mild saturation
condition \ref{it:37.c} suffices) it follows that $p_X:X\to e$ is a
weak equivalence. In fact, one would expect it is even universally so,
i.e., that $p_X$ is an \emph{aspheric} map, as a plausible
generalization of the homotopy lemma \ref{it:48.hlr.c} above (which we
have used already a number of times as our main asphericity criterion
in elementary \Ahat\ldots). The natural idea to prove asphericity of
$p_X$, namely that for any base change $S\to e$, the projection
$X_S=X\times S\to S$ is a weak equivalence, is to apply the precedent
criterion to $X_S$, viewed as an element of $M_{/S}$. As the base
change functor commutes with products, it is clear indeed that for two
maps $f,g:X\rightrightarrows Y$ in $M$, \bI-homotopy of $f$ and $g$
will imply $\bI_S$-homotopy of $f_S$ and $g_S$, with evident
definition of base change for an ``interval''. On the other hand,
asphericity of \bI{} over $e$ implies tautologically asphericity of
$\bI_S$ over $S$, the final object of $M_{/S}$, which is all we need
to care about to get the result that $X_S\to S$ is a weak
equivalence. (The condition \ref{cond:HIb} on homotopy intervals is
definitely misleading here, as it would induce us to put the extra
condition that the base change functor is compatible with initial
elements, which is true indeed if $\varnothing_M$ is strict, but here a
wholly extraneous condition\ldots)

Thus \emph{$X$ contractible implies $X$ aspheric over $e$}. It is
well-known that the converse isn't true, already for $M=\Sssets$, or
$M=\Spaces$ with the usual notions of weak equivalences, taking
aspheric complexes which are not Kan complexes, or ``aspheric'' spaces
(in the sense of singular cohomology) which are not CW-spaces. The
same examples show, too, that even a good honest homotopy interval
need not be contractible,\pspage{108} contrarily to what one would
expect from the intuitive meaning of an ``interval''. In those two
examples, the condition of getting two ``disjoint'' sections of the
aspheric $I$ over $e$ looks kind of trivial -- kind of unrelated to
the question of whether or not $I$ is actually contractible, and not
only aspheric. What comes to mind here is to look for \emph{contractible
  homotopy intervals}, namely homotopy intervals
$(I,\delta_0,\delta_1)$ (including condition \ref{it:37.b} that
$\Ker(\delta_0,\delta_1)=\varnothing_M$) such that $I$ is not only
aspheric over $e$, but even contractible. Existence of contractible
homotopy intervals seems a priori a lot stronger than just existence
of homotopy intervals, but in the cases I've met so far, the two
apparently coincide. As a matter of fact, the first choice of homotopy
interval that comes to mind, in all those examples, \emph{is} indeed a
contractible one -- so much so that till the notes today I was
somewhat confused on this matter, and was under the tacit impression
that homotopy intervals are all contractible, and trivially so! But
probably, as usual with most evidently false impressions, underneath,
something correct should exist, worth being explicitly stated.

First thing that comes to mind is here to look at the simplest case of
an aspheric interval \bI{} which is contractible -- namely when \bI{}
is even \bI-contractible, and more specifically still, when there
exists an \bI-homotopy (an ``elementary'', not a ``composed'' one!)
from $\id_I$ to one of the two constant maps of $I$ into itself,
defined by the two sections $\delta_0,\delta_1$ -- say by $\delta_1$,
to be specific. Such a homotopy is a map
\[ h : I\times I\to I\]
in $M$, having the two properties expressed symbolically by
\[ h(e_0,x)=x, \quad h(e_1,x)=e_1,\]
where $x$ may be viewed as any ``point'' of $I$ with ``values'' in an
arbitrary parameter object $T$ of $M$, i.e., $x:T\to I$, and $e_0,e_1$
are the constant maps $T\to I$ defined by $\delta_0,\delta_1$. If we
view $h$ as a composition law on $I$, these relations just mean that
$e_0$ acts as a left unit, and $e_1$ as a left zero element. Such
composition law in an interval $I$, a would-be homotopy interval as a
matter of fact, has been repetitively used before, and this use
systematized in the ``\hyperref[lem:comparisonlemmaforHI]{comparison lemma
  for homotopy intervals}'' (page \ref{p:60}); there we found that if
$I$ is an object of $M$ with such a composition law, and if there
exists an \emph{aspheric} interval \bJ{} and a morphism of intervals
from \bJ{} into \bI, then $I$ is equally aspheric. This we would now
see as coming from the fact that $I$ being \bI-contractible is
\bJ-contractible, and a fortiori (as \bJ{} is aspheric over $e$)
aspheric over $e$. In any case, we see that when $I$ is endowed with a
composition law as above, then $I$ is aspheric over $e$ if{f} it is
contractible, and equivalently if{f} $\delta_0$ and $\delta_1$ are
homotopic.

In\pspage{109} case when there exists a Lawvere element $L_M$ in $M$,
for instance if $M$ is a topos, this element is automatically endowed
with an idempotent composition law, coming from intersection of
subobjects, and in case $M$ admits a strict initial object, $L_M$ is
endowed with two canonical sections which are indeed ``disjoint'',
corresponding respectively to the ``full'' and the ``empty''
subobjects. Then (as already noticed) previously in the
\hyperref[cor:ofcomparisonlemmaforHI]{corollary} of the comparison
lemma) $M$ admits a homotopy interval if{f} $L_M$ is such an interval,
namely is aspheric over $e_M$. But we can now add that in this case,
there exists even a \emph{contractible} homotopy interval, namely
$L_M$ itself. It is even a ``\emph{strict homotopy interval}'', namely
one admitting a composition $L\times L\to L$ having the properties
above (where $\delta_0$ and $\delta_1$, for the first time, play
asymmetric roles!). Thus when a Lawvere element exists in $M$, there
is an equivalence between the three properties we may expect from a
pair $(M,W)$, with respect to homotopy intervals, namely:
\[ (\exists\text{homotopy int.}) \Leftarrow
(\exists\text{contractible hom.\ int.}) \Leftarrow
(\exists\text{strict hom.\ int.}).\]
It seems that in all cases I've in mind at present, when there exists
a homotopy interval, there exists even a strict one. The only case
besides topoi (where it is so, because of the existence of a Lawvere
element) which I have looked up so far, is the case of \Cat{} and
various full subcategories, all containing $\Simplex_1$ which is indeed
a strict homotopy interval, as it represents the presheaf on \Cat{}
\[ C\mapsto\text{set of all cribles in $C$.}\]
If we take the choice $I=$ two-point category equivalent to final one,
this also is a strict homotopy interval, as it represents the functor
\[ C\mapsto\text{set of all full subcategories of $C$,}\]
hence again an intersection law. The first homotopy interval though is
a lot more important than the second, because the first one is
``characteristic'', namely sufficient for checking the homotopy
relation between any two maps in \Cat, whereas the second isn't. A
``\emph{perfect}'' homotopy interval would be one which is both
\emph{strict} (hence contractible) and \emph{characteristic}.

\bigbreak
\presectionfill\ondate{2.4.}\par

\hangsection[Reflection on some main impressions. The foresight of
\dots]{Reflection on some main impressions. The foresight of an
  ``idyllic picture'' \texorpdfstring{\textup(}{(}of would-be
  ``canonical
  modelizers''\texorpdfstring{\textup)}{)}.}\label{sec:50}%
While writing the notes last time, and afterwards while pondering a
little more about the matter, a few impressions came gradually into
the fore. One was about the interplay of four basic ``homotopy
notions'' which more or less mutually determine each other, namely the
\emph{homotopy relation} between maps, the notion of \emph{homotopy
  interval}, the notion of homotopy equivalences or \emph{homotopisms}
(which has formal analogy to weak equivalences\pspage{110} the was it
is handled), and the notion of \emph{contractible objects}. Another
impression was about the dependence of these notions upon a
preliminary notion of ``weak equivalence'', namely upon
$W \subset \Fl(M)$, being a rather loose one. Thus the construction of
homotopy notions in terms of a given interval \bI{} (including the
category $M_\bI$ with the canonical functor $M\to M_\bI$) is valid for
any interval in any category $M$ with final object and binary products
(instead of binary products, it is even enough that $I$ be
``squarable'', namely all products $X\times I$ exist in $M$). As for
the $W$-homotopy notions, they depend on $W$ via the corresponding
notion of $W$-asphericity over $e$, which is at first sight \emph{the}
natural condition to impose upon an interval \bI, in order for the
corresponding \bI-homotopy notions to fit nicely with $W$ (as
expressed in the homotopy lemma). But then we noticed that a much
weaker condition than asphericity on \bI{} suffices -- namely that the
two sections $\delta_0,\delta_1$ of $I$ over $e$ be $W$-homotopic,
which means essentially that the ``points'' of $I$ they define can be
``joined'' by a finite chain of $W$-aspheric intervals mapping into
$I$. This strongly suggests (in view of the main application we have
in mind, namely to the study of modelizers) that the natural condition
to impose upon intervals, in most contexts of interest to us, will be
merely $0$-connectedness. But this notion is \emph{intrinsic to the
  category} $M$, irrespective again of the choice of any $W$; and
therefore the corresponding homotopy notions in $M$ will turn out (in
the cases at least of greatest interest to us) to be equally intrinsic
to the category $M$. On the other hand (and here fits in the third
main impression that peeled out two day ago), the work carried through
so far in view of the ``observation'' and the (naive) ``question'' of
last week (pages \ref{p:94} and \ref{p:95}) strongly suggests that in
the nicest modelizers (including \Cat{} and the elementary modelizers,
presumably), the notion of weak equivalence $W$ can be described in
terms of the homotopy notions, more specifically in terms of the
notion of contractible objects (when exactly and how should appear in
due course). Thus it will follow that for those modelizers, the
modelizing structure $W$ itself is uniquely determined in terms of the
intrinsic category structure -- thus any equivalence between the
underlying categories of any two such modelizers should automatically
be model-preserving! It will be rather natural to call the modelizers
which fit into this idyllic picture \emph{canonical modelizers}, as
their modelizing structure $W$ is indeed canonically determined by the
category structure. Next thing then would be to try to gain an overall
view of how to get ``all'' canonical modelizers, if possible in as
concrete terms as the overall view we got upon elementary modelizers
\Ahat{} in terms of the corresponding test categories $A$.\pspage{111}

\hangsection{The four basic ``pure'' homotopy notions with
  variations.}\label{sec:51}%
\renewcommand*{\thesubsection}{\Alph{subsection})}%
First thing though I would like to do now, is to elaborate ``from
scratch'' on the four basic homotopy notions and their interplay, much
in the style of a ``fugue with variations'' I guess, and without
interference of a pregiven notion $W$ of weak equivalence --
relationship with a $W$ will be examined only after the intrinsic
homotopy notions and their interrelations are well understood.

We start with a category $M$, without for the time being any specific
assumptions on $M$. The strongest we're going to introduce, I guess,
is existence of finite products, and incidently maybe finite sums and
fiber products. In the cases we have in mind, $M$ is a ``large''
category, therefore it doesn't seem timely here introduce $M\uphat$
and the embedding of $M$ into $M\uphat$.

The most trivial implications between the four basic homotopy notions
are symbolized by the plain arrows in the diagram below, the somewhat
more technical ones by dotted arrows. It is understood these notions
correspond to a given ``homotopy structure'' on $M$, symbolized by the
letter $h$, and which (in the most favorable cases) may be described
at will in terms of any one of the four notions. I'll first describe
separately each of these basic notions, and afterwards the
relationships symbolized by the arrows in the diagram. I recall that
in \emph{interval} in $M$ is just an object $I$, endowed with two
subobjects $e_0,e_1$ which are final objects of $M$, or equivalently,
with two sections $\delta_0,\delta_1$ of $I$ over a fixed final object
$e_M=e$ of $M$. I definitely want to forget entirely for the time
being about any condition of the type $e_0\sand e_1=\varnothing_M$
(initial object of $M$) -- we may later refer to these as
``\emph{separated}'' intervals (namely the endpoints $e_0,e_1$ are
``separated''). We'll denote by $\Int(M)$ the set of all intervals in
$M$, by $\bInt(M)$ the corresponding category (the notion of a
morphism of intervals being the obvious one). Now here's the
\namedlabel{fig:organigram}{organigram}:
\begin{equation}
  \label{eq:51.D}
  \begin{tikzcd}[arrows=Rightarrow]
    \begin{tabular}{@{}c@{}}1) homotopy relation\\
      $R_h\subset\Fl(M)\times\Fl(M)$
    \end{tabular} \ar[d]\ar[r] &
    \begin{tabular}{@{}c@{}}2) homotopism\\
      $W_h\subset\Fl(M)$
    \end{tabular} \ar[l, bend right=10, dashed, swap, "1"]\ar[d] \\
    \begin{tabular}{@{}c@{}}3) homotopy intervals\\
      $\Sigma_h\subset\Int(M)$
    \end{tabular} \ar[u, bend left=25, dashed, "2"] &
    \begin{tabular}{@{}c@{}}4) contractible objects\\
      $C_h\subset\Ob(M)$
    \end{tabular} \ar[l, dashed, swap, "3"]
  \end{tikzcd}
  \tag{D}
\end{equation}

\subsection[Homotopy relation between maps]{Homotopy relation between maps.}
\label{subsec:51.A}
As\pspage{112} a \emph{type of structure}, a homotopy relation between
maps in $M$ is a subset
\[R_h\subset\Fl(M)\times\Fl(M),\]
namely a relation in the set $\Fl(M)$ or arrows of $M$, the
\emph{basic assumption} being that whenever $f$ and $g$ are
``homotopic'' arrows, then they have the same source, and the same
target. Thus, the data $R_h$ is equivalent to giving a ``homotopy
relation'' in any one of the sets $\Hom(X,Y)$, with $X$ and $Y$
objects in $M$. The relevant \emph{saturation condition} is twofold:
\begin{enumerate}[label=\alph*)]
\item\label{it:51.A.a}
  the relation $R_h$ is an equivalence relation, or equivalently, the
  corresponding relations in the sets $\Hom(X,Y)$ are equivalence
  relations;
\item\label{it:51.A.b}
  stability under composition: if $f$ and $g$ are homotopic, then so
  are $vf$ and $vg$, and so are $fu$ and $gu$, for any arrow $v$ or
  $u$ such that the relation makes sense.
\end{enumerate}

When these conditions are satisfied, we'll say we got a \emph{homotopy
  relation} between maps of $M$. This relation between $f$ and $g$
will be denoted by a symbol like
\[ f\hsim g.\]
If the basic assumption is satisfied, but not the saturation
condition, there is an evident way of ``saturating'' the given
relation, getting one $\overline R_h$ which is saturated, i.e., a homotopy
relation in $M$ (in fact, the smallest one containing $R_h$).

Given a homotopy relation $R_h$, we denote by
\[ \Hom(X,Y)_h\]
the corresponding quotient sets of the set $\Hom(X,Y)$, they compose
in an evident way, so as to give rise to a category $M_h$ having the
same objects as $M$, and to a canonical functor
\[ M\to M_h\]
which is the identity on objects, and surjective on arrows. We may
view thus $M_h$ as a \emph{quotient category} of $M$, having the same
objects as $M$. Clearly, $R_h\mapsto M_h$ is a bijective
correspondence between the set of homotopy relations in $M$, and the
set of quotient categories of $M$ satisfying the aforesaid
property. By abuse of language, we may even consider that considering
a homotopy relation in $M$, amounts to the same as giving a functor
$M\to M_h$ from $M$ which is bijective on objects and surjective on
arrows.

When we got a homotopy relation $R_h$ in $M$, we deduce a notion of
\emph{homotopisms}
\[ W_h\subset\Fl(M),\]
namely those arrows in $M$ which become isomorphisms in $M_h$. Also
we\pspage{113} deduce a notion of \emph{homotopy interval}, i.e.,
\[\Sigma_h\subset\Int(M),\]
namely those intervals \bI{} in $M$ such that the two marked sections
be homotopic. (NB\enspace $\Sigma_h$ is non-empty if and only if $M$ has a
final element, in this case $\Sigma_h$ contains all intervals such
that $\delta_0=\delta_1$ -- which we may call \emph{trivial}
intervals, for instance the \emph{final interval} with $I=e$\ldots)
This notion of a homotopy interval is considerably wider than the one
we have worked with so far, however it is clearly the right one in the
context of pure homotopy notions. To avoid any confusion, we better
call this notion by the name of \emph{weak homotopy intervals} --
funnily there won't be any unqualified ``homotopy intervals'' in our
present set-up of ``pure'' homotopy notions!

The two prescriptions above account for two among the plain arrows in
our \ref{fig:organigram}.

We'll often make use of an \emph{accessory assumption} on $R_h$, which
can be expressed by demanding that \emph{the canonical functor $M\to
  M_h$ commute to binary products}, in case we assume already such
products exist in $M$. This can be expressed also by the property that
for two maps
\[f,g: X \rightrightarrows Y_1\times Y_2,\]
$f$ and $g$ are homotopic if{f} so are $f_i$ and $g_i$ ($i\in\{1,2\}$)
(where the ``only if'' part is satisfied beforehand anyhow). This
implies too that in $M_h$ binary products exist, and that $f\sim f'$,
$g\sim g'$ implies $f\times g\sim f'\times g'$. On the other hand, it
is trivial that if $M$ admits a final object, this is equally a final
object of $M_h$ (and hence, under the accessory assumption on $R_h$,
the functor $M\to M_h$ commute to finite products).

\subsection[Homotopisms, and homotopism structures]{Homotopisms.}
\label{subsec:51.B}
As a type of structure on $M$, we got just a subset
\[W_h\subset\Fl(M),\]
without any basic assumption to make. The natural \emph{saturation
  condition} is just the strong saturation for a subset of $\Fl(M)$,
which can be expressed by stating that $W_h$ can be obtained as the
set of arrows made invertible by some functor from $M$ into a category
$M'$, or equivalently, by the localization functor
\[M\to W_h^{-1}M.\]
We may refer to a strongly saturated $W_h$ as a ``\emph{homotopism
  structure}'' (or ``homotopy equivalence structure'') in $M$ -- but
as in the case \ref{subsec:51.A}, we'll have soon enough to make
pretty strong extra assumptions. Maybe we should, at the very least,
demand for the notion of homotopy structure that the canonical functor
above, which is bijective on objects in any case, should be moreover
\emph{surjective on arrows} -- thus I'll take this\pspage{114} as a
\emph{basic assumption} after all. This assumption makes sense of
course independently of any saturation condition. If $W_h$ is not
strongly saturated, then denoting by $\overline W_h$ the subset of $\Fl(M)$
of all arrows made invertible by the canonical functor, this will now
be a strongly saturated set of arrows (in fact the smallest one
containing $W_h$, and giving rise to the same localized category, and
hence satisfying the basic assumption too) -- thus $\overline W_h$ is
indeed a homotopism structure on $M$. When $M$ admits a final object,
this will equally be a final object of $W_h^{-1}M$. We may now define
in terms of $W_h$ the notion of contractible objects in $M$, forming a
subset
\[C_h \subset \Ob(M),\]
as those objects $X$ in $M$ such that the projection $p_X:X\to e$ is
in $W_h$, or equivalently, such that $X$ is a final object in the
localized category $W_h^{-1}M$. This accounts for the third plain
arrow of the \ref{fig:organigram}.

We'll now dwell a little more on the first dotted arrow, namely the
description of a homotopy relation
\[R_h \subset \Fl(M)\times\Fl(M)\]
in terms of $W_h$: the natural choice here is to define $f,g\in\Fl(M)$
to be \emph{homotopic} (or $W_h$-homotopic, if ambiguity may arise)
if{f} their images in the category $W_h^{-1}M$ are equal. This
relation between maps in $M$ clearly satisfies the basic assumption on
source and target, as well as the saturation condition -- it is
therefore a ``homotopy relation'' in $M$, namely the one associated to
$W_h^{-1}$, viewed as a quotient category of $M$. It is clear that we
recover $W_h$ from $R_h$, consequently, by the process described in
\ref{subsec:51.A}.

To make the relationship between the notions 1) and 2) still clearer,
let's denote respectively by
\[ \Hom_1(M),\Hom_2(M)\]
the set of all homotopy relations, resp.\ of homotopism notions, in
$M$. We got maps\footnote{\alsoondate{3.4.} There \emph{is no map} $r_{21}$, only
  $r_{12}$, see correction in \S\ref{sec:52}.}
\[\begin{tikzcd}[cramped]
  \Hom_1(M) \ar[r, shift left, "r_{21}"] &
  \Hom_2{M} \ar[l, shift left, "r_{12}"]
\end{tikzcd},\]
and the relevant fact here is that
\[r_{12}:\Hom_2\to\Hom_1\]
\emph{is injective, and admits $r_{21}$ as a left inverse.} Thus, we
may view $\Hom_2$ as a subset of $\Hom_1$, i.e., the structure of a
``homotopism notion'' on $M$ as a particular case of the structure of
a ``homotopy relation'' on $M$. Namely, a structure of the latter type
can be described in terms of a notion of homotopism in $M$, if{f} the
canonical functor $M\to M_h$ it gives to\pspage{115} is a localization
functor.

For a general $R_h\in\Hom_1(M)$, if we consider the corresponding
$W_h$ ($=r_{21}(R_h)$) in $\Hom_2(M)$, it is clear that the canonical
functor $M\to M_h$ of \ref{subsec:51.A} factors into
\[ M \to W_h^{-1} \to M_h,\]
and $R_h$ ``is in $\Hom_2(M)$'', i.e., $R_h=r_{12}(W_h)$, if{f} the
second functor
\[ W_h^{-1}M \to M_h\]
(which is anyhow bijective on objects and surjective on arrows) is an
isomorphism, or equivalently, faithful. Here is a rather direct
\emph{sufficient} condition on $R_h$ for this to be
so,\footnote{\alsoondate{3.4.} See \S\ref{sec:52} for a
  \scrcommentinline{?} of this rash statement!} namely:
\begin{description}
\item[\namedlabel{cond:51.C.12}{C$_{12}$})] If $f,g:X\to Y$ are homotopic,
there exists a homotopism $X'\to X$, two sections $s_0,s_1$ of $X'$
over $X$, and a map $h:X'\to Y$, such that
\[f_0=hs_0, \quad f_1=hs_1.\]
\end{description}
\begin{remark}
  Intuitively, we are thinking of course of $X'$ as a product $X\times
  I$, where $\bI=(I,e_0,e_1)$ is a weak homotopy interval, and
  $s_0,s_1$ are defined in terms of $e_0,e_1$. In Quillen's somewhat
  different set-up, $X'$ is referred to as a ``cylinder object for
  $X$'', suitable for defining the ``left homotopy relation''
  associated to a given $W_h$. The condition \ref{cond:51.C.12} is not
  autodual, we could state a dual sufficient condition in terms of a
  ``path object for $Y$'', namely a homotopism $Y\to Y'$ endowed with
  two retractions $t_0,t_1$ upon $Y$ -- but we don't have any use for
  this in the present set-up, which (as for as the main emphasis is
  concerned) is by no means autodual, as is Quillen's.
\end{remark}
The condition \ref{cond:51.C.12} above can be viewed equally as a
condition on a $W_h\in\Hom_2(M)$.

We may interpret the set $\Hom_2(M)$ of homotopism notions in $M$ as
the set of all quotient categories $M_h$ of $M$, having the same
objects as $M$, and such that moreover the canonical functor $M\to
M_h$ be a localizing functor. As in \ref{subsec:51.A}, the relevant
``\emph{accessory assumption}'' on $W_h$ (a particular case indeed of
the corresponding one for $R_h$) is that this functor commute to
products. I don't see any simple computational way though to express
this condition directly in terms of $W_h$, as previously in terms of
$R_h$. I would only like to notice here a consequence of this
assumption (I doubt it is equivalent to it), namely that \emph{the
  cartesian product of two homotopisms is again a homotopism} -- which
implies, for instance,\pspage{116} that the product of a finite family
of contractible objects of $M$ is again contractible.

\subsection[Homotopy interval structures]{Weak homotopy intervals.}
\label{subsec:51.C}
We assume $M$ stable under finite products. The type of structure
we've in view is a set of intervals in $M$,
\[\Sigma_h\subset\Int(M),\]
called the ``weak homotopy intervals''. No basic assumption on this
set, it seems; the natural ``\emph{saturation condition}'' is the
following:
\begin{description}
\item[\namedlabel{cond:51.Sat.3}{(Sat 3)}] Any interval
  $\bI=(I,\delta_0,\delta_1)$ in $M$, such that the sections
  $\delta_0,\delta_1$ of $I$ be $\Sigma_h$-homotopic (see below), is
  in $\Sigma_h$.
\end{description}

The assumption on \bI{} means, explicitly, that there exists a finite
chain of sections of $I$
\[s_0=\delta_0, s_1, \dots, s_N=\delta_1,\]
joining $\delta_0$ to $\delta_1$, and for two consecutive
$s_i,s_{i+1}$ an interval \bJ{} in $\Sigma_h$, and a map of intervals
from \bJ{} or $\check\bJ$ to $(I,s_i,s_{i+1})$, i.e., a map $J\to I$,
mapping the two given sections of $J$, one into $s_i$, the other into
$s_{i+1}$ (without specification which is mapped into which).

The significance of this saturation condition becomes clear in terms
of the second dotted arrow of the \ref{fig:organigram}. Namely, in terms of any
subset $\Sigma_h$ of $\Int(M)$, we get a corresponding homotopy
relation between maps, say $R_h$, which is the equivalence relation in
$\Fl(M)$ generated by the ``elementary'' homotopy relation (with
respect to $\Sigma_h$) between maps $f,g$ in $M$, namely the relation $R_0$
\[ f \mathrel{\underset{R_0}{\sim}} g
\;\xLeftrightarrow{\textup{def}}\;
\parbox[t]{0.6\textwidth}{there exists
  \bI{} in $\Sigma_h$, and an \bI-homotopy from $f$ to $g$.}\]
The corresponding equivalence relation $R_h$ in $\Fl(M)$ is already
saturated, namely stable under compositions, moreover it satisfies
condition \ref{cond:51.C.12} above -- thus we may view this homotopy
relation as defined in terms of a homotopisms notion -- thus in fact
the second dotted arrow should go from 3) to 2) rather than from 3) to
1)! Now, if we look at the subset $\Sigma_h$ of $\Int(M)$
defined in terms of $R_h$ as in \ref{subsec:51.A} (namely the
set of ``homotopy intervals with respect to $R_H$''), we get
\[ R_h \subset \overline R_h,\]
and the equality holds if{f} $R_h$ satisfies \ref{cond:51.Sat.3}? At
the same time, in case of arbitrary $R_h$, we get the construction of
its saturation, $\overline R_h$, which may of course be described
alternatively as the smallest saturated subset of $\Int(M)$ containing
$R_h$.

We'll call \emph{weak homotopy interval structures} on $M$, any
set\pspage{117} $\Sigma_h$ of intervals in $M$, satisfying the
saturation condition above. The set of all such structures on $M$ will
be denoted by $\Hom_3(M)$, thus we get two embeddings
\[ \Hom_3(M) \hookrightarrow \Hom_2(M) \hookrightarrow \Hom_1(M),\]
in such a way that a weak homotopy interval structure on $M$ may be
viewed also as a particular case of a homotopism structure on $M$, and
a fortiori as a particular case of a homotopy relation on $M$. Of
course, the homotopy relations or homotopism structures on $M$ we'll
ultimately be interested in, are those stemming from weak homotopy
interval structures on $M$. Recall that $M$ admits finite products,
and these structures satisfy automatically the accessory assumption,
namely commutation of the canonical functor $M\to M_h$ to finite
products.

It is immediate that if we start with a homotopy relation $R_h$, the
corresponding $\Sigma_h$ as defined in \ref{subsec:51.A} is
saturated. Thus, the canonical embedding $r_{13}$ of $\Hom_3$ into
$\Hom_1$ admits a canonical left inverse $r_{31}$, the restriction of
which to $\Hom_2$ is a canonical left inverse $r_{32}$ of the natural
embedding $r_{23}$ of $\Hom_3$ into $\Hom_2$.

\subsection[Contractibility structures]{Contractibility structures.}
\label{subsec:51.D}
(We still assume $M$ admits finite products.) As a type of structure,
it is a set of objects of $M$
\[ C_h\subset \Ob(M),\]
without any ``basic assumption'' on $C_h$ it seems. These objects will
be called the \emph{contractible} objects. Sorry, there \emph{is} a
basic assumption here I just overlooked, namely every $X$ in $C_h$
should have at least one section (thus I better assume beforehand, as
in \ref{subsec:51.C}, that $M$ has a final object $e$). To get the
natural saturation condition on $C_h$, we'll make use of the third
dotted arrow in the \ref{fig:organigram}, by associating to $C_h$ the set
$\Sigma_h$ of ``contractible intervals'', namely intervals
$\bI=(I,\delta_0,\delta_1)$ such that $I$ is in $C_h$. Of course in
general there is no reason that $\Sigma_h$ should be saturated, never
mind -- it defines anyhow (as seen in \ref{subsec:51.C}) a homotopism
notion in $M$, and hence (as seen in \ref{subsec:51.B}) a notion of
contractible objects, i.e., another subset $\overline C_h$ of
$\Ob(M)$. Now it occurs to me that it is by no means clear that the
latter contains $C_h$, which brings near the necessity of a more
stringent \emph{basic assumption} on $C_h$, namely for the very least
\[ C_h \subset \overline C_h\]
(this will imply that any $X$ in $C_h$ has indeed a section over $e$,
as this is automatically the case for $W_h$-contractible objects). The
saturation condition (Sat~4) will of course be equality
\[ C_h = \overline C_h,\]
and\pspage{118} for general $C_h$ (satisfying the basic assumption
$C_h\subset\overline C_h$), $\overline C_h$ can be viewed as the
``saturation'' of $C_h$, namely the smallest saturated subset of
$\Ob(M)$ satisfying the basic assumption, or in other words, the
smallest \emph{contractibility structure} on $M$ such that the objects
in $C_h$ are contractible.

It may be worth while to state more explicitly the basic assumption
here, and the saturation condition on $C_h$.
\begin{description}
\item[\namedlabel{cond:51.Bas.4}{(Bas~4)}]
  For any $X$ in $C_h$, we can find a finite sequence of maps from $X$
  to $X$,
  \[ f_0=\id_X,f_1,\dots,f_N=c_s,\]
  joining the identity map of $X$ to a constant map $c_s$ (defined by
  some section $s$ of $X$), in such a way that two consecutive maps
  $f_i,f_{i+1}$ are $C_h$-homotopic in the strict sense, namely we can
  find $Y_i$ in $C_h$ and two sections $\delta_0^i$ and $\delta_1^i$
  of $Y_i$ over $e$, and a map
  \[ h_i : Y_i\times X\to X,\]
  such that
  \[ h_i\circ(\delta_0^i\times\id_X)=f_i,\quad
  h_i\circ(\delta_1^i\times\id_X)=f_{i+1}.\]
\item[\namedlabel{cond:51.Sat.4}{(Sat~4)}]
  \emph{Any} object $X$ in $M$ satisfying the condition just stated is
  in $C_h$.
\end{description}

The third dotted arrow can be viewed as denoting an embedding of
$\Hom_4(M)$ (the set of all contractibility structures on $M$) into
$\Hom_3(M)$, we finally get a cascade of three inclusions
\[\Hom_4(M) \hookrightarrow
\Hom_3(M) \hookrightarrow
\Hom_2(M) \hookrightarrow
\Hom_1(M),\]
in terms of which a contractibility structure on $M$ can be viewed a a
particular case of any of the three types of homotopy structures on
$M$ considered before.

If we start with a homotopism structure $W_h$ on $M$, and consider the
corresponding set $C_h$ of contractible objects of $M$ (namely objects
$X$ such that $X\to e$ is in $W_h$), it is pretty clear that $C_h$
satisfies the saturation condition \ref{cond:51.Sat.4}, but by no
means clear that it satisfies the basic assumption
\ref{cond:51.Bas.4}, even in the special case when we assume moreover
that $W_h$ comes from a weak homotopy interval structure $\Sigma_h$ on
$M$. The trouble comes from the circumstance that there is no reason
in general that the contractibility of an object $X$ of $M$ can be
described in terms of a sequence of elementary homotopies between maps
$f_i:X\to X$ (joining $\id_X$ to a constant map) involving weak
homotopy intervals $\bI_i$ \emph{which are moreover contractible}. I
doubt this is always so, and there doesn't come either any plausible
extra condition on $\Sigma_h$ which may ensure this, except precisely
that $\Sigma_h$ can be generated (through saturation)\pspage{119} by
the subset $\Sigma_{hc}$ of its contractible elements, which is just
another way of saying that this $\Sigma_h\in\Hom_3(M)$ comes already
from a contractibility structure $C_h\in\Hom_4(M)$! Thus, definitely
the uniformity of formal relationships between successively occurring
notions seems broken here, namely there does not seem to be any
natural retraction $r_{43}$ of $\Hom_3(M)$ onto the subset
$\Hom_4(M)$. For the least, if there is such a retraction, its
definition should be presumably a somewhat more delicate one than the
first that comes to mind. I will not pursue this matter any further
now, as it is not clear if we'll need it later.

It is clear that for any weak homotopy interval structure $\Sigma_h$
on $M$, $\Sigma_h$ is stable under the natural notion of finite
products of intervals (in the sense of the category structure of
$\bInt(M)$). We saw already that this is handy, as the consideration
of products of intervals allows to show that the family of homotopy
relations $\Isim$ in $\Fl(M)$, for variable \bI{} in $\Sigma_h$, is
``filtrant d\'ecroissant'', so we get the relation $\hsim$ as the
filtering direct limit or union of the more elementary relations
$\Isim$. Similarly, if $C_h\subset\Ob(M)$ is a contractibility
structure on $M$, $C_h$ is stable under finite products.
\begin{remark}
  From the way we've been working so far with homotopy notions, it
  would seem that we're only interested here in homotopy notions which
  stem from a structure in $\Hom_3(M)$, namely which can be described
  in terms of a notion of weak homotopy intervals. The focus on
  contractibility has set in only lately, and it is too soon to be
  sure whether we'll be working only with homotopy structures on $M$
  which can be described in terms of a contractibility notion, namely
  which are in $\Hom_4(M)$. In the cases I've had in mind so far, it
  turns out, it seems that the homotopy notions dealt with do come
  from a structure in $\Hom_4(M)$, i.e., from a contractibility
  structure.
\end{remark}

\subsection[Generating sets of homotopy intervals. Two standard ways
of generating multiplicative intervals. Contractibility of
\texorpdfstring{$\Hom(X,Y)$}{Hom(X,Y)}'s]{Generating sets of weak
  homotopy intervals. Contractors.}\label{subsec:51.E}
Let $\Sigma_h$ be a weak homotopy interval structure on
$M$.\footnote{This implies we assume $M$ stable under finite
  products.} A subset $\Sigma_h^0$ is called \emph{generating}, if
$\Sigma_H$ is just its saturation (cf.\ \ref{subsec:51.C} above),
i.e., for any \bI{} in $\Sigma_h$, the two endpoints can be joined by
a chain as in \ref{cond:51.Sat.3}, involving only intervals in $\Sigma_h^0$.
This implies that all homotopy notions dealt with so far can be
checked directly in terms of intervals in $\Sigma_h^0$. We've met a
particular case of this before, when $\Sigma_h^0$ is reduced to just
one element \bI{} -- we then called \bI{} ``characteristic'', but
``\emph{generating weak homotopy interval}'' now would seem the more
appropriate expression. Even when\pspage{120} there should not exist
such a generating interval, the natural next best assumption to make
is the existence of a generating set $\Sigma_h^0$ which is ``small''
(namely an element of the ``universe'' we are working in). The case of
a finite generating set of intervals reduces to the case of a single
one though, by just taking the product of those intervals.

An interesting case is when the generating set $\Sigma_h^0$ consists
of contractible objects of $M$. Such a generating set exists if{f} the
structure considered $\Sigma_h$ comes from a contractibility
structure. About the best we could hope for is the existence of a
single \emph{generating contractible weak homotopy interval} \bI. If
we got any interval \bI{} in $M$, this can be viewed as a generating
contractible weak homotopy interval for a suitable homotopy structure
on $M$ (then necessarily unique) if{f} the identity map of $I$ can be
joined to a constant one by a chain of maps, such that two consecutive
ones are tied by an \bI-homotopy or an $\check\bI$-homotopy. The most
evident way to meet this condition is by a \emph{one-step chain} from
$\id_I$ to the constant map defined by one of the endpoints,
$\delta_1$ say. This brings us back to structure of a composition law
\[ I\times I\to I\]
in $I$, having $e_0$ as a left unit and $e_1$ as a left zero
element. Let's call an interval, endowed with such a composition law,
a \emph{contractor} in $M$. Thus starting from a contractor in $M$ is
about the nicest way to define a homotopy structure in $M$, as a
matter of fact the strongest type of such a structure -- namely a
contractibility structure, admitting a \emph{generating contractible
  weak homotopy interval} (and better still, admitting a
\emph{generating contractor}).

Of course, starting with the weakest kind of homotopy structure on
$M$, namely just a homotopy relation $R_h\in\Hom_1(M)$, if \bI{} is a
homotopy interval which is moreover endowed with a structure of a
contractor, i.e., if it is a contractor such that the end-point
sections $\delta_0,\delta_1$ are homotopic, then $I$ is automatically
contractible (never mind if it is generating or not).

It seems to me that the homotopy structures I've looked at so far
(such as \Cat) and various standard elementary modelizers \Ahat) are
not only contractibility structures, but they all can be defined by a
single contractor each.

Besides the ``\emph{basic contractor}'' $\Simplex_1$ in \Cat, there are
two general ways I've met so far for getting contractors. One has been
made explicit in these notes a number of times, namely the
\emph{Lawvere element} $L_M$ in $M$\pspage{121} if it exists, and if
moreover $M$ has a \emph{strict} initial element,
$\varnothing_M$. Recall that $L_M$ represents the contravariant functor
on $M$
\[X\mapsto\text{set of all subobjects of $X$,}\]
and that the ``full'' and ``empty'' subobjects of $X$, for variable
$X$, define two sections $\delta_0$ and $\delta_1$ of $L_M$. I forgot
to state the extra condition that in $M$ fibered products exist
(intersection of two subobjects would be enough); then the
intersection law endows $L_M$ with a structure of a contractor
$\bL_M$, admitting $\delta_0$ as a unit and $\delta_1$ as a zero
element. Moreover, it is clear that $\bL_M$ as an interval is
separated, i.e., $\Ker(\delta_0,\delta_1)=\varnothing_M$. More precisely
still, $\bL_M$ can be viewed as a \emph{final object of the category
  of all separated intervals in $M$}, namely for any such interval,
there is a \emph{unique} map of intervals
\[\bI\to\bL_M.\]
This implies that if $M$ is endowed with a homotopy structure, such
that there exists a weak homotopy interval which is separated, then
$\bL_M$ is such an interval, and it is moreover contractible. It is
doubtful though, even if we can find a \emph{generating} contractor
for the given homotopy structure on $M$, that the Lawvere contractor
is generating too.

Here now is a second interesting way of getting contractors. We assume
that $M$ admits finite products (as usual). Let $X$ be an objects, and
assume the object $\Hom(X,X)$, representing the functor
\[Y \mapsto \Hom(X\times Y,X) = \Hom_Y(X_Y,X_Y),\]
exists in $M$. (NB\enspace $X_Y$ denotes $X\times Y$, viewed as an object of
$M_{/Y}$.) Composition of endomorphisms of $X_Y$ clearly endow this
functor with an associative composition law, admitting a two-sided
unit, which I call $e_0$. Notice that sections of
\[\bI=\bHom(X,X)\]
can be identified with maps $X\to X$, and the section corresponding to
$\id_X$ is of course the two-sided unit. On the other hand, if $X$
admits sections, i.e., admits ``constant'' endomorphisms, it is
clear that the corresponding sections of \bI{} are left zero
elements. If we choose a section of $X$, \bI{} becomes a
contractor. Its interest lies in the following
\begin{proposition}
  Assume finite products exist in $M$, and $M$ endowed with a homotopy
  structure. Let $X$ be an object of $M$ endowed with a
  section $e_X$, and suppose the object $\bHom(X,X)$ exists, hence a
  contractor \bI{} as seen above. The following two conditions are
  equivalent:
  \begin{enumerate}[label=\alph*),font=\normalfont]
  \item\label{it:51.E.a}
    $X$ is contractible.
  \item\label{it:51.E.b}
    \bI{} is contractible \textup(or, equivalently as seen above,
    \bI{} is a weak homotopy interval, namely the two endpoints are homotopic\textup).
  \end{enumerate}
  Moreover,\pspage{122} this condition implies the following two:
  \begin{enumerate}[label=\alph*),font=\normalfont,resume]
  \item\label{it:51.E.c}
    For any object $Y$ in $M$, if $\bHom(Y,X)$ exists, it is contractible.
  \item\label{it:51.E.d}
    For any $Y$ in $M$, if $\bHom(X,Y)$ exists, the natural map
    \[ Y \mapsto \bHom(X,Y)\]
    \textup(identifying $Y$ to the ``subobject of constant maps from
    $X$ to $Y$''\textup) is a homotopism.
  \end{enumerate}
\end{proposition}

The equivalence of \ref{it:51.E.a} and \ref{it:51.E.b} is just a
tautological translation of contractibility and homotopy relations in
terms of weak homotopy intervals \bJ{} (cf.\ cor.\ \hyperref[cor:51.E.2]{2}
below). That \ref{it:51.E.c} and \ref{it:51.E.d} follow comes from the
fact that the monoid object $I=\bEnd(X)=\bHom(X,X)$ operates on the
left on $\bHom(Y,X)$, on the right on $\bHom(X,Y)$, and the following:
\begin{corollarynum}\label{cor:51.E.1}
  Let \bI{} be a weak homotopy interval, assume the underlying $I$
  ``operates'' on an object $H$, namely we are given a map
  \[ h:I\times H\to H\quad(\text{``operation'' of $I$ on $H$})\]
  satisfying the relations \textup(where $h(u,f)$ is written simply
  $u\cdot f$\textup{):}
  \[ e_0 \cdot f = f, \quad e_1\cdot (e_1 \cdot f) = e_1 \cdot f,\]
  namely $e_0$ acts as the identity and $e_1$ acts as an idempotent
  $p$ on $H$ \textup(a very weak associativity assumption indeed if
  $I$ is a contractor, as $e_1\cdot e_1=e_1$\textup). Assume the image
  of $p$, i.e., $\Ker(\id_H,p)$ exists, let $H_0$ be the corresponding
  subobject of $H$, and
  \[p_0:H\to H_0\]
  the map induced by $p$. Then $p_0$ is a homotopism \textup(and hence
  the inclusion $i:H_0\to H$, which is a section of $p_0$, is a
  homotopism too\textup).
\end{corollarynum}

Because of the saturation property \ref{it:48.cprime} on homotopies, it
is enough to check that $p=ip_0$ is a homotopism (as $p_0i=\id_H$
already is one), and for this it is enough to see it is homotopic to
the identity map of $H$. But a homotopy between the two is realized by
$h$, qed.

The argument for equivalence of \ref{it:51.E.a} and \ref{it:51.E.b}
above can be generalized as follows:
\begin{corollarynum}\label{cor:51.E.2}
  Let $M$ be as before, and $X$ and $Y$ objects such that $H =
  \bHom(X,Y)$ exists in $M$. Let $f,g:X\rightrightarrows Y$ be two
  maps, which we'll identify to the corresponding sections of
  $H$. Then $f$ and $g$ are homotopic maps if{f} they give rise to
  homotopic sections of $H$.
\end{corollarynum}

\subsection{The canonical homotopy structure: preliminaries on
  \texorpdfstring{$\pi_0$}{pi0}.}\label{subsec:51.F}
In order to simplify life, I will in this section make the following
assumptions on $M$ (which presumably, except for the first, could
be\pspage{123} considerably weakened, but these will be sufficient):
\begin{enumerate}[label=\alph*)]
\item\label{it:51.F.a}
  Finite products exist in $M$ (``pour m\'emoire'').
\item\label{it:51.F.b}
  Arbitrary sums exist in $M$, they are ``disjoint'' and ``universal''
  (which implies that $M$ has a \emph{strict} initial object).
\item\label{it:51.F.c}
  Every object in $M$ is isomorphic to a direct sum of $0$-connected ones.
\end{enumerate}

I recall an object is called \emph{$0$-connected} if it is a)
``non-empty'', i.e., non-isomorphic to $\varnothing_M$, and b)
connected, i.e., any decomposition of it into a sum of two subobjects
is trivial (namely, one is ``empty'', the other is ``full''). Also,
under the assumption \ref{it:51.F.b}, well-known standard arguments
show that for any object $X$, a decomposition of $X$ into a direct sum
of $0$-connected components is essentially unique (namely the
corresponding set of subobjects of $X$ is unique), if it exists. The
subobjects occurring in the sum are called the \emph{connected
  components} of $X$, and the set of connected components of $X$ is
denoted as usual by $\pi_0(X)$. We thus get a natural functor
\[ X \mapsto \pi_0(X), \quad M \to \Sets,\]
This can be equally described as a \emph{left adjoint} to the functor
\[ E \mapsto E_M, \quad \Sets \to M,\]
associating to a set $E$ the corresponding ``\emph{constant object}''
$E_M$ of $M$, sometimes also denoted by the product symbol $E\times
e_M = E\times e$ ($e$ the final object of $M$), namely the direct sum
in $M$ of $E$ copies of $e$. The adjunction formula is
\[\Hom_M(X,E_M) \simeq \Hom_\Sets(\pi_0(X),E).\]
The adjunction relation implies that the functor $\pi_0$ commutes with
all direct limits which exist in $M$, and in particular (and trivially
so) to direct sums.

I'll finally assume also, to fix the ideas:
\begin{enumerate}[label=\alph*),resume]
\item\label{it:51.F.d}
  The final object of $M$ is $0$-connected, i.e., $\pi_0(e)=$
  one-point set.
\end{enumerate}

Empty objects of $M$, on the other hand, are of course characterized
by the condition
\[\pi_0(\varnothing_M)=\emptyset.\]

Finally, I'll make in the end a very strong assumption on $M$, which
however is satisfied more or less trivially in the cases we are
interested in, when $M$ is a would-be modelizer:
\begin{enumerate}[label=\alph*),resume]
\item\label{it:51.F.e}
  (Total $0$-asphericity of $M$): the product of two $0$-connected
  objects of $M$ is again $0$-connected.
\end{enumerate}

This is clearly equivalent to the condition
\begin{enumerate}[label=\alph*'),start=5]
\item\label{it:51.F.eprime}
  The functor $\pi_0:M\to\Sets$ commutes to finite products.
\end{enumerate}
\begin{remarks}
  The\pspage{124} crucial assumptions here seem to be \ref{it:51.F.b}
  (which allows definition of a $\pi_0$ functor, and topological
  intuition tied up with connectedness to enter into play), and
  \ref{it:51.F.e}, which implies that with respect to cartesian
  products, the usual intuitive background for connectedness, rooted
  in the example of $M = (\text{topological spaces})$ is indeed
  valid. This condition is clearly stronger than \ref{it:51.F.d},
  which is a mere condition for convenience in itself (otherwise, a
  decomposition of $e$ into connected components would mean a
  corresponding decomposition of $M$ as a product category, and
  everything could be looked at ``componentwise''). As for
  \ref{it:51.F.e} it could probably be dispensed with, by still
  defining $\pi_0(X)$ as a strict pro-set. In the case of modelizers
  anyhow, such generalization seems quite besides the point.

  The condition \ref{it:51.F.e}, however innocent-looking in terms of
  topological intuition, seems to me an extremely strong condition
  indeed. I suspect that in case $M$ is a topos, it is equivalent to
  total asphericity. In any case, if $M$ is the topos associated to a
  locally connected topological space, we've seen time ago that the
  condition \ref{it:51.F.e} implies $X$ is irreducible, and hence
  totally aspheric. In view of this exactingness of \ref{it:51.F.e},
  I'll not use it unless specifically stated.
\end{remarks}

\bigbreak
\presectionfill\ondate{3.4.}\par

\hangsection{Inaccuracies rectified.}\label{sec:52}%
Before pursuing the review of ``pure homotopy notions'' begun in
yesterday's notes, I would like to correct some inaccuracies which
flew in when looking at the relationship between the two first basic
homotopy notions, namely the notion of a homotopy relation, and the
notion of a homotopism structure. As usual, the provisional image I
had in mind was still somewhat vague, while the reasonable
expectations came out more clearly through the process of writing
things down (including factual inaccuracies!).

The two notions clearly correspond to two kinds of ways of
constructing new categories $M'$ in terms of a given one, and a
functor
\[ M \to M'\]
which is bijective on objects, and has the property moreover that for
any category $C$, the corresponding functor
\[\bHom(M',C) \to \bHom(M,C)\]
is a fully faithful embedding in the strict sense, namely injective
on\pspage{125} objects. One way is to take for $M'$ any quotient
category of $M$, by an equivalence relation which is the discrete one
on objects -- thus it corresponds just to an equivalence relation $R$
in $\Fl(M)$, compatible with the source and target maps and with
compositions. The other is to take as $M'$ a localization with respect
to some $W\subset\Fl(M)$, and if we take $W$ strongly saturated we get
a bijective correspondence between these $M'$ and the set of strongly
saturated subsets of $\Fl(M)$.

It wouldn't be any more reasonable to call an arbitrary $R$ as above,
corresponding to an arbitrary quotient category $M'$ with the same
objects as $M$, a ``homotopy relation'' on $M$ (as I did though
yesterday), as it would be to call an arbitrary strongly saturated
$W\subset\Fl(M)$ a ``homotopism structure'' on $M$ (as I nearly did
yesterday, but then rectified in the stride). \emph{The characteristic
  flavor of homotopy theory comes in, when we get an $M'$ which is
  \emph{both} a quotient category and a localization of $M$.} Thus
neither approach, via $R$ or via $W$, is any more contained in the
other, then the converse. We should regard the \emph{homotopy
  structure} on $M$ to be embodied in the basic functor
\[M\to M',\]
which is a description where no choice yet is made between the two
possible descriptions of $M'$, either by an $R$, or by a $W$. If we
describe $M'$ in terms of $R$, the extra assumption to make on $R$ for
calling it a ``\emph{homotopy relation}'', is that the canonical
functor
\[ M \to M/R\]
should be localizing. Alternatively, describing $M'$ by $W$, the extra
assumption to make on $W$ (as we did yesterday) in order to call $W$ a
\emph{homotopism structure} on $M$, is that the canonical localization
functor
\[ M \to W^{-1}M\]
be essentially a passage-to-quotient functor, namely surjective on
arrows (as we know already it is bijective on objects). Thus the set
of all homotopy relations on $M$ is in one-to-one correspondence with
the set of all homotopism structures on $M$, and if we denote these
sets (in accordance with yesterday's provisional notations)
$\Hom_1(M)$ and $\Hom_2(M)$, we get thus a \emph{bijective}
correspondence
\[ \Hom_1(M) \leftrightarrow \Hom_2(M).\]
The set of homotopy structures $\Hom(M)$ on $M$ may be either defined
as the usual quotient set defined by the previous two-member system of
transitive bijections between sets, or more substantially, as a set of
isomorphism classes of categories $M'$ ``under M'', i.e., endowed with
a functor $M\to M'$, and subject to the following two extra
conditions:
\begin{enumerate}[label=\alph*)]
\item\label{it:52.a}
  The\pspage{126} functor $M\to M'$ is bijective on objects,
  surjective on arrows.
\item\label{it:52.b}
  The functor $M\to M'$ is a localization functor (it will be so then,
  in view of \ref{it:52.a}, in the strict sense, namely $M'$ will be
  $M$-isomorphic, not only $M$-equivalent, to a localization $W^{-1}M$).
\end{enumerate}

However, the question arises whether it is possible to define such a
homotopy structure on $M$ in terms of an arbitrary $R$, i.e., an
arbitrary quotient category having the same objects (let $Q(M)=Q$ be
the set of all such $R$'s) or in terms of an arbitrary localization of
$M$, or what amounts to the same, in terms of an arbitrary strongly
saturated $W\subset\Fl(M)$ (let's call $L(M)=L$ the set of all such
$W$'s). The first thing that comes to mind here, is that we got two
natural maps
\begin{equation}
  \label{eq:52.1}
  \begin{tikzcd}[cramped]
    Q \ar[r, bend left, "r"] &
    L \ar[l, bend left, "s"]
  \end{tikzcd}
  \tag{1}
\end{equation}
between $Q$ and $L$, which are defined by the observation that
whenever we have a functor $i:M\to M'$, injective on objects, it
defines both an $R\in Q$ (namely $f\Rsim g$ if{f} $i(f)=i(g)$) and a
$W\in L$ (namely $f\in W$ if{} $i(f)$ is invertible). For defining
$r(R)$ resp.\ $s(W)$, we apply this to the case when $M'=M/R$ resp.\
$W^{-1}M$. Let's look a little at the two cases separately.

Start with $R$ in $Q$, we get $W=r(R)$,
\[ W = \set[\big]{f\in\Fl(M)}{\text{$i(f)$ invertible}},\]
where
\[ i : M\to M/R\]
is the canonical functor, we thus get a canonical functor (compatible
with the structures ``under $M$'')
\begin{equation}
  \label{eq:52.2}
  \alpha_R : M_W \to M_R,
  \tag{2}
\end{equation}
where for simplicity $I$ write
\[ M_R = M/R,\quad M_W = W^{-1}M.\]
We may define $W$ as the largest element in $L$ (for the natural order
relation in $L$, namely inclusion of subsets of $\Fl(M)$) such that a
functor \eqref{eq:52.2} exists (compatible with the functors from
$M$ into both sides) -- such a functor of course is unique (by the
preliminaries on functors $M\to M'$ made at the beginning). In terms of
\eqref{eq:52.1}, we can say that $R$ is actually a homotopy relation
(let's call $Q_0(M)=Q_0$ the subset of $Q$ of all such relations)
if{f} \eqref{eq:52.2} is an isomorphism, or equivalently (as it is
clearly bijective on objects, surjective on arrows) if{f} it is
\emph{injective on arrows}, i.e., \emph{faithful}.

Conversely, start now with $W$ in $L$, we get $R=s(W)$,
\[ R = \set[\big]{(f,g)\in\Fl(M)\times\Fl(M)}{i(f)=i(g)},\]
where\pspage{127} now
\[ i : M\to W^{-1}M=M_W\]
is the canonical functor defined in terms of $W$; we thus get a
canonical functor of categories ``under $M$''
\begin{equation}
  \label{eq:52.3}
  \beta=\beta_W : M_R \to M_W,
  \tag{3}
\end{equation}
as a matter of fact, $R$ is the largest element in $Q$ (for the
inclusion relation of subsets of $\Fl(M)\times\Fl(M)$) for which a
functor \eqref{eq:52.3} exists (then necessarily unique, as
before). We may say that $W$ is a homotopism structure on $M$, i.e.,
$W\in L_0$ (where $L_0$ is the subset of $L$ of all homotopism
structures on $M$) if{f} the functor \eqref{eq:52.3} is an
isomorphisms, or equivalently (as it is clearly bijective on objects,
injective on arrows) if{f} it is \emph{surjective on arrows}.

We may describe $Q_0$ and $L_0$ in a purely set-theoretic way, in
terms of the system $(r,s)$ of maps in \eqref{eq:52.1}, by the formula
(which is just a translation of the definitions of $Q_0$ and $L_0$)
\begin{align*}
  Q_0 &= \set[\big]{q\in Q}{sr(q)=q} \\
  L_0 &= \set[\big]{\ell\in L}{rs(\ell)=\ell} , 
\end{align*}
and we can describe formally the pair of subsets $(Q_0,L_0)$ of $Q,L$
as the largest pair of subsets, such that $r$ and $s$ induce between
$Q_0$ and $L_0$ bijections inverse of each other. In the general
set-theoretic set-up, it is by no means clear, and false in general,
that $r$ maps $Q$ into $L_0$ or $L$ into $Q_0$ (thus both $Q_0$ and
$L_0$ may well be empty, whereas $Q$ and $L$ are not). Thus it is not
clear at all that starting with an arbitrary $R\in Q$, the
corresponding $W=i(R)$ is a homotopism structure, and it is easily
seen that this is \emph{not} so in general, contrarily to what I
hastily stated in yesterday's notes. (Take $M$ which just one object,
therefore defined by a monoid $G$, and $M_R$ corresponds to a quotient
monoid $G'$, we may take $G'=$ unit monoid, thus $W$ is $G$ itself,
and $M_W$ the enveloping group $\overline G$ of $G$ -- the map
$\overline G\to G'$ need not be injective!) In the opposite direction,
starting with an arbitrary $W$ in $L$, the corresponding $R$ need not
be a homotopy relation. (If $M$ is reduced to a point, this amounts to
saying that if we localize a monoid $G$ with respect to a subset
$W\subset G$, by making invertible the elements in $W$, the
corresponding map $G\to \overline G$ need no be surjective!)

What we may say, though, is that if we start with a pair
\[(R,W) \in Q\times L\]
such that the \emph{two} functors \eqref{eq:52.2}, \eqref{eq:52.3}
exist, i.e., such that
\[ W\le r(R)\quad\text{and}\quad R\le s(W),\]
then $(R,W)\in Q_0\times L_0$, i.e., $R$ is a homotopy relation and
$W$ is a homotopism\pspage{128} structure, and the two are
associated. This comes from the fact that both compositions of the
functors \eqref{eq:52.2}, \eqref{eq:52.3} must be the identity
functors (being compatible with the ``under $M$'' structure), hence
$\alpha$ and $\beta$ are isomorphisms, which shows both $R\in Q_0$ and
$W\in L_0$, and that $R$ and $W$ are associated. In terms of the
set-theoretic situation \eqref{eq:52.1}, this may be described by
using the order relations on $Q$ and $L$, and the fact that $r$ and
$s$ are monotone maps, and satisfy moreover
\begin{equation}
  \label{eq:52.star}
  sr(q)\le q, \quad rs(\ell)\le\ell\quad (\text{any $q\in Q, \ell\in L$}),
  \tag{*}
\end{equation}
which implies that the set $C_0$ of pairs $(q,\ell)$ of associated
elements of $Q_0,L_0$ can be described also as
\[ C_0 = \set[\big]{(q,\ell)\in Q\times L}{\ell\le r(q),\; q\le
  s(\ell)}.\]

Thus it doesn't seem evident to get a homotopy structure on $M$, just
starting with an $R\in Q$ or a $W\in L$, without assuming beforehand
that $R$ is a homotopy relation, or $W$ a homotopism structure. The
condition \ref{cond:51.C.12} on page \ref{p:115} may be viewed as a
condition on a pair $(R,W)\in Q\times L$, and it clearly implies
\[ R \le s(W);\]
if we assume moreover $W=r(R)$ we get $R\le s(r(R))$ and hence, in
view of the first inequality \eqref{eq:52.star} above,
\[ R = sr(R), \quad\text{i.e.,}\quad R\in Q_0,\]
i.e., $R$ is a homotopy relation, as asserted somewhat quickly
yesterday.

The only standard way for getting homotopy structures in a general
category $M$ which I can see by now, is in terms of an arbitrary set
$\Sigma_h$ of intervals in $M$ (assuming only that $M$ admits finite
products). As soon as the focus gets upon intervals for describing
homotopy structures, the situation becomes typically non-autodual --
in contrast to Quillen's autodual treatment of homotopy relations (via
``left'' and ``right'' homotopies, involving respectively ``cylinder''
and ``path'' objects). This is in keeping with the highly non-autodual
axiom on universal disjoint sums in $M$, which we finally introduced
by the end of yesterday's reflection.

To come back however upon the relationships between the four basic
``homotopy notions'' introduced in yesterday's notes, I would now
rather symbolize these relations in the following diagram of maps
between sets\pspage{129}
\[\begin{tikzpicture}[commutative diagrams/every diagram]
  \matrix[matrix of math nodes, name=m, commutative diagrams/every cell] {
    Q(M) & L(M) \\
    Q_0=\Hom_1(M) & \Hom_2(M)=L_0 \\
    \Hom_3(M) & \Hom_4(M) \\};
  \path[commutative diagrams/.cd, every arrow, every label]
  (m-1-1) edge[bend left=10] node {$r$} (m-1-2)
  (m-1-2) edge[bend left=10] node {$s$} (m-1-1)
  (m-2-1) edge[commutative diagrams/hook] (m-1-1)
  (m-2-1) edge[<->] node {$\sim$} (m-2-2)
  (m-2-2) edge[commutative diagrams/hook] (m-1-2)
  (m-3-1) edge[commutative diagrams/hook] (m-2-1)
  (m-3-2) edge[commutative diagrams/hook] (m-2-2)
  (m-3-2) edge[commutative diagrams/hook] (m-3-1);
  \path[commutative diagrams/.cd, every arrow, every label]
  (m-1-1) edge[controls=+(190:15mm) and +(150:10mm)] (m-3-1.north west);
  \node (P1) [right=6pt,text width=4cm,align=left] at (m-2-2.east) {("homotopy
    structures" on $M$, defined in terms of homotopy relations, or
    homotopism structures)};
  \draw[decorate,decoration={brace,mirror,amplitude=6pt},yshift=0pt]
    (P1.north west) -- (P1.south west);
  \node (P2) [below,text width=3cm,align=left,xshift=3mm] at (m-3-1.south) {(weak
    homotopy interval structures)};
  \node (P3) [below,text width=3cm,align=left,xshift=5mm] at (m-3-2.south) {(contractibility
    structures)};
  \node[right,yshift=-3mm] at (P3.east) {,};
\end{tikzpicture}\]
which\footnote{D\'eployer le diagramme, trop tass\'e dans les deux
  dimensions. \scrcommentinline{I think I've fixed that; I also reversed the brace for
  clarity.}} in terms of the preceding reflections, and yesterday's, is
self-explanatory. The vertical arrow from $Q(M)$ to $\Hom_3(M)$ is the
canonical retraction -- in terms of the latter, and its composition
with $s:L(M)\to Q(M)$, there \emph{are} ways after all to associate to
any $R$ in $Q(M)$ or $W$ in $L(M)$ a homotopy structure on $M$,
provided only $M$ admits finite products, by using weak homotopy
intervals. If $h$ is the homotopy structure thus defined, we get a
priori a functor
\[ M_h \to M_R=M/R\quad\text{resp.}\quad M_h\to M_W=W^{-1}M,\]
provided we assume $R$ resp.\ $W$ satisfy the ``accessory
assumption'', namely that the corresponding functor $M\to M'$ (where
$M'$ is either $M_R$ or $M_W$) commute with finite products.

The main fact to remember from this whole discussion, it seems to me,
is that there are not really \emph{four}, but only \emph{three}
essentially distinct types of structure (among yesterday's) we may
consider upon $M$ as ``homotopy-flavored'' structures, namely
\[\parbox[c]{0.345\textwidth}{\raggedright
  homotopy structures $\Hom_1(M)\leftrightarrow\Hom_2(M)$}
\Leftarrow
\parbox[c]{0.26\textwidth}{\raggedright
  weak homotopy interval structures $\Hom_3(M)$}
\Leftarrow
\parbox[c]{0.30\textwidth}{\raggedright
  contractibility structures $\Hom_4(M)$.}
\]
It would seem at present that the homotopy structures that naturally
come up in our present ``modelizer story'' are all of the strictest
type, and even describable in terms of just one \emph{generating
  contractible weak homotopy interval} (I would like to drop the
qualification ``weak'', definitely when a contractibility assumption
comes in!), and even a \emph{generating contractor}, with commutative
idempotent composition law.\pspage{130}

\hangsection[Compatibility of a functor $u: M\to N$ with a
homotopy \dots]{Compatibility of a functor \texorpdfstring{$u: M\to N$}{u:
    M->N} with a homotopy structure on
  \texorpdfstring{$M$}{M}.}\label{sec:53}%
Before pursuing yesterday's reflection about the $\pi_0$-functor and
its relation to homotopy structures on $M$, it seems more convenient
to interpolate some more or less obvious ``functorialities'' on the
homotopy notions just developed. They all seem to turn around the
relationship of such notions in $M$ with a more or less arbitrary
functor
\[ u : M\to N,\]
for the time being I am not making any special assumption on $u$. In
terms of the four ways we got for describing a homotopy structure in
$M$, we get four corresponding natural conditions of ``compatibility''
of $u$ with a given homotopy structure $h$, namely:
\begin{description}
\item[\namedlabel{it:53.i}{(i)}]
  If $f \hsim g$ in $M$, i.e. $(f,g)\in R_h$, then $u(f)=u(g)$.
\item[\namedlabel{it:53.iprime}{(i')}]
  If $\bI=(I,\delta_0,\delta_1)\in\Sigma_h$ is a weak homotopy
  interval, then $u(\delta_0)=u(\delta_1)$.
\item[\namedlabel{it:53.ii}{(ii)}]
  If $f\in W_h$, i.e., $f$ is a homotopism, then $u(f)$ is an isomorphism.
\item[\namedlabel{it:53.iiprime}{(ii')}]
  If $X\in C_h$ is a contractible object, then $u(X)$ is a final
  object (NB\enspace here we assume that $u(e_M)$ is a final object in $N$).
\end{description}
(NB\enspace It is understood implicitly, whenever dealing with intervals and
with contractible objects, that $M$ admits finite products.)

The conditions \ref{it:53.i} and \ref{it:53.ii} are clearly
equivalent, and equivalent to the requirement that $u$ factors into
\begin{equation}
  \label{eq:53.1}
  M \to M_h \to N,\tag{1}
\end{equation}
where $M\to M_h$ is the canonical functor of $M$ into the
corresponding \emph{homotopy-types} category. We also have the
tautological implications \ref{it:53.i} $\Rightarrow$
\ref{it:53.iprime} and \ref{it:53.ii} $\Rightarrow$
\ref{it:53.iiprime}. Moreover we have (trivially) \ref{it:53.iiprime}
$\Rightarrow$ \ref{it:53.iprime} whenever the homotopy structure on
$M$ is a weak homotopy interval structure, and \emph{moreover} the
functor $u$ commutes to finite products. All these implications are
summarized in the diagram
\[\begin{tikzcd}[baseline=(O.base),math mode=false,%
  arrows=Rightarrow]
  \ref{it:53.i} \ar[r,Leftrightarrow]\ar[d] &
  \ref{it:53.ii} \ar[d] \\
  \ref{it:53.iprime} \ar[u,bend left,"*\;\;" pos=0.45, "$\Hom_3$"] &
  |[alias=O]| \ref{it:53.iiprime} \ar[l,swap,"$\Hom_4$"]
\end{tikzcd},\]
where the symbol $\Hom_3$ or $\Hom_4$ indicates that the implication
qualified by it is valid provided we assume that the homotopy
structure on $M$ is in $\Hom_3$ (namely is defined in terms of weak
homotopy intervals) resp.\ in $\Hom_4$ (namely is a contractibility
structure), and where (*) denotes the extra assumption of commutation
of $u$ with finite products.

We'll\pspage{131} say the functor $u$ is \emph{compatible with the
  homotopy structure} $h$ on $M$, if it satisfied the equivalent
conditions \ref{it:53.i}, \ref{it:53.ii}, i.e., if it factors as in
\eqref{eq:53.1} above. In case $u$ commutes with finite products, and
if either the homotopy structure $h$ can be described in terms of weak
homotopy intervals, or in terms of contractible objects, the
compatibility of $u$ with $h$ can be checked correspondingly, either
by \ref{it:53.iprime}, or by \ref{it:53.iiprime}.

\hangsection[Compatibility of a homotopy structure with a set $W$ of
\dots]{Compatibility of a homotopy structure with a set
  \texorpdfstring{$W$}{W} of ``weak equivalences''. The homotopy
  structure \texorpdfstring{$h_W$}{h-W}.}\label{sec:54}%
An important particular case is the one when
\[N = W^{-1}M\]
is a localization on $M$ by a set of arrows in $M$
\[W\subset \Fl(M).\]
We'll say that the homotopy structure $h$ on $M$ is \emph{compatible
  with} $W$, if it is with the canonical functor $M\to W^{-1}M$. If
$W$ is \emph{strongly saturated}, this is most readily expressed by
the condition that $W_h\subset W$, i.e., any homotopism is in $W$
(i.e., ``any homotopism is a weak equivalence'', if the elements of
$W$ are named ``weak equivalences''); in case $M$ admits finite
products and the localization functor commutes to these (e.g., the
case $(M,W)$ is a \emph{strict} modelizer), and if moreover $h$ is a
contractibility structure, it is sufficient to check that for any
contractible $X$, the projection $X\to e$ is in $W$.

If we don't assume or know beforehand that $W$ is strongly saturated,
but just saturated say, we may still introduce a more stringent
compatibility condition, by saying that the homotopy structure $h$ and
$W$ are \emph{strictly compatible} if $W_h\subset W$. Using the
saturation condition \ref{it:48.cprime} on $W$, it is easily seen that
in the case when $h$ is a contractibility structure, then $W_h\subset
W$ (strict compatibility) is equivalent to: for contractible $X$, the
projection $X\to e$ is ``universally in $W$'', or (as we'll say)
\emph{$W$-aspheric}. Indeed, to deduce from this that any homotopism
$f:X\to Y$ is in $W$, we are reduced to checking that any endomorphism
of either $X$ or $Y$ which is homotopic to the identity map, is in
$W$. Now this will follow from the assumption, and the following
\begin{proposition}\label{prop:54}
  Assume the homotopy structure $h$ on $M$ can be defined by a
  generating set $\Sigma_h^0$ of weak homotopy intervals
  $\bI=(I,\delta_0,\delta_1)$ which are \emph{$W$-aspheric} \textup(i.e., $I$
  $W$-aspheric over $e$\textup), where $W\subset\Fl(M)$ is any saturated
  subset \textup(in fact, mildly saturated is enough\textup). Then $W$ is the
  inverse image by the canonical functor $M\to M_h=\overline M$ of a
  subset $\overline W\subset\Fl(\overline M)$, i.e., if $f,g$ are
  homotopic arrows in $M$, if one is in $W$, so is the other. Moreover
  \textup(if $W$ is saturated\textup) $W_h\subset W$, i.e., $h$ and
  $W$ are strictly compatible.
\end{proposition}

The\pspage{132} first statement is just the
``\hyperref[lem:hlr]{homotopy lemma}'' part \ref{it:48.hlr.b} (page
\ref{p:99}), the second follows by the argument sketched above.

We're about back now to the context we started with three days ago
(par.\ \ref{sec:48}, page \ref{p:98} and following), where we started
with a $W$ (viewed as a notion of ``weak equivalence''), and in terms
of $W$ constructed various homotopy notions -- namely those, we would
now say, corresponding to the homotopy structure defined by the set of
all intervals \bI{} in $M$ which are $W$-aspheric (i.e., $I$ is
$W$-aspheric over $e$\footnote{\alsoondate{11.4.} By which we mean that $I\to e$ is
  ``universally in $W$''. The terminology used here for
  ``$W$-aspheric'' is highly ambiguous, cf.\ discussion p.\
  \ref{p:181} and following.}). As a matter of fact, we were a little
stricter still, by restricting to intervals which are moreover
``disjoint'' (and which we called ``homotopy intervals'' relative to
$W$), but this restriction now definitely appears as awkward and
artificial. I will henceforth call \emph{homotopy intervals} (with
respect to $W$), any interval (not necessarily a separated one) which
is $W$-aspheric. Let $h_W$ be the corresponding homotopy structure on
$M$, which is a weak homotopy interval structure admitting the set of
all $W$-homotopy intervals as a generating set of weak homotopy
intervals. (Clearly, there will be many weak homotopy intervals for
this structure, which are far from being $W$-aspheric, i.e., far from
being homotopy intervals.) Of course, as stated in the preceding
proposition, $h_W$ and $W$ are strictly compatible, i.e.,
\[W_{h_W} \subset W,\]
i.e., any $h_W$-homotopism is in $W$ (i.e., is a ``weak
equivalence''). As a matter of fact, the definition of homotopy
notions in terms of $W$ we gave in loc.\ sit.\ were just the widest
one we could think of by that time, which would ensure the
``compatibility'' of these notions with $W$, in a sense which wasn't
technically clear (not even definable at that point) as it is now, but
however reasonably clear in terms of mathematical ``bon sens''.
At present though the question arises rather naturally whether the
homotopy structure $h_W$ we selected ``au flair'' by that time is
indeed the best one, namely the widest one, we could get. More
explicitly, this means whether the homotopy structure $h_W$ is the
\emph{widest} (in terms of the natural order relation considered in
the previous paragraph) among all those which are compatible with $W$
in the strict sense $W_{h_W}\subset W$. Now this is certainly not so,
if we are not a little more specific about restricting to homotopy
structure definable in terms of a weak homotopy interval
structure. For instance, if we take for $W$ a homotopism structure on
$M$, compatible with products, corresponding to a homotopy structure
$h$, to say $h_W$ is the ``best'' would imply that $W$ itself can be
described in terms of weak homotopy intervals\pspage{133} which is not
always the case. (Take for instance $M$ to be an abelian category, say
projective complexes of modules and quasi-isomorphism between these;
in this case, more generally whenever $M$ is a ``zero objects'' namely
one which is both initial and final, any interval in $M$ is trivial,
i.e., $\delta_0=\delta_1$, and hence any weak homotopy interval
structure on $M$ is trivial, namely $W_h$ is reduced to
isomorphisms\ldots)

Thus the more reasonable question here is whether any homotopy
structure $h$ on $M$, definable in terms of a weak homotopy interval
structure, and such that $W_h\subset W$, satisfies $h\le
h_W$. Clearly, for such an $h$, any weak homotopy interval \bI{} (for
$h$) satisfies $u(\delta_0)=u(\delta_1)$, where $u:M\to W^{-1}M$ is
the canonical functor (indeed, it is enough for this that
$W_h\subset\overline W$ instead of $W_h\subset W$, where $\overline W$
is the strong saturation of $W$), and conversely, if $W=\overline W$
and if moreover $u$ \emph{commutes to finite products}. On the other
hand, $h\le h_W$ means that any $\bI\in\Sigma_h$ is in $\Sigma_{h_W}$,
which also means that its endpoint sections $\delta_0,\delta_1$ are
$h_W$-homotopic, namely may be joined by a chain of sections, any two
consecutive of which are related by some \bJ-homotopy, where \bJ{} is
a \emph{$W$-aspheric interval}. Thus we get the:
\begin{proposition}
  Let $W$ a saturated set of arrows in $M$ \textup($M$ stable under
  finite products\textup), hence a corresponding homotopy structure
  $h_W$ on $M$, defined in terms of $W$-aspheric
  intervals\footnote{\alsoondate{11.4.} Cf.\ note on preceding page.\label{fn:54.star}} in $M$ as a
  generating set of weak homotopy intervals for $h_W$. Consider the
  following conditions:
  \begin{enumerate}[label=(\roman*),font=\normalfont]
  \item\label{it:54.i}
    $h_W$ is the widest of all homotopy structures $h$ on $M$ which
    are
    \begin{enumerate}[label=(\alph*),font=\normalfont]
    \item\label{it:54.i.a} strictly compatible with $W$, i.e., such
      that $W_h\subset W$ and moreover
    \item\label{it:54.i.b}
      definable in terms of a weak homotopy interval structure.
    \end{enumerate}
  \item\label{it:54.ii}
    For any object $I$ of $M$ and two sections $\delta_0,\delta_1$ of
    $I$ such that $u(\delta_0)=u(\delta_1)$ \textup(where $u:M\to
    W^{-1}M$ is the canonical functor\textup), $\delta_0$ and
    $\delta_1$ are $h_W$-homotopic, namely can be joined by a chain of
    elementary homotopies as above, involving
    \emph{$W$-aspheric}\footref{fn:54.star} intervals \bJ.
  \end{enumerate}
  Then \textup{\ref{it:54.ii}} implies \textup{\ref{it:54.i}}, and conversely if $W$ is
  strongly saturated and moreover $u$ commutes to finite products.
\end{proposition}

In connection with the $\pi_0$-functor, we are going to get pretty
natural conditions in terms of $0$-connectedness for ensuring
\ref{it:54.ii} and hence \ref{it:54.i}, which should apply I guess to
all ``reasonable'' modelizers $(M,W)$. It would thus seem that in
practical terms, the definition of $h_W$ is the best, in all cases of
actual interest to us. Of course, in case $W$ is strongly saturated
and $u$ commutes with finite products, the widest $h$ is the one whose
weak homotopy intervals are triples $(I,\delta_0,\delta_1)$ satisfying
$u(\delta_0)=u(\delta_1)$, which we could have used instead of just
$W$-aspheric intervals, which are\pspage{134} also the intervals such
that $I\to e$ is in $W$ (and hence universally so in terms of the
assumption made of compatibility of $u$ with products).  The trouble
with working with this $h$, rather than with $h_W$ as above, is
twofold though: a)\enspace The assumptions of strong saturation on $W$
and compatibility with products are not so readily verified in the
cases of interest to us, and the second moreover is not always
satisfied, e.g., there are test categories which are not strict, i.e.,
elementary modelizers which are not strict; b)\enspace the condition
$u(\delta_0)=u(\delta_1)$ is not readily verified in terms of $W$
directly, whereas the condition of $W$-asphericity is -- still more so
if $W$ is compatible with products and hence $I\to e$ is $W$-aspheric
just means it is in $W$.

\hangsection{Maps between homotopy structures.}\label{sec:55}%
Let's now look at ``\emph{morphisms}'' between categories endowed with
homotopy structures, $(M_1,h_1)$ and $(M_2,h_2)$ say. The natural
definition here is to take as morphisms between these homotopy
structures the functors $u:M_1\to M_2$ that give rise to a commutative
square of functors
\[\begin{tikzcd}[baseline=(O.base)]
  M_1 \ar[r,"u"]\ar[d] & M_2\ar[d] \\
  \overline M_1 \ar[r,"\overline u"] & |[alias=O]| \overline M_2
\end{tikzcd},\]
where the vertical arrows are the canonical functors into the
respective homotopy-types categories, and $\overline u$ a suitable
functor, necessarily unique.
The existence of $\overline u$ can be expressed at will in terms of
the $\Hom_1$ or $\Hom_2$ structures, namely as
\begin{enumerate}[label=(\roman*)]
\item\label{it:55.i}
  $f \underset{h_1}\sim g$ implies $u(f) \underset{h_2}\sim u(g)$,
\end{enumerate}
or as
\begin{enumerate}[label=(\roman*),resume]
\item\label{it:55.ii}
  $f\in W_{h_1}$ implies $u(f)\in W_{h_2}$.
\end{enumerate}
These conditions, when $M_1$ and $M_2$ have final objects and these
are respected by $u$, imply that $u$ transforms weak homotopy
intervals into weak homotopy intervals, and homotopisms into
homotopisms. Conversely, \emph{if $u$ commutes with finite products},
and if $h_1$ can be defined by a weak homotopy interval structure
(respectively, by a contractibility structure), then for $u$ to be a
morphism of homotopy structures, it is (necessary and) sufficient that
$u$ carry weak homotopy intervals (resp.\ contractible objects) into
same.

In case $h_1$ is described in terms of a generating set
$\Sigma_{h_1}^0$ of weak homotopy intervals, and if \emph{$u$ commutes
  with finite products}, the most economic way often to express that
$u$ is a morphism of homotopy structures, is by the condition that for
any \bI{} in $\Sigma_{h_1}^0$, $u(\bI)$ be a weak homotopy\pspage{135}
interval in $M_2$, namely $u(\delta_0) \underset{h_2}\sim
u(\delta_1)$. If we assume moreover that the intervals \bI{} in
$\Sigma_{h_1}^0$ are contractible, the previous condition is
equivalent to $u(I)$ being a contractible object of $M_2$ for any
\bI{} in $\Sigma_{h_1}^0$. The case I am mainly thinking of, of
course, is the one when $h_1$ can be described by a single generating
weak homotopy interval, possibly even contractible, or even by a
(generating) contractor. In the latter case, because of commutation of
$u$ with finite products, $u(\bI)$ will be equally a contractor -- and
contractible for $h_2$ if{f} $u$ is a morphism of homotopy structures.

In the precedent paragraphs, I forgot to mention the reduction of the
corresponding compatibility conditions (of a functor $u:M\to N$, or of
a saturated $W\subset\Fl(M)$) with a homotopy structure $h$, when the
latter is defined in terms of just a generating set $\Sigma_h^0$ of
weak homotopy intervals, possibly reduced to a single one, and
moreover \emph{$u$ or $W$ is ``compatible with finite products''}. In
the first case, it is enough to check $u(\delta_0)=u(\delta_1)$ for
any \bI{} in $\Sigma_h^0$ -- and if \bI{} is contractible, this
amounts to demanding $u(\bI)$ is a final object of $N$. In the second
case it is enough for strict compatibility, i.e., $W_h\subset W$, to
check that for any \bI{} in $\Sigma_h^0$, $I\to e$ is in $W$ (this
condition is also necessary, if any \bI{} in $\Sigma_h^0$ is
contractible). Even if $W$ is not supposed compatible with finite
products, namely $M\to W^{-1}M$ does not commute with finite products,
it is still sufficient for $W_h\subset W$ that all \bI's in
$\Sigma_h^0$ be \emph{$W$-aspheric} (as stated in
\hyperref[prop:54]{prop}.\ (p.\ \ref{p:131})), and this is necessary
too, if the \bI's are contractible.

\bigbreak

\presectionfill\ondate{4.4.}\par

\hangsection[Another glimpse upon canonical modelizers. Provisional
\dots]{Another glimpse upon canonical modelizers. Provisional working
  plan -- and recollection of some questions.}\label{sec:56}%
This seemingly endless review of generalities on homotopy notions is
getting a little fastidious - and still I am not quite through yet I
feel. One main motivation for embarking on this review was one strong
impression which grew out of the reflections of now just one week ago
(paragraph \ref{sec:48}), namely that the interesting ``test
functors'' from a test category $A$ into a modelizer $(M,W)$ are those
which factor through the full subcategory $M\subc$ of
\emph{contractible objects} of $M$. The presumable extra condition to
put on a functor $A\to M\subc$ to correspond to an actual test functor
$i$ from $A$ to $M$ are strikingly weak, such as $i^*(\bI)$ should be
aspheric over $e_\Ahat$ under the assumption we got a contractible
generating homotopy interval $\bI$ in $M$. In any case, if $M\subc^0$
is any full subcategory of $M\subc$ which gives rise (by taking
intervals in $M\subc^0$) to a family of homotopy intervals which
\emph{generates} the homotopy structure $h_W$ on $M$ associated to
$W$, it should\pspage{136} suffice, if $M$ is a ``canonical
modelizer'', that the asphericity condition on $i^*(\bI)$ should be
verified for any \bI{} in $M\subc^0$; this will presumably turn out in
due course as part of the definition (still ahead) of a ``canonical''
modelizer. Now, the most evident way to meet this condition, is to
take for $i$ a \emph{fully faithful} functor \emph{whose image
  contains} $M\subc^0$, or what amounts essentially to the same, any
full subcategory of $M$ containing $M\subc^0$! As we would like though
$A$ to be ``small'', this will be feasible only if $M\subc^0$ is small
-- hence the significance of the condition of a small generating
family of weak homotopy intervals for $h_W$ (which will imply that we
can find such a family with \emph{contractible} intervals, provided
only the homotopy structure $h_W$ can be described by a
contractibility structure, as we did indeed assume). Just as the
homotopy structure $h_W$ of $M$ was defined in terms of $W$,
conversely the notion of weak equivalences $W$ should be recoverable
in terms of the homotopy structure, and more specifically in terms of
the subcategory $M\subc$ and the small ``generating'' subcategories
$M\subc^0$ of $M\subc$, which we may now as well denote by $A$, by
taking the inclusion functor $i:A\to M$ of such an $A$, hence a
functor
\[i^*:M\to\Ahat,\]
and taking
\[W=(i^*)^{-1}(W_\Ahat).\]
Of course, we'll have still to check under which general conditions
upon a pair $(M,W)$ of a category and a saturated set of arrows $W$,
or rather, upon a pair $(M,M\subc)$ of a category endowed with a
contractibility structure $M\subc=C_h$ (where we think of $h$ as an
$h_W$), is it true that the saturated set of arrows
\[W(A) = (i^*)^{-1}(W_\Ahat)\subset\Fl(M)\]
in $M$ does not depend upon the choice of the full small
homotopy-generating subcategory $A$ of $M\subc$ (if restrictive
conditions are needed indeed). It may be reasonable to play safe, to
restrict at first to subcategories $A$ which are stable under finite
products in $M$, which will ensure that $A$ is a \emph{strict} test
category, i.e., \Ahat{} is a \emph{strict} elementary modelizer,
namely \Ahat{} is totally aspheric. But such restriction -- as well as
to test functors which are fully faithful -- should be a provisional
one, as ultimately we want of course to be able to use test categories
such as $\Simplex$ for ``testing'' rather general (canonical)
modelizers, whereas $\Simplex$ is by no means stable under products, nor
embeddable \emph{faithfully} in modelizers such as \Spaces{} say.

This\pspage{137} expectation of $W$ to be recoverable in terms of the
corresponding homotopy structure $h=h_W$ on $M$ takes its full
meaning, when joined with another one, namely that the latter can be
canonically described in terms of the category structure of $M$ and
the corresponding notion of $0$-connectedness. This latter expectation
is extremely strongly grounded, and I'll come back to it
circumstantially very soon I think (I started on it two days ago, but
then it got too late to take it to the end, and yesterday was spent on
some formal digressions\ldots). The two ``expectations'' put together,
when realized by carefully cutting out the suitable notions, should
imply that the structure of any ``canonical modelizer'' is indeed
determined ``canonically'' in terms of its category structure alone.

To come back to the relationship between test categories and
categories of the type $M\subc^0$, the idea which has been lurking
lately is that possibly, test categories can be viewed as no more, no
less, as categories endowed with a homotopy structure (necessarily
unique) which is a contractibility structure, and for which \emph{all
  objects are contractible}. At any rate, there must be a very close
relationship between the two notions, which I surely want to
understand. But as it would be quite unreasonable to restrict the
notion of a test category (and its weak and strong variants) to
categories admitting finite products, this shows that for a
satisfactory understanding of the above relationship, we should be
able to work with contractibility structures, and presumably too with
weak homotopy interval structures, in categories $A$ where we do not
make the assumption of stability under products. Thus in the outline
of the last few days, I still wasn't general enough it would seem!
This situation reminds me rather strongly of the early stages when
developing the language of sites, and restricting to sites where
fiber-products exist -- this seems by then a very weak and natural
assumption indeed, before it appeared (first to Giraud, I believe)
that it was quite an awkward and artificial restriction indeed, which
had to be overcome in order to work really at ease\ldots

All this now gives a lot of interesting things to look up in the short
run! I'll make a provisional plan of work, as follows:
\begin{enumerate}[label=\alph*)]
\item\label{it:56.a}
  Relation between a homotopy structure and the $\pi_0$ functor, and
  description of the so-called canonical homotopy structures.
\item\label{it:56.b} Write down in the end the ``key result'' on test
  functors $A\to\Cat$ which is overripe since the reflections of four
  days ago (par.\ \ref{sec:47}).  Presumably,\pspage{138} this will
  yield at the same time an axiomatic characterization of $W_\Cat$,
  namely of the notion of weak equivalence for functors between
  categories.
\item\label{it:56.c}
  At this point, we could go on and try and carry through the similar
  characterization for test functors $A\to\Bhat$, where $A$ and $B$
  are both test categories. There are also some generalities to
  develop about ``morphisms'' between test categories, which is ripe
  too for quite a while and cannot be pushed off indefinitely -- here
  would be the right moment surely. If the expected ``key result'' for
  test functors $A\to\Bhat$ carries through nicely it could presumably
  be applied at once in order to study general test functors $A\to M$,
  and thus get the clues for cutting out ``the'' natural notion of a
  canonical modelizer, which ``was in the air'' since the ``naive
  question'' of par.\ \ref{sec:46} (page \ref{p:95}).
\item\label{it:56.d}
  However, there is another approach to canonical modelizers which is
  just appearing, via the idea (described above) of associating
  canonically a notion of ``weak equivalence'' $W$ to a homotopy
  structure of type $\Hom_4$, i.e., to a contractibility structure,
  subject possibly to some restrictions. This ties in, as explained
  above, with a closer look at the relationship between test
  categories, and ``coarse'' contractibility structures (where all
  objects are contractible).
\end{enumerate}

It would seem unreasonable to push off \ref{it:56.a} and \ref{it:56.b}
any longer now -- so I'll begin with these. I am hesitant however
between \ref{it:56.c} and \ref{it:56.d} -- with a feeling that the
later approach \ref{it:56.d} may well turn out to be technically the
most expedient one. Both have to be carried through anyhow, and the
two together should give a rather accurate picture of what canonical
modelizers are about.

On the other hand, there are still quite a bunch of questions which
have been waiting for investigation -- for instance the list of
questions of nearly three weeks ago (page \ref{p:42}). Among the six
questions stated there, four have been settled, or are about to be
settled through the previous program (if it works out), questions
\ref{it:32.4} and \ref{it:32.6} remain, the first one being about
\Ahat{} being a closed model category, and about the homotopy
structure of \Cat. There are a number of more technical questions too,
for instance I did not finish yet my review of the ``standard'' test
categories and never wrote down the proof that \Simplexf{} (simplices
with face operations and no degeneracies) is indeed a weak test
category. But for the time being, all these questions appear as
somewhat marginal with respect to the strong focus the reflection has
been gradually taking nearly since the very\pspage{139} beginning --
namely an investigation of modelizers and, more specifically, the
gradual unraveling of a notion of ``canonical modelizer''. I certainly
feel like carrying this to the end at once, without any digressions
except when felt relevant for the main focus at present. As for
choosing precedence between \ref{it:56.c} and \ref{it:56.d}, it is
still time to decide, when we're through with \ref{it:56.a} and \ref{it:56.b}!

\hangsection[Relation of homotopy structures to $0$-connectedness and
\dots]{Relation of homotopy structures to
  \texorpdfstring{$0$}{0}-connectedness and
  \texorpdfstring{$\pi_0$}{pi-0}. The canonical homotopy structure
  \texorpdfstring{$h_M$}{hM} of a category
  \texorpdfstring{$M$}{M}.}\label{sec:57}%
\textbf{Relation of homotopy structures to $0$-connectedness and to
  \piz.}\enspace Here we're resuming the reflection started in par.\
\ref{sec:51} \ref{subsec:51.F} (page \ref{p:122}). All we did there
was to introduce some conditions on a category $M$, namely
\ref{it:51.F.a} to \ref{it:51.F.e} (page \ref{p:123}), and introduce
the functor
\[\pi_0 : M \to \Sets.\]
As before, we'll assume now $M$ satisfies conditions \ref{it:51.F.a}
to \ref{it:51.F.d} (the first three, or rather \ref{it:51.F.b} and
\ref{it:51.F.c}, are all that is needed for defining the functor
$\pi_0$), and will not assume \ref{it:51.F.e}, or the equivalent
\ref{it:51.F.eprime} of $\pi_0$ commuting to finite products, unless
explicitly specified.

We suppose now, moreover, $M$ endowed with a homotopy structure
$h$. We'll say that $h$ is \emph{$\pi_0$-admissible}, or simply
\emph{$0$-admissible}, if $h$ is compatible with the functor $\pi_0$,
which can be expressed by either one of the following two equivalent
conditions (cf.\ page \ref{p:130}):
\begin{enumerate}[label=(\roman*)]
\item\label{it:57.i}
  $f\hsim g$ implies $\pi_0(f)=\pi_0(g)$,
\end{enumerate}
for any two maps $f,g$ in $M$, or
\begin{enumerate}[label=(\roman*),resume]
\item\label{it:57.ii}
  $f\in W_h$ (i.e., $f$ a homotopism) implies $\pi_0(f)$ bijective.
\end{enumerate}

Another equivalent formulation is that the functor $\pi_0$ factors
through the quotient category $M_h$ of homotopy types
\[M\to M_h\to\Sets,\]
we'll still denote by $\pi_0$ the functor $M_h\to\Sets$ obtained.

If $h$ is $0$-admissible, it satisfies \ref{it:57.iprime} and
\ref{it:57.iiprime} below:
\begin{enumerate}[label=(\roman*')]
\item\label{it:57.iprime}
  For any weak homotopy interval $(I,\delta_0,\delta_1)\in\Sigma_h$, $\pi_0(\delta_0)=\pi_0(\delta_1)$.
\end{enumerate}
As $\pi_0(e_M)$ is a one-point set by the assumption \ref{it:51.F.d}
on $M$, for any section $\delta$ of an object $I$, $\pi_0(\delta)$ may
be described as just an element of $\pi_0(I)$, which is the unique
connected component of $I$ through which factors the given section;
thus \ref{it:57.iprime} can be expressed by saying that for any weak
homotopy interval, the two ``endpoints'' belong to the same connected
component of $I$, a natural condition indeed! It is automatically
satisfied\pspage{140} if $I$ is connected. In fact, \ref{it:57.iprime}
is satisfied if{f} $\Sigma_h$ admits a generating subset $\Sigma_h^0$
made up with \emph{connected} intervals.
\begin{enumerate}[label=(\roman*'),resume]
\item\label{it:57.iiprime}
  Any contractible object $X$ is connected (and hence $0$-connected).
\end{enumerate}

Conversely (cf.\ page \ref{p:130}), if condition \ref{it:51.F.e}
holds, i.e., $\pi_0$ commutes to products, and if moreover $h$ is
definable in terms of $\Sigma_h$, i.e., comes from a weak homotopy
interval structure, then \ref{it:57.iprime} implies
$0$-admissibility. If $M$ admits even a generating family of
contractible weak homotopy intervals, namely if $h$ comes from a
contractibility structure on $M$, then \ref{it:57.iiprime} equally
implies $0$-admissibility.
\begin{remarks}
  1)\enspace We can generalize these converse statements, by dropping
  the condition \ref{it:51.F.e} on $M$, but demanding instead that the
  connected component involved (namely $I$ itself in case
  \ref{it:57.iiprime}) is not only $0$-connected, but even
  ``\emph{$0$-connected over $e_M$}'', which just means that its
  products by any $0$-connected object of $M$ is again
  $0$-connected. In the case \ref{it:57.iiprime}, this condition es
  equally \emph{necessary} for admissibility. (The corresponding
  statements could have been made in the general context of page
  \ref{p:130} of course\ldots)

  2)\enspace The name of $0$-admissibility suggests there may exist
  correspondingly ``higher'' notions of $n$-admissibility for $h$,
  where $n$ is any natural integer. I \emph{do} see a natural
  candidate, namely whenever we have a functor $\pi_n$ from $M$ to
  $n$-truncated homotopy types (as is the case, say, when $M$ is
  either a topos -- we then rather get \emph{pro}homotopy types -- or
  a modelizer). But it would seem that in all cases of geometrical
  significance, and when moreover $h$ is defined by a weak homotopy
  interval structure, that $0$-admissibility implies already
  $n$-admissibility for any $n$.
\end{remarks}

Due to the existence of $\sup$ for an arbitrary subset, in each of the
ordered sets $\Hom_i(M)$ ($i\in\{0,3,4\}$) of homotopy structures in
$M$ (either unqualified, or weak homotopy interval structures, or
contractibility structures), it follows that for each of these three
types of homotopy structures, there is a widest one $h_i$ among all
those of this type which are $0$-admissible. We are interested here,
because of topological motivations, by the case of $\Hom_3$, namely
weak homotopy interval structures. We call the corresponding homotopy
structure $h=h_3$ the \emph{canonical homotopy structure} of $M$. In
case $M$ is totally $0$-connected (by which we mean condition
\ref{it:51.F.e}), the weak homotopy intervals for this structure are
just those intervals for which $\delta_0,\delta_1$ correspond to the
same connected component of $I$ -- and we get a generating set of weak
homotopy intervals, by just taking \emph{all connected intervals}.
\begin{remark}
  If $M$ is not\pspage{141} totally $0$-connected, we still get a
  description of $\Sigma_h$, as those intervals such that for any
  $0$-connected $X$, the corresponding sections of $X\times I$ over
  $X$ correspond to the same connected component of $X\times I$. It is
  enough for this that the common connected component $I_0$ of $I$ for
  $\delta_0,\delta_1$ should be $0$-connected over $e$ -- it is not
  clear to me whether this condition is equally necessary. Anyhow,
  presumably the case $M$ totally $0$-connected will be enough for all
  we'll have to do.
\end{remark}

In case $M$ is totally $0$-connected, the homotopy notions in $M$ for
the canonical homotopy structure are just those which can be described
in terms of ``homotopies'' using connected intervals -- which is
intuitively the first thing that comes to mind indeed, when trying to
mimic most naively, in an abstract categorical context, the familiar
homotopy notions for topological spaces.

Let let
\[ W\subset\Fl(M)\]
be any saturated set of arrows in $M$ (viewed as a notion of ``weak
equivalence'' in $M$). Consider the corresponding homotopy structure
$h_W$ on $M$, defined in terms of $W$-aspheric intervals as a
generating family of weak homotopy intervals. Let $h_M$ be the
canonical homotopy structure on $M$. Thus the condition
\[ h_W \le h_M\]
just means that $h_W$ is $0$-admissible, or equivalently, that
$W$-aspheric objects over $e$ which have a section, are
$0$-connected (for simplicity, I assume from now on $M$ totally
$0$-connected). This looks like a very reasonable condition indeed, if
$W$ should correspond at all to the intuitions associated to the
notion of ``weak equivalence''! As a matter of fact, this condition is
clearly implied by the condition that $W$ itself should be
``\emph{$0$-admissible}'', by which we mean that the functor $\pi_0$
is compatible with $W$, i.e., transforms weak equivalences into
bijections, or equivalently, factors through $M\to W^{-1}M$.

What we are looking for however is conditions on $W$ for the
\emph{equality}
\[ h_W=h_M\]
to hold. When the previous condition (expressing $h_W\le h_M$) is
satisfied, all that remain is to express the opposite inequality,
which is done in the standard way. We thus get:
\begin{proposition}
  Let\pspage{142} $M$ be a category, assume $M$ totally $0$-connected
  \textup(i.e., satisfying conditions
  \textup{\ref{it:51.F.a}\ref{it:51.F.b}\ref{it:51.F.c}\ref{it:51.F.e}}
  of page \ref{p:123}\textup). Let $W\subset\Fl(M)$ be a saturated set
  of arrows in $M$, consider the associated homotopy structure $h_W$
  on $M$, with $W$-aspheric intervals as a generating family of weak
  homotopy intervals. Consider also the canonical homotopy structure
  $h_M$ on $M$, with $0$-connected intervals as a generating family of
  weak homotopy intervals \textup(thus $h_M=h_{W_0}$, where
  $W_0\subset\Fl(M)$ consists of all arrows made invertible by the
  functor $\pi_0:M\to\Sets$\textup). In order for the equality
  $h_W=h_M$ to hold, it is necessary and sufficient that the following
  two conditions be satisfied:
  \begin{enumerate}[label=\alph*),font=\normalfont]
  \item\label{it:57.a}
    Any object $I$ of $M$ which is $W$-aspheric over $e$ and admits a
    section, is connected \textup(it is enough for this that $W$ by
    $0$-admissible, i.e., $f\in W$ imply $\pi_0(f)$ bijective\textup);
  \item\label{it:57.b}
    For any connected object $I$ of $M$ and any two sections
    $\delta_0,\delta_1$ of $I$, these can be ``joined'' by a finite
    chain of sections $s_i$ \textup($0\le i\le n$\textup),
    $s_0=\delta_0$, $s_n=\delta_1$, such that for any two consecutive
    ones, there exists an object $J$, $W$-aspheric over $e$, a map
    $J\to I$ and two sections of $J$ mapped into the sections
    $s_i,s_{i+1}$ of $I$.
  \end{enumerate}
\end{proposition}

The condition \ref{it:57.a} says there are not too many weak
equivalences, whereas \ref{it:57.b} says there are still enough for
``testing'' connectedness in terms of $W$-aspheric intervals. Both
conditions look plausible enough!

Next step one would think of, in this context, is to give conditions
on $W$ (independently of the previous ones) which will allow to
express $W$ in terms of $h_W$. But for this, I should develop first a
description of a notion of ``weak equivalence'' in terms of an
arbitrary homotopy structure $h$ on $M$, as contemplated in the
previous paragraph. I decided however to give precedence to the ``key
result'' still ahead.

Maybe the condition $h_W=h_M$, expressed in the previous proposition,
merits a name -- we'll say that the notion of weak equivalence $W$ is
``\emph{geometric}'', if it satisfies the two conditions \ref{it:57.a}
and \ref{it:57.b} above, or rather the slightly stronger
\namedlabel{it:57.aprime}{a')} in place of a) -- namely
$0$-admissibility of $W$ (plus \ref{it:57.b} of course). The
conditions \ref{it:57.aprime} and \ref{it:57.b} are the explicit ones
for checking -- but for using that $W$ is geometric, the more
conceptual statement $h_W=h_M$ (besides $0$-admissibility) is the
best. Thus, as any \emph{contractor} in $M$ which is connected (hence
a weak homotopy interval for $h_M$) is $h_M$-contractible, it is
$h_W$-contractible, and hence $W$-aspheric over $e$, and conversely of
course. Using this, we get a converse to part \ref{it:51.E.c}
and\pspage{143} \ref{it:51.E.d} of the proposition of page
\ref{p:121}, in the present context, which we may state in a more
complete form as follows:
\begin{proposition}
  Let $M$ be a totally $0$-connected category, $W$ a geometric
  saturated set of arrows in $M$. Assume moreover that for any two
  objects in $M$, the object $\bHom(X,Y)$ in $M$ exists. Let $X$ be an
  object of $M$, then the following conditions \textup{\ref{cond:57.a}} to
  \textup{\ref{cond:57.cdblprime}} are equivalent:
  \begin{description}
  \item[\namedlabel{cond:57.a}{a)}]
    For any object $Y$, $\bHom(Y,X)$ is $0$-connected.
  \item[\namedlabel{cond:57.aprime}{a')}]
    For any object $Y$, $\bHom(Y,X)$ is $W$-aspheric \textup(i.e., its
    map to $e$ is $\in W$\textup).
  \item[\namedlabel{cond:57.adblprime}{a'')}]
    For any object $Y$, $\bHom(Y,X)$ is $W$-aspheric over $e$
    \textup(i.e., its map to $e$ is ``universally in $W$''\textup).
  \item[\namedlabel{cond:57.atplprime}{a''')}]
    For any object $Y$, $\bHom(Y,X)$ is $h_W$-contractible.
  \item[\namedlabel{cond:57.b}{b)}\namedlabel{cond:57.bprime}{b')}%
    \namedlabel{cond:57.bdblprime}{b'')}\namedlabel{cond:57.btplprime}{b''')}]
    Same as above, with $Y$ replaced by $X$.
  \item[\namedlabel{cond:57.c}{c)}]
    For any object $Y$, $Y\to \bHom(X,Y)$ is in $W$.
  \item[\namedlabel{cond:57.cprime}{c')}]
    For any object $Y$, $Y\to \bHom(X,Y)$ is $W$-aspheric.
  \item[\namedlabel{cond:57.cdblprime}{c'')}]
    For any object $Y$, $Y$ as a subobject of $\bHom(X,Y)$ is a
    deformation retract with respect to $h_W$.
  \end{description}
  If moreover $X$ has a section, these conditions are equivalent to:
  \begin{description}
  \item[\namedlabel{cond:57.d}{d)}]
    $X$ is $h_M$-contractible.
  \item[\namedlabel{cond:57.dprime}{d')}]
    $X$ is $h_W$-contractible.
  \end{description}
\end{proposition}
\begin{example}
  Let $M=\Cat$, the assumptions on $M$ are clearly satisfied. If we
  take for $W$ the usual weak equivalences, it is clear too that $W$
  satisfies \ref{it:57.aprime} and \ref{it:57.b} above, i.e., $W$ is
  geometric. Thus the preceding proposition is just an elaboration of
  the result stated after the prop.\ of page \ref{p:97} -- which was
  the moment when it became clear that contractibility was an
  important notion in the context of test categories and test
  functors, and which was the main motivation too for the somewhat
  lengthy trip through generalities on homotopy notions, which is
  coming now (in the long last!) to a provisional end\ldots
\end{example}

To check if the notions developed in this section are handy indeed, I
would still like to try them out in the case $M$ is an elementary
modelizer \Ahat, corresponding to a test category $A$. But it's
getting late this night\ldots

\bigbreak
\presectionfill\ondate{5.4.}\pspage{144}\par

\hangsection[Case of totally $0$-connected category $M$. The category
\dots]{Case of totally \texorpdfstring{$0$}{0}-connected category
  \texorpdfstring{$M$}{M}. The category
  \texorpdfstring{$\overline\Cat$}{(Cat)} of
  \texorpdfstring{\textup(}{(}small\texorpdfstring{\textup)}{)}
  categories and homotopy classes of functors.}\label{sec:58}%
Some comments still on the canonical homotopy structure of a category
$M$, which we assume again totally $0$-connected. Two sections of an
object $K$ are homotopic if{f} they belong to the same connected
component of $X$, thus we get an \emph{injective} canonical map from
homotopy-classes of sections
\[\overline\Gamma(X)=\oHom(e,X) \to \pi_0(X).\]
We are interested in the case when this map is always a bijection, or
what amounts to the same, when $M$ satisfies the extra condition:
\begin{enumerate}[label=\alph*),start=6]
\item\label{it:58.f}
  Every $0$-connected objects of $M$ has a section.
\end{enumerate}
\begin{proposition}
  Let $M$ be a category satisfying conditions \textup{\ref{it:51.F.a}} to
  \textup{\ref{it:51.F.d}} of page \ref{p:123}. If $M$ satisfies
  moreover condition \textup{\ref{it:58.f}} above, it satisfies also
  condition \textup{ref{it:51.F.e}}, i.e., $M$ is totally $0$-connected.
\end{proposition}

Indeed, let $X,Y$ be $0$-connected objects, we must prove $X\times Y$
is $0$-connected. First, it is ``non-empty'', because it got a section
(as $X$ and $Y$ each have a section). All that remains to do is show
that two connected components of $X\times Y$ are equal. By assumption
\ref{it:58.f}, each has a section, say $(s_i,t_i)$ with
$i\in\{1,2\}$. These two sections can be joined by a two-step chain
\[ (s_1,t_1), \quad (s_2,t_1), \quad (s_2,t_2)\]
the first step contained in $X\times t_1$ which is connected, the
second in $s_2\times Y$ which is connected, hence the two sections
belong to the same connected component, qed.

The condition \ref{it:58.f} is \emph{strictly} stronger than
\ref{it:51.F.e} (when \ref{it:51.F.a} to \ref{it:51.F.d} are
satisfied), as we see by taking for $M$ the category of all sheaves on
an irreducible topological space $X$ -- thus $M$ is a totally aspheric
topos, a lot better it would seem than just totally $0$-connected, but
condition \ref{it:58.f} is satisfied if{f} the topos is equivalent to
the final topos or ``one-point topos'', i.e., if{f} the topology of
$X$ is the chaotic one ($X$ and $\emptyset$ are the only open
subsets). We'll say $M$ is \emph{strictly totally $0$-connected}, if
it satisfies the conditions \ref{it:51.F.a} to \ref{it:58.f} (where
\ref{it:51.F.e} is a consequence of the others). Thus, if $M$ is
strictly totally $0$-disconnected, we get for any object $X$ a
canonical bijection, functorial in $X$
\begin{equation}
  \label{eq:58.1}
  \overline\Gamma(X) \eqdef \;\begin{tabular}{@{}c@{}}homotopy classes\\of sections of
  $X$\end{tabular}\; \tosim \pi_0(X).
  \tag{1}
\end{equation}

Assume now, moreover, that $X,Y$ are two objects of $M$ such that the
object $\bHom(X,Y)$ exists. Then, using the observation of
\hyperref[cor:51.E.2]{cor.\ 2} (p.\ \ref{p:122}),\pspage{145} we get
the familiar relationship
\begin{equation}
  \label{eq:58.2}
  \oHom(X,Y) \tosim \pi_0(\bHom(X,Y)),
  \tag{2}
\end{equation}
where $\oHom$ denotes homotopy classes of maps, with respect to the
canonical homotopy structure of $M$.

The first example which I have in mind is the case $M=\Cat$. I was a
little short yesterday about $W_\Cat$ being ``geometric'' -- the
condition \ref{it:57.aprime} of page \ref{p:142} is evident indeed in
terms of the (geometric!) definition of weak equivalence in terms of
non-commutative cohomology (where we need to case only about zero
dimensional cohomology!). But condition \ref{it:57.b} is a consequence
of the fact that $0$-connectedness in \Cat{} can be checked using only
$\Simplex_1$ (i.e., $\Simplex_1$ is a generating contractor for the
canonical homotopy structure of \Cat), and that $\Simplex_1$ is
$W_\Cat$-aspheric over the final category $e$, i.e., that for any $C$
in \Cat, the projection
\[ C\times\Simplex_1\to C\]
is in $W_\Cat$. I want to start being attentive from now on about what
exactly are the formal properties of $W_\Cat$ we are using -- it
really seems they boil down to very few, which we have kept using
without ever having to refer to the ``meaning'' of weak equivalence in
terms of cohomology (except for proving the formal properties we
needed, in terms of the precise definition of weak equivalence we have
been starting with from the very beginning).

In terms of the canonical homotopy structure in \Cat, admitting
$\Simplex_1$ as a generating homotopy interval (a contractor, as a
matter of fact), we now get a notion of two functors from a category
$X$ to a category $Y$ being homotopic, and a corresponding notion of
homotopy classes of functors from $X$ to $Y$, which are in one-to-one
correspondence with the connected components of the category
$\bHom(X,Y)$ of all functors from $X$ to $Y$. An elementary homotopy
from a functor $f$ to another $g$ (with respect of course to the basic
generating interval $\Simplex_1$, which is always understood here) is
nothing else but a morphism from $f$ into $g$. Thus, homotopy classes
of functors are nothing but equivalence classes for the equivalence
relation in $\Hom(X,Y)$ generated by the relation ``there exists a
morphism from $f$ to $g$''. To consider at all this equivalence
relation, rather than the usual one of isomorphy between functors, is
rather far from the spirit in which categories are generally being
used -- as is also the very notion of a functor which is a weak
equivalence, which has been our starting point. Topological motivation
alone, it seems, would induce anyone to introduce such barbaric
looking notions into category theory!

From\pspage{146} our point of view, the main point for paying
attention to the homotopy relation for functors, is of course because
homotopic functors define the same map in the localized category \Hot,
which category for the time being is (together with the modelizers
designed for describing it) our main focus of attention. According to
the general scheme of homotopy theory as reviewed previously, homotopy
classes of functors give rise to a quotient category of \Cat, which is
at the same time a localization of \Cat{} with respect to homotopisms,
and which we'll denote by $\oCat$. Thus we get a factorization of the
canonical functor from \Cat{} to \Hot{}
\[ \Cat \to \oCat \to \Hot.\]
This is just the factorization of the localization functor with
respect to $W_\Cat=W$, through the ``\emph{partial localization}''
with respect to the smaller set $W_h\subset W$ consisting of
homotopisms only. All we have used about $W=W_\Cat$ (from which the
localization $\Hot=W^{-1}\Cat$ is deduced) for getting this
factorization, was that $f\in W$ implies $\pi_0(f)$
bijective,\footnote{\alsoondate{6.4.} As a matter of fact, the condition about
  $\pi_0$ visibly isn't needed!} and that the projections
$C\times\Simplex_1\to C$ are in $W$. (Of course, we assume tacitly that
$W$ is saturated too.) We may view \Hot{} as a localization of \oCat{}
with respect to the set of arrows $\overline W$ corresponding to $W$
-- as a matter of fact, $W$ may be viewed as the inverse image of
$\overline W$ by the canonical functor from \Cat{} to its partial
localization \oCat{} -- the category of ``homotopy types'' relative to
the canonical homotopy structure of \Cat.

Well-known analogies would suggest at this point that we may well be
able to describe \Hot{} in terms of $(\oCat, \overline W)$ by a
\emph{calculus of fractions} -- right fractions presumably, or maybe
either right or left. This may possibly lead to a direct proof of the
notion of weak equivalence we have been working with being strongly
saturated, without having to rely upon Quillen's closed model
theory. But it is not yet the moment to pursue this line of thought,
which would take us off the main focus at present.

\hangsection[Case of the ``next best'' modelizer \Spaces{} -- and
\dots]{Case of the ``next best'' modelizer
  \texorpdfstring{\Spaces}{(Spaces)} -- and need of introducing the
  \texorpdfstring{$\pi_0$}{pi-0}-functor as an extra structure on a
  would-be modelizer \texorpdfstring{$M$}{M}.}\label{sec:59}%
Before pushing ahead, I would like to make still another point about
the work done yesterday -- a point suggested by looking at the case of
the modelizer \Spaces, which after all is the next best ``naive''
modelizer, less close to algebra than \Cat, but still worth being
taken into account! This category satisfies the conditions
\ref{it:51.F.a} to \ref{it:58.f}, \emph{except} the condition
\ref{it:51.F.c} -- which would mean that the connected components of a
space (as defined in terms of usual topology) are open subsets, which
is true (for a space and its open subsets) \emph{only} for locally
connected spaces. The point is that this doesn't (or shouldn't)
really\pspage{147} matter -- the way topological spaces are used as
``homotopy models'' in standard homotopy theory, it is \emph{pathwise
  connected} components that count, and not the topological ones. In
terms of these, there is still a canonical functor
\[\pi_0:\Spaces\to\Sets,\]
this functor however is no longer left adjoint of the functor in
opposite direction, associating to every set $E$, the corresponding
discrete topological space. (To get an adjunction, we should have to
restrict to the category of pathwise locally connected topological
spaces.) It doesn't matter visibly -- all that's being used is that
$\pi_0$ commutes to arbitrary sums, and takes pathwise $0$-connected
spaces into one-point sets.

This suggests that we should generalize the notions around the
``canonical'' homotopy structure on a category $M$, to the case of a
category which need not satisfy the exacting conditions of total
$0$-connectedness, by introducing as an \emph{extra structure} upon
$M$ a given functor
\[\pi_0:M\to\Sets,\]
subject possibly to suitable restrictions. The first which comes to
mind here is commutation with sums -- it doesn't seem though we've had
to use this property so far. All we've used occasionally was existence
of finite products in $M$, and commutation of $\pi_0$ to these.

If we think of $M$ as a would-be modelizer, and therefore endowed with
a hoped-for functor
\[M\to \Hot,\]
there is a natural functor $\pi_0$ indeed on $M$, namely the
composition
\[ M \to \Hot \to \Sets,\]
where the canonical functor
\[\pi_0:\Hot\to\Sets\]
is deduced from the $\pi_0$-functor $\Cat\to\Sets$ considered
previously, by factorization through the localized category \Hot{} of
\Cat. Thus, ``the least we would expect'' from a category $M$ for
being eligible as a modelizer is that there should be a natural
functor $\pi_0$ around, corresponding to the intuition of connected
components. In case of a ``canonical'' modelizer $M$ (maybe we should
say rather: canonical with respect to a given $\pi_0$), there is the
feeling that the functor $M\to\Hot$ we are after could eventually be
squeezed out from just $\pi_0$, and that it could be viewed as
something like a ``total left derived functor'' of the functor
$\pi_0$. But this for the time being is still thin air\ldots

What\pspage{148} we \emph{can} do however at present, in terms of a
given functor $\pi_0$, is to introduce the corresponding notion of
$0$-connectedness (understood: with respect to $\pi_0$), namely
objects $X$ such that $\pi_0(X)$ is a one-point set, the notion of
compatibility of a homotopy structure $h$ on $M$ with $\pi_0$, and the
$\pi_0$-canonical (or simply, ``canonical'') homotopy structure on
$M$, which now should be denoted by $h_\piz$ rather than $h_M$ (unless
we write $h_\bM$ where \bM{} denotes the pair $(M,\piz)$), which is
the widest weak homotopy interval structure on $M$ which is
\piz-admissible, and can be described (assuming \piz{} commutes with
finite products) in terms of all $0$-connected intervals as a
generating family of weak homotopy intervals. The generalities of
par.\ \ref{sec:54} about the relationship of $h_\bM$ with $h_W$ (where
$W\subset\Fl(M)$) should carry over verbatim, as well as those of the
next, provided everywhere $0$-connectedness is understood relative to
the given functor \piz, and ``total $0$-connectedness'' is interpreted
as just meaning that \piz{} commutes to finite products. Thus, our
contact with ``geometry'' via true honest connected components of
objects was of short duration, and back we are to pure algebra with
just a functor given which we call \piz{}, God knows why -- the culprit
for this change of perspective being poor modelizer \Spaces, which was
supposed to represent the tie with so-called ``topology''\ldots

\hangsection[Case of strictly totally aspheric topos. A timid start on
\dots]{Case of strictly totally aspheric topos. A timid start on
  axiomatizing the set \texorpdfstring{\scrW}{W} of weak equivalences
  in \texorpdfstring{\Cat}{(Cat)}.}\label{sec:60}%
I almost forgot I still have to check ``handiness'' of the notions
developed yesterday, on the example of test categories or rather, the
corresponding elementary modelizers \Ahat. As usual, I can't resist
being a little more general, so let's start with an arbitrary topos
\scrA{} first. It always satisfies conditions \ref{it:51.F.a} and
\ref{it:51.F.b} of page \ref{p:123}. Condition \ref{it:51.F.c}, namely
that every object of \scrA{} could be decomposed into a sum of
$0$-connected ones, is equivalent with saying that \scrA{} admits a
generating subcategory $A$ made up with $0$-connected objects. In this
case, \scrA{} is called \emph{locally $0$-connected} or simply,
\emph{locally connected} -- which generalizes the notion known under
this name from topological spaces to topoi. On the other hand,
condition \ref{it:51.F.d} is expressed by saying that the topos
considered is \emph{$0$-connected} -- equally a generalization of the
corresponding notion for spaces. Condition \ref{it:51.F.e}, about the
product of two $0$-connected objects being $0$-connected, is a highly
unusual one in ordinary topology. For a topological space, it means
that the space is irreducible (hence reduced to a point if the space
is Hausdorff). In accordance with the terminology introduced
yesterday, we'll say that \scrA{} is \emph{totally $0$-connected} if
it is locally connected, and if the product of two $0$-connected
objects is again $0$-connected. The standard arguments show that for
this, it is enough\pspage{149} that the product of two elements in $A$
be $0$-connected. The topos is called \emph{strictly totally
  $0$-connected} if it is locally connected, and if moreover every
$0$-connected object admits a section -- which (as we saw earlier
today) implies \scrA{} is totally $0$-connected, as the wording
suggests. It amounts to the same to demand that every ``non-empty''
object have a section -- and for this it is enough that the elements
in $A$ have a section. This latter condition is trivially checked for
all standard test categories I've met so far (they all have a final
object, and there maps of the latter into any other object of $A$). A
noticeable counterexample here is $(\Simplexf)\uphat$ (semisimplicial
face complexes, \emph{without degeneracies}), where the \emph{weak}
test category \Simplexf{} hasn't got a final object ($\Simplex^0$
definitely isn't!) and no $\Simplex_n$ in \Simplexf{} except
$\Simplex_0$ only has got a section.
\begin{remark}
  I wonder, when \scrA{} is totally $0$-connected, and moreover
  modelizing, i.e., the Lawvere element is aspheric over $e$, if this
  implies \scrA{} is totally aspheric, and that every element which is
  ``non-empty'' has a section (i.e., \emph{strict} total $0$-connectedness).
\end{remark}

Next thing is to look at
\[W\subset\Fl(\scrA),\]
the set of weak equivalences (as defined by non-commutative cohomology
of topoi), and see if it is ``geometric'' (page \ref{p:142}). Condition
\ref{it:57.aprime} is clearly satisfied, there remains the condition
\ref{it:57.b}, namely whether $0$-connectedness of an object of
\scrA{} (\scrA{} supposed totally $0$-connected) can be tested, using
``intervals'' which are \emph{aspheric over $e$}. More specifically,
we want to test that two sections of $I$ belong to the same connected
component, using for ``joining'' them intervals that are aspheric over
$e$. The natural idea here is to assume the generating objects in $A$
to be aspheric over $e$ (which implies \scrA{} is \emph{totally
  aspheric}, not only totally $0$-connected), and to use these objects
(endowed with suitable sections) as testing intervals. This goes
through smoothly, indeed, if we assume moreover strict total
zero-connectedness. Thus:
\begin{proposition}
  Let \scrA{} be a topos which is strictly totally aspheric
  \textup(namely totally aspheric, and every ``non-empty'' object has
  a section\textup). Then the set $W\subset\Fl(\scrA)$ of weak
  equivalences in \scrA{} is ``geometric'', and accordingly, the
  homotopy structure $h_W$ defined in terms of aspheric homotopy
  intervals, is the same as the canonical homotopy structure $h_\scrA$
  defined in terms of merely $0$-connected homotopy
  intervals. Moreover, for any set $A\subset\Ob\scrA$ which is
  generating and whose objects are $0$-connected, the set of
  $0$-connected intervals $\bI=(I,\delta_0,\delta_1)$ with $I$ in $A$\kern1pt,
  generate the homotopy structure $h_W$.
\end{proposition}

A\pspage{150} topos \scrA{} as in the proposition (namely strictly
totally aspheric) need not be a modelizer, i.e., the Lawvere element
$L$ need not be aspheric, or what amounts to the same because of
$h_W=h_\scrA$ and $L$ being a contractor, $L$ need not be connected:
take $\scrA=\Sets$! I suspect though this to be the only
counterexample (up to equivalence). For \scrA{} to be a modelizer, we
need only find an object in $A$ which has got two \emph{distinct}
sections (because then they must be disjoint, i.e., $e_0\sand
e_1=\varnothing_\scrA$, because $e$ has only the full and the ``empty''
subobject, as a consequence of every ``non-empty'' object of \scrA{}
having a section), thus getting a ``homotopy interval'' (more
specifically, a separated and relatively aspheric one) as requested
for \scrA{} to be a modelizer. Now for any would-be test category met
with so far (except precisely \Simplexf{} and the like, which are
\emph{not} test categories but only weak ones), this condition that
there are objects in $A$ which have more than just one ``point'' ($=$
section), is trivially verified.

In case of a topos of the type \Ahat, the notion of weak equivalence
in \Ahat{} can be described (independently of cohomological notions)
in terms of the notion of weak equivalence in \Cat, more precisely
\[W_A \eqdef W_\Ahat = i_A^{-1}(\scrW),\]
where
\[i_A:\Ahat\to\Cat, \quad F\mapsto A_{/F}\]
is the canonical functor, and where
\[\scrW=W_\Cat\]
is the set of weak equivalences in \Cat. These of course, for the time
being, are defined in terms of cohomology (including a bit of
non-commutative one in dimension $1$\ldots). We may however start with
any $\scrW\subset\Fl(\Cat)$ and look at which formal properties on
\scrW{} (satisfied for usual weak equivalences) allow our arguments to
go through, in various circumstances. We may make a list of those
which have been used today, and go on this way a little longer, with
the expectation we'll finally wind up with an axiomatic
characterization of weak equivalences, i.e., of \scrW, in terms of the
category \Cat, say.
\begin{enumerate}[label=\alph*)]
\item\label{it:60.a}
  (Pour m\'emoire!) \scrW{} is saturated (cf.\ page \ref{p:101}).
\item\label{it:60.b}
  \scrW{} is $0$-admissible, i.e., if $f:C\to C'$ is in \scrW,
  $\piz(f)$ is bijective.
\item\label{it:60.c}
  $\Simplex_1$ is \scrW-aspheric over $e=\Simplex_0$, i.e., for every $C$
  in \Cat, the projection $C\times\Simplex_1\to C$ is in \scrW.
\item\label{it:60.d}
  Any $C$ in \Cat{} which has a final element is \scrW-aspheric, i.e,
  $C\to e$ is in \scrW.
\end{enumerate}

The\pspage{151} condition \ref{it:60.a} will be tacitly understood
throughout, when taking a \scrW{} to replace usual weak
equivalences. Conditions \ref{it:60.b} and \ref{it:60.c} then were
seen to be enough to imply that $h_\scrW=h_\Cat$. On the other hand,
one sees at once that for the proposition over for a topos \Ahat{}
which is strictly totally aspheric, when we define now
$\scrW_A\subset\Fl(\Ahat)$ as just $i_A^{-1}(\scrW)$, in order to
conclude $h_{\scrW_A}=h_\Ahat$, all we made use of was (besides
saturation of \scrW{} of course, i.e., \ref{it:60.a}) \ref{it:60.b}
and \ref{it:60.d}.

One may object that \ref{it:60.d} isn't expressed in terms of the
category structure of \Cat{} only, but we could express it in terms of
this structure, by the remark that $C$ has a final object if{f} there
exists a $\Simplex^-$-homotopy of $\id_C$ to a constant section of $C$
(this ``section'' will indeed be defined necessarily by a final object
of $C$). As was to be expected, in this formulation, as in
\ref{it:60.c} too, the object $\Simplex_1$ of \Cat{} is playing a
crucial role. But at this point it occurs to me that \ref{it:60.c}
implies \ref{it:60.d}, by the homotopy lemma -- thus for the time
being all we needed was \ref{it:60.a}\ref{it:60.b}\ref{it:60.c}.

\hangsection{Remembering about the promised ``key result'' at last!}\label{sec:61}%
We now in the long last get back to the ``key result'' promised time
ago, and which we kept pushing off. To pay off the trouble of the long
digression in between, maybe it'll come out more smoothly. It shall be
concerned with the functor
\begin{equation}
  \label{eq:61.1}
  i:A\to\Cat,
  \tag{1}
\end{equation}
where in the end $A$ will be (or turn out to be) a strict test
category, and we want to give characterizations for $i$ to be a (weak)
test functor, namely the corresponding functor
\begin{equation}
  \label{eq:61.2}
  i^*:\Cat\to\Ahat\tag{2}
\end{equation}
to induce an equivalence between the localizations, with respect to
``weak equivalences''.\footnote{Plus a little more, see below
  \eqref{eq:61.8}.} We'll now be a little more demanding, and instead
of just assuming it is the usual notion of weak equivalence either in
\Cat{} or in \Ahat, I'll assume that the set
$\scrW_A\subset\Fl(\Ahat)$ is defined in terms of a saturated set
$\scrW\subset\Fl(\Cat)$ of arrows in \Cat, by taking the inverse image
by the functor
\begin{equation}
  \label{eq:61.3}
  i_A:\Ahat\to\Cat\tag{3}
\end{equation}
as above, whereas in the target category \Cat{} of \eqref{eq:61.2},
we'll work with another saturated set $\scrW'$ -- thus besides
\eqref{eq:61.1}, the data are moreover
\begin{equation}
  \label{eq:61.4}
  \scrW,\scrW' \subset\Fl(\Cat),\tag{4}
\end{equation}
two\pspage{152} saturated sets of arrows in \Cat, with no special
assumption otherwise for the time being. We will introduce the
properties we need on these, as well as on $A$ and on $i$, stepwise as
the situation will tell us. We want to derive a set of conditions
ensuring that both \eqref{eq:61.2} and \eqref{eq:61.3} induce
equivalences for the respective localizations, namely
\begin{equation}
  \label{eq:61.5}
  \overline i^*:(\scrW')^{-1}\Cat \tosimeq \scrW_A^{-1}\Ahat
  \quad\text{and}\quad
  \overline i_A:\scrW_A^{-1}\Ahat \tosimeq \scrW^{-1}\Cat.
  \tag{5}
\end{equation}
We may assume beforehand that $A$ is a weak test category (``with
respect to \scrW'') and hence the second functor in \eqref{eq:61.5} is
already an equivalence, in which case the condition that the first
functor in \eqref{eq:61.5} exist and be an equivalence (existence just
meaning the condition
\begin{equation}
  \label{eq:61.star}
  \scrW'\subset (i^*)^{-1}(\scrW_A),\quad\text{i.e.,}\quad
  i^*(\scrW') \subset \scrW_A)\tag{*}
\end{equation}
is equivalent to the corresponding requirement for the composition
\begin{equation}
  \label{eq:61.6}
  f_i:\Cat \xrightarrow{i^*} \Ahat \xrightarrow{i_A} \Cat,\tag{6}
\end{equation}
namely that this induce an equivalence
\begin{equation}
  \label{eq:61.7}
  (\scrW')^{-1}\Cat \tosimeq \scrW^{-1}\Cat.\tag{7}
\end{equation}
As a matter of fact, we are going to be slightly more demanding (in
accordance with the notion of a weak test functor as developed
previously, cf.\ page \ref{p:85}), namely that the inclusion
\eqref{eq:61.star} be in fact an equality
\begin{equation}
  \label{eq:61.8}
  \scrW' = (i^*)^{-1}(\scrW_A),\tag{8}
\end{equation}
the similar requirement for the functor $i_A$ \eqref{eq:61.3} being
satisfied by the very definition of $\scrW_A$ in terms of \scrW{} as
\begin{equation}
  \label{eq:61.9}
  \scrW_A = i_A^{-1}(\scrW).\tag{9}
\end{equation}
In view of this, the extra requirement \eqref{eq:61.8} boils down to
the equivalent requirement in terms of the composition $f_i$
\eqref{eq:61.6}:
\begin{equation}
  \label{eq:61.10}
  \scrW' = f_i^{-1}(\scrW).\tag{10}
\end{equation}

To sum up, we want to at least develop sufficient conditions on the
data $(\scrW,\scrW',A,i)$ for \eqref{eq:61.8} to hold (which allows to
define the first functor in \eqref{eq:61.5}, whereas the second is
always defined), and the functors in \eqref{eq:61.5} to be
equivalences; or equivalently, for \eqref{eq:61.10} to hold, hence a
functor \eqref{eq:61.7}, and for the latter to be an equivalence, and
equally $i_A$ to induce an equivalence for the localizations (i.e.,
the second functor \eqref{eq:61.5} an equivalence). It should be noted
that the latter condition depends only on $(A,\scrW)$, not on $i$ nor
on $\scrW'$ -- it will be satisfied automatically if we assume $A$ to
be a \emph{weak test category} relative to \scrW{} (namely\pspage{153}
$i_A$ and the right adjoint functor
\[j_A:\Cat\to\Ahat\]
to induce quasi-inverse equivalences for the localizations
$\scrW_A^{-1}\Ahat$ and $\scrW^{-1}\Cat$). We have already developed
handy n.s.\ conditions for this in case $\scrW=W_\Cat$ -- and it would
be easy enough to look up which formal properties exactly on $W_\Cat$
have been used in the proof, if need be. At any rate, we know
beforehand that we can find $(A,\scrW)$ such that $A$ be a weak test
category (and even a strict test category!) relative to \scrW. When
$(A,\scrW)$ are chosen this way beforehand, the question just amounts
to finding conditions on $(\scrW',i)$ for \eqref{eq:61.10} to hold and
for \eqref{eq:61.7} to be an equivalence of categories. If we find
conditions which actually can be met, then we get as a byproduct the
formula \eqref{eq:61.10} precisely, which says that there is \emph{just one}
$\scrW'$ satisfying the conditions on $\scrW'$, namely
$f_i^{-1}(\scrW)$! Of course, taking \scrW{} to be just $W_\Cat$, it
will follow surely that $\scrW'$ is just $W_\Cat$ -- i.e., we should
get an axiomatic characterization of weak equivalences.

Let's now go to work, following the idea described in par.\
\ref{sec:57} (pages \ref{p:96}--\ref{p:98}), and expressed mainly in
the basic diagram of canonical maps in \Cat, associated to a given
object $C$ in \Cat:
\begin{equation}
  \label{eq:61.11}
  \begin{tabular}{@{}c@{}}
  \begin{tikzcd}[baseline=(O.base)]
    A_{/C} \ar[r] & A_{\sslash C} & \\
    & A\times C \ar[u]\ar[r] & |[alias=O]| A
  \end{tikzcd},
  \end{tabular}
  \tag{11}
\end{equation}
which will allow to compare $f_i(C) = A_{/C}$ with $C$. When
$\scrW=W_\Cat$, it was seen in loc.\ cit.\ that the two latter among
these three arrows are in \scrW, provided (for the middle one) we
assume that $i$ takes its values in the subcategory of \Cat{} of all
contractible categories. What remained to be done, for getting the
conditions for $f_i$ to be ``weakly equivalent'' to the identity
functor (and hence induce an equivalence for the localizations) was to
write down conditions for the first functor in \eqref{eq:61.11},
$A_{/C}\to A_{\sslash C}$, to be in \scrW. We have moreover to  be
explicit on the conditions to put on a general \scrW, in order for the
two latter maps in \eqref{eq:61.11} to be in \scrW. For the last
functor, this conditions as a matter of fact involves both $A$ and
\scrW, it is clearly equivalent to
\begin{description}
\item[\namedlabel{it:61.W1}{W~1)}]
  $A$ is \scrW-aspheric over $e$.
\end{description}
If we assume that $A$ has a final element, this condition is satisfied
provided \scrW{} satisfies the condition (where there is\pspage{154}
no $A$ anymore!) that any $X$ in \Cat{} with final element is
\scrW-aspheric over $e$ -- a condition which is similar to condition
\ref{it:60.d} above (page \ref{p:150}), but a littler stronger still
(as we want $X\times C\to C$ in \scrW{} for any $C$), it is a
consequence however of condition \ref{it:60.c}, as was seen on page
\ref{p:150} using the remark that $X$ is $\Simplex_1$-contractible. Thus
we get the handy condition
\begin{enumerate}[label=W~\arabic*')]
\item\label{it:61.W1prime}
  \scrW{} satisfies condition \ref{it:60.c} of page \ref{p:150}, i.e.,
  $\Simplex_1$ is \scrW-aspheric\\ over~$e$,
\end{enumerate}
which will allow even to handle the case of an $A$ which is
contractible (for the canonical homotopy structure of \Cat, namely
$\Simplex_1$-contractible), and not only when $A$ has a final object.

To insure that the canonical map
\begin{equation}
  \label{eq:61.star2}
  A\times C\to A_{\sslash C}\tag{*}
\end{equation}
is in \scrW, using the argument on page \ref{p:97}, we'll add one more
condition to the provisional list on page \ref{p:150}, namely:
\begin{enumerate}[label=\alph*),start=5]
\item\label{it:61.e}
  For any cartesian functor $u:F\to G$ of two fibered categories over
  a third one $B$ (everything in \Cat), such that the induced maps on
  the fibers are in \scrW, $u$ is in \scrW.
\end{enumerate}
We are now ready to state the condition we need (stronger than
\ref{it:61.W1prime}):
\begin{description}
\item[\namedlabel{it:61.W2}{W~2)}]
  \scrW{} satisfies conditions \ref{it:60.a} to \ref{it:60.c} (page
  \ref{p:150}) and \ref{it:61.e} above.
\end{description}

As a matter of fact, \ref{it:60.a} to \ref{it:60.c} ensure that
\scrW{} is ``geometric'', i.e., essentially $h_\scrW = h_\Cat$, hence
the proposition page \ref{p:143} applies, to imply that the maps
\[C \to \bHom(i(a), C)\]
are in \scrW{} (they are even $h_\scrW$-homotopisms) and by condition
\ref{it:61.e} this implies that \eqref{eq:61.star2} above is in
\scrW. We don't even need \ref{it:60.b} ($0$-admissibility for \scrW),
as all we care about is $h_\Cat\le h_\scrW$ (not the reverse
inequality), but surely we're going to need \ref{it:60.b} or something
stronger soon enough, as $\scrW=\Fl(\Cat))$ say surely wouldn't do!

Now to the last (namely first) map of our diagram \eqref{eq:61.11},
namely
\begin{equation}
  \label{eq:61.12}
  A_{/C} \to A_{\sslash C}.\tag{12}
\end{equation}
To give sufficient conditions for this to be in \scrW, we want to
mimic the standard asphericity criterion for a map in \Cat, which we
have used constantly before. This leads to the extra
condition
\begin{enumerate}[label=\alph*),resume]
\item\label{it:61.f}
  Let\pspage{155} $u:X'\to X$ be a map in \Cat{} such that for any $a$
  in $X$, the induced category $X'_{/a}$ be \scrW-aspheric, i.e.,
  $X'_{/a}\to e$ is in \scrW{} (or what amounts to the same if we
  assume \ref{it:60.d}, e.g., if we assume the stronger condition
  \ref{it:60.c}, the induced map $X'_{/a}\to X_{/a}$ is in
  \scrW). Then $u$ is in \scrW.
\end{enumerate}

If $u:X'\to X$ satisfies the condition stated above, namely that after
any base change $X_{/a}\to X$, the corresponding map $u_{/a}$ is in
\scrW, we'll say that $u$ is \emph{weakly} \scrW-aspheric (whereas
``\scrW-aspheric'' means that after \emph{any} base change $Y\to X$,
the corresponding $f_Y$ is in \scrW). Thus, condition \ref{it:61.f}
can be stated as saying that \emph{a weakly \scrW-aspheric map in
  \Cat{} is in \scrW}.

For making use of this latter assumption on \scrW, we have to look at
how the induced categories for the functor \eqref{eq:61.12} look like,
which functor (I recall) induced a bijection on objects. These can be
described as pairs $(a,p)$, with $a$ in $A$ and $p$ a map in \Cat
\[p:i(a)\to C.\]
An easy computation shows the
\begin{lemma}
  Let $(a,p)$ as above. The induced category $(A_{/C})_{/(a,p)}$
  \textup(for the functor \textup{\eqref{eq:61.12})} is canonically
  isomorphic to the induced category $A_{/G}$, where $G$ is the
  fibered product in \Ahat{} displayed in the diagram
  \begin{equation}
    \label{eq:61.13}
    \begin{tabular}{@{}c@{}}
    \begin{tikzcd}[baseline=(O.base)]
      G\ar[r]\ar[d] & a\ar[d] \\
      i^*(\bFl(C)) \ar[r] & |[alias=O]|i^*(C)
    \end{tikzcd},
    \end{tabular}
    \tag{13}
  \end{equation}
  where
  \[\bFl(C)\eqdef \bHom(\Simplex_1,C)\]
  and where the second horizontal arrow in \textup{\eqref{eq:61.13}}
  is the $i^*$-transform of the target map in \Cat
  \[\bFl(C) \xrightarrow t C.\]
\end{lemma}
\begin{corollary}
  In order for \textup{\eqref{eq:61.12}} to be weakly \scrW-aspheric,
  it is n.s.\ that the map
  \begin{equation}
    \label{eq:61.14}
    i^*(t) : i^*(\bFl(C)) \to i^*(C)\tag{14}
  \end{equation}
  in \Ahat{} be $\scrW_A$-aspheric \textup(i.e., be ``universally in $\scrW_A$''\textup).
\end{corollary}

To make the meaning of the latter condition clear, it should be noted
that the condition \ref{it:61.f} on \scrW{} guarantees precisely that
for a map $u:F'\to F$ in \Ahat{} ($A$ any category) to be
$\scrW_A$-aspheric, it is n.s.\ that the corresponding map $i_A(u)$ in
\Cat{} be weakly \scrW-aspheric -- the kind of thing we have been
constantly using before of course, when\pspage{156} assuming
$\scrW=W_\Cat$.

It is in the form of \eqref{eq:61.14} that weak \scrW-asphericity of
\eqref{eq:61.12} will actually be checked, whereas it will be
\emph{used} just by the fact that \eqref{eq:61.12} is in \scrW.

\bigbreak
\presectionfill\ondate{6.4.}\par

\hangsection[An embarrassing case of hasty
over-axiomatization. \dots]{An embarrassing case of hasty
  over-axiomatization. The unexpected riches\dots}\label{sec:62}%
I finally stopped with the notes last night, by the time when I
started feeling a little uncomfortable. A few minutes of reflection
then were enough to convince me that definitely I hadn't done quite
enough preliminary scratchwork yet on this ``key result'' business,
and embarked overoptimistically upon a ``mise en \'equation'' of the
situation, with the pressing expectation that a characterization of
weak equivalences should come out at the same time. First thing that
became clear, was that the introduction of two different localizing
sets of arrows $\scrW,\scrW'$ in \Cat{} was rather silly alas, nothing
at all would come out unless supposing from the very start
$\scrW=\scrW'$. Indeed, the crucial step for getting the ``key
result'' on test functors we are out for, goes as follows.

As the target map
\[t: \bFl(C)=\bHom(\Simplex_1,C) \to C\]
in \Cat{} is clearly a homotopy retraction, and $i^*:\Cat\to\Ahat$
commutes with products, we do have a good hold on the condition
\eqref{eq:61.14} of the last corollary, namely that $i^*(t)$ be
$\scrW_A$-aspheric -- e.g., it is enough that the contractor
$i^*(\Simplex_1)$ in \Ahat{} be $\scrW_A$-aspheric over $e_\Ahat$ (for
instance, it is often enough it be $0$-connected!). In view of the
corollary and the condition \ref{it:61.f} on \scrW{}
(\hyperref[p:155]{last page}), we thus get a very good hold upon the
map
\begin{equation}
  \label{eq:62.star}
  A_{/C} \to A_{\sslash C}\tag{*}
\end{equation}
in \Cat{} being in \scrW, and hence on all three maps in the diagram
\eqref{eq:61.11} (page \ref{p:153}) being in \scrW. With this in mind,
the key step can be stated as follows:
\begin{lemma}
  Assume that \scrW{} satisfies the conditions
  \textup{\ref{it:60.a}\ref{it:60.c}\ref{it:61.e}\ref{it:61.f}}
  \textup(pages \ref{p:150}, \ref{p:154}, \ref{p:155}\textup),
  that $A$ is \scrW-aspheric over $e$ \textup(i.e., $A\times C\to C$
  is in \scrW{} for any $C$, which will be satisfied if $A$ is
  $\Simplex_1$-contractible in \Cat, for instance if $A$ has a final or
  initial object\textup), and that the objects $i(a)$ in \Cat{}
  \textup(for any $a$ in $A$\textup) are contractible \textup(for the
  canonical homotopy structure of \Cat, i.e., $\Simplex_1$-contractible,
  or even only for the wider homotopy structure $h_\scrW$ based on
  \scrW-aspheric homotopy intervals). Under\pspage{157} these
  conditions, the following conclusions hold:
  \begin{enumerate}[label=\alph*),font=\normalfont]
  \item\label{it:62.a}
    $\scrW = f_i^{-1}(\scrW)\quad$ \textup(where $f_i=i_Ai^*$ with
    yesterday's notations\textup).
  \item\label{it:62.b}
    The functor $\overline f\!_i$ from $\scrW^{-1}\Cat$ to itself
    induced by $f_i$ \textup(which is defined because of
    \textup{\ref{it:62.a})} is isomorphic \textup(canonically\textup)
    to the identity functor, and hence is an equivalence.
  \end{enumerate}
\end{lemma}

The use we have for the three maps in the diagram \eqref{eq:61.12} is
completely expressed in this lemma. The pretty obvious proof below
would not work at all if in \ref{it:62.a} above, we replace \scrW{} in
the left hand side by a $\scrW'$! We have to prove that for a map
$C\to C'$ in \Cat, this is in \scrW{} if{f} $A_{/C}\to A_{/C'}$
is. Now this is seen from an obvious diagram chasing in the diagram
below, using saturation condition \ref{it:37.b} on \scrW:
\[\begin{tikzcd}[baseline=(O.base)]
  A_{/C} \ar[d]\ar[r] & A_{\sslash C} \ar[d] &
  A\times C \ar[r]\ar[d]\ar[l] & C \ar[d] \\
  A_{/C} \ar[r] & A_{\sslash C}  &
  A\times C \ar[r]\ar[l] & |[alias=O]| C
\end{tikzcd},\]
where all horizontal arrows are already known to be in \scrW{} (the
assumptions in the lemma were designed for just that end). At the same
time, we see that the corresponding statements are equally true for
the functors $C\mapsto A_{\sslash C}$ and $C\mapsto A\times C$, and
that two consecutive among the four functors we got from
$H_\scrW=\scrW^{-1}\Cat$ to itself, deduced by localization by \scrW,
are canonically isomorphic, which proves \ref{it:62.b} by taking the
composition of the three isomorphisms
\[\gamma(A_{/C}) \tosim \gamma(A_{\sslash C}) \tosim \gamma(A\times C)
\tosim \gamma(C),\]
where $\gamma:\Cat\to H_\scrW=\scrW^{-1}\Cat$ is the canonical
functor.

With this lemma, we have everything needed in order to write down the
full closed chain of implications, between various conditions on
$(\scrW, A, i)$, from which to read off the ``key result'' we're
after. Before doing so, I would like still to make some preliminary
comments on the role of \scrW, and on the nature of the conditions we
have been led so far to impose upon \scrW.

A first feature that is striking, is that all conditions needed are in
the nature either of stability conditions (if such and such maps are
in \scrW, so are others deduced from them), or conditions stating that
such and such unqualified maps (the projection $\Simplex_1\times C\to C$
for any $C$, say) are in \scrW. We did not have any use of the only
condition stated so far, namely \ref{it:60.b} (if $f\in\scrW$,
$\piz(f)$ is bijective) of a \emph{restrictive} type on the kind of
arrows allowed in \scrW{} -- which is quite contrary to my
expectations. Thus, all conditions are trivially met if we
take\pspage{158} $\scrW=\Fl(\Cat)$ -- \emph{all} arrows in \Cat! This
circumstance seems tied closely to the fact that, contrarily to quite
unreasonable expectations, \emph{we definitely do \emph{not} get an
  axiomatic characterization of weak equivalences}, in terms of the
type of properties of \scrW{} we have been working with so far. As
soon as one stops for considering the matter without prejudice, this
appears rather obvious. As a matter of fact, using still cohomological
invariants of topoi and categories, there are lots of variants of the
cohomological definition of ``weak equivalence'', which will share all
formal properties of the latter we have been using so far, and
presumably a few more we haven't met yet. For instance, starting with
any ring $k$ (interesting cases would be \bZ, $\bZ/n\bZ$, \bQ), we may
demand on a morphism of topoi
\[f:X\to Y\]
to induce as isomorphism for cohomology with coefficients in $k$, or
with coefficients in any $k$-module, or with any twisted coefficients
which are $k$-modules -- already three candidates for a \scrW,
depending on a given $k$! We may vary still more, by taking, instead
of just one $k$, a whole bunch $(k_i)$ of such, or a bunch of
(constant) commutative groups -- we are thinking of choices such as
all rings $\bZ/n\bZ$, with possibly $n$ being subjected to be prime to
a given set of primes, along the lines of the Artin-Mazur theory of
``localization'' of homotopy types. And we may combine this with an
isomorphism requirement on twisted non-commutative $1$-cohomology, as
for the usual notion of weak equivalence. Also, in all the isomorphism
requirements, we may restrict to cohomology up to a certain dimension
(which will give rise to ``truncated homotopy types''). The impression
that goes with the evocation of all these examples, is that the theory
we have been pursuing, to come to an understanding of ``models for
homotopy types'', while we started with just usual homotopy types in
mind and a corresponding tacit prejudice, is a great deal richer than
what we had in mind. Yesterday's (or rather last night's)
embarrassment of finding out finally I had been very silly, is a
typical illustration of the embarrassment we feel, whenever a
foreboding appears of our sticking to inadequate ideas; still more so
if it is not just mathematics but ideas about ourselves say or about
something in which we are strongly personally implies. This
embarrassment then comes as a rescue, to bar the way to an
unwished-for overwhelming richness dormant in ourselves, ready to wash
away forever those ideas so dear to us\ldots

I\pspage{159} am definitely going to keep from now on a general
\scrW{} and work with this and the corresponding localization, which
in case of ambiguity we better won't denote by \Hot{} any longer (as
we might be thinking of usual homotopy types in terms of usual weak
equivalences), but by $H_\scrW$ or $\Hot_\scrW$, including such
notions as rational homotopy types, etc.\ (for suitable choices of
\scrW). The idea that now comes to mind here is that possibly, the
usual $W_\Cat$ of usual weak equivalences could be characterized as
the \emph{smallest} of all \scrW's, satisfying the conditions we have
been working with so far (tacitly to some extent), and maybe a few
others which are going to turn up in due course -- i.e.\ that the
usual notion of weak equivalence is the \emph{strongest} of all
notions, giving rise to a modelizing theory as we are developing. This
would be rather satisfactory indeed, and would imply that other
categories $H_\scrW$ we are working with \emph{are all localizations
  of \Hot}, with respect to a saturated set of arrows in \Hot,
satisfying some extra conditions which it may be worth while writing
down explicitly, in terms of the internal structures of \Hot{}
directly (if at all possible). All the examples that have been
flashing through my mind a few minutes ago, do correspond indeed to
equivalence notions weaker than so-called ``weak'' equivalence, and
hence to suitable localizations of the usual homotopy category
\Hot. But it is quite conceivable that this is not so for all \scrW's,
namely that the characterization just suggested of $W_\Cat$ is not
valid. This would mean that there are refinements of the usual notion
of homotopy types, which would still however give rise to a homotopy
theory along the lines I have been pursuing lately. There is of course
an immediate association with Whitehead's \emph{simple homotopy types}
-- maybe after all they can be interpreted as elements in a suitable
localization $H_\scrW$ of \Cat{} (and correspondingly, of any one of
the standard modelizers, such as semi-simplicial complexes and the
like)? In any case, sooner or later one should understand what the
smallest of all ``reasonable'' \scrW's looks like, and to which
geometric reality it corresponds. But all these questions are not
quite in the present main line of thought, and it is unlikely I am
going to really enter into it some day\ldots

\hangsection{Review of terminology
  \texorpdfstring{\textup(}{(}provisional\texorpdfstring{\textup)}{)}.}%
\label{sec:63}%
What I should do though immediately, is to put a little order in the
list of conditions for a set \scrW, which came out somewhat
chaotically yesterday. After the notes I still did a little
scratchwork last night, which I want now to write down, before coming
to a formal statement of the ``key result'' -- as this will of course
make\pspage{160} use of some list of conditions on \scrW.

First of all, I feel a review is needed of the few basic notions which
have appeared in our work, relative to a set of arrows
$W\subset\Fl(M)$ in a general category $M$. We will not give to the
maps in $W$ a specific name, such as ``weak equivalences'', as this
may be definitely misleading, in the general axiomatic set-up we want
to develop; here $W_\Cat$ is just one among many possible \scrW's and
correspondingly for a small category $A$, $W_A=W_\Ahat$ is just one
among the many $\scrW_A$'s, associated to the previous \scrW's. When
$M=\Cat$, it will be understood we are working with a fixed set
$\scrW\subset\Fl(\Cat)$, consisting of the basic ``equivalences'', on
which the whole modelizer story hinges. We may call them
\scrW-equivalences -- for the time being there will be no question of
varying \scrW.

Coming back to a general pair $(M, W\subset\Fl(M))$ (not necessarily a
``modelizer''), we may call the maps in $W$ $W$-equivalences. If $M$
has a final object $e$, we get the corresponding notion of
\emph{$W$-aspheric object} of $M$, namely an object $X$ such that the
unique map
\[ X\to e\]
is in $W$, i.e., is a $W$-equivalence. We'll define a
\emph{$W$-aspheric map}
\[f : X\to Y\]
in $M$ as one which is ``universally in $W$'', by which I mean that
for any base-change
\[ Y' \to Y,\]
the fiber-product $X'=X\times_Y Y'$ exists (i.e., $f$ is
``squarable'') and the map
\[ f': X' \to Y'\]
deduced from $f$ by base change is in $W$. The thing to be quite
careful about is that for an object $X$ in $M$, to say that $X$ is
\emph{$W$-aspheric over $e$} (meaning that the \emph{map} $X\to e$ is
$W$-aspheric) implies $X$ is $W$-aspheric \&c., but the converse need
not hold true. This causes a slight psychological uneasiness, due to
the fact I guess that the notion of a $W$-aspheric object has been
defined after all in terms of the \emph{map} $X\to e$, and
consequently may be thought of as meaning is ``in $W$ over
$e$''. Maybe we shouldn't use at all the word ``$W$-aspheric object''
here, not even by qualifying it as ``weakly $W$-aspheric'' to cause a
feeling of caution, but rather refer to this notion as ``$X$ is a
$W$-object'' -- and denote by $M(W)$ the set of all these objects (or
the corresponding full subcategory of $M$, and call $X$ $W$-aspheric
(dropping ``over $e$'')\footnote{But we'll see immediately that this
  conflicts with the standard terminology in topoi -- so no good
  either!} when it is ``universally'' a $W$-object. The terminology we
have been using so far was of course\pspage{161} suggested by the case
when $M$ is a \emph{topos} and $W$ the usual notion of weak
equivalence, but then to call $X$ in $M$ ``$W$-aspheric'' or simply
``aspheric'' does correspond to the usual (absolute) notion of
asphericity for the induced topos $M_{/X}$, only in the case when the
topos $M$ itself is aspheric. This is so in the case I was most
interested in (e.g., $M$ a modelizing topos), but if we want to use it
systematically in the general setting, the term I used of
``$W$-aspheric object'' is definitely misleading. Thus we better
change it now than never, and use the word ``$W$-object'' instead, and
the notation $M(W)$. As for the notion of $W$-aspheric map, in the
present case of a topos with the notion of weak equivalence, it does
correspond to the usual notion of asphericity for the induced morphism
of topoi
\[M_{/X} \to M_{/Y},\]
which is quite satisfactory.

There is still need for caution with the notion of $W$-aspheric maps
($W$-aspheric objects have disappeared in the meanwhile!), when
working in $M=\Cat$ (and the same thing if $M=\Spaces$). Namely, when
we got a map $f:A\to B$ in \Cat, this is viewed for topological
intuition as corresponding to a morphism of topoi
\[\Ahat\to\Bhat.\]
Now, the requirement that $f:A\to B$ should be $W_\Cat$-aspheric is a
lot stronger than the asphericity of the corresponding morphism of
topoi. Indeed, the latter just means that for any base-change in
\Cat{} of the very particular ``localization'' or ``induction'' type,
namely
\[B_{/b} \to B,\]
the corresponding map deduced by base change
\[ f_{/b}: A_{/b} \to B_{/b} \]
is a weak equivalence (or equivalently, that $A_{/b}$ is aspheric),
whereas $W$-asphericity of $f$ means that the same should hold for
\emph{any} base-change $B'\to B$ in \Cat, or equivalently, that for
\emph{any} such base-change, with $B'$ having a final element
moreover, the corresponding category $A'=A\times_B B'$ is aspheric. To
keep this distinction in mind, and because the weaker notion is quite
important and deserves a name definitely, I will refer to this notion
by saying $f$ is \emph{weakly $W$-aspheric} (returning to the case of
a general $\scrW\subset\Fl(\Cat)$) if for any base change of the
particular type $B_{/b}\to B$ above, the corresponding\pspage{162} map
$A_{/b}\to B_{/b}$ is in \scrW. We could express this in terms of the
morphism of topoi $\Ahat\to\Bhat$ by saying that the latter is
$W$-aspheric -- being understood that the choice of an ``absolute''
\scrW{} in \Cat, implies as usual a corresponding notion of
\scrW-equivalence for arbitrary morphisms of arbitrary topoi, in terms
of the corresponding morphism between the corresponding homotopy
types. (This extension to topoi of notions in \Cat{} should be made
quite explicit sooner or later, but visibly we do not need it yet for
the time being.) One relationship we have been constantly using, and
which is nearly tautological, comes from the case of a map
\[u:F\to G\]
in a category \Ahat, hence applying $i_A$ a map in \Cat
\[A_F \to A_G.\]
For this map to be weakly \scrW-aspheric, it is necessary and
sufficient that for any base-change in \Ahat{} of the particular type
\[ G'=A \to G\quad\text{with $a$ in $A$,}\]
the corresponding map
\[ F\times_G a \to a\]
in \Ahat{} be in $\scrW_A$ -- a condition which is satisfied of course
if $u$ is $\scrW_A$-aspheric. This last condition is also necessary,
if we assume that \scrW{} satisfies the standard property we've kept
using all the time in case of \Cat, namely that a map in \Cat{} that
is weakly $W$-aspheric, is in \scrW{} (condition \ref{it:61.f} on page
\ref{p:155}). Thus we get the
\begin{proposition}
  Assume that any map in \Cat{} which is weakly \scrW-aspheric is in
  \scrW, and let $A$ be any small category, $u:F\to G$ a map in
  \Ahat. Then $u$ is $\scrW_A$-aspheric \textup(and hence in
  $\scrW_A$\textup) if{f} $i_A(u) : A_{/F} \to A_{/G}$ is weakly
  \scrW-aspheric.
\end{proposition}

The assumption we just made on \scrW{} is of such constant use, that
we are counting it among those we are making once and for all upon
\scrW{} (which I still have to pass in review).

As for \scrW-aspheric maps in \Cat, this is a very strong notion
indeed when compared with \scrW-equivalence or even with weak
\scrW-asphericity. We did not have any use for it yet, but presumably
this notion will be of importance when it comes to a systematic study
of the internal properties of \Cat{} with respect to \scrW{} (which is
still in our program!).

Coming\pspage{163} back to a general $(M,W)$, we have defined earlier
a canonical homotopy structure $h_W$ on $M$, which we may call
``associated to $W$'' -- this is also the weak homotopy interval
structure on $M$, generated by intervals in $M$ which are $W$-aspheric
over $e$. This makes sense at least, provided in $M$ finite products
exist. If moreover $M$ satisfies the conditions \ref{it:51.F.b} to
\ref{it:51.F.d} of page \ref{p:123} concerning sums and connected
components of objects of $M$, we have defined (independently of $W$)
the \emph{canonical} homotopy structure $h_M$ in $M$, which may be
viewed as the weak homotopy interval structure generated by all
$0$-connected intervals in $M$. It still seems that in all cases we
are going to be interested in, we have the equality
\[h_W = h_M.\]
When speaking of homotopy notions in $M$ (such as $f$ and $g$ being
homotopic maps, written
\[ f\sim g,\]
or a map being a homotopism, or an object being contractible) it will
be understood (unless otherwise stated) that this refers to the
homotopy structure $h_W=h_M$. In case we should not care to impose
otherwise unneeded assumptions which will imply $h_W=h_M$, we'll be
careful when referring to homotopy notions, to say which structure we
are working with.

\hangsection{Review of properties of the ``basic localizer''
  \texorpdfstring{$\scrW_\Cat$}{W(Cat)}.}\label{sec:64}%
We recall that a set of arrows $W\subset\Fl(M)$ is called
\emph{saturated} if it satisfies the conditions:
\begin{enumerate}[label=\alph*)]
\item\label{it:64.a}
  Identities belong to $W$.
\item\label{it:64.b}
  If $f,g$ are maps and $fg$ exists, then if two among $f,g,gf$ are in
  $W$, so is the third.
\item\label{it:64.c}
  If $f:X\to Y$ and $g:Y\to X$ are such that $gf$ and $fg$ are in $W$,
  so are $f$ and $g$.
\end{enumerate}

On the other hand, \emph{strong saturation} for $W$ means that $W$ is
the set of arrows made invertible by the localization functor
\[M\to M_W=W^{-1}M,\]
or equivalently, that $W$ can be described as the set of arrows made
invertible by some functor $M\to M'$. The trouble with strong
saturation is that it is a condition which often is not easy to check
in concrete situations. This is so for instance for the notion of weak
equivalence in \Cat, and the numerous variants defined in terms of
cohomology. Therefore, we surely won't impose the strong saturation
condition\pspage{164} on \scrW{} (which we may call the ``\emph{basic
  localizer}'' in our modelizing story), but rather be happy if we can
prove strong saturation as a consequence of other formal properties of
\scrW, which are of constant use and may be readily checked in the
examples we have in mind.

Let's give finally a provisional list of those properties for a
``localizer'' \scrW.
\begin{description}
\item[\namedlabel{it:64.L1}{L~1)}] \textbf{(Saturation)}
  \scrW{} is saturated, i.e., satisfies conditions
  \ref{it:64.a}\ref{it:64.b}\ref{it:64.c} above.
\item[\namedlabel{it:64.L2}{L~2)}] \textbf{(Homotopy axiom)}
  $\Simplex_1$ is \scrW-aspheric over $e$, i.e., for
  any $C$ in \Cat, the projection
  \[\Simplex_1\times C\to C\]
  is in \scrW.
\item[\namedlabel{it:64.L3}{L~3)}] \textbf{(Final object axiom)}
  Any $C$ in \Cat{} which has a final object is
  in $\Cat(\scrW)$, i.e., $C\to e$ is in \scrW{} (or, as we will still
  say when working in \Cat, $C$ is \scrW-aspheric\footnote{For a
    justification of this terminology,see proposition on p.\
    \ref{p:167} below.}).
\item[\namedlabel{it:64.L3prime}{L~3')}] \textbf{(Interval axiom)}
  $\Simplex_1$ is \scrW-aspheric, i.e., $\Simplex_1\to e$ is in \scrW.
\item[\namedlabel{it:64.L4}{L~4)}] \textbf{(Localization axiom)}
  Any map $u:A\to B$ in \Cat{} which is weakly \scrW-aspheric (i.e.,
  the induced maps $A_{/b} \to B_{/b}$ are in \scrW) is in \scrW.
\item[\namedlabel{it:64.L5}{L~5)}] \textbf{(Fibration axiom)}
  If $f:X\to Y$ is a map in \Cat{} over an object $B$ of \Cat,
  such that $X$ and $Y$ are fiber categories over $B$ and $f$ is
  cartesian, and if moreover for any $b\in B$, the induced map on
  the fibers $f_b:X_b\to Y_b$ is in \scrW, then so is $f$.
\end{description}

These properties are all I have used so far, it seems, in the case
$\scrW=W_\Cat$ we have been working with till now, in order to develop
the theory of test categories and test functors, including ``weak''
and ``strong'' variants, and including too the generalized version of
the ``key result'' which is still waiting for getting into the
typewriter. Let's list at once the implications
\[ \text{\ref{it:64.L2}} \Rightarrow \text{\ref{it:64.L3}} \Rightarrow
\text{\ref{it:64.L3prime},}\]
and
\[\parbox{0.8\textwidth}{if \ref{it:64.L1} and \ref{it:64.L4} hold
  (saturation and localization axioms), then the homotopy axiom
  \ref{it:64.L2} is already implied by the final object axiom
  \ref{it:64.L3}.}\]

Thus, the set of conditions \ref{it:64.L1} to \ref{it:64.L4} (not
including the last one \ref{it:64.L5}, i.e., the fibration axiom) is
equivalent to the conjunction of \ref{it:64.L1} \ref{it:64.L3}
\ref{it:64.L4}. This set of conditions is of such a constant use, that
we'll assume it throughout, whenever there is a \scrW{} around:
\begin{definition}
  A\pspage{165} subset \scrW{} of $\Fl(\Cat)$ is called a \emph{basic localizer},
  if it satisfies the conditions \ref{it:64.L1}, \ref{it:64.L3},
  \ref{it:64.L4} above (saturation, final object and localization
  axioms), and hence also the homotopy axiom \ref{it:64.L2}.
\end{definition}

These conditions are enough, I quickly checked this night, in order to
validify all results developed so far on test categories, weak test
categories, strict test categories, weak test functors and test
functors with values in \Cat{} (cf.\ notably the review in par.\
\ref{sec:44}, pages \ref{p:79}--\ref{p:88}), provided in the case of
test functors we restrict to the case of loc.\ cit.\ when each of the
categories $i(a)$ has a final object. All this I believe is
justification enough for the definition above.

As for the fibration axiom \ref{it:64.L5}, this we have seen to be
needed (at least in the approach we got so far) for handling test
functors $i:A\to\Cat$, while no longer assuming the categories $i(a)$
to have final objects (which was felt to be a significant
generalization to carry through, in view of being able subsequently to
replace \Cat{} by more general modelizers). While still in the nature
of a stability requirement, this fibration axiom looks to me a great
deal stronger than the other axioms. Clearly \ref{it:64.L5}, together
with the very weak ``interval axiom'' \ref{it:64.L3prime}
($\Simplex_1\to e$ is in \scrW) implies the homotopy axiom. It can be
seen too that when joined with \ref{it:64.L1}, it implies the
localization axiom \ref{it:64.L4} (using the standard device of a
mapping-cone for a functor\ldots). Thus, a basic localizer satisfying
the fibration axiom \ref{it:64.L5} can be viewed also as a \scrW{}
satisfying the conditions
\[\text{\ref{it:64.L1} (saturation),\enspace\ref{it:64.L3prime}
  (interval axiom),\enspace\ref{it:64.L5} (fibration axiom).}\]
In the next section, after we will have stated the two key facts about
weak test functors and test functors, which both make use of
\ref{it:64.L5}, we'll presumably, for the rest of the work ahead
towards canonical modelizers, assume the fibration condition on the
basic localizer \scrW.

There are some other properties of weak equivalence $W_\Cat$ and its
manifold variants in terms of cohomology, which have not been listed
yet, and which surely will turn up still sooner or later. Maybe it's
too soon to line them up in a definite order, as their significance is
still somewhat vague and needs closer scrutiny. I'll just list those
which come to my mind, in a provisional order.
\begin{description}
\item[\namedlabel{it:64.La}{L~a)}]
  $0$-admissibility of \scrW, namely $f\in\scrW$ implies $\piz(f)$ bijective.
\end{description}

This condition, together with the homotopy axiom \ref{it:64.L2}, will
imply\pspage{166}
\begin{equation}
  \label{eq:64.1}
  h_\scrW = h_\Cat\quad\parbox[t]{.55\textwidth}{the canonical
    homotopy structure on \Cat{} defined in terms of the generating
    contractor $\Simplex_1$ in \Cat,}
  \tag{1}
\end{equation}
whereas the homotopy axiom alone, I mean without
\hyperref[it:64.La]{a)}, will imply only the inequality
\begin{equation}
  \label{eq:61.1prime}
  h_\Cat\le h_\scrW,\tag{1'}
\end{equation}
which is all we care for at present. The latter implies that two maps
in \Cat{} which are $\Simplex_1$-homotopic (i.e., belong to the same
connected component of $\bHom(X,Y)$ have the same image in the
localization $\scrW^{-1}\Cat=\HotW$, and that any map in \Cat{} which
is a $\Simplex_1$-homotopism is in \scrW, and even is \scrW-aspheric if
it is a ``homotopy retraction'' with respect to the
$\Simplex_1$-structure $h_\Cat$. However, in practical terms, even
without assuming \hyperref[it:64.La]{a)} expressly, we may consider
\eqref{eq:64.1} to be always satisfied. This amounts indeed to the
still weaker condition than \hyperref[it:64.La]{a)}
\begin{description}
\item[\namedlabel{it:64.Laprime}{L~a')}]
  If $C$ in \Cat{} is \scrW-aspheric over $e$, it is $0$-connected,
  i.e., non-empty and connected.
\end{description}

But if it were empty, it would follow that for any $X$ in \Cat,
$\varnothing\to X$ is in \scrW{} and hence \HotW{} is equivalent to the
final category. If $C$ is non-empty and disconnected, choosing two
connected components and one point in each to make $C$ into a weak
homotopy interval for $h_\scrW$, one easily gets that any two maps
$f,g:X\rightrightarrows Y$ in \Cat{} are $h_\scrW$-homotopic, hence
have the same image in \HotW, which again must be the final
category. Thus we get:
\begin{proposition}
  If \scrW{} satisfies \textup{\ref{it:64.L1}},
  \textup{\ref{it:64.L2}} \textup(e.g., \scrW{} a basic
  modelizer\textup), then we have equality \eqref{eq:64.1}, except in
  the case when \HotW{} equivalent to the final category.
\end{proposition}

This latter case isn't too interesting one will agree. Thus, we would
easily assume \eqref{eq:64.1}, i.e.,
\hyperref[it:64.Laprime]{a')}. But the slightly stronger condition
\hyperref[it:64.La]{a)} seems hard to discard; even if we have not
made any use of it so far, one sees hardly of which use a category of
localized homotopy types \HotW{} could possibly be, if one is not even
able to define the \piz-functor on it! Thus, presumably we'll have to
add this condition, and maybe even more, in order to feel \scrW{}
deserves the name of a ``basic localizer''\ldots Among other formal
properties which still need clarification, even in the case of
$W_\Cat$, there is the question of exactness properties of the
canonical functor
\begin{equation}
  \label{eq:64.2}
  \Cat\to\scrW^{-1}\Cat=\HotW.\tag{2}
\end{equation}
I am thinking particularly of the following
\begin{description}
\item[\namedlabel{it:64.Lb}{L~b)}]
  The\pspage{167} functor \eqref{eq:64.2} commutes with finite sums,
\end{description}
possibly even with infinite ones, which should be closely related to
property \hyperref[it:64.La]{a)}, and
\begin{description}
\item[\namedlabel{it:64.Lc}{L~c)}]
  The functor \eqref{eq:64.2} commutes with finite products
\end{description}
(maybe even with infinite ones, under suitable assumptions).

The following property, due to Quillen for weak equivalences, is used
in order to prove for these (and the cohomological analogs) the
fibration axiom \ref{it:64.L5} (what we get directly is the case of
cofibrations, as a matter of fact -- cf.\ prop.\ page \ref{p:97}):
\begin{description}
\item[\namedlabel{it:64.Ld}{L~d)}]
  If $f:C\to C'$ is in \scrW, so is $f\op:C\op\to(C')\op$ for the dual
  categories.
\end{description}

This is about all I have in mind at present, as far as further
properties of a \scrW{} is concerned. The property
\hyperref[it:64.Lc]{c)} however brings to mind the natural (weaker)
condition that the cartesian product of two maps in \scrW{} should
equally be in \scrW. The argument in the beginning of par.\
\ref{sec:40} (p.\ \ref{p:69}) carries over here, and we get:
\begin{proposition}
  Let \scrW{} be a basic localizer, and $C$ in \Cat{} such that $C$ is
  \scrW-aspheric, i.e., $C\to e$ is in \scrW, then $C$ is even
  \scrW-aspheric over $e$, i.e., for any $A$ in $C$, $C\times A\to A$
  is in \scrW.
\end{proposition}

However, the proof given for the more general statement we have in
mind (of proposition on page \ref{p:69}) does not carry over using
only the localization axiom in the form \ref{it:64.L4} it was stated
above, as far as I can see. This suggests a stronger version
\namedlabel{it:64.L4prime}{L~4')} of \ref{it:64.L4} which we may have
to use eventually, relative to a commutative triangle in \Cat{} as on
page \ref{p:70}
\[\begin{tikzcd}[baseline=(O.base),column sep=tiny,row sep=small]
  P' \ar[dr]\ar[rr,"F"] & & P\ar[dl] \\ & |[alias=O]| C &
\end{tikzcd},\]
when assuming that the induced maps (for arbitrary $c$ in $C$)
\[F_{/c}: P'_{/c}\to P_{/c}\]
are in \scrW, to deduce that $F$ itself is in \scrW. However, this
``strong localization axiom'' is a consequence (as is the weaker one
\ref{it:64.L4}) of the fibration axiom \ref{it:64.L5}, which implies
also directly the property we have in mind, namely
\[\text{$f,g\in\scrW$ implies $f\times g\in\scrW$.}\]

To\pspage{168} come to an end of this long terminological and
notational digression, I'll have to say a word still about test
categories, modelizers, and test functors. We surely want to use
freely the terminology introduced so far, while we were working with
ordinary weak equivalences, in the more general setting when a basic
localizer \scrW{} is given beforehand. As long as there is only one
\scrW{} around, which will be used systematically in all our
constructions, we'll just use the previous terminology, being
understood that a ``modelizer'' say will mean a ``\scrW-modelizer'',
namely a category $M$ endowed with a saturated $W\subset\Fl(M)$, such
that $W^{-1}M$ is equivalent some way or other to
$\scrW^{-1}\Cat$. The latter category, however, I dare not just
designate as \Hot, as this notation has been associated to the very
specific situation of just ordinary homotopy types, therefore I'll
always write \HotW{} instead, as a reminder of \scrW{} after all! If
at a later moment it should turn out that we'll have to work with more
than one \scrW{} (for instance, to compare the \scrW-theory to the
ordinary $W_\Cat$-theory), we will of course have to be careful and
reintroduce \scrW{} in our wording, to qualify all notions dependent
on the choice of \scrW.

\bigbreak
\presectionfill\ondate{7.4.}\par

\hangsection[Still another review of the test notions (relative to
given \dots]{Still another review of the test notions
  \texorpdfstring{\textup(}{(}relative to given basic localizer
  \texorpdfstring{\scrW\textup)}{W)}.}\label{sec:65}%
It has been over a week now and about eighty pages typing, since I
realized the need for looking at more general test functors than
before and hit upon how to handle them, that I am grinding stubbornly
through generalities unending on homotopy notions. The grinding is a
way of mine to become familiar with a substance, and at the same time
getting aloof of it climbing up, sweatingly maybe, to earn a birds-eye
view of a landscape and maybe, who knows, in the end start a-flying in
it, wholly at ease\ldots I am not there yet! The least however one
should expect, is that the test story should now go through very
smoothly. As I have been losing contact lately with test categories
and test functors, I feel it'll be worth while to make still another
review of these notions, leading up to the key result I have been
after all that time. It will be a way both to gain perspective, and
check if the grinding has been efficient indeed\ldots

If
\[u:M\to M'\]
is a functor between categories endowed each with a saturated set of
arrows, $W$ and $W'$ say, we'll say $u$ is ``\emph{model preserving}''
(with respect to $(W,W')$ if it satisfies the conditions:
\begin{enumerate}[label=\alph*)]
\item\label{it:65.a}
  $W=u^{-1}(W')$\pspage{169} (hence the functor $\overline
  u:W^{-1}M\to (W')^{-1}M'$ exists),
\item\label{it:65.b}
  the functor $\overline u$ is an equivalence.
\end{enumerate}

We do not assume beforehand that $(M,W)$, $(M',W')$ are modelizers
(with respect to a given basic localizer $\scrW\subset\Fl(\Cat)$), but
in the cases I have in mind, we'll know beforehand at least one of the
pairs to be a modelizer, and it will follow the other is one too.

In all what follows, a basic localizer \scrW{} is given once and for
all (see definition on page \ref{p:165}). For the two main results
below on the mere general test functors, we'll have to assume \scrW{}
satisfies the fibration axiom \ref{it:64.L5} (page \ref{p:164}). We
are going to work with a fixed small category $A$, without any other
preliminary assumptions upon $A$, all assumptions that may be needed
later will be stated in due course. Recall that $A$ can be considered
as an object of \Cat, and we'll say $A$ is \scrW-aspheric if $A\to e$
is in \scrW, which implies that $A\to e$ is even \scrW-aspheric, i.e.,
$A\times C\to C$ is in \scrW{} for any $C$ in \Cat. More generally, if
$F$ is in \Ahat, we'll call $F$ \scrW-aspheric if the category
$A_{/F}$ is \scrW-aspheric. Thus, to say $A$ is \scrW-aspheric just
means that the final object $e_\Ahat$ of \Ahat{} is \scrW-aspheric.

We'll constantly be using the canonical functor
\begin{equation}
  \label{eq:65.1}
  i_A:\Ahat\to\Cat, \quad F\mapsto A_{/F},\tag{1}
\end{equation}
and its right adjoint
\begin{equation}
  \label{eq:65.2}
  j_A = i_A^*:\Cat\to\Ahat, \quad C\mapsto (a\mapsto
  \Hom(A_{/a},C)).\tag{2}
\end{equation}
The category \Ahat{} will always be viewed as endowed with the
saturated set of maps
\begin{equation}
  \label{eq:65.3}
  \scrW_A = i_A^{-1}(\scrW).\tag{3}
\end{equation}
This gives rise to the notions of $\scrW_A$-aspheric map in \Ahat, and
of an object $F$ of \Ahat{} being $\scrW_A$-aspheric over $e_\Ahat$,
namely $F\to e_\Ahat$ being $\scrW_A$-aspheric, which means that
$F\times G\to G$ is in $\scrW_A$ for any $G$ in \Ahat, which by
definition of \scrWA{} means that for any $G$, the map
\[A_{/F\times G} \to A_{/F}\]
in \Cat{} is in \scrW. Using the localization axiom on \scrW, one sees
that it is enough to check this for $G$ and object $a$ of $A$, in
which case $A_{/G}=A_{/a}$ has a final object and hence is
\scrW-aspheric by \ref{it:64.L3}, and the condition amounts to
$A_{/F\times a}$ being \scrW-aspheric, i.e. (with the terminology
introduced above), that $F\times a$ is \scrW-aspheric. Thus, an object
$F$ of \Ahat\pspage{170} is \scrWA-aspheric over $e_\Ahat$ if{f} for
any $a$ in $A$, $F\times a$ is \scrW-aspheric. We should beware that
for general $A$, this does not imply $F$ is \scrW-aspheric, nor is it
implies by it. We should remember that \scrW-asphericity of $F$ is an
``absolute notion'', namely is a property of the induced
\emph{category} $A_{/F}$ or equivalently, of the induced topos
$\Ahat_{/F}\simeq(A_{/F})\uphat$, whereas \scrWA-asphericity of $F$
over $e_\Ahat$ is a relative notion for the \emph{map} of categories
\[A_{/F} \to A\]
or equivalently, for the map of topoi $\Ahat_{/F}\to\Ahat$ (the
localization map with respect to the object $F$ of the topos
\Ahat). More generally, for a map $F\to G$ in \Ahat, the property for
this map of being \scrWA-aspheric is a property for the corresponding
map in \Cat
\[A_{/F} \to A_{/G},\]
namely the property we called \emph{weak \scrW-asphericity} yesterday
(page \ref{p:161}), as we stated then in the prop.\ page
\ref{p:162}. An equivalent way of expressing this is by saying that
for $F\to G$ to be \scrWA-aspheric, i.e., to be ``universally in
\scrWA'', it is enough to check this for base changes $G'\to G$ with
$G'$ an object $a$ in $A$, namely that the corresponding map
\[F\times_Ga\to a\]
in \Ahat{} should be in \scrWA (for any $a$ in $A$ and map $a\to G$),
which amounts to saying that $F\times_Ga$ is \scrW-aspheric for any
$a$ in $A$ and map $a\to G$.

With notations and terminology quite clear in mind, we may start
retelling once again the test category story!

\subsection{Total asphericity.}\label{subsec:65.A}
Before starting, just one important pre-test notion to recall, namely
total asphericity, summed up in the
\begin{propositionnum}\label{prop:65.1}
  The following conditions on $A$ are equivalent:
  \begin{description}
  \item[\namedlabel{it:65.A.i}{(i)}]
    The product in \Ahat{} of any two objects of $A$ is
    \scrW-aspheric.
  \item[\namedlabel{it:65.A.iprime}{(i')}]
    Every object in $A$ is \scrWA-aspheric over the final element $e_\Ahat$.
  \item[\namedlabel{it:65.A.ii}{(ii)}]
    The product of two \scrW-aspheric objects is again \scrW-aspheric.
  \item[\namedlabel{it:65.A.iiprime}{(ii')}]
    Any \scrW-aspheric object of \Ahat{} is \scrWA-aspheric over $e_\Ahat$.
  \end{description}
\end{propositionnum}

This is just a tautology, in terms of what was just said. Condition
\ref{it:65.A.i} is just the old condition \ref{it:31.T2} on test
categories\ldots
\begin{definitionnum}\label{def:65.1}
  If $A$ satisfies these conditions and moreover $A$ is
  \scrW-aspheric, \Ahat is called \emph{totally} \scrW-aspheric.
\end{definitionnum}
\begin{remark}
  In all cases when we have met with totally aspheric \Ahat, this
  condition \ref{it:65.A.i} was checked easily, because we were in one
  of\pspage{171} the two following cases:
  \begin{enumerate}[label=\alph*)]
  \item\label{it:65.A.a} $A$ stable under binary products.
  \item\label{it:65.A.b} The objects of $A$ are \emph{contractible}
    for the homotopy structure $h_\scrWA$ of \Ahat{} associated to \scrWA.
  \end{enumerate}
\end{remark}

In case \ref{it:65.A.b}, in the cases we've met, for checking
contractibility we even could get away with a homotopy interval
$\bI=(I,\delta_0,\delta_1)$ which is in $A$, namely we got
\bI-contractibility for all elements of $A$, and hence for the
products. All we've got to check then, to imply
$h_\scrWA$-contractibility of the objects $a\times b$, and hence their
\scrW-asphericity, is that $I$ itself is \scrWA-aspheric over
$e_\Ahat$, namely the products $I\times a$ are \scrW-aspheric. This
now has to be checked indeed some way or other -- I don't see any
general homotopy trick to reduce the checking still more. In case when
$A=\Simplex$ (standard simplices) say, and while still working with
usual weak equivalences $W_\Cat$, we checked asphericity of the
products $\Simplex_1\times\Simplex_n$ by using a Mayer-Vietories argument,
each product being viewed as obtained by gluing together a bunch of
\emph{representable} subobjects, which are necessarily \scrW-aspheric
therefore. The argument will go through for general \scrW, if we
assume \scrW{} satisfies the following condition, which we add to the
provisional list made yesterday (pages \ref{p:166}--\ref{p:167}) of
extra conditions which we may have to introduce for a basic localizer:
\begin{description}
\item[\namedlabel{it:65.Le}{L~e)}] \textbf{(Mayer-Vietoris axiom)}
  Let $C$ be in \Cat, let $C',C''$ be two full subcategories which are
  cribles (if it contains $a$ in $C$ and if $b\to a$, it contains
  $b$), and such that $\Ob C=\Ob C' \sor \Ob C''$. Assume $C',C''$ and
  $C' \sand C''$ are \scrW-aspheric, then so is $C$.
\end{description}

This condition of course is satisfied whenever \scrW{} is described in
terms of cohomological conditions, as envisioned yesterday (page
\ref{p:158}). We could elaborate on it and develop in this direction a
lot more encompassing conditions (``of \v Cech type'' we could say),
which will be satisfied by all such cohomologically defined basic
localizers. It would be fun to work out a set of ``minimal''
conditions such as \ref{it:65.Le} above, which would be enough to
imply all \v Cech-type conditions on a basic localizer. At first sight,
it isn't even obvious that \ref{it:65.Le} say isn't a consequence of
just the general conditions \ref{it:64.L1} to \ref{it:64.L4} on \scrW,
plus perhaps the fibration axiom \ref{it:64.L5} which looks very
strong. As long as we don't have any other example of basic
localizers than in terms of cohomology, it will be hard to tell!

\subsection{Weak \texorpdfstring{\scrW}{W}-test categories.}
\label{subsec:65.B}

\begin{definitionnum}\label{def:65.2}
  The\pspage{172} category $A$ is a weak \scrW-test category if it
  satisfies the conditions
  \begin{enumerate}[label=\alph*)]
  \item\label{it:65.B.a}
    $A$ is \scrW-aspheric.
  \item\label{it:65.B.b}
    The functor $i_A^*:\Cat\to\Ahat$ is model-preserving, i.e.,
    \begin{enumerate}[label=b\textsubscript{\arabic*})]
    \item\label{it:65.B.b1}
      $\scrW=(i_A^*)^{-1}(\scrWA)$ ($=f_A^{-1}(\scrW)$, where
      $f_A=i_Ai_A^*:\Cat\to\Cat$),
    \item\label{it:65.B.b2}
      the induced functor
      \[ \overline i_A^* : \HotW \eqdef \scrW^{-1}\Cat\to
      \scrWA^{-1}\Ahat\]
      is an equivalence.
    \end{enumerate}
  \end{enumerate}
\end{definitionnum}
\begin{propositionnum}\label{prop:65.2}
  The following conditions on $A$ are equivalent:
  \begin{description}
  \item[\namedlabel{it:65.B.i}{(i)}]
    $A$ is a weak \scrW-test category.
  \item[\namedlabel{it:65.B.ii}{(ii)}]
    The functors $i_A^*$ and $i_A$ are both model-preserving, the
    induced functors
    \[\begin{tikzcd}[cramped]
      \HotW\ar[r,bend left=10] & \scrWA^{-1}\Ahat \ar[l,bend left=10]
    \end{tikzcd}\]
    are equivalences quasi-inverse of each other, with adjunction
    morphisms in \HotW{} and in $\HotA\eqdef\scrWA^{-1}\Ahat$ deduced
    from the adjunction morphisms for the pair of adjoint functors $i_A,i_A^*$.
  \item[\namedlabel{it:65.B.iii}{(iii)}]
    The functor $i_A^*$ transforms maps in \scrW{} into maps in \scrWA
    \textup(i.e., $f_A=i_Ai_A^*$ transforms maps in \scrW{} into maps
    in \scrW\textup), and moreover $A$ is \scrW-aspheric.
  \item[\namedlabel{it:65.B.iiiprime}{(iii')}]
    Same as in \textup{\ref{it:65.B.iii}}, but restricting to maps
    $C\to e$, where $C$ in \Cat{} has a final object.
  \item[\namedlabel{it:65.B.iv}{(iv)}]
    The categories $f_A(C)=A_{/i_A^*(C)}$, where $C$ in \Cat{} has a
    final object, are \scrW-aspheric.
  \end{description}
\end{propositionnum}

The obvious implications are
\[ \text{\ref{it:65.B.ii}}
\Rightarrow \text{\ref{it:65.B.i}}
\Rightarrow \text{\ref{it:65.B.iii}}
\Rightarrow \text{\ref{it:65.B.iiiprime}}
\Rightarrow \text{\ref{it:65.B.iv}}\]
and the proof of \ref{it:65.B.iv} $\Rightarrow$ \ref{it:65.B.ii}
follows from an easy weak asphericity argument and general non-sense
on adjoint functors and localization (cf.\ page \ref{p:35} and prop.\
on page \ref{p:38}).
\begin{remark}
  In case \scrW{} is strongly saturated, and hence $A$ \scrW-aspheric
  just means that its image in \HotW{} is a final object, the
  condition of \scrW-asphericity of $A$ in \ref{it:65.B.iii} or in
  def.\ 2 can be restated, by saying that the endomorphism $\overline
  f_A$ of \HotW{} induced by $f_A$ transforms final object into final
  object -- which is a lot weaker than being an equivalence!
\end{remark}

\subsection{\texorpdfstring{\scrW}{W}-test categories.}
\label{subsec:65.C}
\begin{definitionnum}\label{def:65.3}
  The\pspage{173} category $A$ is a \emph{\scrW-test category} if it
  is a weak \scrW-test category, and if the localized categories
  $A_{/a}$ for $a$ in $A$ are equally weak \scrW-test categories. We
  say $A$ is a \emph{local \scrW-test category} if the localized
  categories $A_{/a}$ are weak \scrW-test categories.
\end{definitionnum}

Clearly, $A$ is a \scrW-test category if{f} if is a local \scrW-test
category, and moreover $A$ is \scrW-aspheric (as the categories
$A_{/a}$ are \scrW-aspheric by \ref{it:64.L3}). Also, $A$ is a local
\scrW-test category if{f} the functors $i_{A_{/a}}^*$ (for $a$ in $A$)
are model preserving.
\begin{propositionnum}\label{prop:65.3}
  The following conditions on $A$ are equivalent:
  \begin{enumerate}[label=(\roman*),font=\normalfont]
  \item\label{it:65.C.i}
    $A$ is a local \scrW-test category.
  \item\label{it:65.C.ii}
    The Lawvere element
    \[L_\Ahat = i_A^*(\Simplex_1)\]
    in \Ahat{} is \scrWA-aspheric over $e_\Ahat$, i.e., the products
    $a\times L_\Ahat$ for $a$ in $A$ are all \scrW-aspheric.
  \item\label{it:65.C.iii}
    There exists a separated interval $\bI=(I,\delta_0,\delta_1)$ in
    \Ahat{} \textup(i.e., an object endowed with two sections such
    that $\Ker(\delta_0,\delta_1)=\varnothing_\Ahat$\textup), such that
    $I$ be \scrWA-aspheric over $e_\Ahat$, i.e., all products $a\times
    I$ are \scrW-aspheric.
  \end{enumerate}
\end{propositionnum}

The obvious implications here are
\[ \text{\ref{it:65.C.i}}
\Rightarrow \text{\ref{it:65.C.ii}}
\Rightarrow \text{\ref{it:65.C.iii}}.\]
on the other hand \ref{it:65.C.iii} $\Rightarrow$ \ref{it:65.C.ii} by
the homotopy interval comparison lemma (p.\ \ref{p:60}), and finally 
\ref{it:65.C.i} $\Rightarrow$ \ref{it:65.C.i} by the criterion for
weak \scrW-test categories of prop.\ \ref{prop:65.2} \ref{it:65.B.iv}, using an
immediate homotopy argument (cf.\ page \ref{p:62}).
\begin{corollarynum}\label{cor:65.1}
  $A$ is a \scrW-test category if{f} it is \scrW-aspheric and
  satisfies \textup{\ref{it:65.C.ii}} or \textup{\ref{it:65.C.iii}} of
  proposition \textup{\ref{prop:65.3}} above.
\end{corollarynum}
\begin{remark}
  In the important case when \Ahat{} is totally \scrW-aspheric (cf.\
  prop.\ \ref{prop:65.1}, the asphericity condition on $L_\Ahat$ or on $I$ in
  prop.\ \ref{prop:65.3} is equivalent to just \scrW-asphericity of $L_\Ahat$ resp.\
  of $I$. In case \Ahat{} is even ``strictly totally \scrW-aspheric'',
  i.e., if moreover every ``non-empty'' object in \Ahat{} admits a
  section, then we've seen that $h_\scrWA=h_\Ahat$ (prop.\ page
  \ref{p:149}, which carries over to a general \scrW{} satisfying
  \ref{it:64.Laprime} of page \ref{p:166}, i.e., provided \HotW{}
  isn't equivalent to the final category, which case we may discard!),
  then condition \ref{it:65.C.ii} just means that the contractor
  $L_\Ahat$ is $0$-connected -- a condition which does not depend upon
  the choice of \scrW.
\end{remark}

\subsection{Strict \texorpdfstring{\scrW}{W}-test categories.}
\label{subsec:65.D}
\begin{propositionnum}\label{prop:65.4}
  The\pspage{174} following conditions on $A$ are equivalent:
  \begin{description}
  \item[\namedlabel{it:65.D.i}{(i)}]
    Both functors $i_A$ and $i_A^*$ are model preserving, moreover
    $i_A$ commutes to finite products ``modulo \scrW''.
  \item[\namedlabel{it:65.D.ii}{(ii)}]
    $A$ is a test category and \Ahat{} is totally \scrW-aspheric.
  \item[\namedlabel{it:65.D.iiprime}{(ii')}]
    $A$ is a weak test category and \Ahat{} is totally \scrW-aspheric.
  \item[\namedlabel{it:65.D.iii}{(iii)}]
    $A$ satisfies conditions \textup{\ref{it:31.T1}}
    \textup{\ref{it:31.T2}} \textup{\ref{it:31.T3}} of page
    \textup{\ref{p:39}}, with ``aspheric'' replaced by ``\scrW-aspheric''.
  \end{description}
\end{propositionnum}

This is not much more than a tautology in terms of what we have seen
before, as we'll get the obvious implications
\[ \text{\ref{it:65.D.i}}
\Rightarrow \text{\ref{it:65.D.iii}}
\Rightarrow \text{\ref{it:65.D.ii}}
\Rightarrow \text{\ref{it:65.D.iiprime}}
\Rightarrow \text{\ref{it:65.D.i}}.\]
\begin{definitionnum}\label{def:65.4}
  If $A$ satisfies the conditions above, it is called a \emph{strict
    \scrW-test category}.
\end{definitionnum}
\begin{remarks}
  \namedlabel{rem:65.D.1}{1})\enspace When we know that the canonical functor from \Cat{} to the
  localization \HotW{} commutes with binary products, then the
  exactness property mod \scrW{} in \ref{it:65.D.i} implies that the
  same holds for the canonical functor from \Ahat{} to its
  localization \HotA, and conversely if \scrW{} is known to be
  saturated.

  \namedlabel{rem:65.D.2}{2})\enspace In the case $\scrW=W_\Cat$ we've seen that condition
  \ref{it:31.T2} implies \ref{it:31.T1}, i.e., the conditions of
  prop.\ \ref{prop:65.1} imply $A$ is \scrW-aspheric, i.e.,
  \Ahat{} is totally \scrW-aspheric. The argument works for any
  \scrW{} defined by cohomological conditions of the type considered
  in yesterday's notes. To have it work for more general \scrW, we
  would have to introduce some \v Cech-type requirement on \scrW,
  compare page \ref{p:171}.

  \namedlabel{rem:65.D.3}{3})\enspace In the statement of the theorem page \ref{p:46}, similar
  to the proposition above, in \ref{it:33.i} no assumption is made on
  $i_A^*=j_A$ -- which I believe was an omission by hastiness -- it is
  by no means clear to me that we could dispense with it, and get away
  with an assumption on $i_A$ alone.
\end{remarks}

\subsection{Weak \texorpdfstring{\scrW}{W}-test functors and
  \texorpdfstring{\scrW}{W}-test functors.}\label{subsec:65.E}

Let\pspage{175}
\[(M,W) , \quad W\subset\Fl(M)\]
be a category endowed with a saturated set of arrows $W$, and
\begin{equation}
  \label{eq:65.4}
  i:A\to M\tag{4}
\end{equation}
a functor, hence a corresponding functor
\begin{equation}
  \label{eq:65.5}
  i^*: M \to \Ahat, \quad
  X\mapsto i^*(X)=(a\mapsto \Hom(i(a),X)).\tag{5}
\end{equation}
\begin{definitionnum}\label{def:65.5}
  The functor $i$ is called a \emph{weak \scrW-test functor}
  (with respect to the given $W\subset\Fl(M)$) if $A$ is
  a weak \scrW-test category and the functor $i^*$ is model-preserving
  (for $W$ and \scrWA), i.e., if $A$ satisfies the three conditions:
  \begin{enumerate}[label=\alph*)]
  \item\label{it:65.E.a}
    $i^*$ is model preserving,
  \item\label{it:65.E.b}
    $i_A^*: \Cat\to\Ahat$ is model preserving,
  \item\label{it:65.E.c}
    $A$ is \scrW-aspheric.
  \end{enumerate}
\end{definitionnum}

The conditions \ref{it:65.E.b} and \ref{it:65.E.c}, namely that $A$ be
a weak \scrW-test category, do not depend of course upon $M$, and it
may seem strange in the definition not to have simply asked beforehand
that $A$ satisfy this preliminary condition -- i.e., reduce to the
case when we start with a weak \scrW-test category $A$. The reason for
not doing so is that we'll find below handy criteria for all three
conditions to hold, without assuming beforehand $A$ to be a weak
\scrW-test category.

As \ref{it:65.E.b} and \ref{it:65.E.c} imply that
\[ i_A: \Ahat\to\Cat\]
is model-preserving too, condition \ref{it:65.E.a} above can be
replaced by the condition
\begin{enumerate}[label=\alph*')]
\item\label{it:65.E.aprime}
  The composition
  \[f_i=i_Ai^*: M\to\Cat\]
  is model-preserving (for the pair $W,\scrW$).
\end{enumerate}

Of course, as conditions \ref{it:65.E.b}, \ref{it:65.E.c} imply that
$(\Ahat,\scrWA)$ is a modelizer (with respect to \scrW), the condition
\ref{it:65.E.a} will imply $(M,W)$ is a modelizer too.

We recall the condition for $i^*$ to be model-preserving decomposes
into two:
\begin{enumerate}[label=a\textsubscript{\arabic*})]
\item\label{it:65.E.a1}
  $W=(i^*)^{-1}(\scrWA)\quad (=f_i^{-1}(\scrW))$,
\item\label{it:65.E.a2}
  The\pspage{176} functor $\overline{i^*}$ induced by $i^*$ on the
  localizations (which exists because of \ref{it:65.E.a1})
  \[W^{-1}M \to \HotA\eqdef \scrWA^{-1}\Ahat\]
  is an equivalence.
\end{enumerate}
\begin{definitionnum}\label{def:65.6}
  The functor $i$ is called a \emph{\scrW-test functor} if this
  functor \emph{and} the induced functors $i_{/a}:A_{/a}\to M$ (for
  $a$ in $A$) are weak \scrW-test functors.
\end{definitionnum}

In view of the definition \ref{def:65.3}, this amounts to the two
conditions:
\begin{enumerate}[label=\alph*)]
\item\label{cond:65.E.a}
  $A$ is a \scrW-test category, i.e., the functors $i_A^*$ and
  $i_{A_{/a}}^*$ are model-preserving and $A$ is \scrW-aspheric,
\item\label{cond:65.E.b}
  the functors $i^*$ and $(i_{/a})^*$ from $M$ into the categories
  \Ahat{} and $(A_{/a})\uphat\to\Ahat_{/a}$ are model-preserving (for
  $W$ and \scrWA{} resp.\ $\scrW_{A_{/a}}$).
\end{enumerate}
\begin{example}
  Consider the canonical functor induced by $i_A$
  \[i_A^0: A\to \Cat, \quad\text{\Cat{} endowed with \scrW,}\]
  this functor is a weak \scrW-test functor (resp.\ a \scrW-test
  functor) if{f} $A$ is a weak \scrW-test category (resp.\ a
  \scrW-test category).
\end{example}

These two definitions are pretty formal indeed. Their justification is
mainly in the two theorems below.

\emph{We assume from now on that the basic localizer \scrW{} satisfies
  the fibration axiom} \ref{it:64.L5} of page \ref{p:164}. Also, we
recall that an object $X$ in \Cat{} is \emph{contractible} (for the
canonical homotopy structure of \Cat) if{f} $X$ is non-empty and the
category $\bHom(X,X)$ is connected -- indeed it is enough even that
$\id_X$ belong to the same connected component as some \emph{constant}
map from $X$ into itself. This condition is satisfied for instance if
$X$ has a final or an initial object.
\begin{theoremnum}\label{thm:65.1}
  We assume that $M=\Cat$, $W=\scrW$, i.e., we've got a functor
  \[i:A\to\Cat, \quad\text{\Cat{} endowed with \scrW,}\]
  and we assume that for any $a$ in $A$\kern1pt, $i(a)$ is contractible
  \textup(cf.\ above\textup), i.e., that $i$ factors through the full
  subcategory $\Cat_{\mathrm{cont}}$ of contractible objects of
  \Cat. The following conditions are equivalent:
  \begin{description}
  \item[\namedlabel{it:65.E.i}{(i)}]
    $i$ is a \scrW-test functor \textup(def.\ \textup{\ref{def:65.6})}.
  \item[\namedlabel{it:65.E.iprime}{(i')}]
    For any $a$ in $A$\kern1pt, the induced functor $i_{/a}:A_{/a}\to\Cat$ is
    a weak \scrW-test functor, and moreover $A$ is \scrW-aspheric.
  \item[\namedlabel{it:65.E.ii}{(ii)}]
    $i^*(\Simplex_1)$ is \scrWA-aspheric over $e_\Ahat$, i.e., the
    products $a\times i^*(\Simplex_1)$ in \Ahat{} are \scrW-aspheric,
    for any $a$ in $A$\kern1pt, and $A$ is \scrW-aspheric.
  \end{description}
\end{theoremnum}

The\pspage{177} obvious implications here are
\[ \text{\ref{it:65.E.i}}
\Rightarrow \text{\ref{it:65.E.iprime}}
\Rightarrow \text{\ref{it:65.E.ii}},\]
for the last implication we only make use, besides $A$ being
\scrW-aspheric, that the functors $(i_{/a})^*$ transform the
projection $\Simplex_1\to e$ in \Cat, which is in \scrW{} by
\ref{it:64.L3prime}, into a map in $\scrW_{A_{/a}}$, i.e., that the
corresponding map in \Cat
\[A_{/a\times i^*(\Simplex_1)} \to A_{/a}\]
be in \scrW, which by the final object axiom implying that $A_{/a}\to
e$ is in \scrW, amounts to demanding that the left-hand side is
\scrW-aspheric, i.e., $a\times i^*(\Simplex_1)$ \scrWA-aspheric.

So we are left with proving \ref{it:65.E.ii} $\Rightarrow$
\ref{it:65.E.i}. By the criterion \ref{it:65.C.iii} of prop.\
\ref{prop:65.3} we know already (assuming \ref{it:65.E.ii}) that $A$
is a local \scrW-test category, hence a \scrW-test category as $A$ is
\scrW-aspheric (cor.\ \ref{cor:65.1}); indeed we can use
$I=i^*(\Simplex_1)$ as a \scrWA-aspheric interval, using the two
canonical sections deduced from the canonical sections of
$\Simplex_1$. The fact that these are disjoint follows from the fact
that $i(a)$ non-empty for any $a$ in $A$ -- we did not yet have to use
the contractibility assumption on the categories $i(a)$. Thus, we are
reduced to proving that $i^*$ is model-preserving -- the same will
then hold for the functors $i_{/a}$ (as required in part
\ref{cond:65.E.b} in def.\ \ref{def:65.6}), as the assumption
\ref{it:65.E.ii} is clearly stable under restriction to the categories
$A_{/a}$. As we know already that $i_A$ is model-preserving (prop.\
\ref{prop:65.2} \ref{it:65.B.i} $\Rightarrow$ \ref{it:65.B.ii}), all
we have to do is to prove the composition $f_i=i_Ai^*$ is
model-preserving. But this was proved yesterday in the key lemma of
page \ref{p:156}. We're through!

\begin{remark}
  The presentation will be maybe a little more elegant, if we
  complement the definition of a \scrW-test functor by the definition
  of a \emph{local \scrW-test functor}, by which we mean that the
  induced functors $i_{/a}:A_{/a}\to M$ are weak test functors, period
  -- which means also that the following conditions hold:
  \begin{enumerate}[label=\alph*)]
  \item\label{it:65.E.rem.a}
    $A$ is a local \scrW-test category (def.\ \ref{def:65.3}), i.e.,
    the functors $(i_{/a})^*: M\to(A_{/a})\uphat$ are all
    model-preserving;
  \item\label{it:65.E.rem.b}
    the functors $(i_{/a})^*$ (for $a$ in $A$) are model-preserving.
  \end{enumerate}
\end{remark}

Thus, it is clear that if $i$ is a \scrW-test functor, it is a local
\scrW-test functor such that moreover $A$ is \scrW-aspheric. The
converse isn't clear in general, because it isn't clear that if $A$ is
a \scrW-test category and moreover all functors $(i_{/a})^*$ are
modelizing, then $i^*$ is equally modelizing. The criterion
\ref{it:65.E.iprime} of theorem \ref{thm:65.1} shows however that this
is so in the case when $(M,W)=(\Cat,\scrW)$, and when we assume
moreover the objects $i(a)$ contractible. We could now reformulate the
theorem\pspage{178} as a twofold statement:
\begin{corollary}
  Under the assumptions of theorem \ref{thm:65.1}, $i$ is a
  \emph{local \scrW-test functor} \textup(i.e., all functors
  $i_{/a}: A_{/a}\to\Cat$ are weak \scrW-test functors, or
  equivalently the functors $(i_{/a})^*$ and $(i_{A_{/a}})^*$ from
  $\Cat \to (A_{/a})\uphat$ are all model-preserving\textup) if{f}
  $i^*(\Simplex_1)$ is \scrWA-aspheric over $e_\Ahat$, i.e., the
  products $a\times i^*(\Simplex_1)$ in \Ahat{} are \scrW-aspheric. When
  this condition is satisfied, in order for $i$ to be a \scrW-test
  functor, namely for $i^*$ to be equally model-preserving, it is
  n.s.\ that $A$ be \scrW-aspheric.
\end{corollary}

\subsection{\texorpdfstring{\scrW}{W}-test functors
  \texorpdfstring{$A\to\Cat$}{A->(Cat)} of strict
  \texorpdfstring{\scrW}{W}-test categories.}
\label{subsec:65.F}
Let again
\[i:A\to\Cat,\quad\text{\Cat{} endowed with \scrW,}\]
be a functor such that the objects $i(a)$ be contractible, we assume
now moreover that \Ahat{} is \emph{totally \scrW-aspheric} (def.\
\ref{def:65.1}), which implies $A$ is \scrW-aspheric. Thus, by the
corollary above $i$ is a test functor if{f} it is a local \scrW-test
functor, and by the criterion \ref{it:65.B.iv} of prop.\
\ref{prop:65.2} (with $C=\Simplex_1$) we see it amounts to the same that
$i$ be a weak \scrW-test functor. (Here we use the assumption of total
\scrW-asphericity, which implies that if $i^*(\Simplex_1)$ is
\scrW-aspheric, it is even \scrWA-aspheric over $e_\Ahat$.) Thus, the
three variants of the test-functor notion coincide in the present
case. With this in mind, we can now state what seems to me the main
result of our reflections so far, at any rate the most suggestive
reformulation of theorem \ref{thm:65.1} in the present case:
\begin{theoremnum}\label{thm:65.2}
  With the assumptions above \textup(\Ahat{} totally \scrW-aspheric
  and the objects $i(a)$ in \Cat{} contractible), the following
  conditions on the functor $i:A\to\Cat$ are equivalent:
  \begin{enumerate}[label=(\roman*),font=\normalfont]
  \item\label{it:65.F.i}
    $i$ is a \scrW-test functor.
  \item\label{it:65.F.ii}
    $i^*: \Cat\to\Ahat$ is model-preserving, i.e., for any map $f$ in
    \Cat, $f$ is in \scrW{} if{f} $i^*(f)$ is in \scrWA{}
    \textup(i.e., if{f} $i_Ai^*(f)$ is in \scrW\textup), and moreover
    the induced functor
    \[\overline{i^*}:\HotW\to\HotA\eqdef\scrWA^{-1}\Ahat\]
    is an equivalence.
  \item\label{it:65.F.iii}
    The functor above exists, i.e., $f$ in \scrW{} implies $i^*(f)$ in
    \scrWA, i.e., $i_Ai^*(f)$ in \scrW.
  \item\label{it:65.F.iv}
    The functor $i^*$ transforms \scrW-aspheric objects into
    \scrWA-objects \textup(i.e., the condition in
    \textup{\ref{it:65.F.iii}} is satisfied for maps $C\to e$ in \scrW\textup).
  \item\label{it:65.F.v}
    The functor $i^*$ transforms contractible objects of \Cat{}
    into\pspage{179} objects of \Ahat, contractible for the homotopy
    structure $h_\scrWA$ associated to \scrWA{} -- or equivalently,
    $i^*$ is a morphism of homotopy structures \textup(cf.\ definition
    on page \textup{\ref{p:134})}.
  \item\label{it:65.F.vi}
    The functor $i^*$ transforms the projection $\Simplex_1\to e$ into a
    map in \scrW, or equivalently \textup(as $A$ is
    \scrW-aspheric\textup) $i^*(\Simplex_1)$ is \scrW-aspheric.
  \end{enumerate}
\end{theoremnum}

We have the trivial implications
\[ \text{\ref{it:65.F.i}}
\Rightarrow \text{\ref{it:65.F.ii}}
\Rightarrow \text{\ref{it:65.F.iii}}
\Rightarrow \text{\ref{it:65.F.iv}}
\Rightarrow \text{\ref{it:65.F.vi}}
\Leftarrow \text{\ref{it:65.F.v}},\]
where the implication \ref{it:65.F.v} $\Rightarrow$ \ref{it:65.F.vi}
is in fact an equivalence, due to the fact that the contractibility
structure on \Cat{} is defined in terms of $\Simplex_1$ as a generating
contractor, and that the assumption $i^*(\Simplex_1)$ \scrW-aspheric
implies already that it is \scrWA-aspheric over $e_\Ahat$ (because $A$
is totally \scrWA-aspheric), and hence contractible as it is a
contractor and $h_\scrWA$ is defined in terms of ``weak homotopy
intervals'' which are \scrWA-aspheric over $e_\Ahat$. Thus, the only
delicate implication is \ref{it:65.F.vi} $\Rightarrow$
\ref{it:65.F.i}, which however follows from theorem \ref{thm:65.1}
\ref{it:65.E.ii} $\Rightarrow$ \ref{it:65.E.i}.

We got the longed-for ``key result'' in the end!

\bigbreak

\presectionfill\ondate{9.4.}\pspage{180}\par

\hangsection{Revising \texorpdfstring{\textup(}{(}and
  fixing?\texorpdfstring{\textup)}{)} terminology again.}
\label{sec:66}%
After writing down nicely, in the end, that long promised key result,
I thought the next thing would be to pull myself up by my bootstraps
getting the similar result first for test functors $A\to\Bhat$ with
values in an elementary modelizer \Bhat, and then for general
``canonical'' modelizers $(M,W)$. So I did a little scratchwork
pondering along those lines, before resuming the typewriter-engined
work. What then turned out, it seems, is that there wasn't any need at
all to pass through the particular case $M=\Cat$ and the somewhat
painstaking analysis of our three-step diagram on page
\ref{p:96}. Finally, the most useful result of all the eighty pages
grinding, since that point, is by no means the so-called key result,
as I anticipated -- the day after I finally wrote it down, it was
already looking rather ``\'etriqu\'e''\scrcomment{``\'etriqu\'e'' can be
  translated as ``narrow-minded'' here, I think.} -- why
all this fuss about the special case of test functors with values in
\Cat! The main result has been finally more psychological than
technical -- namely drawing attention, in the long last, to the key
role of contractible objects and, more specifically, of the
contractibility structure associated to modelizers $(M,W)$, suggesting
that the localizer $W$ should, conversely, be describable in terms of
the homotopy structure $h_W$. This was point \ref{it:56.d} in the
``provisional plan of work'' contemplated earlier this week (page
\ref{p:138}) -- by then I had already the feeling this approach via
\ref{it:56.d} would turn out to be the most ``expedient'' one -- but
it was by then next to impossible for me to keep pushing off still
more the approach \ref{it:56.b} via test functors with values in \Cat,
which I finally carried through. One point which wasn't wholly clear
yet that day, as it is now, is the crucial role played by the
circumstance that for the really nice modelizers (surely for those I'm
going to call ``canonical'' in the end), the associated homotopy
structure is indeed a \emph{contractibility} structure. Here, as so
often in mathematics (and even outside of mathematics\ldots), the main
thing to dig out and discover is where the emphasis belongs -- which
are the really essential facts or notions or features within a given
context, and which are accessory, namely, which will follow suit by
themselves. It took a while before I would listen to what the things I
was in were insistently telling me. It finally got through I feel, and
I believe that from this point on the whole modelizing story is going
to go through extremely smoothly.

Before starting with the work, just some retrospective, somewhat more
technical comments, afterthoughts rather I should say. First of all, I
am not so happy after all with the terminological review a few days
ago\pspage{181} (pages \ref{p:159}--\ref{p:163}), and notably the use
of the word ``aspheric'' in the generalization ``$W$-aspheric'' map
(in a category $M$ endowed with a saturated set of arrows $W$) --
which then practically obliged me, when working in \Cat, to call
``weakly aspheric'' a functor $C'\to C$ which spontaneously I surely
would like to call simply ``aspheric'' -- and as a matter of fact, it
turned out I couldn't force myself to add a ``weakly'' before as I
decreed I should -- or if I did, it was against a very strong feeling
of inappropriateness. That decree precisely is an excellent
illustration of loosing view of where the main emphasis belongs, which
I would like now to make very clear to myself.

In all this work the underlying motivation or inspiration is
geometrico-topological, and expressed technically quite accurately by
the notion of a topos and of maps (or ``morphisms'') of topoi, and the
wealth of geometric and algebraic intuitions which have developed
around these. One main point here is that topoi may be viewed as
\emph{the} natural common generalization of both topological spaces
(the conventional support for so-called topological intuition), and of
(small) categories, where the latter may be viewed as the ideal purely
algebraic objects carrying topological information, including all the
conventional homology and homotopy invariants. This being so, in a
context where working with small categories as ``spaces'', the main
emphasis in choice of terminology should surely be in stressing
throughout, through the very wording, the essential identity between
situations involving categories, and corresponding situations
involving topological spaces or topoi. Thus, it has been about twenty
years now that the needs for developing \'etale cohomology have told
me a handful of basic asphericity and acyclicity properties for a map
of topoi (which apparently have not yet been assimilated by
topologists, in the context of maps of topological spaces\ldots),
including the condition for such a map to be aspheric. This was
recalled earlier (page \ref{p:37}), and the corresponding notion for a
functor $f:C'\to C$ was introduced. The name ``\emph{aspheric map}''
of topoi, or of topological spaces, or of categories, is here a
perfectly suggestive one. As the notion itself is visibly a basic one,
there should be no question whatever to change the name and replace it
say by ``weakly aspheric'', whereas the notion is surely quite a
strong one, and doesn't deserve such minimizing qualification! There
is indeed a stronger notion, which in the context of topological
spaces or \'etale cohomology of schemes reduces to the previous one in
the particular case of a map which is supposed \emph{proper}. This
condition could be expressed by saying that for \emph{any} base-change
$Y'\to Y$\pspage{182} for the map $f:X\to Y$ (at least any base-change
within the given context, namely either spaces or schemes with \'etale
topology), the corresponding map
\[ f': X'=X\times_YY'\to Y'\]
is a weak equivalence, or what amounts to the same, that for any $Y'$
the corresponding map is aspheric. This property, if a name is needed
indeed, would properly be called ``\emph{universally
  aspheric}''. Thus, in \Cat{} a map $f:X\to Y$ will be called either
aspheric, or universally aspheric, when for any base-change of the
special type $Y_{/y}\to Y$, namely ``localization'' in the first case,
or any base-change whatever $Y'\to Y$ in the second case, the
corresponding map $f'$ is a weak equivalence. On the other hand, if
$Y$ is just the final object $e$ of \Cat, it turns out the two notions for
$X$ (being ``aspheric over $e$'' and being ``universally aspheric over
$e$'') coincide, and just mean that $X\to e$ is a weak equivalence. In
accordance with the use which has been prevalent for a long time in
the context of spaces, such an object will be call simply an
\emph{aspheric object} -- which means that the corresponding topos is
aspheric (namely has ``trivial'' cohomology invariants, and hence
trivial homotopy invariants of any kind\ldots).

In case the notion of weak equivalence is replaced by a basic
localizer $\scrW\subset\Fl(\Cat)$, there is no reason whatever to
change anything in this terminology -- except that, if need, we will
add the qualifying \scrW, and speak of \emph{\scrW-aspheric} or
\emph{universally \scrW-aspheric maps} in \Cat, as well of
\emph{\scrW-aspheric objects} of \Cat.

What about terminology for maps and objects within a category \Ahat?
Here the emphasis should be of course perfect coherence with the
terminology just used in \Cat. An object $F$ of \Ahat{} should always
be sensed in terms of the induced topos $\Ahat_{/F}\simeq
(A_{/F})\uphat$, or what amounts to the same, in terms of the
corresponding object $A_{/F}$ in \Cat, which will imply that
``\emph{$F$ is aspheric}'' cannot possibly mean anything else but
$A_{/F}$ is aspheric as an object of \Cat; the same if qualifying by a
\scrW{} -- $F$ is called \emph{\scrW-aspheric} if $A_{/F}$ is a
\scrW-aspheric object of \Cat. Similarly for maps -- thus $f: F\to G$
will be called a weak equivalence, if the corresponding map for the
induced topoi is a weak equivalence, or equivalently, if the
corresponding map in \Cat
\[A_{/F} \to A_{/G}\]
is a weak equivalence. When a \scrW{} is given, we would say instead
(if confusion may arise) that $f$ is a \emph{\scrW-equivalence}. The
map will be called\pspage{183} \emph{aspheric}, or
\emph{\scrW-aspheric}, if the corresponding map in \Cat{} is. It turns
out that (because of the localization axiom on \scrW) this is
equivalent with $f$ being ``universally a \scrW-equivalence'', i.e.,
$f$ being ``universally in \scrWA'', namely that for any base-change
$G'\to G$ in \Ahat, the corresponding map in \Ahat
\[f' : F\times_GG'\to G'\]
be in \scrWA, i.e., be a \scrW-equivalence. Of course, when this
condition is satisfied, then for any base change, $f'$ will be, not
only a \scrW-equivalence, but even \scrW-aspheric -- thus we can say
that $f$ is ``universally \scrW-aspheric'' -- where ``universally''
refers to \emph{base change in} \Ahat. This of course does not mean
(and here one has to be slightly cautious) that the corresponding map
in \Cat{} is universally \scrW-aspheric (which refers to arbitrary
\emph{base change in} \Cat). But this apparent incoherence is of no
practical importance as far as terminology goes, as the work
``\scrW-aspheric map in \Ahat'' is wholly adequate and sufficient for
naming the notion, without any need to replace it by the more
complicated and ambiguous name ``universally \scrW-aspheric'', which
therefore will never be used. We even could rule out the formal
incoherence, by using the words \scrWA-equivalence, \scrWA-aspheric
maps (which are even universally \scrWA-aspheric maps, without any
ambiguity any longer), as well as \scrWA-aspheric objects -- replacing
throughout \scrW{} by \scrWA. In practical terms, I think that when
working consistently with a single given \scrW, we'll soon enough drop
anyhow both \scrW{} and \scrWA{} in the terminology and notations!

A last point which deserves some caution, is that for general $A$,
there is no implication between the two asphericity properties of an
object $F$ of \Ahat, namely of $F$ being \scrW-aspheric (i.e., the
object $A_{/F}$ of \Cat{} being \scrW-aspheric, i.e., the map
\begin{equation}
  \label{eq:66.*}
  A_{/F} \to e\tag{*}
\end{equation}
in \Cat{} being in \scrW), and the property that $F\to e$ be
\scrW-aspheric, namely that map
\begin{equation}
  \label{eq:66.starstar}
  A_{/F} \to A_{/e}=A\tag{**}
\end{equation}
in \Cat{} being aspheric (which also means that the products $F\times
a$ for $a$ in $A$ are \scrW-aspheric objects of \Ahat, i.e., the
categories $A_{/F\times a}$ are \scrW-aspheric, i.e., the maps
\[ A_{/F\times a} \to e\]
in \Cat{} are \scrW-aspheric. A third related notion, weaker than the
last one\pspage{184} is the property that $F\to e$ be a
\scrW-equivalence, which also means that the map
\eqref{eq:66.starstar} in \Cat{} is a \scrW-equivalence, i.e., is in
\scrW. If $A$ is \scrW-aspheric, this third notion however reduces to
the first one, namely $F$ to be \scrW-aspheric.

\emph{These terminological conventions, in the all-important cases of
  \Cat{} and categories of the type \Ahat, should be viewed as the
  basic ones} and there should be no question whatever to requestion
them, because say of the need we are in to devise a terminology,
applicable to the general case of a category $M$ endowed with any
saturated set of maps $W\subset\Fl(M)$ (which are being thought of as
still more general substitutes of ``weak equivalences''). This shows
at once that we will have to renounce to the name of ``$W$-aspheric''
which we have used so far, in order to designate maps which are
``universally in $W$''; indeed, this contradicts the use we are making
of this word, in the case of \Cat. The whole trouble came from this
inappropriate terminology, which slipped in while thinking of the
\Ahat{} analogy, and forgetting about the still more basic \Cat! The
mistake is a course one indeed, and quite easy to correct -- \emph{we
  better refrain altogether from using the word ``aspheric'' in the
  context of a general pair $(M,W)$}, and rather speak of maps which
are ``\emph{universally in $W$}'' or ``\emph{universal
  $W$-equivalences}'', which is indeed more suggestive, and does not
carry any ambiguity. The notion of ``$W$-aspheric map'' should be
reserved to the case when, among all possible base-change maps $Y'\to
Y$ in $M$, we can sort out some which we may view as ``localizing
maps'' -- all maps I'd think in cases $M$ is a topos, and pretty few
ones when $M=\Cat$. As for qualifying objects of $M$, we'll just be
specific in stating properties of the projection $X\to e_M$ -- such as
being a $W$-equivalence, or a universal $W$-equivalence, or a
homotopism for $h_W$ (in which case the name ``contractible object''
is adequate indeed). It may be convenient, when we got a $W$, to
denote by
\[\UW\subset W\subset \Fl(M)\]
the corresponding set of maps which are universally in $W$, a
property which then can be abbreviated into the simple notation
\[ f\in \UW\]
or ``$f$ is in \UW. It should be noted that \UW{} contains all
invertible maps and is stable by composition, but \emph{it need not be
  saturated}, thus $f$ and $fg$ may be in \UW{} without $g$ being so.

This\pspage{185} terminological digression was of a more essential
nature, as a matter of fact, than merely technical. There is still
another correction I want to make with terminology introduced earlier,
namely with the name of a ``contractor'' I used for intervals endowed
with a suitable composition law (page \ref{p:120}). The name in itself
seems to me quite appropriate, however I have now a notion in reserve
which seems to me a lot more important still, a reinforcement it turns
out of the notion of a strict test category -- and which I really
would like to call a \emph{contractor}. I couldn't think of any more
appropriate name -- thus I better change the previous terminology --
sorry! -- and call those nice intervals ``multiplicative intervals'',
thus referring to the composition law as a ``multiplication'' (with
left unit and left zero element). The name which first slipped into
the typewriter, when it occurred that a name was desirable, was not
``contractor'' by the way but ``intersector'', as I was thinking of
the examples I had met so far, when composition laws were defined in
terms of intersections and were idempotent. But this doesn't square
too well with the example of contractors $\bHom(X,X)$, when $X$ is an
object which has a section -- and this example turns out as equally
significant.

One last comment is about the ``\v Cech type'' condition
\ref{it:65.Le} on the basic localizer \scrW, introduced two days ago
(page \ref{p:171}). As giving a ``crible'' in a category amounts to
the same as giving a map
\[ C\to \Simplex_1\]
(by taking the inverse image of the source-object $\{0\}$ of
$\Simplex_1$), and therefore giving two such amounts to a map from $C$
into $\Simplex_1\times\Simplex_1$, we see that the situation when $C$ is
the union of two cribles is expressed equivalently by giving a functor
from $C$ into the subcategory
\[C_0 = \Biggl(
  \begin{tikzcd}[baseline=(O.base),row sep=0pt]
    & b \\ |[alias=O]| a \ar[ur]\ar[dr] & \\ & c
  \end{tikzcd}\Biggr)\]
of $\Simplex_1\times\Simplex_1$ (dual to the barycentric subdivision of
$\Simplex_1$). The asphericity conditions on $C'$, $C''$ and $C' \sand
C''$ then just mean that \emph{this functor is \scrW-aspheric}, which by
the localization axiom \ref{it:64.L4} implies that the functor itself
is a weak equivalence. Thus (by saturation), the conclusion that $C$
should be \scrW-aspheric, just amounts to the following condition,
which is in a way the ``universal'' special case when $C=C_0$, and
$C',C''$ are the two copies of $\Simplex_1$ contained in $C_0$:
\begin{description}
\item[\namedlabel{it:66.Leprime}{L~e')}]
  The category $C_0$ above is \scrW-aspheric.
\end{description}

If\pspage{186} we now look upon the projection map of $C_0$ upon one
factor $\Simplex_1$ (carrying $a$ and $b$ into $\{0\}$ and $c$ into
$\{1\}$), we get a functor which is fibering, and whose fibers are
$\Simplex_1$ and $\Simplex_0$, which are \scrW-aspheric. Hence \emph{the
  fibration axiom \textup{\ref{it:64.L5}} on \scrW{} implies the
  Mayer-Vietoris axiom} \ref{it:65.Le} (page \ref{p:171}). This
argument rather convinces me that the fibration axiom should be strong
enough to imply all \v Cech-type \scrW-asphericity criteria which one
may devise (provided of course they are reasonable, namely hold for
ordinary weak equivalences!). More and more, it seems that the basic 
requirements to make upon a basic localizer, which will imply maybe
all others, are \ref{it:64.L1} (saturation), \ref{it:64.L3prime} (the
``standard interval axiom'', namely $\Simplex_1$ is \scrW-aspheric), and
the powerful fibration axiom \ref{it:64.L5}. This brings to my mind
though the condition \ref{it:64.La} of page \ref{p:165}, namely that
$f\in\scrW$ implies $\piz(f)$ bijective, which wasn't needed really
for the famous ``key result'' I was then after, and which for this
reason I then was looking at almost as something accessory! I now do
feel though that it is quite an essential requirement, even though we
made no formal use of it (except very incidentally, in the slightly
weaker form \ref{it:64.Laprime}, which just means that \HotW{} isn't
equivalent to the final category). I would therefore add it now to the
list of really basic requirements on a ``basic localizer'', and
rebaptize it therefore as \ref{it:66.L6}, namely:
\begin{description}
\item[\namedlabel{it:66.L6}{L~6)}] \textbf{(Connectedness axiom)}
  $f\in\scrW$ implies $\piz(f)$ bijective, i.e., the functor
  $\piz:\Cat\to\Sets$ factors through $\scrW^{-1}\Cat=\HotW$ to give
  rise to a functor
  \[\piz:\HotW\to\Sets.\]
\end{description}
This, as was recalled on page \ref{p:166}, is more than needed to
imply
\[h_\scrW=h_\Cat,\]
namely the homotopy structure in \Cat{} associated to \scrW{} (in
terms of \scrW-aspheric intervals) is just the canonical homotopy
structure (defined in terms of $0$-connected intervals), which is also
the homotopy structure defined by the single ``basic''
(multiplicative) interval $\Simplex_1$.

\bigbreak
\presectionfill\ondate{11.4.}\par

I forget to clear up still another point of terminology -- namely
about ``weak homotopy intervals'' -- it turns out finally we never
quite came around defining what ``homotopy intervals'' which aren't
weak should be! The situation is very silly indeed - so henceforth
I'll just\pspage{187} drop the qualificative ``weak'' -- thus from now
on a ``\emph{homotopy interval}'' (with respect to a given homotopy
structure $h$ in a category $M$) is just an interval whose end-point
sections $\delta_0,\delta_1$ are homotopic. In case $h=h_W$, where $W$
is a given saturated set of arrows in $M$, the notion we get is a lot
wider than the notion of a homotopy interval (with respect to $W$)
introduced earlier (page \ref{p:132}, and which we scarcely ever used
it seems, so much so that I even forgot till this very minute it had
been introduced formally), where we were restricting to intervals for
which $I\to e$ is universally in $W$, i.e., is in \UW{} (we may call
such objects simply \UW-objects). Anyhow, it seems that so far, the
only property of such intervals we kept using from the beginning is
the one shared with all homotopy intervals in the wider sense I am now
advocating. There is just one noteworthy extra property which is
sometimes of importance, especially in the characterization of test
categories, namely the property $\Ker(\delta_0,\delta_1)=\varnothing_M$;
this was referred to earlier by the name ``separated interval'' --
which however may lead to confusion when for objects of $M$ we have
(independently of homotopy notions) a notion of separation. Therefore,
we better speak about \emph{separating intervals} as those for which
$\Ker(\delta_0,\delta_1)=\varnothing_M$ (initial object in $M$), hence a
notion of \emph{separating homotopy interval} (with respect to a given
homotopy structure $h$, or with respect to a given saturated $W$,
giving rise to $h_W$).



\chapter{Asphericity structures and canonical modelizers}
\label{ch:IV}

\presectionfill\alsoondate{11.4.}\pspage{188}\par

\hangsection[Setting out for the asphericity game again: variance of
\dots]{Setting out for the asphericity game again: variance of the
  category \texorpdfstring{\HotA}{(Hot-A)}, for arbitrary small
  category $A$ and aspheric functors.}\label{sec:67}%
A little more pondering and scribbling finally seems to show that the
real key for an understanding of modelizers isn't really the notion of
contractibility, but rather the notion of aspheric objects (besides,
of course, the notion of weak equivalence). At the same time it
appears that the notion of an \emph{aspheric map} in \Cat, more
specifically of a \scrW-aspheric ``map'' (i.e., a functor between
small categories) is a lot more important than being just a highly
expedient technical convenience, as it has been so far -- it is indeed
one of the basic notions of the theory of modelizers we got into. As a
matter of fact, I should have known this for a number of weeks
already, ever since I did some scribbling about the plausible notion
of ``morphism'' between test-categories (as well as their weak and
strong variants), and readily convinced myself that the natural
``morphisms'' here were nothing else but the aspheric functors between
those categories. I kind of forgot about this, as it didn't seem too
urgent to start moving around the category I was working with. If I
had been a little more systematic in grinding through the usual
functorialities, as soon as a significant notion (such as the various
test notions) appears, I presumably would have hit upon the crucial
point about modelizers and so-called ``asphericity structures'' a lot
sooner, without going through the long-winded detour of homotopy
structures, and the still extremely special types of test functors
suggested by the contractibility assumptions. However, I believe that
most of the work I went through, although irrelevant for the
``asphericity story'' itself, will still be useful, especially when it
comes to pinpointing the so-called ``canonical'' modelizers, whose
modelizing structure is intrinsically determined by the category
structure.

First thing now which we've to do is to have a closer look at the
meaning of asphericity for a functor between small categories. There
is no reason whatever to put any restrictions on these categories
besides smallness (namely the cardinals of the sets of objects and
arrows being in the ``universe'' we are implicitly working in, or even
more stringently still, $A$ and $B$ to be in \Cat{} namely to be
objects of that universe). Thus, we will not assume $A$ and $B$ to be
test categories or the like. We will be led to consider, for
\emph{any} small category $A$, the localization of \Ahat{} with
respect to \scrW-equivalences, which I'll denote by \HotAW{} or simply
\HotA:
\begin{equation}
  \label{eq:67.1}
  \HotAW=\HotA=\scrW^{-1}\Ahat.\tag{1}
\end{equation}
These\pspage{189} categories, I suspect, are quite interesting in
themselves, and they merit to be understood. Thus, one of Quillen's
results asserts that (at least for $\scrW=W_\Cat=$ ordinary weak
equivalences, but presumably his arguments will carry over to an
arbitrary \scrW) in case $A$ is a product category $\Simplex\times A_0$,
where $A_0$ is any small category and $\Simplex$ the category of
standard, ordered simplices, then \Ahat{} is a closed model category
admitting as weak equivalences the set \scrWA\footnote{\alsoondate{3.5.}
  This, I quickly became aware, is a misrepresentation of Quillen's
  result -- that ``weak equivalences'' he introduced are a lot
  stronger than \scrWA. I'll have to come back upon this soon
  enough!}; and hence \HotA, the corresponding ``homotopy category'',
admits familiar homotopy constructions, including the two types of
Dold-Puppe exact sequences, tied up with loop- and suspension
functors. It is very hard to believe that this should be a special
feature of the category $\Simplex$ as the multiplying factor -- surely
any test category or strict test category instead should do as
well. As we'll check below, the product of a local test category with
any category $A_0$ is again a local test category, hence a
test-category if both factors are (\scrW-)aspheric. Thus suggests that
maybe for any local test category $A$, the corresponding \Ahat{} is a
closed model category -- but it isn't even clear yet if the same
doesn't hold for \emph{any} small category $A$ whatever! It's surely
something worth looking at.

As we'll see presently, it is a tautology more or less that a functor
\begin{equation}
  \label{eq:67.1again}
  i:A\to B,\tag{1}
\end{equation}
giving rise to a functor
\begin{equation}
  \label{eq:67.2}
  i^*:\Bhat\to\Ahat\tag{2}
\end{equation}
(commuting to all types of direct and inverse limits), induces a
functor on the localizations
\[\overline{i^*}:\HotB\to\HotA,\]
provided $i$ is \scrW-aspheric. Thus,
\[A\mapsto\HotA\]
can be viewed as a functor with respect to $A$, provided we take as
``morphisms'' between ``objects'' $A$ the \emph{aspheric} functors
only -- i.e., it is a functor on the subcategory
$\Cat_{\textup{\scrW-asph}}$ of \Cat, having the same objects as \Cat,
but with maps restricted to be \scrW-aspheric ones.

We'll denote by
\[\HotOf(\scrW) = \scrW^{-1}\Cat\]
the homotopy category defined in terms of the basic localizer
\scrW. For any small category $A$, we get a commutative
diagram\pspage{190}
\begin{equation}
  \label{eq:67.3}
  \begin{tabular}{@{}c@{}}
    \begin{tikzcd}[baseline=(O.base)]
      \Ahat\ar[r,"i_A"]\ar[d,swap,"\gamma_A"] &
      \Cat\ar[d,"\text{$\gamma_\scrW$ or $\gamma$}"] \\
      \HotA \ar[r,"\overline{i_A}"] &
      |[alias=O]| (\HotOf(\scrW))
    \end{tikzcd},
  \end{tabular}
  \tag{3}
\end{equation}
we denote by
\[\varphi_A=\gamma i_A = \overline{i_A}\gamma_A :
\Ahat\to(\HotOf(\scrW))\]
the corresponding composition.

Coming back to the case of a functor $i$ and corresponding $i^*$
(\eqref{eq:67.1again} and \eqref{eq:67.2}), the functors $i^*$, $i_A$,
$i_B$ do not give rise to a commutative triangle, but to a triangle
\emph{with commutation morphism} $\lambda_i$:
\[\begin{tikzcd}[baseline=(O.base),column sep=small]
  \Bhat\ar[rr,"i^*"]\ar[dr,swap,"i_B",""{name=iB,right}] & &
  \Ahat\ar[dl,"i_A",""{name=iA,left}] \\
  & |[alias=O]| \Cat \arrow[bend right=10, from=iA, to=iB, swap, "\lambda_i"] &
\end{tikzcd},\]
i.e., for any $F$ in \Bhat, we get a map
\begin{equation}
  \label{eq:67.4}
  \lambda_i(F): i_Ai^*(F)=A_{/i^*(F)} \to i_B(F)=B_{/F} ,\tag{4}
\end{equation}
the first hand side of \eqref{eq:67.4}, also written simply $A_{/F}$
when there is no ambiguity for $i$, can be interpreted as the category
of pairs
\[(a,p), \quad a\in\Ob A,\quad p:i(a)\to F,\]
where $p$ is a map in \Bhat, $B$ identified as usual to a full
subcategory of \Bhat{} (hence $i(a)$ identified with an object of
\Bhat). The map $\lambda_i(F)$ for fixed $F$ is of course the functor
\[ (a,p) \mapsto (i(a),p).\]
The topological significance of course is clear: interpreting $i$ as
defining a ``map'' or morphisms of the corresponding topoi \Ahat{} and
\Bhat, having $i^*$ as inverse image functor, an object $F$ of \Bhat{}
gives rise to an \emph{induced topos} $\Bhat_{/F} \tosimeq
(B_{/F})\uphat$, and the restriction of the ``topos above'' \Ahat{} to
the induced $\Bhat_{/F}$, or equivalently the result of base change
$\Bhat_{/F}\to\Bhat$, gives rise to the induced morphism of topoi
$\Ahat_{/i^*(F)}\to\Bhat_{/F}$, represented precisely by the map
$\lambda_i(F)$ in \Cat.

The condition of \scrW-asphericity on $i$ may be expressed in manifold
ways, as properties of either one of the three aspects
\[\lambda_i, \quad i^*, \quad i_Ai^*\]
of the situation created by $i$, with respect to the localizing sets
\scrW, \scrWA, $\scrW_B$, or to the notion of aspheric object. As
\scrWA{} is defined in terms of \scrW{} as just the inverse image of
the latter by $i_A$, and the same for aspheric objects, it turns out
that each of the conditions we\pspage{191} are led to express on
$i^*$, can be formulated equivalently in terms of the composition
$i_Ai^*$. I'll restrict to formulate these in terms of $i^*$ only,
which will be the form most adapted to the use we are going to make
later of the notion of a \scrW-aspheric map, when introducing the
so-called ``asphericity structures'' and corresponding ``testing
functors''.

\bigbreak
\presectionfill\ondate{27.5.}\par

\hangsection{Digression on a ``new continent''.}\label{sec:68}%
It has been over six weeks now that I didn't write down any notes. The
reason for this is that I felt the story of asphericity structures and
canonical modelizers was going to come now without any problem, almost
as a matter of routine to write it down with some care -- therefore, I
started doing some scratchwork on a few questions which had been
around but kept in the background since the beginning, and which were
a lot less clear in my mind. Some reflection was needed anyhow, before
it would make much sense to start writing down anything on
these. Finally, it took longer than expected, as usual -- partly
because (as usual too!), a few surprises would turn up on my
way. Also, I finally allowed myself to become distracted by some
reflection on the ``Lego-Teichmüller construction game'', and pretty
much so during last week. The occasion was a series of informal talks
Y.\scrcomment{Y.\ = Yves Ladegaillerie} has been giving in Molino's seminar,
on Thurston's hyperbolic geometry game and his compactification of
Teichm\"uller space. Y.\ was getting interested again in mathematics
after a five year's interruption. He must have heard about my seminar
last year on ``anabelian algebraic geometry'' and the ``Teichmüller
tower'', and suggested I might drop in to get an idea about Thurston's
work. This work indeed appears as closely related in various respects
to my sporadic reflections of the last two years, just with a
diametrically opposed emphasis -- mine being on the algebro-geometric
and arithmetic aspects of ``moduli'' of algebraic curves, his on
hyperbolic riemannian geometry and the simply connected transcendental
Teichmüller spaces (rather than the algebraic modular
varieties). The main intersection appears to be interest in surface
surgery and the relation of this to the Teichmüller modular group. I
took the occasion to try and recollect about the Lego-Teichmüller
game, which I had thought of last year as a plausible, very concrete
way for modelizing and visualizing the whole tower of Teichmüller
groupoids $T_{g,\nu}$ and the main operations among these, especially
the ``cutting'' and ``gluing'' operations. The very informal talk I
gave was mainly intended for Y.\ as a matter of fact, and
it\pspage{192} was an agreeable surprise to notice that the message
this time was getting through. For the fiver or six years since my
attention became attracted by the fascinating melting-pot of key
structures in geometry, topology, arithmetic, discrete and algebraic
groups, intertwining tightly in a kind of very basic
Galois-Teichmüller theory, Y.\ has been the very first person I met
so far to have a feeling for (a not yet dulled instinct I might say,
for sensing) the extraordinary riches opening up here for
investigation. The series of talks I had given in a tentative seminar
last year had turned short, by lack of any active interest and
participation of anyone among the handful of mere listeners. And the
two or three occasions I had the years before to tell about the matter
of two-dimensional maps (``cartes'') and their amazing
algebro-arithmetic implications, to a few highbrow colleagues with
incomparably wider background and know-how than anyone around here, I
met with polite interest, or polite indifference which is the same. As
there was nobody around anyhow to take any interest the these juicy
greenlands, nobody would care to see, because there was no text-book
nor any official seminar notes to prove they existed, after a few
years I finally set off myself for a preliminary voyage.

I thought it was going to take me a week or two to tour it and kind of
recense\scrcomment{``recense'' could be translated as ``survey''}
resources. It took me five months instead of intensive work, and two
impressive heaps of notes (baptized ``La Longue Marche à travers la
théorie de Galois''), to get a first, approximative grasp of some of
the main structures and relationships involved. The main emphasis was
(still is) on an understanding of the action of profinite
Galois-groups (foremost among which \GalQQ{} and the subgroups of
finite index) on non-commutative profinite fundamental groups, and
primarily on fundamental groups of algebraic curves -- increasingly
too on those of modular varieties (more accurately, modular
\emph{multiplicities}) for such curves -- the profinite completions of
the Teichmüller group. The voyage was the most rewarding and exciting
I had in mathematics so far -- and still it became very clear that it
was just like a first glimpse upon a wholly new landscape -- one
landscape surely among countless others of a continent unknown, eager
to be discovered.

This was in the first half of the year 1981 -- just two years ago, it
turns out, but it look almost infinitely remote, because such a lot of
things took place since. Looking back, it turns out there have been
since roughly four main alternating periods of reflection, one period
of reflection on personal matters alternating with one on
mathematics. The next mathematical reflection started with a long
digression on tame topology and the ``déployment''
(``unfolding'')\pspage{193} of stratified structures, as a leading
thread towards a heuristic understanding of the natural stratification
of the Mumford-Deligne compactifications of modular multiplicities
$M_{g,\nu}$ (for curves of genus $g$ endowed with $\nu$ points). This
then led to the ``anabelian seminar'' which turned short, last
Spring. Then a month or two sicknees, intensive meditation for three
or four months, a few more months for settling some important personal
matter; and now, since February, another mathematical reflection
started.

I am unable to tell the meaning of this alternation of periods of
meditation on personal matters and periods of mathematical reflection,
which has been going through my life for the last seven years, more
and more, very much like the unceasing up and down of waves, or like a
steady breathing going through my life, without any attempt any longer
of controlling it one way or the other. One common moving force surely
is the inborn curiosity -- a thirst for getting acquainted with the
juicy things of the inexhaustible world, whether they be the breathing
body of the beloved, or the evasive substance of one's own life, or
the much less evasive substance of mathematical structure and their
delicate interplay. This thirst in itself is of a nature quite
different of the ego's -- it is the thirst of life to know about
itself, a primal creative force which, one suspects, has been around
forever, long before a human ego -- a bundle mainly of fears, of
inhibitions and self-deceptions -- came into being. Still, I am aware
that the ego is strongly involved in the particular way in which the
creative force expresses itself, in my own life or anyone else's (when
this force is allowed to come into play at all\ldots). The motivations
behind any strong energy investment, and more particularly so when it
is an activity attached with any kind of social status or prestige,
are a lot more complex and ego-driven than one generally cares to
admit. True, ambition by itself is powerless for discovering or
understanding or perceiving anything substantial whatever, neither a
mathematical relationship nor the perfume of a flower. In moments of
work and of discovery, in any creative moments in life, as a matter of
fact, ambition is absent; the Artisan is a keen interest, which is
just one of the manifold aspect of love. What however pushes us so
relentlessly to work again, and so often causes our life and passion
gradually to dry off and become insensitive, even to the kindred
passion of a follow-being, or to the unsuspected beauties and
mysteries of the very field we are supposed to be plowing -- this
force is neither love nor keen interest for the things and beings in
this world. It is interesting enough though, and surely deserves a
close look!

I\pspage{194} thought I was starting a retrospective of six weeks of
scratchwork on homotopical algebra -- and it turns out to be a (very)
short retrospective rather of the up-and-down movement of my
mathematical interests and investments during the last
years. Doubtless, the very strongest attraction, the greatest
fascination goes with the ``new world'' of anabelian algebraic
geometry. It may seem strange that instead, I am indulging in this
lengthy digression on homotopical algebra, which is almost wholly
irrelevant I feel for the Galois-Teichmüller story. The reason is
surely an inner reluctance, an unreadiness to embark upon a long-term
voyage, well knowing that it is so enticing that I may well be caught
in this game for a number of years -- not doing anything else day and
night than making love with mathematics, and maybe sleeping and eating
now and then. I have gone through this a number of times, and at times
I thought I was through. Finally, I came to admit and to accept, two
years ago, I was not through yet -- this was during the months of
meditation after the ``long march through Galois theory'' -- which had
been, too, a wholly unexpected fit of mathematical passion, not to say
frenzy. And during the last weeks, just reflecting a little here and
there upon the Teichmüller-Lego game and its arithmetical
implications, I let myself be caught again by this fascination -- it
is becoming kind of clear now that I am going to finish writing up
those notes on algebra, almost like some homework that has got to be
done (anyhow I like to finish when I started something) -- and as soon
as I'm through with the notes, back to geometry in the long last!
Also, the idea is in the air for the last few months -- since I
decided to publish these informal notes on stacks or whatever it'll
turn out to be -- that I may well go on the same way, writing up and
publishing informal notes on other topics, including tame topology and
anabelian algebraic geometry. In contrast to the present notes, I got
heaps of scratchwork done on these in the years before -- in this
respect time is even riper for me to ramble ``publicly'' than on
stacks and homotopy theory!

From Y.\ who looked through a lot of literature on the subject, it
strikes me (agreeably of course) that nobody yet hit upon ``the''
natural presentation of the Teichmüller groupoids, which kind of
imposes itself quite forcibly in the set-up I let myself be guided
by. Technically speaking (and this will rejoice Ronnie Brown I'm
sure!), I suspect one main reason why this is so, is that people are
accustomed to working with fundamental groups and generators and
relations for these and stick to it, even in contexts when this is
wholly inadequate, namely when you get a clear\pspage{195} description
by generators and relations only when working simultaneously with a
whole bunch of base-points chosen with case -- or equivalently,
working in the algebraic context of \emph{groupoids}, rather than
groups. Choosing paths for connecting the base-points natural to the
situation to just one among them, and reducing the groupoid to a
single group, will then hopelessly destroy the structure and inner
symmetries of the situation, and result in a mess of generators and
relations no-one dares to write down, because everyone feels they
won't be of any use whatever, and just confuse the picture rather than
clarify it. I have known such perplexity myself a long time ago,
namely in van Kampen-type situations, whose only understandable
formulation is in terms of (amalgamated sums of) groupoids. Still,
standing habits of thought are very strong, and during the long march
through Galois theory, two years ago, it took me weeks and months
trying to formulate everything in terms of groups or ``exterior
groups'' (i.e., groups ``up to inner automorphism''), and finally
learning the lesson and letting myself be convinced progressively, not
to say reluctantly, that groupoids only would fit nicely. Another
``technical point'' of course is the basic fact (and the wealth of
intuitions accompanying it) that the Teichmüller groups are
fundamental groups indeed -- a fact ignored it seems by most
geometers, because the natural ``spaces'' they are fundamental groups
of are not topological spaces, but the modular ``multiplicities''
$M_{g,\nu}$ -- namely topoi! The ``points'' of these ``spaces'' are
just the structures being investigated (namely algebraic curves of
type $(g,\nu)$), and the (finite) automorphism groups of these
``points'' enter into the picture in a very crucial way. They can be
adequately chosen as part of the system of basic generators for the
Teichmüller groupoid $T_{g,\nu}$. The latter of course is essentially
(up to suitable restriction of base-points) just the fundamental
groupoid of $M_{g,\nu}$. It is through this interpretation of the
Teichmüller groups or groupoids that it becomes clear that the
profinite Galois group \GalQQ{} operates on the profinite completion
of these and of their various variants, and this (it turns out) in a
way respecting the manifold structures and relationships tying them
tightly together.

\bigbreak

\presectionfill\ondate{29.5.}\pspage{196}\par

\hangsection[Digression on six weeks' scratchwork: \emph{derivators},
and \dots]{Digression on six weeks' scratchwork:
  \texorpdfstring{\emph{derivators}}{derivators}, and integration of
  homotopy types.}\label{sec:69}%
Before resuming more technical work again, I would like to have a
short retrospective of the last six weeks' scratchwork, now lying on
my desk as a thickly bunch of scratchnotes, nobody but I could
possibly make any sense of.

The first thing I had on my mind has been there now for nearly twenty
years -- ever since it had become clear, in the SGA~5 seminar on
$L$-functions and apropos the formalism of traces in terms of derived
categories, that Verdier's set-up of derived categories was
insufficient for formulating adequately some rather evident situations
and relationships, such as the addition formula for traces, or the
multiplicative formula for determinants. It then became apparent that
the derived category of an abelian category (say) was too coarse an
object various respects, that it had to be complemented by similar
``triangulated categories'' (such as the derived category of a
suitable category of ``triangles'' of complexes, or the whole bunch of
derived categories of categories of filtered complexes of order $n$
with variable $n$), closely connected to it. Deligne and Illusie had
both set out, independently, to work out some set-up meeting the most
urgent requirements (Illusie's treatment in terms of filtered
complexes was written down and published in his thesis six years later
(Springer Lecture Notes N\textsuperscript{\b o}~239)).%
\scrcomment{\cite{Illusie1971,Illusie1972}}
While adequate for the main tasks then at hand, neither treatment was
really wholly satisfactory to my taste. One main feature I believe
making me feel uncomfortable, was that the extra categories which had
to be introduced, to round up somewhat a stripped-and-naked
triangulated category, were triangulated categories in their own
right, in Verdier's sense, but remaining nearly as stripped by
themselves as the initial triangulated category they were intended to
provide clothing for. In other words, there was a lack of inner
stability in the formalism, making it appear as very much provisional
still. Also, while interested in associating to an abelian category a
handy sequence of ``filtered derived categories'', Illusie made no
attempt to pin down what exactly the inner structure of the object he
had arrived at was -- unlike Verdier, who had introduced, alongside
with the notion of a derived category of an abelian category, a
general notion of triangulated categories, into which these derived
categories would fit. The obvious idea which was in my head by then
for avoiding such shortcomings, was that an abelian category \scrA{}
gave rise, not only to the single usual derived category $D(\scrA)$ of
Verdier, but also, for every type of diagrams,\pspage{197} to the
derived category of the abelian category of all \scrA-valued diagrams
of this type. In precise terms, for any small category $I$, we get the
category $D(\bHom(I,\scrA))$, depending functorially in a
contravariant way on $I$. Rewriting this category $D_\scrA(I)$ say,
the idea was to consider
\begin{equation}
  \label{eq:69.star}
  I\mapsto D_\scrA(I),\tag{*}
\end{equation}
possibly with $I$ suitably restricted (for instance to finite
categories, or to finite ordered sets, corresponding to finite
\emph{commutative} diagrams), as embodying the ``full'' triangulated
structure defined by \scrA. This of course at once raises a number of
questions, such as recovering the usual triangulated structure of
$D(\scrA)=D_\scrA(e)$ ($e$ the final object of \Cat) in terms of
\eqref{eq:69.star}, and pinning down too the relevant formal
properties (and possibly even extra structure) one had to assume on
\eqref{eq:69.star}. I had never so far taken the time to sit down and
play around some and see how this goes through, expecting that surely
someone else would do it some day and I would be informed -- but
apparently in the last eighteen years nobody ever was
interested. Also, it had been rather clear from the start that
Verdier's constructions could be adapted and did make sense for
non-commutative homotopy set-ups, which was also apparent in between
the lines in Gabriel-Zisman's\scrcomment{\cite{GabrielZisman1967}}
book on the foundations of homotopy theory, and a lot more explicitly
in Quillen's axiomatization of homotopical algebra. This
axiomatization I found very appealing indeed -- and right now still
his little book\scrcomment{\cite{Quillen1967}} is my most congenial
and main source of information on foundational matters of homotopical
algebra. I remember though my being a little disappointed at Quillen's
not caring either to pursue the matter of what exactly a
``non-commutative triangulated category structure'' (of the type he
was getting from his model categories) was, just contenting himself to
mumble a few words about existence of ``higher structure'' (then just
the Dold-Puppe sequences), which (he implies) need to be understood. I
felt of course that presumably the variance formalism
\eqref{eq:69.star} should furnish any kind of ``higher'' structure one
was looking for, but it wasn't really my business to check.

It still isn't, however I did some homework on \eqref{eq:69.star} --
it was the first thing indeed I looked at in these six weeks, and some
main features came out very readily indeed. It turns out that the main
formal variance property to demand on \eqref{eq:69.star}, presumably
even the only one, is that for a given map $f:I\to J$ on the indexing
categories of diagram-types $I$ and $J$, the corresponding functor
\[f^*: D(J) \to D(I)\]
\emph{should have both a left and a right adjoint}, say $f_!$ and
$f_*$. In case $J=e$, the two\pspage{198} functors we get from $D(I)$
to $D(e)=\scrD$ (the ``stripped'' triangulated category) can be viewed
as a \emph{substitute for taking, respectively, direct and inverse
  limits in \scrD} (for a system of objects indexed by $I$),
\emph{which in the usual sense don't generally exist in \scrD} (except
just finite sums and products). These operations admit as important
special cases, when $I$ is either one of the two mutually dual
categories
\begin{equation}
  \label{eq:69.starstar}
  \begin{tikzcd}[baseline=(O.base),column sep=small,row sep=tiny]
    & b \\ |[alias=O]| a\ar[ur]\ar[dr] & \\ & c
  \end{tikzcd}
  \quad\text{or}\quad
  \begin{tikzcd}[baseline=(O.base),column sep=small,row sep=tiny]
    & b\ar[dl] \\ |[alias=O]| a & \\ & c\ar[ul]
  \end{tikzcd}\quad,\tag{**}
\end{equation}
the operation of (binary) amalgamated sums or fibered products, and
hence also of taking ``cofibers'' and ``fibers'' of maps, in the sense
introduced by Cartan-Serre in homotopy theory about thirty years
ago. I also checked that the two mutually dual Dold-Puppe sequences
follow quite formally from the set-up. One just has to fit in a
suitable extra axiom to ensure the usual exactness properties for
these sequences.

Except in the commutative case when starting with an abelian category
as above, I did not check however that there is indeed such ``higher
variance structure'' in the usual cases, when a typical ``triangulated
category'' in some sense or other turns up, for instance from a model
category in Quillen's sense. What I did check though in this last
case, under a mild additional assumption which seems verified in all
practical cases is the existence of the operation $f_!=\int_I$
(``integration'') and $f_*=\prod_I$ (cointegration) for the special
case $f:I\to e$, when $I$ is either of the two categories
\eqref{eq:69.starstar} above. I expect that working some more, one
should get under the same assumptions at least the existence of $f_!$
and $f_*$ for any map $f:I\to J$ between finite ordered sets.

My main interest of course at present is in the category \Hot{}
itself, more generally in $\HotOf(\scrW)=\scrW^{-1}\Cat$, where
\scrW{} is a ``basic localizer''. More generally, if $(M,W)$ is any
modelizer say, the natural thing to do, paraphrasing
\eqref{eq:69.star}, is for any indexing category $I$ to endow
$\bHom(I,M)$ with the set of arrows $W_I$ defined by componentwise
belonging to $W$, and to define
\[ D_{(M,W)}(I) = D(I) \eqdef W_I^{-1}\bHom(I,M),\]
with the obvious contravariant dependence on $I$, denoted by $f^*$ for
$f:I\to J$. The question then arises as to the existence of left and
right adjoints, $f_!$ and $f_*$. In case we take $M=\Cat$, the
existence of $f_!$ goes through with amazing smoothness: interpreting
a ``model'' object of $\bHom(I,\Cat)$, namely a functor
\[I\to \Cat\]
in\pspage{199} terms of a cofibered category $X$ over $I$
\[p: X\to I,\]
and assuming for simplicity $f$ cofibering too, $f_!(X)$ is just $X$
itself, the total category of the cofibering, viewed as a (cofibered)
category over $J$ by using the functor $g=f\circ p$! This applies for
instance when $J$ is the final category, and yields the operation of
``integration of homotopy types'' $\int_I$, in terms of the total
category of a cofibered category over $I$. If we want to rid ourselves
from any extra assumption on $f$, we can describe $D(I)$ (up to
equivalence) in terms of the category $\Cat_{/I}$ of categories $X$
over $I$ (not necessarily cofibered over $I$), $W_I$ being replaced by
the corresponding notion of ``\emph{\scrW-equivalences relative to
  $I$}'' for maps $u:X\to Y$ of objects of \Cat{} over $I$, by which
we mean a map $u$ such that the localized maps
\[u_{/i}: X_{/i} \to Y_{/i}\]
are in \scrW, for any $i$ in $I$. Regarding now any category $X$ over
$I$ as a category over $J$ by means of $f\circ p$, this is clearly
compatible with the relative weak equivalences $\scrW_I$ and
$\scrW_J$, and yields by localization the looked-for functor $f_!$.

This amazingly simple construction and interpretation of the basic
$f_!$ and $\int_I$ operations is one main reward, it appears, for
working with the ``basic localizer'' \Cat, which in this occurrence,
as in the whole test- and asphericity story, quite evidently deserves
its name. It has turned out since that in some other respects -- for
instance, paradoxically when it comes to the question of the
relationship between this lofty integration operation, and true honest
amalgamated sums -- the modelizers \Ahat{} associated to test
categories $A$ (namely the so-called ``elementary modelizers'') are
more convenient tools than \Cat. Thus, it appears very doubtful still
that \Cat{} is a ``model category'' in Quillen's sense, in any
reasonable way (with \scrW{} of course as the set of ``weak
equivalences''). I finally got the feeling that a good mastery of the
basic aspects of homotopy types and of basic relationships among
these, will require mainly great ``aisance''\scrcomment{``aisance''
  here could be ``ease'' or ``fluency''} in playing around with a
number of available descriptions of homotopy types by models, no one
among which (not even by models in \Cat, and surely still less by
semisimplicial structures) being adequate for replacing all others.

As for the $f_*$ and cointegration $\prod_I$ operations among the
categories $D(I)$, except in the very special case noted above
(corresponding to fibered products), I did not hit upon any
ready-to-use candidate for it, and I doubt there is any. I do believe
the operations exist indeed, and I even have\pspage{200} in mind a
rather general condition on a pair $(M,W)$ with $W\subset\Fl(M)$, for
both basic operations $f_!$ and $f_*$ to exist between the
corresponding categories $D_{(M,W)}(I)$ -- but to establish this
expectation may require a good amount of work. I'll come back upon
these matters in due course.

There arises of course the question of giving a suitable name to the
structure $I\to D(I)$ I arrived at, which seems to embody at least
some main features of a satisfactory notion of a ``triangulated
category'' (not necessarily commutative), gradually emerging from
darkness. I have thought of calling such a structure a
``\emph{derivator}'', with the implication that its main function is
to furnish us with a somehow ``complete'' bunch (in terms of a
rounded-up self-contained formalism) of categories $D(I)$, which are
being looked at as ``\emph{derived categories}'' in some sense or
other. The only way I know of for constructing such a derivator, is as
above in terms of a pair $(M,W)$, submitted to suitable conditions for
ensuring existence of $f_!$ and $f_*$, at least when $f$ is any map
between finite ordered sets. We may look upon $D(I)$ as a refinement
and substitute for the notion of family of objects of $D(e)=D_0$
indexed by $I$, and the integration and cointegration operations from
$D(I)$ to $D_0$ as substitutes (in terms of these finer objects) of
direct and inverse limits in $D_0$. When tempted to think of these
latter operations (with values in $D_0$) as the basic structures
involved, one cannot help though looking for the same kind of
structure on any one of these subsidiary categories $D(I)$, as these
are being thought of as derived categories in their own right. It then
appears at once that the ``more refined substitutes'' for $J$-indexed
systems of objects of $D(I)$ are just the objects of $D(I\times J)$,
and the corresponding integration and cointegration operations
\[ D(I\times J)\to D(I)\]
are nothing but $p_!$ and $p_*$, where $p:I\times J\to I$ is the
projection. Thus, one is inevitably conducted to look at operations
$f_!$ and $f_*$ instead of merely integration and cointegration --
thus providing for the ``inner stability'' of the structure described,
as I had been looking for from the very start.

The notion of integration of homotopy types appears here as a natural
by-product of an attempt to grasp the ``full structure'' of a
triangulated category. However, I had been feeling the need for such a
notion of integration of homotopy types for about one or two years
already (without any clear idea yet that this operation should be one
out of two\pspage{201} main ingredients of a (by then still very
misty) notion of a triangulated category of sorts). This feeling arose
from my ponderings on stratified structures and the ``screwing
together'' of such structures in terms of simple building blocks
(essentially, various types of ``tubes'' associated to such
structures, related to each other by various proper maps which are
either inclusion or -- in the equisingular case at any rate -- fiber
maps). This ``screwing together operation'' could be expressed as
being a direct limit of a certain finite system of spaces. In the
cases I was most interested in (namely the Mumford-Deligne
compactifications $M\uphat_{g,\nu}$ of the modular multiplicities
$M_{g,\nu}$), these spaces or ``tubes'' have exceedingly simple
homotopy types -- they are just $K(\pi,1)$-spaces, where each $\pi$ is
a Teichmüller-type discrete group (practically, a product of usual
Teichmüller groups). It then occurred to me that the whole homotopy
type of $M\uphat_{g,\nu}$, or of any locally closed union of strata,
or (more generally still) of ``the'' tubular neighborhood of such a
union in any larger one, etc.\ -- that all these homotopy types should
be expressible in terms of the given system of spaces, and more
accurately still, just in terms of the corresponding system of
fundamental groupoids (embodying their homotopy types). In this
situation, what I was mainly out for, was precisely an accurate and
workable description of this direct system of groupoids (which could
be viewed as just one section of the whole ``Teichmüller tower'' of
Teichmüller groupoids\ldots). Thus, it was a rewarding extra feature
of the situation (by then just an expectation, as a matter of fact),
that such a description should at the same time yield a ``purely
algebraic'' description of the homotopy types of all the spaces
(rather, multiplicities, to be wholly accurate) which I could think of
in terms of the natural stratification of $M\uphat_{g,\nu}$. There was
an awareness that this operation on homotopy types could not be
described simply in terms of a functor $I\to\Hot$, where $I$ is the
indexing category, that a functor $i\mapsto X_i: I\to M$ (where $M$ is
some model category such as \Spaces{} or \Cat) should be available in
order to define an ``integrated'' homotopy type $\int X_i$. This
justified feeling got somewhat blurred lately, for a little while, by
the definitely unreasonable expectation that finite limits should
exist in \Hot{} after all, why not! It's enough to have a look though
(which probably I did years ago and then forgot in the meanwhile) to
make sure they don't\ldots

Whether or not this notion of ``integration of homotopy types'' is
more or less well known already under some name or other, isn't quite
clear to me. It isn't familiar to Ronnie Brown visibly, but it seems
he heard about such a kind of thing, without his being specific about
it. It\pspage{202} was the episodic correspondence with him which
finally pushed me last January to sit down for an afternoon and try to
figure out what there actually was, in a lengthy and somewhat rambling
letter to Illusie (who doesn't seem to have heard at all about such
operations). This preliminary reflection proved quite useful lately,
I'll have to come back anyhow upon some of the specific features of
integration of homotopy types later, and there is not much point
dwelling on it any longer at present.

\hangsection[Digression on scratchwork (2): cohomological
\dots]{Digression on scratchwork
  \texorpdfstring{{\normalfont(2)}}{(2)}: cohomological properties of
  maps in \texorpdfstring{\Cat}{(Cat)} and in
  \texorpdfstring{\Ahat}{Ahat}. Does any topos admit a ``dual'' topos?
  Kan fibrations rehabilitated.}\label{sec:70}%
This was a (finally somewhat length!) review of ponderings which
didn't take me more than just a few days, because it was about things
some of which were on my mind for a long time indeed. It took a lot
more work to try to carry through the standard homotopy constructions,
giving rise to the Dold-Puppe sequences, within the basic modelizer
\Cat. Most of the work arose, it now seems to me, out of a block I got
(I couldn't tell why) against the Kan-type condition on complexes, so
I tried hard to get along without anything of the sort. I kind of
fooled myself into believing that I was forced to do so, because I was
working in an axiomatic set-up dependent upon the ``basic localizer''
\scrW, so the Kan condition wouldn't be relevant anyhow. The main
point was to get, for any map $f:X\to Y$ in \Cat, a factorization
\[ X \xrightarrow i Z \xrightarrow p Y,\]
where $p$ has the property that base-change by $p$ transforms weak
equivalences into weak equivalences (visibly a Serre-fibration type
condition), and $i$ satisfies the dual condition with respect to
co-base change; and moreover where either $p$, or $i$ can be assumed
to be in \scrW, i.e., to be a weak equivalence. (This, by the way, is
the extra condition on a pair $(M,W)$ I have been referring to above
(page \ref{p:200}, for (hopefully) getting $f_!$ and $f_*$
operations.) At present, I do not yet know whether such factorization
always exists for a map in \Cat, without even demanding that either
$i$ or $p$ should be in \scrW.

I first devoted a lot of attention to Serre-type conditions on maps in
\Cat, which turned out quite rewarding -- with the impression of
arriving at a coherent and nicely auto-dual picture of cohomology
properties of functors, i.e., maps in \Cat, as far as these were
concerned with \emph{base change behavior} (and \emph{not} co-base
change). Here I was guided by work done long ago to get the étale
cohomology theory off the ground, and where the two main theorems
achieving this aim were precisely two theorems of commutation of
``higher direct images'' $\mathrm R^ig_*$ with respect to base-change
by a map $h$ -- namely, it is OK when either $g$ is \emph{proper}, or
$h$ is \emph{smooth}. It was rather natural then to introduce the
notion of\pspage{203} (cohomological) \emph{smoothness} and (coh.)
\emph{properness} of a map in \Cat, by the obvious base-change
properties. It turned out that these could be readily characterized by
suitable asphericity conditions, which are formally quite similar to
the well-known valuative criterion for a map of finite type between
noetherian schemes to be ``universally open'' (which can be viewed as
a ``purely topological'' variant of the notion of smoothness), resp.\
``universally closed'', or rather, more stringently, proper. These
conditions, moreover, are trivially satisfied when $f: X\to Y$ turns
$X$ into a category over $Y$ which is fibered over $Y$ (i.e.,
definable in terms of a contravariant pseudo-functor $Y\op\to\Cat$),
resp.\ cofibered over $Y$ (namely, definable in terms of a
pseudo-functor $Y\to\Cat$). If we call such maps in \Cat{}
``fibrations'' and ``cofibrations'' (very much in conflict, alas, with
Quillen's neat set-up of fibrations-cofibrations!), it turns out that
\emph{fibrations are smooth}, \emph{cofibrations are proper}. This is
all the more remarkable, as the notions of fibered and cofibered
categories were introduced with a view upon ``large categories'', in
order to pin down some standard properties met in all situations of
``base change'' (and later the dual situation of co-base change) --
the main motivation for this being the need to formulate with a
minimum of precision a set-up for ``descent techniques'' in algebraic
geometry. (These techniques, as well as the cohomological base change
theorems, make visibly a sense too in the context of analytic spaces
say, or of topological spaces, but they don't seem to have been
assimilated yet by geometers outside of algebraic geometry.) That
these typically ``general nonsense'' notions should have such precise
topological implications came as a complete surprise! As a
consequence, a bifibration (namely a map which is both a fibration and
a cofibration) is smooth and proper and hence a (cohomological)
``Serre fibration'', for instance the sheaves $\mathrm R^if_*(F)$,
when $F$ is a constant abelian sheaf on $X$ (i.e., a constant abelian
group object in $X\uphat$), namely $y\mapsto\mathrm H^i(X_{/y},F)$,
are \emph{local systems} on $Y$, i.e., factor through the fundamental
groupoid of $Y$.

A greater surprise still was the \emph{duality between the notions of
  properness and smoothness:} just as a map $f:X\to Y$ in \Cat{} is a
cofibration if{f} the ``dual'' map $f\op:X\op\to Y\op$ is a fibration,
it turns out that \emph{$f$ is proper if{f} $f\op$ is smooth}. This
was a really startling fact, and it caused me to wonder, in the
context of more general topoi than just those of the type $X\uphat$,
whether there wasn't \emph{a notion of duality generalizing the
  relationship between two topoi $X\uphat$ and $X\upvee \eqdef
  (X\op)\uphat$}. Indeed, these two categories of sheaves can be
described intrinsically, one in terms of the other, up to equivalence,
by a natural pairing\pspage{204}
\[X\uphat \times X\upvee \to \Sets\]
commuting componentwise with (small) direct limits, and inducing an
equivalence between either factor with the category of ``co-sheaves''
on the other, namely covariant functors to \Sets{} commuting with
direct limits. But it isn't at all clear, starting with an arbitrary
topos \scrA{} say, whether the category $\scrA'$ of all cosheaves on
\scrA{} is again a topos, and still less whether \scrA{} can be
recovered (up to equivalence) in terms of $\scrA'$ as the category of
all cosheaves on $\scrA'$.

To come back though to the factorization problem raised above (p.\
\ref{p:202}), the main trouble here is that, except the case of an
isomorphism $i: X\tosim Z$, I was unable to pin down a single case of
a map $i$ in \Cat{} such that co-base change by $i$ transforms weak
equivalences (in the usual sense say) into weak equivalences. One
candidate I had in mind, the so-called ``\emph{open immersions}'',
namely functors $i:X\to Z$ inducing an isomorphism between $X$ and a
``sieve'' (or ``crible'') in $Z$ (corresponding to an open sub-topos
of $Z\uphat$), and dually the ``closed immersions'', finally have
turned out delusive -- a disappointment maybe, but still more a big
relief to find out at last how the score was! Almost immediately in
the wake of this negative result in \Cat, and in close connection with
the fairly well understood $\int$ substitute for amalgamated sums in
\Cat, came the big compensation valid in any category $\Ahat$, namely
the fact that for a cocartesian square
\[\begin{tikzcd}
  Y\ar[r,"i"] \ar[d,swap,"g"] & X\ar[d,"f"] \\
  Y'\ar[r,"i'"] & X'
\end{tikzcd}\]
in \Ahat, where $i$ is a monomorphism, if $i$ (resp.\ $g$) is in
\scrWA, so is $i'$ (resp.\ $f$). This implies that \emph{co-base
  change by a monomorphism in \Ahat{} transform weak equivalences into
  weak equivalences}. The common main fact behind these statements is
that for a diagram as above (without assuming $i$ nor $g$ to be in
\scrWA), $X'$ can be interpreted up to weak equivalence as the
``integral'' of the diagram
\[\begin{tikzcd}[baseline=(O.base),sep=small]
  Y\ar[r] \ar[d] & X \\
  |[alias=O]| Y' &
\end{tikzcd},\]
more accurately, the natural map in \Cat
\[\int \begin{tikzcd}[cramped,sep=small]
  A_{/Y}\ar[r]\ar[d] & A_{/X} \\ A_{/Y'} &
\end{tikzcd} \longrightarrow\; A_{/X'}\]
is in \scrW.

The corresponding statement in \Cat{} itself, even if $i$ is
supposed\pspage{205} to be an open immersion in \Cat{} say (or a
closed immersion, which amounts to the same by duality), is definitely
false, in other words: while the modelizer \Cat{} allows for a
remarkably simple description of integration of homotopy types, as
seen in the previous section, in the basic case however of an
integration\scrcomment{this integral is set inline in the typescript}
\[\int \begin{tikzcd}[cramped,sep=small]
  A_{/Y}\ar[r,"i"]\ar[d,"g"] & A_{/X} \\ A_{/Y'} &
\end{tikzcd}\]
(corresponding  to ``amalgamated sums''), and even when $i:Y\to X$ is
an open or closed immersion, this operation does definitely \emph{not}
correspond to the operation of taking the amalgamated sum $X'$ in
\Cat. It does though when $X\to X'$ is smooth resp.\ proper, for
instance if it is a fibration resp.\ a cofibration, and this in fact
implies the positive result in \Ahat{} noted above. This condition
moreover is satisfied if $g: Y\to Y'$ is equally an open resp.\ closed
immersion, in which case the situation is just the one of an ambient
$X'$, and two open resp.\ closed subobjects $X$ and $Y'$, with
intersection $Y$. This is a useful result, but wholly insufficient for
the factorization problem we were after in \Cat, with a view of
performing the standard homotopy constructions in \Cat{} itself. It
may be true still that if $i:Y\to X$ is not only an open or closed
immersion, but a weak equivalence as well, that then $i':Y'\to X'$ is
equally a weak equivalence, or what amounts to the same, that $X'$ can
be identified up to weak equivalence with the homotopy integral (which
indeed, up to weak equivalence, just reduces to $Y'$); but I have been
unable so far to clear up this matter. If true, this would be quite a
useful result, but still insufficient it seems in order to carry out
the standard homotopy constructions in \Cat{} itself.

To sum up, the main drawback of \Cat{} as a modelizer is that, except
in very special cases which are just too restricted, amalgamated sums
in \Cat{} don't have a reasonable meaning in terms of homotopy
operations -- whereas in a category \Ahat, a topos indeed where
therefore amalgamated sums have as good exactness properties as if
working in \Sets, these amalgamated sums (for a two-arrow diagram with
one arrow a monomorphism) \emph{do} have a homotopy-theoretic
meaning. This finally seems to force us, in order to develop some of
the basic structure in $\HotOf(\scrW)$, to leave the haven of the
basic modelizer \Cat, and work in an elementary modelizer \Ahat{}
instead, where $A$ is some \scrW-test category. This then brought me
back finally to the question whether these modelizers are \emph{closed
  model categories} in Quillen's sense, when we take of course for
``weak equivalences'' \scrWA, and moreover as ``cofibrations'' (in the
sense of Quillen's set-up) just the monomorphisms. Relying heavily
upon the result on monomorphisms in \Ahat{} stated above, it seems to
come out that we do get a closed model category indeed -- and even a
simplicial model category, if we are out for this. There is still a
cardinality question to be\pspage{206} settles to get the Quillen
factorizations in the general case, but this should not be too serious
a difficulty I feel. What however makes me still feel a little unhappy
in all this, is rather that I did not get a direct proof for an
elementary modelizer being a closed model category -- I finally have
to make a reduction to the known case of semisimplicial complexes,
settles by Quillen in his notes. This detour looks rather artificial
-- it is the first instance, and presumably the last one, where the
theory I am digging out seems to depend on semi-simplicial techniques,
techniques which moreover I don't really know and am not really eager
to swallow. It's just a prejudice maybe, a block maybe against the
semi-simplicial approach which I never really liked nor assimilated --
but I do have the feeling that the more refined and specific
semi-simplicial techniques and notions (such as minimal fibrations,
used in Quillen's proof, alas!) are irrelevant for an understanding of
the main structures featuring homotopy theory and homotopical
algebra. As for the notion of a Kan complex or a Kan fibration --
namely just a ``fibration'' in Quillen's axiomatic set-up, which I was
finally glad to find, ``ready for use'' -- I came to convince myself
at last that it was a basic notion indeed, and it was no use trying to
bypass it at all price. Thus, I took to the opposite, and tried to pin
down a Quillen-type factorization theorem, and his characteristic
seesaw game between right and left lifting properties, in as great
generality as I could manage.

\bigbreak
\presectionfill\ondate{30.5.}\par

\hangsection[Working program and rambling questions (group objects
\dots]{Working program and rambling questions \texorpdfstring{\textup
    (}{(}group objects as models, Dold-Puppe
  theorem\ldots\texorpdfstring{\textup )}{)}.}\label{sec:71}%
The scratchwork done since last month has of course considerably
cleared up the prospects of my present pre-stacks reflection on
homotopy models, on which I unsuspectingly embarked three months
ago. I would like to sketch a provisional working program for the
notes still ahead.

\namedlabel{it:71.A}{A)} Write down at last the story of asphericity
structures and canonical modelizers, as I was about to when I
interrupted the notes to do my scratchwork.

\namedlabel{it:71.B}{B)} Study the basic modelizer \Cat, and the
common properties of elementary modelizers \Ahat, with a main emphasis
upon base change and co-base change properties, and upon Quillen-type
factorization questions. Here it will be useful to dwell somewhat on
the ``homotopy integral'' variant of taking amalgamated sums in \Cat,
on the analogous constructions for topoi, and how these compare to the
usual amalgamated sums, including the interesting case of topological
spaces. It turns out that the homotopy\pspage{207} integral variant
for amalgamated sums is essentially characterized by a Mayer-Vietoris
type long exact sequence for cohomology, and the cases when the
homotopy construction turns out to be equivalent to usual amalgamated
sums, are just those when such a Mayer-Vietoris sequence exists for
the latter. An interesting and typical case is for topological spaces,
taking the amalgamated sum for a diagram
\[\begin{tikzcd}[baseline=(O.base),sep=small]
  Y\ar[r,hook,"i"] \ar[d,swap,"g"] & X \\
  |[alias=O]| Y' &
\end{tikzcd},\]
when $i$ is a closed immersion and $g$ is proper (which is also the
basic type of amalgamations which occur in the ``unfolding'' of
stratified structures).

In the course of the last weeks' reflections, there has taken place
also a substantial clarification concerning the relevant properties of
a basic localizer \scrW{} and how most of these, including strong
saturation of \scrW, follow from just the first three (a question which
kept turning up like a nuisance throughout the notes!). This should be
among the very first things to write down in this part of the
reflections, as \scrW{} after all is \emph{the} one axiomatic data
upon which the whole set-up depends.

\namedlabel{it:71.C}{C)}\enspace A reflection on the main common features of
the various contexts met with so far having a ``homotopy theory''
flavor, with a hope to work out at least some of the main features of
an all-encompassing new structure, along the lines of Verdier's
(commutative) theory of derived categories and triangulated
categories. The basic idea here, for the time being, seems to be the
notion of a \emph{derivator}, which should account for all the kind of
structure dealt with in Verdier's set-up, as well as in Deligne's and
Illusie's later elaborations. There seems to be however some important
extra features which are not accounted for by the mere derivator, such
as external $\Hom$'s with values in \Hot{} or some closely related
category, and the formalism of basic invariants (such as $\pi_i$,
$\mathrm H_i$ or $\mathrm H^i$), with values in suitable categories
(often abelian ones), which among others allow to check weak
equivalence. Such features seem to be invariably around in all cases I
know of, and they need to be understood I feel.

Coming back to \Hot{} itself and to modelizers $(M,W)$ giving rise to
it, there is the puzzling question about when exactly can we assert
that taking group objects of $M$ and weak equivalences between these
(namely group-object homomorphisms which are also in $W$), we get by
localization a category equivalent to the category of pointed
$0$-connected homotopy types. This is a well-known basic fact when we
take as models semisimplicial complexes or (I guess) topological
spaces -- a fact closely connected to the game of
associating\pspage{208} to any topological group its
``\emph{classifying space}'', defined ``up to homotopy''. I suspect
the same should hold in any elementary modelizer \Ahat, $A$ a test
category, at least in the ``strict'' case, namely when \Ahat{} is
totally aspheric, i.e., the canonical functor $\Ahat\to\Hot$ commutes
to finite products. The corresponding statement for the basic
localizer \Cat{} itself is definitely false. Group objects in \Cat{}
are indeed very interesting and well-known beings (introduced, I
understand from Ronnie Brown, by Henry Whitehead long time ago, under
the somewhat misleading name of ``crossed modules''), yet they embody
not arbitrary pointed $0$-connected homotopy types $X$, but merely those
for which $\pi_i(X)=0$ for $i>2$. Thus, we get only $2$-truncated
homotopy types -- and presumably, starting with $n\mathrm{-Cat}$
instead of \Cat{} as a modelizer, we then should get $(n+1)$-truncated
homotopy types. This ties in with the observation that taking group
objects either in \Cat, or in the full subcategory
$(\mathrm{Groupoids})$ of the latter, amounts to the same -- and
similarly surely for $(n\mathrm{-Cat})$; on the other hand it has been
kind of clear from the very beginning of this reflection that at any
\emph{finite} level, groupoids and $n$-groupoids ($n$ finite) will
only yield \emph{truncated} homotopy types.

A related intriguing question is when exactly does a modelizer $(M,W)$
give rise to a \emph{Dold-Puppe theorem} -- namely when do we get an
actual \emph{equivalence} between the category of \emph{abelian} group
objects of $M$, and the category of chain complexes of abelian groups?
The original statement was in case $M=\Simplexhat=$ semisimplicial
complexes, and doubtlessly it was one main impetus for the sudden
invasion of homotopy and cohomology by semisimplicial calculus -- so
much so it seems that for many people, ``homotopy'' has become
synonymous to ``semisimplicial algebra''. The impression that
semi-simplicial complexes is the God-given ground for doing homotopy
and even cohomology, comes out rather strong also in Quillen's
foundational notes, and in Illusie's thesis. Still, there are too some
cubical theory chaps I heard, who surely must have noticed long ago
that the Dold-Puppe theorem is valid equally for cubical complexes (I
could hardly imagine that it possibly couldn't). Now it turns out that
semisimplicial and cubical complexes are part of a trilogy, together
with so-called ``\emph{hemispherical complexes}'',\scrcomment{I guess
  today we call these ``globular sets''} which look at lot simpler
still, with just two boundary operators and one degeneracy in each
dimension. They can be viewed as embodying the ``primitive structure''
of an \oo-groupoid, the boundary operators being the ``source'' and
``target'' maps, and the degeneracy the map associating to any
$i$-object the corresponding ``identity''. I hit upon this structure
among the first examples of test categories and elementary modelizers,
and have been told since by Ronnie Brown that he has already known for
a while these models, under\pspage{209} the similar name of ``globular
complexes''). Roughly speaking, it can be said that the three types of
complexes correspond to three ``series'' of \emph{regular cellular
  subdivisions} of all spheres $S_n$ where the two-dimensional pieces
are respectively (by increasing order of ``intricacy'') bigons,
triangles, and squares. It shouldn't be hard to show these are the
only series of regular cellular subdivisions of all spheres (one in
each dimension) such that for any cell of such subdivision, the
induced subdivision of the bounding sphere should still be in the
series (up to isomorphism). The existence moreover of suitable
``degeneracy'' maps, which merit a careful general definition in this
context of cellular subdivisions of spheres, is an important common
extra feature of the three basic contexts, whose exact significance I
have not quite understood still. To come back to Dold-Puppe, sure
enough it is still valid in the hemispherical context. Writing down
the equivalence of categories in explicit terms comes out with
baffling simplicity. I wrote it down without even looking for it, in a
letter to Ronnie Brown, while explaining in a PS the ``yoga'' of
associating to a chain complex of abelian groups an \oo-groupoid with
additive structure (a so-called ``Picard category'' but within the
context of \oo-categories or \oo-groupoids, rather than usual
categories).

On the other hand, taking multicomplexes instead of simple ones, which
still can be interpreted as working in a category \Ahat{} for a
suitable test category $A$ (namely, a product category of categories
of the types $\Square$, $\Simplex$, and $\Globe$),
it is clear that Dold-Puppe's theorem as originally stated is no
longer true in these: in such case the category of abelian group
objects of \Ahat{} is equivalent to a category of \emph{multiple}
chain complexes. This shows that definitely, among all possible
elementary modelizers, the three in our trilogy are distinguished
indeed, as giving rise to a Dold-Puppe theorem. A thorough
understanding of this theorem would imply, I feel, an understanding of
which exactly are the modelizers, or elementary modelizers at any rate
giving rise to such a theorem. I wouldn't be too surprised if it
turned out that the three we got are the only ones, up to equivalence.

Another common feature of these modelizers, is that they allow for a
sweeping computational description of cohomology (or homology)
invariants in terms of the so-called ``boundary operations''. This is
visibly connected to the (strongly intuitive) tie between these kinds
of models, and cellular subdivisions of spheres. However, at this
level (unlike the Dold-Puppe story) the regularity feature of the
cellular subdivisions we got, and the fact that we allow for just one
(up to isomorphism) in each dimension,\pspage{210} seems to be
irrelevant. It might be worth while to write down with care what
exactly is needed, in order to define, in terms of a bunch of cellular
structures of spheres, a corresponding test-category (hopefully even a
strict one), and a way of computing in terms of boundary operators the
cohomology of the corresponding models. The more delicate point here
may be that if we really want to get an actual test category, not just
a weak one, we should have ``enough'' degeneracy maps between our
cellular structures, which might well prove an extremely stringent
requirement. Again, it will be interesting if we can meet it otherwise
than just sticking to our trilogy.

\namedlabel{it:71.D}{D)} ``Back to topoi''. They have been my main
intuitive leading thread in the reflections so far, but have remained
somewhat implicit most times. When working with categories $X$ as
models for homotopy types, we have been thinking in reality of the
associated topos $X\uphat$. In the same way, when relativizing over an
object $I$ of \Cat{} the construction of \Hot{} as a derived category,
namely working with categories over $I$ as models for a derived
category $D(I)$, the leading intuition again has been to look at $X$
over $I$ as the topos $X\uphat$ over $I\uphat$. This gives strong
suggestions as to defining $D(I)$ not only in terms of a given
category $I$, but also in terms of an arbitrary topos $T$ as well,
standing for $I\uphat$, and look to what extent the $f^*, f_!, f^*$
formalism of derivators extends to the case when $f$ is a ``map''
between topoi. The definition of some $D(T)$ for a topos $T$ can be
given in a number of ways it seems. As far as I know, the only one
which has been written and used so far consists in taking a suitable
derived category of the category of semisimplicial objects of
$T$. this is done in Illusie's thesis, where there is no mention
though of $f_!$ and $f_*$ functors -- which, maybe, should rather be
written $\mathrm Lf_!$ and $\mathrm Rf_*$, suggesting they are
respectively left and right derived functors of the more familiar
$f_!$ and $f_*$ functors for sheaves (the left and right adjoints of
the inverse image functor $f^*$ for sheaves of sets). One would
expect, at best, $\mathrm Lf_1$ to exist when $f_!$ itself does,
namely when $f^*$ commutes not only with finite inverse limits, but
with infinite products as well. As for $\mathrm Rf_*$, whereas there
is no problem for the existence of $f_*$ itself, already in the case
of a morphism of topoi coming from a map in \Cat, a map say between
finite ordered sets, the existence of $\mathrm Rf_*$ has still to be
established. Thus, presumably I'll content myself with writing down
and comparing a few tentative definitions of $D(T)$ and make some
reasonable guesses as to its variances. Intuitively, the objects of
$D(T)$ may be viewed as ``sheaves of homotopy types over $T$'', or
``relative homotopy types over $T$'', or ``non-commutative chain
complexes over $T$ up to quasi-isomorphism''. As in the case when $T$
is\pspage{211} the final (one point) topos, $D(T)$ is just the
homotopy category \Hot, there must be of course a vast variety of ways
of defining $D(T)$ in terms of model categories, and I would like to
review some which seem significant. In fact, one cannot help but
looking at two mutually dual groups of models categories, giving rise
to (at least) two definitely non-equivalent derived categories $D(T)$
and $D'(T)$ say, which, in case $T=I\uphat$, would correspond to
$D(I)$ and $D(I\op)$. A typical model category for the former is made
up with Illusie's semisimplicial sheaves; a typical model category for
$D'(T)$, on the other hand, should be made up with $1$-stacks on $T$
(in Giraud's sense), for a suitable notion of weak equivalence
between these.

Maybe the most natural models of all, in this context, for ``relative
homotopy types over the topos $T$'', should be \emph{topoi} $X$ over
$T$ (generalizing the categories over an ``indexing category''
$I$). The only trouble with this point of view though is that the best
we can hope for, in term's of Illusie's $D(T)$ say, is that a topos
$X$ over $T$ gives rise to a pro-object of $D(T)$, which needs not
come from an object of $D(T)$ itself, i.e., is not necessarily
``essentially constant''. This brings us back to the simpler and still
more basic question of associating a pro-homotopy type, namely a
pro-object of \Hot, to any topos $X$ -- namely back to the \v
Cech-Verdier-Artin-Mazur construction. This has been handled so far
using the semisimplicial models for \Hot, I suspect though that using
\Cat{} as a modelizer will give a more elegant treatment, as already
suggested earlier. Whichever way we choose to get the basic functor
\[\mathrm{Topoi} \to \Pro\Hot,\]
this functor will allow us, given a basic localizer \scrW, to define
\scrW-equivalences between topoi as maps which become isomorphisms
under the composition of the basic functor above, and the canonical
functor
\[\Pro\Hot \to \Pro(\HotOf(\scrW))\]
deduced from the localization functor
\[\Hot\to(\HotOf(\scrW)).\]
(It turns out, using Quillen's theory, that usual weak equivalences is
indeed the finest of all possible basic localizers \scrW, hence
$\HotOf(\scrW)$ is indeed a localization of \Hot.)

Maybe it is not too unreasonable to expect that all, or most homotopy
constructions, involving a topos $T$, can be expressed replacing
$T$ by its image in $\Pro\Hot$, or in $\Pro(\HotOf(\scrW))$ if the
construction are relative to a given basic localizer \scrW. The very
first example one would like\pspage{212} to look up in this respect,
is surely $D(T)$, defined say à la Illusie.

In the present context, the main point of the property of \emph{local
  asphericity} for a topos (cf.\ section \ref{sec:35}) is that the
corresponding prohomotopy type is essentially constant, i.e., the
topos defines an actual homotopy type. Thus, locally aspheric topoi
and weak equivalences between these should be eligible models for
homotopy types (more accurately, make up a modelizer), just as the
basic modelizer \Cat{} contained in it. The corresponding statements,
when introducing a basic localizer \scrW, should be equally valid. One
might expect, too, a relative variant for the notion of local
asphericity or \scrW-asphericity, in case of a topos $X$ over a given
base topos $T$ (which we may have to suppose already locally
aspheric), with the implication that the corresponding object of
$\Pro(D(T))$ should be again essentially constant.

A last question I would like to mention here is about the meaning of
the notion of a so-called ``modelizing topos'', introduced in a
somewhat formal way in section \ref{sec:35}, as a locally aspheric and
aspheric topos $T$ such that the Lawvere element $L_T$ of $T$ is
aspheric over the final object. (We assumed at first, moreover, that
$T$ be even totally aspheric, but soon after the point of view and
terminology shifted a little and the totally aspheric case was
referred to as a \emph{strictly} modelizing topos, cf.\ page
\ref{p:68}. As made clear there, the expectation suggested by the
terminology is of course that such a topos should indeed be a
modelizer, when endowed with the usual notion of weak
equivalence. This sill makes sense and seem plausible enough, when
usual weak equivalence is replaced by the weaker notion defined in
terms of an arbitrary basic localizer \scrW. A related question is to
get a feeling for how restrictive the conditions put on a modelizing
topos are. It is clear that nearly all topoi met with in practice,
including those associated to the more common topological spaces (such
as locally contractible ones) are locally aspheric -- but what about
the condition on the Lawvere element? For instance, taking a
topological space admitting a finite triangulation and which is
aspheric, i.e., contractible, is the corresponding topos modelizing?

\bigbreak

\presectionfill\ondate{10.6.}\pspage{213}\par

\hangsection{Back to asphericity: criteria for a map in
  \texorpdfstring{\Cat}{(Cat)}.}\label{sec:72}%
In the long last, we'll come back now to the ``asphericity game''!
Let's take up the exposition at the point where we stopped two months
ago (section \ref{sec:67}, p.\ \ref{p:188}). We were then about to
reformulate in various ways the property of asphericity, more
specifically \scrW-asphericity, for a given map
\[ i :A\to B\]
in \Cat. To this end, we introduced the corresponding diagram of maps
in \Cat{} with ``commutation morphism'' $\lambda_i$:
\[\begin{tikzcd}[baseline=(O.base),column sep=small]
  \Bhat\ar[rr,"i^*"]\ar[dr,swap,"i_B",""{name=iB,right}] & &
  \Ahat\ar[dl,"i_A",""{name=iA,left}] \\
  & |[alias=O]| \Cat \arrow[bend right=10, from=iA, to=iB, swap, "\lambda_i"] &
\end{tikzcd}.\]
As already stated, the asphericity condition on $i$ can be expressed
in a variety of ways, as a condition on either of the three
``aspects''
\[\lambda_i,\quad i^*,\quad i_Ai^*\]
of the situation created by $i$, with respect to the localizing sets
\scrW, \scrWA, $\scrW_B$ in the three categories under consideration,
or with respect to the notion of \scrW-aspheric objects in these. As
\scrWA, and aspheric objects in \Ahat, are defined respectively in
terms of \scrW{} and aspheric objects in \Cat{} via the functor $i_A$,
it turns out that the formulations in terms of properties of $i_Ai^*$
reduce trivially to the corresponding formulations in terms of $i^*$,
which will be the most directly useful for our purpose for later work
-- hence well's omit them, and focus attention instead on $\lambda_i$
and $i^*$. Notations are those of loc.\ sit., in particular, we are
working with the categories
\[ \mathrm{Hot}_A^\scrW \quad\text{or} \quad
\mathrm{Hot}_A \eqdef \scrWA^{-1}\Ahat\]
and
\[\HotOf(\scrW) = \scrW^{-1}\Cat,\]
where \scrW{} is a given ``basic localizer'', namely a set of arrows
in \Cat{} satisfying certain conditions.

It seems worthwhile here to be careful to state which exactly are the
properties of \scrW{} we are going to use in the ``asphericity game''
-- therefore I'm going to list them again here, using a labelling
which, hopefully, will not have to be changed again:
\begin{description}
\item[\namedlabel{loc:1}{Loc~1)}]
  ``Mild saturation'' (cf.\ page \ref{p:59}),
\item[\namedlabel{loc:2}{Loc~2)}]
  ``Homotopy condition'': $\Simplex_1\times X\to X$ is in \scrW{} for
  any $X$ in \Cat,
\item[\namedlabel{loc:3}{Loc~3)}]
  ``Localization\pspage{214} condition'': If $X,Y$ are objects in \Cat{} over an
  object $A$,\scrcomment{AG seems to have replaced $A$ with $\mathscr
    S$ here, but it's not clear\ldots} and $u: X\to Y$ an $A$-morphism
  such that for any $a$ in $A$, the induced $u_{/a}:X_{/a}\to Y_{/a}$
  is in \scrW, then so is $u$.
\end{description}

It should be noted that the ``mild saturation'' condition is slightly
weaker than the saturation condition introduced later (p.\
\ref{p:101}, conditions
\ref{it:48.aprime}\ref{it:48.bprime}\ref{it:48.cprime}), namely in
condition \ref{it:48.cprime} (if $f:X\to Y$ and $g:Y\to X$ are such
that $gf,fg\in\scrW$, then $f,g\in\scrW$) we restrict to the case when
$gf$ is the identity, namely $f$ an inclusion and $g$ a retraction
upon the corresponding subobject. On the other hand, in what follows
we are going to use \ref{loc:3} only in case $Y=A$ and $Y\to A$ is the
identity. I stated the condition in greater generality than needed for
the time being, in view of later convenience -- as later it will have
to be used in full strength. It is not clear whether the weaker form
of \ref{loc:3} (plus of course \ref{loc:1} and \ref{loc:2}) implies
already the stronger. We'll see later a number of nice further
properties of \scrW{} implied by these we are going to work with for
the time being -- including \emph{strong} saturation of \scrW, namely
that \scrW{} is the set of arrows made invertible by the localization
functor
\[\Cat\to\scrW^{-1}\Cat=\HotOf(\scrW).\]
\begin{propositionnum}\label{prop:72.1}
  Let as above $i:A\to B$ be a map in \Cat. Consider the following
  conditions on $i$:
  \begin{description}
  \item[\namedlabel{it:72.1.i}{(i)}]
    For any $F$ in \Bhat, $\lambda_i(F): A_{/F} \to B_{/F}$ is in \scrW.
  \item[\namedlabel{it:72.1.iprime}{(i')}]
    For any $F$ as above, $\lambda_i(F)$ is \scrW-aspheric
    \textup(i.e., satisfies the assumption on $u$ in
    \textup{\ref{loc:3}} above, when $Y\to A$ is an identity\textup).
  \item[\namedlabel{it:72.1.idblprime}{(i'')}]
    Same as \textup{\ref{it:72.1.i}}, but restricting to $F=b$ in $B$.
  \item[\namedlabel{it:72.1.ii}{(ii)}]
    For any $F$ in \Bhat, $F$ \scrW-aspheric $\Rightarrow i^*(F)$
    \scrW-aspheric.
  \item[\namedlabel{it:72.1.iiprime}{(ii')}]
    For any $F$ in \Bhat, $F$ \scrW-aspheric $\Leftrightarrow i^*(F)$
    \scrW-aspheric.
  \item[\namedlabel{it:72.1.iii}{(iii)}]
    For any $b$ in $B$, $A_{/b} (\eqdef A_{/i^*(b)})$ is
    \scrW-aspheric, i.e., $i$ is \scrW-aspheric.
  \item[\namedlabel{it:72.1.iv}{(iv)}]
    For any map $f$ in \Bhat, $f\in\scrW_B \Rightarrow i^*(f)\in\scrWA$.
  \item[\namedlabel{it:72.1.ivprime}{(iv')}]
    For any map $f$ as above, $f\in\scrW_B \Leftrightarrow i^*(f)\in\scrWA$.
  \item[\namedlabel{it:72.1.v}{(v)}]
    Condition \textup{\ref{it:72.1.iv}} holds, i.e., $i^*$ induces a
    functor
    \[\overline{i^*} : \mathrm{Hot}_B \to \mathrm{Hot}_A,\]
    and moreover the latter is an \emph{equivalence}.
  \end{description}
  The conditions \textup{\ref{it:72.1.i}} up to
  \textup{\ref{it:72.1.iii}} are all equivalent, call this set of
  conditions \textup{\namedlabel{cond:72.As}{(As)}}
  \textup(\scrW-asphericity\textup). We moreover have the following
  implications between \textup{\ref{cond:72.As}} and the remaining
  conditions \textup{\ref{it:72.1.iv}} to \textup{\ref{it:72.1.v}}:
  \begin{equation}
    \label{eq:72.1.star}
    \begin{tabular}{@{}c@{}}
      \begin{tikzcd}[baseline=(O.base),math mode=false,%
        arrows=Rightarrow,column sep=-3pt]
        \textup{\ref{cond:72.As}}
        \ar[rrrr,dashed,"{if $A,B$ ps.test}"{inner sep=3pt}]\ar[dr] & & & &
        \textup{\ref{it:72.1.v}} \ar[ddll] \\
        & \textup{\ref{it:72.1.ivprime}} \ar[dr] & & \phantom{hello} & \\
        & & |[alias=O]| \textup{\ref{it:72.1.iv}}
        \ar[uull,dashed,bend left,%
        "{\begin{tabular}{@{}l@{}}
            if $A,B$ \\ \scrW-asph.
          \end{tabular}}"{inner sep=-3pt,near end}] & &
      \end{tikzcd},
    \end{tabular}\tag{*}
  \end{equation}
  where\pspage{215} the implication \textup{\ref{it:72.1.iv} $\Rightarrow$
    \ref{cond:72.As}} is subject to $A,B$ begin \scrW-aspheric, and
  \textup{\ref{cond:72.As} $\Rightarrow$ \ref{it:72.1.v}} to $A,B$
  being ``pseudo-test categories'' \textup(for \scrW\textup), namely the canonical
  functors
  \[\overline{i_A}: \mathrm{Hot}_A \to \HotOf(\scrW), \quad
  \overline{i_B} : \mathrm{Hot}_B \to \HotOf(\scrW)\]
  being equivalence.
\end{propositionnum}
\begin{corollary}
  Assume that $A$ and $B$ are pseudo-test categories, and are
  \scrW-aspheric. Then all conditions \textup{\ref{it:72.1.i}} to
  \textup{\ref{it:72.1.v}} of the proposition above are equivalent.
\end{corollary}
\begin{remarknum}
  If we admit strong saturation of \scrW{} (which will be proved
  later), it follows at once that a pseudo-test category is
  necessarily \scrW-aspheric -- hence the conclusion of the corollary
  holds assuming only $A$ and $B$ are pseudo-test categories. Of
  course, it holds a fortiori if $A$ and $B$ are weak test categories,
  or even test categories. Also, strong saturation implies that in
  \eqref{eq:72.1.star} above, we have even the implication:
  \ref{it:72.1.v} $\Rightarrow$ \ref{it:72.1.ivprime}, stronger than
  \ref{it:72.1.v} $\Rightarrow$ \ref{it:72.1.iv}.
\end{remarknum}

\begin{proof}[Proof of proposition]
  It is purely formal -- for the first part, it follows from the
  diagram of tautological implications
  \[\begin{tikzcd}[baseline=(O.base),math mode=false,%
    arrows=Rightarrow,column sep=0pt]
    & \ref{it:72.1.iprime} \ar[d] & & $\phantom{\text{(i')}}$ & \\
    & \ref{it:72.1.i} \ar[dl] \ar[dr] & & & \\
    \ref{it:72.1.iiprime}\ar[d] & &
    \ref{it:72.1.idblprime}\ar[d,Leftrightarrow] & & \\
    \ref{it:72.1.ii}\ar[rr] & & \ref{it:72.1.iii} \ar[rr] & &
    |[alias=O]| \ref{it:72.1.iprime}
  \end{tikzcd},\]
  where the implication \ref{it:72.1.iprime} $\Rightarrow$
  \ref{it:72.1.i} is contained in the assumption
  \hyperref[loc:3]{Loc~3} on \scrW. The implications of the diagram
  \eqref{eq:72.1.star} are about as formal -- there is no point I
  guess writing it out here.
\end{proof}
\begin{remarknum}
  Using the canonical functors
  \[\varphi_A: \Ahat\to\HotOf(\scrW), \quad
  \varphi_B:\Bhat\to\HotOf(\scrW),\]
  there are two other amusing versions still of \ref{it:72.1.i}, which
  look a lot weaker still, and are equivalent however to
  \ref{it:72.1.i}, i.e., to \scrW-asphericity of $i$, namely:
  \begin{description}
  \item[\namedlabel{it:72.1.vi}{(vi)}]
    For any $F$ in \Bhat, $\varphi_B(F)=\gamma(B_{/F}$ and
    $\varphi_A(i^*(F))=\gamma(A_{/F})$ are isomorphic objects of
    $\HotOf(\scrW)$.
  \item[\namedlabel{it:72.1.viprime}{(vi')}]
    Same\pspage{216} as \ref{it:72.1.vi}, with $F$ restricted to be
    object $b$ in $B$.
  \end{description}
  Indeed, we have of course \ref{it:72.1.i} $\Rightarrow$
  \ref{it:72.1.vi} $\Rightarrow$ \ref{it:72.1.viprime}, but also
  \ref{it:72.1.viprime} $\Rightarrow$ \ref{it:72.1.iii} if we admit
  strong saturation of \scrW, which implies that an object $X$ of
  \Cat{} is \scrW-aspheric if{f} its image in $\HotOf(\scrW)$ is a
  final object.
\end{remarknum}
\begin{remarknum}
  If we don't assume $A$ and $B$ to be \scrW-aspheric, the
  implications
  \[ \text{\ref{cond:72.As}} \Rightarrow \text{\ref{it:72.1.ivprime}}
    \Rightarrow \text{\ref{it:72.1.iv}}\]
  are both strict. As an illustration of this point, take for $B$ a
  discrete category (defined in terms of a set of indices $I=\Ob B$),
  thus a category $A$ over $B$ is essentially the same as a family
  $(A_b)_{b\in I}$ of objects of \Cat{} indexed by $I$. In terms of
  this family, we see at once that the three conditions above on
  $i:A\to B$ mean respectively (a)\enspace that all categories $A_b$
  are \scrW-aspheric, (b)\enspace that all categories $A_b$ are
  non-empty, and (c)\enspace condition vacuous. This example brings to
  mind that \emph{condition \textup{\ref{it:72.1.iv}} is a
    Serre-fibration type condition}, we'll come back upon this
  condition when studying homotopy properties of maps in \Cat, with
  special emphasis on base change questions. Likewise, condition
  \ref{it:72.1.ivprime} appears as a strengthening of such Serre-type
  condition, to the effect that the restriction of $A$ over any
  connected component of $B$ should be moreover non-empty. As an
  example, we may take the projection $B\times C\to C$, where $C$ is
  any non-empty object in \Cat.
\end{remarknum}

I would like now to dwell still a little on the case, of special
interest of course for the modelizing story, when $A$ and $B$ are test
categories or something of the kind. The formulation \ref{it:72.1.i}
of the asphericity condition for the functor $i$ can be expressed by
stating that the functor $\overline{i^*}$ between the localizations
$\mathrm{Hot}_B$ and $\mathrm{Hot}_A$ exists, and gives rise (via
$\lambda_i$) to a commutation morphism which is an \emph{isomorphism}
\[\begin{tikzcd}[baseline=(O.base),column sep=small]
  \mathrm{Hot}_B\ar[rr,"\overline{i^*}"]\ar[dr,swap,"\overline{i_B}",""{name=iB,right}] & &
  \mathrm{Hot}_A\ar[dl,"\overline{i_A}",""{name=iA,left}] \\
  & |[alias=O]| \Cat \arrow[bend right=10, from=iA, to=iB, swap,
  "\overline{\lambda_i}"{inner sep=4pt}, "\sim"{inner sep=0pt}] &
\end{tikzcd}.\]
This shows, when $A$ and $B$ are pseudo-test categories for \scrW,
i.e., the functors $\overline{i_A}$ and $\overline{i_B}$ are
equivalences, that the functor $\overline{i^*}$ deduced from the
\scrW-aspheric map $i:A\to B$, \emph{does not depend} (up to canonical
isomorphism) on the choice of $i$, and can be described as the
composition of $\overline{i_B}$ followed by a quasi-inverse for
$\overline{i_A}$. This of course is very much in keeping with the
``inspiring assumption'' (section \ref{sec:28}), which just means that
up to unique isomorphism, there is indeed but \emph{one}\pspage{217}
equivalence from $\mathrm{Hot}_B$ to $\mathrm{Hot}_A$ (both categories
being equivalent to $\HotOf(\scrW)$). Here we are thinking of course
of the extension of the ``assumption'' contemplated earlier, when
usual weak equivalence is replaced by a basic localizer \scrW{} as
above. It seems plausible that the assumption holds true in all cases
-- anyhow we didn't have at any moment to make explicit use of it
(besides at a moment drawing inspiration from it\ldots).
\begin{propositionnum}\label{prop:72.2}
  Let $i:A\to B$ be a map in \Cat.
  \begin{enumerate}[label=\alph*),font=\normalfont]
  \item\label{it:72.2.a}
    If $i$ is \scrW-aspheric, and \Ahat{} is totally \scrW-aspheric,
    then \Bhat{} is totally \scrW-aspheric too.
  \item\label{it:72.2.b}
    Assume $A$ \scrW-aspheric, and that $B$ is a \scrW-test category,
    admitting the separating \scrWB-homotopy interval
    $\bI=(I,\delta_0,\delta_1)$, satisfying the homotopy condition
    \textup{\ref{cond:TH1}} of page \textup{\ref{p:50}}. Consider the following
    conditions:
    \begin{enumerate}[label=\arabic*),font=\normalfont]
    \item\label{it:72.2.b.1}
      $i^*(I)$ is \scrW-aspheric over $e_\Ahat$,
    \item\label{it:72.2.b.2}
      $i$ is \scrW-aspheric,
    \item\label{it:72.2.b.3}
      $i^*(I)$ is \scrW-aspheric.
    \end{enumerate}
    We have the implications
    \[ \textup{\ref{it:72.2.b.1}} \Rightarrow
    \textup{\ref{it:72.2.b.2}} \Rightarrow
    \textup{\ref{it:72.2.b.3}},\]
    hence, if $A$ is totally \scrW-aspheric \textup(hence
    \textup{\ref{it:72.2.b.3}} $\Rightarrow$
    \textup{\ref{it:72.2.b.1})}, all three conditions are equivalent.
  \end{enumerate}
\end{propositionnum}
\begin{proof}
  \ref{it:72.2.a}\enspace We have to prove that if $b,b'$ are in $B$, then
  their product in \Bhat{} is \scrW-aspheric, but by assumption on $i$
  we know that the images of $b,b'$ by $i^*$ are \scrW-aspheric, hence
  (as \Ahat{} is totally aspheric) their product
  \[ i^*(b) \times i^*(b') = i^*(b\times b')\]
  is \scrW-aspheric too, hence so is $b\times b'$ by criterion
  \ref{it:72.1.iiprime} of prop.\ above.

  \ref{it:72.2.b}\enspace The homotopy condition \ref{cond:TH1} referred to
  means that all objects of $B$ are \bI-contractible. As $i^*$
  commutes with finite products, it follows that the objects $i^*(b)$
  of \Ahat, for $b$ in $B$, are $i^*(\bI)$-contractible. When $i^*(I)$
  is \scrW-aspheric over $e_\Ahat$, this implies that so are the
  objects $i^*(b)$, a fortiori they are \scrW-aspheric (as by
  assumption $e_\Ahat$ is \scrW-aspheric). Thus \ref{it:72.2.b.1}
  $\Rightarrow$ \ref{it:72.2.b.2}, and \ref{it:72.2.b.2} $\Rightarrow$
  \ref{it:72.2.b.3} is trivial.
\end{proof}
\begin{corollary}
  Let $i:A\to B$ be a \scrW-aspheric map in \Cat, assume that $A$ is
  totally \scrW-aspheric, and $B$ is a local \scrW-test category. Then
  both $A$ and $B$ are strict \scrW-test categories.
\end{corollary}

Indeed, by \ref{it:72.2.a} above we see that $B$ is totally
\scrW-aspheric, hence $B$ is a strict \scrW-test category. In order to
prove that so is $A$, we only have to show that \Ahat{} admits a
separating homotopy interval for \scrWA. By assumption\pspage{218} on
$B$, there is a separating \scrWB-homotopy interval
$\bI=(I,\delta_0,\delta_1)$ in \Bhat. The exactness properties of $i^*$
imply that $i^*(\bI)$ is a separating interval in \Ahat, the
asphericity condition on $i$ implies that moreover $i^*(I)$ is
\scrW-aspheric, hence \scrW-aspheric over $e_\Ahat$ as \Ahat{} is
totally \scrW-aspheric, qed.

To finish with the more formal properties of the notion  of
\scrW-aspheric maps in \Cat, let's give a list of the standard
stability conditions for this notion, with respect notably to
composition, base change, and cartesian products:

\begin{propositionnum}\label{prop:72.3}
  \textup{\namedlabel{it:72.3.a}{a)}}\enspace
  Consider two maps
  \[ A \xrightarrow i B \xrightarrow j C\]
  in \Cat. Then if $i$ and $j$ are \scrW-aspheric, so is $ji$. If $ji$
  and $i$ are \scrW-aspheric, so is $j$. Any isomorphism in \Cat{} is
  \scrW-aspheric.

  \textup{\namedlabel{it:72.3.b}{b)}}\enspace
  Let
  \[\begin{tikzcd}
    A \ar[d,swap,"i"] & A' \ar[l,swap,"f"]\ar[d,"i'"] \\
    B & B' \ar[l,swap,"g"]
  \end{tikzcd}\]
  be a cartesian square in \Cat, assume $i$ is \scrW-aspheric and $g$
  is fibering \textup(for instance an induction functor $B_{/F}\to B$,
  with $F$ in \Bhat\textup). Then $i'$ is \scrW-aspheric \textup(and,
  of course, $f$ is equally fibering\textup).
  
  \textup{\namedlabel{it:72.3.c}{c)}}\enspace
  Let
  \[i:A\to B, \quad i':A'\to B'\]
  be two \scrW-aspheric maps in \Cat, then
  \[ i\times i' : A\times A' \to B\times B'\]
  is \scrW-aspheric.
\end{propositionnum}
\begin{proof}
  Property \ref{it:72.3.a} is formal, in terms of criterion
  \ref{it:72.1.iiprime} of prop.\ \ref{prop:72.1}. Property
  \ref{it:72.3.c} follows formally from the criterion
  \ref{it:72.1.iii}, and the canonical isomorphism
  \[ (A\times A')_{/b\times b'} \simeq (A_{/b}) \times (A'_{/b'}),\]
  and the fact that a product of two \scrW-aspheric objects of \Cat{}
  is again \scrW-aspheric (cf.\ prop.\ of page \ref{p:167}, making use
  of the localization condition \ref{loc:3} on \scrW{} in its full
  generality). Another proof goes via \ref{it:72.3.a} and
  \ref{it:72.3.b}, by reducing first (using \ref{it:72.3.a}) to the
  case when either $i$ or $i'$ are identities, and using the fact that
  a projection map $C \times B \to B$ in \Cat{} is a fibration.

  We are left with proving property \ref{it:72.3.b}. For this, we note
  that the asphericity condition \ref{it:72.1.iii} on a map $i:A\to B$
  just means that for\pspage{219} any base change of the type
  \[ B_{/b} \to B,\]
  where $b$ is in $B$, the corresponding map
  \[i_{/b}: A\times_B B_{/b} \simeq A_{/b} \to B_{/b}\]
  is in \scrW. Applying this to the case of $i':A'\to B'$, and an
  object $b'$ in $B'$, and denoting by $b$ its image in $B$, using
  transitivity of base change, we get a cartesian square
  \[\begin{tikzcd}[baseline=(O.base)]
    A_{/b} \ar[d,"i_{/b}"] & A'_{/b'}\ar[l]\ar[d,"i'_{/b'}"] \\
    B_{/b} & |[alias=O]| B'_{/b'} \ar[l]
  \end{tikzcd},\]
  and we got to prove $i'_{/b'}$ is in \scrW, i.e., $A'_{/b'}$ is
  \scrW-aspheric, using the fact that we know the same holds for
  $i_{/b}$, i.e., $A_{/b}$ is \scrW-aspheric. Thus, all we have to
  prove is that the first horizontal arrow is in \scrW. But we check
  at once that the condition that $B'\to B$ is fibering implies that
  the induced functor
  \[ B'_{/b'} \to B_{/b}\]
  is fibering too, and moreover has fibers which have final
  objects. Hence by base change, the same properties hold for
  \[ A'_{/b'} \to A_{/b}.\]
  This reminds us of the ``fibration condition''
  \hyperref[it:64.L5]{L~5} (page \ref{p:164}), which should ensure
  that a fibration with \scrW-aspheric fibers is in \scrW. We did not
  include this axiom among the assumptions (recalled above) we want to
  make on \scrW. However, it turns out that the assumptions we do make
  here imply already the fibration condition, as well as the dual
  condition on cofibrations. We'll give a proof later -- in order not
  to diverge at present from our main purpose. There will not be any
  vicious circle, as all we're going to use of prop.\ \ref{prop:72.3}
  for the formalism of asphericity structures and canonical modelizers
  is the first part \ref{it:72.3.a}. I included \ref{it:72.3.b} and
  \ref{it:72.3.b} for the sake of completeness, and because
  \ref{it:72.3.c} is useful for dealing with products of two
  categories, notably of two test categories -- a theme which has been
  long pending, and on which I would like to digress next, before
  getting involved with asphericity structures.
\end{proof}
\begin{remarknum}
  In part \ref{it:72.3.a} of the proposition, if $j$ and $ji$ are
  \scrW-aspheric, we cannot conclude that $i$ is. If $C$ is the final
  object of \Cat, this means that a map between \scrW-aspheric objects
  in \Cat{} need not be \scrW-aspheric.
\end{remarknum}

\bigbreak

\presectionfill\ondate{11.6.}\pspage{220}\par

\hangsection{Asphericity criteria
  \texorpdfstring{\textup(}(continued\texorpdfstring{\textup)}).}%
\label{sec:73}%
I am not quite through yet with generalities on asphericity criteria
for a map in \Cat, it turns out -- it was just getting prohibitively
late last night to go on!

From now on, I'll drop the qualifying \scrW{} when speaking of
asphericity, test categories, modelizers and the like, as by now it is
well understood, I guess, there is a given \scrW{} around in all we
are doing. It'll be enough to be specific in those (presumably rare)
instances when working with more than one basic localizer.

Coming back to the last remark in yesterday's notes, a good
illustration is the case of a functor $i:A\to B$ of aspheric objects
of \Cat, when $A$ is a final object in \Cat, i.e., a one-point
discrete category. Then $i$ \emph{is aspheric if{f} $i(a)$ is an
  initial object of $B$} -- an extremely stringent extra condition
indeed!

In part \ref{it:72.3.a} of proposition \ref{prop:72.3} above, stating
that isomorphisms in \Cat{} are aspheric, it would have been timely to
be more generous -- indeed, \emph{any equivalence of categories in
  \Cat{} is aspheric}. This associates immediately with a map of topoi
which is an equivalence being (trivially so) aspheric (in the case,
say, when \scrW{} is the usual notion of weak equivalence, the only
one for the time being when the notion is extended from maps in \Cat{}
to maps between topoi). In the context of \Cat, the basic modelizer,
we can give a still more general case of aspheric maps, both
instructive and useful:
\addtocounter{propositionnum}{3}
\begin{propositionnum}\label{prop:73.4}
  Let
  \[\begin{tikzcd}
    A \ar[r,shift left,"f"] & B \ar[l,shift left,"g"] B
  \end{tikzcd}\]
  be a pair of adjoint functors between the objects $A,B$ in \Cat,
  with $f$ left and $g$ right adjoint. Then $f$ is aspheric.
\end{propositionnum}

Indeed, by the adjunction formula we immediately get
\[A_{/b} \simeq A_{/g(b)}\]
(as a matter of fact, $f^*(b)=g(b)$), hence this category has a final
object and hence is aspheric, qed.

\addtocounter{remarknum}{4}
\begin{remarknum}
  The conclusion of prop.\ \ref{prop:73.4} is mute about $g$, which
  will rightly strike as unfair. Dualizing, we could say that
  $(g\op,f\op)$ is a pair of adjoint functors between $A\op$ and
  $B\op$, and therefore $g\op$ is aspheric. We will express this fact
  (by lack of a more suggestive name) by saying that $g$ is a
  \emph{coaspheric} map. For a functor $i:A\to B$, in terms of the
  usual criterion for $i\op$, it just\pspage{221} means that for any
  $b$ in $B$, the category
  \[\preslice Ab \eqdef A \times_B(\preslice Bb)
  = \parbox[t]{0.5\textwidth}{category of pairs $(x,p)$, with $x$ in
    $A$ and $p:b\to i(x)$}\]
  is aspheric. (NB\enspace In the case of the functor $g$ from $B$ to
  $A$, the corresponding categories $\preslice Ba$ even have
  \emph{initial objects}.) This conditions comes in here rather
  formally, we'll see later though that it has a quite remarkable
  interpretation, in terms of a very strong property of cohomological
  ``cofinality'' of the functor $i$, implying the usual notion of $i$
  being a ``cofinal'' functor, or $A$ being ``cofinal'' in $B$ (namely
  $i$ giving rise to isomorphisms $\varinjlim_A{} \to \varinjlim_B{}$
  for any direct system $B\to M$ with values in a category $M$
  admitting direct limits), as the ``dimension zero'' shadow of this
  ``all dimensions'' property. To make an analogy which will acquire
  more precise meaning later, the asphericity property for a map $i$
  in \Cat{} can be viewed as a (slightly weakened) version of a
  \emph{proper map with aspheric fibers}, whereas the coasphericity
  property appears as the corresponding version of a \emph{smooth map
    with aspheric fibers}. The qualification ``slightly weakened''
  reflects notably in the fact that the notions of asphericity and
  coasphericity are \emph{not} stable under arbitrary base change --
  but rather, asphericity for a map is stable under base change by
  \emph{fibration functors} (more generally, by smooth maps in \Cat),
  whereas coasphericity is stable under base change by
  \emph{cofibration functors} (more generally, by proper maps in
  \Cat). Thus, whereas a functor $i:A\to B$ which is an equivalence of
  categories is clearly both aspheric and coaspheric, this property is
  not preserved by arbitrary base change, e.g., passage to fibers:
  e.g., some fibers may be empty, and therefore are neither aspheric
  nor coaspheric!
\end{remarknum}

The most comprehensive property for a map in \Cat, implying both
asphericity and coasphericity, is to be in \UW, i.e., a
``\emph{universal weak equivalence}'' -- namely it is in \scrW{} and
remains so after any base change. These maps deserve the name of
``\emph{trivial Serre fibrations}''. They include all maps with
aspheric fibers which are either ``proper'' or ``smooth'', for
instance those which are either cofibering or fibering functors (in
the usual sense of category theory). We'll come back upon these
notions in a systematic way in the next part of the reflections.

\hangsection{Application to products of test categories.}\label{sec:74}%
As announced yesterday, I would like still to make an overdue
digression on products of test categories, before embarking on the
notion of asphericity structure. It will be useful to begin with
some\pspage{222} generalities on presheaves on a product category
$A\times B$, where for the time begin $A,B$ are any two objects of
\Cat. The following notation is often useful, for a pair of presheaves
\[ F\in\Ob\Ahat, \quad G\in\Ob\Bhat,\]
introducing an ``\emph{external product}''
\[ F \boxtimes G\in \Ob(A\times B)\uphat\]
by the formula
\[\mathop{F \boxtimes G}(a,b) = F(a) \times G(b).\]
Introducing the two projections
\[ p_1 : A\times B\to A, \quad p_2:A\times B\to B,\]
and the corresponding inverse image functors $p_i^*$, we have
\[ F \boxtimes G = p_1^*(F) \times p_2^*(G).\]
Of course, $F\boxtimes G$ depends bifunctorially on the pair $(F,G)$,
and it is easily checked, by the way, that the corresponding functor
\[ \Ahat \times \Bhat \to (A\times B)\uphat, \quad
(F,G)\mapsto F\boxtimes G,\]
is fully faithful.

We have a tautological canonical isomorphism
\[ (A\times B)_{/{F \boxtimes G}} \tosim A_{/F} \times B_{/G},\]
and hence the
\begin{propositionnum}\label{prop:74.1}
  \textup{\namedlabel{it:74.1.a}{a)}}\enspace
  If $F$ and $G$ are aspheric objects in \Ahat{} and \Bhat{}
  respectively, then $F\boxtimes G$ is an aspheric object of $(A\times
  B)\uphat$.

  \textup{\namedlabel{it:74.1.b}{b)}}\enspace
  If $u: F\to F'$ and $v: G\to G'$ are aspheric maps in \Ahat{} and
  \Bhat{} respectively, then
  \[ {u\boxtimes v} : F \boxtimes G \to F' \boxtimes G'\]
  is an aspheric map in $(A\times B)\uphat$.
\end{propositionnum}

Indeed, this follows respectively from the fact that a product of two
aspheric objects (resp.\ aspheric maps) in \Cat{} is again aspheric.
\begin{corollarynum}\label{cor:74.1.1}
  Let $F$ in \Ahat{} be aspheric over the final object $e_\Ahat$, then
  $F \boxtimes e_\Bhat = p_1^*(F)$ is aspheric over the final object
  of $(A\times B)\uphat$.
\end{corollarynum}
\begin{corollarynum}\label{cor:74.1.2}
  Assume $A$ and $B$ are totally aspheric, then so is $A\times B$.
\end{corollarynum}

We have to check that the product elements
\[ (a,b) \times (a',b') = (a\times a') \boxtimes (b\times b')\]
are aspheric, which follows from the assumption (namely $a\times a'$
and $b\times b'$ aspheric) and part \ref{it:74.1.a} of the
proposition.

Assume\pspage{223} now $A$ is a local test category, i.e., \Ahat{}
admits a separating interval
\[ \bI=(I,\delta_0,\delta_1)\]
such that $I$ is aspheric over the final object of \Ahat. Then
$p_1^*(\bI)$ is of course a separating interval in $(A\times
B)\uphat$, which by cor.\ \ref{cor:74.1.1} is aspheric over the final
object. Hence
\begin{propositionnum}\label{prop:74.2}
  If $A$ is a local test category, so is $A\times B$ for any $B$ in
  \Cat. If $A$ is a test category, then so is $A\times B$ for any
  aspheric $B$.
\end{propositionnum}

The second assertion follows, remembering that a test category is just
an \emph{aspheric} local test category. Using cor.\ \ref{cor:74.1.2},
we get:
\begin{corollary}
  If $A$ is a strict test category, and $B$ totally aspheric, then
  $A\times B$ is a strict test category. In particular, if $A$ and $B$
  are strict test categories, so is their product.
\end{corollary}

We need only remember that strict test categories are just test
categories that are totally aspheric.

As an illustration, we get the fact that the categories of
multicomplexes of various kinds, which we can even take mixed
(semisimplicial in some variables, cubical in others, and
hemispherical say in others still) are strict modelizers (as generally
granted), which corresponds to the fact that a finite product of
standard semisimplicial, cubical and hemispherical (strict) test
categories $\Simplex$, $\Square$ and $\Globe$, is again a strict test
category.

I would like now to dwell a little upon the comparison of ``homotopy
models'', using respectively two test categories $A$ and $B$, namely
working in \Ahat{} and \Bhat{} respectively. More specifically, we
have two description of $\HotOf$ (short for $\HotOf(\scrW)$ here),
namely as
\[ \HotOf_A = \scrWA^{-1}\Ahat \quad\text{and}\quad \HotOf_B =
\scrWA^{-1}\Bhat,\]
and we want to describe conveniently the tautological equivalence
between these two categories (this equivalence being defined up to
unique isomorphism). The most ``tautological'' way indeed is to use
the basic modelizer \Cat{} and its localization \Hot{} as the
intermediary, i.e., using the diagram of equivalences of categories
\begin{equation}
  \label{eq:74.1}
  \begin{tabular}{@{}c@{}}
    \begin{tikzcd}[baseline=(O.base),column sep=small]
      \HotOf_A\ar[dr,"\equ","\overline{i_A}"'] & &
      \HotOf_B\ar[dl,"\overline{i_B}","\equ"'] \\
      & |[alias=O]| \Hot &
    \end{tikzcd}.
  \end{tabular}
  \tag{1}
\end{equation}
Remember we have a handy quasi-inverse $\overline{j_B}$\pspage{224} to
$\overline{i_B}$, using the canonical functor\scrcomment{In this display, AG
  originally put the definition of $i_B$ instead $j_B$, cf., e.g.,
  \eqref{eq:65.1} and \eqref{eq:65.2} in \S\ref{sec:65} to recall the definitions.}
\[ j_B = i_B^*: \Cat\to\Bhat, \quad X\mapsto(b\mapsto\Hom(B_{/b},X)).\]
Thus, we get a description of an equivalence
\begin{equation}
  \label{eq:74.1prime}
  \overline{j_B}\, \overline{i_A} : \HotOf_A \toequ \HotOf_B,\tag{1'}
\end{equation}
whose quasi-inverse of course is just the similar $\overline{j_A}
\,\overline{i_B}$. We can vary a little this description, admittedly
cumbersome in practice, by replacing $\overline{j_B}$ by the
isomorphic functor $\overline{i^*}$, where $i:B\to\Cat$ is any test
functor from $B$ to \Cat{} (while we have to keep however
$\overline{i_A}$ as it is, without the possibility of replacing $i_A$
by a more amenable test functor).

Another way for comparing $\HotOf_A$ and $\HotOf_B$ arises, as we saw
yesterday, whenever we have an \emph{aspheric} functor
\begin{equation}
  \label{eq:74.2}
  i: A\to B,\tag{2}
\end{equation}
by just taking
\begin{equation}
  \label{eq:74.2prime}
  \overline{i^*}: \HotOf_B \toequ \HotOf_A.\tag{2'}
\end{equation}
This of course is about the simplest way imaginable, all the more as
the functor $i^*:\Bhat\to\Ahat$ commutes to arbitrary direct and
inverse limits, just perfect for comparing constructions in \Bhat{}
and constructions in \Ahat{} -- whereas the functor
$j_Bi_A:\Ahat\to\Bhat$, giving rise to \eqref{eq:74.1} above commutes
just to sums and fibered products, not to amalgamated sums nor to
products, sadly enough! We could add here that if we got \emph{two}
aspheric functors from $A$ to $B$, namely $i$ plus
\[ i':A\to B,\]
then (as immediately checked) any map between these functors
\[u: i\to i'\]
gives rise to an isomorphism between the corresponding equivalences
\[\begin{tikzcd}[baseline=(O.base),column sep=large]
  |[alias=O]| \HotOf_B
  \ar[r,bend left,"\overline{i^*}",""{name=i,below}]
  \ar[r,bend right,"\overline{(i')^*}"',""{name=iprime,above}] &
  \HotOf_A \arrow[from=iprime,to=i,"\overline u"']
\end{tikzcd},\]
which is nothing but the canonical isomorphism referred to yesterday
(both functors being canonically isomorphic to \eqref{eq:74.1prime}),
and yields the most evident way for ``computing'' the latter. Thus, it
turns out that the isomorphism $\overline u$ does not depend upon the
choice of $u$. If we want to ignore this fact and look at the
situation sternly, in a wholly computational spirit, we could present
things by stating that we get a contravariant functor, from aspheric
maps $i:A\to B$ to equivalences $\HotOf_B\to\HotOf_A$:\pspage{225}
\[ \bAsph(A,B)\op \to \bHom(\HotOf_B,\HotOf_A), \quad i\mapsto
\overline{i^*},\]
where $\bAsph$ denotes the full subcategory of the functor category
$\bHom$ made up with aspheric functors. This functor transforms
arbitrary arrows from the left hand side into isomorphisms on the
right, and therefore, it factors through the fundamental groupoid
(i.e., localization of the left hand category with respect to the set
of \emph{all} its arrows):
\[ \bigl(\Pi_1(\bAsph(A,B))\bigr)\op\to\bHom(\HotOf_B,\HotOf_A).\]
If we now remember that we had assumed $A$ and $B$ to be test
categories (otherwise the functors just written would be still
defined, but their values would not necessarily be \emph{equivalences}
between $\HotOf_B$ and $\HotOf_A$, but merely functors between these),
we may hope that this might imply that the category $\bAsph(A,B)$ of
aspheric functors from $A$ to $B$, whenever non-empty, to be
$1$-connected. Whenever this is so, in any case, we get ``a priori''
(namely without any reference to \Hot{} itself) a transitive system of
isomorphisms between the functors $\overline{i^*}$, for $i$ in
$\bAsph(A,B)$, hence a \emph{canonical functor} $\HotOf_B\to\HotOf_A$,
\emph{defined up to unique isomorphism} (and which, in case $A$ and
$B$ are test categories, or more generally pseudo-test categories, is
an equivalence, and the one precisely stemming from the diagram
\eqref{eq:74.1}).

\begin{remark}
  Here the reflection slipped, almost against will, into a related
  one, about comparison of $\HotOf_A$ and $\HotOf_B$ for arbitrary $A$
  and $B$ (not necessarily test categories nor \emph{even} pseudo test
  categories), using aspheric functors $i:A\to B$ to get
  $\overline{i^*}:\HotOf_B\to\HotOf_A$ (not necessarily an
  equivalence). As seen above, this functor depends rather loosely on
  the choice of $i$, and we could develop comprehensive conditions on
  $A$ and $B$, not at all of a test-condition nature, implying that
  just using the isomorphisms $\overline u$ between these functors, we
  get a canonical \emph{transitive} system of isomorphisms between
  them, hence a \emph{canonical} functor $\HotOf_B\to\HotOf_A$, not
  depending on the particular choice of any aspheric functor $i:A\to
  B$. As we are mainly interested in the modelizing case though, I
  don't think I should dwell on this much longer here. Anyhow, in case
  both $A$ and $B$ are pseudo-test categories, and provided only
  $\bAsph(A,B)$ is $0$-connected (not necessarily $1$-connected), it
  follows from comparison with the diagram \eqref{eq:74.1} above that
  the isomorphisms $\overline u$ do give rise indeed to a
  \emph{transitive} system of isomorphisms between the equivalences $\overline{i^*}$.
\end{remark}

In most cases though, such as $A=\Simplex$ and $B=\Square$ say (the
test\pspage{226} categories of standard simplices and standard cubes
respectively), we do not have any aspheric functor $A\to B$ at hand,
and presumably we may well have that $\bAsph(A,B)$ is empty, poor it!
We now assume again that $A$ and $B$ are test categories, hence
$A\times B$ is a test category, and the natural idea for comparison of
$\HotOf_A$ and $\HotOf_B$ is to use the diagram
\[\begin{tikzcd}[baseline=(O.base),column sep=tiny]
  & A\times B\ar[dl,"p_1"']\ar[dr,"p_2"] & \\
  A & & |[alias=O]| B \end{tikzcd},\]
where now $p_1$ and $p_2$ are aspheric (because $B$ and $A$ are
aspheric). Thus, we get a corresponding diagram on the corresponding
modelizers\scrcomment{in the typescript this equation is tagged (2)}
\begin{equation}
  \label{eq:74.3}
  \begin{tabular}{@{}c@{}}
    \begin{tikzcd}[baseline=(O.base),column sep=tiny]
      & \HotOf_{A\times B} & \\
      \HotOf_A\ar[ur,"\overline{p_1^*}","\equ"'] & & |[alias=O]|
      \HotOf_B\ar[ul,"\overline{p_2^*}"',"\equ"]
    \end{tikzcd}.
  \end{tabular}\tag{3}
\end{equation}
These equivalences are compatible with the canonical equivalences with
$\HotOf$ itself, and hence they give rise to an equivalence $\HotOf_A
\toequ \HotOf_B$ which (up to canonical isomorphism) is indeed the
canonical one. This way for comparing $\HotOf_A$ and $\HotOf_B$ looks
a lot more convenient than the first one, as the functors $p_1^*$ and
$p_2^*$ which serve as intermediaries have all desirable exactness
properties, and their very description is the simplest imaginable.

We are very close here to an Eilenberg-Zilber situation, which will
arise more specifically when $A$ and $B$ are both \emph{strict} test
categories, i.e., in the corresponding modelizers $\Ahat,\Bhat$
products of models do correspond to products of the corresponding
homotopy types. As seen above, this implies the same for $A\times
B$. Thus, if $F$ is in \Ahat, $G$ in \Bhat, the objects $F\boxtimes G$
in $(A\times B)\uphat$ is just a $(A\times B)$-model for the homotopy
type described respectively by $F$ (in terms of $A$) and $G$ (in terms
of $B$). As a matter of fact, the relation
\[(A\times B)_{/{F \boxtimes G}} \simeq A_{/F} \times B_{/G},\]
already noticed before, implies that this interpretation holds,
independently of strictness. In the classical statement of
Eilenberg-Zilber (as I recall it), we got $A=B$ (both categories being
just $\Simplex$);\pspage{227} on the one hand the $A\times B$-model is
used in order to get readily the Künneth type relations for homology
and cohomology, whereas we are interested really in the $A$-model
$F\times G$. Using the diagonal map
\[ \delta: A\to A\times A,\]
we have of course
\[F \times G\simeq \delta^*(F\boxtimes G),\]
more generally it is expected that the functor
\[ \delta^*: (A\times A)\uphat \to \Ahat,\]
is modelizing, i.e., gives a means of passing from $(A\times
A)$-models to $A$-models. This essentially translates into $\delta$
being an aspheric map -- something which will not be true for
arbitrary $A$ in \Cat. We get in this respect:
\begin{propositionnum}\label{prop:74.3}
  Let $A$ be an object in \Cat. Then the diagonal map $\delta:A\to
  A\times A$ is aspheric if{f} $A$ is totally aspheric.
\end{propositionnum}

This is just a tautology -- one among many which tell us that the
conceptual set-up is OK indeed! A related tautology:
\begin{corollary}
  Let $i:A\to B$, $i':A\to B'$ be two aspheric functors with same
  source $A$\kern1pt, then the corresponding functor
  \[(i,i'):A \to B\times B'\]
  is aspheric, provided $A$ is totally aspheric.
\end{corollary}

This can be seen either directly, or as a corollary of the
proposition, by viewing $(i,i')$ as a composition
\[A \xrightarrow \delta A\times A \xrightarrow{i\times i'} B\times
B'.\]

To sum up: when we got a bunch of strict test categories, and a
(possibly empty) bunch of aspheric functors between some of these,
using finite products we get a (substantially larger!) bunch of strict
test categories, giving rise to corresponding strict modelizers for
homotopy types; and using projections, diagonal maps, and the given
aspheric functors, we get an impressive lot of aspheric maps between
all these, namely as many ways to ``commute'' from one type of
``homotopy models'' to others. The simplest example: start with just
one strict test category $A$, taking products
\[ A^I\]
where $I$ is any finite set, and the ``simplicial'' maps between
these, expressing contravariance of $A^I$ with respect to $I$. The
case most commonly used is $A=\Simplex$, giving rise to the formalism
of semisimplicial multicomplexes.

To\pspage{228} finish these generalities on aspheric functors, I would
like still to make some comments on $\bAsph(A,B)$, in case $B$ admits
binary products, with $A$ and $B$ otherwise arbitrary in \Cat. We are
interested, for two functors
\[ i,i' : A \rightrightarrows B,\]
in the functor
\[ i\times i': a\mapsto i(a)\times i'(a) : A \to B,\]
which can be viewed indeed as a product object in $\bHom(A,B)$ of $i$
and $i'$. Our interest here is mainly to give conditions ensuring that
$i\times i'$ is aspheric. It turns out that there is no point to this
end assume \emph{both} $i$ and $i'$ to be aspheric, what counts is
that one, say $i$, should be aspheric. The point is made very clearly
in the following
\begin{propositionnum}\label{prop:74.4}
  Let $i:A\to B$ be an aspheric map in \Cat, we assume for simplicity
  that in $B$ binary products exist.
  \begin{enumerate}[label=\alph*),font=\normalfont]
  \item\label{it:74.4.a}
    Let $b_0\in\Ob B$, $i_{b_0}:A\to B$ be the constant functor with
    value $b_0$, consider the product functor
    \[ i\times i_{b_0} : a\mapsto i(a)\times b_0 : A\to B.\]
    In order for this functor to be aspheric, it is necessary and
    sufficient that for any object $y$ in $B$, the object
    $\bHom(b_0,y)$ in \Bhat{} be aspheric \textup(a condition which
    depends only on $b_0$, not upon $i$ nor even upon $A$\textup).
  \item\label{it:74.4.b}
    Assume this condition satisfied for any $b_0$ in $B$, assume
    moreover $A$ totally aspheric. Then for \emph{any} functor
    $i':A\to B$, $i\times i'$ is aspheric.
  \end{enumerate}
\end{propositionnum}
\begin{proof}
  \ref{it:74.4.a}\enspace Asphericity of $i\times i_{b_0}$ means that
  for any object $y$ in $B$, the corresponding presheaf on $A$
  \[ (i\times i_{b_0})^*(y) = (a \mapsto \Hom(i(a)\times b_0, y))\]
  is aspheric. Now the formula defining the presheaf $\bHom(b_0,y)$ on
  $B$ yields
  \[ \Hom(i(a)\times b_0,y) \simeq \Hom(i(a),\bHom(b_0,y)),\]
  hence the presheaf we get on $A$ is nothing but
  \[i^*(\bHom(b_0,y)).\]
  As $i$ is aspheric, the criterion \ref{it:72.1.iiprime} of prop.\
  \ref{prop:72.1} (page \ref{p:214}) implies that this presheaf is
  aspheric if{f} $\bHom(b_0,y)$ is, qed.

  \ref{it:74.4.b}\enspace We may view $i\times i'$ as a composition
  \[ A \xrightarrow\delta A\times A \xrightarrow{i\boxtimes i'} B,\]
  where the ``external product'' $i\boxtimes i': A\times A\to B$ is
  defined by\pspage{229}
  \[{i\boxtimes i'}(a,a') \eqdef i(a)\times i'(a').\]
  By prop.\ \ref{prop:74.3} we know that the diagonal map for $A$ is
  aspheric, thus we are left with proving that $i\boxtimes i'$ is
  aspheric. More generally, we get:
\end{proof}
\setcounter{corollarynum}{0}
\begin{corollarynum}\label{cor:74.4.1}
  Let $B$ in \Cat{} satisfy the conditions of \textup{\ref{it:74.4.a}}
  and \textup{\ref{it:74.4.b}} above, and let $A,A'$ be two objects in
  \Cat, and
  \[i : A\to B, \quad i':A' \to B\]
  two maps with target $B$, hence a map
  \[ i\boxtimes i' : (a,a')\mapsto i(a)\times i'(a') : A\times A'\to
  B.\]
  If $i$ is aspheric, then $i\boxtimes i'$ is aspheric if{f} $A'$ is aspheric.
\end{corollarynum}

We have to express that for any object $b$ in $B$, the presheaf
\[ (a,a') \mapsto \Hom(i(a)\times i'(a'),b) \simeq \Hom(i(a),
\bHom(i'(a'),b))\]
on $A\times A'$ is aspheric. For $a'$ fixed in $A'$, the corresponding
presheaf on $A$ is aspheric, as we saw in \ref{it:74.4.a}. The
conclusion now follows from the useful
\begin{lemma}
  Let $F$ be a presheaf on a product category $A\times A'$, with
  $A,A'$ in \Cat. Assume that for any $a'$ in $A'$, the corresponding
  presheaf $a\mapsto F(a,a')$ on $A$ be aspheric. Then the composition
  \[(A\times A')_{/F} \to A\times A' \xrightarrow{\mathrm{pr}_2} A'\]
  is in \scrW, and hence $F$ is aspheric if{f} $A'$ is aspheric.
\end{lemma}
\begin{proof}[Proof of lemma]
  The composition is fibering, as both factors are. This reminds us of
  the ``fibration condition'' \hyperref[it:64.L5]{L~5} on \scrW{}
  (page \ref{p:164}), as yesterday (page \ref{p:219}), where we stated
  that this condition follows from the conditions \ref{loc:1} to
  \ref{loc:3} reviewed yesterday. This condition asserts that a
  fibration with aspheric fibers is in \scrW{} -- hence the lemma.
\end{proof}
\begin{remarks}
  The results stated in prop.\ \ref{prop:74.4} and its corollary give
  a lot of elbow freedom for getting new aspheric functors in terms of
  old ones, with target category $B$ -- provided $B$ satisfies the two
  assumptions: stability under binary products (a property frequently
  met with, although the standard test categories $\Simplex$,
  $\Square$ and $\Globe$ lack it\ldots), and the asphericity of the
  presheaves $\bHom(b,y)$, for any two objects $b,y$ in $B$. A
  slightly stronger condition (indeed an equivalent one, for fixed
  $b$, when $b$ admits already a section over $e_\Bhat$ and if \Bhat{}
  is totally aspheric, hence an aspheric object is even aspheric over
  $e_\Bhat$\ldots) is \emph{contractibility} of the objects $b$ of
  $B$, for the homotopy interval structure on \Bhat{} defined by
  \scrWB{} (i.e., in terms of homotopy intervals aspheric over
  $e_\Bhat$).\pspage{230} (Compare with propositions on pages
  \ref{p:121} and \ref{p:143}.) Whereas this latter assumption
  admittedly is quite a stringent one, it is however of a type which
  has become familiar to us in relation with test categories, where it
  seems a rather common lot.

  If $B$ satisfies these conditions (as in prop.\ \ref{prop:74.4} and
  if $A$ is totally aspheric, we see from prop.\ \ref{prop:74.4} that
  the category $\bAsph(A,B)$ of aspheric functors from $A$ to $B$ is
  stable under binary products. Now, a non-empty object $C$ of \Cat{}
  stable under binary products gives rise to a category $\Chat$
  which is clearly totally aspheric for any basic localizer \scrW, and
  in particular for the usual one \scrWz{}\scrcomment{later, we'll
    write \scrWoo{} for this instead\ldots} corresponding to usual
  weak equivalence. A fortiori, such a category $C$ is $1$-connected
  (which is easily checked too by down-to-earth direct
  arguments). Thus, the reflections of page \ref{p:225} apply, and
  imply that if $\bAsph(A,B)$ is non-empty, i.e., if there is at least
  \emph{one} aspheric functor $i:A\to B$, then there is a canonical
  transitive system of isomorphisms between all functors
  \begin{equation}
    \label{eq:74.star}
    \HotOf_B \to \HotOf_A\tag{*}
  \end{equation}
  of the type $\overline{i^*}$, and hence there is a canonical functor
  \eqref{eq:74.star}, defined up to unique isomorphism, as announced
  in the remark on page \ref{p:225}.
\end{remarks}

\bigbreak
\presectionfill\ondate{13.6.}\par

\hangsection{Asphericity structures: a bunch of useful tautologies.}%
\label{sec:75}%
The generalities on aspheric maps of the last three sections should be
more than what is needed to develop now the notion of
\emph{asphericity structure} -- which, together with the closely
related notion of \emph{contractibility structure}, tentatively dealt
with before, and the various ``test-notions'' (e.g., \emph{test
  categories} and \emph{test functors}) seems to me the main pay-off
so far of our effort to come to a grasp of a general formalism of
``homotopy models''.

In the case of asphericity structures, just as for the kindred notion
of a contractibility structure, in all instances I could think of a
present, the underlying category $M$ of an asphericity structure is
\emph{not} an object in \Cat{} nor even a ``small category'', but is
``large'' -- namely the cardinality of $\Ob M$ and $\Fl(M)$ are not in
the ``universe'' we are working in, still less is $M$ an object of
\scrU{} -- all we need instead, as usual, is that $M$ be a
\scrU-category, namely that for any two objects of $M$, $\Hom(x,y)$ be
an element of \scrU. Till now, the universe \scrU{} has been present
in our reflections in a very much implicit way, in keeping with the
informal nature of the reflections, which however by and by have
become more formal (as I finally let myself become involved
in\pspage{231} a minimum of technical work, needed for keeping out of
the uneasiness of ``thin air conjecturing''). An attentive reader will
have felt occasionally this implicit presence of \scrU, for instance
in the definition of the basic modelizer \Cat{} (which is, as all
modelizers, a ``large'' category), and in our occasional reference to
various categories as being ``small'' or ``large''. He\scrcomment{or
  she} will have noticed that whereas modelizers are by necessity
large categories (just as \Hot{} itself, whose set of isomorphism
classes of objects is large), test categories are supposed to be small
(and often even to be in \Cat) -- with the effect that \Ahat, the
category of presheaves on $A$, is automatically a \scrU-category
(which would not be the case if $A$ was merely assumed to be a
\scrU-category).

An \emph{asphericity structure} (with respect to the basic localizer
\scrW) on the \scrU-category $M$ consists of a subset
\begin{equation}
  \label{eq:75.1}
  M\subas\subset\Ob M,\tag{1}
\end{equation}
whose elements will be called the \emph{aspheric objects} of $M$ (more
specifically, the \scrW-aspheric objects, if confusion may arise),
this subset being submitted to the following condition:
\begin{description}
\item[\namedlabel{cond:Asstr}{(Asstr)}]
  There exists an object $A$ in \Cat, and a functor $i:A\to M$,
  satisfying the following two conditions:
  \begin{enumerate}[label=(\roman*)]
  \item\label{cond:Asstr.i}
    For any $a$ in $A$, $i(a)\in M\subas$, i.e., $i$ factors through
    the full subcategory (also denoted by $M\subas$) of $M$ defined by
    $M\subas$.
  \item\label{cond:Asstr.ii}
    Let
    \[i^*:M\to\Ahat, \quad x\mapsto i^*(x)=(a\mapsto\Hom(i(a),x)),\]
    then we have
    \begin{equation}
      \label{eq:75.2}
      M\subas = (i^*)^{-1}(\Ahatas),\tag{2}
    \end{equation}
    where $\Ahatas$ is the subset of $\Ob\Ahat$ of all
    (\scrW-)aspheric objects of \Ahat, i.e., the presheaves $F$ on $A$
    such that the object $A_{/F}$ of \Cat{} is \scrW-aspheric.
  \end{enumerate}
\end{description}

In other words, for $x$ in $M$, we have the equivalence
\begin{equation}
  \label{eq:75.2bis}
  x\in M\subas \Leftrightarrow i^*(x)\in\Ahatas,
  \quad\text{i.e.,}\quad
  A_{/x} (\eqdef A_{/i^*(x)}) \in \Cat\subas,\tag{2 bis}
\end{equation}
where $\Cat\subas$ is the subset of $\Ob\Cat$ made up with all
\scrW-aspheric objects of \Cat{} -- i.e., objects $X$ such that
$X\to\Simplex_0$ ($=$ final object) be in \scrW.

Thus, an asphericity structure on $M$ can always be defined by a
functor
\[i:A\to M,\]
with $A$ in \Cat, by the formula \eqref{eq:75.2}; and conversely, any
such functor\pspage{232} defined an asphericity structure $M\subas$ on
$M$, admitting $i$ as a ``\emph{testing functor}'' (i.e., a functor
satisfying \ref{cond:Asstr.i} and \ref{cond:Asstr.ii} above), provided
only we assume that
\begin{equation}
  \label{eq:75.3}
  \text{For any $a$ in $A$, $i^*(i(a))$ is an aspheric object of
    \Ahat.}\tag{3}
\end{equation}
This latter condition is automatically satisfied if $i$ is fully
faithful, for instance, if it is the inclusion functor of a full
subcategory $A$ of $M$. Thus, \emph{any small full subcategory $A$ of
  $M$ defines an asphericity structure on $M$}, and we'll see in a
minute that any asphericity structure on $M$ can be defined this way.

\begin{propositionnum}\label{prop:75.1}
  Let $M\subas$ be an asphericity structure on $M$, and $A,B$ objects
  in \Cat, and
  \[ A \xrightarrow f B \xrightarrow j M\]
  be functors, with $j$ factoring through $M\subas$.
  \begin{enumerate}[label=\alph*),font=\normalfont]
  \item\label{it:75.1.a}
    If $f$ is aspheric, then $j$ is a testing functor if{f} $i=jf$ is.
  \item\label{it:75.1.b}
    If $j$ is fully faithful, and if $i=jf$ is a testing functor, then
    $f$ is aspheric, and $j$ is a testing functor too.
  \end{enumerate}
\end{propositionnum}
\begin{proof}
  \ref{it:75.1.a} follows trivially from
  \[\Bhatas = (f^*)^{-1}(\Ahatas),\]
  which is one of the ways of expressing that $f$ is aspheric
  (criterion \ref{it:72.1.iiprime} of prop.\ \ref{prop:72.1},
  \ref{p:214}). And the first assertion in \ref{it:75.1.b} is a
  trivial consequence of the definition of testing functors, and of
  asphericity of $f$ (by criterion \ref{it:72.1.iii} on p.\
  \ref{p:214}) -- and the second assertion of \ref{it:75.1.b} now
  follows from \ref{it:75.1.a}.
\end{proof}
\begin{corollary}
  Let $i:A\to M$ be a testing functor for $(M,M\subas)$, and let $B$
  any \emph{small} full subcategory of $M\subas$ containing
  $i(A)$. Then the induced functor $A\to B$ is aspheric, and the
  inclusion functor $B\hookrightarrow M$ is a testing functor.
\end{corollary}
\begin{remarks}
  This shows, as announced above, that any asphericity structure on
  $M$ can be defined by a small full subcategory of $M$ (for instance,
  the smallest full subcategory of $M$ containing $i(A)$). We have
  been slightly floppy though, while we defined testing functors by
  insisting that the source should be in Cat{} (which at times will be
  convenient), whereas the $B$ we got here is merely small, namely
  \emph{isomorphic} to an object of \Cat, but not necessarily in
  \Cat{} itself. This visibly is an ``inessential floppiness'', which
  could be straightened out trivially, either by enlarging accordingly
  the notion of testing functor, or by submitted $M$ to the (somewhat
  artificial, admittedly) restriction that all\pspage{233} its small
  subcategories should be in \Cat. This presumably will be satisfied
  by most large categories we are going to consider, and it shouldn't
  be hard moreover to show that any \scrU-category is isomorphic to a
  category $M'$ satisfying the above extra condition.
\end{remarks}

A little more serious maybe is the use we are making here of the name
of a ``testing functor'', which seems to be conflicting with an
earlier use (def.\ \ref{def:65.5}, p.\ \ref{p:175} and def.\
\ref{def:65.6}, p.\ \ref{p:176}), where we insisted for instance that
the source $A$ should be a test category. That's why I have been using
here the name ``testing functor'' rather than ``test functor'', to be
on the safe side formally speaking -- but this is still playing on
words, namely cheating a little. Maybe the name ``\emph{aspherical
  functor}'' instead of ``testing functor'' would be less misleading,
thinking of the case when $M$ is small itself, and endowing $M$ with
the canonical asphericity structure, for which
\[M\subas = M,\]
(admitting the identity functor as a testing functor) -- in which case
the ``testing functors'' $A\to M$ are indeed just the aspherical
functors. The drawback is that when $M$ is small, it may well be
endowed with an asphericity structure $M\subas$ different from the
previous one, in which case the proposed extension of the name
``aspheric functor'' again leads to an ambiguity, unless specified by
``\emph{$M\subas$-aspheric}'' (where, after all, $M\subas$ could be
\emph{any} full subcategory of $M$) -- but then the notion reduces to
the one of an aspheric functor $A\to M\subas$. But the same after all
holds even for large $M$ -- the notion of a ``testing functor'' $A\to
M$ (with respect to an asphericity structure $M\subas$ on $M$) does
not really depend on the \emph{pair} $(M,M\subas)$, but rather on the
(possibly large) category $M\subas$ itself -- namely it is no more no
less than a functor
\[A\to M\subas\]
which satisfies the usual asphericity condition \ref{it:72.1.iii} (of
prop.\ \ref{prop:72.1}, p.\ \ref{p:214}), with the only difference
that $M\subas$ may not be small, and therefore $M\subas\uphat$ may not
be a \scrU-category (and we are therefore reluctant to work with this
latter category at all, unless we first pass to the next larger
universe $\scrU'$\ldots).

This short reflection rather convinces us that the designation of
``testing functors'' as introduced on the page before, by the
alternative name of \emph{$M\subas$-aspheric functors}, or just
\emph{aspheric functors} when no confusion is likely to arise as to
existence and choice of $M\subas$, is satisfactory indeed. I'll use it
tentatively, as a synonym to ``testing functor'', and it will appear
soon enough if this terminology is a good one\pspage{234} or not --
namely if it is suggestive, and not too conducive to
misunderstandings.
\begin{propositionnum}\label{prop:75.2}
  Let $(M,M\subas)$ be an asphericity structure, and let
  \[i:A\to M\]
  be an $M\subas$-aspheric functor, hence $i^*:M\to\Ahat$. Let
  \begin{equation}
    \label{eq:75.4}
    W = (i^*)^{-1}(\scrWA) \subset \Fl(M).\tag{4}
  \end{equation}
  Then $W$ is a mildly saturated subset of $\Fl(M)$, independent of
  the choice of $(A,i)$.
\end{propositionnum}

Mild saturation of $W$ follows from mild saturation of \scrWA, the set
of \scrW-weak equivalences in \Ahat{} (which follows from mild
saturation of \scrW{} and the definition
\[ \text{$\scrWA=i_A^{-1}(\scrW)$ ).}\]
To prove that $W_i$ for $(A,i)$ is the same as $W_{i'}$ for $(A',i')$,
choose a small full subcategory $B$ of $M\subas$ such that $i$ and
$i'$ factor through $B$, let $j:B\to M$ be the inclusion, we only have
to check $W_i=W_j$ (and similarly for $W_{i'}$), which follows from
$i^*=f^*j^*$ (where $f:A\to B$ is the induced functor) and the
relation
\[\scrWA = (f^*)^{-1}(\scrWA),\]
which follows from $f$ being aspheric by prop.\ \ref{prop:75.1}
\ref{it:75.1.b}, and the known property \ref{it:72.1.iv} (prop.\
\ref{prop:72.1}, p. \ref{p:214}) of asphericity, qed.
\begin{remark}
  Once we prove that \scrW{} is even strongly saturated, it will
  follow of course that the sets \scrWA, and hence $W$ above, are
  strongly saturated too.
\end{remark}

We'll call $W$ the set of \emph{weak equivalences} in $M$, for the
given asphericity structure $M\subas$ in $M$. We are interested now in
giving a condition on $M\subas$ ensuring that conversely, $M\subas$ is
known when the corresponding set $W$ of weak equivalences is. We'll
assume for this that $M$ has a final object $e_M$, and we'll call the
asphericity structure $(M,M\subas)$ ``\emph{aspheric}'' if $e_M$ is
aspheric, i.e., $e_M\in M\subas$ (which does not depend of course on
the choice of $e_M$). Now we get the tautology:
\begin{propositionnum}\label{prop:75.3}
  Let $(M,M\subas)$ be an asphericity structure, and $u:x\to y$ a map
  in $M$. If $y$ is aspheric, then $x$ is aspheric if{f} $u$ is
  aspheric.
\end{propositionnum}
\begin{corollary}
  Assume $M$ admits a final object $e_M$, and that $e_M$ is aspheric
  \textup(i.e., $(M,M\subas)$ is aspheric\textup). Then an object $x$
  in $M$ is aspheric if{f} the map $x\to e_M$ is aspheric.
\end{corollary}

Next question then is to state the conditions on a pair $(M,W)$, with
$W\subset\Fl(M)$, for the existence of an aspheric asphericity
structure on $M$,\pspage{235} admitting $W$ as its set of weak
equivalences. We assume beforehand $M$ admits a final object $e_M$. A
n.s.\ condition is the existence of a small category $A$ and a functor
\[i:A\to M\]
satisfying the following two conditions:
\begin{enumerate}[label=(\roman*)]
\item\label{it:75.i}
  for $x$ in $M$ of the form $i(a)$ ($a\in\Ob A$) \emph{or} $e_M$,
  $i^*(x)\in\Ahatas$,
\item\label{it:75.ii}
  $W=(i^*)^{-1}(\scrWA)$.
\end{enumerate}
Another n.s.\ condition, which looks more pleasant I guess, is that
there exist a small full subcategory $B$ of $M$, containing $e_M$,
with inclusion functor $j:B\to M$, such that
\begin{equation}
  \label{eq:75.5}
  W=(j^*)^{-1}(\scrWB).\tag{5}
\end{equation}
In any case, if we want a n.s.\ condition on $W$ for it to be the set
of weak equivalences for some asphericity structure on $M$ (not
necessarily an aspheric one, and therefore maybe not unique), we get:
there should exist a small full subcategory $B$ of $M$, with inclusion
functor $j$, such that \eqref{eq:75.5} above holds. (Here we do not
assume $e_M$ to exist, still less $B$ to contain it.) Similarly, the
corollary of prop.\ \ref{prop:75.1} implies that a subset $M\subas$ of
$\Ob M$ is an asphericity structure on $M$, if{f} there exists a small
full subcategory $B$ in $M$, such that
\begin{equation}
  \label{eq:75.6}
  M\subas = (j^*)^{-1}(\Bhatas);\tag{6}
\end{equation}
moreover, if $M$ admits a final object $e_M$, the asphericity
structure $M\subas$ is aspheric if{f} $B$ can be chosen to contain
$e_M$.
\addtocounter{remarksnum}{1}
\begin{remarksnum}
  The relationship between aspheric asphericity structures $M\subas$
  on $M$, and sets $W \subset \Fl(M)$ of ``weak equivalences'' in $M$
  satisfying the condition above (and which we may call ``weak
  equivalence structures'' on $M$), in case $M$ admits a final object
  $e_M$, is reminiscent of the relationship between ``contractibility
  structures'' $M\subc\subset\Ob M$ on $M$, and those ``homotopism
  structures'' $h_M\subset\Fl(M)$ on $M$ which can be described in
  terms of such a contractibility structure (cf.\ sections
  \ref{sec:51} and \ref{sec:52}). Both pairs can be viewed as giving
  two equivalent ways of expressing one and the same kind of structure
  -- the structure concerned by the first pair being centered on
  asphericity notions, whereas the second is concerned with typical
  homotopy notions rather. We'll see later that any contractibility
  structure defines in an evident way an asphericity structure, and in
  the most interesting cases (e.g., canonical modelizers), it is
  uniquely determined by the latter.
\end{remarksnum}

To\pspage{236} sum up some of the main relationships between the three
asphericity notions just introduced ($M\subas$, $W$,
$M\subas$-aspheric functors), let's state one more tautological
proposition, which is very much a paraphrase of the display given
earlier (prop.\ on p.\ \ref{p:214}) of the manifold aspects of the
notion of an aspheric map between small categories:
\begin{propositionnum}\label{prop:75.4}
  Let $(M,M\subas)$ be an asphericity structure, $A$ a small category,
  $i:A\to M$ a functor, factoring through $M\subas$. Consider the
  following conditions on $i$:
  \begin{description}
  \item[\namedlabel{it:75.4.i}{(i)}]
    $i$ is $M\subas$-aspheric \textup(or ``testing functor'' for
    $M\subas$\textup), i.e.,
    \[ M\subas = (i^*)^{-1}(\Ahatas).\]
  \item[\namedlabel{it:75.4.iprime}{(i')}]
    For $x$ in $M\subas$, $i^*(x)$ is aspheric, i.e.,
    \[ M\subas \subset (i^*)^{-1}(\Ahatas).\]
  \item[\namedlabel{it:75.4.idblprime}{(i'')}]
    \textup(Here, $B$ is a given small full subcategory of $M\subas$
    containing $i(A)$ and which ``\emph{generates}'' the asphericity
    structure $M\subas$, namely such that the inclusion functor
    $j:B\hookrightarrow M$ is $M\subas$-aspheric, i.e.,
    $M\subas=(j^*)^{-1}(\Bhatas)$.\textup)
    For any $x$ in $B$,
    $i^*(x)$ is aspheric, i.e.,
    \[ B \subset (i^*)^{-1}(\Ahatas).\]
  \item[\namedlabel{it:75.4.ii}{(ii)}]
    For any map $u$ in $M$, $u$ is a weak equivalence if{f} $i^*(u)$
    is, i.e.,
    \[ W = (i^*)^{-1}(\scrWA).\]
  \item[\namedlabel{it:75.4.iiprime}{(ii')}]
    If the map $u$ in $M$ is a weak equivalence, so is $i^*(u)$, i.e.,
    \[ W \subset (i^*)^{-1}(\scrWA).\]
  \item[\namedlabel{it:75.4.iidblprime}{(ii'')}]
    \textup(Here, $B$ is given as in \textup{\ref{it:75.4.idblprime}}
    above\textup)
    For any map $u$ in $M$, $i^*(u)$ is a weak equivalence, i.e.,
    \[ \Fl(B) \subset (i^*)^{-1}(\scrWA).\]
  \end{description}
  The conditions
  \textup{\ref{it:75.4.i}\ref{it:75.4.iprime}\ref{it:75.4.idblprime}}
  are equivalent and imply all others, and we have the tautological
  implications \textup{\ref{it:75.4.ii}} $\Rightarrow$
  \textup{\ref{it:75.4.iiprime}} $\Rightarrow$
  \textup{\ref{it:75.4.iidblprime}}. If $M$ admits a final object
  $e_M$ and if $A$ and the asphericity structure $M\subas$ are
  aspheric, then all six conditions \textup(except the last\textup)
  are equivalent; and all six are equivalent if moreover $e_M\in\Ob B$.
\end{propositionnum}
\begin{proof}
  The implications \ref{it:75.4.i} $\Rightarrow$ \ref{it:75.4.iprime}
  $\Rightarrow$ \ref{it:75.4.idblprime} are tautological, on the other
  hand \ref{it:75.4.idblprime} just means that the induced functor $f:
  A\to B$ is aspheric, which by prop.\ \ref{prop:75.1} \ref{it:75.1.a}
  implies that $i$ is aspheric. On the other hand \ref{it:75.4.i}
  $\Rightarrow$ \ref{it:75.4.ii} by the definition of $W$ (cf.\ prop.\
  \ref{prop:75.2}). If $e_M$ exists and $A$ and the asphericity
  structure $M\subas$ are aspheric, and if $B$ contains $e_M$, then
  \ref{it:75.4.iidblprime} implies that the maps $x\to e_M$ for $x$ in
  $B$ are transformed by $i^*$ into a weak equivalence, and as
  $e_\Ahat$ is aspheric, this implies $i^*(x)$ is aspheric,
  i.e., \ref{it:75.4.idblprime}, which proves the last
  statement\pspage{237} of the proposition -- all six conditions are
  equivalent in this case. If no $B$ is given, but still assuming $A$
  and $(M,M\subas)$ aspheric, we can choose a generating subcategory
  for $(M,M\subas)$ large enough in order to contain $i(A)$ and $e_M$,
  and we get that conditions \ref{it:75.4.i} to \ref{it:75.4.iiprime}
  are equivalent, qed.
\end{proof}

\hangsection{Examples. Totally aspheric asphericity structures.}%
\label{sec:76}%
It's time to give some examples.

\namedlabel{ex:76.1}{1)}\enspace
Take $M=\Cat$, $M\subas=$ set of \scrW-aspheric objects in \Cat. We
then got an asphericity structure, as we see by taking any weak test
category $A$ in \Cat{} (def.\ \ref{def:65.2} on page \ref{p:172}), and
the functor
\[i_A: a\mapsto A_{/a} : A \to\Cat,\]
which satisfies indeed \ref{cond:Asstr.i} and \ref{cond:Asstr.ii} of
p.\ \ref{p:231} above.

We may call
\begin{equation}
  \label{eq:76.7}
  (\Cat, \Cat_{\scrW\mathrm{-as}} \eqdef \Cat\subas)\tag{7}
\end{equation}
the ``basic asphericity structure'', giving rise (by taking the
corresponding set $W$ of ``weak equivalences'') to the ``basic
modelizer'' $(\Cat,\scrW)$. Turning attention towards the former
corresponds to a shift in emphasis; whereas previously, our main
emphasis has been dwelling consistently with the notion of ``weak
equivalence'', namely with giving on a category $M$ a bunch of
\emph{arrows} $W$, here we are working rather with notion of
``aspheric \emph{objects}'' as the basic notion. One hint in this
direction comes from prop.\ \ref{prop:72.1} on p.\ \ref{p:214}, when
we saw that for a functor $f:A\to B$ between small categories, giving
rise to $f^*:\Bhat\to\Ahat$ (a map between asphericity structures, as
a matter of fact), asphericity of $f$ can always be expressed in
manifold ways as a property of $f$ relative to the notion of aspheric
\emph{objects} in \Ahat{} and \Bhat, but not as a property relative to
the notion of weak equivalences, unless $A$ and $B$ are assumed to be
aspheric.

Coming back to the case of the ``basic asphericity structures''
\eqref{eq:76.7}, we get more general types of ``aspheric'' functors
$A\to\Cat$ than functors $i_A$, by taking any weak test category $A$
and any weak test functor $i:A\to\Cat$ (cf.\ def.\ \ref{def:65.5}, p.\
\ref{p:175}), provided however the objects $i(a)$ are aspheric. It
would seem though that, for a given $A$, even assuming $A$ to be a
\emph{strict} test category say (and even a ``contractor'' moreover),
and restricting to functors $i:A\to\Cat$ which factor through
$\Cat\subas$ (to make them eligible for being ``aspheric functors''),
the condition for $i$ to be a weak test functor, namely for $i^*$ to
be ``model preserving'', is substantially stronger than mere
``asphericity'': indeed, the latter just means that $i^*$ transforms
weak\pspage{238} equivalences into weak equivalences, i.e., gives rise
to a functor
\[\HotOf{} = \scrW^{-1}\Cat \to \HotOf_A=\scrWA^{-1}\Ahat,\]
whereas the latter insists that this functor moreover should be an
equivalence of categories. Theorem~\ref{thm:65.1} (p.\ \ref{p:176})
gives a hint though that the two conditions may well be equivalent --
this being so at any rate provided the objects $i(a)$ in \Cat{} are,
not only aspheric, but even \emph{contractible}. This reminds us at
once of the ``silly question'' of section \ref{sec:46} (p.\
\ref{p:95}), which was the starting point for the subsequent
reflections leading up to the theorem~\ref{thm:65.1} recalled above;
and, beyond this still somewhat technical result, the ultimate
motivation for the present reflections on asphericity structures. The
main purpose for these, I feel, is to lead up to a comprehensive
answer to the ``silly question''. We'll have to come back to this very
soon!

\namedlabel{ex:76.2}{2)}\enspace
Take $M=\Spaces$, the category of topological spaces, and $M\subas$
the spaces which are weakly equivalent to a point. We get an
asphericity structure, indeed an \emph{aspheric} asphericity structure
(I forgot to make this evident specification in the case \ref{ex:76.1}
above), as we see by taking $A=\Simplex$ for instance, and
\[ i:\Simplex\to\Spaces\]
the ``geometric realization functor'' for simplices, which satisfies
conditions \ref{cond:Asstr.i} and \ref{cond:Asstr.ii} of page
\ref{p:231}, as is well known (cf.\ the book of
Gabriel-Zisman);\scrcomment{\cite{GabrielZisman1967}}
as a matter of fact, $i$ is even a test functor for the modelizer
$(M,W)$ (which is one way of stating the main content of GZ's book).

\namedlabel{ex:76.3}{3)}\enspace
The two examples above are \emph{aspheric} asphericity structures, and
such moreover that $(M,W)$ is (\scrW-)modelizing. These extra features
however are not always present in the next example
\begin{equation}
  \label{eq:76.8}
  M=\Ahat, \quad M\subas=\Ahatas,\quad
  \text{$A$ a small category,}
  \tag{8}
\end{equation}
which is indeed tautologically an asphericity structure, by taking the
canonical inclusion functor (which is fully faithful)
\[ A \to \Ahat,\]
satisfying the conditions \ref{cond:Asstr.i}\ref{cond:Asstr.ii} above
(p.\ \ref{p:231}). This shows moreover that the corresponding notion
of weak equivalence is the usual one, (I forgot to state the similar
fact in example \ref{ex:76.2}, sorry):
\[ W = \scrWA.\]
The asphericity structure is aspheric if{f} $A$ is aspheric. A functor
\[ A'\to A\]
where\pspage{239} $A'$ is another small category, is
$M\subas$-aspheric as a functor from $A'$ to \Ahat, if{f} it is
aspheric.

This last example suggests to call an asphericity structure
$(M,M\subas)$ ``\emph{totally aspheric}'' if $M$ is stable under
finite products, and if the final object of $M$, as well as the
product of any two aspheric objects of $M$, is again aspheric; in
other words, if any finite product in $M$ whose factors are aspheric
is aspheric. We have, in this respect:
\addtocounter{propositionnum}{4}
\begin{propositionnum}\label{prop:76.5}
  Let $(M,M\subas)$ be an asphericity structure, where $M$ is stable
  under finite products. The following conditions are equivalent:
  \begin{enumerate}[label=(\roman*),font=\normalfont]
  \item\label{it:76.5.i}
    $M$ is totally aspheric, i.e., $e_M$ and the product of any two
    aspheric objects of $M$ are aspheric.
  \item\label{it:76.5.ii}
    There exists a small subcategory $B$ of $M$, stable under finite
    products \textup(i.e., containing a final object $e_M$ of $M$ and,
    with any two objects $x$ and $y$, a product $x\times y$ in
    $M$\textup), and which generates the asphericity structure
    \textup(i.e., $j:B\hookrightarrow M$ is
    $M\subas$-aspheric\textup).
  \item\label{it:76.5.iii}
    There exists a small category $A$ such that \Ahat{} is totally
    aspheric \textup(def.\ \ref{def:65.1}, p.\ \ref{p:170}\textup),
    and a $M\subas$-aspheric functor $A\to M$.
  \end{enumerate}
\end{propositionnum}

The proof is immediate. Of course, the asphericity structure in
example \ref{ex:76.3} is totally aspheric if{f} \Ahat{} is totally
aspheric in the usual sense referred to above.

\hangsection{The canonical functor
  \texorpdfstring{$\HotM\to\Hot$}{Hot-M -> (Hot)}.}\label{sec:77}%
Let
\begin{equation}
  \label{eq:77.9}
  \scrM = (M,M\subas)
  \tag{9}
\end{equation}
be an asphericity structure, which will be referred to also as merely
$M$, when no ambiguity concerning $M\subas$ is feared. We'll write
\begin{equation}
  \label{eq:77.10}
  \HotOf_\scrM \text{ (or simply $\HotOf_M$) } = W^{-1}M,
  \tag{10}
\end{equation}
where $W$ of course is the set of weak equivalences in $M$. We are now
going to define a canonical functor
\begin{equation}
  \label{eq:77.11}
  \HotOf_\scrM \to \Hot_\scrW \quad
  (\text{or simply $\Hot\eqdef\scrW^{-1}\Cat$}),
  \tag{11}
\end{equation}
defined at any rate up to canonical isomorphism. For this, take any
aspheric (namely, $M\subas$-aspheric) functor
\[ i : A\to M,\]
and consider the composition
\begin{equation}
  \label{eq:77.12}
  M \xrightarrow{i^*} \Ahat \xrightarrow{i_A} \Cat.
  \tag{12}
\end{equation}
The\pspage{240} three categories in \eqref{eq:77.12} are endowed
respectively with sets of arrows $W$, \scrWA, \scrW, and the two
functors satisfy the conditions
\begin{equation}
  \label{eq:77.13}
  W = (i^*)^{-1}(\scrWA)\quad
  \text{and}\quad
  \scrWA=(i_A)^{-1}(\scrW),\quad
  \tag{13}
\end{equation}
hence $W=(i_{M,i})^{-1}(\scrW)$, where
\[i_{M,i} : M\to \Cat\]
is the composition $i_Ai^*$. Therefore, we get a functor
\begin{equation}
  \label{eq:77.14}
  \overline{i_{M,i}} : \HotOf_M \to \Hot,\tag{14}
\end{equation}
which a priori depends upon the choice of $(A,i)$. If we admit strong
saturation of \scrW, it follows that this functor is ``conservative'',
namely an arrow in $\HotOf_M$ is an isomorphism, provided its image in
\Hot{} is.

To define \eqref{eq:77.11} in terms of \eqref{eq:77.14}, we have to
describe merely a transitive system of isomorphisms between the
functors \eqref{eq:77.14}, for varying pair $(A,i)$. Therefore,
consider two such $(A,i)$ and $(A',i')$, choose a small full
subcategory $B$ of $M\subas$ containing both $i(A)$ and $i'(A')$
(therefore, $B$ is a generating subcategory for the asphericity
structure $M\subas$), and consider the inclusion functor
$j:B\hookrightarrow M$. From prop.\ \ref{prop:75.1} \ref{it:75.1.b} it
follows that the functors
\[f: A\to B, \quad f':A'\to B\]
induced by $i,i'$ are aspheric. Using this, and the criterion
\ref{it:72.1.i} of asphericity (prop.\ \ref{prop:72.1} on p.\
\ref{p:214}) we immediately get two isomorphisms
\[\begin{tikzcd}[baseline=(O.base),column sep=tiny,row sep=small,arrows=phantom]
  & \overline{i_{M,j}} \ar[dl,swap,sloped,"\simeq"{description}]
    \ar[dr,sloped,"\simeq"{description}] & \\
  \overline{i_{M,i}} & & |[alias=O]| \overline{i_{M,i'}}
\end{tikzcd},\]
and hence an isomorphism
\[ \overline{i_{M,i}} \tosim \overline{i_{M,i'}},\]
depending only on the choice of $B$. As a matter of fact, it doesn't
depend on this choice. To see this, we are reduced to comparing the
two isomorphisms arising from a $B$ and a $B'$ such that $B\subset
B'$, in which case the inclusion functor $g:B\to B'$ is aspheric
(prop.\ \ref{prop:75.1}). We leave the details of the verification
(consisting mainly of some diagram chasing and some compatibilities
between the $\lambda_i(F)$-isomorphisms of prop.\ \ref{prop:75.1} on
page \ref{p:214}) to the skeptical reader. As for transitivity for a
triple $(A,i)$, $(A',i')$, $(A'',i'')$, it now follows at once, by
using a $B$ suitable simultaneously\pspage{241} for all three.

Having thus well in hand the basic functor \eqref{eq:77.11} (which in
case of example \ref{ex:76.3} above with $M=\Ahat$, reduces to the
all-important functor $\HotOf_A\to\Hot$), we cannot but define an
asphericity structure to be \emph{modelizing}, as meaning that this
canonical functor is an equivalence of categories. In the case of
$M=\Ahat$ above, this means that $A$ is a pseudo-test category -- a
relatively weak test notion still. It appears as just as small bit
stronger than merely assuming the pair $(M,W)$ to be modelizing, i.e.,
to be a ``\emph{modelizer}'', namely assuming $\HotOf_M$ to be
equivalent (in some way or other\ldots) to \Hot. Maybe we should be a
little more cautious with the use of the word ``modelizing'' though,
and devise a terminology which should reflect very closely the
hierarchy of progressively stronger test notions
\begin{multline*}
  \text{(pseudo-test cat.)}
  \supset
  \text{(weak test cat.)}
  \supset\\
  \text{(test cat.)}
  \supset
  \text{(strict test cat.)}
  \supset
  \text{(contractors)}
\end{multline*}
which gradually has peeled out of our earlier reflections, by
pinpointing corresponding qualifications for an asphericity structure,
as being ``pseudo-modelizing'', ``weakly modelizing'', ``modelizing'',
``strictly modelizing'', and ???. It now appears that a little extra
reflection is needed here -- for today it's getting a little late
though!

\bigbreak
\presectionfill\ondate{18.6} and \ondate{19.6.}\par

\hangsection[Test functors and modelizing asphericity structures: the
\dots]{Test functors and modelizing asphericity structures: the
  outcome \texorpdfstring{\textup(}{(}at
  last!\texorpdfstring{\textup)}{)} of an easy
  ``observation''.}\label{sec:78}%
It is about time now to get a comprehensive treatment, in the context
of asphericity structures, of the relationships suggested time ago in
the ``observation'' and the ``silly question'' of section \ref{sec:46}
(pages \ref{p:94} and \ref{p:95}). The former is concerned with
test-functors from test categories to modelizers, the second more
generally with model-preserving maps between modelizers, having
properties similar to the map $M\to\Ahat$ stemming from a
test-functor. For the time being, as I found but little time for the
``extra reflection'' which seems needed, only the first situation is
by now reasonably clear in my mind.

As no evident series of ``asphericity structure''-notions has
appeared, paralleling the series of test-notions recalled by the end
of last Monday's reflections (see above), I'm going to keep
(provisionally only, maybe) the name of a \emph{modelizing asphericity
  structure} $(M,M\subas)$ as one for which the canonical functor
\begin{equation}
  \label{eq:78.1}
  \HotOf_M = W^{-1}M \to \Hot
  \tag{1}
\end{equation}
(cf.\ \eqref{eq:77.14}, p.\ \ref{p:240}) is an equivalence. This notion
at any rate is\pspage{242} satisfactory for formulating the following
statement, which comes out here rather tautologically, and which
however appears to me as exactly what I had been looking for in the
``observation'' recalled above:
\begin{theoremnum}\label{thm:78.1}
  Let $(M,M\subas)$ be a \emph{modelizing} asphericity structure, $A$
  a pseudo-test category, and
  \[i: A\to M\]
  a functor, factoring through $M\subas$ \textup(i.e., $i(a)$ is
  aspheric for any $a$ in $A$\textup). We assume moreover $M$ has a
  final object $e_M$. Then the six conditions \textup{\ref{it:75.4.i}}
  to \textup{\ref{it:75.4.iidblprime}} of prop.\ \ref{prop:75.4}
  \textup(p.\ \ref{p:236}\textup) \textup(the first of which
  expresses that $i$ is $M\subas$-aspheric\textup) are equivalent, and
  they are equivalent to the following condition:
  \begin{description}
  \item[\namedlabel{it:78.1.iii}{(iii)}]
    The functor $i^*:M\to\Ahat$ gives rise to a functor
    \[ \HotOf_M = W^{-1}M \to \HotOf_A = \scrWA^{-1}\Ahat\]
    \textup(i.e., \textup{\ref{it:75.4.iiprime}} of prop.\
    \ref{prop:75.4} holds, namely $i^*(W)\subset\scrWA$\textup),
    \emph{and} this functor moreover is an equivalence of categories.
  \end{description}
\end{theoremnum}
\begin{proof}
  For the first statement (equivalence of conditions \ref{it:75.4.i}
  to \ref{it:75.4.iidblprime}), by prop.\ \ref{prop:75.4} we need
  only show that the asphericity structure and $A$ as aspheric. But
  this follows from the assumptions and from the
\end{proof}
\begin{lemma}
  Let $(M,M\subas)$ be any modelizing asphericity structure, then a
  final object of $M$ is aspheric \textup(i.e., the structure is
  ``aspheric''\textup). In particular, if $A$ is a pseudo-test
  category \textup(i.e., $(\Ahat,\Ahatas)$ is a modelizing asphericity
  structure\textup), then $A$ is aspheric.
\end{lemma}

We only have to prove the first statement. It is immediate that $e_M$
is a final object of $W^{-1}M = \HotOf_M$ (this is valid whenever we
got a localization $W^{-1}M$ of a category $M$ with final object
$e_M$), hence its image in $\Hot=\scrW^{-1}\Cat$ is a final
object. Thus, we are reduced to proving the following
\begin{corollary}
  Let $(M,M\subas)$ be any asphericity structure, consider the
  composition
  \[ \varphi_M:M\to\HotOf_M=W^{-1} \xrightarrow{\mathrm{can.}} \Hot,\]
  then we get
  \[ M\subas = \set[\Big]{x\in\Ob M}{\text{$\varphi_M(x)$ is a final
      object of \Hot}}.\]
\end{corollary}

Indeed, using the construction of $\HotOf_M\to\Hot$ in terms of a
given $M\subas$-aspheric functor $B\to M$, we are reduced to the same
statement, with $(M,M\subas)$ replaced by $(\Bhat,\Bhatas)$. The
statement then reduced to: the object $B$ in \Cat{} is aspheric if{f}
its image in \Hot{} is a final object (this can be viewed also as the
particular case of the\pspage{243} corollary, when $M=\Cat$,
$M\subas=\Cat\subas$). Now, this follows at once from strong
saturation of \scrW, which (as we announced earlier) followed from
\ref{loc:1} to \ref{loc:3} (and will be proved in part \ref{ch:V} of
the notes). The reader who fears a vicious circle may till then
restrict use of the theorem to the case when we assume beforehand that
$e_M$ and $e_\Ahat$ are aspheric.

It is now clear that the six conditions of prop.\ \ref{prop:75.4} are
equivalent, and they are of course implied by
\ref{it:78.1.iii}. Conversely, they imply \ref{it:78.1.iii}, as
follows from the fact that in the canonical diagram (commutative up to
canonical isomorphism)
\[\begin{tikzcd}[baseline=(O.base),column sep=tiny]
  \HotOf_M\ar[rr]\ar[dr] & & \HotOf_A\ar[dl] \\
  & |[alias=O]| \Hot &
\end{tikzcd},\]
the two downwards arrows are equivalences, qed.

The theorem above seems to me to be exactly \emph{the} ``something
very simple-minded surely'' which I was feeling to get burningly
close, by the end of March, nearly three months ago (p.\ \ref{p:89});
at least, to be ``just it'' as far as the case of \emph{test-functors}
is concerned. We may equally view this theorem as giving the precise
relationship between the notion of a weak test functor or a \emph{test
  functor} (the latest version of which (in the context of
\scrW-notions) appears in section \ref{sec:65} (def.\ \ref{def:65.5}
and \ref{def:65.6}, pp.\ \ref{p:175} and \ref{p:176})), and the notion
of \emph{aspheric functors}, more precisely of $M\subas$-aspheric
functors, introduced lately (p.\ \ref{p:233}). This now is the moment
surely to check if the terminology of test functors and weak ones
introduced earlier, before the relevant notion of asphericity
structures was at hand, is really satisfactory indeed, and if needed,
adjust it slightly.

So let again
\[i:A\to M\]
be a functor, with $A$ small, and $(M,M\subas)$ an asphericity
structure. We don't assume beforehand, neither that $A$ is a
test-category or the like, nor that $(M,M\subas)$ be modelizing. We
now paraphrase def.\ \ref{def:65.5} (p.\ \ref{p:175}) of weak test
functors as follows:\scrcomment{it looks like there are only two
  conditions, but see \ref{it:78.1.1}--\ref{it:78.1.3} below\ldots}
\begin{definitionnum}\label{def:78.1}
  The functor $i$ above is called a \emph{pseudo-test functor}, if it
  satisfies the following three conditions:
  \begin{enumerate}[label=\alph*)]
  \item\label{it:78.1.a}
    The corresponding functor $i^*:M\to\Ahat$ is ``model-preserving'',
    by\pspage{244} which is meant here that
    \[ W = (i^*)^{-1}(\scrWA),\]
    \emph{and} the induced functor
    \[\HotOf_M=W^{-1}M \to \HotOf_A=\scrWA^{-1}\Ahat\]
    is an equivalence.
  \item\label{it:78.1.b}
    The functor $i_A:\Ahat\to\Cat$ is model preserving (for
    $\scrWA,\scrW$), which reduces to the canonical functor
    \[\HotOf_A=\scrWA^{-1}\Ahat\to\Hot=\scrW^{-1}\Cat\]
    being an equivalence (as we know already that
    $\scrWA=(i_A)^{-1}(\scrW)$, by definition of \scrWA).
  \end{enumerate}
\end{definitionnum}

Condition \ref{it:78.1.b} just means that $A$ is a pseudo-test
category, i.e., that the asphericity structure $(\Ahat,\Ahatas)$ it
defines is modelizing. By the corollary of lemma above, it implies
that $A$ is aspheric (which corresponds to condition \ref{it:65.E.c}
of def.\ \ref{def:65.5} in loc.\ cit.). Condition \ref{it:78.1.a}
comes in two parts, the first just meaning that $i$ is
$M\subas$-aspheric -- this translation being valid, at any rate, in
case we assume already $A$ aspheric, and the given asphericity
structure is aspheric, i.e., $M$ admits a final object $e_M$, and
$e_M$ is aspheric. We certainly do want a pseudo-test functor to be
(at the very least) $M\subas$-aspheric, so we should either strengthen
condition \ref{it:78.1.a} to this effect (which however doesn't look
as nice), or assume beforehand $(M,M\subas)$ aspheric. At any rate, if
we use the first variant of the definition, condition \ref{it:78.1.a}
in full then is equivalent (granting \ref{it:78.1.b}) to:
\begin{description}
\item[\namedlabel{it:78.1.aprime}{a')}]
  The functor $i$ is $M\subas$-aspheric, \emph{and} $(M,M\subas)$ is modelizing,
\end{description}
which in turn implies that $(M,M\subas)$ is aspheric. So we may as
well assume asphericity of $(M,M\subas)$ beforehand! At any rate, we
see that the notion we are after can be decomposed into three
conditions, namely:
\begin{enumerate}[label=\arabic*)]
\item\label{it:78.1.1} $A$ is a pseudo-test category, i.e.,
  $(\Ahat,\Ahatas)$ is modelizing.
\item\label{it:78.1.2} $(M,M\subas)$ is modelizing.
\item\label{it:78.1.3} The functor $i$ is $M\subas$-aspheric.
\end{enumerate}

The two first conditions are just conditions on $A$ and on
$(M,M\subas)$ respectively, the third is just the familiar asphericity
condition on $i$.

To get the notion of a \emph{weak test-functor} (def.\ \ref{def:65.5},
p.\ \ref{p:175}) we have to be just one step more specific in
\ref{it:78.1.1}, by demanding that $A$ be even a\pspage{245} weak test
category, namely that the functor
\[i_A^*=j_A:\Cat\to\Ahat\]
be model-preserving (which implies that $i_A$ is too).

Following def.\ \ref{def:65.6} (p.\ \ref{p:176}), we'll say that $i$
is a \emph{test functor} if $i$ \emph{and} the induced functors
$i_{/a}:A_{/a}\to A\to M$ are weak test functors. Using the definition
of a test category (def.\ \ref{def:65.3}, p.\ \ref{p:173}), we see
that this just means that the following conditions hold:
\begin{enumerate}[label=\arabic*')]
\item\label{it:78.1.1prime} $A$ is a test category.
\item\label{it:78.1.2prime} $(M,M\subas)$ is modelizing (i.e., same as
  \ref{it:78.1.2} above).
\item\label{it:78.1.3prime} The functor $i$ and the induced functors
  $i_{/a}:A_{/a}\to M$ are $M\subas$-aspheric.
\end{enumerate}

This last condition merits a name by itself, independently of other
assumptions:

\begin{definitionnum}\label{def:78.2}
  Let $(M,M\subas)$ be any asphericity structure, $A$ a small
  category, and $i:A\to M$ a functor. We'll say that $i$ is
  \emph{totally $M\subas$-aspheric} (or simply \emph{totally
    aspheric}, if no confusion is feared), if $i$ and the induced
  functors $i_{/a}:A_{/a}\to M$ are $M\subas$-aspheric (for any $a$ in
  $A$). We'll say that $i$ is \emph{locally $M\subas$-aspheric} (or
  simply \emph{locally aspheric}) if for any $a$ in $A$, the induced
  functor $i_{/a}:A_{/a}\to M$ is $M\subas$-aspheric.
\end{definitionnum}

Thus, $i$ is totally aspheric if{f} it is aspheric and locally
aspheric. On the other hand, $i$ is a test functor if{f} $A$ is a test
category, $(M,M\subas)$ is modelizing, and $i$ is totally aspheric.

In order for $i$ to be locally aspheric, it is n.s.\ that $i$ factor
through $M\subas$ and for any $x$ in $M\subas$, $i^*(x)$ be
\emph{aspheric over $e_\Ahat$}; if $B$ is any subcategory of $M$
containing $i(A)$ and generating the asphericity structure, it is
enough in this latter condition to take $x$ in $B$.

\begin{remarks}
  Thus, we see that the three gradations for the test-functor notion,
  as suggested by the definitions \ref{def:65.5} and \ref{def:65.6} of
  section \ref{sec:65} and now by the present context of asphericity
  structures, just amount to gradations for the test conditions on the
  category $A$ itself (namely, to be a pseudo-test, a weak test or
  just a plain test -category), and a two-step gradation on the
  asphericity condition for $i$ (namely, that $i$ be either just
  aspheric in the two first cases, or totally aspheric in the third),
  while these two asphericity conditions on $i$ are of significance,
  independently of any specific assumption which we may make on either
  $A$ or $(M,M\subas)$. This\pspage{246} seems to diminish somewhat
  the emphasis I had put formerly upon the notion of a test functor
  and its weak variant, and enhance accordingly the notion of an
  aspheric functor (with respect to a given asphericity structure) and
  the two related asphericity notions for $i$ which spring from it
  (namely, the notions of locally and of totally aspheric functors),
  which now seem to come out as the more relevant and the more general
  ones.
\end{remarks}

To be wholly happy, we still need the relevant reformulation, in terms
of asphericity structures, of the main result of section \ref{sec:65},
namely theorem \ref{thm:65.1} (p.\ \ref{p:177}) characterizing test
functors with values in \Cat, under the assumption that the objects
$i(a)$ be contractible. This theorem, I recall, has been the main
outcome of the ``grinding'' reflections taking their start with the
``observation'' on p.\ \ref{p:94} about ten days earlier. The appreach
then followed, as well as the contractibility assumption for the
objects $i(a)$ made in the theorem, retrospectively look awkward -- it
is clear that the relevant notions of asphericity structures, and of
aspheric functors into these, had been lacking. The theorem stated
last (p.\ \ref{p:242}) looks indeed a lot more satisfactory than the
former, except however in one respect -- namely that the asphericity
condition on $i$, in theorem \ref{thm:65.1} of p.\ \ref{p:177}, can be
expressed by asphericity over $e_\Ahat$ \emph{of just
  $i^*(\Simplex_1)$ alone}, rather than having to take $i^*(x)$ for
all $x$ in some fixed generating subcategory of $(M,M\subas)$
containing $i(A)$. To recover such kind of minimal criterion in a more
general case than \Cat, we'll have to relate the notion of an
asphericity structure with the earlier one of a \emph{contractibility
  structure}; the latter in our present reflections has faded somewhat
into the background, while at an earlier stage the homotopy and
contractibility notions had invaded the picture to an extent as to
overshadow and nearly bring to oblivion the magic of the ``asphericity
game''.

\bigbreak
\presectionfill\ondate{20.6.}\pspage{247}\par

\hangsection[Asphericity structure generated by a contractibility
\dots]{Asphericity structure generated by a contractibility structure:
  the final shape of the ``awkward main result'' on test
  functors.}\label{sec:79}%
I would like now to write out the relationship between asphericity
structures, and contractibility structures as defined in section
\ref{sec:51} \ref{subsec:51.D} (pages \ref{p:117}--\ref{p:119}). First
we'll need to rid ourselves of the smallness assumption for a
generating category of an asphericity structure:
\begin{propositionnum}\label{prop:79.1}
  Let $M$ be a \scrU-category \textup(\scrU{} being our basic
  universe\textup), and $N$ any full subcategory. The following
  conditions are equivalent:
  \begin{enumerate}[label=(\roman*),font=\normalfont]
  \item\label{it:79.1.i}
    There exists an asphericity structure $M\subas$ in $M$ such that
    \textup{\namedlabel{it:79.1.a}{(a)}}\enspace $M\subas\supset N$
    and \textup{\namedlabel{it:79.1.b}{(b)}}\enspace $N$ admits a
    \emph{small} full subcategory $N_0$ generating the asphericity
    structure, i.e., such that the inclusion functor $i_0$ from $N_0$
    to $M$ be $M\subas$-aspheric.
  \item\label{it:79.1.ii}
    There exists a small full subcategory $N_0$ of $N$, such that for
    any $x$ in $N$, $i_0^*(x)$ in $N_0\uphat$ be aspheric, where
    $i_0:N_0\to M$ is the inclusion functor.
  \item\label{it:79.1.iii}
    The couple $(N,N)$ is an asphericity structure.
  \end{enumerate}
  Moreover, when these conditions are satisfied, the asphericity
  structure $M\subas$ in \textup{\ref{it:79.1.i}} is unique.
\end{propositionnum}

Of course, we'll say it is the asphericity structure on $M$
\emph{generated} by the full subcategory $N$ of $M$, and the latter
will be called a \emph{generating subcategory} for the given
asphericity structure.

The proposition is a tautology, in view of the definitions and of
prop.\ \ref{prop:75.1} (p.\ \ref{p:232}) and its corollary. The form
\ref{it:79.1.ii} or \ref{it:79.1.iii} of the condition shows that it
depends only upon the category $N$, not upon the way in which $N$ is
embedded in a larger category $M$. We may call a category $N$
satisfying condition \ref{it:79.1.iii} above an \emph{aspherator}, by
which we would like to express that this category represents a
standard way of generating asphericity structures, through any full
embedding of $N$ into a category $M$. This condition is automatically
satisfied if $N$ is small, more generally it holds if $N$ is
equivalent to a small category. It should be kept in mind that the
condition depends both on the choice of the basic universe \scrU{} and
on the choice of the basic localizer \scrW{} in the corresponding
large category \Cat. It looks pretty sure that the condition is not
always satisfied (I doubt it is for $N=\Sets$ say), but I confess I
didn't sit down to make an explicit example.

Let\pspage{248} now $(M,M\subc)$ be a contractibility structure, such
that there exists a \emph{small} full subcategory $C$ of $M$ which
generates the contractibility structure, i.e., such that (a)\enspace
the objects of $C$ are ``contractible'', i.e., are in $M\subc$ and
(b)\enspace any object of $M\subc$ is $C$-contractible (i.e.,
contractible for the homotopy interval structure admitting the
intervals made up with objects of $C$ as a generating
family). Independently of any smallness assumption upon $C$, we gave
in section \ref{sec:51} \ref{subsec:51.D} (p.\ \ref{p:118}), under the
name of ``basic assumption'' \ref{cond:51.Bas.4}, the n.s.\ condition
on a full subcategory $C$ of a given category $M$, in order that $C$
can be viewed as a generating set of contractible objects, for a
suitable contractibility structure
\[ M\subc\subset \Ob M\]
in $M$ (which is uniquely defined by $C$). It turns out that in case
$C$ contains a final object of $M$, and is stable under binary
products in $M$, the condition \ref{cond:51.Bas.4} depends only upon
the category structure of $C$, and not upon the particular way this
category is embedded in another one $M$.

In loc.\ sit.\ we did not impose, when defining a contractibility
structure, a condition that there should exist a small set of
generators for the structure. From now on, we'll assume that the
(possibly large) categories $M$ we are working with are
\scrU-categories, and that ``contractibility structure'' means
``\scrU-contractibility structure'', namely one such that $M$ admit a
small full subcategory $C$, generating the structure.

We recall too that in the definition of a contractibility structure
$(M,M\subc)$, it has always been understood that $M$ is stable under
finite products.

\begin{propositionnum}\label{prop:79.2}
  With the conventions above, let $(M,M\subc)$ be any contractibility
  structure. Then:
  \begin{enumerate}[label=\alph*),font=\normalfont]
  \item\label{it:79.2.a}
    The full subcategory $M\subc$ of $M$ generates an asphericity
    structure $M\subas$ in $M$ \textup(cf.\ prop.\
    \ref{prop:79.1}\textup).
  \item\label{it:79.2.b}
    Any small full subcategory $C$ of $M$ which generates the
    contractibility structure and such that $\Chat$ be totally
    aspheric, \emph{generates} the asphericity structure $M\subas$,
    i.e., the inclusion functor $i:C\to M$ is $M\subas$-aspheric.
  \item\label{it:79.2.c}
    The asphericity structure $M\subas$ is totally aspheric \textup(cf.\
    prop.\ \ref{prop:76.5}, p.\ \ref{p:239}\textup).
  \end{enumerate}
\end{propositionnum}

The\pspage{249} first statement \ref{it:79.2.a} can be rephrased, by
saying that the category $M\subc$ of contractible objects of $M$ is an
\emph{aspherator}. To prove this, we use the fact that $M\subc$ is
stable under finite products, which implies that we can find a small
full subcategory $C$ of $M\subc$, stable under such products, and
which generates the contractibility structure. The stability condition
upon $C$ implies that \Chat{} is totally aspheric. Therefore,
\ref{it:79.2.a} and \ref{it:79.2.b} will be proved, if we prove that
any small full subcategory $C$ of $M\subc$, such that \Chat{} is
totally aspheric, satisfies the conditions of prof.\ \ref{prop:79.1}
(with $N=M\subc$), namely that for any $x$ in $M\subc$, $i^*(x)$ is an
aspheric object of \Chat, where $i:C\to M$ (or $C\to M\subc$,
equivalently) is the inclusion functor. But from the fact that $i^*$
commutes with finite products, it follows that $i^*(x)$ is
contractible in \Chat, for the homotopy interval structure admitting
as a generating family of homotopy intervals the intervals in \Chat{}
made up with objects of $C$. The assumption of total asphericity upon
\Chat{} implies that these homotopy intervals are aspheric over
$e_\Chat$, and from this follows (as already used a number of times
earlier) that $i^*(x)$ too is aspheric over $e_\Chat$, and hence
aspheric as $e_\Chat$ is aspheric (because of the assumption of total
asphericity). This proves \ref{it:79.2.a} and \ref{it:79.2.b}, and
\ref{it:79.2.c} follows via prop.\ \ref{prop:76.5} (p.\ \ref{p:239}).

We'll call of course the asphericity structure described in prop.\
\ref{prop:79.2} the \emph{asphericity structure generated by the given
  contractibility structure}.

\begin{remarknum}
  Clearly, not every asphericity structure can be generated by a
  contractibility structure, as a necessary condition (presumably not
  a sufficient one) is total asphericity. We don't expect either, in
  case it can, that the generating contractibility structure is
  uniquely defined; however, we do expect in this case that there
  should exist a canonical (largest) choice -- we'll have to come back
  upon this in due course. For the time being, let's only remark that
  all modelizing asphericity structures met with so far, it seems, do
  come from contractibility structures.
\end{remarknum}

It occurs to me that the last statement was a little hasty -- after
all we have met with test categories which are not strict ones, hence
$(\Ahat,\Ahatas)$ is a modelizing asphericity structure (it would be
enough even that $A$ be a pseudo-test category), which isn't totally
aspheric, and a fortiori does not come from a contractibility
structure. Thus,\pspage{250} I better correct the statement, to the
effect that, it seems, all modelizing \emph{totally aspheric}
asphericity structures met with so far are generated by suitable
contractibility structures.

Let $M$ be a category endowed with a contractibility structure
$M\subc$, hence also with an asphericity structure $M\subas$, and let
\[i: A\to M\]
be a functor from a small category $A$ to $M$. We want to give n.s.\
conditions for $i$ to be a \emph{test functor}, in terms of homotopy
notions in $M$ and in \Ahat. To this end, it seems necessary to
refresh memory somewhat and recall some relevant notions which were
developed in part \ref{ch:III} of our notes (sections \ref{sec:54} and
\ref{sec:55}).

It will be convenient to call an object $F$ of \Ahat{} \emph{locally
  aspheric} (resp.\ \emph{totally aspheric}), if its product in
\Ahat{} by any object of $A$, and hence also its product by any
aspheric object of \Ahat, is aspheric (resp.\ and if moreover $F$
itself is aspheric). With this terminology, \Ahat{} is totally
aspheric if{f} every aspheric object of \Ahat{} is totally aspheric,
and if moreover $A$ is aspheric, i.e., $e_\Ahat$ is aspheric. Note
that $F$ is locally aspheric if{f} the map $F\to e_\Ahat$ is aspheric,
or what amounts to the same, if this map is universally in \scrWA, or
equivalently, if the corresponding functor $A_{/F}\to A$ is
aspheric. If $A$ itself is aspheric, and in this case only, this
condition implies already that $F$ is aspheric, i.e., that $F$ is
totally aspheric. We'll denote by
\[\Ahattotas\quad(\text{resp.}\quad\Ahatlocas)\]
the full subcategory of \Ahat{} made up with the totally aspheric
(resp.\ locally aspheric) objects of \Ahat. Thus we get
\[\begin{tikzcd}[baseline=(O.base),column sep=tiny]
  & \Ahattotas = \Ahatlocas \sand \Ahatas
  \ar[dl,hook]\ar[dr,hook] & \\
  \Ahatlocas & & |[alias=O]| \Ahatas
\end{tikzcd}.\]

Now recall that (section \ref{sec:54}) in terms of the set \scrWA{} of
weak equivalences in \Ahat, we constructed a homotopy structure on
\Ahat, more specifically a homotopy interval structure, admitting as a
generating family of homotopy intervals
\[\bI = (I, \delta_0,\delta_1)\]
the set of intervals such that $I$ be a locally aspheric object of
\Ahat, i.e.,
\[ I\in\Ob \Ahatlocas .\]
Let\pspage{251}
\[ h = h_\scrWA\]
be this homotopy structure, hence a corresponding notion of
\emph{$h$-equivalence} or \emph{$h$-homotopy} $\hsim$ for arrows in
\Ahat, a corresponding notion of $h$-homotopisms, i.e., a set of
arrows
\[ W^h\subset\Fl\Ahat, \quad\text{such that}\quad W^h\subset\scrWA,\]
a notion of $h$-homotopy interval (namely an interval such that
$\delta_0$ and $\delta_1$ be $h$-homotopic, for which it is
sufficient, but not necessary, that $I$ be locally aspheric\ldots),
and last not least, a notion of \emph{contractible objects}, making up
a full subcategory
\[ \Ahatc \subset \Ahat, \quad\text{such that}\quad
\Ahatc \subset \Ahatlocas.\]
The latter inclusion, in case $A$ is aspheric, can be equally written
\[\Ahatc \subset \Ahattotas \quad(\text{if $A$ aspheric}).\]

Coming back now to the contractibility structure $(M,M\subc)$, and a
functor $i:A\to M$, we are interested in the corresponding functor
\[u=i^*: M \to M'=\Ahat,\]
where both members will be viewed as being endowed with their
respective homotopy structures -- the one of $M$ being of the most
restrictive type envisioned in section \ref{sec:51}, namely it is
defined in terms of a contractibility structure, whereas the one of
$M'$ is a priori defined in terms of a homotopy interval structure,
but not necessarily in terms of a contractibility structure. Now this
situation has been described in section \ref{sec:53}, as far as
compatibility conditions with homotopy structures are concerned,
independently by the way of any special assumption on $M'$ (such as
being a category of presheaves on some small category $A$), or on the
functor $u$, except for commuting to finite products. Compatibility of
$u$ with the homotopy structures on $M,M'$ can be expressed by either
one of the following six conditions, which are all equivalent:
\begin{description}
\item[\namedlabel{cond:79.H1}{H~1)}]
  $u$ transforms homotopic arrows of $M$ into homotopic arrows in $M'$.
\item[\namedlabel{cond:79.H2}{H~2)}]
  $u$ transforms homotopisms in $M$ into homotopisms in $M'$.
\item[\namedlabel{cond:79.H3}{H~3)}]
  $u$ transforms any homotopy interval $\bI=(I,\delta_0,\delta_1)$ in
  $M$ into a homotopy interval in $M'$ (i.e., if two sections of an
  object $I$ over $e_M$ are homotopic, so are $u(\delta_0)$ and $u(\delta_1)$).
\item[\namedlabel{cond:79.H3prime}{H~3')}]
  Same as \ref{cond:79.H3}, but \bI{} being restricted to be in a
  given family of homotopy intervals, generating for the homotopy
  structure in $M$.
\item[\namedlabel{cond:79.H4}{H~4)}]
  $u$\pspage{252} transforms any contractible object $x$ of $M$ into a
  contractible object of $M'$.
\item[\namedlabel{cond:79.H4prime}{H~4')}]
  Same as \ref{cond:79.H4}, with $x$ being restricted to be in a given
  subcategory $C$ of $M$, generating for the contractibility structure
  of $M$.
\end{description}
\begin{remarknum}\label{rem:79.2}
  It should be noted that the condition upon $C$ stated in
  \ref{cond:79.H4prime}, namely that the (given) $C$ should be
  generating for the contractibility structure $M\subc$ of $M$, means
  exactly two things: (a)\enspace $C\subset M\subc$, and (b)\enspace
  the family of all intervals $\bI=(I,\delta_0,\delta_1)$ made up with
  objects of $C$ (these intervals are necessarily homotopy intervals
  for the homotopy structure of $M$) \emph{generates} the homotopy
  structure of $M$, i.e., two arrows in $M$ are homotopic if{f} they
  can be joined by a chain of arrows, two consecutive among which
  being related by an elementary homotopy involving an interval of
  that family. This reminder being made, it follows that $C$ is
  equally eligible for applying criterion \ref{cond:79.H3prime}, which
  means that we get still another equivalent formulation of
  \ref{cond:79.H4prime}, by demanding merely that for any two sections
  $\delta_0,\delta_1$ of $x$ over $e_M$, the corresponding sections of
  $u(x)$ should be homotopic.
\end{remarknum}

After these preliminaries, we can state at last the following
generalization of the main result of section \ref{sec:65} (th,\
\ref{thm:65.1}, p.\ \ref{p:176}), concerning test functors with values
in \Cat:
\begin{theoremnum}\label{thm:79.1}
  Let $(M,M\subc)$ be a contractibility structure \textup(cf.\ p.\
  \ref{p:248}\textup), $C$ a small full subcategory of $M$ generating
  the contractibility structure \textup(cf.\ remark \ref{rem:79.2}
  above\textup), $A$ a small category, and $i:A\to M$ a functor,
  factoring through $M\subc$. We consider $M$ as endowed equally with
  the asphericity structure $M\subas$ generated by $M\subc$
  \textup(cf.\ prop.\ \ref{prop:79.2}\textup). Then the following
  conditions on $i$ are equivalent:
  \begin{enumerate}[label=(\roman*),font=\normalfont]
  \item\label{it:79.thm1.i}
    The functor $i^*:M\to\Ahat$ deduced from $i$ is compatible with
    the homotopy structures on $M$ and on \Ahat{} \textup(cf.\ pages
    \ref{p:250}--\ref{p:251} for the latter\textup), i.e., $i^*$
    satisfies either one of the six equivalent conditions
    \textup{\ref{cond:79.H1}} to \textup{\ref{cond:79.H4prime}} above.
  \item\label{it:79.thm1.ii}
    The functor $i^*$ is locally $M\subas$-aspheric \textup(def.\
    \ref{def:78.2}, p.\ \ref{p:245}\textup), i.e., for any $x$ in
    $M\subas$, $i^*(x)$ is a locally aspheric object of \Ahat{}
    \textup(cf.\ p.\ \ref{p:250}\textup), i.e., $i^*(x)$ is aspheric
    over $e_\Ahat$.
  \item\label{it:79.thm1.iii}
    Same as \textup{\ref{it:79.thm1.ii}}, but with $x$ restricted to
    be in $C$.
  \end{enumerate}
\end{theoremnum}
\begin{corollarynum}\label{cor:79.1}
  In order for $i$ to be totally $M\subas$-aspheric \textup(def.\
  \ref{def:78.2}, p.\ \ref{p:245}\textup)\pspage{253} it is n.s.\ that
  $A$ be aspheric, and that the equivalent conditions of th.\
  \ref{thm:79.1} be satisfied.
\end{corollarynum}

Indeed, it follows at once from the definitions and from the fact that
the asphericity structure $(M,M\subas)$ is aspheric, i.e., that $e_M$
is aspheric, that $i$ is totally $M\subas$-aspheric if{f} it is
locally $M\subas$-aspheric (i.e., condition \ref{it:79.thm1.ii}), and
if moreover $A$ is aspheric, hence the corollary.

\begin{corollarynum}\label{cor:79.2}
  In order for $i$ to be a test functor, it is n.s.\ that $A$ be a
  test category, that $(M,M\subas)$ be modelizing, and that the
  equivalent conditions of theorem \ref{thm:79.1} be satisfied.
\end{corollarynum}

This follows from cor.\ \ref{cor:79.1} and the reformulation of the
notion of a test-functor, given p.\ \ref{p:245}.

This corollary contains the main result of section \ref{sec:65} as a
particular case, when taking $M=\Cat$ with the usual contractibility
structure, and $C=\{\Simplex_1\}$, except that in loc.\ sit.\ we did
not have to assume beforehand that $A$ be a test category, but only
that $A$ is aspheric: this condition, plus condition
\ref{it:79.thm1.iii} above (namely, $i^*(\Simplex_1)$ locally
aspheric) implies already that $A$ is a test category. In order to get
also this extra result, we state still another corollary:
\begin{corollarynum}\label{cor:79.3}
  Assume
  \[\bI=(I,\delta_0,\delta_1,\mu), \quad\text{with}\quad
  I\in\Ob C,\]
  is a \emph{multiplicative interval} in $M$, i.e., an interval
  endowed with a multiplication $\mu$, admitting $\delta_0$ as a left
  unit and $\delta_1$ as a left zero element \textup(cf.\ section
  \ref{sec:49}, p.\ \ref{p:108} and section \ref{sec:51}, p.\
  \ref{p:120} -- where such an interval was provisionally called a
  ``contractor''\textup). Assume moreover that for any $x$ in
  $M\subc$, the two compositions $x \to $
  \begin{tikzcd}[cramped]
    e_M\ar[r,shift left=1pt,"\ensuremath{\delta_0,\delta_1}"]\ar[r,shift right=2pt] &
    I
  \end{tikzcd} are distinct \textup(which is the case for instance
  if $\Ker(\delta_0,\delta_1)$ exists in $M$ and is a strict initial
  object $\varnothing_M$ of $M$ \textup(i.e., an initial object such
  that any map $x\to\varnothing_M$ in $M$ is an isomorphism\textup), and
  moreover $\varnothing_M\notin M\subc$\textup). Then the conditions of
  th.\ \ref{thm:79.1} imply that $A$ is a local test category, and
  hence a test category provided $A$ is aspheric.
\end{corollarynum}

Indeed, $i^*(\bI)$ is a multiplicative interval in \Ahat{} which is
locally aspheric, and (as follows immediately from the assumptions on
\bI) separating -- hence \Ahat{} is a local test
category.

This\pspage{254} array of immediate corollaries of th.\ \ref{thm:79.1}
do convince me that this statement is indeed ``the'' natural
generalization of the ``awkward'' main result of section
\ref{sec:51}. All we have to do is to prove theorem \ref{thm:79.1}
then.

The three conditions of theorem \ref{thm:79.1} can be rewritten simply
as
\begin{description}
\item[\ref{it:79.thm1.i}]
  $i^*(C) \subset \Ahatc$ (using \hyperref[cond:79.H4prime]{H~4'}),
\item[\ref{it:79.thm1.ii}]
  $i^*(M\subas) \subset \Ahatlocas$,
\item[\ref{it:79.thm1.iii}]
  $i^*(C) \subset \Ahatlocas$,
\end{description}
and because of
\[ C \subset M\subas, \quad \Ahatc\subset\Ahatlocas,\]
it follows tautologically that \ref{it:79.thm1.i} and
\ref{it:79.thm1.ii} both imply \ref{it:79.thm1.iii}. On the other
hand, \ref{it:79.thm1.iii} $\Rightarrow$ \ref{it:79.thm1.i} by
criterion \ref{cond:79.H3prime}, and the definition of the homotopy
interval structure of \Ahat{} in terms of \Ahatlocas. Thus, we are
left with proving \ref{it:79.thm1.i} $\Rightarrow$
\ref{it:79.thm1.ii}. But \ref{it:79.thm1.i} can be rewritten as
\[i^*(M\subc) \subset \Ahatc,\]
which implies
\begin{equation}
  \label{eq:79.star}
  i^*(M\subc) \subset \Ahatlocas.\tag{*}
\end{equation}
That the latter condition implies \ref{it:79.thm1.ii} now follows from
the corollary to the following tautology (with $N=M\subc$), which
should have been stated as a corollary to prop.\ \ref{prop:79.1} above
(p.\ \ref{p:247}):
\begin{lemma}
  Let $(M,M\subas)$ be an asphericity structure, generated by the full
  subcategory $N$ of $M$, let $A$ be a small category, and $i:A\to M$
  a functor \emph{factoring through $N$}. Then $i$ is
  $M\subas$-aspheric if{f} for any $x$ in $N$, $i^*(x)$ is an aspheric
  object of $\Ahat$.
\end{lemma}
\begin{corollary}
  The functor $i$ is locally $M\subas$-aspheric \textup(resp.\ totally
  $M\subas$-aspheric\textup) if{f} $i^*(N)\subset\Ahatlocas$
  \textup(resp.\ $i^*(N)\subset\Ahattotas$\textup).
\end{corollary}
\begin{remarknum}
  Assume\pspage{255} in theorem \ref{thm:79.1} that $A$ is totally aspheric. Then
  condition \ref{it:79.thm1.ii} just means that $i$ is
  $M\subas$-aspheric, and condition \ref{it:79.thm1.iii} that $i^*(x)$
  is aspheric for any $x$ in $C$ (as we got
  $\Ahatlocas=\Ahatas=\Ahattotas$). If moreover $C$ satisfies the
  condition of cor.\ \ref{cor:79.2} above, then these conditions imply
  that $A$ is a strict test category, and if $(M,M\subas)$ is
  modelizing, that $i$ is a test functor as stated in corollary
  \ref{cor:79.2}. These observations sum up the substance of the
  restatement of the main result of section \ref{sec:51}, given in
  theorem \ref{thm:65.2} of p.\ \ref{p:178}, in the present general
  context.
\end{remarknum}

\hangsection{Reminders and questions around canonical
  modelizers.}\label{sec:80}%
In the preceding section, we associated to any contractibility
structure $M\subc$ on a category $M$, an asphericity structure
$M\subas$ ``generated'' by the former in a natural sense. It is this
possibility of associating (in a topologically meaningful way) an
asphericity structure to a given contractibility structure, which
singles out the latter structure type, among the three essential
distinct ``homotopy flavored'' kind of structures developed at length
in sections \ref{sec:51} and \ref{sec:52}, in preference to the two
weaker notions of a homotopy interval structure, and of a homotopism
structure (or, equivalently, or a ``homotopy relation''). The
association
\begin{equation}
  \label{eq:80.star}
  M\subc\mapsto M\subas \quad(\leftrightarrow \text{corresponding
    notion of weak equivalence $W\suba$ in $M$})
  \tag{*}
\end{equation}
finally carried through in the last section, had been foreshadowed
earlier (cf.\ p.\ \ref{p:110} and p.\ \ref{p:142}), but was pushed off
for quite a while, in order to ``give precedence'' to the other
approach in view by then towards more general test functors than
before, leading up finally to the ``awkward main result'' of section
\ref{sec:65}. In \eqref{eq:80.star}, as the asphericity structure
$M\subas$ is totally aspheric and a fortiori aspheric, $M\subas$ and
the corresponding notion of weak equivalence $W\suba=W_{M\subc}$ (which
of course should not be confused with the notion of homotopism
associated to $M\subc$, giving a considerably smaller set of arrows
$W\subc$) determine each other mutually. Till the writing up of
section \ref{sec:51}, it was the aspect ``weak equivalence'' $W\suba$
which was in the fore, whereas the conceptually more relevant aspect
of ``aspheric objects'' did not appear in full light before I was
through with grinding out the ``awkward approach'' (cf.\ p.\
\ref{p:188}).

It is time now to remember the opposite association
\begin{equation}
  \label{eq:80.starstar}
  W\suba\mapsto\text{homotopy structure $h_{W\suba}$,}\tag{**}
\end{equation}
associating\pspage{256} to a notion of weak equivalence in $M$, i.e.,
to any saturated subset
\[W\suba \subset \Fl(M),\]
a corresponding homotopy structure $h_{W\suba}$, a homotopy interval
structure as a matter of fact (section \ref{sec:54}, p.\
\ref{p:131}). Here we are primarily interested of course in the case
when $W\suba$ is associated to a given asphericity structure $M\subas$
in $M$, which we may as well assume to be aspheric, so to be sure that
$W\suba$ and $M\subas$ determine each other. I recall that the weak
interval structure $h_{W\suba}$ can be described by the generating
family of homotopy intervals, consisting of all intervals
\[ \bI=(I,\delta_0,\delta_1) \]
such that $I\to e_M$ be universally in $W\suba$, i.e., such that for
any object $x$ in $M$, the projection
\[ x\times I\to x\]
be in $W\suba$. For pinning down further the exact relationship
between contractibility structures (which may be viewed as just
special types of homotopy interval structures) and asphericity
structures, we are thinking of course more specifically of
\emph{totally aspheric} asphericity structures, in view of prop.\
\ref{prop:79.2} \ref{it:79.2.c} of the preceding section (p.\
\ref{p:248}). It is immediate in this case that for an object $I$ of
$M$, $I\to e_M$ is in $\mathrm UW\suba$ (i.e., is universally in
$W\suba$) if (and only if, of course) $I$ is aspheric, i.e., $I\to
e_M$ is in $W\suba$. The most relevant questions which come up here,
now seem to me the following:

\namedlabel{q:80.1}{1)}\enspace If $M\subas$ is generated by a
contractibility structure $M\subc$, is the homotopy structure
$h_{W\suba}$ associated to $M\subas$ (indeed, a homotopy interval
structure as recalled above) also the one defined by $M\subc$, using
intervals in $M\subc$ as a generating family of homotopy intervals?

\namedlabel{q:80.2}{2)}\enspace Conversely, what extra conditions on a
given asphericity structure $M\subas$ on $M$ are needed (besides total
asphericity) to ensure that the corresponding homotopy structure
$h_{W\suba}$ on $M$ comes from a contractibility structure $M\subc$
(i.e., admits a generating family of homotopy intervals which are
\emph{contractible}),\pspage{257} and that moreover $M\subc$
\emph{generates} $M\subas$?

Before looking up a little these questions, I would like however to
carry through at once the ``idyllic picture'' of canonical modelizers,
foreshadowed in section \ref{sec:50} (p.\ \ref{p:110}), as I feel that
this should be possible at present at no costs. Taking into account
the reflections of the later sections \ref{sec:57} and \ref{sec:59},
we get the following set-up.

Let $M$ be a \scrU-category, stable under finite products, endowed
with a functor
\begin{equation}
  \label{eq:80.1}
  \pi_0^M\text{ or }\piz: M\to\Sets,\tag{1}
\end{equation}
on which we make no assumptions for the time being. We are thinking of
the example when $M$ is a totally $0$-connected category (cf.\ prop.\
on page \ref{p:142} for this notion) and \piz{} is the ``connected
components'' functor, or when $M=\Spaces$ and \piz{} corresponds to
taking sets of arc-wise connected components. According to section
\ref{sec:54} (p.\ \ref{p:131}), we introduce a corresponding homotopy
interval structure $h$ on $M$, admitting the generating family of
homotopy intervals made up with those intervals
$\bI=(I,\delta_0,\delta_1)$ for which $I\to e_M$ is ``universally in
$W_\piz$'', namly
\[\piz(x\times I)\to\piz(x) \quad\text{bijective for any $x$ in
  $M$.}\]
(Under suitable conditions on \piz, this homotopy interval structure
$h$ on $M$ is the widest one ``compatible with \piz'' in the sense of
page \ref{p:130}, namely such that \piz{} transforms homotopisms into
isomorphisms -- cf.\ proposition p.\ \ref{p:133}.) We are interested
in the case when this homotopy structure on $M$ can be described by a
contractibility structure $M\subc$ on $M$, which is then unique of
course, hence well-defined in terms of \piz. Therefore, the
asphericity structure generated by $M\subc$ is equally well defined in
terms of \piz, and likewise the corresponding notion $W\suba$ of
``weak equivalence''. We then get a canonical functor
\begin{equation}
  \label{eq:80.2}
  \HotOf_M = W\suba^{-1}M \to \Hot\tag{2}
\end{equation}
(section \ref{sec:77}, p.\ \ref{p:239}). We'll have to find still a
suitable extra condition on the functor \piz, implying that this
functor is canonically isomorphic deduced from \eqref{eq:80.2} by
composing with the canonical functor\pspage{258}
\begin{equation}
  \label{eq:80.3}
  \Hot \xrightarrow{\piz} \Sets,
  \tag{3}
\end{equation}
which can be defined using a very mild extra condition on the basic
localizer \scrW{} (namely, $f\in\scrW$ implies $\piz(f)$ bijective,
cf.\ condition \ref{it:64.La} on page \ref{p:165}). Thus, there seems
to be a little work ahead after all -- in order to deduce something
like a one to one correspondence, say, between pairs $(M,\piz)$
satisfying suitable conditions, and certain types of asphericity
structures $(M,M\subas)$ (which will have to be assumed totally
aspheric, and presumably a little more still).

The case of special interest to us is the one when the asphericity
structure we get on $M$ in terms of the functor $\pi_0^M$ is
\emph{modelizing}, hence even strictly modelizing (i.e., $(M,W\suba)$
is a strict modelizer), as $(M,M\subas)$ is totally aspheric. If we
assume moreover that the category $M$ is totally $0$-connected and
that \piz{} is just the ``connected components'' functor, then the
modelizing asphericity structure we got on $M$ is canonically
determined by the mere category structure of $M$, and deserves
therefore to be called \emph{the canonical \textup(modelizing\textup)
  asphericity structure on $M$}. A \emph{canonical modelizer} $(M,W)$
is a modelizer which can be obtained from a canonical asphericity
structure $(M,M\subas)$ by $W=W\suba=$ corresponding set of weak
equivalences (for $M\subas$).

\begin{remark}
  The slight sketchy definition we just gave for the canonical
  modelizing asphericity structures, and accordingly for the canonical
  modelizers, is essentially complete. The point however which
  requires clarification is the relationship between the ``connected
  components'' functor on the corresponding (totally $0$-connected)
  category $M$, and the composition of the canonical functors
  \eqref{eq:80.2} and \eqref{eq:80.3}.
\end{remark}

\bigbreak

\presectionfill\ondate{25.6.}\pspage{259}\par

\hangsection[Contractibility as the common expression of homotopy,
\dots]{Contractibility as the common expression of homotopy,
  asphericity and \texorpdfstring{$0$}{0}-connectedness
  notions. \texorpdfstring{\textup(}{(}An overall review of the
  notions met with so far.\texorpdfstring{\textup)}{)}}\label{sec:81}%
I didn't find much time since Monday for mathematical pondering -- the
little I got nonetheless has been enough for convincing myself that
things came out more nicely still than I expected by then. One main
point being that, provided the basic localizer satisfies the mild
extra assumption \ref{loc:4} below, any contractibility structure
$M\subc$ on a category $M$ with finite products can be recovered, in
the simplest imaginable way, in terms of the associated asphericity
structure $M\suba$, or equivalently, in terms of the corresponding set
$W\suba$ of ``weak equivalences''; namely, $M\subc$ is the set of
contractible objects in $M$, for the homotopy interval structure
defined in terms of all intervals made up with objects of
$M\suba$. This implies that the canonical map
\begin{equation}
  \label{eq:81.1}
  \Homtp_4(M) = \Cont(M) \hookrightarrow \WAsph(M),\tag{1}
\end{equation}
from the set of contractibility structures on $M$ to the set of
asphericity structures on $M$ relative to the basic localizer \scrW{}
(or ``\emph{\scrW-asphericity structures}''), is \emph{injective}. In
other words, we may view a contractibility structure (on a category
$M$ stable under finite products), which is an \emph{absolute} notion
(namely independent of the choice of a basic localizer \scrW), as a
``\emph{particular case}'' of a \scrW-asphericity structure (depending
on the choice of \scrW), namely as ``equivalent'' to a
\scrW-asphericity structure, satisfying some extra conditions which
we'll have to write down below.

This pleasant fact associates immediately with two related ones. The
first is just a reminder of our reflections of sections \ref{sec:51}
and \ref{sec:52}, namely that the set of contractibility structures on
$M$ can be viewed as one among four similar sets of ``homotopy
structures'' on $M$
\begin{equation}
  \label{eq:81.2}
  \Homtp_4(M) \hookrightarrow \Homtp_3(M) \hookrightarrow \Homtp_2(M)
  \tosim \Homtp_1(M),\tag{2}
\end{equation}
corresponding to the four basic ``homotopy notions'' met with so far,
namely (besides contractibility structure $\Homtp_4$) the
\emph{homotopy interval structures} ($\Homtp_3$), the \emph{homotopism
  structures} ($\Homtp_2$), and the \emph{homotopy relations} between
maps ($\Homtp_1$). In the sequel, if
\begin{equation}
  \label{eq:81.3}
  M\subc \subset \Ob M\tag{3}
\end{equation}
is a given contractibility structure on $M$, we'll denote by
\begin{equation}
  \label{eq:81.4}
  J\subc \subset \Int(M), \quad
  W\subc \subset \Fl(M), \quad
  R\subc \subset \Fl(M)\times\Fl(M)
  \tag{4}
\end{equation}
the corresponding other three homotopy structures on $M$, where
$\Int(M)$ denotes the set of all ``intervals''
$\bI=(I,\delta_0,\delta_1)$ in $M$, i.e., objects of $M$ endowed with
two sections $\delta_0,\delta_1$ over the final object $e_M$ of $M$.

The\pspage{260} second fact alluded to above is concerned with
behavior of the \scrW-asphericity notions, for varying \scrW, more
specifically for a pair
\begin{equation}
  \label{eq:81.5}
  \scrW \subset\scrW'\subset\Fl(\Cat)\tag{5}
\end{equation}
of two basic localizers, \scrW{} and \scrW', such that \scrW{}
``refines'' \scrW'. It then follows that for any small category $A$,
we have
\begin{equation}
  \label{eq:81.6}
  \scrWA \subset \scrW'_A,\tag{6}
\end{equation}
and accordingly, that \emph{for any \scrW-asphericity structure}
\begin{equation}
  \label{eq:81.7}
  M_\scrW \subset \Ob M\tag{7}
\end{equation}
\emph{on $M$, there exists a unique $\scrW'$-asphericity structure
  $M_{\scrW'}$ on $M$,}
\begin{equation}
  \label{eq:81.8}
  M_\scrW \subset M_{\scrW'}, \tag{8}
\end{equation}
\emph{such that for any small category $A$\kern1pt, a functor $A\to M$ which
  is $M_\scrW$-aspheric \textup(with respect to \scrW\textup) is also
  $M_{\scrW'}$-aspheric \textup(with respect to
  $\scrW'$\textup)}. This is merely a tautology, which we didn't state
earlier, because there was no compelling reason before to look at what
happens when \scrW{} is allowed to vary. Thus, we get a canonical map
\begin{equation}
  \label{eq:81.9}
  \WAsph(M) \to \WprimeAsph(M),\tag{9}
\end{equation}
with the evident transitivity property for a triple
\[\scrW \subset \scrW' \subset \scrW'',\]
in other words we get a functor
\[ \scrW \mapsto \WAsph(M)\]
from the category of all basic localizers (the arrows between
localizers being inclusions \eqref{eq:81.5}) to the category of
sets. The relation of the canonical inclusion \eqref{eq:81.1} with
this functorial dependence of $\WAsph(M)$ on \scrW{} is expressed in
the commutativity of
\[\begin{tikzcd}[baseline=(O.base),column sep=-6pt]
  & \Homtp_4(M) = \Cont(M)
  \ar[dl,hook]\ar[dr,hook] & \\
  \WAsph(M)\ar[rr] & & |[alias=O]| \WprimeAsph(M)
\end{tikzcd}.\]

It is time to write down the ``mild extra assumption'' on \scrW{}
needed to ensure injectivity of \eqref{eq:81.1}, namely the familiar
enough condition:
\begin{description}
\item[\namedlabel{loc:4}{Loc~4)}]
  The set $\scrW\subset\Fl(\Cat)$ of weak equivalences in \Cat{} is
  ``compatible'' with the functor $\piz:\Cat\to\Sets$, i.e.,
  \begin{equation}
    \label{eq:81.10}
    f\in\scrW \Rightarrow \text{$\piz(f)$ bijective.}\tag{10}
  \end{equation}
\end{description}
(Cf.\ pages\pspage{261} \ref{p:213}--\ref{p:214} for the conditions
\ref{loc:1} to \ref{loc:3}.)

Among all basic localizers satisfying this extra condition, there is
one coarsest of all, which we'll call \scrWz, defined by the
condition
\begin{equation}
  \label{eq:81.11}
  f\in\scrWz \Leftrightarrow \text{$\piz(f)$ bijective.}\tag{11}
\end{equation}
for any map $f$ in \Cat. It is clear too that among all basic
localizers there is a finest, which we'll call $\scrWoo$, and which
can be described as
\begin{equation}
  \label{eq:81.12}
  \scrWoo=\bigcap\scrW, \quad\parbox[t]{0.5\textwidth}{intersection
    of the set of all basic localizers in \Cat.}\tag{12}
\end{equation}
We'll see in\scrcomment{see also \cite{Cisinski2004}} part \ref{ch:V}
of the notes that \scrWoo{} is none else than just the usual notion of
weak equivalence we started with, at the very beginning of our
reflections (cf.\ section \ref{sec:17}). Thus, functoriality of
\eqref{eq:81.1} with respect to \scrW{} implies that \eqref{eq:81.1}
for arbitrary \scrW{} can be described in terms of the particular case
\scrWoo, as the composition
\begin{equation}
  \label{eq:81.13}
  \Cont(M) \hookrightarrow \WooAsph(M) \to \WAsph(M).\tag{13}
\end{equation}
On the other hand, the strongest version of injectivity of
\eqref{eq:81.1}, for different \scrW's, is obtained for \scrWz,
i.e., taking the map
\begin{equation}
  \label{eq:81.14}
  \Cont(M) \hookrightarrow \WzAsph(M).\tag{14}
\end{equation}

This last map seems to me of special significance, because the two
sets it relates correspond to ``\emph{absolute}'' notions (not
depending on the choice of some \scrW), and which moreover are both
``\emph{elementary}'', in the sense that they do not depend on
anything like consideration of non-trivial homotopy or (co)homology
invariants of objects of \Cat. As a matter of fact, the notion of a
contractibility structure corresponds to the algebraic translation of
one of the most elementary and intuitive topological notions, namely
contractibility; whereas the notion of a \scrWz-asphericity structure
can be expressed just as ``elementarily'' in terms of the functor
\begin{equation}
  \label{eq:81.15}
  \piz:\Cat\to\Sets,\tag{15}
\end{equation}
which we may call the ``basic'' connected components-functor, which is
nothing but the algebraic counterpart of the basic intuitive notion of
connected components of a space. We'll denote by
\begin{equation}
  \label{eq:81.16}
  M_0=M_\scrWz\subset\Ob M,\quad W_0\subset\Fl(M)\tag{16}
\end{equation}
the set of \scrWz-aspheric objects and the set of \scrWz-weak
equivalences, associated to a given contractibility structure $M\subc$
on $M$. The objects of $M_0$\pspage{262} merit the name of
\emph{$0$-connected objects} of $M$ (with respect to $M\subc$), and
the arrows of $M$ in $W_0$ merit the name of \emph{$0$-connected
  maps}\footnote{\alsoondate{26.6.} This name for \emph{maps} is
  improper though, as it rather suggests the property of being
  ``universally in $W_0$''.} (with respect to $M\subc$). These two
notions of $0$-connectedness determine each other in an evident way
(valid for any \scrW-aspheric \scrW-asphericity structure, for any
basic localizer \scrW\ldots). Explicitly, this can be expressed by
\begin{equation}
  \label{eq:81.17}
  \begin{cases}
    x\in M_0 \Leftrightarrow \bigl( (x\to e_M)\in W_0 \bigr) \\
    (f:x\to y)\in W_0 \Leftrightarrow
    \bigl(\text{$\piz({M_0}_{/x})\to\piz({M_0}_{/y})$
      bijective}\bigr).
  \end{cases}
  \tag{17}
\end{equation}

The injectivity of \eqref{eq:81.1} can be restated by saying that
\emph{any contractibility structure $M\subc$ on a category $M$} (with
finite products) \emph{can be recovered in terms of the corresponding
  notion of $0$-connected objected objects of $M$}, or, equivalently,
in terms of the corresponding notion of $0$-connected maps in
$M$. Still another way of phrasing this result, is in terms of the
canonical functor
\[ M\to W_0^{-1}M=\HotOf_{(M,M_0),\scrWz} \to \HotOf_\scrWz =
\scrWz^{-1}\Cat \toequ \Sets,\]
which is a functor
\begin{equation}
  \label{eq:81.18}
  \piz:M\to\Sets\tag{18}
\end{equation}
canonically associated to the contractibility structure. We may say
that \emph{the contractibility structure $M\subc$ can be recovered in
  terms of the corresponding functor \piz}, more accurately still, in
terms of the isomorphism class of the latter. Indeed, in terms of this
functor \piz, we recover $M_0$ and $W_0$ by the relations:
\begin{equation}
  \label{eq:81.19}
  \begin{split}
    M_0 &= \set[\big]{x\in\Ob M}{\text{$\piz(x)$ is a one-point set}},\\
    W_0 &= \set[\big]{f\in\Fl(M)}{\text{$\piz(f)$ is bijective}}.
  \end{split}\tag{19}
\end{equation}

This shows that the ``nice'' main fact mentioned at the beginning of
today's notes, namely (essentially) injectivity of \eqref{eq:81.1}
(and, moreover an explicit description of a way how to recover an
$M\subc$ in terms of the corresponding \scrW-asphericity structure) is
not really dependent on relatively sophisticated notions such as
``basic localizers'' and corresponding ``asphericity structures'', but
can be viewed as an ``elementary'' result (namely independent of any
consideration of ``higher'' homotopy or homology invariants, apart
from \piz) about the relationship between contractibility structures,
and corresponding $0$-connectedness notions; the latter may at will be
expressed in terms of either one of the three structural data
\begin{equation}
  \label{eq:81.20}
  M_0, \quad W_0, \quad\text{or}\quad \piz.\tag{20}
\end{equation}

This\pspage{263} relationship has been ``in the air'' since section \ref{sec:50}
(p.\ \ref{p:109}--\ref{p:110}), and I kind of turned around it
consistently up to section \ref{sec:60}, without really getting to the
core. One reason for this ``turning around'' has been, I guess, that I
let myself be distracted, not to say hypnotized, by the ``canonical''
\piz{} functor on a category $M$ (which makes really good sense only
when $M$ is ``totally $0$-connected'' as a category, a condition which
should mean, more or less I suppose, that the $0$-connected objects of
$M$, defined in terms of the mere category structure of $M$, define a
\scrWz-asphericity structure on $M$). Even after realizing (in section
\ref{sec:59}) that one should generalize the description of a
contractibility structure in terms of a ``connected components
functor'' \piz, to the case of a functor $\piz: M\to\Sets$ given
beforehand, and satisfying suitable restrictions (which I did not try
to elucidate), I still was under the impression that the
contractibility structures one could get this way must be of an
extremely special nature. In order to become aware of the fact that
this is by no means so, namely that any contractibility structure
could be obtained from a suitable functor \piz, it would have been
necessary to notice that such a structure $M\subc$ defines in a
natural way a functor \piz. There was indeed the realization that
$M\subc$ should allow to define a notion of weak equivalence (cf.\
page \ref{p:136}), but it wasn't clearly realized by then that at the
same time as a notion of weak equivalence $W$ in $M$, we should also
get a canonical functor
\[M \to W^{-1}M \to \Hot,\]
namely something a lot more precise still than a functor with values
in \Sets{} merely! But rather than push ahead in this direction, I
then decided (p.\ \ref{p:138}) that it would be ``unreasonable'' to go
on still longer pushing of investigation of test functors with values
in \Cat, following the approach which had been on my mind for quite a
while by then, and finally sketched (with the promise of a
corresponding generalization of the former ``key result'' on test
functors) in section \ref{sec:47}. Retrospectively, the whole
``grinding'' part \ref{ch:III} of these notes now looks as a rather
heavy and long-winded digression, prompted by this approach to still a
particular case of test functors (namely with values in \Cat, and
more stringently still, factoring through $\Cat\subas$).

This particular case has been of no use in the present part
\ref{ch:IV} of the notes, developing the really relevant notions in
terms of asphericity structures. \emph{Technically speaking, it now
  appears that most of the reflections of part \ref{ch:III} are
  superseded by part \ref{ch:IV}} -- the main exception being the
development of the\pspage{264} various homotopy notions in sections
\ref{sec:51} and \ref{sec:52}. On the other hand, it is clear that the
main ideas which are coming to fruition in part \ref{ch:IV} all
originated during the awkward grinding process in part \ref{ch:III}!

The ``modelizing story'' so far has turned out as the interplay of
three main sets of notions. One is made up with the
``\emph{test-notions}'', centering around the notion of a
\emph{test-category}, as one giving rise to the most elementary type
of ``\emph{modelizers}'', namely the so-called ``\emph{elementary
  modelizers}'' $(\Ahat,\scrWA)$. This was developed in part
\ref{ch:II} (while part \ref{ch:I} was concerned with the initial
motivation of the reflections, namely stacks, forgotten for the time
being!). The second set of notions concerns the so-called
``\emph{homotopy notions}'', developed at some length in sections
\ref{sec:51} to \ref{sec:55}, summarized in the diagram
\eqref{eq:81.2} above. They constitute the main technical content of
part \ref{ch:III} of the notes, with however one major shortcoming:
the relationship between these notions, and $0$-connectedness notions
\eqref{eq:81.20}, was only partially understood in part \ref{ch:III},
namely as a one-way relationship merely, associating to suitable
$0$-connectedness notions in a category $M$, a corresponding homotopy
structure in $M$. The third set of notions may be called
``\emph{asphericity notions}'', they center around the notions of
\emph{aspheric objects} and \emph{aspheric maps} (in a category
endowed with a so-called \emph{asphericity structure}), and more
specifically around the notion of an aspheric map in \Cat, whose
formal properties turn out to be the key for the development of a
theory of asphericity structures. The first and the third set of
notions (namely test notions and asphericity notions) depend on the
choice of a ``\emph{basic localizer}'' \scrW{} in \Cat, whereas the
second set, namely homotopy notions, is ``absolute'', i.e., does not
depend on any such choice, nor on any knowledge of homotopy or
(co)homology invariants.

Whereas the test notions are essentially concerned with modelizers,
namely getting descriptions of the category of homotopy types \Hot{}
in terms of elementary modelizers \Ahat{} (as being
$\scrWA^{-1}\Ahat$), it appears that the homotopy notions, as well as
the asphericity notions, are independent of any modelizing notions and
assumptions. In a deductive presentation of the theory, the test
notions would come last, whereas they came first in these notes -- as
an illustration of the general fact that the deductive approach will
present things roughly in opposite order in which they have been
discovered! The test notions, as an outcome of the attempt to get a
picture of modelizers, have kept acting as a constant guideline in the
whole reflection, even though technically speaking they are
``irrelevant'' for the development of\pspage{265} the main properties
of homotopy and asphericity notions and their interplay.

By the end of part \ref{ch:IV}, there has been some floating in my
mind as to whether which among the two structures, namely
contractibility structures or asphericity structures, should be
considered as ``the'' key structure for an understanding of the
modelizing story. There was a (justified) feeling, expressed first at
the beginning of section \ref{sec:67}, that in some sense, asphericity
structures were ``more general'' than contractibility structures,
which caused me for a while to view them as the more ``basic'' ones. I
would be more tempted at present to hold the opposite view. \emph{The
  notion of a contractibility structure now appears as a kind of hinge
  between the two main sets of notions besides the test notions,
  namely between homotopy notions and asphericity notions.} On the one
hand, as displayed in diagram \eqref{eq:81.2}, the notion of a
contractibility structure appears as the most stringent one among the
four main types of homotopy structures. On the other hand, by
\eqref{eq:81.1} it can be equally viewed as being a special case of an
asphericity structure, and as such it gives rise to (and can be
expressed by) either one of the following four asphericity-flavored
data on a category $M$ (for any given basic localizer \scrW{}
satisfying \ref{loc:4}):\scrcomment{also listed, but scratched out in
  the typescript, was $\piz:M\to\Sets$}
\begin{equation}
  \label{eq:81.21}
  \begin{cases}
    M_\scrW\subset\Ob M,  &\text{$\scrW_{M\subc}$ or simply
      $\scrW_M\subset\Fl(M)$} \\
    \text{$\varphi_M^\scrW$ or $\varphi_M:M\to\HotOf_\scrW$.} &
  \end{cases}\tag{21}
\end{equation}
Restricting to the case when \scrW{} is either \scrWz{} or \scrWoo, we
get the corresponding structures on $M$, namely the three structures
$M_0$, $W_0$, \piz{} of \eqref{eq:81.20} plus the three extra
structures:
\begin{equation}
  \label{eq:81.22}
  \begin{gathered}
  M_\oo = M_\scrWoo \subset\Ob M, \quad
  W_\oo = (\scrWoo)_M\subset\Fl(M), \\
  \varphi_M : M\to\Hot = \scrWoo^{-1}\Cat,
  \end{gathered}\tag{22}
\end{equation}
which can be referred to as \emph{\oo-connected objects} of $M$,
\emph{\oo-connected arrows} of $M$, and the ``\emph{canonical
  functor}'' from $M$ to \Hot. Putting together \eqref{eq:81.4},
\eqref{eq:81.20}, \eqref{eq:81.22}, we see that a contractibility
structure $M\subc$ gives rise to ten different structures (including
$M\subc$ itself) and is determined by each one of these,\pspage{266}
\begin{equation}
  \label{eq:81.23}
  \left\{\begin{aligned}
    &M\subc\subset\Ob M, &&J\subc\subset\Int(M), \!\!
    &&W\subc\subset\Fl(M), &&R\subc\subset\Fl(M)\times\Fl(M) \\
    &M_0\subset\Ob M,    &&
    &&W_0\subset\Fl(M),    &&\piz:M\to\Sets \\
    &M_\oo\subset\Ob M,\!\!\!  &&
    &&W_\oo\subset\Fl(M),\!\!\!  &&\varphi_M:M\to\Hot,
  \end{aligned}\right.\tag{23}
\end{equation}
namely: \emph{contractible objects}, \emph{homotopy intervals},
\emph{homotopisms}, \emph{homotopy relation for maps},
\emph{$0$-connected objects}, \emph{$0$-connected maps},\footnote{This
  name is inadequate, cf.\ note p.\ \ref{p:262}.} the \emph{connected
  components functor}, \emph{\oo-connected} (or ``aspheric'', more
accurately \scrWoo-aspheric) \emph{objects},
\emph{\scrWoo-equivalences} (or simply ``weak equivalences''), and the
(would-be ``modelizing'') \emph{canonical functor from $M$ to homotopy
  types}. Moreover, it should be remembered that, just as the
structures $J\subc$, $W\subc$, $R\subc$ are by no means unrestricted
(the fact that they stem from a contractibility structure being a
substantial restriction), the two asphericity structures (namely, the
$0$-asphericity structure $M_0$ and the \oo-asphericity structure
$M_\oo$) are subject to extra conditions which will be written down
below, implying among others that they are totally aspheric (hence
\piz{} and $\varphi_M$ are compatible with finite products).

We may complement the ten structure data \eqref{eq:81.23} above, by
the following two
\begin{equation}
  \label{eq:81.24}
  R_0, R_\oo \subset \Fl(M)\times\Fl(M),\tag{24}
\end{equation}
defined by
\[(f,g)\in R_0 \Leftrightarrow \bigl(\piz(f)=\piz(g)\bigr), \quad
(f,g)\in R_\oo \Leftrightarrow
\bigl(\varphi_M(f)=\varphi_M(g)\bigr).\]
More generally, for any basic localizer \scrW, we may define
\begin{equation}
  \label{eq:81.25}
  R_\scrW\subset\Fl(M)\times\Fl(M),\quad
  (f,g)\in R_\scrW \Leftrightarrow
  \bigl(\varphi_M^\scrW(f)=\varphi_M^\scrW(g)\bigr).
  \tag{25}
\end{equation}
It is not clear however that $M\subc$, or equivalently $M_\scrW$ or
$\scrW_M$, can be recovered in terms of the equivalence relation
$R_\scrW$ among maps of $M$.

Among the three series of structures appearing in \eqref{eq:81.23}, we
have the tautological relations
\begin{equation}
  \label{eq:81.26}
  \begin{cases}
    M\subc \subset M_\oo \subset M_0 & \\
    W\subc \subset W_\oo \subset W_0 & \\
    R\subc \subset R_\oo \subset R_0 & \\
    \begin{tikzcd}[cramped]
      M \ar[r,"\varphi_M"] \ar[rr,bend right, "\pi_0"] &
      \Hot \ar[r,"\overline{\piz}"] & \Sets
    \end{tikzcd} & \text{(commutative).}
  \end{cases}\tag{26}
\end{equation}
Finally, we may also display some of the main functors defined in
terms of a given contractibility structure $M\subc$:
\begin{equation}
  \label{eq:81.27}
  \begin{tabular}{@{}c@{}}
  \begin{tikzcd}[baseline=(O.base),sep=tiny]
    |[alias=O]| M\ar[r] &
    \overline M \ar[r]\ar[d,equal] &
    \HotOf_M \ar[r]\ar[d,equal] &
    \Hot \ar[r]\ar[d,equal] &
    \HotOf_\scrW \ar[r] &
    \Sets\ar[d,phantom,sloped,"\equ" description] \\
    & W\subc^{-1}M & W_\oo^{-1}M & \HotOf_\scrWoo & & \HotOf_\scrWz
  \end{tikzcd}.
  \end{tabular}\tag{27}
\end{equation}

\bigbreak

\presectionfill\ondate{26.6.}\pspage{267}\par

\hangsection[Proof of injectivity of
$\alpha: \Contr(M)\hookrightarrow\WAsph(M)$. \dots]{Proof of
  injectivity of \texorpdfstring{$\alpha:
    \Contr(M)\hookrightarrow\WAsph(M)$}{alpha:Contr(M)->W-Asph(M)}.
  Application to \texorpdfstring{$\bHom$}{Hom} objects and to products
  of aspheric functors \texorpdfstring{$A\to
    M$}{A->M}.}\label{sec:82}%
Yesterday I have been busy mainly with the readjustment of the overall
perspective on the main notions developed so far, which has sprung
from the new fact stated at the beginning of yesterday's notes: namely
that an arbitrary contractibility structure $M\subc$ (on a category
$M$ stable under finite products) can be recovered in terms of the
associated \scrW-asphericity structure, where \scrW{} is any basic
localizer satisfying (besides the condition \ref{loc:1} to \ref{loc:3}
of p.\ \ref{p:213}--\ref{p:214}) the extra assumption \ref{loc:4} of
p.\ \ref{p:260}. It seems about time now to enter into a little more
technical specifications along the same lines -- and to start with,
give a proof of the ``new fact''! Let's state it again in full:
\begin{theoremnum}\label{thm:82.1}
  Let $M$ be a \scrU-category stable under finite products,
  $M\subc\subset\Ob M$ a contractibility structure on $M$, admitting a
  small full subcategory $C$ which generates the structure. Let
  moreover \scrW{} be any basic localizer satisfying
  \textup{\ref{loc:4}} \textup(compatibility with the \piz-functor
  $\Cat\to\Sets$\textup), and let $M_\scrW\subset\Ob M$ the
  \scrW-asphericity structure generated by $M\subc$\kern1pt,
  $W=\scrW_{M\subc}$ the corresponding set of ``\scrW-equivalences''
  or ``weak equivalences'' in $\Fl(M)$. Consider the homotopy
  structure $h_W$ associated to $W$, i.e.\ \textup(cf.\ section
  \ref{sec:54}\textup), the homotopy structure associated to the
  homotopy interval structure $J$ generated by the set $J_0$ of all
  intervals
  \[\bI=(I,\delta_0,\delta_1)\]
  in $M$ such that
  \[I\in M_\scrW.\]
  Then $h_W$ is the homotopy structure on $M$ associated to the
  contractibility structure $M\subc$\kern1pt, and hence $M\subc$ can be
  described in terms of $M_\scrW$ \textup(or of
  $\scrW_{M\subc}=W$\textup) as the set of objects which are
  contractible for $h_W$, i.e., such that the map $x\to e_M$ is an $h_W$-homotopism.
\end{theoremnum}
\begin{proof}
  Let $M\subc'$ be the set of $h_W$-contractible objects of $M$,
  clearly we have
  \[ M\subc \subset M\subc'\subset M_\scrW \eqdef M\suba.\]
  The theorem amounts to saying that $M\subc$ generates the homotopy
  interval structure $J$ (by which we mean that the set of intervals
  of $M$ made up with objects of $M\subc$ generates the structure
  $J$). Indeed, because of $M\subc\subset M\subc'$, this will imply
  that $J$ is associated to a contractibility structure, namely to
  $M\subc'$. But for an object $x$ of $M$ to be in $M\subc'$, i.e., to
  be contractible for the structure $J$, amounts to be contractible
  for $M\subc$, and hence by saturation of $M\subc$, to be in
  $M\subc$, hence $M\subc'=M\subc$ -- which yields what
  we\pspage{268} want. By the description of $J$ in terms of $J_0$, we
  are now reduced to proving the following
\end{proof}
\begin{lemmanum}\label{lem:82.1}
  Let $I$ be an object of $M\suba=M_\scrW$. Then any two sections of
  $I$ \textup(over $e_M$\textup) are $M\subc$-homotopic.
\end{lemmanum}

Let, as in the previous section, $\scrWz\supset\scrW$ be the largest
of all basic localizers satisfying \ref{loc:4}, i.e.,
\[\scrWz=\set[\big]{f\in\Fl(\Cat)}{\text{$\piz(f)$ a bijection}},\]
therefore, we have
\[ M_\scrW \subset M_\scrWz \eqdef M_0,\]
and we are reduced to proving the lemma for \scrWz{} instead of \scrW,
i.e., for $M_0$ instead of $M_\scrW$. We'll use the small full
subcategory $C$ of $M\subc$ generating the contractibility structure
$M\subc$, we may assume that $C$ is stable under finite
products. Hence \Chat{} is totally \scrWz-connected, i.e., totally
$0$-connected. Moreover, as $C\subset M\subc$, and $e_M$ is in $C$, it
follows that every object of $C$ has a section (over $e_M=e_C$) --
which implies that every non-empty object of \Chat{} has a section,
i.e., \Chat{} is ``strictly totally $0$-connected'' (cf.\ p.\
\ref{p:144} and \ref{p:149}). Note that (by prop.\ \ref{prop:79.2}
\ref{it:79.2.b} of p.\ \ref{p:248}) we have
\[ M_0 = \set[\big]{x\in\Ob M}{\text{$i^*(x)$ is $0$-connected in
    \Chat}},\]
where $i$ is the inclusion functor:
\[i:C\to M.\]
We have to prove that any two sections $\delta_0,\delta_1$ of an
object $I$ of $M_0$ are $M\subc$-homotopic, or what amounts to the
same, $C$-homotopic. This translates readily into the statement that
$i^*(\delta_0)$ and $i^*(\delta_1)$ are $C$-homotopic in \Chat. Thus,
we are reduced to the following lemma (in the case of the topos
\Chat):
\begin{lemmanum}\label{lem:82.2}
  Let \scrC{} be a totally $0$-connected topos such that any non-empty
  object of $C$ has a section, and let $C$ be a small full generating
  subcategory, whose elements are $0$-connected. Then for any
  $0$-connected object $I$ of \scrC, and any two sections
  $\delta_0,\delta_1$ of $I$, these are $C$-homotopic, i.e., they can
  be joined by a finite chain of sections, any two consecutive among
  which can be obtained as the images of two sections $s_i,t_i$ of an
  objects $x_i$ of $C$, by means of a map $h_i:x_i\to I$.
\end{lemmanum}

This lemma is essentially a restatement (cleaned from extraneous
assumptions due to an awkward conceptual background) of the
proposition of page \ref{p:149} (section \ref{sec:60}), and the proof
will be left to the reader.

\subsection*{Application to relation between contractibility and
  objects \texorpdfstring{$\bHom(X,Y)$}{Hom(X,Y)}.}
\label{subsec:82.hom}

I\pspage{269} would like to review here a few things along these lines, which were
somewhat scattered in the notes before (section \ref{sec:51}
\ref{subsec:51.E} p.\ \ref{p:121} and section \ref{sec:57} p.\
\ref{p:143} notably), and at times came out awkwardly because of
inadequate conceptual background. We assume $M$ endowed with a
contractibility structure $M\subc$, and use the notations of the
previous section, especially concerning the subsets of $\Ob M$ and of
$\Fl(M)$
\begin{equation}
  \label{eq:82.1}
  M\subc \subset M_\oo \subset M_\scrW \subset M_0, \quad
  W\subc \subset W_\oo \subset \scrW_M \subset W_0. \tag{1}
\end{equation}
Let $X$ be an object of $M$, such that the object
\begin{equation}
  \label{eq:82.2}
  I = \bHom(X,X)\tag{2}
\end{equation}
exists in $M$ (NB\enspace a priori it is an object in $M\uphat$, we
assume it to be representable). We assume that $X$ is endowed with a
section $c$, hence a section $\delta_1$ of $I$, corresponding to the
constant endomorphism of $X$ with value $c$. We'll denote by
$\delta_0$ the section of $I$ corresponding to the identity of
$X$. Thus,
\begin{equation}
  \label{eq:82.3}
  \bI=(I,\delta_0,\delta_1)\tag{3}
\end{equation}
becomes an interval of $M$, and the composition law of $I=\bHom(X,X)$
turns it into a \emph{multiplicative interval}, admitting respectively
$\delta_0$ and $\delta_1$ as left unit and left zero
element.\scrcomment{in the typescript we have $\pi_0$ and $\pi_1$ as
  the elements; I assume this is a typo?} Then the following
conditions are equivalent:
\begin{description}
\item[\namedlabel{it:82.hom.i}{(i)}]
  $X$ is contractible.
\item[\namedlabel{it:82.hom.ii}{(ii)}]
  $\bHom(X,X)$ is contractible, i.e., $I\in M\subc$.
\item[\namedlabel{it:82.hom.iiprime}{(ii')}]
  $\bHom(X,X)$ is $0$-connected, i.e., $I\in M_0$.
\item[\namedlabel{it:82.hom.iidblprime}{(ii'')}]
  The two sections $\delta_0,\delta_1$ of $I=\bHom(X,X)$ are homotopic.
\item[\namedlabel{it:82.hom.iii}{(iii)}]
  For any object $Y$ in $M$ such that $\bHom(Y,X)$ exists in $M$,
  $\bHom(Y,X)$ is contractible.
\item[\namedlabel{it:82.hom.iiiprime}{(iii')}]
  For any $Y$ as in \ref{it:82.hom.iii}, $\bHom(Y,X)$ is $0$-connected.
\item[\namedlabel{it:82.hom.iv}{(iv)}]
  For any $Y$ in $M$ such that $\bHom(X,Y)$ exists in $M$, the
  canonical map
  \begin{equation}
    \label{eq:82.4}
    Y \to \bHom(X,Y)\tag{4}
  \end{equation}
  is a homotopism (more accurately still, $Y$ as a subobject of
  $\bHom(X,Y)$ is a deformation retract). Moreover, $X$ is $0$-connected.
\item[\namedlabel{it:82.hom.ivprime}{(iv')}]
  For any $Y$ in $M$ as in \ref{it:82.hom.iv}, the map \eqref{eq:82.4}
  induces a bijection
  \[\piz(Y)\to\piz(\bHom(X,Y)),\]
  moreover $X$ is $0$-connected.
\end{description}
NB\enspace Of course, the $0$-connectedness notion and the functor
\piz{} used\pspage{270} here are those associated to the given
structure $M\subc$. The case dealt with in section \ref{sec:57} (p.\
\ref{p:143}) is essentially (it seems) the one when these notions are
the ones canonically defined in terms of the category structure of $M$
alone. In view of the inclusions \eqref{eq:82.1}, we could throw in a
handful more obviously equivalent conditions, using $M_\oo$, $M_\scrW$
or $W_\oo$, $\scrW_M$, instead of $M\subc$, etc.\ but it seems this
would confuse the picture rather than complete it.

Next thing is to look at the canonical map
\begin{equation}
  \label{eq:82.5}
  \overline\Gamma(X) \to \piz(X)\tag{5}
\end{equation}
defined for any object $X$ of $M$, where $\Gamma(X)$ denotes the set
of sections of $X$, and $\overline\Gamma(X)$ denotes the set of
corresponding homotopy classes of sections. (This map was considered
in a slightly different case on page.\ \ref{p:144}) The map
\eqref{eq:82.5} can be viewed as the particular case of
\[ \oHom(Y,X) \to \Hom_{\mathrm{Sets}}(\piz(Y),\piz(X)),\]
obtained when $Y=e_M$, hence $\piz(Y)=$ one-point set. By the standard
description
\begin{equation}
  \label{eq:82.6}
  \piz(X)\simeq \bpiz(C_{/X}) \quad(=\piz(i^*(X)))\tag{6}
\end{equation}
we immediately get that \eqref{eq:82.5} is always surjective, as any
element in $\piz(X)$ is induced by an element $x\to X$ of $C_{/X}$,
hence by a section
\[ e_M\to x\to X,\]
where $e_M\to x$ is a section of the (contractible) object $x$ in
$C$. But \eqref{eq:82.5} is in fact bijective. To see this, let's note
that the map \eqref{eq:82.5} is isomorphic to the corresponding map in
\Sets, with $M$ replaced by \Chat{} and $X$ by $i^*(X)$, $C$ remaining
the same, as follows from \eqref{eq:82.6} and the formula
\[\Gamma(X) \simeq \overline\Gamma(i^*(X)).\]
The latter is a particular case of the
\begin{lemmanum}\label{lem:82.3}
  Let $Y,X$ be two objects of $M$, with $Y$ in $C$. Then the natural
  map
  \begin{equation}
    \label{eq:82.7}
    \oHom(Y,X) \to \oHom(i^*(Y),i^*(X))\tag{7}
  \end{equation}
  between sets of homotopy classes of maps \textup(in $M$ and in
  \Chat{} respectively, the latter being endowed with the
  contractibility structure generated by $C$\textup) is bijective.
\end{lemmanum}

The verification is tautological, due to the fact that for any object
$I$ of $C$ (which is going to play the role of a homotopy interval for
deformations), we got
\[\Hom_M(I\times Y,X) \tosim \Hom_\Chat(I\times Y, i^*(X)).\]

Now,\pspage{271} the map \eqref{eq:82.5} in case of a strictly totally
$0$-connected category (here \Chat) has been dealt with in section
\ref{sec:58} (p.\ \ref{p:114}), it follows easily that the map is
bijective -- hence the
\begin{propositionnum}\label{prop:82.1}
  Let $(M,M\subc)$ be any contractibility structure,
  \[\piz:M\to\Sets,\quad
  \piz(X)\simeq\bpiz({M\subc}_{/X})\]
  the corresponding ``connected components functor'', $X$ any object
  of $M$. Then the canonical map \eqref{eq:82.5} above is bijective.
\end{propositionnum}
\begin{corollary}
  Let $X,Y$ be two objects of $M$, such that $\bHom(X,Y)$ exists in
  $M$. Then the natural map
  \begin{equation}
    \label{eq:82.8}
    \oHom(X,Y) \to \piz(\bHom(X,Y))\tag{8}
  \end{equation}
  is bijective.
\end{corollary}

\bigbreak

To finish these generalities on $\bHom$'s, I would like to generalize
to the context of contractibility structures the results of prop.\
\ref{prop:74.4} p.\ \ref{p:228} (section \ref{sec:74}) on products of
aspheric functors. We start within a context of asphericity
structures:
\begin{propositionnum}\label{prop:82.2}
  Let $(M,M\suba)$ be a \scrW-asphericity structure \textup(for a basic
  localizer \scrW\textup), $C$ a full subcategory of $M$ generating
  this structure, i.e., such that $C\to M$ is $M\suba$-\scrW-aspheric,
  $A$ a small category,
  \[ i : A\to M\]
  a $M\suba$-\scrW-aspheric functor. We assume $M$ stable under binary
  products.
  \begin{enumerate}[label=\alph*),font=\normalfont]
  \item\label{it:82.prop2.a}
    Let $b_0$ in $M$ be such that for any $y$ in $C$, $\bHom(b_0,y)$
    exists in $M$, and for any $x$ in $M\suba$, $x\times b_0$ is in
    $M\suba$. Let $i_{b_0}:A\to M$ be the constant functor with value
    $b_0$, and consider the product functor
    \[ i \times i_{b_0}: a\mapsto i(a)\times b_0 : A\to M.\]
    This functor is $M\suba$-\scrW-aspheric if{f} for any $y$ in $C$,
    $\bHom(b_0,y)$ is in $M\suba$ \textup(a condition which does not
    depend on $i$ nor even on $A$\textup).
  \item\label{it:82.prop2.b}
    Let $B$ be a full subcategory of $M$, such that for any $b_0$ in
    $B$, $x$ in $M\suba$, and $y$ in $C$, we have $x\times b_0\in
    M\suba$ and $\bHom(b_0,y)$ exists in $M$, and is in
    $M\suba$. Assume $A$ totally \scrW-aspheric, let $i':A\to M$ be
    any functor factoring through $B$, then the product functor
    \[ i \times i' : a\mapsto i(a)\times i'(a) : A\to M\]
    is $M\suba$-\scrW-aspheric.
  \end{enumerate}
\end{propositionnum}

The proof is word by word the same as for the analogous statements p.\
\ref{p:228}--\ref{p:229}, and therefore left to the reader. (The
analogous statement to cor.\ \ref{cor:74.4.1} p.\ \ref{p:229} is
equally valid.)

The\pspage{272} case I've in mind is when $M\suba$ is generated by a
contractibility structure $M\subc$, and $C=B=M\subc$. The condition in
\ref{it:82.prop2.b} boils down to existence of $\bHom(b_0,y)$ in $M$,
when $b_0$ and $y$ are both contractible objects of $M$ -- as a matter
of fact, in all cases I'm having in mind, the $\bHom$ exists even
without any contractibility assumption. Thus we get:
\begin{corollary}
  Let $(M,M\subc)$ be a contractibility structure, such that for any
  two contractible objects $x,y$, $\bHom(x,y)$ exists in $M$ \textup(a
  condition which presumably is even
  superfluous\ldots\textup). Consider the \scrW-asphericity structure
  $M_\scrW$ on $M$ generated by $M\subc$. Then the product of a
  $M_\scrW$-\scrW-aspheric functor $i:A\to M$ with a functor $i':A\to
  M$ factoring through $M\subc$ is again $M_\scrW$-\scrW-aspheric. In
  particular, the category of all aspheric functors from $A$ to $M$
  which factor through $M\subc$ is stable under binary products.
\end{corollary}

From this one can deduce as on p.\ \ref{p:230} that if the latter
category is non-empty, i.e, if there does exist an aspheric functor
$A\to M$ through $M\subc$, then one can define a canonical functor
\begin{equation}
  \label{eq:82.9}
  \HotOf_{M,\scrW} \to \HotOf_{A,\scrW},\tag{9}
\end{equation}
defined up to unique isomorphism, via a transitive system of
isomorphisms between the functors deduced from aspheric functors $A\to
M$ factoring through $M\subc$.

\hangsection{Tautologies on \texorpdfstring{$\Imm\alpha$}{Im alpha}, and
  related questions.}\label{sec:83}%
After stating and proving theorem \ref{thm:82.1} of the previous
section, I forgot to give an answer to the most natural question
arising from it -- namely how to characterize those \scrW-asphericity
structures on a category $M$ stable under finite products, which can
be generated by a contractibility structure. The reason for this is
surely that I have nothing better to offer than a tautology: let
\[M\suba \subset \Ob M\]
be the given \scrW-asphericity structure, we assume beforehand that
this structure is totally aspheric (which is a necessary condition for
$M\suba$ to come from a contractibility structure $M\subc$). If
$W\suba$ is the corresponding set of \scrW-equivalences, it follows
that the homotopy structure $h_{W\suba}$ is also the one defined by
the homotopy interval structure $J$ generated by the set $J_0$ of
intervals of $M$ made up with objects of $M\suba$. Let\pspage{273}
\[M\subc \subset \Ob M\]
be the corresponding set of contractible objects (which may not be a
contractibility structure on $M$). The conditions now are the
following:
\begin{enumerate}[label=\alph*)]
\item\label{it:83.a}
  $M\subc$ generates the homotopy interval structure $J$, or
  equivalently, any two sections of an object of $M\suba$ are $M\subc$-homotopic.
\end{enumerate}
This condition is clearly equivalent to saying that $J$ is indeed
generated by a contractibility structure, and the latter is
necessarily $M\subc$.
\begin{enumerate}[label=\alph*),start=2]
\item\label{it:83.b}
  $M\subc$ generates the \scrW-asphericity structure $M\suba$, i.e.,
  there exists a small full subcategory $B$ of $M\subc$ such that the
  inclusion functor $i:B\to M$ be $M\suba$-\scrW-aspheric, i.e., for
  any $x$ in $M\suba$, $i^*(x)$ in \Bhat{} is \scrW-aspheric (i.e.,
  $B_{/x}$ is \scrW-aspheric in \Cat).
\end{enumerate}

One would like too a n.s.\ condition on a functor
\[\varphi:M\to\HotOf_\scrW,\]
for this functor to be isomorphic to the one associated to a
contractibility structure on $M$ -- with special interest in the case
when $\scrW=\scrWz$, i.e., when the functor reduces to a functor
\[\piz:M\to\Sets.\]
Here again, I have nothing to offer except a tautological statement,
which isn't worth the trouble writing down. Nor do I have at present a
compelling feeling that there should exist such a characterization,
under suitable exactness assumptions on $M$ say, and possibly assuming
too that $M$ is stable under the operation $\bHom$. Here also arises
the question \emph{whether a functor $\varphi:M\to\HotOf_\scrW$
  stemming from a contractibility structure $M\subc$ on $M$ can have
  non-trivial automorphisms} -- a question closely connected to the
``inspiring assumption'' of section \ref{sec:28}.

\bigbreak

\presectionfill\ondate{29.6.}\pspage{274}\par

\hangsection[A silly (provisional) answer to the ``silly question'' --
and \dots]{A silly
  \texorpdfstring{\textup(\kernifitalic{1pt}}{(}provisional\texorpdfstring{\textup)}{)}
  answer to the ``silly question'' -- and the new perplexity
  \texorpdfstring{$f_!(M\subas)\subset M'\subas$}{f!(Mas) in
    M'as}?}\label{sec:84}%
Two days have passed without writing any notes. Much of the time I
spent on writing mathematical letters -- one pretty long one to Gerd
Faltings,\scrcomment{\textcite{Grothendieck1983b}} who (on my request)
had sent me preprints of his recent work, notably on the Tate
conjectures for abelian varieties and on the Mordell conjecture, and
had expressed interest hearing about some ideas and conjectures on
``anabelian algebraic geometry''. I had been impressed, from a glance
upon the last of his manuscripts, to see three key conjectures proved
in about forty pages, while they were being considered as quite out of
reach by the people supposed to know. Some ``anabelian'' conjectures
of mine are closely related to the Tate and Mordell conjectures just
proved by Faltings -- Deligne had pointed out to me about two years
ago that a certain fixed-point conjecture (which I like to view as the
basic conjecture at present in the anabelian program) implied
Mordell's, so why loose one's time on it! I have the feeling Faltings
is the kind of chap who may become interested in things which are
supposed to be too far off to be worth looking at, that's why I had
written him a few words, under the moment's inspiration. -- Another
letter, not quite as long, was an answer to a very long and patient
letter of Tim Porter, telling me about a number of things which have
been done by homotopy people, and which I was of course wholly
ignorant of! His letter has been the first echo I got from someone who
read part of the notes on ``Pursuing stacks'', and I was glad he could
make some sense of what he read so far, and conversely -- that not all
he was telling me was going wholly ``above my head''!

Apart from this and the (not unpleasant!) daily routine, I spent a
fair bunch of hours on scratchwork, centering around trying to figure
out the right notion of morphism for asphericity structures. It was a
surprise that the notion should be so reticent for revealing itself --
as a matter of fact, I am not quite sure yet if I got the right
notion, in sufficient generality I mean. Time will tell -- for the
time being, the notion of morphism I have to offer, while maybe too
restrictive, looks really seducing, because it parallels so perfectly
the formalism of aspheric functors $i:A\to M$ and the corresponding
(would-be modelizing) functor $i^*:M\to\Ahat$. Again, it has been a
surprise, right now, that after uncounted hours of unconvincing
efforts today and yesterday, and almost wholly unrelated to these,
this pretty set-up would come out within ten minutes reflection!

\begin{proposition}
  Let\pspage{275} $(M,M\suba)$, $(M',M'\suba)$ be two
  \scrW-asphericity structures \textup(\scrW{} a basic
  localizer\textup), $B\subset M\suba$ a full subcategory of $M$
  generating the asphericity structure $M\suba$, and
  \[f^*:M\to M'\]
  a functor, admitting a left adjoint $f_!$ \textup(hence $f^*$
  commutes to inverse limits\textup). We consider the following six
  conditions, paralleling those of prop.\ \ref{prop:75.4} of section
  \ref{sec:75} \textup(p.\ \ref{p:236}\textup):
  \begin{description}
  \item[\namedlabel{it:84.i}{(i)}]
    $M\suba = (f^*)^{-1}(M'\suba)$,
  \item[\namedlabel{it:84.iprime}{(i')}]
    $M\suba \subset (f^*)^{-1}(M'\suba)$, i.e., $f^*(M\suba) \subset M'\suba$,
  \item[\namedlabel{it:84.idblprime}{(i'')}]
    $B\suba \subset (f^*)^{-1}(M'\suba)$, i.e., $f^*(B) \subset M'\suba$,
  \item[\namedlabel{it:84.ii}{(ii)}]
    $W\suba = (f^*)^{-1}(W'\suba)$,
  \item[\namedlabel{it:84.iiprime}{(ii')}] 
    $W\suba \subset (f^*)^{-1}(W'\suba)$, i.e., $f^*(W\suba)\subset W'\suba$,
  \item[\namedlabel{it:84.iidblprime}{(ii'')}] 
    $\Fl(B)\suba \subset (f^*)^{-1}(W'\suba)$, i.e., $f^*(\Fl(B))\subset W'\suba$,
  \end{description}
  where $W\suba$, $W'\suba$ are the sets of \scrW-equivalences in $M$
  and $M'$ respectively. We'll make the extra assumption \textup(I
  nearly forgot, sorry!\textup):
  \begin{description}
  \item[\namedlabel{cond:84.Awk}{(Awk)}]
    There exists a $M'\suba$-\scrW-aspheric functor
    \[i':A\to M'\]
    \textup($A$ a small category\textup), such that $i=f_!i':A\to M$
    factors through $B$.
  \end{description}
  Under these conditions and with these notations, the following
  holds: The conditions \textup{\ref{it:84.i}, \ref{it:84.iprime},
    \ref{it:84.idblprime}} are equivalent and imply the three others,
  which satisfy the tautological implications \textup{\ref{it:84.ii}}
  $\Rightarrow$ \textup{\ref{it:84.iiprime}} $\Rightarrow$
  \textup{\ref{it:84.iidblprime}}. If $M$ and $M'$ admit final objects
  and their asphericity structures are aspheric, then the fiver first
  conditions \textup{\ref{it:84.i}} to \textup{\ref{it:84.iiprime}}
  are equivalent, and if moreover $e_M\in \Ob B$, all six are
  equivalent.
\end{proposition}
\noindent\emph{Proof:} reduction to loc.\ sit., using the commutative diagram
\begin{equation}
  \label{eq:84.1}
  \begin{tikzcd}[column sep=tiny]
    M\ar[rr,"f^*"]\ar[dr,"i^*"'] & & M'\ar[dl,"{i'}^*"] \\
    & \Ahat &
  \end{tikzcd}\tag{1}
\end{equation}
(up to isomorphism), and the relations
\[M'\suba = ({i'}^*)^{-1}(\Ahata), \quad
W'\suba = ({i'}^*)^{-1}(\scrWA).\]
This proof shows moreover that the conditions \ref{it:84.i},
\ref{it:84.iprime}, \ref{it:84.idblprime} are equivalent each to $i$
being $M\suba$-\scrW-aspheric. Thus, if these conditions are satisfied
(let's say then that $f^*$ is a
\emph{morphism} for the asphericity structures), we may use $i^*$ and
${i'}^*$ for describing the canonical functors from\pspage{276} the
localizations of $M,M'$ to $\HotOf_\scrW$, and therefore from
\eqref{eq:84.1} get a \emph{commutative} diagram (up to canonical
isomorphism)
\begin{equation}
  \label{eq:84.2}
  \begin{tabular}{@{}c@{}}
    \begin{tikzcd}[baseline=(O.base),column sep=tiny]
      \HotOf_{M,\scrW}\ar[rr,"\overline{f^*}"]
      \ar[dr] & & \HotOf_{M',\scrW}\ar[dl] \\
      & |[alias=O]| \HotOf_\scrW &
    \end{tikzcd}.
  \end{tabular}\tag{2}
\end{equation}
From this follows:
\begin{corollary}
  Let $f^*$ be a morphism for the given asphericity structure, and
  assume these to be \emph{modelizing}. Then the induced functor for
  the localizations
  \[\overline{f^*}: W\suba^{-1}M = \HotOf_{M,\scrW} \to
  {W'\suba}^{-1}M' = \HotOf_{M',\scrW}\]
  is an equivalence of categories.
\end{corollary}

This is the answer, at last, of the ``silly question'' of section
\ref{sec:45} (p.\ \ref{p:95})!

We have still to comment though on the restrictive conditions we had
to make, for getting the equivalences stated in the proposition, and
the result stated in the corollary. These conditions are twofold:
a)\enspace existence of a left adjoint $f_!$, which presumably, in
most circumstances we are going to meet, will be equivalent with $f^*$
commuting to inverse limits. It's a pretty restrictive condition, but
of a rather natural kind, often met with in the modelizing situations;
b)\enspace this is the ``awkward'' condition \ref{cond:84.Awk}, which
can be equally stated as follows: there exists a small full
subcategory $B'$ of $M'$ (NB\enspace we may take the full subcategory
defined by $i'(\Ob A)$), \emph{generating} the asphericity structure
$M'\suba$, and such moreover that
\[f_!(B')\subset M\suba,\]
(so that we can choose $C$, a full subcategory of $M$ generating the
asphericity structure, such that $C$ contains $f_!(B')$. Another way
of phrasing this condition on $f^*$ (preliminary to the choice of $C$)
is that the full subcategory of $M'$
\begin{equation}
  \label{eq:84.3}
  \text{$(f_!)^{-1}(M\suba) \sand M'\suba$ generates the asphericity
    structure of $M'$.}\tag{3}
\end{equation}

If this condition, plus the condition \ref{it:84.i} say, which we
express jointly by saying that $f^*$ is a ``morphism of asphericity
structures'', did imply the condition (stronger than
\eqref{eq:84.3})\pspage{277}
\begin{equation}
  \label{eq:84.4}
  f_!(M'\suba) \subset M\suba,\tag{4}
\end{equation}
we would replace \ref{cond:84.Awk} by this condition \eqref{eq:84.4},
which doesn't look awkward any longer, and does not refer to any $A$
or $i'$ whatsoever (and the six conditions \ref{it:84.i} to
\ref{it:84.iidblprime}, with the exception of \ref{it:84.idblprime}
and \ref{it:84.iidblprime}, do not make any reference to any given
aspheric functor). In any case, one might think of defining a
\emph{morphism} $f^*$ of asphericity structures as a functor admitting
a left adjoint satisfying \eqref{eq:84.4}, and such moreover that
\ref{it:84.i} above, or equivalently \ref{it:84.iprime}, namely the
nice symmetric relation to \eqref{eq:84.4}
\begin{equation}
  \label{eq:84.5}
  f^*(M\suba)\subset M'\suba,\tag{5}
\end{equation}
is satisfied. This notion for a morphism, stricter than the one we
adopted provisionally, looks a lot nicer indeed -- the trouble is that
I could not make up my mind if the restriction \eqref{eq:84.4} is not
an unreasonable one. It would be reasonable indeed, if for those cases
which are the most interesting for us, and primarily for any functor
\[f^*=i^*: M\to\Ahat\]
associated to an $M\suba$-\scrW-aspheric functor
\[ i:A\to M\]
for any given ``nice'' asphericity structure $(M,M\suba)$, this
condition is satisfied. Of course, in this case, $i^*$ does admit a
left adjoint
\[ i_!:\Ahat\to M,\]
provided only $M$ is stable under (small) direct limits, which we'll
assume without any reluctance! Thus, we are led to the following
\begin{questionnum}\label{q:84.1}
  Under which conditions is it true, for an $M\suba$-\scrW-aspheric
  functor $i:A\to M$ (where $M$ is endowed with an asphericity
  structure $M\suba$) that
  \begin{equation}
    \label{eq:84.6}
    i_!(\Ahata) \subset M\suba,\tag{6}
  \end{equation}
  i.e., $i_!$ (the canonical extension of $i$ to \Ahat, commuting with
  direct limits) takes aspheric objects into aspheric objects?
\end{questionnum}

Yesterday and today I pondered mainly about the typical case when
$M=\Bhat$, endowed with its canonical asphericity
structure,\pspage{278} where $B$ is a small category, and moreover $i$
is a functor $i:A\to B\subset \Bhat$. This then brings us to the
related
\setcounter{questionnum}{0}
\renewcommand*{\thequestionnum}{\arabic{questionnum}'}
\begin{questionnum}\label{q:84.1prime}
  Let $i:A\to B$ be an aspherical map in \Cat, under which restrictive
  conditions on $A,B,i$ (if any) does
  \[i_!: \Ahat\to\Bhat\]
  take aspheric objects into aspheric ones, i.e., do we
  have\scrcomment{here in the typescript, AG reverts to putting an
    extra ``s'' in the subscript for asphericity structures, but I'm
    sticking to just ``a'' from now on}
  \begin{equation}
    \label{eq:84.7}
    i_!(\Ahata) \subset B\uphat\suba\text{?}\tag{7}
  \end{equation}
\end{questionnum}
\renewcommand*{\thequestionnum}{\arabic{questionnum}}

The question makes sense even without assuming $i$ to be aspheric. If
we drop this asphericity assumption on $i$, we are reduced, for a
given $B$, (by replacing $A$ by $A_{/F}$, where $F$ is a given
aspheric object of \Ahat) to the case when $A$ is aspheric and
$F=e_\Ahat$, i.e., to looking at whether the element
\begin{equation}
  \label{eq:84.8}
  i_!(e_\Ahat) = \varinjlim_A \prescript{}{(\Bhat)}{i(a)}\tag{8}
\end{equation}
in \Bhat{} is aspheric. This element appears as the direct limit in
\Bhat{} of aspheric elements $i(a)$, the indexing category $A$ being
itself aspheric. Thus, it is rather tempting to hope that the limit
might well be aspheric too. Let's assume $\scrW=\scrWoo=$ usual weak
equivalence, it turns out that when $A$ is $1$-connected, so is
\eqref{eq:84.8}, which would seem to give some support to the hope
that \eqref{eq:84.8} is aspheric if $A$ is. However, this is
definitely not so, unless at the very least we assume $B$ to be
equally aspheric (which is however a rather natural assumption, as it
follows from $A$ aspheric if we assume moreover $i$ to be
aspheric). To see this, we take $i$ to be cofibering with
$0$-connected fibers -- the cofibration assumption implies that the
direct limit over $A$ can be computed first fiberwise and then take a
limit over $B$ (the most general version of associativity for direct
limits!), whereas the $0$-connectedness assumption on the fibers
implies that the limit taken over the fiber $A_b$ is just $b$ itself,
hence \eqref{eq:84.8} is isomorphic to
\[ \varinjlim_b \prescript{}{(\Bhat)}{b} = e_\Bhat,\]
which isn't aspheric except precisely when $B$ is; now for topological
reasons it is easy to find an aspheric $A$ cofibered with
$0$-connected fibers over a non-aspheric $B$.

Thus, in the question of asphericity of \eqref{eq:84.8} we'll better
assume both $A$ \emph{and} $B$ aspheric. If either $A$ has a final
object or $B$ has an initial object, asphericity of \eqref{eq:84.8} is
more or less trivial in any case. For certain categories $B$ even
without initial object, \eqref{eq:84.8} is always aspheric when $A$
is, this I checked at any rate for the ordered category
\[ B= \begin{tikzcd}[row sep=-3pt,column sep=small,baseline=(O.base)]
  \alpha\ar[dr] & \\ & |[alias=O]| \gamma \\ \beta\ar[ur] &
\end{tikzcd},\]
giving rise to a rather interesting computation. This suggested, for
arbitrary $B$ again, to look at the dual of the category above as $A$:
\begin{equation}
  \label{eq:84.9}
A= \begin{tikzcd}[row sep=-3pt,column sep=small,baseline=(O.base)]
  \alpha & \\ & |[alias=O]| \gamma\ar[ul]\ar[dl] \\ \beta &
\end{tikzcd}.\tag{9}
\end{equation}
Here the question then amounts to \emph{whether for a diagram}
\[\begin{tikzcd}[row sep=-3pt,column sep=small]
  & a \\ c\ar[ur]\ar[dr] & \\ & b
\end{tikzcd}\]
\emph{in $B$, the amalgamated sum}
\begin{equation}
  \label{eq:84.10}
  a \amalg_c b \tag{10}
\end{equation}
\emph{in \Bhat{} is aspheric.} This I know to be true if either
$c\to a$ or $c\to b$ are monomorphisms (as stated in section
\ref{sec:70}). Now, I don't really expect this to be true in general,
even if $B$ is such an excellent category as $\Simplex$ say, however,
I didn't push through and make a counterexample. If we now want to
remember the asphericity condition on $i$ which we dropped, this
condition, when $A$ has an initial object $\gamma$, just means that
$c=i(\gamma)$ is an initial object of $B$. Even with this extra
condition, I do not really expect \eqref{eq:84.10} to be necessarily
aspheric. Therefore, I do not expect the inclusion \eqref{eq:84.7} to
hold for any aspheric functor $i$, even when $A$ (and hence $B$) are
assumed to be aspheric, without some extra condition, on $A$ or $B$
say.

The condition that \Ahat{} be \emph{totally aspheric} seems in this
context a rather natural one -- for instance, when reducing the
question whether \eqref{eq:84.7} holds to the question of asphericity
of an object \eqref{eq:84.8} (where $A$ stands for $A_{/F}$), and if
we do not want to loose the asphericity assumption for $i$ when taking
the composition $A_{/F}\to A\xrightarrow i B$, we would like that
asphericity of $F$ imply asphericity of $A_{/F}\to A$ -- which means
precisely that $A$ is totally aspheric. The trouble is that in the
would-be counterexample above with $A$ given by \eqref{eq:84.9}, $A$
is stable under Inf, i.e., under binary products, hence \Ahat{}
\emph{is} totally aspheric indeed -- thus it is doubtful that this
extra condition in the ``question \ref{q:84.1prime}\kern1pt'' above is
quite enough. If it does fail indeed, the next best would be to try
the stronger condition ``$A$ is a strict test category'' (NB\enspace
the category \eqref{eq:84.9} isn't a test category!), which leads to
the question whether \eqref{eq:84.8} is aspheric when $A$ is a test
category (no longer a strict one though!) and $i:A\to B$ an aspheric
functor. But I confess I have no idea at present how to handle this
question, and I am dubious there will come out any positive result
along these lines, even when assuming $A$ and $B$ to be both stable
under binary products say, and to be contractors and what-not!

Another\pspage{280} typical case for ``Question \ref{q:84.1}'' (p.\
\ref{p:277}) is the case when
\[ M= \Cat,\]
endowed with the usual asphericity structure, giving rise to the
notion of ``test functor''
\[i:A\to\Cat\]
we have been working with almost from the very start. The most
important case of all is of course the canonical functor
\[ i_A: a\mapsto A_{/a} : A\to\Cat,\]
and we well know that this functor is aspheric (for the natural
asphericity structure of \Cat) if $A$ is a weak test category, and in
this case it is true indeed that not only $i_A^*$, but equally
${i_A}_!$ is modelizing, and transforms aspheric objects into aspheric
objects. But we know too that for more general test functors, for
instance the standard inclusion
\begin{equation}
  \label{eq:84.11}
  i:\Simplex\hookrightarrow\Cat,\tag{11}
\end{equation}
it is no longer true in general that $i_!$ be modelizing -- so is it
at all reasonable to expect
\begin{equation}
  \label{eq:84.12}
  i_!(\Ahata) \subset\Cat\suba \text{?}\tag{12}
\end{equation}
Definitely, I'll have to find out the answer in the typical case
\eqref{eq:84.11}, whether I like it or not!

\bigbreak
\presectionfill\ondate{30.6.}\par

\hangsection[Digression on left exactness properties of $f_!$
functors, \dots]{Digression on left exactness properties of
  \texorpdfstring{$f_!$}{f!} functors, application to the inclusion
  \texorpdfstring{\protect\iSimplexIntoCat}{i:Delta->(Cat)}.}\label{sec:85}%
Just got an impressive heap of reprints and preprints by Tim Porter
(announced in his letter two days ago). Many titles have to do with
``shape theory'', ``coherence'' and ``homotopy limits'' (corresponding
to proobjects in various modelizers, as I understand it from his
letter). The one which most attracts my attention though is ``Cat as a
closed model category'' by
R.W.~Thomason\scrcomment{\textcite{Thomason1980}} -- the title comes
quite as a surprise, as my ponderings on the homotopy structure of
\Cat{} had left me with the definite feeling that there wasn't a
closed model structure on \Cat, with the usual notion of weak
equivalence. I'll have to have a closer look at this paper definitely,
before starting on part \ref{ch:V} of these notes!

It was getting prohibitively late yesterday, so I had to stop. It then
occurred to me that the (admittedly provisional) notion of
\emph{morphism} of asphericity structures I had proposed yesterday
(p.\ \ref{p:275}) is definitely stupid, because there is no reason
whatever that it should be stable under composition! This surely was
the reason for the feeling of uneasiness\pspage{281} caused by the
condition \ref{cond:84.Awk}, which surely deserves its name! This flaw
of course disappears, if we strengthen the condition, as suggested on
p.\ \ref{p:277}, by condition \ref{eq:84.4} -- that $f_!$ should take
aspheric objects into aspheric ones. This makes all the more
imperative the task of finding out whether this condition is
reasonable, and otherwise, what to put in as a substitute. The test
case now seems to be really the one when $f^*$ is the usual nerve
functor
\[f^*=i^*:\Cat\to\Simplexhat,\]
associated to the standard inclusion
\[i:\Simplex\to\Cat.\]
I never before had a closer look upon the corresponding functor
\[i_!:\Simplexhat\to\Cat,\]
left adjoint to the nerve functor, except just rectifying a big
blunder (p.\ \ref{p:22}), and convincing myself that $i_!$ did
\emph{not} take weak equivalences into weak equivalences. But why not
aspheric objects into aspheric ones?

Pondering a bit over this matter today, and trying to see whether
$i_!$ takes \emph{contractible} objects into contractible ones, this
brought up the question whether $i_!$ commutes with finite products --
which will allow then to make use of the generalities on morphisms of
contractibility structures, and will imply that $i_!$ is indeed such a
morphism (as it transforms $\Simplex$, which generates the
contractibility structure of $\Simplexhat$, into contractible elements
of \Cat). As a matter of fact, a number of times in my scratchwork
during the last four months I've met with the question of when a
functor of the type $f_!$ commutes to various types of finite inverse
limits, and this now is the occasion for writing down some useful
general facts in this respect, which had remained somewhat in the air
so far. It seems that the relevant facts can be summed up in two
steps.
\begin{propositionnum}\label{prop:85.1}
  Let $M$ be a category stable under \textup(small\textup) direct
  limits, $A$ a small category, $f:A\to M$ a functor, and
  \[f_!:\Ahat\to M, \quad f^*:M\to\Ahat\]
  the corresponding pair of adjoint functors.
  \begin{enumerate}[label=\alph*),font=\normalfont]
  \item\label{it:85.prop1.a}
    Assume $M$ stable under binary products, and that for any $x$ in
    $M$, the product functor $y\mapsto y\times x$ commutes to
    arbitrary direct limits \textup(with small categories of indices
    of course\textup). Then, in order for $f_!$ to commute with binary
    products,\pspage{282} is it n.s.\ that for any two objects $a,b$
    of $A$\kern1pt, the map
    \begin{equation}
      \label{eq:85.star}
      f_!(a\times b) \to f(a) \times f(b)\tag{*}
    \end{equation}
    in $M$ be an isomorphism.
  \item\label{it:85.prop1.b}
    Assume $M$ stable under fibered products, and that base change in
    $M$ commutes to arbitrary direct limits \textup(with small
    indexing categories of course\textup). Then, in order for $f_!$ to
    commute with fibered products, it is n.s.\ that it does so for any
    diagram
    \[\begin{tikzcd}[baseline=(O.base),column sep=tiny,row sep=small]
      b\ar[dr] & & c\ar[dl] \\ & |[alias=O]| F &
    \end{tikzcd},\]
    with $b,c$ objects in $A$\kern1pt, $F$ in \Ahat.
  \end{enumerate}
\end{propositionnum}
\begin{corollary}
  Assume $M$ stable under finite inverse limits, and that base change
  in $M$ commutes with arbitrary direct limits. The $f_!$ is left
  exact if{f} it satisfies the conditions \textup{\ref{it:85.prop1.a}}
  and \textup{\ref{it:85.prop1.b}} above, and moreover transforms
  $e_\Ahat$ into a final object of $M$ \textup(which also means, if
  $A$ has a final object $e_A$, that $f(e_A)$ is a final object of
  $M$\textup).
\end{corollary}
\noindent\emph{Proof of proposition.}
\textup{\ref{it:85.prop1.a}}\enspace As we have, by definition of
$f_!$, for any $F$ and $G$ in \Ahat
\[f_!(F) = \varinjlim_{A_{/F}} f(a) ,\quad
f_!(G) = \varinjlim_{A_{/G}} f(b) ,\]
deduced from the corresponding relations in \Ahat
\[ F = \varinjlim_{A_{/F}} a , \quad
G = \varinjlim_{A_{/G}} b ,\]
the looked for bijectivity of
\[ f_!(F\times G) \to f_!(F)\times f_!(G)\]
will follow from the assumption \eqref{eq:85.star}, and from the
following lemma (applied in both categories \Ahat{} and $M$), the
proof of which is immediate:
\begin{lemma}
  Let\pspage{283} $M$ be a category stable under direct limits and
  under binary products, and such that the product functors $y\mapsto
  y\times x$ commute to direct limits. Then binary products are
  distributive with respect to direct limits, i.e., if
  \[ u : I\to M, \quad v : J\to M\]
  are two functor of small categories with values in $M$, we have
  \[\varinjlim_{I\times J} u(i)\times v(j) \tosim \bigl( \varinjlim_I
  u(i) \bigr) \times \bigl( \varinjlim_J v(j) \bigr).\]
\end{lemma}

\emph{Proof of} \ref{it:85.prop1.b}. Follows from \ref{it:85.prop1.a},
applying it to the induced functor
\[ f_{/F} : A_{/F} \to M_{/f_!(F)},\]
in order to prove commutation of $f_!$ with a fibered product
corresponding to a diagram
\begin{equation}
  \label{eq:85.starstar}
  \begin{tikzcd}[column sep=tiny,row sep=small]
    G\ar[dr] & & H\ar[dl] \\ &  F &
  \end{tikzcd}\tag{**}
\end{equation}
in \Ahat.

\begin{remarknum}\label{rem:85.1}
  Part \ref{it:85.prop1.b} does not look as nice as part
  \ref{it:85.prop1.a}, as we would like to be able to take $F$ equally
  in $A$ -- what I first expected to get out. If we make only this
  weaker assumption, it will follow at once that $f_!$ commutes to
  fibered products corresponding to diagrams \eqref{eq:85.starstar}
  \emph{with $F$ in $A$}. It doesn't seem finally that from this
  follows commutation of $f_!$ to all fibered products, unless making
  some extra assumptions, such as $f_!$ transforming monomorphisms
  into monomorphisms (which is a necessary condition anyhow). We'll
  have to come back upon this -- for the time being we need only to
  deal with products, for which \ref{it:85.prop1.a} is adequate.
\end{remarknum}

The ``second step'' now is a tautology:
\begin{propositionnum}\label{prop:85.2}
  Under the preliminary assumptions of prop.\ \ref{prop:85.1} above,
  assume moreover that $f$ is fully faithful, and that $f(A)$
  ``generates $M$ by strict epimorphisms'', i.e., that for any object
  $P$ in $M$, we have
  \begin{equation}
    \label{eq:85.starstarstar}
    P \xleftarrow{\sim} \varinjlim_{A_{/P}} f(x).\tag{***}
  \end{equation}
  Then $f_!$ commutes to finite products in \Ahat{} of objects of $A$\kern1pt,
  and to fibered products in \Ahat{} of objects of $A$\kern1pt.
\end{propositionnum}
\noindent\emph{Proof.}\pspage{284} That $f_!$ transforms final object
into final object follows from the formula \eqref{eq:85.starstarstar}
above, by taking $P$ to be the final object of $M$. To get commutation
to binary products, we apply the formula for $P=a\times_M b$, and
notice that
\[ A_{/P} \simeq A_{/P'}, \quad\text{with $P'=a\times_\Ahat b$,}\]
because the inclusion $f:A\hookrightarrow M$ is full (I forgot to say
we may assume $f$ to be the inclusion functor of a full
subcategory). But we have too
\[ f_!(P') = \varinjlim_{A_{/P'}} f(x),\]
hence $f_!(P') \simeq P$. The proof for fibered products is similar,
or can be reduced to the case of products by considering the induced
functor
\[ f_{/a} : A_{/a} \to M_{/a}.\]
\begin{remarknum}\label{rem:85.2}
  The various variants of the notion of a generating subcategory $A$
  of a category $M$ have been dealth with in some detail in
  SGA~4 vol.~1\scrcomment{\textcite{SGA4vol1}}
  exp.~I par.~7 (p.\ 45--60) (Springer LN 269), cf.\ more
  specifically prop.\ 7.2 on page 47, summarizing the main
  relationships. The strongest of all notions considered there is
  generation by strict epimorphisms, which can be expressed by the
  formula \eqref{eq:85.starstarstar} above, and is also equivalent to
  the \emph{functor $f^*$ being fully faithful}.
\end{remarknum}

We come back now to the case $M=\Cat$ -- it is immediate that the
preliminary exactness conditions on $M$ of prop.\ \ref{prop:85.1} are
satisfied. Thus, putting together prop.\ \ref{prop:85.1} and
\ref{prop:85.2}, we get:
\begin{propositionnum}\label{prop:85.3}
  Let $A$ be a full subcategory of \Cat, generating \Cat{} by strict
  epimorphisms, $i:A\hookrightarrow\Cat$ the inclusion functor. Then
  the corresponding functor
  \[i_!:\Ahat\to\Cat\]
  commutes to finite products, and also to fibered products in \Ahat{}
  over an object of \Ahat{} coming from $A$\kern1pt.
\end{propositionnum}

The most familiar case when this applies is indeed the case of the
canonical inclusion of $\Simplex$ into \Cat, as it is well-known that
the corresponding $i^*$, namely the nerve functor, is fully
faithful. From this follows that \emph{a fortiori},\marginpar{hasty!}
for any full subcategory $A$ of \Cat{} containing $\Simplex$, $A$
generates \Cat{} by strict epimorphisms, and hence proposition
\ref{prop:85.3} applies.
\begin{remarknum}\label{rem:85.3}
  I doubt however, even for $A=\Simplex$, the most typical
  case,\pspage{285} that $i_!$ is even left exact, because it looks
  unlikely to me, in view of the explicit description of
  $i_!(K_\bullet)$ for a given ss~complex $K_\bullet$, in terms of
  ``generators and relations'' (cf.\ Gabriel-Zisman's
  book),\scrcomment{\textcite{GabrielZisman1967}} that $i_!$ should
  transform monomorphisms into monomorphisms. For what follows, this
  is irrelevant anyhow.
\end{remarknum}

\hangsection[Bimorphisms of contractibility structure as the (final?)
\dots]{Bimorphisms of contractibility structure as the
  \texorpdfstring{\textup(\kernifitalic{2pt}}{(}final?\texorpdfstring{\textup)}{)}
  answer. Does the notion of a map of asphericity structures exist?}
\label{sec:86}%
After these somewhat painful preliminaries, it seems to me that firm
ground is in sight at last! The main feeling that finally comes out of
these reflections, is that we have a very good hold on checking
whether or not a functor $f_!:\Ahat\to M$ commutes to finite products,
and moreover, that this condition is satisfied in many ``good'' cases,
in more cases at any rate than I suspected. In case of the inclusion
functor of $\Simplex$ into \Cat, more generally of any subcategory $A$
of \Cat{} consisting of contractible objects, containing $\Simplex_1$,
and large enough in order to generate \Cat{} by strict epimorphisms
(or, as we'll say, to \emph{generate strictly} the category \Cat), it
follows that we have the inclusion
\begin{equation}
  \label{eq:86.1}
  i_!(\Ahatc) \subset \Cat\subc,\tag{1}
\end{equation}
i.e., $i_!$ transforms \emph{contractible} objects into contractible
objects, or, equivalently, is a \emph{morphism for the homotopy
  structures on \Ahat{} and \Cat}. (Note that the two first
assumptions on the subcategory $A$ of \Cat{} imply that the homotopy
structure on \Ahat{} generated by the ``intervals'' made up with
elements of $A$, is a contractibility structure, admitting $A$ as a
generating set of contractible objects. The situation is really nice
only when assuming that \Ahat{} is totally aspheric though, which will
imply that $A$ is a test category (even a strict one) and that the
\scrW-asphericity structure on \Ahat{} associated to the
contractibility structure just described is the usual one (cf.\ prop.\
\ref{prop:79.2} \ref{it:79.2.b} p.\ \ref{p:248}, where we take $C$ to
be $A$).)

The inclusion \eqref{eq:86.1} of course strongly resembles the
problematic inclusion \eqref{eq:84.4} (p.\ \ref{p:277}) of yesterday's
notes, namely, here to $i_!$ taking aspheric objects into aspheric
ones. I doubt though this latter is true even for the standard case
$A=\Simplex$? But before trying at all costs to see whether it holds
or not, the point I wish to make is that for the purpose we have in
view, namely defining suitable conditions, stable under composition,
on a pair $(f_!,f^*)$ of adjoint functors between two asphericity
structures $M$ and $M'$, -- for this purpose, a relation of the type
\eqref{eq:86.1}, namely\pspage{286}
\begin{equation}
  \label{eq:86.2}
  f_!(M'\subc) \subset M\subc,\tag{2}
\end{equation}
is just as good as the relation \eqref{eq:84.4} on p.\ \ref{p:277},
provided we make the evident assumption needed for \eqref{eq:86.2} to
make sense, namely \emph{the given asphericity structures on $M,M'$ to
  be associated to contractibility structures}. This condition
\eqref{eq:86.2} is (considerably!) weaker than
\begin{equation}
  \label{eq:86.2prime}
  f_!(M'\suba)\subset M\suba,\tag{2'}
\end{equation}
due to the inclusion $M'\subc\subset M'\suba$ and to the fact that
(with the notations and using the results of section \ref{sec:81}) the
inclusion \eqref{eq:86.2} is implied by the apparently weaker
inclusion (under the assumption \ref{loc:4}):
\[f_!(M'\subc)\subset M_0,\]
where $M\subc\subset M\suba\subset M_0$. On the other hand,
\eqref{eq:86.2} \emph{is evidently stable under composition}.

To be specific about ``as good as'', let's come back to the statement
of the ``pretty'' proposition on p.\ \ref{p:275}, somewhat marred by
the ``awkward'' condition \ref{cond:84.Awk} we had to throw in, which
looks so ugly because of its lack of stability under
compositions. Let's replace this condition by the slightly stronger
one \eqref{eq:86.2}, which is stable under composition. Condition
\eqref{eq:86.2} does imply \ref{cond:84.Awk} indeed, as we see by
taking for $A$ a small full subcategory of $M\subc$, generating the
asphericity structure of $M$, and for $B$ any full subcategory of $M$
contained in $M\suba$ and big enough in order a)\enspace to contain
$f_!(A)$ and b)\enspace to generate the asphericity structure of
$M$. We may take for instance for $A$ any subcategory of $M'\subc$
closed under finite products in $M'$ and generating the
contractibility structure, and accordingly for $B$, so that the pair
$(A,B)$ will match for any choice of the basic localizer \scrW. Thus,
the six conditions \ref{it:84.i} to \ref{it:84.iidblprime} of the
proposition are equivalent, provided in the last
\ref{it:84.iidblprime} we assume moreover that $B$ contain a final
object $e_M$ of $M$. Moreover, it is clear now that the conditions can
be viewed also as a property of the functor
\begin{equation}
  \label{eq:86.3}
  {f_!}\subc : M'\subc \to M\subc\tag{3}
\end{equation}
induced by $f_!$, which we may adequately express by saying that this
functor \eqref{eq:86.3} is \emph{\scrW-aspheric} (where the two sides
of \eqref{eq:86.3} are categories which need not be small, but which
are at any rate \emph{\scrW-aspherators} in the sense of section
\ref{sec:79}, p.\ \ref{p:247}).

Stated\pspage{287} this way, the condition just obtained on the pair
of adjoint functors $(f_!,f^*)$, between the two categories $M,M'$
endowed with contractibility structures, namely \eqref{eq:86.2} and
\begin{equation}
  \label{eq:86.4}
  f^*(M\suba)\subset M'\suba,\tag{4}
\end{equation}
are not quite symmetric, as the condition on $f_!$ is expressed in
terms of contractible objects, whereas the condition on $f^*$ is
expressed in terms of aspheric ones. However, assuming that \scrW{}
satisfies \ref{loc:4} and using theorem \ref{thm:79.1} p.\
\ref{p:252}, we see that \eqref{eq:86.4} is equivalent to the
condition
\begin{equation}
  \label{eq:86.5}
  f^*(M\subc)\subset M'\subc,\tag{5}
\end{equation}
which does not depend any longer on the choice of \scrW!

Finally, we are led to a notion of pure homotopy theory, in terms of
contractibility structures alone, without the intrusion of the choice
of a basic localizer \scrW{} and corresponding asphericity notions. We
may call the pair a ``\emph{bimorphism}'' of contractibility
structures, and introduce it via the following summing-up statement:
\begin{scholie}
  Let\scrcomment{a
    ``\href{https://en.wikipedia.org/wiki/Scholia}{scholium}'' is a
    critical or explanatory comment extracted from preexisting
    propositions} $(M,M\subc)$ and $(M',M'\subc)$ be two
  contractibility structures, each admitting a small generating
  subcategory for the contractibility structure, say $C$ and $C'$
  respectively. Let $(f_!,f^*)$ be a pair of adjoint functors
  \[f_!:M'\to M, \quad f^*:M\to M'.\]
  Let's consider the following inclusion conditions
  \begin{description}
  \item[\namedlabel{it:86.bang}{($!$)}]
    $f_!(M'\subc)\subset M\subc$,
  \item[\namedlabel{it:86.bangprime}{($!$')}] 
    $f_!(C')\subset M_0$,
  \item[\namedlabel{it:86.star}{($*$)}]
    $f^*(M\subc)\subset M'\subc$,
  \item[\namedlabel{it:86.starprime}{($*$')}] 
    $f^*(C)\subset M'_0$,
  \end{description}
  where $M_0$ and $M'_0$ are \textup(as in section
  \ref{sec:81}\textup) the sets of ``$0$-connected'' objects in $M$
  and $M'$ respectively, for the given contractibility structures.
  \begin{enumerate}[label=\alph*),font=\normalfont]
  \item\label{it:86.a}
    Clearly, in view of
    \[C \subset M\subc \subset M_0, \quad C'\subset M'\subc\subset
    M'_0,\]
    condition \textup{\ref{it:86.bang}} implies
    \textup{\ref{it:86.bangprime}} and condition
    \textup{\ref{it:86.star}} implies
    \textup{\ref{it:86.starprime}}. Moreover, it is even true that conditions \textup{\ref{it:86.star}} and
    \textup{\ref{it:86.starprime}} are equivalent, and the same holds
    for \textup{\ref{it:86.bang}} and \textup{\ref{it:86.bangprime}}
    provided $f_!$ commutes with finite products. We'll say that
    $(f_!,f^*)$ is a \emph{bimorphism} for the given contractibility
    structures, if the inclusions \textup{\ref{it:86.bang}} and
    \textup{\ref{it:86.star}} do hold.
  \item\label{it:86.b} Assume\pspage{288} we got a bimorphism
    $(f_!,f^*)$, i.e., \textup{\ref{it:86.bang}} and
    \textup{\ref{it:86.star}} above hold. Let \scrW{} be a basic
    localizer, hence the sets $M_\scrW$ and $M'_\scrW$ of
    \scrW-aspheric objects in $M$ and $M'$ respectively, and the sets
    $\scrW_M$ and $\scrW_{M'}$ of \scrW-equivalences in $M$ and
    $M'$. Then the following relations hold
    \begin{description}
    \item[\namedlabel{it:86.W}{(\scrW)}] $M_\scrW =
      (f^*)^{-1}(M'_\scrW)$,
    \item[\namedlabel{it:86.Wprime}{(\scrW')}] $\scrW_M =
      (f^*)^{-1}(\scrW_{M'})$.      
    \end{description}
  \item\label{it:86.c}
    Under the assumption of \textup{\ref{it:86.b}}, hence
    \textup{\ref{it:86.Wprime}} holds and $f^*$ gives rise to a
    functor $\overline{f^*}$ between the localizations
    \[\HotOf_{M,\scrW} = \scrW_M^{-1}M \quad\text{and}\quad
    \HotOf_{M',\scrW}=\scrW_{M'}^{-1}M',\]
    the following diagram is commutative up to canonical isomorphism:
    \[\begin{tikzcd}[baseline=(O.base),column sep=tiny]
      \HotOf_{M,\scrW}\ar[rr,"\overline{f^*}"]
      \ar[dr] & & \HotOf_{M',\scrW}\ar[dl] \\
      & |[alias=O]| \HotOf_\scrW &
    \end{tikzcd},\]
    where the vertical functors are the canonical functors of section
    \ref{sec:77}. In particular, if the latter are equivalences
    \textup(i.e., $M$ and $M'$ are \scrW-modelizing\textup), then so
    is $\overline{f^*}$.
  \item\label{it:86.d}
    Assume merely that the inclusion \textup{\ref{it:86.bang}} holds,
    and let \scrW{} as in \textup{\ref{it:86.b}}\textup{\ref{it:86.c}}
    be a basic localizer, satisfying moreover
    \textup{\ref{loc:4}}. Then the following conditions on the pair
    $(f_!,f^*)$ are equivalent:
    \begin{enumerate}[label=(\roman*),font=\normalfont]
    \item\label{it:86.d.i}
      The pair is a bimorphism, i.e., \textup{\ref{it:86.star}}
      \textup(or equivalently, \textup{\ref{it:86.starprime})} holds.
    \item\label{it:86.d.ii}
      The inclusion $\subset$ in \textup{\ref{it:86.W}} above holds.
    \item\label{it:86.d.iii}
      The inclusion $\subset$ in \textup{\ref{it:86.Wprime}} above holds.
    \item\label{it:86.d.iv}
      The functor induced by $f_!$\kern1pt,
      \[{f_!}\subc : M'\subc \to M\subc\]
      is a \scrW-aspheric functor between the aspherators $M'\subc,M\subc$.
    \item\label{it:86.d.v}
      \textup(\kern1pt If $A$ is a given small category, and
      \[i':A\to M'\]
      a given functor, \emph{factoring through $M'\subc$}, and
      $M'_\scrW$-\scrW-aspheric.\textup) The composition
      \[i = f_!\,i' : A\to M\]
      is $M_\scrW$-\scrW-aspheric.
    \end{enumerate}
  \end{enumerate}
\end{scholie}
\begin{comments}
  In terms of what is known to us since sections \ref{sec:81} and
  \ref{sec:82}, and notably that any two sections of a $0$-connected
  object of $M$ or of $M'$ are homotopic, the Scholie is just a long
  tautology. I have taken great case however in stating it, so as to
  get as clear a view as possible of exactly what the relevant
  relationships are. In \ref{it:86.a} the conditions \ref{it:86.star}
  and \ref{it:86.starprime} are just different ways of stating that
  $f^*$ is a morphism of contractibility structures, which can still
  be expressed in various other ways, compare p.\
  \ref{p:251}--\ref{p:252}. Similarly, if we assume that $f_!$
  commutes with finite products, then \ref{it:86.bang} or
  \ref{it:86.bangprime} can be viewed as expressing that $f_!$ is a
  morphism of homotopy structures (in opposite direction), which again
  could be expressed in various other ways, for instance in terms of
  the homotopy relation for maps, or in terms of homotopisms. This
  implies that in any case $f^*$ induces a functor between the strict
  localizations
  \[\overline{f^*}^{\mathrm c} : W\subc^{-1}M=\overline M \to
  {W'\subc}^{-1}M' = \overline{M'},\]
  and similarly for $f_!$, when $f_!$ commutes to finite products. I
  doubt however that the latter localized map is of geometric
  relevance, except maybe in the cases when $f_!$ gives rise to
  relations similar to \ref{it:86.W} and \ref{it:86.Wprime} in
  \ref{it:86.b} above, and both $f_!$ and $f^*$ are model-preserving
  with respect to \scrW-equivalences, and define quasi-inverse
  equivalences between the localizations $\HotOf_{M,\scrW}$ and
  $\HotOf_{M',\scrW}$ -- a rather exceptional case indeed.
\end{comments}

After stating the Scholie, there is scarcely a doubt left in my mind
about the notion of bimorphism, which finally peeled out of
reflections, being a relevant one. This Scholie, rather than the
``seducing'' proposition of section \ref{sec:84} (p.\ \ref{p:275}),
now seems to me the adequate ``answer'' of the not-so-silly-after-all
question of section \ref{sec:46} (cf.\ page \ref{p:95} -- nearly 200
pages ago!). The only minor uncertainty remaining in my mind is
whether or not in the notion of a bimorphism of contractibility
structures, we should insist that $f_!$ should commute to finite
products. It seems that it will be hard to check condition
\ref{it:86.bang}, except via an apparently weaker form such as
\ref{it:86.bangprime}, and using commutation of $f_!$ to finite
products. But on the other hand, the Scholie makes good sense without
assuming such commutation property. If the notion of a bimorphism is
going to be useful, time only will tell which terminology use is the
best.

More puzzling however is this facit,\scrcomment{``facit'' comes from
  the Latin verb ``\emph{facere}'', to do, so it means the result. Or
  it could be a typo for ``fact'', which essentially means the
  same\ldots} that we still have not been able to give a satisfactory
definition of a \emph{morphism} of \emph{asphericity} structures,
independently of any restrictive assumption on these, such\pspage{290}
as being generated by a contractibility structure. Even when making
such an assumption, it is not wholly clear yet that there isn't a good
notion of a morphism $f^*$ of contractibility structures, without
having to assume there exists a left adjoint $f_!$. Maybe a little
more pondering on the situation would be useful. If it should become
clear that (except for the obvious notion of \emph{equivalence} of
asphericity structures) there was \emph{not} any reasonable notion of
a morphism of asphericity structures, this would probably mean that
the notion of an asphericity structure should not be viewed as a main
structure type in homotopy theory in its own right, but rather,
mainly, as an important by-product of a contractibility structure. At
any rate, since section \ref{sec:81}, I feel that the main emphasis
has definitely been shifting towards contractibility structures, which
seems now \emph{the} main type of structure dominating in the
modelizing story, as this ``story'' is gradually emerging into light.



\chapter{Abelianization I}
\label{ch:V}

\presectionfill\ondate{1.7.}\pspage{291}\par

\hangsection[Comments on Thomason's paper on closed model structure
\dots]{Comments on Thomason's paper on closed model structure on
  \texorpdfstring{\Cat}{(Cat)}.}\label{sec:87}%
I couldn't resist last night and had to look through Thomason's
preprint on the closed model structure of
\Cat.\scrcomment{\textcite{Thomason1980}} The paper is really pleasant
reading -- and it gives exactly what had been lacking me in my
reflections lately on the homotopy theory of \Cat{} -- namely, a class
of neat monomorphisms $Y\to X$ which all have the property that cobase
change by these preserves weak equivalence, the so-called \emph{Dwyer
  maps}.\scrcomment{These Dwyer maps are not closed under retracts, so
  \textcite{Cisinski1999} introduced the
  \href{http://ncatlab.org/nlab/show/Dwyer+map}{\emph{pseudo-Dwyer
      maps}} which are, and \emph{they} are now called Dwyer maps.} I
had hoped for a while that ``open immersions'' and their duals, the
``closed immersions'' (namely sieve and cosieve maps, in Thomason's
wording), have this property, and when it turned out they hadn't, I
had been at a loss of what stronger property to put instead, wide
enough however to allow for the standard factorization statements for
a map to go through. The definition of a Dwyer map is an extremely
pretty one -- it is an open immersion $Y\to X$ such that the induced
map from $Y$ into its closure $\overline Y$ should have a right
adjoint. Now this implies that $Y\to\overline Y$ is aspheric, and I
suspect that this extra condition on an open immersion $Y\to X$ should
be sufficient to imply that it has the ``cofibration property'' above
with respect to weak equivalence. It would mean in a sense that the
given open immersion is very close to being a closed immersion too,
without however being a direct summand necessarily. The dual notion is
that of a closed immersion such that the corresponding interior
$\mathring Y$ of $Y$ in $X$ gives a map $\mathring Y \to Y$ which has
a left adjoint, or only which is ``coaspheric''. By Quillen's duality
principle, if one notion works well for pushouts, so does the
other. -- With this notion in hands, it shouldn't be difficult now to
get a closed model structure on \Cat{} a lot simpler than
Thomason's. Visibly, he was hampered by the standing reflex: homotopy
= semi-simplicial algebra, which caused him to pass by the detour of
the category \Simplexhat{} of semi-simplicial complexes, rather than
just working in \Cat{} itself. I'll have to come back upon this in
part \ref{ch:V} of the notes, where I intend to investigate the
homotopy properties of \Cat{} and elementary modelizers \Ahat,
including the existence of closed model structures.

Some comments of Thomason's at the end of his preprint, about
application to algebraic K-theory, seem to indicate that the notion of
``integration'' and ``cointegration'' of homotopy types I have been
interested in, has been studied (under the name holim and hocolim) in
the context of closed model categories by Anderson (his paper appeared
in 1978).\scrcomment{\textcite{Anderson1978}} As I am going to develop
some ideas along these lines in part \ref{ch:VI} on derivators, I
should have a look at what Anderson does, notably\pspage{292} what his
assumptions on the indexing categories are. Thomason seems to believe
that the closed model structure of \Cat{} is essential for being able
to take homotopy limits -- whereas it is clear a priori to me that the
notion depends only on the notion of weak equivalence. Indeed, he
seems to consider the possibility of taking homotopy limits in \Cat{}
as the main application of his theorem, and in order to apply
Anderson's results, believes it is necessary to be able to give
concrete characterizations for ``fibrations'' and ``cofibrations'' of
his closed model structure. Now, it turns out that the case he is
interested in (for proving ``Lichtenbaum's
conjecture'')\scrcomment{aka the
  \href{https://en.wikipedia.org/wiki/Quillen-Lichtenbaum_conjecture}{Quillen–Lichtenbaum
    conjecture}} is a typical
case of \emph{direct} homotopy limits, namely ``integration'' -- which
can be described directly in \Cat{} in such an amazingly simple way
(as sketched in section \ref{sec:69}, p.\
\ref{p:198}--\ref{p:199}). Thus, I feel for this application, the
closed model structure is wholly irrelevant. As for cointegration, I
do not expect that there is a comparably simple construction of this
operation within the modelizer \Cat, but presumably there is when
taking \oo-Gr-stacks as models (as suggested by the ``geometric''
approach to cohomology invariants, via stacks, where the operation of
``direct image'', namely cohomology precisely, is the obvious one,
whereas inverse images are more delicate to define, by an adjunction
property with respect to direct images\ldots). When working in \Cat,
cointegration of homotopy types should be no more nor less involved
than in any closed model category say, and involve intensive recourse
to ``fibrations'' (in the sense of the closed model structure, or more
intrinsically, in the sense that base change by these should preserve
weak equivalence). Now the latter have become quite familiar to me
during my long scratchwork on cohomology properties of maps in \Cat,
and I'll have to try to put it down nicely in part \ref{ch:V} of the
notes.

\hangsection[Review of pending questions and topics (questions 1) to
\dots]{Review of pending questions and topics
  \texorpdfstring{\textup(questions \textup{1)} to \textup{5)},
    including characterizing canonical modelizers\textup)}{(questions
    1) to 5), including characterizing canonical
    modelizers)}.}\label{sec:88}%
The present ``part \ref{ch:IV}'' on asphericity structures (and their
relations to contractibility structures) turns out a lot longer than I
anticipated, and the end is not yet quite in sight! Therefore, before
pursuing, I would like to make a review of the questions along these
lines which seem to require elucidation, and then decide which I'm
going to deal with, before going over to part \ref{ch:V}.

\namedlabel{q:88.1}{1)}\enspace Whereas the relevant notions of
``morphisms'' and ``bimorphisms'' for contractibility structures seem
to me well understood, there remains\pspage{293} a certain feeling of
uneasiness with respect to asphericity structures, which haven't got
yet a reasonable notion of morphism. Thus, I have left unanswered the
questions raised in section \ref{sec:84} around the inclusion
condition
\[ f_!(M\suba) \subset M'\suba,\]
so it is not impossible that, while relying on the mere feeling that
the inclusion is just not reasonable and that the answer to the
specific questions are presumably negative, I am about to miss some
unexpected important fact! It seems I developed kind of a block
against checking -- maybe the answers are well-known and Tim Porter
will tell me\ldots Maybe I better leave the question for a later
moment, as it will ripen by itself if I leave it alone\ldots

\namedlabel{q:88.2}{2)}\enspace I should at last introduce
contractors, and morphisms of such. When I set out on part \ref{ch:IV}
of the notes, I expected that the notion of a contractor would be one
main notion, alongside with the notion of a canonical modelizer -- it
turns out that so far I didn't have much use yet for one or the
other. Contractors can be viewed as categories generating
contractibility structures, just as aspherators are there for
generating asphericity structures. However, whereas any small category
is an aspherator, the same is definitely not so for contractors, as we
demand that every object in $C$ should be contractible, for the
homotopy structure in \Chat{} generated by $C$ itself. If we except
the case of a contractor equivalent to the final category (a so-called
\emph{trivial} contractor), a contractor is a strict test category --
thus the notion appears somewhat as a hinge notion between the test
notions, and the ``pure'' homotopy notions and more specifically,
contractibility structures. Writing up some scratchnotes I got should
be a pure routine matter.

\namedlabel{q:88.3}{3)}\enspace A lot more interesting seems to me to
try and resolve a persistent feeling of uneasiness which has been
floating, throughout the long-winded reflections on homotopy and
asphericity structures in parts \ref{ch:III} and \ref{ch:IV} of the
notes. This is tied up with this fact, that my treatment of the main
notions, namely contractibility and asphericity, has been consistently
\emph{non-autodual}. More specifically, when a category $M$ is
endowed with either a contractibility or an asphericity structure, it
does \emph{not} follow that the opposite category $M\op$ is too in a
natural way. When defining homotopy relations and homotopism
structures (section \ref{sec:51}, \ref{subsec:51.A} and
\ref{subsec:51.B}), these were autodual notions, but the notion of a
homotopy interval structure, which we used in order to pass from a
contractibility\pspage{294} structure to the corresponding notions of
homotopy equivalence between maps and of homotopisms, is highly
non-autodual too. It breaks down altogether when $M$ is a ``pointed''
category, namely contains an object which is both initial and final --
in this case, for any homotopy interval structure on $M$, any two maps
in $M$ are homotopic, and any map is a homotopism, hence any object is
contractible!

Our initial motivation, namely understanding ``modelizers'' for
ordinary homotopy types, made it very natural to get involved in
non-autodual situations, as the homotopy category \Hot{} itself, and
the usual model categories for it, displays strongly non-autodual
features. (Thus, whereas the usual test categories all have a final
object, it is easy to see that a test category cannot possibly have an
initial object.) However, the homotopy and asphericity notions we then
came to develop make sense and are familiar indeed not only in the
``modelizing story'', but in any situation whatever which turns up,
giving rise to anything like a ``homology'' or a ``homotopy''
theory. To give one specific example, starting with an abelian
category \scrA, the corresponding homology theory is concerned with the
category of \scrA-valued complexes, say $\mathrm
K^\bullet(\scrA)$. The most basic notions here are the three homotopy
notions and two asphericity notions, namely: homotopy equivalence
between maps, homotopism, contractible objects, and quasi-isomorphisms
(= ``weak equivalences''), and acyclic (= ``aspheric'') objects. The
two first homotopy notions determine each other in the usual way, and
define the third, namely an object is contractible (or null-homotopic)
if{f} the map $0\to X$, or equivalently $X\to 0$, is a homotopism. On
the other hand, if we use the mapping-cylinder construction for a map,
the set of contractible objects determines the set of homotopisms, as
the maps whose mapping cylinder is contractible. Likewise, weak
equivalences determine aspheric objects, and conversely if the mapping
cylinder construction is given. The question that now comes to mind
immediately is whether the two sets of notions, the three pure
homotopy notions (determining each other), and the two ``asphericity
notions'' (determining each other too), mutually determine each
other, as in the non-commutative set-up we have worked in so far. We
can also remark that the functor
\[ \mathrm H^0 : \mathrm K^\bullet(\scrA)\to\scrA\]
visibly plays the part of the functor \piz{} in the non-commutative
set-up, it gives rise moreover to the $\mathrm H^i$ functors (any
integer $i$) by composing with the iterated shift functor (where the
shift of $X$ is just the mapping\pspage{295} cylinder of $0\to X$),
thus the set of aspheric objects formally from the functor $\mathrm
H^0$ (as the objects $X$ such that $\mathrm H^i(X)=0$ for all $i$),
much as in the non-commutative set-up the functor \piz{} for a
contractibility structure determines the latter, and hence also the
corresponding asphericity structure.

Thus, the question arises of formulating the basic structures, namely
contractibility and asphericity structures, in an autodual way,
applying both to the autodual situation just described, and to the
non-autodual one we have been working out in the notes -- and if
possible even, in all situations met with so far where a homology or
homotopy theory of some kind of other has turned up. Of course, it is
by no means sure a priori that we can do so, by keeping first nicely
apart the two sets of notions (contractibility and asphericity),
namely defining them separately, and then showing that a
contractibility structures determines an asphericity structures, and
is determined by the latter. Maybe we'll have to define from the
outset a richer kind of structure, where both ``pure'' homotopy
notions and asphericity notions are involved. Also, the familiar
generalization of mapping cylinders, namely homotopy fibers and
cofibers, and the corresponding long exact sequences, will evidently
play an important role in the structure to be described. Now, this
again ties in with the corresponding structure of a derivator, as
contemplated in section \ref{sec:69}, namely ``integration'' and
``cointegration'' of diagrams in a given category.

Definitely, this reflection is going to lead well beyond the scope of
the present part \ref{ch:IV}\scrcomment{I guess it is because this
  question is already growing in prominence that AG later decided to include
  the present section in part \ref{ch:V}!} -- it relates rather to
part \ref{it:71.C} of the working program envisioned by the end of May
(section \ref{sec:71}, p.\ \ref{p:207}--\ref{p:210}), and rather
belong to part \ref{ch:VI} of these notes, which presumably will
center around the notion of a derivator. There is however a more
technical question, and of more limited scope, which deserves some
thought and goes somewhat in the same direction, namely: how to define
(for a given contractibility or asphericity structure on $M$) an
\emph{induced} structure on a category $M_{/a}$, where $a$ is in $M$?
There is a little perplexity in my mind, even when $M$ is of the type
\Chat{} say, with $C$ a contractor and $a$ in $C$, taking the
canonical structure on \Chat{} -- because with the most evident choice
of an ``induced'' asphericity structure on $(\Chat)_{/a} \simeq
(C_{/a})\uphat$, namely the usual notion of aspheric objects, this
structure\pspage{296} will practically never be totally aspheric
(unless we take $a$ to be a final object of $C$), hence will not be
associated to a contractibility structure -- whereas we expect that
the contractibility structure of \Chat{} should induce one on
$\Chat_{/a}$. Presumably, the ``correct'' notions of induced
structure, in the case of asphericity structure and the corresponding
notion of weak equivalence, should be considerably stronger than the
one I just envisioned, and correspond to the intuition of
``\emph{fiberwise} homotopy types'' over the object $a$ (visualized as
a space-like object). A careful description of such induced structures
seems to me to be needed, and the natural place to be the present part
\ref{ch:IV} of the notes.

\namedlabel{q:88.4}{4)}\enspace A little reflection on semi-simplicial
homotopy notions (and their analogons when $\Simplex$ is replaced by a
general test category $\Delta$) seems needed, in order to situate the
following fact: ss homotopy notions, namely for ss objects in any
category $A$, behave well with respect to \emph{any} functor
\begin{equation}
  \label{eq:88.1}
  A\to B,\tag{1}
\end{equation}
without having to assume that this functor commutes with finite
products, whereas in the context of homotopy structures, when have a
functor between categories endowed say with homotopy interval
structures (for instance, with contractibility structures), such a
functor
\begin{equation}
  \label{eq:88.2}
  M\to N\tag{2}
\end{equation}
behaves well with respect to homotopy notions only, it would seem, if
we assume beforehand it commutes with finite products (plus, of
course, that it transforms a given generating family of homotopy
intervals of $M$ into homotopy intervals of $N$). In
case
\[ M = \bHom(\Delta\op,A), \quad N=\bHom(\Delta\op,B),\]
and \eqref{eq:88.2} comes from a functor \eqref{eq:88.1} which does
not commutes to finite products, neither does \eqref{eq:88.2} -- and
still \eqref{eq:88.2} is well-behaved with respect to ss homotopy
notions! It should be noted of course that the semi-simplicial
homotopy notions in $M$ can be defined, even without assuming that in
$A$ finite products exist, namely in situations when $M$ does not
admit finite products -- and hence, strictly speaking, the set-up of
section \ref{sec:51} does not apply. All this causes a slightly
awkward feeling, which I would like to clarify and see what's going
on. I suspect it should be simple enough to do it here and now.

First, assume that $A$ is stable under finite products, and under
direct sums (with small indexing set say -- for what we want to do
with $\Simplex$, finite direct sums even would be enough). We'll use
sums only with summands\pspage{297} equal to a chosen final object $e$
of $A$, in order to get a functor
\[ I\mapsto I_A : \Sets\to A,\]
where $I_A$ is the ``constant'' object of $A$ with value $I$, namely a
sum of $I$ copies of $e$ (sometimes also written $I\times e$). Using
this functor, we get a functor
\[\Delta\!\uphat = \bHom(\Delta\op,\Sets) \to M =
\bHom(\Delta\op,A),\]
which I denote by
\[K\mapsto K_A,\]
associating to any ss~set the corresponding ``constant'' (relative to
$A$) ss~object of $A$.\scrcomment{I think it's a bit unclear here
  whether the test category $\Delta$ is actually assumed to be
  $\Simplex$ after all\ldots} On the other hand, because $A$ admits
finite products, so does $M$, which enables us to make use of the
homotopy notions developed in sections \ref{sec:51} etc. Thus, if
\[\bI=(I,\delta_0,\delta_1)\]
is any interval in $\Delta\!\uphat$, considering the corresponding
``$A$-constant'' interval $\bI_A$ in $M$, we get homotopy notions in
$M$, which we may refer to as \emph{\bI-homotopy} (dropping the
subscript $A$). They can all be deduced from the elementary
\bI-homotopy between maps in $M$, which is expressed in the known way,
in terms of a map in $M$
\[ h : \bI_A \times X\to Y,\]
where $X$ and $Y$ are the source and target in $M$ of the two
considered maps, between which we want to find an elementary
\bI-homotopy $h$. This map $h$ decomposes componentwise into
\begin{equation}
  \label{eq:88.star}
  h_n : (I_n)_A \times X_n\to Y_n,\tag{*}
\end{equation}
and each $h_n$ can be interpreted, in view of the definition of $I_n$,
as a map
\begin{equation}
  \label{eq:88.starprime}
  h_n': I_n \to \Hom_A(X_n,Y_n),\tag{*'}
\end{equation}
provided we assume that taking products in $A$ is distributive with
respect to the sums we are taking, whence
\[ (I_n)_A \times X_n \simeq (I_n)_{X_n} = \text{direct sum of $I_n$
  copies of $X_n$.}\]
Now to give $h$, or equivalently a sequence of maps $h_n$ in
\eqref{eq:88.star} ``functorial in $n$ for variable $n$'', amounts to
giving a sequence of maps \eqref{eq:88.starprime}, satisfying a
corresponding compatibility relation for variable $n$. The point of
course is that (for any $I$ in $\Delta\!\uphat$ and $X,Y$ in $M$)
\emph{the set of data \eqref{eq:88.starprime} plus the compatibility
  condition make sense, formally, independently of any exactness
  assumptions on $A$}. Thus, it can be taken\pspage{298} as the formal
ingredient of a definition of ``elementary \bI-homotopy'' between two
maps in $M$, without any assumptions whatever on the category $A$ we
start with. The standard case is the one when $I=\Simplex_1$, the
``unit interval'', but never mind. The definition works just as well,
when $\Simplex$ is replaced by any (let's say small) category $\Delta$
whatever. On the other hand, it is immediate that for a functor $M\to
N$ as above, induced by a functor $A\to B$, for two maps in $M$, any
elementary homotopy between them gives rise to an elementary homotopy
of their images in $N$ -- which is just the well-known fact (in case
$\Delta=\Simplex$, $I=\Simplex_1$) that $M\to N$ is compatible with
simplicial homotopy notions.

In order to fit this into the general framework of section
\ref{sec:51}, let's remark that if $A$ is a full subcategory of a
category $A'$, then for a pair of maps in $M$, the elementary
\bI-homotopies between these are the same as when considering the
given maps as maps in $M'=\bHom(\Delta\op,A')$, in which $M$ is embedded
as a full subcategory. Now, any (small, say) category $A$ can be
embedded canonically into $A'=\Ahat$ as a full subcategory, and any
functor
\[f:A\to B\]
embeds in the corresponding functor
\[f_!:\Ahat\to\Bhat,\]
hence the functor
\[\varphi:M\to N\]
embeds in the corresponding functor
\[\varphi':M'\to N', \quad M'=\bHom(\Delta\op,\Ahat),
N'=\bHom(\Delta\op,\Bhat).\]
As \Ahat, \Bhat{} satisfy the required exactness properties, it
follows that the \bI-homotopy notions in $M',N'$ can be interpreted in
terms of the notions of section \ref{sec:51}, with respect to $\bI_A$
and $\bI_B$, defined now as (componentwise) constant objects of
\Ahat{} and \Bhat{} respectively. Still, $f_!$ commutes to finite
products only if $f$ does, so we are still left with explaining why
$M'\to N'$ is well-behaved with respect to \bI-homotopy
notions. Equivalently, we need only see this in the case of
\eqref{eq:88.2} $M\to N$, when $A$ and $B$ are supposed to have the
required exactness properties to allow for the interpretation given
above of the \bI-homotopy notions in terms of the formalism of section
\ref{sec:51}, and when moreover $f$ (as $f_!$ above)\pspage{299}
commutes with sums. This now is readily expressed by the relations
\begin{align*}
  &\varphi'(\bI_A) \simeq \bI_B , \\
  &\varphi'(\bI_A\times X) \tosim \varphi'(\bI_A)\times \varphi'(X),
\end{align*}
i.e., while $\varphi'$ does \emph{not} commute to finite products in
general, however it \emph{does} commute to the products which enter in
the description of elementary homotopies (as these products can be
expressed in terms of direct sums in $A,B$, and $f$ commutes to
these).

These reflections suggest that the notions of homotopy interval
structures and contractibility structures may be generalized, in a way
that the underlying category need no longer be stable under finite
products nor even admit a final object; and likewise, the notion of a
morphism of such structures may be generalized, without assuming that
the underlying functor should commute with finite products. The
thought that this kind of generalization may be needed had already
occurred before in these notes, in connection with the corresponding
situation about twenty five years ago, when the notion of a site was
developed. But at present, the extension doesn't seem urgent yet, and
I better stop here this long digression!

The remaining questions possibly to deal with in part \ref{ch:IV} are
all concerned with modelizers. I'll try to be brief!

\namedlabel{q:88.5}{5)}\enspace Consider an ``algebraic structure
type'', and the category $M$ of its set-theoretic realizations. I am
looking for a comprehensive set of sufficient conditions on $M$ to
ensure that $M$ is a ``canonical modelizer''. It seems natural to
assume beforehand that in $M$ (where at any rate small direct and
inverse limits must exist) internal $\bHom$'s exist, and more
generally, for $X,Y$ two objects over an object $S$ of $M$, that
$\bHom_S(X,Y)$ -- this implies that base change $S'\to S$ in $M$
commutes with small direct limits and a fortiori, that direct sums are
universal -- we may as well suppose them disjoint too. One feel quite
willing too to throw in the total $0$-connectedness assumption (cf.\
section \ref{sec:58}), and that every non-empty object has a section
over the final object. This preliminary set of conditions on an
algebraic structure species is of course highly unusual, however it is
satisfied for must ``elementary'' algebraic structures (by which I
mean $M\equeq \Ahat$ for some small category $A$), as well as for
$n$-stacks or \oo-Gr-stacks, for any $n$ between $0$ and \oo. The hope
now is, in terms of these assumptions, to give a necessary and
sufficient condition in order that a)\enspace the ``canonical''
homotopy structure on $M$ be a contractibility structure, and moreover
b)\enspace the latter structure be ``modelizing'', by which we mean
that the\pspage{300} associated \scrWoo-asphericity structure (\scrWoo{}
= usual weak equivalences) be modelizing, which will imply that for
\emph{any} basic localizer \scrW, the corresponding \scrW-asphericity
structure is modelizing.

Even if I don't look into this question now, it'll turn up soon enough
in a similar shape, when it comes to prove modelizing properties for
categories of stacks of various kinds. The best we could hope for
would be a statement in terms of the category structure of $M$ alone,
with no assumption that $M$ be defined in terms of an algebraic
structure type. If I try to formulate anything by way of wishful
testing conjecture, what comes to mind is: is it enough that there
should exist a separating contractible interval? So the first I would
try to get an idea, is to see how to make a counterexample to
this\ldots

\bigbreak
\presectionfill\ondate{3.7.}\par

\hangsection[Digression (continued) on left exactness properties of
$f_!$ \dots]{Digression
  \texorpdfstring{\textup(continued\textup)}{(continued)} on left
  exactness properties of \texorpdfstring{$f_!$}{f!}
  functors.}\label{sec:89}%
In connection with the left exactness properties of a $f_!$ functor,
considered three days ago (section \ref{sec:85}), I have been befallen
by some doubts whether any subcategory of \Cat{} containing the
subcategory of standard simplices is strictly generating. I wrote
there (p.\ \ref{p:284}) that as this is true for $\Simplex$ itself, it
``follows a fortiori'' for any subcategory $A$ of \Cat{} containing
$\Simplex$. Assuming $A$ to be full and denoting by $i$ the inclusion
functor, this is known to be equivalent (cf.\ remark \ref{rem:85.2}
same page) to $i^*:\Cat\to\Ahat$ being fully faithful, and in this
form, it doesn't look so obvious that when this is true for one full
subcategory, $A_0$ say, it should be true for any larger one $A$. This
thought had been lingering for a second while writing the ``a
fortiori'' and I then brushed it aside, because of the formulation of
being generating in terms of strict epimorphisms. Only the next day
did it occur to me that it is by no means clear that if a family of
maps $X_i\to X$ in a category $M$ is strictly epimorphic, any larger
family with same target $X$ should be ``a fortiori'' strictly
epimorphic too -- the ``a fortiori'' is known to apply only in the
case of the similar notions of epimorphic, or universally strictly
epimorphic, families of maps. After a little perplexity, I found the
situation was saved, in the case I was interested in, through the fact
that it was known from Giraud's article on
descent\scrcomment{\textcite{Giraud1964}} (Bull.\ Soc.\ Math.\ France,
Mémoire 2, 1964, prop.\ 2.5, p.\ 28) that $\Simplex$ and even the
smaller subcategory of simplices of dimension $\le 2$, is
even\pspage{301} generating by ``\emph{universally} strict
epimorphisms'', a notion which is stable under enlargement of the
family of maps, as recalled above. Thus, the statement made on p.\
\ref{p:284} does hold true. And I just checked today that, while this
stability property by enlargement is surely not always true for a
family of maps which is strictly epimorphic, however, it \emph{is}
true that if a full subcategory $A_0$ of a category $M$ is generating
by strict epimorphisms (or, as we'll say, is ``strictly generating''),
then so is any larger full subcategory $A$. This is seen by an easy
direct argument, in terms of the initial definition, as meaning that
for any object $X$ in $M$, the family of maps $a_i\to X$ with target
$X$ and source in the given subcategory ($A$ say) should be strictly
epimorphic. (For\scrcomment{\textcite{SGA4vol1}} the definition of
common variants of the notion of epimorphism, see the ``Glossaire'' at
the end of chapter 1, SGA~4, vol.~1.)

It occurred also to me that (as suspected in remark \ref{rem:85.3},
loc.\ cit.) the functor
\[i_!: \Simplexhat\to\Cat\]
coming from the inclusion functor $i:\Simplex\to\Cat$ is \emph{not}
left exact (for another reason though than first contemplated), namely
because \Cat{} is known \emph{not to be a topos} (for instance, an
epimorphism need not be \emph{strict} (or, what amounts here to the
same, \emph{effective}) -- as stated in the cited result of
Giraud). Indeed, we have the following
\begin{proposition}[which should belong to section \ref{sec:85}!]
  Let $M$ be a \scrU-category stable under small direct limits, $A$ a
  small full subcategory, $i:A\to M$ the inclusion functor, hence a
  functor
  \[i_!:\Ahat\to M.\]
  If $A$ is strictly generating \textup(i.e., $i^*:M\to\Ahat$ fully
  faithful\textup), then $i_!$ is left exact if{f} $M$ is a topos.
\end{proposition}

Indeed, the inclusion functor into \Ahat{} of $M'$, the essential
image of $i^*$ in \Ahat, admits a left adjoint ($i_!$ essentially). By
the criterion of Giraud, left exactness of this adjoint, or
equivalently of $i_!$, means that $M'$ is the category of sheaves on
$A$ for a suitable site structure on $A$, qed.
\begin{corollary}
  If $M$ is \emph{not} a topos, then $i_!$ does \emph{not} commute to
  fibered products in \Ahat{} of diagrams of the type
  \begin{equation}
    \label{eq:89.star}
    \begin{tabular}{@{}c@{}}
      \begin{tikzcd}[baseline=(O.base),column sep=tiny,row sep=small]
        b\ar[dr] & & c\ar[dl] \\ & |[alias=O]| F &
      \end{tikzcd},
    \end{tabular}\tag{*}
  \end{equation}
  with $b,c$ in $A$ and $F$ in \Ahat,\pspage{302} while it does
  commute to finite products, and to fibered products of any two
  objects of \Ahat{} over an object of $A$\kern1pt.
\end{corollary}

The ``while'' comes from prop.\ \ref{prop:85.1} and prop.\
\ref{prop:85.2} of section \ref{sec:85}, which imply too that, if $M$ is
strictly generating and whether or not $M$ is a topos, left exactness
of $i_!$ is equivalent with commutation to fibered products of the
diagrams \eqref{eq:89.star}. Hence the corollary.

This corollary answers also the perplexity raised in remark
\ref{rem:85.1} (p.\ \ref{p:283}), as to a hypothetical sharper version
of part \ref{it:85.prop1.b}, concerning fibered products. As
anticipated there, it turns out that this sharper version is not
valid, -- not without additional assumptions at any rate.

\hangsection[Review of questions (continued): 6) Existence of test
\dots]{Review of questions
  \texorpdfstring{\textup(continued\textup)}{(continued)}:
  \texorpdfstring{\textup{6)}}{6)} Existence of test functors and
  related questions. Digression on strictly generating
  subcategories.}\label{sec:90}%
After this digression on exactness properties of $f_!$ functors, let's
come back to the review of those questions not yet dealt with, which
seem more or less to belong to the present part \ref{ch:IV} of the
notes. We had stopped two days ago with the question \ref{q:88.5} of
finding some simple characterization of canonical modelizers,
comparable maybe in simplicity to the characterization we found for
test categories (in part \ref{ch:II}). This question may well turn out
to be related to the following one.

\namedlabel{q:90.6}{6)}\enspace This is the question of finding handy
existence theorems for test functors, whereas so far our attention to
test functors had been turned towards a thorough understanding of the
very notion of a test functor and its variants. I have the feeling
that, after the reflections of sections \ref{sec:78} and \ref{sec:86}
notably, the notion in itself is about understood now, so that time is
getting ripe for asking for existence theorems. As all modelizers we
have been meeting so far were associated to asphericity structures, it
seems reasonable to restrict to these, namely to the case of a given
modelizing asphericity structure
\[(M,M\suba),\]
and, if need be, even restrict to the case when this structure is
associated to a contractibility structure $M\subc$. We suppose given
moreover a test category $A$, which we may (if needed) assume to be
strict even, or even a contractor (i.e., the objects of $A$ in \Ahat{}
are moreover contractible, for the homotopy interval structure in
\Ahat{} defined by all intervals coming from $A$). The question then
is whether there exists a test functor
\[A\to M.\]
This (under the assumptions made) just reduces to the existence
of\pspage{303} a functor which be $M\suba$-\scrW-aspheric. Here,
\scrW{} is a given basic localizer, with respect to which we got an
asphericity structure. The most important case for us surely is the
one when $\scrW=\scrWoo$, namely usual weak equivalence. It is
immediate indeed that an $M\suba$-\scrW-aspheric functor is equally
aspheric for the corresponding $\scrW'$-asphericity structure of $M$,
for any basic localizer $\scrW'\supset\scrW$. Thus, if we get an
aspheric functor for \scrWoo, the finest basic localizer of all, we
get ipso facto an aspheric functor for any basic localizer
\scrW. (Note also that if an asphericity structure is modelizing for a
given \scrW, the corresponding $\scrW'$-structure is modelizing too,
for any $\scrW'\supset\scrW$; and the analogous fact holds for the
notion of a test category -- namely a \scrW-test category is also a
$\scrW'$-test category, and similarly for total asphericity of \Ahat{}
and hence for the condition of being a \emph{strict} test category.)

In case $M$ is even endowed with a contractibility structure, we will
be interested, more specifically still, in aspheric functors factoring
not only through $M\suba$, but even through $M\subc$:
\[ i:A\to M\subc,\]
while replacing the asphericity requirement on this functor, by the
stronger one that for any $x$ in $M\subc$, the object $i^*(x)$ in
\Ahat{} be \emph{contractible} (for the homotopy structure in \Ahat{}
defined by homotopy intervals coming from objects in $A$, say). In
other words, we are interested in the question of existence of
\emph{bimorphisms} of contractibility structures (in the sense of
section \ref{sec:86}) from $(M,M\subc)$ to $(\Ahat,\Ahatc)$. It may be
noted that in both cases (working with asphericity structures or with
the contractibility structures instead), in this existence question,
we may altogether forget $M$ itself, and consider it as an existence
question for functors from $A$ into either $M\suba$, or $M\subc$, with
the property that for any object $x$ in the target category $M\suba$
or $M\subc$, the object $i^*(x)$ in \Ahat{} be either aspheric, or
contractible. In the second case, we may even restrict $x$ to be in
any given subcategory $C$ of $M\subc$ generating the contractibility
structure -- and in the cases met with so far, we can find such a $C$
reduced to just one object $I$. In the case of asphericity structures,
the same holds when taking for $C$ a subcategory generating the
asphericity structure, provided however $C$ contains the image of $A$
by $i$ (which gives little hope to have $C$ restricted to just one
element!)

The interest of finding criteria for existence of \scrW-aspheric
or\pspage{304} or more stringently still of ``c-\emph{aspheric
  functor}'' (as we may call them) is rather evident, as it gives a
way, via $i^*$, for any homotopy type described in terms of a
``model'' $x$ in $M$, to find a corresponding model $i^*(x)$ in
\Ahat. The situation would be more satisfactory still if we could find
the test functor $i$ such that the corresponding functor
\begin{equation}
  \label{eq:90.1}
  i_!:\Ahat\to M\tag{1}
\end{equation}
be modelizing too (assuming $M$ to be stable under small direct
limits, so that $i_!$ is defined as the left adjoint of
\begin{equation}
  \label{eq:90.2}
  i^*:M\to\Ahat\quad\text{.)}\tag{2}
\end{equation}
In this case, for a homotopy type described by a model $K$ in \Ahat,
$i_!(K)$ gives a description of the same by a model in $M$.

Maybe we should remember though that even if we do not know about any
test functor from $A$ to $M$, still we always can find in three steps
a modelizing functor
\[M\to\Ahat,\]
namely a composition
\begin{equation}
  \label{eq:90.star}
  M \xrightarrow{j^*} \Bhat \xrightarrow{i_B}
  \Cat\xrightarrow{j_A=i_A^*} \Ahat,\tag{*}
\end{equation}
where $j:B\to M$ is an $M\suba$-aspheric functor from an auxiliary
small category, which we may assume to be a test category, by a mild
extra assumption on $M$ (cf.\ cor.\ \ref{cor:79.3} p.\
\ref{p:253}). The modelizing functor we thus get has the disadvantage
of not being left exact, whereas the looked-for functor $i^*$ commutes
to small inverse limits. Still, the composition \eqref{eq:90.star} is
pretty near to being left exact, it commutes to fibered products
(because $i_B$ does) which is the next best -- we can view it as a
left-exact functor from $M$ to $\Ahat_{/E}$, where $E$ is the image in
\Ahat{} of the final object of $M$ (assuming $e_M$ exists).

There is another advantage still of having a test functor $i:A\to M$,
rather than merely using \eqref{eq:90.star}, namely it allows us to
``enrich'' the category structure of $M$, in such a way as to get
``external $\Hom$'s'' of objects of $M$, with ``values in \Ahat'', by
defining, for $x,y$ in $M$, the object $\bHom(A)(x,y)$ of \Ahat{} as
\begin{equation}
  \label{eq:90.3}
  \bHom(A)(x,y) = \bigl\{ a\mapsto\Hom_M(i(a)\times x,y)\bigr\}.
  \tag{3}
\end{equation}
Such enriched structure, when $A=\Simplex$, plays an important part in
the second part of Quillen's treatment of homotopical algebra, under
the name of (semi-)simplicial categories, especially with the notion
of (semi-)simplicial model categories, which looks quite handy
indeed.\pspage{305} We should of course define composition of the
$\bHom(A)$'s, as required too in Quillen's set-up. This is done by
relating the $\bHom(A)$'s to the well-known internal $\bHom$'s in
$M\uphat$ -- which will show at the same time that for formula
\eqref{eq:90.3} to make sense, we do not really have to assume $M$ be
stable under binary products, as we can interpret the products
$i(a)\times x$ as being taken in $M\uphat$, as well as the $\Hom$, so
as to get
\[\Hom_{M\uphat}(i(a)\times x,y) \simeq \Hom_{M\uphat}(i(a),
\bHom_{M\uphat}(x,y)),\]
hence
\begin{equation}
  \label{eq:90.4}
  \bHom(A)(x,y)\simeq i^*(\bHom_{M\uphat}(x,y)),\tag{4}
\end{equation}
where $i^*$ in the right hand side is interpreted as a functor
\[ i^*:M\uphat\to\Ahat,\]
rather than $M\to\Ahat$. (I leave to the reader the task of enlarging
the basic universe, as need may be\ldots) As $i^*$ commutes to
products, the evident composition of the internal $\bHom$'s in
$M\uphat$ gives rise to the looked-for composition of the
$\bHom(A)$'s, with the required associativity properties. Of course,
in case $M$ is stable under binary products and $\bHom$'s, which
apparently is going to be the case in all modelizing situations, there
is no need in the interpretation \eqref{eq:90.4} to introduce the
prohibitively large $M\uphat$, and we can work in $M$ throughout.

There is an important relation though on the external $\bHom(A)$'s
which we would like to be true for a satisfactory formalism, namely
\begin{equation}
  \label{eq:90.5}
  \Gamma_\Ahat(\bHom(A)(x,y)) \fromsim \Hom_M(x,y),\tag{5}
\end{equation}
where $\Gamma_\Ahat$ just means $\Hom_\Ahat(e_\Ahat,\dots)$. This is
equivalent to the requirement
\begin{equation}
  \label{eq:90.6}
  i_!(e_\Ahat)\simeq e_M\tag{6}
\end{equation}
(assuming a final object $e_M$ in $M$ to exist), or equivalently
\begin{equation}
  \label{eq:90.7}
  i(e_A)\simeq e_M\tag{7}
\end{equation}
if we assume moreover $e_A$ to exist. Thus, it will be natural to ask
for test functors satisfying the extra condition \eqref{eq:90.5} or
\eqref{eq:90.6} -- and when trying to construct test functors in
various situations (even without being aware of constructing test
functors, as Mr~Jourdain\scrcomment{the reference is of course to
  Molière's \emph{comédie-ballet},
  \href{https://en.wikipedia.org/wiki/Le_Bourgeois_gentilhomme}{Le
    Bourgeois gentilhomme}} was ``doing prose without knowing
it''\ldots), the very first thing everybody has been doing
instinctively was to write down formula \eqref{eq:90.7},\pspage{306} I
would bet!

The motivation for wanting to find test functors being reasonably
clear by now, what kind of existence theorems may we hope for? When
$A$ is such a beautiful test category as $\Simplex$, $\Square$ or
$\Globe$, I would expect that for practically any $M$ endowed with a
modelizing contractibility structure say, under mild restrictions
(such as the exactness assumptions which are natural in the modelizing
story), there should exist a test functor indeed. What I feel less
definite about is whether it is reasonable to expect we can find $i$
even such that $i_!$ be modelizing too, in which case we would expect
of course that the pair of equivalences of categories
\begin{equation}
  \label{eq:90.8}
  \begin{tikzcd}[cramped,sep=small]
    \HotOf_A \ar[r,shift left] & \HotOf_M \ar[l,shift left]
  \end{tikzcd}\tag{8}
\end{equation}
defined in terms of $i_!$ and $i^*$ should be quasi-inverse to each
other, and the adjunction maps deduced from those between the functors
$i_!$ and $i^*$ themselves. This in turn is equivalent with the
adjunction morphism
\begin{equation}
  \label{eq:90.9}
  F\to i^*i_!(F)\tag{9}
\end{equation}
being a weak equivalence, for any object $F$ in \Ahat. A test functor
satisfying this exacting extra property merits a name of its own, we
may call it a \emph{perfect test functor} (or a \emph{perfect aspheric
  functor}, when not making any modelizing assumptions on $A$ or
$M$). Thus, the existence problem of finding test functors can be
sharpened to the one of finding perfect ones. Remember though that the
most familiar test functor of all (besides the geometric realization
functor $\Simplex\to\Spaces$, namely the inclusion
\begin{equation}
  \label{eq:90.10}
  i:\Simplex\to\Cat\tag{10}
\end{equation}
giving rise to the nerve functor (introduced for the first time, I
believe, in a Bourbaki talk of mine, on passage to quotient by a
preequivalence relation in the category of schemes\ldots), is
\emph{not} perfect. The most natural perfect test functor from
$\Simplex$ into \Cat, more generally from any weak test category $A$
into \Cat, is of course $i_A$ -- the functor indeed which has been
dominating the whole modelizing picture in our reflections from the
start. In the case of $A=\Simplex$, Thomason discovered another
perfect test functor, conceptually less simple surely, namely
$i_!\Sd^2$, where $\Sd$ is the ``barycentric subdivision functor''.
I suspect there must be an impressive bunch of perfect test functors
from $\Simplex$ with values in more or less any given modelizer, not
only the basic one -- and the question here is to get a clear picture
of\pspage{307} how to get them, and the same of course for test
functors which need not be perfect, including \eqref{eq:90.10}.

Next question then would be to see whether the existence theorems we
may get for $\Simplex$, or its siblings $\Square$ and $\Globe$ and the
like, still hold true for a more or less arbitrary test category, or
contractor. If so, this would be a very strong confirmation of the
feeling which has been prompting the reflections in part \ref{ch:II},
namely that for the purpose of having ``all-purpose''-models for
homotopy types (insofar as this is feasible), any strict test
category, or any contractor at any rate, is just as good as simplices
or cubes, which people have kept working with for the last twenty-five
years. If not, it will be quite interesting indeed to come a grasp of
what the relevant extra features of $\Simplex$ and the like are, and
how restrictive they are.

I doubt I will dive into these questions, still less come to a clear
picture, in the present part of the notes. Still, before leaving the
topic now, I would like to write down some hints I came upon while
doing my scratchwork on homotopy properties of \Ahat{}
categories. When looking for functors
\[A\to M\]
having some specified properties (such as being a test functor, or a
perfect one, etc.), we may view this question as meaning that we are
looking for an object with specified properties in the category
\[M^A =\bHom(A,M).\]
Presumably, this category is endowed with an asphericity or
contractibility structure if $M$ is (as we assume), presumably even a
modelizing one. This reminds me that as far as the notion of weak
equivalence goes, there may be even several non-equivalent ways of
finding such structure on $M^A$, one being modelizing, whereas
another, more useful one in some respects, is not. Thus, if $M$ is of
the type \Bhat, we may rewrite $M^A$ as
\[ M^A \equeq (A\op\times B)\uphat \equeq P\uphat;\]
hence we get the weak equivalence notion coming from $P\uphat$,
disregarding its product structure, which is modelizing indeed if $B$
is a test category and $A$ aspheric, hence $A\op\times B$ a test
category. The structure which should be of more relevance though for
our purpose should be a considerably finer one (namely with a smaller
set of weak equivalences), which we may\pspage{308} visualize best
maybe by writing
\[M^A=(\bHom(A\op,M\op))\op,\]
i.e., interpreting the dual of $M^A$ as the category of $A$-objects of
$M\op$, for instance (if $A=\Simplex$) as the dual of the category of
ss~objects of $M\op$. Now, Quillen has given handy conditions, in
terms of projectives of $M\op=N$, namely in terms of injectives of
$M$, for the category $\bHom(\Simplexop,N)$ of ss~objects of $N$ to
be a closed model category -- hence the dual category $M^A$ will turn
out as a closed model category too, under suitable conditions
involving existence of injective objects in $M$. These conditions are
satisfied for instance when $M$ is a topos, and notably when $M$ is of
the type \Bhat{} -- quite an interesting particular case indeed!
Assuming that $M$ is stable under both types of limits, so is $M^A$,
hence there is an initial and final object, and according to Quillen's
factorization axiom, the map from the former to the latter can be
factored through an object
\[\text{$F$ in $M^A$,}\quad\text{i.e.,}\quad F:A\to M,\]
which is \emph{cofibering}, and such that $F\to e$ is a \emph{trivial
  fibration}. The idea is that these conditions mean more or less, at
any rate imply, that $F$ is a test functor.

I hit upon this ``way out'' while trying to construct test functors
from $\Simplex$ to any elementary modelizer \Bhat, in order to try and
check that \Bhat{} is a (semi)simplicial model category in the sense
of Quillen. The intuitive idea of constructing inductively the
components $F_n$ of $F$ was simple enough, still I got stuck in some
messiness and did not try to push through this way, all the less as
this naive approach had no chance of generalizing to the case of a
more or less general test category $A$. Of course, for the time being
Quillen's theorems, about certain categories $\bHom(A\op,N)$ being
closed model categories, is equally restricted to the case when
$A=\Simplex$, which looks as usual like a rather arbitrary
assumption. Thus, to ``test'' whether the feeling about \emph{any}
test category more or less being ``just as good'' as $\Simplex$, a
second point would be to see whether Quillen's theorems extend, which
presumably is going to be very close to the first point I raised.

I take this occasion to raise a third point -- where there is no
reason to restrict to an $M$ which is modelizing (neither was there
such reason before, when phrasing everything in terms of aspheric or
c-aspheric functors, rather than test functors\ldots). Namely,
assuming as\pspage{309} above that for any test category or contractor
$A$, $M^A$ or $M^{A\op}$ can be endowed with an asphericity structure
or a closed model structure, or at any rate with a set of weak
equivalences, hence a localization or corresponding ``homotopy
category''. What one would expect now is that up to (canonical)
equivalence, the latter does not depend upon the choice of $A$, and
hence is the same for arbitrary $A$ as when using $\Simplex$, i.e.,
simplices. This should be true at any rate for $M^{A\op}$ and when $M$
is a topos -- which means that \emph{Illusie's derived category
  $\mathrm D_\bullet(X)$ of the category of semisimplicial sheaves on
  a topos $X$, could be constructed by using, instead of ss~objects,
  $A$-objects of the category of sheaves on $X$, where $A$ is any test
  category.} This should be one of the main points to settle in part
\ref{ch:VII} of the notes.

There is a slight discrepancy though between the first point, about
existence of test functors with values in $M$, depending on a given
asphericity or contractibility structure of $M$, and the second and
the third, which seem to depend only on the category structure of
$M$. This is further evidence that the set of questions raised here is
still far from being clear in my mind yet. Stating them now, however
confusingly, is a first step towards clarification!

\bigbreak
\presectionfill\ondate{4.7.}\par

Just still two comments about the existence questions for test
functors, before going over to the next questions in our present
review. One is that for given $A$, to prove that for rather general
modelizing $M$ there exists a test functor from $A$ to $M$, we are
reduced to the case when $M$ is of the type \Bhat, where $B$ is a test
category -- namely, it is enough to take a $B$ such that there exists
a test functor $B\to M$. If we can even find a perfect test functor
from some $B$ to $M$, then likewise the existence question for perfect
test functors from $A$ to $M$ is reduced to the case when
$M=\Bhat$. These comments may be useful for applying to the situation
Quillen's model theory, as envisioned on the previous page -- as his
criteria for $\bHom(\Simplexop,N)$ to be a closed model category
apply when $N$ is the dual of a topos, for instance the dual of
\Bhat. The second comment is about Thomason's result concerning the
standard inclusion
\[i:\Simplex\hookrightarrow \Cat,\]
which can be expressed by saying that, although $i$ itself is not a
\emph{perfect} test functor, however, for any integer $n\ge2$, the
composition\pspage{310}
\[i_n=i_!\Sd^n \alpha: \Simplex\to\Simplexhat\to\Simplexhat\to\Cat\]
is a perfect test functor, where $\alpha:\Simplex\to\Simplexhat$ is
the canonical inclusion. It is tempting to surmise that this result is
not special to $i$ alone, but that it holds for a large class, if not
all, test functors from $\Simplex$ to asphericity or contractibility
modelizers. Here, $\Sd^n$ denotes the $n$'th iterate of the
barycentric subdivision functor $\Sd$ (following now the notation in
Thomason's paper, which presumably is standard, while I have been
using ``Bar'' in part \ref{ch:II} of the notes). Presumably, functors
analogous to $\Sd$ can be defined in any elementary modelizer \Ahat,
as suggested by the natural constructions arising in connection with
the factorization property for a closed model structure on
\Ahat. Thus, possibly there is a general method in view for deducing
perfect test functors from ordinary ones. However, it definitely seems
to me that the natural place for these existence questions is in part
\ref{ch:V} of the notes, as they seem intimately related to an
understanding of the homotopy structures of elementary modelizers, and
more specifically to the closed model structures to which such
modelizers give rise in various ways.

\starsbreak

Before proceeding, I would like to state still another afterthought to
the reflections of section \ref{sec:89}, about strictly generating
subcategories of a category $M$. I recall that a family of arrows in
$M$ with same target $X$
\[u_i:X_i\to X\]
is called \emph{strictly epimorphic}, if for every object $Y$ of $M$,
the corresponding map
\begin{equation}
  \label{eq:90.starbis}
  \Hom(X,Y) \to \prod_i\Hom(X_i,Y)\tag{*}
\end{equation}
is injective (which is expressed by saying that the family $(u_i)$ is
\emph{epimorphic}), and if \emph{moreover} the following obviously
necessary condition for an element $(f_i)$ of the product set of
\eqref{eq:90.starbis} to be in the image of the map \eqref{eq:90.starbis},
is also sufficient:
\begin{description}
\item[\namedlabel{cond:90.Comp}{(Comp)}]
  For any two indices $i,j$ (possibly equal) and any commutative
  square
  \[\begin{tikzcd}[sep=tiny]
    & T\ar[dl,"v_i"']\ar[dr,"v_j"] & \\
    X_i\ar[dr,"u_i"'] & & X_j\ar[dl,"u_j"] \\
    & X &
  \end{tikzcd}\]
  in $M$, the relation $f_iv_i=f_jv_j$ holds.
\end{description}
It is immediate that the condition for $(u_i)$ to be strictly
epimorphic\pspage{311} depends only on the \emph{sieve} $X_0$ of $X$
in $M$ (namely, the subobject of $X$, viewed as an object of
$M\uphat$) generated by the $u_i$'s -- we'll say also that this
\emph{sieve is strictly epimorphic}. One should beware that this does
not mean of course that $X_0\to X$ is epimorphic in $M\uphat$ (which
would imply $X_0=X$, i.e., that one of the $u_i$'s admits a section,
i.e., a right inverse); nor is it true that if a sieve is strictly
epimorphic, a large one should be so too -- which means that when
adding more arrows to a strictly epimorphic family, the family need
not stay str.\ ep.

We'll say that the family $(u_i)$ is \emph{universally strictly
  epimorphic} if it is strictly epimorphic, i.e., the corresponding
sieve $X_0$ is, and if the latter remains so by arbitrary base change
$X'\to X$ in $M$, i.e., if the corresponding sieve $X'_0$ of $X'$ is
strictly epimorphic too. If the fibered products
\[X'_i = X_i\times_X X'\]
exist in $M$, this condition also means that the corresponding family
of maps
\[ u'_i:X'_i\to X'\]
is strictly epimorphic. The condition that $(u_i)$ be univ.\ str.\
epimorphic again depends only on the generated sieve, it is moreover
\emph{stable under base change}, and equally \emph{stable under adding
  new arrows}, i.e., replacing a sieve in $X$ by a larger one.

It should be noted that if the fibered products $X_i\times_X X_j$
exist in $M$, then the compatibility condition \ref{cond:90.Comp}
above is equivalent to the one obtained by restricting to
\[T = X_i\times_X X_j,\]
with $v_i$ and $v_j$ the two projections.

Assume the indices $i$ are objects of a category $I$, and the $X_i$
are the values of a functor
\[I\to M,\]
and that the family of arrows $(u_i)$ turns $X$ into the direct limit
in $M$ of the $X_i$:
\[X=\varinjlim_i X_i,\]
then it follows immediately that the family $(u_j)$ is strictly
epimorphic.

After these terminological preliminaries, we're ready to give
the\pspage{312} following useful statement, which is lacking in
SGA~4\scrcomment{\textcite{SGA4vol1}} 
Chap.~I (compare loc.\ cit.\ prop.\ 7.2, page 47, giving part of the
story):
\begin{proposition}
  Let $M$ be a \scrU-category, $A$ a small full subcategory, $i:A\to
  M$ the inclusion functor, hence a functor
  \[i^*:M\to \Ahat.\]
  For any object $X$ of $M$, we consider the family $F_X$ of all
  arrows in $M$ with target $X$, source in $A$\kern1pt. The following
  conditions are equivalent:
  \begin{enumerate}[label=(\roman*),font=\normalfont]
  \item\label{it:90.i}
    For any $X$ in $M$, the family $F_X$ is strictly epimorphic.
  \item\label{it:90.ii}
    For any $X$ in $M$, the family $F_X$ is universally strictly
    epimorphic.
  \item\label{it:90.iii}
    For any $X$ in $M$, $F_X$ turns $X$ into a direct limit of the
    composition functor $A_{/X}\to A\to M$, i.e.,
    \[X \fromsim \varinjlim_{A_{/X}} a.\]
  \item\label{it:90.iv}
    The functor $i^*$ is fully faithful.
  \end{enumerate}
\end{proposition}

Proof left to the reader (who may consult loc.\ cit.\ for
\ref{it:90.i} $\Rightarrow$ \ref{it:90.iii} $\Leftrightarrow$
\ref{it:90.iv}, so that only \ref{it:90.i} $\Rightarrow$
\ref{it:90.ii} is left to prove).
\begin{definition}
  When the equivalent conditions above are satisfied, we'll say that
  $A$ is a \emph{strictly generating} subcategory of $M$.
\end{definition}

NB\enspace The notion makes sense too without assuming $A$ to be
small, nor $M$ to be a \scrU-category, by passing to a larger universe
(it is immediate for \ref{it:90.i} or \ref{it:90.iii} that the
condition does not depend on the choice of the universe.)
\begin{corollary}
  If $A$ is strictly generating in $M$, then so is any larger full
  subcategory $B$.
\end{corollary}

This is clear by criterion \ref{it:90.ii} (whereas it isn't by any one
of the other two criteria!).

\hangsection[Review of questions (continued): 7) Homotopy types of
\dots]{Review of questions \texorpdfstring{\textup(continued\textup):
    \textup{7)}}{(continued): 7)} Homotopy types of finite type,
  \texorpdfstring{\textup{8)}}{8)} test categories with boundary
  operations, \texorpdfstring{\textup{9)}}{9)}
  miscellaneous.}\label{sec:91}%
I see three more questions to review -- presumably they will be a lot
shorter than the last!

\namedlabel{q:91.7}{7)}\enspace\textbf{Description of homotopy types
  ``of finite type'',} in terms of an elementary modelizer \Ahat. In
terms of the modelizer \Spaces, a natural finiteness condition on a
homotopy type is that it may be described (up\pspage{313} to
isomorphism) as the homotopy type of a \emph{space admitting a finite
  triangulation}. In terms of ss~sets, i.e., of the modelizer
\Simplexhat, the natural finiteness condition, suggested by the
algebraic formalism, is that the homotopy type be isomorphic to one
defined by an object \emph{``of finite presentation'' in} \Simplexhat,
namely one which is a direct limit of a \emph{finite} diagram in
\Simplexhat, made up with simplices, i.e., coming from a diagram in
$\Simplex$. It is clear that the first finiteness condition implies
the second, by using a total order on the set of vertices of the
triangulation. The converse shouldn't be hard either, using an
induction argument on the number of simplices occurring in the
diagram, and using the fact that any quotient object in \Simplexhat{}
of a simplex is again a simplex, hence also a subobject of a simplex
is a union of subsimplices; from this should follow by induction that
the geometric realization of a ss~set of finite presentation is
endowed with a natural compact \emph{piecewise linear structure}, and
hence can be finitely triangulated. Presumably, one can even find a
canonical triangulation, using twofold barycentric subdivision $\Sd^2$
(again!) on any simplex. All this is surely standard knowledge, and I
don't feel like diving into technicalities on this matter, unless I am
forced to.

If we start with an arbitrary test category $A$, the notion of an
object of finite presentation in \Ahat{} still makes sense. Indeed, in
any category $M$, stable under filtering small direct limits, we may
define objects of finite presentation as those for which the
corresponding \emph{covariant} functor
\[Y\mapsto \Hom(X,Y)\]
commutes to filtering direct limits. If $M$ is stable under some type
of finite direct limits, say under any finite direct limits, then so
is the full subcategory $M\subfp$ of objects of finite presentation of
$M$. In the case when $M=\Ahat$, $A$ any small category, it is obvious
that objects of $A$, and hence finite direct limits of such, are of
finite presentation, and it is not difficult to show that the converse
equally holds. As a matter of fact $\Ahatfp$ \emph{can be viewed as
  the solution of the $2$-universal problem of ``adding finite direct
  limits to $A$''}.

Coming back to the case when $A$ is a test category, and hence \Ahat{}
is modelizing, one may ask for conditions upon $A$ which ensure that
the homotopy types of finite presentation are exactly those isomorphic
to the homotopy types defined by objects of \Ahatfp. One expects that
some stringent\pspage{314} extra condition is needed on $A$ to ensure
this. To see this, let's take a finite group $G$, and an aspherical
topological space $E_G$ upon which $G$ operates freely, with quotient
$B_G$, a classifying space of $G$. If $G\ne1$, the homotopy type of
$B_G$ isn't of finite type, because $B_G$ has non-vanishing cohomology
groups in arbitrarily high dimensions, as well known. We could
transpose the following construction in either \Cat{} or \Simplexhat{}
say, but we may as well work in the modelizer $M=\Spaces$, and take
any small full subcategory $A$ of $M$, containing the unit interval
but not the empty space, and stable under finite products -- which
implies that $A$ is a strict test category. We'll take $A$ large
enough to contain $E_G$, and small enough to be made up with aspheric
spaces, hence the inclusion functor
\[i: A\to\Spaces\]
is a test functor. Now consider the quotient object \emph{in} \Ahat
\[F = E_G / G,\]
i.e., the presheaf on $A$
\[ F : T\mapsto \Hom(T,E_G)/G \simeq \Hom(T,B_G),\]
where the last isomorphism comes from the fact that the object $T$ of
$A$ is aspheric and hence $1$-connected. Thus, we get
\[F\simeq i^*(B_G),\]
hence the homotopy type defined by $F$ is the homotopy type of $B_G$,
which is not of finite type, despite the fact that $F$ is of finite
presentation.

In order to ensure that the homotopy type defined by any object in
\Ahatfp{} be of finite type, it may be useful perhaps to make on $A$
the assumption that any quotient in \Ahat{} of an object in $A$ is
isomorphic to an object in $A$, and that the set of all subobjects in
\Ahat{} of an object $a$ in $A$ (i.e., the set of all sieves on $a$)
is finite -- possibly too that for any two objects $a$ and $b$ of $A$,
$\Hom(a,b)$ is finite. As for the opposite inclusion, namely that any
homotopy type of finite type can be described by an object of \Ahatfp,
this would follow from the existence of a \emph{perfect} test functor
from $\Simplex$ into \Ahat, factoring through\pspage{315}
\Ahatfp. Thus, the present question about finiteness conditions, seems
to be related (possibly) to the previous one about existence of
various types of test functors.

The condition for a strict test category $A$ we are looking at is
surely satisfied, besides $\Simplex$, by the cubical and hemispherical
test categories $\Square$ and $\Globe$, and surely also by any finite
products of these. I add this comment, of course, in order to ``push
through'' the point that not any more with respect to finiteness
conditions on homotopy types, than (presumably at least, for the time
being)in any other essential respect concerning the ability for
expressing basic situations and facts in homotopy or cohomology
theory, the category of simplices stands singled out by itself from
all other test categories. Nor does it seem that the ``trinity''
\[ \Simplex, \Square, \Globe\]
has this property, with the only exception so far, possibly, of the
Dold-Puppe theorem (as no other test category except these is known to
me for which a Dold-Puppe theorem in its strict form holds true
(compare reflections section \ref{sec:71})).

\namedlabel{q:91.8}{8)}\enspace
I could make the same point in favor of more general test categories
than the trinity above, when it comes to the existence of an algorithm
for computing homology and cohomology groups, using suitable
\emph{boundary operators}. What is meant by these is clear of course
for the three types of complexes, but then it extends in an obvious
way to multicomplexes too -- which means that for the test categories
deduced from the trinity by taking finite products, we still get an
algorithm for cohomology via boundary operators. Of course, for any
test category $A$, using a test functor (if we can find one) from one
of the three above (say) into \Ahat{} will allow us to reduce
``computation'' of homology and cohomology invariants in terms of a
model in \Ahat, to the case of the corresponding type of complexes --
hence again an algorithm (similar to the familiar one of computing the
cohomology of an object of \Cat{} semisimplicially, via the
nerve). But this is cheating of course! The question I want to raise
here is about existence of ``boundary operations'' \emph{in} $A$,
similar to the familiar ones used for the three basic types of
complexes, and allowing to compute the homology and cohomology groups
of an object of \Ahat{} in the usual way, involving suitable
\emph{signs} $+$ or $-$ associated to the various boundary
operations.\pspage{316} It shouldn't be hard, I feel, to pin down
exactly what is needed for getting such a formalism. The intuitive
idea behind it (suggested by the example of standard complexes and
multicomplexes) is that \emph{such a formalism should be associated to
  cellular decompositions of $n$-cells for variable $n$}, such that
the interior of each $n$-cell should be an open cell of the
subdivision. There may of course be several $n$-cells for the same
$n$, which are not combinatorially isomorphic. When trying to express
this idea in a precise way, we are led to assume, as an extra
structure on the would-be test category $A$, a functor
\begin{equation}
  \label{eq:91.1}
  i:A\to\Ord\tag{1}
\end{equation}
of $A$ into the category of ordered sets, such that for any $a$ in
$A$, $i(a)$ be a \emph{finite} ordered set, whose geometric
realization (cf.\ section \ref{sec:22}) is an $n$-cell for suitable
$n\eqdef\dim(a)$. We assume moreover that $i(a)$ has a largest element
$e(a)$, and that the geometric realization of $i(a)\setminus\{e(a)\}$
is the bounding $(n-1)$-sphere of the $n$-cell $\abs{i(a)}$:
\begin{equation}
  \label{eq:91.2}
  \abs{i(a)^*} \simeq \mathrm S^{n-1}, \quad\text{where}\quad
  i(a)^* \eqdef i(a)\setminus\{e(a)\},\tag{2}
\end{equation}
which will imply the precedent condition, namely
\begin{equation}
  \label{eq:91.3}
  \abs{i(a)}\simeq \mathrm B^n,\tag{3}
\end{equation}
as $\mathrm B^n$ can be identified to the cone over $\mathrm
S^{n-1}$. As another condition, we need that
\medbreak
\noindent(\namedlabel{eq:91.4}{4})\hfill%
\parbox[t]{0.9\textwidth}{For any $a$ in $A$ and $x\in i(a)$, there
  exists $b$ in $A$ and an isomorphism
  \[i(b)\tosim i(a)_{/x} \eqdef\set[\big]{y\in i(a)}{y\le x}
  \hookrightarrow i(a),\]
  induced by a map
  \[b\to a\]
  in $A$.}\par
\medbreak
\noindent It is enough to make this assumption for $x$ of codimension
$1$ in $i(a)$, which will imply that it is true for any $x$. It seems
reasonable on the other hand to make the assumption that for a given
$x$, the object $b$ in \eqref{eq:91.4}, viewed as an object of
$A_{/a}$, is determined up to a unique isomorphism, we may call it
$a_x$, and $\partial_x$ the canonical map of $b$ into $a$
\begin{equation}
  \label{eq:91.5}
  \partial_x:a_x\to a \qquad
  \begin{tabular}[t]{@{}c@{}}
    ($x\in i(a)$, of codim.\ $1$ in $i(a)$, \\
    i.e., $\dim(x)=\dim(a)-1$).
  \end{tabular}
  \tag{5}
\end{equation}
As an extra structure, we need for any $a$ in $A$
\begin{equation}
  \label{eq:91.6}
  \omega_a,\quad\text{an \emph{orientation} of the $n$-cell
    $i(a)$}\quad(n=\dim(a)).\tag{6}
\end{equation}
This allows us, for any $x$ as in \eqref{eq:91.5}, to define a
signature
\begin{equation}
  \label{eq:91.7}
  \varepsilon(x)\quad\text{or}\quad \varepsilon_a(x)\in\{+1,-1\},\tag{7}
\end{equation}
which\pspage{317} will be $+1$ or $-1$, depending on whether in the
inclusion
\[\abs{\partial_x}:\abs{i(a_x)}\to\abs{i(a)},\]
the orientation $\omega(a_x)$ is induced ``à la Stokes'' by the
orientation $\omega(a)$ of the ambient $n$-cell, or not. Having the
boundary operations \eqref{eq:91.5} with their signatures
\eqref{eq:91.7}, and the decomposition
\begin{equation}
  \label{eq:91.8}
  A = \coprod_{n\ge0} A_n,\quad\text{where $A_n=\set[\big]{a\in\Ob A}{\dim
      a=n}$,}\tag{8}
\end{equation}
we get in the usual way, for any contravariant functor $K_\bullet$
from $A$ with values in an additive category, a corresponding chain
complex in this category, with components
\begin{equation}
  \label{eq:91.9}
  K_n = \coprod_{a\in A_n} K_\bullet(a),\tag{9}
\end{equation}
and boundary operators defined in the usual way via \eqref{eq:91.5}
and \eqref{eq:91.7}. Applying this to the case of the category \Ab{}
of abelian groups, or to its dual, and to the abelianization of an
object $X$ of \Ahat, we obtain a tentative way for computing the
homology and cohomology groups of the homotopy type of $A_{/X}$, and
similarly for any system of twisted coefficients on $X$. The question
which arises here is to write down a set of natural extra conditions
on the data, which will ensure that we do get a canonical isomorphism
between the ``homology'' and ``cohomology'' groups thus constructed,
and the usual homology and cohomology invariants of the object
$A_{/X}$ of \Cat. Moreover, we would like too to have conditions to
ensure that $A$ is a test category, or even a strict one.

One difficulty here, if one really wants a test category and not just
a weak one (which may not be without any problem either), is that
presumably for this, we'll need suitable \emph{degeneracy operations},
which may well turn out a very exacting condition indeed! The
skeptical reader may wonder, as I am just doing myself, whether there
will be any example within the set-up I propose, which does not reduce
to a finite product of test categories in our trinity.\pspage{318}

I just spent a while trying to find some convincing example, by using
a suitable full subcategory $A_0$ of \Ord, made up with finite sets
satisfying the assumption \eqref{eq:91.2} above for $i(a)$, under some
additional assumption on $A_0$ such as stability under finite products
and under passage from $a$ to an object $a_{/x}$, and that $A_0$
contain the ordered set
\[I =
\begin{tikzcd}[baseline=(O.base),cramped,row sep=-3pt,column sep=small]
  \bullet\ar[dr] & \\ & |[alias=O]| \bullet \\ \bullet\ar[ur] &
\end{tikzcd},\]
whose geometric realization is the segment $\mathrm B^1$ with its
usual cellular decomposition. In terms of $A_0$ and introducing
``orientations'' of object of $A_0$, the idea was to define another
category $A$ (of pairs $(a,\omega)$, with $a$ in $A_0$ and $\omega$ an
orientation of $\abs a$), \emph{stable under finite products} so that
\Ahat{} is totally aspheric, together with a functor
$i:A\to\Ord\hookrightarrow\Cat$ such that $i^*(I)$ should be
representable, and hence furnish the homotopy interval needed to
ensure that $A$ is a test category. The first idea that comes to mind,
namely define a map from $(a,\omega)$ to $(a',\omega')$ as merely a
map from $a$ to $a'$ in $A_0$, is nonsense unfortunately, as in the
data \eqref{eq:91.6}, the orientations will not be stable under
isomorphisms, a condition which I forgot to state before, and which is
visibly needed in order to be able to define the differential between
the $K_n$'s. If we try to define $A$ taking into account this
compatibility condition, we loose existence of products, anyhow
$i^*(I)$ isn't representable anymore, so why should it be aspheric
over the final object, so why should the functor $i$ be a test functor?

The difficulty I find in carrying through any explicit example for a
``test category with boundary operations'', except those which stem
from our trinity, is rather intriguing I feel. The question is whether
maybe in this direction, one might get at an intrinsic description of
the trinity, in terms of the rather natural structure species of a
``test category with boundary operations''. This is the second
instance where the thought arises that the three standard test
categories $\Simplex$, $\Square$ and $\Globe$ may be distinguished in
some respects -- the first instance was in relation to the Dold-Puppe
theorem.

\namedlabel{q:91.9}{9)}\enspace\textbf{Miscellaneous residual
  questions from part \ref{ch:II}.} One of these was about the category
$\Simplexf$ of simplices without degeneracies being a weak test
category (cf.\ section \ref{sec:43}) -- while it is definitely
\emph{not} a test category. It seems worth while to write down a proof
for this, maybe too for the analogous statements for
$\Square^{\mathrm f}$ and $\Globe^{\mathrm f}$. This reminds
me\pspage{319} too that I never got around to introducing formally the
hemispherical test category, which presumably will be very useful when
it comes to studying stacks -- this too could be done in part
\ref{ch:IV}, as well as proving of course that $\Globe$ is a strict
test category indeed, or better still, a contractor. It may be fun too
constructing test functors from any one of the three basic test
categories in the trinity, to the category of complexes defined by the
two others -- six cases altogether to consider! But as I am not in the
process of writing the ``Elements d'Algèbre
Homotopique'', maybe I will skip this!

In the same section \ref{sec:43}, I raised the question as to whether
the ordered set of all non-empty finite subsets of a given infinite
set, viewed as a category in the usual way for an ordered set, was a
weak test category (on page \ref{p:78} it was seen not to be a test
category, and it is immediate then that it is not totally aspheric
either). One interesting application, as noticed there, would be to
the effect that \Ord, the subcategory of \Cat{} defined by ordered
sets, is a modelizer (for the induced notion of weak
equivalences). Now the question arises moreover whether this
modelizing structure comes from an asphericity structure, or even from
a contractibility structure -- and the same question arises in the
more general situation described in the proposition of page
\ref{p:74}.

A last question along these lines I would like to clear up, is the
relation of total \scrW-asphericity for an asphericity structure, for
variable \scrW, when $\scrW\subset\scrW'$. Assuming the localizers
satisfy the condition \ref{loc:4}, is it true that total
\scrW-asphericity is equivalent to total $\scrW'$-asphericity -- or
equivalently, is it equivalent to total $0$-connectedness?

\bigbreak

\presectionfill\ondate{5.7.}\pspage{320}\par

\hangsection[Short range working program, and an afterthought on
\dots]{Short range working program, and an afterthought on
  abelianization of homotopy types: a handful of questions around the
  Whitehead and Dold-Puppe theorems.}\label{sec:92}%
The review on ``pending questions and topics'' related to part
\ref{ch:IV} of the notes has taken pretty much longer than
expected. It was quite useful though, to get a clearer view of what
those questions are about, and to get a feeling for what to include
and develop, and where. As I do not intend to spend my life on the
task, not even one year, it is becoming clear that I am not going to
get the whole picture of all the questions touched -- and some
presumably I am going to leave just aside, as they do not seem
indispensable for a comprehensive overall picture of what I'm
after. This seems to me to be the case for the questions \ref{q:91.7}
and \ref{q:91.8}, concerned with finiteness conditions for homotopy
types in terms of models, and with test categories with boundary
operations. At the opposite side, it seems that the questions
\ref{q:88.2}, part of \ref{q:88.3}, and \ref{q:91.9}, about the notion
of contractor, induced asphericity and contractibility structures on a
category $M_{/a}$, and ``miscellaneous'' left-overs from part
\ref{ch:II}, should be dealt with in part \ref{ch:IV} -- whose end now
is in sight after all!  On the other hand, questions \ref{q:88.1},
another part of \ref{q:88.3}, \ref{q:88.5} and \ref{q:90.6}, about
morphisms of asphericity structures and related problems, about an
autodual treatment of asphericity and contractibility notions, about a
handy criterion for canonical modelizers, and about existence theorems
for various kinds of test functors or aspheric functors, while I feel
that I should come at least to a considerably clear understanding of
these matters than now, the adequate place for developing such
reflection is definitely \emph{not} in the present part \ref{ch:IV},
but belong to one or the other of the three parts still ahead in our
overall reflection on the modelizing story.

During our review, we came a number of times upon situations when the
question arose as to whether one point I like to make, namely that a
more or less arbitrary (strict) test category ``is just as good'' as
the sacrosanct test category $\Simplex$, or its twin brothers
$\Square$ and $\Globe$, is a valid one or not. I would like to list
here these situations, with a view of coming back to it later:
\begin{enumerate}[label=\alph*)]
\item\label{it:92.a}
  Existence theorems for test functors (cf.\ section \ref{sec:90}).
\item\label{it:92.b}
  \Ahat{} and various other categories constructed in terms of $A$,
  such as $\bHom(A\op,M)$, are closed model categories (under suitable
  assumptions\ldots).
\item\label{it:92.c}
  Independence of the derived category of $\bHom(A\op,M)$ on the
  choice of test category, notably when $M$ is a topos or the dual of
  a\pspage{321} topos (with suitable assumptions on $A,M$\ldots).
\item\label{it:92.d}
  Possibility of expressing finite type of a homotopy type in terms of
  \Ahatfp, for suitable test categories $A$.
\item\label{it:92.e}
  Possibility of defining boundary operations within a test category
  -- and/or getting Dold-Puppe type relations.
\end{enumerate}

\starsbreak

Before resuming more technical work with the matters left over for
part \ref{ch:IV}, I would like still to write down some afterthoughts,
concerning the question of boundary operations in a test category
(question \ref{q:91.8} in our review). It occurred to me that perhaps
it isn't such a good idea, to try at all costs to subordinate this
question to a question of cellular decompositions of spheres, however
natural this idea may be in view of the examples of the standard types
of complexes and multicomplexes. In this connection, I remember that
among my first thoughts when starting unwittingly on the modelizing
story, was that a ``test category'' $A$ (namely one such that \Ahat{}
should be ``modelizing'') should more or less correspond to such
decompositions. Soon after it came as a big surprise that so little
was needed in fact for $A$ to merit the name of a test category -- and
that the relevant conditions had nothing to do with cellular
decompositions of this or that. The same may well turn out, when
looking for a generalization of the standard simplicial, cubical or
hemispherical chain complexes, giving rise to the homology and
cohomology invariants of a given ``complex''. The kind of set-up I
proposed in yesterday's notes, for a formalism of boundary operations
in a test category $A$, now looks to me in some respects
somewhat\scrcomment{``étriqué'' can again be translated as
  ``narrow-minded''} ``étriqué'', and I'll try another start in a
different spirit.

In order not to get involved in irrelevant technicalities, I assume
that the basic localizer \scrW{} is $\scrWoo=$ usual weak
equivalence. It seems that one basic fact for writing down
a relationship between homotopy types and ``homology types'', is the
existence of a canonical ``abelianization functor''
\begin{equation}
  \label{eq:92.1}
  \Hot\to\D_\bullet(\Ab) \quad (\eqdef \HotabOf),\tag{1}
\end{equation}
where \Ab{} is the abelian category of abelian groups, and
$\D_\bullet$ designates the ``derived category'' of the category
$\Ch_\bullet(\Ab)$ of chain complexes of abelian groups, namely its
localization with respect to ``weak equivalences'', i.e.,
quasi-isomorphisms:\pspage{322}
\begin{equation}
  \label{eq:92.2}
  \D_\bullet\Ab = W^{-1}\Ch_\bullet\Ab,\tag{2}
\end{equation}
where $W$ means ``quasi-isomorphisms'', i.e., maps inducing
isomorphisms for all homology groups. The most common way for defining
the canonical functor \eqref{eq:92.1}, where as usual here \Hot{} is
defined as $\scrW^{-1}\Cat$, is via the test category $\Simplex$, as
the composition in the bottom row of
\begin{equation}
  \label{eq:92.3}
  \begin{tabular}{@{}c@{}}
    \begin{tikzcd}[baseline=(O.base)]
      \Cat \ar[r,"i^*"]\ar[d] &
      \Simplexhat \ar[r,"\Wh_{\Simplex}"]\ar[d] &
      \Simplexhatab \ar[r,"\equ\supDP"]\ar[d] &
      \Ch_\bullet\Ab \ar[d] \\
      \Hot \ar[r,"\equ"] &
      \HotOf_{\Simplex} \ar[r] &
      \HotabOf_{\Simplex} \ar[r,"\equ"] &
      |[alias=O]| \D_\bullet\Ab
    \end{tikzcd},
  \end{tabular}\tag{3}
\end{equation}
which is deduced from the top row by passing to the localized
categories. The subscript $\mathrm{ab}$ in \Simplexhatab{} denotes the
category of abelian group objects in \Simplexhat, the functor
\begin{equation}
  \label{eq:92.4}
  \Wh_{\Simplex} : \Simplexhat\to\Simplexhatab\tag{4}
\end{equation}
is the ``\emph{abelianization functor}'' obtained by composing a
presheaf $\Simplexop\to\Sets$ with the abelianization functor
\[\Sets\to\Ab, \quad X\mapsto\bZ^{(X)}.\]
We call this functor $\Wh_\Simplex$ also the ``Whitehead functor'', as
its main property is expressed in \emph{Whitehead's theorem}, namely
that it is compatible with weak equivalences (where weak equivalences
in \Simplexhatab{} are defined in terms of the underlying
semisimplicial sets, forgetting the addition laws). The localized
category of \Simplexhatab{} with respect to the latter notion of weak
equivalence is denoted by $\HotabOf_\Simplex$, the functor
\begin{equation}
  \label{eq:92.4prime}
  \HotOf_\Simplex \to \HotabOf_\Simplex\tag{4'}
\end{equation}
induced by $\Wh_\Simplex$ may be equally (and still more validly) be
designated by $\Wh_\Simplex$. The two subscripts $\Simplex$ (in $\Wh$
and in $\HotabOf$) refer to the fact that the notions make still a
sense when $\Simplex$ is replaced by an arbitrary small category $A$,
cf.\ below.

The functor $\DP$ in the top row is the well-known \emph{Dold-Puppe
  functor}, which is an equivalence of categories. As for $i^*$, it is
defined in terms of an arbitrary test functor
\[i:\Simplex\to\Cat,\]
which may be either the standard inclusion (which is the more commonly
used one) or the canonical functor $a\mapsto\Simplex_{/a}$, called
$i_\Simplex$. The functors $i^*$ corresponding to different choices of
$i$ are of course in general non-isomorphic, however (as follows from
section \ref{sec:77})\pspage{323} the corresponding functor
\[\overline{i^*}:\Hot\to\HotOf_\Simplex\]
is independent of such choice, up to canonical isomorphism.

Instead of $\Simplex$, we could have worked with $\Square$ or $\Globe$
instead, as these give rise to a Dold-Puppe functor (which is still an
equivalence), and (almost certainly, see below) to a corresponding
variant of the ``Whitehead theorem''. We thus get two other ways for
defining a canonical ``abelianization functor'' \eqref{eq:92.1} for
homotopy types, and it should be an easy and pleasant exercise to show
these three functors are canonically isomorphic, using the fact that a
product of test categories is again a test category.

\begin{remark}
  Of course, when concerned mainly with defining a functor
  \eqref{eq:92.1} we don't really need Whitehead and Dold-Puppe
  theorems -- indeed, instead of taking the functor
  \[\DP\circ\Wh_\Simplex:\Simplexhat\to\Ch_\bullet\Ab,\]
  we could have taken directly (using the standard boundary operations
  between the components of a semisimplicial abelian group) the
  standard chain complex structure of $\bZ^{(K_\bullet)}$ (for
  $K_\bullet$ in \Simplexhat), without taking the trouble and
  normalizing it à la Dold-Puppe -- and it is a lot more trivial than
  either Whitehead's or Dold-Puppe's theorem, that the latter functor
  transforms weak equivalences into quasi-isomorphisms; moreover, it
  gives rise to the same functor
  \begin{equation}
    \label{eq:92.5}
    \HotOf_\Simplex\to\D_\bullet\Ab\quad(=\HotabOf)\tag{5}
  \end{equation}
  as $\DP\circ\Wh_\Simplex$.
\end{remark}

Let now $A$ be \emph{any} small category, we are interested in the
functor
\begin{equation}
  \label{eq:92.6}
  \HotOf_A \to \D_\bullet\Ab\quad (=\HotabOf)\tag{6}
\end{equation}
obtained as the composition
\begin{equation}
  \label{eq:92.7}
  \HotOf_A\to\Hot \to\D_\bullet\Ab,\tag{7}
\end{equation}
where the second functor is the abelianization functor
\eqref{eq:92.1}, and the first is the canonical functor, deduced by
localization from
\begin{equation}
  \label{eq:92.8}
  i_A:\Ahat\to\Cat,\quad a\mapsto A_{/a}.\tag{8}
\end{equation}
We see immediately that for $A=\Simplex$, the functor \eqref{eq:92.6}
reduces to \eqref{eq:92.5} up to canonical isomorphism -- and the same
of course when\pspage{324} $A$ is either $\Square$ or $\Globe$. In
these three cases, the functor \eqref{eq:92.6} can be factorized in a
natural way through the category
\[\HotabOf_A = W^{-1}\Ahatab,\]
where now $W$ stands for the set of ``weak equivalences'' in \Ahatab,
defined in the same way as above in the case $A=\Simplex$. The
question then arises, for any given $A$, as to whether such a
factorization can be still obtained, and how exactly.

This formulation is inspired by the description of abelianization of
homotopy types via the (slightly sophisticated diagram)
\eqref{eq:92.3}. When following the more naive approach of the remark
above, this leads us to the closely related question of defining
\eqref{eq:92.6} via a composition
\begin{equation}
  \label{eq:92.9}
  \Ahat\xrightarrow{\Wh_A} \Ahatab \xrightarrow L \Ch_\bullet\Ab,\tag{9}
\end{equation}
(for a suitable functor $L$, cf.\ below), by passing to localizations.

Both approaches seem to me of interest. The first one, to make sense
at all as stated, relies on the existence of a canonical functor
\begin{equation}
  \label{eq:92.9prime}
  \HotOf_A\to\HotabOf_A,\tag{9'}
\end{equation}
induced by the abelianization functor
\[\Wh_A:\Ahat\to\Ahatab,\quad X\mapsto\bZ^{(X)},\]
i.e., on the \emph{validity of Whitehead's theorem, with $\Simplex$
  replaced by $A$}. This looks like an interesting question, whose
answer should be in the affirmative. At any rate, if we can find an
aspheric functor
\begin{equation}
  \label{eq:92.10}
  j:\Simplex\to\Ahat\tag{10}
\end{equation}
(with respect to the standard asphericity structure of $A$), then the
answer is affirmative, as we are immediately reduced to the known case
$A$ is replaced by $\Simplex$. Thus, the answer is possibly tied up
with the question of existence of test functors, which we'll deal with
presumably in part \ref{ch:V}. It should be noted though that if such
a functor \eqref{eq:92.10} exists, then necessarily $A$ is aspheric,
and even totally aspheric -- a substantial restriction indeed.

More generally, let
\begin{equation}
  \label{eq:92.10prime}
  j:B\to\Ahat\tag{10'}
\end{equation}
any aspheric functor with respect to the standard asphericity
structure of \Ahat, where $B$ is any small category. The corresponding
functor\pspage{325}
\begin{equation}
  \label{eq:92.11}
  j^*:\Ahat\to\Bhat\tag{11}
\end{equation}
then satisfies
\[(j^*)^{-1}(\scrWB)=\scrWA,\]
and the corresponding functor for the localizations gives rise to a
commutative diagram
\[\begin{tikzcd}[baseline=(O.base),column sep=tiny]
  \HotOf_A \ar[rr,"\overline{j^*}"]\ar[dr,"\overline{i_A}"'] & &
  \HotOf_B \ar[dl,"\overline{i_B}"] \\
  & |[alias=O]| \Hot &
\end{tikzcd},\]
and hence the corresponding diagram
\begin{equation}
  \label{eq:92.12}
  \begin{tikzcd}[column sep=tiny]
    \HotOf_A \ar[rr,"\overline{j^*}"]\ar[dr] & &
    \HotOf_B \ar[dl] \\
    & \D_\bullet\Ab=\HotabOf &
  \end{tikzcd}
  \tag{12}
\end{equation}
is commutative, where the vertical arrows are the canonical functors
\eqref{eq:92.6}. Coming back to the base $B=\Simplex$ say, this shows
that \eqref{eq:92.6} can be viewed as the composition
\[ \HotOf_A \xrightarrow{\overline{j^*}} \HotOf_\Simplex
\xrightarrow{\overline{\Wh_\Simplex}} \HotabOf_\Simplex
\xrightarrow[\equ]{\overline{\DP}} \HotabOf,\]
and hence it can be inserted in the commutative diagram
\begin{equation}
  \label{eq:92.13}
  \begin{tabular}{@{}c@{}}
    \begin{tikzcd}[baseline=(O.base)]
      \HotOf_A\ar[r]\ar[d] & \HotOf_\Simplex\ar[d] & \\
      \HotabOf_A\ar[r] & \HotabOf_\Simplex \ar[r,"\equ"] &
      |[alias=O]| \HotabOf
    \end{tikzcd},
  \end{tabular}\tag{13}
\end{equation}
where the functor
\[\HotabOf_A \to \HotabOf_\Simplex\]
is induced by $j^*\subab$. Thus, we get the wished for factorization
of \eqref{eq:92.6} via $\HotabOf_A$, \emph{provided we can find an
  aspheric functor} \eqref{eq:92.10}. It should not be hard moreover
to see that the factorizing functor obtained from \eqref{eq:92.13},
namely
\begin{equation}
  \label{eq:92.14}
  \HotabOf_A\to\HotabOf,\tag{14}
\end{equation}
does not depend up to canonical isomorphism on the choice of $j$, at
least in the case when $A$ is a contractor, using the end remarks of
section \ref{sec:82} (p.\ \ref{p:272}) concerning products of aspheric
functors.

The question remains whether we can define \eqref{eq:92.14} for (more
or less) any small category $A$, without having to rely upon the
existence of\pspage{326} an aspheric functor \eqref{eq:92.10}, in such
a way that it factors the canonical functor \eqref{eq:92.6} (granting
Whitehead's theorem holds for $A$), and moreover that for an aspheric
functor \eqref{eq:92.10prime} $j:B\to\Ahat$, giving rise to
\[\overline{j^*\subab}:\HotabOf_A\to\HotabOf_B,\]
then corresponding diagram
\begin{equation}
  \label{eq:92.15}
  \begin{tikzcd}[column sep=tiny]
    \HotabOf_A\ar[rr]\ar[dr] & & \HotabOf_B\ar[dl] \\ & \HotabOf &
  \end{tikzcd}\tag{15}
\end{equation}
should commute, where the vertical arrows are the functors
\eqref{eq:92.14}.

For defining \eqref{eq:92.14} in this general case, recalling that
\[\HotabOf \fromequ \HotabOf_\Simplex,\]
we can't help it and have to use diagram \eqref{eq:92.3} and the
functor
\[i_A:\Ahat\to\Cat,\]
or rather (\Cat{} serving only as an intermediary) the functor
\[u \eqdef i^*i_A : \Ahat\to\Simplexhat,\]
so that \eqref{eq:92.6} can be viewed as deduced by localization of
\Ahat{} from the composition
\begin{equation}
  \label{eq:92.star}
  \Ahat\xrightarrow
  u\Simplexhat\xrightarrow{\Wh_\Simplex}\Simplexhatab\to\HotabOf_\Simplex
  \quad (\toequ\HotabOf).\tag{*}
\end{equation}
One difficulty here is that $u$ does not commute to finite products,
and hence doesn't induce a functor
\[\Ahatab\to\Simplexhatab,\]
it would seem. Now this difficulty, I just noticed, can be overcome,
using the fact that $i_A$ and hence also $u$ \emph{commutes to fibered
  products}, or, what amounts to the same, induces an \emph{exact}
functor
\[u_0 : \Ahat\to\Simplexhat_{/F}\toequ(\Simplex_{/F})\uphat,\]
where
\[F=u(e_\Ahat)=i^*(A)\]
is a suitable object in \Simplexhat. A fortiori, $u_0$ commutes to
finite products, hence transforms abelian group objects into same,
i.e., induces
\[u_{0\mathrm{ab}} : \Ahatab\to(\Simplex_{/F})\uphat\subab,\]
on\pspage{327} the other hand, we do have too a natural functor
\begin{equation}
  \label{eq:92.16}
  \alpha^{\mathrm{ab}}_!:(\Simplex_{/F})\uphat\subab\to\Simplexhatab,\tag{16}
\end{equation}
defined as the left adjoint of the evident functor
\[\alpha^*\subab:\Simplexhatab\to(\Simplex_{/F})\uphat\subab\]
induced by the left exact functor $\alpha^*$, where $\alpha$ is the
``localization morphism of topoi''
\[\alpha:\Simplex_{/F}\to \Simplex_F\]
defined by the object $F$ of \Simplexhat. We thus get a diagram
\begin{equation}
  \label{eq:92.17}
  \begin{tabular}{@{}c@{}}
    \begin{tikzcd}[baseline=(O.base)]
      \Ahat\ar[r,"\Wh_A"] \ar[d,"u_0"'] \ar[dd,bend right=50,"u"'] &
      \Ahatab\ar[d,"u_{0\mathrm{ab}}"]  & & \\
      (\Simplex_{/F})\uphat \ar[r,"\Wh_{\Simplex_{/F}}"]
      \ar[d,"\alpha_!"'] &
      (\Simplex_{/F})\uphat\subab \ar[d,"\alpha_!^{\mathrm{ab}}"] & & \\
      \Simplexhat\ar[r,"\Wh_\Simplex"] & \Simplexhatab \ar[r] &
      \HotabOf_\Simplex \ar[r] & |[alias=O]| \HotabOf
    \end{tikzcd},
  \end{tabular}
  \tag{17}
\end{equation}
containing \eqref{eq:92.star} above as the composition of maps in the
left-hand vertical column and in the bottom row. The lower square in
this diagram commutes (up to natural isomorphism) -- which is a
general fact surely for morphisms of topoi $f:X\to Y$ such that $f_!$
exists, which allows to define too a functor $f_!^{\mathrm{ab}}$ as
\eqref{eq:92.16} above. The upper square though doesn't look at all
commutative, too bad! The only hope left now is that the natural
compatibility arrow for this square (there is bound to be one, isn't
there!), when composed with the lower square (so as to give a
compatibility map for the composite rectangle) \emph{should give rise
  to a weak equivalence in} \Simplexhatab, for any choice of an object
$X$ in \Ahat.

It would seem to me that a reasonable functor \eqref{eq:92.14} will
exist, without an existence assumption of a test functor
\eqref{eq:92.10}, exactly in those cases when the composite rectangle
in \eqref{eq:92.17} is ``commutative up to weak equivalence''. I have
no idea whether or not this is true for any small category $A$, not
even (I confess) when there \emph{is} an aspheric functor
\eqref{eq:92.10} -- as a matter of fact, I don't feel like going any
further now in this direction, and trying to check anything
whatsoever.

I was a little rash in the definite statement I made about the
``exact'' assumption to make for a ``reasonable'' functor
\eqref{eq:92.14} to exist; another seems needed still, namely that the
functor\pspage{328}
\begin{equation}
  \label{eq:92.18}
  \Ahatab\to\Simplexhatab\tag{18}
\end{equation}
we obtain by composing the two arrows in the right hand vertical
column, transforms weak equivalences into same, which is needed in
order to deduce \eqref{eq:92.14} by passing to the localized
categories. When the two assumptions are satisfied, then the functor
\eqref{eq:92.14} obtained from \eqref{eq:92.18} does factorize
\eqref{eq:92.6} as required, and it should be clear too that it
satisfies the compatibility \eqref{eq:92.15}.

Thus, it seems there are good prospects for getting canonical functors
\[\HotOf_A \to \HotabOf_A \to \HotabOf,\]
whose composition is \eqref{eq:92.6}, i.e., inserting into the
commutative diagram
\begin{equation}
  \label{eq:92.19}
  \begin{tabular}{@{}c@{}}
    \begin{tikzcd}[baseline=(O.base)]
      \HotOf_A\ar[r]\ar[d] & \HotabOf_A\ar[d] \\
      \HotOf \ar[r] & |[alias=O]| \HotabOf
    \end{tikzcd}.
  \end{tabular}
  \tag{19}
\end{equation}
The next question then which arises is \emph{whether the second
  vertical arrow} (namely \eqref{eq:92.14}) \emph{is an equivalence,
  whenever the first one is}, i.e., when $A$ is a pseudo-test
category, or whether this is true if we make some familiar extra
assumption on $A$, such as being an actual test category say.

I feel I am getting gradually back into thin air conjecturing, I
wouldn't go on too long this way! This whole $\HotabOf_A$ business was
just a digression, which then took me longer than expected, it doesn't
seem to have much to do with what I have been out for in this section,
namely afterthoughts about \emph{boundary operations in a test
  category}, which are designed to gave a ``computational''
description of the canonical functor \eqref{eq:92.6}, the latter being
defined without any restriction nor difficulty for \emph{any} small
category $A$. More accurately still, we want to describe the
composition of \eqref{eq:92.6} with the canonical functor
$\Ahat\to\HotOf_A$, namely $\Ahat\to\HotabOf$, via a suitable functor
$\Ahat\to\Ch_\bullet\Ab$, with the expectation that the latter should
factor through \Ahatab{} via the abelianization functor $\Wh_A$. In
other words, we are looking for commutative diagrams
\begin{equation}
  \label{eq:92.20}
  \begin{tabular}{@{}c@{}}
    \begin{tikzcd}[baseline=(O.base),sep=tiny]
      & \Ahat\ar[dl]\ar[rr]\ar[dd,"K"] & & \HotOf\ar[dd] \\
      \Ahatab\ar[dr,"L"'] & & & \\
      & \Ch_\bullet\Ab\ar[rr] & & |[alias=O]| \HotabOf
    \end{tikzcd},
  \end{tabular}
  \tag{20}
\end{equation}
where\pspage{329} all functors in the diagram, except for $K$ and $L$,
are the canonical ones familiar to us. The question then is how to
define a suitable $L$, such that the corresponding square (where
$K=L\circ\Wh_A$) should commute up to (canonical?) isomorphism.

There \emph{is} such an $L$, whenever we have an aspheric functor
\eqref{eq:92.10} (where $\Simplex$ may be replaced by one of its
twins), using corresponding semisimplicial chain complexes -- but, as
already remarked yesterday, taking things this way is ``cheating''!
Conceivably too, there are quite general theorems asserting that a
functor from a category \Ahat{} say to a derived additive category
such as $\HotabOf$ can be lifted to the category of models (here
abelian chain complexes) it comes from, and possibly even in a way
factoring through \Ahatab? I don't intend to dive into these questions
either, but rather, make a comment on a general method for
constructing certain functors
\begin{equation}
  \label{eq:92.21}
  L:\Ahatab\to\Ch_\bullet\Ab,\tag{21}
\end{equation}
(maybe not in a way to give rise to a commutative diagram
\eqref{eq:92.20}), as suggested by the standard chain complexes
associated to the three types of complexes, using simplices, cubes or
hemispheres, or multicomplexes (using products of the standard test
categories). Writing
\[\Ahatab\simeq \bHom(A\op,\Ab),\]
we remark that the standard constructions of chain complexes
associated to semisimplicial (say) complexes of abelian groups, makes
sense not only for such abelian complexes, but more generally for
\emph{complexes with values in \emph{any} additive category}, $M$
say. This induces us to look more generally, for any such $M$, for a
functor
\begin{equation}
  \label{eq:92.22}
  L_M:\bHom(A\op,M)\to\Ch_\bullet(M),\tag{22}
\end{equation}
in a way \emph{compatible with additive functors}
\[M\to M'\]
(in the obvious sense of the word).

Now, for any category $B$ (here $A\op$) one can define an
``\emph{enveloping additive category}'' $\Add(B)$, together with a
canonical functor
\begin{equation}
  \label{eq:92.23}
  B\to \Add(B),\tag{23}
\end{equation}
which is ``$2$-universal'' for all possible functors of $B$ into any
additive category $M$. More specifically, for any such $M$, the
corresponding\pspage{330} functor ``composition with
\eqref{eq:92.23}'' is an \emph{equivalence}
\begin{equation}
  \label{eq:92.24}
  \mathbf{Homadd}(\Add(B),M) \toequ \bHom(B,M).\tag{24}
\end{equation}
This condition defines \eqref{eq:92.22} ``up to canonical equivalence''
-- but we'll give an explicit description in a minute. Before doing
so, let's just remark that the universal property of \eqref{eq:92.23}
implies that to give a system of functors $L_M$ as above, ``amounts to
the same'' as giving a chain complex $L_\bullet$ in $\Add(B)$. More
accurately, the category of all systems $L_M$ (where maps are defined
in an evident way) is \emph{equivalent} to the category
$\Ch_\bullet(\Add(B))$, where $B=A\op$. The functors $L$
\eqref{eq:92.21} we are specifically interested in, are those which
are associated to some chain complex in $\Add(A\op)$,
\begin{equation}
  \label{eq:92.25}
  L_\bullet\in\Ob(\Ch_\bullet(\Add(A\op))),\tag{25}
\end{equation}
by the formula
\begin{equation}
  \label{eq:92.26}
  L(X) = \widetilde X(L_\bullet)\quad\text{for any $X$ in
    $\Ahatab\simeq\bHom(A\op,\mathrm{Ab})$}\tag{26}
\end{equation}
where
\[\widetilde X: \Add(A\op)\to\Ab\]
is the \emph{additive} functor corresponding to $X$.

We are thus led to the question: if $A$ is any small category, does
there exist a chain complex $L_\bullet$ in $\Add(A\op)$, the additive
envelope of $A\op$, giving rise to a functor \eqref{eq:92.21} via
\eqref{eq:92.26} and hence to a diagram \eqref{eq:92.20}, such that
the square in \eqref{eq:92.20} commutes up to isomorphism? And when
this is so, what kind of unicity statement, if any, can be made for
$L_\bullet$ (such as being unique up to chain homotopy say), and what
about the structure of the category of all pairs
$(L_\bullet,\lambda)$, where $\lambda$ is a compatibility isomorphism
making the square in \eqref{eq:92.20} commute?

We are far here from the rather narrow set-up in yesterday's notes,
and as far as existence goes, if no extra conditions are put upon
$L_\bullet$, it seems likely that for a rather large class of small
categories $A$ (if not all) it should hold true. At any rate, the
functor $L$ obtained from a functor \eqref{eq:92.10}, i.e., ``by
cheating'', is visibly associated to an $L_\bullet$. Sorry, we have to
assume that the functor $j$ factors even through $A$ itself, i.e., is
just an aspherical functor between the small categories $\Simplex$ and
$A$, a much more stringent condition on $A$ to be sure -- and which
implies that there is an induced\pspage{331} functor
\begin{equation}
  \label{eq:92.27}
  \Add(j\op):\Add(\Simplexop)\to\Add(A\op),\tag{27}
\end{equation}
hence we get an $L_\bullet$ as the image of the canonical chain
complex $L_\bullet^\Simplex$ we got in $\Add(\Simplexop)$.

\bigbreak
\presectionfill\ondate{6.7.}\par

\hangsection[The afterthought continued: abelianizators, and
\dots]{The afterthought continued: abelianizators, and ``standard''
  abelianizators for categories with boundary
  operators.}\label{sec:93}%
It is time to give the promised construction of $\Add(B)$, the
additive envelope of $B$, for any given category $B$. The obvious idea
is to enlarge the sets $\Hom(a,b)$, for $a$ and $b$ in $B$, by taking
linear combinations with coefficients in \bZ, i.e., writing
\begin{equation}
  \label{eq:93.1}
  \Hom_{\Add(B)}(a,b) = \bZ^{(\Hom(a,b))},\tag{1}
\end{equation}
and composing these $\Hom_{\Add}$ in the obvious way. This is not
quite enough though, as we still have to add new objects, namely
direct sums of objects in $B$. The most convenient way for doing so
seems by defining an object of $\Add(B)$ to be defined by a finite set
$I$ (in the given universe), namely the indexing set for taking the
direct sum, and a map
\[ I\to\Ob B,\]
in other words, the new objects are just \emph{families} of objects of
$B$
\[(b_i)_{i\in I},\]
indexed by finite sets. We'll however denote by
\begin{equation}
  \label{eq:93.2}
  \bigoplus_{i\in I}b_i\tag{2}
\end{equation}
the corresponding object of $\Add(B)$, as this will turn out to be the
direct sum indeed of the images of the $b_i$'s in $\Add(B)$ -- but of
course we'll ignore the possible existence of direct sums in $B$
itself, when they exist, and not confuse \eqref{eq:93.2} with a direct
sum taken in $B$. Writing $\Homadd$ instead of $\Hom_{\Add(B)}$ for
the sake of abbreviation, the maps between objects \eqref{eq:93.2} are
defined by matrices in the obvious way
\begin{multline}
  \label{eq:93.3}
  \Homadd((a_i)_{i\in I},(b_j)_{j\in J}) = \\
  \set[\big]{(u_{ij})_{(i,j)\in I\times
      J}}{u_{ij}\in\Homadd(a_i,b_j)=\bZ^{\Hom(a_i,b_j)}},\tag{3}
\end{multline}
while composition of maps is defined by the composition of
matrices. We thus get a new category $\Add(B)$ and a functor
\begin{equation}
  \label{eq:93.4}
  B\to\Add(B),\tag{4}
\end{equation}
it\pspage{332} is immediately checked that $\Add(B)$ is an additive
category and that the functor \eqref{eq:93.4} has the $2$-universal
property for functors from $B$ into any additive category, stated in
yesterday's notes (p.\ \ref{p:330}).

\begin{remark}
  The same construction essentially applies when considering the
  universal problem of mapping $B$ into any $k$-\emph{additive}
  category $M$ (where $k$ is any commutative ring with unit), i.e.,
  an additive category $M$ endowed with a ring homomorphism
  \[k\to\End(\id_M),\]
  replacing $\bZ$ by $k$ in formulas \eqref{eq:93.1} and
  \eqref{eq:93.3}. The ``abelianization'' questions touched at in
  yesterday's notes still make sense in terms of ``$k$-linearization''
  -- a notion much in the spirit of our introduction of a general
  basic localizer \scrW, as the very notion of $k$-linearization will
  give rise to a corresponding basic localizer $\scrW_k$\ldots
\end{remark}

Let's come back to the case when $B=A\op$, $A$ being a small category,
and to our question about chain complexes
\[L_\bullet \quad\text{in}\quad \Ch_\bullet(\Add(A\op))\]
giving rise to a commutative diagram \eqref{eq:92.20} (p.\
\ref{p:328}), up to isomorphism. A pair
\[(L_\bullet,\lambda),\]
where $L_\bullet$ is a chain complex as above, and $\lambda$ a
compatibility isomorphism for the square in \eqref{eq:92.20}, could be
suggestively called an \emph{abelianizator} for the small category
$A$. The question of existence, and uniqueness up to homotopy say, of
an abelianizator for $A$ seems especially relevant when $A$ is a test
category say, and hence \Ahat{} modelizes homotopy types. In any case,
in terms of an abelianizator we get an additive functor
\begin{equation}
  \label{eq:93.5}
  L:\Ahatab\to\Ch_\bullet\Ab,\tag{5}
\end{equation}
and the question arises whether this is compatible with weak
equivalences and quasi-isomorphisms; maybe even if this is not
automatic, we should insist it holds when defining the notion of an
abelianizator. When this is OK, then by passing to localizations we
deduce from \eqref{eq:93.5} a functor
\begin{equation}
  \label{eq:93.6}
  \HotabOf_A\to\HotabOf,\tag{6}
\end{equation}
i.e., a functor \eqref{eq:92.14} as looked for in yesterday's notes,
giving rise to the commutative diagram \eqref{eq:92.19} (p.\
\ref{p:328}) -- whereas commutativity of diagrams of the type
\eqref{eq:92.15} (p.\ \ref{p:326}) looks less obvious.

When\pspage{333} $A$ is a finite product of copies taken from among
the three standard test categories $\Simplex$, $\Square$, $\Globe$,
the standard chain complex structure on multicomplexes does furnish us
with a \emph{``canonical'' abelianizator} for $A$, which we may denote
by $L_\bullet^A$ (as we did yesterday for $A=\Simplex$). This
``standard'' abelianizator has some very remarkable extra features
which I would like to pin down, which had caused our rather narrow
focus in the notes of two days ago (section \ref{sec:91}).

\namedlabel{cond:93.a}{a)}\enspace
There is a ``dimension map''
\begin{equation}
  \label{eq:93.7}
  \dim:\Ob A\to \bN.\tag{7}
\end{equation}
It can be described (in the particular case above at any rate) in
terms of the intrinsic category structure of $A$, by associating to
every $a$ in $A$ the ordered set
\begin{equation}
  \label{eq:93.8}
  i(a)=\text{set of subobjects of $a$ in $A$.}\tag{8}
\end{equation}
(NB\enspace not to be confused with subobjects of $a$ in \Ahat, namely
sieves in $a$). This is an ordered set with a largest object (namely
$a$ itself), and which turns out to be finite (in the particular cases
considered), hence of finite combinatorial dimension (equal to the
dimension of the geometrical realization $\abs{i(a)}$), and we have
\begin{equation}
  \label{eq:93.9}
  \dim(a)=\dim i(a).\tag{9}
\end{equation}

\namedlabel{cond:93.b}{b)}\enspace
The $n$'th component $L_n$ of $L_\bullet^A=L_\bullet$ is given by
\begin{equation}
  \label{eq:93.10}
  L_n = \bigoplus_{\dim(a)=n} a\tag{10}
\end{equation}
(where the direct sum of course is taken in $\Add(A\op)$ as in
\eqref{eq:93.2} above), which makes sense when we assume (as is the
case in our example) that the map \eqref{eq:93.7} is ``finite'', i.e.,
has finite fibers. For instance, in all our ``standard'' cases, there
is just \emph{one} object of $A$ which is of dimension $0$, and this
is also the final object.

\namedlabel{cond:93.c}{c)}\enspace
The differential operator
\begin{equation}
  \label{eq:93.11}
  d_n:L_n\to L_{n-1}\tag{11}
\end{equation}
can be obtained in the following way. We have only to describe $d_n$
on each summand $a$ of $L_n$, i.e., by \eqref{eq:93.10} on each $a$ in
\begin{equation}
  \label{eq:93.12}
  A_n=\set[\big]{a\in\Ob A}{\dim a=n}.\tag{12}
\end{equation}
In view of \eqref{eq:93.3}, this restriction $d_n\mid a$ can be
described as a linear combination of elements in the disjoint sum of
the sets
\[\Hom(b,a),\quad\text{with}\quad b\in A_{n-1}\quad (\text{$a\in A_n$
  fixed}).\]
This\pspage{334} being clear, the non-zero coefficients which occur in
this linear combination are all $\pm1$, and moreover the maps which
target $a$
\[b\to a\quad(b\in A_{n-1})\]
which occur with non-zero coefficient are exactly \emph{all
  monomorphisms} from objects $b$ in $A_{n-1}$ into $a$. Thus, the
differential operators are known, when we know, for all monomorphisms
in $A$
\begin{equation}
  \label{eq:93.13}
  \partial:b\to a,\quad\text{with}\quad \dim a=\dim b+1\tag{13}
\end{equation}
(the so-called \emph{``boundary maps''}), the corresponding
coefficients
\begin{equation}
  \label{eq:93.14}
  \varepsilon_{\partial}\in\{\pm1\}.\tag{14}
\end{equation}
Instinct tells us, at this point, that we may get into trouble, when
trying to define (in a more or less general case) boundary operations
in such a way, because of the ambiguity in the definition of
subobjects, namely, because of possible \emph{existence of
  isomorphisms which may not be identities}. But precisely, in the
standard cases we are copying from, \emph{any isomorphism is an
  identity}!

\namedlabel{cond:93.d}{d)}\enspace
For describing the ``signatures'' \eqref{eq:93.14}, in one of the
``standard'' cases, we still need to remark that for any $a$ in $A$,
we have
\begin{equation}
  \label{eq:93.15}
  \text{$i(a)$ is an $n$-cell, with $n=\dim(a)$,}\tag{15}
\end{equation}
and the choice of the signatures will be determined by a choice of
orientations
\begin{equation}
  \label{eq:93.16}
  \text{$\omega_a$ an \emph{orientation} of the $n$-cell
    $\abs{i(a)}$,}\tag{16}
\end{equation}
(a notion which could be given a purely combinatorial definition, by
induction on the dimension of a given ordered set whose geometrical
realization is a variety\ldots). We then get a the \emph{``Stokes
  rule''}
\begin{flalign}\label{eq:93.17}
  &&&\parbox[t]{0.9\textwidth}{For a boundary map $\partial:b\to a$,
  $\varepsilon_{\partial}=+1$ if{f} $\omega_b$ is ``induced'' à la
  Stokes by $\omega_a$, via the induced orientation on the boundary of
  $i(a)$ (which is the union of the images of all $i(b)$'s, for all
  boundary operations with target $a$).}
  \tag{17}
\end{flalign}

Whether or not we are in a ``standard'' case, if $A$ is any category
such that for any object $a$ of $A$, the ordered set $i(a)$ of its
subobjects in $A$ is finite, and its geometrical realization is an
$n$-cell (call $n$ the ``dimension'' of $a$), and if moreover for a
given $n$, the set $A_n$ of objects with dimension $n$ is finite, and
also (to be on the safe side!) assuming that all isomorphisms are
identities, then for \emph{any} choice of orientations
\eqref{eq:93.16}, giving rise to a system of\pspage{335} signatures
\eqref{eq:93.13} by the ``Stokes rule'' \eqref{eq:93.17}, the
corresponding operators \eqref{eq:93.11} do turn the family $(L_n)$
into a chain complex, namely we have the relations
\begin{equation}
  \label{eq:93.18}
  d_{n-1}d_n=0.\tag{18}
\end{equation}
This follows immediately from the well-known anti-commutativity
property of (twofold) induction of orientation on boundaries.

Things are a little more delicate if we don't assume that isomorphisms
are identities, even if (by compensation) we should insist that two
distinct objects are never isomorphic. To define $d_n$, we then must
\emph{choose}, for any subobject $b$ of dimension $n-1$ of an object
$a$ of dimension $n$, just \emph{one} representative monomorphism
\eqref{eq:93.13} of $b$. The coherence condition then needed in order
to get \eqref{eq:93.18} is that any square diagram
\begin{equation}
  \label{eq:93.19}
  \begin{tikzcd}[sep=tiny]
    & c\ar[dl]\ar[dr] & \\
    b\ar[dr] && b'\ar[dl] \\
    & a &
  \end{tikzcd}
  \tag{19}
\end{equation}
made up with such restricted boundary maps, should commute -- a
somewhat delicate condition, presumably hard to ensure, for the
choices involved for defining the ``strict'' boundary maps in $A$.

In one case as in the other, we are very close of course to the set-up
envisioned in section \ref{sec:91} -- it wouldn't be hard even to fit
the case considered here into this set-up, if we make the slight extra
assumption that any map in $A$ factors into an
epimorphism-with-section, followed by a monomorphism (which is true
indeed in the ``standard'' cases), which will ensure that for varying
$a$, $i(a)$ is indeed a functor with values in \Ord, as stated in
loc.\ cit. But from the point of view of construction of
abelianizators, it would seem that the existence of the functor
\begin{equation}
  \label{eq:93.20}
  i:A\to\Ord\tag{20}
\end{equation}
is irrelevant.

Our main question now, of course, is about \emph{the chain complex
  $L_\bullet$ being an abelianizator or not}. The question is
interesting even in the standard cases, by choosing the orientations
\eqref{eq:93.16} in a way different from the standard one. Are the
corresponding chain complexes in $\Add(A\op)$ necessarily
chain-homotopic?

It just occurs to me that indeed, between the chain
complexes\pspage{336}
\begin{equation}
  \label{eq:93.22}
  L_\bullet^\omega = ((L_n)_{n\in\bN},(d_n^\omega)_{n\in\bN})\tag{22}
\end{equation}
associated to all possible systems of orientations
\begin{equation}
  \label{eq:93.23}
  \omega=(\omega_a)_{a\in\Ob(A)}\tag{23}
\end{equation}
of the various cells $i(a)$, there is a canonical \emph{transitive
  system} of isomorphisms, by defining the isomorphism
\begin{equation}
  \label{eq:93.24}
  u_{\omega,\omega'}:L_\bullet \tosim L_\bullet'\tag{24}
\end{equation}
for two different choices $\omega,\omega'$ of systems of orientations,
by
\begin{equation}
  \label{eq:93.25}
  u_{\omega,\omega'} \mid a=\varepsilon_a^{\omega,\omega'}\,
  \id_a,\quad \varepsilon_a^{\omega,\omega'}\in\{\pm1\},\tag{25}
\end{equation}
where the sign $\varepsilon_a^{\omega,\omega'}$ is equal to $+1$ or
$-!$, \emph{according to whether $\omega_a$ and $\omega_a'$ are equal
  or not}. It is immediate that \eqref{eq:93.24} then is an
isomorphism componentwise, respecting degrees, and commuting to the
respective differential operators. Transitivity of the isomorphisms
\eqref{eq:93.24} for a triple $(\omega,\omega',\omega'')$ is equally
immediate. This implies that by this transitive system of
isomorphisms, \emph{we may identify all the chain complexes}
$L_\bullet^\omega$ in $\Add(A\op)$ to a single chain complex,
canonically isomorphic to each $L_\bullet^\omega$, and which we may
just designate by $L_\bullet$. This chain complex now is defined
intrinsically in terms of the category structure of $A$ (up to
canonical isomorphism), in the ``safe'' case at any rate when every
isomorphism in $A$ is an identity, so that in the construction of
$L_\bullet^\omega$ there enters no other choice besides
$\omega$. Otherwise as seen above (precedent page), we must still
suitably choose the so-called ``strict'' boundary operators
\eqref{eq:93.13}, among all monomorphisms $b\hookrightarrow a$ in $A$
such that $\dim a=\dim b+1$.

In the first case say (isomorphisms being identities), all conditions
considered for $A$ are stable under finite products, that's why in
terms of the three standard cases of $\Simplex$, $\Square$ and
$\Globe$, we could construct others by taking finite products. The
three standard test categories may be viewed as particularly
``economic'' of skillful ways of ``cutting out'' a suitable bunch of
cellular decompositions, and of eliminating automorphisms (by total
ordering of vertices and the like\ldots), so as to ensure: a)\enspace
that isomorphisms in $A$ are identities, b)\enspace the canonical
chain complex $L_\bullet$ in $\Add(A\op)$ is an abelianizator, and
c)\enspace $A$ moreover is a strict test category, and even a
contractor. On the other hand, as all these conditions (plus the
condition \eqref{eq:93.15} of course about the $i(a)$'s representing
$n$-cells)\pspage{337} are stable under taking products (of finite
non-empty families of categories $A_i$), hence in terms of the three
standard cases, the possibility of satisfying them too by the
``multistandard'' test categories. I wonder if there are any other
ways (up to equivalence). If we take categories such as
$\widetilde\Simplex$ (non-ordered simplices), we still get
contractors, but objects have non-trivial automorphisms, and if we
take categories such as \Simplexf{} (ordered simplices without
degeneracy operations, only boundary maps), it is true that
isomorphisms are identities, but the category is no longer a test
category but only a weak one.

If we do not insist on the rigidity assumption (isomorphisms are
identities), but on suitable choice of so-called ``\emph{strict}
boundary operations'' within $A$, then it would seem after all that we
do have a lot more elbow freedom than it seemed by the end of our
reflections on that matter two days ago (cf.\ p.\ \ref{p:318}), where
the picture of the relevant data and corresponding construction of
chain complexes was still a little confused. Let now $A$ be the
category called $A_0$ in loc.\ cit.  We don't have to modify it in
order to introduce orientations of cells $i(a)$ as extra structure and
take account of this in defining a new notion of maps. Therefore, it
is clear that $A$ \emph{just as it is, is a strict test category}
(presumably \emph{not} a contractor though). There is problem of
course of isomorphisms which are not identities, and particularly of
non-trivial automorphisms -- for instance the object $I$ (playing the
part of the unit segment) has a non-trivial automorphism, the
elimination of which does not look so trivial! However, there
\emph{is} a rather evident way of cutting out \emph{strict} boundary
maps, in a way as to satisfy the transitivity condition of p.\
\ref{p:335} -- namely by taking boundary maps \eqref{eq:93.13}
$\partial:b\to a$ which are \emph{inclusions} in the strict sense,
namely the inclusion map of a \emph{subset} of $a$, endowed with the
induced order relation.

Thus, there are many other cases still than just multi-standard test
categories for getting a canonical chain complex $L_\bullet$ in
$\Add(A\op)$, and for which now the question makes sense as to whether
$L_\bullet$ is an abelianizator. In the construction above, we were
careful to assume that the full subcategory $A$ of \Ord, besides
containing $I$, was stable under finite products, so as to make sure
it comes out as a test category. The silly thing is that this
condition is not satisfied by any one among the standard test
categories -- thus, it seems reasonable to try and\pspage{338} replace
it by a suitable substitute, such as the existence, of any two objects
$a$ and $b$ in $A$, of a cellular subdivision of
$\abs{i(a)}\times\abs{i(b)}$, made up with cells of the type
$\abs{i(c)}$, and inducing on the latter the given cellular structure
of $\abs{i(c)}$. The problem is now (besides getting or not an
abelianizator $L_\bullet$) \emph{whether $A$ is at any rate a weak
  test category} (in view of the example \Simplexf, we can't expect
now of course to get an actual test category). Maybe I'll come back to
this later, when writing down a proof for \Simplexf{} being a weak
test category, i.e., a more general result along these lines should
come out alongside.

\begin{remarks}
  I feel the canonical chain complex $L_\bullet$ in $\Add(A\op)$
  constructed in this section, under suitable assumptions on the small
  category $A$, merits a name of its own. We may call it the
  \emph{standard abelianizator} of $A$ -- but this is reasonable only
  if it turns out that in all cases when it can be constructed, it is
  an abelianizator indeed. Another convenient name may be the
  \emph{Dold-Puppe chain complex}, as in the three standard cases, the
  standard Dold-Puppe construction of the ``normalization'' of an
  abelian complex (ss say) can be viewed as being performed in the
  ``universal'' case, namely for $A\op\to\Add(A\op)$, and the
  corresponding ``full'' chain complex, namely $L_\bullet$ -- with
  this grain of salt though that we still have to enlarge $\Add(A\op)$
  slightly, so as to make stable under taking direct summands
  corresponding to projectors. But then it occurs to me that the name
  of Dold-Puppe chain complex is much more suitable for the
  \emph{result of normalization} applied to $L_\bullet$, which (if I
  got it right) is the ``new'' complex discovered by Dold-Puppe,
  together with the inverse construction, whereas $L_\bullet$ had
  already been known for ages (even if not under its universal
  disguise\ldots).
\end{remarks}

\hangsection[Afterthought (continued): retrospective on the ``De Rham
\dots]{Afterthought
  \texorpdfstring{\textup(continued\textup)}{(continued)}:
  retrospective on the ``De Rham complex with divided powers'' and on
  some wishful thinking about linearization of homotopy types and
  arbitrary ground-ring extension in homotopy types.}\label{sec:94}%
In the last section, as in the two preceding days, our emphasis with
abelianization of homotopy types has been to look at it in terms of
more or less arbitrary test categories and the corresponding
elementary modelizers, and even in terms of arbitrary small
categories. This has causes as spinning a kind of dream for a while,
with the Whitehead and Dold-Puppe theorems and generalized boundary
maps as our main thread. Now this reminds me of a rather different
line of thoughts tied up with abelianization, quite independently of
playing around with variable modelizers -- a question which has been
intriguing me for a very long time now, ever since I got acquainted a
little with the very notion of homotopy types, and the corresponding
homology and\pspage{339} cohomology invariants. This is the question
of \emph{how far a homotopy type can be expressed in terms of homology
  or cohomology invariants} (or both together), plus \emph{some
  relevant extra structure}, the most important surely being
cup-products in cohomology (or, dually, ``interior'' operation of
cohomology on homology).\scrcomment{I wonder whether AG knew of the
  Steenrod algebra\ldots yes, see below} Once the notion of derived
categories of various kinds had become familiar, in the early sixties,
the question would appear as expressing, or recovering, a homotopy
type, namely an object in the (highly non-abelian) ``derived
category'' \Hot, in terms of its abelianization in
$\HotabOf=\D_\bullet\Ab$, \emph{endowed with suitable extra
  structure}. It was about clear that this extra structure had to
include, as its main non-commutative item, the fundamental group
$\pi$, so as to allow for description of homology and cohomology
invariants with twisted coefficients. The most natural candidate for
expressing this would be the chain complex associated to the universal
covering, viewed as an object in the derived category
\begin{equation}
  \label{eq:94.26}
  \D_\bullet(\bZ(\pi))\tag{26}
\end{equation}
of chain complexes of modules over the group ring $\bZ(\pi)$. Another
important structural item, giving rise to all cup-products with
non-twisted coefficients, is the diagonal map for the abelianization
\[L_\bullet\quad\text{in}\quad \HotabOf=\D_\bullet\Ab,\]
namely a map
\begin{equation}
  \label{eq:94.27}
  L_\bullet \to L_\bullet \Lotimes L_\bullet,\tag{27}
\end{equation}
where $\Lotimes$\scrcomment{aka the $\Tor$ functor} is the ``total''
left derived functor of tensor product. This map is subjected to
suitable conditions, concerning mainly commutativity and
associativity. In case of a non-$1$-connected space, i.e., $\pi\ne1$,
it shouldn't be hard combining the two structural items so as to get a
structure embodying at any rate cup-products with arbitrary twisted
coefficients. One key question in my mind, which I never really looked
into, was whether these two structures were enough in order to
reconstruct entirely (up to canonical isomorphism) the (pointed,
$0$-connected) homotopy type giving rise to it, and hence also any
other homotopy invariants, such as ``operations'' on cohomology and
the like, K-invariants, etc.

If I got it right, it has been known now for quite a while that even
for a $1$-connected homotopy type, so that the relevant structure
reduces to \eqref{eq:94.27}, that this is \emph{not} quite enough for
recovering the homotopy type, maybe not even the rational homotopy
type. I believe I first got this from Sullivan, namely
that\pspage{340} what was needed for recovering a $1$-connected
rational homotopy type was not merely \eqref{eq:94.27} (where now
$L_\bullet$ is an object of $\D_\bullet(\bQ)$ rather than of
$\D_\bullet(\bZ)=\D_\bullet\Ab$), which reduces more or less (under
suitable finiteness assumption) to knowing the rational cohomology
ring, but an anti-commutative and associative \emph{differential graded
algebra} over $\bQ$ (giving rise to \eqref{eq:94.27} by
duality). Thus, $1$-connected rational homotopy types are expressible
as objects of the derived category defined in terms of such algebras,
and the obvious notion of quasi-isomorphism for these. To any space or
ss~set, Sullivan associates a corresponding \emph{``De Rham complex''}
with rational coefficients, in order to get a functor from rational
homotopy types to the derived category obtained from those algebras --
and (if I remember it right) this is an equivalence of categories,
provided one restricts to $1$-connected homotopy types, and
correspondingly to $1$-connected algebras. Probably somebody must have
explained to me by then (it was in 1976 more or less) why not every
eligible differential algebra could be recovered (up to isomorphism in
the derived category) by the corresponding cohomology algebra, namely
why it was not necessarily isomorphic to the latter, endowed with zero
differential operator; I am afraid I forgot it since! Also, it was
well-known by the informed people (as I was told too) that there where
obstructions against expressing the multiplicative structure in
cohomology with (say) integer coefficients, in terms of an
anti-commutative differential graded $\bZ$-algebra; so there was no
hope, I was informed, for defining something like a ``De Rham complex
with integer coefficients'' for an arbitrary topological space.

All this was very interesting indeed -- still, I found it hard to
believe that, while succeeding in constructing De Rham complexes with
rational coefficients for arbitrary spaces, by looking at the
algebraic De Rham complex on the enveloping affine space for the
various singular simplices of a space, that the same could not be
achieved with integral coefficients. Of course, the basic Poincaré
lemma for algebraic differential forms was no longer true, however
this reminded me strongly of a similar difficulty met with in
algebraic geometry, and which is overcome by working with suitable
``divided power structures'' -- as Poincaré's lemma becomes valid when
replacing usual polynomials (as coefficients for differential forms)
by ``polynomials with divided powers''. Then I got quite excited and
involved in a formalism of De Rham complexes with divided powers
for\pspage{341} arbitrary semisimplicial sets, which took me a few
weeks to work out and alongside getting back into homotopy and
cohomology formalism again. I had the feeling that this structure, or
the technically more adequate dual ``coalgebra'' structure, might well
turn out to be the more refined version of \eqref{eq:94.27} needed for
recovering homotopy types -- or at any rate $1$-connected ones. I gave
a talk about the matter at IHES while things were still hot in my mind
-- but it doesn't seem it went really through. It doesn't seem this
structure (which was worked out independently by someone else too, I
understand) has become a familiar notion to
topologists.\scrcomment{I'm guessing the reference is to \textcite{Cartan1976} and
  \textcite{Miller1978}\ldots anyhow, there has been many developments since
  concerning Witt vectors, crystalline cohomology, etc.\ldots} Maybe
one reason is that most topologists and homotopy theorists never
really got acquainted with the formalism of derived categories -- and
it seems that moreover, by the mid-seventies, it had even become
altogether unfashionable and ``mal vu'' to make any mention of them,
let alone work with them, also among some of the people who during
some time had been helping develop it. Now one of the main points I
was making in that talk was a somewhat delicate property of derived
categories of abelian categories, with respect to binomial
coefficients -- too bad!

I have not heard since about any work done in this direction I am
reflecting about now (somewhat retrospectively) -- namely recovering
homotopy types from their abelianization, \emph{plus} extra
structure. For all I know, \emph{the} relevant structure may well be
the differential algebra with divided powers structure embodied by the
De Rham complex (with a bigraduation however instead of just a
graduation), or its coalgebra version -- viewed as defining an object
of a suitable derived category. (Of course, when there is a
fundamental group $\pi$ around, one will have to look at a slightly
more complex structure still, involving operations of $\pi$, by
looking at the De Rham complex of the universal covering.) If it is
just the matter of describing homotopy types in terms of other models
than semisimplicial complexes, it must be admitted that the new models
are of incomparably more intricate description than the complexes!
There \emph{is} however one feature of it which greatly struck me by
that time, and still seems to me quite intriguing, namely \emph{that
  this structure, although definitely not ``abelian'' anymore} (due to
multiplication as well as to divided power structure), \emph{makes a
  sense over any commutative ground ring} (or even scheme, etc.). When
this ring is $\bQ$, the ``models'' we get modelize rational homotopy
types, which was the starting point of my reflections\pspage{342}
about seven years ago. Replacing $\bQ$ by a more general ring, this
suggests that \emph{there might exist a notion of ``homotopy types''
  over any ground ring $k$} -- and a corresponding notion of ground
ring extension for homotopy types. For abelianizations of homotopy
types, this is particularly ``obvious'', as being just the functorial
dependence of the derived category $\D_\bullet(k)$ with respect to the
ground ring $k$, corresponding to ring extension in a chain
complex. For a week or two I played around with this idea, which on
the semisimplicial level tied in with expressing homotopy types of
some simple spaces (such as standard $K(\pi,n)$ spaces and fibrations
between these) in terms of some simple \emph{semisimplicial schemes}
(affine and of finite type over $\Spec(\bZ)$), by taking $\bZ$-valued
points of these; ring extension $\bZ\to k$ was interpreted in the
scheme-theoretic sense.

I didn't go on very long, as soon after I was taken by personal
matters and never took up the matter later -- and maybe it was an
altogether unrealistic or silly attempt. If I remember it right, the
idea lurking was something of this kind, that \emph{there was a
  functor from \Hot{} to} (if not an equivalence of \Hot{} with\ldots)
\emph{a suitable derived category of some category of semisimplicial
  schemes over} $\Spec(\bZ)$, and that the base change intuition, as
suggested by the abelianized theory or by the subtler ``divided power
De Rham theory'', would reflect in naive base change $\bZ\to k$ for
schemes.

I was then looking mainly at $1$-connected structures, but there was
an idea too that nilpotent fundamental groups might fit into the
picture, with the hope that such a group (under suitable restrictions,
finite presentation and torsion freeness say) could be expressed in a
canonical way in terms of an affine nilpotent group scheme of finite
type over $\Spec(\bZ)$, by taking the integral points of the
latter. It seems (if I remember right) that this is not quite true
though -- that one couldn't hope for much better than getting a
nilpotent algebraic group scheme \emph{over $\bQ$} -- and that one
would recover the discrete group one started with only ``up to
commensurability''. Possibly, there may be \emph{an equivalence
  between localization of the category of nilpotent groups as above
  \textup(with respect to monomorphisms with image of finite
  index\textup) and affine nilpotent connected algebraic group schemes
  over $\bQ$}, or equivalently, group schemes whose underlying scheme
is isomorphic to standard affine space.

\bigbreak

\presectionfill\ondate{7.7.}\pspage{343}\par

\hangsection{Contractors}\label{sec:95}%
After this cascade of ``afterthoughts'' on abelianization of homotopy
types, it is time now to resume some more technical work, and get
through with this unending part \ref{ch:IV}, in accordance with the
short range working program I had come to four days ago (section
\ref{sec:92}, p.\ \ref{p:320}). I'll take up the three topics stated
there -- namely contractors, induced structures, and ``miscellaneous''
-- in that order, as reviewed previously. Thus, we'll start with
contractors. I have in mind now mainly the definition of contractors,
and a few basic facts following easily from what is already known to
us.

The first thought that comes to my mind is to define a contractor as a
category $A$ such that the set $\Ob(A)$ of all objects of $A$ is a
contractibility structure on $A$, i.e., that there exists a
contractibility structure on $A$ for which every object in $A$ is
contractible. The trouble with this definition is that it makes the
implicit assumption that $A$ is stable under finite products -- as the
notion of a contractibility structure was defined only in a category
satisfying this extra assumption (cf.\ section \ref{sec:51},
\ref{subsec:51.D}). Now, this assumption is \emph{not} satisfied by
the three standard test categories, including $\Simplex$, which surely
we do want to consider as contractors! The next thought then,
suggested by this reflection, is to embed $A$ into \Ahat{} to supply
the products which may be lacking in $A$, and demand there exist a
contractibility structure on \Ahat, such that the objects in $A$ be
contractible and moreover generate; or, what amounts to the same, that
for the homotopy interval structure on \Ahat{} defined by intervals
coming from $A$ as a generating family, the objects of $A$ are
contractible (which implies that this structure ``is'' indeed a
contractibility structure). This condition (in the more general case,
when $A$ appears as a full subcategory of any larger category $M$) has
been restated in wholly explicit terms as the \emph{``basic
  assumption''} \ref{cond:51.Bas.4} on a set of objects, in order that
it generate a contractibility structure (section \ref{sec:51}, p.\
\ref{p:118}). It is immediate that in the case when $A$ itself is
stable under finite products in the ambient category, that this
condition is intrinsic to $A$ and just amounts to the first definition
we had in mind.

Still,\pspage{344} we will call a category $A$ satisfying the
condition \ref{cond:51.Bas.4} with respect to the embedding
\begin{equation}
  \label{eq:95.1}
  A\hookrightarrow\Ahat
  \tag{1}
\end{equation}
a \emph{precontractor}, as we'll expect something more still from a
contractor, which will be automatically satisfied in the particular
case when $A$ is stable under finite products. Roughly speaking, we
want to have a satisfactory relation between contractibility and
asphericity in \Ahat{} -- we'll make this more precise below. For the
time being, let's dwell just a little more on the notion of a
precontractor.

A second thought about contractors, coming alongside with the first,
is that for any full embedding of $A$ into a larger category
\begin{equation}
  \label{eq:95.2}
  f:A\to M, \quad\text{$M$ stable under finite products,}\tag{2}
\end{equation}
$f(A)$ should generate in $M$ a contractibility structure. In the
particular case when $A$ is stable under finite products (and hence
the notion of a precontractor, already defined, coincides with the
notion of a contractor), this is indeed so provided $A$ is a
(pre)contractor, and moreover $f$ commutes to finite products. When
$A$ is just assumed to be a precontractor (without an assumption about
stability of $A$ under products), we'll assume in compensation that
$M$ is stable under small direct limits, which allows to take the
canonical extension $f_!$ of $f$ to \Ahat{} in a way commuting to
direct limits
\begin{equation}
  \label{eq:95.3}
  f_!:\Ahat\to M,\tag{3}
\end{equation}
and we can now state: \emph{if $f_!$ commutes to finite products}
(cf.\ prop.\ \ref{prop:85.1}, \ref{it:85.prop1.a}, p.\ \ref{p:281}),
\emph{then $f(A)$ generates a contractibility structure in $M$.} This
statement is true even without assuming that the functor $f$ is fully
faithful (and follows immediately from the criterion
\ref{cond:51.Bas.4} of p.\ \ref{p:118}); however, in the particular
case when $f$ is fully faithful, we have a handy criterion (prop.\
\ref{prop:85.2}, p.\ \ref{p:283}) for $f_!$ to commute to finite
products, namely that $f(A)$ be a strictly generating subcategory of
$M$, or equivalently, that the functor
\begin{equation}
  \label{eq:95.4}
  f^*:M\to\Ahat\tag{4}
\end{equation}
(right adjoint to $f_!$) be fully faithful. In this case, we may
identify $M$ (up to equivalence) to a full subcategory of \Ahat{}
containing $A$, and the fact that $f(A)$ generates a contractibility
structure in $M$ follows immediately directly (without having to rely
on existence of direct limits in $M$, nor even existence of $f_!$). To
sum up:
\begin{propositionnum}\label{prop:95.1}
  Let\pspage{345} $A$ be a small category, $M$ a category stable under
  finite products, $f:A\to M$ a functor, we assume $A$ is a
  \emph{precontractor}. Then $f(A)$ generates a contractibility
  structure in $M$ in each of the following three cases:
  \begin{enumerate}[label=\alph*),font=\normalfont]
  \item\label{it:95.prop1.a}
    $A$ stable under finite products, and $f$ commutes to these.
  \item\label{it:95.prop1.b}
    There exists a functor $f_!:\Ahat\to M$ extending $f$, and
    commuting to final object and binary products in \Ahat{} of
    objects in $A$\kern1pt.
  \item\label{it:95.prop1.c}
    The functor $f$ is fully faithful and strictly generating.
  \end{enumerate}
\end{propositionnum}

Of course, the validity of the conclusion in either case
\ref{it:95.prop1.b} or \ref{it:95.prop1.c}, for fixed $A$ and variable
$M$ and $f$, \emph{characterizes} the property for $A$ of being a
precontractor, and the same for \ref{it:95.prop1.a} if we assume
beforehand that $A$ is stable under finite products. Thus, we may view
the proposition \ref{prop:95.1} as the most comprehensive statement of
the meaning of this property.
\begin{propositionnum}\label{prop:95.2}
  Let $A$ be a precontractor. Let \Ahatc{} be the set of contractible
  objects in \Ahat{} for the contractibility structure generated by
  the subcategory $A$\kern1pt, \Ahatas{} \textup(resp.\ \Ahatlocas\textup) the
  set of aspheric \textup(resp.\ locally aspheric -- cf.\ p.\
  \ref{p:250}\textup) objects of \Ahat, $h$ the homotopy structure on
  \Ahat{} associated to the contractibility structure \Ahatc, i.e.,
  generated by the intervals in \Ahat{} coming from $A$\kern1pt. As usual,
  \scrWA{} denotes the set of weak equivalences in \Ahat{} -- it is
  understood here that the basic localizer $\scrW\subset\Fl\Cat$ is
  $\scrWoo$ = usual weak equivalence. The following conditions on $A$
  are equivalent:
  \begin{description}
  \item[\namedlabel{it:95.prop2.i}{(i)}]
    \Ahat{} is totally aspheric \textup(i.e., $\Ob A\subset\Ahatlocas$\textup).
  \item[\namedlabel{it:95.prop2.ii}{(ii)}]
    The asphericity structure \Ahatas{} on \Ahat{} is generated by the
    contractibility structure \Ahatc.
  \item[\namedlabel{it:95.prop2.iiprime}{(ii')}]
    $\Ahatc\subset\Ahatas$.
  \item[\namedlabel{it:95.prop2.iidblprime}{(ii'')}]
    $\Ahatc\subset\Ahatlocas$.
  \item[\namedlabel{it:95.prop2.iii}{(iii)}]
    Any $h$-homotopism is in \scrWA{} \textup(i.e., \scrWA{} is
    ``strictly compatible'' with the homotopy structure $h$,
    \textup(cf.\ section \ref{sec:54}\textup), i.e., $h\le
    h'=h_{\scrWA}$\textup).
  \item[\namedlabel{it:95.prop2.iv}{(iv)}]
    The homotopy structure $h$ is equal to the homotopy structure
    $h'=h_W$ associated to $W=\scrWA$ \textup(cf.\ section
    \ref{sec:54}\textup). 
  \end{description}
\end{propositionnum}
\begin{proof}[Proof of proposition]
  Immediate from what is known to us, via
  \[\text{\ref{it:95.prop2.i}} \Rightarrow
  \text{\ref{it:95.prop2.ii}} \Rightarrow
  \text{\ref{it:95.prop2.iiprime}} \Leftrightarrow
  \text{\ref{it:95.prop2.iidblprime}} \Rightarrow
  \text{\ref{it:95.prop2.i}}
  \quad\text{and}\quad
  \text{\ref{it:95.prop2.ii}} \Rightarrow
  \text{\ref{it:95.prop2.iv}} \Rightarrow
  \text{\ref{it:95.prop2.iii}} \Rightarrow
  \text{\ref{it:95.prop2.i}.}\]
\end{proof}
\begin{definitionnum}\label{def:95.1}
  A\pspage{346} small category $A$ is called a \emph{contractor} if it
  is a precontractor, and if moreover it is totally aspheric or,
  equivalently, satisfies one of the equivalent condition
  \ref{it:95.prop2.i} to \ref{it:95.prop2.iv} of prop.\
  \ref{prop:95.2}.
\end{definitionnum}

Equivalently, this also means that
\begin{equation}
  \label{eq:95.5}
  \Ob A \subset \Ahatc\sand \Ahatlocas,\tag{5}
\end{equation}
i.e., every object in $A$ is contractible and locally aspheric, where
the set \Ahatc{} of ``contractible'' objects of \Ahat{} is defined in
terms of the homotopy interval structure $h$ generated by all
intervals in \Ahat{} coming from objects in $A$. (Thus structure is
not necessarily a contractibility structure, but it is when $A$ is a
precontractor, namely $\Ob A\subset\Ahatc$.)

The most trivial example of a contractor is the final category
$\Simplex_0$, and more generally, any category equivalent to it
(NB\enspace a category equivalent to a precontractor resp.\ to a
contractor is again a precontractor resp.\ a contractor). Such a
contractor will be called \emph{trivial}. For a trivial contractor
$A$, we get an equivalence
\[\Ahat\equeq\Sets.\]

If
\[ (M,M\subc)\]
is a contractibility structure, and $A\subset M$ any small full
subcategory of $M$ generating the contractibility structure, then $A$
is a precontractor, hence a contractor if{f} $A$ is totally aspheric,
which will be the case if $A$ is stable in $M$ under binary products,
a fortiori if it is stable under finite products, i.e., contains
moreover a final object of $M$. Thus, \emph{the contractibility
  structure of $M$ can always be generated by a full subcategory $A$
  of $M$ which is a contractor}.

Apart from these two examples, the most interesting examples of
contractors are of course the three standard test categories
$\Simplex$, $\Square$ and $\Globe$, and also their finite
products. Note that the notion of a precontractor or of a contractor
is clearly stable under finite products.
\begin{propositionnum}\label{prop:95.3}
  Let $A$ be a precontractor, assume $A$ non-trivial, i.e.,
  non-equivalent to the final category. Then $A$ contains a separating
  interval, and hence it is a \textup(strict\textup) test category if
  $A$ is tot.\ asph., i.e., is a contractor.
\end{propositionnum}
\begin{proof}
  As any object of $A$ has a section (over $e_\Ahat$), it follows
  immediately that any non-empty object of \Ahat{} has a section too,
  hence \emph{any non-empty subobject of $e_\Ahat$ is equal to
    $e_\Ahat$}. This implies that for any\pspage{347} interval
  \[\bI=(I,\delta_0,\delta_1)\]
  in \Ahat, either \bI{} is separating, i.e.,
  $\Ker(\delta_0,\delta_1)$ is the empty object of \Ahat, or
  $\delta_0=\delta_1$. If no interval coming from $A$ was separating,
  this then would just mean that any two sections (over $e_\Ahat$) of
  an object $I$ of $A$ are equal. By the definition of the homotopy
  structure in \Ahat{} generated by these intervals, this would imply
  that any homotopism in \Ahat{} is an isomorphism, and hence that any
  contractible object for this structure is isomorphic to the final
  object $e_\Ahat$. As by assumption on $A$ all objects of $A$ are
  contractible, this would mean that $A$ is trivial, which is against
  our assumptions, qed.
\end{proof}
\begin{corollarynum}\label{cor:95.1}
  Let $(M,M\subc)$ be a contractibility structure.
  \begin{enumerate}[label=\alph*),font=\normalfont]
  \item\label{it:95.cor1.a}
    The following conditions are equivalent \textup(and will be
    expressed by saying that this given contractibility structure is
    \emph{trivial}\textup):
    \begin{enumerate}[label=(\roman*),font=\normalfont]
    \item\label{it:95.cor1.a.i}
      Two maps in $M$ which are homotopic are equal.
    \item\label{it:95.cor1.a.ii}
      Any homotopism in $M$ is an isomorphism.
    \item\label{it:95.cor1.a.iii}
      Any homotopy interval in $M$ is ``trivial'', i.e., any two
      homotopic sections of an object of $M$ are equal.
    \item\label{it:95.cor1.a.iv}
      Any contractible object of $M$ is a final object, i.e., $M\subc$
      is just the set of all final objects of $M$.
    \item\label{it:95.cor1.a.v}
      Any two sections of a contractible object are equal.
    \end{enumerate}
  \item\label{it:95.cor1.b}
    Assume\pspage{348} the contractibility structure $M\subc$
    non-trivial, i.e., there exists an interval
    \[\bI=(I,\delta_0,\delta_1)\quad\text{with}\quad I\in M\subc,
    \delta_0\ne\delta_1.\]
    Then for any small category $A$ and any functor
    \[i:A\to M\]
    factoring through $M\subc$, the interval $i^*(\bI)$ in \Ahat{} is
    separating. Hence if \textup(for a given basic localizer
    \scrW\textup) $i$ is totally \scrW-aspheric \textup(cf.\ theorem
    \ref{thm:79.1} cor.\ \ref{cor:79.1} p.\ \ref{p:252}\textup), hence
    $i^*(I)$ is totally aspheric in \Ahat, then $A$ is a \scrW-test
    category. In particular, if $A$ is totally \scrW-aspheric and $i$
    is $M\suba$-\scrW-aspheric, then $A$ is a strict \scrW-test category.
  \end{enumerate}
\end{corollarynum}
\begin{proof}
  Part \ref{it:95.cor1.a} is a tautology in terms of section
  \ref{sec:51}. For part \ref{it:95.cor1.b}, to prove that $i^*(\bI)$
  is separating, we only have to check that for any $a$ in $A$, the
  two compositions
  \[\begin{tikzcd}[cramped,sep=huge]
    a \to i^*(e_M) \simeq e_\Ahat
    \ar[r,shift left=1pt,"\ensuremath{i^*(\delta_0),i^*(\delta_1)}"]\ar[r,shift right=2pt] &
    i^*(I)
  \end{tikzcd}\]
  are distinct, or what amounts to the same by the definition of
  $i^*$, that the compositions
  \[\begin{tikzcd}[cramped,sep=small]
    i_!(a) \eqdef x \to e_M
    \ar[r,shift left=2pt]\ar[r,shift right=2pt] &
    I
  \end{tikzcd}\]
  are distinct. As $x$ is in $M\subc$, it has a section over $e_M$, so
  it is enough to check that the compositions with $e_M\to x$ are
  distinct, which just means that $\delta_0\ne\delta_1$, qed.
\end{proof}
\begin{remarks}
  \namedlabel{rem:95.1}{1)}\enspace Part \ref{it:95.cor1.b} of the
  corollary replaces cor.\ \ref{cor:79.3} on page \ref{p:253}, which
  is a little monster of incongruity (as I just discovered) -- namely,
  two of the assumption on \bI{} made there (namely that \bI{} be a
  multiplicative interval, and $I\in\Ob C$) are useless if we assume
  just $I$ contractible, moreover the awkward separation assumption
  made there just reduces by the trivial argument above to the
  assumption $\delta_0\ne\delta_1$!

  \namedlabel{rem:95.2}{2)}\enspace It should be noted that the
  homotopy structure on \Ahat{} envisioned in th.\ \ref{thm:79.1} of
  section \ref{sec:79} (p.\ \ref{p:252}) is \emph{not} defined as in
  the present section, in terms of intervals in \Ahat{} coming from
  $A$ (call this structure $h$), but as
  \[ h' = h_\scrWA\]
  defined in terms of intervals in $\Ahatlocas$; this depends a priori
  on the choice of \scrW, as it has to because th.\ \ref{thm:79.1}
  gives a criterion for the functor $i$ to be $M\suba$-\scrW-aspheric
  which does depend on \scrW. (It surely won't be the same if we take
  $\scrW=\scrWoo=$ usual weak equivalence, or $\scrW=\Fl(\Cat)$ hence
  $\scrWA=\Fl(A)$ and\pspage{349} the condition that $i$ be
  $M\suba$-\scrW-aspheric is always satisfied!) However, let's assume
  $A$ to be totally \scrW-aspheric and every object of $A$ has a
  section (over the final object $e_\Ahat$ of \Ahat) or, what amounts
  to the same, every ``non-empty'' object of \Ahat{} has a section --
  we'll say in this case $A$ is \emph{``strictly totally
    \scrW-aspheric''} (compare section \ref{sec:60}, p.\ \ref{p:149},
  in the particular case $\scrW=\scrWoo$, and with \Ahat{} replaced by
  an arbitrary topos). Let's assume moreover that \scrW{} satisfied
  \ref{loc:4}. In this case, the homotopy structure $h'=h_\scrWA$ does
  not depend on the choice of \scrW, namely it is the so-called
  ``canonical homotopy structure''
  \[h'' = h_\Ahat\]
  of the (strictly totally $0$-connected) category \Ahat{} (cf.\
  section \ref{sec:57}), which in the special case of a category
  \Ahat{} can also be defined as the homotopy structure $h_\scrWz$
  associated to ${\scrWz}_A$, where \scrWz{} is the coarsest basic
  localizer satisfying \ref{loc:4}, i.e.,
  \[\scrWz=\set[\big]{f\in\Fl\Cat}{\text{$\piz(f)$ bijective}}.\]
  The proof of this fact $h'=h''$, i.e.,
  \[h_\scrWA = h_\Ahat \quad ({}=h_{\scrWz})\]
  is essentially the same as for the similar prop.\ (section
  \ref{sec:60}, p.\ \ref{p:149}). The condition \ref{loc:4} on \scrW,
  i.e., $\scrW\subset\scrWz$ clearly implies
  \[h_\scrWA\subset h_{{\scrWz}_A},\]
  and to get the opposite inequality, for which we'll use the
  assumption on $A$, we only have to prove that for any $0$-connected
  object $K$ of \Ahat, any two sections are ($h'=h_\scrWA$)-homotopic,
  a fortiori (as $h\le h'$ by the assumption of total
  \scrW-asphericity of $A$) it is enough to prove they are
  $h$-homotopic. Now this follows from lemma \ref{lem:82.2}, p.\
  \ref{p:268}, applied to $\scrC=\Ahat$, $C=A$.
\end{remarks}

The\pspage{350} preceding reflections thus prove the following
afterthought to theorem \ref{thm:79.1} of section \ref{sec:79}:
\begin{propositionnum}\label{prop:95.4}
  Let $(M,M\subc)$ be a contractibility structure, $A$ a small
  category, $i:A\to M$ a functor factoring through $M\subc$. Let
  moreover \scrW{} be a basic localizer satisfying
  \textup{\ref{loc:4}}. We assume $A$ \emph{strictly totally
    \scrW-aspheric}, i.e., totally \scrW-aspheric and moreover any
  object of $A$ has a section \textup(over $e_\Ahat$\textup).
  \begin{enumerate}[label=\alph*),font=\normalfont]
  \item\label{it:95.prop4.a}
    The homotopy structure $h_\scrWA$ on \Ahat{} is equal to the
    canonical homotopy structure $h_\Ahat$ defined by $0$-connected
    intervals, and equal also to the homotopy structure $h$ defined by
    intervals coming from $A$:
    \begin{equation}
      \label{eq:95.6}
      h=h_\scrWA=h_\Ahat.\tag{6}
    \end{equation}
  \end{enumerate}
  In what follows, we assume \Ahat{} endowed with this homotopy
  structure, and denote by \Ahatc{} the set of all contractible object
  in \Ahat. We equally endow \Ahat{} with its canonical
  \scrW-asphericity structure, and $M$ with the \scrW-asphericity
  structure associated to its contractibility structure $M\subc$. With
  these conventions:
  \begin{enumerate}[label=\alph*),font=\normalfont,start=2]
  \item\label{it:95.prop4.b}
    The following conditions on $i$ are equivalent, where
    \[i^*:M\to\Ahat\]
    is the functor defined as usual in terms of $i$:
    \begin{enumerate}[label=(\roman*),font=\normalfont]
    \item\label{it:95.prop4.b.i}
      $i^*$ is compatible with the homotopy structures \textup(cf.\
      criteria on pages \ref{p:251}--\ref{p:252}\textup), which can be
      expressed also by
      \[i^*(M\subc)\subset\Ahatc\]
      \textup(a condition independent from \scrW, in view of
      \textup{\ref{it:95.prop4.a})}.
    \item\label{it:95.prop4.b.ii}
      $i$ is \scrW-aspheric, i.e.,
      \[M_\scrW \subset (i^*)^{-1}(A\uphat_\scrW)\]
      \textup(where $M_\scrW$ and $A\uphat_\scrW$ are the sets of
      \scrW-aspheric objects in $M$ and in \Ahat\textup).
    \item\label{it:95.prop4.b.iii}
      \textup(For a given full subcategory $C$ of $M$ generating the
      contractibility structure $M\subc$\textup):
      \[i^*(C)\subset A\uphat_\scrWz = \text{set of $0$-connected
        objects of \Ahat.}\]
    \end{enumerate}
  \item\label{it:95.prop4.c}
    Assume these conditions hold, and moreover that the
    contractibility structure of $M$ is non-trivial. Then $A$ is a
    strict \scrW-test category.
  \end{enumerate}
\end{propositionnum}
\begin{proof}
  Part\pspage{351} \ref{it:95.prop4.a} has been proved in remark
  \ref{rem:95.2} above, and in view of th.\ \ref{thm:79.1}, p.\
  \ref{p:252}, the equivalence of \ref{it:95.prop4.b.i} and
  \ref{it:95.prop4.b.ii} is clear, hence also the equivalence with
  \ref{it:95.prop4.b.iii} by applying loc.\ cit.\ to \scrWz{} instead
  of \scrW. Part \ref{it:95.prop4.c} now follows from prop.\
  \ref{prop:95.3} cor.\ \ref{cor:95.1} \ref{it:95.cor1.b}.
\end{proof}
\begin{corollary}
  Under the conditions of \textup{\ref{it:95.prop4.c}} above, if
  $M\subc$ is \scrW-modelizing, then $i$ is a \scrW-test functor, and
  induces an \emph{equivalence}
  \[\HotOf_{M,\scrW} \eqdef \scrW_M^{-1} M \tosimeq \HotOf_{A,\scrW}
  \eqdef \scrWA^{-1}\Ahat.\]
\end{corollary}

\bigbreak

\presectionfill\ondate{8.7.}\par

\hangsection{Vertical and horizontal topoi\dots
  \texorpdfstring{\textup(afterthought on
    terminology\textup)}{(afterthought on
    terminology)}.}\label{sec:96}%
Yesterday's notes have proceeded very falteringly, to my surprise,
while everything seemed ready for smooth sailing. A number of times,
after going on for a page or two ``following my nose'' (as they say in
German),\scrcomment{I'm pretty sure that ``to follow one's nose'' is an
  English idiom in the sense of following one's instinct, while the
  meaning of going straight ahead is shared with the German ``der Nase nach
  gehen''} or for half a page, it turned out it just wasn't right that
way and I would feel quite stupid and put the silly pages away as
scratchpaper and have another start. There wouldn't have been any
point dragging the poor reader (if there is still one left\ldots)
along on my stumbling path, where it was a matter merely of getting
some technical adjustments right. Maybe it is just that attention was
distracted, perhaps precisely through this (partly mistaken, and
anyhow not too inspiring) feeling that everything was kind of cooked
already, and what was left to do was just swallow! What came out in
the process was that finally things were not so clear yet in my mind
as I thought they were. It is a frequent experience that whenever one
wants to go ahead too quickly, one finds oneself dispersing stupidly a
hell of a lot of energy\ldots

There occurred to me some inadequacies with terminology. One is about
the property of certain categories (contractor or precontractors for
instance) that every object of $A$ has a section (over $e_\Ahat$),
which can be viewed also as a property of the topos $\Ahat=\scrA$,
namely that any ``non-empty'' object of the topos has a section. This
is immediately seen (for any given topos \scrA) to imply the property
that the final object $e_\scrA$ has only the two trivial subobjects,
the ``empty'' and the ``full'' one -- or equivalently, that any
subtopos of the topos is either the empty of the full one, -- a
property, too of obvious geometric significance. In case of a topos of
the type \Ahat, one immediately sees the two properties are equivalent
-- but this is not\pspage{352} true for an arbitrary topos: for
instance the classifying topos $B_G$ of a discrete group $G$ has the
second property, but visibly not the first unless $G$ is the unit
group. (I recall that the category of sheaves on $B_G$ is the category
\Gsets{} of sets on which $G$ operates.) I feel both properties for a
topos merit a name. The first (every sheaf has a section) can be
viewed as the strongest conceivable (I would think) global asphericity
property for a topos, as far as $\mathrm H^1$ goes at any rate, as the
$\mathrm H^1$ of $X$ with coefficients in \emph{any} group object will
be zero. (But I confess I didn't try and look if any precontractor,
say, is aspheric\ldots) The second property (every subtopos is
trivial) comes with a rather different flavor, it suggests the image
of just one ``point'' -- and as a matter of fact, the étale topos of a
scheme, say has this property if{f} it is reduced to a point. Such a
topos may called ``punctual'' (not to be confused though with some
other meanings suggested by this word, such as being equivalent to the
topos defined by a one-point topological space, namely \scrA{} being
equivalent to \Sets) or ``atomic'' (which has rather unpleasant
connotations though nowadays!), or maybe ``vertical'' (this image is
suggested by the $B_G$ above) -- the ``base'', i.e., the final object
of \scrA{} being very ``small'' (in terms of harboring subobjects), so
the inner structure is expressed like a kind of tower, related (in the
case of $B_G$) to the ``Galois tower'' of subgroups of $G$\ldots The
corresponding notion of a ``horizontal'' topos is visibly the one when
\scrA{} admits the subobjects of the final object as a generating
family. In terms of these definitions, a topos is horizontal and
vertical if{f} it is either the ``empty'' or the ``final'' (or
``one-point'') topos. This brings to mind that in the notion of
verticality, we should exclude the ``empty'' topos (which formally
satisfies the condition -- every sheaf has a section). This brings to
my attention too that I certainly do not want to consider an empty
category $A$ (defining the ``empty'' topos \Ahat) as a precontractor,
although formally (in terms of yesterday's definition) it is. Thus, I
suggest I'll introduce the following
\begin{definition}
  A topos is called \emph{vertical} if it is not an ``empty'' topos
  (i.e., the category of sheaves on it is not equivalent to the final
  category $\Simplex_0$), and if moreover any open subtopos is either
  the ``empty'' or the ``full'' one (hence the same for any subtopos,
  whether open or not). A small category $A$ is called vertical, if
  the associated topos (with category of sheaves \Ahat)
  is,\pspage{353} or equivalently, if $A$ is non-empty and any object
  of $A$ has a section (over $e_\Ahat$). A topos is called
  \emph{horizontal} if the family of all subobjects of the final
  object in the category of sheaves \scrA{} is generating.
\end{definition}

For instance, the topos associated to a topological space is
horizontal -- in particular, an ``empty'' topos is horizontal. A topos
is both horizontal and vertical if{f} it is a ($2$-)final topos, i.e.,
equivalent to the topos defined by a one-point topological space
(i.e., the category of sheaves is equivalent to \Sets).

The property of verticality, I feel, is of interest in its own right,
as exemplified notably by lemma \ref{lem:82.2} p.\ \ref{p:268} (which
we used yesterday), and the related proposition of section
\ref{sec:60} (p.\ \ref{p:149}). It does not seem at all subordinated
to notions such as total asphericity or total $0$-connectedness, and
goes in an entirely different direction -- thus, the terminology
``\emph{strictly} totally aspheric'' (or totally $0$-connected), which
I still used yesterday (hesitatingly, I should say), is definitely
inadequate. I would rather say ``totally aspheric (or totally
$0$-connected) \emph{and} vertical''.

Another point is about the terminology of \emph{totally aspheric} and
\emph{locally aspheric} objects in a category \Ahat{} (with respect to
a given basic localizer \scrW), introduced in section \ref{sec:79}
(p.\ \ref{p:250}), and still used yesterday. This terminology does not
seem inadequate by itself, I introduced it because it struck me as
suggestive (and the notions it refers to do deserve a name, in order
to be at ease). The trouble here is that it conflicts with another
possible meaning, in accordance with the principle insisted upon
forcefully in the reflections of section \ref{sec:66} -- namely that
for objects or arrows in \Cat, or within a category \Ahat, the
terminology used for naming properties for these should be in
accordance with the terminology used for the corresponding topoi or
maps of topoi. Now, we do have already the notions of a locally
aspheric and totally aspheric topos, which therefore should imply
automatically the meaning of these notions for an object of \Cat{}
(which was done satisfactorily months ago), or for an object of a
category \Ahat. But in the latter case, there is definitely conflict
with the terminology introduced on p.\ \ref{p:252}. This conflict has
not manifested itself yet in any concrete situation, while the
unorthodox terminology has been used quite satisfactorily a number of
times. Therefore, I would like to keep it, as long as I am not forced
otherwise.

\starsbreak

We\pspage{354} were faced yesterday with three different homotopy
structures $h,h',h''$ on a category \Ahat, for a given small category
$A$, which make sense for any $A$, and which in case $A$ is a
contractor all coincide. The exact relationship between these
structures in more general cases has remained somewhat confused, and
in order to dispel the resulting feeling of uneasiness, I took finally
the trouble today to write it out with some case. One of these
structures, $h'$, depends on the choice of a basic localizer \scrW,
whereas the two others don't\ldots

\bigbreak
\presectionfill\ondate{12.7.}\par

\hangsection[``Projective'' topoi. Morphisms and bimorphisms of
\dots]{``Projective'' topoi. Morphisms and bimorphisms of
  contractors.}\label{sec:97}%
I was interrupted in my notes by visiting friends arriving in close
succession -- then since yesterday I have been busy mainly with letter
writing. Now, I am ready to take up the thread where I left it --
namely some afterthoughts to the reflections of section \ref{sec:95}
on contractors.

First an afterthought to the afterthoughts! I had introduced the name
``\emph{vertical} topos'' for a topos admitting only the two trivial
open subtopoi (page \ref{p:352}), whereas the stronger property that
every ``non-empty'' sheaf has a section remained unnamed (which is no
real drawback as long as we are restricting to topoi of the type
\Ahat, where indeed the two notions coincide). Now, the latter
property can be viewed as the property that every sheaf $F$ \emph{such
  that $F\to e$ be epimorphic}, should admit a section. It is this
last property which does merit to ``be viewed as the strongest
conceivable asphericity property for a topos'' as I commented on it
last Friday (p.\ \ref{p:352}). After I had written this down as a kind
of selfevidence, a doubt turned up though and I qualified the comment
by added ``as far as $\mathrm H^1$ goes at any rate, as the $\mathrm
H^1$ with coefficients in \emph{any} group object will be zero''. I
didn't pause then to see if the doubt was founded -- quite evidently
it isn't, except for $\mathrm H^0$, as it is clear by the usual shift
argument, using embedding of an abelian sheaf into an injective one,
that if for given $k$ (here $k=1$) $\mathrm H^k(X,F)=0$ for any
abelian sheaf $F$, then the same holds for $\mathrm H^n$ with any
$n\ge k$ -- i.e., the global cohomological dimension of $X$ is
$<k$. This implies that any small vertical category (a fortiori any
\emph{precontractor}) is aspheric, provided it is $0$-connected,
indeed its cohomology variants with values in \emph{any} sheaf of
coefficients (not necessarily commutative as far as $\mathrm H^1$
goes) are trivial.

The\pspage{355} property for a topos $X$, with category of sheaves
\scrA, that any object $F$ in \scrA{} covering the final object
$e_\scrA$ should have a section, can be expressed by saying that the
latter is a \emph{projective} object in the category \scrA. Following
the principle to use the same names for properties of a topos, and
corresponding properties of the final sheaf on it, we may call a topos
with the above property a \emph{projective topos}. Thus, the
``non-empty'' topoi such that every ``non-empty'' sheaf has a section,
are exactly the topoi which are both vertical and projective.

\starsbreak

Here is the promised ``exact relationship'' between the three standard
homotopy structures $h,h',h''$ on \Ahat, where $A$ is any small
category (cf.\ end of section \ref{sec:96}, p.\ \ref{p:354}):
\begin{equation}
  \label{eq:97.star}
  \begin{tabular}{@{}c@{}}
    \begin{tikzcd}[baseline=(O.base),row sep=small,column sep=large]
      h \ar[d,equal]\ar[r,phantom,"\le"{description}]
      \ar[r,invisible,"(\text{$A$ tot. \scrW-asph.})"{inner sep=1.2ex}]
      &
      h' \ar[d,equal]\ar[r,phantom,"\le"{description}]
      \ar[r,invisible,"(\scrW\subset\scrWz)"{inner sep=1.2ex}]
      &
      h'' \ar[d,equal]\ar[r,phantom,"\le"{description}]
      \ar[r,invisible,"(\text{$A$ vertical})"{inner sep=1.2ex}]
      & h \\
      |[alias=O]| h_{\Ahat\!,\;A} & h_\scrWA & h_{{\scrWz}_A}=h_\Ahat &
    \end{tikzcd},
  \end{tabular}\tag{*}
\end{equation}
where \scrW{} is a given basic localizer. Above each one of the three
conditional inequalities between homotopy structures $h,h',h''$ I
wrote the natural assumption on $A$ or \scrW{} validating it, and in
the diagram I have recalled the definition of the three homotopy
structures. Apropos the description $h=h_{\Ahat\!,\;A}$, the notation
used here is $h_{M,A}$ when $M$ is a category stable under finite
products and $A$ a full subcategory, for designating the homotopy
structure on $M$ generated by intervals in $M$ coming from
$A$. Apropos $h''=h_\Ahat$, I recall the notation $h_M$ for
designating the canonical homotopy structure on a category $M$
satisfying suitable conditions (section \ref{sec:57}). Also, I recall
\[\scrWz =\set[\big]{f\in\Fl\Cat}{\text{$\piz(f)$ bijective}}.\]
The diagram implies that if $\scrW\subset\scrWz$, i.e., \scrW{}
satisfies \ref{loc:4}, and if moreover $A$ is vertical and totally
\scrW-aspheric, then all three homotopy structures coincide. Also,
taking $\scrW=\scrWz$, we see that $h=h''$ if $A$ is vertical and
totally $0$-connected, which is lemma \ref{lem:82.2} of p.\
\ref{p:268} for $\Ahat$, $A$.

From \eqref{eq:97.star} it follows of course that if $A$ is a
contractor, then (for any \scrW{} satisfying \ref{loc:4}) the three
homotopy structures $h, h', h''$ on $\Ahat$ coincide. In case the
contractor $A$ is not trivial, hence $A$ is\pspage{356} a strict test
category and \Ahat{} is \scrW-modelizing, it follows that \Ahat{} is
even a \emph{canonical modelizer} (with respect to \scrW), i.e.,
defined in terms of the \scrW-asphericity structure associated to the
``canonical'' homotopy structure $h''=h_\Ahat$ on \Ahat{} (cf.\ prop.\
\ref{prop:95.2} \ref{it:95.prop2.ii} p.\ \ref{p:345}). These, for the
time being, together with the modelizers \Cat{} and \Spaces, are the
main examples we got of canonical modelizers. Presumably, stacks
should give another sizable bunch of canonical modelizers, not of the
type \Ahat.

\starsbreak

We still have to say a word about morphisms between contractors $A,
B$. The first thing that comes to my mind is that this should be a
functor
\begin{equation}
  \label{eq:97.0}
  f:A\to B\tag{0}
\end{equation}
such that the corresponding functor
\begin{equation}
  \label{eq:97.1}
  f^*:\Bhat\to\Ahat\tag{1}
\end{equation}
should be compatible with the homotopy structures, which can be
expressed, as we know, in manifold ways, the most natural one here
being the following two
\begin{equation}
  \label{eq:97.2}
  f^*(B) \subset \Ahatc\tag{2}
\end{equation}
or
\begin{equation}
  \label{eq:97.3}
  f^*(\Bhatc) \subset \Ahatc,\tag{3}
\end{equation}
which are both implied by the apparently weaker one
\begin{equation}
  \label{eq:97.4}
  f^*(B) \subset A\uphat_{\scrWz} \eqdef \text{set of $0$-connected
    objects of \Ahat,}\tag{4}
\end{equation}
and equivalently still, as $\Ahatc\subset A\uphat_\scrW\subset
A\uphat_\scrWz$ (where \scrW{} is a basic localizer satisfying
\ref{loc:4}), to the condition
\begin{equation}
  \label{eq:97.5}
  f^*(B)\subset A\uphat_\scrW \quad ( \eqdef \text{set of
    \scrW-aspheric objects of \Ahat}).\tag{5}
\end{equation}
Thus, the condition for $f$ to be a ``morphism of contractors'' just
boils down to the long familiar \emph{\scrW-asphericity} of $f$, and
implies the following relation, apparently stronger than
\eqref{eq:97.5}:
\begin{equation}
  \label{eq:97.5prime}
  B\uphat_\scrW=(f^*)^{-1}(A\uphat_\scrW).\tag{5'}
\end{equation}
This comes almost as a surprise (after a four day interruption in
contact with the stuff!) -- but it occurs to me now that we got
already a more general statement with prop.\ \ref{prop:95.4} of
section \ref{sec:95} (p.\ \ref{p:350}), which includes the
situation\pspage{357} when instead of $f:A\to B$, we got a functor
\begin{equation}
  \label{eq:97.6}
  f:A\to \Bhat, \quad\text{factoring through \Bhatc}\tag{6}
\end{equation}
(which need not factor through $B$), or equivalently a functor
\begin{equation}
  \label{eq:97.7}
  f_!:\Ahat\to\Bhat\tag{7}
\end{equation}
commuting with small direct limits, or equivalently still, a functor
$f^*$ in opposite direction, commuting with small inverse limits, but
in the last two cases with the extra condition that $f_!(A)\subset
\Bhatc$. We may want to extend the notion of morphism of contractors
to include this situation, hence expressed by the two following
conditions on a functor $f$ \eqref{eq:97.6} or $f_!$ \eqref{eq:97.7},
or on the pair $(f_!,f^*)$ of adjoint functors
\begin{equation}
  \label{eq:97.8}
  \begin{cases}
    f_!(A)\subset\Bhatc &\\
    f^*(B)\subset\Ahatc &\quad.
  \end{cases}\tag{8}
\end{equation}
However, in order for this notion to be stable under composition, we
should strengthen the first of the relations \eqref{eq:97.8} into
\begin{equation}
  \label{eq:97.9}
  f_!(\Ahatc)\subset\Bhatc,\tag{9}
\end{equation}
which follows automatically whenever $f_!$ commutes to finite products
(cf.\ section \ref{sec:85}), but may not follow from \eqref{eq:97.8}
in general, even in the case when $f$ factors through $B$, i.e., in
the case we start with a functor \eqref{eq:97.0} $f:A\to B$. Thus, we
get \emph{two} plausible notions of a morphism of contractors, neither
of which implies the other, and I feel unable to predict which one
will prove the more useful. As far as terminology goes, it seems
reasonable to reserve the name ``morphism of contractors'' to the
first notion, as the second is adequately characterized as a
\emph{bimorphism} between the contractibility structures
$(\Ahat,\Ahatc)$ and $(\Bhat,\Bhatc)$ defined by the contractors $A$
and $B$ (cf.\ section \ref{sec:86}).

\bigbreak
\presectionfill\ondate{16.7.}\par

\hangsection[Sketch of proof of \protect\smashSimplexf{} being a weak test
category -- and \dots]{Sketch of proof of
  \texorpdfstring{\protect\Simplexf{}}{Delta-f} being a weak test
  category -- and perplexities about its being
  aspheric!}\label{sec:98}%
Next point on my provisional program is induced structures
(asphericity of contractibility structures) on a category $M_{/a}$,
when one is given on $M$ -- but finally I decided to skip this, as
there was no urgent need for clarifying this and I am not writing a
treatise, thanks Gods! I felt more interested writing down the proof
that the category \Simplexf{} of standard ordered simplices without
degeneracies is a weak test category, as announced months ago, in
section \ref{sec:43}. It then\pspage{358} seemed to come out rather
simply, but I didn't keep notes of the proof I thought I found, which
caused me spending now a day or two feeling a little stupid, as the
stuff was resisting while I felt it shouldn't! It did come out in the
end I guess -- and still I feel a little stupid, with the impression
of having bypassed definitely some very simple argument which had
presented itself as a matter of evidence by the end of March. On the
other hand, I was led to reflect on some other noteworthy features of
the situation, so I don't feel I altogether have been loosing my time.

There are four main variants of categories of standard simplices,
inserting into a diagram
\begin{equation}
  \label{eq:98.1}
  \begin{tabular}{@{}c@{}}
    \begin{tikzcd}[baseline=(O.base)]
      \Simplex\ar[r,"\beta"] & \Simplextilde \\
      \Simplexf\ar[u,"\alpha"]\ar[r,"\beta^{\mathrm f}"] &
      |[alias=O]| \Simplextildef\ar[u,"\widetilde\alpha"] 
    \end{tikzcd},
  \end{tabular}\tag{1}
\end{equation}
where \Simplextilde{} denotes the category of non-ordered standard
simplices $\Simplex_n$, and where the exponent $f$ in \Simplexf{} and
\Simplextildef{} denotes restriction to maps which are injective,
namely compositions of boundary maps (plus symmetric in the
non-ordered case). I recall that $\Simplex$ and \Simplextilde{} are
contractors, wheres \Simplexf{} and \Simplextildef{} are not even test
categories. We'll see however that \Simplexf{} is a \emph{weak} test
category, and presumably the same kind of argument should apply to
prove that \Simplextildef{} is a weak test category too. On the other
hand, from the point of view of the modelizing story, the main common
property of the four functors in \eqref{eq:98.1} should be
asphericity. However, I checked this for $\alpha$ and $\beta$ only, as
this was all I needed for getting the desired result on \Simplexf. As
a matter of fact, $\beta$ is even better than being aspheric, it is a
\emph{morphism of contractors}; more precisely still, for any object
$E$ in \Simplextilde, namely essentially a finite non-empty set,
choosing one point $a$ in $E$, one easily constructs an elementary
homotopy for $\beta^*(E)$ from the identity map to the constant map
defined by $\beta^*(a)$. I do not know on the other hand whether
$\beta$ defines a bimorphism between the canonical contractibility
structures on \Simplexhat{} and \Simplextildehat, namely whether
$\beta_!$ transforms contractible elements into contractible ones
(which would be clear if we knew that $\beta_!$ commutes to finite
products).

We'll come back upon proof of asphericity of $\alpha$, and show at
once how this implies that \Simplexf{} is a weak test category, or
equivalently, that the canonical functor\pspage{359}
\[i_\Simplexf : \Simplexf\to\Cat\]
is aspheric, for the canonical asphericity structure on \Cat. (I
should have noted that the asphericity statements are meant here in
the strongest possible sense, namely with respect to $\scrW=\scrWoo=$
usual weak equivalence.) Now, for any $\Simplex_n$ in \Simplexf, the
category
\[ i_\Simplexf(\Simplex_n) \eqdef \Simplexf_{/\Simplex_n}\]
is canonically isomorphic to the category associated to the ordered
set of all non-empty subsets of $\Simplex_n$, hence we get a canonical
isomorphism
\begin{equation}
  \label{eq:98.2}
  i_\Simplexf \simeq \widetilde i(\beta\alpha),\tag{2}
\end{equation}
where
\begin{equation}
  \label{eq:98.3}
  \widetilde i:\Simplextilde\to\Cat\tag{3}
\end{equation}
is the standard test functor, associating to any non-ordered simplex
$E$ the category associated to the ordered set of all non-empty
subsets of $E$. We know already (section \ref{sec:34}) that
$\widetilde i$ is aspheric, i.e., $\widetilde i^*$ transforms aspheric
objects of \Cat{} into aspheric ones, hence the same holds for its
composition with the aspheric functor $\beta\alpha$, hence also for
$i_\Simplexf$, qed.

Thus, we are left with proving that $\alpha$ is aspheric, i.e., that
the categories
\begin{equation}
  \label{eq:98.star}
  \Simplexf_{/\alpha^*(\Simplex_n)}\tag{*}
\end{equation}
are aspheric. Now, let's denote by \Simplexprimef{} the category
deduced from \Simplexf{} by adding an initial object $\varnothing$
(which we may view as being the empty simplex), which defines an
``open subcategory'' $U$, namely as sieve in \Simplexprimef, in such a
way that \Simplexf{} appears as the ``closed subcategory'', i.e.,
cosieve in \Simplexprimef{} complementary to $U$. One immediately
checks that the category \eqref{eq:98.star} is canonically isomorphic
to
\begin{equation}
  \label{eq:98.starstar}
  (\Simplexprimef)^{n+1} \setminus U^{n+1},\tag{**}
\end{equation}
where $U^{n+1}$ is the open subcategory defined by the initial object
of the ambient category $(\Simplexprimef)^{n+1}$. Now, asphericity of
\eqref{eq:98.starstar} and hence of \eqref{eq:98.star} follows from
the following two lemmas:
\begin{lemmanum}\label{lem:98.1}
  The category \Simplexf{} is aspheric.
\end{lemmanum}
\begin{lemmanum}\label{lem:98.2}
  Let $(X_i,U_i)_{i\in I}$ be a finite non-empty family of pairs
  $(X_i,U_i)$, where $X_i$ is a small category, $U_i$ an open
  subcategory. We assume that for any $i$ in $I$, $X_i$ and the closed
  complement $Y_i$ of $U_i$ in $X_i$ are aspheric. Let $X$ be the
  product of the $X_i$'s, $U$ the products of the $U_i$'s,\pspage{360}
  then \textup($X$ and\textup) $X\setminus U=Y$ are aspheric too.
\end{lemmanum}
\noindent\emph{Proof of lemma \ref{lem:98.2}:} by an immediate
induction, we are reduced to the case when $I$ has just two elements,
$I=\{1,2\}$, but then $X\setminus U$ can be viewed as the union of the
two closed subcategories $X_1\times Y_2$ and $X_2\times Y_1$, whose
intersection is $Y_1\times Y_2$. As all three categories are aspheric
(being products of aspheric categories), it follows by the well-known
Mayer-Vietoris argument that so is $X\times U$, qed.

Thus, we are left with proving that \Simplexf{} is aspheric. Somewhat
surprisingly, that's where I spent a number of hours not getting
anywhere and feeling foolish! There is a very simple heuristic
argument though involving the standard calculation of the cohomology
invariants of any semisimplicial ``complex'', in terms of the standard
boundary operations: if we admit that the same calculations are valid
when working with semisimplicial ``face complexes'', i.e., objects of
\Simplexfhat, then it is enough to apply this to the final object of
\Simplexfhat{} (including for computation of the non-commutative
$\mathrm H^1$ with constant coefficients) to get asphericity of
\Simplexf. As a matter of fact, this argument would give directly
asphericity of $\alpha$, bypassing altogether the categories
\eqref{eq:98.starstar} and lemmas \ref{lem:98.1} and \ref{lem:98.2}.
Apparently, I got a block against the down-to-earth computational
approach to cohomology via semisimplicial calculations, and have been
trying to bypass it at all price -- and not succeeding! Then,
curiously enough, when finding no other way out than look at those
boundary operations and try to understand what they meant (something I
remember vaguely have been doing once ages ago!), this brought me back
again to the abelianization story of sections \ref{sec:92} and
\ref{sec:93}, and to a more comprehensive way for looking at
``abelianizators'', and get an existence and unicity statement for
these. (At any rate, for a suitably strengthened version of these.)
This seems to me of independent interest, and worth being written down
with some care.

\hangsection{Afterthoughts on abelianization IV: Integrators.}\label{sec:99}%
When writing down (in sections \ref{sec:92} and \ref{sec:93}) some
rambling reflections about ``abelianization'' and ``abelianizators'',
there has been a persistent feeling of uneasiness, which I kept
pushing aside, as I didn't want to spend too much thought on this
``digression''. This uneasiness had surely something to do with the
way abelianization (of an object $X$ say of an elementary modelizer
\Ahat) was handled, so that it was designed in a more or less
exclusive way for embodying information about the ``homology
structure'' of the homotopy type modelized by $X$, or\pspage{361}
equivalently, to describe its cohomology invariants with arbitrary
\emph{constant} coefficients. Now, among the strongest reflexes I
acquired in the past while working with cohomology, was systematically
to look at coefficients which are arbitrary sheaves (abelian say), and
to view constant or locally constant coefficients as being just
particular cases. This reflex has been remaining idle, not to say
repressed, during nearly all of the reflections of the last four
months, due to the fact that in the whole modelizing story woven
around weak equivalence, there was a rather exclusive emphasis on
constant and locally constant coefficients, disregarding any other
coefficients throughout. Probably, while reflecting on abelianization,
a more or less underground reminiscence must have been around of the
semisimplicial boundary operators having a meaning for computing
cohomology of ``something'', with coefficients in arbitrary sheaves --
and also that to get it straight, one had to be careful not to get
mixed up in the variances. But I just didn't want to dive into all
this again if I could help it -- and now it is getting clear, after a
day or two of feeling silly, that it can't be helped, and I'll have to
write things down at last, however ``well-known'' they may be.

Let $A$ be a small category. In section \ref{sec:93} we defined an
``abelianizator'' for $A$ to be a chain complex $L_\bullet$ in the
additive envelope $\Add(A\op)$ of the category $A\op$ opposite to $A$,
satisfying a suitable condition of commutativity (in the diagram
\eqref{eq:92.20} of p.\ \ref{p:328}), and endowed with a mild extra
structure $\lambda$, expressing this commutativity. The function of an
abelianizator in loc.\ cit.\ was essentially to allow for a
simultaneous ``computational'' description of the homology structure
of the homotopy types stemming from a variable object $X$ in \Ahat, or
equivalently, to describe cohomology of such $X$ (as an object of a
suitable derived category say, to get it at strongest) with
coefficients in any (\emph{constant}) ring or abelian
group. Introducing by an independent symbol the opposite category
\begin{equation}
  \label{eq:99.1}
  B=A\op,\tag{1}
\end{equation}
I want now to establish a relationship between this property or
function of a chain complex $L_\bullet$ in $\Add(B)$, involving
objects in \Ahat{} and their abelianizations in \Ahatab, with an
apparently different one, in terms of a variable object of \Bhatab{}
(\emph{not} \Ahatab{} this time!), namely expressing \emph{cohomology
  of $B$} (i.e., of the topos \Bhat{} defined by $B$) with
coefficients in an \emph{arbitrary} abelian presheaf $F$ on $B$, i.e.,
an arbitrary object\pspage{362} of \Bhatab. I will first describe this
property of (possible) function of a chain complex in $\Add(B)$,
forgetting for the time being the category $A=B\op$ and the homotopy
types defined by objects $X$ in \Ahat. Once this property is well
understood, it will be time to show it implies the previous one
relative to $A$ and objects of \Ahat, and presumably is even
equivalent with it.

First, we'll have to interpret the category $\Add(B)$, which was
constructed somewhat ``abstractly'' in section \ref{sec:93} (as the
solution of a universal problem stated in section \ref{sec:92}), as a
full subcategory of the category \Bhatab{} of abelian presheaves on
$B$. It will be useful to keep in mind the following diagram of
canonical functors
\begin{equation}
  \label{eq:99.2}
  \begin{tabular}{@{}c@{}}
    \begin{tikzcd}[baseline=(O.base)]
      B\ar[r,"\alpha_B"] \ar[d]\ar[dr,"\beta_B"] & \Bhat
      \ar[d,"\Wh_B"] \\
      \Add(B) \ar[r,"\gamma_B"'] & |[alias=O]| \Bhatab
    \end{tikzcd},
  \end{tabular}\tag{2}
\end{equation}
where $\alpha_B$ is the canonical inclusion, $\Wh_B$ is the
abelianization functor, $\beta_B$ the composition of the two, and
$\gamma_B$ the \emph{additive} functor factoring $\beta_B$, in virtue
of the universal property of $\Add(B)$. This functor is defined up to
canonical isomorphism, and the lower triangle of \eqref{eq:99.2} is
commutative, up to a given commutativity isomorphism. Also, we'll use
the composition of the following sequence of canonical equivalences of
categories:
\begin{multline*}
  \Bhatab = \bHom(B\op,\Ab) \tosim \bHom(B,\Ab\op)\op \\
  \tosim \bHomadd(\Add(B),\Ab\op)\op
  \tosim \bHomadd(\Add(B)\op,\Ab),
\end{multline*}
i.e., a canonical equivalence of category
\begin{equation}
  \label{eq:99.3}
  \Bhatab\tosim\bHomadd(\Add(B)\op,\Ab),\quad
  F\mapsto\widetilde F,\tag{3}
\end{equation}
which is a particular case of
\begin{equation}
  \label{eq:99.3prime}
  \bHom(B\op,M)\tosim \bHom(\Add(B)\op,M),\tag{3'}
\end{equation}
where $M$ is any additive category. If $F$ is an abelian presheaf on
$B$, i.e., an object in the left-hand side of \eqref{eq:99.3}, we'll
denote by
\begin{equation}
  \label{eq:99.4}
  \widetilde F: \Add(B)\op \to \Ab\tag{4}
\end{equation}
the corresponding additive functor. Now, this functor can be
interpreted very nicely in terms of the functor $\gamma_B$ in
\eqref{eq:99.2}, by the canonical isomorphism of abelian groups
\begin{equation}
  \label{eq:99.5}
  \widetilde F(L)\tosim \Hom_\Bhat(\gamma_B(L),F),\tag{5}
\end{equation}
functorial\pspage{363} with respect to the pair $(F,L)$ in
$\Bhatab\times\Add(B)\op$. This formula in turn implies easily that
the functor $\gamma_B$ is \emph{fully faithful}. Thus, we can
interpret $\Add(B)$ as the full subcategory of \Bhatab{} whose objects
are all finite direct sums (in \Bhatab) of objects of the type
$\Wh_B(b)$, with $b$ in $B$. In terms of this interpretation,
$\gamma_B$ is just an inclusion functor, and on the other hand, for
$F$ in the ambient category \Bhatab, $\widetilde F$ is just the
restriction to the subcategory $\Add(B)$ of the contravariant functor
on \Bhatab{} represented by $F$.

This situation is the exact ``additive'' analogon of the situation of
$B$ embedded in \Bhat{} as a full subcategory, the functor on $B$
defined by an object $F$ of \Bhat{} being the restriction to $B$ of
the contravariant functor on \Bhat{} represented by $F$, i.e., an
object of $F(b)$ or $\widetilde F(b)$ can (often advantageously) be
interpreted as a map in \Bhat, $b\mapsto F$. Moreover, the fact that
\[\alpha_B^*: F\mapsto \widetilde F: \Bhat \tosim \bHom(b\op,\Sets)\]
is an equivalence (in fact, an isomorphism even), is paralleled by the
equivalence \eqref{eq:99.3}, which can likewise be interpreted as
$\gamma_B^*$, or more accurately as the canonical factorization of the
purely set-theoretic $\gamma_B^*:\Bhatab\to\bHom(\Add(B)\op,\Sets)$
through $\bHomadd(\Add(B)\op,\Ab)$\ldots

The objects $L$ of the full subcategory $\Add(B)$ of \Bhatab{} have a
very strong common property, namely they are \emph{projectives}, and
they are of \emph{finite presentation} (``small'' in Quillen's
terminology), namely for variable $F$ in \Bhatab, the functor
\[F\mapsto \Hom_\Bhatab(L,F)\]
commutes with filtering direct limits. Both properties are immediate,
and they nearly characterize the objects in $\Add(B)$ -- more
accurately, it is immediately checked that the projectives of finite
presentation in \Bhatab{} are exactly those which are isomorphic to
\emph{direct factors} of objects in $\Add(B)$. It shouldn't be hard to
check that the full subcategory of \Bhatab{} made up with the
projectives of finite presentation can be identified up to equivalence
to the ``Karoubi envelope'' of the category $\Add(B)$ (obtained by
formally adding images of projectors), or equivalently, can be
described as the solution of the $2$-universal problem defined by
sending $B$ into categories which are, not only additive, but moreover
stable under taking images of projectors (i.e., endomorphisms $u$ of
objects, such that $u^2=u$).

We'll\pspage{364} henceforth identify $\Add(B)$ to a full subcategory
of \Bhatab{} (by replacing the solution of the universal problem,
constructed in section \ref{sec:95}, by the essential image in
\Bhatab{} say), and rewrite \eqref{eq:99.5} simply as
\begin{equation}
  \label{eq:99.5prime}
  \widetilde F(L)\simeq \Hom(L,F),\tag{5'}
\end{equation}
the $\Hom$ being taken in \Bhatab, category of abelian presheaves on
$B$. Accordingly, if $L_\bullet$ is a chain complex in $\Add(B)$,
hence in \Bhatab, the corresponding \emph{cochain} complex $\widetilde
F(L_\bullet)$ in \Ab{} can be interpreted as
\begin{equation}
  \label{eq:99.6}
  \widetilde F(L_\bullet) \simeq \Hom^\bullet(L_\bullet,F),\tag{6}
\end{equation}
where the symbol $\Hom^\bullet$ means taking $\Hom$'s componentwise.

What we're after here is to find a \emph{fixed} chain complex
$L_\bullet$ in $\Add(B)$, such that for \emph{any} abelian presheaf
$F$ on $B$, the cochain complex \eqref{eq:99.6} in \Ab{} should be
isomorphic (in the derived category $\D^\bullet\Ab$ of cochain
complexes in \Ab{} with respect to quasi-isomorphisms) to the
``integration'' of $F$ over the topos \Bhat, i.e., to $\mathrm
R\Gamma_B(F)$:
\begin{equation}
  \label{eq:99.star}
  \Hom^\bullet(L_\bullet,F)\simeq \mathrm R\Gamma_B(F)\text{
    ?}\quad(\text{isom.\ in $\D^\bullet\Ab$}),\tag{*}
\end{equation}
namely to the total right derived functor $\mathrm R\Gamma_B$ (taken
for the argument $F$) of the ``sections'' functor
\begin{equation}
  \label{eq:99.7}
  \Gamma_B(F) \eqdef \varprojlim_{B\op} F.\tag{7}
\end{equation}
Now, using the fact that the components of the chain complex
$L_\bullet$ are projective, hence $\Ext^i(L_n,F)=0$ for $i>0$ (any
$n$, any $F$), we get at any rate a canonical isomorphism in
$\D^\bullet\Ab$, or in $\D\Ab$:
\begin{equation}
  \label{eq:99.8}
  \Hom^\bullet(L_\bullet,F)\simeq \mathrm R\Hom(L_\bullet,F),\tag{8}
\end{equation}
i.e., an interpretation of \eqref{eq:99.6} as a ``hyperext''. Now,
let's remember that $\mathrm R\Gamma_B(F)$ (as on any topos) can be
interpreted equally as
\begin{equation}
  \label{eq:99.9}
  \mathrm R\Gamma_B(F)\simeq\mathrm R\Hom(\bZ_B,F),\tag{9}
\end{equation}
where $\bZ_B$ denotes the constant presheaf on $B$ with value
$\bZ$. Thus, the wished-for isomorphism \eqref{eq:99.star} will follow
most readily from a corresponding isomorphism in $\D_\bullet(\Bhatab)$
between $L_\bullet$ and $\bZ_B$. But using again the fact that the
components of $L_\bullet$ are projective, we see that to give a map in
the derived category of $L_\bullet$ into $\bZ_B$ amounts to the same
as to give an \emph{augmentation}
\begin{equation}
  \label{eq:99.10}
  L_\bullet\to\bZ_B,\tag{10}
\end{equation}
and the map is an isomorphism in $\D_\bullet(\Bhatab)$
if{f}\pspage{365} the augmentation \eqref{eq:99.10} turns $L_\bullet$
into a (projective) \emph{resolution} of $\bZ_B$.

We now begin to feel in known territory again! Let's call
``\emph{integrator}'' on $B$ any projective resolution of $\bZ_B$, and
let's call the integrator ``\emph{special}'' (by lack of a more
suggestive name) if its components $L_n$ are in $\Add(B)$, or what
amounts to the same, if it can be viewed as a chain complex in
$\Add(B)$, endowed with the extra structure \eqref{eq:99.10}. Of
course, $\bZ_B$ is no longer in $\Add(B)$ in general, and therefore
the data \eqref{eq:99.10} has to be interpreted as a map $L_0\to\bZ_B$
external to $\Add(B)$, or equivalently (via \eqref{eq:99.5prime}) as
an object
\begin{equation}
  \label{eq:99.11}
  \lambda\quad\text{in}\quad \widetilde{\bZ}_B(L_0) = \bZ^{(I_0)},\tag{11}
\end{equation}
where $I_0$ is the set of indices used in order to express $L_0$ as
the direct sum in \Bhatab{} of elements of $B$. We know, by the
general principles of homological algebra, that any two integrators
must be chain homotopic, hence, if they are special, as $\Add(B)$ is a
full additive subcategory of \Bhatab, they must be chain homotopic in
$\Add(B)$.

As for existence of integrators, it follows equally from general
principles, as we know that \Bhatab{} has ``enough projectives''
(which is a very special feature indeed of \Bhatab, coming from the
fact that the topos \Bhat{} has enough projectives, namely the objects
of $B$\ldots). It isn't clear though that there exists a special
integrator, because when trying inductively to construct the
resolution $L_\bullet$ of $\bZ_B$ with components in $\Add(B)$, it
isn't clear that the kernel of $L_n\to L_{n-1}$ is ``of finite type'',
namely is isomorphic to a quotient of an object of $\Add(B)$ (or,
equivalently, is a quotient of a projective of finite
presentation). If we take for instance $B$ to be the one-object
groupoid defined by a group $G$, an integrator on $B$ is just a
resolution of the constant $G$-module $\bZ$ by projective
$\bZ[G]$-modules, and the integrator is special of{f} the components
are even free modules of finite type -- I doubt such a resolution
exists unless $G$ itself is finite. This example seems to indicate
that the existence of a special integrator for $B$ is a very strong
condition on $B$, of the nature of a (homological) finiteness
condition. Maybe this condition, more than most others, singles out
the three standard test categories and their finite products, from
arbitrary test categories (even strict ones and contractors\ldots).

Even in case a strict integrator doesn't exist for $B$, there is
a rather evident way out to get ``the next best'' in terms of
computations, namely replacing the very much finitely restricted
category $\Add(B)$ by a larger category\pspage{366}
\begin{equation}
  \label{eq:99.12}
  \Addinf(B)\hookrightarrow \Bhatab\tag{12}
\end{equation}
deduced from $B$ by adding, not merely finite direct sums (and linear
combinations of maps), but equally infinite ones. The construction can
be given ``formally'' as in section \ref{sec:93}, and it can be
checked that this category satisfies the obvious $2$-universal
property with respect to functors
\[f:B\to M\]
from $B$ to infinitely additive categories $M$ (namely additive
categories stable under direct sums), and functors $M\to M'$ which are
not merely additive, but commute to small direct sums. Moreover, it is
checked that the category $\Addinf(B)$ thus constructed embeds by a
fully faithful functor into \Bhatab{} as indicated in
\eqref{eq:99.12}, and hence can be identified up to equivalence to a
full subcategory of \Bhatab. The formulas \eqref{eq:99.5} and
\eqref{eq:99.5prime} are still valid, when $L$ is in $\Addinf(B)$ only
instead of $\Add(B)$. The objects of $\Addinf(B)$ in \Bhatab{} are
still projective (as direct sums of projectives), but of course no
longer of finite presentation. In compensation, any element in
\Bhatab{} is now a quotient of an object in $\Addinf(B)$. As a
consequence, the projectives in \Bhatab{} can be characterized as the
direct factor of objects of $\Addinf(B)$, and presumably the full
subcategory of \Bhatab{} made up with all projectives can again be
described (up to equivalence) as the Karoubi envelope of $\Addinf(B)$,
or equivalently, as the solution of the $2$-universal problem of
sending $B$ into infinitely additive karoubian categories (karoubian =
every projective has an image, i.e., corresponds to a direct sum
decomposition). We may call an integrator $L_\bullet$ for $B$ with
components in $\Addinf(B)$ ``\emph{quasi-special}''. We did just what
was needed in order to be sure now that there always exist
quasi-special integrators; moreover, these integrators are unique up
to chain homotopy in $\Addinf(B)$. The interpretation \eqref{eq:99.11}
of the augmentation structure \eqref{eq:99.10} on $L_\bullet$ is still
valid in the quasi-special case, with the only difference that now the
indexing set $I_0$ need not be finite anymore.

\bigbreak

\presectionfill\ondate{17.6.}\pspage{367}\par

\hangsection{Abelianization V: Homology versus cohomology.}\label{sec:100}%
Yesterday I introduced the notion of an \emph{integrator} for any
small category $B$, to be just a projective resolution of $\bZ_B$ in
the category \Bhatab{} of all abelian presheaves on $B$, where $\bZ_B$
denotes the constant presheaf with value $\bZ$. Such an object in
$\Ch_\bullet(\Bhatab)$ exists, due to the existence of sufficiently
many projectives in \Bhatab, and it is unique up to homotopism of
augmented chain complexes, which encourages us to denote it by a
canonizing symbol, namely
\begin{equation}
  \label{eq:100.1}
  L_\bullet^B\to \bZ_B.\tag{1}
\end{equation}
As will become clear in the sequel, $L_\bullet^B$ can be viewed as
embodying \emph{homology properties} of $B$, i.e., of the topos
associated to $B$ (whose category of sheaves of sets is \Bhat). The
way we hit upon it though was in order to obtain a ``computational''
way for computing \emph{cohomology} of $B$ (i.e., of the associated
topos) with coefficients in any abelian presheaf $F$ in \Bhatab, by a
canonical isomorphism
\begin{equation}
  \label{eq:100.2}
  \mathrm R \Gamma_B(F) \tosim \Hom^\bullet(L_\bullet^B,F)\tag{2}
\end{equation}
in the derived category $\D^\bullet\Ab$, where $\Hom^\bullet$ denotes
the cochain complex obtained by applying $\Hom$ componentwise. Passing
to the cohomology groups of both members, this gives rise to
\begin{equation}
  \label{eq:100.3}
  \mathrm H^i(B,F) \simeq \mathrm H^i \Hom^\bullet(L_\bullet^B,F).\tag{3}
\end{equation}
The designation ``computational'' takes a rather concrete meaning,
when we choose $L_\bullet^B$ to have its components in the infinitely
additive envelope $\Addinf(B)$ of $B$, which (as we saw yesterday) can
be viewed as a full subcategory of \Bhatab, made up with projectives,
and such that any object in \Bhat{} is quotient of an object coming
from $\Addinf(B)$; this ensures that there exist indeed integrators
which are ``\emph{quasi-special}'', i.e., are made up with objects of
$\Addinf(B)$, and hence can be interpreted as chain complexes of this
additive category. Thus, any component $L_n$ can now be written, in an
essentially canonical way, as
\begin{equation}
  \label{eq:100.4}
  L_n = \bigoplus_{i\in I_n} \bZ^{(b_i)} ,\tag{4}
\end{equation}
where
\begin{equation}
  \label{eq:100.5}
  (b_i)_{i\in I_n} \tag{5}
\end{equation}
is a family of objects of $B$ indexed by $I_n$ (NB\enspace for
simplicity of notations, we assume the $I_n$'s mutually disjoint,
otherwise we should write the\pspage{368} general object in the family
\eqref{eq:100.5} $b_i^n$ rather than $b_i$). Thus, the $n$'th
component of the cochain complex of the second member of
\eqref{eq:100.2} can be explicitly written as
\begin{equation}
  \label{eq:100.6}
  \Hom^n(L_\bullet^B,F) = \Hom(L_n^B,F) \simeq \bigoplus_{i\in I_n} F(b_i),\tag{6}
\end{equation}
and the coboundary operators between these components can be made
explicit in a similar way, by means of (possibly infinite) matrices,
whose entries are $\bZ$-linear combinations of maps from some $b_i^n$
to some $b_j^{n-1}$ ($i\in I_n$, $j\in I_{n-1}$). We feel a little
happier still when the direct sums \eqref{eq:100.4} yielding the
components $L_n$ are finite, i.e., the sets $I_n$ are finite, which
also means that $L_\bullet^B$ can be interpreted as a chain complex in
the additive envelope $\Add(B)$ of $B$, as contemplated in the first
place -- in which case the integrator will be called
``\emph{special}''.

The formula \eqref{eq:100.2} immediately generalizes when $F$ is
replaced by a complex of presheaves $F^\bullet$, with degrees bounded
from below (NB\enspace as the notation indicates, the differential
operator is of degree $+1$), to
\begin{equation}
  \label{eq:100.7}
  \mathrm R\Gamma_B(F^\bullet) \tosim \Hom^{\bullet\bullet}(L_\bullet^B,F^\bullet),\tag{7}
\end{equation}
where now the left-hand side designates hypercohomology of $B$ (i.e.,
of the corresponding topos), viewed as an objects of the right derived
category $\D^+\Ab$ of the category of abelian groups, and where
$\Hom^{\bullet\bullet}$ designates the double complex obtained by
taking $\Hom$'s componentwise, or more accurately, the object in
$\D^+\Ab$ defined by the associated simple complex.

An interesting special case of \eqref{eq:100.7} is obtained when
starting with a complex of abelian groups $K^\bullet$ bounded from
below, i.e., defining an object of the right derived category
$\D^+\Ab$, and taking
\begin{equation}
  \label{eq:100.8}
  F^\bullet=K_B^\bullet=p_B^*(K^\bullet),\tag{8}
\end{equation}
the corresponding \emph{constant complex of presheaves} on $B$, which
may be viewed equally as the inverse image of $K^\bullet$ by the
projection
\begin{equation}
  \label{eq:100.9}
  p_B:B\to\Simplex_0\quad(\text{the final category}),\tag{9}
\end{equation}
which geometrically interprets as the canonical morphism of the topos
associated to $B$ to the final (or ``one-point'') topos. The second
member of \eqref{eq:100.7} can be rewritten componentwise, using the
adjunction formula for the pair $(p_!\supab,p^*)$ (where the
qualifying $B$ is omitted now in the notation $p$):
\[\Hom(L_n,p^*(K^m))\simeq\Hom(p_!\supab(L_n),K^M),\]
so that \eqref{eq:100.7} can be rewritten as\pspage{369}
\begin{equation}
  \label{eq:100.10}
  \mathrm R \Gamma_B(K_B^\bullet) \simeq
  \Hom^{\bullet\bullet}(p_{B!}\supab(L_\bullet^B),K^\bullet),\tag{10}
\end{equation}
where this time the $\Hom$'s in the right-hand side of
\eqref{eq:100.10} are taken in \Ab, not in \Bhatab.

This formula very strongly suggests to view the chain complex of
abelian groups
\begin{equation}
  \label{eq:100.11}
  p_{B!}\supab(L_\bullet^B),\tag{11}
\end{equation}
which is in fact a complex of projective (hence free) abelian groups
defined up to chain homotopy, as embodying the global homology
structure of $B$ (or of the corresponding topos), more accurately
still, as embodying the homology structure of the corresponding
homotopy type. It is easily seen that the corresponding object of
$\D_\bullet\Ab$ depends covariantly on $B$ when $B$ varies in the
category \Cat, so that we get a functor
\[\Cat\to \D_\bullet\Ab \eqdef \Hotab,\]
which in view of \eqref{eq:100.10} (an isomorphism functorial not only
with respect to $K^\bullet$, but equally with respect to $B$) factors
through the localization \Hot{} of \Cat, thus yielding a canonical
functor
\begin{equation}
  \label{eq:100.12}
  \Hot\to\Hotab,\tag{12}
\end{equation}
which deserves to be called the \emph{abelianization functor}, from
homotopy types to ``abelian homotopy types''. This cannot be of course
anything else (up to canonical isomorphism) but the functor
\eqref{eq:92.1} of section \ref{sec:92} (p.\ \ref{p:321}), but
obtained here in a wholly ``intrinsic'' way, without having to pass
through the particular properties of a particular test category such
as $\Simplex$ or one of its twins. One possible way to check this
identity would be by proving that an isomorphism \eqref{eq:100.10} is
valid when replacing (for a given $B$ in \Cat{} and $K^\bullet$ in
$\D^+\Ab$) the chain complex \eqref{eq:100.11} by the corresponding
one deduced from the map \eqref{eq:92.1} defined p.\ \ref{p:321} (via
the diagram \eqref{eq:92.3} on p.\ \ref{p:322}), and checking moreover
that an object $\ell_\bullet$ of $\D^-\Ab$ is known up to canonical
isomorphism, when we know the corresponding functor
\begin{equation}
  \label{eq:100.13}
  K^\bullet\mapsto \Hom_{\D\Ab}(\ell_\bullet,K^\bullet)\tag{13}
\end{equation}
on $\D^+\Ab$. Presumably, this latter statement holds when replacing
\Ab{} by any abelian category, but I confess I didn't sit down to
check it, nor do I remember having seen it stated somewhere -- as I
don't remember either having seen anywhere a comprehensive treatment
about the relationship between homology and cohomology. So maybe my
present reflections do fill a gap, or at any rate give some
indications as to how to fill it\ldots

I played around some yesterday and today with the formalism of
integrators, notably with respect to maps
\[ f:B'\to B\]
between small categories, and the corresponding \emph{integration
  functor}
\[ f_!\supab: {B'}\subab\uphat \to \Bhatab,\]
and its left derived functor $\mathrm Lf_!\supab$. Thus, the chain
complex in \Bhatab
\begin{equation}
  \label{eq:100.14}
  L_\bullet^{B'/B}\quad\text{or}\quad
  L_\bullet^f \eqdef f_!\supab(L_\bullet^{B'}),\tag{14}
\end{equation}
which has projective components (and even is a chain complex in
$\Addinf(B)$ resp.\ in $\Add(B)$, if $L_\bullet^{B'}$ is quasi-special
resp.\ is special), and is defined up to chain homotopism, embodies
the relative homology properties of $B'$ over $B$, i.e., of $f$, in
much the same way as \eqref{eq:100.11} embodies the global homology
properties of $B$ (i.e., of $B$ over one point). When the functor $f$
is ``coaspheric'', i.e., the functor
\[f\op:{B'}\op \to B\op\]
between the opposite categories is aspheric, then $L_\bullet^{B'/B}$
is again an integrator on $B$, and the converse should hold too
provided we take the meaning of ``coaspheric'' and ``aspheric'' with
respect to a suitable basic localizer $\scrW=\scrW_\oo^{\bZ}$ --
presumably, we'll come back upon this in part \ref{ch:V} or part
\ref{ch:VI} of the notes. For the time being, it seems more
interesting to give now the precise relationship between the notion of
an \emph{integrator} for $B$, and the notion of an
\emph{abelianizator} for the dual category $A=B\op$, introduced in
section \ref{sec:93}.

\medbreak

\noindent\textbf{Remarks.} \namedlabel{rem:100.1}{1})\enspace It is a familiar fact that when
working in \v Cech-flavored contexts, such as general topoi, or étale
topoi for schemes and the like, one has throughout and from the start
a good hold upon \emph{cohomology} notions, whereas it is a lot more
subtle to squeeze out adequate homology notions, which (to my
knowledge) can be carried through only indirectly via cohomology, and
using suitable finiteness and duality statements within the cohomology
formalism. Historically however, homology was introduced before
cohomology via cellular decompositions of spaces, with a more direct
appeal to geometric intuition. This preference for homology rather
than cohomology seems to be still prevalent among most homotopy
theorists, who have a tendency to view a topological space (however
wild it may be) as being no more no less than its singular complex. A
comprehensive statement establishing, in a suitable wide enough
context, essential equivalence between the two viewpoints, seems to be
still lacking, as far as I\pspage{371} know -- although a fair number
of partly overlapping results in this direction are known, among the
oldest being the relevant ``universal coefficients formulæ'' relating
homology and cohomology (reducing all to a formula of the type
\eqref{eq:100.2} or \eqref{eq:100.10} above), or Cartan's old seminar
on Leray's sheaf theory, introducing singular homology with
coefficients in a sheaf and proving that on a topological variety,
this was (up to dimension shift and twist by the twisted integers)
essentially the same as singular cohomology (with coefficients in
sheaves too). It is not sure that an all-inclusive statement of
equivalence between homology and cohomology (in those situations when
such equivalence is felt hold indeed) does at all exist -- at any
rate, according to what kind of coefficients one wants to consider,
and what kind of extra structures one is interested in when dealing
with homology and cohomology invariants, it seems that each of the two
points of view has an originality and advantages of its own and cannot
be entirely superseded by the other. From the contexts I have been
mainly working in, there definitely was no choice, namely cohomology
(including non-commutative one) was the basic data, while sheaves and
their generalizations (such complexes of sheaves, or stacks) were the
coefficients. I don't remember of any moment where I would have paused
and asked myself \emph{why} in most contexts where I was working in
(whose common denominator was topoi), there wasn't any direct hold on
anything like homology invariants. The reason for this inertness of
mine, probably, is that the cohomology formalisms I hit upon were
self-contained enough, so as to leave no regret for the absence of a
homology formalism, or at any rate of a more or less direct
description of it independently of cohomology. Another reason, surely,
is that I didn't have too much contact with topologists and
homotopists and their everyday tools, such as Steenrod operations,
homology of the symmetric group, and the like. This question of ``why
this reluctance of homology to show forth'' has finally surfaced only
during these very last days, when the answer for it (or one possible
answer at any rate) is becoming evident: namely, that \emph{for a
  general topos}, embodied by a category of sheaves (of sets) \scrA,
\emph{there are not enough projectives in \scrA, and not even enough
  projectives in $\scrAab$}, the category of abelian sheaves. It
is becoming apparent (what surely everybody has known ages) that in
technical terms, \emph{doing ``homology'' is working with projectives,
  while doing ``cohomology'' is working with injectives}. As there are
enough injectives in $\scrAab$ but not enough projectives,
cohomology is around and homology not, period!

There is however a rather interesting class of topoi admitting
sufficiently\pspage{372} many projective sheaves of sets, and hence
sufficiently many projective abelian sheaves -- namely the topoi
\Bhat{} defined in terms of small categories $B$. They include the
topoi which can be described in terms of semisimplicial complexes and
the like, and can be viewed equally as the topoi which are ``closest
to algebra'' or ``purely algebraic'' in a suitable sense -- for
instance, definable directly in terms of arbitrary presheaves, without
any reference to the notion of site and of localization. (The
intuition of localization remaining however and indispensable guide
even in the so-called ``algebraic'' set-up.) Moreover, the morphisms
which arise most naturally among such topoi, namely those associated
to maps
\[f:B'\to B\]
in \Cat, besides the traditional adjoint pair $(f^*,f_*)$ of functors
between sheaves of sets, gives rise equally to a functor
\begin{equation}
  \label{eq:100.15}
  f_!:{B'}\uphat\to\Bhat\tag{15}
\end{equation}
left adjoint to $f^*$ (i.e., $f^*$ commutes to small inverse limits,
not only to small direct limits and to finite inverse limits),
inserting in a triple of mutually adjoint functors (from left to
right)
\begin{equation}
  \label{eq:100.16}
  (f_!,f^*,f_*).\tag{16}
\end{equation}
The functors $f^*$ and $f_*$ induce corresponding adjoint functors on
abelian sheaves (due to the fact that they commute to finite
products), $f^*\subab$ and $f_*\supab$, whereas $f_!$ does not in
general transform group objects into group objects; however, as
$f^*\subab$ commutes to small inverse limits, it does admit again a
left adjoint $f_!\supab$, so as to give again a triple
\begin{equation}
  \label{eq:100.17}
  (f_!\supab,f^*\subab,f_*\supab)\tag{17}
\end{equation}
of mutually adjoint functors. Now, whereas the derived functors
\[\text{$f^*$ or $\mathrm L f^*\subab$,} \quad
\text{$\mathrm R f_*$ or $\mathrm R f_*\supab$}\]
of $f^*\subab$ and $f_*\supab$ have been extensively used in the
every-day cohomology formalism of topoi, the existence in certain
cases (such as the one we are interested in here) of a functor
$f_!\supab$ and of its left derived functor
\begin{equation}
  \label{eq:100.18}
  \text{$\mathrm L f_!$ or $\mathrm Lf_!\supab: \D^-({B'}\uphat\subab)
    \to \D^-(\Bhatab),$}\tag{18}
\end{equation}
seems to me to have been widely overlooked so far, except in extremely
particular cases such as inclusion of an open subtopos; at any rate, I
have been overlooking it till lately, when it came to my attention
through\pspage{373} the writing of these notes. (Namely, first in
connection with my reflections on derivators (cf.\ section
\ref{sec:69}), and now in connection with the reflections on
abelianization.) In view of my reflections on derivators, I would like
to view the functor \eqref{eq:100.18} as an operation of
``integration'', whereas the traditional functor
\begin{equation}
  \label{eq:100.19}
  \mathrm Rf_* : \D^+({B'}\uphat\subab) \to \D^+(\Bhatab)\tag{19}
\end{equation}
is viewed as ``cointegration'' (which I prefer to my former way of
calling it an ``integration''). The first should be viewed as
expressing \emph{homology} properties of the map $f$ in \Cat{} (or
between the corresponding topoi), just as the latter expresses
\emph{cohomology} properties of $f$. This does check with the
corresponding qualifications ``integration'' -- ``cointegration'' --
as well as with the intuition, when $B$ is reduced to a point,
identifying the first to a kind of direct sum (= integration), whereas
the latter is viewed as a kind of direct product (=
\emph{co}integration). The idea behind the terminology will go through
maybe when looking at the particular case when $B'$ is a sum of copies
of $B$, namely a product of $B$ by a discrete category $I$, and
\[f : B'=B\times I\to B\]
the projection.

The point I want to make here, mainly to myself, is that in the
present context \emph{when \eqref{eq:100.18}, namely integration,
  exists, this operation presumably is by no means less meaningful and
  important than the familiar $\mathrm Rf_*$ or cointegration} -- or
\emph{equivalently stated, that the homology properties of $f$ are
  just as meaningful and deserving close attention, as the cohomology
  properties}, which so far have been the only ones I have been
looking at. Presumably, when following this recommendation, a few
unexpected facts and relationships should come out, such as various
``duality'' relationships between homology properties of $f$, and
cohomology properties of the map $f\op$ between the opposite
categories. (This is suggested by some of the scratchwork I made on
derivators and cohomology properties of maps in \Cat.) The only
trouble is that such change or broadening of emphasis as I am now
suggesting will require a certain amount of extra attention, which I
am not too sure to be willing to invest in the subject, namely
algebraic topology. Thus presumably, my main emphasis will remain with
cohomology, rather than homology. I am no longer convinced though that
this point of view is technically more adequate than the dual one.

\namedlabel{rem:100.2}{2})\enspace All the reflections of yesterday's
notes as well as today's can be extended, when replacing throughout
abelian presheaves by presheaves of\pspage{374} $k$-modules, and
additive envelopes by $k$-linear ones, where $k$ is any given
commutative ring. Of course, the category \Ab{} and its various
derived categories will have to be replaced accordingly by the
category \kMod{} of $k$-modules etc. The same holds for the
relationship I am going to write down between integrators for $B$ and
abelianizators for $A=B\op$. For simplicity of notations, I am going
to keep the exposition in the \Ab-framework I have started with, and
leave the necessary adjustments to the reader.

\bigbreak
\presectionfill\ondate{18.7.}\par

\hangsection[Abelianization VI: The abelian integration operation
\dots]{Abelianization VI: The abelian integration operation
  \texorpdfstring{$\mathrm Lf_!\supab$}{Lf!ab} defined by a map
  \texorpdfstring{$f$ in \Cat{} \textup(}{f in (Cat) (}versus abelian
  cointegration
  \texorpdfstring{$\mathrm Rf_*$\textup)}{Rf*)}.}\label{sec:101}%
Finally with yesterday's non-technical reflections on homology versus
cohomology, it was getting prohibitively late, and there could be no
question to deal with the relationship between integrators (for $B$)
and abelianizators (for $B\op=A$). Also, I feel I should give some
``computational'' details about the functor $f_!\supab$ associated to
a map in \Cat
\[f:B'\to B,\]
namely
\[f_!\supab:{B'}\uphat\subab \to \Bhatab,\]
which is a lot less familiar to me than its right adjoint and
biadjoint $f^*$ and $f_*$. One way to get a ``computational hold''
upon it is by noting that $f_!\supab$ commuting to small direct limits
and a fortiori being right exact, and moreover any object $F'$ in
${B'}\uphat\subab$ being a cokernel of a map between ``\emph{special
  projectives}'' in ${B'}\uphat\subab$, i.e., between objects in
$\Addinf(B')$, namely inserting into an exact sequence
\[L_1' \xrightarrow d L_0' \to F' \to 0 \quad\text{with $L_0',L_1'$ in
  $\Addinf(B')$,}\]
the functor $f_!\supab$ (via its values on any $F'$ say) is
essentially known, when we know its restriction to the subcategory
$\Addinf(B')$, as we'll get a corresponding exact sequence in \Bhatab
\[f_!\supab(L_1') \to f_!\supab(L_0') \to f_!\supab(F') \to 0,\]
describing $f_!\supab(F')$ as a cokernel of a map $f_!\supab(d)$
corresponding to a map in $\Addinf(B')$. The relevant fact now is that
we have a commutative diagram of functors (up to can.\ isomorphism as
usual)
\begin{equation}
  \label{eq:101.1}
  \begin{tabular}{@{}c@{}}
    \begin{tikzcd}[baseline=(O.base)]
      \Addinf(B')\ar[r]\ar[d,"\Addinf(f)"'] &
      {B'}\uphat\subab\ar[d,"f_!\supab"] \\
      \Addinf(B)\ar[r] & |[alias=O]| \Bhatab
    \end{tikzcd},
  \end{tabular}\tag{1}
\end{equation}
where\pspage{375} the horizontal arrows are the canonical inclusion
functors, and $\Addinf(f)$ is the ``tautological'' extension of $f:B'
\to B$ to the infinitely additive envelopes, defined computationally
as
\begin{equation}
  \label{eq:101.2}
  \Addinf(f)(L') \simeq \bigoplus_{i\in I}\bZ^{f(b_i')}\tag{2}
\end{equation}
for an object of $\Addinf(B')$ written canonically as
\begin{equation}
  \label{eq:101.3}
  L'=\bigoplus_{i\in I}\bZ^{(b_i')}.\tag{3}
\end{equation}
Here, for an object $b$ in a small category $B$, we denote by the more
suggestive symbol $\bZ^{(F)}$ the abelianization $\Wh_B(F)$ of an
object $F$ of \Bhat, and accordingly of $F$ is an object $b$ in
$B$. The fact that \eqref{eq:101.2} is equally an expression for
$f_!\supab$ follows immediately from commutation of $f_!\supab$ to
small direct sums, and from the canonical isomorphism
\begin{equation}
  \label{eq:101.4}
  f_!\supab(\bZ^{F'}) \simeq \bZ^{(f_!(F'))},\tag{4}
\end{equation}
i.e., commutation up to canonical isomorphism of the diagram
\begin{equation}
  \label{eq:101.5}
  \begin{tabular}{@{}c@{}}
    \begin{tikzcd}[baseline=(O.base)]
      {B'}\uphat \ar[r,"\Wh_{B'}"] \ar[d,"f_!"'] &
      {B'}\uphat\subab \ar[d,"f_!\supab"] \\
      \Bhat\ar[r,"Wh_B"] & |[alias=O]| \Bhatab
    \end{tikzcd},
  \end{tabular}\tag{5}
\end{equation}
the verification of which is immediate. (For a generalization to
sheaves endowed with arbitrary ``algebraic structures'' and taking
free objects, see\scrcomment{\textcite{SGA4vol1}} SGA~4 I~5.8.3,
p.~30.)

Of course, \eqref{eq:101.1} and \eqref{eq:101.2} imply that
$f_!\supab$ maps $\Add(B')$ into $\Add(B)$, and induces the
tautological extension $\Add(f)$ of $f$ to the additive
envelopes. Thus, \eqref{eq:101.1} and \eqref{eq:101.5} can be inserted
into a beautiful commutative diagram (up to canonical isomorphism)
\begin{equation}
  \label{eq:101.6}
  \begin{tabular}{@{}c@{}}
    \begin{tikzcd}[baseline=(O.base),column sep=small]
      B' \ar[rr,hook] \ar[d,"f"'] \ar[ddrrr,hook,%
      dash pattern=on 33 pt off 30pt on 60pt] & &
      \Add(B')\ar[rr,hook] \ar[d,"\Add(f)"] & &
      \Addinf(B')\ar[rr,hook] \ar[d,"\Addinf(f)"'] & &
      {B'}\uphat\subab \ar[d,"f_!\supab"] \\
      B \ar[rr,hook] \ar[ddrrr,hook] & &
      \Add(B) \ar[rr,hook] & &
      \Addinf(B)\ar[rr,hook] & &
      \Bhatab \\
      & & & {B'}\uphat \ar[d,"f_!"] \ar[uurrr,"\Wh_{B'}" pos=0.15,%
      dash pattern=on 40pt off 30pt on 60pt] & & & \\
      & & & |[alias=O]| \Bhat \ar[uurrr,"\Wh_B" pos=0.4] & & &
    \end{tikzcd}.
  \end{tabular}\tag{6}
\end{equation}

The formula \eqref{eq:101.1} (or equivalently, \eqref{eq:101.2}) can
be viewed as giving a computational description of the left derived
functor
\begin{equation}
  \label{eq:101.7}
  \mathrm Lf_!\supab:\D^-({B'}\uphat\subab)\to \D^-(\Bhatab).\tag{7}
\end{equation}
Indeed,\pspage{376} by general principles of homological algebra, for
any small category $B$, from the fact that $\Addinf(B)$ is made up
with projective objects of \Bhatab{} and that any object in \Bhatab{}
is isomorphic to a quotient of an object in this subcategory, it
follows that
\begin{equation}
  \label{eq:101.8}
  \D^-(\Bhatab)\equeq W_B^{-1}\Comp^-(\Addinf(B)),\tag{8}
\end{equation}
i.e., the left derived category $\D^-(\Bhatab)$ is equivalent with the
category obtained by localizing, with respect to the set $W_B$ of
homotopy equivalences, the category $\Comp^-(\ldots)$ of differential
complexes in the additive category $\Addinf(B)$, with degrees bounded
from above (the differential operator being of degree $+1$, according
to my preference for cohomology notation, sorry!). An object of
$\D^-(\Bhatab)$ may thus be viewed as being essentially the same as a
differential complex in $\Addinf(B)$ with degrees bounded from above,
and given ``up to homotopism''.  The similar description holds for
$\D^-({B'}\uphat\subab)$, and in terms of these descriptions, the
``integration functor'' (in the abelian context) \eqref{eq:101.7} can
be described by
\begin{equation}
  \label{eq:101.9}
  \mathrm Lf_!\supab(L_\bullet') = \Addinf(L_\bullet'),\tag{9}
\end{equation}
i.e., by applying componentwise the tautological extension
$\Addinf(f)$ of $f$ to the differential complexes in
$\Addinf(B')$. This very concrete description applies notably to the
complex
\begin{equation}
  \label{eq:101.10}
  L_\bullet^{B'/B}\quad\text{or}\quad L_\bullet^f\eqdef f_!\supab(L_\bullet^{B'})\tag{10}
\end{equation}
introduced yesterday, whenever a (quasi-special) integrator
$L_\bullet^{B'}$ for $B'$ has been chosen. Applying this to the case
of the map
\[p_B:B\to\Simplex_0,\]
we get (for a given integrator $L_\bullet^B$ for $B$) an explicit
description of the abelianization of the homotopy type of $B$ in terms
of the chain complex $p_{B!}\supab(L_\bullet^B)$ in \Ab, with the
$n$'th component given by
\begin{equation}
  \label{eq:101.11}
  \bigl(L_\bullet^{B/\mathrm{pt}}\bigr)_n = \bZ^{(I_n)},\tag{11}
\end{equation}
where $I_n$ is the set of indices used for describing $L_n^B$ as the
direct sum of objects of the type $\bZ^{(b_i)}$.

Returning to the case of a general map $f:B'\to B$, maybe I should
still write down the formula generalizing \eqref{eq:100.2} or
\eqref{eq:100.10} of yesterday's notes (pages \ref{p:367} and
\ref{p:369}), relating $L_\bullet^{B'/B}$ to the cohomology properties
of the map $f$, i.e., to cointegration relative to $f$. The
formula\pspage{377} expresses cointegration $\mathrm Rf_*$ with
coefficients coming from downstairs, namely $f^*(K^\bullet)$, where
$K^\bullet$ is any differential complex in \Bhatab{} with degrees
bounded from below (thus defining an object in $\D^+(\Bhatab)$). The
relevant formula is\scrcomment{see section~\ref{sec:139}, bottom of
  p.~\ref{p:588}, for corrections to this formula and
  \eqref{eq:101.12prime}, \eqref{eq:101.13} below\dots}
\begin{equation}
  \label{eq:101.12}
  \mathrm Rf_*(f^*(K^\bullet)) \simeq \bHom^{\bullet\bullet}(L_\bullet^{B'/B},K^\bullet),\tag{12}
\end{equation}
an isomorphism in $\D^+(\Bhatab)$, where $\bHom^{\bullet\bullet}$
designates the double complex in \Bhatab{} obtained by applying
$\bHom$ componentwise, more accurately the associated simple complex,
and where $\bHom$ is the internal $\bHom$ in the category \Bhatab,
namely the (pre)sheaf of additive homomorphisms of a given abelian
(pre)sheaf ($L_n$ say) into another ($K^m$ say). The proof of
\eqref{eq:101.12} is essentially trivial, it is just the computational
interpretation, in terms of using projective resolutions, of the
adjunction formula ``localized on $B$
\begin{equation}
  \label{eq:101.12prime}
  \mathrm Rf_*(\mathrm Lf^*(K^\bullet)) \simeq \mathrm R\bHom(\mathrm
  Lf_!\supab(\bZ_{B'}), K^\bullet),\tag{12'}
\end{equation}
which is a particular case of the more general ``adjunction formula''
\begin{equation}
  \label{eq:101.13}
  \mathrm Rf_*(\mathrm R\bHom(F_\bullet', \mathrm Lf^*(K^\bullet))
  \simeq \mathrm R\bHom(\mathrm Lf_!\supab(F'_\bullet),K^\bullet),\tag{13}
\end{equation}
valid for
\[ \text{$F_\bullet'$ in $\D^-({B'}\uphat\subab)$,}\quad
\text{$K^\bullet$ in $\D^+(\Bhatab)$,}\]
\eqref{eq:101.12prime} following from \eqref{eq:101.13} by taking
$F_\bullet'=\bZ_{B'}$.

\begin{remarks}
  We may view \eqref{eq:101.13}, and its particular case
  \eqref{eq:101.12} or \eqref{eq:101.12prime}, as the main formula
  relating the \emph{homology} and \emph{cohomology} invariants for a
  map $f$ in \Cat, or equivalently, the (abelian) \emph{integration}
  and \emph{cointegration} operations defined by $f$. It now occurs to
  me that this formula, and the variance formalism in which it
  inserts, is valid more generally whenever we have a map $f$ between
  two ringed topoi, such that $f_!$ exists for sheaves of sets, hence
  there exists too a corresponding functor $f_!^{\mathrm{mod}}$ for
  sheaves of modules. The fact that we have been restricting to the
  case of the constant sheaves of rings defined by $\bZ$ isn't
  relevant, and (in the case of topoi defined by objects in \Cat,
  hence with sufficiently many projective sheaves of sets) the
  formalism of the subcategories $\Add(B)$ and $\Addinf(B)$ in
  \Bhatab{} can be generalized equally to arbitrary sheaves of rings
  on $B$. At present, I don't see though any striking particular case
  where this generalization would seem useful.\pspage{378}
\end{remarks}

\hangsection[Abelianization VII: Integrators (for $A\op$) are
\dots]{Abelianization VII: Integrators
  \texorpdfstring{\textup(\kernifitalic{2pt}}{(}for
  \texorpdfstring{$A\op$\textup)}{Aop)} are abelianizators
  \texorpdfstring{\textup(\kernifitalic{2pt}}{(}for
  \texorpdfstring{$A$\textup)}{A)}.}\label{sec:102}%
We now focus attention upon the pair of mutually dual small categories
\begin{equation}
  \label{eq:102.1}
  (A,B), \quad\text{with $B=A\op$, i.e., $A=B\op$,}\tag{1}
\end{equation}
and recall the equivalence of section \ref{sec:93} following from the
universal property of $\Add(B)$
\begin{equation}
  \label{eq:102.2}
  \Ahatab = \bHom(A\op,\Ab) \equeq \bHomadd(\Add(A\op),\Ab),\tag{2}
\end{equation}
which we parallel with the formula \eqref{eq:99.3} of section
\ref{sec:99} (p.\ \ref{p:362}), which reads when replacing in it $B$
by $A$
\[\Ahatab\simeq\bHomadd(\Add(A)\op,\Ab)\text{;}\]
this immediately suggests a canonical equivalence of categories
\begin{equation}
  \label{eq:102.3}
  \Add(A\op)\equeq \Add(A)\op,\tag{3}
\end{equation}
following immediately indeed from the $2$-universal properties of
these categories. We complement \eqref{eq:102.2} by the similar
formula
\begin{equation}
  \label{eq:102.4}
  F\mapsto \widetilde F: \quad
  \Ahatab\toequ\bHomaddinf(\Addinf(B),\Ab), \quad B=A\op,\tag{4}
\end{equation}
where $\bHomaddinf$ denotes the category of infinitely additive
functors between two infinitely additive categories. In view of the
emphasis lately on chain complexes in $\Addinf(B)$ rather than in
$\Add(B)$, in order to reconstruct say the derived category
$\D_\bullet(\Bhatab)$ of chain complexes in \Bhatab, and get existence
of ``integrators'' with components in $\Addinf(B)$ (whereas there may
be none with components in $\Add(B)$), it is formula \eqref{eq:102.4}
rather than \eqref{eq:102.2} which is going to be relevant for our
homology formalism. Using \eqref{eq:102.4}, we get a canonical
biadditive pairing
\begin{equation}
  \label{eq:102.star}
  \Ahatab \times \Addinf(B) \to \Ab,\tag{*}
\end{equation}
which visibly is exact with respect to the first factor, and which we
may equally interpret as a functor
\begin{equation}
  \label{eq:102.5}
  L\mapsto \widetilde L: \Addinf(B) \to \bHomex(\Ahatab,\Ab),\tag{5}
\end{equation}
where $\bHomex$ denotes the category of \emph{exact} (hence additive)
functors from an abelian category to another one.

It can be shown that the pairing \eqref{eq:102.star} can be extended
canonically to a pairing
\begin{equation}
  \label{eq:102.6}
  \Ahatab \times \Bhatab \to \Ab\tag{6}
\end{equation}
commuting to small direct limits in each variable, and identifying
(up\pspage{379} to equivalence) each left hand factor to the category
of functors from the other functor to \Ab{} which commute with small
direct limits (much in the same way as the corresponding relationship
between \Ahat{} and \Bhat, with \Ab{} being replaced by \Sets), and
the accordingly the functor \eqref{eq:102.5} is equally fully
faithful, and extends to a fully faithful functor from \Bhatab{} to
$\bHomadd(\Ahatab,\Ab)$, inducing in fact an equivalence between
\Bhatab{} and the full subcategory $\bHom_!(\Ahatab,\Ab)$ of
$\bHom(\Ahatab,\Ab)$ made up by all functors $\Ahatab\to\Ab$ which
commute to small direct limits. But for what we have in mind at
present, these niceties are not too relevant yet it seems -- all what
matters is that an object $L$ of $\Addinf(B)$ defines an exact functor
\[\widetilde L :\Ahatab\to\Ab,\]
depending functorially on $L$, in an infinitely additive way. Thus, as
noted in section \ref{sec:92} (but where $\Addinf(B)$ was replaced by
the smaller category $\Add(B)$, which has turned out insufficient for
our purposes), whenever we have a chain complex $L_\bullet$ in
$\Addinf(B)$, we get a corresponding functor
\begin{equation}
  \label{eq:102.7}
  \widetilde L_\bullet:\Ahatab\to\Ch_\bullet\Ab\tag{7}
\end{equation}
from \Ahatab{} to the category of chain complexes of \Ab, which is
moreover an \emph{exact} functor. Generalizing slightly the
terminology introduced in section \ref{sec:93}, where we restricted to
chain complexes with components in $\Add(B)$ rather than in $\Add(B)$,
we'll say that $\widetilde L_\bullet$ is an \emph{abelianizator for
  $A$}, if the following diagram commutes up to isomorphism:
\begin{equation}
  \label{eq:102.8}
  \begin{tabular}{@{}c@{}}
    \begin{tikzcd}[baseline=(O.base)]
      \Ahat\ar[r]\ar[d,"{%
        \begin{array}{@{}c@{}}
          \Wh_A \\ \text{(abelianization)}
        \end{array}}"'] & \HotOf_A \ar[r] & \Hot \ar[d,"{%
        \begin{array}{@{}c@{}}
          \text{``absolute''} \\ \text{abelianization} \\ \text{functor}
        \end{array}}" inner sep=1em] \\
      \Ahatab\ar[r,"\widetilde L_\bullet"] &
      \Ch_\bullet\Ab\ar[r] & |[alias=O]| \Hotab
    \end{tikzcd}.
  \end{tabular}\tag{8}
\end{equation}
More accurately, an abelianizator is a pair $(\widetilde L_\bullet,\lambda)$,
where $\lambda$ is an isomorphism of functors $\Ahat\to\Hotab$ making
the diagram commute. Here, I like to view the abelianization functor
\begin{equation}
  \label{eq:102.9}
  \Hot\to\Hotab\eqdef\D_\bullet\Ab\tag{9}
\end{equation}
as the one described directly in section \ref{sec:100} via integrators
of arbitrary modelizing objects in \Cat, without any reference to an
auxiliary test category such as $\Simplex$ or the like.

The point of \eqref{eq:102.8} is that via an ``abelianizator'' for
$A$, we want to be able to give a)\enspace a \emph{simultaneous} handy
expression, in terms of ``computable'' chain complexes in \Ab,
of\pspage{380} abelianization of homotopy types modelized by a
variable object $X$ in \Ahat, and b)\enspace we want that the chain
complex in \Ab{} expressing abelianization of $X$, should be
expressible in terms of the ``tautological abelianization'' $\Wh_A(X)
= \bZ^{(X)}$ of $X$ itself, by a formula moreover which should make
sense functorially with respect to an arbitrary abelian presheaf,
i.e., an object $F$ in \Ahatab.

The main fact I have in view here is that whenever the chain complex
$L_\bullet$ in $\Addinf(B)$ is endowed with an augmentation
\begin{equation}
  \label{eq:102.10}
  L_\bullet\to\bZ_B\tag{10}
\end{equation}
turning it into a resolution of $\bZ_B$, i.e., into a (quasi-special)
\emph{integrator} for $B$, then ipso facto $L_\bullet$ is an
abelianizator for $A$, the commutation isomorphism $\lambda$ being
canonically defined by the augmentation \eqref{eq:102.10}.

Some comments, before proceeding to a proof. Presumably, the converse
of our statement holds too -- namely that the natural functor we'll
get from quasi-special integrators for $B$ to abelianizators
$(L_\bullet,\lambda)$ for $B$ is an equivalence (even an isomorphism!)
between the relevant categories. I don't feel like pursuing this --
the more relevant fact here, whether or not a converse as contemplated
holds, is that we can pin down at any rate a special class of
abelianizators for $A$, namely those which come from (quasi-special)
integrators for $B$, and these abelianizators are defined up to chain
homotopism in $\Addinf(B)$. In this sense, \emph{we get an existence
  and unicity statement for abelianizators in $A$}, as strong as we
possibly could hope for. In practical terms, it would seem, \emph{an
  abelianizator for $A$ will be no more no less than just a
  \textup(quasi-special\textup) integrator for $B$, namely a
  projective resolution of $\bZ_B$ in \Bhatab, whose components
  satisfy a mild extra assumption besides being projective}.

Here, I am struck by a slight discrepancy in terminology, as we would
rather have a correspondence
\[
\begin{cases}
  &\text{integrators for $B$ $\to$ abelianizators for $A$} \\
  &\text{quasi-special int.s for $B$ $\to$ quasi-special abelian.s for $A$,}
\end{cases}\]
and the same for ``special'' integrators and abelianizators. As I
still feel that the general appellation of an ``integrator'' for
\emph{any} projective resolution of $\bZ_B$ is adequate (without
insisting that the components should be in $\Addinf(B)$), this kind of
forces us to extend accordingly still the notion of an abelianizator
for $A$. This does make sense, using the pairing \eqref{eq:102.6}
(which we had dismissed as an ``irrelevant nicety for the time
being''!), and the corresponding equivalence\pspage{381}
\begin{equation}
  \label{eq:102.11}
  \Bhatab \toequ \bHom_!(\Ahatab, \Ab),\tag{11}
\end{equation}
where the index $!$ denotes the full subcategory of $\bHom$ made up
with functors commuting to small direct limits. It is immediate that
projective objects in \Bhatab{} give rise to objects in $\bHom_!$
which are \emph{exact} functors from \Ahatab{} to \Ab, and I'll have to
check that the converse also holds. If so, a chain complex in
\Bhatab{} with projective components can be interpreted as being just
an \emph{arbitrary exact functor commuting to small sums}
\[ \Ahatab \to \Ch_\bullet\Ab,\]
(never minding whether or not it can be described ``computationally''
in terms of objects in $\Add(B)$ or in $\Addinf(B)$!) -- which is all
that is needed in order to complete the diagram \eqref{eq:102.8}, and
wonder if it commutes up to isomorphism! And the most natural
statement here is that this is indeed so whenever this functor, viewed
as a chain complex in the abelian category
\begin{equation}
  \label{eq:102.12}
  \bHom_!(\Ahatab,\Ab),\tag{12}
\end{equation}
is a (projective) resolution of the canonical object $\widetilde\bZ_B$
of the category \eqref{eq:102.12}, coming from the object $\bZ_B$ of
the left-hand side of \eqref{eq:102.11}. Now, this functor is just the
familiar ``direct limit'' functor
\begin{equation}
  \label{eq:102.13}
  \widetilde\bZ_B \simeq \varinjlim_B : \Ahatab \eqdef \bHom(B,\Ab)
  \to \Ab,\tag{13}
\end{equation}
which can be equally interpreted as
\begin{equation}
  \label{eq:102.14}
  \widetilde \bZ_B \simeq p_{A!}\supab : \Ahatab \to \Ab,\tag{14}
\end{equation}
namely (abelian) ``integration'' with respect to the map in \Cat
\[p_A:A\to\Simplex_0.\]
Thus, ultimately, \emph{abelianizators for $A$} (or what we may call
``standard abelianizators'', if there should turn out to be any
others, and that they are worth looking at) \emph{turn out to be no
  more, no less than just a projective resolution, in the category
  \eqref{eq:102.12} of functors from \Ahatab{} to \Ab{} commuting with
  small direct limits, of the most interesting object in the category,
  namely the functor}
\begin{equation}
  \label{eq:102.15}
  p_{A!}\supab \simeq \varinjlim_B : \Ahatab\to\Ab.\tag{15}
\end{equation}
We are far indeed from the faltering reflections of section
\ref{sec:91}, about computing homology and cohomology of homotopy
models described in terms of test categories deduced some way or other
from cellular decompositions of spheres!

\bigbreak

\presectionfill\ondate{11.8.}\pspage{382}\par

\hangsection{Integrators versus cointegrators.}\label{sec:103}%
It has been over three weeks now I haven't been working on the
notes. Most part of this time was spent wandering in the Pyrenees with
some friends (a kind of thing I hadn't been doing since I was a boy),
and touring some other friends living the simple life around there, in
the mountains. I was glad to meet them and happy to wander and breathe
the fresher air of the mountains -- and very happy too after two weeks
to be back in the familiar surroundings of my home amidst the gentle
hills covered with vineyards\ldots Yesterday I resumed mathematical
work -- I had to spend the day doing scratchwork in order to get back
into it, now I feel ready to go on with the notes. I'll have to finish
in the long last with that abelianization story I got into
unpremeditatedly -- which turns out to be essentially the same thing
as some systematics about (commutative) cohomology and homology, in
the context of ``models'' in \Cat, or in a category \Ahat{} (with $A$
is \Cat). We were out for proving a statement about ``integrators''
for a small category $B$ being ``abelianizators'' for the dual
category $A=B\op$. The proof I had in mind for this is somewhat
indirect via cohomology, and follows the proof I gave myself a very
long time ago (in case $A\Simplex$), that the usual semi-simplicial
boundary operations do give the correct (topos-theoretic) cohomology
invariants for any object $X$ in \Ahat{} (i.e., any semisimplicial
set), for any locally constant coefficients on $A_{/X}$. The idea was
to replace $A_{/X}$ by the dual category
$(A_{/X})\op=\preslice{A\op}X$ (which, according to a nice result of
Quillen, has a homotopy type canonically isomorphic to the one defined
by $A_{/X}$), and use the canonical functor
\[f = (p_X)\op : (A_{/X})\op \to A\op\eqdef B,\]
which is a cofibration with discrete fibers, and hence gives rise, for
any abelian presheaf $F$ on the category $C$ upstairs, to an
isomorphism
\[\mathrm R \Gamma(C,F) \simeq \mathrm R\Gamma(B,f_*(F))\]
(due to $\mathrm Rf_*(F) \fromsim f_*(F)$, as $f_*$ is exact,
due to the fact that $f$ is a cofibration with discrete fibers). We'll
get Quillen's result about the isomorphism $C\simeq C\op$ in \Hot, for
any object $C$ in \Cat, very smoothly in part \ref{ch:VI}, as a result
of the asphericity story of part \ref{ch:IV}. However, I now realize
that the proof of the fact about abelianizators via Quillen's result
and cohomology is rather awkward, as what we're after now is typically
a result on homology, not cohomology -- and I was really turning it
upside down in order to fit it at all costs into the\pspage{383} more
familiar (to me) cohomology pot! Therefore, I'm not going to write out
this proof, as ``the'' natural proof is going to come out by itself,
once we got a good conceptual understanding of homology, cohomology
and abelianization, in the context of ``spaces'' embodies by objects
of \Cat. Thus, I feel what is mainly needed now is an overall review
of the relevant notions and facts along these lines -- most of which
we've come in touch with before, be it only ``en passant''.

Before starting, just an afterthought on terminology. It occurred to
me that the name of an ``integrator'' (for $A$), for a projective
resolution of the constant abelian presheaf $\bZ_A$ in \Ahatab, is
inaccurate -- as it was meant to suggest that its main use is for
allowing computation, for an arbitrary abelian presheaf (or complex of
such presheaves) $F$ on $A$, of $\mathrm R\Gamma(A,F)$, which we were
thinking of by that time as the ``integration'' of $F$ over $A$ (or
over the associated topos). But it has turned out that for the sake of
coherence with a broader use of the notions of ``integration'' and
``cointegration'' (compare section \ref{sec:69}), the appropriate
designation of $\mathrm R\Gamma(A,F)$ is ``\emph{co}integration'' of
$F$ over $A$, not integration. Therefore, the appropriate designation
for a projective resolution $L_\bullet^A$ of $\bZ_A$, allowing
computational expression of cointegration, is ``cointegrator'' (for
$A$) rather than ``integrator''. On the other hand, in terms of the
dual category $B=A\op$, it turns out that such $L_\bullet^A$ allows
computational expression of \emph{integration} (i.e., homology) over
$B$, and therefore it seems adequate to call $L_\bullet^A$ also an
\emph{``integrator'' for $B$}. Moreover, it turns out that such an
integrator for $B$ is equally an ``abelianizator'' for $B$, i.e., it
allows simultaneous computation of the ``abelianizations'' of the
homotopy types defined by arbitrary objects $X$ in \Ahat, in terms of
the abelianization $\Wh_A(X)=\bZ^{(X)}$ of $X$ (cf.\ sections
\ref{sec:93} and \ref{sec:102}) -- and possibly the converse holds
too. Whether this is so or not, there doesn't seem at present much
sense to bother about abelianizators which do not come from
integrators, while the latter have the invaluable advantage (besides
mere existence) of being unique up to homotopism. Thus, in practical
terms, it would seem that abelianizators (for a given small category
$B$) are no more no less than just integrators (for the same $B$,
i.e., cointegrators for $A=B\op$) -- and I would therefore suggest to
simply drop the designation ``abelianizator'' for the benefit of the
synonym ``integrator'', which fits more suggestively into the pair of
dual notions integrator---cointegrator.\pspage{384}

\hangsection{Overall review on abelianization
  \texorpdfstring{\textup{(1)}}{(1)}: Case of
  pseudo-topoi.}\label{sec:104}%
I'll have after all to give a certain amount of functorial ``general
non-sense'' which I've tried to bypass so far.

\subsection{Pseudo-topoi and adjunction equivalences.}
\label{subsec:104.A}
In what follows, ordinary capital letters as $A,B,\ldots$ will
generally denote small categories (mostly objects in \Cat), whereas
round capital letters $\scrA,\scrB,\scrM$ will denote \scrU-categories
which may be ``large'', for instance $\scrA=\Ahat$, $\scrB=\Bhat$,
etc.  For two such categories $\scrA,\scrB$, we denote by
\begin{equation}
  \label{eq:104.1}
  \bHom_!(\scrA,\scrB), \quad \bHom^!(\scrA,\scrB)\tag{1}
\end{equation}
the full subcategories of the functor category $\bHom(\scrA,\scrB)$,
made up with all functors which commute with small direct or inverse
limits respectively. This notation is useful mainly in case \scrA{}
and \scrB{} are stable under small direct resp.\ inverse limits, in
which case the same holds true for the corresponding category
\eqref{eq:104.1}, because as a full subcategory of
$\bHom(\scrA,\scrB)$ (where direct resp.\ inverse limits exist and are
computed componentwise) it is stable under direct resp.\ inverse
limits. Thus, the inclusion functors
\begin{equation}
  \label{eq:104.2}
  \bHom_!(\scrA,\scrB) \to \bHom(\scrA,\scrB),\quad
  \bHom^!(\scrA,\scrB) \to \bHom(\scrA,\scrB)\tag{2}
\end{equation}
commute with direct resp.\ inverse limits, i.e., those limits in the
categories \eqref{eq:104.1} are computed equally componentwise.

The canonical inclusion
\[\scrA \hookrightarrow \bHom(\scrA\op,\Sets)\]
factors into a fully faithful inclusion functor
\begin{equation}
  \label{eq:104.3}
  \scrA \hookrightarrow \bHom^!(\scrA\op,\Sets).\tag{3}
\end{equation}
Let's recall the non-trivial useful result:
\begin{propositionnum}\label{prop:104.1}
  Assume the \scrU-category \scrA{} is stable under small direct
  limits, and admits a small full subcategory $A$ which is
  ``generating for monomorphisms'', i.e., any monomorphism $i:X\to Y$
  in \scrA{} such that $\Hom(Z,i):\Hom(Z,X)\to\Hom(Z,Y)$ is bijective
  for any $Z$ in $C$, is an isomorphism. Then the fully faithful
  functor \eqref{eq:104.3} is an equivalence, i.e., any functor
  \[\scrA\op\to\Sets\]
  that commutes with small inverse limits is representable.
\end{propositionnum}

For a proof,\scrcomment{\textcite{SGA4vol1}} see SGA~4~I~8.12.7.

\begin{corollarynum}\label{cor:104.prop1.1}
  If \scrA{} satisfies the assumptions above, then \scrA{} is equally
  stable under small \emph{inverse} limits.
\end{corollarynum}

For\pspage{385} the sake of brevity, we'll say that a \scrU-category
satisfying the assumptions of prop.\ \ref{prop:104.1} is a
\emph{pseudo-topos} (as these conditions are satisfied for any
topos). We get at once the
\begin{corollarynum}\label{cor:104.prop1.2}
  Let $\scrA,\scrB$ be two pseudo-topoi. Then a functor from $\scrA$
  to $\scrB\op$ \textup(resp.\ from $\scrB\op$ to \scrA\textup) has a
  right adjoint \textup(resp.\ a left adjoint\textup) if{f} it
  commutes to small direct limits \textup(resp.\ to small inverse
  limits\textup). Thus, taking right and left adjoints we get two
  equivalences of categories, quasi-inverse to each other
  \begin{equation}
    \label{eq:104.4}
    \bHom_!(\scrA,\scrB\op)\leftrightarrows\bHom^!(\scrB\op,\scrA),\tag{4}
  \end{equation}
  and the two members of \eqref{eq:104.4} are canonically equivalent
  to the category
  \begin{equation}
    \label{eq:104.5}
    \bHom^{!!}(\scrA\op,\scrB\op; \Sets)\tag{5}
  \end{equation}
  of functors
  \[\scrA\op\times\scrB\op\to\Sets\]
  which commute with small inverse limits with respect to either
  variable \textup(the other being fixed\kern1pt\textup), \eqref{eq:104.5}
  being viewed as a full subcategory of $\bHom(\scrA\op\times\scrB\op,
  \Sets)$.
\end{corollarynum}
\begin{remarks}
  \namedlabel{rem:104.1}{1})\enspace The first-hand side of
  \eqref{eq:104.4} is tautologically isomorphic to the category
  $\bHom^!(\scrA\op,\scrB)$ (as for any two \scrU-categories we have
  the tautological isomorphism
  \begin{equation}
    \label{eq:104.6}
    (\bHom_!(\mathscr P,\mathscr Q))\op \simeq \bHom^!(\mathscr P\op,
    \mathscr Q\op)\quad\text{),}\tag{6}
  \end{equation}
  thus, the equivalence \eqref{eq:104.4} can be seen more
  symmetrically as an equivalence
  \begin{equation}
    \label{eq:104.4prime}
    \bHom^!(\scrA\op,\scrB)\simeq\bHom^!(\scrB\op,scrA),\tag{4'}
  \end{equation}
  both categories being equivalent to \eqref{eq:104.5} using the
  equivalence \eqref{eq:104.3} for the second, and the corresponding
  equivalence for \scrB{} for the first, plus the tautological
  isomorphism
  \begin{equation}
    \label{eq:104.7}
    \bHom^!(\scrP, \bHom^!(\scrQ,\scrM)) \simeq
    \bHom^{!!}(\scrP, \scrQ; \scrM),\tag{7}
  \end{equation}
  for any three \scrU-categories $\scrP,\scrQ,\scrM$.

  \namedlabel{rem:104.2}{2})\enspace When \scrA{} is a topos, \scrB{}
  any \scrU-category stable under small inverse limits, then we may
  interpret the category $\bHom^!(\scrA\op,\scrB)$ as the category of
  \emph{\scrB-valued sheaves} on the topos (defined by) \scrA. When
  \scrA{} and \scrB{} are both topoi, then the equivalence
  \eqref{eq:104.4} states that \scrB-valued sheaves on \scrA{} can be
  identified with \scrA-valued sheaves on \scrB, and both may be
  identified with set-valued ``bi-sheaves'' on $\scrA\times\scrB$. In
  case the topoi \scrA, \scrB{} are defined respectively by
  \scrU-sites $A$, $B$ (not necessarily small ones),\pspage{386} these
  bisheaves can be interpreted in a rather evident way as bisheaves on
  $A\times B$, namely functors
  \[A\op\times B\op\to\Sets\]
  which are sheaves with respect to each variable (the other being
  fixed). It is easy to check that the category of all such bisheaves
  is again a topos, and that the latter is a $2$-product of the two
  topoi \scrA, \scrB{} in the $2$-category of all topoi -- it plays
  exactly the same geometrical role as the usual product for two
  topological spaces\ldots
\end{remarks}

\subsection{Abelianization of a pseudo-topos.}
\label{subsec:104.B}
Let \scrA{} be a pseudo-topos, and let's denote by
\begin{equation}
  \label{eq:104.8}
  \scrAab\tag{8}
\end{equation}
the category of abelian group-objects in \scrA. It is immediate that
the forgetful functor
\begin{equation}
  \label{eq:104.9}
  \scrAab\to\scrA\tag{9}
\end{equation}
commutes with small direct limits (and that such limits exist in
$\scrAab$, whereas they exist in \scrA{} by cor.\ \ref{cor:104.prop1.1}
above) -- thus, we may expect that this functor admits a left
adjoint. When so, this will be denoted by
\begin{equation}
  \label{eq:104.10}
  \Wh_\scrA : \scrA\to\scrAab,\tag{10}
\end{equation}
we'll write also
\[\Wh(X)=\bZ^{(X)}\]
when no confusion may arise. The abelianization functor exists for
instance when \scrA{} is a topos, in this case it is well-known that
$\scrAab$ is not only an additive category, but an \emph{abelian}
category with small filtering direct limits which are exact, and a
small generating subcategory. This in turn ensures, as well-known too,
that any object of $\scrAab$ can be embedded into an injective one,
and from this follows (cf.\ SGA~4
I~7.12)\scrcomment{\textcite{SGA4vol1}} that $\scrAab$ admits also a
small full subcategory which is \emph{co}generating with respect to
epimorphisms, in other words that $\scrAab$ is not only an abelian
pseudo-topos, but that the dual category $(\scrAab)\op$ is a
pseudo-topos too. Conversely (kind of), without assuming \scrA{} to be
a topos, if we know some way or other (but this may be hard to check
directly\ldots) that $(\scrAab)\op$ is a pseudo-topos, then it follows
from cor.\ \ref{cor:104.prop1.2} above that the abelianization functor
$\Wh_\scrA$ exists, and this in turn implies that $\scrAab$ is a
pseudo-topos, i.e., admits a small full subcategory which is
generating with respect to monomorphisms (as we see by taking such a
full subcategory\pspage{387} $A$ in \scrA, and the full subcategory in
$\scrAab$ generated by $\Wh_\scrA(A)$).

Let now \scrM{} be an \emph{additive} \scrU-category, which is
moreover a \emph{pseudo-cotopos}, i.e., the dual category $\scrM\op$
is a pseudo-topos. Using twice the corollary \ref{cor:104.prop1.2} above,
for the pair of pseudotopoi $(\scrA,\scrM\op)$ and
$(\scrAab,\scrM\op)$, we get the sequence of equivalences of
categories
\[\bHom_!(\scrA,\scrM)\equeq\bHom^!(\scrM,\scrA)\op \fromequ
\bHom^!(\scrM,\scrAab)\op \equeq \bHom_!(\scrAab,\scrM),\]
where the second equivalence of categories comes from the fact that
any functor
\[f:\scrM\to\scrA\]
from an \emph{additive} category \scrM{} to a category \scrA, which
commutes with finite products, factors canonically through
$\scrAab\to\scrA$. We are interested now in the composite
equivalence
\begin{equation}
  \label{eq:104.11}
  \bHom_!(\scrA,\scrM)\equeq\bHom_!(\scrAab,\scrM),\tag{11}
\end{equation}
defined under the only assumption that \scrM{} is additive and the
categories \scrA, $\scrAab$ and $\scrM\op$ are pseudotopoi
(without having to assume the existence of the abelianization functor
$\Wh_\scrA$). This equivalence is functorial for variable additive
pseudo-cotopos \scrM, when we take as ``maps'' $\scrM\to\scrM'$
functors which commute to small direct limits (a fortiori, these are
right exact and hence additive). In case $\scrAab$ itself is among
the eligible \scrM's, i.e., is a pseudo-cotopos (not only
pseudotopos), we may say that \emph{$\scrAab$ $2$-represents the
  $2$-functor $\scrM\mapsto\bHom_!(\scrA,\scrM)$} on the $2$-category
of all additive pseudo-cotopoi and functors between these commuting to
small direct limits. As we noticed above, the assumption just made
implies that $\Wh_\scrA$ exists. On the other hand, assuming merely
existence of $\Wh_\scrA$ (besides \scrA{} being a pseudo-topos), which
implies that $\scrAab$ is equally a pseudotopos as we say above,
it is readily checked that the equivalence \eqref{eq:104.11} can be
described as
\begin{equation}
  \label{eq:104.12}
  F\mapsto F\circ\Wh_\scrA : \bHom_!(\scrAab,\scrM) \toequ
  \bHom_!(\scrA,\scrM).\tag{12} 
\end{equation}
Thus, we get the
\begin{propositionnum}\label{prop:104.2}
  Let \scrA{} be a pseudotopos such that the abelianization functor
  \eqref{eq:104.10} exists \textup(\kern2pt for instance \scrA{} a
  topos\textup). Then $\scrAab$ is a pseudotopos and an additive
  category. Moreover, for any additive category \scrM{} which is a
  pseudo-cotopos, the functor \eqref{eq:104.12} is an equivalence of
  categories. 
\end{propositionnum}
\setcounter{corollarynum}{0}
\begin{corollarynum}\label{cor:104.prop2.1}
  Let\pspage{388} \scrA{} be a pseudotopos such that $\scrAab$ is
  a pseudotopos \textup(\kern2pt for instance, \scrA{} is a
  topos\textup). Then the abelianization functor $\Wh_\scrA$ exists,
  and this functor is $2$-universal for functors from \scrA{} into
  \scrU-categories which are both additive and are pseudotopoi
  \textup(maps between these being functors which commute with small
  direct limits\textup).
\end{corollarynum}

By duality, using \eqref{eq:104.6}, we can restate the equivalence
\eqref{eq:104.12} as
\begin{equation}
  \label{eq:104.13}
  F\mapsto F\circ\Wh\op : \bHom^!(\scrAab\op,\scrN) \toequ
  \bHom^!(\scrA\op,\scrN),\tag{13}
\end{equation}
valid provided \scrA{} is a pseudo-topos, $\Wh_\scrA$ exists, and
\scrN{} is an additive category which is moreover a pseudotopos.

Take for instance $\scrN=\Ab$, the category of abelian groups, i.e.,
\[\scrN=\scrB\subab, \quad\text{where $\scrB=\Sets$,}\]
then as already noticed above the left hand side of \eqref{eq:104.13}
is canonically equivalent with
$\bHom^!(\scrAab\op,\scrB)=\bHom^!(\scrAab\op,\Sets)$, as
$\scrAab\op$ is additive; on the other hand, by prop.\
\ref{prop:104.1} applied to the pseudotopos $\scrAab$ we get
\[\scrAab\toequ \bHom^!(\scrAab\op,\Sets) \quad\bigl(\fromequ
\bHom^!(\scrAab\op,\Ab)\bigr),\]
and hence an equivalence
\begin{equation}
  \label{eq:104.14}
  \scrAab \toequ\bHom^!(\scrA\op,\Ab), \quad
  F\mapsto\bigl(X\mapsto\Hom(X,F)\bigr),\tag{14}
\end{equation}
valid whenever \scrA{} is a pseudotopos such that $\Wh_\scrA$
exists. When \scrA{} is a topos, this corresponds to the familiar fact
that an abelian group object in the category of sheaves (of sets) on a
topos, can be interpreted equally as a sheaf on the topos with values
in the category \Ab{} of abelian groups.

\subsection{Interior and exterior operations
  \texorpdfstring{$\otimes_{\bZ}$ and $\Hom_{\bZ}$}{tensorZ and
    HomZ}.}\label{subsec:104.C}
In the first place, I want to emphasize the basic tensor product
operation
\begin{equation}
  \label{eq:104.15}
  (F,G)\mapsto F\otimes_{\bZ} G : \scrAab\times\scrAab \to \scrAab\tag{15}
\end{equation}
between abelian group objects of the pseudotopos \scrA, defined as
usual argumentwise as the solution of the universal problem, expressed
by the ``Cartan isomorphism''
\begin{equation}
  \label{eq:104.16}
  \Hom_\scrAab(F\otimes G, H) \simeq \Bil_{\bZ}(F,G ; H),\tag{16}
\end{equation}
where we dropped the subscript $\bZ$ in the tensor product, and where
$\Bil_{\bZ}$ or simply $\Bil$ denotes the set of maps $F\times G\to H$
which are ``biadditive'' in the usual sense of the word. Thus, the
existence of \eqref{eq:104.15}\pspage{389} just means, by definition,
that for any pair $(F,G)$ of objects in \scrAab, the functor in
$H$
\[H \mapsto \Bil(F,G;H)\]
is representable. It is clear that this functor commutes with small
inverse limits, hence by prop.\ \ref{prop:104.1} it is representable,
provided we know that \scrAab{} is a pseudo-cotopos (for instance,
when \scrA{} is a topos, in which case the existence of tensor
products is anyhow a familiar fact). The familiar Bourbaki
construction of a tensor product amounts on the other hand to viewing
$F\otimes G$ as a quotient of $\Wh_\scrA(F\times G)=\bZ^{(F\times G)}$ by
suitable ``relations'', i.e., as the cokernel of a map in
\scrAab
\[L_1\to L_0=\Wh_\scrA(F\times G),\]
where, as a matter fact, $L_1$ can be described as
\[L_1=\Wh_\scrA(F\times F\times G) \times \Wh_\scrA(F\times G\times
G).\]
Thus, if we know beforehand that cokernels exist in \scrAab{}
(which would indeed follow from \scrAab{} being a pseudotopos, but
may be checked more readily in terms of suitable exactness properties
of \scrA{} directly), plus the existence of course of $\Wh_\scrA$, the
tensor product functor \eqref{eq:104.15} exists. (On the other hand,
no use is made here of the assumption that \scrA{} be a pseudotopos.)

Let's assume the tensor product functor \eqref{eq:104.15} exists. Then
it is readily checked it is associative and commutative up to
canonical isomorphisms, giving rise to the usual compatibilities. If
moreover $\Wh_\scrA$ exists, we readily get the canonical isomorphism
\begin{equation}
  \label{eq:104.17}
  \Wh_\scrA(X\times Y)\eqdef\bZ^{(X\times Y)} \fromsim
  \Wh_\scrA(X)\otimes \Wh_\scrA(Y)\quad \bigl(\eqdef
  \bZ^{(X)}\otimes\bZ^{(Y)}\bigr),\tag{17} 
\end{equation}
compatible of course with the commutativity and associativity
isomorphisms for the operations $\times$ and $\otimes$. If on the
other hand \scrA{} admits moreover a final object (as it does if
\scrA{} is a pseudotopos and hence stable under small direct limits),
then
\begin{equation}
  \label{eq:104.18}
  \bZ_\scrA \eqdef \Wh_\scrA(e) = \bZ^{(e)}\tag{18}
\end{equation}
is a two-sided unit for the tensor product operation.

In what follows, we are interested in categories of the type
\begin{equation}
  \label{eq:104.19}
  \scrA^\scrM \eqdef\bHom_!(\scrA,\scrM), \quad\scrA_\scrN
  \eqdef\bHom^!(\scrA\op,\scrN)\tag{19} 
\end{equation}
where now \scrA{} is assumed to be a fixed pseudotopos, and \scrM{}
and \scrN{} are\pspage{390} \emph{additive} categories, \scrM{} being
moreover a pseudo-cotopos, \scrN{} a pseudotopos. We assume moreover
that $\Wh_A$ exists, and hence the categories \eqref{eq:104.19} can be
interpreted up to equivalence, via \eqref{eq:104.12} and
\eqref{eq:104.13}, as
\begin{equation}
  \label{eq:104.19prime}
  \bHom_!(\scrAab, \scrM), \quad \bHom^!(\scrAab\op,\scrN).\tag{19'}
\end{equation}
Let's remark that the dual of a category of one of the types
\eqref{eq:104.19} (or equivalently, \eqref{eq:104.19prime}) is
isomorphic to a category of the other type, more accurately, by
\eqref{eq:104.6} we get
\begin{equation}
  \label{eq:104.20}
  \bHom_!(\scrA,\scrM)\op\tosim\bHom^!(\scrA\op,\scrN),\tag{20}
\end{equation}
i.e., $(\scrA^\scrM)\op\tosim\scrA_\scrN$, with $\scrN=\scrM\op$. In
case \scrA{} is a topos, the objects of the second category
$\scrA_\scrN$ in \eqref{eq:104.19} (or equivalently, in
\eqref{eq:104.19prime}) can be interpreted as \emph{\scrN-valued
  sheaves} on the topos \scrA, whereas the object of the first,
$\scrA^\scrM$, may be called, correspondingly, \emph{cosheaves on
  \scrA{} with values in \scrM}. Thus, in virtue of \eqref{eq:104.20},
\scrM-valued cosheaves on \scrA{} can be interpreted as
$scrN=\scrM\op$-valued sheaves on \scrA, the corresponding categories
of cosheaves and sheaves being however \emph{dual} to each other. In
the next subsection \ref{subsec:105.D}, when $\scrA=\Ahat$, we'll
interpret moreover \scrM-valued cosheaves on \scrA{} (or on $A$, as
we'll call them equivalently) as \scrM-valued \emph{sheaves} on the
(topos associated to the) dual category $B=A\op$, and in this context
the difference between the categories of cosheaves and of sheaves
(which for the time being appear as categories \emph{dual} to each
other) will disappear altogether, provided we allow the ground topos
\scrA{} to change (from $\scrA=\Ahat$ to the ``dual'' topos
$\scrB=\Bhat$).

In terms of the expressions \eqref{eq:104.19prime} of the category of
``cosheaves'' and ``sheaves'' we are interested in, we want now to
define an external operation of the fixed category \scrAab{} on those
categories, using the tensor product operation \eqref{eq:104.15} in
\scrAab. What is needed visibly for this end is that for fixed $F$,
the functor
\begin{equation}
  \label{eq:104.21}
  G\mapsto F\otimes G\tag{21}
\end{equation}
from \scrAab{} to itself should commute with small direct limits --
hence composing it with a functor $L$ in $\bHom_!(\scrAab,\scrM)$ will
yield a functor in the same category, which will be the looked-for
external tensor product $F\oast L$, i.e.,
\begin{equation}
  \label{eq:104.22}
  \mathop{F\oast L}(G) = L(F\otimes G),\tag{22}
\end{equation}
which we may write more suggestively as
\begin{equation}
  \label{eq:104.22prime}
  G * (F \oast L) = (G \otimes F) * L,\tag{22'}
\end{equation}
with the notation
\begin{equation}
  \label{eq:104.23}
  H * L \eqdef L(H), \quad\text{for $H$ in \scrAab, $L$ in
    $\scrA^\scrM=\bHom_!(\scrA,\scrM)$,}\tag{23}
\end{equation}
which\pspage{391} will be convenient mainly in the context of the next
subsection (when \scrA{} is of the type \Ahat).

The exactness property needed for the functor \eqref{eq:104.21} is
equivalent with the property that for any object $H$ in \scrAab, the
functor
\[G\mapsto \Hom(F\otimes G,H) \simeq \Bil(F,G;H)\]
from $\scrAab\op$ to \Sets{} commute with small inverse limits. As
\scrAab{} is a pseudotopos (prop.\ \ref{prop:104.2}), this is
equivalent by prop.\ \ref{prop:104.1} with this functor being
representable. By definition of $\Bil$, we get
\[\Bil(F,G;H)\simeq\Hom_{\scrA\subab\uphat}(G, \bHom_{\bZ}(F,H)),\]
where the object
\begin{equation}
  \label{eq:104.24}
  \bHom_{\bZ}(F,H)\tag{24}
\end{equation}
is taken in the category of presheaves $\scrA\subab\uphat$ (cheating a
little with universes here\ldots). To sum up, the condition we want
for \eqref{eq:104.21} just boils down to the representability of the
abelian group objects \eqref{eq:104.24} in $\scrA\subab\uphat$, for
any two objects $F,G$ in \scrAab, i.e., essentially to the existence
of ``internal $\bHom$'s'' in the category \scrAab{} (endowed with the
tensor product $\otimes_{\bZ}$), satisfying the familiar Cartan
isomorphism formula
\begin{equation}
  \label{eq:104.25}
  \Hom(F\otimes G,H) \simeq\Hom(G,\bHom(F,H)).\tag{25}
\end{equation}

To sum up, what is needed for a nice formalism of interior and
exterior tensor products and $\Hom$'s for ``sheaves'' and
``presheaves'' on the pseudotopos \scrA, are the following assumptions
on \scrA:
\begin{enumerate}[label=\arabic*)]
\item\label{it:104.C.1}
  Tensor products and corresponding internal $\bHom$'s exist in
  \scrAab, and
\item\label{it:104.C.2}
  the abelianization functor $\Wh_\scrA$ \eqref{eq:104.15} exists,
\end{enumerate}
the latter assumption being needed in order to feel at ease with the
equivalence between the categories \eqref{eq:104.19} of sheaves and
cosheaves, and their ``abelianized'' interpretations
\eqref{eq:104.19prime}.

Under these assumptions, we define the exterior tensor product
operation of \scrAab{} upon a category of cosheaves $\scrA^\scrM$
\begin{equation}
  \label{eq:104.26}
  \scrAab\times\scrA^\scrM \to \scrA^\scrM, \quad
  (F,L)\mapsto F\oast L,\tag{26}
\end{equation}
by formula \eqref{eq:104.22} (which may be written also under the form
\eqref{eq:104.22prime}). This operation has the obvious associativity
property
\begin{equation}
  \label{eq:104.27}
  F\oast(G\oast L)\simeq (F\otimes G)\oast L,\tag{27}
\end{equation}
and moreover the unit $\bZ_\scrA$ for the internal tensor product in
\scrAab{} operators on $\scrA^\scrM$ as the identity functors,
i.e.,\pspage{392}
\begin{equation}
  \label{eq:104.28}
  \bZ_\scrA\oast L\simeq L.\tag{28}
\end{equation}

Using the tautological duality relation \eqref{eq:104.20} between
categories of cosheaves (with values in \scrM) and categories of
sheaves (with values in $\scrN=\scrM\op$), we deduce accordingly an
associative and unitary operation of \scrAab{} on any category of the
type $\scrA_\scrN$, namely \scrN-valued sheaves on \scrA. This
operation is most conveniently denoted by the $\bHom$ symbol
\begin{equation}
  \label{eq:104.29}
  (F,K)\mapsto \bHom(F,K): \scrAab\op\times\scrA_\scrN\to\scrA_\scrN,\tag{29}
\end{equation}
its explicit description in terms of $K$, viewed as a functor
\[K:\scrAab\op\to\scrN\]
is by
\begin{equation}
  \label{eq:104.30}
  \bHom(F,L)(G)=K(G\otimes F)\tag{30}
\end{equation}
for $K$ in $\scrA_\scrN$, $F$ and $G$ in \scrAab. This may be written
more suggestively as
\begin{equation}
  \label{eq:104.30prime}
  \Hom(G, \bHom(F,L)) \simeq \Hom(G\otimes F, K),\tag{30'}
\end{equation}
similar to \eqref{eq:104.22prime}, where we use the notation $\Hom$
(non-bold!)\scrcomment{In the typescript this parenthetical remark
  says ``non-underlined!'', as internal $\bHom$'s there are underlined
  rather than in boldface.} in analogy to \eqref{eq:104.23} for
denoting $K(H)$, namely
\begin{equation}
  \label{eq:104.31}
  \Hom(H,K)\eqdef K(H).\tag{31}
\end{equation}
If confusion is feared, we may put a subscript $\bZ$ in all $\Hom$'s
and $\bHom$'s just introduced, as well as in the internal and external
tensor product operations $\otimes$, $*$, $\oast$.

\begin{comments}
  Take for instance
  \[\scrN=\Ab,\quad\text{hence $\scrA_\scrN\simeq\scrAab$}\]
  by prop.\ \ref{prop:104.1} (compare \eqref{eq:104.14}). By this
  equivalence, the $\Hom$ in \eqref{eq:104.31} is just the usual
  $\Hom$ set corresponding to the category structure of \scrAab{} (the
  $\Hom$ endowed moreover with the structure of abelian group, coming
  from the fact that \scrAab{} is an additive category), whereas
  formula \eqref{eq:104.30prime} shows that the ``external'' $\bHom$
  in this case is nothing but the usual internal $\Hom$ as
  contemplated in \eqref{eq:104.24}. This I hope will convince the
  reader of the adequacy of the notation used (in case of exterior
  operation of abelian sheaves on \scrN-valued sheaves), and of the
  convention \eqref{eq:104.31}. We would like to give a similar
  justification for the notations \eqref{eq:104.22} \eqref{eq:104.23}
  used in connection with operation of abelian sheaves on \scrM-valued
  emph{co}sheaves, by interpreting this as the internal tensor product
  operation in \scrAab, for suitable choice of \scrN. This I am afraid
  cannot be done for\pspage{393} an arbitrary \scrA{} satisfying our
  assumptions, even when \scrA{} is a topos and even when it is of the
  special type \Ahat, as I do not know of any \scrM{} such that
  \begin{equation}
    \label{eq:104.star}
    \scrA^\scrM\simeq\scrAab.\tag{*}
  \end{equation}
  However, in case $\scrA=\Ahat$, introducing the dual category
  $B=A\op$ and $\scrB=\Bhat$, we get (see \ref{subsec:105.D} below)
  \begin{equation}
    \label{eq:104.starstar}
    \scrA^\scrM \equeq \scrB_\scrM \equeq
    \scrB\subab\quad\text{if}\quad \scrM=\Ab,\tag{**}
  \end{equation}
  hence we do get a canonical isomorphism \eqref{eq:104.star} provided
  $A=B$ say and hence $\scrA=\scrB$. Let's look at any rate at the
  simplest case, namely when $A$ is a final object in \Cat, hence
  \scrA{} can be identified with \Sets, and \scrAab{} with \Ab, the
  identification between the categories \Ab{} and $\bHom_!(\Ab,\Ab)$
  being obtained (we hope!) by associating to every object $L$ in \Ab,
  the functor
  \[F\mapsto L\otimes F: \Ab\to\Ab.\]
  This being so, the external operation \eqref{eq:104.22} of \Ab{} on
  $\bHom_!\equeq\Ab$ can be interpreted (using this identification) as
  the interior tensor product operation in \Ab. On the other hand, the
  operation $*$ of \eqref{eq:104.23} is equally interpreted as nothing
  but the tensor product in \Ab, which justifies the notation
  suggesting a tensor product. It would be nice checking corresponding
  compatibilities for a general object $A$ in \Cat{} satisfying
  $A=A\op$, namely a direct sum of one-object categories $A_i$ defined
  each in terms of a \emph{commutative} monoid $M_i$ -- I didn't work
  it out myself, sorry!

  It should be noted that the relationship between the two exterior
  $\Hom$'s in \eqref{eq:104.29} and \eqref{eq:104.31} is essentially
  the same as between the two exterior tensor-type operations
  \eqref{eq:104.26} and \eqref{eq:104.23}, $\Hom(F,K)$ designating an
  object in \scrN{} and $\bHom(F,K)$ an \scrN-valued sheaf on \scrA,
  just as $F*L$ designates an object in \scrM{} and $F\oast L$ an
  \scrM-valued cosheaf on \scrA; the graphical device of
  \emph{bold-facing}\scrcomment{in the typescript: \emph{underlining}}
  the symbol $\Hom$ (used for sheaves) corresponds to the device of
  \emph{circling} the symbol $*$ (used for cosheaves). With this
  luxury of explanations, I hope the notations introduced here are
  getting through\ldots
\end{comments}

\bigbreak

\presectionfill\ondate{12.8.}\pspage{394} and \ondate{13.8.}\par

\hangsection[Review (2): duality equivalences for ``algebraic'' topoi
and \dots]{Review \texorpdfstring{\textup{(2)}}{(2)}: duality
  equivalences for ``algebraic'' topoi and abelian topoi.}\label{sec:105}%
Let's go on with the overall review of abelianization.
\addtocounter{subsection}{3}

\subsection{Duality for topoi of the type
  \texorpdfstring{\Ahat}{Ahat}, and tentative generalizations.}
\label{subsec:105.D}
The main fact, it seems, which will give rise to duality statements
for topoi of the type \Ahat{} is the following, rather familiar one:
\addtocounter{propositionnum}{2}
\begin{propositionnum}\label{prop:105.3}
  Let $A$ be a small category, \scrM{} a \scrU-category stable under
  small direct limits,
  \begin{equation}
    \label{eq:105.32}
    \varepsilon_A:A\hookrightarrow\Ahat\tag{32}
  \end{equation}
  the canonical inclusion functor. Then the following functor is an
  equivalence of categories:
  \begin{equation}
    \label{eq:105.33}
    F\mapsto F\circ\varepsilon_A: \bHom_!(\Ahat,\scrM) \to \bHom(A,\scrM).\tag{33}
  \end{equation}
\end{propositionnum}

As I am at a loss to give a reference for this standard fact of
category theory, I'll give in guise of a proof the indication that a
quasi-inverse functor for \eqref{eq:105.33} is given by the familiar
construction
\begin{equation}
  \label{eq:105.34}
  i\mapsto i_! : \bHom(A,\scrM)\to\bHom_!(\Ahat,\scrM),\tag{34}
\end{equation}
where for any functor
\[i:A\to\scrM,\]
the functor
\[i_!:\Ahat\to\scrM\]
is defined by the formula
\begin{equation}
  \label{eq:105.35}
  i_!(F) = \varinjlim_{\text{$a$ in $A_{/F}$}} i(a).\tag{35}
\end{equation}
It is readily checked (and we have already used a number of times)
that this functor admits a right adjoint
\begin{equation}
  \label{eq:105.36}
  i^*:\scrM\to\Ahat, \quad i^*(x)=(a \mapsto \Hom(i(a),x)),\tag{36}
\end{equation}
and hence $i_!$ commutes to small direct limits, i.e., is in
$\bHom_!(\Ahat,\scrM)$, hence \eqref{eq:105.34}.
\begin{remark}
  If \scrM{} is a pseudotopos, then by prop.\ \ref{prop:105.3} above
  and by corollary \ref{cor:104.prop1.2} of prop.\ \ref{prop:104.1}
  (p.\ \ref{p:385}) the functor
  \begin{equation}
    \label{eq:105.37}
    i\mapsto i^*: \bHom(A,\scrM) \to \bHom^!(\scrM,\Ahat)\tag{37}
  \end{equation}
  (which in any case is fully faithful) is equally an equivalence of
  categories. 
\end{remark}

The\pspage{395} equivalence \eqref{eq:105.33} of prop.\
\ref{prop:105.3} can be interpreted by saying that the canonical
functor \eqref{eq:105.32} from $A$ to \Ahat{} is \emph{$2$-universal},
for functors of $A$ into \scrU-categories stable under small direct
limits, taking as ``maps'' between such categories functors which
commute to small direct limits. Thus, the \scrU-category \Ahat{} may
be viewed as ``the'' category deduced from $A$ by adding arbitrary
direct limits (disregarding the direct limits which may perchance
already exist in $A$\ldots).

Taking the duals of the two members of \eqref{eq:105.33}, we get an
equivalent statement of prop.\ \ref{prop:105.3}:
\setcounter{corollarynum}{0}
\begin{corollarynum}\label{cor:105.prop3.1}
  Let \scrN{} be a \scrU-category stable under small inverse
  limits. Then the functor
  \begin{equation}
    \label{eq:105.38}
    F\mapsto \varepsilon_A\op : \bHom^!(\Ahat\op,\scrN)\to \bHom(A\op,\scrN)\tag{38}
  \end{equation}
  is an equivalence of categories.
\end{corollarynum}

In terms of the topos
\[\scrA=\Ahat,\]
we may interpret the left-hand side of \eqref{eq:105.33} as the
category of \scrM-valued \emph{cosheaves} on this topos, which by
prop.\ \ref{prop:105.3} can be interpreted (up to equivalence) as the
category of functors from $A$ to \scrM. Dually, the left-hand side of
\eqref{eq:105.38} can be viewed as the category of \scrN-valued
\emph{sheaves} on the topos \scrA, which (via the right-hand side) can
be interpreted up to equivalence as the category of functors
$A\op\to\scrN$, i.e., as the category of \scrN-valued
\emph{presheaves} on $A$. As $A$ endowed with the coarsest
(``chaotic'') site structure is a generating site for the topos \scrA,
the equivalence \eqref{eq:105.38} may be viewed as a particular case
of the familiar fact, according to which (up to equivalence) the
category of \scrN-valued sheaves on a topos can be constructed in
terms of \scrN-valued sheaves on any \scrU-site defining this
topos. (When the site structure is chaotic, then those sheaves are
just arbitrary \scrN-valued presheaves.)

Assume now that the \scrU-category \scrM{} is stable under \emph{both}
types of small limits (direct and inverse). Then applying
\eqref{eq:105.33} for $(A,\scrM)$ and \eqref{eq:105.38} for
$(B,\scrM)$ where $B=A\op$, we get an equivalence
\begin{equation}
  \label{eq:105.39}
  \delta_A^\scrM:\bHom_!(\Ahat,\scrM) \toequ\bHom^!(\Bhat\op,\scrM),\tag{39}
\end{equation}
i.e.\ (as announced in yesterday's notes, p.\ \ref{p:390}), we get:
\begin{corollarynum}\label{cor:105.prop3.2}
  Let $A$ be a small category, $B=A\op$ the dual category, \scrM{} any
  \scrU-category stable under small direct and inverse limits. Then
  the category of \scrM-valued \emph{cosheaves} on the topos \Ahat{}
  is equivalent to\pspage{396} the category of \scrM-valued
  \emph{sheaves} on the topos \Bhat.
\end{corollarynum}

This equivalence, defined up to unique isomorphism, is deduced from
the diagram of canonical equivalences
\begin{equation}
  \label{eq:105.40}
  \begin{tabular}{@{}c@{}} 
    \begin{tikzpicture}[commutative diagrams/every diagram,baseline=(O.base)]
      \node (A) at (-3.5cm,0) {\makebox[6ex][l]{$\bHom_!(\Ahat,\scrM)\eqdef(\Ahat)^\scrM$}};
      \node (B) at (3.5cm,0) {\makebox[6ex][r]{$(\Bhat)_\scrM\eqdef\bHom^!(\Bhat\op,\scrM)$}};
      \node (O) at (-1.2cm,-1.2cm) {$\bHom(A,\scrM)$};
      \node (P) at (1.2cm,-1.2cm) {$\bHom(B\op,\scrM)$};

      \path[commutative diagrams/.cd, every arrow, every label]
      (A) edge node {$\equ$} (O)
      (B) edge node[swap] {$\equ$} (P);
      \path[commutative diagrams/.cd, every arrow, every label, equal]
      (O) to (P);
    \end{tikzpicture},
  \end{tabular}\tag{40}
\end{equation}
and depends upon the choice of a quasi-inverse of the second vertical
equivalence in \eqref{eq:105.40} (which choice can be made, via the
dual of \eqref{eq:105.35}, via the choice of small inverse limits in
\scrM).

\begin{remark}
  It is felt that the way we got the equivalence \eqref{eq:105.39} via
  \eqref{eq:105.40}, the role of $A$ and $B$ in it is symmetric. To
  give a more precise statement, consider the equivalence deduced from
  \eqref{eq:105.39} by passing to the dual categories of the two
  members -- using the tautological isomorphisms \eqref{eq:104.6} of
  page \ref{p:385}, we get an equivalence
  \begin{equation}
    \label{eq:105.39prime}
    (\delta_A^\scrM)' : \bHom^!((\Ahat)\op,\scrN) \toequ
    \bHom_!(\Bhat,\scrN), \quad \text{where
      $\scrN=\scrM\op$;}\tag{39'}
  \end{equation}
  this equivalence is \emph{canonically quasi-inverse to the
    equivalence $\delta_B^\scrN$ in opposite direction}, associated to
  the pair $(B,\scrN=\scrM\op)$ instead of $(A,\scrM)$.
\end{remark}

In the rest of this subsection \ref{subsec:105.D}, we'll elaborate on
some particular cases of the equivalence \eqref{eq:105.39} between
cosheaves and sheaves.

\textbf{Case \namedlabel{case:105.1}{1\textsuperscript{\b o})}.}%
\enspace Assume $\scrM=\Sets$,
then by prop.\ \ref{prop:104.1} (p.\ \ref{p:384}) the right-hand side
of \eqref{eq:105.39} is canonically equivalent to \Bhat{} itself,
hence we get an equivalence
\begin{equation}
  \label{eq:105.41}
  \Bhat \toequ \bHom_!(\Ahat,\Sets).\tag{41}
\end{equation}
If we want to keep track of the symmetry aspect described in the
remark above, we may consider the functor \eqref{eq:105.41} as being
deduced from a canonical ``pairing'' between the categories \Ahat{}
and \Bhat
\begin{equation}
  \label{eq:105.42}
  \delta_A : \Ahat\times\Bhat\to\Sets,\tag{42}
\end{equation}
which is an object in
\[\bHom_{!!}(\Ahat,\Bhat; \Sets),\]
i.e., which commutes to small direct limits in each variable. (For
this interpretation, compare the dual statement contained in formula
\eqref{eq:104.7} of page \ref{p:385} -- and note that
\eqref{eq:105.41}, being an equivalence, commutes to small direct
limits, i.e., is in a category $\bHom_!(\Bhat, \bHom_!(\Ahat,\Sets))$.)
This\pspage{397} pairing gives rise, in a symmetric way, to the
functor \eqref{eq:105.41} (which is an equivalence) \emph{and} to a
functor
\begin{equation}
  \label{eq:105.41prime}
  \Ahat\toequ \bHom_!(\Bhat,\Sets)\tag{41'}
\end{equation}
which (it turns out) is none other (up to canonical isomorphism) than
\eqref{eq:105.41} with $B$ replaced by $A$ (and hence $A$ replaced by
$B$), and therefore is equally and equivalence. \emph{Thus, the
  pairing \eqref{eq:105.42} between the two topoi \Ahat{} and \Bhat{}
  has the remarkable property that it defines an equivalence of each
  of these topoi with the category of \textup(set-valued\textup)
  cosheaves on the other.} I do not know of any other example of a
pair of topoi related in such a remarkable way, which we may express
by saying that the two topoi are ``\emph{dual}'' to each other.

We still have to give an explicit expression for the pairing
\eqref{eq:105.42}, plus a convenient notation. I'll write
\begin{equation}
  \label{eq:105.43}
  \delta_A(F,G) \eqdef F * G \quad \text{for $F$ in \Ahat, $G$ in \Bhat,}\tag{43}
\end{equation}
and I'll use the canonical equivalence (valid for any pair of small
categories $A,B$ -- not necessarily dual to each other -- and any
\scrU-category stable under direct limits), deduced by twofold
application of prop.\ \ref{prop:105.3}, plus the tautological
isomorphism similar to \eqref{eq:104.7} p.\ \ref{p:385}:
\begin{equation}
  \label{eq:105.44}
  \bHom_{!!}(\Ahat,\Bhat;\scrM) \toequ \bHom(A\times B,\scrM), \quad
  F\mapsto F\circ(\varepsilon_A\times\varepsilon_B),\tag{44}
\end{equation}
which shows that \eqref{eq:105.42} is known up to canonical
isomorphism when we know its restriction to the full subcategory
$A\times B=A\times A\op$, identifying as usual an object $a$ of $A$ to
its image in \Ahat, and similarly for $B$. If $b$ is an object of $A$,
we'll denote by $b\op$ the same object viewed as an object of
$A\op=B$. With these conventions (including \eqref{eq:105.43}) we get
the nice formula
\begin{equation}
  \label{eq:105.45}
  a * b\op = \Hom_A(b,a) \quad(\; = \Hom_B(a\op,b\op)), \quad
  \text{for $a,b$ in $A$,}\tag{45}
\end{equation}
which has the required symmetry property -- which, for general objects
$F$ in \Ahat{} and $G$ in \Bhat, can be stated as a bifunctorial
isomorphism
\begin{equation}
  \label{eq:105.46}
  F * G \simeq G * F,\tag{46}
\end{equation}
where the operation $*$ in the first member refers to the pair $(A,B)$,
and in the second to the pair $(B,A)$.

From \eqref{eq:105.45} we easily deduce the more general formula for
$F*G$, when either $F$ or $G$ is in $A$ resp.\ $B$, namely\pspage{398}
\begin{equation}
  \label{eq:105.47}
  a * G \simeq G(a), \quad F * b\op \simeq F(b).\tag{47}
\end{equation}
\begin{remarks}
  We are mainly interested here in abelianization and (commutative)
  homology and cohomology, and hence in sheaves and cosheaves with
  values in additive (even abelian) categories, we are not going to
  use for the time being the relationship between \Ahat{} and \Bhat{}
  just touched upon. We could elaborate a great deal more on it, for
  instance introducing a canonical pairing (more accurately, a
  bi-sheaf) with opposite variance to \eqref{eq:105.42}
  \begin{equation}
    \label{eq:105.48}
    \Ahat\op\times\Bhat\op\to\Sets\tag{48}
  \end{equation}
  (or what amounts to the same, canonical functors adjoint to each
  other
  \[\Ahat\op\to\Bhat,\quad\Bhat\op\to\Ahat\quad),\]
  deduced (via the equivalence dual to \eqref{eq:105.44}
  \[\bHom^{!!}(\Ahat\op,\Bhat\op;\scrN) \toequ \bHom(A\op\times
  B\op,\scrN)\quad)\]
  from the co-pairing
  \[A\op\times B\op=B\times A\to\Sets\]
  given by
  \[(b\op,a)\mapsto\Hom(b,a)\quad(\text{for $a,b$ in $A$}).\]
  The two pairings \eqref{eq:105.42} and \eqref{eq:105.48} can be
  given a common interpretation as $\Hom$-sets in a suitable category
  \[\scrE=\scrE(A),\]
  which is the union of the two full subcategories \Ahat{} and
  $\Bhat\op$ (which may be interpreted as deduced from $A$ by
  adjoining respectively small direct and small inverse limits to it)
  intersecting in the common subcategory $A$, the $\Hom$-sets in the
  two directions between an object $F$ of \Ahat{} and an object $G\op$
  of $\Bhat\op$ (corresponding to an object $G$ of \Bhat) being given
  respectively by the pairings \eqref{eq:105.48} and
  \eqref{eq:105.42}.
  The full relationship between these pairings is most conveniently
  expressed, it seems, by the composition law of maps in $\scrE(A)$,
  and associativity for this law. The symmetry of the situation with
  respect to the pair $(A,B)$ is expressed by the canonical
  isomorphism of categories
  \[\scrE(A)\op\simeq\scrE(B)\quad (\text{where $B=A\op$}).\]
\end{remarks}

\textbf{Case \namedlabel{case:105.2}{2\textsuperscript{\b o})}.}%
\enspace Of direct relevance for the abelianization story
is the particular case of the equivalence \eqref{eq:105.39}, obtained
by taking
\[\scrM=\Ab.\]
Using formula \eqref{eq:104.12} (page \ref{p:387}) for the pair
$(\Ahat,\Ab)$ and formula \eqref{eq:104.14} for \Bhat{} in guise of
\scrA, we get the canonical equivalence\pspage{399}
\begin{equation}
  \label{eq:105.49}
  \Bhatab \toequ \bHom_!(\Ahatab,\Ab),\tag{49}
\end{equation}
which should be viewed as the ``abelian'' analogon of the equivalence
\eqref{eq:105.41} above (corresponding to the case
$\scrM=\Sets$). This equivalence again may be viewed as described (in
analogy to \eqref{eq:105.45}) by a canonical pairing
\begin{equation}
  \label{eq:105.50}
  \Ahatab\times\Bhatab\to\Ab,\tag{50}
\end{equation}
which commutes to small direct limits in each variable, i.e., can be
viewed as an object in the category of ``abelian bi-cosheaves''
\begin{equation}
  \label{eq:105.star}
  \bHom_{!!}(\Ahatab,\Bhatab;\Ab),\tag{*}
\end{equation}
and gives rise simultaneously to the equivalence \eqref{eq:105.49},
and to the symmetric equivalence
\begin{equation}
  \label{eq:105.49prime}
  \Ahatab\toequ\bHom_!(\Bhatab,\Ab)\tag{49'}
\end{equation}
of the category of abelian sheaves on \Ahat{} with the category of
abelian cosheaves on \Bhat{} (which is just \eqref{eq:105.49} with $A$
replaced by $B$, up to canonical isomorphism at any rate).

Using the equivalence
\begin{equation}
  \label{eq:105.51}
  \bHom_{!!}(\Ahatab,\Bhatab;\Ab) \toequ \bHom(A\times B,\Ab)\tag{51}
\end{equation}
(which is a particular case of the evident abelian analogon of the
equivalence \eqref{eq:105.44}), we see that the pairing
\eqref{eq:105.50} is described, up to canonical isomorphism, by its
composition with
\[(a,b)\mapsto (\Wh_\Ahat(a),\Wh_\Bhat(b)): A\times B\to
\Ahatab\times\Bhatab,\]
and the latter, as is readily checked, is given by
\begin{equation}
  \label{eq:105.52}
  \Wh_A(a) *_{\bZ} \Wh_B(b\op) \simeq \bZ^{(\Hom(b,a))} \quad\text{for
    $a,b$ in $A$,}\tag{52}
\end{equation}
where we have written $\Wh_A$ instead of $\Wh_\Ahat$ for brevity and
accordingly for $B$, and where the pairing \eqref{eq:105.50} is
denoted by the symbol $*_{\bZ}$, in analogy with the notation $*$ in
\eqref{eq:105.43}, the index $\bZ$ being added in order to avoid
confusion with the non-abelian case \eqref{eq:105.43} (and the index
being dropped when no such confusion is to be feared). The formula
\eqref{eq:105.52} can be written, with different notations
\begin{equation}
  \label{eq:105.52prime}
  \bZ^{(a)} *_{\bZ} \bZ^{(b\op)} \simeq \bZ^{(\Hom(b,a))}
  \quad\text{for $a,b$ in $A$.}\tag{52'}
\end{equation}
Comparing with the similar formula \eqref{eq:105.45}, this suggests
the generalization
\begin{equation}
  \label{eq:105.53}
  \Wh_A(F) *_{\bZ} \Wh_B(G) \simeq \bZ^{(F*G)},\tag{53}
\end{equation}
or with the exponential notation
\begin{equation}
  \label{eq:105.53prime}
  \bZ^{(F)} *_{\bZ} \bZ^{(G)} \simeq \bZ^{(F*G)},\tag{53'}
\end{equation}
valid\pspage{400} for $F$ in \Ahat{} and $G$ in \Bhat. As both members
of \eqref{eq:105.53} commute with small direct limits in each
variable, the formula \eqref{eq:105.53} follows from the particular
case \eqref{eq:105.52}, in view of the equivalence of categories
\eqref{eq:105.51}.

\medbreak

\noindent\textbf{Remarks.}\enspace\namedlabel{rem:105.1}{1})%
\enspace In order to appreciate the significance of the pairing
\eqref{eq:105.50}, we may forget altogether about the non-additive
categories \Ahat{} and \Bhat, and view \eqref{eq:105.50} as a
remarkable ``duality'' relationship between two additive
\scrU-categories, stable under small direct limits, say \scrP{} and
\scrQ, endowed with a ``pairing''
\[\scrP\times\scrQ\to\Ab\quad\text{in
  $\bHom_{!!}(\scrP,\scrQ;\Ab)$,}\]
giving rise to two functors which are \emph{equivalences}
\[ \scrQ\toequ \bHom_!(\scrP,\Ab),\quad \scrP\toequ
\bHom_!(\scrQ,\Ab),\]
identifying each of \scrP, \scrQ{} to the category of ``abelian
cosheaves'' on the other. In the particular case \eqref{eq:105.50},
$\scrP=\Ahatab$ and $\scrQ=\Bhatab$ with $B=A\op$, each of these
categories is even an ``abelian topos'' by which I mean an abelian
category \scrP{} stable under small filtering direct limits, with the
latter being \emph{exact}, and moreover \scrP{} admitting a small
generating subcategory. (These categories are sometimes called,
somewhat misleadingly, ``Grothendieck categories''. Of course, an
``abelian topos'' is by no means a category which is a topos, besides
being abelian!) There are many other examples of dual pairs of abelian
topoi. One evident generalization is by taking
\[\scrP=A_k\uphat,\quad\scrQ=B_k\uphat,\quad\text{where again
  $B=A\op$,}\]
where $k$ is any commutative ring, and $A_k\uphat\eqdef A_\kMod\uphat$
is the category of presheaves of $k$-modules on $A$, or equivalently,
of objects in \Ahat{} endowed with a structure of a $k$-module -- and
accordingly for the notation $B_k\uphat$. Indeed, the generalities
\ref{subsec:104.B} and \ref{subsec:104.C} in yesterday's notes about
abelian sheaves and cosheaves, as well as today's, could be developed
replacing throughout abelian group objects and additive categories by
$k$-module objects and $k$-linear categories. In case $k$ is not
supposed commutative, one still should get a duality pairing
\[A_{k\op}\uphat\times B_k\uphat\to\Ab,\]
where $k\op$ denotes the ring opposite to $k$ (i.e., a duality pairing
between presheaves on $A$ of right $k$-modules, and copresheaves on
$A$ of left $k$-modules), given in terms of $*_{\bZ}$ in
\eqref{eq:105.50} by the formula\pspage{401}
\[M *_{k} N = (M*_{\bZ} N)^\natural,\]
where in the right-hand side $P=M*_{\bZ}N$ is viewed as a
bi-$k$-module via the right and left $k$-module structures on $M$ and
$N$ respectively and bifunctoriality of $*_{\bZ}$, and where for any
bimodule $P$, we write
\[P^\natural \eqdef P /
\parbox[t]{0.65\textwidth}{sub-$\bZ$-module generated by elements of
  the type $s.x-x.s$ for $x$ in $P$ and $s$ in
  $k$.}\]
I confess I didn't do the checking that this does give rise indeed to
a duality pairing as desired. When $A=B=$ final category, then the
pairing above is just the pairing given by usual tensor product
\[(M,N)\mapsto M\otimes_k N\]
between right and left $k$-modules, which is immediately checked to be
dualizing indeed. More generally, posing
\[\scrP=(k\op\textup{-Mod}),\quad Q=\kMod,\]
it is immediately checked that for any additive category \scrM{}
stable under small direct limits, we get a canonical
equivalence\scrcomment{Clearly, the ``d'' in $k\subd$ is for
  ``dexter'' (on the right), and the ``s'' below is for ``sinister''
  (on the left).}
\[\bHom_!(\scrP,\scrM) \toequ (k\textup{-}\scrM), \quad
F\mapsto F(k\subd),\]
where $(k\textup{-}\scrM)$ denotes the category of objects $L$ of
\scrM{} endowed with a structure of a ``left $k$-module in \scrM'',
i.e., a ring homomorphism $k\to\End_\scrM(L)$, $k\subd$ denotes $k$
viewed as a right $k$-module, and $F(k\subd)$ is viewed as an object
of $(k\textup{-}\scrM)$ via the operations of $k$ on it coming from
left multiplication of $k$ upon $k\subd$. A quasi-inverse equivalence
is obtained by associating to an object $L$ in $(k\textup{-}\scrM)$
the functor
\[M\mapsto M\otimes_k L: (k\op\textup{-Mod})\to\scrM.\]
Dually, we get an equivalence (if \scrN{} stable under small inv.\
limits)
\[\bHom^!(\scrQ,\scrN) \toequ (k\textup{-}\scrN), \quad
F\mapsto F(k\subs),\]
where $k\subs$ denotes $k$ viewed as a left $k$-module, so that the
\emph{contravariant} functor $F$ transforms the endomorphisms of
$k\subs$ (obtained by \emph{right} operation of $k$ on $k\subs$ via
right multiplication) into a \emph{left} operation of $k$ on
$F(k\subs)$; the quasi-inverse is given by the familiar $\Hom_k$
operation, it associates to the left $k$-module $L$ in \scrM{} the
functor
\[M\mapsto \Hom_k(M,L): \kMod\to\scrN.\]

Comparing the two pairs of equivalences, we get the
equivalence\pspage{402}
\begin{equation}
  \label{eq:105.starbis}
  \bHom_!(scrP,\scrM) \equeq \bHom^!(\scrQ\op,\scrM),\tag{*}
\end{equation}
valid when \scrM{} is stable under both (small) direct and inverse
limits, and which should be viewed as an abelian analogon of the
equivalence \eqref{eq:105.39}.

\namedlabel{rem:105.2}{2})%
\enspace It is well-known that an abelian topos \scrP{} is
equivalent to a category \kMod, for a suitable ring $k$ (not
necessarily commutative) if{f} it admits an object $L$ which is
\namedlabel{it:105.rem2.a}{a)}\enspace generating, and
\namedlabel{it:105.rem2.b}{b)}\enspace ``\emph{ultraprojective}'',
i.e., the functor $X\mapsto\Hom(L,X)$ commutes to small direct
limits. The condition \ref{it:105.rem2.b} (for an object of an
\emph{abelian} category stable under small direct limits or
equivalently, under small direct sums) is equivalent with $L$ being:
\namedlabel{it:105.rem2.b1}{b\textsubscript{1})}\enspace projective,
and \namedlabel{it:105.rem2.b2}{b\textsubscript{2})}\enspace of
``\emph{finite presentation}'', i.e., the functor $X\mapsto\Hom(L,X)$
commutes with small filtering direct limits. This observation suggests
one common feature of all the examples of abelian duality pairings
considered so far, namely that the abelian topoi under consideration
in the pairing \emph{have a small set of ultraprojective
  generators}. I don't know if a structure theory of such categories
(which are the abelian analogons for topoi equivalent to topoi of the
special type \Ahat, with $A$ in \Cat) has been worked out yet. I
didn't do it at any rate -- but the natural thing to expect is that
these abelian topoi \scrP{} (which we may call ``\emph{algebraic}''
ones, just as an ordinary topos equivalent to one of the type \Ahat{}
may be called ``algebraic'', which equally means that the set of
ultraprojective objects in it is generating\ldots) are exactly those
equivalent to a category of the type
\[\Homadd(P\op,\Ab),\]
where $P$ is any \emph{small additive} category, and where $\Homadd$
denotes the category of additive functors from one additive category
to another. Instead of assuming $P$ small, we may as well take $P$
merely ``essentially small'', i.e., equivalent to a small category,
with the benefit that for a given \scrP, there is a canonical choice
of an additive category $P$ together with an equivalence
\[\scrP\toequ\Homadd(P\op,\Ab),\]
namely by taking
\[P =
\begin{tabular}[t]{@{}l@{}}
  full subcategory of \scrP{} made up with all
  ultraprojective\\objects in \scrP.
\end{tabular}\]
As we saw earlier, in case $\scrP=\Ahatab$, $P$ is nothing but the
abelian Karoubi envelope of the category $A$, namely the Karoubi
envelope of the additive category $\Add(A)$ (cf.\ sections
\ref{sec:93} and \ref{sec:99}). Another choice for $P$ in this case
would be just $\Add(A)$ itself, whose objects are more amenable to
computations.\pspage{403}

Associating to any small category $P$ the
algebraic abelian topos $\bHomadd(P\op,\Ab)$ should be viewed of
course as the abelian analogon of $A\mapsto\Ahat$, associating to a
small category $A$ the corresponding algebraic topos. It merits a
notation of its own, say
\[P\supamp\eqdef \bHomadd(P\op,\Ab),\]
and as in the non-additive case, we get a canonical inclusion functor
\[\varepsilon_P:P\to P\supamp\]
which is additive. (Its composition with the canonical functor
$P\supamp\to P\uphat$ is the canonical inclusion
previously denoted by $\varepsilon_P$ too from $P$ to $P\uphat$.) Next
thing we'll expect, in analogy to prop.\ \ref{prop:105.3}, is that for
any additive category \scrM{} stable under small direct limits, the
following canonical functor is an equivalence of categories:
\begin{equation}
  \label{eq:105.starstar}
  \bHom_!(P\supamp,\scrM) \toequ \bHomadd(P,\scrM), \quad
  F\mapsto F\circ\varepsilon_P.\tag{**}
\end{equation}
The proof, via construction of a quasi-inverse functor, should be
about the same as for prop.\ \ref{prop:105.3}, which should go through
once we get the abelian analogon of the well-known fact in \Ahat, that
any object $F$ in \Ahat{} can be recovered as a direct limit in
\Ahat{} of objects of $A$, according to $A_{/F}$ as an indexing
category -- which makes us expect that we get too:
\[F \fromsim \varinjlim_{\text{$a$ in $P\supamp_{/F}$}} a\quad
(\text{direct limit in $P\supamp$}).\]
From \eqref{eq:105.starstar} we get as in cor.\ \ref{cor:105.prop3.1},
passing to the dual categories, the dual equivalence
\[\bHom^!({P\supamp}\op,\scrN) \toequ \bHomadd(P\op,\scrN),\]
valid if \scrN{} is an additive category stable under small inverse
limits. Hence, if \scrM{} is additive and stable under both types of
limits, the equivalence
\begin{equation}
  \label{eq:105.starstarstar}
  \bHom_!(P\supamp,\scrM)\toequ \bHom^!({Q\supamp}\op,\scrM), \quad
  \text{with $Q=P\op$,}\tag{***}
\end{equation}
between abelian cosheaves on $P\supamp$ and abelian sheaves in
$Q\supamp$, with values in the same additive category (in analogy to
\eqref{eq:105.39}). In the particular case $\scrM=\Ab$, this then
gives rise to the equivalence
\[Q\supamp \equeq \bHom_!(P\supamp,\Ab)\]
and to the corresponding pairing
\[P\supamp\times Q\supamp\to\Ab\]
which is a duality, namely induces an equivalence between each abelian
topos\pspage{404} $P\supamp$, $Q\supamp$ and the category of abelian
cosheaves on the other.

We may call an abelian topos ``\emph{reflexive}'' if it can be
inserted in a pair $(\scrP,\scrQ)$ of dually paired abelian topoi --
where \scrQ, or the ``dual'' of \scrP, is defined up to equivalence in
terms of \scrP{} as $\bHom_!(\scrP,\Ab)$, the category of cosheaves on
\scrP{} with values in \Ab. Thus, it seems that a sufficient condition
for reflexivity is ``algebraicity'' of \scrP, namely the existence of
a small generating family made up with ultraprojective
objects. (NB\enspace In the non-abelian case, it is well-known that a
topos \scrA{} is ``algebraic'', i.e., equivalent to a topos \Ahat,
if{f} it admits such a generating family -- and as we saw in
\ref{case:105.1}, as a consequence of \eqref{eq:105.39} for
$\scrM=\Sets$, such a topos is indeed ``reflexive''.) I wouldn't be
too surprised if this sufficient condition for reflexivity turned out
to be necessary too, at any rate if we want a property stronger still
than reflexivity, namely validity of a duality equivalence
\eqref{eq:105.starstarstar} for sheaves and cosheaves with values in
an arbitrary additive \scrU-category stable under small direct and
inverse limits, satisfying (for varying \scrM) suitable compatibility
assumptions.

\namedlabel{rem:105.3}{3})%
\enspace With respect to this duality equivalence
\eqref{eq:105.starstarstar}, I am a little unhappy still, as I do not
see how to get (for a general dual pair \scrP, \scrQ{} of abelian
topoi) a functor i one direction or the other between the two
categories
\[\bHom_!(\scrP,\scrM),\quad\bHom^!(\scrQ\op,\scrM),\]
in terms of just the duality pairing. The same perplexity holds in the
non-abelian case. This is one of the reasons that make me feel that I
haven't yet a thorough understanding of the duality formalism I am
developing here, except in the ``algebraic'' case (granting for the
latter that the tentative theory just outlined for algebraic
\emph{abelian} topoi is indeed correct).

\namedlabel{rem:105.4}{4})%
\enspace To finish with the comments on the (pre-homological)
duality formalism for algebraic topoi and algebraic abelian topoi, I
still would like to add that the category $\scrE(A)$ (union of \Ahat{}
and $\Bhat\op$, with $B=A\op$) introduced in \ref{case:105.1} (cf.\
remark on page \ref{p:398}) admits also an abelian analogon. In the
non-abelian case still, the simplest way to construct the category
$\scrE(A)$, is via an equivalent category canonically embedded in the
category $\bHom(\Bhat,\Sets)$ as a strictly full subcategory
(``strict'' referring to the fact that with any object it contains all
isomorphic ones), namely the union $\overline{\scrE}(A)$ of the
(strictly full) subcategories\pspage{405} $\bHom_!(\Bhat,\Sets)$
(equivalent to \Ahat) and the subcategory of representable functors
(equivalent to $\Bhat\op$). The intersection of these two categories
contains of course $A$ (embedded in $\bHom(\Bhat,\Sets)$ by
associating to $a$ in $A$ the functor $G\mapsto G(a)$ from
$\Bhat=\bHom(A,\Sets)$ to \Sets), but in general need not be quite
equivalent to $A$ -- it turns out to be the ``Karoubi envelope'' of
$A$, obtained by adjoining to $A$ formally images (or equivalently,
coimages) of projectors in $A$. The more immediate interpretation of
this intersection, is that it is equivalent to the dual category of
the category of ultraprojective objects in \Bhat{} (and the latter can
be viewed as $\Kar(B)$, but formation of the Karoubi envelope up to
equivalence commutes to taking dual categories\ldots). All these
constructions immediately extend to the abelian set-up, starting with
a small \emph{additive} category $P$, instead of $A$.

After this endless procession of remarks, which are really digressions
for what we're after (namely abelianization and duality in the context
of small categories as homotopy models), it is time to resume our main
line of thought in this subsection, namely looking at interesting
particular cases for the general duality relation \eqref{eq:105.39}.

\textbf{Case \namedlabel{case:105.3}{3\textsuperscript{\b o})}.}%
\enspace This is the case when \scrM{} is an \emph{additive} category,
stable under both types of small limits. If we assume moreover that
\scrM{} and $\scrM\op$ are both pseudotopoi, using the equivalences
\eqref{eq:104.12} and \eqref{eq:104.13} (p.\ \ref{p:387}),
\eqref{eq:105.39} may be interpreted as an equivalence
\begin{equation}
  \label{eq:105.54}
  \bHom_!(\Ahatab,\scrM)\equeq \bHom^!(\Bhatab,\scrM),\tag{54}
\end{equation}
interpreting \scrM-valued cosheaves on the abelian topos \Ahatab{} in
terms of \scrM-valued sheaves on the dual abelian topos \Bhatab, as
anticipated in a more general situation in the remark above (cf.\
formula \eqref{eq:105.starstarstar} on page \ref{p:403}). There,
however, the assumption that \scrM{} and/or $\scrM\op$ should be
pseudotopoi didn't seem to come in at all, so we expect this condition
to be irrelevant indeed. This will of course follow, if the same holds
for \eqref{eq:104.12} (hence by duality for \eqref{eq:104.13}), namely
that the canonical functor
\begin{equation}
  \label{eq:105.55}
  \bHom_!(\Ahatab,\scrM)\toequ\bHom_!(\Ahat,\scrM), \quad
  F\mapsto F\circ\Wh_A,\tag{55}
\end{equation}
is an equivalence, under the only assumption that the additive
category \scrM{} is stable under small direct limits (without assuming
that \scrM{} be a pseudotopos). The line of thought of the remark
\ref{rem:105.2} above suggests a way for proving this, via an
equivalence\pspage{406}
\begin{equation}
  \label{eq:105.56}
  \bHom_!(\Ahatab,\scrM)\toequ\bHomadd(\Add(A),\scrM), \quad
  F\mapsto F\circ j_A,\tag{56}
\end{equation}
where
\[j_A: \Add(A)\to\Ahatab\]
is the canonical inclusion functor (cf.\ section \ref{sec:97}). This,
and the dual equivalence (deduced from \eqref{eq:105.55}, taking
$\scrN=\scrM\op$)
\begin{equation}
  \label{eq:105.57}
  \bHom^!(\Ahatab\op,\scrN)
  \toequ\bHomadd(\Add(A)\op,\scrN),\tag{57}
\end{equation}
valid for any additive category stable under small inverse limits,
will immediately imply an equivalence \eqref{eq:105.54} by a direct
argument as in the remark above, without passing through the
non-abelian case \eqref{eq:105.39}. At any rate, \eqref{eq:105.56}
implies that \eqref{eq:105.55} is an equivalence, as is seem by
looking at the commutative diagram
\begin{equation}
  \label{eq:105.58}
  \begin{tabular}{@{}c@{}}
    \begin{tikzcd}[baseline=(O.base),row sep=small,column sep=-1em]
      & \bHom_!(\Ahat,\scrM)\ar[dr] & \\
      \bHom_!(\Ahatab,\scrM)\ar[dr]\ar[ur] & &
      \bHom(A,\scrM) \\
      & |[alias=O]| \bHomadd(\Add(A),\scrM)\ar[ur] &
    \end{tikzcd},
  \end{tabular}\tag{58}
\end{equation}
where the two right-hand arrows are equivalences, which implies that
one of the two left-hand arrows is an equivalence if{f} the other is.

Thus, for getting \eqref{eq:105.54} and \eqref{eq:105.55} without
extraneous assumptions on \scrM, we are left with proving
\eqref{eq:105.56}. Now, writing
\[\scrP=\Ahatab,\quad P=\Add(A),\]
we do have indeed an equivalence
\[\scrP\equeq P\supamp \eqdef \bHomadd(P\op,\Ab),\]
as seen from \eqref{eq:105.57} taking $\scrN=\Ab$ (which satisfies the
extra assumptions). So we may as well prove \eqref{eq:105.56} in the
more general case when $P$ is any small additive category and \scrP{}
is defined as $P\supamp$, namely prove the equivalence
\eqref{eq:105.starstar} of page \ref{p:403} above, as I don't expect
the particular case at hand here to be any simpler. The suggestion for
a proof there there seems convincing, I guess I should check it works
during some in-between scratchwork\ldots

\bigbreak

\presectionfill\ondate{14.8.}\pspage{407} and \ondate{15.8}\par

\hangsection[Review (3): A formulaire for the basic integration and
\dots]{Review \texorpdfstring{\textup{(3)}}{(3)}: A formulaire for the
  basic integration and cointegration operations \texorpdfstring{$*$
    and $\Hom$}{* and Hom}.}\label{sec:106}%
\phantomsection\addcontentsline{toc}{subsection}{\numberline {E)}A
  formulaire around the basic operations \texorpdfstring{$*$ and
    $\Hom$}{* and Hom}.}%
\textbf{\namedlabel{subsec:106.E}{E)}\enspace A formulaire around the
  basic operations $*$ and $\Hom$.}\enspace I would like to dwell a
little more still on the duality formalism weaving around formula
\eqref{eq:105.39} (p.\ \ref{p:395}), stating that for two given small
categories $A$ and $B$ dual to each other
\[B=A\op, \quad A=B\op,\]
and any \scrU-category \scrM{} stable under both types of small
limits, \scrM-valued \emph{cosheaves} on \Ahat{} may be interpreted as
\scrM-valued \emph{sheaves} on the dual topos \Bhat. This
identification preserves variance (i.e., \eqref{eq:105.39} is an
\emph{equivalence} of categories, not an antiequivalence:
\begin{equation}
  \label{eq:106.59}
  (\Ahat)^\scrM = \bHom_!(\Ahat,\scrM) \toequ
  \bHom^!((\Bhat)\op,\scrM) = \BhatM\quad\text{).}\tag{59}
\end{equation}
It should not be confused with the tautological interpretation of
\scrM-valued cosheaves on \Ahat{} as $\scrM\op$-valued sheaves on the
\emph{same} topos, an identification \emph{reversing variance}, as
expressed by the canonical antiequivalence between the corresponding
categories -- an anti-isomorphism even (reflecting its tautological
nature):
\begin{equation}
  \label{eq:106.60}
  \bHom_!(\Ahat,\scrM)\op \tosim \bHom^!((\Ahat)\op,\scrM\op),
  \text{ i.e., }((\Ahat)^\scrM)\op \tosim \AhatM.\tag{60}
\end{equation}
In the latter formula, the basic topos \Ahat{} remains the same in
both sides, it is the category of values that changes from \scrM{} to
the dual one $\scrM\op$, whereas in formula \eqref{eq:106.59} =
\eqref{eq:105.39}, it is the opposite. In terms of the tautological
formula \eqref{eq:106.60} (a particular case of formula
\eqref{eq:104.6} p.\ \ref{p:385}), the not-so-tautological formula
relating cosheaves and sheaves can be reformulated as a formula in
terms of sheaves only (due to our preference for sheaves rather than
cosheaves\ldots):
\begin{equation}
  \label{eq:106.61}
  \bHom^!(\Ahat,\scrM)\op \equeq \bHom^!(\Bhat,\scrM\op),\quad
  \text{i.e.,}\quad (\AhatM)\op \equeq B\uphat_{\scrM\op},\tag{61} 
\end{equation}
namely \scrM-valued sheaves on the topos \Ahat{} can be interpreted as
sheaves on the dual topos with values in the dual category $\scrM\op$,
this interpretation \emph{reversing variances}. In the homology and
cohomology formalism which is to follow, due to habits of long
standing, I prefer systematically to take as coefficients
\emph{sheaves} rather than cosheaves -- hence rule out cosheaves in
favor of sheaves via \eqref{eq:106.60}. From this point of view the
relevant basic duality statement is \eqref{eq:106.61} rather than
\eqref{eq:106.59}.

On the \emph{cosheaves} side, yesterday's diagram \eqref{eq:105.58} of
equivalences gives a fourfold description of cosheaves on the topos
\Ahat{} with values in an \emph{abelian} category \scrM{} stable under
small direct limits. We could still enlarge this diagram, by including
in it a fifth category equivalent\pspage{408} to the four others,
namely
\[\bHomaddinf(\Addinf(A),\scrM),\]
the category of infinitely additive functors from the infinitely
additive envelope of $A$ into $M$ (cf.\ section \ref{sec:99} p.\
\ref{p:366} for description of the category $\Addinf(A)$). Rather than
writing down the larger diagram here, I'll write down the dual
enlarged one, for the dual topos \Bhat{} and for various expressions
of the category of \emph{sheaves} on this topos, with values in a
category \scrM{} stable this time under small inverse limits:
\begin{equation}
  \label{eq:106.62}
  \begin{tabular}{@{}c@{}}
    \begin{tikzpicture}[commutative diagrams/every diagram,baseline=(E.base)]
      \node (A) at (0,1.2cm) {$\bHom^!((\Bhat)\op,\scrM)$};
      \node (B) at (-3.2cm,0) {$\bHom^!((\Bhatab)\op,\scrM)$};
      \node (C) at (2.2cm,0) {$\bHom(B\op,\scrM)$};
      \node (D) at (-3.2cm,-1.2cm) {$\bHommultinf(\Addinf(B)\op,\scrM)$};
      \node (E) at (1.5cm,-1.2cm) {$\bHomadd(\Add(B)\op,\scrM)$};

      \path[commutative diagrams/.cd, every arrow, every label]
      (A) edge node {$\equ$} (C)
      (B) edge node {$\equ$} (A)
      (B) edge node[swap] {$\equ$} (D)
      (D) edge node {$\equ$} (E)
      (E) edge node {$\equ$} (C);
    \end{tikzpicture},
  \end{tabular}\tag{62}
\end{equation}
where $\bHommultinf$ denotes the category of ``infinitely
multiplicative'' functors from one additive category stable under
infinite products to another. Recalling for the extreme right term of
\eqref{eq:106.62} that $B\op=A$, we see that this term is
\emph{identical} to the corresponding term in the diagram (even the
enlarged one) \eqref{eq:105.58} -- hence, if \scrM{} is stable under
both types of limits, the ten categories occurring altogether in the
two diagrams are mutually equivalent (as a matter of fact, there are
nine only which are mutually different). It may be noted that there is
still another pair of corresponding terms in the two diagrams for
which the equivalence between them may be viewed as tautological,
namely
\[\bHomadd(\Add(A),\scrM) \equeq \bHomadd(\Add(B)\op,\scrM),\]
due to the tautological equivalence of categories
\[\Add(B)\op \equeq \Add(B\op) = \Add(A).\]

As emphasized in yesterday's notes, the canonical pairing (deduced
from \eqref{eq:106.59} by taking $\scrM=\Ab$)
\begin{equation}
  \label{eq:106.63}
  \Ahatab\times\Bhatab\to\Ab,\quad
  (F,G)\mapsto F*_{\bZ} G,\tag{63}
\end{equation}
deserves special attention, giving rise to an equivalence between each
of the mutually dual abelian topoi \Ahatab, \Bhatab{} with the
category of abelian cosheaves on the other
\begin{equation}
  \label{eq:106.64}
  \Bhatab \toequ \bHom_!(\Ahatab,\Ab), \quad
  \Ahatab \toequ \bHom_!(\Bhatab,\Ab)\tag{64}
\end{equation}
(cf.\ \eqref{eq:105.49} and \eqref{eq:105.49prime} page
\ref{p:399}). It should be kept in mind that besides this duality
pairing between two abelian topoi, there is important
extra\pspage{409} structure in this abelianized duality context,
embodied by the \emph{tensor product structure} on both abelian topoi
\Ahatab{} and \Bhatab, as contemplated in section \ref{sec:104}
\ref{subsec:104.C}, in a somewhat more general context. Corresponding
to this extra structure on \Ahatab{} say, we saw that this category of
abelian sheaves ``operates'' covariantly (by an operation denoted by
$\oast_{\bZ}$ or simply $\oast$) on any category of \scrM-valued
cosheaves on \Ahatab, and contravariantly (by an operation denoted by
$\bHom_{\bZ}$ or simply $\bHom$, if no confusion may arise) on any
category of \scrN-valued sheaves on \Ahatab, where \scrM, \scrN{} are
additive categories, stable under small direct resp.\ inverse
limits. The latter operation (cf.\ page \ref{p:392})
\begin{equation}
  \label{eq:106.65}
  (L,K) \mapsto \bHom_{\bZ}(L,K) : (\Ahatab)\op \times \AhatN \to
  \AhatN,\tag{65}
\end{equation}
involving sheaves, will be used in the sequel ``tel
quel'',\scrcomment{``tel quel'' = ``as is''} its
definition is of a tautological character, independent of duality. As
for the former operation involving cosheaves, we may view it via the
duality relation \eqref{eq:106.59} as an operation of \Ahatab{} on
\scrM-valued \emph{sheaves} on the dual abelian topos \Bhatab, and
this operation will be denoted by the same symbol $\oast_{\bZ}$:
\begin{equation}
  \label{eq:106.66}
  (L,M') \mapsto L\oast_{\bZ} M' : \Ahatab \times \BhatM
  \to \BhatM.\tag{66} 
\end{equation}
Replacing $A$ by $B$ in \eqref{eq:106.66}, we get an operation of
\Bhatab{} upon \AhatM,
\begin{equation}
  \label{eq:106.67}
  (M,L') \mapsto M \oast_{\bZ} L' :
  \AhatM\times\Bhatab\to\AhatM.\tag{67}
\end{equation}
Whenever convenient, we'll allow ourselves to write $L'\oast M$
instead of $M\oast L'$ (which doesn't seem to lead to any trouble),
and will henceforth (unless special need should arise) drop the
subscripts $\bZ$.

Thus, for a given category $\AhatM$ of \scrM-valued sheaves on
\Ahatab, there is a twofold operation on this category, namely
\Ahatab{} itself operates (the operation defined by $L$ in \Ahatab{}
depending contravariantly on $L$) as well as the dual abelian topos
\Bhatab{} (the operation defined by $L'$ in \Bhatab{} depending
covariantly on $L'$), \scrM{} being any additive category stable under
small inverse and direct limits (in order to ensure existence of both
types of operations). I would like to dwell a little more on this
twofold structure, as I don't feel to have understood it thoroughly
yet. It is this second operation mainly which hasn't become really
familiar yet, still less its relationship to the first, more
familiar\pspage{410} operation is understood. I'll have to play around
a little more with it for being really at ease. It's worth the while,
as the $\bHom$ and corresponding $\Hom$ operation is the key operation
for expressing cohomology of \Ahat{} (with coefficients in $K$), where the
$\oast$ and corresponding $*$ operation is the key for expressing
homology of \Ahat{} (with coefficients in $M$, where $K$ and $M$ are
the sheaves occurring in \eqref{eq:106.65} and \eqref{eq:106.67}
respectively).

A typical special case is the one when $A$ is the final category,
hence $B=A$ and $\AhatM\simeq \scrM$, in which case
\eqref{eq:106.65} and \eqref{eq:106.67} are the two familiar exterior
operations of \Ab{} on any additive category stable under the two
types of small limits
\[ (L,X)\mapsto \bHom_{\bZ}(L,X): \Ab\op\times\scrM\to\scrM,\]
and
\[ (X,L)\mapsto X\otimes_{\bZ} L:\Ab\times\scrM\to\scrM,\]
each one of these two operations being deducible from the other by the
usual device of replacing \scrM{} by the dual category
$\scrM\op$. When $\scrM=\Ab$, these are just the usual internal
$\Hom=\bHom$ and tensor product operations. This very particular case
shows at once that we shouldn't expect in general the operations
$\bHom_{\bZ}(L,{-})$ of \Ahatab{} and ${-} * L'$ of \Bhatab{} upon
\AhatM{} to commute up to isomorphism -- we shouldn't expect, for a
given $L$ in \Ahatab{} or a given $L'$ in \Bhatab{} the commutation
relation to hold, except when this object is ``projective of finite
presentation'', i.e., is a direct factor of an object of $\Add(A)$
resp.\ of $\Add(B)\equeq\Add(A)\op$. Another fact becoming evident by
this particular case, is that whereas it is true that in the
equivalence \eqref{eq:106.64}
\[\Bhatab \equeq \bHom_!(\Ahatab,\Ab)\]
a projective object $L'$ in \Bhatab{} gives rise to a functor
$\Ahatab\to\Ab$ which is \emph{exact} (besides commuting to small
direct limits) -- and even to a functor commuting to small inverse
limits if $L'$ is ultraprojective, i.e., projective and of finite
presentation -- the converse to this (as contemplated on page
\ref{p:381}) does \emph{not} hold true. Indeed, in the particular case
$A=\Simplex_0$, when $L'$ is just an object in \Ab, the exactness
property envisioned, i.e, exactness of the functor $M\mapsto
M\otimes_{\bZ} L'$ from \Ab{} to itself, only means that $L'$ is a
\emph{flat} $\bZ$-module (i.e., torsion-free), which does \emph{not}
imply that it is projective (i.e., free). 

The\pspage{411} feeling I had earlier today, that the familiar looking
operation \eqref{eq:106.65} $\bHom(L,K)$ of abelian sheaves on \Ahat{}
upon \scrM-valued ones (\scrM{} an additive category stable under
small inverse limits) was well-understood, whereas the less familiar
one $M*L'$ in \eqref{eq:106.67} was not, turns out to be mistaken. In
computational terms, and writing $F$ for $K$ in \eqref{eq:106.65} and
$M$ in \eqref{eq:106.67}, the three basic data
\[\text{$F$ in \AhatM,} \quad
\text{$L$ in \Ahatab,} \quad
\text{$L'$ in \Bhatab}\]
should surely be interpreted as just functor
\begin{equation}
  \label{eq:106.68}
  F:A\op\to\scrM, \quad
  L:A\op\to\Ab,\quad
  L':B\op=A\to\Ab,\tag{68}
\end{equation}
and the practical question of ``computing'' $\bHom(L,F)$ or $F\oast
L'$ thus amounts to describing, directly in terms of these data, the
corresponding objects in \AhatM{} as again a functor
\[\bHom(L,F)\quad\text{or}\quad F\oast L':A\op\to\scrM.\]
It would seem that neither of the two can be expressed in simplistic
computational terms, via the data \eqref{eq:106.68}. I feel I have to
come to terms with this fact and get as close as I can to an explicit
expression of both. The point I want to make first, is that this
question of expressing $\bHom(L,F)$, or of expressing the operation
$F\oast L'$, is essentially the same, via replacement of \scrM{} by
$\scrM\op$, and mere interchange of $A$ and $B$. More accurately,
passing from $F$ to the corresponding functor $F\op$ between the dual
categories, we may view the data \eqref{eq:106.68} as being functors
\begin{multline}
  \label{eq:106.68prime}
  F\op:B\op\to\scrN,\quad
  L':B\op\to\Ab,\quad
  L:A\op\to\Ab\\
  \text{ (where $\scrN=\scrM\op$),}\tag{68'}
\end{multline}
i.e., a set of data like \eqref{eq:106.68}, with $(A,\scrM)$ replaced
by the dual pair $(B,\scrN)$. This being understood, we have the
tautological isomorphisms
\begin{equation}
  \label{eq:106.69}
  \begin{cases}
    \bHom(L,F)\op \simeq F\op \oast L & \\
    (F\oast L')\op \simeq \bHom(L',F\op) &\text{,}
  \end{cases}\tag{69}
\end{equation}
where the first members involve operations relative to the pair
$(A,\scrM)$, the second members operations relative to the dual pair
$(B,\scrN)$. This makes very clear, it seems to me, that the
operations $\bHom$ \eqref{eq:106.65} and $\oast$ \eqref{eq:106.67} may
be viewed as the same type of operation, simply viewed with two
different pairs of spectacles -- one being $(A,\scrM)$, the other the
dual pair $(B,\scrN)$. Thus, if we got a good understanding of one of
the two operations, embodied by a comprehensive
formulaire\scrcomment{I'm leaving in ``formulaire'' (form), even
  though ``formula'' seems to work better\ldots} for it, by just
dualizing we should get just as good a formulaire and corresponding
comprehension for the dual operation.

Now, it is clear indeed that it is the operations
\eqref{eq:106.65}\pspage{412} which is closer to my experience, it
makes sense however, independently of any duality statements, in the
vastly more general context of topoi \scrA{} (or even only pseudotopoi
satisfying some mild extra conditions, cf.\ section \ref{sec:104}
\ref{subsec:104.C}) instead of just \Ahat, provided we make on \scrM{}
the mild extra assumption of being a pseudotopos, needed in this more
general context in order to ensure equivalence between the category
$\scrA_\scrM$ of \scrM-valued sheaves, and the category
$\bHom^!({\scrA\subab}\op,\scrM)$ (compare \eqref{eq:104.13} p.\
\ref{p:388}). What I should do then is, first to write down a basic
``formulaire'' for the $\bHom$ operation in this general and familiar
context, then see how it can be used for clarifying the computational
puzzle raised on the previous page, in the case of the $\bHom$
operation, and finally dualize the formulaire and computational
insight, for getting a hold on the dual operation $\oast$.

We'll need too the $\Hom$ operation (non-bold-face) of
\eqref{eq:104.31} p.\ \ref{p:392}
\begin{equation}
  \label{eq:106.70}
  (L,F)\mapsto \Hom(L,F) : {\scrA\subab}\op \times \scrA_\scrM \to
  \scrM,\tag{70} 
\end{equation}
with values in \scrM, not $\scrA_\scrM$, where $\Hom(L,F)$ denotes the
value on $L$ of the functor
\[\widetilde F : {\scrA\subab}\op\to\scrM\]
defined by $F$, $F$ being viewed for the time being as an object in
$\bHom^!(\scrA\op,\scrM)$:
\begin{equation}
  \label{eq:106.71}
  F:\scrA\op\to\scrM.\tag{71}
\end{equation}
(In case $\scrA=\Ahat$, we get the description \eqref{eq:106.68} of
$F$ by taking the restriction of \eqref{eq:106.71} to the subcategory
$A$ of \Ahat.)

In the following formulaire, $L$, $L'$ are objects in $\scrA\subab$,
$F$ is an object in $\scrA_\scrM$, where \scrM{} is an additive
category stable under small inverse limits, which moreover is assumed
to be a pseudotopos (i.e., admits a small set of objects generating
with respect to monomorphisms, and is stable under small direct
limits) in case the topos \scrA{} is not equivalent to a category
\Ahat. We denote by
\[X\mapsto\bZ^{(X)}, \quad \scrA\to\scrA\subab\]
the abelianization functor, which in the case $\scrA=\Ahat$ is just
``componentwise abelianization'', i.e.,
\[\bZ^{(X)}(a)=\bZ^{(X(a))}\quad\text{for $a$ in $A$.}\]
We recall that the constant sheaf $\bZ_\scrA$ on \scrA{} with value
$\bZ$ can be also described as
\begin{equation}
  \label{eq:106.72}
  \bZ_\scrA=\bZ^{(e)},\tag{72}
\end{equation}
where $e$ is the final object of \scrA, and that the \emph{sections}
functor on \scrA\pspage{413} is defined as
\begin{equation}
  \label{eq:106.73}
  \Gamma_\scrA(F) \eqdef F(e);\tag{73}
\end{equation}
in case $\scrA=\Ahat$, this can equally be interpreted as the inverse
limit functor for the functor $A\op\to\scrM$ defined by $F$:
\begin{equation}
  \label{eq:106.74}
  \Gamma_\scrA(F) \simeq \varprojlim_{A\op} F(a).\tag{74}
\end{equation}
We are now ready to give a basic formulaire for the operations $\bHom$
and $\Hom$, and their relations to the abelianization functor and to
the sections functor (i.e., to inverse limits, in case $\scrA=\Ahat$).
\begin{equation}
  \label{eq:106.75}
  \left\{%
    \renewcommand*{\arraystretch}{1.1}%
    \begin{array}{@{}r@{}rl@{}}
      \left\{\rule{0pt}{3ex}\right. &
      \begin{tabular}{@{}r@{}}
        \textcircled{a} \\ 
        a')
      \end{tabular} &
      \begin{array}{@{}l@{}}
        \Hom(\bZ^{(X)},F)\simeq F(X) \\
        \bHom(\bZ^{(X)},F)\simeq (Y\mapsto F(X\times Y):\scrA\op\to\scrM)
      \end{array} \\
      \left\{\rule{0pt}{3ex}\right. &
      \begin{tabular}{@{}r@{}}
        \textcircled{b} \\
        b')
      \end{tabular} &
      \begin{array}{@{}l@{}}
        \Hom(L',\bHom(L,F)) \simeq \Hom(L'\otimes L,F) \\
        \bHom(L',\bHom(L,F)) \simeq \bHom(L'\otimes L,F)
      \end{array} \\
      \left\{\rule{0pt}{3ex}\right. &
      \begin{tabular}{@{}r@{}}
        c) \\
        c')
      \end{tabular} &
      \begin{array}{@{}l@{}}
        \Hom(\bZ_\scrA,F)\simeq\Gamma_\scrA(F) \\
        \bHom(\bZ_\scrA,F)\simeq F
      \end{array} \\
      & \text{d)} & \Hom(L,F)\simeq\Gamma_\scrA\bHom(L,F) \\
      \left\{\rule{0pt}{6ex}\right. &
      \begin{tabular}{@{}r@{}}
        \textcircled{e} \\
        $\phantom{()}$ \\
        e') \\
        $\phantom{()}$
      \end{tabular} &
      \begin{tabular}{@{}l@{}}
        The functor $L\mapsto\Hom(L,F):{\scrA\subab}\op\to\scrM$ \\
        commutes to small inverse limits \\
        Similar statement as e) for \\
        $L\mapsto\bHom(L,F):{\scrA\subab}\op\to\scrA_\scrM$.
      \end{tabular}
    \end{array}
  \right.\tag{75}
\end{equation}
\textbf{Comments on the formulaire \eqref{eq:106.75}.}\enspace I have
limited myself to chose canonical isomorphisms
(\hyperref[eq:106.75]{a)} to \hyperref[eq:106.75]{d)}) and exactness
properties (\hyperref[eq:106.75]{e)} and \hyperref[eq:106.75]{e')})
which seem to me the most relevant for what follows. Other exactness
and variance properties are commutation of the functors
\[F\mapsto\Hom(L,F):\scrA_\scrM\to\scrM\quad\text{and}\quad
F\mapsto\bHom(L,F):\scrA_\scrM\to\scrA_\scrM\]
to small inverse limits, and compatibility of formation of $\Hom(L,F)$
and $\bHom(L,F)$ with functors
\[u:\scrM\to\scrM'\]
commuting to small inverse limits. As for formulæ for varying topos
\scrA, corresponding to a morphism of topoi, we'll come back upon this
in a later section, in relation with the homology and cohomology
invariants of maps in \Cat. Also, I am completely disregarding here
compatibilities between canonical isomorphisms (surely the reader
won't complain about this). All this as far as omissions are
concerned.

As\pspage{414} for the formulas included in \eqref{eq:106.75}, the
three basic ones, including all others in a more or less formal way,
are the circled ones \hyperref[eq:106.75]{a)},
\hyperref[eq:106.75]{e)} and \hyperref[eq:106.75]{b)}. The properties
\hyperref[eq:106.75]{a)} and \hyperref[eq:106.75]{e)} jointly can be
viewed as the characterization up to canonical isomorphism, for fixed
$F$, of the operation $\Hom(L,F)$, i.e., of the functor
\[\widetilde F: L\mapsto\Hom(L,F): {\scrA\subab}\op\to\scrM,\]
factoring the functor
\[F:\scrA\op\to\scrM\]
via the abelianization functor $\Wh_A:X\mapsto \bZ^{(X)}$. In terms of
a), the formula b) can be viewed as essentially the definition of
$\bHom(L,F)$ via $\Hom({-},F)$, more specifically we get
\begin{equation}
  \label{eq:106.76}
  \bHom(L,F)(X) \simeq \Hom(\bZ^{(X)},\bHom(L,F)) \simeq
  \Hom(\bZ^{(X)}\otimes L, F).\tag{76}
\end{equation}
Taking $L=\bZ^{(X)}$ and using
\[\bZ^{(X)}\otimes \bZ^{(Y)} \simeq\bZ^{(X\times Y)}\]
(\eqref{eq:104.17} page \ref{p:389}), \eqref{eq:106.76} gives
\hyperref[eq:106.75]{a')}, whereas \hyperref[eq:106.75]{b')} follows
via \eqref{eq:106.76} applied to both members, from associativity of the
operation $\otimes$. Formula \hyperref[eq:106.75]{c)} is the
particular case of \hyperref[eq:106.75]{a)} for $X=e$, in the same way
\hyperref[eq:106.75]{c')} follows from
\hyperref[eq:106.75]{a')}. Formula \hyperref[eq:106.75]{d)} follows
from \hyperref[eq:106.75]{c)}, \hyperref[eq:106.75]{b)} and the
relation
\[\bZ_\scrA\otimes L\simeq L.\]
The exactness property \hyperref[eq:106.75]{e')} is equivalent to the
similar exactness statement for the functors
\[L\mapsto\bHom(L,F)(X)\simeq\Hom(\bZ^{(X)}\otimes L,F),\]
for $X$ in \scrA, and thus reduces to \hyperref[eq:106.75]{e)} with
$L$ replaced by $\bZ^{(X)}\otimes L$.

I would like now to come back to the question of ``computation'' of
$\Hom(L,F)$ and $\bHom(L,F)$. We may for this end assume \scrA{} to be
described by a \emph{site} $A$ -- which, in case the ``topology'' on
$A$ defining the site structure is the chaotic one, brings us back to
the particular case $\scrA=\Ahat$ we are mainly interested in at
present. Accordingly, we'll consider the objects $F$ in $\scrA_\scrM$
as being functors
\begin{equation}
  \label{eq:106.77}
  F:A\op\to\scrM\tag{77}
\end{equation}
satisfying the standard exactness properties for sheaves (with respect
to the given site structure on $A$). In terms of \eqref{eq:106.71},
this is just the composition of the functor \eqref{eq:106.71} with the
canonical functor
\begin{equation}
  \label{eq:106.78}
  A\to\scrA=A^{\sim},\tag{78}
\end{equation}
associating to an object $a$ in $A$ the presheaf represented by it, in
the\pspage{415} most common case when this presheaf is a sheaf for any
choice of $a$, otherwise we take the sheaf associated to it. In the
first case (which we may reduce to if we prefer, by suitable choice of
the site $A$ for given topos \scrA) the functor \eqref{eq:106.78} is
fully faithful and moreover and embedding, therefore, we'll identify
an object $a$ in $A$ with the corresponding object in \scrA. Thus, the
description \eqref{eq:106.76} of $\bHom$ in terms of $\Hom$ may be
interpreted, from this point of view, as a formula with $X=a$ in $A$,
i.e., as describing the sheaf $\bHom(L,F)$ as a functor on
$A\op$. Accordingly, the question of describing the sheaf $\bHom(L,F)$
is reduced to the question of describing the objects $\Hom(L',F)$ in
\scrM, for $L'=\bZ^{(X)}\otimes L$. Thus, the main question here is to
give a ``computational'' description of the object $\Hom(L,F)$ in
\scrM, for $L$ in $\scrA\subab$ and $F$ in $\scrA_\scrM$, i.e., $F$
and $L$ being sheaves on $A$
\begin{equation}
  \label{eq:106.79}
  F:A\op\to\scrM,\quad
  L:A\op\to\Ab.\tag{79}
\end{equation}
The rule of the game here is to do so, using just
\hyperref[eq:106.75]{a)} in case of $X=a$ in $A$, and the exactness
property \hyperref[eq:106.75]{e)}.

It seems most convenient here to introduce again the additive envelope
$\Add(A)$ of the category $A$, which we'll assume to be small in what
follows, and the canonical additive functor
\begin{equation}
  \label{eq:106.80}
  \varepsilon\subab:\Add(A)\to\scrA\subab,\tag{80}
\end{equation}
extending the functor
\[a\mapsto\bZ^{(a)}:A\to\scrA\subab.\]
For a given $F$ \eqref{eq:106.77}, it follows from formula
\hyperref[eq:106.75]{(75 a))} that the composition
\[\widetilde F\circ\varepsilon\subab\op:\Add(A)\op
\xrightarrow{\varepsilon\subab\op} {\scrA\subab}\op
\xrightarrow{\widetilde F} \scrM\]
is just the canonical extension $\Add(F)$ of $F$ to $\Add(A)\op$,
whose value on the general object
\[x = \bigoplus_{i\in I} \Wh_A(a_i)\quad
\text{($I$ a finite indexing set)}\]
of $\Add(A)$ (where $\Wh_A(a)=\bZ^{(a)}$ as an object in
$\Add(A)\subset\Ahatab$) is just
\begin{equation}
  \label{eq:106.81}
  \Add(F)(x)=\prod_{i\in I}F(a_i).\tag{81}
\end{equation}
Now, it is easily checked that for any object $L$ in $\scrA\subab$,
i.e., any sheaf $L:A\op\to\Ab$, we have a canonical isomorphism in
$\scrA\subab$
\begin{equation}
  \label{eq:106.82}
  L\fromsim \varinjlim_{\text{$(x,u)$ in $\Add(A)_{/L}$}} \varepsilon\subab(x)\tag{82}
\end{equation}
(compare with the similar isomorphism on page \ref{p:403}). Using
\hyperref[eq:106.75]{(75~e)}, we\pspage{416} deduce from this the
expression
\begin{equation}
  \label{eq:106.83}
  \Hom(L,F)=\widetilde F(L)\simeq \varinjlim_{\text{$(x,u)$ in
      $\Add(A)_{/L}$}} \Add(F)(x),\tag{83}
\end{equation}
which in an evident way is functorial in $L$ for variable $L$.

This is about the best which can be done in general, it seems to me,
by way of ``computational'' expression of $\Hom(L,F)$ in terms of $F$
and $L$ given as in \eqref{eq:106.79}. Of course, the symbol
$\Add(A)_{/L}$ is relative to the canonical functor
$\varepsilon\subab$ \eqref{eq:106.80}, which is a full embedding in
case the site structure on $A$ is the chaotic one, i.e.,
$\scrA=\Ahat$. In computational terms, this category is rather
explicit, an object of the category is just a pair
\[(x,u)=\bigl((a_i)_{i\in I}, (u_i)_{i\in I}\bigr)\]
where $I$ is a finite indexing set, $(a_i)_{i\in I}$ a family of
objects of $A$, and for $i$ in $I$, $u_i$ is an element of $L(a_i)$ --
I'll leave to the reader the description of maps between such
objects. The value of $\Add(F)(x)$ is given by \eqref{eq:106.81}
above.

\begin{remark}
  The expression \eqref{eq:106.83} of $\Hom(L,F)=\widetilde F(L)$
  makes sense, provided only the additive category \scrM{} is stable
  under small inverse limits, without having to assume that \scrM{} be
  a pseudotopos. This makes us suspect that the functor
  \[\scrP\mapsto \scrP\circ\Wh :
  \bHom^!(\scrA\subab\op,\scrM)\to\scrA_\scrM\]
  is an equivalence (\eqref{eq:104.13} p.\ \ref{p:388}) without this
  extra assumption, provided \scrA{} is an actual topos (not only a
  pseudotopos as in loc.\ cit.). Indeed, we get a reasonable candidate
  for a quasi-inverse functor
  \[F\mapsto\widetilde
  F:\scrA_\scrM\to\bHom^!(\scrA\subab\op,\scrM).\]
  The only point still to check, with \eqref{eq:106.83} defining
  $\widetilde F$ for given $F$ in $\scrA_\scrM$, is that we get a
  functorial isomorphism
  \[\widetilde F(\bZ^{(a)}) \simeq F(a)\]
  for $a$ in $A$. In case $\scrA=\Ahat$, this follows from the fact
  that \eqref{eq:106.80} is fully faithful, hence $\Add(A)_{/L}$ for
  $L=\bZ^{(a)}$ admits a final object -- hence the limit
  \eqref{eq:106.83} is the value of $\Add(F)$ on the latter, namely
  $F(a)$.
\end{remark}

I feel the little program on the $\bHom$ and $\Hom$ operations, as
contemplated on page \ref{p:412}, is by now completed; all we've got
to do still is to dualize to get corresponding results for $*$ and
$\oast$. It's just a matter of essentially copying the formulaire
\eqref{eq:106.75}, which I'll do\pspage{417} for the sake of getting
more familiar with the more unusual operations $*$ and $\oast$. Now
of course, we'll have to restrict to the case $\scrA=\Ahat$, and use
the interpretation \eqref{eq:106.68} of the data $F$, $L'$ as functors
on $A\op$ and on $A$ with values in \scrM{} and \Ab{} respectively,
where now \scrM{} is an additive category stables under small
\emph{direct} limits. By duality, the ``sections'' or ``inverse
limits'' functor $\varprojlim_{A\op}$ (or ``cointegration'') is
replaced by the direct limit functor $\varinjlim_{A\op}$ (or
``integration''). With this in mind, we get the following
transcription of \eqref{eq:106.75}:
\begin{equation}
  \label{eq:106.84}
  \left\{%
    \renewcommand*{\arraystretch}{1.15}%
    \begin{array}{@{}r@{}rl@{}}
      \left\{\rule{0pt}{3.2ex}\right. &
      \begin{tabular}{@{}r@{}}
        a) \\ 
        a')
      \end{tabular} &
      \begin{array}{@{}l@{}}
        \smash{F*\bZ^{(b\op)}\simeq F(b)\quad\text{for any $b$ in $A$, hence
        $b\op$ in $B$}} \\
        \smash{F\oast\bZ^{(b\op)}\simeq(a\mapsto F(a\lor b)\simeq
        \varinjlim_{\text{$(x,\alpha)$ in $\preslice A{a\lor
        b}$}} F(x))}
      \end{array} \\
      \left\{\rule{0pt}{3.2ex}\right. &
      \begin{tabular}{@{}r@{}}
        b) \\
        b')
      \end{tabular} &
      \begin{array}{@{}l@{}}
        \smash{(F*L')*L'' \simeq F*(L'\otimes L'')} \\
        \smash{(F\oast L')*L'' \simeq F\oast(L'\otimes L'')}
      \end{array} \\
      \left\{\rule{0pt}{3.2ex}\right. &
      \begin{tabular}{@{}r@{}}
        c) \\
        c')
      \end{tabular} &
      \begin{array}{@{}l@{}}
        \smash{F*\bZ_\Bhat \simeq \varinjlim_{A\op} F} \\
        \smash{F \oast \bZ_\Bhat \simeq F}
      \end{array} \\
      & \text{d)} &\smash{F*L' \simeq \varinjlim_{A\op} F\oast L'} \\
      \left\{\rule{0pt}{6.2ex}\right. &
      \begin{tabular}{@{}r@{}}
        e) \\
        $\phantom{()}$ \\
        e') \\
        $\phantom{()}$
      \end{tabular} &
      \begin{tabular}{@{}l@{}}
        The functor $\smash{L'\mapsto F*L':\Bhatab\to\scrM}$ \\
        commutes to small direct limits \\
        Similar statement as e) for \\
        $\smash{L'\mapsto F\oast L':\Bhatab\to\AhatM}$.
      \end{tabular}
    \end{array}
  \right.\tag{84}
\end{equation}
\textbf{Comments.}\enspace This formulaire doesn't look wholly
symmetric to \eqref{eq:106.75}, due to the fact that we gave
\eqref{eq:106.75} in a somewhat more general context than topoi of the
type \Ahat{} only. This accounts for the letter $X$ or $Y$ in
\eqref{eq:106.75} (designating there an arbitrary object of \Ahat)
being replaced by a small letter $a$ or $b$ (designating objects in
$A$), which allows the dualization to be done. A slight trouble then
occurs when $A$ is not stable under binary products $a\times b$, these
products are only in \Ahat{} not in $A$, which accounts for the
slightly more complicated formula \hyperref[eq:106.84]{a')} of
\eqref{eq:106.84} in comparison to \hyperref[eq:106.75]{(75~a'))},
whose more explicit form, in the present context of data as in
\eqref{eq:106.68}, would be
\begin{equation}
  \label{eq:106.85}
  \bHom(\bZ^{(a)},F)\simeq\bigl(b\mapsto F(a\times
  b)\simeq\varprojlim_{\text{$(x,\alpha)$ in $A_{/a\times b}$}}
  F(x)\bigr).\tag{85} 
\end{equation}
Accordingly, the symbol $a\lor b$ (``sum'') in
\hyperref[eq:106.84]{(84~a'))} denotes the element $(a\op\times
b\op)\op$ of $(\Bhat)\op$ and can be identified with the \emph{sum} of
$a$ and $b$ in the category $A$\pspage{418} whenever the sum exists in
$A$. Accordingly, the category $\preslice A{a\lor b}$, dual to
$B_{/a\op\times b\op}$, can be described as
\begin{equation}
  \label{eq:106.86}
  \preslice A{a\lor b}=
  \begin{tabular}[t]{@{}l@{}}
    category of all triples $(x,u,v)$, with $x$ in $A$ and \\
    $u:a\to x$, $v:b\to x$ maps in $A$,
  \end{tabular}\tag{86}
\end{equation}
the maps in this category from $(x,u,v)$ to $(x',u',v')$ being just
maps $x\to x'$ ``compatible'' with the pairs $\alpha=(u,v)$ and
$\alpha'=(u',v')$ in the obvious way.

\begin{remarks}
  \namedlabel{rem:106.1}{1})\enspace
  An interesting particular case (although admittedly a little strange
  looking in our modelizing context!) is the one when $A$ is an
  additive category, hence stable under both binary sum and product
  operation, the two operations being canonically isomorphic, and
  written as $a\oplus b$. In this case, comparison of
  \eqref{eq:106.85} and \eqref{eq:106.86} shows that for a given
  object $a$ in $A$, hence $a\op$ in $B$, the operation
  $\Hom(\bZ^{(a)},{-})$ on \AhatM{} is canonically isomorphic to the
  operation ${-}\oast\bZ^{(a\op)}$. This immediately extends to a
  canonical isomorphism
  \begin{equation}
    \label{eq:106.87}
    \bHom(L,F)\simeq F\oast\check L \quad
    \text{for $L$ in $\Add(A)\subset\Ahatab$, 
      $F$ in \AhatM,}
    \tag{87}
  \end{equation}
  where we have denoted by
  \begin{equation}
    \label{eq:106.88}
    L\mapsto\check L: \Add(A)\op\toequ \Add(A\op)=\Add(B)\tag{88}
  \end{equation}
  the canonical antiequivalence between $\Add(A)$ and $\Add(B)$. In
  case $A$ is the final category, namely an additive category reduced
  to the zero object, and if we take moreover $\scrM=\Ab$,
  \eqref{eq:106.87} is the familiar formula of linear algebra, valid
  when $L$ is a free $\bZ$-module of finite type. It should be noted
  that $A$ being stable under binary products, it follows that
  $\Add(A)$ is stable under tensor products, and similarly for
  $\Add(B)$, and that the equivalence \eqref{eq:106.88} is compatible
  with tensor products. The relation \eqref{eq:106.87} is about the
  only relationship I could think of between the two types of
  operations upon a given category \AhatM.

  2)\enspace There are still two other, more trivial operations on a
  category \Ahat, of a similar nature to the two operations $\Hom$
  and $*$ considered so far. The more familiar one is componentwise
  tensor product
  \begin{equation}
    \label{eq:106.89}
    (L,F)\mapsto L\otimes F: \Ahatab\times\AhatM\to\AhatM,\tag{89}
  \end{equation}
  defined by
  \[{L\otimes F}\,(a) = L(a)\otimes F(a),\]
  where the second member denotes external tensor product of the
  abelian group $L(a)$ with the object $F(a)$ of \scrM{} (defined when
  \scrM{} is additive and stable under small direct limits). The
  other, deduced from \eqref{eq:106.89} by duality\pspage{419}
  \begin{equation}
    \label{eq:106.90}
    (L',F)\mapsto\bHom(L',F) : \Bhatab\times\AhatN\to\AhatN\tag{90}
  \end{equation}
  is defined when the additive category \scrN{} is stable under small
  inverse limits, and can be equally described as taking external
  $\Hom$'s componentwise
  \[\bHom(L',F)(a) = \Hom(L'(a),F(a)).\]
  These operations make sense too when \Ahat{} is replaced by an
  arbitrary topos \scrA, \Bhatab{} being replaced by the category of
  abelian cosheaves on \scrA. It doesn't seem worthwhile here to dwell
  on them, as they don't seem to be so relevant for the homology and
  cohomology formalism we want to develop in the next sections. I like
  to point out, though, that in the cohomology formalism of ringed
  topoi the tensor product operation \eqref{eq:106.89} and the derived
  operation $\overset{\mathrm L}{\otimes}$ on the relevant derived
  categories $\D_\bullet$ play an important role, and it is likely
  therefore that in a more extensive development of the homology and
  cohomology formalism within the context of topoi \Ahat{} and maps in
  \Cat, the same will hold for the dual operation \eqref{eq:106.90}
  too.

  The reader who may feel confused by the manifold use of the symbol
  $\bHom$ should notice that there is no possibility of confusion
  reasonably between \eqref{eq:106.90} and \eqref{eq:106.65} (p.\
  \ref{p:409}), as the argument $L'$ in \eqref{eq:106.90} is in
  \Bhatab, whereas the argument $L$ in \eqref{eq:106.65} is in an
  altogether different category \Ahatab. In the case when $A$ is the
  final category say, hence $A=B$, and a confusion might arise, the
  two operations turn out to be actually the same (up to canonical
  isomorphism). A similar remark applies to the fear of confusion
  between the kindred operations $\otimes$ and $\oast$. I daresay I
  devoted a considerable amount of attention on terminology and
  notation around the abelianization story -- and it does seem that a
  pretty coherent formalism is emerging indeed.
\end{remarks}

\bigbreak
\presectionfill\ondate{17.8.}\par

\hangsection{Review \texorpdfstring{\textup{(4)}}{(4)}: Case of
  general ground ring \texorpdfstring{$k$}{k}.}\label{sec:107}%
\phantomsection\addcontentsline{toc}{subsection}{\numberline
  {F)}Extension of ground ring from \texorpdfstring{$\bZ$ to $k$
  \textup($k$}{Z to k (k}-linearization\texorpdfstring{\textup)}{)}.}%
\textbf{\namedlabel{subsec:107.F}{F)}\enspace Extension of ground ring
  from $\bZ$ to $k$ ($k$-linearization).}\enspace%
I would like still to make a quick review of the main facts and
formulas of the last two sections, replacing throughout the ground
ring $\bZ$ by an arbitrary \emph{commutative} ring $k$, and additive
categories \scrM{} and additive functors between these, by
$k$-additive categories and $k$-additive functors. This will allow us
to check that the conceptual and notational set-up we got so far
extends smoothly to $k$-linearization.

Let's\pspage{420} recall that a \emph{$k$-additive category} \scrM{}
is an additive category endowed with the extra structure given by a
homomorphism of commutative rings
\begin{equation}
  \label{eq:107.91}
  k\to\End(\id_\scrM),\tag{91}
\end{equation}
where the second member denotes the (commutative) ring of all
endomorphisms of the identity functor of \scrM{} to itself. Defining
accordingly the notion of \emph{$k$-additive functor} between two
$k$-additive categories $\scrM$, $\scrM'$, we'll denote by
\begin{equation}
  \label{eq:107.92}
  \bHom_k(\scrM,\scrM')\subset\bHom(\scrM,\scrM')\tag{92}
\end{equation}
the full subcategory of $\bHom(\scrM,\scrM')$ made up with such
functors. Thus, we get a canonical fully faithful inclusion
\begin{equation}
  \label{eq:107.93}
  \bHom_k(\scrM,\scrM') \hookrightarrow \bHom_{\bZ}(\scrM,\scrM')
  \eqdef \Homadd(\scrM,\scrM').\tag{93}
\end{equation}
Wed defined accordingly the categories $\bHom_{k!}$, $\bHom_k^!$ as
full subcategories of \eqref{eq:107.92}, and the category
\[\bHom_k(\scrP,\scrQ;\scrM) \subset \bBiadd(\scrP,\scrQ;\scrM)\]
the full subcategory of $\bHom(\scrP\times\scrQ,\scrM)$ made up with
\emph{$k$-bilinear} functors, namely functors $k$-additive in each
argument (in case $k=\bZ$,\scrcomment{Actually, it was previously
  denoted by just $\bHom(\scrP,\scrQ;\scrM)$\ldots} this is the
category denoted previously by $\bBiadd$), and similarly for the
notations $\bHom_{k!!}$ and $\bHom_k^{!!}$.

It should be noted that for a given additive category \scrM, there is
a ``best'' choice for endowed it with a $k$-linear structure, in such
a way that any $k'$-linear structure just corresponds to ``ground ring
restriction'' with respect to suitable (well-defined) ring
homomorphism
\[k'\to k;\]
we just take the ``tautological'' linear structure with
\[k=\End(\id_\scrM),\]
and \eqref{eq:107.91} the identity.

If $A$ is any small category, we'll denote by
\begin{equation}
  \label{eq:107.94}
  \Ahatk = A\uphat_{k\textup{-Mod}} \simeq \bHom(A\op,\kMod)\tag{94}
\end{equation}
the category of objects in \Ahat{} endowed with a structure of
$k$-module, i.e., the category of presheaves on $A$ with values in the
category \kMod{} of $k$-modules (in the given basic universe
\scrU). This is of course a $k$-additive category, which for $k=\bZ$
reduces to the category\pspage{421} of additive presheaves on $A$:
\[ A\uphat_{\bZ} \eqdef \Ahatab.\]
We have, for a homomorphism of commutative rings
\[k\to k',\]
a corresponding functor between additive topoi
\begin{equation}
  \label{eq:107.95}
  \Ahatk \to A\uphat_{k'}, \quad
  F\mapsto F\otimes_k k' = (a\mapsto F(a)\otimes_k k'),\tag{95}
\end{equation}
by which we may interpret if we wish, in a rather evident way the
$k'$-linear topos $A\uphat_{k'}$ as deduced from the $k$-linear one
\Ahatk{} by ``ground ring extension'' $k\to k'$, namely as the
solution of a $2$-universal problem with respect to categories
$\bHom_{k!}(\Ahatk,\scrM)$, where \scrM{} is a $k'$-additive category
stable under small direct limits. The $k$-abelianization functor
\begin{equation}
  \label{eq:107.96}
  \Ahat\to\Ahatk, \quad X\mapsto k^{(X)}\quad\bigl(\;\simeq (a \mapsto
  k^{(X(a))})\bigr)\tag{96} 
\end{equation}
or $\Wh_{\Ahat,k}$, is defined as the composition
\[\Ahat\to A\uphat_{\bZ} = \Ahatab\to \Ahatk,\]
where the first functor is the familiar abelianization
$X\mapsto \bZ^{(X)}$, and the second is ground ring extension for
$\bZ\to k$. If \scrM{} is any $k$-additive category stable under small
direct limits, \eqref{eq:107.96} gives rise to a functor which is an
equivalence of categories $F\mapsto(X\mapsto F(k^{(X)}))$
\begin{equation}
  \label{eq:107.97}
  \bHom_{k!}(\Ahatk,\scrM)\toequ\bHom_!(\Ahat,\scrM) \quad (\;\toequ
  \bHom(A,\scrM)),\tag{97} 
\end{equation}
where the second equivalence is the familiar one of prop.\
\ref{prop:105.3} (p.\ \ref{p:394}), independent of any abelian
assumptions. Dually, we get an equivalence
\begin{equation}
  \label{eq:107.98}
  \bHom_k^!((\Ahatk)\op,\scrM)\toequ \bHom((\Ahat)\op,\scrM)
  \quad(\;\toequ \bHom(A\op,\scrM)),\tag{98}
\end{equation}
where \scrM{} is any $k$-additive category stable under small inverse
limits. From \eqref{eq:107.97} \eqref{eq:107.98} and replacing in
\eqref{eq:107.98} $A$ by the dual category $B=A\op$, and assuming the
$k$-additive category \scrM{} is stable under both types of small
limits, we get the duality equivalence
\begin{equation}
  \label{eq:107.99}
  \bHom_{k!}(\Ahatk,\scrM)\equ \bHom_k^!(\Bhatk,\scrM) \quad (\;\equ
  \bHom(A,\scrM)).\tag{99} 
\end{equation}
This may be viewed as giving two alternative descriptions, by the two
members of \eqref{eq:107.99}, of the category
\[\BhatM = \bHom(B\op=A,\scrM)\]
of \scrM-valued presheaves on $B$ (defined without any use of the
$k$-additive structure of \scrM). The left-hand side interpretation
\eqref{eq:107.99}, via \scrM-valued\pspage{422} $k$-additive cosheaves
on the $k$-additive topos \Ahatk, gives rise to the operations $*_k$
and $\oast_k$ of \Ahatk{} upon \BhatM{} (operations previously denoted
by $*$ and $\oast$ when $k=\bZ$ and no confusion would arise from
dropping subscripts), and similarly the interpretation by right-hand
side of \eqref{eq:107.99}, via \scrM-valued $k$-additive sheaves on
the $k$-additive topos \Bhatk, gives rise to the operations $\Hom_k$
and $\bHom_k$ of \Bhatk{} upon \BhatM. Replacing in this comment $A$
by $B$, hence \Bhat{} by \Ahat, namely in terms of operations upon the
category of \scrM-valued sheaves on the topos \Ahat{} (or \scrM-valued
presheaves on $A$), we get the mutually dual pair of operations
\begin{equation}
  \label{eq:107.100}
  \begin{aligned}
    (F,L')&\mapsto F*_kL' : \AhatM\times\Bhatk\to\scrM, \\
    (F,L')&\mapsto F\oast_kL' : \AhatM\times\Bhatk\to\AhatM
  \end{aligned}
  \tag{100}
\end{equation}
and
\begin{equation}
  \label{eq:107.101}
  \begin{aligned}
    (L,F)&\mapsto \Hom_k(L,F) : \Ahatk\times\AhatM\to\scrM, \\
    (L,F)&\mapsto \bHom_k(L,F) : \Ahatk\times\AhatM\to\AhatM.
  \end{aligned}
  \tag{101}
\end{equation}
The operations \eqref{eq:107.100} are ruled by formulaire
\eqref{eq:106.84} (with subscripts $k$ added), whereas the operations
\eqref{eq:107.101} are ruled by formulaire \eqref{eq:106.75} with
subscripts (see moreover for the latter comments on page \ref{p:417},
and formula \eqref{eq:106.85} for \hyperref[eq:106.75]{(75~a'))}; they
are valid provided the additive category is stable under small direct
resp.\ inverse limits. Moreover, we get a ``computational'' expression
of $\Hom_k(L,F)$ by a formula extending \eqref{eq:106.83} which we'll
still have to write down, and correspondingly for $F*_kL'$ (by a dual
formula, which we forgot to include in the previous section). To do
so, we have to introduce still
\begin{equation}
  \label{eq:107.102}
  \Add_k(A)\subset\Ahatk,\tag{102}
\end{equation}
the $k$-additive envelope of $A$, which may be described (beside by
the familiar $2$-universal property in the context of $k$-additive
categories and functors from $A$ into these) as the full subcategory
of \Ahatk{} generated by finite sums of objects of the type $k^{(a)}$
with $a$ in $A$ -- i.e., the general object of $\Add_k(A)$ may be
written
\[\bigoplus_{i\in I}k^{(a_i)},\]
where $(a_i)_{i\in I}$ is any finite family of objects of $A$. When
the finiteness condition on $I$ is dropped, we get a larger full
subcategory
\begin{equation}
  \label{eq:107.103}
  \Addinf_k(A)\subset\Ahatk,\tag{103}
\end{equation}
which may also be interpreted as ``the'' solution of the $2$-universal
problem of sending $A$ into categories which are $k$-additive and
moreover infinitely additive, i.e., stable under small direct
sums. Enlarging the subcategories \eqref{eq:107.102} and
\eqref{eq:107.103} of \Ahatk{} by adjoining all objects of
\Ahatk\pspage{423} isomorphic to direct factors of objects in the
considered subcategory, we get to (strictly) full subcategories of
\Ahatk{} containing the latter, which may be interpreted as being just
the subcategory $\Proj(\Ahatk)$ of \emph{projective} objects of
\Ahatk{} when starting with \eqref{eq:107.103}, and as the subcategory
$\UlProj(\Ahatk)$ of \emph{ultraprojective} objects, namely objects
projective and of finite presentation, when starting with
\eqref{eq:107.102}. These may be equally interpreted as the abstract
Karoubi envelops of the categories \eqref{eq:107.103} and
\eqref{eq:107.102}, deduced from these formally by adjoining images
(=coimages) of projectors (or equivalently, as $2$-universal solutions
of the $2$-universal problem of sending the given category
\eqref{eq:107.103} or \eqref{eq:107.102} into ``karoubian
categories'', namely categories stable under images (=coimages) of
projectors, with maps between these being functors commuting to those
images or coimages of projectors):
\begin{equation}
  \label{eq:107.104}
  \Proj(\Ahatk)\equeq\Kar(\Addinf_k(A)), \quad
  \UlProj(\Ahatk)\equeq\Kar(\Add_k(A)).\tag{104}
\end{equation}
Accordingly, these two categories may be equally described, directly
in terms of $A$, as the solutions of the two $2$-universal problems,
obtained from mapping $A$ into $k$-additive karoubian categories,
which in the first case (corresponding to $\Proj(\Ahatk)$) are
moreover assume infinitely additive.

To sum up the situation, we get in \Ahatk{} a diagram of four
remarkable full subcategories \eqref{eq:107.102}, \eqref{eq:107.103},
\eqref{eq:107.104}, which may be interpreted (as well as \Ahatk{}
itself) as the solutions of five corresponding ``$k$-additive''
$2$-universal problems, in terms of sending $A$ into $k$-additive
categories satisfying suitable extra exactness assumptions (namely
being karoubian for the two categories in \eqref{eq:107.104}, being
infinitely additive for the two categories $\Addinf_k(A)$ and its
Karoubi envelope $\Proj(\Ahatk)$, and being stable for small direct
limits in case of \Ahatk). Including equally the non-additive
categories $A$ and \Ahat{} and the functors $A\to\Add_k(A)$,
$\Ahatk\to\Ahat$, we get a seven term diagram of canonical functors
between categories of presheaves upon $A$:
\begin{equation}
  \label{eq:107.105}
  \left\{
    \begin{tabular}{@{}c@{}}
      \begin{tikzcd}[baseline=(O.base),sep=small]
        A\ar[d] & & \\
        \Add_k(A)\ar[r,hook]\ar[d,hook] &
        \UlProj(\Ahatk)\equeq\KarAdd_k(A) \ar[d,hook,shift right=3.5em] & \\
        \Addinf_k(A)\ar[r,hook] &
        \Proj(\Ahatk)\equeq\KarAddinf_k(A)\ar[r,hook] & \Ahatk\ar[d] \\
        & & |[alias=O]| \Ahat
      \end{tikzcd},
    \end{tabular}\right.
  \tag{105}
\end{equation}
where the five categories in the two intermediate lines are
$k$-additive as well as all functors between them in the diagram,
which are moreover fully faithful. For any $k$-additive\pspage{424}
category \scrM{} stable under small direct limits, taking cosheaves on
\Ahat{} with values in \scrM, and their restrictions to the six other
categories in the diagram \eqref{eq:107.105}, we get a transposed
seven term diagram as follows, part of which reduces to the four term
diagram \eqref{eq:105.58} (p.\ \ref{p:406}) in case $k=\bZ$:
\begin{widematter}
  \begin{equation}
  \label{eq:107.106}
  \left\{
    \begin{tabular}{@{}c@{}}
      \begin{tikzcd}[baseline=(O.base),sep=small]
        \bHom_!(\Ahat,\scrM) \ar[d,"\equ"'] & & \\
        \bHom_{k!}(\Ahatk,\scrM) \ar[r,"\equ"] &
        \bHomaddinfkar_k(\Proj(\Ahatk),\scrM) \ar[r,dash,"\equ"]
        \ar[d,"\equ"] &
        \bHomaddinf_k(\Addinf_k(A),\scrM) \ar[d,"\equ"] \\
        & \bHomaddkar_k(\UlProj(\Ahatk),\scrM) \ar[r,dash,"\equ"] &
        \bHomadd_k(\Add_k(A),\scrM) \ar[d,"\equ"] \\
        & & |[alias=O]| \bHom(A,\scrM)
      \end{tikzcd},
    \end{tabular}\right.
  \tag{106}
\end{equation}
\end{widematter}
where the meaning of the symbols used (such as index $k$, suffixes
``add'' or ``addinf'' and ``kar'') for qualifying $\bHom$ and denoting
various full subcategories of $\bHom$ categories, is clear from the
explanations given previously. Replacing $A$ by $B$ and
$\bHom_!(\Ahat,\scrM)$ by $\bHom^!(\Bhat,\scrM)$, we get a diagram
``dual'' to \eqref{eq:107.106} (containing the five-term diagram
\eqref{eq:106.62} (p.\ \ref{p:408}) in case $k=\bZ$), which we'll not
write out here, valid for any $k$-additive category \scrM{} stable
under small inverse limits. When \scrM{} is a $k$-additive category
stable under both types of small limits, then the last term of the
diagram \eqref{eq:107.106} is equal to the last term of the dual one,
hence a system of fourteen mutually equivalent categories (compare p.\
\ref{p:408}, when we considered ten among them only!), expressing as
many ways for interpreting the notion of an \scrM-valued copresheaf on
$A$, i.e., an object of $\bHom(A,\scrM)$ (which is one among the
fourteen\ldots).

Let's comment a little on the significance of the various five
$k$-additive categories appearing in \eqref{eq:107.105}. The largest
one \Ahatk{} is there precisely as the all-encompassing category of
$k$-additive presheaves, where to carry through all kinds of
$k$-linear constructions between presheaves on $A$. The significance
of the (second largest) subcategory $\Proj(\Ahatk)$, made up with all
projective objects of \Ahatk, comes mainly from homological algebra
and emphasis upon replacing objects of \Ahatk{} by projective
resolutions; these are chain complexes in $\Proj(\Ahatk)$, which may
be viewed as being defined (by any given object in \Ahatk) ``up to
chain homotopy''. More sweepingly still, we get from general
principles the canonical equivalence of categories
\begin{equation}
  \label{eq:107.star}
  \D^-(\Ahatk) \fromequ \mathrm K^-(\Proj(\Ahatk)),\tag{*}
\end{equation}
where $\D^-$ designates the ``derived category bounded from
above''\pspage{425} of a given \emph{abelian} category (defined in
terms of differential operators with degree $+1$, and
quasi-isomorphisms between complexes with degrees bounded from above),
whereas $\mathrm K^-$ designates localization of the category of
differential complexes with degrees bounded from above of a given
\emph{additive} category, localization being taken with respect to
homotopisms.

As any object of $\Proj(\Ahatk)$ is a direct factor of an object in
$\Addinf_k(A)$, and hence, any object in \Ahatk{} is isomorphic to a
quotient of an object in $\Addinf_k(A)$, it follows again from general
principles that the categories in \eqref{eq:107.star} are equally
equivalent to $\mathrm K^-(\Addinf_k(A))$, hence
\begin{equation}
  \label{eq:107.107}
  \D^-(\Ahatk) \fromequ \mathrm K^-(\Proj(\Ahatk)) \fromequ
  \mathrm K^-(\Addinf_k(A)).\tag{107}
\end{equation}
The advantage of $\Addinf_k(A)$ over $\Proj(\Ahatk)$ is that its
objects, and maps between objects, are more readily described in
computational terms, just working with small direct sums of objects of
the type $k^{(a)}$ (with $a$ in $A$), and corresponding matrices, with
entries in free $k$-modules $k^{(\Hom(a,b))}$. Thus, if we call
\emph{cointegrator} (with coefficients in $k$) \emph{for} $A$ any
projective resolution of the constant presheaf $k_A$ with value $k$,
and denote such object by $L_k^A$, we may view $L_k^A$ as an object
determined up to unique isomorphism, either in $\mathrm
K^-(\Proj(\Ahatk))$, or in $\mathrm K^-(\Addinf_k(A))$ -- and it is
the latter interpretation which looks the most convenient. Objects in
the first category, namely complexes with degrees bounded from above
and projective components, which happen to be in the first (i.e.,
components are in $\Addinf_k(A)$, i.e., are direct sums of objects of
the type $k^{(a)}$) may be called ``\emph{quasi-special}'' (extending
the terminology previously used for cointegrators and integrators, in
case $k=\bZ$). We'll call them \emph{special} if the components are
even in $\Add_k(A)$. The category $\Add_k(A)$ and its Karoubi envelope
$\UlProj(\Ahatk)$ may be viewed both as embodying \emph{finiteness
  conditions}, and similarly for the two corresponding $\mathrm K^-$
categories, which are of course equivalent:
\begin{equation}
  \label{eq:107.108}
  \mathrm K^-(\UlProj(\Ahatk)) \fromequ \mathrm K^-(\Add_k(A)),\tag{108}
\end{equation}
and presumably the canonical functor from \eqref{eq:107.108} to
\eqref{eq:107.107} is fully faithful, under suitable coherence
conditions at any rate\ldots

\bigbreak
\presectionfill\ondate{22.8.}\pspage{426}\par

\hangsection{Review \texorpdfstring{\textup{(5)}}{(5)}: Homology and
  cohomology \texorpdfstring{\textup(absolute case\textup)}{(absolute
    case)}.}\label{sec:108}%
Since last Monday, namely for about one week, I have been mainly taken
by a rather dense sequence of encounters and events, the center of
which has been the unexpected news of my granddaughter Ella's death at
the age of nine, by a so-called health accident. I resumed some
mathematical pondering last night. Today, I got a short letter from
Ronnie Brown, mainly with the announcement of the loss of his son
Gabriel, twenty years old, which occurred about the same time by a
climbing accident. It is a good thing that Ronnie felt like telling me
in a few words about this, while we have never yet seen each other and
our letters so far have been restricted to mathematics, with maybe
sometimes some personal comments about his or my own involvement in
mathematics. It is through these, surely, that a mutual sympathy has
come into being, not merely motivated by a common interest in
mathematics -- and this sympathy I feel has been the main force giving
life to our correspondence while mathematically speaking more than
once it has been rather a ``dialogue de
sourds''.\scrcomment{``dialogue of the deaf''} (This is due mainly to
my illiteracy homotopy in theory, and to my reluctance to get really
involved in any ``technical'' matters, until I am really forced to by
what I am just doing.)

I want now to go on with the overall review on ``abelianization'' and
its relation to the homology and cohomology formalism for small
categories, serving as models for homotopy types.

\addtocounter{subsection}{6}
\subsection{Homology and cohomology (absolute case).}
\label{subsec:108.G}
My aim is to give a perfectly dual treatment of cohomology and
homology, which is one main reason why I have to take as coefficients
for both, not merely usual abelian presheaves on a given small
category $A$, or sheaves of $k$-modules for a given ring $k$, but more
generally sheaves with values in any abelian category \scrM, stable
under small direct or inverse limits (according as to whether we are
interesting in taking homology, or cohomology invariants). It will
then turn out that homology of $A$ for \scrM-valued presheaves (or
complexes of such) is ``the same'' as cohomology of the dual category
$B$, with coefficients in the corresponding $\scrM\op$-valued ones.

As I am a lot more familiar with cohomology, it is by this I'll begin
again. Here, as in the case of an arbitrary topos \scrX, the
cohomology invariants $\mathrm H^i(\scrX,F)$ with values in an abelian
sheaf $F$ may be\pspage{427} viewed as being just the invariants
$\Ext^i(\bZ_\scrX,F)$ in the category of all abelian sheaves, where
$\bZ_\scrX$ is the constant sheaf on \scrX{} with value $\bZ$. The
similar fact holds when $F$ is any sheaf of modules over a sheaf of
rings $\scrO_\scrX$ on \scrX, with $\bZ_\scrX$ being replaced by
$\scrO_\scrX$ in the interpretation above:
\[\mathrm H^i(\scrX,F) \simeq \Ext^i_{\scrO_\scrX}(\scrO_\scrX, F),\]
which is often quite useful in the cohomology formalism. We are going
to restrict here to the case of a constant sheaf of rings, defined by
a fixed commutative ring $k$, which will allow us to play around as
announced with the duality relation between $A$ and $B=A\op$, provided
moreover we take more general coefficients still, namely (pre)sheaves
with value in a given $k$-additive category \scrM{} stable under the
relevant limits. (Presumably, the case of a locally constant
commutative sheaf of rings could be dealt with too, but we'll not dive
into this here!) Another important (and by now familiar?) conceptual
point is that, rather than the $\Ext^i$'s which give only partial
information, we are interested in the object they come from (as the
``cohomology objects''), namely the objects $\mathrm
R\Hom_\scrOX(L,F)$ in a suitable derived category. In the present
case when \scrX{} is the topos associated to the small category $A$,
hence the category of \scrOX-modules has sufficiently many projective
(namely direct sums of sheaves of the type $\scrO_\scrX^{(a)}$ with
$a$ in $A$), the $\mathrm R\Hom_\scrOX(L,F)$ may be computed, taking a
projective resolution $L_\bullet$ of $L$, by the formula
\[\mathrm R\Hom_\scrOX(L,F) \simeq \Hom_\scrOX^\bullet(L_\bullet,F)\]
(an isomorphism in $\D^+\kMod$ say), and similarly when replacing $L$
and $F$ by arguments $L_\bullet$ and $F^\bullet$ in $\D^-$ and $\D^+$
of the category of \scrOX-modules. As a result, we get a pairing,
computable here using projective resolutions of the argument
$L_\bullet$:
\begin{equation}
  \label{eq:108.star}
  (L_\bullet,F^\bullet)\mapsto\mathrm
  R\Hom_\scrOX(L_\bullet,F^\bullet) : \D^-(\scrOX)\times\D^+(\scrOX)
  \to \D^+\kMod,\tag{*}
\end{equation}
where $k$ is a commutative ring and \scrOX{} is endowed with a
structure of $k$-algebra. Using $\bHom_\scrOX$ and its total derived
functor, we get likewise
\begin{equation}
  \label{eq:108.starstar}
  (L_\bullet,F^\bullet)\mapsto\mathrm
  R\bHom_\scrOX(L_\bullet,F^\bullet) : \D^-(\scrOX)\times\D^+(\scrOX)
  \to \D^+(\scrOX),\tag{**}
\end{equation}
with
\[\mathrm R\bHom_\scrOX(L_\bullet,F^\bullet) \simeq
\bHom_\scrX^{\bullet\bullet}(\scrL_\bullet,F^\bullet),\]
where $\scrL_\bullet$ is a projective resolution of $L_\bullet$, and
$\bHom^{\bullet\bullet}$ stands for the simple complex associated to
the double complex obtained by taking\pspage{428} $\bHom$'s
componentwise (and we have the similar formula of course for the
$\Hom$'s and $\mathrm R\Hom$'s non-bold-faced).

Taking $L_\bullet=\scrOX$ (or any resolution of \scrOX), the $\mathrm
R\Hom$ invariant \eqref{eq:108.star} reduces to the total derived
functor $\mathrm R\Gamma$ of the sections functor
\begin{equation}
  \label{eq:108.i}
  \mathrm R\Hom_\scrOX(\scrOX,F^\bullet)\simeq\mathrm R\Gamma_\scrX(F^\bullet)\tag{i}
\end{equation}
(whereas $\mathrm R\bHom_\scrOX(\scrOX,F^\bullet)\simeq F^\bullet$ of
course), which in turn allows to give the following familiar
expression of $\mathrm R\Hom$ in terms of $\mathrm R\bHom$:
\begin{equation}
  \label{eq:108.ii}
  \mathrm R\Hom_\scrOX(L_\bullet,F^\bullet) \simeq \mathrm
  R\Gamma_\scrX(\mathrm R\bHom_\scrOX(L_\bullet,F^\bullet)),\tag{ii}
\end{equation}
coming from the similar isomorphism
$\Hom_\scrOX\simeq\Gamma_\scrX\bHom_\scrOX$. All this is standard
cohomology formalism, valid on an arbitrary ringed topos
$(\scrX,\scrOX)$, except for the possibility of computing $\RHom$ and
$\RbHom$ by taking \emph{projective resolutions of the first argument}
(rather than injective ones of the second), which is special to the
case when $\scrX=\Ahat$, to which we'll now restrict.

Let now $A$ be a fixed small category, $k$ a fixed commutative ring,
\scrM{} a $k$-additive category, stable under small inverse limits. We
want to define a total derived functor of the functor
\begin{equation}
  \label{eq:108.109}
  (L,F)\mapsto \Hom_k(L,F):\Ahatk \times \AhatM\to\scrM,\tag{109}
\end{equation}
which should be a functor
\begin{equation}
  \label{eq:108.110}
  (L_\bullet,F^\bullet)\mapsto \RHom_k(L_\bullet,F^\bullet) :
  \D^-(\Ahatk)\times \D^+(\AhatM) \to \D^+(\scrM),\tag{110}
\end{equation}
and similarly
\begin{equation}
  \label{eq:108.111}
  (L_\bullet,F^\bullet)\mapsto\RbHom_k(L_\bullet,F^\bullet) :
  \D^-(\Ahatk) \times \D^+(\AhatM) \to \D^+(\AhatM).\tag{111}
\end{equation}
For this, in order for $\D^+(\AhatM)$ to be defined, we better assume
\scrM{} to be an \emph{abelian} category, hence \AhatM{} is abelian
too. Of course, we'll write
\begin{equation}
  \label{eq:108.112}
  \begin{aligned}
    \Ext_k^i(L_\bullet,F^\bullet) &=
    \mathrm H^i(\RHom_k(L_\bullet,F^\bullet)), \\
    \bExt_k^i(L_\bullet,F^\bullet) &=
    \mathrm H^i(\RbHom_k(L_\bullet,F^\bullet)),
  \end{aligned}\tag{112}
\end{equation}
these global and local $\Ext^i$ may be viewed as ``external''
$\Ext^i$'s, as contrarily to the familiar case, the components of the
two arguments $L_\bullet$ and $K^\bullet$ are not in the same category
-- just as the $\Hom_k$ in \eqref{eq:108.109} and the corresponding
\begin{equation}
  \label{eq:108.109prime}
  \bHom_k:\Ahatk\times\AhatM\to\AhatM\tag{109'}
\end{equation}
has arguments in the two different categories \Ahatk{} and \AhatM.

As\pspage{429} we don't know about the existence of enough injective
in \AhatM, the only way for defining the pairings \eqref{eq:108.110},
\eqref{eq:108.111} is now by using projective resolutions of the first
argument, writing
\begin{equation}
  \label{eq:108.113}
  \begin{aligned}
    \RHom_k(L_\bullet,F^\bullet) &=
    \Hom_k^{\bullet\bullet}(\scrL_\bullet,F^\bullet), \\
    \RbHom_k(L_\bullet,F^\bullet) &=
    \bHom_k^{\bullet\bullet}(\scrL_\bullet,F^\bullet) 
  \end{aligned}\tag{113}
\end{equation}
where $\scrL_\bullet$ is a projective resolution of $L_\bullet$ in
\Ahatk. As the latter is defined up to chain homotopy, it follows that
for fixed $L_\bullet$ and $F^\bullet$, the second members of
\eqref{eq:108.113} are defined up to chain homotopy, i.e., they may be
viewed as objects in $\mathrm K^+(\scrM)$ and $\mathrm K^+(\AhatM)$
respectively. They are defined as such, even without assuming \scrM{}
to be abelian and hence $\D^+(\scrM)$ and $\D^+(\AhatM)$ to be
defined. When we make this assumption, in order to check that the
formulæ \eqref{eq:108.113} do define pairings as in \eqref{eq:108.110}
and \eqref{eq:108.111}, we still have to check that for a
quasi-isomorphism
\[F^\bullet\to (F')\bullet\]
in \AhatM, the corresponding maps between $\RHom_k$ and $\RbHom_k$ are
quasi-isomorphisms too. This will follow immediately, provided we
check that for fixed \emph{projective} $L$ in \Ahatk{} and variable
$F$ in \AhatM, the functors
\[F\mapsto\Hom_k(L,F)\quad\text{and}\quad
F\mapsto\bHom_k(L,F)\]
from \AhatM{} to \scrM{} and \AhatM{} respectively are exact. Now,
this is clear for $\Hom_k$ when $L$ is of the type $k^{(a)}$, hence
$\Hom(L,F)\simeq F(a)$, hence it follows when $L$ is a small direct
sum of objects $k^{(a_i)}$, hence
\[\Hom_k(L,F) = \prod_i F(a_i),\]
provided we make on \scrM{} the mild extra assumption that \emph{a
  small direct product of epimorphisms is again an epimorphism}. As
any projective object of \Ahatk{} is a direct factor of a small direct
sum as above, the exactness result we want then follows, hence the
looked-for pairing \eqref{eq:108.110}. The same argument will hold for
$\bHom_k$, provided we check exactness of the functors
\[F\mapsto \bHom_k(k^{(a)},F)=(b\mapsto F(a\times b)).\]
Here it seems we get into trouble when $A$ is not stable under binary
products -- in this case there is little chance that the functor above
be exact, even when restricting to the case $\scrM=\Ab_k \eqdef\kMod$,
hence $\AhatM=\Ahatk$ and $L$, $F$ have values in the same category
(namely presheaves of $k$-modules). This may seem strange, as we know
(and recalled above) that in this standard case there is no problem
for defining a\pspage{430} pairing \eqref{eq:108.111} $\RbHom_k$. The
point here is that, whereas a reasonable $\RbHom$ can be defined
indeed, it \emph{cannot} be computed in terms of a projective
resolution of the first argument as in \eqref{eq:108.113}; or
equivalently, that for projective $L$ it is not necessarily true that
\[\bExt_k^i(L,F)=0\quad\text{for $i>0$;}\]
this in turn relates to the observation that, contrarily to what
happens for the notion of injective sheaves of modules (on an
arbitrary topos), it is not true that the property for a sheaf of
modules to be projective is stable under localization (even for
a constant sheaf of rings $k$ on a topos \Ahat). Indeed, the
localization of $k^{(a)}$ with respect to $A_{/b}$ (with $a$ and $b$
in $A$) is $k^{(a')}$ where $a'$ is $a\times b$ viewed as an object in
$A_{/b}$, and for any sheaf of $k$-modules $F$ on $A_{/b}$ we have
\[\Ext_{A_{/a}}^i(k^{(a')},F) = \mathrm H^i(A_{/a'=a\times b},F),\]
which need not be zero for $i>0$. If it was, this would imply that
$a\times b$ is $k$-acyclic (rather, that its connected components
are), a rather strong property indeed when $a\times b$ is not in
$A$\ldots

Thus, when $A$ is not stable under binary products, it doesn't seem
that there exists a pairing \eqref{eq:108.111} as I expected, except
(possibly) when there are enough injectives in \AhatM{} -- a case I do
not wish to examine for the time being, as I am mainly interested now
in a formalism using projective resolutions instead of injective
ones. Anyhow, for the purpose of subsuming the cohomology functor
$\RGamma_A$ under the $\RHom_k$ formalism, by formula
\begin{equation}
  \label{eq:108.114}
  \RGamma_A(F^\bullet) = \RHom_k(k_A,F^\bullet) \simeq
  \Hom_k^{\bullet\bullet}(L_\bullet^A,F^\bullet),\tag{114}
\end{equation}
where
\begin{equation}
  \label{eq:108.115}
  L_\bullet^A\to k_A\tag{115}
\end{equation}
is a projective resolution of $k_A$, it is the pairing $\RHom_k$ and
not $\RbHom_k$ which is the relevant one. Let's recall that a
projective resolution \eqref{eq:108.115} is called a
\emph{cointegrator} (for the category $A$, with coefficients in $k$),
as by formula \eqref{eq:108.114} it allows indeed to express
``cointegration'' of any \scrM-valued presheaf or complex of such
presheaves (with degrees bounded from below).

Thus, for the time being we just got the pairing $\RHom_k$ in
\eqref{eq:108.110}, and the corresponding functor
\begin{equation}
  \label{eq:108.116}
  F^\bullet\mapsto \RGamma_A(F^\bullet) : \D^+(\AhatM) \to \D^+(\scrM),\tag{116}
\end{equation}
and\pspage{431} not the pairing $\RbHom_k$ in \eqref{eq:108.111}, and
hence no formula \eqref{eq:108.ii} (p.\ \ref{p:428}) relating the two
-- which makes me feel a little silly! I'll have to come back upon
this later. At present, let's dualize what we got, assuming now that
\scrN{} is a $k$-additive abelian category stable under small
\emph{direct} limits, and such that \emph{a small direct sum of
  monomorphisms in \scrN{} is again a monomorphism}. We then get a
pairing
\begin{equation}
  \label{eq:108.117}
  (F_\bullet,L'_\bullet) \mapsto F_\bullet \Last_k L'_\bullet :
  \D^-(\AhatN)\times\D^-(\Bhatk) \to \D^-(\scrN),\tag{117}
\end{equation}
defined by the formula (dual to \eqref{eq:108.112})
\begin{equation}
  \label{eq:108.118}
  F_\bullet \Last_k L'_\bullet \simeq F_\bullet *_k \scrL'_\bullet,\tag{118}
\end{equation}
where in the second member $\scrL'_\bullet$ is a projective resolution
of $L_\bullet$ in \Bhatk{} ($B=A\op$ being of course the dual category
of $A$), and the $*_k$ denotes the simple complex associated to the
double complex obtained by applying $*_k$ componentwise. Using the
composite equivalence
\begin{multline}
  \label{eq:108.119}
  F_\bullet\mapsto(F_\bullet)\op: (\D^-(\AhatN))\op\equeq
  \D^+((\AhatN)\op)\equeq \D^+(\BhatM),\\
  \text{ with $\scrM=\scrN\op$,}\tag{119}
\end{multline}
we get the tautological duality isomorphism
\begin{equation}
  \label{eq:108.120}
  \bigl(F_\bullet\Last_k L'_\bullet)\op \simeq \RHom_k(L'_\bullet,
  (F_\bullet)\op),\tag{120} 
\end{equation}
where the expression $\Last_k$ in the first member is relative to the
pair $(A,\scrN)$, whereas the expression $\RHom_k$ in the second is
relative to the dual pair $(B,\scrM)$. Symmetrically, we get
\begin{equation}
  \label{eq:108.120prime}
  \bigl(\RHom_k(L_\bullet,F^\bullet)\bigr)\op \simeq L_\bullet \Last_k
  (F_\bullet)\op,\tag{120'} 
\end{equation}
which is essentially the inverse isomorphism of \eqref{eq:108.120},
but for the pair $(B,\scrM)$ instead of $(A,\scrN)$.

We still should dualize the functor $\RGamma_A$ \eqref{eq:108.116}
(defined by \eqref{eq:108.114}), which we do, recalling that
$\Gamma_A$ is just the inverse limit functor $\varprojlim_{A\op}$,
which is dual to the direct limit functor $\varinjlim_{A\op}$, thus,
``integration'' of \scrN-valued presheaves is just (at least morally)
the total left derived functor of the latter, and may be denoted by
\begin{equation}
  \label{eq:108.121}
  \mathrm L{\varinjlim_{A\op}},\tag{121}
\end{equation}
while using for $\RGamma_A$ the equivalent notation, dual to
\eqref{eq:108.121}
\begin{equation}
  \label{eq:108.122}
  \mathrm R{\varprojlim_{A\op}} = \RGamma_A.\tag{122}
\end{equation}
I am not wholly happy, though, with the purely algebraic flavor of
these notations, not really suggestive of the manifold geometric
intuitions\pspage{432} surrounding the familiar homology and
cohomology notations $\mathrm H_\bullet$ and $\mathrm H^\bullet$. This
flavor is at least partially preserved, it seems to me, in the
notation $\RGamma_A$ (because of the geometric intuition tied with the
sections functor), whereas there is not yet a familiar geometric
notion of a ``cosections functor''. As we would like to have the
duality symmetry reflected as perfectly as possible in the notation, I
am going to use the notations
\begin{equation}
  \label{eq:108.123}
  \left\{
    \begin{aligned}
      \RH^\bullet(A,F^\bullet) &= \mathrm
      R{\varprojlim_{A\op}}(F^\bullet)
      \;(\;=\RGamma_A(F^\bullet))
      : \D^+(\AhatM)\to\scrM \\
      \LH_\bullet(A,F_\bullet) &= \mathrm
      L{\varinjlim_{A\op}}(F_\bullet) : \D^-(\AhatN)\to\scrN.
    \end{aligned}
  \right.\tag{123}
\end{equation}
With these notations, the duality isomorphisms
(\ref{eq:108.120},\ref{eq:108.120prime}) take the form (as announced):
\begin{equation}
  \label{eq:108.124}
  \left\{
    \begin{aligned}
      (\RH^\bullet(A,F^\bullet))\op &\simeq
      \LH_\bullet(B,(F^\bullet)\op) \\
      (\LH_\bullet(A,F_\bullet))\op &\simeq
      \RH^\bullet(B,(F_\bullet)\op)\quad,
    \end{aligned}
  \right.\tag{124}
\end{equation}
where the first members are defined in terms of cohomology resp.\
homology invariants with respect to the pair $(A,\scrM)$ resp.\
$(A,\scrN)$, whereas the second members denote homology resp.\
cohomology invariants with respect to the dual pairs $(B,\scrN)$
resp.\ $(B,\scrM)$.
\begin{remarks}
  This perfect symmetry, or rather essential identity, between
  ``homology'' or ``integration'' and ``cohomology'' or
  ``cointegration'', is obtained here at the price of working with
  presheaves with values in rather general abelian categories,
  subjected to some simple exactness properties. It should be
  remembered moreover that for the time being, $\RH^\bullet$ has not
  been defined as the total right derived functor of the sections of
  inverse limits functor, therefore the notations \eqref{eq:108.121}
  and \eqref{eq:108.122} are somewhat misleading. To feel really at
  ease, we should still work out conditions that ensure that \AhatM{}
  has enough injectives and that $\RHom_k$ can be defined also using
  such resolutions -- in which case we'll expect too to have a
  satisfactory formalism for the $\RbHom_k$ functor.
\end{remarks}

\bigbreak
\presectionfill\ondate{23.8.}\pspage{433}\par

\hangsection[Review (6): A further step in linearization: coalgebra
\dots]{Review \texorpdfstring{\textup{(6)}}{(6)}: A further step in
  linearization: coalgebra structures
  \texorpdfstring{$P\to P\otimes_k P$ in \Cat}{P->PkP in
    (Cat)}.}\label{sec:109}%
I still did a little scratchwork last night, about the question of
existence of enough injectives or projectives in a category
\[\AhatN = \bHom(A\op,\scrN),\]
where \scrN{} is a $k$-additive abelian category. Introducing the
small $k$-additive category
\[P=\Add_k(A),\]
and remembering the canonical equivalence
\[\AhatN=\bHom(A\op,\scrN)\equeq \bHom_k(P\op,\scrN),\]
the question just stated may be viewed as a particular case of the
same question for a category of the type
\begin{equation}
  \label{eq:109.125}
  \PampN \eqdef \bHom_k(P\op,\scrN),\tag{125}
\end{equation}
where now $P$ is \emph{any} small $k$-additive category. (Compare with
the reflections on pages \ref{p:403}, \ref{p:404}.) It is immediate
that in \PampN{} exist all types of (direct or inverse) limits
which exist in \scrN, and they are computed ``componentwise'' for each
argument $a$ in $P$ -- from this follows that if \scrN{} is abelian,
so is \PampN.
\addtocounter{propositionnum}{3}
\begin{propositionnum}\label{prop:109.4}
  Assume the $k$-additive category \scrN{} is stable under small
  direct limits, and is abelian, and that any object of \scrN{} is
  isomorphic to the quotient of a projective object. Then the same
  holds for \PampN. Assume moreover that any projective object $x$ of
  \scrN{} is $k$-flat, i.e., the functor
  \[U\mapsto U\otimes_k x : \AbOf_k\to\scrN\]
  is exact, i.e., transform monomorphisms into monomorphisms. Then for
  any projective object $F$ in \PampN, the functor
  \[L'\mapsto F *_k L' : Q\supamp\to\scrN\quad (\text{where
    $Q=P\op$})\]
  is exact, i.e.\ \textup(as it is known to commute to small direct
  limits\textup), it transforms monomorphisms into monomorphisms.
\end{propositionnum}
\begin{comments}
  Here, the operation $*_k$ (similar to a tensor product) is defined
  as in the non-additive set-up (with \PampN, $Q\supamp$ being
  replaced by \AhatN, \Bhatk) reviewed in section \ref{sec:107}, and
  follows from the canonical equivalence of categories
  \begin{equation}
    \label{eq:109.star}
    \PampN \fromequ \bHom_{k!}(Q\supamp, \scrN)\tag{*}
  \end{equation}
  (this\pspage{434} is formula \eqref{eq:105.starstar} of page
  \ref{p:403} with $P,\scrN$ replaced by $Q,\scrM$). It should be
  noted that the assumptions made in prop.\ \ref{prop:109.4} are the
  weakest possible for the conclusions to hold (for any $k$-additive
  small category $P$), as these conclusions, in case $P =$ final
  category, just reduce to the assumptions.
\end{comments}

Here is the outline of a proof of prop.\ \ref{prop:109.4}. Using only
stability of \scrN{} under small direct limits (besides
$k$-additivity) we define a canonical $k$-biadditive pairing
\begin{equation}
  \label{eq:109.126}
  P\supamp \times\scrN \to \PampN ,
  \quad
  (L,x)\mapsto L\otimes x\eqdef
  (a\mapsto L(a)\otimes_k x)\tag{126}
\end{equation}
(NB\enspace I recall that $P\supamp$ is defined as
\[P\supamp = \bHom_{\bZ}(P\op,\AbOf)\fromequ\bHom_k(P\op,\AbOf_k),\]
here we interpret an object of $P\supamp$ is a $k$-additive functor
\[L:P\op\to\AbOf_k\;(\;\eqdef\kMod)\quad\text{.)}\]
The relevant fact here for objects of \PampN{} of the type $L\otimes_k
x$ is
\begin{equation}
  \label{eq:109.127}
  \begin{aligned}
    \Hom_\PampN(L\otimes_kx,F) &\simeq \Hom_\scrN(x,\Hom_k(L,F)) \\
    &\simeq \Hom_{P\supamp}(L, \Hom(x,F)),
  \end{aligned}\tag{127}
\end{equation}
where in the second term,
\[\Hom_k(L,F)\in\Ob\scrN\]
is defined in a way dual to $F*_k L'$ (cf.\ comments above), using the
equivalence (dual to \eqref{eq:109.star} above)
\begin{equation}
  \label{eq:109.starprime}
  \PampN\simeq\bHom_k^!((P\supamp)\op,\scrN),\tag{*'}
\end{equation}
which is defined only, however, when \scrN{} is stable under small
inverse limits (hence the first isomorphism in \eqref{eq:109.127}
makes sense only under this extra assumption); on the other hand, in
the third term in \eqref{eq:109.127}
\[\Hom(x,F) \eqdef (a\mapsto\Hom_\scrN(x,F(a))\quad\text{in
  $P\supamp$,}\]
and the isomorphism between the first and third term in
\eqref{eq:109.127} makes sense and is defined without any extra
assumption on \scrN.

We leave to the reader to check \eqref{eq:109.127} (where one is
readily reduced to the case when $L$ is an object $a$ in $P$, using
the commutation of the three functors obtained
$(P\supamp)\op\to\AbOf_k$ with small inverse limits). It follows, when
\scrN{} is abelian:
\begin{equation}
  \label{eq:109.128}
  \text{$L$ projective in $P\supamp$, $x$ projective in \scrN}
  \Rightarrow
  \text{$L\otimes_kx$ proj.\ in \PampN.}\tag{128}
\end{equation}

Assume now that any object of \scrN{} is quotient of a projective
one, and let $F$ be any object in \PampN. Formula \eqref{eq:109.127}
for $L=a$ in $P$ reduces\pspage{435} to the down-to-earth formula
\begin{equation}
  \label{eq:109.127prime}
  \Hom_\PampN(a\otimes_kx,F)\simeq\Hom(x,F(a)).\tag{127'}
\end{equation}
Now, let for any $a$ in $P$
\[x_a\to F(a)\]
be an epimorphism in \scrN, with $x_a$ projective. From
\eqref{eq:109.127prime} we get a map
\[a\otimes_kx_a\to F\]
in \PampN, hence a map
\[\scrF = \bigoplus_{\text{$a$ in $P$}} a\otimes_kx_a \to F,\]
it is easily seen that this is epimorphic (because the maps $x_a\to
F(a)$ are), and \scrF{} is projective as a direct sum of projective
objects. This proves the first statement in prop.\
\ref{prop:109.4}. For the second statement, we'll use the formula
\begin{equation}
  \label{eq:109.129}
  (L\otimes_k) *_k L' \simeq (L *_k L') *_k x\tag{129}
\end{equation}
for $L$ in $P\supamp$, $L'$ in $Q\supamp$, and $x$ in \scrN{} -- for
the proof, we may reduce to the case when $L$ is in $P$, $L'$ in
$Q=P\op$, say $L=a$ and $L'=b\op$, in which case both members identify
with $\Hom_P(b,a)\otimes_k x$. To prove that for $F$ projective,
$L'\mapsto F*_kL'$ takes monomorphisms into monomorphisms, using that
$F$ is a direct factor of objects of the type $a\otimes_k x$ with $a$
in $P$ and $x$ in \scrN, we are reduced to the case $F=a\otimes_kx$,
in which case by \eqref{eq:109.129} the functor reduces to
\[L'\mapsto L'(a)\otimes_kx,\]
which is again exact by the assumption that any projective object in
\scrN{} (and hence $x$) is $k$-flat.
\begin{remark}
  It is not automatic that a projective object in a $k$-additive
  abelian category be $k$-flat -- take for instance $k=\bZ$ and
  $\scrN=\AbOf_{\bF_p}$, where $\bF_p$ is a finite prime field, then
  all objects in \scrN{} are projective, whereas only the zero objects
  are \bZ-flat.
\end{remark}

We leave to the reader to write down the dual statement of prop.\
\ref{prop:109.4}, concerning injectives in a category \PampM, where
now \scrM{} is a $k$-additive abelian category stable under small
inverse limits, and possessing sufficiently many injectives (hence the
same holds in \PampM), and assuming eventually that these injectives $x$
are ``$k$-coflat'', i.e.,
\[U\mapsto\Hom_k(U,x) : \AbOf_k\op\to\scrM\]
is exact (i.e., transforms monomorphisms in $\AbOf_k$ into
epimorphisms in\pspage{436} \scrM), which implies that for $F$
injective in $P\supamp$, the functor
\[L\mapsto\Hom_k(L,F) : P\supamp\to\scrM\]
is exact.

To sum up, we get the
\begin{corollarynum}\label{cor:109.prop4.1}
  Let $P$ be any small $k$-additive category, and let \scrM{} be a
  $k$-additive category which satisfies the following assumptions:
  \begin{enumerate}[label=\alph*),font=\normalfont]
  \item\label{it:109.cor1.a}
    \scrM{} is abelian, and stable under small inverse limits,
  \item\label{it:109.cor1.b}
    \scrM{} has ``sufficiently many injectives'',
  \item\label{it:109.cor1.c}
    injective objects of \scrM{} are $k$-coflat,
  \item\label{it:109.cor1.d}
    any product \textup(with small indexing family\textup) of
    epimorphisms in \scrM{} is an epimorphism.
  \end{enumerate}
  Consider the $k$-biadditive pairing
  \begin{equation}
    \label{eq:109.130}
    (L,F)\mapsto\Hom_k(L,F) : P\supamp \times\PampM\to\scrM.\tag{130}
  \end{equation}
  This pairing admits a total right derived functor
  \begin{equation}
    \label{eq:109.131}
    (L_\bullet,F^\bullet)\mapsto\RHom_k(L_\bullet,F^\bullet) :
    \D^-(P\supamp) \times \D^+(\PampM) \to \D^+(\scrM),\tag{131}
  \end{equation}
  which can be computed using either projective resolutions of
  $L_\bullet$, or injective resolutions of $F^\bullet$, or both
  simultaneously.
\end{corollarynum}

We have a dual statement, concerning the pairing
\begin{equation}
  \label{eq:109.130prime}
  (F,L')\mapsto F*_kL': \PampN\times Q\supamp \to\scrN,
  \quad\text{with $Q=P\op$,}\tag{130'}
\end{equation}
giving rise to a total left derived functor
\begin{equation}
  \label{eq:109.131prime}
  (F_\bullet,L_\bullet') \mapsto F_\bullet\Last_k L_\bullet' :
  \D^-(\PampN) \times \D^-(Q\supamp) \to \D^-(\scrN),\tag{131'}
\end{equation}
using projective resolutions of either $F_\bullet$, or $L_\bullet'$,
or both. Here, \scrN{} is a $k$-additive category satisfying the
properties dual to \ref{it:109.cor1.a} to \ref{it:109.cor1.d} above,
i.e., such that $\scrM=\scrN\op$ satisfies the properties stated in
the corollary. We have the evident duality relations between the two
kinds of operations $\RHom_k$ and $\Last_k$, embodied by the formulæ
\eqref{eq:108.120} and \eqref{eq:108.120prime} of page \ref{p:431},
where the categories \Ahatk, \AhatM, \AhatN, etc.\ are replaced by
$P\supamp$, \PampM, \PampN, etc.\ (where the etc.'s refer to
replacement of $A$ by $B=A\op$ and of $P$ by $Q=P\op$).

Next question is to extend the $\RHom_k$ formalism to a $\RbHom_k$
formalism (and similarly from $\Last_k$ to $\Loast_k$), as envisioned
yesterday. To do so, in the wholly $k$-additive set-up we are now
working in, we still\pspage{437} need (in case of $\RbHom_k$ an
(``interior'') tensor product structure on $P\supamp$ (and dually for
$\Loast$, requiring a tensor product structure on $Q\supamp$), so as
to give rise to a $k$-biadditive pairing
\begin{equation}
  \label{eq:109.132}
  (L,F)\mapsto\bHom_k(L,F): P\supamp\times\PampM\to\PampM,\tag{132}
\end{equation}
(and dually,
\begin{equation}
  \label{eq:109.132prime}
  (F,L')\mapsto F\oast_kL' : \PampN\times Q\supamp\to\PampN
  \text{ ),}\tag{132'}
\end{equation}
for which we want to take the total right derived functor $\RbHom_k$
(or dually, the total left derived functor $\Loast_k$). As we say
yesterday, though, taking projective resolutions of $L$ will not do
(except in very special cases, such as $P=\Add_k(A)$ with $A$ in
\Cat{} stable under binary products), because for $L$ projective in
$P\supamp$, the functor
\[F\mapsto\bHom_k(L,F):\PampM\to\PampM\]
has little chance to be exact. Taking injective resolutions of
$F_\bullet$, however, we expected \emph{would} do in ``reasonable''
cases -- here the question is whether for $F$ injective, the functor
\[L\mapsto\bHom_k(L,F) \eqdef (a\mapsto\Hom_k(L\otimes_ka,F)) :
(P\supamp)\op\to\scrM\]
is exact (where $L\otimes_ka$ denotes the given tensor product within
$P\supamp$), i.e., whether the functors (for $F$ in \PampM, $a$ in
$P$)
\begin{equation}
  \label{eq:109.starbis}
  L\mapsto\Hom_k(L\otimes_ka,F) : (P\supamp)\op\to\scrM\tag{*}
\end{equation}
transform monomorphisms of $P\supamp$ into epimorphisms of \scrM. When
\scrM{} satisfies the conditions dual to
\ref{it:109.cor1.a}\ref{it:109.cor1.b}\ref{it:109.cor1.c} in the
corollary, hence $L\mapsto\Hom_k(L,F)$ is exact when $F$ in \PampM{}
is injective, exactness of \eqref{eq:109.starbis} above will follow
from exactness of
\[L\mapsto L\otimes_ka : P\supamp\to P\supamp.\]
Now, this latter exactness holds in the case we are interesting in
mainly, when $P=\Add_k(A)$ and hence
$P\supamp\equeq\bHom(A\op,\AbOf_k)$ and when tensor product for
presheaves on $A$ is defined as usual, componentwise -- then the
objects of $\Add_k(A)=P$ correspond to $k$-flat presheaves, hence
tensor product by these is exact. Thus:
\begin{corollarynum}\label{cor:109.prop4.2}
  Assume \scrM{} satisfies the conditions of corollary
  \ref{cor:109.prop4.1} above, let $A$ be any small category and
  $P=\Add_k(A)$. Then the pairing
  \[\bHom_k:\Ahatk\times\AhatM\to\AhatM\]
  \textup(cf.\ \eqref{eq:108.109prime} page \ref{p:428}\textup) admits
  a total right derived functor $\RbHom_k$
  \eqref{eq:108.111},\pspage{438} which may be computed using
  injective resolutions of the second argument $F^\bullet$ in
  $\RbHom_k(L_\bullet,F^\bullet)$ \textup(but not, in general, by
  using projective resolutions of $L_\bullet$\textup), i.e.,
  \begin{equation}
    \label{eq:109.133}
    \RbHom_k(L_\bullet,F^\bullet) \simeq
    \bHom_k^{\bullet\bullet}(L_\bullet,\scrF^\bullet),\tag{133}
  \end{equation}
  where $\scrF^\bullet$ is an injective resolution of $F^\bullet$
  \textup(i.e., a complex in \AhatM{} with degrees bounded from below
  and injective components, endowed with a quasi-isomorphism
  $F^\bullet\simeq\scrF^\bullet$\textup).
\end{corollarynum}

Applying now $\RGamma_A$ to both members of \eqref{eq:109.133}, and
using the similar isomorphism for $\RHom_k$ (valid by cor.\
\ref{cor:109.prop4.2}, we get the familiar formula
\begin{equation}
  \label{eq:109.134}
  \RHom_k(L_\bullet,F^\bullet)\simeq\RH^\bullet(\RbHom_k(L_\bullet,F^\bullet)),\tag{134}
\end{equation}
(where in accordance with \eqref{eq:108.124}, we wrote $\RH^\bullet$
instead of $\RGamma_A$), where however $F^\bullet$ is now a complex of
presheaves with values in a $k$-additive category \scrM{} (satisfying
the conditions dual to \ref{it:109.cor1.a} to \ref{it:109.cor1.d} in
cor.\ \ref{cor:109.prop4.1}), not just a presheaf of $k$-modules.

Replacing \scrM{} by a category \scrN{} satisfying the assumptions of
cor.\ \ref{cor:109.prop4.1}, we get likewise a total left derived
functor
\begin{equation}
  \label{eq:109.135}
  (F_\bullet,L_\bullet')\mapsto F\Loast_k L_\bullet' : \D^-(\AhatN)
  \times \D^-(\Bhatk) \to \D^-(\AhatN),\tag{135}
\end{equation}
which can be defined using projective resolutions of the first
argument $F_\bullet$ (but not by using projective resolutions of
$L_\bullet'$), giving rise to the isomorphism (dual to
\eqref{eq:109.134})
\begin{equation}
  \label{eq:109.134prime}
  F_\bullet \Last_k L_\bullet' \simeq \LH_\bullet(F_\bullet \Loast_k
  L_\bullet'). \tag{134'}
\end{equation}
The duality relationship between $\RbHom_k$ and $\Loast_k$ can be
expressed by two obvious formulæ, similar to \eqref{eq:108.120} and
\eqref{eq:108.120prime} for $\RHom_k$ and $\Last_k$, which we leave to
the reader.

\begin{remarks}
  \namedlabel{rem:109.1}{1)}\enspace If we want to consider
  ``multiplicative structure'' in the purely ``$k$-linear'' set-up,
  where the data is a small $k$-additive category $P$, rather than a
  small category $A$ (giving rise to $P=\Add_k(A)$), and corresponding
  cohomology and homology operations $\RbHom_k$ and $\Loast_k$, not
  only $\RHom_k$ and $\Last_k$, the natural thing to do, it seems, is
  to introduce a ``\emph{diagonal map}''
  \begin{equation}
    \label{eq:109.starbisbis}
    P \to P\otimes_k P,\tag{*}
  \end{equation}
  where the tensor product in the second member can be defined in a
  rather evident way (as solution of the obvious $2$-universal problem
  in terms of $k$-biadditive functors on $P\times P$), which will give
  rise in the ``usual'' way to a tensor product operation in both
  $P\supamp$ and $Q\supamp$ (where\pspage{439} $Q=P\op$). We'll come
  back upon this later, I expect. This structure
  \eqref{eq:109.starbisbis} will be the $k$-linear analogon of the
  usual diagonal map
  \begin{equation}
    \label{eq:109.starstar}
    A\to A\times A\tag{**}
  \end{equation}
  for a small category $A$, giving rise by $k$-linearization to
  \[\Add_k(A)\to\Add_k(A\times A)\equeq \Add_k(A)\otimes_k\Add_k(A),\]
  namely a structure of type \eqref{eq:109.starbisbis}. It just
  occurred to me, through the reflections of these last days, that the
  ``coalgebra structure'' \eqref{eq:109.starbisbis} may well turn out
  (taking $k=\bZ$) to be the more sophisticated structure than a usual
  coalgebra structure (cf.\ p.\ \ref{p:339} \eqref{eq:94.27}), needed
  in order to grasp ``in linear terms'' the notion of a homotopy type,
  possibly under restrictions such as $1$-connectedness, as pondered
  about in section \ref{sec:94}. This looks at any rate a more
  ``natural'' object than the De~Rham complex with divided powers,
  referred to in loc.\ cit., and is more evidently adapted to our
  point of view of using small categories as models for homotopy
  types. The greater sophistication, in comparison to De~Rham type
  complexes, lies in this, that here the objects serving as models
  (whether small categories, or small additive categories endowed with
  a diagonal map) are objects in a \emph{$2$-category}, whereas
  De~Rham complexes and the like are just objects in ordinary
  categories, without any question of taking ``maps between
  maps''. This feature implies, ``as usual'' (or in duality rather to
  familiar situations with tensor product functors\ldots) that the
  (anti)commutative and associative \emph{axioms} familiar from linear
  algebra (in the case of usual $k$-algebras or $k$-coalgebras),
  should be replaced by commutativity and associativity \emph{data},
  namely given isomorphisms (not identities) between two natural
  functors
  \[P\to P\otimes_k P, \quad P\to P\otimes_k P \otimes_k P\]
  deduced from \eqref{eq:109.starbisbis}. The axioms now will be more
  sophisticated, they will express ``coherence conditions'' on these
  data -- one place maybe where this is developed somewhat, in the
  context of diagonal maps \eqref{eq:109.starbisbis}, might be
  Saavedra's thesis.\scrcomment{\textcite{Saavedra1972}; see also in
    particular \textcite{DeligneMilne1982,Deligne1990,Deligne2002}}
  (The more familiar case, when starting with a tensor product
  operation on a category, together with associativity and/or
  commutativity data, has been done with care by various
  mathematicians, including Mac~Lane, Bénabou, Mme~Sinh Hoang Xuan,
  and presumably it should be enough to ``reverse arrows'' in order to
  get ``the'' natural set of coherence axioms for a diagonal map
  \eqref{eq:109.starbisbis}). The ``intriguing feature'' with the
  would-be De~Rham models for homotopy types (cf.\ p.\ \ref{p:341},
  \ref{p:342}), namely that the latter make sense over any\pspage{440}
  commutative ground ring $k$, not only \bZ, with corresponding notion
  of ring extension $k\to k'$, carries over to structures of the type
  \eqref{eq:109.starbisbis}. Indeed, for any $k$-additive category
  $P$, it is easy to define a $k'$-additive category
  \[P\otimes_kk',\quad \
  \text{for given homomorphism $k\to k'$,}\]
  for instance as the solution of the obvious $2$-universal problem
  corresponding to mapping $P$ $k$-additively into $k'$-additive
  categories, or more evidently by taking the same objects as for $P$,
  but with
  \[\Hom_{P'}(a,b) = \Hom_P(a,b)\otimes_k k'.\]
  Thus, any ``coalgebra structure in \Cat'' \eqref{eq:109.starbisbis}
  over the ground ring $k$, gives rise to a similar structure over
  ground ring $k'$.

  Of course, among the relevant axioms for the diagonal functor
  \eqref{eq:109.starbisbis}, is the existence of unit objects in
  $P\supamp$ and $Q\supamp$, which may be viewed equally as
  $k$-additive functors (defined up to unique isomorphism)
  \[P\op\to\AbOf_k, \quad P\to\AbOf_k,\]
  playing the role I would think of ``augmentation'' and
  ``coaugmentation'' in the more familiar set-up of ordinary
  coalgebras. Denoting these objects by $k_P$ and $k_Q$ respectively
  (in analogy to the constant presheaves $k_A$ and $k_B$ on $A$ and
  $B$), $\RHom(k_P,{-})$ now allows expression of cohomology or
  cointegration, and ${-}\Last_k(k_Q)$ allows expression of homology
  or integration (for complexes in \PampM{} say). ``Constant
  coefficients'' on $P$, i.e., in $P\supamp$ may now be defined, as
  objects in $P\supamp$ of the type
  \[U\otimes_k k_P,\]
  where $U$ is in $\AbOf_k$, i.e., is any $k$-module, and hence we get
  homology and cohomology invariants with coefficients in any such $U$
  (or complexes of such), and surely too cup and cap products\ldots
  Also, quasi-isomorphisms of structure \eqref{eq:109.starbisbis}
  (with units) can now be defined in an evident way, hence a derived
  category which merits to be understood, when $k=\bZ$, in terms of
  the homotopy category \Hot. I wouldn't expect of course that for any
  small category $A$, the abelianization $\Add(A)$ together with its
  diagonal map allows to recover the homotopy type, unless $A$ is
  $1$-connected. As was the case visibly for De~Rham complexes, if we
  hope to recover general homotopy types (not only $1$-connected
  ones), we should work with slightly more sophisticated structures
  still,\pspage{441} involving a group (or better still, a groupoid)
  and an operation of it on a structure of type
  \eqref{eq:109.starbisbis} (embodying a universal covering\ldots).

  Here I am getting, though, into thin air again, and I don't expect
  I'll ponder much more in this direction and see what comes out. The
  striking fact, however, here, is that quite unexpectedly, we get
  further hold and food for this thin-air intuition (which came up
  first in relation to De~Rham structures with divided powers), that
  there may be a reasonable (and essentially just one such) notion of
  a ``homotopy type over the ground ring $k$'' for any commutative
  ring $k$, reducing for $k=\bZ$ to usual homotopy types, and giving
  rise to base change functors
  \[\HotOf(k)\to\HotOf(k')\]
  for any ring homomorphism $k\to k'$. And I wonder whether this might
  not come out in some very simplistic way, in the general spirit of
  our ``modelizing story'', without having to work out in full a
  description of homotopy types by such sophisticated models as
  De~Rham complexes with divided powers, or coalgebra structures in
  \Cat, and looking up maybe the relations between these. (How by all
  means hope to recover a De~Rham structure from a stupid structure
  \eqref{eq:109.starbisbis}???)
  
  \namedlabel{rem:109.2}{2)}\enspace The condition \ref{it:109.cor1.d}
  in corollary \ref{cor:109.prop4.1} is needed in order to ensure that
  a derived functor $\RHom_k(L_\bullet,F^\bullet)$ may be defined
  using \emph{projective} resolutions of $L_\bullet$, whereas
  conditions \ref{it:109.cor1.b}, \ref{it:109.cor1.c} ensure that a
  functor $\RHom_k$ may be defined using \emph{injective} resolutions
  of $F^\bullet$. It is a well-known standard fact of homological
  algebra that in case both methods work (namely here, when all four
  assumptions are satisfied) that the two methods yield the same
  result, which may equally be described by resolving simultaneously
  the two arguments. (NB\enspace condition \ref{it:109.cor1.a} is
  needed anyhow for $\Hom_k$ to be defined and for \PampM{} being an
  abelian category, which allows to define $\D^+(\PampM)$.) Our
  preference goes to the first method, which in case $P=\Add_k(A)$ and
  $L=k_A$, conduces to computation of cohomology
  $\RGamma_A(F^\bullet)$ in terms of a ``cointegrator'' $L_\bullet^A$
  on $A$. However, when it comes to introducing the variant
  $\RbHom_k$, this method breaks down, as we saw, it is the other one
  which works. We thus get a satisfactory formalism of $\RHom_k$ and
  $\RbHom_k$ (including formula \eqref{eq:109.134} relating them via
  $\RH^\bullet$) using only assumptions
  \ref{it:109.cor1.a}\ref{it:109.cor1.b}\ref{it:109.cor1.c}. 

  \namedlabel{rem:109.3}{3)}\enspace If we want to extend the
  $\RbHom_k$ formalism to the set-up when the data $A$ is replaced by
  a $k$-additive category $P$ endowed with a diagonal map as in remark
  \ref{rem:109.1}, the proof on page \ref{p:437} shows that what is
  needed\pspage{442} is exactness of the functor $L\mapsto L\otimes a$
  from $P\supamp$ to $P\supamp$, for any $a$ in $P$ -- which is a
  ``flatness'' condition on $a$. It is easily checked that his
  condition is satisfied provided $\Hom_P(b,a)$ is a flat $k$-module,
  for any $b$ in $P$ (more generally, an object $M$ in $P\supamp$ is
  flat for the tensor product structure in $P\supamp$, provided $M(b)$
  is flat for any $b$ in $P$). Thus, there will be a satisfactory
  $\RbHom_k$ theory provided the $k$-modules $\Hom_P(b,a)$, for $a,b$
  in $P$, are $k$-flat. It is immediately checked that this also means
  that any projective object $L$ in $P\supamp$ is ``$k$-flat'' (with
  respect to external tensor product $U\mapsto U\otimes_k L :
  \AbOf_k\to P\supamp$, as in prop.\ \ref{prop:109.4}), or
  equivalently still, that this holds when $L$ is any object $a$ in
  $P$. In case $P=\Add_k(A)$, when any object of $P$ is a finite sum
  of objects $k^{(x)}$ with $x$ in $A$, the modules $\Hom(b,a)$ for
  $b,a$ in $L$ are finite sums of modules of the type
  $\Hom(k^{(y)},k^{(x)})=k^{(\Hom(y,x))}$, and hence are projective,
  not only flat. It would seem that in the general case of $P$ endowed
  with a diagonal map, the ``natural'' assumption to make in order to
  have everything come out just as nicely as when $P$ comes from an
  $A$, is that the $k$-modules $\Hom_P(b,a)$ (for $a,b$ in $P$) should
  be projective, not only flat. Flatness however seems to be all that
  is needed in order to ensure that when $L$ in $P\supamp$ is
  projective (hence $L(a)$ is flat for any $a$ in $P$) and $F$ in
  \PampM{} is injective, then the objects $\Hom_k(L,F)$ and
  $\bHom_k(L,F)$ in \scrM{} and \PampM{} respectively are
  injective. This implies that for a $k$-additive functor
  \[u:\scrM\to\scrM'\]
  between categories \scrM, $\scrM'$ satisfying the assumptions of
  cor.\ \ref{cor:109.prop4.1}, and $u$ commuting moreover to small
  inverse limits (and hence to formation of $\Hom_k(L,F)$), we get a
  canonical isomorphism
  \begin{equation}
    \label{eq:109.136}
    \mathrm R u (\RHom_k(L_\bullet,F^\bullet)) \simeq
    \RHom_k(L_\bullet,\mathrm Ru^P(F^\bullet)),\tag{136}
  \end{equation}
  where
  \[u^P:\PampM\to P\supamp_{\scrM'}\]
  denotes the extension of $u$, and $\mathrm Ru$, $\mathrm Ru^P$ are
  the right derived functors. When $u$ is exact, we may replace
  $\mathrm Ru$, $\mathrm Ru^P$ by $u$, $u^P$ (applied componentwise to
  complexes, without any need to take an injective resolution
  first). There is a formula as \eqref{eq:109.136} with $\RHom_k$
  replaced by $\RbHom_k$, which I skip, as well as the dual formulas.
\end{remarks}



\chapter{Schematization}
\label{ch:VI}

\presectionfill\ondate{24.8.}\pspage{443}\par

\hangsection[More wishful thinking on ``schematization'' of homotopy
\dots]{More wishful thinking on ``schematization'' of homotopy
  types.}\label{sec:110}%
I pondered some more about homotopy types over a ground ring $k$, just
enough to become familiar again with the idea, and more or less
convinced that that there should exist such a thing, which should
amount, kind of, to putting a ``continuous'' structures (namely the
very rich structure of a scheme) upon something usually visualized as
something ``discrete'' -- namely a homotopy type. The basic analogy
here is free \bZ-modules $M$ of finite type -- a typical case of a
``\emph{discrete}'' structure. It gives rise, though, to a
\emph{vector bundle} $W(M)$ over the absolute base $S_0=\Spec(\bZ)$,
whose \bZ-module of sections is $M$, and the functor
\[M\mapsto W(M)\]
from free \bZ-modules of finite type to vector bundles over $S_0$ is
fully faithful. When $M$ is an arbitrary \bZ-module, i.e., an abelian
group, $W(M)$ still makes sense, namely as a functor
\[k\mapsto M\otimes_\bZ k\]
on the category of all commutative \bZ-algebra (i.e., just commutative
rings); it is no longer representable by a scheme over $S_0$ (except
precisely when $M$ is free of finite type), but it is very close
still, intuitively and technically too, to a usual vector bundle (the
``vector'' structure coming from the $k$-module structure on
$W(M)(k)=M\otimes_\bZ k$). Again, the functor $M\mapsto W(M)$ from
$\AbOf$ to the category of ``generalized vector bundles'' over $S_0$
is fully faithful. Working with semisimplicial \bZ-modules (say)
rather than just \bZ-modules, and more specifically with those
corresponding to $K(\pi,n)$ types, and using Postnikov
``dévissage''\scrcomment{``dévissage'' = ``decomposition''}
of a general homotopy type, one may hope to ``represent'' this type,
in a more or less canonical way in terms of the successive
semisimplicial Postnikov fibrations, by a semi-simplicial object in
the topos (say) of all functors from \bZ-algebras to sets which are
``sheaves'' for a suitable site structure on the dual category (namely
the category of affine schemes over \bZ) -- the so-called ``flat''
topology seems OK. (NB\enspace To eliminate logical difficulties, we
may have to restrict somewhat the rings $k$ used as arguments, for
instance take them to be of finite type over \bZ{} -- never mind such
technicalities!). This approach may possibly work, when restricting to
$1$-connected homotopy types, or at any rate to the case when the
fundamental groups of the connected components are abelian. If we wish
a ``schematization'' of arbitrary homotopy types, we may think of
going about it by keeping the fundamental group of groupoid
``discrete'', and ``schematize'' the Postnikov\pspage{444} truncation
involving only the higher homotopy groups ($\pi_i$ with $i\ge2$). This
suggests that for any integer $n\ge1$, there may be a ``schematization
above level $n$'' for a given homotopy type, leaving the Cartan-Serre
truncation of level $n$ discrete, and ``schematizing'' the
corresponding total Postnikov fiber (involving homotopy groups $\pi_i$
for $i\ge n$). Maybe such a schematization can be constructed equally
for level $0$, even without assuming the fundamental groupoid to be
abelian, only nilpotent -- but then we may have to change from ground
ring \bZ{} to the considerably coarser one \bQ{} (compare comments at
the end of section \ref{sec:94}). In any case, the key step in this
approach would consist in checking that, after schematization has been
carried through successfully up to a certain level in the successive
elementary Postnikov fibrations, the next elementary fibration
(described by a cohomology class in $\mathrm H^{n+2}(X_n,\pi_{n+1})$)
comes from a ``schematic'' one, and that the latter is essentially
unique; in other words, that the canonical map from the ``schematic''
$\mathrm H^{n+2}$ (with ``quasi-coherent'' coefficients) to the usual
``discrete'' one is bijective. Maybe this hope is wholly unrealistic
though. One fact which calls for some skepticism about this approach,
comes in when looking at the case of an ``abelian'' homotopy type,
described by a semisimplicial abelian group $X_\bullet$, in which case
we expect that base change $\bZ\to k$ should be just the usual base
change
\[X_\bullet \mapsto X_\bullet \otimes_\bZ k\]
(if $X_\bullet$ has torsion-free, i.e., flat components, at any
rate). But when $X_\bullet$ has homology torsion, the universal
coefficients formula shows us that the homology (= homotopy) groups of
$X_\bullet\otimes_\bZ k$ are \emph{not} just the groups $\mathrm
H_i(X_\bullet) \otimes_\bZ k=\pi_i(X_\bullet)\otimes_\bZ k$, as we
implicitly were assuming it seems in the approach sketched above, when
schematizing the homotopy groups $\pi_i$ one by one via
$W(\pi_i)$. Thus, maybe Postnikov dévissage isn't a possible approach
towards schematization of homotopy types, and one will have to work
out rather a comprehensive yoga of reconstructing a homotopy type from
one kind or other of ``abelianization'' or ``linearization'' of
homotopy types, endowed with suitable extra structure embodying
``multiplicative'' features of the homology and cohomology
structure. At any rate, I did not hit upon any ``simplistic way'' to
define homotopy types over any ground ring, and I have some doubts
there is any, in terms of the general non-sense we did so far.

Besides\pspage{445} this, I spent hours to try and put some order into
the mess of all $\Hom$ and tensor product type operations between
categories \Ahatk, \AhatM, \Bhatk, \BhatM{} (or their $\&$-style
generalizations), and the duality and Cartan-type isomorphisms between
these. There are a few more still than the fair bunch met with in
these notes so far -- I finally renounced to get really through and
work out a wholly satisfactory set of notations, taking into account
all symmetries in the situation. I realized this might well take days
of work, while at present there is no real need yet for it. I
sometimes find it difficult to find a proper balance in these notes
between the need of working on reasonably firm ground, and working out
suggestive terminology and notations for what is coming up, and on the
other hand my resolve not to get caught again by the ``Eléments de
Géométrie Algébrique'' style of work, when it was understood that
everything had to be worked out in complete detail and in greatest
generality, for the benefit of generations of
``usagers''\scrcomment{I'm leaving in ``usagers'' for color, though
  ``users'' would surely do\ldots} (besides my own till the end of my
life!\kern.35\fontdimen3\font.\kern.35\fontdimen3\font.). As a matter
of fact, this whole ``abelianization story'', going on now for well
over a hundred pages and nothing really startling coming out -- just
things I feel I should have known for ages, has been won (so to say)
over an inner reluctance against these ``digressions'' in the main
line of thought, the reluctance of one who is in a hurry to get
through. I know well this old reluctance, feeling silly whenever
working out ``trivial details'' with utmost care; as I know too that
through this work only would come to a thorough understanding of what
is going on, and new intuitions or relationships would flash up
sometimes and open up unexpected landscapes and provide fresh
impetus. The same has happened innumerable times too within the last
seven years, when ``meditating'' on personal matters -- constantly
``the-one-in-a-hurry'' has turned out to be just the servant of the
inner resistances against renewal, against a fresh, innocent look upon
things familiar, and consistently ignores as ``irrelevant''. It
doesn't seem the-one-in-a-hurry gets at all discouraged for not
getting his way many times -- he seems to be just as stubborn as the
one who likes to take his time and look up things thoroughly!

\bigbreak

\presectionfill\ondate{25.8.}\pspage{446}\par

\hangsection[Complexes of ``unipotent bundles'' as models, and
\dots]{Complexes of ``unipotent bundles'' as models, and ``schematic''
  linearization.}\label{sec:111}%
Still about ``schematization'' of homotopy types! Here is a tentative
approach, without any explicit use of Postnikov fibrations nor
abelianization, although both are involved implicitly. If $n$ is any
natural integer, I'll denote by
\[\HotOf_n\]
the full subcategory of the pointed homotopy category
$\HotOf^\bullet$, made up with $n$-connected homotopy types, with the
extra assumption for $n=0$ that the fundamental group be abelian. For
any (commutative) ring $k$, I want to define a category
\begin{equation}
  \label{eq:111.a}
  \HotOf_n(k),\tag{a}
\end{equation}
depending covariantly on $k$, in such a way that we have an
equivalence
\begin{equation}
  \label{eq:111.b}
  \HotOf_n(\bZ)\toequ\HotOf_n,\tag{b}
\end{equation}
which should come from a canonical functor
\begin{equation}
  \label{eq:111.c}
  \HotOf_n(k) \to \HotOf_n\tag{c}
\end{equation}
defined for any $k$, and which should be viewed as ground-ring
restriction from $k$ to \bZ{} -- more generally, for any ring
homomorphism
\[k'\to k\]
we expect a restriction functor, beside the ring extension functor
\begin{equation}
  \label{eq:111.d}
  \HotOf_n(k')\to\HotOf_n(k)\quad\text{and}\quad
  \HotOf_n(k)\to\HotOf_n(k').\tag{d}
\end{equation}
Among other important features to expect, is that for any object $X$
in $\HotOf_n(k)$, the homotopy groups $\pi_i(X)$ (defined via
\eqref{eq:111.c}) should be naturally endowed with structures of
$k$-modules, and the ring extension and restriction functors
\eqref{eq:111.d} should be compatible with these.

Here is an idea for getting such a theory. For given $k$, we first
define an auxiliary category
\begin{equation}
  \label{eq:111.e}
  U(k),\tag{e}
\end{equation}
whose objects may be called ``unipotent bundles over $k$''. These
``bundles'' will not be quite schemes over $k$, they will be defined
as functors
\begin{equation}
  \label{eq:111.f}
  \Alg_{/k}\to\Sets\tag{f}
\end{equation}
where $\Alg_{/k}$ is the category of (commutative) $k$-algebras (in
the basic universe \scrU). The opposite category may be identified
with the category
\[\Aff_{/k}\]
of\pspage{447} affine schemes over $k$, thus, we'll be working in the
category of functors (or presheaves over $\Aff_{/k}$)
\begin{equation}
  \label{eq:111.fprime}
  (\Aff_{/k})\op \to \Sets\tag{f'}
\end{equation}
more specifically, $U(k)$ will be a full subcategory of this category
of functors. We'll endow $\Aff_{/k}$ with one of the standard site
structures, the most convenient one here is the
fpqc\scrcomment{``fidèlement plat et quasi-compact''} topology
(faithfully flat quasi compact topology), and work in the category of
\emph{sheaves} in the latter. In terms of the interpretation
\eqref{eq:111.f} as covariant functors on $\Alg_{/k}$, this just means
that we are restricting to functors $X$ which 1)\enspace commute to
finite products and 2)\enspace are ``compatible with faithfully flat
descent'', i.e., for any map
\[k'\to k''\]
in $\Alg_{/k}$ such that $k''$ becomes a faithfully flat algebra over
$k'$, the following diagram in \Sets
\[X(k') \to X(k'') \rightrightarrows X(k'' \otimes_{k'} k'')\]
is exact. Thus, $U(k)$ will be defined as a full subcategory of the
category of such functors, or ``sheaves''.

One way for defining $U(k)$ is to present it as the union of a
sequence of subcategories $U_m(k)$ ($m$ a natural integer). We'll take
$U_0(k)$ to be just reduced to the final functors (i.e., $X(k')$ is a
one-point set for any $k'$ in $\Alg_{/k}$), and define inductively
$U_{m+1}(k)$ in terms of $U_m(k)$ as follows. For any $k$-module $M$,
let $W(M)$ be the corresponding ``vector bundle'', defined by
\begin{equation}
  \label{eq:111.g}
  W(M)(k')=M\otimes_k k',\tag{g}
\end{equation}
then an object $X$ is in $U_{m+1}(k)$ if{f} there exists an object $X_m$
in $U_m(k)$, and a $k$-module $M=M_{m+1}$, in such a way that
$X=X_{m+1}$ should be isomorphic to a ``torsor'' over $X_m$, with
group $W(M)$. I am not too sure, here, whether we should view the
objects of $U_m(k)$ as endowed with the extra structure consisting in
giving the modules $M_1,\dots,M_m$ used in the inductive construction,
and moreover the successive fibrations -- if so, then of course the
categories $U_m(k)$, and their union $U(k)$, will no longer be
interpreted as a mere subcategory of the category of\pspage{448}
sheaves (of sets) just described. Possibly, both approaches are of
interest and will yield non-equivalent notions of schematization. On
the other hand, although definitely $X$ is not representable by a
usual scheme over $k$ unless the $k$-modules $M_i$ are projective of
finite type (in which case $X$ will be an affine scheme, and even
isomorphic, at least locally over $\Spec(k)$, to standard affine space
$E_k^d$ for suitable $d$), it is felt that $X$, as far as cohomology
properties go, should be very close to being an affine scheme, and
that presumably its cohomology groups $\mathrm H^i$ with coefficients
in ``quasi-coherent'' sheaves such as $W(M)$ should vanish for $i>0$;
consequently, presumably the torsors used for the inductive
construction of $X$ are trivial, which means that $X$ is in fact
isomorphic to the product of all $W(M_i)$'s. In the case when we
disregard the successive fibration structure, this means that the
objects of $U(k)$ are just sheaves of sets which are isomorphic to
some $W(M)$ (where morally $M$ is the direct sum of the modules $M_i$
which have been used in our inductive definition). This gives then a
rather trivial description of the objects of $U(k)$ (and all
$U_m(k)$'s are already equal to $U(k)$, for $m\ge1$!), it should be
remembered, however, that maps in $U(k)$ from a $W(M)$ to a $W(M')$
are a lot more general than just $k$-linear maps $M\to M'$ (they may
be viewed as ``polynomial maps from $M$ to $M'$'').

The category $U(k)$ will be endowed with the \emph{sections functor}
\begin{equation}
  \label{eq:111.h}
  X\mapsto X(k) : U(k) \to \Sets.\tag{h}
\end{equation}
Now, let $A$ be any test category, for instance $A=\Simplex$, and
consider the functor
\begin{equation}
  \label{eq:111.hprime}
  \bHom(A\op,U(k)) \to \bHom(A\op,\Sets)\tag{h'}
\end{equation}
induced by \eqref{eq:111.h}. The second member modelizes homotopy
types, which therefore allows us to define homotopy invariants for
objects in the first member, and hence to define the property of
$n$-connectedness, and (if $n=0$) of $0$-connectedness with abelian
fundamental group. As a matter of fact, we would like to define a
subcategory $\scrM_n(k)$ of the first member, so that is should become
clear that for an object $X_*$ in it, its homotopy groups are endowed
with structures of $k$-modules, so that for $n=0$ the abelian
restriction on $\pi_1$ is superfluous, because automatic. To be
specific, we better restrict maybe to $A=\Simplex$, more familiar to
us, so that $X_*$ may be viewed as a ``complex'' ($n\mapsto X_n$) with
components in $U(k)$. The kind of restriction I am thinking of for
defining the\pspage{449} model category $\scrM_n(k)$ is:
\begin{equation}
  \label{eq:111.i}
  X_i=e \quad\text{for}\quad i\le n,\tag{i}
\end{equation}
where $e$ is the final object of $U(k)$, and also maybe, if we adopt
the more refined version of $U(k)$ as a strictly increasing union of
subcategories $U_m(k)$,
\begin{equation}
  \label{eq:111.iprime}
  \text{$X_i$ is in $U_{i-n}$, for any $i\ge n$.}\tag{i'}
\end{equation}
One may have to play around some more to get ``the correct''
description of the model category, which I tentatively propose to
define simply as a suitable full subcategory
\begin{equation}
  \label{eq:111.j}
  \scrM_n(k)\subset \bHom(A\op,U(k)).\tag{j}
\end{equation}
The functor \eqref{eq:111.hprime} allows to define a notion of weak
equivalence in $\scrM_n(k)$, hence a localized category $\HotOf_n(k)$,
and a functor \eqref{eq:111.c} from this category to $\HotOf_n$. The
ring extension and restriction functors \eqref{eq:111.d} are equally
defined in an evident way, via corresponding functors on the model
categories (with the task, however, to check that these are compatible
with weak equivalences). The key point here is to check that for
$k=\bZ$, the functor \eqref{eq:111.c} (namely \eqref{eq:111.b}) is
indeed an equivalence of categories. Thus, the main task seems to cut
out carefully a description of a model category $\scrM_n(k)$, in terms
of semisimplicial objects say, in a category such as $U(k)$, in such a
way as to give rise to an equivalence of categories \eqref{eq:111.b}.

One point which is still somewhat misty in this (admittedly overall
misty!) picture, is how to get, for an object $X$ in $\HotOf_n(k)$,
the promised operation of $k$ on the homotopy groups $\pi_i(X)$. I was
thinking about this when suggesting the conditions \eqref{eq:111.i}
and \eqref{eq:111.iprime} above on $k$-models for homotopy types --
but I really doubt these are enough. On the other hand, it seems hard
to imagine there be a good notion of homotopy types over $k$, without
the homotopy groups to be $k$-modules over $k$, not just abelian
groups. Even more still, there should be moreover a ``linearization
functor''
\begin{equation}
  \label{eq:111.k}
  \HotOf_n(k) \to \D_\bullet(\AbOf_k),\tag{k}
\end{equation}
with values in the derived category of the category of chain complexes
in $\AbOf_k= k\textup{-Mod}$ of $k$-modules, presumably coming by
localization from a functor
\begin{equation}
  \label{eq:111.kprime}
  \scrM_n(k)\to\Ch_\bullet(\AbOf_k),\tag{k'}
\end{equation}
and giving rise to a commutative diagram\scrcomment{footnote ``only for $k=\bZ$''}\pspage{450}
\begin{equation}
  \label{eq:111.l}
  \begin{tabular}{@{}c@{}}
    \begin{tikzcd}[baseline=(O.base)]
      \HotOf_n(k)\ar[r]\ar[d] & \D_\bullet(\AbOf_k)\ar[d] \\
      \HotOf_n\ar[r] & |[alias=O]| \D_\bullet\Ab      
    \end{tikzcd},
  \end{tabular}\tag{l}
\end{equation}
where the second horizontal arrow is the usual abelianization functor
for homotopy types, and the second vertical one comes from the ring
restriction functor $\AbOf_k\to\AbOf=\AbOf_\bZ$. In
\eqref{eq:111.kprime}, $\Ch_\bullet$ denotes the category of chain
complexes, and it looks rather mysterious again how to get such a
functor \eqref{eq:111.kprime}. We may of course think of the trivial
abelianization
\[X_\bullet \mapsto k^{(X_\bullet)} \eqdef (n\mapsto k^{(X_n)}),\]
where for an object $X$ in $U(k)$, or more generally any sheaf on
$\Aff_{/k}$, $k^{(X)}$ defines a trivial $k$-linearization of this
sheaf, in the sense of the topos of all such sheaves. Anyhow,
$k^{(X)}$ is a sheaf, not a $k$-module, so we should still take
sections to get what we want -- but this functor looks not only
prohibitively large and inaccessible, but just silly! A much better
choice for $k$-linearizing objects of $U(k)$ specifically seems the
following. Disregarding the fibration structures, such an object $X$
is isomorphic to an object $W(M)$, $M$ some $k$-module. We look for a
$k$-linearization
\begin{equation}
  \label{eq:111.m}
  X\;(=W(M)) \to W(L(X)),\tag{m}
\end{equation}
where $L(X)$ is a suitable $k$-module. Now, among all maps
\[X\to W(N)\]
of $X$ into sheaves \emph{of the type $W(N)$}, there is a universal
one, which in terms of $M$ can be described as
\[N=\Gamma_k(M),\]
where $\Gamma_k$ denotes the ``algebra with divided powers generated
by $M$'', the canonical map $M\to \Gamma_k(M)$ or rather
\begin{equation}
  \label{eq:111.n}
  W(M)\to W(\Gamma_k(M)), \quad
  x\mapsto \exp(x)=\sum_{i\ge0} x^{(i)}\tag{n}
\end{equation}
being \emph{the} ``\emph{universal} polynomial map'' of $M$ with
values in a module $N$ (or rather, of $W(M)$ into $W(N)$). Here,
$x^{(i)}$ denotes the $i$'th divided power of $x$, which is an element
of $\Gamma_k^i(M)$. It just occurs to me that this expression of
$\exp(x)$, the universal map, is infinite, thus, it doesn't take its
values in $W(N)$ actually, but in a suitable completion of it -- this
doesn't seem too serious a drawback, though! The point I wish to make
here, is that for given $X$ in $U(k)$, defining\pspage{451}
\begin{equation}
  \label{eq:111.o}
  L_k(X) = \Gamma\uphat_k(M) \eqdef \prod_{i\ge0} \Gamma_k^i(M),
  \quad\text{($k$-linearization of $X$),}\tag{o}
\end{equation}
where $M$ is any $k$-module endowed with an isomorphism
\[u: X\simeq W(M),\]
the $k$-module $L_k(X)$ \emph{does not depend, up to unique
  isomorphism, on the choice of a pair $(M,u)$}, because for two
$k$-modules $M,M'$, any morphism of \emph{sheaves of sets}
\[v:W(M)\to W(M')\]
induces a homomorphism of $k$-modules
\[\Gamma\uphat_k(v): \Gamma\uphat_k(M) \to \Gamma\uphat_k(M')\]
(compatible not with multiplications, but with diagonal maps\dots),
which will be an isomorphism if $v$ is.

Thus, we do have, it seems to me, a good candidate for
$k$-linearization. To check it is suitable indeed, the main point
seems to check that the corresponding diagram \eqref{eq:111.l}
commutes up to canonical isomorphism, the crucial case of course being
$k=\bZ$. This now looks like a rather down-to-earth question, which
seems to me a pretty good test, whether the intuition of
schematization of homotopy types is a sound one. Let's rephrase it
here. For this, let's first restate the description of the category
$U(k)$ (coarse version) in the more down-to-earth terms of linear
algebra. Objects may be viewed as just $k$-modules $M$, whereas
(non-additive) ``maps'' from $M$ to $M'$ (defined previously as maps
$W(M)\to W(M')$ of sheaves of sets) are described as just continuous
$k$-linear maps
\begin{equation}
  \label{eq:111.p}
  f:\Gamma\uphat_k(M)\to\Gamma\uphat_k(M'),\tag{p}
\end{equation}
which are moreover compatible with the natural augmentations to $k$,
and with the natural diagonal maps:
\[\varepsilon:\Gamma\uphat_k(M)\to k, \quad
\Delta:\Gamma\uphat_k(M)\to\Gamma\uphat_k(M)\hatotimes_k\Gamma\uphat_k(M)
\;(\simeq\Gamma\uphat_k(M\times M)),\]
(the latter deduced from the usual linear diagonal map $M\to M\times
M$). When $M$ is looked at as being embedded in $\Gamma\uphat_k(M)$ by
the exponential map \eqref{eq:111.n}, it is identified (if I remember
it right) to the set of elements $\xi$ in $\Gamma\uphat_k(M)$
satisfying the relations
\begin{equation}
  \label{eq:111.q}
  \varepsilon(\xi)=1, \quad \Delta(\xi)=\xi\otimes\xi,\tag{q}
\end{equation}
where $\varepsilon$ is the augmentation and $\Delta$ the diagonal map,
hence \eqref{eq:111.p} induces a map (in general not
additive)\pspage{452}
\begin{equation}
  \label{eq:111.pprime}
  \Gamma(f):M\to M'\tag{p'}
\end{equation}
(corresponding to the action of $f$, viewed as a map $W(M)\to W(M')$,
on \emph{sections} of $W(M)$) -- and likewise after any ring extension
$k\to k'$, defining a map
\[\Gamma(f)_{k'}: M\otimes_k k' \to M'\otimes_k k'\]
from $W(M)$ to $W(M')$ -- which is the description of the map of
sheaves $W(M)\to W(M')$ associated to a map \eqref{eq:111.p}. We have
thus a description, in terms of linear algebra, of a category $U(k)$,
and of a ``sections'' functor
\begin{equation}
  \label{eq:111.r}
  \Gamma: U(k)\to\Sets, \quad
  X\mapsto X(k),\tag{r}
\end{equation}
which is essentially the functor \eqref{eq:111.h} above, viewed in a
different light. Now to our

\begin{questionnum}\label{q:111.1}
  The category $U(k)$ and the functor \eqref{eq:111.r} being defined
  as above in terms of linear algebra over the ground ring $k$, let
  $X_*=(n\mapsto X_n)$ be a semisimplicial object in $U(k)$, consider
  the corresponding semi-simplicial set $X_*(k)=(n\mapsto X_n(k))$,
  and the semi-simplicial $k$-module $L(X_*)=(n\mapsto L(X_n))$,
  where, for an object $M$ of $U(k)$, $L(M)$ is defined as
  \begin{equation}
    \label{eq:111.s}
    L(M)=\Gamma\uphat_k(M),\tag{s}
  \end{equation}
  which depends functorially on $M$ in $U(k)$. Then in the derived
  category of \Ab, is there a canonical isomorphism between $L(X_*)$
  (with ground ring restricted from $k$ to \bZ) and the abelianization
  $\bZ^{(X_*(k))}$ of $X_*(k)$?\footnote{\scrcommentinline{\dots
      unreadable\dots} $k=\bZ$, and components $X_n$ projective.}
\end{questionnum}

We may have to throw in some extra assumption on $X_*$, at any rate
\begin{equation}
  \label{eq:111.t}
  X_0=e\quad\text{(final object of $U(k)$),}\tag{t}
\end{equation}
giving rise to $L(X_0)\simeq k$. Also, we may have to restrict to
$k=\bZ$, or otherwise correct the obvious drawback that the two chain
complexes don't have isomorphic $\mathrm H_0$ (one is $k$ I guess, a
$k$-module in any case, the other is \bZ), by truncating accordingly
the two chain complexes (``killing'' their $\mathrm H_0$). There is a
natural candidate for a map
\begin{equation}
  \label{eq:111.u}
  \bZ^{(X_*(k))} \to L(X_*),\tag{u}
\end{equation}
by using the functorial map
\begin{equation}
  \label{eq:111.uprime}
  \bZ^{(M)} \to L(M)=\Gamma\uphat_k(M),\tag{u'}
\end{equation}
deduced from the inclusion \eqref{eq:111.n}
\[M\hookrightarrow L(M),\]
and we may still specify the question above, by asking \emph{whether
  \eqref{eq:111.u} induces an isomorphism for homology groups} in
dimension $i>0$.

I\pspage{453} am not too sure whether all this isn't just complete
nonsense -- it is worth getting it clear whether it is or not, at any
rate! There is one case of special interest, the ``simplest'' one in a
way, namely when the simplicial maps between the $X_n$'s, each
represented by a $k$-module $M_n$, are in fact $k$-linear, in other
words, when $X_*$ comes from a semisimplicial $k$-module $M_*$ -- more
specifically still, when this is a $K(\pi,n)$ type, say the nicest
semisimplicial model of this, using the Kan-Dold-Puppe functor for the
chain complex of $k$-modules, having $\pi$ in degree $n$ and zero
elsewhere. Then the left-hand side of \eqref{eq:111.u} gives rise to
the Eilenberg-Mac~Lane homology groups
\begin{equation}
  \label{eq:111.v}
  \mathrm H_i(\pi,n;\bZ),\tag{v}
\end{equation}
which I guess should be $k$-modules a priori, because of the
operations of $k$ upon $\pi$, and the question then arises whether
these can be computed using the right-hand complex $L(X_*)$. Maybe
such a thing is even a familiar fact for people in the know? If it
turned out to be false even for $k=\bZ$, my faith in schematization of
homotopy types would be seriously shaken I confess\dots

\starsbreak

After a little break for dinner, just one more afterthought. Working
with the completions $\Gamma\uphat_k$ may seem a little forbidding,
and all the more so if it used for computing \emph{homology}
invariants, not cohomology (the latter more likely to involve infinite
products in the corresponding (cochain) complexes\dots). On the other
hand, as was explicitly stated from the beginning, the natural context
here seems to be \emph{pointed} homotopy types, and hence ``pointed''
algebraic paradigms for these -- an aspect we lost sight of, when
looking for a suitable description of some category $U(k)$ of
``unipotent bundles'' over $k$. It would seem that in the ``question''
above, we should therefore insist that $X_*$ should be a
semi-simplicial object of the category $U(k)^\bullet$ of ``pointed''
objects of $U(k)$, namely objects $X$ endowed with a section over the
final object $e$ (i.e., with an element in $X(k)$). This will be
automatic at any rate in terms of the condition \eqref{eq:111.t},
$X_0=e$. The point I wish to make is that the category of pointed
objects of $U(k)$ admits a somewhat simpler description (by choosing
the marked point as the ``origin'' for parametrization of the given
object $X$ of $U(k)$ by a $k$-module $M$), by model-objects which are
still arbitrary $k$-modules $M$, but the ``maps'' now being
\emph{$k$-linear continuous maps}
\begin{equation}
  \label{eq:111.pprimebis}
  f:\Gamma_k(M)\to\Gamma_k(M')\tag{p'}
\end{equation}
between\pspage{454} the $k$-modules $\Gamma_k$, without having to pass
to completions, satisfying compatibility with augmentations and
diagonal maps, and the extra condition (expressing that
$\Gamma(f)(0)=0$ in exponential notation):
\[f(1)=1,\]
i.e., $f$ reduced to component $\Gamma_k^0(M)\simeq k$ of degree zero
is just the identity of $k$ with $k\simeq\Gamma_k^0(M')$. Accordingly,
we have a less awkward $k$-linearization functor than $L$ in
\eqref{eq:111.s}, namely ``pointed linearization'' $L\subpt$:
\begin{equation}
  \label{eq:111.sprime}
  M\mapsto L\subpt(M)\eqdef \Gamma_k(M) : U(k)^\bullet\to\AbOf_k,\tag{s'}
\end{equation}
which seems to me the better candidate for describing
linearization. Thus, we better rephrase now the ``question'' above in
terms of \eqref{eq:111.sprime} rather than \eqref{eq:111.s}. One
trouble however is that the comparison map \eqref{eq:111.u} takes
values in $L(X_*)$, not $L\subpt(X_*)$, therefore, we may still have
to use the ``prohibitive'' $L(X_*)$ as an intermediary for comparing
the complexes $\bZ^{(X_*(k))}$ and $L\subpt(X_*)$. It may be noted now
that, while the first chain complex embodies Eilenberg-Mac~Lane
homology \eqref{eq:111.v} (in the special case considered above), the
second one $L\subpt(X_*)$ (in that same case) describes the value of
the total derived functor of the familiar $\Gamma_k$ functor, on the
``argument'' $\pi$ placed in degree $n$, and the statement that the
two are ``the same'' does sound like some standard Dold-Puppe type
result which everybody is supposed to know from the cradle -- sorry!

\starsbreak

After another break (visit, tentative nap), still another
afterthought. The final shape we arrived at for the ``question''
above, when working in the ``pointed'' category $U(k)^\bullet$ of
``pointed unipotent bundles over $k$'', was whether for any
semisimplicial object $X_*$ in $U(k)^\bullet$ satisfying
\eqref{eq:111.t} above, i.e., $X_0=e$ (final object of
$U(k)^\bullet$), the two canonical maps of chain complexes (in fact,
semisimplicial abelian groups)
\begin{equation}
  \label{eq:111.w}
  \bZ^{(X_*(k))} \to L(X_*) \hookleftarrow L\subpt(X_*)\tag{w}
\end{equation}
are quasi-isomorphisms (for $\mathrm H_i$ with $i>0$). I am not too
sure yet if some extra conditions on $X_*$ are not required for this
to be reasonable -- I want to review two that came into my mind.

As the maps \eqref{eq:111.w} are functorial for varying $X_*$, it
would follow\pspage{455} from a positive answer that whenever
\begin{equation}
  \label{eq:111.x}
  X_*\to X_*'\tag{x}
\end{equation}
is a map of semisimplicial objects in $U(k)^\bullet$ satisfying
condition \eqref{eq:111.t}, and such that the corresponding map
\begin{equation}
  \label{eq:111.xprime}
  X_*(k) \to X_*'(k)\tag{x'}
\end{equation}
is a weak equivalence, and hence the map between the
\bZ-abelianizations is a weak equivalence too, i.e., a
quasi-isomorphism, that the same holds for the corresponding map
\begin{equation}
  \label{eq:111.y}
  L\subpt(X_*)\to L\subpt(X_*')\tag{y}
\end{equation}
for the ``schematic'' $k$-linearizations. Now, this is far from being
an evident fact by itself, except of course in the case when the map
\eqref{eq:111.x} above is a homotopism. Take for instance the case
when we start with a map of chain complexes in $\AbOf_k$
\begin{equation}
  \label{eq:111.xdblprime}
  M_\bullet \to M_\bullet',\tag{x''}
\end{equation}
hence a map between the associated semisimplicial $k$-modules
\begin{equation}
  \label{eq:111.xtplprime}
  M_*\to M_*',\tag{x'''}
\end{equation}
which may be viewed as giving rise to a (componentwise linear) map
between the associated semisimplicial objects $X_*$, $X_*'$ in
$U(k)^\bullet$ via the canonical functor
\[\AbOf_k\to U(k)^\bullet;\]
the corresponding map \eqref{eq:111.y} is then just the componentwise
extension of \eqref{eq:111.xtplprime} to the enveloping algebras with
divided powers
\begin{equation}
  \label{eq:111.yprime}
  \Gamma_k(M_*) \to \Gamma_k(M_*').\tag{y'}
\end{equation}
The map \eqref{eq:111.xprime} can now be identified with
\eqref{eq:111.xtplprime}, hence it is a weak equivalence if{f}
\eqref{eq:111.xtplprime} is, i.e., if{f} \eqref{eq:111.xtplprime} is a
quasi-isomorphism. If we assume moreover the components of
$M_\bullet$, $M_\bullet'$ to be projective objects in $\AbOf_k$, then
from the assumption that \eqref{eq:111.xdblprime} is a
quasi-isomorphism it does follow that it is a chain homotopism, hence
by Kan-Dold-Puppe the map \eqref{eq:111.xtplprime} is a semisimplicial
homotopism, and hence the same holds for \eqref{eq:111.yprime}, and
\eqref{eq:111.yprime} therefore is indeed a quasi-equivalence. But
without the assumption that components are projective, it is surely
false that the mere fact that \eqref{eq:111.xdblprime} is a
quasi-isomorphism, implies that \eqref{eq:111.yprime} is -- otherwise
this would mean that in order to compute the left derived functor of
the non-additive functor\pspage{456}
\[\Gamma_k:\AbOf_k\to\AbOf_k,\]
it is enough, for getting its value on a chain complex $M_\bullet'$
say, to replace $M_\bullet'$ by $M_*'$ and apply the functor
$\Gamma_k$ componentwise, without first having to take a projective
resolution $M_\bullet$ of $M_\bullet'$ -- something rather absurd
indeed! Thus, the statement made on page \ref{p:454}, that when taking
for $M_\bullet'$ the $k$-module $\pi$ placed in degree $n$ and zero in
all other degrees, the corresponding $L\subpt(X_*')=L_k(M_*')$
embodies the value of the left derived functor $\mathrm L\Gamma_k$ on
$M_\bullet'$, is visibly incorrect if we don't assume moreover that
$\pi$ is projective (flat, presumably, would be
enough\dots). Otherwise, we should first replace $\pi$ by a projective
(or at any rate flat) resolution, which we shift by $n$ to get
$M_\bullet$, and then take $\Gamma_k(M_*)$ to get the correct value of
$\mathrm L\Gamma_k(M_\bullet')$.

This convinces me that in the question as to whether the maps in
\eqref{eq:111.w} are quasi-isomorphisms (the more crucial one of
course being the first of the two), \emph{we should assume moreover
  that the components of $X_*$ are described in terms of
  \emph{projective} $k$-modules $M_n$}, or at any rate $k$-modules
that are flat. Accordingly, we should make the same restriction on the
semisimplicial schematized model $X_*$, in order for the description
we gave of ``$k$-linearization'' as $L\subpt(X_*)$ (or $L(X_*)$, never
mind which) to be topologically meaningful. Very probably, in the
whole schematization set-up, namely in the very definition we gave of
$U(k)$ and $U(k)^\bullet$, we should stick to the same restriction. If
I insisted first (with some inner reluctance, I admit) on taking
$k$-modules $M$ unrestricted, this was because I was thinking of $M$,
more specifically of the $M_i$'s occurring in the inductive
``dévissage'' of an object of $U_m(k)$ (when thinking of the more
refined version of $U(k)$), as essentially the homotopy groups of the
homotopy type we want to modelize, or rather, as the components of the
corresponding semisimplicial $k$-modules (denoted $M_*'$ some minutes
ago). I was still thinking of course, be it implicitly, in terms of
Postnikov dévissage, despite yesterday's remark that to use such
dévissage literally may cause trouble (p.\ \ref{p:444}). Thus, the
feeling which gets into the fore now is that \emph{we should kind of
  forget Postnikov, and work with semisimplicial ``schematic'' models
  built up with $k$-modules which are \emph{projective}}, or at any
rate flat (namely torsion free, if $k=\bZ$).

It\pspage{457} may be remarked that if $M$ is any $k$-module, then the
property that $M$ be projective, or flat, can be described in terms of
the isomorphism class of the corresponding object $X$ in
$U(k)^\bullet$, or equivalently, of the functor $W(M)$ on $\Alg_{/k}$,
with values in the category of pointed sets. Indeed, the isomorphism
class of the $k$-module $\Gamma_k(M)$ depends only on the class of
$X$, and it is easily seen that $M$ is projective, resp.\ flat, if{f}
$\Gamma_k(M)$ is. (The ``only if'' is standard knowledge of
commutative algebra, the ``if'' comes from the fact that $M$ is a
direct factor of $\Gamma_k(M)$, hence projective resp.\ flat if the
latter is.) Using this, one even checks that it is enough to know the
isomorphism class of $X$ in $U(k)$ -- because two objects of
$U(k)^\bullet$ are isomorphic in $U(k)^\bullet$ if{f} they are in
$U(k)$.

\starsbreak

The second afterthought is that in the question on page \ref{p:452},
we should definitely assume $k=\bZ$. Already when asserting hastily
(p.\ \ref{p:453}) that en Eilenberg-Mac~Lane homology group $\mathrm
H_i(\pi,n;\bZ)$ should automatically inherit a structure of a
$k$-module, whenever $\pi$ had one, I was feeling uncomfortable,
because after all the dependence of this group upon variable $\pi$ is
not at all additive, thus the operation of $k$ upon this group
stemming from its operations on $\pi$ was no too likely to come out
additive! To take one example, take $n=1$, i.e., we just take ordinary
group homology for $\pi$, and assume that $k$ is free of finite type
over \bZ, and $\pi=M$ free of finite type over $k$ (for instance
$M=k$), hence free of finite type over \bZ{} too. Then it is
well-known that
\[\mathrm H_*(M; \bZ) \simeq \bigwedge_\bZ M,\]
the exterior algebra of $M$ over \bZ, which surely is not endowed with
a structure of a $k$-module in any natural way! This, if there was any
such structure (natural or not) on the highest non-zero term
(corresponding to the rank $d$ of $M$ over \bZ), it would follow that
we get a ring homomorphism from $k$ to $\bZ\simeq\End_\bZ(\mathrm
H_d\simeq\bZ)$, and we may choose $k$ in such a way that there is no
such homomorphism, for instance $k=\bZ[T]/(T^2+1)$. In this case,
there cannot be \emph{any} isomorphism between the $\mathrm H_d$'s of
the two members of \eqref{eq:111.u}! Another point, which I hit upon
first, is that when $X_*$ is defined in terms of a semisimplicial
object $M_*$ of $\AbOf_k$, then the functorial dependence of the first
term in \eqref{eq:111.w} with respect to\pspage{458} varying $M_*$, is
that to a direct sum corresponds the componentwise tensor product
\emph{over \bZ}, whereas for the last term of \eqref{eq:111.w} we have
to take tensor products \emph{over $k$} (in the middle term,
\emph{completed} tensor products over $k$) -- the two variances are
clearly at odds with each other.

Thus, when working with would-be ``$k$-homotopy types'' as defined
here via semisimplicial objects in $U(k)^\bullet$, \emph{we should
  altogether drop the idea that the homology groups of the
  corresponding semisimplicial set $X_*(k)$ are $k$-modules}. I
wouldn't really look at them as being ``the'' homology groups of the
$k$-homotopy type $X_*$, these should be rather given via
$k$-linearization $L\subpt$ and they are $k$-modules, indeed, they
come from a canonical object of $\D_\bullet(\AbOf_k)$, namely
$L\subpt(X_*)$. In the example just looked at, presumably we should
get the exterior algebra of $M$ \emph{over $k$} (not \bZ!). This makes
me suspect even that, except in the case $k=\bZ$, this semisimplicial
set $X_*(k)$ doesn't make much sense, namely its homology invariants
(and presumably its homotopy groups too) are not really relevant for
the $k$-homotopy type $X_*$, \emph{which has invariants of its own
  which are completely different. Thus, I am not at all convinced any
  more that the homotopy groups $\pi_i(X_*)$ carry $k$-module
  structures} (as expected at the beginning, p.\ \ref{p:446}) -- but
to clear our mind on that matter, we should take off from the
simplistic example when $X_*$ comes from a semisimplicial $k$-module
$M_*$, in which case of course we have by Whitehead's isomorphism
\[\pi_i(X_*(k)) \simeq \mathrm H_i(M_\bullet) \; (\simeq\pi_i(M_*)),\]
where $M_\bullet$ is the chain complex associated to $M_*$. Also, it
is now becoming obvious whether weak equivalences for objects $X_*$
should be defined (as we did) via $X_*(k)$, as suggested by the
``simplistic'' example above. If we define homology of $X_*$ by the
formula
\begin{equation}
  \label{eq:111.z}
  \LH_\bullet(X_*) \eqdef L\subpt(X_*)\quad
  \text{viewed as an object of $\D_\bullet(\AbOf_k)$,}\tag{z}
\end{equation}
and accordingly, the homology modules
\begin{equation}
  \label{eq:111.zprime}
  \mathrm H_i(X_*) \eqdef \pi_i(L\subpt(X_*)) = \mathrm
  H_i(\LH_\bullet(X_*)),\tag{z'} 
\end{equation}
maybe the better idea for defining weak equivalences $X_*\to X_*'$ is
by demanding that \emph{they should be transformed into
  quasi-isomorphisms by the total homology functor $\LH_\bullet$}, or
equivalently, induce isomorphisms for the homology modules
\eqref{eq:111.zprime}. If the answer to our crucial question is
affirmative (with the corrections made, including $k=\bZ$),
then\pspage{459} in case $k=\bZ$, the new definition just given for
weak equivalences is equivalent to the old one in terms of $X_*(\bZ)$
(taking \bZ-valued points), provided at any rate we admit or rather
assume that $X_*(\bZ)$ and $X_*'(\bZ)$ are simply connected, which
will be automatic if we work in the category of (schematized) models
$\scrM_1(k)$, the condition \eqref{eq:111.t} above ($X_0=e$) being
replaced by \eqref{eq:111.i} with $n=1$, i.e., by
\begin{equation}
  \label{eq:111.alpha}
  X_0=X_1=e,\tag{$\alpha$}
\end{equation}
to be on the safe side! Under this extra assumption at any rate, I
feel definitely more confident with the new definition of weak
equivalence, via the homology invariants \eqref{eq:111.zprime}, rather
than the old one. At any rate, the question of the two definitions
being equivalent or not should be cleared up, namely:
\begin{questionnum}\label{q:111.2}
  Let $k$ be any ring, define the category $U(k)^\bullet$ of ``pointed
  unipotent bundles over $k$'' in terms of projective $k$-modules,
  with maps defined as in \eqref{eq:111.pprimebis} page
  \ref{p:453}. (This is equally the correct set-up for
  question~\ref{q:111.1} on page \ref{p:452}, besides the extra
  condition $k=\bZ$, as we saw before.) Let
  \[u:X_*\to X_*'\]
  be a map of semisimplicial objects in $U(k)$, satisfying both the
  extra assumptions \eqref{eq:111.alpha} above. Then is it true that
  the corresponding map
  \begin{equation}
    \label{eq:111.beta}
    X_*(k)\to X_*'(k)\tag{$\beta$}
  \end{equation}
  of semisimplicial sets is a weak equivalence, if{f} the map
  \[L\subpt(X_*)\to L\subpt(X_*')\]
  is, i.e., if{f} $u$ induces an isomorphism for the homology
  invariants $\mathrm H_i$ defined in \eqref{eq:111.zprime} above (via
  the abelianization functor $L\subpt=\Gamma_k$).
\end{questionnum}

If instead of \eqref{eq:111.alpha} we only assume $X_0=X_0'=e$, it
seems that we may have to throw in some other extra condition on $X_*$
and $X_*'$, in order for the definition of weak equivalence in terms
of the mere homology invariants to be reasonable -- a condition which,
at the very least, and in case $k=\bZ$ say, should ensure that the
homotopy types defined by the two terms in \eqref{eq:111.beta} should
have \emph{abelian fundamental groups} (which doesn't look at first
sight to be automatic). At any rate, the conditions
\eqref{eq:111.alpha} above, which should be viewed as a
$1$-connectedness assumption, are natural enough, and it is natural
too to try first to push through a theory of schematization of
homotopy types, under this assumption.

\bigbreak

\presectionfill\ondate{26.8.}\pspage{460}\par

\hangsection[Postnikov dévissage and Kan condition for schematic
\dots]{Postnikov dévissage and Kan condition for schematic
  complexes.}\label{sec:112}%
I am continuing the ``wishful thinking'' about schematization of
homotopy types -- a welcome break in the ``overall review'' on
linearization (in the context of the modelizer \Cat), which had been
getting a little fastidious lately!

I'll admit, as one firm hold in all the wishfulness, that in the
``Eilenberg-Mac~Lane case'' of p.~\ref{p:453}, when moreover $k=\bZ$
and the components $X_i$ of the semisimplicial unipotent bundle $X_*$
are projective, the two maps \eqref{eq:111.w} of page \ref{p:454} are
indeed quasi-isomorphisms. From this should follow the similar
statement, when $X_*$ comes from a chain complex of $k$-modules
$M_\bullet$ with projective coefficients, by reducing to the case when
only a finite number of components of $M_\bullet$ are non-zero (by
suitable passage to the limit), and then by induction on the number of
these components, using the fact that the three terms in
\eqref{eq:111.w} depend on $X_*$ in a ``multiplicative'' way, namely
direct sums being transformed into tensor products. (NB\enspace Under
the assumptions of projectivity made, we may as well express the
quasi-isomorphisms \eqref{eq:111.w} we start with as being
semisimplicial \emph{homotopisms}, and remark that componentwise
tensor-product of such homotopisms is again one.) From this, using the
relevant spectral sequences in homology, should follow that the maps
\eqref{eq:111.w} are still quasi-isomorphisms, whenever $X_*$ can be
``unscrewed'' (``dévissé'') as a finite successive fibering with
fibers of the type $M_*$ as above. Another passage to the limit will
yield the same result for an infinite dévissage, provided the fibers
$M(i)_*$ ($i=1,2,\ldots$) are ``way-out'', i.e., for given $n$, only a
finite number of components $M(i)_n$ are non-zero (it amounts to the
same to demand that the sequence $M(i)_\bullet$ of corresponding chain
complexes be ``way out''). This will give already, it seems, a fair
number of cases when \eqref{eq:111.w} are
quasi-isomorphisms. (Admittedly, working this out will involve a fair
amount of work, especially for getting the relevant properties of
$L_\bullet$ and $L\subpt$, which should mimic very closely the known
ones for usual linearization, including spectral sequences or, more
neatly, transitivity isomorphisms in the relevant derived
categories\dots) The main point here is that those special types of
$X_*$'s (we may call them \emph{Kan-Postnikov complexes} in $U(k)$)
are enough in order to modelize, via the corresponding semisimplicial
sets $X_*(\bZ)$, arbitrary pointed homotopy types with abelian
$\pi_1$. This is seen of course using Postnikov\pspage{461} dévissage
of a given homotopy type, and replacing every homotopy group $\pi_i$
by a shifted projective resolution (a two-step resolution will do
here) $M(i)_\bullet$, as indicated on page \ref{p:456}, and using the
corresponding semisimplicial \bZ-modules $M(i)_*$ as fibers, in the
successive fiberings. To see that what is being done on the
``discrete'' level, working with semisimplicial \emph{sets}, can be
``followed'' in an essentially unique way on the ``schematic'' level,
we hit now of course upon the key difficulty, pointed out on page
\ref{p:444}, about the Postnikov cohomology group $\mathrm
H^{n+2}(X(n+1),\pi_{n+1})$ being isomorphic to the corresponding
``schematic'' group. Now, that this is so indeed should follow from
the homology isomorphisms \eqref{eq:111.w} I hope, just dualizing the
result to cohomology. All this seems to sound kind of reasonable, it
seems, even that for a given homotopy type in $\HotOf_0$, we should be
able to squeeze out this way a unique \emph{isomorphism type}, at any
rate, of semisimplicial unipotent bundles -- but to see whether it
does work, or if there is some major blunder which turns the whole
into nonsense, will come out only from careful, down-to-earth work,
which I am not prepared to dive into.

It occurred to me that the ``Kan-Postnikov'' complexes in $U(k)$ have
some special features among all possible complexes with $X_0=e$, and
also that some extra feature \emph{are} needed, if we want the maps in
\eqref{eq:111.w} to be quasi-isomorphisms. I want to dwell upon this a
little. First of all, the condition $X_0=e$ is indeed essential, as we
see by taking a constant complex with value $X_0$, then the homology
of the three chain complexes \eqref{eq:111.w} reduces to degree zero,
and the $\mathrm H_0$'s are respectively
\[\bZ^{(X_0(\bZ))}=\bZ^{(M_0)}, \quad
L_0(X_0)=\Gamma_bZ\uphat(M_0), \quad
L\subpt(X_0)=\Gamma_\bZ(M_0),\]
where $M_0$ is the \bZ-module giving rise to $X_0$ -- and none of the
two maps is an isomorphism, unless $M_0=0$.

Take now the next simplest case, when $X_*$ comes from a monoid object
$G$ in $U(k)$ in the usual way; then what we are after, in dual terms
of cohomology rather than homology (taking the dual complexes of those
in \eqref{eq:111.w}), amounts essentially to asking whether the usual
discrete cohomology of the discrete monoid $G(\bZ)$ can be computed,
using \emph{polynomial cochains} rather than arbitrary ones. Now, this
we did admit as ``well-known'' in the most evident case of all, when
$G$ is being represented as an object of $U(\bZ)$ by a projective
(hence free) \bZ-module $M$, the multiplication law is just usual
addition. It still looks reasonable enough when the monoid $G$ is a
\emph{group}, with $M$ of finite type say. In this case, the Borel
theory of algebraic\pspage{462} affine groups over a field (here, the
field of fractions \bQ{} of \bZ) tells us that $G_\bQ$ is a
\emph{nilpotent} algebraic group, and that therefore it admits a
composition series with factors isomorphic to the additive group
$\bG_{\mathrm a \bQ}$; presumably, the same dévissage then can
be obtained over the base \bZ, and using induction on the length of
the composition series, and the Hochschild-Serre type of relations
(traditionally expressed by a spectral sequence) between group
cohomology of a group, quotient group and corresponding subgroup, we
should get the wished for quasi-isomorphisms \eqref{eq:111.w}.

Take now, however, the simplest case of a monoid which isn't a group,
namely the multiplicative law on the affine line, given by the
polynomial law
\[W(\bZ) \times W(\bZ) \to W(\bZ) : (x,y) \mapsto xy.\]
the corresponding ``discrete'' monoid is just $\bZ^{(\times)}$, namely
the integers with multiplication, its $\mathrm H^1$ with coefficients
in \bZ{} is just the group of all homomorphisms
\begin{equation}
  \label{eq:112.star}
  \bZ^{(\times)} \to \bZ^{(+)} = \bZ,\tag{*}
\end{equation}
and denoting by \bP{} the set of all primes and using the prime
decomposition of integers, we find that
\[\Hom(\bZ^{(\times)}, \bZ^{(+)}) \simeq \bZ^\bP,\]
i.e., a family $(n_p)_{p\in\bP}$ of integers being associated the
homomorphism
\[\pm\prod_{p\in\bP} p^{\alpha_p} \mapsto \sum_{p\in\bP}\alpha_p
n_p.\]
On the other hand, the schematic $\mathrm H^1$ consists of all
homomorphisms \eqref{eq:112.star} that can be expressed by a
polynomial, hence induce a homomorphism of algebraic group schemes
$\bG\subm\to\bG\suba$, and it is well-known (and immediately checked)
that there is only the zero homomorphism!

Thus, it turns out that the assumptions made yesterday on $X_*$, in
order for the ``linearization theorem'' (!) to hold, namely the maps
\eqref{eq:111.w} to be quasi-isomorphisms, are definitely not strong
enough yet! One may think of throwing in the extra condition $X_1=e$,
so as to rule out monoids altogether (and even groups, too bad!), but
I don't think this helps at all (didn't try though to make a
counterexample). On the other hand, just restricting to Kan-Postnikov
complexes seems rather awkward, we definitely don't want to drag along
Postnikov fibrations as a compulsory ingredient of the complexes we
work with. The idea which comes up here is just to ``drop Postnikov
and keep Kan'' -- namely \emph{introduce a Kan type condition on
  semisimplicial complexes in $U(k)$}. If we mimic\pspage{463}
formally the usual ``discrete'' Kan condition, we get that (for given
pair of integers $k,n$ with $0\le k\le n$) a certain map from $X_n$,
to a certain finite projective limit defined in terms of the boundary
maps $X_{n-1}\to X_{n-2}$, should be epimorphic. Now, clearly $U(k)$
is by no means stable under fiber products, except under very special
assumptions (including differential transversality conditions, at any
rate), and on the other hand one feels that the notion of
``epimorphism'' one will have to work with in $U(k)$ will have to be a
lot more exacting than the map $X(k)\to Y(k)$ on sections being
surjective, or the usual categorical meaning within $U(k)$, which
looks kind of silly here. Even the most exacting surjectivity
condition on $X\to Y$, namely that it admit a section doesn't quite
satisfy me -- what I really want is that $X$ should be a trivial
bundle over $Y$, more specifically that $X$ is isomorphic to a product
$Y\times Z$, in such a way that the given map $X\to Y$ identifies with
the projection $Y\times Z\to Y$. Maybe this is too exacting a
condition, however, and hard to check in computational terms sometimes
(?), maybe we should be content with demanding only that $X\to Y$ has
a section, and moreover is ``smooth'', i.e., has everywhere a
surjective tangent map (which may be expressed on the corresponding
sheaf on $\Aff_{/k}$ by the familiar condition of ``formal
smoothness'', namely possibility of lifting sections over arbitrary
infinitesimal neighborhoods\dots). We'll have to choose at any rate
some such strong ``surjectivity'' notion in $U(k)$, which we'll call
``\emph{submersions}'' say. Thus, I feel a ``Kan complex'' in $U(k)$
should have boundary maps which are submersions. What we should do, is
to pin down some simple ``Kan condition'' on a complex $X_*$, in terms
of ``submersions'', in such a way as to ensure, at any rate
\begin{enumerate}[label=\alph*)]
\item\label{it:112.a}
  that for any pair $(n,k)$ with $0\le k\le n$, the object $X_*(n,k)$
  of ``horns of type $(n,k)$ of $X_*$'', expressed by the suitable
  finite inverse limits (in terms of boundary maps $X_{n-1}\to
  X_{n-2}$) is representable in $U(k)$, and
\item\label{it:112.b}
  the canonical map $X_n\to X_*(n,k)$ is a submersion,
\end{enumerate}
and such of course that all Kan-Postnikov complexes should satisfy
this Kan condition, at the very least.

The first non-trivial case, in view of $X_0=e$, is $n=2$, in which
case $X_*(2,k)$ is trivially representable by $X_1\times X_1$, and the
condition we get is that the three natural maps coming from boundary
maps
\begin{equation}
  \label{eq:112.starstar}
  \begin{tikzcd}[sep=small]
    X_2 \ar[r]\ar[r,shift left=1.5]\ar[r,shift right=1.5] & X_1\times X_1
  \end{tikzcd}\tag{**}
\end{equation}
should\pspage{464} be submersions for a ``\emph{Kan complex}''. In
case $X_*$ is defined by a monoid object $G$ as above, this clearly
implies that $G$ is a group -- which rules out the counterexample
above!

Of course, the very first thing we'll expect from a ``good notion'' of
Kan complexes in $U(k)$, is that for $k=\bZ$, it should make the
linearization theorem work, namely that maps in \eqref{eq:111.w}
p.~\ref{p:454} are quasi-isomorphisms. The next thing, very close to
this one but for arbitrary ground ring $k$ now, is that a map
\[X_* \to X_*'\quad \text{(with $X_0=X_0'=X_1=X_1'=e$)}\]
of Kan complexes in $U(k)$ is a homotopism if{f} it induces an
isomorphism on the homology modules \eqref{eq:111.zprime}
(p.~\ref{p:458}) -- which sounds reasonable precisely because we are
working with Kan complexes. If this is so, the homotopy category
$\HotOf_1(k)$ of $1$-connected homotopy types over $k$ may be
identified with a category of Kan complexes ``up to homotopy'', as
usual (but working now with complexes of unipotent bundles over
$k$). Third thing, still over arbitrary ground ring, would be a
development of the usual homotopy formalism in the unipotent context,
including (one hopes) homotopy fibers of maps, and Postnikov
dévissage. Again, it is hard to imagine how to get such dévissage,
without getting hold inductively of homotopy invariants $\pi_i$ which
are $k$-modules. This should come out if we are able to define
homotopy fibers as for an $(n-1)$-connected $k$-homotopy type (defined
here as one whose homology invariants $\mathrm H_i$ are zero for $i\le
n-1$), $\pi_n$ should be no more, no less than $\mathrm H_n$, which is
indeed a $k$-module. Coming back to $k=\bZ$ again, this should imply
that for $n\ge1$ at any rate, the canonical functor
\begin{equation}
  \label{eq:112.starstarstar}
  \HotOf_n(\bZ)\to\HotOf_n\tag{***}
\end{equation}
induces a bijection between isomorphism classes of schematic and
ordinary $n$-connected homotopy types -- and it will be hard to
believe this can be so, without this functor being actually an
equivalence of categories -- the expected apotheosis of the theory!
Maybe to this end, one may even be able to introduce reasonable
internal $\bHom$'s within the category of schematic Kan complexes, in
a way compatible with the familiar notion in the discrete set-up.

If\pspage{465} it is possible indeed to construct Postnikov dévissage
of a schematic Kan complex over any ground ring $k$, it is clear that
this is compatible with restriction of ground ring, hence it would
seem that formation of the homotopy invariants $\pi_i$ is compatible
with restriction of rings (whereas, as we noticed yesterday, the same
does definitely \emph{not} hold for the homology invariants $\mathrm
H_i$). Taking restriction to the ground ring \bZ, this shows that the
canonical functor \eqref{eq:112.starstarstar} from schematic to
discrete homotopy types is compatible with taking homotopy groups (but
not with homology) -- thus, the relation between $X_*$ and the complex
of sections $X_*(k)$ seems to be a rather close one, via the homotopy
groups, which are the same (and thus, the homotopy groups of $X_*(k)$
seem to turn out to be $k$-modules after all!). By the way, speaking
of ``restriction of ground ring'' for Kan complexes was a little
hasty, in view of the projectivity condition on the components, which
a priori seems to oblige us to assume $k$ to be a projective
$k_0$-module (for a given ring homomorphism
\[k_0\to k\quad\text{).}\]
Still, the remark about the sections functor $X_*\to X_*(k)$ makes
sense, without having to assume $k$ to be a projective \bZ-module!
Also, we feel that, by analogy of what can be done in the linear
set-up, when we define a total derived functor
\[\D_\bullet(\AbOf_k) \to \D_\bullet(\AbOf_{k_0})\]
without any assumption on the ring homomorphism $k_0\to k$, a notion
of ring restriction for schematic homotopy types should make sense
without any restriction, as was surmised yesterday. As in the linear
case, we should allow ourselves to work with schematic complexes which
are \emph{not} projective, but be prepared to take ``resolutions'' (in
some sense) of such general complexes by the more restricted ones
(with projective components).

There is no such difficulty in the case of the ring extension functor,
which transforms projective bundles over $k_0$ into projective bundles
over $k$. The reflections above suggest that, whereas ring extension
is compatible with taking total homology invariants $\LH_\bullet$, via
the corresponding functor
\[\D_\bullet(\AbOf_{k_0}) \to \D_\bullet(\AbOf_k),\]
it is compatible too with taking homotopy invariants $\pi_i$
separately.

\bigbreak

\presectionfill\ondate{27.8.}\pspage{466}\par

\hangsection[``Soft'' versus ``hard'' Postnikov dévissage, $\pi_1$ as
a group \dots]{``Soft'' versus ``hard'' Postnikov dévissage,
  \texorpdfstring{$\pi_1$}{pi-1} as a group scheme.}\label{sec:113}%
What I was thinking of last night (see last sentence) is that whereas
for total homology (not for the separate $\mathrm H_i$'s) we have the
comprehensive formula
\begin{equation}
  \label{eq:113.A}
  \LH_\bullet(X_* \otimes_k k') \simeq \LH_\bullet(X_*) \Lotimes_k k'
  \quad\text{(in $\D_\bullet(k')$),}\tag{A}
\end{equation}
(where $\Lotimes_k$ denotes the left derived functor of the ring
extension functor for $k\to k'$, and $X_*'=X_*\otimes_k k'$ denotes
ring extension for the semisimplicial unipotent bundle $X_*$), for
homotopy modules we should have the term-by-term isomorphisms
\begin{equation}
  \label{eq:113.B}
  \pi_i(X_*\otimes_k k')\leftarrow \pi_i(X_*) \otimes_k k'.\tag{B}
\end{equation}
This however was pretty rash indeed (it was time to go to sleep I
guess!). Whereas the map on sections
\[X_*(k)\to X_*'(k')\]
does induce a map
\begin{equation}
  \label{eq:113.Bprime}
  \pi_i(X_*) \simeq \pi_i(X_*(k)) \to \pi_i(X_*')\simeq\pi_i(X_*'(k'))
  \tag{B'} 
\end{equation}
which surely is $k$-linear, and hence induces a map \eqref{eq:113.B},
this map is certainly \emph{not} an isomorphism without some flatness
restriction either on $k'$ over $k$, or on the $k$-modules
$\pi_j(X_*)$ for $j<i$, as we had noted already three days ago when
looking at the case when $X_*$ comes from a chain complex $M_\bullet$
in $(\AbOf_k)$ (with projective coefficients say), and hence $X_*'$
comes from
\[M_\bullet' = M_\bullet \otimes_k k'.\]
If we look at the description of the \emph{$k$-module}
$\pi_i(X_*)=\pi_i$ in terms of a Postnikov dévissage of $X_*$, we
should recall that the semisimplicial group object $M(i)_*$ which
enters into the picture as the $i$'th step fiber is \emph{not} the one
defined \emph{directly} (via the Kan-Dold-Puppe functor) by $\pi_i$
placed in degree $i$, but rather by the chain complex $M(i)_\bullet$
\emph{with projective components}, obtained by taking first a shifted
projective resolution of $\pi_i$. Thus, by ring extension we get from
this dévissage of $X_*$ another one of $X_*'$, whose successive fibers
are
\[M(i)_*' = M(i)_* \otimes_k k',\]
corresponding to the chain complexes
\[M(i)_\bullet' = M(i)_\bullet \otimes_k k'.\]
The\pspage{467} latter has $\pi_i\otimes_k k'$ as homology module in
degree $i$ (and zero homology in degree $j<i$), but the homology
modules in degree $j>i$ need not be zero. In other words, there is a
canonical \emph{augmentation}
\begin{equation}
  \label{eq:113.Bdblprime}
  M(i)_*' \to K(\pi_i\otimes_k k', i),\tag{B''}
\end{equation}
(where the second member is the semisimplicial $k'$-module defined by
an $i$-shifted projective resolution of $\pi_i(X)\otimes_k k'$), but
this augmentation need not be a quasi-isomorphism, unless we make the
relevant flatness assumptions. To sum up, \emph{the dévissage of}
\[X_*' = X_*\otimes_k k'\]
\emph{deduced from a Postnikov dévissage of $X_*$ is \emph{not} a
  Postnikov dévissage of $X_*'$}, unless we assume either $k'$ flat
over $k$, or the $k$-modules $\pi_j(X_*)$ flat. This at the same time
solves the puzzle raised on page \ref{p:444}, and points towards a
\emph{serious shortcoming of the \textup(usual\textup) Postnikov
  dévissage -- namely that it is not compatible with ground ring
  extension}, or, as we would say in the language of algebraic
geometry, that this construction is not ``geometric'' -- a harsh thing
to say indeed!

At this point the idea comes up that we may define another dévissage,
a lot more natural in the spirit of a theory of ``abelianization'' of
homotopy types it would seem, and which is ``geometric'', namely
compatible with ring extension. Here, we'll have to work, though, with
the ``prohibitive'' abelianization functor
\[U(k)\to\AbOf_k, \quad
X\mapsto L(X)\; (=\Gamma\uphat_k(M))\]
(where $M$ is a $k$-module ``representing'' the object $X$ in $U(k)$),
as we'll need the functorial embedding\scrcomment{see also page
  \ref{p:474}}
\begin{equation}
  \label{eq:113.C}
  X\to W(L(X))\tag{C}
\end{equation}
(where
\[W : \AbOf_k\to U(k)\]
denotes the canonical functor from $k$-modules to unipotent
$k$-bundles); this is licit anyhow, if we admit that the canonical map
\[L\subpt(X_*) \to L(X_*)\]
(cf.\ page \ref{p:454}~\eqref{eq:111.w}) is a quasi-isomorphism, i.e.,
a weak equivalence, for the complexes $X_*$ we are working with. (This
of course should hold over an arbitrary ground ring, not just \bZ.)
Applying \eqref{eq:113.C} componentwise, we get for a complex of
bundles $X_*$ a canonical map into its\pspage{468} abelianization
\begin{equation}
  \label{eq:113.Cprime}
  X_*\to W(L(X_*)).\tag{C'}
\end{equation}
Postnikov's construction, for $(n-1)$-connected $X_*$, consists in
composing this map with the ``augmentation''
\[W(L(X_*)) \to K(\pi_n,n)\]
(where $\pi_n\simeq\mathrm H_n$ is the first possibly non-trivial
homology module of the chain complex corresponding to $L(X_*)$), and
after this only take the homotopy fiber, and iterate (the homotopy
fiber will be $n$-connected now). This process, in a way, breaks the
natural abelianization into pieces, a brutal thing to do one will
admit, all the more so as we start with a beautiful complex with
projective components, and kind of destroy its unmarred harmony by
tearing out of breaking off the most showy part, $\mathrm H_n$ to name
it, which now looks so lost and awkward we really can't just leave it
as it is, we first have to take a projective resolution of it,
choosing it as we may\dots But we won't do all this, will we, and
rather keep abelianization and \eqref{eq:113.Cprime} as God gave them
to us, and take the homotopy fiber (we hope God will give this
too\dots), and repeat the process, without even having to care at any
stage which $\pi_i$'s vanish and which not. Let's call this the
``\emph{soft Postnikov dévissage}'', in contrast to the ``brutal
one''. In describing the process (which of course makes sense in the
discrete context as well as in the schematic one), I implicitly
admitted that the $X_*$ we start with is $1$-connected, or for the
very least has abelian $\pi_1$ (a notion we'll have to come back to,
in the schematic set-up). But we may as well apply it to a (discrete)
$K(G,1)$ type, $G$ any discrete group, then it amounts to taking the
descending filtration of $G$ by iterated commutator groups, which is a
finite filtration if{f} $G$ is solvable. Maybe it would be more
natural still to take the similar descending filtration, suitable for
the study of \emph{nilpotent} groups rather than solvable ones,
with\scrcomment{this second superscript in this equation is rather
  hard to read in the typescript\dots}
\begin{equation}
  \label{eq:113.D}
  G^{(n+1)} = [G, G^{(n-1)}],\tag{D}
\end{equation}
where $[A,B]$ denotes the subgroup of $G$ generated by commutators
\[(a,b)=aba^{-1}b^{-1},\]
with $a$ in $A$, $b$ in $B$. It doesn't seem there is a similar
distinction to make in case we start with a $1$-connected $X_*$, more
specifically if $X_0=X_1=e$. There may be some extra caution needed,
however, when we assume only $X_0=e$ without assuming
$1$-connectedness, even when $\pi_1$ is abelian,\pspage{469} because
of the possibility of operation of $\pi_1$ upon the $\pi_i$'s. Maybe,
when trying to modelize usual homotopy types by complexes of unipotent
bundles over \bZ, we should restrict to homotopy types which are not
only $0$-connected and have abelian $\pi_1$, but moreover with $\pi_1$
operating trivially (or for the very least, in a unipotent way) upon
the higher $\pi_i$'s. At any rate, as soon as $\pi_1$ operates
non-trivially (on itself, or on the higher $\pi_i$'s) there will
presumably be two non-equivalent ways for defining soft Postnikov
dévissage, corresponding to the two standard descending commutator
group series in a discrete group $G$ The more relevant in view of
unipotent schematization would seem to be the ``nilpotent'' one.

Restricting for simplicity to the $1$-connected case $X_0=X_1=e$, I
would expect soft Postnikov dévissage to be the key for an
understanding as well of the behavior of the $\pi_i(X_*)$ modules with
respect to ring extension, as of the full relationship between these
invariants, and the homology invariants $\mathrm H_i(X_*)$.

\starsbreak

I still should have a look upon complexes $X_*$ in $U(k)^\bullet$
satisfying (as always in this game) $X_0=e$, but not necessarily
$X_1=e$. Even when we make the Kan assumption (plus
``\emph{smoothness}'' of the components, by which I mean that their
$k$-linearization is projective), I don't feel too sure yet if they
fit into a good formalism, for instance (when $k=\bZ$) if they satisfy
the ``linearization theorem'' (quasi-isomorphy for the two maps in
\eqref{eq:111.w} p.~\ref{p:454}). If we start for instance with a
group object $G$ of $U(k)$ and let $X_*$ be the corresponding
semisimplicial complex in $U(k)^\bullet$, then we get an isomorphism
\begin{equation}
  \label{eq:113.E}
  \pi_1(X_*) \eqdef \pi_1(X_*(k)) \simeq G(k),\tag{E}
\end{equation}
which shows us that $\pi_1(X_*)$ need not be abelian even when
$k=\bZ$. If we assume that $X_1=G$ is ``of finite presentation''
(namely the projective $k$-module which describes $X_1$ is of finite
presentation), or what amounts to the same, representable by an actual
(group) scheme, it is true, however, that $G$ and hence $\pi_1=G(k)$
is \emph{nilpotent} (this holds for any $k$). It looks an intriguing
question whether $\pi_1$ is nilpotent under the only assumption that
$X_*$ is a smooth Kan complex with $X_1$ a scheme (without assuming
anymore $X$ comes from a group object). At any rate, it follows
from\pspage{470} the Kan condition that $\pi_1$ may be interpreted as
a quotient set of $E \eqdef X_1(k)$ (without having to pass to the
full free group generated by this set), with a set of relations
\begin{equation}
  \label{eq:113.F}
  z_i = x_iy_i, \quad \text{$i$ in $X_2(k)=I$},\tag{F}
\end{equation}
indexed by $X_2(k)$, where
\[i\mapsto x_i, \quad i\mapsto y_i, \quad i\mapsto z_i\]
are the three boundary maps, remembering moreover the Kan condition
that the three maps
\[\begin{tikzcd}[cramped,sep=small]
  I\ar[r]\ar[r,shift left=1.5]\ar[r,shift right=1.5] & E\times E,
\end{tikzcd}\quad
i\mapsto(x_i,y_i), \quad i\mapsto(x_i,z_i),\quad i\mapsto(y_i,z_i)\]
are \emph{surjective}, which implies indeed that any element of the
group $\pi$ described by the set of generates $E$ and relations
\eqref{eq:113.F} comes from an element in $E$. Replacing $k$ by any
$k$-algebra $k'$, we see that we have a presheaf
\[k' \mapsto \pi_1(X_*(k')) = \pi_1(X_*\otimes_k k')\]
on the category $\Aff_{/k}$ of affine schemes over $k$, with values in
the category of groups, which may be viewed (as a presheaf of sets) as
a quotient presheaf of the presheaf on $\Aff_{/k}$ defined by
$X_1$. We feel that this presheaf will fit into a reasonable
``schematic'' set-up, only if it turns out to be a \emph{sheaf}, and
more exactingly still, if this sheaf is isomorphic (as a sheaf of
sets) to one stemming from an object of $U(k)$, i.e., if it is
isomorphic to a sheaf $W(M)$, for suitable $k$-module $M$ (not
necessarily a projective one). If we denote by $G$ this object of
$U(k)$, it will be endowed with a group structure, and \emph{it is
  this group object of $U(k)$, rather than just the set-theoretic
  group of its sections}, i.e., \emph{of $k$-valued ``points'', which
  merits to be viewed as the ``true'' $\pi_1(X_*)$}. To say it
differently, whereas the higher $\pi_i(X_*)\eqdef \pi_i(X_*(k))$ (for
$i\ge2$), in the cases considered so far, should be viewed as being
not mere abelian groups, but moreover endowed with a natural
$k$-module structure, in the case when $i=1$, i.e., for the
fundamental group $\pi_1(X_*)$, the natural structure to expect on
this (possibly non-commutative) group is a ``unipotent schematic''
structure, namely essentially a pointed ``parametrization'' of this
group by elements of a suitable $k$-module $M_1$, in such a way that
the composition law is expressed in terms of a polynomial law, making
sense therefore not only for $k$-valued points, i.e., for elements of
$M$, but for $k'$-valued points as well (for any $k$-algebra $k'$),
namely defining a group law on\pspage{471} $W(M)(k') = M\otimes_k
k'$. If henceforth we denote by $\pi_1(X_*)$ this group object of
$U(k)$, the relevant formula now is
\begin{equation}
  \label{eq:113.G}
  \pi_1(X_*(k')) \simeq \pi_1(X_*)(k'),\tag{G}
\end{equation}
a group isomorphism functorial with respect to variable $k$-algebra
$k'$, which will imply the corresponding isomorphism
\begin{equation}
  \label{eq:113.Gprime}
  \pi_1(X_*\otimes_k k') \simeq \pi_1(X_*)\otimes_k k'\tag{G'}
\end{equation}
of groups objects in $U(k')$, i.e., formation of the ``schematic''
$\pi_1$ is compatible with ground ring extension $k\to k'$ (provided
$\pi_1$ exists, for a given $X_*$).

We will expect the map of passage to quotient
\begin{equation}
  \label{eq:113.H}
  X_1\to \pi_1(X_*)=G\tag{H}
\end{equation}
to be ``epimorphic'' in a very strong sense, stronger even than just
in the sense of presheaves, the first thought that comes to mind here
is that it should be a ``\emph{submersion}'', in the sense suggested in
yesterday's reflections in connection with the description of the Kan
condition. If, however, we want to be able to get for
$\pi_1(X_*(k))=\pi_1(X_*)(k)$\scrcomment{in the typescript there is
  here a tautological equation, but I think this is what was
  meant\dots} \emph{any} abelian group beforehand, in the case $k=\bZ$
say, without demanding that it be a projective $k$-module, and still
get it via an $X_*$ with \emph{smooth} components, this shows that
when defining a notion of ``submersion'' for objects of $U(k)$ which
may not be smooth, we should not be quite as demanding as suggested
yesterday (cf.\ page~\ref{p:463}), but find a definition which will
include also any map $X\to Y$ coming from an epimorphism $M\to N$ of
$k$-modules (which will allow us to take $X_1$ as associated to a
projective $k$-module admitting the given $\pi_1$ as its
quotient). One idea that comes to mind here, is to take this property
as the \emph{definition} of a submersion, as an arrow in $U(k)$ which
is isomorphic to one obtained from an epimorphism in $\AbOf_k$. This,
of the three definitions that have come to my mind so far for this
notion, is the one which looks the most convincing to me. I wouldn't
expect too much from a complex $X_*$, even a smooth one and satisfying
the Kan condition, unless (in terms of the three boundary maps from
$X_2$ to $X_1$) it gives rise, as just explained, to a group object
$G=\pi_1(X_*)$ in $U(k)$, together with a \emph{submersion}
\eqref{eq:113.H}. Thus, definitely, when defining a schematic model
category $\scrM_n(k)$ of $n$-connected ss~complexes of unipotent
bundles over $k$, I feel like insisting in case $n=0$ at least upon
this extra condition (plus of course $X_0=e$). The still more
stringent condition one may think of, in order to have schematic
models as\pspage{472} close as one may wish to $k$-modules, is to
demand that moreover $G$ is isomorphic to the object of $U(k)$ defined
by a $k$-module $M_1$, the group law moreover coming from the addition
law in $M$. This condition is stronger still than merely demanding
that $G$ be commutative, even when $M_1$ is free of rank one, because
one knows that over a non-perfect field $k$ there may be ``forms'' of
the additive group $\bG\suba$ which are not isomorphic to $\bG\suba$;
presumably, there should be similar examples over \bZ{} too, with rank
larger than one, however.

We feel, however, that the case when $G=\pi_1(X_*)$ is a non-linear or
even a non-abelian group object of $U(k)$ is still worthy of
interest. The first test it would seem, to check if we do have a good
notion indeed, is to see if it does satisfy to the ``linearization
theorem'' in case $k=\bZ$, i.e., the maps \eqref{eq:111.w} on page
\ref{p:454} are quasi-isomorphisms. Another key test, which now makes
sense for arbitrary $k$, is whether for a smooth Kan complex in
$\scrM_0(k)$ (i.e., satisfying the extra assumption involving $G$),
$X_*$ is homotopic to a bundle over $K(G,1)$, with a $1$-connected
fiber, or more specifically, a fiber $Y_*$ satisfying
$Y_0=Y_1=e$. Among other features to expect is a natural operation of
the group object $G$ on the $k$-modules $\pi_i(X_*)$, as well as
$\mathrm H_i(X_*)$, If however we wish, for $k=\bZ$, to use models in
$\scrM_0(\bZ)$ for describing possibly homotopy types with nilpotent
$\pi_1$ say,\footnote{and moreover unipotent action on the $\pi_i$'s}
and devise a corresponding equivalence between suitable homotopy
categories, we should first investigate the question of the
relationship between nilpotent discrete groups, and group objects of
$U(k)$ -- a question already touched on earlier in our reflection on
linearization (see end of section \ref{sec:94}), and of separate
interest.

\hangsection{Outline of a program.}\label{sec:114}%
During these four days of reflection on schematization of homotopy
types, a relatively coherent picture has gradually been emerging from
darkness. How far this image reflects substantial reality, not just
daydreaming, I would be at a loss to tell now. Maybe some substantial
corrections will have to be made still, besides getting in other ideas
for a more complete picture -- I would be amazed at any rate if
everything should turn out as just nonsense! If it doesn't, there is
surely a lot of work ahead to get everything straightened out and
ready-to-use. I will leave it at that I suppose, for the time being --
maybe just finish this digression by a quick review of the set-up, and
of some main questions which have come out.

For\pspage{473} a given ground ring $k$, the basic category we'll use
of ``schematic'' objects over $k$ is the category of \emph{unipotent
  bundles} over $k$, which may be defined as the category of functors
from $\Alg_{/k}$ to \Sets{} isomorphic to functors of the type
\[W(M) = (k' \mapsto M\otimes_k k'),\]
where $M$ is any $k$-module. We do not restrict, here, $M$ to be
projective or flat, as we definitely want to have, for a ring
homomorphism $k\to k'$, a problemless functor ``\emph{restriction of
  scalars}''
\[U(k')\to U(k),\]
inserting in the commutative diagram
\begin{equation}
  \label{eq:114.I}
  \begin{tabular}{@{}c@{}}
    \begin{tikzcd}[baseline=(O.base)]
      \AbOf_{k'}\ar[r,"\text{restr.}"]\ar[d,"W_{k'}"'] &
      \AbOf_k\ar[d,"W_k"] \\
      U(k')\ar[r,"\text{restr.}"] &
      |[alias=O]| U(k)
    \end{tikzcd}.
  \end{tabular}\tag{I}
\end{equation}
Another reason is that we want that the $k$-modules of the type
$\pi_i(X_*)$ which will come out should be eligible for defining
objects in $U(k)$. We are more specifically interested though in
``\emph{smooth}'' objects of $U(k)$, namely those that correspond to
projective $k$-modules. (We prefer to call them ``smooth'' rather than
``projective'', in order to avoid confusion with the notion of a
projective object in the usual categorical sense for $U(k)$.) Another
relevant notion is the notion of a \emph{submersion}, namely a map in
$U(k)$ isomorphic to one coming from an epimorphism $M\to N$ in
$\AbOf_k$. (If the latter can be chosen to have a projective kernel,
we may speak of a smooth submersion.) The ring restriction functor
transforms submersions into submersions, and also smooth objects into
smooth ones provided $k'$ is projective as a module over $k$. We also
have a \emph{ring extension functor} from $U(k)$ to $U(k')$, giving
rise to a diagram (I') similar to \eqref{eq:114.I} above, it
transforms submersions into submersions, smooth objects into smooth
ones.

The smoothness condition is likely to come in in two ways, one is via
\emph{flatness} (we may call an object of $U(k)$ ``flat'' when it
isomorphic to some $W(M)$ with $M$ a flat $k$-module), whereas
projectivity is needed in order to ensure that in certain cases, weak
equivalences are homotopisms. Flatness is the kind of condition which
ensure the validity of ``naive'' universal coefficients formulæ for
homotopy or homology objects, whereas projectivity may be needed in
case of such formulæ for cohomology rather than homology.

The\pspage{474} description I just recalled of $U(k)$ is the one most
intuitive to my mind, other people may prefer the more computational
one on page \ref{p:451} in terms of $\Gamma\uphat_k(M)$ (endowed with
its augmentation to $k$ and its diagonal map), which is of importance
in its own right. It shows the existence of a canonical
\emph{$k$-linearization functor}
\begin{equation}
  \label{eq:114.J}
  L: U(k)\to\AbOf_k,\tag{J}
\end{equation}
giving rise to the commutative diagram (up to can.\ isom.)
\begin{equation}
  \label{eq:114.Jprime}
  \begin{tabular}{@{}c@{}}
    \begin{tikzcd}[baseline=(O.base),column sep=tiny]
      \AbOf_k\ar[rr,"W"]\ar[dr,"\Gamma\uphat_k"'] & &
      U(k)\ar[dl,"L"] \\
      & |[alias=O]| \AbOf_k
    \end{tikzcd},
  \end{tabular}\tag{J'}
\end{equation}
where
\[\Gamma\uphat_k(M) = \prod_{i\ge0} \Gamma^i_k(M).\]
This linearization is not quite compatible with ring extension, it
becomes so only when we view it as a functor with values, not just in
the category $\AbOf_k$ of $k$-modules, but of separated and complete
linearly topologized $k$-modules, the ring extension functor for these
being the \emph{completed} tensor product. This is a little (or big?)
technical drawback for this notion of linearization. We have a
canonical embedding
\begin{equation}
  \label{eq:114.K}
  x\mapsto\exp(x) : X\to W\uphat(L(X))\quad \text{(cf.\ p.\
    \ref{p:450}, \eqref{eq:111.n})}\tag{K}
\end{equation}
functorial in $X$, where for a topological $k$-module $M$ as above,
described as a filtering inverse limit of discrete ones $M_i$, we
define
\begin{equation}
  \label{eq:114.Kprime}
  W\uphat(M) = (k'\mapsto\varprojlim_i M_i\otimes_k k').\tag{K'}
\end{equation}
The map \eqref{eq:114.K} has a universal property with respect to all
possible maps $X\to W\uphat(M)$ with $M$ a linearly topologized
separated and complete $k$-module, which accounts for its role as
``linearization''. It should be noted here that the map
\eqref{eq:113.C} on page \ref{p:467} doesn't quite exist, we have
corrected this point here -- definitely we cannot in \eqref{eq:114.K}
replace $W\uphat$ by $W$. Of course, linearization $L$ (or its variant
$L\subpt$) doesn't commute in any sense whatever to restriction of
ground ring.

The image of $X$ in $W\uphat(L(X))$ is characterized by the simple
formulæ \eqref{eq:111.q} p.\ \ref{p:451}. Maps from $X$ to $Y$ may be
described as just continuous $k$-linear maps from $L(X)$ to $L(Y)$,
compatible with augmentations and diagonal maps.

We'll\pspage{475} more specifically work in the category
$U(k)^\bullet$ of \emph{pointed} objects of $U(k)$, namely objects
endowed with a section over the final object $e$, the so-called
\emph{pointed unipotent bundles}. We now have a functor
\begin{equation}
  \label{eq:114.L}
  W^\bullet : \AbOf_k \to U(k)^\bullet\tag{L}
\end{equation}
deduced from $W$ using the fact that $W(0)=e$, and a ``\emph{pointed
  linearization functor}''
\begin{equation}
  \label{eq:114.M}
  L^\bullet\text{ or }L\subpt : U(k)^\bullet \to \AbOf_k,\tag{M}
\end{equation}
giving rise to a commutative diagram similar to \eqref{eq:114.Jprime}
\begin{equation}
  \label{eq:114.Mprime}
  \begin{tabular}{@{}c@{}}
    \begin{tikzcd}[baseline=(O.base),column sep=tiny]
      \AbOf_k\ar[rr,"W^\bullet"]\ar[dr,"\Gamma_k"'] & &
      U(k)^\bullet\ar[dl,"L^\bullet"] \\
      & |[alias=O]| \AbOf_k
    \end{tikzcd};
  \end{tabular}\tag{M'}
\end{equation}
the notation $L^\bullet$ seems here the most coherent one, but may
bring about confusion with the similar notation for some cochain
complex say, therefore we had first used the alternative notation
$L\subpt$, to which one may still come back if needed. This time the
functor $L^\bullet$ commutes to ring extension without any grain of
salt. We have of course a canonical embedding
\begin{equation}
  \label{eq:114.N}
  L^\bullet(X)\hookrightarrow L(X)\tag{N}
\end{equation}
defined via the corresponding embedding for an object $M$ of $\AbOf_k$
\begin{equation}
  \label{eq:114.Nprime}
  \Gamma_k(M)\hookrightarrow \Gamma\uphat_k(M),\tag{N'}
\end{equation}
by which $L(X)$ may be viewed as the completion of $L^\bullet(X)$ with
respect to the topology it induces on it, which is a \emph{canonical
  topology} on $L^\bullet(X)$. Maps in $U(k)^\bullet$ correspond to
\emph{$k$-linear maps}
\[L^\bullet(X)\to L^\bullet(Y)\]
which are moreover \emph{continuous, and compatible with
  coaugmentation} (i.e., transforms $1$ into $1$) \emph{as well as
  with augmentations and diagonal maps}. I wonder if there is any
simple characterization of submersions in $U(k)^\bullet$ in terms of
the corresponding map between the linearizations. At any rate, an
object $X$ of $U(k)^\bullet$ is smooth resp.\ flat if{f}
$L^\bullet(X)$ is a projective resp.\ a flat $k$-module.

For any natural integer, we want now to define a model category
\begin{equation}
  \label{eq:114.O}
  \scrM_n(k)\subset\bHom(\Simplexop, U(k)^\bullet),\tag{O}
\end{equation}
which should be a full subcategory of the category of
semisimplicial\pspage{476} objects in $U(k)^\bullet$. We'll get a
functor
\begin{equation}
  \label{eq:114.P}
  X_* \mapsto X_*(k) :
  \scrM_n(k)\to\bHom(\Simplexop,\mathrm{Sets}^\bullet) \tag{P} 
\end{equation}
from this category to the category of semisimplicial pointed sets. For
$n\ge1$, the only condition, it seems, to impose upon $X_*$ in the
second member of \eqref{eq:114.O}, i.e., upon a semisimplicial pointed
unipotent bundle over $k$, in order to belong to $\scrM_n(k)$, is
\begin{equation}
  \label{eq:114.Qn}
  X_i = e\quad\text{for}\quad i\le n.\tag{Q\textsubscript{n}}
\end{equation}
This, for $n=0$, reduces to the common condition
\begin{equation}
  \label{eq:114.Qzero}
  X_0=e,\tag{Q\textsubscript{0}}
\end{equation}
which definitely is not enough, though, to get a category of
``models'' $\scrM_0(k)$ whose objects should have the kind of
properties we are after. There are various kinds of extra restrictions
one may want to impose, according to the type of situations one wants
to describe, some hints along these lines are given on pages
\ref{p:469}--\ref{p:472}. For a preliminary study, the case $n\ge1$,
and more specifically, the case $n=1$, is quite enough, the latter
corresponding to the restrictions
\begin{equation}
  \label{eq:114.Qone}
  X_0=X_1=e.\tag{Q\textsubscript{1}}
\end{equation}

From $P$ we get a functor
\begin{equation}
  \label{eq:114.Pprime}
  \scrM_n(k)\to\HotOf_n^\bullet =
  \begin{tabular}[t]{@{}l@{}}
    category of pointed $n$-connected\\
    homotopy types,
  \end{tabular}\tag{P'}
\end{equation}
we define a map in $\scrM_n(k)$ to be a \emph{weak equivalence} if its
image by \eqref{eq:114.P} is, i.e., its image by \eqref{eq:114.Pprime}
is an isomorphism, and localizing by weak equivalences we get the
category
\begin{equation}
  \label{eq:114.R}
  \HotOf_n(k)\tag{R}
\end{equation}
of ``$n$-connected schematic homotopy types over $k$'', together with
a ``sections functor'' induced by \eqref{eq:114.P}
\begin{equation}
  \label{eq:114.Rprime}
  \HotOf_n(k)\to\HotOf_n.\tag{R'}
\end{equation}
One main point in our definitions is that \emph{we hope this functor
  to be an equivalence of categories, in the case when $k=\bZ$}, and
of course $n\ge1$.

The description just given of categories $\HotOf_n(k)$ is suitable for
defining functors of restriction of ground ring for $k\to k'$
\begin{equation}
  \label{eq:114.S}
  \HotOf_n(k')\to\HotOf_n(k),\tag{S}
\end{equation}
compatible with the sections functor \eqref{eq:114.Rprime} for $k$ and
$k'$. It isn't directly suited, though, for describing ring extension
-- as a matter of\pspage{477} fact, ring extension for homotopy types
(an operation of greater interest than ring restriction surely) is
\emph{not} expressed, in general, by just performing the trivial ring
extension operation
\[X_*\mapsto X_*\otimes_k k'\]
on models in $\scrM_n(k)$, unless we assume $k'$ to be flat over $k$
say -- but even in this case it is by no means clear a priori that the
operation above transforms weak equivalences into weak
equivalences. This is very clearly shown by the linear analogon, the
categories $\scrM_n(k)$ being replaced by the categories of chain
complexes in $\AbOf_k$ say, or by $\Comp^-(\AbOf_k)$ or the like. In
order to correctly describe ground ring extension on homotopy types,
we'll have first to take a suitable ``resolution'' of $X_*$, namely
replace $X_*$ by some $K_*$ say, endowed with a weak equivalence
\[K_*\to X_*,\]
and $K_*$ satisfying some extra assumptions. Maybe flatness of the
components would be enough here. For other purposes, we may have to
use resolutions which are even smooth (componentwise), or which
satisfy a suitable \emph{Kan condition} (or a type outlined on page
\ref{p:463}), or both. Our expectation is that, when we restrict to
the subcategory
\[\mathrm{sK}\scrM_n(k)\]
of the model category $\scrM_n(k)$ made up with smooth Kan complexes,
\emph{that the category $\HotOf_n(k)$ may be described simply in terms
  of such sK-complexes ``up to homotopy''}, as usual. If this is so,
the ground ring extension functor follows trivially from a
corresponding functor on the sK-model categories
\begin{equation}
  \label{eq:114.T}
  \sKM_n(k)\to\sKM_n(k'), \quad X_*\mapsto X_*\otimes_k k',\tag{T}
\end{equation}
hence
\begin{equation}
  \label{eq:114.Tprime}
  \HotOf_n(k)\to\HotOf_n(k').\tag{T'}
\end{equation}

From the sections functor \eqref{eq:114.Rprime} we get homotopy
invariants $\pi_i$ for an object in $\HotOf_n(k)$, but the relevant
$k$-module structure on these is not apparent on this definition. We
have a better hold, via linear algebra over $k$, upon homology
invariants of $X_*$, which are $k$-linear objects, and are definitely
distinct (unless $k=\bZ$) from the homology invariants of $X_*(k)$,
which for general $k$ are definitely of little interest it seems. The
definition of homology goes via the pointed linearization functor
\eqref{eq:114.M}
\begin{equation}
  \label{eq:114.U}
  \left\{\begin{tabular}{@{}l@{}}
    $\LH_\bullet(X_*) \eqdef L\subpt(X_*)$, viewed as an object
      in $\D_\bullet(\AbOf_k)$, \\
    $\mathrm H_i(X_*) \eqdef \mathrm H_i(\LH_\bullet(X_*)) =
      \pi_i(L\subpt(X_*))$ in $\AbOf_k$,
  \end{tabular}\right.\tag{U}
\end{equation}
where\pspage{478} the $\mathrm L$ in $\LH_\bullet$ suggests that we
are taking something similar to a total left derived functor, and
where definitely in the right-hand member we had to write $L\subpt$
and not $L^\bullet$, in order not to get sunk into a morass of
confusion! In the formulæ \eqref{eq:114.U} we should assume however
that $X_*$ is a smooth Kan complex, which will imply (if indeed
$\HotOf_n(k)$ may be described in terms of $\sKM_n(k)$ as said above)
that $\LH_\bullet$ may be viewed as a functor
\begin{equation}
  \label{eq:114.Uprime}
  \HotOf_n(k)\to\D_\bullet(\AbOf_k),\tag{U'}
\end{equation}
and likewise the $\mathrm H_i$'s are functors from $\HotOf_n(k)$ to
$\AbOf_k$. In order to compute these homology invariants for an
arbitrary complex in $\scrM_n(k)$, we'll first have to resolve it by a
$\mathrm{sK}$ complex, and then apply \eqref{eq:114.U}.

We expect \emph{that a map $X_*\to Y_*$ in $\scrM_n(k)$ is a weak
  equivalence if{f} the corresponding map for $\LH_\bullet$ is a
  quasi-isomorphism}, in other words we expect the functor
\eqref{eq:114.Uprime} to be ``conservative'': a map in the first
category is an isomorphism if{f} its image in the second one is. A
second main feature we expect from linearization, is that in the case
$k=\bZ$ it corresponds to the usual abelianization of homotopy
types. This statement, when made more specific as in \eqref{eq:111.w}
page \ref{p:454}, decomposes into two distinct ones. One is of
significance over an arbitrary ring $k$, and states that \emph{for a
  $\mathrm{sK}$ complex $X_*$, the inclusion} (coming from
\eqref{eq:114.N})
\begin{equation}
  \label{eq:114.V}
  L\subpt(X_*)\to L(X_*)\tag{V}
\end{equation}
\emph{from $L\subpt$ into its completion, when viewed as a map of
  chain complexes in $\AbOf_k$} (using the simplicial differential
operator, or passing to the corresponding ``normalized'' chain
complexes first) \emph{is a quasi-isomorphism}. Whether this is always
so or not, or whether noetherian conditions on $k$ or some finiteness
conditions for the components of $X_*$ are needed, looks like a rather
standard question of linear homological algebra! On the other hand,
using the exponential embedding \eqref{eq:114.K} for sections, we get
another map of semisimplicial $k$-modules
\begin{equation}
  \label{eq:114.Vprime}
  k^{(X_*(k))} \to L(X_*),\tag{V'}
\end{equation}
and here \emph{the question again} (the expectation I might say?)
\emph{is whether this is a quasi-isomorphism}. This would just mean
(if coupled with quasi-isomorphy of \eqref{eq:114.V}) that the
homology invariants \eqref{eq:114.U} are just the usual homology
invariants of the discrete homotopy type modeled by $X_*(k)$, but with
coefficients not in \bZ, but in $k$. We certainly do expect this to be
true for $k=\bZ$ -- which was the content of ``question
\ref{q:111.1}''\pspage{479} on page \ref{p:452} (taken up again on
page \ref{p:454} and following). Of course, in case we don't assume
$k=\bZ$, writing $\bZ^{(X_*(k))}$ instead of $k^{(X_*(k))}$ as we did
in \eqref{eq:111.w} p.~\ref{p:454} now looks kind of silly, and the
idea in this $k$-linear context to take $k$-valued homology of
$X_*(k)$ rather than \bZ-valued one is evident enough! However, I was
confused by the misconception that the \emph{internal} homology of
$X_*(k)$ should carry $k$-linear structure, as this was what I
expected too for the invariants $\pi_i(X_*(k))$ (which seems to turn
out to be correct). This misconception was corrected a few pages later
(p.~\ref{p:457}) but still I kept dragging along the silly first
member of \eqref{eq:111.w}. Anyhow, it now just occurs to me that
except in case $k=\bZ$, it is definitely \emph{false} that
\eqref{eq:114.Vprime} is a quasi-isomorphism, except in some wholly
trivial cases. Indeed, let $\mathrm H_i$ be the first non-vanishing
homology invariant \eqref{eq:114.U} of $X_*$ (or more safely still,
take $i=n+1$), then we definitely expect to have a canonical
isomorphism of $k$-modules
\begin{equation}
  \label{eq:114.Vbis}
  \pi_i(X_*(k)) \simeq \mathrm H_i(X_*) = \mathrm H_i\tag{V}
\end{equation}
but we equally have by Hopf's theorem, as the lower $\pi_j$'s of
$X_*(k)$ are zero
\[\pi_i(X_*(k)) = \mathrm H_i(X_*(k),\bZ),\]
hence
\begin{equation}
  \label{eq:114.Vprimebis}
  \mathrm H_i \simeq \mathrm H_i(S,\bZ), \quad\text{where
    $S=X_*(k)$,}\tag{V'}
\end{equation}
which is not compatible with the guess that $\mathrm H_i\simeq\mathrm
H_i(S,k)$ ($\simeq\mathrm H_i\otimes_\bZ k$). Presumably, the
isomorphism \eqref{eq:114.Vprimebis} above is induced by the first map
in \eqref{eq:111.w} above, but (except for $k=\bZ$) we should expect
in \eqref{eq:111.w} to have an isomorphism only for the lowest
dimensional homology groups which are occurring in the two first
members. Anyhow, it appears after all that this map in
\eqref{eq:111.w} is the more reasonable one compared to
\eqref{eq:114.Vprime} above, as it yields an isomorphism on homology
in the key dimension $n+1$, whereas \eqref{eq:114.Vprime} apparently
will practically \emph{never} give an isomorphism.

What is mainly lacking still in this review of the expected main
features of schematization of homotopy types, is description of the
$k$-module structure on the homotopy groups
\begin{equation}
  \label{eq:114.W}
  \pi_i(X_*) \eqdef \pi_i(X_*(k)),\tag{W}
\end{equation}
or preferably still, a direct description of those invariants as
$k$-modules, working within the model category $\scrM_n(k)$. This, as
suggested in yesterday's notes (p.~\ref{p:464}), may be achieved
\emph{by developing a theory of Postnikov dévissage within
  $\scrM_n(k)$ and using \eqref{eq:114.V}} in order to pull ourselves
by the bootstraps, defining homotopy finally in terms of
homology.\pspage{480} At this point it should be noted that the
dévissage we'll have to use here is the ``brutal'' one, which we
frowned upon earlier today! To develop such a formalism, it seems
essential to work with smooth Kan complexes and projective resolutions
of the $k$-modules $\pi_i$ as they appears one by one. Whether we want
to describe ``hard'' or ``soft'' Postnikov dévissage (see
p.~\ref{p:468} for the latter), one common key step is the
\emph{linearization map} coming from the exponential map
\eqref{eq:114.K} applied componentwise
\begin{equation}
  \label{eq:114.Wbis}
  X_* \to W\uphat(L(X_*)),\tag{W}
\end{equation}
which we would like to look upon as defining a homotopy class of maps
in $\scrM_n(k)$
\begin{equation}
  \label{eq:114.Wq}
  X_* \to W(L\subpt(X_*)),\tag{W?}
\end{equation}
where the second member moreover is endowed with its natural abelian
group structure (its components are abelian group objects of $U(k)$
and the simplicial maps are additive). To pass from
\eqref{eq:114.Wbis} to \eqref{eq:114.Wq}, it is felt that the
essential step is that \eqref{eq:114.V} above be a quasi-isomorphism,
hence, applying the functor $W$, we should get a weak equivalence,
hence an isomorphism in the derived category $\HotOf_n(k)$, hence
\eqref{eq:114.Wbis} implies \eqref{eq:114.Wq}. The main flaw in this
``argument'' comes from the $W\uphat$ in the second member of
\eqref{eq:114.Wbis}, which isn't quite the same as $W$
definitely. Thus, some further amount of work will be needed,
presumably, to get \eqref{eq:114.Wq} from \eqref{eq:114.Wbis}. Of
course, we can't possibly just keep \eqref{eq:114.Wbis} as it is, as
for getting dévissage we need a map in $\scrM_n(k)$, whereas the map
\eqref{eq:114.Wbis} is just a map of semisimplicial sheaves on
$\Aff_{/k}$, where the second member is \emph{not} in $\scrM_n(k)$,
i.e., its components are not in $U(k)$. Once we got \eqref{eq:114.Wq}
factoring \eqref{eq:114.Wbis} up to homotopy (NB\enspace of course we
assume $X_*$ to be an $\mathrm{sK}$ complex in all this), we still
need a reasonable notion of homotopy fibers of maps in $\scrM_n(k)$,
in order to push through the inductive step.

Thus, a large part of the weight of the work ahead may well lie upon
\emph{developing the standard homotopy constructions within the model
  category $\scrM_n(k)$}, as contemplated on page \ref{p:464}. This
should be fun, if it can be done indeed! One difficulty here seems to
be that Quillen's standard machines won't work, not ``telles quelles''
at any rate, because of the category $U(k)$ failing to be stable under
finite limits -- it doesn't even have fiber products. But I think I'll
stop my ponderings on schematization here\dots

\bigbreak

\presectionfill\ondate{28.8.}\pspage{481}\par

\hangsection{\texorpdfstring{$L(X)$}{L(X)} as the pro-quasicoherent
  substitute for
  \texorpdfstring{$\scrO_k(X)$}{Ok(X)}.}\label{sec:115}%
For the last four days, while reflecting on ``schematization'', each
time I think I am going to be through with that unforeseen green apple
within an hour or two, and get back to ``l'ordre du jour'' -- and
overnight something else still appears I feel I should still look into
just a little; and there I am again, sure enough, with some extra
reflection on ``schematic linearization'' which I hadn't quite
understood yet, it appears to me now. These last days I had given up
numbering formulas as usual by Arab ciphers (1), (2), etc., as I
didn't want to ``cut'' the numbering of that unending ``review'' of
section \ref{sec:104} to \ref{sec:109} which wasn't quite finished
yet, got it only till formula \eqref{eq:109.136}. But now I will stop
this nonsense with numberings (a), (b) and (A), (B), after all even if
there are in-between ``Arab'' formulas now, this doesn't prevent me,
when it comes to it, to start a ``review'' section with formula (137)
and go on till (1000) if I like\dots  And now to the schematic
linearization functor again, for unipotent bundles!

When writing up that schematization program yesterday, some technical
difficulties appeared at the end (see page before) for a proper
understanding of the relationship between the two linearizations
$L\subpt$ and $L$, in order to define, in a suitable derived category,
a map
\[X_*\to W(L\subpt(X_*))\]
using the canonical term-by-term exponential map
\[X_*\to W\uphat L(X_*).\]
It seems to me that the exact significance of the objects $L(X_*)$ or
$W\uphat L(X_*)$ isn't quite understood yet, and that the confusion
which occurred between which usual kind of linearization we should
compare this with, whether $X\mapsto \bZ^{(X)}$ as I did first, or
$X\mapsto k^{(X)}$ as it occurred to me yesterday (pages
\ref{p:478}--\ref{p:479}) , is quite typical of this lack of
understanding. It now occurred to me that neither term, for a general
ground ring $k$ (namely, not assuming $k=\bZ$), is reasonable, whereas
the reasonable ``usual'' kind of linearization comparing with $L(X)$
(when $X$ is a unipotent bundle or a ss~complex of such) is
\begin{equation}
  \label{eq:115.1}
  X\mapsto \scrO_k^{(X)},\tag{1}
\end{equation}
where $\scrO_k$ is the basic \emph{quasicoherent} sheaf of rings over
$k$, i.e., over $\Aff_{/k}$, given by the tautological functor
\begin{equation}
  \label{eq:115.2}
  \scrO_k: (\Aff_{/k})\op \equeq \Alg_{/k} \to \Rings, \quad
  k'\mapsto k',\tag{2}
\end{equation}
associating\pspage{482} to any affine scheme $S=\Spec(k')$ over $k$
the ring of sections of its usual Zariski structure sheaf, which ring
is canonically isomorphic to $k'$ itself! The operation
\eqref{eq:115.1} is the usual linearization operation with respect to
this sheaf of rings, working in the topos of fpqc sheaves of sets over
$k$ which we described at some length in section \ref{sec:111}
(p.~\ref{p:447}). As I was fearing that working in such a thing would
cause anguish to a number of prospective readers, I took pains to
translate unipotent bundles from the geometric language which is the
suggestive one, to the language of commutative algebra which is more
liable to hide than to disclose geometrical meaning; so much so that
in the process I myself lost contact somewhat with the geometric
flavor, and more specifically still with this basic fact, that in our
context of unipotent bundles and complexes of such, the ``natural''
coefficients for cohomology (such as the $\mathrm
H^{n+2}(X(n)_*,\pi_{n+1})$ groups occurring in Postnikov dévissage)
are by no means ``\emph{discrete}'' ones such as \bZ{} or $k$, but
\emph{quasi-coherent} ones, namely provided by quasi-coherent sheaves
of $\scrO_k$-modules or complexes of such. Thus, in the above
Postnikov obstruction group, $\pi_{n+1}$ does \emph{not} stand as a
constant group of coefficients (if it was, this would drag us into the
niceties and difficulties of étale cohomology for the components
$X(n)_i$ of the semisimplicial unipotent bundle $X(n)_*$); but using
the $k$-module structure of $\pi_{n+1}$ for defining a
\emph{quasi-coherent} sheaf of modules $W(\pi_{n+1})$ ``over $k$'',
i.e., over $\Aff_{/k}$, it is this ``continuous'' sheaf (or ``vector
bundle'') over $k$, lifted of course to the various components
$X(n)_i$, which yields the correct answer. This was kind of clear in
my mind the very first day when I started reflecting on
schematization, even before introducing formally unipotent bundles
(pages \ref{p:443}--\ref{p:444}), but this instinctive understanding
later became dulled somewhat, largely due, it seems to me, to the
concession I had made to algebra, giving up to some extent the
language of geometry.

Let's recall that the operation \eqref{eq:115.1} may be defined as the
solution of a universal problem, namely sending the non-linear object
$X$ into a ``linear'' one, namely into a sheaf of $\scrO_k$-modules
(or a $\scrO_k$-module, as we'll simply say). This is expressed by a
canonical map of sheaves of sets
\begin{equation}
  \label{eq:115.3}
  X \to \scrO_k^{(X)}\tag{3}
\end{equation}
(which I am tempted to call the \emph{``exponential'' map} for $X$,
and denote by a corresponding symbol such as $\exp_X$), giving rise,
for every module (over $\scrO_k$) to a corresponding map which is
\emph{bijective}\pspage{483}
\begin{equation}
  \label{eq:115.4}
  \Hom_{\scrO_k}(\scrO_k^{(X)},F) \simeq
  \Hom(X,F)\quad(\simeq\Gamma(X,F_X)),\tag{4} 
\end{equation}
where in the last member (included as a more geometric interpretation
of the second) $F_X$ denotes the restriction of $F$ to the object $X$,
more accurately to the topos (or site) induced on $X$ by the ambient
topos (or site) we are working in. Thus, we may indeed view
\eqref{eq:115.1} as the most perfect notion of linearization, as far
as generality goes -- it makes sense of course in any ringed topos
(without even a commutativity assumption!). The only trouble is that,
even for such a down-to-earth $X$ as a unipotent bundle, the standard
affine line $E_k^1$ say, the sheaf $\scrO_k^{(X)}$ in \eqref{eq:115.1}
is not quasi-coherent and therefore not too amenable it seems to
computations -- thus, we get easily from \eqref{eq:115.4} a canonical
map (for general $X$)
\[k^{(X(k))} \to \Gamma(k,\scrO_k^{(X)}) \quad(=\scrO_k^{(X)}(k))\]
(where the $\Gamma$ in the second member denotes sections over $k$,
i.e., value of a functor on $\Alg_{/k}$ on the initial object $k$, and
remembering in the first member that the ring of sections of $\scrO_k$
is $k$), but I would be at a loss to make a guess as for reasonable
conditions for this map to be an isomorphism! This may seem a
prohibitive ``contra'' against using at all such huge sheaves as
$\scrO_k^{(X)}$, the point though is that in most questions where such
linearizations are introduced (mainly questions where interest lies in
computing cohomology invariants), one is practically never interesting
in taking the groups of sections of these, but rather in looking at
their maps into sheaves of modules $F$ precisely, which is achieved by
\eqref{eq:115.4}, or taking more generally their $\Ext^i$ with such an
$F$, which is achieved by the similar familiar formula
\begin{equation}
  \label{eq:115.5}
  \Ext^i_{\scrO_k}(\scrO_k^{(X)},F) \simeq \mathrm H^i(X,F_X),\tag{5}
\end{equation}
more neatly
\begin{equation}
  \label{eq:115.5prime}
  \RHom_{\scrO_k}(\scrO_k^{(X)},F) \simeq \RGamma_X(F_X),\tag{5'}
\end{equation}
valid of course again for any ringed topos. In the present context
however, the ``coefficients'' $F$ we are interested in, as was just
pointed out, are not arbitrary $\scrO_k$-modules, but rather
\emph{quasi-coherent} ones. Thus, if we get a variant of
\eqref{eq:115.1}
\begin{equation}
  \label{eq:115.6}
  X\to L(X)\tag{6}
\end{equation}
with $L(X)$ some \emph{quasi-coherent} module, giving rise to
\eqref{eq:115.4}, or even to \eqref{eq:115.5} and
\eqref{eq:115.5prime}, this would be for us a perfectly good
substitute for \eqref{eq:115.3}, which would deserve the name of a
``\emph{quasicoherent envelope}'' of $X$. Of\pspage{484} course, this
module $L(X)$ would be unique up to unique isomorphism, as the
solution of a universal problem embodied by \eqref{eq:115.4}, namely
as the quasi-coherent module representing the functor
\begin{equation}
  \label{eq:115.7}
  F\mapsto\Hom(X,F)\simeq\Gamma(X,F_X)\tag{7}
\end{equation}
on the category this time of all \emph{quasicoherent}
$\scrO_k$-modules.

For the unipotent schematization story, we are more specifically
interested in the case when $X$ comes from a quasicoherent module
itself, by forgetting its module structure. Now, as well-known, the
functor
\begin{equation}
  \label{eq:115.8}
  M\mapsto W(M):\AbOf_k=\kMod \to
  \begin{tabular}[t]{@{}l@{}}
    category of quasicoherent\\
    $\scrO_k$-modules, $\Qucoh(k)$ say
  \end{tabular}\tag{8}
\end{equation}
is an equivalence of categories. Thus, for $X$ defined by such an $M$,
the question of representability of \eqref{eq:115.7} within the
category of quasicoherent modules, amounts to the similar question in
$\AbOf_k$ for the functor
\begin{equation}
  \label{eq:115.7prime}
  N\mapsto\Hom(W(M),W(N)),\tag{7'}
\end{equation}
where the $\Hom$ denotes homomorphisms of sheaves \emph{of sets} of
course. Now, as suggested first, somewhat vaguely still, in section
\ref{sec:111} (page \ref{p:450}), we have an alternative expression of
this functor, via
\begin{multline}
  \label{eq:115.9}
  \Hom(W(M),W(N))\simeq\Homcont_k(\Gamma\uphat_k(M),N)\\
  \simeq \varinjlim_i \Hom_k(\Gamma_k(M)(i),N),\tag{9}
\end{multline}
where in the second member, $\Homcont_k$ denotes the set of
$k$-ho\-mo\-mor\-phisms which are continuous on
\begin{equation}
  \label{eq:115.10}
  \Gamma\uphat_k(M) = \prod_{i\ge0}\Gamma^i_k(M)\tag{10}
\end{equation}
(endowed with the product of discrete topologies), and in the third we
have written
\begin{equation}
  \label{eq:115.10prime}
  \Gamma_k(M)(i) = \prod_{j\le i}\Gamma_k^j(M)\tag{10'}
\end{equation}
for the product of the $i$ first factors occurring in
\eqref{eq:115.10}. The map \eqref{eq:115.9} is deduced in the evident
way from the exponential map
\begin{equation}
  \label{eq:115.11}
  M\to W\uphat\Gamma\uphat_k(M) \eqdef \varprojlim_i W(\Gamma_k(M)(i)).\tag{11}
\end{equation}
(NB\enspace The relation between the description \eqref{eq:115.9} of
maps $W(M)\to W(N)$ with the description given p.~\ref{p:451} in terms
of maps \eqref{eq:111.p} from $\Gamma\uphat_k(M)$ to
$\Gamma\uphat_k(N)$, is by associating to such a map $f$ its
composition with the projection of the target upon its factor
$N$\dots). An incorrect way of expressing \eqref{eq:115.9}, which I
slipped into in section \ref{sec:111} and kind of remained in till
now, is by pretending that the $k$-module $\Gamma\uphat_k(M)$
represents the\pspage{485} functor \eqref{eq:115.7prime}, this is
clearly false, as we do not have any canonical map from $W(M)$ into
$W(\Gamma\uphat_k(M))$, only into $W\uphat\Gamma\uphat_k(M)$ -- we
have an embedding
\begin{equation}
  \label{eq:115.12}
  W(\Gamma\uphat_k(M))\hookrightarrow W\uphat\Gamma\uphat_k(M),\tag{12}
\end{equation}
but it is clear that in general, the exponential map \eqref{eq:115.11}
does not factor through the first term in \eqref{eq:115.12}. (It does
of course when we look at sections over $k$ only, but when we go over
to a general $k'$, we hit into the trouble that formation of inverse
limits does not commute with ring extension $\otimes_k k'$!) We may
however express \eqref{eq:115.9} by stating that the functor
\eqref{eq:115.7prime} is ``\emph{prorepresentable}'' by the
\emph{pro-object}
\begin{equation}
  \label{eq:115.13}
  \Pro \Gamma_k(M) \eqdef (\Gamma_k(M)(i))_{i\ge0} \quad
  \text{in $\Pro(\AbOf_k)$,}\tag{13}
\end{equation}
this is even a \emph{strict} pro-object (the transition morphisms are
epimorphisms), which implies that the functor it prorepresents is
representable if{f} this projective system is ``essentially constant''
in the most trivial sense, which means here
\[\Gamma^i_k(M)=0\quad\text{for large $i$,}\]
a condition which presumably is satisfied only for $M=0$! Thus, the
``correct'' interpretation of non-pointed quasi-coherent linearization
seems to me to be the corresponding functor, which I would like now to
call $L_k$ or simply $L$ as before but with slightly different
meaning:
\begin{equation}
  \label{eq:115.14}
  L\text{ or }L_k : U(k)\to \Pro(\AbOf_k) \equeq \Pro(\Qucoh(k)),\tag{14}
\end{equation}
where $\Qucoh(k)$ is defined in \eqref{eq:115.8}. In computational
terms, I would like to view $L(X)$ (for a unipotent bundle $X$) to be
a pro-$k$-module, but in terms of geometric intuition, I would see it
rather as a pro-$\scrO_k$-module, i.e., essentially as an inverse
system of quasicoherent modules. It is in these latter terms that the
construction we just gave generalizes to unipotent bundles over
arbitrary ground schemes, not necessarily affine ones. As for the
``pointed'' quasicoherent linearization functor
\begin{equation}
  \label{eq:115.15}
  L\subpt \text{ or } {L_k}\subpt: U(k)^\bullet \to (\AbOf_k) \equeq \Qucoh(k),\tag{15}
\end{equation}
which I like best to view as taking values $L\subpt(X)$ which are
quasicoherent sheaves, it maps into $L$ by
\begin{equation}
  \label{eq:115.16}
  L\subpt(X)\hookrightarrow L(X),\tag{16}
\end{equation}
interpreting objects of a category as special cases of
pro-objects. We'll denote by
\begin{equation}
  \label{eq:115.17}
  WL(X) \in \Pro(U(k))\tag{17}
\end{equation}
the\pspage{486} pro-unipotent bundle defined in terms of $L(X)$ via
the canonical extension
\[\Pro(W)\text{ or simply } W:\Pro(\AbOf_k)\to\Pro(U(k))\]
of $W$ (cf.~\eqref{eq:115.8}) to pro-objects. Thus, instead of the map
\eqref{eq:115.6} which doesn't quite exist, we get a canonical
``exponential'' map
\begin{equation}
  \label{eq:115.18}
  X\to WL(X)\tag{18}
\end{equation}
in $\Pro U(k)$, which has of course little chance to factor through
\begin{equation}
  \label{eq:115.16prime}
  WL\subpt(X)\to WL(X)\tag{16'}
\end{equation}
deduced from \eqref{eq:115.16} by applying $W$. It is via this map
\eqref{eq:115.18} that we may declare that $L(X)$ prorepresents the
functor \eqref{eq:115.7prime} -- it may be viewed as the universal map
of the type
\[X\to W(N),\]
where now $N$ is (not just a $k$-module, but) a variable object in
$\Pro \AbOf_k$. Whereas the pro-object $L(X)$ is of a ``$k$-linear''
nature and may be viewed as the (\emph{quasi-coherent})
\emph{$k$-linearization of the unipotent bundle $X$}, the pro-object
$WL(X)$ of $U(k)$ has lost its $k$-linear nature, we would rather view
it as the canonical ``abelianization'' of $X$, retaining mainly its
additive structure (plus maybe operation of $k$ on it, which is a lot
weaker, though, than structure of an $\scrO_k$-module\dots).

I would like now to examine if the quasi-coherent pro-object $L(X)$,
which has been obtained as the suitable quasi-coherent substitute for
$\scrO_k^{(X)}$ in order to get the basic isomorphism \eqref{eq:115.4}
for quasicoherent $F$, may serve the same purpose for
$\RHom_{\scrO_k}$, in analogy to \eqref{eq:115.5},
\eqref{eq:115.5prime}. Quite generally, if
\[\Gamma=(\Gamma_\alpha)\]
is any pro-$\scrO_k$-module, let's define for any module $F$
\begin{equation}
  \label{eq:115.19}
  \RHom_{\scrO_k}(\Gamma,F) \eqdef \Hom_{\scrO_k}(\Gamma,C^\bullet(F)),\tag{19}
\end{equation}
where $C^\bullet(F)$ is an injective resolution of $F$ -- thus, the
definition extends the usual one when $\Gamma$ is just an
$\scrO_k$-module. Our expectation now would be
\begin{equation}
  \label{eq:115.20}
  \RHom_\scrOk(L(X),F) \tosim \RGamma_X(F_X),\tag{20}
\end{equation}
giving rise to
\begin{equation}
  \label{eq:115.20prime}
  \Ext_\scrOk^n(L(X),F) \tosim \mathrm H^n(X,F_X),\tag{20'}
\end{equation}
for\pspage{487} any unipotent bundle $X$ over $k$, and any
quasicoherent sheaf
\[F=W(N),\]
where $N$ is any $k$-module. Of course, \eqref{eq:115.19} yields (for
general $\Gamma$)
\begin{equation}
  \label{eq:115.19prime}
  \Ext_\scrOk^i(\Gamma,F)\simeq \varinjlim_\alpha\Ext_\scrOk^i(\Gamma_\alpha,F),\tag{19'}
\end{equation}
so that \eqref{eq:115.20prime} may be rewritten more explicitly, if
$X\simeq W(M)$, as
\begin{equation}
  \label{eq:115.21}
  \mathrm H^n(X,W(N)_X) \simeq \bigoplus_{i\ge0}\Ext_k^n(\Gamma_k^i(M),N).\tag{21}
\end{equation}
At any rate, we have a canonical map \eqref{eq:115.20} in
$\D^+(\AbOf_k)$, hence maps \eqref{eq:115.20prime}, in view of the
isomorphism \eqref{eq:115.5} and the canonical map
\begin{equation}
  \label{eq:115.star}
  \scrO_k^{(X)} \to L(X)\tag{*}
\end{equation}
deduced from \eqref{eq:115.18}, and the question now is whether these
are isomorphisms. We may of course assume in \eqref{eq:115.21}
\[X=W(M),\quad n\ge1.\]
If $M$ is projective, so are the modules $\Gamma_k^i(M)$, and hence
the second member in \eqref{eq:115.21} is zero, so we should check the
first member is too. This is clear when $M$ is of finite type, hence
$X$ representable by an affine scheme, whose quasicoherent cohomology
is well-known therefore to vanish in dim.~$n>0$. The general case
should be a consequence of this, representing $X$ as the filtering
direct limit of its submodules which are projective of finite type --
this should work at any rate when $M$ is free with a basis which is at
most countable, using the standard so-called ``Mittag-Leffler''
argument for passage to limit. Thus, in case $M$ projective,
\eqref{eq:115.21} and hence \eqref{eq:115.20} seems OK indeed. When
$M$ is \emph{not} projective, however, there must be some $k$-module
$N$ such that $\Ext_k^1(M,N)\ne0$, and hence the second member of
\eqref{eq:115.19} is non-zero for $n=1$, which should imply rather
unexpectedly
\[\mathrm H^1(X,F)=\mathrm H^1(W(M),W(N)_{W(M)})\ne0,\]
whereas till this very moment I had been under the impression that
quasi-coherent cohomology of unipotent bundles should be zero, just as
for affine schemes! Maybe it has been familiar to Larry Breen for a
long time that this is \emph{not} so? Maybe also for what we want to
do it isn't really basic to find out whether \eqref{eq:115.21} is true
in full generality, as for the purpose of studying Postnikov type
dévissage, the\pspage{488} unipotent bundles $X$ we are going to work
with will be smooth, i.e., $M$ projective (and we may even get away
with free $M$'s, if we need so). The natural idea here for getting
\eqref{eq:115.20} via \eqref{eq:115.21} in full generality, is to use
a projective resolution of $M$ (even a free one), but I'll not try to
work this out now. The main impression which remains is that for the
more relevant cases (involving cohomology groups of a \emph{smooth}
unipotent bundle $X$ at any rate, with quasicoherent coefficients),
the quasicoherent ``pro''-linearization $L(X)$ is just as good for
computing cohomology invariants, as the forbidding $\scrO_k^{(X)}$
modules we were shrinking from.

If now we take an $X_*$ instead of just $X$, namely a ss~complex in
$U(k)$, and assuming the components $X_n$ to be smooth (to be safe),
the isomorphisms \eqref{eq:115.20prime} should give rise to
isomorphisms (for $F=W(N)$)
\begin{equation}
  \label{eq:115.22}
  \mathrm H^n(X_*,F) \simeq \Ext_\scrOk^n(L(X_*),F) \simeq
  \Ext_k^n(L(X_*),N),\tag{22} 
\end{equation}
where the $\Ext^n$ should be viewed as hyperext functors (not
term-by-term), and where in the last member $L(X_*)$ may be
interpreted as a chain complex in $\Pro(\AbOf_k)$. The first member of
\eqref{eq:115.22} is the kind of group occurring as obstruction group
in the Postnikov-type dévissage of $X_*$ into linear structures
$W(M(i)_*)$. The chain-pro-complex $L(X_*)$ may still look a little
forbidding, our hope, though, now is that in the ``pointed'' case we
are really interested in, we may replace $L(X_*)$ by $L\subpt(X_*)$,
which is just a true honest chain complex in $\AbOf_k$. Now, from
\eqref{eq:115.16} we get indeed a canonical map
\begin{equation}
  \label{eq:115.23}
  \Ext_k^n(L(X_*),N) \to \Ext_k^n(L\subpt(X_*),N)\tag{23}
\end{equation}
and we hope that this is an isomorphisms, under suitable assumptions
on $X_*$, the most basic one I can think of now being
\[X_0=e.\]
The map \eqref{eq:115.23} was defined as the transposed of a map of
chain complexes in $\Pro \AbOf_k$
\begin{equation}
  \label{eq:115.24}
  L\subpt(X_*)\to L(X_*),\tag{24}
\end{equation}
deduced from \eqref{eq:115.16} by applying it componentwise. We
recognize here, but with a different interpretation (which seems to me
``the correct'' one), the second map in the often referred-to diagram
\eqref{eq:111.w} of page \ref{p:454}, or \eqref{eq:114.V} in
yesterday's reflections (p.~\ref{p:478}). To say that it gives rise to
isomorphisms \eqref{eq:115.23}, for any $n$ and any module $N$, should
be equivalent to saying that \eqref{eq:115.24} is a quasi-isomorphism
-- but to make sure I should demand a little work on foundations
matters on pro-complexes I guess; also, to see\pspage{489} if the
assumption that \eqref{eq:115.24} is a quasi-isomorphism should imply
the same statement with $L(X_*)$ replaced by its componentwise
projective limit -- namely that \eqref{eq:114.V} on p.~\ref{p:478} is,
which we'll need of course in case $k=\bZ$ for the so-called
``linearization theorem''. Thus, we get three isomorphism or
quasi-isomorphism statements, concerning \eqref{eq:115.23},
\eqref{eq:115.24} and \eqref{eq:114.V} in yesterday's notes, which are
at any rate closely related, and which one hopes to be true, because
this seems needed for a schematization theory of homotopy types to
work. But I should confess I have not tried even to get any clue as to
why this should be true, under the only assumptions, say, that the
components of $X_*$ should be smooth (and possibly the Kan
assumption?), plus $X_0=e$ say.

Now to the second ingredient of the looked-for ``linearization
theorem'', which previously was the first map in \eqref{eq:111.w}
p.~\ref{p:454}, or \eqref{eq:114.Vprime} on page \ref{p:478},
involving maps of
\[\text{either $\bZ^{(X_*(k))}$ or $k^{(X_*(k))}$}\]
into what was previously called $L(X_*)$, and which we would now
rather denote by
\[\varprojlim L(X_*)\quad (\simeq\Gamma\uphat_k(M_*)\text{ if
  $X_*=W(M_*)$}).\]
We made sure that, unless $k=\bZ$, none of the two had any chance to
be a quasi-isomorphism. The only positive thing that came out in this
direction was that the first one of these maps would induce an
isomorphism on $\mathrm H_i$ in the critical dimension (namely $i=n+1$
in the $n$-connected case). We now understand why, for $k\ne\bZ$, we
would not get any actual quasi-isomorphism -- namely, the ``correct''
naive linearization which compares reasonably with $L(X_*)$ should not
be relative to a \emph{constant} ring such as \bZ{} or $k$, but
relative to $\scrO_k$, via the map \eqref{eq:115.star}
(p.~\ref{p:487}) giving rise now to
\begin{equation}
  \label{eq:115.25}
  \scrO_k^{(X_*)}\to L(X_*).\tag{25}
\end{equation}
The more reasonable question now, making good sense really for any
ground ring $k$, is whether this map (under the usual assumptions say
on $X_*$) is a quasi-isomorphism. I wouldn't really but it is, as I
have some doubts as to whether the homology sheaves of the first
member (both members of course being viewed as chain complexes of
\scrOk-modules or ``pro'' such) are quasi-coherent -- but for the time
being I am not sure either if those of the second member are
essentially constant pro-objects! But even if \eqref{eq:115.25} isn't
a quasi-isomorphism, it does behave like one for all practical
purposes of computing quasi-coherent cohomology it would seem, as this
boils down indeed to the isomorphisms \eqref{eq:115.20} or
\eqref{eq:115.20prime}.

\bigbreak

\presectionfill\ondate{7.9.}\pspage{490}

\hangsection{The need and the drive.}\label{sec:116}%
For ten days I haven't written any notes, and the time when I stopped
looks a lot more remote still. For two days still after I last wrote
on the notes, I kept pondering about schematization of homotopy types
-- it were rather lively days, first day I found the amazingly simple
description of the homotopy groups in the schematic set-up, which got
me quite excited; next day, from a phone call to Illusie, it turned
out that the key assumption in all my ponderings on schematization,
namely the ``abelianization theorem'' (sic) asserting isomorphism
between discrete and schematic (namely quasi-coherent) homology (or
equivalently, cohomology) invariants, was definitely false:
consequently, the canonical functor from schematic to discrete
homotopy types turns out definitely not to be an equivalence. This
completely overthrows the idyllic picture in my mind about the
relationship between schematic and discrete homotopy types -- but the
reflection on schematic homotopy types ``in their own right'' had by
then proceeded far enough, the very day before, so that my faith in
the relevance of schematic homotopy types wasn't seriously shaken --
rather, I got excited at drawing a systematic
``bilan''\scrcomment{``bilan'' translates as ``assessment'',
  ``results'', ``balance sheet'', or even ``death toll''\dots} from
the evidence now at hand, about the prospects of developing a theory
of ``schematic homotopy types'' satisfying some basic formal
properties, whether or not such theory be based on semisimplicial
unipotent bundles as models, or on any other kind of models making
sense over arbitrary ground rings. Before reverting to a review of the
more formal properties of abelianization in the context of the basic
modelizer \Cat, I would like still to write down with some care what
had thus come up Monday and Tuesday last week\dots

The next days I felt a great fatigue in all my body, and I then
stopped (till yesterday) any involvement in mathematical reflection. I
am glad I followed this time the hint that had come to me through my
body, rather than brush it aside and go on rushing ahead with the work
I was so intensely involved in, as had been a rule in my life for many
years. This time, I understood that the reluctance of the body to
follow that forward rush, even though I was taking good care of myself
with sleep and food, had strong reasons, which had nothing to do with
neither sleep or food nor with my general way of life. Rather, during
the weeks before and also during those very days, a number of things
had occurred in my life, not all visibly related and of differing
weight and magnitude, to none of which I had really devoted serious
reflection, nor even a minimum of time and attention needed for giving
me a\pspage{491} chance to let these things and their meaning
``enter''. In lack of this, there was little chance my response to
current events would be any better than purely mechanical, and my
interaction with some of the people I love would be in any way
creative. There was this need for being attentive, an urgent need
springing from life itself and which I was about to ignore -- and
there was this drive, this impatience driving me recklessly ahead,
with no look left nor right. Of course, I did know about the need,
``somewhere'' -- and in my head too I kind of knew, but the head was
prejudiced as usual and would take no notice, not of the need and not
of the conflict between an urgent need and a powerful, ego-invested
passion. The head was prejudiced and foolish -- so it was my body
finally which told me: now you stop this nonsense and you take care of
what you well know you better take care \emph{first}, and \emph{now}!
And its language was strong and simple enough and cause me to listen.

Thus, the main work I was involved in for these last seven days was to
let a number of things ``enter'' -- mainly things that were being
revealed through the death of my granddaughter Ella. It surely
\emph{was} ``work'', taking up the largest part of my nights and my
days, -- so much so that I can't really say there is any less fatigue
now than seven days ago. It seems to me, though, it isn't quite the
same fatigue -- this time it is the fatigue coming from work done, not
from work shunned. The ``work done'' wasn't really done by me, I feel,
rather work taking place within me, and ``my'' main contribution has
been to allow it to take place, by providing the necessary time and
quietness; and of course, also, to allow the outcome of this work to
become conscious knowledge, rather than burying it away in some dark
corner of the mind. The rough material, as well as the outcome of this
work, have not been this time new facts or new insights; rather,
things which I had come to perceive and notice, for some time already,
over the last two years, without granting them the proper weight and
perspective -- somehow as if I didn't quite believe what I was
unmistakingly perceiving, or didn't take it quite seriously. Such a
thing, I noticed, happens quite often, not only with me, and takes
care of making even the most lively perceptions innocuous, by
disconnecting them at all price from the image of reality and of
ourselves we are carrying around with us, that the image remains
static, unaffected by any kind of ``information'' flowing around or
through us.

\bigbreak

\presectionfill\ondate{10.9.}\pspage{492}

\hangsection[``Schematic'' versus ``formal'' homology and cohomology
\dots]{``Schematic'' versus ``formal'' homology and cohomology
  invariants.}\label{sec:117}%
Maybe the best will be to write up (and possibly develop some) my
reflections (of the two days after I stopped with the notes) roughly
in the order as they occurred.

There were some somewhat technical afterthoughts. One was about the
logical difficulty coming from the site $\Aff_{/k}$ on page
\ref{p:446} not being a \scrU-site (where \scrU{} is the universe we
are working in), hence strictly speaking, the category of all sheaves
on this site is possibly not even a \scrU-category (i.e., the $\Hom$
objects need not be small, i.e., with cardinal in \scrU), still less a
topos, and hence the standard panoply of notions and constructions in
a topos does not apply. This doesn't look really serious, though, one
way out is to limit beforehand the ``size'' of the unipotent bundles
we are allowing, i.e., of the $k$-modules describing them, in terms
say of cardinality of a family of generators of the latter -- and then
restrict accordingly the size of the $k$-algebras $k'$ taken as
``arguments'' for our sheaves, i.e., as objects of the basic site we
are working on. For instance, when working with unipotent bundles of
finite type only (i.e., corresponding to $k$-modules of finite type --
a rather interesting and natural finiteness condition anyhow on the
components of a schematic model $X_*$), it is appropriate to work on
the ``fppf site'',\scrcomment{``fidèlement plat de présentation
  finie''} where the arguments $k'$ are $k$-algebras of \emph{finite
  presentation}. If we should be unwilling to be limited by a fixed
size restriction on the unipotent bundles we are working with (and
hence also on the corresponding homotopy types), we may have to work
with a hierarchy of size restrictions and passage from one to any
other less stringent one -- a technical nuisance to be sure, if we
don't find a more elegant way out, but surely not a substantial
difficulty. At the present heuristic stage of reflections, it doesn't
seem worth while really to dwell on such questions any longer.

Another afterthought is about the functor
\[M\mapsto W(M): \AbOf_k = \kMod\to(\scrOk\textup{-Mod})\]
from $k$-modules to \scrOk-modules -- a fully faithful functor we
know, whose essential image by definition consists of the so-called
\emph{quasi-coherent} sheaves. A little caution is needed, as this
functor is right exact, but \emph{not exact}, i.e., it does not
commute with formation of kernels, because for a $k$-algebra $k'$
which isn't flat, the functor
\[M\mapsto M\otimes_k k'\]
doesn't.\pspage{493} When we wrote down a formula such as
\eqref{eq:115.21} on page \ref{p:487}, we were implicitly making use
of the assumption that for $k$-modules $M,N$ we have a canonical
isomorphism
\begin{equation}
  \label{eq:117.1}
  \RHom_k(M,N) \fromsim \RHom_\scrOk(W(M),W(N)),\tag{1}
\end{equation}
at any rate such a formula is needed if we want to view
\eqref{eq:115.21} as a more explicit way of writing \eqref{eq:115.20}
or \eqref{eq:115.20prime}. Using a free resolution $L_\bullet$ of $M$,
and the fact that
\[\Ext_\scrOk^i(\scrOk, W(N)) \simeq \mathrm H^i(\Spec(k),W(N)) = 0
\quad\text{if $i>0$,}\]
stemming from the cohomology properties of flat descent, we easily get
a map \eqref{eq:117.1}, but I did not check that this map is an
isomorphism, the difficulty coming from the fact that $W(L_\bullet)$
need not be a resolution of $W(M)$, unless $M$ is flat. This
perplexity already arises in the fppf context -- surely Larry Breen
should know the answer. For what we are after here it doesn't seem to
matter too much, as the computations of loc.\ cit.\ were of interest
mainly (maybe exclusively) in the case when $M$ is projective or at
any rate flat, i.e., when working with flat unipotent bundles -- in
which case \eqref{eq:117.1} is indeed an isomorphism.

The interpretation of polynomial maps between $k$-modules $M$, $N$ in
terms of the topological augmented coalgebras $\Gamma\uphat_k$
associated to these may seem a little forbidding to some readers. In
the all-important and typical case when $M$ and $N$ are projective and
of finite type, things become quite evident, though, by just dualizing
the more familiar concepts around polynomial functions and
homomorphisms between rings of such. Thus, $W(M)$ is just a usual
vector-bundle, hence also a true honest affine scheme over $k$, whose
affine ring is the ring of polynomial functions on $M$, which can be
identified with $\Sym_k^*(M\op)$, the symmetric algebra on the dual
module:\scrcomment{To make sense of these formulæ I have added some
  $\op$s to denote dual modules. I hope this matches AGs
  intentions\dots}
\begin{equation}
  \label{eq:117.2}
  W(M) \simeq \Spec(\Sym_k^*(M\op)),\tag{2}
\end{equation}
and similarly for $W(N)$. Polynomial maps from $M$ to $N$, i.e., maps
from the $k$-scheme $W(M)$ to the $k$-scheme $W(N)$, just correspond
to $k$-algebra homomorphisms
\begin{equation}
  \label{eq:117.3}
  \Sym_k^*(N\op) \to \Sym_k^*(M\op)\tag{3}
\end{equation}
(irrespective of the graded structures). As each of these algebras, as
a $k$-module, is a filtering direct limit of its projective submodules
of finite type such as
\begin{equation}
  \label{eq:117.4}
  \Sym_k(M\op)(i) = \bigoplus_{j\le i} \Sym_k^j(M),\tag{4}
\end{equation}
the\pspage{494} dual module
\[\Gamma_k\uphat(M)\simeq (\Sym_k^*(M))\upvee \simeq
\prod_j(\Sym_k^j(M) \simeq \Gamma_k^j(M))\]
may be viewed as the inverse limit of the duals of those submodules
(compare p.~\ref{p:484} \eqref{eq:115.11}), and topologized
accordingly; linear maps \eqref{eq:117.3} may be interpreted in terms
of continuous maps between the dual structures (or, equivalently,
between the corresponding pro-objects)
\begin{equation}
  \label{eq:117.3prime}
  \Gamma\uphat_k(M) \to \Gamma\uphat_k(N),\tag{3'}
\end{equation}
and compatibility of \eqref{eq:117.3} with multiplication and units is
expressed by compatibility of \eqref{eq:117.3prime} with
``comultiplication'', i.e., diagonal maps, and with augmentations. On
the other hand, maps $W(M)\to W(N)$ respecting the ``pointed
structures'' coming from zero sections, correspond to map
\eqref{eq:117.3prime} transforming $1$ into $1$, besides the other
requirements. Such maps, it turns out, automatically induce a map
between the submodules
\begin{equation}
  \label{eq:117.4bis}
  \Gamma_k(M)\to\Gamma_k(N)\tag{4}
\end{equation}
which is rather evident indeed, if we remind ourselves of the fact
that the submodule $\Gamma_k(M)$ may be viewed as the (topological)
dual of $\Sym_k^*(M)$, topologized by the powers of the augmentation
ideal
\[\Sym_k^+(M)=\bigoplus_{i>0} \Sym_k^i(M),\]
or equivalently of the corresponding adic completion
\begin{equation}
  \label{eq:117.5}
  {\Sym_k^*}\uphat(M) = \varprojlim_i \Sym_k^*(M)(i) \simeq
  \prod_{j\ge0}\Sym_k^j(M),\tag{5} 
\end{equation}
and correspondingly for $N$. The ``pointed'' assumption on a map
$W(M)\to W(N)$ in terms of the corresponding homomorphism of
$k$-algebras \eqref{eq:117.3}, just translates into compatibility with
the augmentations, or equivalently, with the corresponding ideals
$\Sym_k^+$, which implies that it induces a homomorphism of the
corresponding adic rings, and hence by duality a homomorphism
\eqref{eq:117.4} on their duals.

These reminders bring near that working with the (discrete)
$k$-algebras $\Sym_k^*(M)$, or equivalently, with the
\emph{topological} coalgebras $\Gamma_k\uphat(M)$, amounts to working
with ``unipotent bundles'' (projective and of finite type), which are
just usual schemes (of a rather particular structure of course),
whereas working with the \emph{topological} $k$-algebras
${\Sym_k^*}\uphat(M)$, or equivalently, with the discrete coalgebras
$\Gamma_k(M)$, amounts to working with \emph{formal schemes}, namely
essentially, with the formal completions of the former along the zero
sections. (Of course, the topological objects just
considered\pspage{495} may be equally viewed as being
\emph{pro-modules} endowed with suitable extra structure.)
Correspondingly, we will expect the cohomology invariants constructed
in terms of (the apparently more forbidding) $\Gamma_k\uphat$ to
express quasi-coherent cohomology of the corresponding \emph{schemes},
or semisimplicial systems of such; whereas, working with the
(apparently more anodyne!) $\Gamma_k$ will lead to cohomology
invariants of \emph{formal schemes} and semisimplicial systems of
such. Both types of invariants are of interest it would seem, the one
however which looks the more relevant in connection with studying
ordinary homotopy types in terms of schematic ones, is surely the
first. On the other hand, there isn't any reason whatever to believe
that under fairly general conditions, these two types of invariants
are going to be isomorphic, by the evident map from ``schematic'' to
``formal'' quasicoherent cohomology. To say it differently, I do no
longer expect that under reasonably wide assumptions, the map
\eqref{eq:115.24} of p.~\ref{p:488}
\begin{equation}
  \label{eq:117.6}
  L\subpt(X_*)\to L(X_*),\tag{6}  
\end{equation}
is a quasi-isomorphism, nor behaves like one with respect to taking
$\Ext^i$'s with values in a quasicoherent module, as I was hastily
surmising for about one week, while loosing track of the geometric
meaning of the algebraic objects I was playing around with. To give
just one example, take $X_*$ to be the standard semisimplicial
unipotent bundle associated to the group-object $\bG\suba$ in
$\scrU(k)$, namely just the usual affine line with addition law. The
$\Ext^i$'s of the two members of \eqref{eq:117.6} with values in the
$k$-module $k$ may be interpreted as either schematic or formal
Hochschild cohomology of the additive group, with coefficients in
$\scrO_k=\bG\suba$. The map of the former into the latter is not
always an isomorphism, already in dimension $1$, where the two groups
to be compared are just the groups of endomorphisms of $\bG\suba$, and
of the corresponding formal group. If $k$ is of char.\ $p>0$, a prime,
then the latter group can be described as the group of all formal
power series of the type
\[F(t)=\sum_{i\ge0} c_i t^{(p^i)},\]
whereas for the first group we must restrict to those $F$ which are
polynomials, i.e., only a finite number of the coefficients $c_i$ are
non-zero.

One may of course object to this example, because the $X_*$ we are
working with is not simply connected, and because the example does not
apply over a ring such as \bZ, which is the one we are interested in
most of all. I am convinced now, however, that even when assuming
$X_1=X_0=e$ and $k=\bZ$, \eqref{eq:117.6} is very far from being a
quasi-isomorphism, even for such\pspage{496} basic structures as
$K(\bZ,n)$ with $n\ge2$. At any rate, the ``way-cut'' argument I have
finally been thinking of, in order to check \eqref{eq:117.6} is a
quasi-isomorphism, rests on vanishing assumptions which (as I was
informed by Illusie the next day) are wholly unrealistic. This finally
clears up, it would seem, a tenacious misconception which has been
sticking to my first heuristic ponderings about the homology and
cohomology formalism for schematic homotopy types: \emph{one should be
  very careful not to substitute the ``pointed'' linearization
  $L\subpt(W(M)) \simeq \Gamma_k(M)$ for the non-pointed one
  $L(W(M))\simeq\Gamma_k\uphat(M)$, in computing homology and
  cohomology invariants of schematic homotopy types}. To say it
differently, in order to be able to compute (or just define) the
``schematic'' homology and cohomology invariants, we do need as a
model a full-fledged semisimplicial unipotent bundle, not just the
corresponding formal one, giving rise to invariants of it's own,
namely ``formal'' homology and cohomology, which are definitely
distinct from the former.

\hangsection[The homotopy groups $\pi_i$ as derived functors of the
``Lie \dots]{The homotopy groups \texorpdfstring{$\pi_i$}{pi-i} as
  derived functors of the ``Lie functors''. Lack of satisfactory
  models for \texorpdfstring{$S^2$ and $S^3$}{S2 and S3}.}\label{sec:118}%
It is all the more remarkable, in view of the preceding findings, that
the homotopy invariants
\[\pi_i(X_*), \quad i\ge0,\]
of a pointed semisimplicial unipotent bundle $X_*$ (still assuming the
components $X_n$ to be smooth, i.e., to correspond to projective
modules) turn out to be invariants of the corresponding ``formal''
object, and, more startling still, of the corresponding
``infinitesimal'' object of order $1$. More specifically, consider the
``Lie functor'' or ``tangent space at the origin''
\begin{equation}
  \label{eq:118.7}
  \Lie: U(k)^\bullet \to \AbOf_k, \quad X=W(M)\mapsto M,\tag{7}
\end{equation}
which we'll need only for the time being for smooth $X$, when the
geometric meaning of it is clear. This functor transforms
semisimplicial pointed bundles into semisimplicial $k$-modules
$\Lie(X_*)$, thus we should get, besides abelianization, another
remarkable functor, from $\scrM_1(k)$ say to $\D_\bullet(\AbOf_k)$:
\begin{equation}
  \label{eq:118.8}
  \Lie: \scrM_1(k)\to\D_\bullet(\AbOf_k), \quad \text{via $X_*\mapsto
    \Lie(X_*)$,}\tag{8} 
\end{equation}
granting that $X_*\mapsto\Lie(X_*)$ transforms quasi-isomorphisms into
quasi-isomorphisms. Now, that this must be so follows from the really
startling formula
\begin{equation}
  \label{eq:118.9}
  \pi_i(X_*)\simeq\pi_i(\Lie(X_*)) \quad
  \begin{tabular}[t]{@{}c@{}}
    ($\simeq\mathrm{H}_i$ of
    the associated chain\\
    complex in $\AbOf_k$),
  \end{tabular}\tag{9}
\end{equation}
where the left-hand side, I recall, is defined as
\begin{equation}
  \label{eq:118.10}
  \pi_i(X_*)=\pi_i(X_*(k)).\tag{10}
\end{equation}
I\pspage{497} don't have, I must confess, any direct description of
such an isomorphism \eqref{eq:118.9}, valid for any semisimplicial
bundle $X_*$, say, satisfying the assumptions
\[X_0=X_1=e, \quad \text{$X_n$ smooth for any $n$ (maybe flat is
    enough),}\]
plus possibly (if needed) a Kan type condition. However, we have such
isomorphisms \eqref{eq:118.9} in a tautological way, when $X_*$ comes
in the usual way from a chain complex in $\AbOf_k$ with projective
components, hence also when $X_*$ admits a Postnikov-type dévissage
into ``abelian'' pieces as above. If we admit that any $X_*$
satisfying the assumptions is homotopic to one admitting such a
dévissage, the isomorphisms \eqref{eq:118.9} should follow, except of
course that extra work would be still needed to get naturality of
\eqref{eq:118.9}.

I have the feeling however that, besides the specific abelianization
functor in the schematic context, \emph{formula \eqref{eq:118.9}
  should be made a cornerstone of a theory of schematic homotopy
  types, and serve as ``the'' natural definition of the homotopy
  invariants} of a model $X_*$, within the context of schematic models
and without any need a priori to tie them up with, let one subordinate
their study to, invariants of the corresponding discrete homotopy type
$X_*(k)$. Accordingly, weak equivalences should be defined (for
semisimplicial bundles satisfying the suitable assumptions at any
rate) as maps inducing isomorphisms for the $\pi_i$ invariants, namely
inducing quasi-isomorphisms for the corresponding Lie chain
complexes. It is immediately checked that this implies that the
corresponding map for ``formal homology'', namely
\[L\subpt(X_*)\to L\subpt(X_*)\]
is then a quasi-isomorphism too, when viewed as a map of chain
complexes and hopefully the same should hold for ``schematic
homology''
\[L(X_*)\to L(X_*),\]
and of course one would expect converse statements to hold too.

I would like to comment a little on the significance of formula
\eqref{eq:118.9}. As far as I know, this is the only fairly general
formula, not reducing to an ``abelian'' case, where the homotopy
groups $\pi_i$ appear as just the $\mathrm H_i$ invariants of a
suitable chain complex, defined up to unique isomorphism in the
relevant derived category. This chain complex comes here, moreover,
with an amazingly simple description, of immediate geometrical
significance, and suggestive of relationships with\pspage{498} the
homology invariants a lot more precise, presumably, than those
currently used so far. Of course, the significance of \eqref{eq:118.9}
for the study of the usual, ``discrete'' homotopy types, will be
subordinated to how difficult it will turn out for such a homotopy
type to be a)\enspace realizable by a schematic over \bZ, and
b)\enspace to get hold of a more or less explicit description of such
a schematic homotopy type, via say a semisimplicial bundle as a
model. One is of course thinking more specifically of the case of the
spheres $S^n$, the first case (besides the trivial $S^1$ case) being
$S^2$, the sequence of homotopy groups of which (as far as I know) it
not understood yet. Viewing the spheres $S^n$ as successive
suspensions of $S^1$, where $S^1$ is fitting nicely into the formalism
of schematic homotopy types as a $K(\bZ,1)$ (except that the
$1$-connectedness condition $X_1=e$ is not satisfied), this brings
near the question of defining the suspension operation in a relevant
derived category $\scrM_0(k)$ or $\scrM_1(k)$ (whereas before we had
met with the ``dual'' question of constructing homotopy fibers of
maps, such as loop spaces). Thus it would seem that \emph{a
  breakthrough in getting hold of the standard homotopy constructions
  within the schematic context}, assuming that these constructions do
still make a sense, \emph{may well mean a significant advance for the
  understanding of the homotopy groups of spheres}. This looks like a
very strong motivation for trying to carry through those constructions
(possibly even construct a corresponding ``derivator'' embodying any
kinds of finite ``integration'' and ``cointegration'' operations on
schematic homotopy types) and at the same time warns us that such
work, if at all feasible, will most probably be a highly non-trivial
one.

The question of ``schematizing'' the homotopy types of $S^2$ and $S^3$
reminds me of a fact which struck me a long time ago (maybe
J.~P.~Serre or someone else pointed it out to me first). Namely, in
some respects there doesn't seem to be really satisfactory algebraic
models for these homotopy types, taking into account the basic
relationship between the two, namely: the $2$-sphere (or,
equivalently, the projective complex line) is a \emph{homogeneous
  space} under the quaternionic \emph{group} $S^3$ (or, equivalently,
under the complex linear group $\mathrm{SL}(2,\bJ)$).\scrcomment{I'm
  guessing that $\mathrm{SL}(2,\bJ)$ is another notation for
  $\mathrm{SU}(2)$\dots} This relationship, and its manifold
``avatars'' in the realms of discrete groups, Lie groups, algebraic
groups or group schemes etc., is one of the few key situations met
with, and of equal basic significance, in the most diverse quarters in
mathematics, from topology to arithmetic. Thus, $\mathrm{SL}(2)$ as a
simple Lie (or algebraic) group of minimum rank $1$, plays the role
of\pspage{499} \emph{the} basic building block for building up the
most general semisimple groups, whereas $P^1$ may be viewed as being
the most significant homogeneous space under this group, namely the
first and most elementary case of flag manifolds. In view of this
significance of $S^2$ and $S^3$, it is all the more reasonable that no
simple, non-plethoric semisimplicial model say in
\scrcommentinline{unreadable}, in terms say of a semisimplicial group
having the homotopy type of
\[S^3 \sim \mathrm{SL}(2,\bC)=G,\]
and a subgroup having the homotopy type of
\[S^1 \sim \bC^\times = \mathrm{GL}(1,\bC) \sim K(\bZ,1)\]
and playing the part of a Borel subgroup or a maximal torus, in such a
way that the quotient will have the homotopy type of
\[S^2= P^1_\bC \simeq S^3/S^1 \simeq G/B.\]
It would be tempting now to try and construct such a model of the
situation in terms of semisimplicial unipotent bundles over $\bZ$ --
which would at the same time display the homotopy groups of $S^2$ and
$S^3$ (not much of a difference!) via formula \eqref{eq:118.9}. All
the more tempting of course, as it is felt that the geometric objects
and their relationship, the homotopy shadow of which we want to
modelize schematically, are themselves already, basically, most
beautiful schemes over $\Spec(\bZ)$!

Another way of getting a display of the homotopy groups of $S^2$ and
$S^3$ would be in terms of a (discrete) model of the situation above,
in terms of ``hemispherical complexes'' rather than semisimplicial
ones. On the other hand, there is no reason why a theory of schematic
homotopy types could not be carried through as well, using
hemispherical complexes rather than semisimplicial ones. The latter
kind of complexes have the advantage that they have
become thoroughly familiar through constant use
by topologists and homotopy people for thirty years or so -- the
former however are newcomers, have the advantage of still greater
formal simplicity (just two boundary operations, and just one
degeneracy), and more importantly still, of allowing for a direct
computational description of the homotopy invariants, in the discrete
set-up. When working with hemispherical complexes of unipotent bundles
as models for schematic homotopy types, we'll get then \emph{two}
highly different descriptions of the homotopy invariants $\pi_i$, one
by the ``infinitesimal'' formula \eqref{eq:118.9} interpreting them as
derived functors of the Lie functor, the other one via the\pspage{500}
hemispherical set $X_*(k)$, handled as if it came from an actual
\oo-groupoid (by taking its source and target operations etc.) even if
it does not. (I confess I did not check that this process does
correctly describe the homotopy groups of a hemispherical set, even
without assuming it comes from an \oo-groupoid or an \oo-stack, but
never mind for the time being\dots)\enspace Among the interesting
things still ahead (once we get a little accustomed to working with
hemispherical complexes) is to try and understand how these two
descriptions relate to each other, which may be one means for a better
understanding of the basic formula \eqref{eq:118.9}, in the context of
hemispheric schematic models.

\bigbreak

\presectionfill\ondate{26.9.}\par

\hangsection[Breakdown of an idyllic picture -- and a tentative next
\dots]{Breakdown of an idyllic picture -- and a tentative next best
  ``binomial'' version of the ``comparison theorem'' for schematic
  versus discrete linearization.}\label{sec:119}%
After the last notes (of September 10) I was a little sick for a few
days, then I was taken by current tasks from professional and family
life, which left little leisure for mathematical reflection, except
once or twice for a couple of hours, by way of recreation. It would
seem now that in the days and weeks ahead, there will be more time to
go on with the notes, and I feel eager indeed to push ahead. Also, I
more or less promised the publisher, Pierre Bérès, that a first volume
would be ready for the printer by the end of this calendar year, and I
would like to keep this promise.

I still have to tie in with the reflections and happenings of the end
of last month, as I started upon with the last notes (of
Sept.~10). Next thing then to report upon is the ``coup de
théâtre''\scrcomment{a sudden or unexpected event in a play\dots}
occurring through the phone call to Luc Illusie. When I told him about
what by then still looked to me as the key assumption for a theory of
schematization of homotopy types, namely that the homology of
$K(\pi,n)$ should be computable in terms of derived functors of the
``divided power algebra''-functor $\Gamma_\bZ$, he at once felt rather
skeptical, and later he called me back to tell me it was definitely
false. He could not give me an explicit counterexample for $\mathrm
H_i(n,\bZ;\bZ)$, say, with given $n$ and $i$, rather he said that when
suitably ``stabilizing'' the assumption I had in mind, it went against
results of Larry Breen on $\Ext^i$ functors of $\bG\suba$ with itself
over prime fields $\bF_p$. I don't know if I am going some day to give
into Illusie's argument and into Larry Breen's results -- however,
even before I got Illusie's confirmation that definitely my assumption
was wrong, I convinced myself that at any rate it was false for
$n=1$. This is a non-simply connected case and hence not\pspage{501}
entirely conclusive maybe, but still it was enough to shake my
confidence that the assumption was OK. The counterexample is in terms
of cohomology
\[\mathrm H^2(1,\bZ;\bZ)=\mathrm H^2(\bZ,\bZ) = \mathrm
  H^2(B_\bZ,\bZ)\]
rather than homology, as usual. As the classifying space $B_\bZ=S^1$
of \bZ{} is one dimensional, its cohomology is zero in dimensions
$i\ge2$. On the other hand, $\mathrm H^2$ classifies central
extensions of \bZ{} by \bZ, and an immediate direct argument shows
indeed that such extensions (indeed, any extension of \bZ{} by any
group) split. If we take the schematic $\mathrm H^2$, defined by
Hochschild cochains which are polynomial functors, we get the
classification of central extensions of $\bG\suba$ by $\bG\suba$, as
group schemes over the ring of integers. Now, it is easy to find such
an extension which does \emph{not} split, the first one one may think
of being the group scheme representing the functor
\[k\mapsto W_2(k) =
  \begin{tabular}[t]{@{}l@{}}
    group of truncated power series $1+at+bt^2$ \\
    in $k[t]/(t^3)$,
  \end{tabular}\]
where $a,b$ are parameters in the commutative ring $k$, the group
structure being multiplication. These parameters define in an evident
way a structure of an extension of $\bG\suba$ by $\bG\suba$ upon
$W_2$, a splitting of which would correspond to a group homomorphism
\[\bG\suba \to W_2, \quad a\mapsto 1+at+P(a)t^2,\]
with
\[P\in \bZ[t].\]
Expressing compatibility with the group laws gives the condition
\[P(a+a') = P(a)+P(a')+aa',\]
which has, as unique solutions in $\bQ[t]$, expressions
\[P(t)=t(t-1)/2 + ct\]
with $c$ in \bQ, none of which has coefficients in \bZ. This argument
shows in fact that for given ring $k$, $(W_2)_k$ is a split extension
if{f} $2$ is invertible in $k$, in which case a splitting is given by
\[a\mapsto 1+at+(a(a-1)/2)t^2.\]

This example brings near one plausible ``reason'' why the expected
comparison statement about discrete and schematic linearization could
not reasonably hold true, and in particular why we shouldn't expect
discrete and schematic Hochschild cohomology (for group schemes over
\bZ{} such as $\bG\suba$\pspage{502} or successive extensions of such)
to give the same result. Namely, the latter is computed in terms of
cochains which are polynomial functions \emph{with coefficients in
  \bZ}, whereas there exist polynomial functions \emph{with
  coefficients in \bQ} (not in \bZ) which, however, give rise to
integer-valued functions on the group of integer-valued points. Such
are the binomial expressions
\[P_n(t) = t(t-1)\dots(t-n+1)/n!\quad\text{(for $n\in\bN$).}\]
These (in the case of just one variable $t$) are known to form a basis
of the \bZ-module of all integer-valued functions on \bZ, and these is
a corresponding basis for integer-valued functions on $\bZ^r$, for any
natural integer $r$. Thus, the hope still remains that a sweeping
comparison theorem for discrete versus ``schematic'' linearization
might hold true, provided it is expressed in such a way that the
``schematic models'' we are working with should be built up with
``schemes'' (of sorts) described in terms of spectra not of polynomial
algebras $\bZ[t]$ and tensor powers of these, but rather of ``binomial
algebras'' $\bZ{\angled t}$ built up with the binomial expressions
above, and tensor powers of such. If we want to develop a
corresponding notion of homotopy types over a general ground ring $k$,
we should then require upon $k$ an extra structure of a ``binomial
ring'' (as introduced in the Riemann-Roch Seminar SGA~6 in some talks
of Berthelot),\scrcomment{\textcite{SGA6}} namely a ring
endowed with operations
\[x\mapsto \binom xn : k\to k \quad (n\in\bN),\]
satisfying the formal properties of the binomial functions $x\mapsto
P_n(x)$ in the case $k=\bZ$ or $\bQ$. Whereas linearization of
homotopy types via De~Rham complexes with divided powers relies on a
``commutative algebra with divided powers'' (which was developed
extensively by Berthelot and others for the needs of crystalline
cohomology), linearization via unipotent bundles (assuming it can be
done in such a way as to ensure that any discrete homotopy type can be
``schematized'' in an essentially unique way) might well rely on the
development of a ``binomial commutative algebra'' and a corresponding
notion of ``binomial schemes''. There should be a lot of fun ahead
developing the necessary algebraic machinery, which may prove of
interest in its own right.\footnote{See comments next section
  p.~\ref{p:506}--\ref{p:507}.} It should be realized, however, that
for a ring $k$ to admit a binomial structure is a rather strong
restriction -- thus, for a given prime $p$, no field of char.~$p$
(except possibly the prime field?) admits such a structure. This
remark may temper somewhat the enthusiasm for pushing in this
direction, even granting that a ``binomial comparison theorem'' for
discrete versus ``binomial'' linearization holds true.

Maybe\pspage{503} it is worthwhile to give a down-to-earth formulation
of such a comparison statement. For any free \bZ-module $M$ of finite
type, let
\[\Symbin_\bZ(M) \subset \Sym_\bQ^*(M_\bQ)\quad\text{(where
    $M_\bQ=M\otimes_\bZ \bQ$)}\]
be the subalgebra of the algebra of polynomial functions on
$M_\bQ\upvee \simeq (M\upvee)_\bQ$ which are integral-valued on
$M\upvee$ = dual module of $M$. Now let $L_*$ be any semisimplicial
\bZ-module whose components are free of finite type, and consider the
canonical map of cosemisimplicial \bZ-modules
\[\Symbin_\bZ(L_*\upvee) \to \Maps(L_*,\bZ),\]
described componentwise in an obvious way. The question is whether
this is a weak equivalence, i.e., induces a quasi-isomorphism for the
associated cochain complexes, under the extra assumption that $L_*$ is
$0$-connected, i.e., the associated chain complex has zero $\mathrm
H_0$ (and possibly, if needed, assuming even $1$-connectedness, i.e.,
$\mathrm H_0$ \emph{and} $\mathrm H_1$ of the associated chain complex
are both zero). Presumably, by easy dévissage arguments one should be
able to reduce to the case when $L_*$ is a $K(\bZ,n)$ type, and more
specifically still, that it is the semisimplicial abelian complex
associated to the chain complex reduced to \bZ{} placed in degree $n$
(where $n\ge1$). Thus, the question is whether Eilenberg-Mac~Lane
cohomology (with coefficients in \bZ) for $K(\pi,n)$ types (or more
specifically, $K(\bZ,n)$ types) can be expressed in terms of derived
functors of the $\Symbin(M\upvee)$ functor. At any rate, whether
$\Symbin$ is just the right functor to fit in or not, it looks like an
interesting question whether Eilenberg-Mac~Lane cohomology (or, more
relevantly still, homology) can be expressed in terms of the derived
functor of a suitable non-additive contravariant (resp.\ covariant)
functor $B$ from \Ab{} (or from the subcategory of free \bZ-modules)
to itself. If so, there should be a way of defining a (possibly
somewhat sophisticated) notion of ``$B$-schematic homotopy types''
(over a ground ring $k$ endowed with suitable extra structure, such as
a binomial structure), in terms of ``unipotent $B$-bundles'', in such
a way that any ``discrete'', namely usual homotopy type, satisfying a
suitable $1$-connectedness restriction, admits an essentially unique
``$B$-schematization''.

I don't feel like pursuing these questions here, which would take me
too far off the main line of investigation I've been out for. At any
rate, whether or not Eilenberg-Mac~Lane homology may be expressed in
terms of the total left derived functor of a suitable functor from
\Ab{} to itself, it would seem that the somewhat naive approach
towards schematic homotopy types we have been following, valid over an
arbitrary (commutative) ground ring $k$ without any extra structure
needed on it, is worthwhile\pspage{504} pursuing even for the mere
sake of studying ordinary homotopy types. The main reason for feeling
this way is the amazingly simple description of the homotopy modules
$\pi_i$ of a homotopy type defined in terms of a semisimplicial (or
hemispherical) unipotent bundle, as derived functors (so to say) of
the Lie functor (cf.\ previous section \ref{sec:118}). The main test
for deciding whether there is indeed a rewarding new tool to be dug
out, is to see whether or not in the model categories $\scrM_n(k)$ we
have been working with so far, the standard homotopy constructions
(around loop spaces and suspensions) make sense, and in such a way of
course that the canonical functor from schematic to discrete homotopy
types should commute to these operations. It may well turn out that to
get a handy formalism, one will have still to modify more or less the
conceptual set-up of unipotent bundles I've been tentatively working
with so far. I already lately hit upon suggestions of such
modifications, and presumably I'm going to discuss this still, before
leaving the topic of schematizations.

Another reason which makes me feel that there should exist a notion of
homotopy types over more general ground rings $k$ than \bZ, is that
for a number of rings, such a notion has been known for quite a
while. If I got it right, already in the late sixties (even before I
withdrew from the mathematical milieu) I heard about such things as
homotopy types over residue class rings $\bZ/n\bZ$, or over rings
(such as \bQ) which are localizations of \bZ, or over rings such as
$\bZ\uphat$ or $\bZ_p$ ($p$-adic integers) which are completions of
\bZ{} with respect to a suitable linear topology. Last week, which was
the first time I was at a university after the Summer vacations, I
took from the library the
Bousfield-Kan\scrcomment{\textcite{BousfieldKan1972}} Lecture Notes
book on homotopy limits (which had been pointed out to me by Tim
Porter in June, when he had taken the trouble to tell me about ``shape
theory'' and its relations to (filtering) homotopy limits), and while
glancing through it, I noticed there is a systematic treatment of such
homotopy types.

At this very minute I had a closer look upon the introduction of part
I, it turns out that Bousfield and Kan are working with an arbitrary
(commutative) ground ring $k$, and they are defining corresponding
$k$-completion $k_\oo X$ of a homotopy type $X$, rather than a notion
of ``homotopy type over $k$''. But the two kinds of notions are surely
closely related, the $k$-completion of BK presumably should have more
or less the meaning of ground ring extension $\bZ\to
k$.\footnote{definitely not, in general!} At any rate, for
$1$-connected spaces and $k=\bZ$ the completion operation seems to be
no more no less than just the identity, thus it would seem that the
implicit notion of ``homotopy type over \bZ'' should be just the
ordinary ``discrete'' notion of homotopy types -- unlike the notion
of\pspage{505} schematic homotopy types (over $k=\bZ$) defined via
semisimplicial unipotent bundles. Definitely, an understanding of
schematic homotopy types will have to include the (by now classical)
Bousfield-Kan ideas, and these are also relevant for my reflections on
``integration'' and ``cointegration'' operations (in connection with
the notion of a derivator (section \ref{sec:69})), called in their
book ``homotopy direct limits'' and ``homotopy inverse limits'' (in
the special case of the derivator associated to ordinary homotopy
types, if I got it right). It came as a surprise that in their book,
these operations are developed mainly as technical tools for
developing their theory of $k$-completions, whereas in my own
reflection they appeared from the start as ``the'' main operations in
homotopical as well as homological algebra. There had been quite a
similar surprise when Tim Porter had sent me a reprint (in July, just
before I stopped with my notes for a month or so) of Don Anderson's
beautiful paper ``Fibrations and Geometric
Realizations''\scrcomment{\textcite{Anderson1978}} (Bulletin of
the Amer.\ Math.\ Soc.\ September 1978), where a very general and (as
I feel) quite basic existence theorem for integration and
cointegration operations in the set-up of closed model categories of
Quillen of the type precisely I was after, is barely alluded to at the
end of the introduction, and comes more or less as just a by-product
of work done in view of a result on geometric realizations which (to
an outsider like me at any rate) looks highly technical and not
inspiring in the least!

It is becoming clear that I cannot put off much longer getting
acquainted with the main ideas and results of Bousfield-Kan's book,
which definitely looks like one of the few basic texts on foundational
matters in homotopy theory. Still, before doing some basic reading, I
would like to write down the sporadic reflections on schematic
homotopy types I went into during the last weeks, while they are still
fresh in my mind!

\bigbreak

\presectionfill\ondate{28.9.}\pspage{506}\par

\hangsection[Digression on the Lazare ``analyzers'' for ``binomial''
\dots]{Digression on the Lazare ``analyzers'' for ``binomial''
  commutative algebra, \texorpdfstring{$\lambda$}{lambda}-commutative
  algebra, etc.}\label{sec:120}%
Yesterday (prompted by the reflections from the day before, cf.\
section \ref{sec:119}), I pondered a little on the common features of
the various set-ups for ``commutative algebra'' (possibly, too, for
corresponding notions of ``schemes'') one gets when introducing extra
operations on commutative rings or algebras, such as divided power
structure (on a suitable ideal) or a $\lambda$-structure with
operations $\lambda^i$ paraphrasing exterior powers, or a binomial
structure with operations $x\mapsto\binom xn$ paraphrasing binomial
coefficients, or an $S$-structure with operations $S^i$ of Adams' type
paraphrasing sums of $i$\textsuperscript{th} powers of roots of a
polynomial (a weakened version of a $\lambda$-structure). It seems
that the unifying notion here is the notion of an ``analyzer''
(analyseur) of Lazare, ``containing'' the Lazare analyzer for
commutative rings (not necessarily with unit), so that the components
$\Omega_n$ ($n\ge-1$) are commutative rings (not necessarily with
unit). In case the extra operations we want to introduce on
commutative rings are to be defined on rings \emph{with units} (not
just on a suitable ideal of such a ring, as is the case for the
divided power structure), and they all can be defined in terms of
operations involving just one argument, the reasonable extra axiom on
the corresponding analyzer (as suggested by the examples at hand) is
that for $n\ge1$, $\Omega_n$ can be recovered in terms of $\Omega_0$
and $\Omega_{-1}=k_0$ (the latter acting as a ground ring for the
theory) as the $(n+1)$-fold tensor power of $\Omega_0$ over
$k_0$. Thus, the whole structure of the analyzer may be thought of as
embodied in the system $\Omega=(k_0,\Omega_0)$, where $k_0$ is a
commutative ring (with unit, now), $\Omega_0$ a commutative
$k_0$-algebra with unit, endowed moreover with a composition operation
$(F,G)\mapsto F\circ G$, satisfying a bunch of simple axioms I don't
feel like writing down here. The simplest case of all of course
(corresponding to usual commutative algebra ``over $k_0$'' as a ground
ring, with no extra structure on commutative algebras with unit over
$k_0$) is $\Omega_0=k[T]$, with the usual composition of polynomials,
$T$ acting as the two-sided unit for composition. In the general case,
$\Omega_0$ and its tensor powers $\Omega_n$ over $k_0$ are going to
play the part played by polynomial rings in ordinary commutative
algebra. There should be a ready generalization, in this spirit, of
taking the symmetric algebra of a $k_0$-module (which, for a module
free and of finite type will yield an $\Omega$-structure isomorphic to
one of the $\Omega_n$'s). An $\Omega$-structure on a set $k$ amounts
to giving a structure of a commutative $k_0$-algebra with unit on $k$,
plus a map
\[\Omega_0\to\Maps(k,k), \quad F\mapsto (x\mapsto F(x))\]
compatible\pspage{507} with the structures of $k_0$-algebras as well
as composition operations, and satisfying moreover two conditions for
\[F(x+y), \quad\text{resp.~$F(xy)$}\]
in terms of two diagonal maps
\[
  \begin{tikzcd}[cramped,sep=small]
    \Omega_0 \ar[r,shift left,"\Delta\suba"]
    \ar[r,shift right,"\Delta\subm"'] &
    \Omega_1 = \Omega_0 \otimes_{k_0} \Omega_0
  \end{tikzcd}
\]
(which may be described in terms of the composition structure on
$\Omega_0$, as expressing the compositions $F\circ (G'+G'')$ resp.\
$F\circ(G'G'')$), and moreover one trivial condition for $F(\lambda)$
when $\lambda$ in $k$ comes from $k_0$, namely compatibility of the
map $k_0\to k$ with the operations of $\Omega_0$ on both $k_0$ and
$k$. (NB\enspace the operation of $\Omega_0$ on $k_0$ is defined by
\[F(\lambda)=F\circ\lambda,\]
where $k_0$ is identified to a subring of $\Omega_0$.)

I didn't pursue much further these ponderings, just one digression
among many in the main line of investigation! I also read through the
preprint of David W.\ Jones on
Poly-$T$-complexes,\scrcomment{\textcite{Jones1983}} which Ronnie
Brown (acting as David Jones' supervisor) had sent me a while
ago. There he develops a notion of polyhedral cells, with a view of
using these instead of simplices or cubes for doing combinatorial
homotopy theory. As I had pondered a little along this direction (cf.\
sections \ref{sec:91}, topic \ref{q:91.8}, and section \ref{sec:93}),
I was hoping that some of the perplexities I had been meeting would be
solved in David Jones' notes -- for instance that there would be handy
criteria for a category $M$ made up with such polyhedral cells to be a
weak test category, namely that objects of $M\upvee$ may be used as
models for homotopy types; also, that the ``standard'' chain complex
constructed in section \ref{sec:93} is indeed an ``abelianizator'' for
$M$, i.e., may be used for computing homology of objects of
$M\upvee$. David Jones' emphasis, however, is a rather different one
-- he seems mainly interested in generalizing the theory of ``thin''
structures of M.~K.~Dakin\scrcomment{\textcite{Dakin1977}} from the
simplicial to the more general polyhedral set-up and prove a
corresponding equivalence of categories. Thus, my perplexities remain
-- they are admittedly rather marginal in the main line of thought,
and I doubt I'll stop to try and solve them.

\hangsection{The basic pair of adjoint functors
  \texorpdfstring{$\smash{\widetilde K:\Hotabz \leftrightarrows \Hotz:
    \LtH_*}$}{Hotab0<->Hot0}.}\label{sec:121}%
The\pspage{508} tentative approach towards defining and studying
``schematic'' homotopy types I have been following lately relies
heavily on a suitable notion of ``linearization'' of such homotopy
types. One can imagine that many different approaches (for instance
via De~Rham complexes with divided powers, or additive small
categories with diagonal maps) may be devised for ``schematic''
homotopy types, but in any case it seems likely that a suitable notion
of linearization will play an important role. It may be worthwhile
therefore to try and pin down the wished-for main features of such a
theory, with the hope maybe of getting an axiomatic description for
it, with a corresponding unicity statement. Before doing so, the first
thing to do seems to review some main formal features of linearization
for ordinary (``discrete'') homotopy types.

Recall the definition
\begin{equation}
  \label{eq:121.1}
  \HotabOf \eqdef \D_\bullet\Ab =
  \begin{tabular}[t]{@{}l@{}}
    derived category of the category of\\
    abelian chain complexes, with\\
    respect to quasi-isomorphisms,
  \end{tabular}
  \tag{1}
\end{equation}
and the two canonical functors
\begin{equation}
  \label{eq:121.2}
  \begin{tikzcd}[baseline=(O.base),cramped]
    \HotOf \ar[r,shift left,"\LH_\bullet"] &
    |[alias=O]| \HotabOf \ar[l,shift left,"K_\pi"]
  \end{tikzcd},\tag{2}
\end{equation}
Where $\LH_\bullet$ is the ``abelianization functor'', and $K_\pi$ is
defined via the Kan-Dold-Puppe functor, associating to a chain complex
the corresponding semisimplicial abelian group. The diagram
\eqref{eq:121.2} may be viewed (up to equivalence) as deduced from the
corresponding diagram
\begin{equation}
  \label{eq:121.2prime}
  \begin{tikzcd}[baseline=(O.base),cramped]
    \Simplexhat \ar[r,shift left,"W"] &
    |[alias=O]| \Simplexhatab \ar[l,shift left,"\DP"]
  \end{tikzcd}\quad(\equeq\Ch_\bullet\Ab)
  \tag{2'}
\end{equation}
by passing to the suitable localized categories $\HotOf$ and
$\HotabOf$. In the diagram \eqref{eq:121.2prime} $W$ is left adjoint
to $\DP$, which is now just the forgetful functor. One main fact about
\eqref{eq:121.2} is
\begin{equation}
  \label{eq:121.3}
  \text{$K_\pi$ is right adjoint to $\LH_\bullet$,}\tag{3}
\end{equation}
which is just a neater way for expressing the familiar fact that for
given abelian group $\pi$ and natural integer $n$, the object
$K(\pi,n)$ in $\HotOf$ (namely the image by $K_\pi$ of the chain
complex $\pi[n]$ reduced to $\pi$ in degree $n$) represents the
cohomology functor
\[X\mapsto \mathrm H^n(X,\pi) \simeq
  \Hom_{\HotabOf}(\LH_\bullet(X),\pi[n]).\]
The functor $K_\pi$ may be called the ``\emph{Eilenberg-Mac~Lane
  functor}'', as its values are immediate generalizations of the
Eilenberg-Mac~Lane objects\pspage{509} $K(\pi,n)$. As any object in
$\HotabOf$ is isomorphic to a product of objects $\pi[n]$, it follows
that in order to check the adjunction formula between $\LH_\bullet$ and
$K_\pi$ it is enough to do so for objects in $\HotabOf$ of the type
$\pi[n]$, which is the ``familiar fact'' just recalled. The notation
$K_\pi$ in \eqref{eq:121.2} is meant to suggest the Eilenberg-Mac~Lane
$K(\pi,n)$ object generalized by the objects $K_\pi(L_\bullet)$, and
also to recall that we recover the homology invariants of $L_\bullet$
from $K_\pi(L_\bullet)$ via the $\pi_i$ invariants, by the formula
\begin{equation}
  \label{eq:121.4}
  \pi_i(K_\pi(L_\bullet))\simeq\mathrm H_i(L_\bullet),\tag{4}
\end{equation}
which implies by the way that the functor $K_\pi$ is ``conservative'',
i.e., a map in $\HotabOf$ which is transformed into an isomorphism is
an isomorphism. This should not be confused with the stronger property
of being fully faithful, or equivalently of the left adjoint
$\LH_\bullet$ being a localization functor, or equivalently still, the
adjunction morphism
\begin{equation}
  \label{eq:121.5}
  \LH_\bullet(K_\pi(L_\bullet))\to L_\bullet\tag{5}
\end{equation}
being an isomorphism in $\HotabOf$, which is definitely false!

The other adjunction morphism
\begin{equation}
  \label{eq:121.6}
  X\to K_\pi(\LH_\bullet(X)) \eqdef X\subab\tag{6}
\end{equation}
is still more interesting, its effect on the homotopy invariants
$\pi_i$ are the Hurewicz homomorphisms
\begin{equation}
  \label{eq:121.7}
  \pi_i(X)\to \mathrm H_i(X) \eqdef \pi_i(X\subab) \simeq \mathrm
  H_i(\LH_\bullet(X));\tag{7} 
\end{equation}
introducing the homotopy fiber of \eqref{eq:121.6} (in the case of
pointed homotopy types) and denoting by $\gamma_i(X)$ its homotopy
invariants, we get the exact sequences of J.~H.~C.~Whitehead (as
recalled in a letter from R.~Brown I just got)
\begin{equation}
  \label{eq:121.7prime}
  \cdots \to\gamma_i(X) \to \pi_i(X) \to \tH_i(X) \to
  \gamma_{i-1}(X) \to \cdots,\tag{7'}
\end{equation}
where $\tH_i$ denotes the ``reduced'' homology group of a
\emph{pointed} homotopy type, equal to $\mathrm H_i$ for $i\ne0$ and
to $\Coker(\mathrm H_0(\mathrm{pt}) \to \mathrm H_0(X))\simeq \mathrm
H_0(X)/\bZ$ for $i=0$.

The case of pointed homotopy types seems of importance for schematic
homotopy types, and deserves some extra mention and care. We may
factor diagram \eqref{eq:121.2} into
\begin{equation}
  \label{eq:121.8}
  \begin{tikzcd}[baseline=(O.base)]
    \Simplexhat \ar[r,shift left,"\alpha"] &
    {\Simplexhat}^\bullet \ar[l, shift left,"\beta"]
    \ar[r,shift left,"\widetilde W"] &
    |[alias=O]| \Simplexhatab \ar[l, shift left, "\widetilde\DP"]
  \end{tikzcd},\tag{8}
\end{equation}
where $\Simplexhat^\bullet$ is the category of \emph{pointed}
semisimplicial complexes, $\beta$ the forgetful functor from these to
non-pointed complexes, and $\alpha$ its left adjoint, which may be
interpreted as\pspage{510}
\begin{equation}
  \label{eq:121.9}
  \alpha(X_*) = X_* \amalg e_*,\tag{9}
\end{equation}
where $e_*$ is the final object of \Simplexhat, and the second member
is pointed by its summand $e_*$. The functor $\widetilde\DP$ comes
from applying componentwise the obvious functor from $\Ab$ to
$\pSets$ (pointed sets), $\widetilde W$ is its left adjoint. Passing to
the suitable localized categories, we get from \eqref{eq:121.8}
\begin{equation}
  \label{eq:121.10}
  \begin{tikzcd}[baseline=(O.base)]
    \HotOf \ar[r,shift left,"\alpha"] \ar[rr,bend left,"\LH_\bullet"] &
    \HotOf^\bullet \ar[l,shift left,"\beta"]
    \ar[r,shift left,"\LtH_\bullet"] &
    |[alias=O]| \HotabOf\ar[l,shift left,"\widetilde K_\pi"]
    \ar[ll,bend left,"K_\pi"]
  \end{tikzcd},\tag{10}
\end{equation}
factoring \eqref{eq:121.2}, where $\alpha$ is now defined by the
formula similar to \eqref{eq:121.9}
\begin{equation}
  \label{eq:121.9prime}
  \alpha(X) = X\amalg e,\tag{9'}
\end{equation}
where $e$ denotes the final object of $\HotOf$, and is used for
defining the pointed structure of the second member. The functor
$\widetilde W$ in \eqref{eq:121.8} can be described (as is seen
componentwise) as
\begin{equation}
  \label{eq:121.11}
  \widetilde W(X_*) = W(\beta(X_*)) / \Imm W(e_*),\tag{11}
\end{equation}
where
\[W(e_*) \to W(X_*), \quad\text{i.e.,}\quad
  \bZ^{(e)} \to \bZ^{(X)}\]
is deduced from the pointing map $e_*\to X_*$. Accordingly, we get an
expression
\[\LtH_\bullet(X) \simeq \LH_\bullet(X)/\LH_\bullet(e),\]
more accurately, an exact triangle
\begin{equation}
  \label{eq:121.12}
  \begin{tabular}{@{}c@{}}
    \begin{tikzcd}[baseline=(O.base),column sep=small]
      & \LtH_\bullet(X) \ar[dl] & \\
      \LH_\bullet(e)=\bZ[0] \ar[rr] & &
      |[alias=O]| \LH_\bullet(X)\ar[ul]
    \end{tikzcd},
  \end{tabular}\tag{12}
\end{equation}
where $X$ is any object in $\HotOf^\bullet$ and $\LH_\bullet(X)$ is
short for $\LH_\bullet(\beta(X))$. Of course, the functors $\LtH$ and
$\widetilde K_\pi$ are still adjoint (one hopes!), hence for any pointed
homotopy type $X$ an adjunction map in $\HotOf^\bullet$
\begin{equation}
  \label{eq:121.6prime}
  X\to\widetilde K_\pi(\LH_\bullet(X)),\tag{6'}
\end{equation}
and \eqref{eq:121.7} and \eqref{eq:121.7prime} are deduced from
\eqref{eq:121.6prime} and its homotopy fiber, rather than from
\eqref{eq:121.6} where it doesn't really make sense because of lack of
canonical base points for taking $\gamma_i$'s and homotopy fibers.

\bigbreak

\presectionfill\ondate{2.10.}\pspage{511}\par

\hangsection[Rambling reflections on $\LtH_*$, Postnikov
invariants, \dots]{Rambling reflections on
  \texorpdfstring{$\LtH_*$}{LH}, Postnikov invariants,
  \texorpdfstring{$S(H,n)$}{S(H,b)}'s -- and on the non-existence of a
  ``total homotopy''-object \texorpdfstring{$\mathrm L\pi_*(X)$}{Lpi(X)} for
  ordinary homotopy types.}\label{sec:122}%
Since the last notes, I have been doing three days' scratchwork
(including today's) on various questions around abelianization in
general (for discrete homotopy types) and on Postnikov dévissage, in
connection with the review on some main formal properties of
abelianization, started upon in the previous section~\ref{sec:121}. I
didn't get anything really new for me, rather it was just part of the
necessary rubbing against the things, in order to get a better feeling
of what they are like, or what they are likely to be like -- what is
likely to be true, and what not. The most interesting, maybe, is that
I got an inkling of a fairly general version of a Kan-Dold-Puppe kind
of relationship, in terms of derived categories, valid presumably for
any local test category, and in particular for categories like
$\Simplex_{/X}$, with $X$ in \Simplexhat. It would be untimely,
though, to build up still more the (already pretty high) tower of
digressions, and for the time being I'll stick to what is relevant
strictly to my immediate purpose -- namely getting through with the
wishful thinking about schematization! Thus, I'll be content to work
with the category \Simplex{} in order to construct models for homotopy
types and perform constructions with them (such as abelianization),
without getting involved at present in looking up how much is going
over (and how) to the case of more general small categories $A$\dots

It occurred to me that the variant $\widetilde W$ or $\LtH_\bullet$
for the abelianization functors $W$ and $\LH_\bullet$, introduced in
the previous section using a pointed structure for the argument $X$ in
\Simplexhat{} or $\HotOf^\bullet$, could be advantageously defined
without this extra structure. The construction for $\widetilde W$
which follows goes through indeed in any topos whatever (not only
\Simplexhat). Let $X$ be an object in \Simplexhat, then there is a
canonical augmentation map
\[\varepsilon : W(X) =\bZ^{(X)} \to \bZ_\Simplex\]
towards the constant semisimplicial group $\bZ_\Simplex$,
corresponding to the constant map
\[X\to \bZ_\Simplex\]
with value $1$. We thus get a functorial exact sequence
\begin{equation}
  \label{eq:122.1}
  0 \to \widetilde W(X) \to W(X) \xrightarrow\varepsilon{} \bZ \to
  0,\tag{1} 
\end{equation}
where $\widetilde W(X)$ denotes the kernel of the augmentation above,
hence an extension of \bZ{} by $\widetilde W(X)$, or (what amounts
essentially to the same) a torsor under the group object $\widetilde
W(X)$, which we'll denote by
\begin{equation}
  \label{eq:122.2}
  W(X)(1) = \varepsilon^{-1}(1).\tag{2}
\end{equation}
Moreover,\pspage{512} the canonical map
\[X\to W(X)=\bZ^{(X)}\]
factors (by construction) through
\begin{equation}
  \label{eq:122.3}
  X\to W(X)(1)\quad(\hookrightarrow W(X)),\tag{3}
\end{equation}
and this map is universal for all maps of $X$ into torsors under
abelian group objects of \Simplexhat. It is an immediate consequence
of Whitehead's theorem that a weak equivalence $X\to Y$ in
\Simplexhat{} induces weak equivalences in \Simplexhat
\begin{equation}
  \label{eq:122.4}
  \widetilde W(X)\to\widetilde W(Y), \quad
  W(X)(1) \to W(Y)(1)\tag{4}
\end{equation}
(and also $W(X)(n)\to W(Y)(n)$ for any $n\in\bZ$), and accordingly
that the map between the normalized chain complexes corresponding to
$\widetilde W(X)$ and $\widetilde W(Y)$, namely by definition
\begin{equation}
  \label{eq:122.4prime}
  \LtH_\bullet(X) \to \LtH_\bullet(Y)\tag{4'}
\end{equation}
is a quasi-isomorphism. Thus, we get a functor
\[\LtH_\bullet: \HotOf\to\HotabOf,\]
whose composition $\LtH_\bullet\circ \beta$ with the forgetful functor
\[\beta:\HotOf^\bullet\to\HotOf\]
is canonically isomorphic to the functor $\LtH_\bullet$ of
p.~\ref{p:510} (formula \ref{eq:121.9}). More specifically, when $X$
is in $\Simplexhat^\bullet$, then the map
\[\bZ=\bZ^{(e)} \to \bZ^{(X)}\quad\text{deduced from $e\to X$}\]
defines a splitting of the extension \eqref{eq:122.1}, i.e., an
isomorphism
\begin{equation}
  \label{eq:122.5}
  \bZ^{(X)} = W(X) \simeq \widetilde W(X)\oplus \bZ,\tag{5}
\end{equation}
and accordingly, $\widetilde W(X)$ may be interpreted, at will, as a
quotient group of $W(X)$ (which we did on p.~\ref{p:510}), or as a
subobject of $W(X)$, which is the better choice, because it makes good
sense without using any pointed structure. Again, without using the
pointed structure, we get a canonical exact triangle in $\D\Ab$
interpreting \eqref{eq:122.1}
\begin{equation}
  \label{eq:122.6}
  \begin{tikzcd}[column sep=small]
    & \bZ\ar[dl] & \\
    \LtH_\bullet(X) \ar[rr] & &
    \LH_\bullet(X) \ar[ul]
  \end{tikzcd}\tag{6}
\end{equation}
replacing (\eqref{eq:121.12} p.~\ref{p:510}), which (in the case
considered there, namely $X$ pointed) should be replaced by the more
precise relationship
\begin{equation}
  \label{eq:122.5prime}
  \LH_\bullet(X)\simeq \LtH_\bullet(X) \oplus\bZ,\tag{5'}
\end{equation}
an isomorphism functorial for $X$ in $\HotOf^\bullet$.

Truth\pspage{513} to tell, these generalities are more interesting
still for an argument $X$ in a category like $\Simplexhat_{/Y}\equeq
(\Simplex_{/Y})\uphat$, where $Y$ is an arbitrary object in
\Simplexhat, rather than in \Simplexhat{} itself. In the latter case,
$X$ always admits a pointed structure, i.e., a section over the final
object (provided $X$ is non-empty), hence there always exists a
splitting for \eqref{eq:122.1}, and hence one for \eqref{eq:122.6};
whereas for an object $X$ over $Y$, there does not always exist a
section over $Y$, and accordingly it may well happen that the
extension similar to \eqref{eq:122.1} (but taken ``over $Y$'') does
not split, i.e., that the torsor $W_{/Y}(X)(1)$ is not trivial. The
class of this torsor is an element
\begin{equation}
  \label{eq:122.7}
  c(X/Y) \in \mathrm H^1(Y, \widetilde W_{/Y}(X)),\tag{7}
\end{equation}
a very interesting invariant indeed, giving rise to the Postnikov
invariants when $Y$ is one of the $X_n$'s occurring in the Postnikov
dévissage of $X$. It was while trying to understand the precise
relationship between a cohomology group as in \eqref{eq:122.7}, with
``continuous'' coefficients (by which I mean to suggest that
$\widetilde W_{/Y}(X)$ is viewed intuitively as a fiberspace over the
``space'' $X$, whose fibers are topological abelian groups which are
by no means discrete, but got a bunch of non-vanishing $\pi_i$'s!),
and a cohomology group with ``discrete'' coefficients, more accurately
with coefficients in a complex of chains in the category of abelian
sheaves over $X$, i.e., over $\Simplex_{/X}$ whose homology sheaves
are definitely of a ``discrete'' nature (for instance, they are
locally constant if $X\to Y$ is a Kan fibration, and their fibers are
\bZ-modules of finite type provided we make a mild finite-type
assumption on the fibers of $X\to Y$), that I got involved in a more
general reflection on a Kan-Dold-Puppe type of relationship. I hope to
come back to this when getting back to the general review of
linearization begun in part \ref{ch:V} of these notes, and abruptly
interrupted after section~\ref{sec:109}, when getting caught
unwittingly by the enticing mystery of schematization!

Another afterthought to the previous notes is that I am going to
denote by $K, \widetilde K$, the functors
\begin{equation}
  \label{eq:122.7bis}
  K:\HotabOf\to\HotOf, \quad
  \widetilde K:\HotabOf\to\HotOf^\bullet, \quad
  K=\beta\circ\widetilde K,\tag{7}
\end{equation}
denoted by $K_\pi$, $\widetilde K_\pi$ in section~\ref{sec:121}. This
gives a formula such as
\begin{equation}
  \label{eq:122.8}
  K(\pi[n]) = K(\pi,n),\tag{8}
\end{equation}
where $\pi[n]$ denotes the chain complex in \Ab{} which is $\pi$
placed in degree $n$. I was intending to complement the former
notation $K_\pi$ by a similar notation $K_{\mathrm H}$ (where $\mathrm
H$ stands for ``homology'', as $\pi$ was standing\pspage{514} for
``homotopy'', and $K$ means ``complex''), but I finally found the
notation $S$ (suggestive of ``spheres'') more congenial, giving rise
to
\begin{equation}
  \label{eq:122.8prime}
  S(\mathrm H[n]) = S(\mathrm H,n) =
  \begin{tabular}[t]{@{}l@{}}
    sphere-like homotopy type whose $\LtH_\bullet$ \\
    is isomorphic to $\mathrm H[n]$.
  \end{tabular}\tag{8'}
\end{equation}
But I am anticipating somewhat on some of the wishful thinking still
ahead, involving the would-be description (preferably in the schematic
set-up, but more elementarily maybe in the discrete one of ordinary
homotopy types) of a functor
\[S: \HotabOf = \D_\bullet\Ab \to \HotOf,\]
whose most important formal property should be
\[\LH_\bullet(S(L_\bullet)) \simeq L_\bullet,\]
compare with formula \eqref{eq:122.8prime} for
$L_\bullet=\pi[n]$.\scrcomment{in \eqref{eq:122.8prime}, the $\mathrm
  H[n]$ was originally written $\pi[n]$. I'm not sure what to make of
  this\dots} However, such a formula cannot hold for any $L_\bullet$,
i.e., an arbitrary chain complex in \Ab{} is not quasi-isomorphic to
an object $\LH_\bullet(X)$ with $X$ in \Hot, as a necessary condition
for this is that $\mathrm H_0(L_\bullet)$ should be free. Thus, we
better restrict to the $0$-connected case, and take $L_\bullet$
subject to
\[\mathrm H_0(L_\bullet) = 0,\quad\text{i.e., $L_\bullet$ in
    $\D_{\ge1}\Ab \eqdef \HotabOfzc$,}\]
and describe $S$, or better still $\widetilde S$, as a functor
\begin{equation}
  \label{eq:122.9}
  \widetilde S: \HotabOfzc \to \HotOfzc \hookrightarrow \HotOf^\bullet,\tag{9}
\end{equation}
subject to the condition
\begin{equation}
  \label{eq:122.10}
  \LtH_\bullet(\widetilde S(L_\bullet)) \simeq L_\bullet,\tag{10}
\end{equation}
giving rise to
\begin{equation}
  \label{eq:122.10prime}
  \tH_i(\widetilde S(L_\bullet)) \simeq \mathrm H_i(L_\bullet),\tag{10'}
\end{equation}
which should be compared to the familiar formula
\begin{equation}
  \label{eq:122.11}
  \pi_i(\widetilde K(L_\bullet)) \simeq \mathrm H_i(L_\bullet).\tag{11}
\end{equation}

Maybe formula \eqref{eq:122.10} is not quite enough for characterizing
the functor $\widetilde S$ up to unique isomorphism, and I confess I
didn't try to construct in some canonical way a functor
\eqref{eq:122.9} satisfying \eqref{eq:122.10} -- and I am not quite
sure even whether such a functor exists. One main evidence for
existence of such a functor would be the existence and unicity (up to
unique isomorphism) of the homotopy types $S(H,n)$ in
\eqref{eq:122.8prime} -- and I have a vague remembrance of having
looked through a paper, fifteen or twenty years ago, where such spaces
were indeed introduced and studied. Thus, a functor $\widetilde S$ is
perhaps nowadays more or less standard knowledge in\pspage{515}
homotopy theory. I was led to postulate the existence of a functor
$\widetilde S$ in the context of \emph{schematic} homotopy types, by
reasons of symmetry in the $\LH_\bullet$ and $\mathrm L\pi_\bullet$
formalism, which I'll try to get through a little later. It should be
clear from the outset, though, that in the schematic set-up, even the
mere existence of (schematic) homotopy types $S(\bZ,n)$ (corresponding
to ordinary spheres), and even for $n=2$ or $3$ only, is very far from
trivial, and as a matter of fact isn't even known (in the set-up of
semisimplicial unipotent bundles). This shows at the same time that if
such a functor $\widetilde S$ can actually be constructed in the
schematic case, it is likely to be a lot more interesting still than
in the discrete case, as it will presumably give information on
homotopy groups of spheres (via the description of the $\pi_i$'s of a
semisimplicial unipotent bundle via the Lie functor (cf.\
section~\ref{sec:118})).

There is a fourth tentative functor in between the categories
$\HotabOf$ and $\HotOf$, or more accurately, from a suitable
subcategory $\HotOf^\bullet(0)$ of $\HotOf^\bullet$ (corresponding to
$0$-connected pointed homotopy types whose $\pi_1$ is abelian and acts
trivially on the higher $\pi_i$'s) to $\HotabOfzc=\HotabOf(0)$,
strongly suggested by the schematic case of section~\ref{sec:118},
namely a functor
\begin{equation}
  \label{eq:122.star}
  \Lpi_\bullet: \HotOf^\bullet(0) \to \HotabOf^\bullet(0)\text{???,}\tag{*}
\end{equation}
whose main functorial property should be
\begin{equation}
  \label{eq:122.12}
  \mathrm H_i(\Lpi_\bullet(X)) \simeq \pi_i(X).\tag{12}
\end{equation}
To play really safe, one might hope such a functor to be defined at
any rate for $1$-connected pointed homotopy types. Another equally
important property, suggested by the schematic set-up as well as by
\eqref{eq:122.11} and \eqref{eq:122.12}, is the formula
\begin{equation}
  \label{eq:122.13}
  \Lpi_\bullet(\widetilde K(L_\bullet)) \simeq L_\bullet,\tag{13}
\end{equation}
similar to formula \eqref{eq:122.10} (with the pair $(\widetilde
S,\LtH_\bullet)$ replaced by the pair $(\widetilde
K,\Lpi_\bullet)$). Of course, for $X$ of the type $\widetilde
K(L_\bullet)$, \eqref{eq:122.12} follows from \eqref{eq:122.13}

\bigbreak

\presectionfill\ondate{3.10.}\par

We'll see, though, that a functor $\Lpi_\bullet$ as in
\eqref{eq:122.star}, satisfying the properties above, \emph{does not
  exist}. Indeed, applying $\Lpi_\bullet$ to the adjunction morphism
\[X\to \widetilde K(\LtH_\bullet(X))\]
(cf.\ p.~\ref{p:509}, \eqref{eq:121.6}), and using \eqref{eq:122.13},
we should get a Hurewicz homomorphism\pspage{516}
\begin{equation}
  \label{eq:122.14}
  \Lpi_\bullet(X)\to\LtH_\bullet(X)\tag{14}
\end{equation}
more precise than the separate homomorphisms
\[\pi_i(X)\to\tH_i(X),\]
and applying this to an object $\widetilde K(L_\bullet)$ and applying
\eqref{eq:122.13} again, we should get a functorial homomorphism
\begin{equation}
  \label{eq:122.15}
  L_\bullet \to \LtH_\bullet(\widetilde K(L_\bullet))\tag{15}
\end{equation}
in opposite direction from the adjunction morphism
\begin{equation}
  \label{eq:122.16}
  \LtH_\bullet(\widetilde K(L_\bullet)) \to L_\bullet\tag{16}
\end{equation}
(p.~\ref{p:509}, \eqref{eq:121.5}), the composition of the two being
the identity in $L_\bullet$. In other words, \eqref{eq:122.15} should
be a canonical splitting of the natural adjunction morphism
\eqref{eq:122.16}. Take for instance $L_\bullet=\pi[n]$, then the
first member of \eqref{eq:122.16} is the Eilenberg-Mac~Lane homology
$\LtH_\bullet(\pi,n)$, whose first non-vanishing homology group is
$\pi$, and \eqref{eq:122.16} is just the canonical augmentation
\[\LtH_\bullet(\pi,n)\to\pi[n],\]
and taking the $\Ext_\bZ^i$ of both members with $M[n]$ ($M$ any
object in \Ab) yields the transposed canonical homomorphism
\begin{equation}
  \label{eq:122.17}
  \Ext_\bZ^i(\pi,M)\to\mathrm H^{n+i}(\pi,n;M),\tag{17}
\end{equation}
which is an isomorphism for $i=0$ (whereas for $i<0$ both members are
zero). For $i=1$, we get an exact sequence (``universal
coefficients'')
\begin{equation}
  \label{eq:122.18}
  0 \to \Ext_\bZ^i(\pi,M) \to \mathrm H^{i+1}(\pi,n;M) \to
  \Hom_\bZ(\mathrm H_{i+1}(\pi,n), M) \to 0.\tag{18}
\end{equation}
A canonical splitting of \eqref{eq:122.15} would yield a canonical
splitting of this exact sequence, which definitely looks somehow as
``against nature''! Surely, it shouldn't be hard to find, for any
given integer $n\ge1$, suitable abelian groups $\pi$ and $M$, such
that there does not exist a splitting of \eqref{eq:122.18} stable
under the action of the group
\[G = \Aut(\pi) \times \Aut(M)\]
(acting by ``transport de structure'' on the three terms of
\eqref{eq:122.18}); presumably even, looking at $\Aut(\pi)$ should be
enough. As I am not at all familiar with Eilenberg-Mac~Lane
cohomology, except a little in the case $n=1$, i.e., for usual group
cohomology, I didn't try to construct an example for any $n$, only one
for $n=1$, i.e., in a non-simply connected case, which again is a
little less convincing as if it was one for $n=2$\dots

Thus, let's take $n=1$, a vector space of finite dimension over the
prime field $\bF_2$, and $M=\bF_2$, in this case standard
calculations\pspage{517} show that the exact sequence
\eqref{eq:122.18} can be identified with the familiar exact sequence
\begin{equation}
  \label{eq:122.18prime}
  0\to V'\to \Sym^2(V')\to \bigwedge^2 V'\to 0,\tag{18'}
\end{equation}
where the first arrow associates to every linear form on $V$ its
square, i.e., the same form (!) but viewed as being a quadratic form
on $V$; and the second arrow associates to every quadratic form on $V$
the associated bilinear form, which is alternate because of
char.~$2$. It is well-known I guess that for $\dim V\ge2$, there does
not exist a splitting of \eqref{eq:122.18prime} which is stable under
$\Aut(V)\simeq\Aut(V')$.

\begin{remark}
  If we admit that a similar example can be found for non-splitting of
  \eqref{eq:122.18}, with $n\ge2$ (which shouldn't be hard I guess for
  someone familiar with Eilenberg-Mac~Lane homology), this shows that
  there does \emph{not} exist (as was contemplated at the beginning,
  cf.\ section~\ref{sec:111}) an equivalence of categories between the
  schematic $1$-connected homotopy types (defined via the model
  category $\scrM_1(\bZ)$ of semisimplicial unipotent bundles over
  \bZ{} satisfying $X_0=X_1=e$) and $1$-connected pointed usual
  homotopy types, as this would imply the existence of a functor
  $\Lpi_\bullet$, hence of a canonical splitting of
  \eqref{eq:122.18}. The same argument will show that the once
  hoped-for comparison theorem between usual Eilenberg-Mac~Lane
  homology (or cohomology), and the schematic one, is false, because
  in the schematic set-up $\Lpi_\bullet$ does exist, and the direct
  construction of \eqref{eq:122.15} (a functorial splitting of
  \eqref{eq:122.16}) is then anyhow a tautology. Thus after all, we
  don't have to rely on delicate results of Breen's on
  $\Ext^i(\bG\suba,\bG\suba)$ over the prime fields $\bF_p$ (as
  suggested by Illusie), in order to get this ``negative'' result,
  causing unreasonable expectations to crash\dots
\end{remark}

\bigbreak

\presectionfill\ondate{4.10.}\pspage{518}\par

\hangsection[The hypothetical complexes ${}^n\Pi_*=\Lpi_*(S^n)$, and
\dots]{The hypothetical complexes
  \texorpdfstring{${}^n\Pi_*=\Lpi_*(S^n)$}{nPi=Lpi(Sn)}, and comments
  on homotopy groups of spheres.}\label{sec:123}%
For the last few days -- since I resumed the daily mathematical
ponderings (interrupted more or less for nearly one month), the wind
in my sails has been rather slack I feel. It has been nearly six weeks
from now that I started pondering on schematization and on schematic
homotopy types -- in the process I got rid of some misplaced
expectations, very well. Still, the unpleasant feeling remains of not
having really any hold yet on those would-be schematic homotopy types,
due mainly to my incapacity so far of performing (in the model
category of semisimplicial unipotent bundles say) the basic homotopy
operations of taking homotopy fibers and cofibers of maps. In lack of
this, I am not (``morally'') sure yet that there does exist indeed a
substantial reality of the kind I have been trying to foreshadow. I am
not even wholly clear yet of how to define the notion of ``weak
equivalence'' in t he model category $\scrM_0(k)$ (or in the smaller
one $\scrM_1(k)$, to play safe, as even the definition of $\scrM_0(k)$
isn't too clear yet) -- there are three ways to define weak
equivalences, using $\LH_\bullet$, or $\Lpi_\bullet$, or the sections
functor from $\scrM_1(k)$ to ordinary semisimplicial complexes, and it
isn't clear yet whether these are indeed equivalent. We may of course
call ``weak equivalence'' in $\scrM_1(k)$ a map which becomes an
isomorphism by any of these three functors. But I wouldn't say I am
wholly confident yet that there exists a localization $\HotOf_1(k)$ of
$\scrM_1(k)$, with respect to this notion of weak equivalence or some
finer one, in such a way that in $\HotOf_1(k)$ one may introduce the
two types of ``fibration'' and ``cofibration'' sequences with the
usual properties, and the sections functor
\[\HotOf_1(k)\to\HotOf_1\]
should respect this structure, and moreover give rise to functorial
isomorphisms
\[\mathrm H_i(\Lpi_\bullet(X_*)) \simeq \pi_i(X_*(k)),\]
nor even do I feel wholly confident that a theory of this kind can be
developed, possibly with different kinds of models for $k$-homotopy
types from the ones I have been using. At any rate, it seems to me
that the two kinds of conditions or features I have just been stating,
are indeed the crucial ones, plus (I should add) existence in
$\HotOf_1(k)$ of an object
\[S^2_k = S(2,k),\]
standing for the $2$-sphere, with
\[\LtH_\bullet(S^2_k) \simeq \bZ[2],\]
whose\pspage{519} image in $\HotOf_1$ (in case $k=\bZ$ say) should be
the (homotopy type of the) ordinary $2$-sphere. Taking suspensions,
we'll get from this the $n$-spheres $S^n_k=S(n,k)$ over $k$, for any
$n\ge2$. (NB\enspace In case the relevant homotopy theoretic
structures can be introduced in a suitable larger $\HotOf_0(k)$,
containing the objects $\widetilde K(L_\bullet)$ for $L_\bullet$ a
$0$-connected chain complex of $k$-modules, $\HotOf_0(k)$ contains an
object $S(1,k)=\widetilde K(1,k)$, and $S^2_k$ may be simply described
as its suspension.)\enspace One main motivation for trying to push
through a theory of schematic homotopy types may be the hope that this
may provide new insights into the homotopy groups of spheres. A first
interesting consequence would be that for any given sphere $S^n$ (in
the set-up now of ordinary homotopy theory), the series of all its
homotopy groups $\pi_i(S^n)=\Hom_{\HotOf}(S^i,S^n)$ may be viewed as
being the homology modules of an object of $\D_\bullet\Ab$, namely
$\Lpi_\bullet(S^n_\bZ)$, canonically associated to $S^n$. As was seen
by an example in the previous section, the existence of such a
\emph{canonical} representation is by no means a trivial or innocuous
fact, it is indeed definitely false for arbitrary homotopy types --
here then it would turn out as a rather special feature of the full
subcategory of $\HotOf^\bullet$ made up with the spheres $S^n$
($n\ge2$). Thus, for any $\alpha\in\pi_i(S^n)$, i.e., for any map
\[\alpha:S^i\to S^n\]
in $\HotOf^\bullet$, there should be an induced map
\[\Lpi_\bullet(\alpha):\Lpi_\bullet(S^i)\to\Lpi_\bullet(S^n)\]
with associativity condition for a composition
\[S^i \xrightarrow\alpha{} S^n \xrightarrow\beta{} S^m.\]
If we write
\[{}^n\Pi_\bullet \eqdef \Lpi_\bullet(S^n),\]
this defines on the set of chain complexes (more accurately, abelian
homotopy types) ${}^n\Pi_\bullet$ ($n\ge2$) quite a specific
structure, which merits to be investigated a priori (under the
assumption, of course, of the existence of a reasonable theory of
schematic homotopy types, satisfying the criteria above). Maybe after
all, it will turn out that the homotopy groups of spheres
\emph{cannot} be fitted into such an encompassing structure -- so much
the better, as this will show that the kind of theory I started trying
to dig out doesn't exist, which will clear up the situation a great
deal! But if there does exist such a canonical structure, it surely
shouldn't be ignored, and it should be pleasant work to try and pin
down exactly what extra information on homotopy groups of spheres
is\pspage{520} involved in such extra structure. It may be noted for
instance that, when taking $i=n$ hence $\pi_i(S^n)=\bZ$, we get an
operation of the multiplicative monoid $\bZ^{\mathrm{mult}}$ on
${}^n\Pi_\bullet$, whose action on the homology groups of
${}^n\Pi_\bullet$, i.e., on the homotopy groups $\pi_j(S^n)$, is
surely an important extra structure on these, which I hope has been
studied by the homotopy people extensively. Surely it has been known,
too, for a long time (at any rate since Artin-Mazur's foundational
work on profinite homotopy types) that when $\pi_j(S^n)$ is finite
(i.e., practically in all cases except $j=n$) this action comes from a
continuous action of the multiplicative monoid
\[ M={\bZ\uphat}^{\mathrm{mult}},\]
where $\bZ\uphat$ is the completion of \bZ{} with respect to the
topology defined by subgroups of finite index:
\[\bZ\uphat = \lim \bZ/n\bZ \simeq \prod_p \bZ_p,\]
where in the last member the product is taken over all primes, and
$\bZ_p$ denotes $p$-adic integers. Taking homotopy types over the ring
$k=\bZ\uphat$, we should get that this monoid $M$ operates on the
object
\[{}^n\Pi_\bullet \biggl(\bigoplus^l{}_\bZ \bZ\uphat\biggr) \quad\text{in
    $\D_\bullet(\AbOf_{\bZ\uphat})$,}\]
as this may be interpreted as $\Lpi_\bullet(S^n_k)$, and
$k^{\mathrm{mult}}$ operates on $S^n_k=S(n,k)$ for any ring $k$ (if
indeed $S(n,\pi)$ depends functorially on the variable $k$-module
$\pi$). At any rate, it is surely well-known that
$M={\bZ\uphat}^{\mathrm{mult}}$ operates on the profinite completions
of all spheres. For the odd-dimensional spheres, and taking the
restrictions of this operation to the largest subgroup
\[M^\times = {\bZ\uphat}^\times \simeq \prod_p \bZ_p^\times\]
of $M$, this group may be viewed as the profinite Galois group of the
maximal cyclotomic extension of the prime field \bQ{} (deduced by
adjoining all roots of unity), and its action on the profinite
completion ${S^n}\uphat$ of $S^n$ may be viewed as the canonical
Galois action, when $S^{n=2m-1}$ is interpreted as the homotopy type
of affine (complex) $m$-space minus the origin (which makes sense as a
scheme over the prime field \bQ). This was the interpretation in my
mind, when stating that the action of $\bZ^{\mathrm{mult}}$ and of its
completion ${\bZ\uphat}^{\mathrm{mult}}$ on the homotopy groups of
spheres was in important extra structure on these.

\hangsection[Outline of a program (second version): an autodual
\dots]{Outline of a program \texorpdfstring{\textup(}{(}second
version\texorpdfstring{\textup)}{)}: an autodual formulaire for the
basic ``four functors'' \texorpdfstring{$\LtH_*$, $\Lpi_*$,
  $\widetilde K$, $\widetilde S$}{LH, Lpi, K, S}.}\label{sec:124}%
I\pspage{521} guess I'll skip giving a more or (rather!) less complete
axiomatic description of $\HotOf_0^\bullet$ or of
$\HotOf_0^\bullet(k)$, in terms of the pair of adjoint functors
$\LtH_\bullet$ and $\widetilde K$ and Postnikov dévissage, although I
did go through this exercise lately -- it doesn't really shed new
light on the approach we are following here towards ``schematic''
homotopy types over arbitrary ground rings. In contrast with the two
other approaches I have been thinking of before (namely De~Rham
complexes with divided powers, and small categories with diagonal
maps, cf.\ sections \ref{sec:94} and \ref{sec:109} respectively), the
characteristic feature of this approach seems to be that it takes into
account the existence of the canonical functor
$\widetilde K: \D_\bullet\Ab\to\Hot$, paraphrased by a functor
\begin{equation}
  \label{eq:124.1}
  \widetilde K:\D_\bullet(\AbOf_k)\zc\to\HotOf_0(k) \tag{1}
\end{equation}
compatible with ring extension $k\to k'$; while all three approaches
have in common that they yield a paraphrase
\begin{equation}
  \label{eq:124.2}
  \LtH_\bullet:\HotOf_0(k)\to\D_\bullet(\AbOf_k)\zc\tag{2}
\end{equation}
of the (total) homology functor, in a way again compatible with
extension of ground rings. In the two earlier approaches this latter
functor comes out in a wholly tautological way, whereas in the present
approach via semisimplicial schemes, it isn't quite so tautological
indeed. If we want to define a functor $\widetilde K$ in the context
say of De~Rham algebras with divided powers, one may think of
associating first to a chain complex in $\AbOf_k$ the corresponding
semisimplicial $k$-module, view this as a semisimplicial set, and take
its De~Rham complex with divided powers and coefficients in $k$ (as we
want an object over $k$). But this is visibly silly, except possibly
for $k=\bZ$, as the result doesn't depend on the $k$-module structure
of the chain complex we start with. It is doubtful, anyhow, that in
this context, the tautological functor \eqref{eq:124.2} admits a right
adjoint (which we then would call $\widetilde K$ of course). --\enspace
Another noteworthy difference between the present approach towards
defining ``schematic'' homotopy types, and the two earlier ones
(``earlier'' in these notes, at any rate!) is that for $k=\bZ$ say,
the category of ``homotopy types over \bZ'' we get here maps into the
category of ordinary homotopy types, whereas it was the opposite with
the two other approaches.

What, however, I still would like to do, is a little more daydreaming
about the expected formal properties of the basic functors in between
the two categories appearing in \eqref{eq:124.1} and \eqref{eq:124.2}
-- namely, essentially, between ``schematic homotopy types'' and
``abelian'' ones. It will be\pspage{522} convenient to denote the
category of the latter by
\begin{equation}
  \label{eq:124.3}
  \HotabOf_0(k)\tag{3}
\end{equation}
rather than $\D_\bullet\Ab\zc$, the subscript $0$ denoting the
$0$-connectedness condition. The reason for this change in notation is
that, according to the choices made in the basic definitions, the
``$k$-linear'' algebraic interpretation of this category may vary
still. Strictly speaking, in the approach we have been following in
terms of semisimplicial unipotent bundles, the category
\eqref{eq:124.3} cannot really be described as the subcategory of
$0$-connected objects of $\D_\bullet(\AbOf_k)$, we have seen, rather,
that in order to define a functor $\LtH_\bullet$, we have to work with
chain complexes in $\Pro(\AbOf_k)$ rather than in $(\AbOf_k)$; cf.\
section \ref{sec:115}, where this point is made rather
forcefully. (Recall that taking projective limits instead wouldn't
help, because then the linearization functor $\LtH_\bullet$ would no
longer commute to ground ring extension!)\enspace The trouble then is
that in order to define
\[\widetilde K:\HotabOf_0(k) \to \HotOf_0(k),\]
we are obliged to enlarge accordingly the category of models used for
describing $\HotOf_0(k)$, namely take semisimplicial
\emph{pro}-unipotent bundles, rather than just semisimplicial
unipotent bundles. This sudden invasion of the picture by pro-objects
may appear forbidding -- but maybe it will appear less so, or at any
rate kind of natural and inescapable, if we recall that when it will
come to working with homotopy types of general topoi, these are anyhow
``prohomotopy types'' rather, and they can be described only by
working systematically in terms of pro-objects of various
categories. We'll see, however, that there may be a way out of the
``pro''-mess, by using a slightly different, more or less dual
approach towards the idea of ``unipotent bundles''. Of course, one may
also think of using the ``pointed linearization'' $L\subpt$ for
$\LH_\bullet$, instead of $L$, but as was pointed out in section
\ref{sec:117}, this will lead to ``formal'' homology and cohomology
invariants rather than to ``schematic'' ones, and the relation between
these and corresponding invariants for ordinary homotopy types will be
looser still; at any rate, it is suited for describing Postnikov
dévissage in the set-up of ``formal schematic'' homotopy types only,
not for schematic homotopy types. Still, the theory of the former may
have an interest in its own right, even though its relation to
ordinary homotopy types isn't so clear at present. Thus, we may state
that there are around at present three or four different candidates
for possibly fitting the daydreaming I want to write down.\pspage{523}
With this in mind, we shouldn't be too specific (for the time being)
about the exact meaning to be given to the symbols $\HotabOf_0(k)$ and
$\HotOf_0(k)$. It may be safer to replace these by $\HotabOf_1(k)$ and
$\HotOf_1(k)$ (where the subscript $1$ means ``$1$-connected''), all
the more so as we haven't been able yet, in the context of unipotent
bundles, to give even a tentative precise definition of $\HotOf_0(k)$,
i.e., of $\scrM_0(k)$ (mimicking the condition of abelian $\pi_1$ with
trivial action on the higher $\pi_i$'s), except by expressing things
via the associated ordinary homotopy type (passing over to $X_*(k)$),
which looks kind of stupid indeed! But nothing is ``safe'' here
anyhow, and the subscript $0$ looks definitely more natural here than
subscript $1$, so we may as well keep it!

\bigbreak

\presectionfill\ondate{5.10.}\par

What we are mainly interested in, for expressing the relationships
between the non-additive category $\HotOf_0(k)$ and its additive
counterpart $\HotabOf_0(k)$, is an array of four functors in between
these, two of which being $\LtH_\bullet$ ((total) ``homology'', or
``linearization'') and $\widetilde K$ (Eilenberg-Mac~Lane functor) in
\eqref{eq:124.1} and \eqref{eq:124.2}, the two remaining ones being
$\Lpi_\bullet$ (``total homotopy'') and $\widetilde S$ (the
``spherical functor'', compare p.~\ref{p:514} \eqref{eq:122.9}),
fitting into the following diagram
\begin{equation}
  \label{eq:124.4}
  \begin{tabular}{@{}c@{}}
    \begin{tikzcd}[baseline=(O.base),sep=large]
      \HotabOf_0(k) \ar[r,shift left,"\widetilde K"]
      \ar[r,shift right,"\widetilde S"'] &
      \HotOf_0(k) \ar[l,bend right,"\LtH_\bullet"{swap}]
      \ar[l,bend left,"\Lpi_\bullet"{name=O}] \\
    \end{tikzcd}.
  \end{tabular}\tag{4}
\end{equation}
I would now like to discuss the main formal properties to be expected
from this set of functors.

\subsection[Adjunction properties]{Adjunction properties:}
\label{subsec:124.A}
\begin{equation}
  \label{eq:124.5}
  \left\{
    \begin{tabular}{@{}l@{}}
      a)\enspace $(\LH_\bullet,\widetilde K)$ is a pair of adjoint functors,
    \\
      b)\enspace $(\widetilde S, \Lpi_\bullet)$ is a pair of adjoint functors.
    \end{tabular}
  \right.\tag{5}
\end{equation}

Thus, we get adjunction homomorphisms (for $L_\bullet$ in
$\HotabOf_0(k)$ and $X$ in $\HotOf_0(k)$)
\begin{equation}
  \label{eq:124.6a}
  \begin{cases}
    X\to \widetilde K(\LtH_\bullet(X)) &
    \text{($\eqdef \widetilde W(X)$ (``pointed linearization''))} \\
    \LtH_\bullet(\widetilde K(L_\bullet)) \to L_\bullet &
    \text{,}
  \end{cases}\tag{6a}
\end{equation}
and\pspage{524}
\begin{equation}
  \label{eq:124.6b}
  \begin{cases}
    \widetilde S(\Lpi_\bullet(X)) \to X & \text{} \\
    L_\bullet \to \Lpi_\bullet(\widetilde S(L_\bullet)) & \text{.}
  \end{cases}\tag{6b}
\end{equation}

The adjunction formulas for the two pairs in \eqref{eq:124.5} are
\begin{equation}
  \label{eq:124.7}
  \left\{
    \begin{tabular}{@{}lll@{}}
      a) & $\Hom(\LtH_\bullet(X),L_\bullet) \simeq \Hom(X,\widetilde
           K(L_\bullet))$ &\\
      b) & $\Hom(L_\bullet, \Lpi_\bullet(X)) \simeq \Hom(\widetilde
           S(L_\bullet), X)$ & .
    \end{tabular}
  \right.\tag{7}
\end{equation}

The first-hand side of \hyperref[eq:124.7]{a)} should be viewed as
\emph{cohomology} of $X$ with coefficients in $L_\bullet$. More
specifically, as $\widetilde{\mathrm H}^0(X,L_\bullet)$, and in case
$L_\bullet = \pi[n]$, posing
\begin{equation}
  \label{eq:124.8a}
  K(\pi[n]) \eqdef K(\pi,n),\tag{8a}
\end{equation}
we get that the Eilenberg-Mac~Lane objects $K(\pi,n)$ represents the
cohomology functor
\[ X\mapsto \widetilde{\mathrm H}^0(X,\pi[n]) = \mathrm
H^n(X,\pi).\]

Symmetrically, the first-hand side of \hyperref[eq:124.7]{(7b)} should
be viewed as a ``mixed homotopy module'' of $X$ (relative to the
``co-coefficient'' $L_\bullet$), I am tempted to denote it as
\[\pi_0(X,L_\bullet),\quad\text{resp.\ $\pi_n(X,H)$ if
  $L_\bullet=H[n]$;}\]
thus, posing
\begin{equation}
  \label{eq:124.8b}
  \widetilde S(H[n]) \eqdef S(H,n) \quad\text{(sphere-like object for
    $H,n$),}\tag{8b}
\end{equation}
we get that the $H$-sphere $S(H,n)$ over $k$ represents the
(covariant) functor on $\HotOf_0(k)$
\[X\mapsto \pi_n(X,H).\]
In case $H=k$, we get
\[\pi_n(X,k) = \Hom(k[n],\Lpi_\bullet(X)) = \mathrm
H_n(\Lpi_\bullet(X)) \eqdef \pi_n(X),\]
and we get that the ``usual'' homotopy module functor
\[X\mapsto \pi_n(X)\]
is represented by the ``usual'' $n$-sphere (over $k$ however!)
$S(n,k)$, as it should. (This is of course the main justification why
we expect $(\widetilde S,\Lpi_\bullet)$ to be a pair of adjoint
functors, whereas adjunction for the pair $(\LtH_\bullet,\widetilde
K)$ is already fairly familiar from the set-up of ordinary homotopy
types.)

The adjunction formulæ \eqref{eq:124.7} show that the objects
$K(L_\bullet)$ in $\HotOf_0(k)$ are group objects, whereas the objects
$\widetilde S(L_\bullet)$ are co-group objects. This gives some inner
justification for calling $\widetilde W(X)=\widetilde
K(\LtH_\bullet(X))$ the\pspage{525} (pointed) ``linearization'' of the
schematic homotopy type $X$, which maps into the former by the first
adjunction morphism in \eqref{eq:124.6a}. Dually, we may call
$\widetilde S(\Lpi_\bullet(X))$ the (pointed) ``co-linearization'' of
$X$, it maps into $X$ by the first adjunction morphism in
\eqref{eq:124.6b}.

The source $\LtH_\bullet(\widetilde K(L_\bullet))$ in the second
adjunction map \eqref{eq:124.6a} may be viewed as Eilenberg-Mac~Lane
type (total) homology, corresponding to the chain complex
$L_\bullet=\pi[n]$ (but in the schematic sense of course!), whereas
the target $\Lpi_\bullet(\widetilde S(L_\bullet))$ in the second
adjunction map \eqref{eq:124.6b} may be viewed as ``total homotopy''
of a sphere-type space, reducing in case $L_\bullet=\pi[n]$ to the
total homotopy of the standard $n$-sphere over $k$, $S(n,\pi)$. In
case $k=\bZ$, we expect of course the homology groups of this chain
complex to be the ordinary homotopy groups of the usual sphere $S^n$
-- cf.~\ref{subsec:124.G} below.

\subsection[Inversion formulæ]{Inversion formulæ:}
\label{subsec:124.B}%
two functorial isomorphisms
\begin{equation}
  \label{eq:124.9}
  \left\{
    \begin{tabular}{@{}lll@{}}
      a) & $\Lpi_\bullet(\widetilde K(L_\bullet)) \tosim L_\bullet$ & , \\
      b) & $\LtH_\bullet(\widetilde S(L_\bullet)) \tosim L_\bullet$ & .
    \end{tabular}
  \right.\tag{9}
\end{equation}

Maybe we should have begun with these, as they kind of fix the meaning
of the two functors $\widetilde K, \widetilde S$ in terms of
$\LtH_\bullet, \Lpi_\bullet$, which may be viewed as embodying
respectively the two main sets of invariants of a homotopy type,
namely homology and homotopy.

\subsection[Hurewicz map]{Hurewicz map:}
\label{subsec:124.C}%
\begin{equation}
  \label{eq:124.10}
  \Lpi_\bullet(X) \to \LtH_\bullet(X).\tag{10}
\end{equation}

We get such a map by applying $\Lpi_\bullet$ to the linearization map
in \eqref{eq:124.6a}, and using \hyperref[eq:124.9]{(9a)};
symmetrically, we may apply $\LtH_\bullet$ to the colinearization map
\eqref{eq:124.6b}, and use \hyperref[eq:124.9]{(9b)}. We get a priori
two maps \eqref{eq:124.10}, and the statement is that these two maps
are the same. Moreover, we have the all-important \emph{Hurewicz
  theorem}: The first non-vanishing homology objects for
$\Lpi_\bullet(X)$ and $\LtH_\bullet(X)$ occur in the same dimension,
$n$ say, and \eqref{eq:124.10} induces an isomorphism for $\mathrm
H_n$, an epimorphism for $\mathrm H_{n+1}$.

\subsection[Exactness properties]{Exactness properties:}
\label{subsec:124.D}%
they can be stated shortly by saying that in each one of the two pairs
of adjoint functors \eqref{eq:124.5}, the left adjoint one respects
cofibration sequences, whereas the right adjoint respects fibration
sequences. This may be detailed as four exactness statements, one for
each one of the four basic functors.

Thus,\pspage{526} $\LtH_\bullet$ respects cofibration sequences, which
means essentially that for such a sequence in $\HotOf_0(k)$
\[Y\to X\to Z,\quad\text{$Z$ the homotopy cofiber of $Y\to X$,}\]
we get a corresponding exact triangle in $\HotabOf_0(k)$
\begin{equation}
  \label{eq:124.11b}
  \begin{tabular}{@{}c@{}}
    \begin{tikzcd}[baseline=(O.base),column sep=small]
      & \LtH_\bullet(Z) \ar[dl,"i"'] & \\
      \LtH_\bullet(Y) \ar[rr] & &
      |[alias=O]| \LtH_\bullet(X) \ar[ul]
    \end{tikzcd},
  \end{tabular}\tag{11b}
\end{equation}
which is the most complete and elegant way, I guess, for expressing
behaviour of homology with respect to homotopy cofibers and
suspensions. Dually, $\Lpi_\bullet$ respects fibration sequences,
which means that for such a sequence in $\HotOf_0(k)$
\[Z\to X\to Y,\quad\text{$Z$ the (connected) homotopy fiber of $X\to Y$,}\]
we get a corresponding exact triangle in $\HotabOf_0(k)$
\begin{equation}
  \label{eq:124.11a}
  \begin{tabular}{@{}c@{}}
    \begin{tikzcd}[baseline=(O.base),column sep=small]
      & \Lpi_\bullet(Z) \ar[dl,"i"'] & \\
      \Lpi_\bullet(X) \ar[rr] & &
      |[alias=O]| \Lpi_\bullet(Y) \ar[ul]
    \end{tikzcd},
  \end{tabular}\tag{11a}
\end{equation}
which is a more complete and elegant way of stating the exact homotopy
sequence of a fibration, in terms of total homotopy.

Formula \eqref{eq:124.11b} implies that $\LtH_\bullet$ commutes with
suspension functors $\Sigma$, the latter in $\HotabOf_0(k)$
(visualized as a derived chain complex category) is just the shift
functor
\[\Sigma\subab: L_\bullet\mapsto L_\bullet[1]\quad
(L[1]_n = L_{n-1},\]
which is a \emph{fully faithful} functor. Dually, formula
\eqref{eq:124.11a} implies that $\Lpi_\bullet$ commutes to the
loop-space functors $\Omega^0$; the latter in $\HotabOf_0(k)$
cannot be described as just a shift in opposite direction
\[L_\bullet \mapsto L_\bullet[-1],\]
as this will get us out of $0$-connected chain complexes, we have to
truncate, moreover, at dimension $1$ (afterwards, or at dimension $2$
beforehand). This functor is not fully faithful therefore -- we loose
something when passing from $L_\bullet$ to
$\Omega\subab^0(L_\bullet)$.

As for the two functors $\widetilde K,\widetilde S$ in the opposite
direction, from $\HotabOf_0(k)$ (i.e., essentially chain complexes) to
$\HotOf_0(k)$, that the first one respects fibration sequences is
surely quite familiar a fact (and kind of tautological) in the set-up
of ordinary homotopy types. That the second, less familiar functor
$\widetilde S$ should respect cofibration sequences\pspage{527} should
be useful in order to give a more or less explicit construction of
$\widetilde S(L_\bullet)$, for given $L_\bullet$, in terms of the
``spheres over $k$'' $S(H,n)$.

\begin{remarks}
  1)\enspace The superscript $0$ for the loop functor $\Omega^0$ or
  $\Omega\subab^0$ should not be confused with the subscript $n$ in
  the iterated loop-space functor $\Omega_n$, it is added here to
  suggest that we are taking the \emph{neutral} connected component of
  the ``true'' full loop-space, this being imposed by the restriction
  of working throughout with $0$-connected objects; likewise, the
  ``homotopy fiber'' operation in the present context should be viewed
  as meaning ``connected component at the marked point of the full
  homotopy fiber''. These necessary readjustments of the usual notions
  is being felt as an unwelcome feature (of which I have become aware
  only at this very moment, I confess, through the writing of the
  notes). Intuitively, the restriction to $0$-connected (pointed)
  homotopy types doesn't seem to imperative, technically, however,
  when working with semisimplicial schematic models, we had felt like
  introducing the condition $X_0=e$, which may be viewed as a strong
  form of a $0$-connectedness condition. Truth to tell, it isn't so
  clear that this condition is going to be of great utility -- anyhow,
  it isn't enough to ensure what we're really after at present (namely
  abelian $\pi_1$ and trivial action on the higher $\pi_i$'s), a
  condition which is anyhow independent of any $0$-connectedness type
  assumption, and is moreover (it would seem) stable under the basic
  fiber and cofiber operations. To sum up, it may well turn out that
  we better replace the categories of $0$-connected homotopy types
  $\HotabOf_0(k)$ and $\HotOf_0(k)$, by slightly larger ones, so as to
  get rid of the $0$-connectedness restriction. This, however, is at
  present a relatively minor point, and therefore we'll leave the
  notations as they are.

  2)\enspace Behaviour of $\LtH_\bullet$ with respect to
  \emph{fibration} sequences (instead of cofibration sequences) is a
  relatively delicate matter, it is governed by the Leray spectral
  sequence, whose initial term is the homology of the base $Y$ with
  ``coefficients'' in the homology of the fiber. I wonder if there is
  anything similar for the behavior of $\Lpi_\bullet$ with respect to
  cofibration sequences?

  3)\enspace In the display of the main expected properties of the
  four basic functors, there is a striking \emph{symmetry}, which we
  tried to stress by the way of presented the main formulæ. One way to
  express this symmetry is to say that things look as if there was an
  auto-duality in the pair of categories $(\HotabOf_0(k),
  \HotOf_0(k))$, namely a pair of contravariant involutive functors
  $(D\subab,D)$, each of which interchanges\pspage{528} fibration and
  cofibration sequences, i.e., transforms fibers into cofibers and
  vice versa, and the pair of functors interchanging $\LtH_\bullet$
  and $\Lpi_\bullet$ on the one hand, $\widetilde K$ and $\widetilde
  S$ on the other. This heuristic formulation is compatible with all
  we have stated so far, except one thing -- namely, the suspension
  functor $\Sigma\subab$ in $\HotabOf_0(k)$ is fully faithful, whereas
  the (supposedly ``dual'') loop-space functor $\Omega\subab^0$ is
  not. Thus, it will be more accurate to say that, if we view the
  formulaire developed so far as being the description of a certain
  structure type, whose basic ingredients are two categories
  $\scrH\subab$ and \scrH, endowed with fiber and cofiber operations,
  tied by four functors as above satisfying a bunch of properties,
  then the axioms are autodual in an obvious sense; namely if they are
  satisfied for a system $(\scrH\subab,\scrH,\LtH_\bullet,\widetilde
  K,\widetilde S,\Lpi_\bullet)$, they are equally satisfied by the
  system $(\scrH\subab\op,\scrH\op,\Lpi_\bullet,\widetilde
  S,\widetilde K,\LtH_\bullet)$. (Of course, among the properties of
  the structure type, we are \emph{not} going to include that the
  suspension functor in $\scrH\subab$ is fully faithful!)\enspace As
  already stated before, the main reason for introducing the fourth
  functor $\widetilde S$ in the picture, was because it was felt that
  this was lacking in order to round it up. Thus for any kind of
  notion or statement in this set-up, suggested by some kind of
  geometric insight, it becomes automatic to look at its dual and see
  whether it makes sense.

  As for formal autoduality, we already noticed there is none, even in
  $\HotabOf_0(k)$ just by itself. Thinking of this category as
  $\D_\bullet(\AbOf_k)$, namely as a full subcategory of the derived
  category $\D(\AbOf_k)$, we cannot help, though, but thinking of the
  standard ``dualizing'' functor
  \[L_\bullet \mapsto \RHom(L_\bullet,k): \D^-(\AbOf_k)\to
  \D^+(\AbOf_k),\]
  inducing a perfect duality within the category of ``perfect''
  objects in $\D(\AbOf_k)$, namely objects which are isomorphic to
  those which may be described by complexes in $\AbOf_k$ having a
  bounded span of degrees, and all of whose components are projective
  of finite type. But this autoduality of course transforms
  \emph{chain} complexes into \emph{cochain} complexes -- we have to
  shift these in order to get chain complexes again. This suggests
  that maybe we'll hit upon an actual autoduality, provided we go over
  from $\HotOf_0(k)$ to a suitable ``stabilized'' category, the
  suspension category say (deduced from the initial one by introducing
  formally a quasi-inverse for the suspension functor). The homology
  functor $\LtH_\bullet$ extends to a functor from the suspension
  category to the corresponding category for $\HotabOf_0(k)$, say
  $\D^-(\AbOf_k)$. The $\Lpi_\bullet$\pspage{529} functor, though, is
  lost on our way -- I am afraid I am confusing the kind of duality I
  am after here, with a rather different type of duality, kin to
  Poincaré duality, and discovered I believe by J.~H.~C.~Whitehead, in
  the context of spaces having the homotopy type of a finite
  complex. For such a complex, the idea (as far as I remember) is to
  embed $X$ into a large-dimensional sphere, and to take the
  complement $X\upvee$ of an open tubular neighborhood of $X$. Up to
  suspension, the homotopy type of $X\upvee$ does not depend on the
  choices made and (shifting back by $n$) we thus get a canonical
  object in the suspension category, depending contravariantly on
  $X$. The functor $X\mapsto X\upvee$ (if I got it right) is an
  autoduality of the relevant full subcategory of the suspension
  category, and the (reduced) homology functor $\LtH_\bullet$ maps
  this subcategory into perfect complexes, and commutes to
  autodualities. The whole story seems tailored towards a study of
  duality relations for the functor $\LtH_\bullet$ exclusively --
  without any reference to the homotopy invariants $\pi_i$. I don't
  know if one can devise a similar story for $\Lpi_\bullet$, by
  stabilizing with respect to the loop space functor (or is this just
  nonsense?). At any rate, there doesn't seem any autoduality in view,
  exchanging homology and homotopy invariants\dots
\end{remarks}

\subsection[Conservativity properties]{Conservativity properties:}
\label{subsec:124.E}%
The functors $\LtH_\bullet$ and $\Lpi_\bullet$ are both conservative,
i.e., a map in $\HotOf_0(k)$ which by either of these functors becomes
an isomorphism, is an isomorphism. (Here of course the
$0$-connectedness assumption for the homotopy types we are working
with is essential, as far as the functor $\Lpi_\bullet$ is concerned
at any rate.)

\subsection[Base change properties]{Base change properties:}
\label{subsec:124.F}%
All four functors, and the adjunction and inversion maps
\eqref{eq:124.7} and \eqref{eq:124.9}, are compatible with ring
extension $k\to k'$, it being understood that such a ring homomorphism
defines functors
\begin{equation}
  \label{eq:124.12}
  \HotabOf_0(k) \to \HotabOf_0(k'), \quad
  \HotOf_0(k) \to \HotOf_0(k')\tag{12}
\end{equation}
compatible with the fibration and cofibration structures.

In opposite direction, there should be, too, a ``restriction'' functor
(less important, though, I feel), I didn't try to find out what its
formal properties should be with respect to the formalism developed
here.

\subsection[Comparison with ordinary homotopy types]{Comparison with
  ordinary homotopy types:}
\label{subsec:124.G}%
We got functors
\begin{equation}
  \label{eq:124.13}
  \HotabOf_0(k) \to \HotabOf = \D_\bullet\Ab, \quad
  \HotOf_0(k) \to \HotOf_0,\tag{13}
\end{equation}
(the first being interpreted as ``forget $k$'' functor, the second as
a sections functor on semisimplicial schematic models). This pair
of\pspage{530} functors is compatible with the functors $\widetilde
K$, but (even if $k=\bZ$) definitely \emph{not} with their left
adjoints $\LtH_\bullet$, namely with homology. Despite this fact, we
hope if $k=\bZ$ that the functors \eqref{eq:124.13} are compatible
with the functors $\widetilde S$ (assuming that $\widetilde S$ can
actually be constructed in the discrete set-up too), so that spheres
are transformed into spheres. Compatibility with $\Lpi_\bullet$
doesn't have a meaning strictly speaking, as this functor is not
defined on ordinary homotopy types, however, the functors $\pi_i$ are
defined in both contexts, and the sections functor should commute to
these (for arbitrary $k$). Thus, the only serious incompatibility
trouble concerns the homology functor $\LH_\bullet$ and should not
arise, however, for sphere-like objects $S(H,n)$. I forgot to state
from the outset that the two functors \eqref{eq:124.13} are expected
to respect fibration and cofibration sequence structures, of course.

\bigbreak

\presectionfill\ondate{13.10.}\par

\hangsection[Digression on Baues' cofibration and fibration
\dots]{Digression on Baues' cofibration and fibration categories, and
  on weak equivalences alone as the basic data for a model
  category.}\label{sec:125}%
During the last eight days, I have been busy with a number of things,
which left little leisure for mathematical ponderings. I got a letter
from H.~J.~Baues, who had seen my notes, which induced him to send me
his preprint ``On\scrcomment{I'm not sure what this refers to; perhaps
  \textcite{Baues1985}\dots} the homotopy classification problem'' (Chapters~I
to~V + chapter~Ext). I spent two evenings looking through parts of
these, where Baues carries as far as possible (namely quite far
indeed) the homotopy formalism in the context of his so-called
``cofibration (or fibration) categories'', using as his leading
thread his ideas on ``obstruction theory''. He makes the point that he
tried by his basic notion to pinpoint the weakest axiomatic set-up,
sufficient however for developing all the major familiar (and even
some not at all familiar!) notions, operations and statements of usual
homotopy theory. In his letter, he suggested that maybe in any
``universe'' where homotopy constructions make sense, one or the other
of his two mutually dual set-ups should be around. Such suggestion was
of course quite interesting for my present reflections, as I do have
the hope indeed that there exists a ``universe'' of schematic homotopy
types, which may be described in terms of the models (namely
semisimplicial unipotent bundles) I have been using so far, or at any
rate by closely related kinds of models. More specifically, I do hope
for the two kinds of operations to make sense, namely ``integration''
(including homotopy cofibers and suspensions) and ``cointegration''
(including homotopy fibers and loop objects), which in Baues' set-up
should correspond on the category of models to a structure of a
cofibration category and a fibration category, respectively. As the
kind of models I have been working with\pspage{531} don't allow for an
amalgamated sum (or ``push-out'') construction, except in the trivial
case when one of the two arrows to be amalgamated is an isomorphism,
it is clear that we cannot hope to get with these models a cofibration
category. There may be, however, a fibration category structure, and I
started playing around a little with a possible notion of Kan
fibration -- without coming to a definite conclusion yet, however. I
feel I shouldn't dwell much longer still, for the time being, on
getting off the ground a theory of schematic homotopy types, maybe
I'll come back upon it later, after next chapter, when (hopefully)
I'll be a little more ``in the game'' of Kan type conditions, closed
model category structures and Baues' ``halved'' variants of these.

At any rate, whether or not a homotopy cofiber construction can be
carried through for schematic homotopy types, it seems rather clear to
me that Baues' suggestion or expectation, about the set-up of
cofibration and fibration categories he ended up with, is not quite
justified. Already by the time (in the late sixties) when I first
heard from Quillen about his approach (and the same applied to Baues'
ones), when applied say to semisimplicial objects in some category $A$
as ``models'', is pretty strongly relying on the existence of
``enough'' projectives (or dually, of enough injectives, if working
with cosemisimplicial objects) in $A$. When $A$ is a topos say, then
the category $\bHom(\Simplexop,A)$ of semisimplicial objects of $A$
doesn't carry, it seems (except in very special cases, say of a
totally disconnected topos), a reasonable structure of a Quillen model
category, nor even (I would think) of a cofibration or a fibration
category in the sense of Baues; however, I am pretty sure that the
derived category (with respect to the notion of weak equivalence
introduced by Illusie\scrcomment{AG wanted to insert a reference, but
  didn't, probably to \textcite{Illusie1971,Illusie1972}\dots}) has
important geometric meaning and is indeed a ``universe for homotopy
types'' -- and the first steps in developing such theory have been
taken by Illusie already in Chapter~I of his thesis. Quite similarly,
if $A$ is an abelian category, Verdier's derived category $\D(A)$,
deduced from the category $\mathrm C(A)$ of complexes in $A$ by
localizing with respect to the set $W$ of quasi-isomorphisms, does
allow for a homotopy formalism (as does any ``triangulated category''
in the sense of Verdier), however, it doesn't seem that this formalism
may be deduced from a structure of a cofibration category or a
fibration category on $\mathrm C(A)$. The same I guess holds for the
subcategories $\mathrm C^-(A)$ and $\mathrm C^+(A)$, giving rise to
$\D^-(A)$ and $\D^+(A)$, which are more important still than $\D(A)$
in the everyday cohomology formalism of ``spaces'' of all kinds
(namely, essentially, of ringed topoi). When\pspage{532} $A$ admits
enough projectives or enough injectives, respectively, these
categories (as pointed out by Quillen) are associated to closed model
structures on $\mathrm C^-(A)$ and $\mathrm C^+(A)$, respectively, but
(it seems) not otherwise.

It strikes me that nobody apparently so far has tried to develop
homotopy theory, starting as basic data with a category of models $M$
together with a set $W$ of weak equivalences in $M$, satisfying
suitable assumptions, without giving moreover a notion of
``fibration'' or of ``cofibrations'' in $M$. Various examples,
including the all-important case (from my point of view) of \Cat,
suggest that the choice of the notion of fibration or cofibration
isn't really so imperative, that it is to a certain extent arbitrary,
different choices (compatible with the basic data $W$) leading to
essentially ``the same'' homotopy theory. This is due to the fact that
they lead to the same notions of ``integration'' and ``cointegration''
of homotopy types, which depends indeed only on $W$ and which, in my
eyes, are \emph{the} two main operations of homotopy theory (compare
section~\ref{sec:69}). They seem to me the key for defining such a
thing as a specific ``homotopy theory'', independently of any
particular choice of a category of models (+ extra structure on it,
and notably, weak equivalences) used for describing it. The precise
technical notion achieving this is of course the notion of a
``derivator'' -- and I do hope that it shouldn't be too awfully hard
to construct, for instance, the ``derivator of schematic homotopy
types'', and maybe even characterize it axiomatically (as well as the
derivator of usual homotopy types), without having to make any mention
of models for such characterization. One may say that, after the one
major step in the foundations of homological algebra, consisting in
introducing the derived category of an abelian category (and
systematically working with derived categories for stating the main
facts about cohomology of all kind of ``spaces'', namely topoi, such
as usual topological spaces, schemes and the like\dots), the work on
foundations more or less stopped short, while the next step to take
was to come to a grasp of the full structure involved in derived
categories, namely the structure of a derivator. And it turns out that
as well the step of taking derived categories (namely localizations of
suitable model categories), as the next, namely introducing the
corresponding structure of a derivator, make good sense also in
non-abelian set-ups, namely for doing ``homotopy theory'', thus viewed
as a non-commutative version of homology theory.\pspage{533}

\hangsection[Dissymmetry of $\LtH_*$ versus $\Lpi_*$ (no
filtration dual to \dots]{Dissymmetry of
  \texorpdfstring{$\LtH_*$}{LH} versus
  \texorpdfstring{$\Lpi_*$}{Lpi} (no filtration dual to Postnikov
  dévissage).}\label{sec:126}%
Let's come back, though, to schematic homotopy types. Last thing we
looked at was the four basic functors between ``schematic homotopy
types'' and ``linear schematic homotopy types'' over a given ground
ring $k$, making up the two basic categories
\[\HotabOf_0(k)\quad\text{and}\quad
\HotOf_0(k).\]
We have been dwelling somewhat on the remarkable formal symmetry to be
expected for these relations. It is tempting, then, to try and dualize
any kind of basic notion or construction which makes sense in terms of
the basic data, namely the four functors and the adjunction and
inversion maps relating them. Maybe the very first thing which forces
attention is the Postnikov dévissage of an object of $\HotOf_0(k)$,
which had been (together with abelianization) the very starting point
for our approach towards schematic homotopy types. One main ingredient
of this dévissage is the ``tower'' of the Cartan-Serre type functors
\begin{equation}
  \label{eq:126.1}
  X\mapsto X_n\tag{1}
\end{equation}
and maps
\begin{equation}
  \label{eq:126.2}
  X\to X_n\tag{2}
\end{equation}
(with $n$ a natural integer), where $X_n$ is deduced from $X$ by
``killing its $\pi_i$'s for $i>n$''. We may call an object $Y$ of
$\HotOf_0(k)$ \emph{$n$-co-connected} (a notion in a way symmetric to
$n$-connectedness) if $\pi_i(X)=0$ for $i>n$, and denote by
\begin{equation}
  \label{eq:126.3}
  \HotOf_0(k)_{/n}\tag{3}
\end{equation}
the full subcategory of $\HotOf_0(k)$ made up with these objects,
which is therefore the inverse image by the functor
\[\Lpi_i : \HotOf_0(k)\to\HotabOf_0(k)\]
of the corresponding full subcategory
\begin{equation}
  \label{eq:126.4}
  \HotabOf_0(k)_{\le n}\tag{4}
\end{equation}
of $\HotabOf_0(k)$. The most natural way for defining $X_n$ in terms
of $X$ (in an axiomatic set-up with the four functors as basic data),
together with the canonical ``fibration'' \eqref{eq:126.2}, is by
describing \eqref{eq:126.2} as the ``universal'' map (in
$\HotOf_0(k)$) of $X$ into an $n$-co-connected object. In other words,
we are surmising that the inclusion functor
\[\HotOf_0(k)_{\le n} \hookrightarrow \HotOf_0(k)\]
admits a left adjoint, and the latter is denoted by $X\mapsto
X_n$. This\pspage{534} description gives rise at once to the familiar
``tower'' structure for variable $n$
\begin{equation}
  \label{eq:126.5}
  X_{n+1} \to X_n \to \dots \to X_1 \to X_0 (=e).\tag{5}
\end{equation}
Intuitively, the maps in \eqref{eq:126.5} are viewed as being
(surjective) fibrations between (connected) ``spaces'', $X$ being
viewed as a kind of inverse limit of the $X_n$'s. Dually, we would
expect to get a sequence of inclusions
\begin{equation}
  \label{eq:126.star}
  {}_0X \hookrightarrow {}_1X \hookrightarrow \cdots \hookrightarrow
  {}_nX \hookrightarrow {}_{n+1}X \hookrightarrow \cdots,\tag{*}
\end{equation}
with $X$ appearing as a kind of direct limit. In the set-up of
ordinary homotopy types, modelized by semisimplicial sets, one will
think at once of the filtration by skeleta -- which however isn't
quite the right thing surely, because if we dualize the familiar
Cartan-Serre requirement on \eqref{eq:126.2} (namely that it induces
an isomorphism for $\pi_i$ for $i\le n$), we see that the
``inclusion''
\[{}_nX \hookrightarrow X\]
should induce an isomorphism for $\mathrm H_i$ for $i\le n$, which
isn't quite true for the skeletal filtration (it is OK for $i\le n-1$
only); the condition
\[\mathrm H_i({}_nX) = 0\quad\text{for $i>n$}\]
is OK, though, for this filtration. So the next idea would be to
modify a little the $n$-skeleton ${}_nX$, to straighten things out. I
played around some along these lines, and after some initial optimism,
came to the feeling that there does not exist (in the discrete nor in
the schematic set-up) such an increasing canonical filtration of a
homotopy type. I didn't make any formal statement and proof for this
(in the set-up of ordinary homotopy types, say), however, in the
process of playing around in became soon clear that in various ways,
\emph{there is some essential dissymmetry} between the seemingly
``dual'' situations, when trying to get the two types of
``filtrations'' of the object $X$. One dissymmetry occurs already in
the very definition of the subcategory \eqref{eq:126.3}, and of the
corresponding ``dual'' subcategory
\begin{equation}
  \label{eq:126.6}
  {}_{(\le n)}\HotOf_0(k)\tag{6}
\end{equation}
of objects satisfying
\[\mathrm H_i(Y) = 0\quad\text{for $i>n$.}\]
Namely, both subcategories \eqref{eq:126.3} and \eqref{eq:126.6} are
defined in terms of the \emph{same} subcategory of $\HotabOf_0(k)$,
namely \eqref{eq:126.4}, as the inverse image of the\pspage{535}
latter by either $\Lpi_\bullet$, or $\LH_\bullet$. Now, the point is
that the properties of this subcategory, with respect to the inner
structure of $\HotabOf_0(k)$ (involving the ``left'' and the ``right''
homotopy operations, notably the suspension and the loop functors,
respectively) are by no means autosymmetric. The main dissymmetry, it
would seem, turns out in this, that \eqref{eq:126.4} is stable under
the loop functor, and by no means under the suspension functor. This
is the reason why, even if we assume that (in analogy to what happens
for the categories \eqref{eq:126.3}) the inclusion functors from the
categories \eqref{eq:126.5} into $\HotOf_0(k)$ do admit the relevant
(namely right) adjoints (which I greatly doubt anyhow\dots), and thus
give rise for any object $X$ to an increasing filtration
\eqref{eq:126.star}, the Postnikov-type relationships between $X_n$
and $X_{n+1}$ cannot be quite dualized to a similar relationship
between ${}_nX$ and ${}_{n+1}X$. We may think next, of course, of
defining an increasing sequence of subcategories \eqref{eq:126.5} of
$\HotOf_0(k)$ in terms of $\LH_\bullet$ and a corresponding sequence
of subcategories of $\HotabOf_0(k)$, different from the categories
\eqref{eq:126.4}, and stable under suspensions. But there doesn't seem
to be anything reasonable around along these lines.

To sum up, it doesn't seem one should overemphasize the somewhat
startling symmetry which appeared in section~\ref{sec:124} between
``homotopy'' embodied in $\Lpi_\bullet$, and ``homology'' embodied in
$\LH_\bullet$ -- in some essential respects, it would seem that the
corresponding two functors do have non-mutually symmetry properties. I
guess I have to apologize for having taken that long for coming to a
conclusion which, presumably, must be felt as a kind of self-evidence
by all homotopy people!

\hangsection[Schematic homotopy types and Illusie's derived category
\dots]{Schematic homotopy types and Illusie's derived category for a
  topos.}\label{sec:127}%
Before leaving (for the time being) the topic of schematic homotopy
types and schematization, I would like still to add a few comments,
about various possibilities for working with different kinds of models
for defining schematic homotopy types. My point here is not in
replacing the basic test-category we are working with, here
$\Simplex$, by some other (say the category $\Globe$ of standard
hemispheres) -- this choice I feel should be more or less irrelevant,
so we may as well keep $\Simplex$. Thus, we are going to work with
semisimplicial objects, and the main question then is (for a given
ground ring $k$) to say precisely what kind of objects we are allowing
(or imposing!) as components for our complexes. At any rate, they
should be ``objects over $k$'', and the most\pspage{536} encompassing
choice for such objects seems to be sheaves on the category of all
schemes over $k$ (or equivalently, of all affine schemes over $k$),
for a suitable topology such as the fpqc topology (compare
section~\ref{sec:111}, pages \ref{p:446}--\ref{p:447}). Apart from the
technico-logical nuisance of this not being a \scrU-site, where
\scrU{} is our basic universe (cf.\ p.~\ref{p:492}), which we'll
ignore here (as it isn't really too serious a difficulty), when
working the semisimplicial sheaves, the embryo of foundational work of
Illusie's on the derived category of a topos via semisimplicial
sheaves\scrcomment{\textcite{Illusie1971}, chapter 1} becomes
available. Thus, for a map of semisimplicial sheaves over $k$,
\begin{equation}
  \label{eq:127.1}
  X_*\to Y_*,\tag{1}
\end{equation}
we know already what it means to be a quasi-isomorphism (= weak
equivalence). It seems likely that, even when drastically restricting
the sheaves allowed as components for our semisimplicial models, the
notion of quasi-isomorphism relevant for these models should be the
same as Illusie's for the encompassing topos. This is one point we
have been neglecting so far. I am not going here to recall Illusie's
definition (in terms of the homotopy sheaves
$\boldsymbol\pi_i(X_*,s)$\scrcomment{in the typescript, the ``$\pi$''
  is underlined\dots} where $s$ is a section of $X_*$), or translate
it into cohomological terms. One point I want to make here, is that it
is by no means automatic that if \eqref{eq:127.1} is a weak
equivalence, the same holds for the corresponding map of
semisimplicial sets
\begin{equation}
  \label{eq:127.2}
  X_*(k) \to Y_*(k).\tag{2}
\end{equation}
Thus, there is no more-or-less tautological ``sections'' functor from
Illusie's derived category to the category \Hot{} of usual homotopy
types. (There surely \emph{is} a canonical functor, though, in the
case of a general topos $T$ and Illusie's corresponding derived
category, namely what we would like to call the \emph{``cointegration
  over $T$'' functor}, which should come out of the formalism of
stacks we haven't begun to develop yet.)\enspace
Presumably, however, when working with semisimplicial models whose
components are restricted to be unipotent or something pretty close to
these (see below for examples), whenever \eqref{eq:127.1} is a weak
equivalence, the same will hold for \eqref{eq:127.2}. As the choice of
objects (let's call them simply the ``bundles'') we are allowing
should clearly be stable under ring extension, it will then follow
that more generally, for any algebra $k'$ over $k$, the corresponding
map
\begin{equation}
  \label{eq:127.2prime}
  X_*(k')\to Y_*(k')\tag{2'}
\end{equation}
is\pspage{537} a weak equivalence, too. Thus, a ``schematic homotopy
type over $k$'' should define a functor $\Alg_{/k} \to \HotOf$.

For any choice of a ``section''
\begin{equation}
  \label{eq:127.3}
  s\in X_0(k)\tag{3}
\end{equation}
of $X_*$, Illusie's constructions yield sheaves
\begin{equation}
  \label{eq:127.4}
  \bpi_i(X_*,s)\quad
  (i\ge0),\tag{4}
\end{equation}
where \bpiz{} is a sheaf of sets, $\bpi_1$ a sheaf of groups, acting
on the higher $\bpi_i$'s which are abelian sheaves. We will be mainly
interested of course in the case when \bpiz{} is the final sheaf
(we'll say that $X_*$ is \emph{relatively $0$-connected over $k$}),
and moreover $\bpi_1$ is abelian and its action on the higher $\bpi_i$
is trivial (let's say in this case that the relative homotopy type
defined by the semisimplicial sheaf $X_*$ is
\emph{pseudo-abelian}). In this case, up to canonical isomorphisms,
the abelian sheaves $\bpi_i(X_*,s)$ do not depend on the choice of
$s$, and (provided the pseudoabelian notion is defined locally) they
make sense, independently even of the existence of a section $s$. Our
hope is, by suitably restricting our notion of ``bundle'', to get on
the abelian sheaves $\bpi_i$ for $i\ge2$, or even on all $\bpi_i$
($i\ge1$), a natural structure of an $\scrO_k$-module, and one
moreover which in many ``good'' cases turns it into a quasi-coherent
$\scrO_k$-module. Also, the natural maps
\begin{equation}
  \label{eq:127.5}
  \pi_i(X_*(k)) \to \bpi_i(X_*)(k)\tag{5}
\end{equation}
should be isomorphisms, i.e., taking $\pi_i$ and $\bpi_i$ should
commute to taking sections (which will imply indeed that if
\eqref{eq:127.1} is a quasi-isomorphism, then so is
\eqref{eq:127.2}). Thus, a first basic question here is to find a
suitable notion of a ``bundle'' over $k$, in such a way that for
semisimplicial bundles satisfying some mild extra assumption (such as
$X_0=X_1=e$) implying that $X_*$ is pseudoabelian, the maps
\eqref{eq:127.5} should be isomorphisms. I expect this to be true for
the notion we have been working with so far, namely for ``unipotent
bundles'' (cf.\ section~\ref{sec:111}), but I haven't made any attempt
yet to prove this. But even granting this property for a given notion
of ``bundles'', the module structure on the $\bpi_i$'s at this point
remains still a mystery, as long as we don't tie them in with the Lie
functor (compare section~\ref{sec:118})\dots

\begin{remark}
  Maybe, in the set-up of a general topos $T$ and taking sections on
  the latter, the maps \eqref{eq:127.5} are isomorphisms, whenever
  the\pspage{538} homotopy sheaves $\bpi_i$ satisfy the relations
  \begin{equation}
    \label{eq:127.6}
    \mathrm H^i(T, \bpi_j(X_*))=0 \quad
    \text{for $i>0$,}\tag{6}
  \end{equation}
  as a consequence, maybe, of some spectral sequence whose abutment is
  the graded homotopy of the cointegration of $X_*$ over $T$. If so,
  then \eqref{eq:127.5} are quasi-isomorphisms whenever the sheaves
  $\bpi_i$ can be endowed with the structure of a quasi-coherent
  $\scrO_k$-module.
\end{remark}

\hangsection[Looking for the right notion of ``bundles''; $V$-bundles
\dots]{Looking for the right notion of ``bundles'';
  \texorpdfstring{$V$}{V}-bundles versus
  \texorpdfstring{$W$}{W}-bundles.}\label{sec:128}%
Next requirement about the notion of ``bundles'' we are going to work
with, is about existence of a ``linearization functor''; associating
to every bundle $X$ another bundle, namely its ``linearization''
$L(X)$, endowed moreover with the structure of an
$\scrO_k$-module. \emph{We insist that for a given bundle $X$, $L(X)$
  should be at any rate a bundle too}. More specifically, if $U(k)$
denotes the category of all bundles over $k$, we should also introduce
a corresponding ``$k$-linear'' category, or more accurately an
$\scrO_k$-linear category, $U\subab(k)$, whose objects should be
objects of $U(k)$ endowed with the extra structure of an
$\scrO_k$-module, possibly subject to some restrictions -- thus, we'll
get a forgetful functor
\begin{equation}
  \label{eq:128.7}
  K: U\subab(k)\to U(k).\tag{7}
\end{equation}
When we take $U(k)$ to be unipotent bundles in the sense of
section~\ref{sec:111}, the evident choice for $U\subab(k)$ is to take
the category of quasicoherent $\scrO_k$-modules, equivalent to the
category of $k$-modules, $\AbOf_k$. Dually, we may take $U(k)$ to be
the category of sheaves over $k$ isomorphic to the underlying sheaf of
sets of a vector bundle $V(M)$ over $k$ associated to a $k$-module
$M$, by the requirement
\begin{equation}
  \label{eq:128.8}
  V(M)(k') = \Hom_\kMod(M,k').\tag{8}
\end{equation}
The evident corresponding choice for $U\subab(k)$ is to take the
category of all vector bundles over $k$, which is equivalent to the
category $\AbOf_k\op$ \emph{opposite} to the category of $k$-modules,
as the functor $M\mapsto V(M)$ is contravariant in $M$. Maybe we
should distinguish between these two choices of bundles by different
notations, namely
\begin{equation}
  \label{eq:128.9}
  \text{$U_W(k)$ and $U_V(k)$,}\tag{9}
\end{equation}
where the subscripts $W$ and $V$ are meant to suggest the standard
descriptions of objects via (the underlying sheaves of sets of the
$\scrO_k$-modules)
\[\text{$W(M)$ and $V(M)$}\]
respectively, where (we recall)
\begin{equation}
  \label{eq:128.10}
  W(M)(k')=M\otimes_k k'.\tag{10}
\end{equation}

Reverting\pspage{539} to a general notion of ``bundles'' $U(k)$, we
define now a linearization functor
\begin{equation}
  \label{eq:128.11}
  L: U(k)\to U\subab(k)\tag{11}
\end{equation}
to be a functor left adjoint to the forgetful functor
\eqref{eq:128.7}, i.e., giving rise to an adjunction isomorphism
\begin{equation}
  \label{eq:128.12}
  \Hom_{U(k)}(X,K(L)) \simeq \Hom_{U\subab(k)}(L(X),L),\tag{12}
\end{equation}
where $X$ is in $U(k)$, $L$ in $U\subab(k)$. Passing over to
semisimplicial objects and the corresponding derived categories, the
functors $K$ and $L$ in \eqref{eq:128.7} and \eqref{eq:128.11} should
give rise to the functors $\widetilde K$ and $\LH_\bullet$ of
section~\ref{sec:124}.

When we take $U(k)=U_W(k)$, there is a drawback, though, because (as
we saw in section~\ref{sec:115}) the functor $K$ \emph{does not} admit
a left adjoint, only a \emph{proadjoint}, associating to an object $X$
a \emph{proobject} $L(X)$ of $U\subab(k)$. As pointed out in
section~\ref{sec:124} (p.~\ref{p:522}), if we want a nice pair $(K,L)$
of adjoint functors, this forces us to work with pro-unipotent bundles
instead of just unipotent ones, thus getting out of the haven of
sheaves over $k$ and into the somewhat dubious sea of prosheaves and
semisimplicial prosheaves, which have not been provided for in
Illusie's foundational ponderings! If we do stick to the $W$-approach,
this promises us a fair amount of extra sweat, putting in proobjects
everywhere, not too enticing a prospect, is it?

It would seem that we are better off with the $V$-approach, in which
case the bundles (more accurately, $V$-bundles) we are working with
are actual schemes, indeed affine schemes over $k$, as we get
\begin{equation}
  \label{eq:128.13}
  V(M)\simeq \Spec(\Sym_k(M)).\tag{13}
\end{equation}
Now, let more generally $X$ be any scheme, and let's look at maps from
$X$ into any vector bundle $V(M)$, we get
\begin{equation}
  \label{eq:128.star}
  \begin{split}
    \Hom(X,V(M)) &\eqdef
    \Hom_{\scrO_X}(p^*(\widetilde M), \scrO_X) \\
    &\qquad\qquad\overset{\textup{adjunction}}{\simeq}
    \Hom_{\scrO_S}(\widetilde M, p_*(\scrO_X)),     
  \end{split}\tag{*} 
\end{equation}
where
\[p:X\to S=\Spec k\]
is the structural map of $X$, and $\widetilde M$ the restriction of
$W(M)$ to the usual Zariski site of $S=\Spec(k)$. It is well-known
that under a rather mild restriction on $X$ (namely $X$ quasi-compact
and quasi-separated), always satisfied when $X$ is affine, $p_*$ takes
quasi-coherent sheaves on $X$ (for the usual small Zariski site) into
quasi-coherent sheaves on $S=\Spec k$; thus, if\pspage{540} $A$ is the
$k$-module (a $k$-algebra as a matter of fact) such that
\[\widetilde A\simeq p_*(\scrO_X),\quad
\text{i.e., $A=\Gamma(S,p_*(\scrO_X))\simeq\Gamma(X,\scrO_X)$,}\]
the last member of \eqref{eq:128.star} is
\[\Hom_{\scrO_X}(\widetilde M,\widetilde A) \simeq \Hom_k(M,A) \simeq
\Hom_{\scrO_k}(V(A),V(M)),\]
and finally we get
\begin{equation}
  \label{eq:128.14}
  \Hom(X,V(M)) \simeq \Hom_{\scrO_k}(V(A),V(M)).\tag{14}
\end{equation}
This shows that the forgetful functor from the category ${U\subab}_V$
of all vector bundles over $k$ to the category of all $k$-schemes
which are quasi-compact and quasi-separated, admits a left adjoint
$L_V$, where
\begin{equation}
  \label{eq:128.15}
  L_V(X)=V(A), \quad
  \text{where $A=\Gamma(S,{p_X}_*(\scrO_X)) \simeq
    \Gamma(X,\scrO_X)$.}\tag{15} 
\end{equation}
In case $X$ is affine, $A$ is just the affine ring of $X$, which is a
$k$-algebra, from which we retain only (in formula \eqref{eq:128.15})
the $k$-module structure.

Thus, as far as the notion of $k$-linearization goes, $V$-bundles
behave in a considerably nicer way than $W$-bundles, without any need
to go over to proobjects. Thus, it may be preferable to work with
$V$-bundles rather than $W$-bundles. We may wonder at this point why
not admit, then, as ``bundles'' any $k$-scheme $X$ which is
quasi-compact and quasi-separated, or at any rate any affine
$k$-scheme, as these can be quite conveniently described in terms of
$k$-algebras. Thus, semisimplicial affine schemes over $k$ just
correspond to \emph{co}-semisimplicial commutative algebras over $k$,
and likewise for maps -- and linearization just corresponds to
forgetting the algebra structure in this co-semisimplicial object, and
retain the structure of a co-semisimplicial $k$-module, which
corresponds dually to a semisimplicial vector bundle. This is a
perfectly simple relationship -- why bother about restricting the
``bundles'' from arbitrary affine $k$-schemes to those which are
isomorphic to a vector bundle?

We should remember here our initial motivation, which was, in case
when $k=\bZ$, to get a category of ``schematic'' homotopy types as
close as possible to the usual one. One plausible way of achieving
this is by restricting the notion of a bundle the more we can, so as
to get still the possibility of ``schematization'' for a very sizable
bunch of ordinary homotopy types. Postnikov dévissage then suggested
to work with so-called ``unipotent bundles'', and it was
almost\pspage{541} a matter of chance 50/50 that we took first the
choice of using $W$-bundles, rather than $V$-bundles which are the
dual choice, giving rise in some respects to a simpler algebraic
formalism (and to a less satisfactory one in some others\dots). In
both cases, an instinct of ``economy'' is leading us. It isn't always
clear that instinct isn't misleading at times -- after all, it would
be nice too to have a so-called ``schematic'' homotopy type (and hence
a usual one) associated to rather general types of semisimplicial
schemes, say. But here already, if we want to get a usual homotopy
type just by taking sections, we've seen that this isn't so automatic,
that this is tied up with the expectation that the maps
\eqref{eq:127.5} should be isomorphisms, which presumably will not be
true unless we make rather drastic extra assumptions on the
``bundles'' we are working with.

Thus, a first main test whether the choice of $W$-bundles or of
$V$-bundles is a workable one, is to see whether this condition on
\eqref{eq:127.5} is satisfied, possibly under a suitable extra
assumption, such as $X_0=e$ or $X_0=X_1=e$. The next test, presumably
a deeper one, is whether these choices allow for a description of the
homotopy sheaves $\bpi_i(X_*)$ in terms of the Lie functor, under the
natural flatness assumption on the components of $X_*$ (compare
section~\ref{sec:118}). The only clue so far for such a relationship
comes from Postnikov dévissage, and this relationship isn't proved
either, except in the more or less tautological ``Postnikov case'',
and without naturality. In a way, it appears as a rather strange kind
of relationship, one which implies that the homotopy type of a
semisimplicial bundle is very strongly dominated (almost determined,
one might say) by its formal completion along a section -- and even by
the first-order infinitesimal neighborhood already. Instinct again
tells us that reducing, as it were, a scheme to a tangent space at one
of its points might make sens when the scheme is isomorphic to affine
$n$-space, and that it is surely nonsense for general affine schemes
(such as an elliptic curve minus a point say!). To say it differently,
we feel that such a thing may be reasonable only when the given
schemes $X_n$ may be thought of as ``homotopically trivial'' in some
sense or other, which for an algebraic curve, say, over an
algebraically closed field $k$ is surely \emph{not} the case, unless
precisely the curve is isomorphic to the affine line.

\bigbreak

\presectionfill\ondate{14.10.}\pspage{542}\par

\hangsection[Quasi-coherent homological quasi-isomorphisms, versus
\dots]{Quasi-coherent homological quasi-isomorphisms, versus weak
  equivalences.}\label{sec:129}%
I have not been clear enough in yesterday's notes, when introducing
the linearization functor \eqref{eq:128.11}
\[L: U(k)\to U\subab(k)\]
for a suitable notion of ``bundles'' and ``linear bundles'', that it
is by no means automatic that such a functor (even when its existence
is granted) will induce a linearization functor
\begin{equation}
  \label{eq:129.16}
  \LH_\bullet:\HotOf_0(k)\to \HotabOf_0(k)\tag{16}
\end{equation}
for the corresponding derived categories. In other words, it is by no
means clear (and we haven't even tried yet to prove in the $U_W$
set-up of ``unipotent bundles'') that if
\begin{equation}
  \label{eq:129.17}
  X_*\to Y_*\tag{17}
\end{equation}
is a weak equivalence of semisimplicial objects of $U(k)$, that the
corresponding map
\begin{equation}
  \label{eq:129.17prime}
  L(X_*)\to L(Y_*)\tag{17'}
\end{equation}
of semisimplicial objects in $U\subab(k)$ will again be a weak
equivalence. At any rate, for this statement even to make sense, we'll
have to make clear what notion of weak equivalence we are taking for
maps between semisimplicial objects of $U\subab(k)$, so as to get a
localized category $\HotabOf_0(k)$. The first obvious choice, of
course, that comes to mind is to take the notion of weak equivalence
for the corresponding semisimplicial sheaves of sets, which presumably
is going to be the right one. Thus, the existence of a ``total
homology'' or ``linearization'' functor \eqref{eq:129.16} may be
viewed as the ``schematic'' analog of W.~H.~C.~Whitehead's theorem in
the discrete set-up, which we had been puzzling about already in
section~\ref{sec:92} (in connection with replacing the test category
$\Simplex$ by a more or less arbitrary small category). In the case of
semisimplicial $W$-bundles or $V$-bundles, under some extra
assumption, presumably, such as $X_0=e$ or $X_0=X_1=e$ on the
semisimplicial objects we are working with, I do expect that the
linearization functor transforms weak equivalences into weak
equivalences, and that the converse holds true, too. This seems to me
to be the third main test (besides the two considered on previous
page~\ref{p:541}) about a notion of ``bundle'' being suited for
developing a theory of schematic homotopy types.

This suggests that when working with more general semisimplicial
sheaves over $k$, such as semisimplicial schemes, say, it may be
useful\pspage{543} to introduce a notion of ``quasi-coherent
homological quasi-isomorphism'' between such objects, as a map
\eqref{eq:129.17} such that the corresponding map
\eqref{eq:129.17prime} should be a quasi-isomorphism, which presumably
may be interpreted also in cohomological terms, as usual. As just
noticed, it is doubtful that this is implied by \eqref{eq:129.17}
being a weak equivalence, and even if it should be implied, it looks
considerably weaker in a way -- just as in the set-up of usual
homotopy types, a homology equivalence is a considerably weaker notion
than homotopy equivalence, unless we make a $1$-connectedness
assumption. We could reinforce, of course, this notion of
(quasi-coherent) ``homological'' or ``cohomological''
quasi-isomorphism by taking to the cohomological version of it, and
instead of quasi-coherent coefficients \emph{coming from the base}
$S=\Spec(k)$, admit equally the analog of ``\emph{twisted
  coefficients}'', which would amount here to taking as coefficients
quasi-coherent sheaves $F_*$ on $Y_*$, such that for any structural
morphism
\[\varphi:Y_m\to Y_n\]
associated to a map in $\Simplex$, the corresponding map of
quasi-coherent sheaves on $Y_m$
\[\varphi^*(F_n)\to F_m\]
should be an isomorphism. (These ``twisted coefficients'' should
correspond to quasi-coherent sheaves on $S=\Spec(k)$, \emph{on which
  the sheaf of groups $\bpi_1(Y_*)$ operates}, in case at least when
$Y_*$ is relatively $0$-connected and endowed with a section, so that
$\bpi_1(Y_*)$ makes sense.)\enspace My point here is that it may be
interesting to take the derived category of a suitable category of
semisimplicial schemes (submitted to some very mild conditions, such
as quasi-compactness and quasi-separation, say) with respect to this
notion of q.c.h.~quasi-isomorphism -- with the hope that the set of
maps with respect to which we are now localizing is wide enough, so as
to get \emph{the same} derived category as when working only with
``bundles'', namely getting just schematic homotopy types. This would
give a very strong link between more or less arbitrary semisimplicial
schemes over $k$ (and stronger still when $k=\bZ$), and ordinary
homotopy types, of a much subtler nature than the known link with
ordinary pro-homotopy types via étale cohomology with discrete
coefficients. But maybe the daydreaming is getting here so much out of
reach or maybe simply crazy, that I better stop along this
line!\pspage{544} 

\hangsection{A case for non-connected bundles.}\label{sec:130}%
In a more down-to-earth line, and reverting to the ``unipotent''
approach still (in either $W$- or $V$-version for ``unipotent
bundles''), I would still like to point out one rather mild extension
of the set-up as contemplated initially. The suggestion here is to
admit as components for our semisimplicial models not merely sheaves
of sets which are ``unipotent bundles'', but equally direct sums of
such. After making such extension, the linearization functor $L$
(\eqref{eq:128.11} p.~\ref{p:539}) still makes sense, provided
$U\subab(k)$ is stable under infinite sums (otherwise, we'd have to
restrict to a finite number of connected components for our
``bundles''). Presumably, once we get into this, we will have to admit
``twisted'' finite direct sums as well, in order to have basic notions
compatible with descent -- never mind such technicalities at the
present stage of reflections! Thus, applying componentwise (i.e., to
each component $X_n$) the \piz-functor (``connected components''), a
semisimplicial ``bundle'' $X_*$ in the wider sense gives rise to an
associated usual semisimplicial set, let's call it $\pi_{0/k}(X_*)$,
together with a map
\begin{equation}
  \label{eq:130.18}
  X_*\to \pi_{0/k}(X_*),\tag{18}
\end{equation}
where the second-hand side is interpreted as a semisimplicial constant
object over $k$. (For simplicity, we have assumed here $\Spec(k)$ to
be $0$-connected, and that the direct sums involved in the $X_n$'s are
not twisted\dots)\enspace Intuitively, we may interpret
\eqref{eq:130.18} as defining $X_*$ as a (``strict'', namely
componentwise connected) \emph{schematic} homotopy type, lying
``over'' the \emph{discrete} (or ``constant'') homotopy type
$\pi_{0/k}(X_*)=E_*$. The latter introduces homotopy invariants of its
own, which strictly speaking shouldn't be viewed as being of a
``schematic'' nature. Thus, when our exclusive emphasis is on studying
the ``strict'' schematic homotopy types, we'll restrict our models
$X_*$ by demanding that the associated discrete homotopy type
$\pi_{0/k}(X_*)$ should be \emph{aspheric}, i.e., isomorphic in
$\HotOf$ to a one-point space. Under this restriction, presumably,
working with those slightly more general models should give (up to
equivalence) the same derived category $\HotOf_0(k)$, as when working
with (connected) unipotent bundles. Allowing connected components may
prove useful for giving a little more ``elbow freedom'' in working
with models, as it allows for instance anodyne operations such as
taking direct sums. Thus, all constant semisimplicial objects will be
allowed, or at any rate those which correspond to aspheric
semisimplicial sets -- which includes notably the objects of
\Simplexhat{} represented\pspage{545} by the standard simplices
$\Simplex_n$; and it is useful indeed to be able to have these among
our models.

My motivation for suggesting to allow for connected components for our
bundles, comes from an attempt to perform the ``integration''
operation for an indexes family of schematic homotopy types, under
suitable assumptions -- one among these being that the indexing small
category $I$ should be aspheric (in order to meet the above
asphericity requirement on $E_*$, for the semisimplicial bundle
obtained). When just paraphrasing the construction of integration for
ordinary homotopy types (as indicated in section~\ref{sec:69} -- we're
going to come back upon this in a later chapter), we just cannot help
but run into a bunch of connected components, for the components of
the semisimplicial sheaf expressing the integrated type. The general
idea here, of course, is to perform the basic homotopy operations
(essentially, integration and cointegration operations) \emph{in the
  context of semisimplicial sheaves}, and then look and see whether
(or when) this doesn't take us out of the realm semisimplicial
``bundles''. I am sorry I didn't work out any clear-cut result along
these lines yet, and I am going to leave it at that for the time
being.

One serious drawback, however, when allowing for connected components
of the components $X_n$ or our semisimplicial bundles, is that we can
hardly expect anymore that the homotopy sheaves of $X_*$ may be
expressed in terms of tangent sheaves, as contemplated in
section~\ref{sec:118}. At any rate, this extension of our category of
``models'' would seem a reasonable one only if we are able to show
that we get the same derived category (up to equivalence) as when
working with the more restricted models, using as components $X_n$
only (connected) unipotent bundles. I didn't make any attempt either
to try and prove such a thing.

\presectionfill\ondate{16.9.}\par

\hangsection[Tentative description of the spherical functor
$\smash{\widetilde S}$, and \dots]{Tentative description of the spherical
  functor \texorpdfstring{$\widetilde S$}{S}, and ``infinitesimal''
  extension of the basic notion of ``bundles''.}\label{sec:131}%
Yesterday again, I didn't do any mathematics -- instead, I have been
writing a ten pages typed report on the preparation and use of kimchi,
the traditional Korean basic food of fermented vegetables, which I
have been practicing now for over six years. Very often friends ask me
for instructions for preparing kimchi, and a few times already I
promised to put it down in writing, which is done now. Besides this, I
wrote to Larry Breen to tell him a few words about my present
ponderings as he is the one person I would think of for whom my
rambling\pspage{546} reflections on schematization and on schematic
homotopy types may make sense.

Definitely, my suggestion in the last notes, to work with
non-connected ``bundles'', isn't much more than the reflection of my
inability so far to make a breakthrough and get the ``left'' homotopy
constructions in terms of semisimplicial (connected) unipotent bundles
alone. Besides the serious drawback already pointed out at the end of
section~\ref{sec:130}, another one occurred to me -- namely that with
the suggested extension, one is losing track of the extra condition
$X_0=e$ (or even $X_0=X_1=e$) on our models, which is often remaining
implicit in the notes, and which, however, is an important
restriction, actually needed for the kind of things we want to do. For
one thing, when we drop it, and even when we otherwise restrict to
components $X_n$ which are standard affine spaces $E_k^{d_n}$, say,
there is no hope of getting a Lie-type description for the homotopy
sheaves $\bpi_i(X_*)$ -- for instance, the Lie-type invariants we get
by using different sections of $X_0$ over $S=\Spec(k)$ are by no means
related (as the $\bpi_i$ should) by a transitive system of canonical
isomorphisms; this is easily seen already when taking the ``trivial''
semisimplicial bundle described by $X_n=X_0$ for all $n$ (the constant
functor with value $X_0$ from $\Simplexop$ to $U(k)$)! Of course, even
when allowing connected components for the ``bundles'' $X_n$, we may
formally still throw in the condition $X_0=e$ -- but this is cheating
and no use at all, if we remember that the main motivation for
allowing connected components was in order to be able (in suitable
cases) to perform an integration operation of schematic homotopy
types. But when following the standard construction, for the resulting
$X_*$ the component $X_0$ as well as all others will have lots of
connected components, and hence the condition $X_0=e$ will not hold
true. The same remark applies also to the ``constant'' semisimplicial
bundles (namely with components which are ``constant'' schemes over
$S$) defined by the standard simplices $\Simplex_n$. To sum up, while
it \emph{is} certainly quite useful to view semisimplicial
``bundles'', used for describing the ``schematic homotopy types'' we
are after, as particular cases of more general semisimplicial sheaves
on the fpqc site of $S=\Spec(k)$, we will probably have to be quite
careful in keeping the ``bundles'' we are working with restricted
enough, and not confuse our models for schematic homotopy types
(allowing for a nice description of the basic $\Lpi_\bullet$ and
$\LH_\bullet$ invariants) with more general semisimplicial sheaves
which may enter the picture in various ways.

I\pspage{547} would like, however, to suggest still another extension
of the notion of a ``bundle'', which maybe will prove something better
than just a random way out of embarrassment! It has to do with an
attempt to come to a description of one among the four ``basic
functors'' of section~\ref{sec:124}, namely the ``spherical'' functor
\begin{equation}
  \label{eq:131.1}
  \widetilde S: \HotabOf_0(k)\to\HotOf_0(k),\tag{1}
\end{equation}
which for the time being is remaining hypothetical, due notably to my
inability so far to carry out the suspension operation for schematic
homotopy types. There are two basic formal properties giving us some
clues about this functor, namely, it should be left adjoint to the
``total homotopy'' or ``Lie'' functor $\Lpi_\bullet$, and it should be
right inverse to the (reduced) total homology functor $\LtH_\bullet$
(cf.\ pages~\ref{p:523}--\ref{p:525}). Let's work for the sake of
definiteness, for describing $\HotOf_0(k)$, with semisimplicial
$V$-bundles $X_*$ satisfying the extra assumption $X_0=e$ (and
presumably a little more), corresponding therefore to
co-semisimplicial algebra $A^*$ satisfying $A^0=k$, such that the
components $A^n$ be isomorphic to symmetric algebras
\begin{equation}
  \label{eq:131.2}
  A^n\simeq \Sym_k(M^n),\tag{2}
\end{equation}
where we'll assume the modules $M^n$ (or equivalently, the algebras
$A^n$) \emph{flat} over $k$. The condition $A^0=k$ implies that the
algebras $A^n$ are $k$-augmented, and if $\mathfrak m^{(n)}$ is the
augmentation ideal, we get
\begin{equation}
  \label{eq:131.3}
  \mathfrak m^{(n)} / (\mathfrak m^{(n)})^2 \simeq M^n,\tag{3}
\end{equation}
thus, the modules $M^n$ may be viewed as the components of a
co-semisimplicial $k$-module $M^*$, giving rise to a semisimplicial
vector bundle $V(M^*)$ over $k$, and we may write
\begin{equation}
  \label{eq:131.4}
  \Lpi_\bullet(X_*) \simeq V(M^*),\tag{4}
\end{equation}
where the second member is viewed as a chain complex (rather than as a
semisimplicial module) in the usual way. Consider now an object $L_*$
in $\HotabOf_0(k)$, it is described in terms of a co-semisimplicial
$k$-module $L^*$ by the formula
\begin{equation}
  \label{eq:131.5}
  L_*=V(L^*),\tag{5}
\end{equation}
and we get now
\begin{equation}
  \label{eq:131.star}
  \Hom_{\scrO_k}(L_*,\Lpi_\bullet(X_*)) \simeq
  \Hom_{\scrO_k}(V(L^*),V(M^*)) \simeq
  \Hom_{\scrO_k}(M^*,L^*).\tag{*}
\end{equation}
Now, if $M$ is any $k$-module, let
\begin{equation}
  \label{eq:131.6}
  D(M) = \Sym_k(M)(1) \simeq k\oplus M\tag{6}
\end{equation}
be\pspage{548} the corresponding augmented $k$-algebra admitting $M$
as a square-zero augmentation ideal, and let
\begin{equation}
  \label{eq:131.7}
  I(M) = \Spec(D(M))\tag{7}
\end{equation}
its spectrum, which is a first-order infinitesimal neighborhood of
$S=\Spec(k)$. For variable $M$, $D(M)$ depends covariantly, $I(M)$
contravariantly on $M$. Thus,
\begin{equation}
  \label{eq:131.8}
  D(L^*)\tag{8}
\end{equation}
is a co-semisimplicial algebra, depending covariantly on $L^*$, hence
contravariantly on $L_*=V(L^*)$, and accordingly
\begin{equation}
  \label{eq:131.9}
  I(L^*) = \Spec(D(L^*))\tag{9}
\end{equation}
(where in both members $I$ and $\Spec$ are applied again
componentwise) is a semisimplicial affine scheme, indeed a first-order
infinitesimal one, depending contravariantly on $L^*$, hence
covariantly on $L_*$. Now, the last member of \eqref{eq:131.star} may
be interpreted non-linearly as
\begin{equation}
  \label{eq:131.starstar}
  \Hom_{\textup{co-ss $k$-alg}}(A^*,D(L^*))
  \simeq
  \Hom(I(L^*), X_*),\tag{**}
\end{equation}
where in the second term the $\Hom$ means homomorphisms of
semisimplicial schemes. Let's now write (with an obvious
afterthought!)
\begin{equation}
  \label{eq:131.10}
  \widetilde S(L_*) = I(L^*)\quad
  \text{whenever $L_*=V(L^*)$,}\tag{10}
\end{equation}
the sequence of isomorphisms \eqref{eq:131.star} and
\eqref{eq:131.starstar} may be summed up by
\begin{equation}
  \label{eq:131.11}
  \Hom(L_*,\Lpi_\bullet(X_*)) \simeq
  \Hom(\widetilde S(L_*),X_*),\tag{11} 
\end{equation}
where the first $\Hom$ means maps of semisimplicial $\scrO_k$-modules,
whereas the second is a $\Hom$ of semisimplicial schemes. It looks
very much like the adjunction formula (p.~\ref{p:524},
\hyperref[eq:124.7]{(7b)}) we are after, with however two big grains
of salt. The smaller one is that this formula doesn't take place in
derived categories, but rather in the categories of would-be
models. The considerably bigger grain of salt is that $\widetilde
S(L_*)$ isn't at all in our model category, its components are purely
infinitesimal, first-order schemes, and far from being vector bundles!
Essentially, what we have been using for getting the (admittedly quite
tautological) adjunction formula \eqref{eq:131.11} is that the Lie
functor on schemes-with-section over $S$ is representable in an
obvious way, namely by the scheme-with-section
\[I(k)=\Spec(k[T]/(T^2)),\]
which is a first-order infinitesimal scheme-with-section over $k$.

Our\pspage{549} tentative $\widetilde S$ functor in \eqref{eq:131.10}
has been constructed in the most evident way, in order to satisfy an
adjunction formula \eqref{eq:131.11}, valid on the level of
semisimplicial objects (and carrying over, hopefully, to the similar
adjunction formula for suitable derived categories). Next question is
then, what about the inversion formula
\begin{equation}
  \label{eq:131.12}
  \LtH_\bullet( \widetilde S(L_*)) \simeq L_*\text{?}\tag{12}
\end{equation}
The question makes sense, as $\LtH_\bullet$ is defined for any
semisimplicial affine scheme-with-section over $k$, or equivalently
for any co-semisimplicial augmented $k$-algebra $A^*$ over $k$, by
taking the augmentation ideal $\mathfrak m^*$ in the latter and
retaining only its linear (co-semisimplicial) structure. Keeping this
in mind, formula \eqref{eq:131.12} comes out indeed a tautology again!

The tentative description we just got is indeed of a most seducing
simplicity, as seducing indeed as the description of homotopy in terms
of the Lie functor, and closely related to the latter. It gives as a
particular case an exceedingly simple description of the sought-for
``spheres over $k$'' $S(k,n)$. But it is clear that this description
is liable to makes sense only at the price of suitably extending the
notion of a ``bundle'' we are working with, in a rather different
direction, I would say, from adding (or allowing) connected
components, as suggested in the previous section. Maybe we might view
it, though, as a kindred, but somewhat subtler extension of our
initial bundles, namely that we are now allowing, not a discrete
non-trivial ``set'' or $k$-scheme of connected components, but rather,
an infinitesimal one. More specifically, the suggestion which comes to
mind here, is to call now ``bundle'' over $k$ any scheme $X$ over $k$
admitting a subscheme
\[X_*\subset X\]
in such a way that $X_0$ should be a $V$-bundle (namely isomorphic to
a vector bundle) over $k$, and $X$ should be an infinitesimal
neighborhood of $X_0$, i.e., $X_0$ should be definable by a
quasi-coherent ideal on $X$ which is nilpotent. Equivalently, in terms
of the affine ring $A$ of $X$, we are demanding that $A$ should admit
a nilpotent ideal $J$ (which is of course not part of its structure),
such that $A/J$ should be isomorphic to a symmetric algebra over $k$
(with respect to some $k$-module $M$). Possibly, we may have moreover
to impose further flatness restrictions.

When working with this extended notion of bundles, there is no problem
for describing for the corresponding semisimplicial
models\pspage{550} the three functors $\LtH_\bullet$, $\widetilde K$,
$\widetilde S$. Indeed, as we just recalled, the first of the three
functors is well-defined and has an evident description for all
semisimplicial affine schemes over $k$. As for $\widetilde K$ and
$\widetilde S$, they are obtained in terms of a variable
co-semisimplicial $k$-module $L^*$ (representing the semisimplicial
vector bundle $L_*=V(L^*)$) by applying componentwise the functor
$\Sym_k$ and the first-order truncation $\Sym_k({-})(1)$, respectively
-- one may hardly imagine something simpler! This brings to my
attention that in terms of this description, we get a canonical
functorial map
\begin{equation}
  \label{eq:131.13}
  \widetilde S(L_*) \hookrightarrow \widetilde K(L_*)\tag{13}
\end{equation}
when working with the semisimplicial models, and hence presumably a
corresponding map for the functors between the relevant derived
categories $\HotabOf_0(k)$ and $\HotOf_0(k)$. Working either in the
model or in the derived categories, this map, as a matter of fact, may
be deduced from the basic formulaire of section~\ref{sec:124}, where
it had by then escaped my attention. Indeed, $\widetilde S$ is a left
adjoint of $\Lpi_\bullet$, and $\widetilde K$ a right adjoint of
$\LtH_\bullet$, to give such a map \eqref{eq:131.13} is equivalent
with giving either one of two maps
\begin{equation}
  \label{eq:131.14}
  \left\{
    \begin{tabular}{@{}lll@{}}
      a) & $L_*\to \Lpi_\bullet(\widetilde K(L_*))$ & \\
      b) & $\LtH_\bullet(\widetilde S(L_*)) \to L_*$ & ,
    \end{tabular}
  \right.\tag{14}
\end{equation}
and the formulaire provides for two such maps, namely the ``inversion
isomorphisms'' (\eqref{eq:124.9}, p.~\ref{p:525}). Thus, there is an
extra property which was forgotten in the formulaire, namely that the
two maps \eqref{eq:131.13} associated to the two inversion
isomorphisms should be the same. A nicer way, then, to state the
formulaire is to consider the map \eqref{eq:131.13} as a basic datum,
and say that the two maps in \eqref{eq:131.14} deduced from it by the
adjunction property should be isomorphisms. The situation is
reminiscent of the two ways by which we could obtain the Hurewicz map
(p.~\ref{p:525}, \ref{subsec:124.C}) -- presumably, the basic data for
the formulaire of section~\ref{sec:124} should be the functors
$\widetilde S$ and $\widetilde K$ and the map \eqref{eq:131.13}
between them, with the property that the relevant adjoint functors
$\Lpi_\bullet$ and $\LtH_\bullet$ exist, and that the corresponding
maps in \eqref{eq:131.14} should be isomorphisms, which then will
allow to define a \emph{unique} Hurewicz map
$\Lpi_\bullet\to\LtH_\bullet$.

As long as we are sticking to the purely formal aspect, and even when
working in the larger context of semisimplicial affine schemes over
$k$ satisfying merely $X_0=e$, or more generally still, dropping the
last restriction and taking ``$k$-pointed'' semisimplicial affine
schemes instead, the whole ``four functors formalism'' (including even
$\Lpi_\bullet$) as contemplated in section~\ref{sec:124} (and
with\pspage{551} the extra feature \eqref{eq:131.13} above as just
notice) goes over very smoothly, in an essentially tautological
way. As recalled on p.~\ref{p:547}, the $\Lpi_\bullet$-functor, when
interpreted on the co-semisimplicial side of the dualizing functor,
appears as a quotient of the $\LtH_\bullet$-functor, the latter
identified to the functor obtained by taking augmentation ideals of
co-semisimplicial algebras -- the quotient being obtained by dividing
out by the squares of the latter ideals. Dually, we get the Hurewicz
map for semisimplicial vector bundles, which is always an
inclusion. Again, imagine something simpler! The only trouble (but an
extremely serious one indeed!) is that in this general set-up, the
relation of the so-called $\Lpi_\bullet$-functor to homotopy
groups-or-sheaves becomes a very dim one. Definitely, the only firm
hope here is that the relationship between the two is OK (as
contemplated in section~\ref{sec:118}) whenever the components $X_n$
are actual flat vector bundles, satisfying moreover $X_0=e$ (at the
very least) -- plus possibly even some extra Kan type conditions
(sorry for the vagueness of even this one ``firm hope''!). If we take
already the ``next best'' set of assumptions, namely essentially that
the $X_n$ be flat ``bundles'' in the sense above (not necessarily
vector bundles, though), then the hoped-for relationship again seems
to vanish. The first case of interest, of course, is the case when
$X_*$ is of the form $\widetilde S(L_*)$, which includes (if our
$\widetilde S$ functor is ``the right one'' indeed) the $n$-spheres
over $k$. We get in this case (namely when $X_*$ is a first-order
neighborhood of the marked section) the trivial, and really stupid
relation
\[\Lpi_\bullet(\widetilde S(L_*)) \simeq L_*\quad\text{(!!!),}\]
which translates into: the homotopy groups of a sphere, computed in
the most naive ``Lie'' way, are canonically isomorphic to its homology
groups! Not much of a success\dots

\hangsection[A crazy tentative wrong-quadrant (bi)complex for \dots]{A
  crazy tentative wrong-quadrant
  \texorpdfstring{\textup(bi\textup)}{(bi)}complex for the homotopy
  groups of a sphere.}\label{sec:132}%
This makes it very clear that, while the functors $\widetilde K$,
$\widetilde S$, $\LtH_\bullet$ in our new context of semisimplicial
``bundles'' make perfectly good sense as they are, the $\Lpi_\bullet$
functor computed naively (taking tangent spaces) definitely doesn't,
except when actually working with flat (hence, essentially ``smooth'')
vector bundles as components of our semisimplicial models. This, after
all, shouldn't be too much of a surprise, if we remember the way
differentials and tangent spaces fit into a sweeping homology or
cohomology formalism. It has become quite familiar to people ``in the
know'' that taking the sheaf of $1$-differentials, say, or its dual,
or a sheaf of $1$-differentials or a tangent sheaf along a section,
behaves\pspage{552} as ``the'' good object in terms of homological
algebra and obstruction theory in various geometric situations,
\emph{only} in the case when the relative scheme ($X$ say) we are
working with is \emph{smooth} over the base scheme $S$ -- which in the
present case amounts to saying (when $X=X_n$ is a component of a
semisimplicial ``bundle'') that $X$ is indeed a flat vector bundle
over $S$. In more general cases, the work of
André-Quillen-Illusie\scrcomment{\textcite{Andre1967,Andre1974,Quillen1970,Illusie1971};
  see also \textcite{SGA6}\dots} tells us that the relevant object
which replaces $\Omega^1_{X/S}$ or its dual $\mathrm T_{X/S}$ is the
relative tangent or cotangent \emph{complex} $\mathrm L_\bullet^{X/S}$
or $\mathrm L^\bullet_{X/S}$, the second being the dual of the other
\begin{equation}
  \label{eq:132.15}
  \mathrm L^\bullet_{X/S} = \RbHom_{\scrO_X}(\mathrm L^{X/S}_\bullet,
  \scrO_X),\tag{15} 
\end{equation}
these objects being viewed, respectively, as objects in the derived
categories $\D_\bullet(\scrO_X)$ and $\D^\bullet(\scrO_X)$ (deduced
from chain and cochain complexes of $\scrO_X$-modules). When $X$ is
endowed with a section over $S$, the naive differentials and
codifferentials along this section should in the same way be replaced
by the co-Lie and Lie complexes
\begin{equation}
  \label{eq:132.16}
  \ell_\bullet(X/S,s) \quad\text{and}\quad
  \ell^\bullet(X/S,s)\simeq
  \RbHom(\ell_\bullet(X/S,s), \scrO_S),\tag{16}
\end{equation}
where $s$ is the given section, obtained from the previous complexes
by taking its inverse images $\mathrm Ls^*$ by $s$. As a matter of
fact, the chain complex $\mathrm L_\bullet^{X/S}$ can be realized
canonically, up to unique isomorphism, via a semisimplicial module
\[\mathrm L_*^{X/S}\]
on $X$, whose components are free $\scrO_X$-modules. Accordingly, we
get \eqref{eq:132.16} in terms of a well-defined semisimplicial
$\scrO_S$-module,
\begin{equation}
  \label{eq:132.17}
  \ell_*(X/S,s)\tag{17}
\end{equation}
whose components are free -- and as $S=\Spec(k)$, we may interpret
this more simply as a semisimplicial $k$-module with free
components. When we apply this to the components $X_n$ of a
semisimplicial bundle $X_*$, we get however the co-Lie invariants; in
order to get the relevant Lie invariants we'll have to take the duals
\begin{equation}
  \label{eq:132.17prime}
  \ell^*(X/S,s) = \Hom_k(\ell_*(X/S,s),k),\tag{17'}
\end{equation}
where $X$ is any one among the $X_n$'s, and $s$ its marked
section. Thus, the ``corrected'' description of $\Lpi_\bullet$, by
using the André-Quillen-Illusie version of the ``Lie-functor along a
section'', would seem to be
\begin{equation}
  \label{eq:132.18}
  \Lpi_\bullet(X_*) = \ell^*(X_*/k,e_*)\quad
  \text{(?),} \tag{18}
\end{equation}
where\pspage{553} now the second member appears as a mixed complex of
$k$-modules
\[(n,p) \mapsto \ell^p(X_n/k,e_n) : \Simplexop\times\Simplex \to
\AbOf_k,\]
contravariant with respect to the index $n$, covariant with respect to
$p$. Translating this via Kan-Dold-Puppe, we get a bicomplex of
$k$-modules, which we'll write in cohomological notation (with the two
partial differential operators of degree $+1$)
\[(C^{n,p}) = (C^{n,p}(X_*))\]
situated in the ``quadrant''
\[n\le 0, p\ge 0.\]
As we finally want an object of the derived category $\D(\AbOf_k)$ of
the category of $k$-modules, and even an object in the subcategory
$\D_\bullet(\AbOf_k)$, the evident thing that seems to be done now is
to take the associated simple complex, which hopefully may prove to be
the ``correct'' expression of the looked-for $\Lpi_\bullet$ -- and
this (if any) should be the precise meaning of \eqref{eq:132.18}.

The associations for getting \eqref{eq:132.18} are very tempting,
indeed, the expression we got makes us feel a little uneasy,
though. The main point is that the quadrant where our bicomplex lies
in is one of the two ``awkward'' ones, which implies that a)\enspace
for a given total degree, there are an infinity of summands occurring
(and one has to be careful, therefore, if these should be ``summands''
indeed, or rather ``factors'', namely if we should take an infinite
direct sum, or an infinite product instead); and b)\enspace the total
complex will have components of any degree both positive and negative,
and it isn't clear at all that it should be (as an object of
$\D(\AbOf_k)$) of the nature of a chain complex, namely that its
cohomology modules vanish for (total) degree $d>0$. If it should turn
out that this is not so (I didn't yet check any particular case), this
would imply for the least that \eqref{eq:132.18} should be corrected,
by taking the relevant truncation of the second-hand side.

Working with the $\mathrm L^*_{X/S}$ and $\ell^*(X/S,s)$ invariants
brings in a slightly awkward feature of its own which we have been
silent about, namely (except under suitable finiteness conditions) it
brings in non-quasi-coherent $\scrO_X$ or $\scrO_S$-modules. This may
encourage us to dualize \eqref{eq:132.18}, which will amount to
working with the co-semisimplicial algebra $A^*$ expressing $X_*$, and
taking componentwise the reduced (via augmentations) André-Quillen
complexes (rather than their duals). At any rate,\pspage{554} the
would-be expression of ``total co-homotopy'' of $X_*$ we'll get this
way isn't so much more appealing than \eqref{eq:132.18} -- it lies
still in one of the wrong quadrants, which definitely makes us feel
uncomfortable.

In principle, the tentative formula \eqref{eq:132.18}, when applied
say to an object such as
\[S(n,k) \eqdef \widetilde S(k[n]),\]
gives a rather explicit (but for the time being highly hypothetical!)
expression of the homotopy modules of the $n$-sphere over $k$, which
in case $k=\bZ$ are hoped to be just the homotopy groups of the
ordinary $n$-sphere. To test whether this makes at all sense, we'll
have to understand first the structure of the André-Quillen ``Lie
complex'' of an algebra \eqref{eq:131.6} of the type $D(M)$, for
variable $k$-module $M$. I haven't started looking into this yet, and
I doubt I am going to do it presently.

At any rate, whether or not the formula \eqref{eq:132.18} we ended up
with is essentially correct, in order (among other things) to get a
method for computing the total homotopy $\Lpi_\bullet$ for a
semisimplicial bundle $X_*$ that isn't a flat vector bundle, we'll
have to find out some more or less explicit means of replacing $X_*$
by some $X_*'$ whose components \emph{are} flat vector bundles, and
which is isomorphic to $X_*$ in the relevant derived category
$\HotOf_0(k)$. When $k$ isn't a field, the question arises already
even when the components $X_*$ \emph{are} vector bundles, when these
are not flat. The first idea that comes to mind, from the
André-Quillen theory precisely, in terms of the co-semisimplicial
algebra $A^*$ expressing $X_*$, is to take ``projective resolutions''
of the various components $A^n$ by polynomial algebras. Again we end
up with a mixed complex, this time a complex of augmented algebras
depending on two indices $n,p$, covariant in $n$ and contravariant in
$p$, or the reverse if we replace those algebras $A^n_p$ by their
spectra $X_n^p$. In any case, we again end up in a ``wrong quadrant''!

The hesitating question that comes to mind now is whether it is at all
feasible to work with a category of models which isn't a category of
\emph{semisimplicial} bundles say, but one of such \emph{mixed
  wrong-quadrant} bundles; namely, use these as ``models'' for getting
hold of a reasonable derived category of ``schematic homotopy types''?
I never heard of anything such yet, and I confess that at this point
my (anyhow rather poor!) formal intuition of the situation breaks down
completely -- maybe the suggestion is complete nonsense, for some
wholly trivial reason! Maybe Larry Breen could tell me at once -- or
someone else who has more feeling than I for semisimplicial and
cosemisimplicial models?



\chapter{Abelianization II}
\label{ch:VII}

\presectionfill\ondate{22.10.}\pspage{555}\par

\hangsection{Birth of Suleyman.}\label{sec:133}%
Again nearly a week has passed by without writing any notes -- the
tasks and surprises of life took up almost entirely my attention and
my energy. It were days again rich in manifold events -- most
auspicious one surely being the birth, three days ago, of a little
boy, Suleyman, by my daughter. The birth took place at ten in the
evening, in the house of a common friend, in a nearby village where my
daughter had been awaiting the event in quietness. It came while
everyone in the house was in bed, the nearly five years old girl
sleeping next to her mother giving birth. The girl awoke just after
the boy had come out, and then ran to tell Y (their hostess) she got a
little brother. When I came half an hour later, the little girl was
radiant with joy and wonder, while telling me in whispers, sitting
next to her mother and to her newborn brother, what had just
happened. As usual with grownups, I listened only distractedly,
anxious as I was to be useful the best I could. There wasn't too much
by then I could do, though, as the mother knew well what was to be
done and how to do it herself -- presently, tie the chord twice, and
then cut it in-between. ``Now you are on your own, boy!'' she told him
with a smile when it was done. A little later she and Y helped the
baby to a warm bath -- and still later, while the girl was asleep next
to her brother and Y in her room, the mother took a warm bath herself
in the same basin, to help finishing with the labors. When her mother
arrived a couple of hours later, everybody in the house was asleep
except me in the room underneath, awaiting her arrival while taking
care of some fallen fruit which had been gathered that very day.

This birth I feel has been beneficial, a blessing I might say, in a
number of ways. The very first one I strongly perceived, was about
Suleyman's sister, who had so strongly participated in his birth. This
girl has been marked by conflict, and her being is in a state of
division, surely as strongly as any other child of her age. However,
after this experience, if ever her time should come to bear a child
and give birth, she will do so with joy and with confidence, with no
secret fear or refusal interfering with her labors and with the act of
giving birth -- this most extraordinary act of all, maybe, which a
human being is allowed to accomplish; this unique privilege of woman
over man, this blessing, carried by so many like a burden and like a
curse\dots\enspace And there are other blessing too in this birth
which I perceive more or less clearly, and surely others still which
escape my conscious attention.\pspage{556}

I suspect that potentially, in every single instance anew, the act of
giving birth, and the sudden arrival and presence of a newborn, are a
blessing, carry tremendous power. In many cases, however, the greatest
part (if not all) of this beneficial power is lost, dispersed through
the action of crispation and fear (including the compulsory medical
``mise en scène'', from which it gets over harder to get rid in our
well-to-do countries\dots). A great deal could be said on this matter
in general, while my own leanings at present would be, rather, to
ponder about the manifold personal aspects and meaning of the
particular, manifold event I have just been involved in. But to yield
fully now to this leaning of mine would mean to stop short with the
mathematical investigation and with these notes. The drive carrying
this investigation is very much alive, though, and I have been feeling
it was time getting back to work -- to this kind of work, I mean.

This brings me to somewhat more mundane matters -- such as the
beginning of school and teaching duties. This year I am in charge of
preparing the ``concours d'agrégation'',\scrcomment{these are the
  civil service competitive examinations for high school (lycée)
  teachers\dots} something which, I was afraid, would be rather
dull. Rather surprisingly, the first session of common work on a
``problème d'agrègue'' wasn't dull at all -- the problem looked
interesting, and two of the three students who turned up so far were
interested indeed, and the atmosphere relaxed and friendly. It looks
as though I was going again to learn some mathematics through
teaching, or at any rate to apply things the way I have known and
understood them (in a somewhat ``highbrow'' way, maybe) to the more
down-to-earth vision going with a particular curriculum (here the
``programme d'agrègue'') and corresponding virtuosity tests --
something rather far, of course, from my own relation to mathematical
substance! Besides these arpeggios, we started a microseminar with
three participants besides me, on the Teichmüller groupoids; I expect
one of the participants is going to take a really active interest in
the stuff. I have been feeling somewhat reluctant to start this
seminar while still involved with the homotopy story, which is going
to keep me busy easily for the next six months still, maybe even
longer. Sure enough, the little I told us about some of the structures
and operations to be investigated, while progressively gaining view of
them again after a long oblivion -- and as discovering them again
hesitatingly, while pulling them out of the mist by bribes and bits --
this little was enough to revive the special fascination of these
structures, and all that goes with them.\pspage{557} It seems to me
there hasn't been a single thing in mathematics, including
motives,\scrcomment{unreadable margin note}
which has exerted such a fascination upon me. It will be hard, I'm
afraid, to carry on a seminar on such stuff, and go on and carry to
their (hopefully happy!) end those ponderings and notes on homotopical
algebra, so-called -- and a somewhat crazy one too at times, I am
afraid! We'll see what comes out of all this! The night after the
seminar session at any rate, and already during the two hours drive
home, the Teichmüller stuff was brewing anew in my head. Maybe it
would have gone on for days, but next day already several things came
up demanding special attention, last not least being the birth of
Suleyman\dots 

\bigbreak\scrcomment{the typescript says ``26.8.'' but it
  must be a misprint!}

\presectionfill\ondate{26.10.}\par

\hangsection[Back to linearization in the modelizer \Cat: another
\dots]{Back to linearization in the modelizer
  \texorpdfstring{\Cat}{(Cat)}: another handful of questions around
  Kan-Dold-Puppe.}\label{sec:134}%
It is time to come back to the ``review'' of linearization of
(ordinary) homotopy types, and the homology and cohomology formalism
in the context of the modelizer \Cat, which has been pushed aside now
for over two months, for the sake of that endless digression on
schematization of homotopy types. The strong tie between these two
strains of reflection has been the equal importance in both of the
linearization process. Technically speaking though, linearization, as
finally handled in the previous chapter (via the so-called
``integrators''), looks pretty much different from the similar
operation in the schematic set-up, due (at least partly) to the choice
we made of models for expressing schematic homotopy types, namely
taking semisimplicial scheme-like objects; which means, notably,
working with $\Simplex$ rather than more general test categories, and
relying heavily on the Kan-Dold-Puppe relationship. In the discrete
set-up, it turned out (somewhat unexpectedly) that the latter can be
replaced, when working with more general categories $A$ than
$\Simplex$ (which need not even be test categories or anything of the
kind), by the $\mathrm Lp_!\supab$ operation (when $p:A\to e$ is the
map from $A$ to the final object $e$ in \Cat), computable in terms of
a choice of an ``integrator'' for $A$. More specifically, recalling
that
\begin{equation}
  \label{eq:134.1}
  p_!\supab : \Ahatab \to \Ab\tag{1}
\end{equation}
is defined as the left adjoint of the inverse image functor
\[p^*:\Ab\to\Ahatab,\]
associating to every abelian group the corresponding ``context''
abelian\pspage{558} presheaf on $A$, and
\begin{equation}
  \label{eq:134.2}
  \LH_\bullet = \mathrm Lp_!\supab : \D^-(\Ahatab) \to \D^-\Ab, \quad
  \text{inducing $\D_\bullet(\Ahatab)\to\D_\bullet\Ab$},\tag{2}
\end{equation}
is its left derived functor (computed using projective resolution of
complexes bounded above in \Ahatab), composing \eqref{eq:134.2} with
the tautological inclusion functor
\[\Ahatab\hookrightarrow\D_\bullet(\Ahatab)\hookrightarrow
  \D^-(\Ahatab),\]
we get a canonical functor
\begin{equation}
  \label{eq:134.3}
  \mathrm Lp_!\supab \restrto \Ahatab:\Ahatab\to \D_\bullet\Ab,\tag{3}
\end{equation}
which may be expressed in terms of an integrator $L_\bullet^B$ for
$B=A\op$, i.e., a projective resolution of the constant presheaf
$\bZ_B$ in \Bhatab, as the composition
\begin{equation}
  \label{eq:134.4}
  \Ahatab\to\Ch_\bullet\Ab\to\D_\bullet\Ab,\tag{4}
\end{equation}
where the first arrow is
\begin{equation}
  \label{eq:134.5}
  F\mapsto F *_\bZ L_\bullet^B : \Ahatab\to\Ch_\bullet\Ab,\tag{5}
\end{equation}
and the second is the canonical localization functor. The functors
\eqref{eq:134.2} and \eqref{eq:134.3} may be viewed as the (total)
homology functors of $A$, with coefficients in complexes of abelian
presheaves, resp.\ in abelian presheaves simply. When we focus
attention on the latter, we may introduce in \Ahatab{} the set of
arrows which become isomorphisms under the total homology functor
\eqref{eq:134.3}, let's call them ``\emph{abelian weak equivalences}''
in \Ahatab, not to be confused with the notion of quasi-isomorphism
for a map between complexes in \Ahatab. Let's denote by
\begin{equation}
  \label{eq:134.6}
  \HotabOf_A = (W_A\supab)^{-1}\Ahatab\tag{6}
\end{equation}
the corresponding localized category of \Ahatab, where $W_A\supab$
denotes the set of abelian weak equivalences in \Ahatab. Thus, the
choice of an integrator $L_\bullet^B$ for $A$ (i.e., a cointegrator
for $B$) gives rise to a commutative diagram of functors
\begin{equation}
  \label{eq:134.7}
  \begin{tabular}{@{}c@{}}
    \begin{tikzcd}[baseline=(O.base)]
      \Ahatab \ar[r]\ar[d] & \Ch_\bullet\Ab\ar[d] \\
      \HotabOf_A \ar[r] &
      |[alias=O]| {\ooalign{$\D_\bullet\Ab$\cr$\D_\bullet\Ab\eqdef\HotabOf$\hidewidth\cr}}
    \end{tikzcd}\qquad\qquad,
  \end{tabular}\tag{7}
\end{equation}
where the lower horizontal arrow is defined via \eqref{eq:134.3}
(independently of the choice of $L_\bullet^B$), the vertical arrows
being the localization functors. Beware that even when $A=\Simplex$
and $L_\bullet^B$ is the usual, ``standard''\pspage{559} integrator
for $\Simplex$, the upper horizontal arrow in \eqref{eq:134.7} is
\emph{not} the Kan-Dold-Puppe equivalence of categories, it has to be
followed still by the ``normalization'' operation. Thus, \emph{we
  certainly should not expect in any case the functor \eqref{eq:134.5}
  to be an equivalence -- however, we suspect that when $A$ is a test
  category} (maybe even a weak test category would do it), \emph{then
  the lower horizontal arrow in \eqref{eq:134.7}}
\begin{equation}
  \label{eq:134.8}
  \HotabOf_A \to \D_\bullet\Ab=\HotabOf \tag{8}
\end{equation}
\emph{is an equivalence of categories}. Whenever true, for a given
category $A$, this statement looks like a reasonable substitute (on
the level of the relevant derived categories) for the Kan-Dold-Puppe
theorem, known in the two cases $\Simplex$ and $\Globe$.

There are however still two important extra features in the case
$A=\Simplex$, which deserve to be understood in the case of more
general $A$. One is that a map $u:F\to G$ in \Ahatab{} is in
$W_A\supab$ (i.e., is an ``abelian weak equivalence'') if{f} it is in
$W_A$, i.e., if{f} it is a weak equivalence when forgetting the
abelian structures. In terms of a final object $e$ in $A$ (when such
object exists indeed),\footnote{it is enough that $A$ be $1$-connected
  instead of having a final object} viewing the categories $A_{/F}$
(for $F$ in \Ahatab) pointed by the zero map $e\to F$, which yields
the necessary base-point $e_F$ for defining homotopy invariants
$\pi_i(F)$ for $i\ge0$, the relationship just considered between
abelian weak equivalence and weak equivalence, will follow of course
whenever we have functorial isomorphisms
\begin{equation}
  \label{eq:134.9}
  \pi_i(F) \eqdef \pi_i(A_{/F},e_F) \simeq \mathrm H_i(A,F),\tag{9}
\end{equation}
which are known to exist indeed in the case $A=\Simplex$. It should be
noted that in section~\ref{sec:92}, when starting (in a somewhat
casual way) with some reflections on ``abelianization'', we introduced
a category (denoted by $\HotabOf_A$) by localizing \Ahatab{} with
respect to the maps which are weak equivalences (when forgetting the
abelian structures), whereas it has now become clear that, in case the
latter should not coincide with the ``abelian weak equivalences''
defined in terms of linearization or homology, it is the category
\eqref{eq:134.6} definitely which is the right one. Still, the
question of defining isomorphisms \eqref{eq:134.9} when $A$ has a
final object, and whether the equality
\begin{equation}
  \label{eq:134.10}
  W_A\supab = \forg_A^{-1}(W_A)\tag{10}
\end{equation}
holds, where
\[\forg_A:\Ahatab\to\Ahat\]
is the ``forgetful functor'', should be settled for general $A$.

The\pspage{560} other ``extra feature'' is about the relationship of
$W_A\supab$ with the notion of homotopism. In case $A=\Simplex$, a map
in \Ahatab{} is a weak equivalence (or equivalently, an abelian weak
equivalence) if{f} it is a homotopism when forgetting the abelian
structures -- this follows from the well-known fact that
semisimplicial groups are Kan complexes. There is of course also a
notion of homotopism in the stronger abelian sense -- a particular
case of the notion of a homotopism between semisimplicial objects in
an arbitrary category (here in \Ab). The Kan-Dold-Puppe theory implies
that if $F$ and $G$ in \Simplexhatab{} have as values projective
\bZ-modules, then a map $F\to G$ in \Simplexhatab{} is a weak
equivalence if{f} it is an ``abelian'' homotopism; and likewise, two
maps $u,v: F\rightrightarrows G$ are equal in $\HotabOf_\Simplex$ if{f}
they are homotopic (in the strict, \emph{abelian} sense of the
word). The corresponding statements still make sense and are true,
when replacing \bZ{} by any ground ring $k$, working in
$\Simplex_k\adjuphat$ rather than $\Simplex_\bZ\adjuphat =
\Simplexhatab$ -- and more generally still, when working in
\[\Simplex\adjuphat_\scrM = \bHom(\Simplexop, \scrM),\]
where \scrM{} is any abelian category. Now, replacing $\Simplex$ by an
arbitrary object $A$ in \Cat, these statements still make sense, it
would seem, provided we got on \Ahat{} a suitable ``homotopy
structure'', more specifically, a suitable ``homotopy interval
structure'' (in the sense of section~\ref{sec:51}). In
section~\ref{sec:97} (p.~\ref{p:355}) we reviewed the three standard
homotopy interval structures which may be introduced on a category
\Ahat, and their relationships -- the main impression remaining was
that, in case $A$ is a \emph{contractor} (cf.\ section~\ref{sec:95}),
those three structures coincide, and may be defined also in terms of a
\emph{contractibility structure} on \Ahat{} (section~\ref{sec:51}),
the latter giving rise (via the general construction of
section~\ref{sec:79}) to the usual notion of weak equivalence $W_A$ in
\Ahat. Thus, we may hope that the feature just mentioned for
$\Simplexhat$ may be valid too for \Ahat, whenever $A$ is a
contractor. In case $A$ is only a test category, then the most
reasonable homotopy interval structure on \Ahat, I would think, which
possibly may still yield the desired ``extra feature'', is the one
defined in terms of the set $W_A$ of weak equivalences in \Ahat{} as
in section~\ref{sec:54} (namely $h_{W_A}$).

\hangsection[Proof of ``integrators are abelianizators'' (block
against \dots]{Proof of ``integrators are abelianizators''
  \texorpdfstring{\textup(}{(}block against homology
  receding?\texorpdfstring{\textup)}{)}.}\label{sec:135}%
This handful of questions (some of which we met with before) is mainly
a way of coming into touch again with the abelianization story, which
has been becoming somewhat remote during the previous two months. I am
not sure I am going to\pspage{561} make any attempt, now or later, to
come to an answer. The last one, anyhow, seems closely related to the
formalism of closed model structure on categories \Ahat, and the
proper place for dealing with it would seem to be rather next chapter
VII,\scrcomment{I guess AG thought this was still
  Part~\ref{ch:VI}. Since there was never a second volume of PS, we
  refer to \textcite{Cisinski2006} for a discussion of closed model
  structures on presheaf topoi, now known as \emph{Cisinski model
    structures}\dots} where such structures are going to be
studied. What I would like to do, however, here-and-now, is to
establish at last the long promised relationship ``integrators are
abelianizators'', which we have kept turning around and postponing
ever since section~\ref{sec:99}, when those integrators were finally
introduced, mainly for this purpose (of furnishing us with
``abelianizators'').

First of all, I should be more outspoken than I have been before, in
defining the ``\emph{abelianization functor}'' (or ``\emph{absolute
  Whitehead functor}''):
\begin{equation}
  \label{eq:135.1}
  \bWh: \Hot\eqdef\scrW^{-1}\Cat \to \HotabOf
  \eqdef\D_\bullet\Ab,\tag{1} 
\end{equation}
without any use of the semisimplicial machinery which, at the
beginning of our reflections, had rather obscured the picture
(section~\ref{sec:92}). Defining such a functor amounts to defining a
functor
\begin{equation}
  \label{eq:135.2}
  \LH_\bullet:\Cat\to\D_\bullet\Ab,\tag{2}
\end{equation}
or ``total homology functor'', which should take weak equivalences
into isomorphisms. For an object $A$ in \Cat, we define
\begin{equation}
  \label{eq:135.3}
  \LH_\bullet(A) = \LH_\bullet(A,\bZ_A) =
  \mathrm L{p_A}_!\supab(\bZ_A),\tag{3} 
\end{equation}
where $\bZ_A$ is the constant abelian presheaf on $A$ with value \bZ,
and
\[p_A : A\to e\]
is the map to the final object of \Cat. We have to define the
functorial dependence on $A$. More generally, for pairs
\[(A,F), \quad\text{with $A$ in \Cat, $F$ in \Ahatab,}\]
the expression
\[\LH_\bullet(A,F)\quad\text{in $\D_\bullet\Ab$}\]
is functorial with respect to the pair $(A,F)$, where a map
\[(A,F) \to (A',F')\]
is defined to be a pair $(f,u)$, where
\begin{equation}
  \label{eq:135.4}
  f: A\to A', \quad
  u: F\to f^*(F'),\tag{4}
\end{equation}
the composition of maps being the obvious one. To see that such a pair
defines a map
\begin{equation}
  \label{eq:135.5}
  \LH_\bullet(f,u): \LH_\bullet(A,F) \to \LH_\bullet(A',F'),\tag{5}
\end{equation}
we\pspage{562} use
\[p_A = p_{A'} \circ f,\]
hence
\begin{equation}
  \label{eq:135.star}
  (p_A)_!\supab \simeq (p_{A'})_!\supab \circ f_!\supab,
  \quad\text{and}\quad
  \mathrm L(p_A)_!\supab \simeq \mathrm L(p_{A'})_!\supab
  \circ \mathrm Lf_!\supab,\tag{*}
\end{equation}
taking into account that $f_!\supab$ maps projectives to
projectives. Hence, we get
\[\LH_\bullet(A,F)\eqdef \mathrm L(p_A)_!\supab(F) \simeq
  \mathrm L(p_{A'})_!\supab(\mathrm Lf_!\supab(F)) \eqdef
  \LH_\bullet(A',\mathrm Lf_!\supab(F)).\]
This, in order to get \eqref{eq:135.5}, we need only define a map in
$\D^-({A'}\uphat\subab)$
\begin{equation}
  \label{eq:135.6}
  \mathrm Lf_!\supab(F)\to F',\tag{6}
\end{equation}
which will be obtained as the composition
\begin{equation}
  \label{eq:135.6prime}
  \mathrm Lf_!\supab(F)\to f_!\supab(F) \to F',\tag{6'}
\end{equation}
where the first map in \eqref{eq:135.6prime} is the canonical
augmentation maps towards the $\mathrm H_0$ object, and where the
second corresponds to $u$ in \eqref{eq:135.4} by adjunction. This
defines the map \eqref{eq:135.5}, and compatibility with compositions
should be a tautology. Hence the functor \eqref{eq:135.2}. To get
\eqref{eq:135.1}, we still have to check that when
\[f:A\to A'\]
is a weak equivalence, then
\[\LH_\bullet(f): \LH_\bullet(A) \to \LH_\bullet(A')\]
is an isomorphism in $\D_\bullet\Ab$. It amounts to the same to check
that for any object $K^\bullet$ in $\D^+\Ab$, the corresponding map
between the $\RHom$'s with values in $K^\bullet$ is an isomorphism in
$\D^+\Ab$. But the latter map can be identified with the map for
\emph{cohomology}
\[\RH^\bullet(A', K_{A'}^\bullet) \to \RH^\bullet(A, K_A^\bullet),\]
with coefficients in the constant complex of presheaves defined by
$K^\bullet$ on $A'$ and on $A$, which is an isomorphism, by the very
definition of weak equivalences in \Cat{} via cohomology.

Now, the statement ``an integrator is an abelianizator'' may be
rephrased rather evidently, without any reference to a given
integrator, as merely the commutativity, up to canonical isomorphism,
of the following diagram for a given $A$ in \Cat:\pspage{563}
\begin{equation}
  \label{eq:135.7}
  \begin{tabular}{@{}c@{}}
    \begin{tikzcd}[baseline=(O.base),column sep=large]
      \Ahat \ar[r,"\varphi_A"]\ar[d,"\Wh_A"'] &
      \Hot\ar[d,"{\bWh}"] \\
      \Ahatab \ar[r,"{\LH_\bullet(A,{-})}"'] &
      |[alias=O]| {\ooalign{$\Hotab$\cr$\Hotab\eqdef\D_\bullet\Ab$\hidewidth\cr}}
    \end{tikzcd}\qquad\qquad,
  \end{tabular}\tag{7}
\end{equation}
or equivalently, of the corresponding diagram where the categories
\Ahat, \Ahatab{} are replaced by the relevant localizations:
\begin{equation}
  \label{eq:135.7prime}
  \begin{tabular}{@{}c@{}}
    \begin{tikzcd}[baseline=(O.base)]
      \HotOf_A \ar[r]\ar[d] &
      \HotOf \ar[d] \\
      \HotabOf_A \ar[r] &
      |[alias=O]| \HotabOf
    \end{tikzcd}.
  \end{tabular}\tag{7'}
\end{equation}
The left vertical arrow in \eqref{eq:135.7} is of course the trivial
abelianization functor in \Ahat:
\[\Wh_A:\Ahat\to\Ahatab, \quad
  X\mapsto \bZ^{(X)}=\bigl(a\mapsto \bZ^{(X(a))}\bigr).\]
Going back to the definitions, the commutativity of \eqref{eq:135.7}
up to isomorphism, means that for $X$ in \Ahat, there is a canonical
isomorphism
\begin{equation}
  \label{eq:135.8}
  \LH_\bullet(A_{/X},\bZ) \simeq\LH_\bullet(A, \bZ^{(X)}),\tag{8}
\end{equation}
functorial with respect to $X$. To define \eqref{eq:135.8}, let
\[f : A'\eqdef A_{/X} \to A\]
be the canonical functor, then we get (using \eqref{eq:135.star} of
the previous page, with the roles of $A$ and $A'$ reversed)
\[\LH_\bullet(A',\bZ_{A'}) \simeq \LH_\bullet(A, \mathrm
  Lf_!\supab(\bZ_{A'})),\]
and the relation \eqref{eq:135.8} follows from the more precise
relation
\begin{equation}
  \label{eq:135.9}
  \mathrm Lf_!\supab(\bZ_{A'}) \tosim \bZ^{(X)}.\tag{9}
\end{equation}
To get \eqref{eq:135.9}, we remark that we have the tautological
relation
\begin{equation}
  \label{eq:135.9prime}
  f_!\supab(\bZ_{A'}) \simeq \bZ^{(X)},\tag{9'}
\end{equation}
hence the map \eqref{eq:135.9}. To prove it is an isomorphism amounts
to proving
\[\mathrm L_if_!\supab(\bZ_{A'}) = 0\quad\text{for $i>0$,}\]
but we have indeed
\begin{equation}
  \label{eq:135.10}
  \mathrm L_if_!\supab = 0\quad
  \text{for $i>0$, i.e., $f_!\supab$ is \emph{exact},}
  \tag{10}
\end{equation}
not only right exact, a rather special feature, valid for a
localization\pspage{564} functor like $f:A_{/X}\to A$, namely for a
functor which is fibering (in the sense of the theory of ``fibered
categories'') and has discrete fibers. It comes from the explicit
description of $f_!\supab$, as
\begin{equation}
  \label{eq:135.11}
  f_!\supab(F) = \Bigl( a\mapsto \bigoplus_{u\in X(a)} F(a)_u\Bigr),\tag{11}
\end{equation}
where, for a group object $F$ in $\Ahat_{/X}\simeq(A_{/X})\uphat$, and
$u$ in $X(a)$ (defining therefore an object of $A_{/X}$) $F(a)_u = $
fiber of $F(a)$ at $u\in X(a)$ is the corresponding abelian
group. I'll leave the proof of \eqref{eq:135.11} to the reader, it
should be more or less tautological.

Once the whole proof is written down, it looks so simple that I feel
rather stupid and can't quite understand why I have turned around it
for so long, rather than writing it down right away more than three
months ago! The reason surely is that I have been accustomed so
strongly to expressing everything via cohomology rather than homology,
that there has been something like a block against doing work
homologically, when it is homology which is involved. This block has
remained even after I took the trouble of telling myself quite
outspokenly (in section~\ref{sec:100}) that homology was just as
important and meaningful as cohomology, and more specifically still
(in section~\ref{sec:103}) that the proof I had in mind first, via
Quillen's result about $A\simeq A\op$ in \Hot{} and via cohomology,
was an ``awkward'', an ``upside-down'', one. The review on
abelianization I went into just after made things rather worse in a
sense, as there I took great pains to make the point that homology and
cohomology were just one and the same thing (so why bother about
homology!). Still, I \emph{did} develop some typically homologically
flavored formalism with the $*_k$ operation, and I hope that at the
end that block of mine is going to recede\dots

\bigbreak

\presectionfill\ondate{27.10.}\pspage{565}\par

\hangsection[Preliminary perplexities about a full-fledged ``six
\dots]{Preliminary perplexities about a full-fledged ``six operations
  duality formalism'' within
  \texorpdfstring{\Cat}{(Cat)}.}\label{sec:136}%
I just spent a couple of hours, after reading the notes of last night,
trying to get a better feeling of the basic homology operation in the
context of the basic localizer \Cat, namely taking the left derived
functor $\mathrm Lf_!\supab$ of the ``unusual'' direct image functor
for abelian presheaves, associated to a map
\[f: A \to B\]
in \Cat. This led me to read again the notes of section~\ref{sec:92},
when I unsuspectingly started on an ``afterthought, later gradually
turning into a systematic reflection on abelianization. With a
distance of nearly four months, what strikes me most now in these
notes is awkwardness of the approach followed at start, when yielding
to the reflex of laziness of describing abelianization of homotopy
types via the semisimplicial grindmill. The uneasiness in these notes
is obvious throughout -- I kind of knew ``au fond''\scrcomment{``au
  fond'' -- at the bottom -- deep down} that dragging in
the particular test category $\Simplex$ was rather silly. In
section~\ref{sec:100} only, does it get clear that the best
description for abelianization, with the modelizer \Cat, is via the
unusual direct image $p_!\supab$ corresponding to the projection
\[p=p_A:A\to e\]
of the ``model'' $A$ in \Cat{} to the final object (formulæ
\eqref{eq:100.11} and \eqref{eq:100.12} page~\ref{p:359}), by applying
$p_!\supab$ to a ``cointegrator'' $L_\bullet^A$ for $A$, namely to a
projective resolution of $\bZ_A$ (where $A$ is written $B$ by the way,
as I had been led before to replace a given $A$ in \Cat{} by its
``dual'' or opposite $B=A\op$, bound as I was for interpreting
``integrators'' for $A$ in terms of ``cointegrators'' for $B$\dots);
and in the next section the step is finally taken (against the
``block''!) to write
\[p_!\supab(L_\bullet^B)=\mathrm Lp_!\supab(\bZ_B),\]
which inserts abelianization into the familiar formalism of (left)
derived functors of standard functors. The reasonable next thing to do
was of course what I finally did only yesterday, namely check the
commutativity of the diagram \eqref{eq:135.7} of p.~\ref{p:563},
namely compatibility of this notion of abelianization (or homology)
with Whitehead's abelianization, when working with models coming from
\Ahat, $A$ any object in \Cat. This by the way, when applied to the
case $A=\Simplex$, gives at once the equivalence of the intrinsic
definition of abelianization, with the semisimplicial one we started
with -- provided we remark that in this case, $\mathrm Lp_!\supab$
(rather, its restriction to $\Simplexhatab$) may be equally
interpreted as the Kan-Dold-Puppe\pspage{566} functor (more
accurately, the composition of the latter with the localization
functor $\Ch_\bullet\Ab\to\D_\bullet\Ab=\HotabOf$). The equivalence of
the two definitions of abelianization is mentioned on the same
p.~\ref{p:369}, somewhat as a chore I didn't really feel then to dive
into. Besides the ``block'' against homology, the picture was being
obscured, too, by the ``computational'' idea I kept in mind, of
expressing $\mathrm Lp_!\supab(F)$ for $F$ in \Ahatab{} in terms of an
``integrator'' for $A$, i.e., as $F *_\bZ L_\bullet^B$ -- whereas I
should have known best myself that for establishing formal properties
relating various derived functors, the particular approaches used for
``computing'' them more or less elegantly are wholly irrelevant\dots

One teaching I am getting out of all this, is that when expressing
abelianization, or presumably any other kind of notion or operation of
significance for homotopy types, one should be careful, for any
modelizer one chooses to work in, to dig out the description which
fits smoothly those particular models. Clearly, when working with
semisimplicial models, the description via tautological abelianization
$\Wh_\Simplex$ and using Kan-Dold-Puppe is the best. When working
within the modelizer \Cat, though, making the detour through
\Simplexhat{} is awkward and makes us just miss the relevant
facts. Once we got this, it gets clear, too, what to do when
$\Simplex$ is replaced by any other object $A$ in \Cat{} when taking
models in \Ahat{} (never minding even whether \Ahat{} is indeed a
``modelizer''): namely, apply $\Wh_A$, and take total homology!

This brings to my mind another example, very similar indeed. Some time
after I got across Thomason's nice paper, showing that \Cat{} is a
closed model category (see comments in section~\ref{sec:87}), I got
from Tim Porter another reprint of
Thomason's,\scrcomment{\textcite{Thomason1979}} where he gives the
description of ``homotopy colimits'' (or ``integration'', as I call
it) within the modelizer \Cat, in terms of the total category
associated to a fibered category (compare section~\ref{sec:69}). There
he grinds through a tedious, highly technical proof, whereas the
direct proof when describing ``integration'' as the left adjoint
functor to the tautological ``inverse image'' functor, is more or less
tautological, too. The reason for this awkwardness is again that,
rather than being content to work with the models as they are,
Thomason refers to the
Bousfield-Kan\scrcomment{\textcite{BousfieldKan1972}} description of
colimits in the semisimplicial set-up which he takes as his definition
for colimits. I suspect that \Cat{} is the one modelizer most ideally
suited for expressing the ``integration'' operation, and that the
Bousfield-Kan description is just the obvious, not-quite-as-simple one
which can be deduced\pspage{567} from the former, using the relevant
two functors between \Cat{} and \Simplexhat{} which allow to pass from
one type of models to another. (Sooner or later I should check in
Bousfield-Kan's book whether this is so or not\dots)\enspace Thomason,
however, did the opposite, and it is quite natural that he had to pay
for it by a fair amount of sweat! (Reference of the paper: Homotopy
colimits in the category of small categories, Math.\ Proc.\ Cambridge
Philos.\ Soc.\ (1979), 85, p.~91--109.)

The proof written down yesterday for compatibility of abelianization
of homotopy types with the Whitehead abelianization functor within a
category \Ahat, still goes through when replacing abelianization by
$k$-linearization, with respect to an arbitrary ground ring $k$ (not
even commutative). It wasn't really worth while, though, to introduce
a ring $k$, as the general result should follow at once from the case
$k=\bZ$, by the formula
\begin{equation}
  \label{eq:136.1}
  \mathrm Lp_!^k(k_A) \simeq \mathrm Lp_!^\bZ(\bZ_A) \Lotimes_\bZ
  k,\tag{1} 
\end{equation}
where in the left-hand side we are taking the left derived functor for
the functor
\[p_!^k : A\uphat_k \to (\AbOf_k)\]
generalizing $p_!\supab = p_!^\bZ$ (with $p: A\to e$ as above), and
where in the right-hand side we are using the ring extension functor
\[\Lotimes_\bZ : \D_\bullet\Ab\to\D_\bullet(\AbOf_k)\]
for the relevant derived categories. This reminds me of the need of
developing a more or less exhaustive formulaire around the basic
operations
\begin{equation}
  \label{eq:136.2}
  \mathrm Lf_!, \quad f^*, \quad \mathrm Rf_*,\tag{2}
\end{equation}
including the familiar one for the two latter, valid more generally
for maps between ringed topoi, and including also a ``\emph{projection
  formula}'' generalizing \eqref{eq:136.1}. Such a formula will be no
surprise, surely, to a reader familiar with a duality context (such as
étale cohomology, or ``coherent'' cohomology of noetherian schemes,
say), where a formalism of the ``four variances'' $f_!$, $f^*$, $f_*$,
$f^!$ and the ``two internal operations'' $\Lotimes$ and $\RbHom$ can
be developed -- it would seem that the formal properties of the triple
\eqref{eq:136.2}, together with the two internal operations just
referred to, are very close (for the least) to those of the slightly
richer one in duality set-ups, including equally an ``unusual inverse
image'' $f^!$, right adjoint to $\mathrm Rf_*$ (denoted sometimes
simply by $f_*$). This similarity\pspage{568} is a matter of course,
as far as the two last among the functors \eqref{eq:136.2}, together
with the two internal operations, are concerned, as in both contexts
(homology and cohomology formalism within \Cat{} on the one hand, and
the ``sweeping duality formalism'' on the other) the formulaire
concerning these four operations
\[f^*, \quad \mathrm Rf_*, \quad {\Lotimes}, \quad \RbHom\]
is no more, no less than just the relevant formulaire in the context
of arbitrary (commutatively) ringed topoi, and maps between such. The
common notation $f_!$ (occurring in $\mathrm Lf_!$ in the \Cat{}
context, in $\mathrm Rf_!$ in the ``duality'' context) is a very
suggestive one, for the least, and I am rather confident that most
reflexes (concerning formal behavior of $f_!$ with respect to the
other operations) coming from one context, should be OK too in the
other. Whether this is just a mere formal analogy, or whether there is
a deeper relationship between the two kinds of contexts, I am at a
loss at present to say. It doesn't seem at all unlikely to me that
among arbitrary maps in \Cat, one can single out some, by some kind of
``finite type'' condition, for which a functor
\begin{equation}
  \label{eq:136.star}
  f^! : \D^+(B\uphat_k) \to \D^+(A\uphat_k)\tag{*}
\end{equation}
can be defined (where $f:A\to B$), right adjoint to the familiar
$\mathrm Rf_*$ functor, so that \eqref{eq:136.2} can be completed to a
sequence of \emph{four} functors
\begin{equation}
  \label{eq:136.3}
  \mathrm Lf_!, \quad f^*, \quad \mathrm Rf_*, \quad f^!,\tag{3}
\end{equation}
forming a sequence of mutually adjoint functors between derived
categories, in the usual sense (the functor immediately to the right
of another being its right adjoint). I faintly remember that Verdier
worked out such a formalism within the context of discrete or
profinite groups, or both, in a Bourbaki talk he gave, this being
inspired by the similar work he did within the context of usual
topological spaces. In the latter, the finiteness condition required
for a map
\[f:X\to Y\]
of topological spaces to give rise to a functor $f^!$ (between, say,
the derived categories of the categories of abelian sheaves on $Y$ and
$X$) is mainly that $X$ should be locally embeddable in a product
$Y\times \bR^d$ -- a rather natural condition indeed! As we are using
\Cat{} as a kind of algebraic paradigm for the category of topological
spaces, this last example for ``sweeping duality'' makes it rather
plausible that something of the same kind should exist indeed in
\Cat{} -- and likewise for the context of groups (which we may view as
just particular cases of\pspage{569} models in \Cat).

I am sorry I was a bit confused, when describing $f^!$ as a right
adjoint to $\mathrm Rf_*$ -- I was thinking of the analogy with a map
of schemes, or spaces, which is not only ``of finite type'' in a
suitable sense, but moreover \emph{proper} -- in which case, in those
duality contexts, $\mathrm Rf_*$ is canonically isomorphic with the
functor denoted by $f_!$ or $\mathrm Rf_!$. Otherwise, the
(``non-trivial'') ``duality theorem'' will assert, rather, that the
pair
\[\mathrm Rf_!,\quad f^!\]
is a pair of adjoint functors, just as is the pair
\[f^*,\quad \mathrm Rf_*\quad\text{(or simply $f_*$).}\]
But even when $f$ is assumed to be proper, the sequence
\eqref{eq:136.3} isn't a sequence of adjoint functors in the standard
duality contexts, namely $\mathrm Rf_!$ is by no means left adjoint to
$f^*$, i.e., $f^*$ isn't isomorphic to $f^!$ (except in extremely
special cases, practically I would think only étale maps are OK, which
in the context of \Cat{} would correspond to maps in \Cat{} isomorphic
to a map $A_{/X}\to A$ for $X$ in \Ahat, namely maps which are
fibering with discrete fibers). This \emph{does} make an important
discrepancy indeed, between the two kinds of contexts -- and increases
my perplexity, as to whether or not one should expect a ``four
variance duality formalism'' to make sense in \Cat. If so, presumably
the $f_!$ or $\mathrm Rf_!$ it will involve (perhaps via a suitable
notion of ``proper'' maps in \Cat, as already referred to earlier
(section~\ref{sec:70})) will be different after all from the $\mathrm
Lf_!$ we have been working with lately, embodying homology
properties. But so does $\mathrm Rf_!$ too, in a rather strong sense,
via the ``duality theorem''!

In the various duality contexts, a basic part is played by the three
particular classes of maps: proper maps, smooth maps, and immersions,
and factorizations of maps into an immersion followed by either a
proper, or a smooth map. In the context of Cat, there is a very
natural way indeed to define the three classes of maps, as we'll see
in the next chapter, presumably -- so natural indeed, that it is hard
to believe that there may be any other reasonable choice! One very
striking feature (already mentioned in section~\ref{sec:70}) is that
the two first notions are ``dual'' to each other in the rather
tautological sense, namely that a map $f:A\to B$ in \Cat{} is proper
if{f} the corresponding map $f\op:A\op\to B\op$ is smooth (whereas the
notion of an immersion is\pspage{570} autodual). How does this fit
with the expectation of developing a ``four variance'' duality
formalism within \Cat? It rather heightens perplexity at first sight!
Proper maps include cofibrations (in the sense of category theory, not
in Kan-Quillen's sense!); dually, smooth maps include
fibrations. Consequently, maps which are moth smooth and proper
include bifibrations, and hence are not too uncommon. Now, how strong
a restriction is it for a map $f:A\to B$ in \Cat{} to factor into
\[f=p\circ i,\]
where $i$ is an immersion $A\to A'$ (namely, a functor identifying $A$
to a full subcategory of $A'$, whose essential image includes, with
any two objects $x,y$, any other $z$ which is ``in between'': $x\to
z\to y$), and $p$ is both proper and smooth (a bifibration, say)?

I start feeling like a battle horse scenting gunpowder again -- still,
I don't think I'll run into it. Surely, there is something to be
cleared up, and perhaps once again a beautiful duality formalism with
the six operations and all will emerge out of darkness -- but this
time I will not do the pulling. Maybe someone else will -- if he isn't
discouraged beforehand, because the big-shots all seem kind of blasé
with ``big duality'', derived categories and all that. As for my
present understanding, I feel that the question isn't really about
homotopy models, or about foundations of homotopy and cohomology
formalism -- at any rate, that I definitely don't need this kind of
stuff, for the program I have been out for. I shouldn't refrain, of
course, to pause on the way every now and then and have a look at the
landscape, however remote or misty -- but I am not going to forget I
am bound for a journey, and that the journey should not be an unending
one\dots

\bigbreak

\presectionfill\alsoondate{27.10.}\par

\hangsection[Looking for the relevant notions of properness and
\dots]{Looking for the relevant notions of properness and smoothness
  for maps in \texorpdfstring{\Cat}{(Cat)}. Case of ordered sets as a
  paradigm for cohomology theory of conically stratified
  spaces.}\label{sec:137}%
It occurred to me that I have been a little rash yesterday, when
asserting that the notions of ``smoothness'' and ``properness'' for
maps in \Cat{} which I hit upon last Spring is the only ``reasonable''
one. Initially, I referred to these notions by the names
``cohomologically smooth'', ``cohomologically proper'', as a measure
of caution -- they were defined by properties of commutation of base
change to formation of the Leray sheaves $\mathrm R^if_*$ (i.e.,
essentially, to ``cointegration''), which were familiar to me for
smooth resp.\ proper maps in the context of schemes, or ordinary
topological spaces. These\pspage{571} cohomological counterparts of
smoothness and properness fit very neatly into the homology and
cohomology formalism, and I played around enough with them, last
Spring as well as more than twenty years ago when developing étale
cohomology of schemes, for there being no doubt left in my mind that
these notions are relevant indeed. However, I was rather rash
yesterday, while forgetting that these cohomological versions of
smoothness and properness are considerably weaker than the usual
notions. Thus, in the context of schemes over a ground field, the
product of any two schemes is cohomologically smooth over its factors
-- or equivalently, any scheme over a field $k$ is cohomologically
smooth over $k$! Similarly, in the context of \Cat, any object in
\Cat, namely any small category, is both cohomologically smooth and
proper over the final object $e$ (as it is trivially ``bifibered''
over $e$). On the other hand, it isn't reasonable, of course, to
expect any kind of Poincaré-like duality to hold for the cohomology
(with twisted coefficients, say) of an arbitrary object $A$ in
\Cat. To be more specific, it is easy to see that in many cases, the
functor
\[\mathrm Rf_* : \D^+(\Ahatab)\to\D^+\Ab\]
does not admit a right adjoint (which we would call $f^!$). For
instance, when $A$ is discrete, then $f^!$ exists (and may then be
identified with $f^*$) if{f} $A$ is moreover \emph{finite} -- a rather
natural condition indeed, when we keep in mind the topological
significance of the usual notion of properness! This immediately
brings to mind some further properties besides base change properties,
which go with the intuitions around properness -- for instance, we
would expect for proper $f$, the functors $f_*$ and $\mathrm Rf_*$ to
commute to filtering direct limits, and the same expectation goes with
the assumption that $\mathrm Rf_*$ should admit a right adjoint. This
exactness property is not satisfied, of course, when $A$ is discrete
infinite. We now may think (still in case of target category equal to
$e$) to impose the drastic condition that the category $A$ is
finite. Such restriction however looks in some respects too weak, in
others too strong. Thus, it will include categories defined by finite
groups, which goes against the rather natural expectation that
properness + smoothness, or any kind of Poincaré duality, should go
with \emph{finite cohomological dimension}. On the other hand, there
are beautiful infinite groups (such as the fundamental group of a
compact surface, or of any other compact variety that is a $K(\pi,1)$
space\dots) which satisfy Poincaré duality.

These\pspage{572} reflections make it quite clear that the notions of
properness and of smoothness for maps in \Cat, relevant for a duality
formalism, have still to be worked out. Two basic requirements to be
kept in mind are the following: \namedlabel{rem:137.1}{1)}\enspace for
a proper map $f: A\to B$, and any ring of coefficients $k$, the
functor
\[\mathrm Rf_* : \D^+(A\uphat_k)\to\D^+(B\uphat_k)\]
should admit a right adjoint $f^!$, and
\namedlabel{rem:137.2}{2)}\enspace for a smooth map $f$ factored as
$f=g\circ i$, with $g$ proper and $i$ an ``open immersion'', the
composition $f^!\eqdef i^*\circ g^!$ should be expressible as
\begin{equation}
  \label{eq:137.5}
  f^!: K^\bullet \mapsto f^*(K)\Lotimes_k T_f(k)[d_f],\tag{5}
\end{equation}
where $T_f(k)$ is a presheaf of $k$-modules on $A$ locally isomorphic
to the constant presheaf $k_A$ ($T_f(k)$ may be called the
\emph{orientation sheaf} for $f$, with coefficients in $k$), and $d_f$
is a natural integer (which may be called the \emph{relative
  dimension} of $A$ over $B$, or of $f$). (For simplicity, I assume in
\ref{rem:137.2} that $A$ is connected, otherwise $d_f$ should be
viewed as a function on the set of connected components of
$A$.)\enspace This again should give the correct relationship, for $f$
as above, between the (for the time being hypothetical) $\mathrm Rf_!$
($\eqdef \mathrm Rg_*\circ\mathrm Lf_!$) and our anodyne $\mathrm
Lf_!$, for an argument $L_\bullet$ in $\D^\bullet(A\uphat_k)$ say:
\begin{equation}
  \label{eq:137.6}
  \mathrm Rf_!(L_\bullet) \simeq \mathrm Lf_!(L_\bullet \otimes
  T_f^{-1})[-d_f],\tag{6} 
\end{equation}
where the left-hand side is just $\mathrm Rf_*(L_\bullet)$, if we
assume moreover $f$ to be proper.

This precise relationship between the two possible versions of an
$f_!$ operation between derived categories, namely $\mathrm Lf_!$
embodying homology, defined for any map $f$ in \Cat, and
$\mathrm Rf_!$ embodying ``cohomology with proper supports'', defined
for a map that may be factored as $g\circ i$ with $g$ ``proper'' and
$i$ an (open, if we like) immersion, relation valid if $f$ is moreover
assumed to be ``smooth'', greatly relaxes yesterday's perplexity,
coming from a partial confusion in my mind between the operations
$\mathrm Rf_!$ and $\mathrm Lf_!$. (Beware the notation $\mathrm Rf_!$
is an abuse, as it doesn't mean at all anything like the right derived
functor of the functor $f_!$!)\enspace At the same time, I feel a lot
less dubious now about the existence of a ``six operations'' duality
formalism in the \Cat{} context -- I am pretty much convinced, now,
that such a formalism exists indeed. The main specific work ahead is
to get hold of the relevant notions of proper and smooth maps. The
demands we have on\pspage{573} these, besides the relevant base change
properties, are so precise, one feels, that they may almost be taken
as a definition! Maybe even the ``almost'' could be dropped -- namely
that a comprehensive axiomatic set-up for the duality formalism could
be worked out, in a way applicable to the known instances as well as
to the presently still unknown one of \Cat, by going a little further
still than Deligne's exposition in
SGA~5\scrcomment{\textcite[Dualité]{SGA4andhalf}} (where the notions
of ``smooth'' and ``proper'' maps were supposed to be given
beforehand, satisfying suitable properties). Before diving into such
axiomatization game, one should get a better feeling, though, through
a fair number of examples (not all with $e$ as the target category
moreover), of how the proper, the smooth and the proper-and-smooth
maps in \Cat{} actually look like. Here, presumably, Verdier's work in
the context of discrete infinite groups should give useful clues.

Other important clues should come from the opposite side so to say --
namely ordered sets. Such a set $I$, besides defining in the usual way
a small category and hence a topos, may equally be viewed as a
topological space, more accurately, the topos it defines may be viewed
as being associated to a topological space, admitting $I$ as its
underlying set (cf.\ section~\ref{sec:22}, p.~\ref{p:18}). This
topological space is noetherian if{f} the ordered set $I$ is -- for
instance if $I$ is finite. In such a case, an old algebraic geometer
like me will feel in known territory, which maybe is a delusion,
however -- at any rate, I doubt the duality formalism for topological
spaces (using factorizations of maps $X\to Y$ via embeddings in spaces
$Y\times\bR^d$) makes much sense for such non-separated
spaces. However, as we saw in section~\ref{sec:22}, when $I$ satisfies
some mild ``local finiteness'' requirement (for instance when $I$ is
finite), we may associate to it a \emph{geometrical realization} $\abs
I$, which is a locally compact space (a compact one indeed if $I$ is
finite) endowed with a ``conical subdivision'' (index by the opposite
ordered set $I\op$), hence canonically triangulated via the
``barycentric subdivision''. The homotopy type of this space is
canonically isomorphic to the homotopy type of $I$, viewed as a
``model'' in \Cat. What is more important here, is that a (pre)sheaf
of sets (say) on the category $I$ may be interpreted as being
essentially the same as a sheaf of sets on the geometric realization
$\abs I$ \emph{which is locally constant} (\emph{and hence constant})
\emph{on each of the \emph{``open''} strata or ``cones'' of $I$}. This
description then carries over to sheaves of $k$-modules. The ``clue''
I had in mind is that \emph{within the context of locally finite
  ordered sets, the looked-for ``six operations duality formalism''
  should be no more, no less than the accurate reflection of the same
  formalism within the context\pspage{574} of \textup(locally
  compact\textup) topological spaces}, as worked out by Verdier in one
of his Bourbaki talks -- it being understood that when applying the
latter formalism to spaces endowed with conical stratifications, maps
between these which are compatible with the stratifications (in a
suitable sense which should still be pinned down), and to sheaves of
modules which are compatible with the stratifications in the sense
above, these will give rise (via the ``six operations'') to sheaves
satisfying the same compatibility. (NB\enspace when speaking of
``sheaves'', I really mean complexes of sheaves $K^\bullet$, and the
compatibility condition should be understood as a condition for the
ordinary sheaves of modules $\mathrm H^i(K^\bullet)$.)\enspace This
remark should allow to work out quite explicitly, in purely algebraic
(or ``combinatorial'') terms, the ``six operations'' in the context of
locally finite ordered sets, at any rate.

This interpretation suggests that an ordered set $I$ should be viewed
as a ``proper'' object of \Cat{} if{f} $I$ is finite. In the same
vein, whenever $I$ is ``locally proper'', namely locally finite, and
moreover its topological realization $\abs I$ is a ``smooth''
topological space in the usual sense, namely is a topological
\emph{variety} (for which it is enough that for every $x$ in $I$, the
topological realization
\begin{equation}
  \label{eq:137.7}
  \abs{I_{>x}}\tag{7}
\end{equation}
of the set of elements $y$ such that $y>x$ should be a sphere), we
would consider $I$ as a ``smooth'' object in \Cat. If $A$ is any
object in \Cat{} and $I$ is an ordered set which is finite resp.\
``smooth'' in the sense above, we will surely expect $A\times I\to A$
to be ``proper'' resp.\ ``smooth'' for the duality formalism we wish
to develop in \Cat.

It should be kept in mind that for the algebraic interpretation above
to hold, for sheaves on a conically stratified space $\abs I$ in terms
of an ordered set $I$, we had to take on the indexing set $I$ for the
strata the order relation \emph{opposite} to the inclusion relation
between (closed) strata -- otherwise, the correct interpretation of
sheaves constant on the open strata is via \emph{covariant} functors
$I\to\Sets$, i.e., \emph{presheaves on $I\op$} (not $I$). At any rate,
as there is a canonical homeomorphism
\begin{equation}
  \label{eq:137.8}
  \abs I\simeq\abs{I\op}\tag{8}
\end{equation}
respecting the canonical barycentric subdivisions of both sides,
notions for $I$ (such as properness, or smoothness) which are
expressed as intrinsic properties of the corresponding topological
space $\abs I$ (independently of its subdivision) are autodual -- they
hold for $I$ if{f} they do for $I\op$.\pspage{575} This is in sharp
contrast with the more naive notions of ``cohomological'' properness
and smoothness via base change operations, which are interchanged by
duality (which was part of yesterday's perplexities, now straightening
out\dots).

I feel I should be a little more outspoken about the relevant notion
of ``proper maps'' between ordered sets, which should be the algebraic
expression of the geometric notion alluded to above, of a map between
conically stratified spaces being ``compatible with the conical
stratifications''. In order for the corresponding direct image functor
for sheaves to take sheaves compatible with the stratification above
to sheaves of same type below, we'll have to insist that the image by
$f$ of a strata above should be a whole strata below. Now, it is clear
that any map
\begin{equation}
  \label{eq:137.9}
  f:I\to J\tag{9}
\end{equation}
between ordered sets takes flags into flags, and hence induces a map
\begin{equation}
  \label{eq:137.10}
  \abs f:\abs I \to \abs J\tag{10}
\end{equation}
between the geometric realizations, compatible with the barycentric
triangulations and thus taking simplices into simplices. But even when
$I$ and $J$ are finite (hence ``proper''), the latter map does not
always satisfy the condition above. Thus, when $I$ is reduced to just
one point, hence \eqref{eq:137.9} is defined by the image $j\in J$ of
the latter, the corresponding map \eqref{eq:137.10} maps the unique
point of $\abs I$ to the barycenter of the stratum $C_j$ in $\abs J$,
which is a stratum of $\abs J$ if{f} $C_j$ is a minimal stratum, i.e.,
$j$ is a \emph{maximal} element in $J$. A natural algebraic condition
to impose upon $f$, in order to ensure that $\abs f$ satisfy the
simple geometric condition above, is that for any $x$ in $I$, and any
$y'$ in $J$ such that
\[y'>y \eqdef f(x),\]
there should exist an $x'$ in $I$ satisfying
\begin{equation}
  \label{eq:137.11}
  x'>x \quad\text{and}\quad f(x')=y'.\tag{11}
\end{equation}
This condition strongly reminds us of the valuative criterion for
properness in the context of schemes, where the relation $y'\ge y$ or
$y\to y'$, say, should be interpreted as meaning that $y'$ is a
\emph{specialization} of $y$. However, in the valuative criterion for
properness (for a map of preschemes of finite type over a noetherian
prescheme say), if one wants actual properness\pspage{576} indeed
(including separation of $f$, not just that $f$ is universally
closed), one has to insist that the $x'$ above should be unique: every
specialization $y'$ of $y=f(x)$ lifts to a \emph{unique}
specialization $x'$ of $x$. If we applied this literally in the
present context, this would translate into the condition that $\abs f$
should map \emph{injectively} each closed stratum $C_x$ of $\abs I$ --
which would exclude such basic maps as the projection $I\to e$ to just
one point!

One may wonder why trouble about the analogy with algebraic geometry
and any extra condition on $f$ besides the one we got. The point,
however, is that we would like the map
\begin{equation}
  \label{eq:137.12}
  C_x=\abs{I_{\ge x}} \to C_y=\abs{I_{\ge y}}\tag{12}
\end{equation}
between corresponding strata induced by $\abs f$ to be
``cohomologically trivial'' in a suitable sense, not only surjective
-- in analogy,\footnote{analogy: aspheric fibers! \scrcommentinline{I
    \emph{think} that's what the footnote says; it's hard to
    read\dots}} say, with the usual notion of maps between
triangulated spaces; if we don't have some condition of this type, we
will have no control over the structure of the direct image and the
higher direct images $\mathrm R^i\abs f_*$ of an ``admissible'' sheaf
upstairs. I didn't really analyze the situation carefully, in terms of
what we are after here in the context of those stratified geometric
realizations. I feel pretty sure, though, that there the general
notion of ``cohomological properness'' which I worked out last Spring
fits in just right, to give the correct answer. The criterion I
obtained (necessary and sufficient for the relevant compatibility
property of $\mathrm R^if_*$ with arbitrary base change $J'\to J$ in
\Cat) reads as follows: let
\begin{equation}
  \label{eq:137.13}
  I(f,x,x') =\set{x'\in I}{\text{$x'>x$ and $f(x')=y'$}}\tag{13}
\end{equation}
be the subset of $I$ satisfying the conditions \eqref{eq:137.11}
above. Instead of demanding only that this set be non-empty, or going
as far as demanding that it should be reduced to just one point, we'll
demand that this set should be \emph{aspheric}. (When we are concerned
with cohomology with commutative coefficients only, presumably it
should be enough to demand only that this set be \emph{acyclic} --
which should be enough for the sake of a mere duality
formalism\dots)\enspace In more geometric terms, this should mean, I
guess, that the inverse image, by the map \eqref{eq:137.12} above
between closed strata, of any closed stratum below, should be
aspheric. I doubt this condition holds under the mere assumption that
the sets \eqref{eq:137.13} be non-empty, i.e., \eqref{eq:137.12} be
surjective -- I didn't sit down, though, to try and make an example.

Thus,\pspage{577} we see that when working with the notion of
properness of maps, even for such simple gadgets as finite ordered
sets, which should be viewed as ``proper'' (or ``compact'') objects by
themselves, this notion is far from being a wholly trivial one -- for
instance, it does not hold true that any map between such ``proper''
objects is again ``proper''. This now seems to me just a mathematical
``fact of life'', which we may not disregard when working with finite
ordered sets, say, in view of expressing in algebraic terms some
standard operations in the cohomology theory of sheaves on topological
spaces, endowed with suitable conical stratifications. The fault,
surely, is not with the notion of conical stratification itself, which
may be felt by some as being ad hoc, awkward and what not. I know the
notion is just right -- but at any rate, even when working with strata
which are perfect topological cells (so that nobody could possibly
object to them), exactly the same facts of life are there -- not every
map between such cellular decompositions, mapping cells \emph{onto}
cells, will fit into a ``combinatorial'' description when it comes to
describing the standard operations of the cohomology of sheaves, for
sheaves compatible with the stratifications\dots

To sum up, it seems to me that definitely, there is a very rich
experimental material available already at present, in order to come
to a feeling of what duality is like in the context of \Cat, and for
developing some of the basic intuitions needed for working out,
hopefully, ``the'' full-fledged duality theory which should hold in
\Cat. Before leaving this question, I would like to point out still
one other property connected with the intuitions around ``properness''
-- more generally, around maps ``of finite type'' in a suitable sense,
which sometimes may translate into: factorizable as a composition
$f=g\circ i$, with $g$ proper and $i$ an immersion. This is about
stability of ``constructibility'' or ``finiteness'' conditions for
(complexes of) sheaves of $k$-modules, with respect to the standard
operations $\mathrm Rf_*$, $\mathrm Rf_!$, $f^!$ (stability by $f^*$
being a tautology in any case). Such stability of course, whenever it
holds, is an important feature, for instance for making ``virtual''
calculations in suitable ``Grothendieck groups'' (where Euler-Poincaré
type invariants may be defined). It should be recalled, however, that
the six fundamental operations in the duality formulaire, as well as
nearly all of the formulaire itself, make sense (and formulæ hold
true) without any finiteness conditions on the complexes of sheaves we
work with, except just boundedness conditions on the degrees of those
complexes.

\begin{remarks}
  1)\enspace One\pspage{578} of the (rather few) instances in the
  duality formulaire where finiteness conditions are clearly needed,
  is the so-called ``biduality theorem'', when taking $\RbHom$'s with
  values in a so-called ``dualizing complex'', which here should be an
  object
  \[R^\bullet\quad\text{in $\D^{\mathrm b}(A\uphat_k)$}\]
  (for given $A$ in \Cat{} and given ring of coefficients $k$, for
  instance $k=\bZ$). The question arises here, for any given $A$ and
  $k$, a)\enspace whether there exists a dualizing complex $R^\bullet$
  (which, as usual, will be unique up to dimension shift and ``twist''
  by an invertible sheaf of $k$-modules on $A$), and b)\enspace can
  such dualizing complex be obtained as
  \[R^\bullet=p^!(k_e),\]
  where $p:A\to e$ is the projection to the final object of \Cat? It
  should be easy enough, once the basic duality formalism is written
  down for ordered sets, as contemplated above, to show that when $A$
  is a finite ordered set (or only locally finite of finite
  combinatorial dimension), then $p^!(e_k)$ is indeed a dualizing
  complex. To refer to something more ``en vogue'' at present than
  those poor dualizing complexes, it is clear, from what I heard from
  Illusie and Mebkhout about the ``complexe d'intersection'' for
  stratified spaces, that this complex can be described also in the
  context of locally finite ordered sets (if I got the stepwise
  construction of this complex right). As the construction here
  corresponds to stratifications where the strata are by no means
  even-dimensional, I am not too sure, though, if the complex obtained
  this way is really relevant -- it isn't a topological invariant of
  the topological space at any rate, independently of its
  stratification -- a bad point indeed. Too bad!

  2)\enspace Here is a rather evident example showing that the
  asphericity condition on \eqref{eq:137.13} is not automatic. Take
  $I$ with a smallest element $x$ (hence $C_x=\abs I$), and
  $J=\Simplex_1=(0\to 1)$. To give a map $I\to J$ amount to give
  $I_0=f^{-1}(0)\subset I$, which is any open subset of $I$ (i.e.,
  such that $a\in I_0$, $b\le a$ implies $b\in I_0$). If we take
  $I_0=\{x\}$, then properness of $f$ is equivalent with
  $I_1=I\setminus I_0=I\setminus\{x\}$ being aspheric, while the
  weaker condition contemplated first (the sets \eqref{eq:137.13}
  non-empty) means only that $I_1\ne\emptyset$. Now, $I_1$ may of
  course be taken to be any (finite say) ordered set -- the
  construction made amounting to taking the cone over the map
  $\abs{I_1}\to e$. In this examples, all fibers of $\abs f:\abs I\to
  \abs J$ are (canonically) homeomorphic to $I_1$, except the fiber at
  the ``barycenter'' $1$, which\pspage{579} is reduced to a point --
  visibly not a very ``cellular'' behavior when $\abs{I_1}$ isn't
  aspheric! At any rate, as the sheaves $\mathrm R^i\abs f_*(F)$ (for
  $F$ a constant sheaf above, say) may be computed fiberwise, we see
  that if $I_1$ isn't aspheric, these sheaves (which are constant on
  $\abs J\setminus\{1\}$) are \emph{not} going to be constant on $\abs
  J\setminus\{0\}$, as they should if we want an algebraic paradigm of
  operations like $\mathrm R\abs f_*$ in terms of sheaves on finite
  ordered sets.

  In this example we could take $\abs I$ to be a perfect $n$-cell
  ($n\ge1$), hence $I_1$ is an $(n-1)$-sphere, whereas $\abs J$ isn't
  really a (combinatorial) $1$-cell, as its boundary has just
  \emph{one} point $0$, instead of two. The $1$-cell structure
  corresponds to the ordered set (\emph{opposite} to the ordered set
  formed by the two endpoints and the dimension $1$ stratum)
  \[J =
    \begin{tikzcd}[row sep=tiny,column sep=small,cramped]
      & y\ar[dl]\ar[dr] & \\
      0 & & 1
    \end{tikzcd}.\]
  For any ordered set $I$, to give a map $f:I\to J$ amounts to giving
  the two subsets
  \[I_0=f^{-1}(0), \quad I_1=f^{-1}(1),\]
  subject to the only condition of being open and disjoint. In terms
  of strata of $\abs I$, this means that we give two sets $I_0,I_1$ of
  strata, containing with any stratum any smaller one, and having no
  stratum in common. The condition that the sets \eqref{eq:137.13}
  should be non-empty says that any point in $I$ which is neither in
  $I_0$ nor $I_1$ admits majorants in both -- or geometrically, any
  stratum which isn't in $I_0$ nor in $I_1$ meets both $\abs{I_0}$ and
  $\abs{I_1}$, i.e., contains strata which are in $I_0$ and strata
  which are in $I_1$. Even when $\abs I$ is a combinatorial $2$-cell,
  i.e., a polygonal disc, this condition does not imply asphericity of
  the sets \eqref{eq:137.13} (not even $0$-connectedness). To see
  this, we take the set \eqref{eq:137.13} where $x$ is the dimension
  $2$ stratum, i.e., the smallest element of $I$, mapped to $y$ (the
  smallest in $J$), and $y'$ either $0$ or $1$. The reader will easily
  figure out on a drawing the structure of the map $\abs f$, as I just
  did myself: the fibers at the endpoints of the segment $\abs J$ are,
  as given, discrete with cardinal $m$, the fiber at a point different
  from the endpoints and from their barycenter are disjoint sums of
  $m$ segments (hence homotopic to the former fibers), whereas the
  fiber at the barycenter is the union of $m$ segments meeting in
  their\pspage{580} common middle, hence is contractible. Thus, the
  $\mathrm R^0\abs f_*$ of a constant sheaf on the disc $\abs I$ is by
  no means constant on the open, dimension one stratum of $\abs J$.

  These examples bring to my mind that for any map $f:I\to J$ between
  finite ordered sets, possibly submitted to the mild restriction that
  the sets \eqref{eq:137.13} should be non-empty, the homotopy types
  of those ordered sets \eqref{eq:137.13} (for the order relation
  induced by $I$) should be exactly the homotopy types of the fibers
  of the maps \eqref{eq:137.12} between strata. Thus, asphericity of
  these ordered sets should express no more, no less than the
  contractibility of those fibers. This latter condition is exactly
  what is needed in order to ensure stability by $\mathrm R\abs f_*$
  of the notion of complexes of sheaves compatible with the
  stratifications. 
\end{remarks}

\bigbreak

\presectionfill\ondate{28.10.}\par

\hangsection[Niceties and oddities: $\mathrm Rf_!$ commutes to
\dots]{Niceties and oddities: \texorpdfstring{$\mathrm Rf_!$}{Rf!}
  commutes to colocalization, not localization.}\label{sec:138}%
Yesterday and the day before, I got involved in that unforeseen
digression around the foreboding of a ``six operations duality
formalism'' in \Cat, and suitable notions for smoothness and (more
important still) of properness for a map in \Cat. This digression
wholly convinced me that the usual duality formalism should hold in
\Cat{} too. Working this out in full detail should be a most pleasant
task indeed, and presumably the best, or even the only way for gaining
complete mastery of the cohomology formalism within \Cat{} or, what
more or less amounts to the same, for topoi which admit sufficiently
many projective objects. In the previous two sections, I referred to
such duality formalism as one concerned with sheaves of $k$-modules,
for any fixed ring $k$ -- but from the example of étale duality for
schemes, say, it is likely that instead of fixing a ring $k$, we may
as well take objects $A$ in \Cat{} endowed with an arbitrary sheaf of
rings $\scrO_A$ (which we'll only have to suppose commutative when
concerned with the two internal operations $\Lotimes$ and $\RbHom$),
and taking maps of such ringed objects as the basic maps. In the
present context, the usual ``six operations'' in duality theory will
be enriched, however, by still another one, namely $\mathrm Lf_!$,
defined for any map $f$ between such ringed objects (not to be
confused with the $\mathrm Rf_!$ operation, defined only under
suitable finiteness assumptions, such as ``properness'', for the
underlying map in \Cat), whose relationship to the other operations
should be understood and added to the standard duality
formulaire. One\pspage{581} such formula, namely the precise relation
between $\mathrm Lf_!$ and $\mathrm Rf_!$ for a smooth $f$, was given
yesterday (p.~\ref{p:572} \eqref{eq:137.6}).

The day before, I was out for trying to get a better understanding of
the $\mathrm Rf_!$ operation (including the case of non-constant
sheaves of rings for the modules we work with). This operation still
remains unfamiliar to me, very unlike my old friend $\mathrm Rf_*$ --
there is a number of things about it which are not quite clear yet in
my mind, even when just taking the functor $f_!$ between modules,
before taking a left derived functor. For instance, for general rings
of operators $\scrO_A$ and $\scrO_B$, when $f$ is a ringed map
\begin{equation}
  \label{eq:138.1}
  f: (A,\scrO_A) \to (B,\scrO_B),\tag{1}
\end{equation}
it doesn't seem that formation of $f_!$ commutes to restriction of
rings of operators (to the constant rings $\bZ_A$ and $\bZ_B$, say),
namely that for a given $\scrO_A$-module $F$, $f_!(F)$ may be
interpreted as just ${f_0}_!\supab(F)$ with suitable operations of
$\scrO_B$ on the latter, where
\begin{equation}
  \label{eq:138.2}
  f_0:A\to B\tag{2}
\end{equation}
is the map in \Cat{} underlying $f$. When reducing to a suitable
``universal'' case, this may be expressed by saying that for given map
$f_0$ and given abelian sheaf $F$ on $A$, if we define
\[\scrO_A = \bEnd_\bZ(F), \quad \scrO_B=f_*(\scrO_A),\]
there doesn't seem to be a natural operation of $\scrO_B$ upon
${f_0}_!\supab(F)$; all we can say, it seems, is that the ring of
global sections of $\scrO_B$ operates on ${f_0}_!\supab(F)$. More
generally, reverting to the general case of a map $f$ of ringed
objects in \Cat, the ring of global sections of $\scrO_B$ operates on
${f_0}_!\supab(F)$. When $\scrO_B$ is a constant sheaf of rings $k_B$,
this implies of course that $k$ operates on this abelian sheaf on $B$,
from which will follow by an obvious argument that with this structure
of a sheaf of $k$-modules, ${f_0}_!\supab(F)$ may indeed be identified
with $f_!(F)$.

Going over to $\mathrm Lf_!$ which we would like to express via
$\mathrm L{f_0}_!\supab$, the situation is worse, as we still have to
check (granting $\scrO_B$ to be constant) that for $F$ a projective
$\scrO_A$-module, we got
\begin{equation}
  \label{eq:138.3}
  \mathrm L_i{f_0}_!\supab(F)=0\quad\text{for $i>0$.}\tag{3}
\end{equation}
This isn't always true, even when $B$ is the final object in \Cat, and
$A$ has a final object, hence $\scrO_A$ is a projective module
over\pspage{582} itself, and \eqref{eq:138.3} reads
\[\mathrm H_i(A,\scrO_A) =0 \quad\text{for $i>0$,}\]
which isn't always true. (NB\enspace if it were for any commutative
ring $\scrO_A$ on $A$, it would be too for any abelian sheaf $M$ on
$A$, as we see by taking $\scrO_A=\bZ_A\otimes M$, hence $A$ would be
homological dimension $0$, a drastic restriction, indeed, even when
$A$ has a final object.)

When however $\scrO_A$ is equally a constant sheaf of rings, say
$\scrO_A=k'_A$, then the relation \eqref{eq:138.3} holds for any
projective module on $A$. We need only check it for
\[F=\scrO_A^{(a)} = {k'}^{(a)},\]
with $a$ in $A$, then \eqref{eq:138.3} follows from
\begin{equation}
  \label{eq:138.4}
  \mathrm L_i {f_0}_!\supab( M^{(a)} )=0 \quad
  \text{for $i>0$, $M$ in \Ab,}\tag{4}
\end{equation}
To check \eqref{eq:138.4}, we take a projective resolution of
$M^{(a)}$, by using a projective resolution $L_\bullet$ of $M$ in
\Ab{} and taking $L_\bullet^{(a)}$ (using the fact that the functor
\[L\mapsto L^{(a)}:\Ab\to\Ahatab\]
is exact). We then get
\[\mathrm L{f_0}_!(M^{(a)}) = {f_0}_! (L_\bullet^{(a)}) =
  L_\bullet^{(b)} \simeq M^{(b)},\]
where $b=f_0(a)$, and where the last equality stems from exactness of
$L\mapsto L^{(b)}$. (NB\enspace the relation \eqref{eq:138.4}
generalizes a standard acyclicity criterion in the homology theory of
discrete groups\dots)

Thus, we get finally that in case both rings $\scrO_A$ and $\scrO_B$
are constant, that formation of $\mathrm Lf_!$ commutes to restriction
of operator rings (provided the rings to which we are restricting are
constant too -- say they are just the absolute $\bZ_A$ and
$\bZ_B$). Presumably, a little more care should show the similar
result for locally constant rings.

Reverting to the case \eqref{eq:138.1} of a general map between ringed
objects in \Cat, our inability, for a given $\scrO_A$-module $F$ on
$A$, to define an operation of $\scrO_B$ upon ${f_0}_!\supab$ (while
there is an operation of the ring $\Gamma(B,\scrO_B)$ upon it), is
tied up with this difficulty, that formation of ${f_0}_!\supab$, and a
fortiori of $f_!$ for arbitrary rings $\scrO_A$ and $\scrO_B$, does
not commute to ``localization'' (as $f_*$ and $\mathrm Rf_*$ does),
namely to base change of the type\pspage{583}
\begin{equation}
  \label{eq:138.5}
  B_{/b}\to B,\tag{5}
\end{equation}
where $b$ is a given object in $B$, and $B_{/b}$ designates as usual
the category of all ``objects over $b$'' in $B$, i.e., of all arrows
in $B$ with target $b$. This gives rise to the cartesian square in
\Cat{}
\begin{equation}
  \label{eq:138.6}
  \begin{tabular}{@{}c@{}}
    \begin{tikzcd}[baseline=(O.base)]
      A_{/b}\ar[d]\ar[r] & A\ar[d] \\
      B_{/b}\ar[r] & |[alias=O]| B
    \end{tikzcd},
  \end{tabular}\tag{6}
\end{equation}
where $A_{/b}$ is the category of all pairs
\[\text{$(a,u)$ with $u:f_0(a)\to b$,}\]
which may be identified equally with the category $A_{/f_0^*(b)}$. The
commutation property we have in mind is a tautology for $f_*$, and it
follows for $\mathrm Rf_*$, because the inverse image by $A_{/b}\to A$
of an injective module on $A$ is an injective module on $A_{/b}$. The
latter fact is true, more generally, for any ``localization map'', of
the type $A_{/X}\to A$, with $X$ in \Ahat, i.e., any map which is
fibering with discrete fibers. Thus, the commutation property for
$\mathrm Rf_*$ is valid more generally for any base change of the type
\[B_{/Y}\to B,\]
with $Y$ in \Bhat{} (not necessarily in $B$). More generally still, it
can be shown to hold for any map
\[B'\to B\]
which is fibering (not necessarily with discrete fibers).

On the other side of the mirror, when taking $f_!$ and its left
derived functor, already the former definitely does \emph{not} commute
to base change of the type \eqref{eq:138.5}, i.e., to ``localization
on the base'' -- something a little hard to get accustomed to! The
base changes which will do here, are those of the dual type
\begin{equation}
  \label{eq:138.5prime}
  \preslice Bb\to B,\tag{5'}
\end{equation}
we may call them maps of ``colocalization'' on the base $B$. It gives
rise to a cartesian diagram in \Cat{} dual to \eqref{eq:138.6}
\begin{equation}
  \label{eq:138.6prime}
  \begin{tabular}{@{}c@{}}
    \begin{tikzcd}[baseline=(O.base)]
      \preslice Ab\ar[d]\ar[r] & A\ar[d] \\
      \preslice Bb\ar[r] & |[alias=O]| B
    \end{tikzcd},
  \end{tabular}\tag{6'}
\end{equation}
where\pspage{584} now $\preslice Ab$ is the category of pairs
\[\text{$(a,u)$ with $u:b\to f_0(a)$.}\]
We'll have to assume now that both shaves of rings $\scrO_A$,
$\scrO_B$ are constant, and correspond to the same ring $k$. Thus,
modules on $A$ and $B$ are just contravariant functors from these
categories to the category $\AbOf_k$ of $k$-modules, and accordingly,
$f_!$ (left adjoint to the composition functor $f^*$) may be computed
by a well-known formula, involving direct limits on the categories
$\preslice Ab$:
\begin{equation}
  \label{eq:138.7}
  f_!(F)(b) \simeq\varinjlim_{\preslice Ab} F(a),\tag{7}
\end{equation}
where the limit in the second member is relative to the
composition\scrcomment{It seems to me that we're forgetting
  that $F$ is contravariant. Should be easily fixed by inserting ${}\op$'s,
  though\dots}
\[\preslice Ab\to A\xrightarrow F{} \AbOf_k.\]
(NB\enspace This formula is dual to the formula for $f_*$, the right
adjoint of $f^*$, closer to intuition -- to mine at any rate --
because the $\varprojlim$ in
\[f_*(F)(b) \simeq \varprojlim_{A_{/b}} F(a)\]
may be ``visualized'' as the set of ``sections'' of $F$ over
$A_{/b}$.)\enspace From this formula \eqref{eq:138.7} follows at once
commutation of $f_!$ with colocalization. To get the corresponding
result for $\mathrm Lf_!$, we have only to use the fact that the
relevant inverse image functor (corresponding to $\preslice Ab\to A$)
transforms projective modules into projective modules. As the map
$\preslice Bb\to B$ is a \emph{cofibering functor with discrete
  fibers}, so is the map
\[h:A'=\preslice Ab\to A\]
deduced by base change. Now, for any such functor $h:A'\to A$ between
small categories, the inverse image functor
\[h^*:\Ahat \to {A'}\uphat, \quad\text{or}\quad
  h_k^*:A\uphat_k \to {A'}\uphat_k\]
carries indeed projectives into projectives. This statement is dual
formally to the corresponding statement for injectives, valid when we
make on $h$ the dual assumption of being \emph{fibering with discrete
  fibers} -- in the latter case the (well-known) proof comes out
formally from the fact that the left adjoint functor $h_!$ or $h_!^k$
carries monomorphisms into monomorphisms -- a fact that we used in
section~\ref{sec:135} in the form $\mathrm Lh_!^k=h_!^k$ (in case
$k=\bZ$). In the present case, the proof is essentially the dual one
-- as a matter of fact, as there are enough\pspage{585} projectives,
the statement about $h^*$ or $h_k^*$ taking projectives into
projectives is equivalent with the right adjoint $h_*$ or $h_*^k$ (for
sheaves of sets, resp.\ sheaves of $k$-modules) transforming
epimorphisms into epimorphisms, which in the case of $k$-modules can
be written equally under the equivalent form
\begin{equation}
  \label{eq:138.8}
  \mathrm R^ih_*=0 \quad\text{for $i>0$.}\tag{8}
\end{equation}
Now, this exactness property for $h_*$, in the case when $h$ is
cofibering with discrete fibers, is I guess well known (it is
well-known to me at any rate), and comes from the specific computation
of $h_*(F)$ for $F$ in ${A'}\uphat$, valid whenever $h$ is
\emph{cofibering} (with arbitrary fibers)
\begin{equation}
  \label{eq:138.9}
  h_*(F)(x) \simeq \Gamma(A'_x, F\restrto A'_X)\quad
  \text{for $x$ in $A$,}\tag{9}
\end{equation}
where $A'_X$ is the fiber of $A'$ over $x$ (a category not to be
confused with $A'_{/x}$, the two being closely related,
however\dots). In case the fibers of $h$ are discrete, the right-hand
side of \eqref{eq:138.9} may be written as a product, hence the
formula
\begin{equation}
  \label{eq:138.10}
  h_*(F)(x) \simeq \prod_{\text{$x'$ in $A'_X$}} F(x'),\tag{10}
\end{equation}
(which may be viewed as the dual of the formula \eqref{eq:135.11}
p.~\ref{p:564}). As in the category of sets (and hence also in
$\AbOf_k$) a product of epimorphisms is again an epimorphism, the
result we want follows indeed.

\begin{remarks}
  1)\enspace The results just given, as well as their proofs,
  concerning inverse images of injectives or projectives, are valid
  not only in the case of a common constant sheaf of rings on $A$ and
  $A'$, but more generally for any sheaf of rings $\scrO_A$ on $A$,
  when taking on $A'$ the ``induced'' sheaf of rings
  \begin{equation}
    \label{eq:138.11}
    \scrO_{A'} = h^*(\scrO_A).\tag{11}
  \end{equation}
  However, it doesn't seem that the result about commutation of
  $\mathrm Lf_!$ to colocalization is valid under the corresponding
  assumption $\scrO_A=f^*(\scrO_B)$, without assuming moreover
  $\scrO_B$ to be constant (hence $\scrO_A$ too), because already for
  the functor $f_!$ itself for modules it doesn't seem that
  commutation will hold.

  2)\enspace I should correct as silly mistake I made at the very
  beginning of this section, when rashly stating that the functor
  $\mathrm Lf_!$ may be defined for \emph{any} map \eqref{eq:138.1}
  between ringed objects in \Cat. I was thinking of the fact that for
  any ringed object $(A,\scrO_A)$ in \Cat, there are indeed enough
  projectives in the category of $\scrO_A$-modules --\pspage{586}
  thus, the modules
  \begin{equation}
    \label{eq:138.12}
    \scrO_A^{(a)}, \quad\text{for $a$ in $A$,}\tag{12}
  \end{equation}
  are clearly projective, and there are ``sufficiently many''. Thus,
  any additive functor from $\Mod(\scrO_A)$ to an abelian category
  admits a total left derived functor. However, it is not always true,
  for a map \eqref{eq:138.1} of small ringed categories, that the
  corresponding inverse image functor for modules
  \[G\mapsto f^*(G) = {f_0}^*(G) \otimes_{{f_0}^*(\scrO_B)} \scrO_A\]
  admits a left adjoint $f_!$, or what amounts to the same, that this
  functor (which is right exact) is left exact and commutes to small
  products. It isn't even true, necessarily, when we assume $A$ and
  $B$ to be the final category! Left exactness of $f^*$ just means
  flatness of $\scrO_A$ over $\scrO_B$, i.e., that for any $a$ in $A$,
  $\scrO_A(a)$ is flat as a module over $\scrO_B(b)$, where
  $b=f(a)$. As for commutation to small products, it amounts (together
  with the first condition) to the still more exacting condition that
  $\scrO_A(a)$ should be a \emph{projective module of finite type over
    $\scrO_B(b)$}, for any $a$ in $A$. This condition is so close to
  the condition $\scrO_A={f_0}^*(\scrO_B)$ (already considered in
  \eqref{eq:138.11} above), that for a bird's eye view as we are
  aiming at here, we may as well assume this slightly stronger
  condition! Anyhow, as noticed before, to get commutation of $f_!$
  with restriction of scalars and with colocalization, even this
  assumption isn't enough, apparently, and we'll have to assume
  moreover that $\scrO_B$ (hence also $\scrO_A$) is constant, or for
  the least, locally constant.

  This brings us back to the case when a fixed ring $k$ is given, and
  when we are working with categories of presheaves of $k$-modules --
  a situation studied at length in chapter~\ref{ch:V}. In case the
  target category $B$ is the final category, hence $\mathrm Lf_!$ is
  just (absolute) total homology of $A$, this then may be computed
  nicely, using an ``integrator'' for $A$, namely a projective
  resolution of $k_{A\op}$ in $(A\op)\uphat_k$. This is more or less
  where we ended up by the end of chapter~\ref{ch:V}, when developing
  a nicely autodual homology-cohomology set-up, replacing the category
  of $k$-modules $\AbOf_k$ by a more-or-less arbitrary abelian
  category. I was about to go on and carry through a similar treatment
  in the ``relative'' case, namely for an arbitrary map $f$ in \Cat{}
  (but then I got caught unsuspectingly by that unending digression on
  schematic homotopy types, finally making up a whole chapter by
  itself). Maybe it is still worthwhile to come back to
  this\pspage{587} without necessarily grinding through a complete
  formulaire for the five main operations we got so far (namely
  $\mathrm Lf_!$, $f^*$, $\mathrm Rf_*$, $\Lotimes$, $\RbHom$). Not
  later than two pages ago or so, we were faced again with two visible
  dual statements, one about $\mathrm Lf_!$, the other about $\mathrm
  Rf_*$ -- and feeling silly not to be able to merely deduce one from
  the other!
\end{remarks}

\bigbreak

\presectionfill\ondate{29.10.}\par

\hangsection[Retrospective on ponderings on abelianization, and on
\dots]{Retrospective on ponderings on abelianization, and on coalgebra
  structures in \texorpdfstring{\Cat}{(Cat)}.}\label{sec:139}%
Last night I still pondered a little about the $\mathrm Lf_!$
operation, and did some more reading in the notes of
chapter~\ref{ch:V} of about three months ago, which had been getting a
little distant in my mind. In those notes a great deal of emphasis
goes with the notions of integrators and cointegrators -- as a matter
of fact, that whole chapter sprung from an attempt to understand the
meaning of certain ``standard complexes'' associated with standard
test categories such as $\Simplex$, which then led us to the notions
of an integrator (via the intermediate one of an
``abelianizator''). This in turn brought us to the $*_k$ formalism,
expressing most conveniently the relationship between categories such
as $A\uphat_k$ and $(A\op)\uphat_k$ ($k$ any commutative ring), and
its various avatars. I was so pleased with this formalism and its
``computational'' flavor, that in my enthusiasm I subsumed under it
the dual treatment of homology and cohomology for an object $A$ of
\Cat{} (with coefficients in a complex of \scrM-valued presheaves,
\scrM{} being an abelian category satisfying some mild conditions), in
section~\ref{sec:108}, as a particular case of the total derived
functors $F_\bullet \Last_k L_\bullet'$ and
$\RHom_k(L_\bullet,K^\bullet)$ (with values in the derived category of
\scrM). The main point here was using \emph{projective resolutions} of
the argument $L_\bullet$ in $B\uphat_k$ or $A\uphat_k$ (where
$B=A\op$), which in the most important case was just the constant
sheaf of rings $k_B$ or $k_A$ -- rather than resolve the argument
$F_\bullet$ (projectively) or $K^\bullet$ (injectively), namely the
coefficients for homology or cohomology. It was the enthusiasm of the
adept of a game he just discovered -- I was going to try it out for
the next step, namely relative homology and cohomology $\mathrm Lf_!$
and $\mathrm Rf_*$, with $f$ any map in \Cat{} -- but then I got
caught by the more fascinating schematization game. Coming back now
upon the rather routine matter of looking up a comprehensive mutually
dual treatment of $\mathrm Lf_!$ and $\mathrm Rf_*$, it doesn't seem
that the formalism of integrators and cointegrators is going to be of
much help. To be more specific, in order to compute (or simply define)
$\mathrm Lf_!(F_\bullet)$\pspage{588} or $\mathrm Rf_*(K^\bullet)$,
for a general map
\[f:A\to B\]
in \Cat, and $F_\bullet\in\D^-(A\uphat_k)$,
$K^\bullet\in\D^+(A\uphat_k)$, I do not see any means of bypassing
projective resolutions of $L_\bullet$, injective ones of
$K^\bullet$. Taking the more familiar case of $\mathrm
Rf_*(K^\bullet)$, one natural idea of course would be to take a
projective resolution $L_\bullet^A$ of $k_A$ (i.e., a cointegrator for
$A$), and write tentatively
\begin{equation}
  \label{eq:139.1}
  \mathrm Rf_*(K^\bullet) \overset{\text{?}}{\simeq}
  f_*(\bHom_k^{\bullet\bullet}(L_\bullet^A,K^\bullet).\tag{1}
\end{equation}
The (misleading) reflex inducing us to write down this formula, is
that this formula looks as if it were to boil down to the similar
(correct) formula for the maps $A_{/b}\to e$, when taking the
localizations on $B$,
\[B_{/b}\to B.\]
For this intuition to be correct, it should be true that the
restriction (or ``localization'') of $L_\bullet^A$ to $A_{/b}$ (which
is of course a resolution of the constant sheaf $k$ on $A_{/b}$) is
indeed a cointegrator on $A_{/b}$, namely that its components are
still projective. This, however, we suspect, will hold true only under
very special assumptions -- as in general, it is inverse image by
\emph{colocalization} (not by localization) that takes projectives
into projectives. Thus, I don't expect a relation \eqref{eq:139.1} to
hold, except under most exacting conditions on $f$ and $K^\bullet$,
which I didn't try to pin down. To take an example, assume $A$ has a
final object, hence $k_A$ is projective and we may take
$L_\bullet^A=k_A$, then \eqref{eq:139.1} reads (when $K^\bullet=K$ is
reduced to degree zero)
\[\mathrm Rf_*(K) \simeq f_*(K), \quad
  \text{i.e., $\mathrm R^if_*(K)=0$ for $i>0$,}\]
which need not hold true even if $K$ is constant ($=k_A$ say). For
instance, we may start with an arbitrary object $A_0$ of \Cat{} and
add a final object $e_1$ to get $A$ (intuitively, it is the ``cone''
over $A_0$), which is mapped into the cone over $e$,
$B=\Simplex_1=(0\to 1)$, in the obvious way. Then for any sheaf $K$ on
$A$, with restriction $K_0$ to $A_0$, we get
\[\mathrm R^if_*(K)_0\text{ (fiber at $0$) } = \mathrm H^i(A_0,K_0),\]
which needs not be $0$ for $i>0$.

There is a big blunder at the end of section~\ref{sec:101}, where the
formulæ \eqref{eq:101.12}, \eqref{eq:101.12prime}, \eqref{eq:101.13}
(p.~\ref{p:377}) (supposed to be ``essentially trivial''), are false
for essentially the same kind of reason. The formulæ ran into the
typewriter as a matter of course, as\pspage{589} they looked just the
same as familiar ones from the standard duality formulaire (with
$\mathrm Lf_!$, $f^*$ replaced by $\mathrm Rf_!$, $f^!$). The first
one reads, in case as above when $A$ has a final object $\varepsilon$
\begin{equation}
  \label{eq:139.2}
  \mathrm Rf_*(f^*(K)) \overset{\text{?}}{\simeq}
  \bHom_k(f_!(k_A),K),\tag{2}
\end{equation}
But we have 
\[f_!(k_A)=f_!(k^{(\varepsilon)}) = k^{(b)} , \quad
  \text{where $b=f_!(\varepsilon)=f(\varepsilon)$,}\]
and if $f$ takes final object into final object, we thus get
$f_!(k_A)=k_B$, and the right-hand side of \eqref{eq:139.2} is just
$K$, and hence \eqref{eq:139.2} implies
\[\mathrm R^if_*(f^*(K)) = 0 \quad \text{for $i>0$,}\]
which, however, needs not be true even for $K=k_B$, as we saw with the
previous example.

To come back to the general $\mathrm Lf_!$, $\mathrm Rf_*$ formalism,
it seems it can't be helped, we'll have to do the usual silly thing
and just resolve the argument involved, projectively in one case,
injectively in the other. Even in the ``absolute case'' when $B=e$, we
couldn't help it, either taking such resolutions, when it comes to
working with the internal operations $\RbHom_k(L_\bullet,K^\bullet)$
or $F_\bullet\Lotimes_k L_\bullet'$, as we saw already in
section~\ref{sec:108} (which should have tuned down a little my
committedness to integrators, but it didn't!).

Once this got clear, in order to turn the homological algebra mill,
all we still need is a handy criterion for existence of ``enough''
projectives or injectives in categories of the type
\[ \AhatM = \bHom(A\op, \scrM),\]
with \scrM{} an abelian category. In section~\ref{sec:109} we got such
a criterion (prop.~\ref{prop:109.4} p.~\ref{p:433}) -- the simplest
one can imagine: it is sufficient that (besides the stability under
direct or inverse limits, needed anyhow in order for the functor $f_!$
resp.\ $f_*$ from \AhatM{} to \BhatM\pspage{590} to exist) that in
\scrM{} there should be enough projectives resp.\ injectives! This
gives what is needed, surely, in order to grind through a mutually
dual treatment of relative homology $\mathrm Lf_!$ and cohomology
$\mathrm Rf_*$. I have the feeling that the little work ahead, for
defining the basic operations and working out the relevant formulaire
(including ``cap products''), is a matter of mere routine, and I don't
really expect any surprise may come up. Therefore, I don't feel like
grinding through this, and rather will feel free to use whenever
needed the most evident formulæ, as the adjunction formulæ between the
three functors
\[ \mathrm Lf_!, \quad f^*, \quad \mathrm Rf_*,\]
transitivity isomorphisms for a composition of maps, possibly also
various projection formulæ -- trying however to be careful with base
change questions, notably for $\mathrm Lf_!$, and not repeat the same
blunders!

In retrospect, the main role for me of the reflections of
chapter~\ref{ch:V} on abelianization has been to become a little more
familiar with the \emph{homology} formalism in \Cat, namely
essentially with the $\mathrm Lf_!$ operation, which has been more or
less a white spot in my former experience, centered rather upon
cohomology. An interesting, still somewhat routine byproduct of these
ponderings has been the careful formulation of the duality
relationship around the pairings between sheaves and co-sheaves,
namely the operations $*_k$ and $\oast_k$. Granting this, the game
with integrators and cointegrators boils down to the standard reflex,
of taking projective resolutions of the ``unit'' sheaf or cosheaf,
$k_A$ or $k_{A\op}$ -- which are the most obvious objects at hand of
all, in our coefficient categories.

The one idea which seems to me of wider scope and significance in this
whole reflection on abelianization, is the ``further step in
linearization'', whereby the models in \Cat{} are replaced by their
$k$-additive envelope, endowed moreover with their natural
\emph{diagonal map} (section~\ref{sec:109}). The psychological effect
of this discovery has been an immediate one -- it triggered at once
the reflection on schematization of chapter~\ref{ch:VI}. This
reflection, apparently, took me into a rather different direction from
those ``$k$-coalgebra structures in \Cat'', which had reawakened and
made more acute the feeling that homotopy types should make sense
``over any ground ring''. Coming back now to the homology and
cohomology formalism within \Cat, it remains\pspage{591} a very
striking fact indeed for me, that as far as I can see at present, all
basic operations in the (commutative) homology and cohomology
formalism in \Cat, and their basic properties and interrelations,
should make sense for these linearized objects, which therefore may be
viewed as more perfect carriers still than small categories for
embodying the relevant formalism. To what extent this feeling is
indeed justified, cannot of course be decided beforehand -- only
experience can tell. For instance, does the ``six operations duality
formalism'' contemplated lately for \Cat{} (sections~\ref{sec:136},
\ref{sec:137}), which has still to be worked out, carry over to this
wider, linearized set-up? This will become clearer when the relevant
notions in \Cat{} are understood, so that it will become a meaningful
question whether for a map in \Cat, the property of being proper, or
smooth, or an immersion, may be read off in terms of properties of the
corresponding map of coalgebra structures in \Cat. And what about
subtler types of cohomology operations, such as the Steenrod
operations (which, I remember, may be defined in the context of
cohomology with coefficients in general sheaves of $\bF_p$-modules)?

When concerned with homology and cohomology formalism in \Cat,
involving general sheaves of coefficients, not merely locally constant
ones, we are leaving, strictly speaking, the waters of a reflection on
``homotopy models''. The objects of \Cat{} now are no longer viewed as
mere models for homotopy types, but rather, each one as defining a
topos, with the manifold riches it carries; a richness similar to the
one of a topological space, almost all of which is being stripped off
when looking at the mere homotopy type -- including even such basic
properties as dimension, compactness, smoothness, cardinality and the
like. When passing from an object in \Cat{} to the corresponding
coalgebra structure in \Cat, much of this richness is preserved --
maybe everything, indeed, which can be expressed in terms of
\emph{commutative} sheaves of coefficients. As for the realm of
non-commutative cohomology formalism (which is supposed to be the main
theme of these notes, with overall title ``Pursuing Stacks''!), it
doesn't seem, not at first sight at any rate, that much of this could
be read off the enveloping coalgebra $P=\Add_k(A)$ of a given object
$A$ in \Cat{} -- except of course in the case when $A$ can be
recovered in terms of $P$, maybe as the category of ``exponential''
pairs $(x,u)$, where $x$ is in $P$ and
\[ u:\delta(x) \simeq x \otimes_k x\]
an isomorphism.

\bigbreak

\presectionfill\ondate{4.11.}\pspage{592}\par

\hangsection{The meal and the guest.}\label{sec:140}%
Life keeps pushing open the doors of that well-tempered hothouse of my
mathematical reflections, as a fresh wind and often an impetuous one,
sweeping off the serene quietness of abstraction, -- a breath rich
with the manifold fragrance of the world we live in. This is the world
of conflict, weaving around each birth and each death and around the
lovers' play alike -- it is the world we all have been born into
without our choosing. I used to see it as a stage -- the stage set for
our acting. Our freedom (rarely used indeed!) includes choosing the
role we are playing, possibly changing roles -- but not choosing or
changing the stage. It doesn't seem the stage ever changes during the
history of mankind -- only the decors kept changing. More and more,
however, over the last years, I have been feeling this world I am
living in, the world of conflict, somehow as a \emph{meal} -- a meal
of inexhaustible richness. Maybe the ultimate fruit and meaning of all
my acting on that stage, is that parts of that meal, of that richness,
be actually eaten, digested, assimilated -- that they become part of
the flesh and bones of my own being. Maybe the ultimate purpose of
conflict, so deeply rooted in every human being, is to be the raw
material, to be eaten and digested and changed into understanding
about conflict. Not a collective ``understanding'' (I doubt there is
such a thing!), written down in textbooks or sacred books or whatever,
nor even something expressed or expressible in words necessarily --
but the kind of immediate knowledge, rather, the walker has about
walking, the swimmer about swimming, or the suckling about milk and
mother's breasts. My business is to be a learner, not a teacher --
namely to allow this process to take place in my being, letting the
world of conflict, of suffering and of joy, of violence and of
tenderness, enter and be digested and become knowledge about myself.

I am not out, though, to write a ``journal intime'' or meditation
notes, so I guess I better get back to the thread of mathematical
reflection where I left it, rather than write allusively about the
events of these last days, telling me about life and about myself
through one of my children.

\bigbreak

\presectionfill\ondate{12.11.}\par

Maybe at times I like to give the impression, to myself and hence
others, that I am the easy learner of things of life, wholly
relaxed,\pspage{593} ``cool'' and all that -- just keen for learning,
for eating the meal and welcome smilingly whatever comes with its
message, frustration and sorrow and destructiveness and the softer
dishes alike. This of course is just humbug, an image
d'Épinal\scrcomment{the Épinal prints being proverbial for a naive
  depiction, showing only good aspects\dots} which at whiles I'll kid
myself into believing I am like. Truth is that I am a hard learner,
maybe as hard and reluctant as anyone. At any rate, the inbuilt
mechanisms causing rejection of the dishes unpalatable to my wholly
conditioned, wholly ego-controlled taste, are as much present in me
now as they have ever been in my life, and as much as in anyone else I
know. This interferes a lot with the learning -- it causes a
tremendous amount of friction and energy dispersion (wholly unlike
what happens when I am learning mathematics, say, namely discovering
things about any kind of substance which my own person and ego is not
part of\dots). If there has been anything new appearing in my life, it
is surely not the end of this process of dispersion, or the end of
inertia, closely related to dispersion. It is something, rather, which
causes learning to take place all the same -- be it the hard way, as
it often happens, very much like the troubled digestion of one who
took a substantial meal reluctantly or in a state of nervousness, of
crispation. Once one is through with the digestion, though, the food
one ate is transformed into flesh and muscles, blood and bones and the
like, just as good and genuine as if the meal had been taken relaxedly
and with eager appetite, as it deserved. What really counts for the
process of assimilation to be able to take place, is that in a certain
sense the food, palatable or not, be \emph{accepted} -- not vomited,
or just kept in the bowels like a foreign body, sometimes for
decennaries.\scrcomment{decennary = decade} The remarkable fact I come
to know through experience, is that even after having been kept thus
inertly for a lifetime, a process of digestion and assimilation may
still come into being and transform the obtrusive stuff into living
substance.

During the last week I have been sick for a few days -- a
grippe\scrcomment{grippe = the flu} I might say, but surely a case of
troubled digestion too. It seems though I am through now -- till the
next case at any rate! I daresay life has been generous with me for
these last three months, while I haven't even taken the trouble to
stop with the mathematical nonsense for any more than a week or
two. This week, too, I still did some mathematical scratchwork, still
along the lines of abelianization, which keeps showing a lot richer
than suspected.


\backmatter

\raggedright
\printbibliography

@Article{Giraud1964,
  author	= {Giraud, Jean},
  title		= {M\'ethode de la descente},
  journal	= {Bull. Soc. Math. France M\'em.},
  fjournal	= {Soci\'et\'e Math\'ematique de France. Bulletin.
		  Suppl\'ement. M\'emoire},
  volume	= {2},
  year		= {1964},
  pages		= {viii+150},
  mrclass	= {14.55},
  mrnumber	= {0190142 (32 \#7556)},
  mrreviewer	= {R. Deheuvels},
  eprint        = {94540},
  eprinttype    = {eudml}
}

@Book{FGA,
  author	= {Grothendieck, A.},
  title		= {Fondements de la Géométrie Algébrique},
  note		= {Multigraphié, Secrétariat mathématique 11, rue Pierre Curie, Paris},
  publisher	= {Séminaire Bourbaki},
  year		= {1957-1962}
}

@Book{Giraud1971,
  author	= {Giraud, Jean},
  title		= {Cohomologie non ab\'elienne},
  note		= {Die Grundlehren der mathematischen Wissenschaften, Band
		  179},
  publisher	= {Springer-Verlag, Berlin-New York},
  year		= {1971},
  pages		= {ix+467},
  mrclass	= {14F20 (14L99 18D30 18F20 55B99 55F65)},
  mrnumber	= {0344253 (49 \#8992)},
  mrreviewer	= {R. T. Hoobler}
}

@Book{Illusie1971,
  author	= {Illusie, Luc},
  title		= {Complexe cotangent et d\'eformations. {I}},
  series	= {Lecture Notes in Mathematics, Vol. 239},
  publisher	= {Springer-Verlag, Berlin-New York},
  year		= {1971},
  pages		= {xv+355},
  mrclass	= {14B10 (14D15 14F20 55D20)},
  mrnumber	= {0491680 (58 \#10886a)},
  mrreviewer	= {Larry Breen},
  doi           = {10.1007/BFb0059052},
}

@Book{Illusie1972,
  author	= {Illusie, Luc},
  title		= {Complexe cotangent et d\'eformations. {II}},
  series	= {Lecture Notes in Mathematics, Vol. 283},
  publisher	= {Springer-Verlag, Berlin-New York},
  year		= {1972},
  pages		= {vii+304},
  mrclass	= {14B10 (14D15 14F20 55D20)},
  mrnumber	= {0491681 (58 \#10886b)},
  mrreviewer	= {Larry Breen},
  doi           = {10.1007/BFb0059573}
}

@Book{Saavedra1972,
  author	= {Saavedra Rivano, Neantro},
  journal	= {Bulletin de la Société Mathématique de France},
  pages		= {4170--430},
  publisher	= {Société mathématique de France},
  title		= {Catégories tannakiennes},
  eprint        = {87193},
  eprinttype    = {eudml},
  volume	= {100},
  year		= {1972}
}

@InCollection{Benabou1967,
  author	= {B{\'e}nabou, Jean},
  title		= {Introduction to bicategories},
  booktitle	= {Reports of the {M}idwest {C}ategory {S}eminar},
  series	= {Lecture Notes in Mathematics, Vol. 47},
  pages		= {1--77},
  publisher	= {Springer, Berlin},
  year		= {1967},
  mrclass	= {18.10},
  mrnumber	= {0220789 (36 \#3841)},
  mrreviewer	= {J. R. Isbell},
  doi           = {10.1007/BFb0074299}
}

@Book{GabrielZisman1967,
  author	= {Gabriel, Peter and Zisman, Michel},
  title		= {Calculus of fractions and homotopy theory},
  series	= {Ergebnisse der Mathematik und ihrer Grenzgebiete, Band
		  35},
  publisher	= {Springer-Verlag New York, Inc., New York},
  year		= {1967},
  pages		= {x+168},
  mrclass	= {55.40 (18.00)},
  mrnumber	= {0210125 (35 \#1019)},
  mrreviewer	= {A. K. Bousfield},
  doi           = {10.1007/978-3-642-85844-4}
}

@Book{Andre1967,
  author	= {Andr{\'e}, Michel},
  title		= {M\'ethode simpliciale en alg\`ebre homologique et
		  alg\`ebre commutative},
  series	= {Lecture Notes in Mathematics, Vol. 32},
  publisher	= {Springer-Verlag, Berlin-New York},
  year		= {1967},
  pages		= {iii+122},
  mrclass	= {18.20 (13.00)},
  mrnumber	= {0214644 (35 \#5493)},
  mrreviewer	= {D. W. Knudson}
}

@Book{Quillen1967,
  author	= {Quillen, Daniel G.},
  title		= {Homotopical algebra},
  series	= {Lecture Notes in Mathematics, Vol. 43},
  publisher	= {Springer-Verlag, Berlin-New York},
  year		= {1967},
  pages		= {iv+156 pp. (not consecutively paged)},
  mrclass	= {18.20 (55.00)},
  mrnumber	= {0223432 (36 \#6480)},
  mrreviewer	= {A. K. Bousfield},
  doi           = {10.1007/BFb0097438}
}

@InCollection{Quillen1970,
  author	= {Quillen, Daniel G.},
  title		= {On the (co-)homology of commutative rings},
  booktitle	= {Applications of {C}ategorical {A}lgebra ({P}roc. {S}ympos.
		  {P}ure {M}ath., {V}ol. {XVII}, {N}ew {Y}ork, 1968)},
  pages		= {65--87},
  publisher	= {Amer. Math. Soc., Providence, R.I.},
  year		= {1970},
  mrclass	= {13.90 (18.00)},
  mrnumber	= {0257068 (41 \#1722)},
  mrreviewer	= {S. Yuan},
  url           = {http://www.ams.org/books/pspum/017/0257068}
}

@Book{SGA4vol1,
  shorthand     = {{SGA}~4.1},
  editor        = {Artin, Michael and Grothendieck, Alexander and Verdier, Jean-Louis},
  title		= {Th\'eorie des topos et cohomologie \'etale des sch\'emas.
		  {T}ome 1: {T}h\'eorie des topos},
  series	= {Lecture Notes in Mathematics, Vol. 269},
  note		= {S{\'e}minaire de G{\'e}om{\'e}trie Alg{\'e}brique du
		  Bois-Marie 1963--1964 (SGA~4)},
  publisher	= {Springer-Verlag, Berlin-New York},
  year		= {1972},
  pages		= {xix+525},
  mrclass	= {14-06},
  mrnumber	= {0354652 (50 \#7130)},
  doi           = {10.1007/BFb0081551}
}

@Book{SGA4andhalf,
  shorthand     = {{SGA}~4$\frac12$},
  author	= {Deligne, Pierre},
  title		= {Cohomologie \'etale},
  series	= {Lecture Notes in Mathematics, Vol. 569},
  note		= {S{\'e}minaire de G{\'e}om{\'e}trie Alg{\'e}brique du
		  Bois-Marie SGA~4$\frac12$},
  publisher	= {Springer-Verlag, Berlin-New York},
  year		= {1977},
  pages		= {iv+312pp},
  mrclass	= {14F20},
  mrnumber	= {0463174 (57 \#3132)},
  mrreviewer	= {J. S. Milne},
  doi           = {10.1007/BFb0091516}
}

@Book{SGA2,
  shorthand     = {{SGA}~2},
  Title		= {Cohomologie locale des faisceaux cohérents et théorèmes de Lefschetz locaux et globaux},
  note		= {S{\'e}minaire de G{\'e}ometrie Alg{\'e}brique du
		  Bois-Marie (SGA 2), Laszlo, Yves (ed.)},
  publisher	= {Amsterdam: North-Holland Publishing Company},
  year		= {1968}
}

@Book{SGA6,
  shorthand     = {{SGA}~6},
  editor        = {Berthelot, Pierre and Grothendieck, Alexander and Illusie, Luc},
  title		= {Th\'eorie des intersections et th\'eor\`eme de
		  {R}iemann-{R}och},
  series	= {Lecture Notes in Mathematics, Vol. 225},
  note		= {S{\'e}minaire de G{\'e}om{\'e}trie Alg{\'e}brique du
		  Bois-Marie 1966--1967 (SGA~6)},
  publisher	= {Springer-Verlag, Berlin-New York},
  year		= {1971},
  pages		= {xii+700},
  mrclass	= {14-06},
  mrnumber	= {0354655 (50 \#7133)},
  doi           = {10.1007/BFb0066283}
}

@Book{Andre1974,
  author	= {Andr{\'e}, Michel},
  title		= {Homologie des alg\`ebres commutatives},
  note		= {Die Grundlehren der mathematischen Wissenschaften, Band
		  206},
  publisher	= {Springer-Verlag, Berlin-New York},
  year		= {1974},
  pages		= {xv+341},
  mrclass	= {18H20 (14D15 13DXX 12GXX)},
  mrnumber	= {0352220 (50 \#4707)},
  mrreviewer	= {R. M. Fossum},
  doi           = {10.1007/978-3-642-51449-4}
}

@Book{BousfieldKan1972,
  author	= {Bousfield, Aldrigde K. and Kan, Daniel M.},
  title		= {Homotopy limits, completions and localizations},
  series	= {Lecture Notes in Mathematics, Vol. 304},
  publisher	= {Springer-Verlag, Berlin-New York},
  year		= {1972},
  pages		= {v+348},
  mrclass	= {55J05},
  mrnumber	= {0365573 (51 \#1825)},
  mrreviewer	= {Harold Hastings},
  doi           = {10.1007/978-3-540-38117-4}
}

@Book{Raynaud1975,
  author	= {Raynaud, Mich{\`e}le},
  title		= {Th\'eor\`emes de {L}efschetz en cohomologie coh\'erente et
		  en cohomologie \'etale},
  note		= {Bull. Soc. Math. France, M{\'e}m. No. 41, Suppl{\'e}ment
		  au Bull. Soc. Math. France, Tome 103},
  publisher	= {Soci\'et\'e Math\'ematique de France, Paris},
  year		= {1975},
  pages		= {176},
  mrclass	= {14F20},
  mrnumber	= {0407021 (53 \#10804)},
  mrreviewer	= {G. Horrocks},
  eprint        = {94714},
  eprinttype    = {eudml}
}

@Book{GCS,
  author	= {Sinh, Hoang Xuan},
  title		= {$\Gr$-catégories},
  year		= {1975},
  url		= "https://agrothendieck.github.io/divers/GCSscan.pdf"
}

@Book{JOU69,
  author	= {Jouanolou, Jean-Pierre},
  title		= {Catégories dérivées et cohomologie $\ell$-adique},
  year		= {1969}
}

@Letter{Grothendieck1983b,
  author	= "Alexander Grothendieck",
  year		= "June 27, 1983",
  note		= "Letter to Gerd Faltings, describing ``anabelian algebraic
		  geometry''. Published in \textcite{LochakSchneps2000}",
  url		= "https://agrothendieck.github.io/divers/LGF2783.pdf"
}

@Article{Cartan1976,
  author	= {Cartan, Henri},
  title		= {Th\'eories cohomologiques},
  journal	= {Invent. Math.},
  fjournal	= {Inventiones Mathematicae},
  volume	= {35},
  year		= {1976},
  pages		= {261--271},
  issn		= {0020-9910},
  mrclass	= {55B20},
  mrnumber	= {0431137 (55 \#4139)},
  mrreviewer	= {Jean-Michel Lemaire},
  doi           = {10.1007/BF01390140},
}

@Article{Miller1978,
  author	= {Miller, Edward Y.},
  title		= {de {R}ham cohomology with arbitrary coefficients},
  journal	= {Topology},
  fjournal	= {Topology. An International Journal of Mathematics},
  volume	= {17},
  year		= {1978},
  number	= {2},
  pages		= {193--203},
  issn		= {0040-9383},
  mrclass	= {55B30 (58A10)},
  mrnumber	= {0467722 (57 \#7575)},
  mrreviewer	= {I. Vaisman},
  doi           = {10.1016/S0040-9383(78)90025-3}
}

@Article{Anderson1978,
  author	= {Anderson, Donald W.},
  title		= {Fibrations and geometric realizations},
  journal	= {Bull. Amer. Math. Soc.},
  fjournal	= {Bulletin of the American Mathematical Society},
  volume	= {84},
  year		= {1978},
  number	= {5},
  pages		= {765--788},
  issn		= {0002-9904},
  mrclass	= {55D99},
  mrnumber	= {0500935 (58 \#18426)},
  mrreviewer	= {J. P. May},
  doi           = {10.1090/S0002-9904-1978-14512-1},
}

@Article{Thomason1979,
  author	= {Thomason, Robert W.},
  title		= {Homotopy colimits in the category of small categories},
  journal	= {Math. Proc. Cambridge Philos. Soc.},
  fjournal	= {Mathematical Proceedings of the Cambridge Philosophical
		  Society},
  volume	= {85},
  year		= {1979},
  number	= {1},
  pages		= {91--109},
  issn		= {0305-0041},
  coden		= {MPCPCO},
  mrclass	= {18F25},
  mrnumber	= {510404 (80b:18015)},
  mrreviewer	= {Daniel R. Grayson},
  doi		= {10.1017/S0305004100055535},
}

@Article{Thomason1980,
  author	= {Thomason, Robert W.},
  title		= {Cat as a closed model category},
  journal	= {Cahiers Topologie G\'eom. Diff\'erentielle},
  fjournal	= {Cahiers de Topologie et G\'eom\'etrie Diff\'erentielle},
  volume	= {21},
  year		= {1980},
  number	= {3},
  pages		= {305--324},
  issn		= {0008-0004},
  mrclass	= {18A99 (55U35)},
  mrnumber	= {591388 (82b:18005)},
  mrreviewer	= {Harold Hastings},
  eprint        = {91233},
  eprinttype    = {eudml}
}

@InBook{DeligneMilne1982,
  author        = {Deligne, Pierre and Milne, James S.},
  title         = {Tannakian Categories},
  booktitle	= {Hodge cycles, motives, and {S}himura varieties},
  bookauthor	= {Deligne, Pierre and Milne, James S. and Ogus, Arthur and
		  Shih, Kuang-yen},
  series	= {Lecture Notes in Mathematics, Vol. 900},
  publisher	= {Springer-Verlag, Berlin-New York},
  year		= {1982},
  pages		= {101--228},
  isbn		= {3-540-11174-3},
  mrclass	= {14Kxx (10D25 12A67 14A20 14F30 14K22)},
  mrnumber	= {654325 (84m:14046)},
  url           = {http://www.jmilne.org/math/xnotes/tc.html}
}

@Thesis{Dakin1977,
  author	= {Dakin, M. Keith},
  title		= {Kan complexes and multiple groupoid structures},
  note          = {Published in Esquisses Math\'ematiques, Vol. 32 (1983)},
  type          = {Ph.D. Thesis},
  year          = {1977},
  institution   = {University College of North Wales, Bangor},
  url           = {http://ehres.pagesperso-orange.fr/Cahiers/dakinEM32.pdf},
}

@Book{Jones1983,
  author	= {Jones, David W.},
  title		= {Poly-{$T$}-complexes},
  series	= {U.C.N.W. Pure Maths Preprint},
  volume	= {84},
  note		= {Dissertation, University College of North Wales, Bangor,
		  1983},
  publisher	= {University College of North Wales, Bangor},
  year		= {1983},
  pages		= {171},
  mrclass	= {18G55 (20L15 55U10)},
  mrnumber	= {812869 (87a:18022b)},
  mrreviewer	= {Johannes Huebschmann}
}

@InCollection{Baues1985,
  author	= {Baues, Hans Joachim},
  title		= {On homotopy classification problems of {J}.\ {H}.\ {C}.\
		  {W}hitehead},
  booktitle	= {Algebraic topology, {G}\"ottingen 1984},
  series	= {Lecture Notes in Math.},
  volume	= {1172},
  pages		= {17--55},
  publisher	= {Springer, Berlin},
  year		= {1985},
  mrclass	= {55P15 (55S37)},
  mrnumber	= {825772 (87f:55006)},
  mrreviewer	= {Guido Mislin},
  doi		= {10.1007/BFb0074422},
}

@InCollection{Deligne1990,
  author	= {Deligne, Pierre},
  title		= {Cat\'egories tannakiennes},
  booktitle	= {The {G}rothendieck {F}estschrift, {V}ol.\ {II}},
  series	= {Progr. Math.},
  volume	= {87},
  pages		= {111--195},
  publisher	= {Birkh\"auser Boston, Boston, MA},
  year		= {1990},
  mrclass	= {14A99 (12H05 18A99)},
  mrnumber	= {1106898 (92d:14002)},
  mrreviewer	= {James Milne},
  doi           = {10.1007/978-0-8176-4575-5_3}
}

@Book{LochakSchneps2000,
  title		= {Geometric {G}alois actions. 1},
  series	= {London Mathematical Society Lecture Note Series},
  volume	= {242},
  editor	= {Schneps, Leila and Lochak, Pierre},
  note		= {Around Grothendieck's ``Esquisse d'un programme''},
  publisher	= {Cambridge University Press, Cambridge},
  year		= {1997},
  pages		= {iv+293},
  isbn		= {0-521-59642-4},
  mrclass	= {14-06},
  mrnumber	= {1483106 (98e:14003)},
  doi		= {10.1017/CBO9780511666124}
}

@Article{Cisinski1999,
  author	= {Cisinski, Denis-Charles},
  title		= {La classe des morphismes de {D}wyer n'est pas stable par
		  retractes},
  journal	= {Cahiers Topologie G\'eom. Diff\'erentielle Cat\'eg.},
  fjournal	= {Cahiers de Topologie et G\'eom\'etrie Diff\'erentielle
		  Cat\'egoriques},
  volume	= {40},
  year		= {1999},
  number	= {3},
  pages		= {227--231},
  issn		= {0008-0004},
  mrclass	= {18G30 (18G35 55U35)},
  mrnumber	= {1716777 (2000j:18008)},
  mrreviewer	= {Pilar C. Carrasco},
  eprint        = {91620},
  eprinttype    = {eudml}
}

@Article{Deligne2002,
  author	= {Deligne, Pierre},
  title		= {Cat\'egories tensorielles},
  note		= {Dedicated to Yuri I. Manin on the occasion of his 65th
		  birthday},
  journal	= {Mosc. Math. J.},
  fjournal	= {Moscow Mathematical Journal},
  volume	= {2},
  year		= {2002},
  number	= {2},
  pages		= {227--248},
  issn		= {1609-3321},
  mrclass	= {18D10 (18E30 20C99)},
  mrnumber	= {1944506 (2003k:18010)},
  mrreviewer	= {Viktor Ostrik},
  url           = {http://www.ams.org/distribution/mmj/vol2-2-2002/abst2-2-2002.html#deligne_abstract}
}

@Article{Cisinski2004,
  author	= {Cisinski, Denis-Charles},
  title		= {Le localisateur fondamental minimal},
  journal	= {Cah. Topol. G\'eom. Diff\'er. Cat\'eg.},
  fjournal	= {Cahiers de Topologie et G\'eom\'etrie Diff\'erentielle
		  Cat\'egoriques},
  volume	= {45},
  year		= {2004},
  number	= {2},
  pages		= {109--140},
  issn		= {1245-530X},
  mrclass	= {18E35 (18E15 55U35)},
  mrnumber	= {2072934 (2005f:18009)},
  mrreviewer	= {Manuel Bullejos Lorenzo},
  eprint        = {91678},
  eprinttype    = {eudml}
}

@Article{Cisinski2006,
  author	= {Cisinski, Denis-Charles},
  title		= {Les pr\'efaisceaux comme mod\`eles des types d'homotopie},
  journal	= {Ast\'erisque},
  fjournal	= {Ast\'erisque},
  number	= {308},
  year		= {2006},
  pages		= {xxiv+390},
  issn		= {0303-1179},
  isbn		= {978-2-85629-225-9},
  mrclass	= {55-02 (18F20 18G50 55P60 55U35)},
  mrnumber	= {2294028 (2007k:55002)},
  mrreviewer	= {Philippe Gaucher},
  url           = {http://www.math.univ-toulouse.fr/~dcisinsk/ast.pdf}
}

@Online{Kunzer2015,
  editor	= {K{\"u}nzer, Matthias},
  title		= {Correspondance Alexandre Grothendieck -- Ronald Brown},
  note		= {Avec la collaboration de R. Brown et G. Maltsiniotis},
  date		= {2015},
  url		= {http://webusers.imj-prg.fr/~georges.maltsiniotis/ps/agrb_web.pdf}		  
}

@Online{LGD7874,
  editor	= {K{\"u}nzer, Matthias},
  title		= {Letter to P Deligne, dated 7.8.74},
  date		= {1974},
  url		= {https://agrothendieck.github.io/divers/LGD7874.pdf}		  
}

@Online{tapisQ,
  title		= {Tapis de Quillen},
  note		= {Transcription in progress},
  date		= {1968},
  url		= {https://agrothendieck.github.io/divers/tapisQscan.pdf}		  
}

@Online{AG76IHES,
  title		= {De Rham complex with divided powers},
  note		= {IHES},
  date		= {1976},
  url		= {https://agrothendieck.github.io/divers/AG75IHES.pdf}		  
}

@Online{Gletters,
  title		= {Letters (collection)},
  note		= {Transcribed by Mateo Carmona},
  url		= {https://agrothendieck.github.io/divers/letters.pdf}		  
}

\end{document}
